\documentclass[10pt,reqno]{amsart}
\usepackage{amsfonts,amsrefs,latexsym,amsmath,amssymb,mathrsfs,verbatim,cancel}
\usepackage{url,color}
\usepackage{pifont}
\usepackage{upgreek}
\usepackage{fancyhdr}
\usepackage{hyperref}
\usepackage{kpfonts}
\usepackage{marvosym}
\usepackage[percent]{overpic}
\usepackage{caption}
\usepackage{subcaption}
\usepackage{pict2e}
\usepackage[hypcap=false]{caption}
\usepackage{xcolor}
\usepackage{soul}
\usepackage{wasysym}
\usepackage{scalerel}
\usepackage{accents}
\usepackage{slashed}
\usepackage{enumitem}
\usepackage{tabularx}
\usepackage{float}
\setlist[enumerate]{label*=\arabic*.}
\usepackage{array}
\newcolumntype{C}[1]{>{\centering\arraybackslash}m{#1}} 
\usepackage{multirow}
\usepackage{venturis2}

\topmargin - .4 in

\newtheorem{theorem}{Theorem}[section]
\newtheorem{proposition}[theorem]{Proposition}
\newtheorem{lemma}[theorem]{Lemma}
\newtheorem{corollary}[theorem]{Corollary}

\theoremstyle{definition}
\newtheorem{definition}[theorem]{Definition}
\newtheorem{remark}[theorem]{Remark}
\newtheorem{assumption}[theorem]{Assumption}
\newtheorem{notation}{Notation}[section]

\newcommand*\eqdef{\overset{\mbox{\tiny{def}}}{=}}
\newcommand{\R}{\mathbb{R}}
\newcommand{\T}{\mathbb{T}}
\newcommand{\p}{\partial}
\newcommand{\mr}{\mathring} 
\newcommand{\rmD}{\mathrm{d}}

\newcommand{\mytr}{{\mbox{\upshape tr}}} 

\newcommand{\otimesarray}{\overset{\diamond}{\otimes}}

\newcommand{\tander}{\mathcal{P}}
\newcommand{\tandersmall}{\mathcal{P}_*}
\newcommand{\comder}{\mathcal{Z}}
\newcommand{\comdersmall}{\mathcal{Z}_*}
\newcommand{\comderdoublesmall}{\mathcal{Z}_{**}}
\newcommand{\tanderY}{\mathcal{Y}}

\DeclareMathOperator{\diag}{diag}

\DeclareMathOperator*{\esssup}{ess\,sup}

\newcommand{\RRiemannPS}{\mathcal{R}_{(+)}^{\textnormal{PS}}}
\newcommand{\LRiemannPS}{\mathcal{R}_{(-)}^{\textnormal{PS}}}
\newcommand{\dataRRiemannPS}{\mathring{\mathcal{R}}_{(+)}^{\textnormal{PS}}}
\newcommand{\dataLRiemannPS}{\mathring{\mathcal{R}}_{(-)}^{\textnormal{PS}}}

\newcommand{\velocityPS}{v^{\textnormal{PS}}}
\newcommand{\LogDensityPS}{\LogDensity^{\textnormal{PS}}}
\newcommand{\SpeedPS}{c^{\textnormal{PS}}}
\newcommand{\SpeedPSdot}{\dot{c}^{\textnormal{PS}}}
\newcommand{\backgroundSpeedPS}{\overline{c^{\textnormal{PS}}}}
\newcommand{\backgroundSpeedprimePS}{\overline{\dot{c}^{\textnormal{PS}}}}

\newcommand{\LunitPS}{L^{\textnormal{PS}}}
\newcommand{\uLunitPS}{\underline{L}^{\textnormal{PS}}}

\newcommand{\muXPS}{\breve X^{\textnormal{PS}}}
\newcommand{\XPS}{X^{\textnormal{PS}}}

\newcommand{\muPS}{\upmu^{\textnormal{PS}}}
\newcommand{\LsmallPS}{L_{(Small)}^{\textnormal{PS}}}

\newcommand{\creasePS}{\partial_- \mathcal{B}^{\textnormal{PS}}}

\newcommand{\twoarghypPS}[2]{{^{(#2)}\widetilde{\Sigma}_{#1}^{\textnormal{PS}}}}
\newcommand{\threearghypPS}[3]{{^{(#3)}\widetilde{\Sigma}_{#1}^{#2;PS}}}

\newcommand{\CSregion}{\mathfrak{CS}}

\newcommand{\FPS}{F^{\textnormal{PS}}}
\newcommand{\FPSdot}{\dot{F}^{\textnormal{PS}}}
\newcommand{\FPStwodots}{\ddot{F}^{\textnormal{PS}}}

\newcommand{\PSdataalpha}{\mathring{\upalpha}^{\textnormal{PS}}}
\newcommand{\blowupdeltaPS}{\mathring{\updelta}_*^{\textnormal{PS}}}
\newcommand{\transversalsizedeltaPS}{\mathring{\updelta}^{\textnormal{PS}}}
\newcommand{\blowuptimePS}{T_{\textnormal{Shock}}^{\textnormal{PS}}}
\newcommand{\PSmutransversalHessiansize}{M_2^{\textnormal{PS}}}

\newcommand{\initialsmall}{\mathring{\upepsilon}}
\newcommand{\fundbootsmall}{\varepsilonup}
\newcommand{\auxbootsmall}{\varepsilonup^{1/2}}

\newcommand{\RRiemannpertinitial}{\mathring{\mathcal{R}}_{(+)}}
\newcommand{\LRiemannpertinitial}{\mathring{\mathcal{R}}_{(-)}}
\newcommand{\vtwopertinitial}{\mathring{v}^2} 
\newcommand{\vthreepertinitial}{\mathring{v}^3}
\newcommand{\vortrenormalizedpertinitial}{\mathring{\vortrenormalized}}
\newcommand{\GradEntpertinitial}{\mathring{\GradEnt}}
\newcommand{\VortVortpertinitial}{\mathring{\VortVort}}
\newcommand{\DivGradEntpertinitial}{\mathring{\DivGradEnt}}

\newcommand{\Lsmall}{L_{(\textnormal{Small})}}
\newcommand{\Xsmall}{X_{(\textnormal{Small})}}
 
\newcommand{\Yvf}[1]{Y_{(#1)}}
\newcommand{\Yvfsmall}[1]{Y_{(#1;\textnormal{Small})}}

\newcommand{\controlvars}{\upgamma}
\newcommand{\badcontrolvars}{\underline{\upgamma}}
\newcommand{\spertinitial}{\mathring{s}}
\newcommand{\wavearray}{\vec{\Psi}}
\newcommand{\wavearraypartial}{{\vec{\Psi}}_{\textnormal{(Partial)}}}
\newcommand{\Wtrans}{\breve{W}}
\newcommand{\Wtransarg}[1]{{ ^{(#1)} \mkern-2mu \breve{W}}}
\newcommand{\Wtranstwoarg}[2]{{ ^{(#1)} \mkern-2mu {\breve{W}^{#2}}}}

\newcommand{\Rtransarg}[1]{{ ^{(#1)}  \mkern-2mu \breve{R}}}

\newcommand{\Rtransunitarg}[1]{{ ^{(#1)}  \mkern-2mu \hat{R}}}

\newcommand{\hypnormalarg}[1]{{^{(#1)}  \mkern-2mu \widetilde{N}}}

\newcommand{\hypunitnormalarg}[1]{{ ^{(#1)} \mkern-2mu \hat{N}}}

\newcommand{\hypg}{\widetilde{g}}
\newcommand{\hypginverse}{\widetilde{g}\, ^{-1}}
\newcommand{\gtorusroughfirstfund}{\widetilde{\gtorus}}
\newcommand{\gtorusroughinversefirstfund}{\widetilde{\gtorus} \,^{-1}}
\newcommand{\gtorusCOV}{\mathbf{A}}

\newcommand{\timefunction}{\uptau}
\newcommand{\timefunctionarg}[1]{{^{(#1)} \mkern-1mu \uptau}}
\newcommand{\PStimefunctionarg}[1]{{^{(#1)} \mkern-1mu \uptau_{\textnormal{PS}}}}
\newcommand{\timefunctionboot}{\timefunction_{\textnormal{Boot}}}

\newcommand{\mupositive}{\mulevelsetvalue_0}
\newcommand{\boringregionmupositive}{\mulevelsetvalue_1}
\newcommand{\awayfromsmallneightborhoodnmupositive}{\upmu_2}
\newcommand{\Ntop}{N_{\textnormal{top}}}

\newcommand{\newtimefunction}{\timefunctionarg{\textnormal{Interesting}}} 
\newcommand{\tisafunctionalonglevelsetsofnewtimefunctionarg}[1]{{^{(\textnormal{Interesting})}\mathfrak{t}_{#1}}} 
\newcommand{\levelsetgeneratornewtimefunction}{\breve{H}} 
\newcommand{\partialderivativewithrespecttonewtimefunction}{\breve{G}}

\newcommand{\HarmlessWave}[1]{\textnormal{Harmless}_{(\textnormal{Wave})}^{[1,#1]}}

\newcommand{\zetatan}{\upzeta^{(\textnormal{Tan--}\wavearray)}} 
\newcommand{\zetatrans}{\upzeta^{(\textnormal{Trans--}\wavearray)}} 
\newcommand{\angktan}{\angk^{(\textnormal{Tan--}\wavearray)}} 
\newcommand{\angktrans}{\angk^{(\textnormal{Trans--}\wavearray)}}

\newcommand{\hypersurfacecontrolwave}{\mathbb{Q}} 
\newcommand{\hypersurfacecontrolwavepartial}{\mathbb{Q}^{(\textnormal{Partial})}}
\newcommand{\spacetimeintegralcontrolwave}{\mathbb{K}}
\newcommand{\spacetimeintegralcontrolwavepartial}{\mathbb{K}^{(\textnormal{Partial})}} 
\newcommand{\totalcontrolwave}{\mathbb{W}}
\newcommand{\totalcontrolwavepartial}{\mathbb{W}^{(\textnormal{Partial})}}
\newcommand{\hypersurfacecontrolVort}{\mathbb{V}}
\newcommand{\hypersurfacecontrolGradEnt}{\mathbb{S}}
\newcommand{\hypersurfacecontrolVortVort}{\mathbb{C}}
\newcommand{\hypersurfacecontrolDivGradEnt}{\mathbb{D}}
\newcommand{\toricontrolVort}{\mathbb{V}^{(\textnormal{Rough Tori})}}
\newcommand{\toricontrolGradEnt}{\mathbb{S}^{(\textnormal{Rough Tori})}}
\newcommand{\toricontrolVortVort}{\mathbb{C}^{(\textnormal{Rough Tori})}}
\newcommand{\toricontrolDivGradEnt}{\mathbb{D}^{(\textnormal{Rough Tori})}}

\newcommand{\fullymodquant}[1]{{ ^{(#1)} \mkern-1mu \mathscr{X}}}
\newcommand{\partialmodquant}[1]{{ ^{(#1)} \mkern-1mu \widetilde{\mathscr{X}}}}
\newcommand{\fullymodquantinhom}{\mathfrak{X}}
\newcommand{\partialmodquantinhom}[1]{{^{(#1)} \mkern-1mu \widetilde{\mathfrak{X}}}}

\newcommand{\hypthreearg}[3]{{^{(#3)} \mkern-2mu {\widetilde{\Sigma}_{#1}^{#2}}}}

\newcommand{\Mrough}[1]{\mathcal{M}_{#1}}
\newcommand{\twoargMrough}[2]{{^{(#2)} \mkern-3mu {\mathcal{M}_{#1}}}}
\newcommand{\PStwoargMrough}[2]{{^{(#2)} \mkern-3mu {\mathcal{M}_{#1}^{\textnormal{PS}}}}}
\newcommand{\roughtori}[1]{\widetilde{\ell}_{#1}}
\newcommand{\twoargroughtori}[2]{{^{(#2)} \mkern-2mu {\widetilde{\ell}_{#1}}}}
\newcommand{\PStwoargroughtori}[2]{{^{(#2)} \mkern-2mu {\widetilde{\ell}_{#1}^{\textnormal{PS}}}}}

\newcommand{\extendedtwoargMrough}[2]{{^{(#2)} \mkern-2mu {\widetilde{\mathcal{M}}_{#1}}}}
\newcommand{\newregionextendedtwoargMrough}[2]{{^{(#2)} \mkern-2mu {\widetilde{\mathcal{M}}_{#1}^{(\textnormal{New region})}}}}

\newcommand{\COVframe}{\mathscr{O}}
\newcommand{\COVL}{\uplambda}

\newcommand{\volRoughHypersurface}{\mathrm{d} \underline{\varpi}}
\newcommand{\volPuRoughCoordinates}{\mathrm{d} \overline{\varpi}}
\newcommand{\volroughtorus}{\mathrm{d} \varpi_{\, \, \gtorusroughfirstfund}}
\newcommand{\volroughtorusarg}[1]{\mathrm{d} \varpi_{\, \, #1}}
\newcommand{\volMRoughCoordinates}{\mathrm{d} \boldsymbol{\varpi}}
\newcommand{\volcanonical}{\mathrm{d} \textnormal{\upshape vol}}

\newcommand{\geop}[1]{\frac{\partial}{\partial #1}}
\newcommand{\roughgeop}[1]{\frac{\widetilde{\partial}}{\widetilde{\partial} #1}}

\newcommand{\Flatdiv}{\mbox{\upshape div}\mkern 1mu}
\newcommand{\Flatcurl}{\mbox{\upshape curl}\mkern 1mu}

\newcommand{\RRiemann}{\mathcal{R}_{(+)}}
\newcommand{\LRiemann}{\mathcal{R}_{(-)}}
\newcommand{\almostRiemann}{\mathcal{R}}

\newcommand{\Ent}{s}
\newcommand{\GradEnt}{S}

\newcommand{\Density}{\varrho}

\newcommand{\LogDensity}{\uprho}

\newcommand{\vortrenormalized}{\Omega}

\newcommand{\Speed}{c}

\newcommand{\VortVort}{\mathcal{C}}
\newcommand{\DivGradEnt}{\mathcal{D}}

\newcommand{\Error}{\mathrm{Error}}
\DeclareMathOperator{\ErrorIBP}{\mbox{\upshape Error}}
\newcommand{\Errortop}{\mbox{\upshape Error}_{N}^{(\textnormal{Top})}} 

\newcommand{\Errortoparg}[1]{\mbox{\upshape Error}_{#1}^{(\mbox{\upshape \footnotesize Top})}} 
\newcommand{\Errorsubcriticalarg}[1]{\mbox{\upshape Error}^{(\mbox{\upshape \footnotesize Sub-critical})}_{#1}}

\newcommand{\enmomem}{\mathbf{Q}}

\newcommand{\Ricfour}{\mathbf{Ric}}
\newcommand{\Riemfour}{\mathbf{Riem}}
\newcommand{\Riemtorus}{\widetilde{\mathfrak{Riem}}}
\newcommand{\Rictorus}{\widetilde{\mathfrak{Ric}}}
\newcommand{\Scalartorus}{\widetilde{\mathfrak{R}}}
\newcommand{\Gausstorus}{\widetilde{\mathfrak{K}}}

\newcommand{\Transport}{\mathbf{B}}

\newcommand{\muX}{\breve X}

\newcommand{\Jen}[1]{^{(#1)} \mkern-1mu \mathbf{J}}
\newcommand{\Jenarg}[2]{{^{(#1)} \mkern-1mu {\mathbf{J}^{#2}}}}

\newcommand{\deform}[1]{{^{(#1)} \mkern-1mu \pmb{\pi}}}
\newcommand{\deformarg}[3]{{^{(#1)} {\pmb{\pi}_{#2 #3}}}}

\newcommand{\angdeform}[1]{{{^{(#1)} \mkern-1mu \pi \mkern-8mu / } \,}}

\newcommand{\gfour}{\mathbf{g}}
\newcommand{\hfour}{\mathbf{h}}

\newcommand{\gtorus}{g \mkern-7.5mu / }
\newcommand{\gtorusdoublearg}[2]{{g \mkern-7.5mu / }_{#1 #2}}

\newcommand{\Chfour}{\pmb{\Gamma}}

\newcommand{\Lunit}{L}
\newcommand{\uLunit}{\underline{L}}
\newcommand{\Lgeo}{L_{(\textnormal{Geo})}}

\newcommand{\smoothtorusproject}{{\Pi \mkern-12mu / } \,\, }
\newcommand{\roughtorusproject}{\widetilde{\smoothtorusproject}}
\newcommand{\Sigmatproject}{\Pi}

\newcommand{\Dfour}{\mathbf{D}}

\newcommand{\angD}{ {\nabla \mkern-11mu / \, \,} }
\newcommand{\roughangrmD}{\widetilde{\rmD} \mkern-8mu / \,}
\newcommand{\roughangD}{\widetilde{\angD}}
\newcommand{\roughangDarg}[1]{\roughangD_{\mkern -3mu #1}}
\newcommand{\roughangDuparg}[1]{\roughangD^{\mkern -1mu #1}}
\newcommand{\angDarg}[1]{{\angD_{\mkern-3mu #1}}}
\newcommand{\angDuparg}[1]{{\angD^{\mkern-1mu #1}}}

\newcommand{\angdiv}{\mbox{\upshape{div} $\mkern-17mu /$\,\,}}
\newcommand{\roughangdiv}{\widetilde{\angdiv}}

\newcommand{\angLap}{ {\Delta \mkern-9mu / \, } }
\newcommand{\angrmD}{\rmD \mkern-8mu / \,}
\newcommand{\argangrmD}[1]{{\rmD \mkern-8mu / \,}_{#1}}
\newcommand{\newangD}{ {\nabla \mkern-10mu / \,} }
\newcommand{\newangDarg}[1]{ {\nabla \mkern-10mu / \,}_{#1} }
\newcommand{\newangDsquared}{ {\nabla \mkern-10mu / \,}^2 }
\newcommand{\newangDsquareddoublearg}[2]{ {\nabla \mkern-10mu / \,}_{#1 #2}^2 }

\newcommand{\angV}{ { {V \mkern-14mu /} \, } }
\newcommand{\angVarg}[1]{ {{V \mkern-14mu /}^{\mkern7mu #1} \, } }
\newcommand{\angk}{ { {k \mkern-7mu /} \, } }

\newcommand{\angG}{{\vec{G} \mkern -10mu /} \,\,}

\newcommand{\roughangxi}{ { {\widetilde{\upxi} \mkern-8mu /} \, } }
\newcommand{\roughangxitwouparg}[2]{ { {\widetilde{\upxi \mkern-8mu /}^{#1 #2} } \, } }

\newcommand{\weight}{\mathscr{W}}

\newcommand{\Lie}{\mathcal{L}}
\newcommand{\SigmatLie}{\underline{\mathcal{L}}}
\newcommand{\angLie}{ { \mathcal{L} \mkern-8mu / } }

\newcommand{\smoothfunction}{\mathrm{f}}

\newcommand{\Roughtoritangentvectorfieldarg}[1]{{^{(#1)} \mkern-2mu U}}

\newcommand{\CurrentboundaryerrorperfectRderivative}{\mathfrak{P}}

\newcommand{\Currentboundaryerrorhavetocontrolprincipal}{\mathfrak{E}_{(\textnormal{Principal})}}
\newcommand{\Currentboundaryerrorhavetocontrollowerorder}{\mathfrak{E}_{(\textnormal{Lower-order})}}

\newcommand{\EllipticHyperbolicCurrentIntegralIdentityTotalSpacetimeErrorTerm}{\mathfrak{M}}

\newcommand{\Nullhypersurfaceproject}{\overline{\Pi}}
\newcommand{\Nullhypersurfacemetric}{\overline{e}}
\newcommand{\Nullhypersurfaceinversemetric}{\overline{E}}

\newcommand{\SigmatTan}{V}

\newcommand{\ehcurrent}{\mathscr{J}}

\newcommand{\ellipticCoerciveQuadratic}{\mathscr{Q}}

\newcommand{\Singletan}{P}

\newcommand{\InterestingCHOV}{{^{(\textnormal{Interesting})} \mkern-1mu \mathscr{T}}}  
\newcommand{\InverseInterestingCHOV}{{^{(\textnormal{Interesting})} \mkern-1mu \mathscr{T}^{-1}}}  

\newcommand{\MLeft}{\mathcal{M}_{\textnormal{Left}}}
\newcommand{\MRight}{\mathcal{M}_{\textnormal{Right}}}
\newcommand{\MInteresting}{\mathcal{M}_{\textnormal{Interesting}}}
\newcommand{\MSingular}{\mathcal{M}_{\textnormal{Singular}}}

\newcommand{\smallneighborhoodofcreasetwoarg}[2]{{^{(#2)} \mkern-2mu {\mathcal{N}_{#1}}}}

\newcommand{\domainofgraphofmulevelsetinsingularregion}[2]{\mathscr{D}_{#2}^{#1}}

\newcommand{\mulevelsetvalue}{\mathfrak{m}}
\newcommand{\muxmulevelsetvalue}{\mathfrak{n}}

\newcommand{\tisafunctiononlevelsetsofmu}[1]{T_{#1}}

\newcommand{\apriorimain}{\mathscr{F}}
\newcommand{\aprioripartial}{\mathscr{G}}
\newcommand{\apriorilower}{\mathscr{H}}

\newcommand{\nullgeneratorofsingularboundary}{Q}
\newcommand{\embeddingofsingularboundaryintogeometriccoordinatespace}{\mathscr{S}}

\newcommand{\muderivativevectorfield}{J}
\newcommand{\muxmuderivativevectorfield}{K}
\newcommand{\derivativevectorfieldatfixemuandmuxmu}[1]{P_{(#1)}}

\newcommand{\PSLmunottoonegativeparameter}{\mathfrak{p}}

\newcommand{\PSCHOVJacobianroughtomumuxmu}[1]{{^{(#1)} \mkern-1mu \mathbf{J}_{\textnormal{PS}}}}
\newcommand{\PSInverseCHOVJacobianroughtomumuxmu}[1]{{^{(#1)} \mkern-1mu {\mathbf{J}_{\textnormal{PS}}^{-1}}}}

\newcommand{\PSCHOVgeotorough}[1]{{^{(#1)} \mkern-1mu \mathscr{T}_{\textnormal{PS}}}}
\newcommand{\PSInverseCHOVgeotorough}[1]{{^{(#1)} \mkern-1mu {\mathscr{T}_{\textnormal{PS}}^{-1}}}}

\newcommand{\PSCartesiantisafunctiononlevelsetsofroughtimefunctionarg}[2]{\mathfrak{t}_{#1,#2}^{\textnormal{PS}}}

\newcommand{\PSthirdorderTaylorremaindercoefficientfunction}{\mathcal{F}}

\newcommand{\PSUpsilon}{{\Upsilon_{\textnormal{PS}}}}
\newcommand{\PSInverseUpsilon}{{\Upsilon_{\textnormal{PS}}^{-1}}}

\newcommand{\PSboringregionmupositive}{\mulevelsetvalue_1^{\textnormal{PS}}}
\newcommand{\PSWtransarg}[1]{{ ^{(#1)} \mkern-2mu \breve{W}^{\textnormal{PS}}}}

\newcommand{\PSdataamplitude}{\mathfrak{a}}
\newcommand{\PSdatamuHessianTaylorcoefficient}{\mathfrak{b}}

\newcommand{\PSBigDelta}{\Delta^{\textnormal{PS}}}

\newcommand{\PSinterestingusmallmultipleofamplitude}{\upzeta}

\newcommand{\PSinitialkeymufunctioncoefficients}{F}

\newcommand{\PSLmusourcetermfunction}{H}
\newcommand{\antiderivativePSLmusourcetermfunction}{\mathfrak{H}}

\newcommand{\embeddingofmuadapatedtorinCartesianspace}[1]{{^{(#1)} \mkern-2mu \mathfrak{I}}}

\newcommand{\Rtransnormsmallfactorarg}[1]{{^{(#1)} \mkern-2mu r}}

\newcommand{\almostRiemannfunction}{F}

\newcommand{\ambient}{\mathscr{A}}

\newcommand{\twoargmumuxtorus}[2]{\breve{\mathbf{T}}_{#1,#2}}

\newcommand{\Cartesiantisafunctiononlevelsetsofroughtimefunctionarg}[2]{\mathfrak{t}_{#1,#2}}
\newcommand{\Cartesiantisafunctiononmumxtoriarg}[2]{\mathfrak{T}_{#1,#2}}
\newcommand{\Eikonalisafunctiononmumuxtoriarg}[2]{\mathfrak{U}_{#1,#2}}

\newcommand{\embeddatahypersurfacearg}[1]{{^{(#1)}E}}
\newcommand{\extendedembeddatahypersurface}{\mathscr{E}}

\newcommand{\inthyp}[2]{{^{(\textnormal{Interesting})}\Sigma_{#1}^{#2}}}

\newcommand{\datahypfortimefunctionarg}[1]{\breve{\mathbb{X}}_{#1}}
\newcommand{\datahypfortimefunctiontwoarg}[2]{\breve{\mathbb{X}}_{#1}^{#2}}
\newcommand{\PSdatahypfortimefunctiontwoarg}[2]{{^{\textnormal{PS}} \mkern-1mu \breve{\mathbb{X}}_{#1}^{#2}}}

\newcommand{\mulevelsetarg}[1]{\breve{\mathbb{M}}_{#1}}
\newcommand{\mulevelsettwoarg}[2]{{\breve{\mathbb{M}}_{#1}^{#2}}}

\newcommand{\domainforembeddingdatahypfortimefunctiontwoarg}[2]{{^{(#1)} \mkern-2mu \mathscr{H}_{#2}}}

\newcommand{\embeddingdatahypfortimefunctionarg}[1]{{^{(#1)} \mkern-2mu H}}

\newcommand{\scalarembeddingdatahypfortimefunctionarg}[1]{{^{(#1)} \mkern-2mu h}}

\newcommand{\flowmapWtransargtwoarg}[2]{{^{(#1)} \mkern-2mu \iota_{#2}}}

\newcommand{\composedflowmapdiffeoarg}[1]{{^{(#1)}\mkern-2mu F}}
\newcommand{\domaincomposedflowmapdiffeoarg}[1]{{^{(#1)} \mkern-2mu \mathscr{F}}}
\newcommand{\flowmapWtransMatrixarg}[1]{{^{(#1)} \mkern-2mu M}}

\newcommand{\nullhyparg}[1]{\mathcal{P}_{#1}}

\newcommand{\nullhypthreearg}[3]{{^{(#1)} \mkern-2mu \mathcal{P}_{#2}^{#3}}}
 
\newcommand{\multipliervectorfield}{\breve{T}}

\newcommand{\Fullset}{\mathscr{Z}}
\newcommand{\Tanset}{\mathscr{P}}
\newcommand{\Angularset}{\mathscr{Y}}

\newcommand{\CHOVgeotorough}[1]{{^{(#1)} \mkern-2mu \mathscr{T}}}
\newcommand{\InverseCHOVgeotorough}[1]{{^{(#1)} \mkern-2mu {\mathscr{T}^{-1}}}}
\newcommand{\CHOVroughtomumuxmu}[1]{{^{(#1)} \mkern-1mu \Phi}}
\newcommand{\InverseCHOVroughtomumuxmu}[1]{{^{(#1)} \mkern-1mu {\Phi^{-1}}}}
\newcommand{\CHOVJacobianroughtomumuxmu}[1]{{^{(\CHOVroughtomumuxmu{#1})} \mkern-2mu \mathbf{J}}}
\newcommand{\InverseCHOVJacobianroughtomumuxmu}[1]{{^{(\CHOVroughtomumuxmu{#1})} \mkern-2mu {\mathbf{J}^{-1}}}}
\newcommand{\CHOVgeotomumuxmu}{\breve{\mathscr{M}}}
\newcommand{\InverseCHOVgeotomumuxmu}{\breve{\mathscr{M}}^{-1}}
\newcommand{\PSCHOVgeotomumuxmu}{\breve{\mathscr{M}}_{\textnormal{PS}}}

\newcommand{\CHOVJacobiangeotomumuxmu}{{^{(\CHOVgeotomumuxmu)} \mkern-2mu \mathbf{J}}}
\newcommand{\InverseCHOVJacobiangeotomumuxmu}{{^{(\CHOVgeotomumuxmu)} \mkern-2mu {\mathbf{J}^{-1}}}}
\newcommand{\PSCHOVJacobiangeotomumuxmu}{{^{(\PSCHOVgeotomumuxmu)} \mkern-2mu \mathbf{J}}}
\newcommand{\PSInverseCHOVJacobiangeotomumuxmu}{{^{(\PSCHOVgeotomumuxmu)} \mkern-2mu {\mathbf{J}^{-1}}}}
\newcommand{\nullform}{\mathfrak{Q}}

\newcommand{\mainnullform}{\mathfrak{M}}

\newcommand{\upmuboot}{\mulevelsetvalue_{\textnormal{Boot}}}

\newcommand{\FlowmapLrougharg}[1]{{^{(#1)} \mkern-2mu \widetilde{\Lambda}}}
\newcommand{\InverseFlowmapLrougharg}[1]{{^{(#1)} \mkern-2mu {\widetilde{\Lambda}^{-1}}}}
\newcommand{\FlowmapLroughtwoarg}[2]{{^{(#1)} \mkern-2mu \widetilde{\Lambda}^{#2}}}

\newcommand{\argLrough}[1]{{^{(#1)} \mkern-2mu \widetilde{L}}}
\newcommand{\twoargLrough}[2]{{^{(#1)} \mkern-2mu {\widetilde{L}^{#2}}}}

\newcommand{\mydiam}{{\mkern-1mu \scaleobj{.75}{\blacklozenge}}}

\newcommand{\leftu}{U_2}
\newcommand{\interestingu}{U_{\mbox{\tiny \Radioactivity}}}
\newcommand{\rightu}{U_1}
\newcommand{\farrightu}{U_0}

\newcommand{\secondtransversalderivativemulowerbound}{M_2}
\newcommand{\secondtransversalderivativemulowerboundPS}{M_2^{\textnormal{PS}}}

\newcommand{\csspacelikehypersurface}{\underline{\mathcal{S}}}

\numberwithin{equation}{section}
 
\begin{document}
\title{The emergence of the singular boundary from the crease in $3D$ compressible Euler flow}
\author[LA,JS]{Leonardo Abbrescia$^{* \dagger}$ and Jared Speck$^{** \dagger\dagger}$}

\thanks{$^{*}$Vanderbilt University, Nashville, TN, USA.
\texttt{leonardo.abbrescia@vanderbilt.edu}}

\thanks{$^{**}$Vanderbilt University, Nashville, TN, USA.
\texttt{jared.speck@vanderbilt.edu}}

\thanks{$^\dagger$ LA gratefully acknowledges support from an NSF Postdoctoral Fellowship.}

\thanks{$^{\dagger\dagger}$JS gratefully acknowledges support from NSF grant \# DMS-2054184 and NSF CAREER grant \# DMS-1914537.
}

\begin{abstract}
We study the Cauchy problem for the $3D$ compressible Euler equations under an
arbitrary equation of state with positive speed of sound, 
aside from that of a Chaplygin gas.
For open sets of smooth initial data with non-trivial vorticity and entropy, 
our main results yield a constructive proof of the formation, structure, and stability of the singular boundary,
which is the set of points where the solution forms a shock singularity, i.e., where some 
first-order Cartesian coordinate partial derivatives of the velocity and density blow up.
We prove that in the solution regime under study,
the singular boundary has the structure of a degenerate, acoustically null 
$3D$ submanifold-with-boundary. Our approach yields the full structure of a neighborhood 
of a connected component of the crease, which is a $2D$  
acoustically spacelike
submanifold equal to the past boundary of the singular boundary.
In the study of shocks, the crease plays the role of the ``true initial singularity''
from which the singular boundary emerges,
and it is a crucial ingredient for setting up the shock development problem.
These are the first results revealing the totality of these structures
without symmetry, irrotationality, or isentropicity assumptions. Moreover,
even within the sub-class of irrotational and isentropic solutions,
these are the first constructive results revealing these structures
without a strict convexity assumption on the shape of the singular boundary.
Our proof relies on a new method: the construction of rough foliations
of spacetime, dynamically adapted to the exact shape of the singular boundary
and crease, where the latter is provably two degrees less differentiable than the fluid.
Our results also set the stage for our forthcoming paper,
in which we will prove the emergence and stability of a Cauchy horizon,
which emanates from the crease and ``evolves'' in a direction that is ``opposite''
the singular boundary in a sense determined by the intrinsic acoustic geometry of the flow.

\bigskip

\bigskip

\end{abstract}

\maketitle

 \tableofcontents

\newpage

\section{Introduction} \label{S:INTRO}
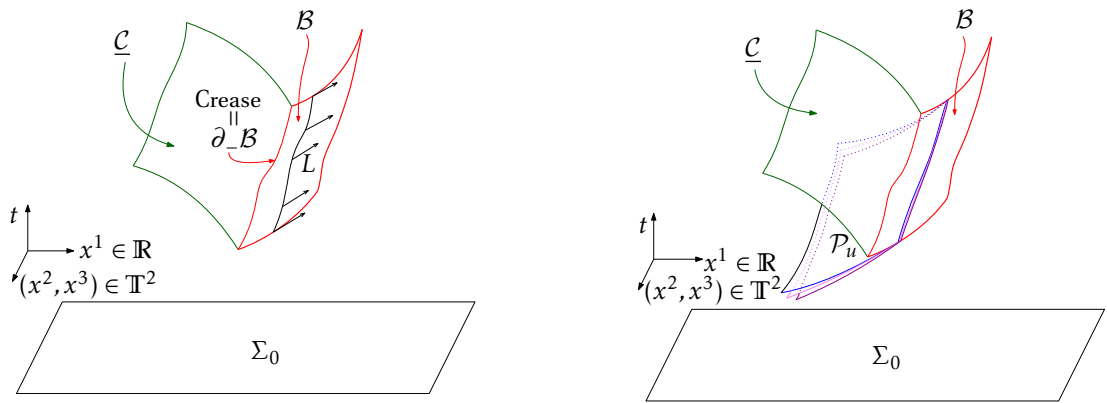
\begin{figure}[ht]
\centering
\begin{subfigure}{.5\textwidth}
  \centering
  	\begin{overpic}[scale=.36, grid = false, tics=5, trim=-.5cm 0cm -1cm -.5cm, clip]{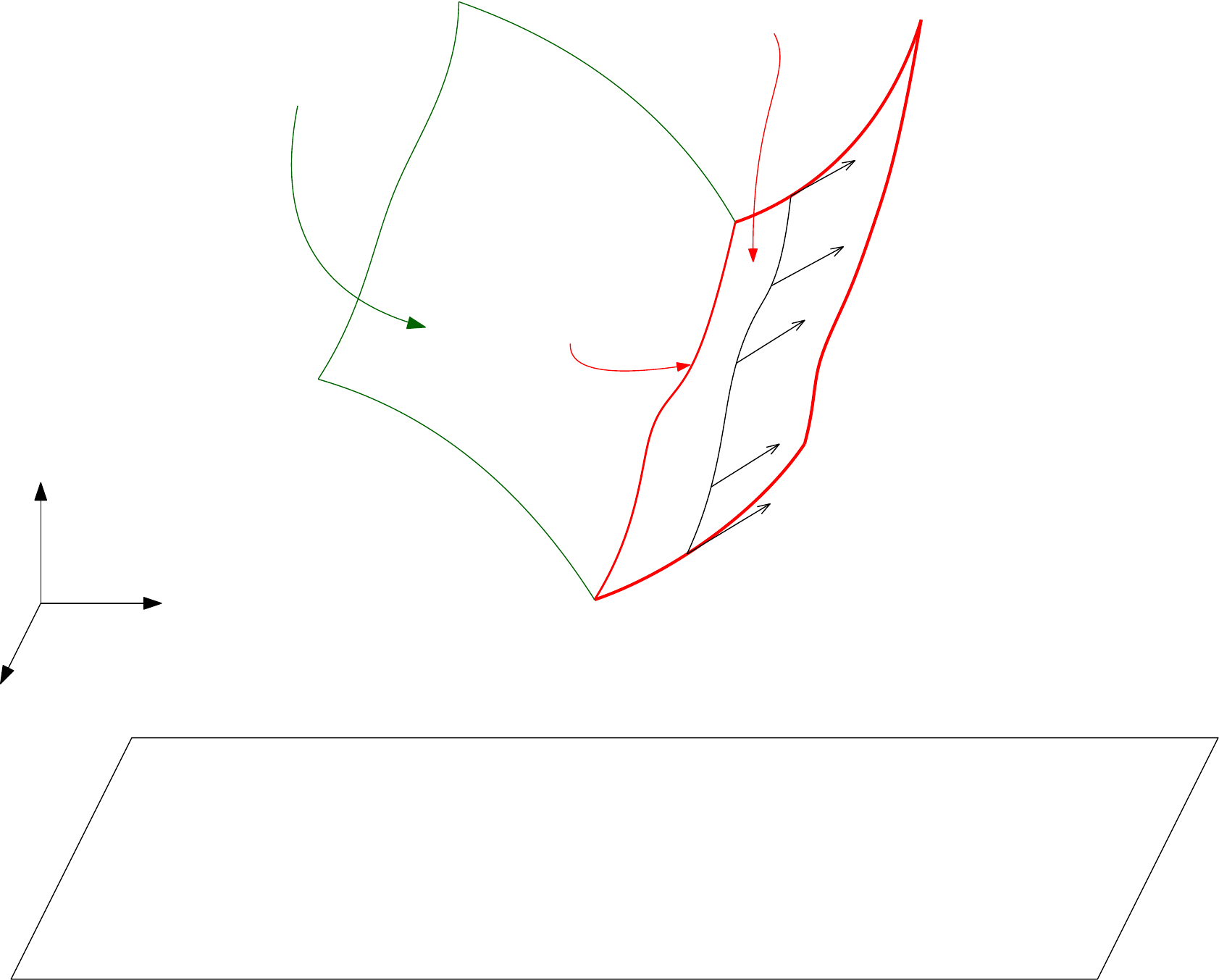}
			\put (50,7) {$\Sigma_0$}
			\put (4,20) {$(x^2,x^3) \in \mathbb{T}^2$}
			\put (2,34) {$t$}
			\put (16,27) {$x^1 \in \mathbb{R}$}
			\put (39,57) {$\mbox{\upshape Crease}$}
			\put (46,53) {\rotatebox{90}{$=$}}
			\put (42,49) {$\partial_- \mathcal{B}$}
			\put (60,44) {$\Lunit$}
			\put (59,72.5) {$\mathcal{B}$}
			\put (23,68.5) {$\underline{\mathcal{C}}$}
		\end{overpic}
		\caption{A localized subset of the maximal classical development in Cartesian coordinate space}
		\label{F:MAXDEVELOPMENTINCARTESIAN}
\end{subfigure}%
\begin{subfigure}{.5\textwidth}
 \centering
\begin{overpic}[scale=.36,grid=false]{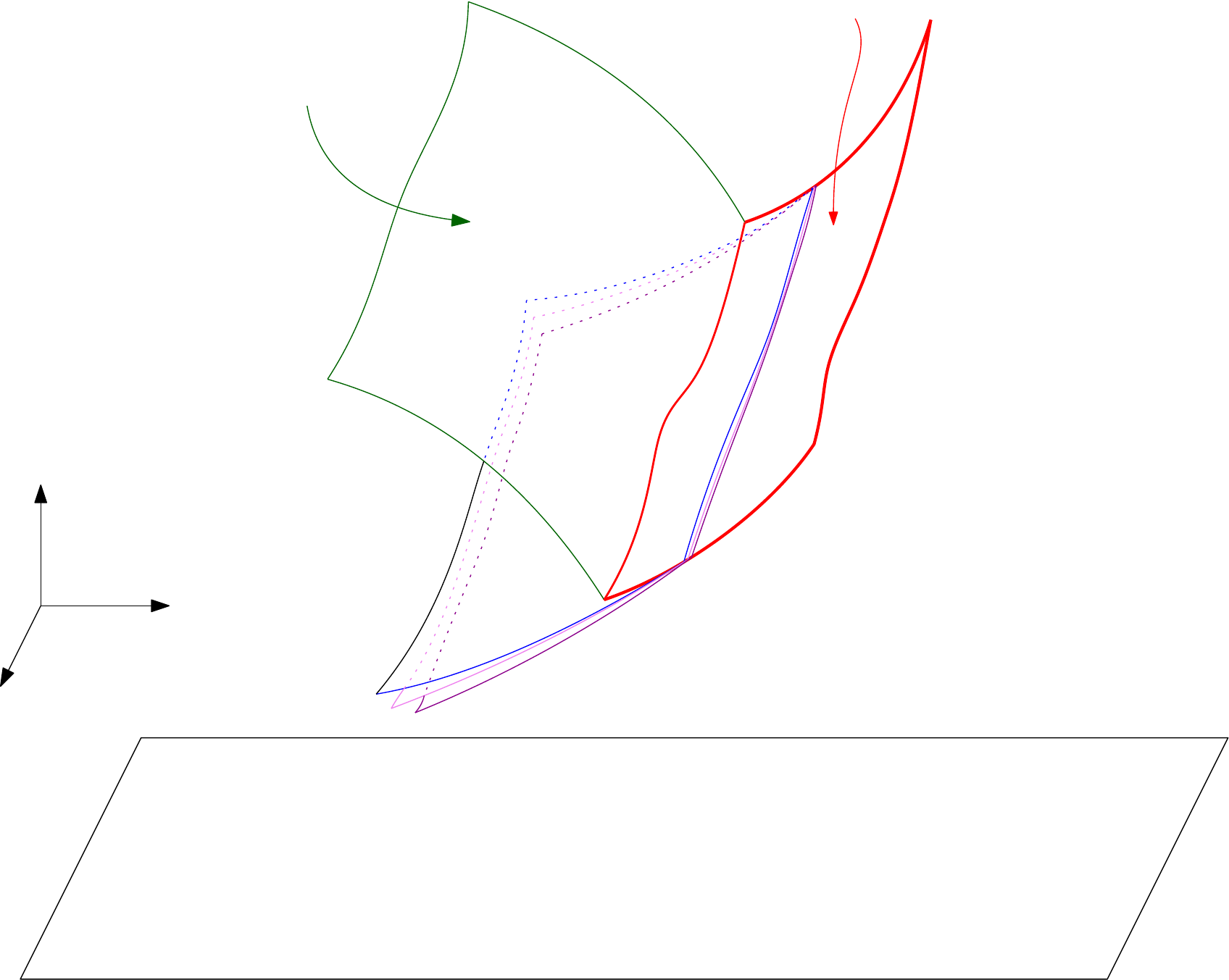} 
			\put (50,8) {$\Sigma_0$}
			\put (1.2,22) {$(x^2,x^3) \in \mathbb{T}^2$}
			\put (0,36) {$t$}
			\put (41,32) {$\nullhyparg{u}$}
			\put (14,28) {$x^1 \in \mathbb{R}$}
			\put (68,80) {$\mathcal{B}$}
			\put (23,74) {$\underline{\mathcal{C}}$}
\end{overpic}
  \caption{Infinite density of the characteristics $\nullhyparg{u}$ on $\mathcal{B}$}
  \label{F:INFININTEDENSITYOFCHARACTERISTICSONSINGULARBOUNDARY}
\end{subfigure}
\caption{Cartesian coordinate space illustrations of the main results}
\label{F:INTROPICTURESOFMAINRESULTS}
\end{figure}

	This is the first of two papers in which we construct a large (though compact) portion of
	the maximal classical globally hyperbolic development 
	(which we refer to as the ``maximal development'' for short from now on)
	-- up to the boundary --
	of the initial data for open sets of initially smooth, \emph{shock-forming} solutions 
	to the $3D$ compressible Euler equations with vorticity and dynamic entropy,
	and without symmetry assumptions.
	Roughly speaking, the maximal classical development
	is the largest possible classical solution 
	+ 
	region that is determined by the initial data.
	Our main results are the constructions of the (singularity-forming) solution
	$+$ localized region depicted in Fig.\,\ref{F:INTROPICTURESOFMAINRESULTS} 
	(see Remark~\ref{R:UPSILONNOTATIION} for comments on our notation in the figure)
	and a proof of their stability under Sobolev-class perturbations of the initial data on 
	the Cauchy hypersurface $\Sigma_0 \eqdef \lbrace t = 0 \rbrace$.
	\begin{quote}
	Our papers provide the first results that 
	construct and fully justify Fig.\,\ref{F:INTROPICTURESOFMAINRESULTS} for open sets of initial data
	without symmetry, irrotationality, or isentropicity assumptions.
	\end{quote}
	
	The present paper concerns the analysis up to the \emph{singular boundary}, 
	which we denote by\footnote{More precisely,
	we follow the solution up to a compact portion of the singular boundary that we denote by
	$\mathcal{B}^{[0,\muxmulevelsetvalue_0]}$ in our main theorems. \label{FN:SINGULARBOUNDARYPORTION}} 
	``$\mathcal{B}$,''
	while our companion work \cite{lAjS20XX} concerns the analysis up to the \emph{Cauchy horizon}, 
	which we denote by ``$\underline{\mathcal{C}}$.'' 
	Roughly, $\mathcal{B}$ is the submanifold-with-boundary\footnote{It is not obvious that the set of blowup-points 
	has the structure of a submanifold-with-boundary in Cartesian coordinate space. 
	Indeed, uncovering this structure is one of the main results of 
	the paper. While this structure holds for open sets of solutions, including the solutions we handle
	in this paper, for other solutions, the set of blowup-points might fail to have the structure of a 
	submanifold-with-boundary.
	\label{FN:MANIFOLDNOTOBVIOUS}} 
	of points where the fluid's first derivatives blowup
	(though the fluid variables themselves remain bounded),
	while $\underline{\mathcal{C}}$ is a future boundary that ``feels the influence'' of
	the singularity, even though the solution remains classical along $\underline{\mathcal{C}} \backslash \partial_- \mathcal{B}$.
	Readers can jump to Theorem~\ref{T:ABBREVIATEDSTATEMENTOFMAINRESULTS}
	for an abbreviated version of the main results of this paper, 
	and Theorems~\ref{T:EXISTENCEUPTOTHESINGULARBOUNDARYATFIXEDKAPPA} and \ref{T:DEVELOPMENTANDSTRUCTUREOFSINGULARBOUNDARY}
	for precise, extended statements.
	As we describe in Sect.\,\ref{SS:THETWOOPENPROBLEMS}, the papers resolve several open problems 
	and allow one to properly set up the \emph{shock development problem},
	which is the problem of (locally) describing the transition of the solution
	from classical to weak, past the ``initial singularity,''
	which is the acoustically\footnote{The word ``acoustic'' refers to  
	the acoustical metric $\gfour$ of definition~\eqref{E:ACOUSTICALMETRIC}, 
	which is the solution-dependent Lorentzian metric that dictates the geometry of sound waves.
	Throughout the paper, all Lorentzian geometric notions such as spacelike, 
	timelike, null, etc. are with respect to $\gfour$. \label{FN:ACOUSTICALMETRIC}} 
	spacelike co-dimension $2$ 
	submanifold that we denote by ``$\partial_- \mathcal{B}$'' 
	(the past boundary of $\mathcal{B}$) in the figure. 
	In full generality, the shock development problem is open, though there has been
	inspiring progress, which we describe in Sect.\,\ref{SSS:PROGRESSONSHOCKDEVELOPMENT}.
	It is crucially important, for example in setting up the shock development problem, 
	that our approach yields a complete description of the
	initial singularity. We highlight the following key point:	
	\begin{quote}
		The ``initial singularity'' should \underline{not} be thought
		of as a point in spacetime, but rather the set
		$\partial_- \mathcal{B}$ depicted in Fig.\,\ref{F:MAXDEVELOPMENTINCARTESIAN}, which we refer to 
		as \emph{the crease}, and which we prove has the structure of a $2$-dimensional acoustically spacelike
		submanifold in the solution regime under study.
	\end{quote}
	
\renewcommand{\thesubfigure}{\Alph{subfigure}} 

While we treat in detail a specific regime in our two papers 
-- perturbations of simple isentropic plane-symmetric solutions --
the methods we develop are robust and could be applied to other regimes,
such as perturbations of non-vacuum steady state solutions in $\mathbb{R}^{1+3}$.
We already stress that our results apply to situations in
which $\mathcal{B}$ fails to be strictly convex 
(that is, whenever it fails to be strictly concave up),
and that substantial new ideas are needed in this case to follow the solution up to $\mathcal{B}$
and to understand its structure; 
in Fig.\,\ref{F:INTROPICTURESOFMAINRESULTS} we have depicted
a $\mathcal{B}$ that, while ``convex in the $x^1$-direction,'' it
fails to be convex ``in the $(x^2,x^3)$-directions.''

\subsection{Outline of the remainder of the introduction}
	\label{SS:INTROOUTLINE}
	In Sect.\,\ref{SS:FIRSTVERSIONOFEQUATIONS}, we introduce the compressible Euler equations,
	though not in the form we use to prove our main results.
	In Sect.\,\ref{SS:OVERVIEWOFDIFFICULTIES}, we highlight some the main 
	challenges in the proofs of our main results and 
	discuss a few of the most important new ideas we use to overcome them.
	In Sect.\,\ref{SS:THETWOOPENPROBLEMS}, we provide an overview of the two open problems
	that we resolve in the paper, and we describe how they are connected to the shock development problem.
	In Sect.\,\ref{SS:ABBREVIATEDSTATEMENTOFMAINRESULTS}, we state Theorem~\ref{T:ABBREVIATEDSTATEMENTOFMAINRESULTS},
	which is a first, somewhat informal, abbreviated version of our main results; see
	Theorems~\ref{T:EXISTENCEUPTOTHESINGULARBOUNDARYATFIXEDKAPPA} and \ref{T:DEVELOPMENTANDSTRUCTUREOFSINGULARBOUNDARY}
	for the extended, precise statements.
	In Sect.\,\ref{SS:REMARKSONRESULTSANDMETHODS}, 
	having stated Theorem~\ref{T:ABBREVIATEDSTATEMENTOFMAINRESULTS},
	we provide some extended remarks on the ideas and methods we use in the proof.
	In Sect.\,\ref{SS:WORKSBUILDINGTOWARDSMAINTHEOREM}, 
	we describe the most relevant precursor results to this paper,
	focusing on those works that yielded methods that we use here.
	In Sect.\,\ref{SS:BACKGROUNDONSHOCKS}, we discuss the history of the study of shock formation
	and highlight some important developments in the subject.
	In Sect.\,\ref{SS:IDEASBEHINDPROOF}, we provide an overview of 
	the main ideas of the proofs of our main results, in particular highlighting
	various technical issues that were not discussed earlier in the introduction.
	In Sect.\,\ref{SS:PAPEROUTLINE}, we provide an outline of the remainder of the paper.

\subsection{First version of the equations and local well-posedness}		
\label{SS:FIRSTVERSIONOFEQUATIONS}
	Our main results concern the Cauchy problem for the $3D$ compressible Euler equations,
	which can be formulated as a quasilinear hyperbolic PDE system
	in the velocity $v : \mathbb{R} \times \Sigma \rightarrow \mathbb{R}^3$,
	the density $\varrho : \mathbb{R} \times \Sigma \rightarrow [0,\infty)$,
	and the entropy $\Ent : \mathbb{R} \times \Sigma \rightarrow \mathbb{R}$.
	In this paper, $\Sigma \eqdef \mathbb{R} \times \mathbb{T}^2$ denotes the ``space manifold,''
	that is, we assume that ``space'' is diffeomorphic to $\mathbb{R} \times \mathbb{T}^2$.
	In our setup, the Cartesian coordinates on spacetime are $(t,x^1,x^2,x^3)$,
	where $t$ is the standard Cartesian time function,
	$x^1$ denotes the standard Cartesian coordinate on $\mathbb{R}$, and $(x^2,x^3)$ denote
	standard Cartesian coordinates on $\mathbb{T}^2 \eqdef [-\pi,\pi]^2$ (with the endpoints identified).
	The spatial topology $\mathbb{R} \times \mathbb{T}^2$ 
	allows us to simplify various aspects of our approach, leading to a cleaner
	presentation of the analysis. However, our analysis
	is local in spacetime and, with modest additional effort, 
	all our results could be extended to other spatial topologies such
	as $\mathbb{R}^3$, without cut-off functions or partitions of unity.
	Relative to the Cartesian coordinates $(t,x^1,x^2,x^3)$, 
	the $3D$ compressible
	equations can be expressed as follows\footnote{Throughout, if $\mathbf{V}$ is a vectorfield and $f$ is a scalar function,
	then $\mathbf{V} f \eqdef \mathbf{V}^{\alpha} \partial_{\alpha}f$ denotes the derivative of $f$
	in the direction of $\mathbf{V}$. \label{FN:VECTORFIELDDERIVATIVEOFSCALARFUNCTION}} 
	(see \cite{dCsM2014} for background on the equations):
\begin{subequations}
\begin{align}
	\Transport v^i 
	& = 
	- \frac{\partial_i p}{\varrho},
	&& (i=1,2,3),
	\label{E:INTROTRANSPORTVI}
		\\
		 \label{E:INTROTRANSPORTDENSITY}
	\Transport \varrho
	& = - \varrho \Flatdiv v,
	&&
		\\
	\Transport \Ent
	& = 0,
	&&
\label{E:INTROBS}
\end{align}
\end{subequations}
where 
$p$ is the pressure,
$\Transport$ denotes the \emph{material vectorfield}:
\begin{align}\label{E:MATERIALDERIVATIVEVECOTRFIELD}
\Transport 
& \eqdef \p_t + v^a \p_a,
\end{align}
and $\Flatdiv$ is the standard Euclidean divergence operator (see Def.\,\ref{D:EUCLIDEANDIVERGENCEANDCURL}).
To close the equations, we assume an \emph{equation of state} $p = p(\varrho,\Ent)$.
Our results apply for \emph{any} sufficiently smooth equation of state 
-- except for that of the Chaplygin gas (see Sect.\,\ref{SSS:LOGDENSITYASSUMPTIONSONEOSANDNORMALIZATIONS}) --
with positive sound speed $\Speed$ defined by:
	\begin{align} \label{E:SOUNDSPEED}
			\Speed
			& \eqdef \sqrt{p_{;\varrho}},
\end{align}
where $p;_{\varrho}$ is the partial derivative of the equation of state with respect
to the density at fixed entropy.

We assume that smooth initial data 
$(v,\varrho,\Ent)|_{\Sigma_0}$
for \eqref{E:INTROTRANSPORTVI}--\eqref{E:INTROBS} 
are prescribed along the spacelike hypersurface
$\Sigma_0 \eqdef \lbrace t = 0 \rbrace = \lbrace 0 \rbrace \times \mathbb{R} \times \mathbb{T}^2$.
We consider only initial data such that $\varrho|_{\Sigma_0} > 0$, thereby avoiding the severe
degeneracies that can occur at fluid-vacuum boundaries.
It is well-known that the equations \eqref{E:INTROTRANSPORTVI}--\eqref{E:INTROBS}
are locally well-posed for
non-vacuum initial data on $\Sigma_0 \eqdef \lbrace t= 0 \rbrace$
such that $(v,\varrho,\Ent)|_{\Sigma_0} \in H^3(\Sigma_0)$.
To follow the solution to the singular boundary, we assume that the data belong
to a sufficiently high order Sobolev space, where different solution variables
have distinct, directionally dependent amounts of regularity;
see Sect.\,\ref{S:ASSUMPTIONSONTHEDATA} for our detailed assumptions.
All of the known approaches to studying shocks away from symmetry rely on
the assumption that the data belong to a high order
Sobolev space. This is due to possibly singular energy estimates at the high
derivative levels -- even in the ``good'' geometric
coordinate system $(t,u,x^2,x^3)$, described below, which ``hides'' the singularity
at the low to mid derivative levels; see Sects.\,\ref{SSS:REGULARESTIMATESONROUGHFOLIATIONS} 
and \eqref{SSS:GEOMETRICENERGYESTIMATESONTHEROUGHFOLIATIONS} 
for further discussion of these fundamental technical issues.

\begin{quote}
We stress up front that our analysis crucially relies on a 
geometric reformulation of \eqref{E:INTROTRANSPORTVI}--\eqref{E:INTROBS} 
as a system of covariant wave equations coupled to transport-div-curl equations,
where the nonlinear terms exhibit remarkable null and regularity properties. 
The reformulation was derived in \cite{jS2019c} (see also the precursor \cite{jLjS2020a}),
and we recall it in Theorem~\ref{T:GEOMETRICWAVETRANSPORTSYSTEM}.
\end{quote}

\subsection{An overview of the degeneracies and difficulties in the problem}
\label{SS:OVERVIEWOFDIFFICULTIES}
The analysis needed to fully justify Fig.\,\ref{F:INTROPICTURESOFMAINRESULTS} is fraught with degeneracies and difficulties.
While many prior works on shocks have constructed the solution in \emph{strict subsets}
of the region in Fig.\,\ref{F:MAXDEVELOPMENTINCARTESIAN}, our papers are the first 
to fully grapple with the degeneracies and 
construct the solution in the entire region.
To handle solutions with non-zero vorticity and dynamic entropy,
we rely on an arsenal of geometric and analytic techniques, developed in
earlier works \cites{dC2007,jLjS2018,jS2019c,jLjS2020a,lAjS2020},
combined with key new ideas that we describe below.
Here, we highlight some of the main challenges in the analysis and
mention some of the methods we use to overcome them.
\begin{itemize} 
	\item (\textbf{Singularities}). The solution forms a shock along
		the submanifold-with-boundary $\mathcal{B}$, i.e., some of its first-order partial derivatives with respect
		to the Cartesian coordinates blow up, though the solution itself remains bounded.
		The blowup-dynamics are extremely rich: some quantities exhibit blowup of their derivatives
		-- but only derivatives in certain directions\footnote{Roughly, along $\mathcal{B}$, many quantities' 
		derivatives in directions \emph{transversal} to the characteristics blow up,
		while their tangential derivatives remain bounded, much like in the simple case of 
		Burgers' equation in $1D$: $\partial_t \Psi + \Psi \partial_x \Psi = 0$.} -- while other quantities 
		\underline{and} their derivatives in all directions remain bounded.
		We refer to Theorem~\ref{T:DEVELOPMENTANDSTRUCTUREOFSINGULARBOUNDARY} for the details.
	\item (\textbf{Nonlinear geometric optics and geometric coordinates}) As in many prior works on shock formation,
		to follow the solution up to $\mathcal{B}$ and to obtain a precise understanding of the singularity formation,
		we cannot rely on the Cartesian coordinates, which are not adapted to the singularity. 
		Instead, we rely on nonlinear geometric optics.
		Specifically, nonlinear geometric optics yields a ``geometric coordinate system'' 
		(see Def.\,\ref{D:GEOMETRICCOORDIANTESANDPARTIALDERIVATIVEVECTORFIELDS}) 
		-- one that is globally homeomorphic (in the compact region under study) 
		to the Cartesian coordinate system but not diffeomorphic to it along $\mathcal{B}$ --
		relative to which the solution remains rather smooth; see Sect.\,\ref{SS:REMARKSONRESULTSANDMETHODS}.
		The key ingredient in implementing nonlinear geometric optics
		is an \emph{eikonal function} $u$ solving the acoustic \emph{eikonal equation}
		$(\gfour^{-1})^{\alpha \beta} \partial_{\alpha} u \partial_{\beta} u = 0$,
		where the \emph{acoustical metric} $\gfour$ is the solution-dependent Lorentzian metric (see \eqref{E:ACOUSTICALMETRIC}) 
		that captures the intrinsic geometry of sound waves.
		The level sets of $u$, which we denote by $\nullhyparg{u}$, are characteristic for the compressible Euler equations;
		roughly, the $\nullhyparg{u}$ represent surfaces along which sound waves can propagate.
		In the present paper, our \emph{geometric coordinate system} is $(t,u,x^2,x^3)$, where $u$ is the 
		eikonal function and $t,x^2,x^3$ are standard Cartesian coordinate functions.
	\item (\textbf{Degeneracies in the acoustic geometry}) 	
		The fluid singularity along $\mathcal{B}$ is intimately tied to degeneracies in the acoustic geometry; 
		roughly, along $\mathcal{B}$, the characteristic hypersurfaces $\nullhyparg{u}$
		(which we also refer to as ``characteristics,''
		``null hypersurfaces,'' ``acoustically null hypersurfaces,'' or ``$\gfour$-null hypersurfaces''), 
		develop \emph{infinite density} along $\mathcal{B}$. In the present paper, 
		the infinite density of the $\nullhyparg{u}$ is characterized
		by the vanishing of a function $\upmu$, the inverse foliation density, which we describe below in detail;
		$\upmu$ is positive in the maximal development, except along $\mathcal{B}$.
		We depict this infinite density in Fig.\,\ref{F:INFININTEDENSITYOFCHARACTERISTICSONSINGULARBOUNDARY}, where we show
		three distinct null hypersurfaces ``piling up'' along $\mathcal{B}$.
		It is important to appreciate that in the solution regime we are studying, 
		\emph{distinct characteristic hypersurfaces, viewed as submanifolds in Cartesian coordinate space,
		never actually intersect}\footnote{This fact follows from the homeomorphism property
		of the map $\Upsilon$, as described in Theorem~\ref{T:ABBREVIATEDSTATEMENTOFMAINRESULTS}.} 
		on $\mathcal{B}$, even though their density becomes infinite.
		This phenomenon is crucial for properly setting up the shock development problem.
		
	\begin{center}
	\begin{figure}  
		\begin{overpic}[scale=.6, grid = false, tics=5, trim=-.5cm -1cm -1cm -.5cm, clip]{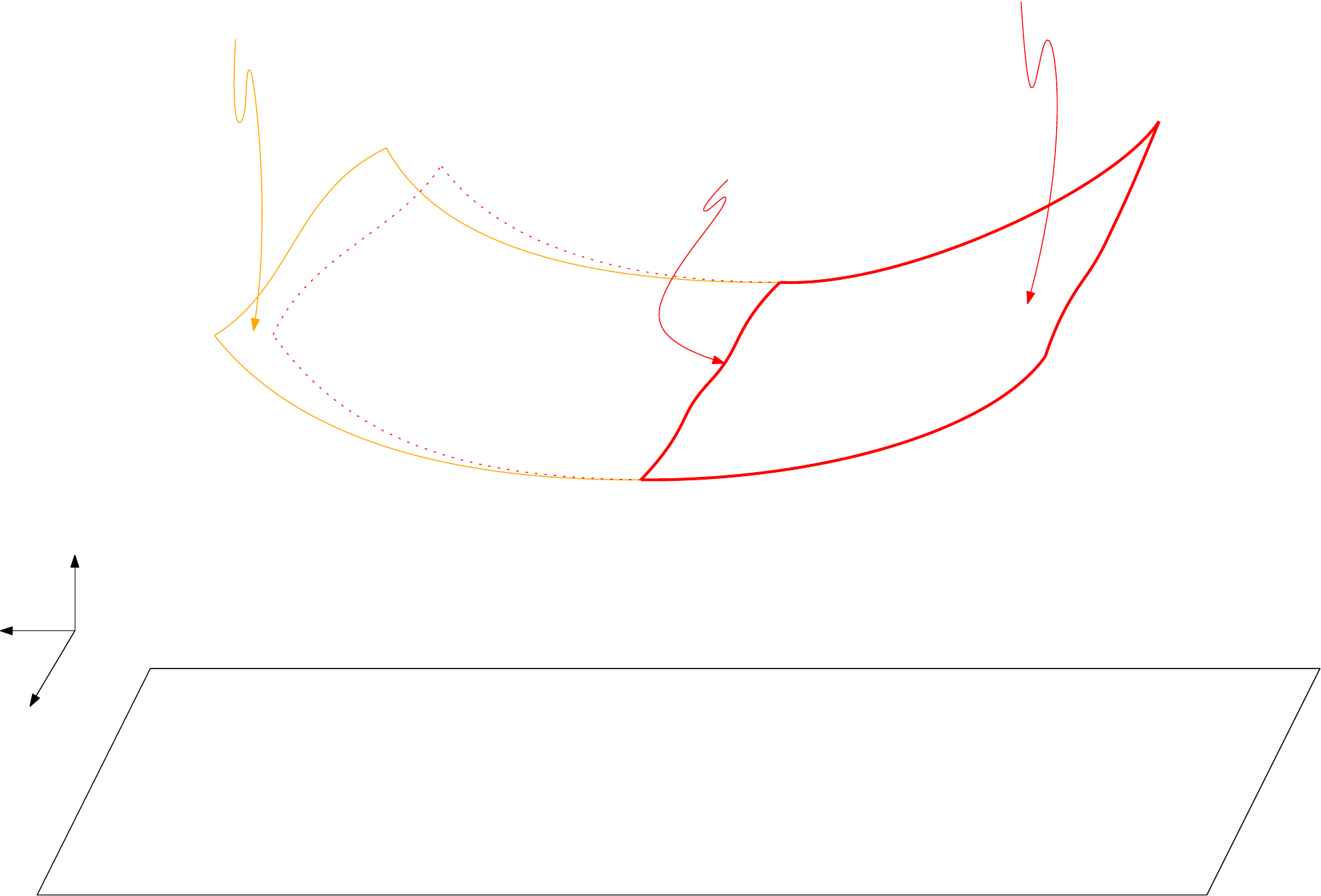}
			\put (73,68) {$\mathcal{B}$}
			\put (52,55.5) {$\partial_-\mathcal{B}$}
			\put (18,66) {$\underline{\mathcal{C}}$}
			\put (48,10) {$\Sigma_0$}
			\put (-3,15.5) {$(x^2,x^3) \in \mathbb{T}^2$}
			\put (6,27) {$t$}
			\put (-5,22.5) {$u \in \mathbb{R}$}
		\end{overpic}
	\caption{The singular boundary and Cauchy horizon in geometric coordinates}
\label{F:SINGULARBOUNDRYANDCAUCHYHORIZONINGEOMETRIC}
	\end{figure}
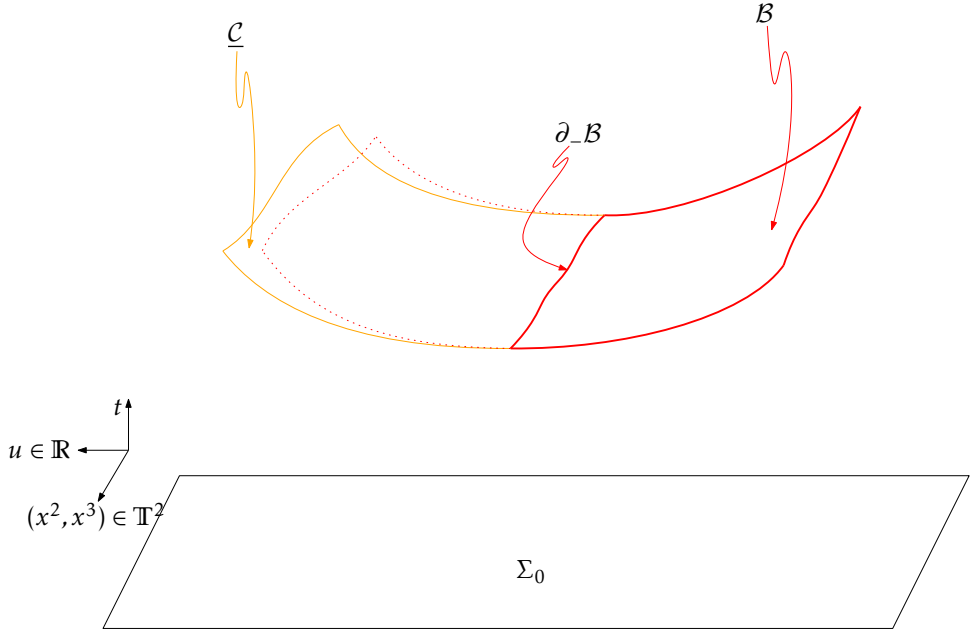
\end{center}

		The singular boundary $\mathcal{B}$, viewed as a subset of geometric coordinate space 
		(see Fig.\,\ref{F:SINGULARBOUNDRYANDCAUCHYHORIZONINGEOMETRIC}), is:
		\begin{align} \label{E:INTROSINGBOUNDARY}
			\mathcal{B}
			& = \lbrace \upmu = 0 \rbrace
				\cap
				\lbrace \muX \upmu \leq 0 \rbrace,
		\end{align}	
		where the vectorfield $\muX$, described later on 
		(see Fig.\,\ref{F:COMMUTATORVECTORFIELDSINCARTESIANCOORDINATES} for a depiction of $\muX$ in Cartesian coordinate space), 
		is transversal to the $\nullhyparg{u}$ and satisfies $\muX u = 1$.
		The complementary set  
		$
		\lbrace \upmu = 0 \rbrace
				\cap
				\lbrace \muX \upmu > 0 \rbrace
		$
		is not part of the singular boundary.
		In the context of Fig.\,\ref{F:SINGULARBOUNDRYANDCAUCHYHORIZONINGEOMETRIC} in geometric coordinate space,
		the irrelevance of
		$
		\lbrace \upmu = 0 \rbrace
				\cap
				\lbrace \muX \upmu > 0 \rbrace
		$
		for $\mathcal{B}$ and the maximal development can be understood as follows: 
		the singular boundary portion
		$\mathcal{B}$ in the figure cannot be extended to the left into the region where
		$\muX \upmu$ would be positive
		because before that region has a chance to dynamically develop,  
		it will be cut off by a Cauchy horizon 
		emanating from the crease $\partial_- \mathcal{B}$, 
		which is the left boundary of $\mathcal{B}$ in the figure.
		The crease $\partial_- \mathcal{B}$ is characterized by:
		\begin{align} \label{E:INTROCREASE}
			\partial_- \mathcal{B}
			& = \lbrace \upmu = 0 \rbrace
				\cap
				\lbrace \muX \upmu = 0 \rbrace.
		\end{align}	
		From \eqref{E:INTROSINGBOUNDARY}--\eqref{E:INTROCREASE}, the status of $\partial_- \mathcal{B}$
		as a boundary of $\mathcal{B}$ is clear.
		In the regime under study, \emph{the two level sets on RHS~\eqref{E:INTROCREASE}
		intersect transversally}, which is what gives the crease the structure of a $2D$ submanifold. 
		The transversality of the intersection is a consequence of \emph{acoustical transversal convexity},
		which is mild condition satisfied by open sets of solutions and which we describe below in more detail.
		By doing \emph{formal} Taylor expansions starting from the crease, 
		one can check that the unattainable portion $\{\upmu = 0\} \cap \{ \muX \upmu >0\}$ 
		would indeed lie to the causal future of the Cauchy horizon $\underline{\mathcal{C}}$;
		see the dotted portion in Fig.\,\ref{F:SINGULARBOUNDRYANDCAUCHYHORIZONINGEOMETRIC},
		which formally depicts $\{\upmu = 0\} \cap \{ \muX \upmu >0\}$.
	\item (\textbf{Null hypersurfaces and PDE energy degeneracies}). In Fig.\,\ref{F:MAXDEVELOPMENTINCARTESIAN},
		$\mathcal{B}$ and $\underline{\mathcal{C}}$ are acoustically null hypersurfaces emanating from the crease,
		and in particular, $\mathcal{B}$ is ruled by acoustically null curves whose tangent vectors
		are denoted by $\Lunit$ in the figure; see Prop.\,\ref{P:DESCRIPTIONOFSINGULARBOUNDARYINCARTESIANSPACE} for a
		proof of these properties of $\mathcal{B}$, and see Remarks~\ref{R:ACOUSTICALMETRICDEGENERACIESALONGSINGULARBOUNDARY} and
		Remark~\ref{R:NONUNIQUENESSOFINTEGRALCURVESOFLUNIT}
		for a discussion of some degeneracies that occur along $\mathcal{B}$.
		As is well-known, any $L^2$-type energy that one uses to control solutions 
		necessarily degenerates along null hypersurfaces,
		becoming only positive semi-definite instead of positive definite.
		This is a particularly challenging issue in the present context, 
		where singularities are forming along all of $\mathcal{B}$.
	\item  (\textbf{Regularity and rough foliations}).
		Many objects in the construction 
		(in particular, the crease $\partial_- \mathcal{B}$)
		are \underline{less regular} than the fluid solution, which leads to difficult regularity theory for the problem.
		To ``detect'' these rough objects as they emerge in the course of the evolution, 
		we rely on a new family of \emph{rough foliations}
		given by the level sets of \emph{rough time functions}, described below.
		The word ``rough'' refers to the fact that the rough time functions are also less regular than the fluid solution.
		There is another way in which the problem of shock formation can be viewed as a low regularity problem:
		the piling up of the characteristics is tied to the blowup of 
		a Euclidean-unit-length derivative of various fluid variables in directions \emph{transversal}
		to the characteristics, even though the solution remains rather smooth in directions \emph{tangent}
		to the characteristics. In particular,
		we are forced to close the estimates knowing that with respect
		to the Cartesian differential structure, there will be \emph{no differentiability}
		in transversal directions at the end of the classical evolution. 
		As we already mentioned, the geometric coordinates $(t,u,x^2,x^3)$
		partially ameliorate this difficulty in the sense
		that the solution remains rather smooth with respect to them.
		Nonetheless, as in other works on shock formation,
		our high order geometric energies can still become singular;
		see Sect.\,\ref{SSS:GEOMETRICENERGYESTIMATESONTHEROUGHFOLIATIONS}.
		This is one of the main technical challenges in the PDE analysis
		since singular high order energy estimates make it difficult for us to prove that
		the solution's partial derivatives with respect to the geometric coordinates 
		remain bounded at the lower derivative levels.
	\item (\textbf{One submanifold of $\mathcal{B}$ at a time via a family of rough time functions}).
		Our approach to constructing $\mathcal{B}$ is to show that it can be foliated by
		a family of $\muxmulevelsetvalue$-parameterized submanifolds $\twoargmumuxtorus{0}{-\muxmulevelsetvalue}$
		with $\muxmulevelsetvalue \geq 0$ a real parameter,
		and to construct each $\twoargmumuxtorus{0}{-\muxmulevelsetvalue}$, ``one $\muxmulevelsetvalue$ at a time;''
		see Fig.\,\ref{F:INTROCARTESIANROUGHFOLIATIONS}. Our construction is such that the crease $\partial_- \mathcal{B}$ 
		coincides with $\twoargmumuxtorus{0}{0}$.
		One might wonder why we didn't try to derive all of $\mathcal{B}$ ``at the same time.''
		From the discussion three points above, we see that 
		that approach would have effectively required us to work with foliations of spacetime that contain 
		or are asymptotic to level sets of $\upmu$ in regions where $\upmu$ is small
		(recall that $\mathcal{B}$ is a portion of $\lbrace \upmu = 0 \rbrace$).
		The difficulty is that for real numbers $\mulevelsetvalue$ small and positive, near the crease,
		the level sets $\lbrace \upmu = \mulevelsetvalue \rbrace$ have 
		a $\gfour$-timelike portion,
		along which top-order $L^2$ estimates for the solution are not available.
		This can formally be understood as the statement that
		in the solution regime under study, along the surface $\lbrace \upmu = 0 \rbrace$,
		RHS~\eqref{E:MUWEGHTEDSPACETIMEGRADMUNORMSQUARED} would be positive\footnote{In the regime under study,
		the term $\Lunit \upmu$ on RHS~\eqref{E:MUWEGHTEDSPACETIMEGRADMUNORMSQUARED} is strictly negative
		and the term $|\angD \upmu|_{\gtorus}^2$ is of negligible size.} 
		in the regions 
		where $\muX \upmu > 0$. 
		See Fig.\,\ref{F:SINGULARBOUNDARYANDOTHERLEVELSETSOFMUGEOMETRICCOORDS}, 
		where these level sets become $\gfour$-timelike in the region near their intersection with 
		the Cauchy horizon $\underline{\mathcal{C}}$.
		In fact, for $\mulevelsetvalue > 0$ small, 
		any surface that agreed with the level set $\lbrace \upmu = \mulevelsetvalue \rbrace$ up to
		second order along the surface $\lbrace \muX \upmu = 0 \rbrace$ 
		would suffer from the same difficulty:
		it would necessarily contain a $\gfour$-timelike portion,
		along which $L^2$ estimates for the solution are not available.
		
		\begin{remark}[Impossibility of $C^2$ spacelike foliations and the limited regularity of $\newtimefunction$]
		\label{R:INTROLIMITEDREGULARITYOFNEWTIMEFUNCTION}
		The upshot is that in the solution regime under study, 
		it is impossible to detect the entire singular boundary $\mathcal{B}$
		by deriving estimates on $C^2$ (relative to the geometric coordinates) 
		$\gfour$-spacelike foliations of spacetime 
		such that $\mathcal{B}$ (including its past boundary $\partial_- \mathcal{B}$) 
		is contained in the interior of one of the leaves of the foliation.
		This difficulty is connected to the following issue: in our main theorem,
		namely Theorem~\ref{T:DEVELOPMENTANDSTRUCTUREOFSINGULARBOUNDARY}, 
		the region $\MInteresting$ 
		that we study (which contains $\mathcal{B}$) is foliated by the level sets of a $C^{1,1}$ 
		(relative to the geometric coordinates) time function $\newtimefunction$, 
		and \emph{this $C^{1,1}$ regularity is optimal} given the shape of region;
		see Remark~\ref{R:NEWTIMEFUNCTIONISC11ANDNOTBETTERANDCONNECTIONTOCAUSALSTRUCTUREOFMUZEROLEVELSET}.
		In our approach to the PDE analysis, the $C^{1,1}$ regularity of $\newtimefunction$
		would be \emph{insufficient} for our proofs of some of our estimates
		(e.g., the co-dimension two Gauss curvature estimates we derive in Lemma~\ref{L:GAUSSCURVATUREOFROUGHTORILINFINITYESTIMATE}).
		As we explain below, we avoid these difficulties by avoiding working directly with 
		$\newtimefunction$; we instead work with a one-parameter family of time functions
		$\timefunctionarg{\muxmulevelsetvalue}$, which are more regular than $\newtimefunction$.
		\end{remark}
		
		\begin{center}
		\begin{figure} 
		\begin{overpic}[scale=.6, grid = false, tics=5, trim=-.5cm -1cm -1cm -.5cm, clip]{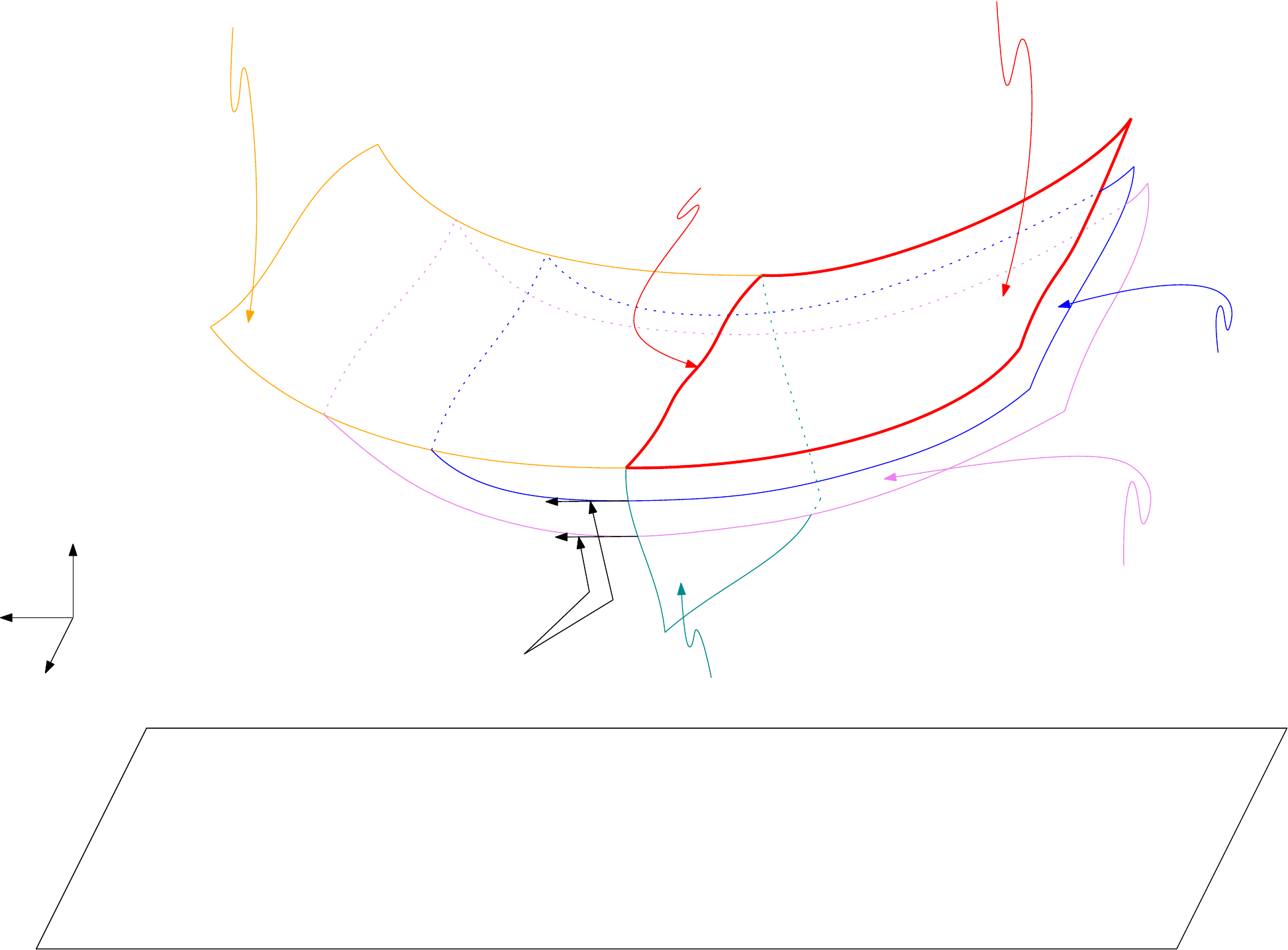}
			\put (73,74) {$\mathcal{B}$}
			\put (51.5,60.5) {$\partial_-\mathcal{B}$}
			\put (18.5,72.5) {$\underline{\mathcal{C}}$}
			\put (48,10) {$\Sigma_0$}
			\put (1,22) {$(x^2,x^3) \in \mathbb{T}^2$}
			\put (6,32) {$t$}
			\put (87,45) {$\{\upmu = \mulevelsetvalue\}$}
			\put (80,30.5) {$\{\upmu = 2 \mulevelsetvalue\}$}
			\put (49,22) {$\{\muX \upmu = 0\}$}
			\put (62,31) {$\{\muX \upmu < 0\}$}
			\put (30,31) {$\{\muX \upmu > 0\}$}
			\put (38,24) {$\muX$}
			\put (-5,28) {$u \in \mathbb{R}$}
	  \end{overpic}
		  \caption{The singular boundary and the level sets $\{\upmu = \mulevelsetvalue\}$ 
			and $\{ \upmu = 2\mulevelsetvalue\}$ for a small $\mulevelsetvalue > 0$, in geometric coordinates}
			\label{F:SINGULARBOUNDARYANDOTHERLEVELSETSOFMUGEOMETRICCOORDS}
		  \end{figure}
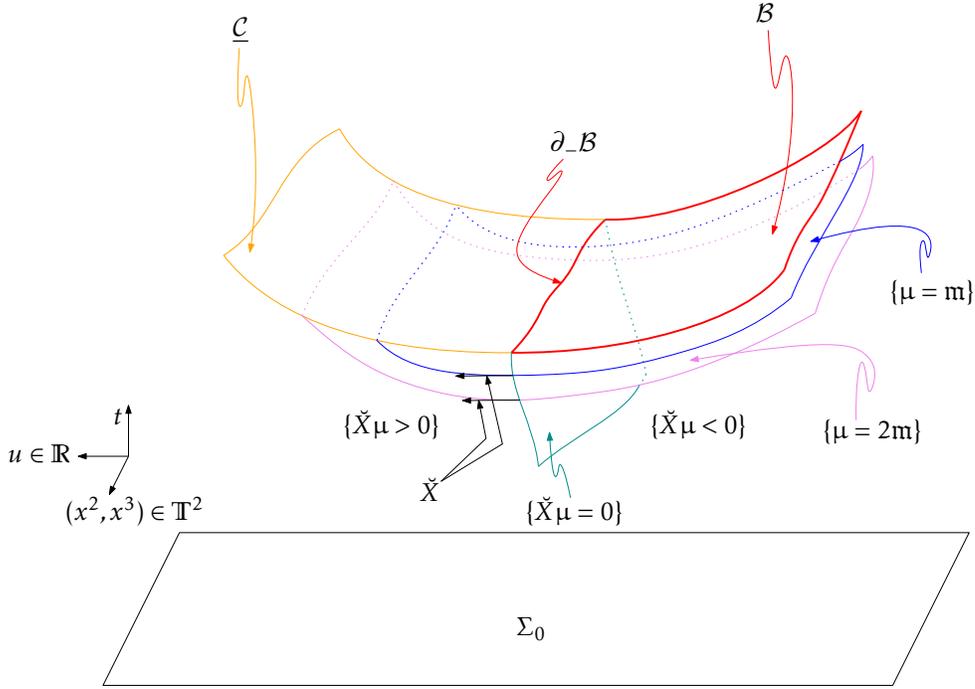
		\end{center}
		
		To construct the submanifolds $\twoargmumuxtorus{0}{-\muxmulevelsetvalue} \subset \mathcal{B}$ and circumvent
		the difficulties noted in Remark~\ref{R:INTROLIMITEDREGULARITYOFNEWTIMEFUNCTION},
		we proceed as follows.
		For each $\muxmulevelsetvalue \in [0,\muxmulevelsetvalue_0]$, where $\muxmulevelsetvalue_0 > 0$ is a small, data-dependent constant,
		we construct a \emph{rough time function} $\timefunctionarg{\muxmulevelsetvalue}$,
		which is $C^{2,1}$ relative to the geometric coordinates; 
		see Sect.\,\ref{S:ROUGHTIMEFUNCTIONANDROUGHSUBSETS}.
		Note that $\timefunctionarg{\muxmulevelsetvalue}$ is one degree more differentiable
		than the time function $\newtimefunction$ mentioned above, which is sufficient
		for avoiding the difficulties we highlighted in Remark~\ref{R:INTROLIMITEDREGULARITYOFNEWTIMEFUNCTION}.
		Our construction of $\timefunctionarg{\muxmulevelsetvalue}$
		depends on the eikonal function $u$, i.e., our construction of
		$\timefunctionarg{\muxmulevelsetvalue}$ relies on nonlinear geometric optics.
		For each $\muxmulevelsetvalue \in [0,\muxmulevelsetvalue_0]$, the level sets of $\timefunctionarg{\muxmulevelsetvalue}$
		yield a foliation of spacetime by $\gfour$-spacelike hypersurfaces, and 
		$\twoargmumuxtorus{0}{-\muxmulevelsetvalue} \subset \lbrace \timefunctionarg{\muxmulevelsetvalue} = 0 \rbrace$.
		Hence, to construct $\twoargmumuxtorus{0}{-\muxmulevelsetvalue}$, it suffices to control the fluid
		up to the level set $\lbrace \timefunctionarg{\muxmulevelsetvalue} = 0 \rbrace$; the vast majority of
		our efforts in this paper are dedicated towards that task.
		The union $\bigcup_{\muxmulevelsetvalue \in [0,\muxmulevelsetvalue_0]} \twoargmumuxtorus{0}{-\muxmulevelsetvalue}$ is the portion of $\mathcal{B}$
		that we construct in our main theorem. 
		The portion of $\mathcal{B}$ that formally corresponds to $\muxmulevelsetvalue < 0$
		never has a chance to emerge in the maximal classical development because it is cut off by the 
		Cauchy horizon $\underline{\mathcal{C}}$. In Fig.\,\ref{F:SINGULARBOUNDRYANDCAUCHYHORIZONINGEOMETRIC},
		the ``irrelevant portion'' of $\mathcal{B}$, corresponding to $\muxmulevelsetvalue < 0$,
		is formally delineated by dotted curves 
		(we cannot actually construct this portion, and we have displayed it only to illustrate that 
		$\underline{\mathcal{C}}$ lies below it).
		We could have extended our results to handle a larger range of $\muxmulevelsetvalue$ values, i.e., $\muxmulevelsetvalue > \muxmulevelsetvalue_0$;
		we avoided this because we would have had to modify our construction of the $\timefunctionarg{\muxmulevelsetvalue}$ for large 
		$\muxmulevelsetvalue$, which would have lengthened the paper.
		We refer Sect.\,\ref{SS:THETWOOPENPROBLEMS} for a
		more detailed overview of our construction of the rough time functions and the $\twoargmumuxtorus{0}{-\muxmulevelsetvalue}$.
		
		\begin{remark}[The terminology ``rough time function'']
		\label{R:TERMINOLOGYROUGHTIMEFUNCTION}
			The word ``rough'' in ``rough time function''
			refers to the fact that the elements of
			$\lbrace \timefunctionarg{\muxmulevelsetvalue} \rbrace_{\muxmulevelsetvalue \in [0,\muxmulevelsetvalue_0]}$ 
			are less regular than the fluid; see Sect.\,\ref{SSS:INTROROUGHTIMEFUNCTIONS}
			for further details.
			In the present paper, each $\timefunctionarg{\muxmulevelsetvalue}$ is $C^{2,1}$ 
			(with respect to the geometric coordinates)
			because the fluid variables are $C^{3,1}$.
			The fluid variables are $C^{3,1}$ because
			we have only assumed limited differentiability on the fluid data in directions
			transversal to the characteristics $\nullhyparg{u}$, 
			even though we assumed they are much smoother in the $\nullhyparg{u}$-tangential directions;
			see Sect.\,\ref{S:ASSUMPTIONSONTHEDATA} for our data assumptions.
			Despite the terminology ``rough time function,'' if 
			we had instead assumed that the fluid data were $C^{\infty}$, then 
			each $\timefunctionarg{\muxmulevelsetvalue}$ would also have been $C^{\infty}$.
			In contrast, even with $C^{\infty}$ fluid data, the time function $\newtimefunction$ from
			Remark~\ref{R:INTROLIMITEDREGULARITYOFNEWTIMEFUNCTION} would have had only $C^{1,1}$ regularity.
		\end{remark}
	\item (\textbf{Lack of strict convexity}).
		Observe that in Fig.\,\ref{F:MAXDEVELOPMENTINCARTESIAN},
		there are points $q \in \mathcal{B}$ such that the tangent plane
		to $\mathcal{B}$ at $q$, which we denote by $T_q \mathcal{B}$, 
		does not lie below $\mathcal{B}$.
		Moreover, there are points $q \in \partial_- \mathcal{B}$ such that
		every (three-dimensional) causal plane containing the two-dimensional
		subspace $T_q \partial_- \mathcal{B}$ 
		(which appears to be one-dimensional in Fig.\,\ref{F:MAXDEVELOPMENTINCARTESIAN}, due to our 
		suppression of a spatial dimension), fails to lie below $\mathcal{B}$
		\emph{even locally near $q$}.
		Let us informally refer to these phenomena as ``absence of strict convexity,''
		where here, ``convexity'' informally refers to ``upwards bending,''
		and ``strict convexity'' 
		-- though not featured Fig.\,\ref{F:MAXDEVELOPMENTINCARTESIAN} -- would refer to ``upwards bending in all directions.''
		We clarify that these notions implicitly refer to the $1+3$-dimensional Cartesian coordinate space 
		since the tangent planes referred to above are understood to be standard Cartesian-flat hyperplanes
		and Fig.\,\ref{F:MAXDEVELOPMENTINCARTESIAN} is a picture in Cartesian coordinate space.
		The absence of strict convexity
		poses serious technical difficulties.
		\begin{quote}
			In the absence of symmetry and strict convexity, the entirety of the crease $\partial_- \mathcal{B}$ 
			has never before been fully constructed for any open set of shock-forming solutions to any hyperbolic PDE.
			In particular, this aspect of our main results is new even in the case of irrotational and isentropic solutions.
		\end{quote}
		The rough time function $\timefunctionarg{\muxmulevelsetvalue}$ and corresponding rough foliations allow us to derive, 
		through a fully constructive approach,  
		the structure of the singular boundary,
		even if it is not strictly convex. 
		Instead of strict convexity, we rely on
		\emph{acoustical transversal convexity} 
		(which we refer to as transversal convexity for short),
		which allows us to handle, for example, perturbations of symmetric solutions, 
		where strict convexity of the singular boundary
		can fail due to the approximate symmetry. 
		Roughly, transversal convexity is a form of convexity only in 
		a particular direction, specifically in a direction that is transversal to the level sets of the eikonal function $u$.
		Note that transversal convexity refers to the structure of $\mathcal{B}$ viewed as an embedded submanifold
		of geometric coordinate space (as opposed to Cartesian coordinate space).
		In the solution regime under study, transversal convexity is captured by our 
		data assumption \eqref{E:DATATASSUMPTIONMUTRANSVERSALCONVEXITY} on the inverse foliation density $\upmu$,
		which we are able to propagate throughout the evolution (see \eqref{E:MUTRANSVERSALCONVEXITY}).
		The singular boundary $\mathcal{B}$ in Fig.\,\ref{F:SINGULARBOUNDRYANDCAUCHYHORIZONINGEOMETRIC} 
		enjoys transversal convexity 
		(roughly, it is ``parabolic in the $u$ direction'' at fixed $(x^2,x^3)$ near $\partial_- \mathcal{B}$). 
		One could check that in the solution regime under study, 
		the transversal convexity of $\mathcal{B}$ in geometric coordinates also implies,
		in the Cartesian coordinate picture,
		the convexity of the $x^1$-parameterized curves in $\mathcal{B}$ along which $(x^2,x^3)$ are fixed;
		while we do not directly need this ``Cartesian transversal convexity'' in our analysis, 
		we have exhibited the ``upwards bending\footnote{In $1+1$ dimensions, 
		under transversal convexity,
		the embedding of the singular boundary in Cartesian coordinate space
		can be modeled by the $u$-parameterized curve $t = u^2$, 
		$x^1 = u^2 + u^3$, for $u \leq 0$. 
		Near the origin (which models the crease) in $(t,x^1)$-space, 
		this singular boundary-modeling curve is asymptotic to the graph of 
		$t = (x^1)^2 + (x^1)^{3/2}$ for $x^1$ small and positive, which bends upwards 
		(c.f.\ the singular boundary in Fig.\,\ref{F:MAXDEVELOPMENTINCARTESIAN})
		and has regularity $C^{1,1/2}$ (c.f.\ the 
		regularity of $\Upsilon(\mathcal{B}^{[0,\muxmulevelsetvalue_0]})$ stated in
		Theorem~\ref{T:ABBREVIATEDSTATEMENTOFMAINRESULTS}).
		\label{FN:REGULARITYANDCONCAVITYOFSINGULARBOUNDARYINCARTESIANCOORDINATESPACE}} 
		of $\mathcal{B}$ in the $x^1$-direction'' 
		in Fig.\,\ref{F:MAXDEVELOPMENTINCARTESIAN}.
		Our assumption of transversal convexity is close to optimal in the sense that without it,
		the qualitative character of the singular boundary can dramatically change,
		even for plane-symmetric solutions; see Sect.\,\ref{SS:REMARKSONRESULTSANDMETHODS} for
		further discussion.
		We also highlight that the assumption of transversal convexity played a crucial role in Christodoulou's resolution
		\cite{dC2019} of the restricted shock development problem, which we describe below.
	\item  (\textbf{Degenerate wave and elliptic estimates on curved domains}).
		To close the high order $L^2$ estimates, we must
		adapt a variety of hyperbolic energy estimates and
		``top-order'' elliptic estimates for the vorticity, entropy, and geometry to the 
		precise shape of 
		$\partial_- \mathcal{B}$,
		$\mathcal{B}$, 
		and $\underline{\mathcal{C}}$,
		\emph{which are not known in advance}. The shock singularity introduces degeneracies into these estimates,
		and when we control the top-order derivatives of the vorticity and entropy using elliptic estimates, 
		our handling of these degeneracies 
		requires our observation of special cancellations
		within delicately constructed ``elliptic-hyperbolic'' identities.
		To exhibit the cancellations, we rely on the full nonlinear structure of the 
		geometric formulation of compressible Euler flow provided by Theorem~\ref{T:GEOMETRICWAVETRANSPORTSYSTEM}.
		In constructing the elliptic-hyperbolic identities (see Sect.\,\ref{S:ELLIPTICHYPERBOLICIDENTITIES}), 
		we rely on the framework we developed in \cite{lAjS2020}. However, 
		for the purposes of the present paper, we had to substantially upgrade that framework to accommodate the structure of the singularity.
		The key new object that we use to derive the elliptic-hyperbolic identities is a well-constructed \emph{characteristic current},
		defined in Def.\,\ref{D:PUTANGENTELLIPTICHYPERBOLICCURRENT}.
\end{itemize}

	\subsection{The two open problems that we resolve and connection to the shock development problem}
	\label{SS:THETWOOPENPROBLEMS}
	In his breakthrough 2007 monograph \cite{dC2007},
	Christodoulou gave a sharp description of the stable formation of shock singularities, 
	starting from open sets of smooth initial data, 
	in solutions to the $3D$ irrotational and isentropic relativistic Euler equations.
	Together with Miao, he later extended his results to the $3D$ compressible Euler equations \cite{dCsM2014},
	again for irrotational and isentropic solutions.
	These results revealed a large subset of the maximal classical development,
	including a portion of the boundary.
	
	In the wake of \cite{dC2007}, there have been many exciting developments on 
	the formation of shock singularities and the subsequent evolution of the solution as a weak solution,
	after the shock; see Sect.\,\ref{SS:BACKGROUNDONSHOCKS}
	for further discussion. However, two fundamental problems have remained open:
	
	\begin{enumerate}
		\item (\textbf{The full structure of bounded portions of the maximal classical development}).
			As is explained on \cite{dC2007}*{Pages 929, 968--969}, Christodoulou's approach
			yields a union of developments of the initial data, where each of his developments
			can be foliated by portions of \emph{Cartesian-flat}\footnote{By a ``Cartesian-flat'' hypersurface, we mean a plane with
			respect to the standard rectangular coordinates on $\mathbb{R}^{1+3}$.} 
			spacelike hypersurfaces and portions of characteristic hypersurfaces.
			By varying the ``angle of tilt'' of these Cartesian-flat hypersurfaces and varying the initial data of the characteristic
			hypersurfaces,
			one obtains (see \cite{dC2007}*{Pages 929, 968--969}) ``a larger part'' of the maximal development. 
			While this approach yields a sharp description of some portion of the maximal development,
			the precise portion that it reveals is not made explicit through the construction. 
			Moreover, from Fig.\,\ref{F:MAXDEVELOPMENTINCARTESIAN}, one can infer that
			for some solutions, there are portions of the boundary (in particular, portions of $\partial_- \mathcal{B}$)
			that are \emph{not accessible} 
			through Cartesian-flat spacelike foliations.
			This is connected to the lack of strict convexity of the singular boundary, as we discussed in
			Sect.\,\ref{SS:OVERVIEWOFDIFFICULTIES};
			see also Fig.\,\ref{F:INTERESTINGREGIONMAINRESULTS} and, in Sect.\,\ref{SS:REMARKSONRESULTSANDMETHODS},
			our discussion of the points $b_1$ and $b_2$ featured in the figure.
			Hence, the following problem is glaring:
			\begin{quote}
				Can one derive the \emph{full structure} of the maximal development, 
				at least in some region of spacetime that includes a neighborhood of the boundary
				that contains the crease?
			\end{quote}
			This problem, while mathematically rich in itself, is important for two other fundamental reasons:
			\textbf{I)} The breakthrough result \cite{fEhRjS2019} 
			shows that in general, one cannot ensure uniqueness of the maximal classical development 
			until one constructs it and proves that it enjoys some crucial structural qualitative properties.
			For quasilinear hyperbolic PDEs, the question of uniqueness of the maximal classical development is global in nature.
			However, the results we derive in the present article and the companion \cite{lAjS20XX}
			exhibit a localized version of the crucial property that the maximal development ``lies on one side of its boundary.''
			In \cite{fEhRjS2019}, the authors showed, roughly speaking, that if this property holds globally, 
			then the maximal classical development is unique.
			\textbf{II)} As is explained in \cite{dC2019},
			the full structure of a neighborhood of the boundary of the maximal development
			is an essential ingredient for properly setting up the
			aforementioned shock development problem.
		\item (\textbf{Removing the irrotationality and isentropicity assumptions}).
			Assuming a positive resolution to the first problem, 
			a second one of clear physical importance stands out:
			\begin{quote}
				Can one extend the result away from the irrotational and isentropic 
				class of solutions, that is, to handle solutions with non-zero vorticity and dynamic entropy?
			\end{quote}
			It is of fundamental importance to understand such ``general solutions'' because vorticity and entropy
			will form in the weak solution (see, e.g., \cite{dC2019}), i.e., ``after the first shock,'' even if the original
			initial data are irrotational, isentropic, and $C^{\infty}$. That is, if one aims towards developing
			a global-in-space-and-time theory that accommodates the formation of shocks and their subsequent interactions,
			\emph{vorticity and entropy are an unavoidable aspect of the dynamics.}
	\end{enumerate}
	
	In the present paper and its companion \cite{lAjS20XX}, we resolve both problems
	for the $3D$ compressible Euler equations \eqref{E:INTROTRANSPORTVI}--\eqref{E:INTROBS}.
	In the notation of Fig.~\ref{F:MAXDEVELOPMENTINCARTESIAN},
	starting from smooth data on $\Sigma_0$, we construct the classical solution
	in the region lying in between $\Sigma_0$ and\footnote{This union is not disjoint, for the crease $\partial_- \mathcal{B}$
	is a past boundary of both of the closed sets $\underline{\mathcal{C}}$ and $\mathcal{B}$.} 
	$\underline{\mathcal{C}} \cup \partial_- \mathcal{B} \cup \mathcal{B}$. More precisely, our results apply to 
	open sets of smooth data satisfying assumptions described below
	(see, for example, Theorem~\ref{T:ABBREVIATEDSTATEMENTOFMAINRESULTS}), 
	and we control the solution in bounded regions
	of spacetime that contain \emph{all}\footnote{More precisely, we construct an entire connected component of $\partial_- \mathcal{B}$.} 
	of $\partial_- \mathcal{B}$
	and a full neighborhood of it in $\underline{\mathcal{C}}$ and $\mathcal{B}$.
	For convenience, in the present paper, we have studied the solution only in 
	a single, ``spatially local'' region.
	However, our approach could be used as a building block to study the solution
	across space, at least in regions where the solution exhibits the property of
	acoustical transversal convexity, mentioned in Sect.\,\ref{SS:OVERVIEWOFDIFFICULTIES}.
	
\begin{figure}[ht]
\centering
\begin{subfigure}{.5\textwidth}
  \centering
  	\begin{overpic}[scale=.36, grid = false,trim=-.5cm -1cm -1cm -1.9cm, clip=true]{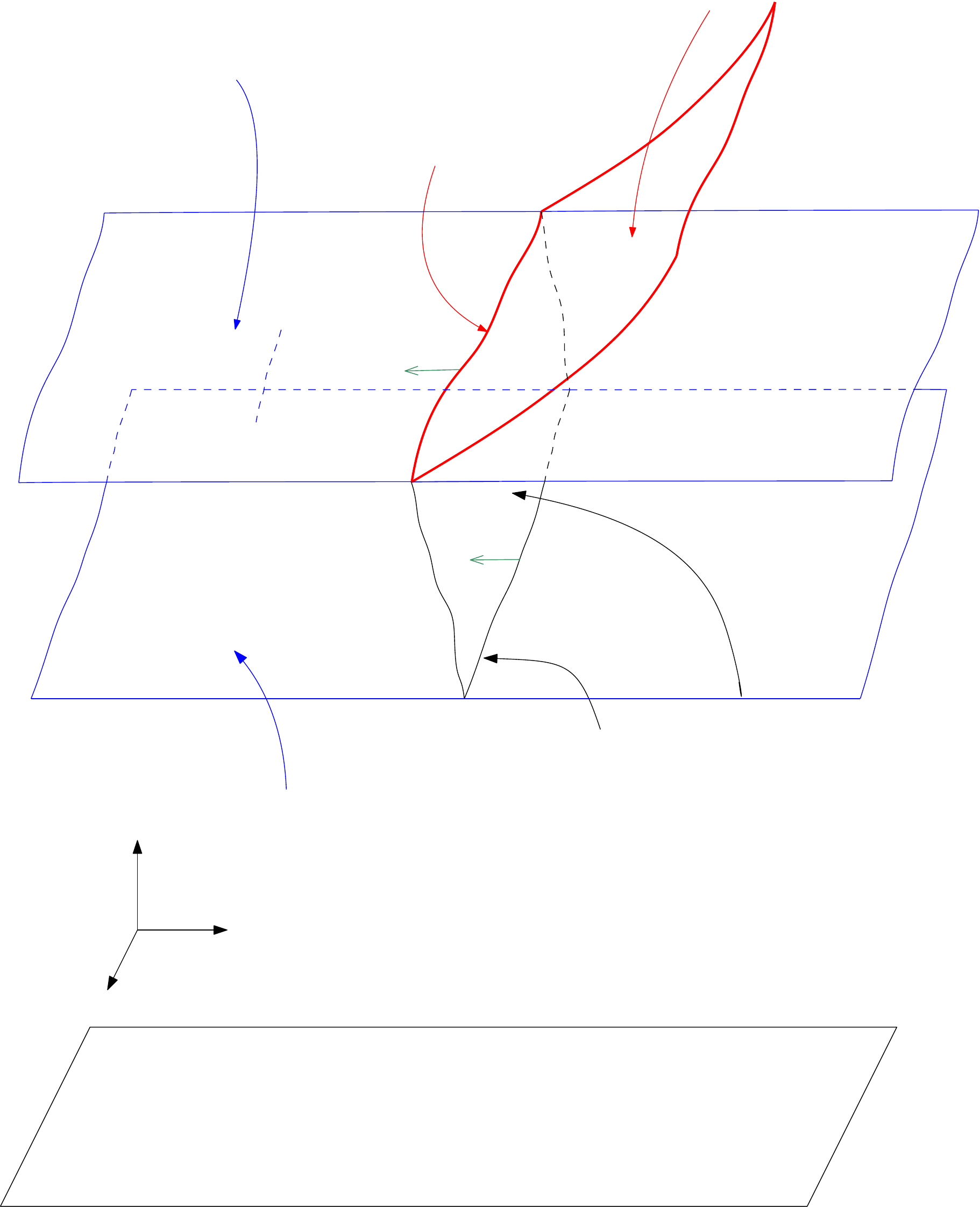}
			\put (36,11) {$\Sigma_0$}
			\put (11,19) {$(x^2,x^3) \in \mathbb{T}^2$}
			\put (10,29) {$t$}
			\put (20,23.5) {$x^1 \in \mathbb{R}$}
			\put (31,90) {$\twoargmumuxtorus{0}{0}$}
			\put (34,85) {\rotatebox{90}{$=$}}
			\put (31,82) {$\partial_- \mathcal{B}$}
			\put (54,93) {$\mathcal{B}$}
			\put (12,89) {$\hypthreearg{0}{[-\rightu,\leftu]}{0}$}
			\put (43,36.5) {$\twoargmumuxtorus{-\timefunction}{0}$}
			\put (55,37) {$\datahypfortimefunctionarg{0}$}
			\put (18,30) {$\hypthreearg{\timefunction}{[- \rightu,\leftu]}{0}$}
			\put (24.5,65) {$\Wtransarg{0}$}
			\put (31,52) {$\Wtransarg{0}$}
		\end{overpic}
		\caption{Rough foliations adapted to the crease}
		\label{F:CARTESIANROUGHFOLIATIONCREASE}
\end{subfigure}%
\begin{subfigure}{.5\textwidth}
 \centering
\begin{overpic}[scale=.36,grid=false,trim=-.5cm -1cm -1cm -.5cm, clip=true]{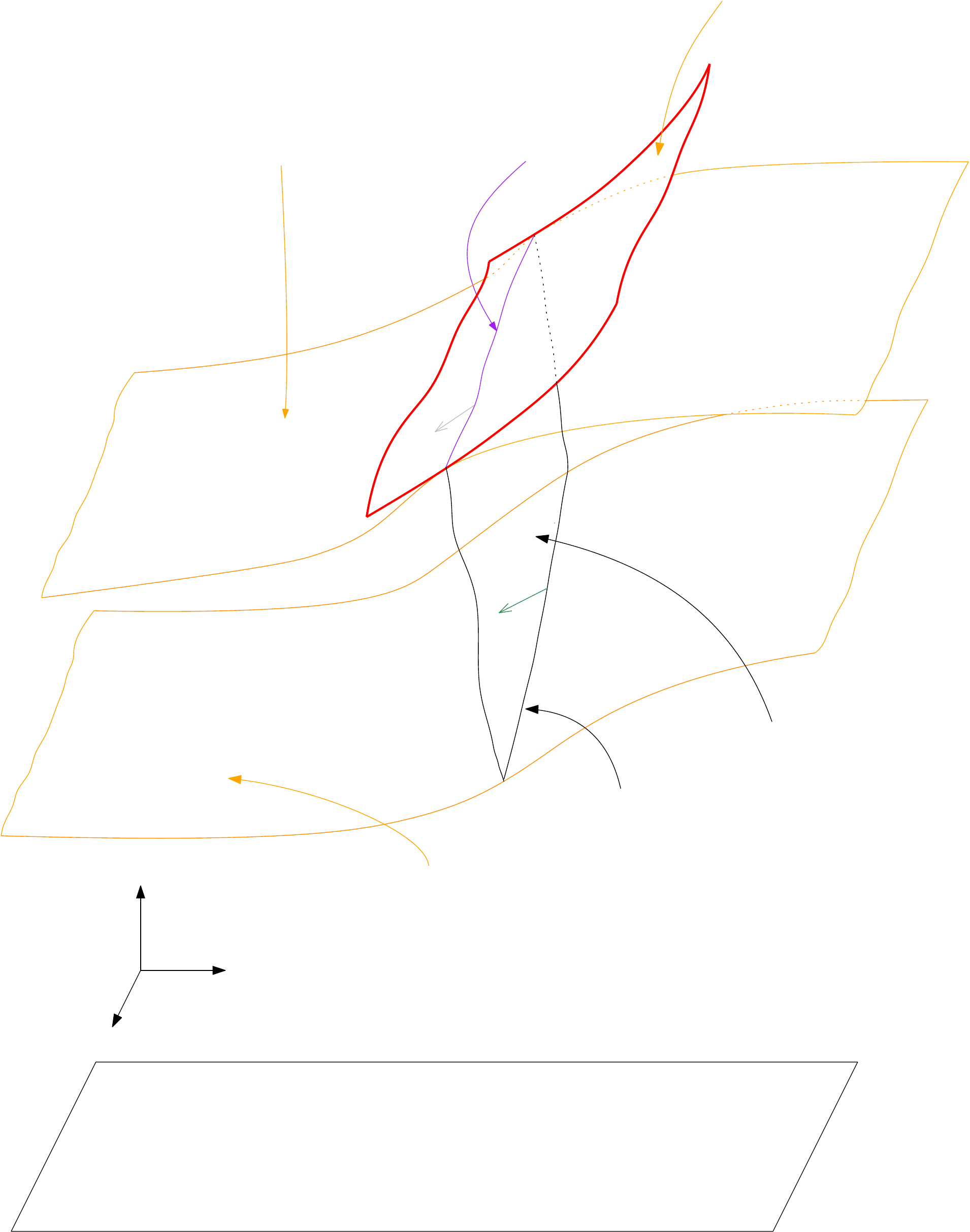} 
			\put (36,11) {$\Sigma_0$}
			\put (11.5,19) {$(x^2,x^3) \in \mathbb{T}^2$}
			\put (10,29) {$t$}
			\put (20,23.5) {$x^1 \in \mathbb{R}$}
			\put (39,88) {$\twoargmumuxtorus{0}{-\muxmulevelsetvalue}$}
			\put (58,98) {$\mathcal{B}$}
			\put (46,35) {$\twoargmumuxtorus{-\timefunction}{-\muxmulevelsetvalue}$}
			\put (15,90) {$\hypthreearg{0}{[-\rightu,\leftu]}{\muxmulevelsetvalue}$}
			\put (59,39) {$\datahypfortimefunctionarg{-\muxmulevelsetvalue}$}
			\put (30,27.5) {$\hypthreearg{\timefunction}{[- \rightu,\leftu]}{\muxmulevelsetvalue}$}
			\put (29,66) {$\Wtransarg{\muxmulevelsetvalue}$}
			\put (34,52.5) {$\Wtransarg{\muxmulevelsetvalue}$}
\end{overpic}
  \caption{Rough foliations adapted to a non-crease torus $\twoargmumuxtorus{0}{-\muxmulevelsetvalue} \subset \mathcal{B}$}
	\label{F:CARTESIANROUGHFOLIATIONFUTUREOFCREASE}
\end{subfigure}
\caption{Rough foliations adapted to the singular boundary, depicted in Cartesian coordinate space}
\label{F:INTROCARTESIANROUGHFOLIATIONS}
\end{figure}
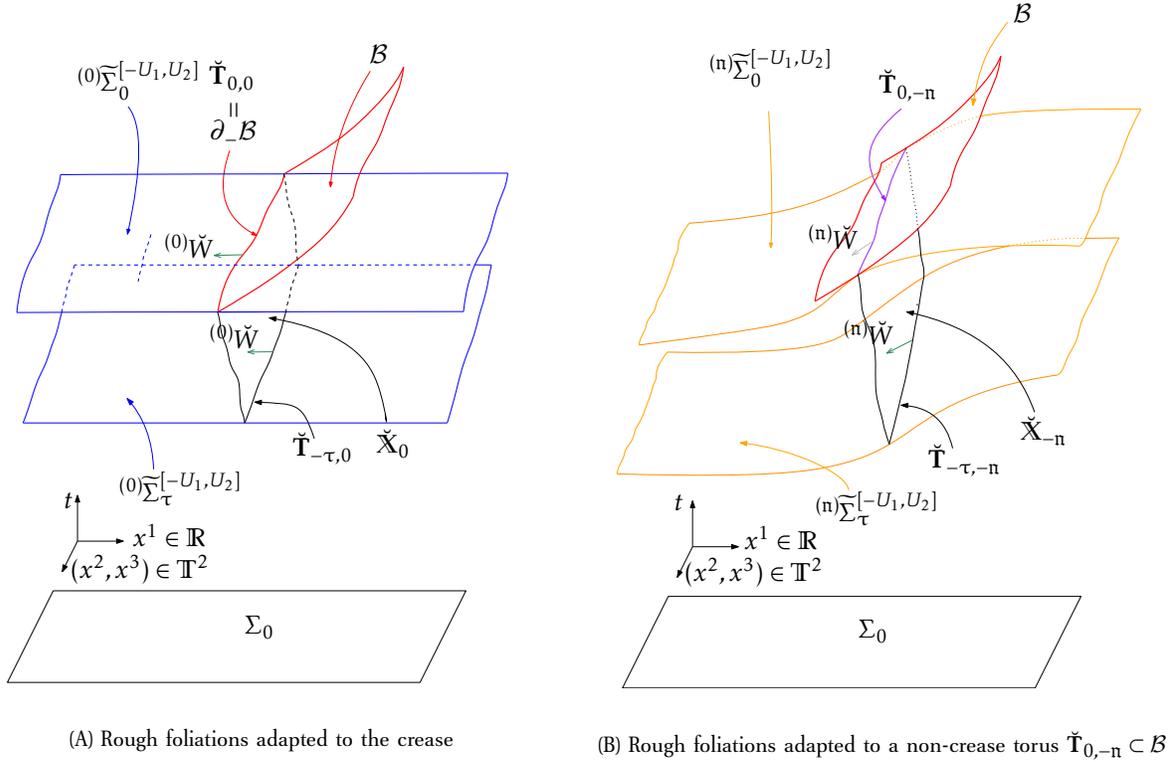
	
	In the present paper, we follow the solution up 
	to the singular portion of the boundary, denoted by $\mathcal{B}$
	in Figs.~\ref{F:MAXDEVELOPMENTINCARTESIAN} and \ref{F:CARTESIANROUGHFOLIATIONCREASE},
	where some first derivatives of the density and velocity blow up.
	We again emphasize that $\mathcal{B}$ contains its past boundary, $\partial_- \mathcal{B}$,
	which is the aforementioned crease.
	In \cite{lAjS20XX}, we construct the Cauchy horizon, denoted by 
	$\underline{\mathcal{C}}$ in Fig.~\ref{F:MAXDEVELOPMENTINCARTESIAN},
	which is an acoustically null hypersurface such that no fluid singularity forms along it,
	but it is nonetheless a boundary of the maximal classical development
	because its causal past intersects the singularity.
	As we mentioned earlier, our analysis in the present paper fundamentally relies on 
	a new family of acoustically spacelike \emph{rough foliations} that are \underline{precisely} and dynamically adapted
	to the shape of $\mathcal{B}$, which is not known in advance. 
	The foliations are level sets of a one-parameter family of \emph{rough time functions}
	$\lbrace \timefunctionarg{\muxmulevelsetvalue} \rbrace_{\muxmulevelsetvalue \in [0,\muxmulevelsetvalue_0]}$,
	where $\muxmulevelsetvalue_0 > 0$ is a constant depending on the initial data.
	Each $\timefunctionarg{\muxmulevelsetvalue}$ is defined on a portion of the maximal classical development union its boundary
	and has a range $[\timefunction_0,0]$ for some constant $\timefunction_0  < 0$ depending on the initial data.
	Moreover, $\timefunctionarg{\muxmulevelsetvalue}$ has the crucial property that its zero level set $\lbrace \timefunctionarg{\muxmulevelsetvalue} = 0 \rbrace$
	is tangent to $\mathcal{B}$ and intersects it in a submanifold. More precisely,
	$\mathcal{B} \cap \lbrace \timefunctionarg{\muxmulevelsetvalue} = 0 \rbrace = \twoargmumuxtorus{0}{-\muxmulevelsetvalue}$
	is a torus with spacetime co-dimension $2$ such that 
	the tori $\twoargmumuxtorus{0}{-\muxmulevelsetvalue}$ foliate a neighborhood of $\partial_- \mathcal{B}$ in $\mathcal{B}$,
	i.e., $\bigcup_{\muxmulevelsetvalue \in [0,\muxmulevelsetvalue_0]} \twoargmumuxtorus{0}{-\muxmulevelsetvalue}$
	is a neighborhood of $\partial_- \mathcal{B}$ in $\mathcal{B}$.
	In particular, the crease $\partial_- \mathcal{B}$ is equal to
	$\twoargmumuxtorus{0}{0}$.
	In Fig.\,\ref{F:CARTESIANROUGHFOLIATIONCREASE}, we exhibit two level sets of $\timefunctionarg{0}$,
	where the top one contains the crease.
	In Fig.\,\ref{F:CARTESIANROUGHFOLIATIONFUTUREOFCREASE}, for some $\muxmulevelsetvalue > 0$,
	we exhibit two level sets of $\timefunctionarg{\muxmulevelsetvalue}$,
	where the top one contains the torus $\twoargmumuxtorus{0}{-\muxmulevelsetvalue}$.
	
	We already highlight that the hypersurfaces
	$\datahypfortimefunctionarg{0}$
	and
	$\datahypfortimefunctionarg{-\muxmulevelsetvalue}$
	in Fig.\,\ref{F:CARTESIANROUGHFOLIATIONFUTUREOFCREASE}
	play a crucial role in our construction of 
	$\timefunctionarg{0}$
	and
	$\timefunctionarg{\muxmulevelsetvalue}$ respectively.
	In particular, by construction,
	$\timefunctionarg{0}$
	and
	$\timefunctionarg{\muxmulevelsetvalue}$
	solve transport equations with initial data 
	given on 	
	$\datahypfortimefunctionarg{0}$
	and
	$\datahypfortimefunctionarg{-\muxmulevelsetvalue}$ respectively;
	see Sect.\,\ref{SSS:ROUGHTIMEFUNCTIONSANDINITIALROUGHSLICE} 
	for the details.
	The surfaces 
	$\datahypfortimefunctionarg{0}$
	and
	$\datahypfortimefunctionarg{-\muxmulevelsetvalue}$
	are transversal to $\mathcal{B}$
	and intersect $\mathcal{B}$ in the tori
	$\twoargmumuxtorus{0}{0}$
	and
	$\twoargmumuxtorus{0}{-\muxmulevelsetvalue}$
	respectively.
	While the transversality of $\datahypfortimefunctionarg{0}$
	and
	$\datahypfortimefunctionarg{-\muxmulevelsetvalue}$
	to $\mathcal{B}$ is crucial for our analysis,
	the causal structure of 
	$\datahypfortimefunctionarg{0}$
	and
	$\datahypfortimefunctionarg{-\muxmulevelsetvalue}$
	(i.e., whether they are acoustically timelike, spacelike or null)
	is \emph{not} important because
	\emph{we do not have to derive any energy estimates for the fluid}
	along
	$\datahypfortimefunctionarg{0}$
	or
	$\datahypfortimefunctionarg{-\muxmulevelsetvalue}$.
	In particular,
	$\datahypfortimefunctionarg{0}$
	and
	$\datahypfortimefunctionarg{-\muxmulevelsetvalue}$
	can be acoustically timelike, spacelike, or null.
	
	In both papers, our approach relies on giving a complete description of the dynamics that in particular shows 
	what blows up and what remains regular and yields a sharp description of the structure and
	regularity of $\mathcal{B}$, $\partial_- \mathcal{B}$, and $\underline{\mathcal{C}}$.
	To handle the presence of vorticity and entropy,
	we develop modified versions of the integral identities that we discovered in \cite{lAjS2020},
	adapted here to the precise structure of $\mathcal{B}$.
	More precisely, in the present paper,
	we rely on a new family of ``elliptic-hyperbolic''
	integral identities that are adapted to the rough foliations and the structure of $\mathcal{B}$,
	and we develop new analytic techniques to handle the following difficulty, which permeates the paper:
	\emph{the singular boundary is ruled by curves that have acoustically null tangent vectors,\footnote{In the context of non-degenerate
	Lorentzian metrics,
	hypersurfaces that are ruled by curves with null tangent vectors are usually called ``null hypersurfaces,'' 
	and as is well-known, those tangent vectors are also normals to the surface. However, we avoid referring to the singular boundary
	as a ``null hypersurface'' because of the following severe degeneracy, which is tied to the many analytical difficulties
	we must handle in the PDE analysis: the singular boundary is
	a portion of the zero level set of a function ``$\upmu$'' whose gradient with respect to the Cartesian coordinates 
	\emph{blows up} along its zero level set; below we discuss $\upmu$ in detail.
	Hence, in the Cartesian differential structure, the gradient one-form of $\upmu$, 
	which in non-degenerate contexts
	can be viewed as the metric dual to its normal, is ill-defined. \label{FN:CHARACTEROFSINGULARBOUNDARY}}
	which leads to severe degeneracies in the estimates}, 
	especially in the top-order energy estimates for the vorticity and entropy.
	See Prop.\,\ref{P:DESCRIPTIONOFSINGULARBOUNDARYINCARTESIANSPACE} and Remark~\ref{R:NONUNIQUENESSOFINTEGRALCURVESOFLUNIT}
	for detailed information
	on the structure of the singular boundary, including its causal structure and its
	differential-topological properties as a subset of Cartesian coordinate space.
	
	\begin{center}
	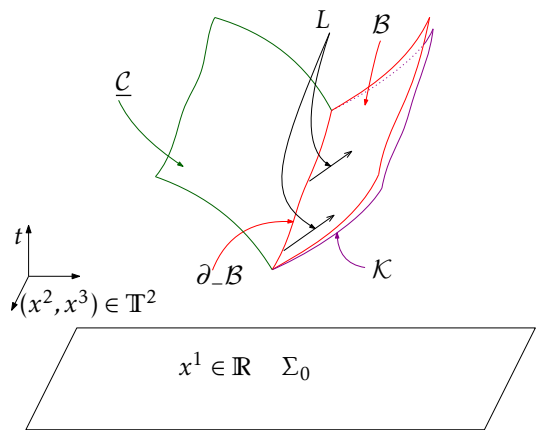
\begin{figure}  
		\begin{overpic}[scale=.4, grid = false, tics=5, trim=-.5cm -1cm -1cm -.5cm, clip]{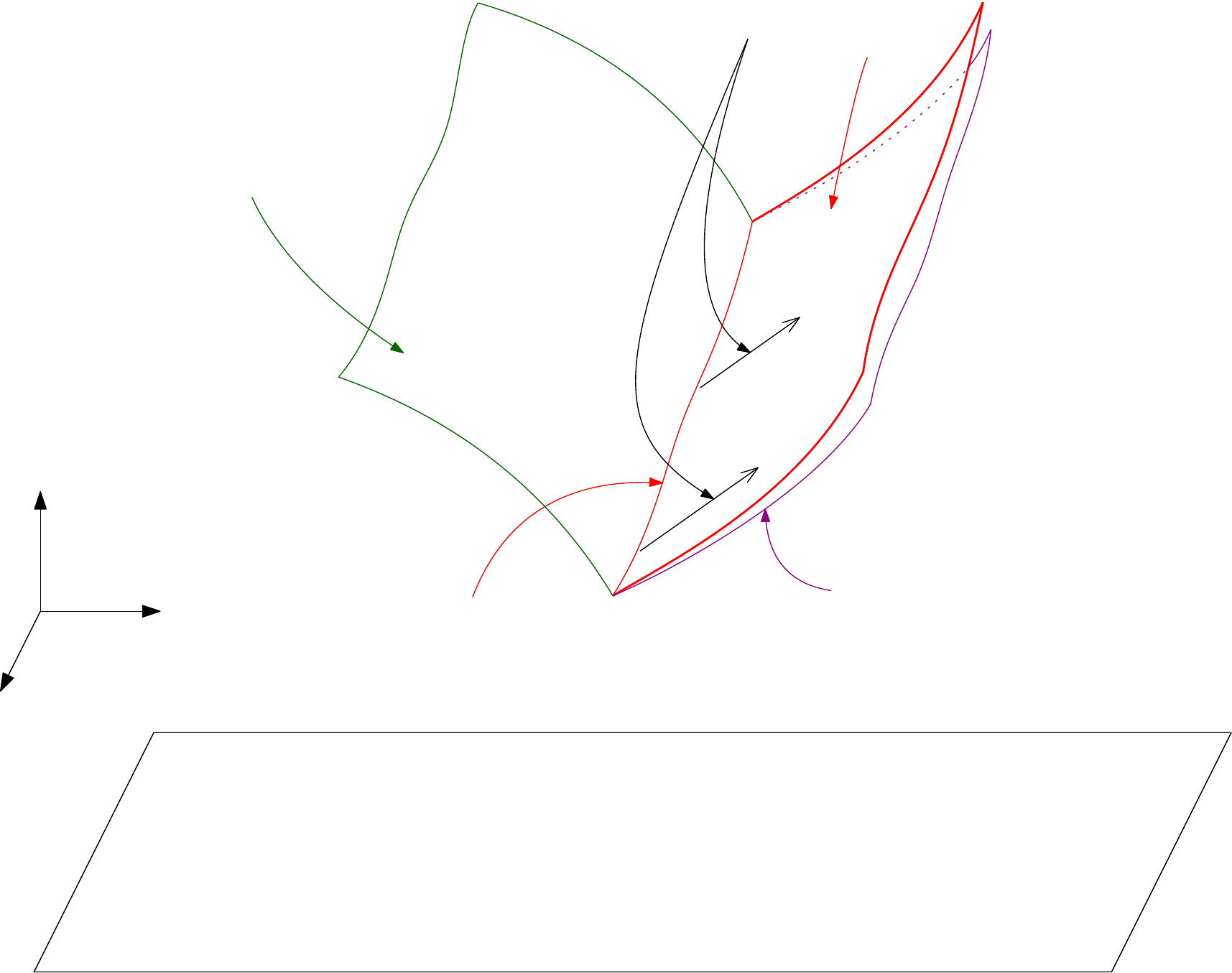}
			\put (50,15) {$\Sigma_0$}
			\put (4,26) {$(x^2,x^3) \in \mathbb{T}^2$}
			\put (3,38) {$t$}
			\put (32,15) {$x^1 \in\mathbb{R}$}
			\put (35,31) {$\partial_- \mathcal{B}$}
			\put (66,75) {$\mathcal{B}$}
			\put (66,32) {$\mathcal{K}$}
			\put (21,65) {$\underline{\mathcal{C}}$}
			\put (56,76) {$\Lunit$}
		\end{overpic}
		\caption{A localized subset of the maximal classical development and the shock hypersurface in Cartesian coordinate space}
	\label{F:MAXDEVELOPMENTWITHSHOCKHYPERSURFACEINCARTESIAN}
	\end{figure}
\end{center}	

	In the solution regime that we treat in our main results, 
	the crease $\partial_- \mathcal{B}$ has the structure of a co-dimension $2$
	spacelike submanifold; see Fig.\,\ref{F:MAXDEVELOPMENTWITHSHOCKHYPERSURFACEINCARTESIAN}.
	This co-dimension $2$ spacelike submanifold structure 
	is essential for properly setting up the shock development problem using
	the known techniques \cite{dC2019}.
	In fact, without this structure, it is not even clear whether the shock development problem is well-posed.
	In particular, the shock hypersurface (with boundary), denoted by ``$\mathcal{K}$'' in
	Fig.\,\ref{F:MAXDEVELOPMENTWITHSHOCKHYPERSURFACEINCARTESIAN}, emanates from the crease
	and starts off tangent to the singular boundary. Roughly, $\mathcal{K}$ is the hypersurface
	of discontinuity for the weak solution that develops to the future of the crease. 
	$\mathcal{K}$ is not part of the constructions we provide in this paper; 
	instead, its construction would be part of the resolution of the shock development problem.
	Let us again highlight one of the key difficulties in the problem: 
	$\mathcal{B}$ is ruled by acoustically null curves with tangent vectors
	denoted by ``$\Lunit$'' in Fig.\,\ref{F:MAXDEVELOPMENTWITHSHOCKHYPERSURFACEINCARTESIAN}.
	This leads to rather severe degeneracies
	in the analysis because, as is well-known, the coercive quantities that can be used to control the solution
	become degenerate along null hypersurfaces. 
	Here, the notion of null is with respect to the \emph{acoustical metric}
	$\gfour$, the solution-dependent Lorentzian metric (see \eqref{E:ACOUSTICALMETRIC}) 
	that captures the intrinsic geometry of sound waves in
	the flow. See below for extended discussion of these fundamental issues.

\subsection{Abbreviated statement of the main results}	
	\label{SS:ABBREVIATEDSTATEMENTOFMAINRESULTS}
The full statement of our main results is quite lengthy, 
due to the intricate geometric structures
and the highly tensorial nature of the singularity.
In Theorem~\ref{T:ABBREVIATEDSTATEMENTOFMAINRESULTS}
we provide an abbreviated, slightly informal statement of the main results.
In Theorems~\ref{T:EXISTENCEUPTOTHESINGULARBOUNDARYATFIXEDKAPPA} and \ref{T:DEVELOPMENTANDSTRUCTUREOFSINGULARBOUNDARY},
we provide full statements of the main results.

\begin{center}
	\begin{figure}[ht]  
		\begin{overpic}[scale=.36, grid = false, tics=5, trim=-.5cm -1cm -1cm -.5cm, clip]{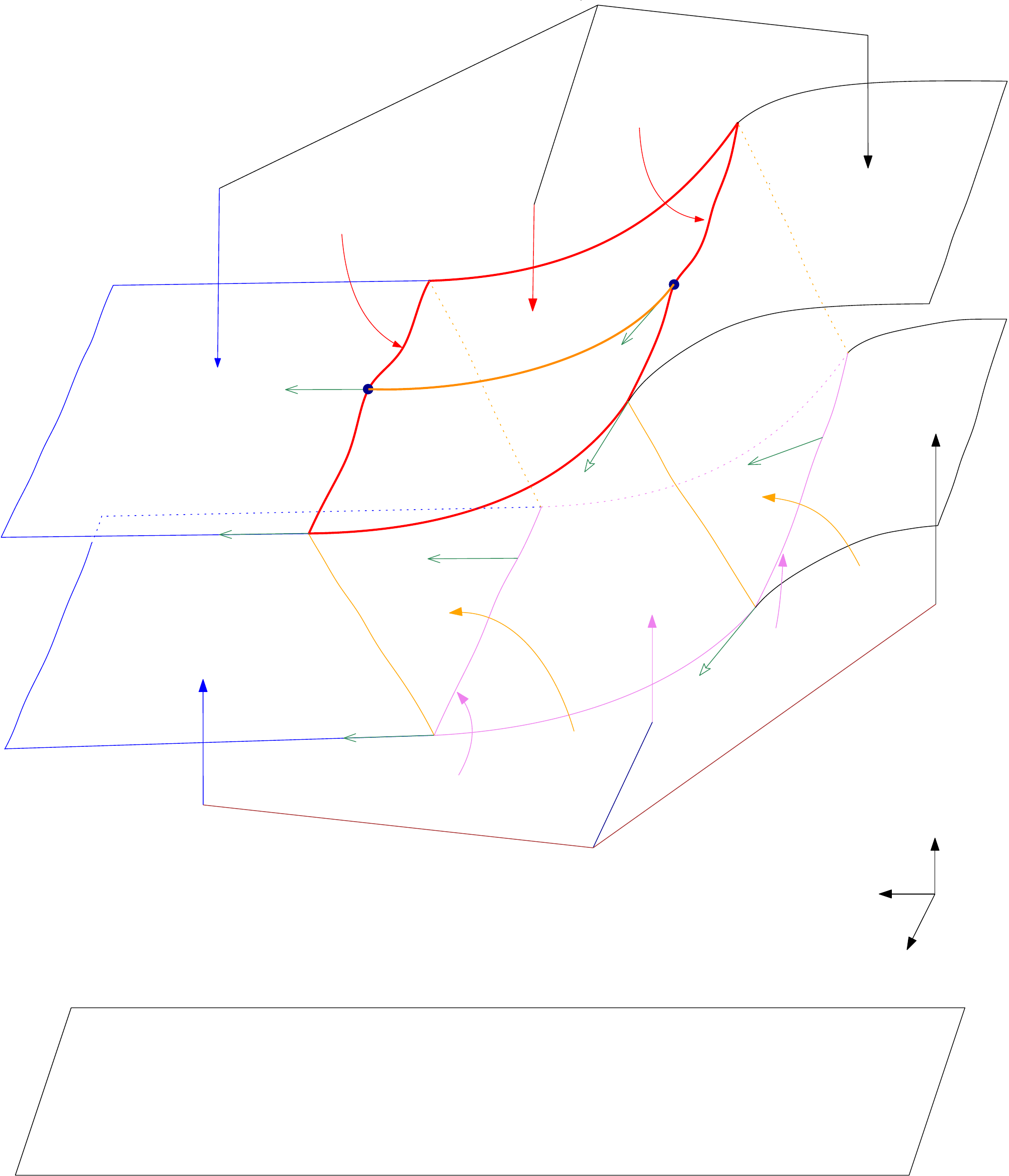}
			\put (40,10) {$\Sigma_0$}
			\put (61,18.5) {$(x^2,x^3) \in \mathbb{T}^2$}
			\put (78,28) {$t$}
			\put (62.5,25) {$u \in\mathbb{R}$}
			\put (15,50) {$\MLeft$}
			\put (45,53) {$\MSingular$}
			\put (65,70) {$\MRight$}
			\put (30,25) {$\inthyp{\timefunction_0}{[- \rightu,\leftu]}$}
			\put (30,100) {$\inthyp{0}{[- \rightu,\leftu]}$}
			\put (32,35) {$\twoargmumuxtorus{-\timefunction}{0}$}
			\put (61,45) {$\twoargmumuxtorus{-\timefunction}{-\muxmulevelsetvalue_0}$}
			\put (26,81) {$\twoargmumuxtorus{0}{0}$}
			\put (49,89.5) {$\twoargmumuxtorus{0}{-\muxmulevelsetvalue_0}$}
			\put (34,80) {$\mathcal{B}^{[0,\muxmulevelsetvalue_0]}$}
			\put (42,35) {$\datahypfortimefunctiontwoarg{0}{[\timefunction_0,0]}$}
			\put (66,49) {$\datahypfortimefunctiontwoarg{-\muxmulevelsetvalue_0}{[\timefunction_0,0]}$}
			\put (17,65) {$\Wtransarg{0}$}
			\put (28,52) {$\Wtransarg{0}$}
			\put (15,53) {$\Wtransarg{0}$}
			\put (25,40) {$\Wtransarg{0}$}
			\put (43,66) {$\Wtransarg{\muxmulevelsetvalue_0}$}
			\put (53,59) {$\Wtransarg{\muxmulevelsetvalue_0}$}
			\put (54,41) {$\Wtransarg{\muxmulevelsetvalue_0}$}
			\put (45,56) {$\Wtransarg{\muxmulevelsetvalue_0}$}
			\put (29,69) {$b_1$}
			\put (53,77) {$b_2$}
	\end{overpic}
		\caption{The region $\MInteresting = \MLeft \cup \MSingular \cup \MRight$ in
		geometric coordinate space}
	\label{F:INTERESTINGREGIONMAINRESULTS}
	\end{figure}
\end{center}

\begin{center}
	\begin{figure}[ht]  
		\begin{overpic}[scale=.6, grid = false, tics=5, trim=-.5cm -1cm -1cm -.5cm, clip]{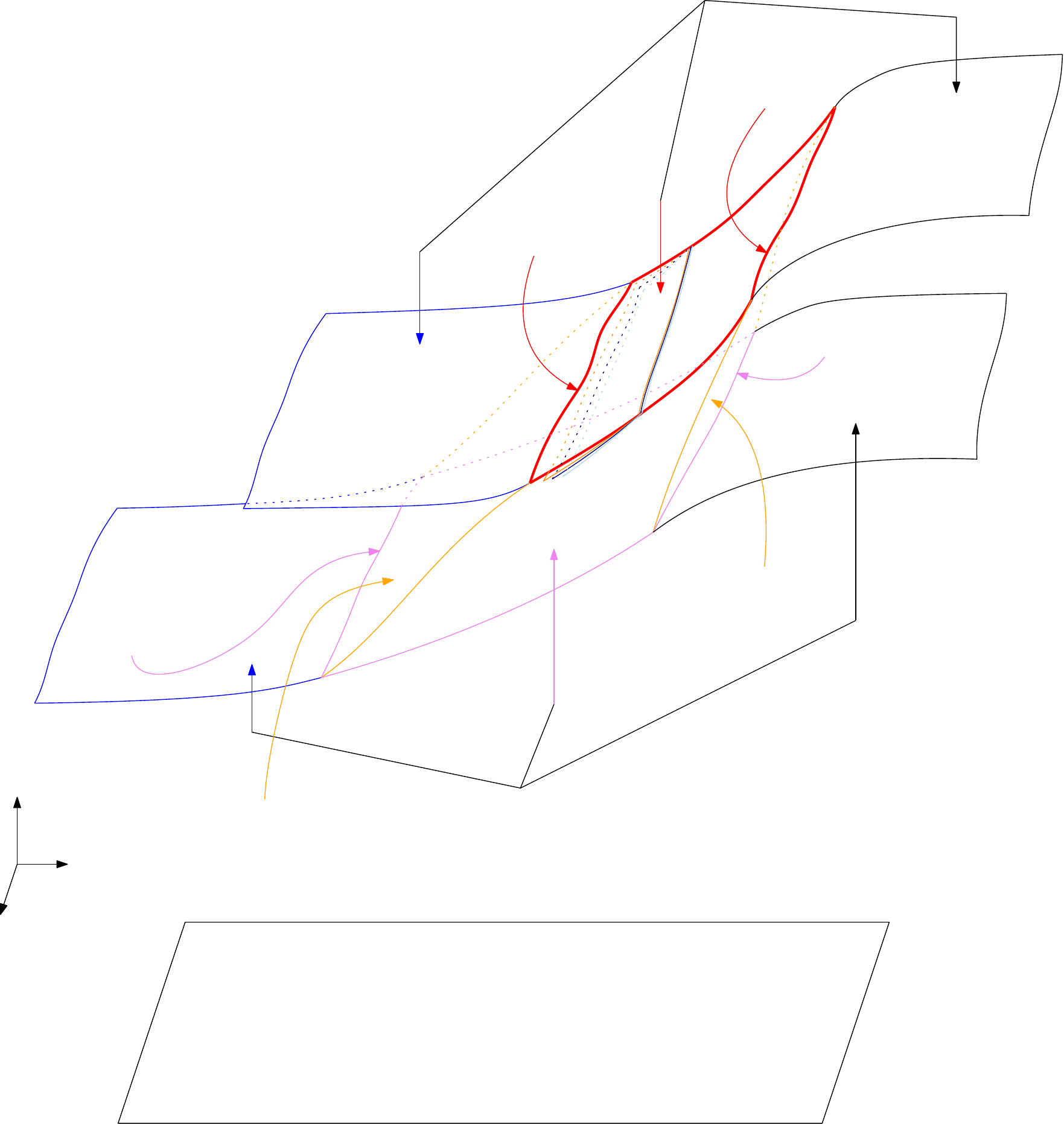}
			\put (40,10) {$\Sigma_0$}
			\put (35,24) {\textrm{Cauchy stability region}}
			\put (-2,20) {$(x^2,x^3) \in \mathbb{T}^2$}
			\put (2,29) {$t$}
			\put (9,26) {$x^1 \in\mathbb{R}$}
			\put (15,50) {$\Upsilon(\MLeft)$}
			\put (65,69) {$\Upsilon(\twoargmumuxtorus{\timefunction_0}{-\muxmulevelsetvalue_0})$}
			\put (30,28) {$\Upsilon(\inthyp{\timefunction_0}{[- \rightu,\leftu]})$}
			\put (45,99) {$\Upsilon(\inthyp{0}{[- \rightu,\leftu]})$}
			\put (75,75) {$\Upsilon(\MRight)$}
			\put (7,45) {$\Upsilon(\twoargmumuxtorus{-\timefunction_0}{0})$}
			\put (39,78) {$\Upsilon(\twoargmumuxtorus{0}{0})$}
			\put (59,89) {$\Upsilon(\twoargmumuxtorus{0}{-\muxmulevelsetvalue_0})$}
			\put (18,28) {$\Upsilon(\datahypfortimefunctiontwoarg{0}{[\timefunction_0,0]})$}
			\put (61,47) {$\Upsilon(\datahypfortimefunctiontwoarg{-\muxmulevelsetvalue_0}{[\timefunction_0,0]})$}
			\put (38,54) {$\Upsilon(\MSingular)$}
			\put (45,81) {$\Upsilon(\mathcal{B}^{[0,\muxmulevelsetvalue_0]})$}
	\end{overpic}
		\caption{The region $\Upsilon(\MInteresting)$ from Fig.\,\ref{F:INTERESTINGREGIONMAINRESULTS}, 
		i.e., $\MInteresting$ mapped into Cartesian coordinate space under $\Upsilon$}
	\label{F:INTERESTINGREGIONMAINRESULTSCARTESIAN}
	\end{figure}
\end{center}

We first provide two pictures illustrating the region that we study in the theorem.
In Fig.\,\ref{F:INTERESTINGREGIONMAINRESULTS}, we display the region
in geometric coordinates $(t,u,x^2,x^3)$. In Fig.\,\ref{F:INTERESTINGREGIONMAINRESULTSCARTESIAN},
we display the region in Cartesian coordinates $(t,x^1,x^2,x^3)$,
where in the labels, $\Upsilon(t,u,x^2,x^3) = (t,x^1,x^2,x^3)$
is the change of variables map between the two coordinate systems.

\begin{remark}[Notation involving $\Upsilon$]	
	\label{R:UPSILONNOTATIION}
	In most of the article, we consider sets such as the singular boundary $\mathcal{B}$
	to be subsets of geometric coordinate space $\mathbb{R}_t \times \mathbb{R}_u \times \mathbb{T}^2$,
	and we denote the image of these sets in Cartesian coordinate space $\mathbb{R}_t \times \mathbb{R}_{x^1} \times \mathbb{T}^2$
	by explicitly indicating the change of variables map $\Upsilon$, e.g., by
	$\Upsilon(\mathcal{B})$.
	We use this notation in particular in Fig.\,\ref{F:INTERESTINGREGIONMAINRESULTSCARTESIAN}.
	However, in many of the other figures that depict regions in Cartesian coordinate space, 
	such as Figs.\,\ref{F:INTROPICTURESOFMAINRESULTS} and \ref{F:INTROCARTESIANROUGHFOLIATIONS},
	we have suppressed the map $\Upsilon$ so as to not clutter the figure.
\end{remark}

		\begin{theorem}[Abbreviated statement of the main results] \label{T:ABBREVIATEDSTATEMENTOFMAINRESULTS}
		\hfill
		\begin{itemize}
		\item Fix any of the ``admissible'' background (shock-forming) simple isentropic plane-symmetric 
			solutions that we construct in Appendix~\ref{A:PS},
			where we define ``admissible'' in Def.\,\ref{AD:ADMISSIBLEBACKGROUND};
			for each background solution, only a single Riemann invariant\footnote{By our definition 
			\eqref{AE:RIEMANNINVARIANTS} of Riemann invariants,
			the fluid state $\RRiemannPS = 0$ corresponds to a fluid with a non-zero constant density 
			$\bar{\varrho} > 0$.}  
			is non-vanishing, and we denote it by $\RRiemannPS$.
			For each such background solution, there exist numbers 
			$\timefunction_0 < 0$, $\muxmulevelsetvalue_0 > 0$, 
			$\rightu > 0$, 
			$\farrightu > \rightu$,
			and $\leftu > 0$ 
			such that the data of $\RRiemannPS$ is
			compactly supported\footnote{We made the assumption that the initial data of the background Riemann invariant
			$\RRiemannPS$ 
			are compactly supported only for convenience; 
			this assumption could be eliminated without much additional effort.} 
			in the Cauchy hypersurface portion 
			$
			\Sigma_0^{[-\rightu,\leftu]} 
			=
			\left\lbrace 
		(0,-x^1,x^2,x^3) \ | \ 
		-x^1 \in [-\rightu,\leftu], 
			\,
		(x^2,x^3) \in \mathbb{T}^2
	\right\rbrace
$
in Cartesian coordinate space $\mathbb{R}_t \times \mathbb{R}_{x^1} \times \mathbb{T}^2$.
\item Let $\Ntop \geq 24$ be an integer.	
\item Assume that the data-norm $\mathring{\Delta}_{\Sigma_0^{[-\farrightu,\leftu]}}^{\Ntop+1}$
	defined in \eqref{E:PERTURBATIONSMALLNESSINCARTESIANDIFFERENTIALSTRUCTURE},
	which measures the $H^{\Ntop+1}(\Sigma_0^{[-\rightu,\leftu]})$-closeness of the perturbed 
	data to the background data, is sufficiently small. Note that the perturbed data do not have to be compactly supported.
\end{itemize}

Then the following conclusions hold, where the role of $\timefunction_0$ and $\muxmulevelsetvalue_0$ is described below.

\noindent \underline{\textbf{Classical existence in the Cauchy stability region}}.
		The perturbed solution exists classically with respect to the Cartesian coordinates
		and remains uniformly bounded in the ``Cauchy stability region'' 
		trapped in between the flat hypersurface portion $\Sigma_0^{[-\rightu,\leftu]}$ and 
		the curved hypersurface portion $\inthyp{\timefunction_0}{[- \rightu,\leftu]}$,
		where $\inthyp{\timefunction_0}{[- \rightu,\leftu]}$ is a portion of the $\timefunction_0$-level set of the
		time function $\newtimefunction$, described below.
		We will not further discuss the behavior of the
		solution in the Cauchy stability region, 
		noting only that it is depicted in Fig.\,\ref{F:INTERESTINGREGIONMAINRESULTSCARTESIAN}.
	
		\medskip
		
		\noindent \underline{\textbf{Acoustic geometry and regular behavior in the geometric coordinates}}.
		There exists an \textbf{eikonal function} $u$ solving the eikonal equation
		$(\gfour^{-1})^{\alpha \beta} \partial_{\alpha} u \partial_{\beta} u = 0$
		such that the \textbf{geometric coordinates} $(t,u,x^2,x^3)$ form a global coordinate system on 
		the compact subset $\MInteresting = \MLeft \cup \MSingular \cup \MRight$
		depicted in Fig.\,\ref{F:INTERESTINGREGIONMAINRESULTS}. 
		The figure is a subset of geometric coordinate space
		$\mathbb{R}_t \times \mathbb{R}_u \times \mathbb{T}^2$,
		and 
		$\MInteresting 
			\subset 
			\lbrace 
				(t,u,x^2,x^3) 
				\in 
				\mathbb{R} \times \mathbb{R} \times \mathbb{T}^2
				\ | \ - \rightu \leq u \leq \leftu 
			\rbrace$.
		The \textbf{acoustical metric}
		$\gfour$ is the fluid-solution-dependent Lorentzian metric (see \eqref{E:ACOUSTICALMETRIC}) 
		that captures the intrinsic geometry of sound waves in the flow.
		The level sets of $u$, which we denote by $\nullhyparg{u}$ and which are depicted as subsets of Cartesian coordinate space
		in Fig.\,\ref{F:INFININTEDENSITYOFCHARACTERISTICSONSINGULARBOUNDARY},
		are characteristic surfaces for the system \eqref{E:INTROTRANSPORTVI}--\eqref{E:INTROBS}.
		The fluid blowup is \underline{not visible} in the geometric coordinates:
		the fluid solution and its up-to-mid-order partial derivatives in the geometric coordinate system
		remain bounded on $\MInteresting$.
		
		\medskip
		
		\noindent \underline{\textbf{The singularity is caused by the infinite density of the characteristics}}.
		With respect to the Cartesian coordinates, there is singularity formation, 
		described below, that coincides with the vanishing of the \textbf{inverse foliation density}
		$\upmu \eqdef - \frac{1}{(\gfour^{-1})^{\alpha \beta}\partial_{\alpha} t \partial_{\beta} u}$.
		The vanishing of $\upmu$ signifies the infinite density of the level sets of $u$ 
		(viewed as a function of the Cartesian coordinates). Within $\MInteresting$, 
		this singular behavior occurs
		precisely along the singular boundary portion $\mathcal{B}^{[0,\muxmulevelsetvalue_0]}$
		depicted in Fig.\,\ref{F:INTERESTINGREGIONMAINRESULTS}.
		
		\medskip
		
		\noindent \underline{\textbf{Behavior of the change of variables map}}.
		The change of variables map $\Upsilon(t,u,x^2,x^3) = (t,x^1,x^2,x^3)$
		is a \textbf{homeomorphism} from $\MInteresting$ onto its image in Cartesian coordinate space.
		It Jacobian determinant is $\approx - \upmu$, and since $\upmu > 0$ on
		$\MInteresting \backslash \mathcal{B}^{[0,\muxmulevelsetvalue_0]}$,
		the map is a diffeomorphism on $\MInteresting \backslash \mathcal{B}^{[0,\muxmulevelsetvalue_0]}$.
		
		\medskip
		
		\noindent \underline{\textbf{The shock singularity in the Cartesian differential structure}}.
		There exists a $\Sigma_t$-tangent vectorfield $X$ 
		of Euclidean length $\sqrt{\sum_{a=1}^3 (X^a)^2} \approx 1$ such that the following occurs.
		Within the compact subset
		$\Upsilon(\MInteresting) = \Upsilon(\MLeft \cup \MSingular \cup \MRight)$ 
		of Cartesian coordinate space
		depicted in Fig.\,\ref{F:INTERESTINGREGIONMAINRESULTSCARTESIAN},
		with respect to the Cartesian coordinates,
		$|X^a \partial_a v^1|$
		and
		$|X^a \partial_a \varrho|$
		blow up precisely on the \textbf{singular boundary} portion
		$\Upsilon(\mathcal{B}^{[0,\muxmulevelsetvalue_0]})$, described below.
		The solution is smooth in $\Upsilon\left(\MInteresting \backslash \mathcal{B}^{[0,\muxmulevelsetvalue_0]} \right)$,
		and for $\alpha = 0,1,2,3$, $i=1,2,3$, and $A=2,3$,
		the fluid quantities 
		$\varrho$, 
		$v^i$, 
		$\Ent$, 
		$\Flatcurl v^i$, 
		$\partial_i \Ent$, 
		and 
		$\gfour_{ab} \Yvf{A}^a \partial_{\alpha} v^b$
		remain bounded on all of $\Upsilon(\MInteresting)$,
		where $\Yvf{2}$ and $\Yvf{3}$ are the vectorfields from \eqref{E:INTROTANGENTCOMMUTATORS}, 
		and they are tangent to $\nullhyparg{u} \cap \Sigma_t$.
		Moreover, the up-to-mid-order derivatives of these quantities
		with respect to the $\nullhyparg{u}$-tangent vectorfields
		$\lbrace \Lunit, \Yvf{2}, \Yvf{3} \rbrace$
		from \eqref{E:INTROTANGENTCOMMUTATORS}
		also remain bounded on all of $\Upsilon(\MInteresting)$.
		
		\medskip
		
		\noindent \underline{\textbf{The structure of $\mathcal{B}^{[0,\muxmulevelsetvalue_0]}$ and $\partial_- \mathcal{B}^{[0,\muxmulevelsetvalue_0]}$}}.
		The singular boundary $\mathcal{B}^{[0,\muxmulevelsetvalue_0]}$, 
		viewed as a subset of geometric coordinate space
		$\mathbb{R}_t \times \mathbb{R}_u \times \mathbb{T}^2$,
		is contained in the top boundary of $\MInteresting$
		and is a subset of the level set $\lbrace \upmu = 0 \rbrace$, specifically, a portion that 
		can be realized as the limit as $\mulevelsetvalue \downarrow 0$ of portions of the level sets
		$\lbrace (t,u,x^2,x^3) \ | \ \upmu(t,u,x^2,x^3) = \mulevelsetvalue \rbrace$ that are
		either null or spacelike with respect to the acoustical metric $\gfour$.
		$\mathcal{B}^{[0,\muxmulevelsetvalue_0]}$ is
		an embedded $3$-dimensional submanifold-with-boundary of geometric coordinate space
		$\mathbb{R}_t \times \mathbb{R}_u \times \mathbb{T}^2$.
		Moreover, we can decompose 
		$\mathcal{B}^{[0,\muxmulevelsetvalue_0]} = \bigcup_{\muxmulevelsetvalue \in [0,\muxmulevelsetvalue_0]} \twoargmumuxtorus{0}{-\muxmulevelsetvalue}$,
		where each $\twoargmumuxtorus{0}{-\muxmulevelsetvalue}$
		is a $C^{1,1}$, $2$-dimensional, spacelike submanifold of geometric coordinate space
		$\mathbb{R}_t \times \mathbb{R}_u \times \mathbb{T}^2$
		that is a graph over $\mathbb{T}^2$
		such that $\upmu X \upmu \equiv - \muxmulevelsetvalue$ along $\twoargmumuxtorus{0}{-\muxmulevelsetvalue}$.
		The torus $\twoargmumuxtorus{0}{0}$,
		which we refer to as the \textbf{crease} and also denote by $\partial_- \mathcal{B}^{[0,\muxmulevelsetvalue_0]}$, 
		is the past boundary of $\mathcal{B}^{[0,\muxmulevelsetvalue_0]}$.
		
		\medskip
		
		\noindent \underline{\textbf{The structure of $\Upsilon(\mathcal{B}^{[0,\muxmulevelsetvalue_0]})$ and 
		$\Upsilon(\partial_- \mathcal{B}^{[0,\muxmulevelsetvalue_0]})$}}.
		The change of variables map $\Upsilon(t,u,x^2,x^3) = (t,x^1,x^2,x^3)$
		is a \textbf{homeomorphism} from the singular boundary portion 
		$\mathcal{B}^{[0,\muxmulevelsetvalue_0]}$ 
		onto its image $\Upsilon(\mathcal{B}^{[0,\muxmulevelsetvalue_0]})$ 
		in Cartesian coordinate space. Moreover, $\Upsilon$ is a diffeomorphism from 
		$\mathcal{B}^{[0,\muxmulevelsetvalue_0]} \backslash \partial_- \mathcal{B}^{[0,\muxmulevelsetvalue_0]}$
		onto its image 
		$\Upsilon(\mathcal{B}^{[0,\muxmulevelsetvalue_0]} \backslash \partial_- \mathcal{B}^{[0,\muxmulevelsetvalue_0]})$.
		The image set $\Upsilon(\mathcal{B}^{[0,\muxmulevelsetvalue_0]})$ 
		is an embedded $3$-dimensional submanifold-with-boundary in Cartesian coordinate space
		that has regularity\footnote{The $C^{1,1/2}$ embedding of $\Upsilon(\mathcal{B}^{[0,\muxmulevelsetvalue_0]})$ is
		provided by Prop.\,\ref{P:DESCRIPTIONOFSINGULARBOUNDARYINCARTESIANSPACE}; it is the map
		$(z,x^2,x^3) \rightarrow \Upsilon \circ \embeddingofsingularboundaryintogeometriccoordinatespace(\sqrt{z},x^2,x^3)$
		from the proposition.
		\label{FN:REGULARITYOFSINGULARBOUNDARYINCARTESIANSPACE}.} 
		$C^{1,1/2}$ with respect to the Cartesian coordinates.
		In addition, $\Upsilon(\mathcal{B}^{[0,\muxmulevelsetvalue_0]} \backslash \partial_- \mathcal{B}^{[0,\muxmulevelsetvalue_0]})$
		is a null hypersurface with respect to the acoustical metric $\gfour$ on $\Upsilon(\MInteresting)$,
		and 	
		$\Upsilon(\mathcal{B}^{[0,\muxmulevelsetvalue_0]} \backslash \partial_- \mathcal{B}^{[0,\muxmulevelsetvalue_0]})$ is ruled,
		in a degenerate sense explained in Remark~\ref{R:NONUNIQUENESSOFINTEGRALCURVESOFLUNIT}, 
		by integral curves of the $\gfour$-null vectorfield $\Lunit = \Lunit^{\alpha} \partial_{\alpha}$
		on $\Upsilon(\MInteresting)$.
		Furthermore, 
		for $\muxmulevelsetvalue \in [0,\muxmulevelsetvalue_0]$, 
		the images $\Upsilon(\twoargmumuxtorus{0}{-\muxmulevelsetvalue})$ of the tori $\twoargmumuxtorus{0}{-\muxmulevelsetvalue}$,
		including the crease $\partial_- \mathcal{B}^{[0,\muxmulevelsetvalue_0]} = \twoargmumuxtorus{0}{0}$,
		are $2$-dimensional, $C^{1,1}$, $\gfour$-spacelike submanifolds of Cartesian coordinate space
		$\mathbb{R}_t \times \mathbb{R}_{x^1} \times \mathbb{T}^2$
		that are graphs over $\mathbb{T}^2$.
		
		\medskip
	
		\noindent \underline{\textbf{Rough time functions reveal the tori foliating $\mathcal{B}^{[0,\muxmulevelsetvalue_0]}$}}.
		There exists a one-parameter family of \textbf{rough time functions}
		$\lbrace \timefunctionarg{\muxmulevelsetvalue} \rbrace_{\muxmulevelsetvalue \in [0,\muxmulevelsetvalue_0]}$,
		each with range $[\timefunction_0,0]$,
		such that the level sets $\lbrace \timefunctionarg{\muxmulevelsetvalue} = \timefunction \rbrace$
		with $\timefunction \in [\timefunction_0,0)$
		do not intersect $\mathcal{B}^{[0,\muxmulevelsetvalue_0]}$,
		while $\lbrace \timefunctionarg{\muxmulevelsetvalue} = 0 \rbrace \cap \mathcal{B}^{[0,\muxmulevelsetvalue_0]} = \twoargmumuxtorus{0}{-\muxmulevelsetvalue}$.
		That is, each rough time function reveals the structure of the two-dimensional torus $\twoargmumuxtorus{0}{-\muxmulevelsetvalue}$
		in $\mathcal{B}^{[0,\muxmulevelsetvalue_0]}$, 
		and $\twoargmumuxtorus{0}{-\muxmulevelsetvalue} \subset \lbrace \timefunctionarg{\muxmulevelsetvalue} = 0 \rbrace$.
		Moreover, each $\timefunctionarg{\muxmulevelsetvalue}$ is one degree less differentiable with respect
		to the geometric coordinates than the fluid solution,
		and the tori $\twoargmumuxtorus{0}{-\muxmulevelsetvalue}$ are two degrees
		less differentiable than the fluid solution.
		
		Moreover, there exists a rough time function $\newtimefunction$ with range $[\timefunction_0,0]$ that foliates
		$\MInteresting$.
		$\newtimefunction$ is
		$C^{1,1}$ with respect to the geometric coordinates - \textbf{and not more regular};
		see Remark~\ref{R:NEWTIMEFUNCTIONISC11ANDNOTBETTERANDCONNECTIONTOCAUSALSTRUCTUREOFMUZEROLEVELSET}.
		Two of its level sets intersected with $\lbrace u \in [- \rightu,\leftu] \rbrace$, 
		namely $\inthyp{0}{[- \rightu,\leftu]}$ and
		$\inthyp{\timefunction_0}{[- \rightu,\leftu]}$, are depicted in Fig.\,\ref{F:INTERESTINGREGIONMAINRESULTS}.
		Finally, $\mathcal{B}^{[0,\muxmulevelsetvalue_0]} \subset \lbrace \newtimefunction = 0 \rbrace$.
		
		\end{theorem}
		
		\begin{remark}[Remarks on our initial data assumptions]
		\label{R:REMARKSONINITIALDATAASSUMPTIONS}
			In our main results,
			we have chosen to assume the smallness of the norm $\mathring{\Delta}_{\Sigma_0^{[-\farrightu,\leftu]}}^{\Ntop+1}$
			defined in \eqref{E:PERTURBATIONSMALLNESSINCARTESIANDIFFERENTIALSTRUCTURE}
			because this immediately allows us to conclude that 
			our results hold for open sets of initial data that are $H^{\Ntop+1}$-close 
			to the data of a background solution.
			However, the main assumptions on the initial data that we use in our PDE analysis
			are actually stated in Sect.\,\ref{S:ASSUMPTIONSONTHEDATA},
			in terms of data-size parameters that satisfy 
			assumptions stated in Sect.\,\ref{S:PARAMETERSANDSIZEASSUMPTIONSANDCONVENTIONSFORCONSTANTS}.
			All these assumptions follow as consequences of
			the smallness of the norm $\mathring{\Delta}_{\Sigma_0^{[-\farrightu,\leftu]}}^{\Ntop+1}$.
			That is, the proof of our main results would go through under only 
			the assumptions stated in Sect.\,\ref{S:ASSUMPTIONSONTHEDATA}
			and the parameter-size assumptions stated in Sect.\,\ref{S:PARAMETERSANDSIZEASSUMPTIONSANDCONVENTIONSFORCONSTANTS},
			modulo the remarks we make at the beginning of Sect.\,\ref{SS:CONTROLOFDATAFORROUGHTORIENERGYESTIMATES}.
		\end{remark}

\subsection{Remarks on the main results and methods}
\label{SS:REMARKSONRESULTSANDMETHODS}
Before proceeding, we make a series of remarks about our main results
and our proof framework,
with an emphasis on the new methods we introduce in this paper,
how they connect to methods developed in other papers,
and how they connect to open problems.

\begin{itemize}
		\item (\textbf{Use of acoustic geometry to detect the singularity and to ``hide it'' by unfolding the characteristics}).
		As in Christodoulou's breakthrough work \cite{dC2007} on irrotational and isentropic shock formation, 
		our analysis fundamentally relies 
		on nonlinear geometric optics (which we also loosely refer to as ``the acoustic geometry''), 
		implemented via an eikonal function $u$, which solves
		the eikonal equation $(\gfour^{-1})^{\alpha \beta} \partial_{\alpha} u \partial_{\beta} u = 0$.
		The \emph{acoustical metric} $\gfour$ is a Lorentzian metric 
		whose Cartesian components $\gfour_{\alpha \beta}$
		are functions of the fluid variables;
		see \eqref{E:ACOUSTICALMETRIC} for the precise formula.
		The level sets $\nullhyparg{u}$ of $u$ are null hypersurfaces (also known as ``characteristics'' or ``characteristic hypersurfaces''),
		and they correspond to the propagation of sound waves.
		Infinite density of the level sets of $u$ (viewed as a function of the Cartesian coordinates) 
		signifies the formation of a shock,
		and in the regime under study, it coincides with the blowup of first-order partial derivatives of $v$ and $\varrho$
		in directions \emph{transversal} to the $\nullhyparg{u}$.
 		As is shown in Fig.\,\ref{F:INFININTEDENSITYOFCHARACTERISTICSONSINGULARBOUNDARY},
		\emph{these phenomena occur along the entire singular boundary.}
		Moreover, with the help of $u$, we can construct a ``geometric coordinate system''
		$(t,u,x^2,x^3)$, relative to which the solution remains rather smooth.
		This is crucial for our derivation of PDE estimates up to top-order.
	\item (\textbf{The structure of the singular boundary and the crease}).
		We prove that in the solution regime under study,
		relative to a differential structure on spacetime tied to the eikonal function,
		the singular boundary $\mathcal{B}$ has the structure of a $3D$ submanifold-with-boundary;
		see Fig.\,\ref{F:MAXDEVELOPMENTWITHSHOCKHYPERSURFACEINCARTESIAN}.
		Its past boundary $\partial_- \mathcal{B}$ is the crease, 
		which we prove is a $2D$ acoustically spacelike submanifold,
		where our notion of ``spacelike'' is 
		relative to $\gfour$.
		These structures are stable, and their availability is fundamental for our approach.
		These structures also have important implications for the shock development problem, 
		described below.
		\item (\textbf{Prior works}). 
		There are many prior works on shock formation for compressible Euler solutions,
		including Riemann's famous work \cite{bR1860} in one spatial dimension.
		Although much less is known in multi-dimensions, there has been dramatic progress over the last several decades,
		starting with Alinhac's foundational works on quasilinear wave equations
		\cites{sA1999a,sA1999b,sA2002} and
		Christodoulou's breakthrough monograph on the irrotational and isentropic 
		relativistic Euler equations \cite{dC2007};
		see Sect.\,\ref{SS:BACKGROUNDONSHOCKS} for a discussion of some additional key developments in the history of the subject.
		While prior results have led to a revolution in our understanding of multi-dimensional
		shock formation, they all were limited in one or more of the following ways:
		\textbf{i)} they treated only irrotational and isentropic solutions;
		\textbf{ii)} the methods allowed one only to follow the solution to the constant-Cartesian-time 
			hypersurface of first blowup;
		\textbf{iii)} the methods applied only to fully non-degenerate singularities (see below for their definition)
			which, as it turns out, corresponds to understanding the blowup at the unique first (relative to Cartesian time)
			point that is contained inside a strictly convex crease;
		\textbf{iv)} the methods yielded a description only of some implicit portion of the 
			boundary of the maximal development of the data
			and in particular yielded only the portion of the crease
			that is tangent to some \emph{flat} spacelike hypersurface $\Sigma$
			such that a neighborhood of the crease lies in its future.
			A key point is that without strict convexity, 
			there can be points in the crease such that
			\underline{no open neighborhood of them is accessible}
			through such an approach. 
			In particular, for the solutions
			we treat in our main results, the full structure of the crease
			is not accessible through any of these approaches;
			see Fig.\,\ref{F:MAXDEVELOPMENTWITHSHOCKHYPERSURFACEINCARTESIAN}.
	\item (\textbf{Relevance for the shock development problem}).
		The crease and its sub-manifold structure are crucial ingredients 
		needed to set up the \emph{shock development problem}, 
		which, as we mentioned earlier,
		is the problem of describing the transition of the solution
		from classical to weak, past the initial singularity.
		It turns out that the crease is part of the maximal classical development
		\emph{and}, once one constructs the weak solution, it will also be part of the weak solution;
		see  \cite{dC2019}*{Section~1.5} for a detailed discussion of the connection between
		the maximal classical development and the weak solution.
		Away from symmetry, the shock development problem is an outstanding open problem.
		In Christodoulou's approach to studying shock developments 
		under the simplifying assumption that the vorticity and entropy are vanishing\footnote{The irrotational and isentropic weak solutions 		
		constructed by Christodoulou in \cite{dC2019} are not solutions to the compressible Euler equations because they do not 
		respect the jump in entropy and vorticity that must occur across the shock hypersurface. Instead, they are solutions
		to a closely related hyperbolic PDE system that is equivalent to the compressible Euler equations for \emph{classical}
		isentropic and irrotational classical solutions. \label{FN:RESTRICTEDSHOCKDEVELOPMENTPROBLEM}}
		across the shock hypersurface \cite{dC2019}, 
		the rest of the singular boundary is also important because it
		plays the role of a mathematical \emph{barrier} in the construction
		of the weak solution. However, aside from the crease, the singular boundary
		is not part of the weak solution. More precisely, in general, 
		there is a portion of the maximal classical development that does not
		agree with the weak solution because the classical development does not account for the shock hypersurface.
		The shock hypersurface, which we denote by $\mathcal{K}$ in Fig.\,\ref{F:MAXDEVELOPMENTWITHSHOCKHYPERSURFACEINCARTESIAN},
		is not part of the constructions of this paper, and we mention it to highlight that the 
		weak solution -- once it is constructed -- will disagree with the classical solution in the region
		in between $\mathcal{B}$ and $\mathcal{K}$.
		We also highlight that $\mathcal{K}$ is supersonic relative to the acoustical metric in the 
		maximal classical development, 
		and that $\mathcal{K}$ \emph{is ``allowed'' to emerge from the crease only if one imposes a weak
		formulation of the flow}, i.e., there is no such thing as ``$\mathcal{K}$'' in the classical formulation.
	\item (\textbf{The Cauchy horizon}).
		In our forthcoming paper, we will describe the emergence of a Cauchy horizon
		from the crease. The Cauchy horizon, denoted by $\underline{\mathcal{C}}$ in
		Fig.\,\ref{F:MAXDEVELOPMENTWITHSHOCKHYPERSURFACEINCARTESIAN},
		is a $\gfour$-null hypersurface with past boundary equal to the crease,
		and it is also a crucial ingredient in the setting up the shock development problem.
		Like the crease, the Cauchy horizon is part of the maximal classical development
		and the weak solution (once one constructs it). However, in that paper, we will show that
		the solution remains smooth up to the Cauchy horizon (away from the crease, that is),
		which is in stark contrast to what happens along the singular boundary.
	\item (\textbf{Solution regimes other than perturbations of simple isentropic plane-waves}).
		Simple isentropic plane-waves, mentioned already in 
		Theorem~\ref{T:ABBREVIATEDSTATEMENTOFMAINRESULTS}, are compressible Euler solutions such that 
		$\RRiemann = \RRiemann(t,x^1)$
	and
	$\LRiemann = v^2 = v^3 \equiv 0$,
	where the ``almost\footnote{For isentropic plane-symmetric solutions,
	$\almostRiemann_{(\pm)}$ \emph{are} Riemann invariants.} 
	Riemann invariants'' $\almostRiemann_{(\pm)}$
	are defined in Def.\,\ref{D:ALMOSTRIEMANNINVARIANTS}.
	Although for definiteness we have focused on general small perturbations of such solutions,
	the techniques we develop here are robust and can be applied to other solution regimes of physical interest,
	such as the regime corresponding to small, spatially-decaying perturbations\footnote{In this regime, the dispersive tendency
	of sound waves competes against the transport phenomena associated to vorticity and entropy, 
	possibly leading to exceptionally complicated long-time dynamics.}
	of non-vacuum constant fluid states in $\mathbb{R}^{1+3}$,
	where dispersive effects play a fundamental role in the dynamics.
	\item 	(\textbf{Building blocks that can accommodate degeneracies}).
	We have focused our attention on perturbations of
	simple isentropic plane-symmetric waves because such solutions are of physical interest
	and, from the point of view of analysis, they are challenging to study because
	they exhibit degeneracies tied to lack of strict convexity (as we mentioned above)
	of their singular boundaries. This is the first paper to fully grapple with 
	these degeneracies and as such, our results are new even in the sub-class of irrotational
	and isentropic solutions.
	The presence of these degeneracies forced us 
	to develop techniques that we anticipate can be used, 
	in view of finite speed of propagation, 
	as building blocks to study 
	the global structure of much more general shock-forming solutions across space.
	\item (\textbf{New, solution-dependent foliations are needed}). 
		To access the entire crease/singular boundary, in general, one cannot exclusively rely on arguments based
		on analysis on regions bounded by the characteristic surfaces $\nullhyparg{u}$,
		surfaces $\Sigma_t$ of constant Cartesian time,
		or, for that matter, any other family of surfaces that are ``pre-specified'' in the sense that they are
		explicitly parameterized with respect to the Cartesian coordinates
		(e.g., one can not generally use ``tilted'' spacelike hypersurfaces that are planes with respect to the Cartesian coordinates).
		In particular, the crease is typically not contained in a fixed $\nullhyparg{u}$
		or $\Sigma_t$; see Fig.\,\ref{F:MAXDEVELOPMENTINCARTESIAN}.
		Hence, one of our key new ideas in the paper is to replace the surfaces $\Sigma_t$ with better ones.
		That is, we construct foliations by spacelike hypersurfaces 
		that are precisely adapted to the shape of the singular boundary
		(though we also fundamentally rely on foliations by the $\nullhyparg{u}$). 
		We construct these foliations by ``flowing out'' (see below), along
		the integral curves of a ``generating vectorfield'' that is transversal to $\nullhyparg{u}$,
		from co-dimension $2$ topological tori that are adapted to the anticipated shape of the singularity.
		In Fig.\,\ref{F:CARTESIANROUGHFOLIATIONCREASE}, the generating vectorfield is denoted by $\Wtransarg{0}$, and
		crucially, it is also transversal to the hypersurface $\datahypfortimefunctionarg{0}$,
		which is foliated by the tori, two of which are denoted by
		$\twoargmumuxtorus{0}{0}$
		and
		$\twoargmumuxtorus{-\timefunction}{0}$.
		Similarly, 
		in Fig.\,\ref{F:CARTESIANROUGHFOLIATIONFUTUREOFCREASE},
		the generating vectorfield is denoted by $\Wtransarg{\muxmulevelsetvalue}$, and
		it is also transversal to the hypersurface $\datahypfortimefunctionarg{- \muxmulevelsetvalue}$,
		which is foliated by the tori, two of which are denoted by
		$\twoargmumuxtorus{0}{-\muxmulevelsetvalue}$
		and
		$\twoargmumuxtorus{-\timefunction}{-\muxmulevelsetvalue}$.
		While there have been many 
		other works on the formation of shocks without symmetry assumptions, described in 
		Sect.\,\ref{SS:WORKSBUILDINGTOWARDSMAINTHEOREM},
		those works have relied on a ``background'' time function (typically the Cartesian time function $t$)
		and the corresponding foliations.
		Note that the PDE analysis of compressible Euler flow, 
		which inevitably involves some kind of energy estimates with respect to the spacelike foliations,
		\emph{necessarily halts} when the spacelike foliations first intersect the
		singular boundary, because singularities form at the intersection points.
		This limits the portion of the singular
		boundary that can be detected through background foliations.
	\item (\textbf{Remarks on the simpler fully non-degenerate sub-regime}).
		Despite the previous point, we note that there is one sub-regime in which background time functions
		\emph{can} be used to detect the structure of the crease
		(though perhaps not a full neighborhood of it in the singular boundary):
		the regime in which the crease and singular boundary are acoustically strictly convex
		(see below and Sect.\,\ref{SS:OVERVIEWOFDIFFICULTIES}). 
		In Fig.\,\ref{F:STRICTLYCONVEX}, we depict an acoustically strictly convex singular boundary.
			
			\renewcommand{\thesubfigure}{\Alph{subfigure}} 
\begin{figure}[ht]
\centering
\begin{subfigure}{.5\textwidth}
 \centering
\begin{overpic}[scale=.36,grid=false,tics=5]{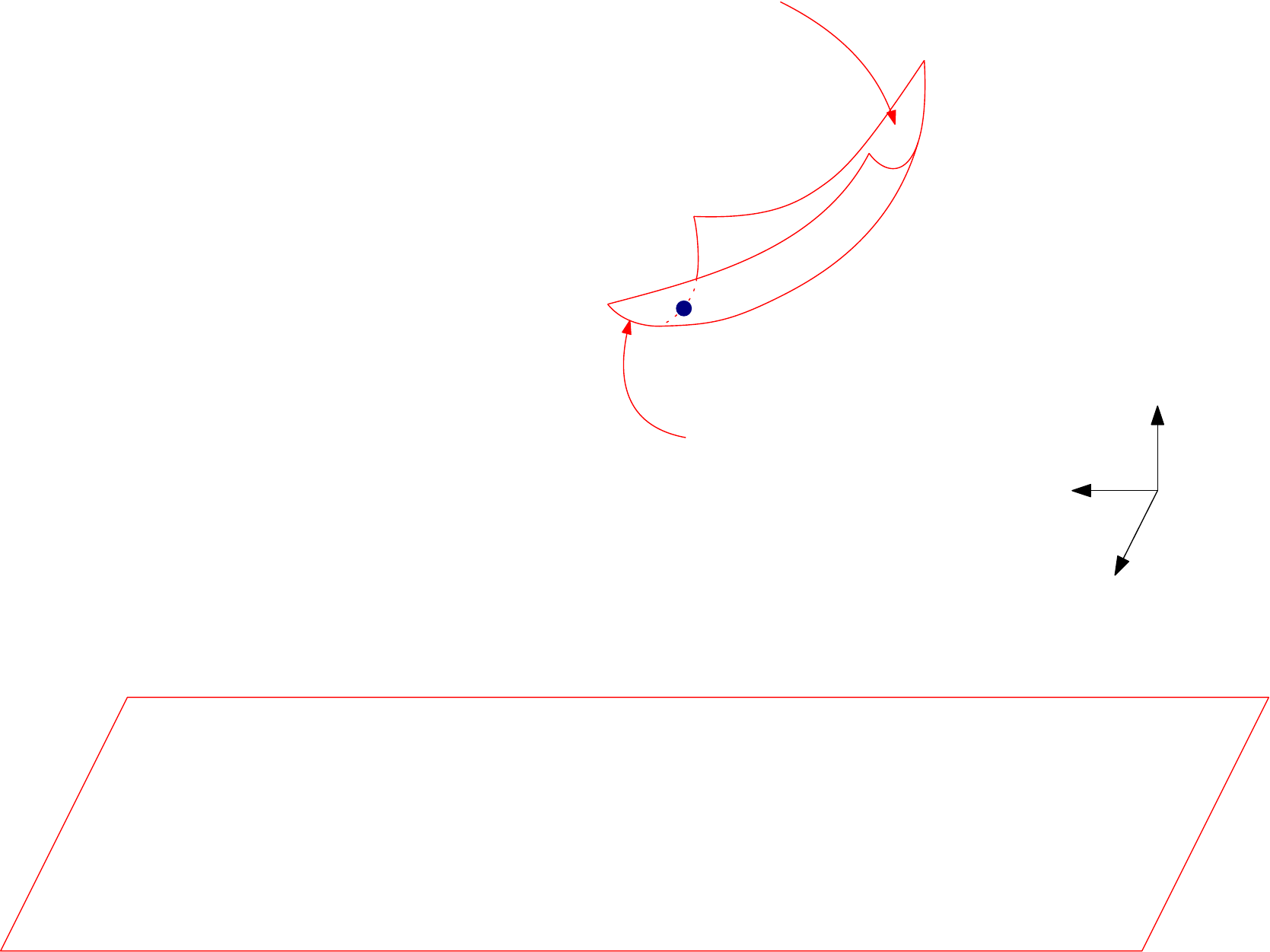} 
			\put (50,10) {$\Sigma_0$}
			\put (67,24) {$(x^2,x^3) \in \mathbb{T}^2$}
			\put (92,38) {$t$}
			\put (70,34) {$u \in\mathbb{R}$}
			\put (48,35) {$\partial_- \mathcal{B}$}
			\put (60,76) {$\mathcal{B}$}
			\put (52,44) {$b_*$}
\end{overpic}
  \caption{An acoustically strictly convex crease and singular boundary in geometric coordinate space}
  \label{F:STRICTLYCONVEXGEOMETRIC}
\end{subfigure}%
\begin{subfigure}{.5\textwidth}
  \centering
  	\begin{overpic}[scale=.36, grid = false, tics=5, trim= 0cm 0cm 0cm -1.5cm, clip]{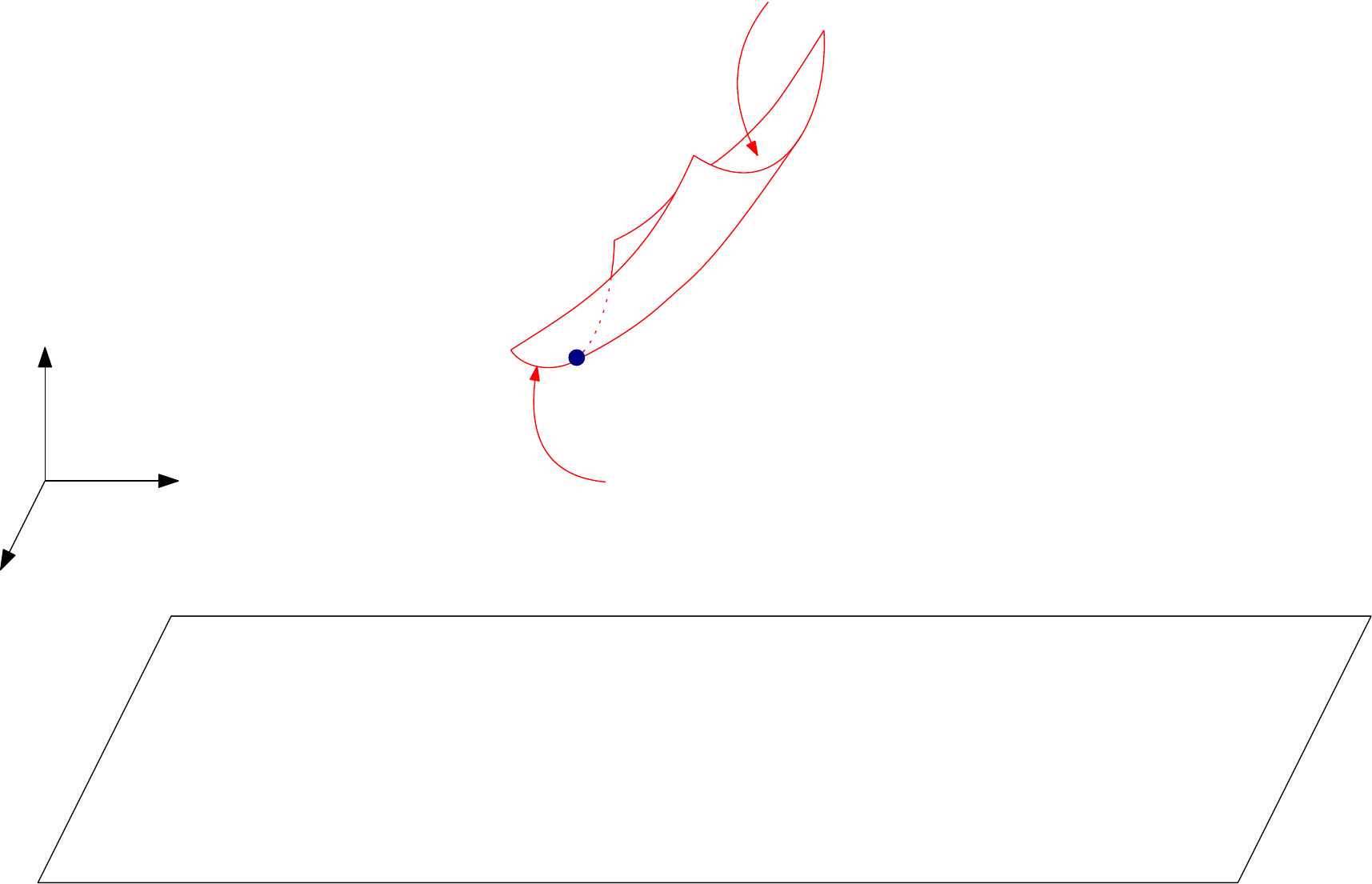}
			\put (50,10) {$\Sigma_0$}
			\put (1,22) {$(x^2,x^3) \in \mathbb{T}^2$}
			\put (.5,35) {$t$}
			\put (14,28) {$x^1 \in \mathbb{R}$}
			\put (39,24) {$\partial_- \mathcal{B}$}
			\put (54,65) {$\mathcal{B}$}
			\put (40.5,33) {$b_*$}
		\end{overpic}
		\caption{The same crease and singular boundary in Cartesian coordinate space}
		\label{F:STRICTLYCONVEXCARTESIAN}
\end{subfigure}%
\caption{Acoustically strictly convex crease and singular boundary}
\label{F:STRICTLYCONVEX}
\end{figure}
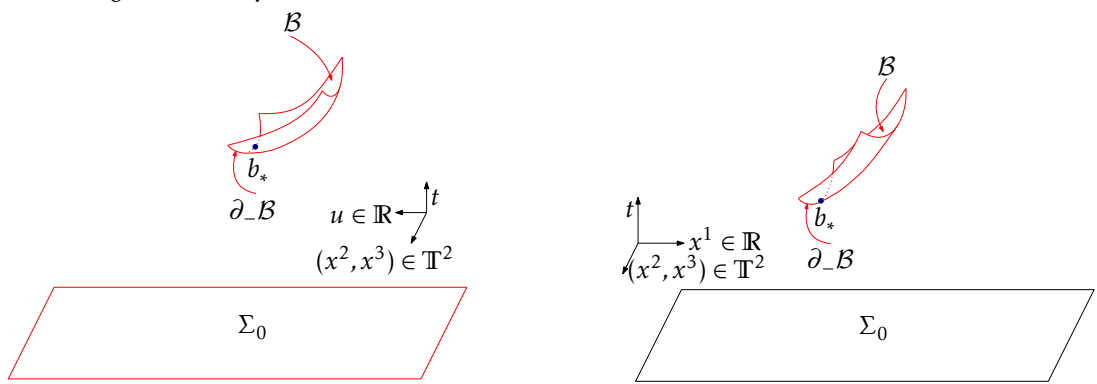
		We emphasize that symmetric solutions and their general perturbations fall outside of 
		the acoustically strictly convex regime.
		Away from symmetry, the acoustically strictly convex regime was first studied by Alinhac 
		\cites{sA1999a,sA1999b,sA2002}
		in the context of shock-forming solutions
		to quasilinear wave equations, where he followed the solution to the constant-Cartesian-time hypersurface 
		of first blowup.
		In particular, due to the strict convexity of the crease, the singular set 
		within the constant-Cartesian-time hypersurface of first blowup is an \emph{isolated point},
		denoted by $b_*$ in Fig.\,\ref{F:STRICTLYCONVEX},
		and the isolated nature of this point was fundamental for his approach. 
		We refer to such singularities as \emph{fully non-degenerate};
		see Sect.\,\ref{SSS:ALINHACSBLOWUPRESULTS} for further discussion.
	\item  (\textbf{Acoustical transversal convexity})
	Here and throughout, by an ``acoustically strictly convex'' 
	hypersurface $\mathcal{B}$, we mean that relative to a system of geometric
	coordinates $(t,u,x^2,x^3)$ tied to an acoustic eikonal function $u$, 
the tangent planes to $\mathcal{B}$ lie below $\mathcal{B}$;
see Fig.\,\ref{F:STRICTLYCONVEXGEOMETRIC} for an example of such a surface.
The key point is that acoustical strict convexity 
is \emph{absent}
in the singular boundaries of open sets of physically relevant shock-forming solutions,
including perturbations of isentropic plane-symmetric simple waves.
To prove our main results, we rely on an \emph{acoustical transversal convexity} 
(we will refer to it as ``transversal convexity'' for short),
which is substantially weaker than strict convexity 
and it is close to optimal in the following sense: 
without transversal convexity, 
the crease can fail to have the structure of a $2D$ sub-manifold,
which would drastically alter the qualitative character of the singular boundary.
We highlight that transversal convexity can be ensured along a portion of the singular boundary
containing the crease by making appropriate
open assumption on the initial data. In the present paper, we  
ensure transversal convexity via the data assumptions stated in \eqref{E:DATATASSUMPTIONMUTRANSVERSALCONVEXITY}.
Roughly, by an ``acoustically transversally convex'' 
hypersurface $\mathcal{B}$,
we mean that relative to the geometric 
coordinates $(t,u,x^2,x^3)$, there exists 
a family of \emph{lines} that are tangent to $\mathcal{B}$  -- as opposed to the full tangent plane --
such that
\textbf{i)} the lines lie below\footnote{Our approach to studying the singular boundary is local,
and thus, in our main results, we only need the tangent lines to locally lie below it.} 
$\mathcal{B}$
and \textbf{ii)} the lines are \emph{transversal} to the level sets of $u$.
In the context of Fig.\,\ref{F:INTERESTINGREGIONMAINRESULTS},
\emph{lack} of \underline{strict convexity} is exhibited by the observation that at $b_1$,
the straight line that is tangent to $\twoargmumuxtorus{0}{0}$ (this line is not drawn in the figure)
lies \emph{above} $\mathcal{B}^{[0,\muxmulevelsetvalue_0]}$,
and similarly, at $b_2$,
the straight line that is tangent to $\twoargmumuxtorus{0}{-\muxmulevelsetvalue_0}$
lies above $\mathcal{B}^{[0,\muxmulevelsetvalue_0]}$.
In Fig.\,\ref{F:INTERESTINGREGIONMAINRESULTS},
\underline{transversal convexity} is exhibited, for example,
by the fact that along the curve in $\mathcal{B}^{[0,\muxmulevelsetvalue_0]}$ joining $b_1$ and $b_2$,
the tangent lines (which are not drawn) lie \emph{below} $\mathcal{B}^{[0,\muxmulevelsetvalue_0]}$.
\item (\textbf{Rough time functions}).
			Our new foliations are level sets of a $\muxmulevelsetvalue$-indexed family of
		``rough time functions,'' denoted by $\timefunctionarg{\muxmulevelsetvalue}$,
		where $\muxmulevelsetvalue$ is a non-negative real parameter. 
		The $\timefunctionarg{\muxmulevelsetvalue}$
		solve a $\muxmulevelsetvalue$-dependent transport equation (see \eqref{E:TRANSPORTEQUATIONFORROUGHTIMEFUNCTION})
		with coefficients 
		and initial data (see \eqref{E:INITIALCONDITIONFORROUGHTIMEFUNCTION})
		determined by the acoustic geometry
		and the up-to-first-order derivatives of the solution.
		Moreover, we choose the initial 
		data surface for $\timefunctionarg{\muxmulevelsetvalue}$ 
		to be a level set of a first derivative of the acoustic geometry,
		more precisely, the $-\muxmulevelsetvalue$ level set of a first derivative of the null lapse in a direction transversal 
		to the characteristics. More precisely, each level set of
		$\timefunctionarg{\muxmulevelsetvalue}$ is the solution a transport equation
		with (constant) initial data given on a topological torus that itself
		is equal to the intersection of a level set of the null lapse with a level set
		of one of its transversal derivatives. 
		In Fig.\,\ref{F:INTROCARTESIANROUGHFOLIATIONS}, we denote the initial data hypersurfaces for the transport equations
		by $\datahypfortimefunctionarg{0}$
		and $\datahypfortimefunctionarg{-\muxmulevelsetvalue}$,
		we denote some tori in the data hypersurfaces by
		$\twoargmumuxtorus{-\timefunction}{0}$,
		$\twoargmumuxtorus{0}{0}$, 
		$\twoargmumuxtorus{-\timefunction}{-\muxmulevelsetvalue}$,
		and
		$\twoargmumuxtorus{0}{-\muxmulevelsetvalue}$,
		and we denote the vectorfields that transport $\timefunctionarg{0}$ and $\timefunctionarg{\muxmulevelsetvalue}$
		respectively by $\Wtransarg{0}$ and $\Wtransarg{\muxmulevelsetvalue}$.
		These vectorfields are tangent to the level sets of
		$\timefunctionarg{0}$ and $\timefunctionarg{\muxmulevelsetvalue}$ respectively,
		and they are also tangent to the singular boundary along
		the tori $\twoargmumuxtorus{0}{0}$ and $\twoargmumuxtorus{0}{-\muxmulevelsetvalue}$
		respectively.
		We emphasize that this construction is quite delicate and that
		\textbf{it is crucial that we consider only the case $\muxmulevelsetvalue \geq 0$ in order to ensure that
		$\timefunctionarg{\muxmulevelsetvalue}$ has a $\gfour$-timelike normal or $\gfour$-null normal}.
		In total, the construction is ``dynamic'' and ``solution-adapted,'' 
		and we achieve it through a bootstrap argument.
		The singular boundary is \emph{provably one degree less differentiable} than the fluid solution,
		and the same is true for our rough time functions. Moreover, the tori (such as $\twoargmumuxtorus{-\timefunction}{-\muxmulevelsetvalue}$)
		are two degrees less differentiable.
		To handle this fundamental difficulty, 
		we work with multiple coordinate systems, as we describe in the next point.
	\item (\textbf{Three distinct coordinate systems}).
		There is no known approach to studying multi-dimensional shock formation 
		that relies only on commuting the equations
		with the Cartesian coordinate partial derivative vectorfields.
		Because of the roughness of the time functions $\timefunctionarg{\muxmulevelsetvalue}$, 
		we cannot close our top-order energy estimates using only commutation vectorfields 
		that are adapted\footnote{However, for some crucial low-derivative-level estimates, 
		we \emph{do} use commutators that are adapted to the level sets of $\timefunctionarg{\muxmulevelsetvalue}$.} 
		to their level sets; 
		such an approach would lead to the loss of a derivative.
		Hence, as in other works on shocks, we construct appropriate commutators with the help of
		the eikonal function $u$. In total, our study of the flow relies on understanding the behavior of the
		fluid in three coordinate systems as well as carefully controlling the relationships 
		-- which degenerate near the singular boundary --
		between the coordinate systems:
			\begin{enumerate}
			\item The standard Cartesian coordinates $(t,x^1,x^2,x^3)$, relative to which the 
				compressible Euler equations are initially formulated
				and relative to which the singularity is visible.
			\item The geometric coordinates $(t,u,x^2,x^3)$ (where $u$ is the eikonal function), 
			for constructing suitable commutator
			and multiplier vectorfields that are adapted to the singularity and have \emph{sufficient regularity},
			and relative to which the solution remains rather smooth; 
			see Sect.\,\eqref{SSS:GEOMETRICENERGYESTIMATESONTHEROUGHFOLIATIONS}.
			\item The rough adapted coordinates $(\timefunctionarg{\muxmulevelsetvalue},u,x^2,x^3)$, where
			the level sets of $\timefunctionarg{\muxmulevelsetvalue}$
			are \emph{precisely adapted to the shape} of the singular boundary.
		\end{enumerate}
	\item (\textbf{Elliptic-hyperbolic estimates for the vorticity and entropy}).
		It is difficult to control the top-order derivatives of the vorticity and entropy
		on regions that are adapted to the shape of the singular boundary.
		In particular, the analytic framework we use, which is based on the
		formulation of compressible Euler flow provided by Theorem~\ref{T:GEOMETRICWAVETRANSPORTSYSTEM},
		requires that we prove that the vorticity and entropy gradient are exactly as differentiable
		as the velocity and density, all the way up to the singular boundary.
		To this end, we develop an upgraded version of the ``elliptic-hyperbolic'' integral identities 
		from \cite{lAjS2020}. The upgraded version, which we implement with the help of a new \emph{characteristic current},
		allows us to handle severe degeneracies that arise near the crease.
		The degeneracies are tied to the fact that the level sets of our time functions
		$\timefunctionarg{\muxmulevelsetvalue}$ \emph{become asymptotically null} the near crease
		for $\muxmulevelsetvalue$ small and positive; see Figs.\,\ref{F:MAXDEVELOPMENTINCARTESIAN} 
		and 
		\ref{F:CARTESIANROUGHFOLIATIONFUTUREOFCREASE}.
		We refer to Sect.\,\ref{S:ELLIPTICHYPERBOLICIDENTITIES} for our derivation of the elliptic-hyperbolic identities.
	\item	(\textbf{Transversal convexity}).
	As we have mentioned, our analysis relies on an open data-assumption that we call 
	\emph{acoustical transversal convexity}
	(which we refer to as transversal convexity for short);
	we refer to \eqref{E:DATATASSUMPTIONMUTRANSVERSALCONVEXITY} for the precise, technical 
	data estimates that capture transversal convexity.
	In studying the flow, 
	we are able to propagate the transversal convexity all the way up to the singular boundary,
	which in particular ensures that the level sets $\lbrace \upmu = 0 \rbrace$ 
	and $\lbrace \muX \upmu = 0 \rbrace$
	on RHS~\eqref{E:INTROCREASE} (which defines the crease) intersect transversally,
	Our assumption of transversal convexity is \underline{weaker} than the assumption that the graph
	of $\lbrace \upmu = 0 \rbrace$, viewed as a subset of geometric coordinate space, 
	is strictly convex, and our weaker assumption is what allows us, for example,
	to handle perturbations of symmetric solutions. 
	Equivalently, under our approach, transversal convexity allows 
	us to understand the structure of the crease and the singular boundary for singularities 
	that are more general than the fully non-degenerate ones featured in Fig.\,\ref{F:STRICTLYCONVEX},
	including the perturbations of simple isentropic plane-waves that we treat in detail in our main theorems,
	as is depicted in Fig.\,\ref{F:MAXDEVELOPMENTINCARTESIAN}.
	
	Our transversal convexity assumption is close to optimal in the sense that without it,
	\emph{the qualitative character of the singular boundary can dramatically change},
	e.g., without transversal convexity,
	the crease could fail to have the structure of a $2D$ submanifold,
	and then one could not even properly set up the shock development problem 
	using the known approach \cite{dC2019}; see the next point.
	The qualitative change of the structure of the crease in the absence of transversal convexity
	can readily be seen in the context of the simple isentropic plane-symmetric solutions
	that we study in Appendix~\ref{A:PS}. More precisely,
	using the explicit formulas provided by Cor.\,\ref{AC:PSEXPLICITEXPRESSIONSFORSOLUTION},
	one can construct simple isentropic plane-symmetric fluid initial data 
	such that transversal convexity fails and such that 
	the crease $\lbrace \upmu = 0 \rbrace \cap \lbrace \muX \upmu = 0 \rbrace$,
	viewed as subset of $1+3$-dimensional geometric coordinate space,
	fails to have the structure of a $2D$ submanifold. For example, 
	one can any consider smooth initial data such that the data-term
	$
		\frac{d}{du}
	\antiderivativePSLmusourcetermfunction[\dataRRiemannPS(u)]
	$
	on RHS~\eqref{AE:SPEEDTIMESMUPLANESYMMETRYWITHEXPLICITSOURCE} 
	achieves its negative minimum along a closed interval $I$ of $u$-values.
	It is then easy to see, 
	using \eqref{AE:SPEEDTIMESMUPLANESYMMETRYWITHEXPLICITSOURCE},
	\eqref{AE:MUXMUPLANESYMMETRYWITHEXPLICITSOURCE},
	and the fact that $\muX = \geop{u}$ for simple isentropic plane-symmetric solutions,
	that $\muX \upmu = 0$ at all points in $\Sigma_{T_{\textnormal{Shock}}}$
	where $\upmu$ vanishes, where $T_{\textnormal{Shock}}$ is the Cartesian time of first blowup.
	In particular, the level sets $\lbrace \upmu = 0 \rbrace$ 
	and $\lbrace \muX \upmu = 0 \rbrace$ \emph{coincide} on the interval of $u$ values $I$ within
	$\Sigma_{T_{\textnormal{Shock}}}$, signifying a dramatic failure of
	transversal convexity. Moreover, using the characterization~\eqref{E:INTROCREASE}
	of the crease, we see that the crease,
	viewed as a subset of $1+3$-dimensional geometric coordinate space,
	is the following set:
	$
	\left\lbrace (t,u,x^2,x^3) \ | \ t = T_{\textnormal{Shock}}, u \in I, (x^2,x^3) \in \mathbb{T}^2 \right\rbrace$,
	which is a $3D$ submanifold-with-boundary (in particular, it is not a $2D$ submanifold).
	\item (\textbf{Homeomorphism property of the change of variables map}).
	Finally, we highlight that  -- thanks to the transversal convexity -- we are able to prove
	that the change of variables map 
	$\Upsilon(t,u,x^2,x^3) \eqdef (t,x^1,x^2,x^3)$
	from geometric to Cartesian coordinates
	is a homeomorphism all the way up to the singular boundary, 
	and it is a diffeomorphism away from the singular boundary;
	see Prop.\,\ref{P:CHOVGEOMETRICTOCARTESIANISINJECTIVEONREGIONWECAREABOUT}.
	In particular, $\Upsilon$ is a bijection
	on the region $\MInteresting$ depicted in Fig.\,\ref{F:INTERESTINGREGIONMAINRESULTS},
	and \emph{the level sets of $u$ never actually intersect} 
	$\MInteresting$ (even along the singular boundary portion $\mathcal{B}^{[0,\muxmulevelsetvalue_0]}$!),
	despite the fact that the density of the level sets of $u$ in Cartesian coordinate space becomes infinite along the singular boundary. 
	Under transversal convexity, 
	the map $u \rightarrow - u^3$ provides a crude scalar caricature of 
	the degenerate behavior of the $x^1$ Cartesian coordinate 
	along a constant Cartesian time slice $\Sigma_t$ that happens to be tangent to the crease.
	\textbf{This not just a mathematical curiosity}; the bijective property of $\Upsilon$ up to and including the singular boundary
	is crucial for formulating the shock development problem,
	and this property can dramatically fail without transversal convexity.
	As in the previous point, this difficulty
	can readily be seen in the context of the simple isentropic plane-symmetric solutions
	that we study in Appendix~\ref{A:PS}. For example, in the example discussed in the previous point,
	using that $\geop{u} x^1 = - \upmu$ in plane symmetry, one finds that
	$x^1 = x^1(t,u,x^2,x^3)$ is constant along the entire crease
	$
	\left\lbrace (t,u,x^2,x^3) \ | \ t = T_{\textnormal{Shock}}, u \in I, (x^2,x^3) \in \mathbb{T}^2 \right\rbrace
	$,
	thereby exhibiting a dramatic failure of the injectivity of $\Upsilon$.
	In the Cartesian coordinate space picture in $1+1$ dimensions (i.e., ignoring the $x^2$ and $x^3$ directions), 
	the shock formation in this example corresponds to a continuum of characteristic curves 
	all intersecting in a single point in Cartesian coordinate space.
	This is in stark contrast to the behavior of the characteristic hypersurfaces 
	on the singular boundary in the transversally convex regime,
	as we show in Prop.\,\ref{P:DESCRIPTIONOFSINGULARBOUNDARYINCARTESIANSPACE}
	and Fig.\,\ref{F:NONUNIQUENESSOFINTEGRALCURVESOFL}.
\end{itemize}

\subsection{Works building up toward Theorem~\ref{T:ABBREVIATEDSTATEMENTOFMAINRESULTS}}	
\label{SS:WORKSBUILDINGTOWARDSMAINTHEOREM}
In this section, we describe some key prior works that developed some of the technology
we use here. While we have already mentioned some of them, 
here we provide additional details.
We refer to Sect.\,\ref{SS:BACKGROUNDONSHOCKS} for a more extensive (though far from comprehensive)
discussion of the history of works on shock formation.

\subsubsection{Christodoulou's framework from \cite{dC2007}}
\label{SSS:CHRISTODOULOUIRROTATIONALFRAMEWORK}
Our proof of Theorem~\ref{T:ABBREVIATEDSTATEMENTOFMAINRESULTS} relies on the framework of nonlinear geometric optics
developed by Christodoulou \cite{dC2007} in his groundbreaking proof of stable shock formation
for the irrotational and isentropic relativistic Euler equations in three spatial dimensions.
In this setting, the equations reduce to a scalar quasilinear wave equation for a potential function
$\Phi$:
\begin{align} \label{E:QLW}
	(\gfour^{-1})^{\alpha \beta}(\pmb{\partial} \Phi) \partial_{\alpha} \partial_{\beta} \Phi 
	& = 0.
\end{align}
In \eqref{E:QLW}, $\gfour_{\alpha \beta}$ is the \emph{acoustical metric},
a Lorentzian metric whose rectangular component functions $\gfour_{\alpha \beta}$
are nonlinear functions of $\pmb{\partial} \Phi \eqdef (\partial_t \Phi, \partial_1 \Phi, \partial_2 \Phi, \partial_3 \Phi)$ 
depending on the equation of state:
$\gfour_{\alpha \beta} = \gfour_{\alpha \beta}(\pmb{\partial} \Phi)$.
To implement nonlinear geometric optics, Christodoulou constructed an
\emph{eikonal function} $u$, that is a solution to the \emph{eikonal equation}
\begin{align} \label{E:INTROEIKONAL}
	(\gfour^{-1})^{\alpha \beta}(\pmb{\partial} \Phi) \partial_{\alpha} u \partial_{\beta} u 
	& = 0.
\end{align}
The level sets of $u$, denoted by $\nullhyparg{u}$, are characteristic for equation \eqref{E:QLW}.
Notice that the principal coefficients on LHS~\eqref{E:INTROEIKONAL} depend on $\pmb{\partial} \Phi$
and thus the evolution of $u$ is coupled to that of $\Phi$, signifying dynamic nature of the geometry 
and the quasilinear nature of the flow.
In \cite{dC2007}, Christodoulou used $u$ to construct a system of \emph{geometric coordinates}
$(t,u,\vartheta^1,\vartheta^2)$ on spacetime such that the solution remains rather smooth relative to these
coordinates. The formation of the singularity can be recovered as a degeneracy between the 
geometric coordinates and the standard ones, $(t,x^1,x^2,x^3)$.
The degeneracy is signified by the vanishing of the \emph{inverse foliation density} $\upmu$, defined by:
\begin{align} \label{E:INTROINVERSEFOLIATION}
	\upmu 
	& 
	\eqdef - \frac{1}{(\gfour^{-1})^{\alpha \beta}(\pmb{\partial} \Phi) \partial_{\alpha} t \partial_{\beta} u}.
\end{align}
In the region of classical existence, one has $\upmu > 0$,
and when $\upmu \to 0$, the density of the $\nullhyparg{u}$ in Cartesian coordinate space
becomes infinite, signifying the ``piling up'' of the characteristics.
In Figs.\,\ref{F:INFININTEDENSITYOFCHARACTERISTICSONSINGULARBOUNDARY} and \ref{F:INTERESTINGREGIONMAINRESULTSCARTESIAN},
we show three distinct characteristic hypersurfaces piling up along the singular boundary,
along which $\upmu$ vanishes.
The hard part of the proof is to derive energy estimates in regions where $\upmu$ is small. 
We refer readers to Sect.\,\ref{SS:BACKGROUNDONSHOCKS} for further discussion on Christodoulou's framework.

\subsubsection{Wave equations beyond fluid mechanics}
The wave equations studied by Christodoulou in \cite{dC2007}
were Euler--Lagrange equations that were invariant with respect to the Poincar\'{e} group.
In \cite{jS2016b}, we extended the results of \cite{dC2007}
to apply to all wave equations of type \eqref{E:INTROEIKONAL}
and of type: 
\begin{align} \label{E:INTROCOVARIANTWAVE}
	\square_{\gfour(\Psi)} \Psi 
	& = \nullform(\pmb{\partial} \Psi, \pmb{\partial} \Psi)
\end{align}
that fail to satisfy Klainerman's null condition \cite{sK1984},
where $\square_{\gfour(\Psi)}$ is the covariant wave operator of 
$\gfour(\Psi)$ (see Def.\,\ref{D:COVWAVEOP}) 
and $\nullform$ is a null form relative to $\gfour$ (see Def.\,\ref{D:STANDARDNULLFORMS}).

\subsubsection{Nearly simple plane-symmetric waves}
In \cite{jSgHjLwW2016}, the second author and his collaborators extended the methods of \cite{dC2007} and \cite{jS2016b} 
to prove stable shock formation
for a large class of quasilinear wave equations on the spacetime $\mathbb{R} \times \Sigma$,
where $\Sigma = \mathbb{R} \times \mathbb{T}$ was the two-dimensional spatial manifold.
The initial data we treated were analogs of the data from Theorem~\ref{T:ABBREVIATEDSTATEMENTOFMAINRESULTS}.
More precisely, the data were (asymmetric) perturbations of simple plane-symmetric waves,
which are solutions that depend only on $(t,x^1) \in \mathbb{R} \times \mathbb{R}$
and which feature a wave moving only in one direction\footnote{In plane symmetry, one can study
the flow by constructing Riemann invariants, in which case simple plane-symmetric waves would feature only a single non-zero Riemann invariant. See Appendix~\ref{A:PS}, in which we use Riemann invariants to construct the background
solutions whose perturbations we study in our main results.} 
(say to the right).
As in the relativistic case, for irrotational and isentropic solutions to the $2D$ compressible Euler equations,
the dynamics reduces to a quasilinear wave equation of type \eqref{E:INTROEIKONAL}.
Hence, as a special case, the results of
\cite{jSgHjLwW2016}
yielded stable shock formation 
for nearly simple and isentropic plane-symmetric solutions
to the $2D$ irrotational and isentropic compressible Euler equations.

\subsubsection{A new formulation of the flow and stable shock formation in the presence of vorticity}
In \cite{jLjS2020a}, the second author and Luk 
developed a new formulation of barotropic\footnote{Barotropic equations of state are such that the pressure
$p$ can be expressed as a function of the density $\varrho$, i.e., with no dependence on $\Ent$.} 
compressible Euler flow with remarkable 
regularity and null structures, which in many regimes, allows one to study the flow as if it was
a perturbation of the quasilinear wave equation \eqref{E:INTROCOVARIANTWAVE}.
In \cite{jS2019c}, the second author 
derived a similar new formulation for all equations of state
in which the pressure is a function of the density and entropy,
thus allowing one to incorporate thermodynamic effects into the framework of \cite{dC2007}.
In particular, the equations of \cite{jS2019c} include a system of
transport-div-curl equations for the vorticity and entropy, which allows one to propagate a gain of one
derivative for these quantities relative to standard estimates. 
In this article, we use the equations of \cite{jS2019c} to prove our main results,
and the gain in regularity is crucial for our approach.
In Theorem~\ref{T:GEOMETRICWAVETRANSPORTSYSTEM}, we recall the new formulation of the flow derived in \cite{jS2019c}.
In \cite{mDjS2019}, we derived a similar new formulation for the relativistic Euler equations.

In \cite{jLjS2018}, the second author and Luk used the equations of
\cite{jLjS2020a} and the technology of \cites{dC2007,jS2016b}
to prove the first stable shock formation result for the $2D$ compressible Euler equations.
The authors treated open sets of initial data with vorticity that are close to the 
data of an irrotational simple plane-wave solution.
The main theorem yielded the full structure of the set of blowup-points
within the constant-time hypersurface $\Sigma_{T_{\textnormal{Shock}}}$ of first blowup;
in the context of Fig.\,\ref{F:MAXDEVELOPMENTINCARTESIAN},
the authors understood the structure of
$\Sigma_{T_{\textnormal{Shock}}} \cap \mathcal{B}$.

In \cite{jLjS2021}, the second author and Luk used the equations of \cite{jS2019c}
to extend the results of \cite{jLjS2020a}
to the $3D$ case in the presence of vorticity and entropy.
The proof was much more difficult than the $2D$ case because
our regularity theory for the vorticity and entropy relied
on elliptic estimates across space, 
which are difficult to derive near the shock.
In particular, the authors used the elliptic estimates to handle the vorticity stretching term in the equations,
which, as is well-known, vanishes for $2D$ barotropic solutions.
These elliptic estimates relied on the
transport-div-curl equations for the vorticity and entropy that we derived in \cite{jS2019c},
and our approach yielded elliptic estimates only along complete, flat hypersurfaces of constant time.
As in the $2D$ case, this approach yielded the full structure of the set of blowup-points
within the constant-time hypersurface $\Sigma_{T_{\textnormal{Shock}}}$ of first blowup,
i.e., 
in the context of Fig.\,\ref{F:MAXDEVELOPMENTINCARTESIAN},
the singular set portion $\Sigma_{T_{\textnormal{Shock}}} \cap \mathcal{B}$.
As we have already mentioned, in order to study the flow beyond the 
Cartesian-flat hypersurface $\Sigma_{T_{\textnormal{Shock}}}$,
one needs additional ingredients,
including the identities that we describe 
in Sect.\,\ref{SSS:REMARKABLELOCALIZED}.

\subsubsection{Remarkable localized integral identities}
\label{SSS:REMARKABLELOCALIZED}
In \cite{lAjS2020}, we derived new localized integral identities for $3D$ compressible Euler flow
that allow us to extend the elliptic estimates for the vorticity and entropy to yield
estimates on \emph{arbitrary} spacetime regions that are bounded by 
spacelike or null hypersurfaces.
We use this crucial ingredient in the present paper to derive 
top-order vorticity and entropy estimates up to the singular boundary.

\subsection{Additional history and results tied to shocks and singularities}	
\label{SS:BACKGROUNDONSHOCKS}
A fundamental issue, discovered by Riemann \cite{bR1860} in the context of one spatial dimension, is that initially smooth 
compressible Euler solutions typically develop shock singularities in finite time. 
Roughly, shocks are singularities such that $\rho$, $v$, $\Ent$ remain bounded
but some first derivative of $\rho$ and $v$ blows up. Riemann's analysis relied on his discovery of Riemann invariants for isentropic solutions in $1D$. Relative to Riemann invariants, the equations reduce to a quasilinear system for two transport equations with distinct characteristic directions; see Appendix~\ref{A:PS}. The proof of the blowup of the solution's first derivatives
then follows from differentiating the equations to obtain a Riccati-type structure, much like in the model case of Burgers' equation $\partial_t \Psi + \Psi \partial_x \Psi = 0$.

\subsubsection{The advanced state of the $1D$ theory and the key difference with multi-dimensions}
\label{SSS:ADVANCEDTHEORYIN1D}
In $1D$, for large class of hyperbolic quasilinear systems (such as strictly hyperbolic systems),
there is an advanced theory capable of describing the global behavior of solutions,
including the formation of shocks and the subsequent interactions of the shock waves.
We refer readers to the compendium \cite{cD2010} for a history of the subject
and a comprehensive introduction to the main techniques.

A principal reason for the advanced state of the $1D$ theory is that the
equations are well-posed for initial data in appropriate bounded variation (BV) spaces.
The state of affairs is dramatically different in multi-dimensions;
Rauch's fundamental work \cite{jR1986} showed that in multi-dimensions,
quasilinear hyperbolic systems are typically \emph{ill-posed} for data in BV spaces.
In fact, the only known well-posedness results in multi-dimensions are for initial
data in $L^2$-type Sobolev spaces. For this reason, in multi-dimensions,
one is forced to derive energy estimates, which can be incredibly difficult in regions
containing singularities; this is the main reason why the theory 
of multi-dimensional shock waves is so much less developed compared to the $1D$ case.

\subsubsection{The first blowup-result in multi-dimensions without symmetry: proof by contradiction}
\label{SSS:SIDERIANMETHODS}
Providing a \emph{constructive} proof of shock formation in higher dimensions turns out to be a very hard problem. 
In a nutshell, the reason is that away from $1D$, it seems necessary to carefully track the evolution
of characteristic hypersurfaces, which is much more difficult compared to the $1D$ case.
The characteristic geometry (e.g., the eikonal function $u$ in the context of the present article)
is much more difficult to construct and control, and, crucially, all proofs of even local well-posedness 
rely on energy estimates in $L^2$-based Sobolev spaces, which are difficult to derive near singularities.
In \cite{tS1985}, Sideris proved an influential, \emph{non-constructive} 
stable blowup-result in $3D$, the first one for multi-dimensional
compressible Euler flow. He assumed that the equation of state is 
barotropic\footnote{Barotropic equations of state are such that $p = p(\varrho)$.}
and that it satisfies a convexity assumption. 
His proof applied to a large set of data, but it did not reveal the nature of the blowup;
his arguments relied on virial-type identities, 
and he showed blowup through a contradiction-argument.

\subsubsection{Alinhac's constructive proof of the formation of isolated singularities}
\label{SSS:ALINHACSBLOWUPRESULTS}
The first constructive results on shock formation in multidimensions were by Alinhac 
\cites{sA1999a,sA1999b,sA2002}, who used nonlinear geometric optics (i.e., he constructed characteristic surfaces) 
and Nash--Moser estimates to prove stable shock formation in $2D$
for a class of quasilinear wave equations
of the form:
\begin{align} \label{E:ALINHACSWAVEEQUATIONS}
	(\gfour^{-1})^{\alpha \beta}(\pmb{\partial} \Phi) \partial_{\alpha} \partial_{\beta} \Phi 
	& = 0,
\end{align}
whenever they fail to satisfy Klainerman's null condition \cite{sK1984}.
He used the Nash--Moser estimates to avoid derivative loss in his control of the regularity of the characteristic surfaces.
His proof applied to open sets of \emph{non-degenerate} initial data,
which he showed lead to the formation of \emph{fully non-degenerate} singularities
in which $\pmb{\partial}^2 \Phi$ blows up while $\Phi$ and $\pmb{\partial} \Phi$ remain bounded.
Roughly, Alinhac's framework
applied to initial data such the singular boundary has the strictly convex
structure depicted in Fig.\,\ref{F:STRICTLYCONVEX},
and it allowed him to follow the solution up to the lowest point
on the singular boundary, but not further.
In Fig.\,\ref{F:STRICTLYCONVEX},
we denote this lowest point by ``$b_*$.''
Alinhac's framework has also been applied to other wave equations;
see, for example, \cites{bDiWhY2012,bDiWhY2015,dBiWyH2016}.

\subsubsection{Christodoulou's breakthrough on irrotational, isentropic shock formation}
\label{SSS:CHRISTODOULOUBREAKTHROUGHONSHOCKFORMATION}
In \cite{dC2007}, Christodoulou proved a stable shock formation result
for open sets of irrotational and isentropic initial solutions to the relativistic
Euler equations in three spatial dimensions. 
In this context, one can study the flow with the help of a potential function $\Phi$,
and relativistic Euler flow reduces to a quasilinear wave equation of type
\eqref{E:ALINHACSWAVEEQUATIONS}. Aside from a single exceptional equation of state corresponding
to the graph of a timelike minimal surface in Minkowski space, all wave equations
of irroational and isentropic relativistic fluid mechanics fail to satisfy the null condition,
i.e., the basic mechanism driving shock formation is present all equations but one.
The data that Christodoulou treated were compact perturbations of non-vacuum constant state data,
and his shock formation results revealed the instability of these states under irrotational and isentropic perturbations.
To study the solution, Christodoulou
used a refined version of nonlinear geometric optics,
based on techniques that he co-developed with Klainerman in their proof of the stability of the Minkowski space \cite{dC1993}.
Specifically, Christodoulou constructed an eikonal function $u$, i.e., a solution to the eikonal equation
$(\gfour^{-1})^{\alpha \beta}(\pmb{\partial} \Phi) \partial_{\alpha} u \partial_{\beta} u = 0$,
and he showed that the wave equation solution's first Cartesian derivatives, 
$\pmb{\partial} \Phi$, remain quite smooth relative to the geometric
coordinates. The level sets $\nullhyparg{u}$ of $u$ are characteristic hypersurfaces for the wave equation.
The surfaces depend on the solution itself, reflecting the quasilinear nature of the flow.

As we mentioned already in Sect.\,\ref{SSS:CHRISTODOULOUIRROTATIONALFRAMEWORK},
as in the present paper, in Christodoulou's framework, 
the formation of the shock corresponds to the vanishing of the inverse foliation density $\upmu$.
In the region of classical existence, one has $\upmu > 0$,
and when $\upmu \to 0$, the density of the $\nullhyparg{u}$ in Cartesian coordinate space
becomes infinite, signifying the piling up of the characteristics.
The blowup of $\pmb{\partial}^2 \Phi$ is a consequence of the degeneracy 
between the Cartesian and a system\footnote{In \cite{dC2007}, the ``angular'' coordinate functions $\vartheta^1$
and $\vartheta^2$ were constructed so as to be constant along the integral curves of the null generator $\Lunit$.
In contrast, in the present paper, in the role of the ``angular'' coordinate functions, we use the standard
Cartesian coordinates $x^2$ and $x^3$, which are not typically constant along the integral curves of $\Lunit$.
This minor difference turns out to have no substantial effect on the analysis. 
\label{FN:OURANGULARCOORDINATESAREDIFFERENTTHATCHRISTODOULOUS}} 
of ``geometric coordinates'' $(t,u,\vartheta^1,\vartheta^2)$,
where we schematically represent the degeneracy caused by the vanishing of $\upmu$ 
as follows: $\pmb{\partial} \sim \frac{1}{\upmu} \frac{\partial}{\partial u}$;
the blowup of $\pmb{\partial}^2 \Phi$ 
then follows from this relation, 
from proving that $\upmu \to 0$ in finite time, and
from proving a lower of the form 
$|\frac{\partial}{\partial u} \partial \Phi| \gtrsim 1$.

Importantly, Christodoulou's proof relied on geometric energy estimates relative to a foliation of spacetime
by the $\nullhyparg{u}$, 
rather than Nash--Moser estimates.
To control the geometry without derivative loss, he relied on techniques and renormalized quantities
that have their roots in \cite{dCsK1993,sKiR2003}.
Crucially, Christodoulou's geometric approach allowed him to prove shock formation for 
a larger class of singularities than the fully non-degenerate ones treated by Alinhac.
He was able to show that blowup occurs at one or more points
even for solutions whose singular boundaries do not have to enjoy the strict convexity displayed in
Fig.\,\ref{F:STRICTLYCONVEX}, though in the absence of strict convexity, the full structure of the maximal
development was not revealed. Importantly, the work \cite{dC2007}
also yielded a sharp conditional global existence result, 
which showed that irrotational, isentropic near-constant-state solutions are global
unless a shock forms.

\subsubsection{Shock-formation works that built upon Christodoulou's framework}
\label{SSS:SHOCKFORMATIONBUILDUPONCHRSITODOULOU}
Many authors have used Christodoulou's approach to prove stable shock formation
for various quasilinear hyperbolic PDEs. 
For example, there are stable shock formation results for:
\begin{itemize}
	\item A larger class of wave equations \cites{sMpY2017,jS2016b}.
	\item Solution regimes that are different than the small, compactly supported data regime: 
		\cites{gHsKjSwW2016,sMpY2017}.
	\item Solutions that exist classically precisely on a past-infinite half-slab \cite{sM2018}.
	\item Various systems involving multiple speeds of propagation, 
		some with symmetry \cites{dCdRP2016,xAhCsY2020},
		and some without \cites{jS2018b,jS2019a}.
\end{itemize}

\subsubsection{Shock formation in the presence of vorticity and entropy}
\label{SSS:SHOCKFORMATIONWITHVORTICITYANDENTROPY}
The aforementioned paper \cite{jLjS2018} was to first to prove stable shock formation for open
sets of $2D$ compressible Euler solutions with vorticity under an arbitrary\footnote{As in our main results,
the Chaplygin gas equation of state is exceptional and is not known to lead to shock formation.} 
barotropic equation of state. The initial data were not required to be fully non-degenerate,
i.e., transversal convexity did not play a role
in the proof of blowup.
The main theorem followed the solution to the constant-time hypersurface of first blowup,
and it gave a complete description of what blows up and what does not, as well as a precise description of 
the set of singular points at the time of first blowup,
i.e., in the context of Fig.\,\ref{F:MAXDEVELOPMENTINCARTESIAN},
a complete description of $\Sigma_{T_{\textnormal{Shock}}} \cap \mathcal{B}$. 
Since the geometric setup relied on the methods of \cite{dC2007}, the
singularities were allowed to be more general than the 
fully non-degenerate ones described in Sect.\,\ref{SSS:ALINHACSBLOWUPRESULTS}.

The recent work \cite{jLjS2021} extended \cite{jLjS2018} 
to the case of the $3D$ compressible Euler equations under an arbitrary (non-Chaplygin gas)
equation of state with vorticity and entropy. 
As in \cite{jLjS2018}, stable shock formation was proved in \cite{jLjS2021} without transversal convexity,
though in a sub-regime of solutions with transversal convexity, additional information
on the H\"{o}lder regularity of the solution with respect to the Cartesian coordinates was derived.\footnote{With
transversal convexity, the solutions was shown to enjoy $C^{1/3}$-H\"{o}lder regularity with
respect to the Cartesian coordinates up to the singularity, while not enjoying
$C^{(1/3})^+$-H\"{o}lder regularity. Providing a more detailed treatment of this regime 
was inspired by H\"{o}lder regularity results derived  
in \cites{tBsSvV2019a,tBsSvV2019b,tBsSvV2020} for some fully non-degenerate solutions.}

\subsubsection{A new approach to proving the formation of fully-non-degenerate shock singularities via self-similarity}
\label{SSS:NEWSHOCKAPPROACHBASEDONSELFSIMILARITY}
In \cites{tBsSvV2019b,tBsSvV2020}, the authors provided a philosophically interesting new approach
for proving the formation of shock singularities in $3D$ compressible Euler solutions
under an adiabatic equation of state without symmetry, 
and with vorticity (and also entropy in \cite{tBsSvV2020}).
The approach allows one to follow the solution to the time of first blowup,
and the singularities produced are isolated within the constant-time hypersurfaces
of first blowup. That is, in the context of
Fig.\,\ref{F:STRICTLYCONVEX}, the approach allows one
to follow the solution up to the point $b_*$.
 Such singularities are analogs of the non-degenerate singularities that
Alinhac studied \cites{sA1999a,sA1999b,sA2002}
in the case of quasilinear wave equations.
The framework of \cites{tBsSvV2019b,tBsSvV2020} 
relied on modulation parameters to show that
for open sets of smooth data, a singularity develops
in finite time, and it is a perturbation of a self-similar Burgers'-type shock.

The aforementioned $3D$ works were preceded by the work \cite{tBsSvV2019a} in $2D$ azimuthal symmetry with vorticity.
In the recent work \cite{tBsI2020}, in the same symmetry class,
the authors constructed shock-forming solutions whose cusp-like spatial behavior (with respect to the standard coordinates)
is \emph{non-generic}; such solutions are unstable. 
In the language of the present paper, 
such solutions do not exhibit the quantitative transversal convexity \eqref{E:MUTRANSVERSALCONVEXITY}
that we use in proving our main results.

\subsubsection{Self-similar blowup for non-hyperbolic PDEs}
\label{SSS:SELFSIMILARBLOWUPFORNONHYPERBOLICPDES}
Notably, there are \emph{non-hyperbolic} PDES with solutions that exhibit self-similar blowup 
modeled on a self-similar Burgers shock. Examples include the Euler--Poison equations
\cite{qYlZ2021},
Burgers' equation with transverse viscosity \cite{cCteGnM2018}, 
the Burgers--Hilbert equations \cite{rY2020}, the fractal Burgers equation \cite{krCrcMVgP2021}, 
and various dispersive or dissipative perturbations of the Burgers equation \cite{sjOfP2021};
see also \cites{cCteGsInM2018,cCteGnM2019}.
It would be interesting to investigate the extent to which the set of singular points
can be understood, i.e., to try to derive results for these equations that are analogs of the results of the present paper.

\subsubsection{Implosion singularities}
The recent breakthrough work \cite{fMpRiRjS2019b} in spherical symmetry showed that for
adiabatic equations of state $p = \varrho^{\upgamma}$ with $\upgamma > 1$, 
there exist $C^{\infty}$ initial data such that the corresponding
solution's density and velocity blow up at the center of symmetry in finite time. 
In fact, there are infinitely many such singularity-forming solutions,
collectively exhibiting a discrete sequence of blowup-rates.
These ``implosion singularities'' are much more severe than shocks.
The methods of \cite{fMpRiRjS2019b} suggest that the implosion might enjoy co-dimension stability
under perturbations of the initial data without symmetry, 
though possibly not full stability corresponding to 
open sets of data.

\subsubsection{Inviscid limits}
\label{SSS:INVISCIDLIMITS}
A physically important and mathematically interesting problem is to study
the relationship between the formation of shock singularities in classical solutions to hyperbolic PDEs 
and the behavior of classical solutions when a small amount of viscosity is added to the equations.
Of particular interest is to understand the zero viscosity limiting behavior of classical solutions.
The recent paper \cite{sCgC2022}, which concerns Burgers' equation in $1D$ and its viscous analog,
is the only work to date that yields information about the zero viscosity limit 
of classical solutions all the way up until the time of first singularity formation for the inviscid solution.
The main results show that the $1D$ viscous Burgers' equation 
solution can be decomposed into a singular piece and a smoother piece,
and that the viscous solution converges to the singular piece in $\| \cdot \|_{L^{\infty}}$ as the viscosity vanishes
(the $L^{\infty}$ norm is shown to be bounded from above by the viscosity parameter to a positive power).
Here, the $\| \cdot \|_{L^{\infty}}$ norm is over the entire region of classical existence of the inviscid solution,
i.e., the $L^{\infty}$-convergence holds all the way until the time of first blowup for the solution to the inviscid Burgers' equation. The results of \cite{sCgC2022} 
apply for non-degenerate initial data, where the notion of non-degeneracy in \cite{sCgC2022} 
is essentially equivalent to the transversal convexity assumption satisfied by the solutions featured in our main theorems.
An important open problem is to extend the results of \cite{sCgC2022} to the compressible Euler equations
(where the corresponding viscous model is the compressible Navier--Stokes equations) 
and to multiple spatial dimensions.

\subsubsection{Progress on the shock development problem}
\label{SSS:PROGRESSONSHOCKDEVELOPMENT}
In Majda's celebrated works \cites{aM1981,aM1983a}, he proved linear stability and local well-posedness  
results for a class of weak solutions to the $3D$ compressible Euler equations arising from a set of
discontinuous initial data.
More precisely, he assumed that the data on $\mathbb{R}^3$ were piecewise smooth and jumped across
a smooth two-dimensional hypersurface such that the jumps are consistent with  
well-known Rankine--Hugoniot jump conditions. He also assumed the data satisfy a certain ``stability condition.'' 
His main result was the construction of a local, unique weak solution in a subset of $\mathbb{R}^{1+3}$,
and a corresponding three-dimensional shock hypersurface, 
across which the solution jumps in accordance with the Rankine--Hugoniot jump conditions.
This is known as the \emph{shock front problem}. 

An important problem, distinct from the shock front problem,
is to determine whether/how Majda's discontinuous initial conditions
can develop from an initially smooth solution. That is, one would like to describe the \emph{transition} of compressible Euler solutions
from being smooth, to developing a ``first singularity'' -- which in the language of the present paper is the crease -- 
and finally to becoming a weak solution that develops a shock hypersurface (emanating from the crease),
across which the solution jumps. This is known as the \emph{shock development problem}. 
In Fig.\,\ref{F:MAXDEVELOPMENTWITHSHOCKHYPERSURFACEINCARTESIAN}, we denote the shock hypersurface by ``$\mathcal{K}$,''
and we show its emergence from the crease $\partial_- \mathcal{B}$. We stress 
that Majda's works did not study the flow in a neighborhood of $\partial_- \mathcal{B}$, but
rather started from initial conditions on a flat hypersurface $\Sigma_{t_0}$ 
such that the co-dimension surface $\Sigma_{t_0} \cap \mathcal{K}$, across which the data jump, 
is already assumed to exist.
A crucial issue is that the jump conditions
imply that initially irrotational and isentropic smooth solutions will develop dynamic entropy and non-trivial vorticity as the solution
jumps across the shock hypersurface. This means that it is not possible to solve the true shock development problem in the class of irrotational and isentropic solutions. The full problem in $3D$ without symmetry assumptions remains open, but there has been 
inspiring progress in recent years:
\begin{itemize}
	\item In \cite{dCaL2016}, Christodoulou--Lisibach used Riemann invariants 
		and a pair of eikonal functions (one ``ingoing'' and the other ``outgoing'')
		to solve the shock development problem
		for spherically symmetric, barotropic solutions to the relativistic Euler equations.
		Mathematically, the problem comprised a system of quasilinear transport equations
		coupled to the free boundary problem of tracking the location of the 
		shock hypersurface. The problem was solved through an iteration scheme that relied
		on the assumption that the crease exhibits transversal convexity.
		Relative to geometric coordinates analogous to the ones that we use in our main results,
		the data are assumed to be smooth along the crease and Cauchy horizon, and this allows one to
		use Taylor expansions to approximate the expected location
		of the shock hypersurface and the behavior of the solution along it; 
		these Taylor expansions form a crucial ingredient in controlling the iterates.
	\item In his breakthrough monograph \cite{dC2019}, Christodoulou extended the methods from \cite{dCaL2016}
		and solved the ``restricted'' shock development problem	
		\emph{without symmetry assumptions} for the compressible Euler equations and relativistic Euler equations
		in an arbitrary number of spatial dimensions. The word ``restricted'' means that he studied only irrotational and isentropic solutions,
		and that he ignored the jump in entropy and vorticity across the shock hypersurface, thereby producing a weak solution to 
		a hyperbolic PDE system that approximates the real one.
		Christodoulou assumed that along the crease, the solution satisfies 
		the same kind of transversal convexity that the solutions produced by our main results enjoy. 
		In particular, this class of singularities is more general than the fully non-degenerate ones
		treated by Alinhac in Sect.\,\ref{SSS:ALINHACSBLOWUPRESULTS}.
		The main new difficulty in \cite{dC2019} is deriving energy estimates for the solution
		and the acoustic geometry. In particular, relative to a system of geometric coordinates,
		the high order energies exhibit degenerate behavior, much like the high order energies in the present paper.
	\item In \cite{hcYzL2021}, 
			Huicheng--Lu solved the shock development problem
			for a class of first-order, scalar, divergence-form hyperbolic equations in $2D$,
			starting from initial singularities that are fully non-degenerate,
			i.e., such that the crease is strictly convex,
			as in Fig.\,\ref{F:STRICTLYCONVEX}.
	\item In \cite{tBsI2020},
			Buckmaster--Drivas--Shkoller--Vicol extended the methods from \cites{tBsSvV2019a,tBsSvV2019b,tBsSvV2020}
			to solve the shock development problem
			for the $2D$ compressible Euler equations 
			in azimuthal symmetry with vorticity and entropy under the
			adiabatic equations of state $p = \frac{1}{\upgamma} \varrho^{\upgamma} \exp(\Ent)$
			for initial singularities that are fully non-degenerate
			with respect to variations in the radial and angular variable.
			In particular, \cite{tBsI2020} provides the first solution to a compressible Euler
			shock development problem with vorticity.
			As in \cite{dCaL2016}, the problem comprised transport equations coupled to the issue of tracking the
			location of the free boundary, and it was solved via an iteration scheme based on Taylor expansions
			relative to a coordinate system in which the solution is rather smooth.
			A new feature compared to \cite{dCaL2016} is that there is a third characteristic 
			direction in \cite{tBsI2020}, corresponding to the transporting of vorticity and entropy by
			the material derivative vectorfield.
\end{itemize}

We also highlight that recently, there have been other interesting works on weak solutions to the compressible
Euler equations in $1D$ without entropy.
In \cite{aL2021}, Lisibach proved local existence for the shock reflection (off of a wall) 
problem in plane-symmetry.
In \cite{aL2022}, he studied the interaction of two shocks in plane-symmetry
and proved local existence of a weak solution near the interaction point.

\subsection{Ideas behind the proof of the main results}
\label{SS:IDEASBEHINDPROOF}
In this section, we provide a more detailed overview of 
some key ideas behind the proof of our main results.

\subsubsection{Almost Riemann invariants}
\label{SSS:INTROALMOSTRIEMANN}
To study perturbations of simple isentropic plane-symmetric waves,
we find it convenient to replace the scalar functions $\LogDensity$
(where $\LogDensity$ is the logarithmic density from Def.\,\ref{D:LOGDENS})
and $v^1$ with the ``almost\footnote{For plane-symmetric isentropic solutions, 
$\RRiemann$, $\LRiemann$ are Riemann invariants, but away from symmetry, they are not.} 
Riemann invariants'' (see Def.\,\ref{D:ALMOSTRIEMANNINVARIANTS}) 
$\RRiemann$, $\LRiemann$ which carry the same information.

Here and in the rest of the paper, 
$\wavearray \eqdef (\RRiemann,\LRiemann,v^2,v^3,\Ent)$
are the fluid ``wave variables'' 
and $\gfour_{\alpha \beta} = \gfour_{\alpha \beta}(\wavearray)$ denotes 
the acoustical metric (see definition~\eqref{E:ACOUSTICALMETRIC}), 
which drives the propagation of sound waves.

\subsubsection{``Late-time'' assumptions on the data for perturbations of simple isentropic plane-waves}
\label{SSS:PERTURBATIONSOFSIMPLEWAVESANDCAUCHYSTABILITY}
We consider the ``bona fide'' initial data to be the data specified along
the flat spacelike hypersurface $\Sigma_0$ and a portion of a null hypersurface,
denoted by $\nullhyparg{- \rightu}^{[0,\frac{4}{\mathring{\updelta}_*}]}$ in Fig.\,\ref{F:LATETIMEDATAHYPERSURFACEANDNULLDATAHYPERSURFACE}
(i.e., we are studying the Cauchy problem for spacelike-characteristic data),
where the parameter $\mathring{\updelta}_* > 0$ (see \eqref{E:DELTASTARDEF}) depends on the data and
is such that the Cartesian time of first blowup is approximately
$\mathring{\updelta}_*^{-1}$. 
However, the structures we need to detect the singular boundary appear only late
in the evolution, and we therefore find it convenient to 
state our data-assumptions on a ``late'' hypersurface portion of constant rough time,
denoted by $\hypthreearg{\timefunction_0}{[- \rightu,\leftu]}{\muxmulevelsetvalue}$ 
in Fig.\,\ref{F:LATETIMEDATAHYPERSURFACEANDNULLDATAHYPERSURFACE},
that will end up being close to the singular boundary $\mathcal{B}^{[0,\muxmulevelsetvalue_0]}$.
That is, we find it convenient 
to state the assumptions on the spacelike level sets 
$
\hypthreearg{\timefunction_0}{[- \rightu,\leftu]}{\muxmulevelsetvalue}
\eqdef
\lbrace \timefunctionarg{\muxmulevelsetvalue} = \timefunction_0 \rbrace \cap \lbrace - \rightu \leq u \leq \leftu \rbrace$
of the rough time functions,
which we describe in detail in Sect.\,\ref{SSS:INTROROUGHTIMEFUNCTIONS}.
In Sect.\,\ref{S:ASSUMPTIONSONTHEDATA},
we state all of these data-assumptions.
In Appendix~\ref{A:OPENSETOFDATAEXISTS}, we use
Cauchy stability arguments to show that these assumptions
are satisfied by perturbations of the bona fide data on $\Sigma_0$
corresponding to a class of shock-forming simple isentropic plane-wave solutions.
In Appendix~\ref{A:PS}, we use standard arguments to construct these
simple isentropic plane-wave solutions.

\begin{center}
	\begin{figure}  
		\begin{overpic}[scale=.4, grid = false, tics=5]{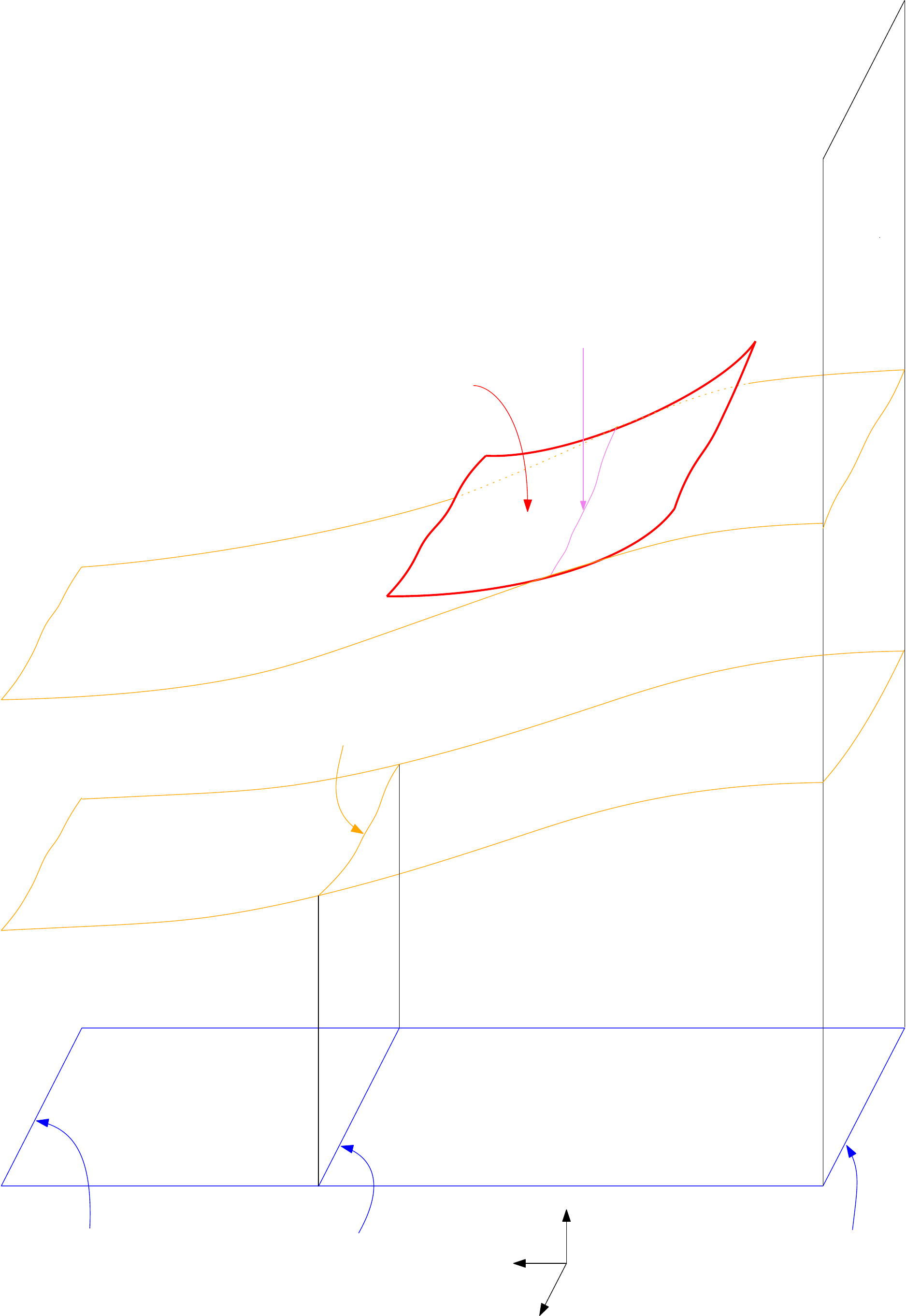}
			\put (42,76) {$\twoargmumuxtorus{0}{-\muxmulevelsetvalue}$}
			\put (32,71) {$\mathcal{B}^{[0,\muxmulevelsetvalue_0]}$}
			\put (38,15) {$\Sigma_0$}
			\put (4,5) {$\ell_{0,\leftu}$}
			\put (24,5) {$\ell_{0,u}$}
			\put (61,5) {$\ell_{0,-\rightu}$}
			\put (26,25) {$\nullhyparg{u}$ (truncated)}
			\put (21,45) {$\twoargroughtori{\timefunction_0,u}{\muxmulevelsetvalue}$}
			\put (32,-2.5) {$(x^2,x^3) \in \mathbb{T}^2$}
			\put (41,6) {$t$}
			\put (31,3) {$u \in\mathbb{R}$}
			\put (10,33) {$\hypthreearg{\timefunction_0}{[- \rightu,\leftu]}{\muxmulevelsetvalue}$}
			\put (10,50) {$\hypthreearg{0}{[- \rightu,\leftu]}{\muxmulevelsetvalue}$}
			\put (64,80) {$\nullhyparg{- \rightu}^{[0,\frac{4}{\mathring{\updelta}_*}]}$}
		\end{overpic}
		\caption{Data hypersurfaces in geometric coordinate space}
	\label{F:LATETIMEDATAHYPERSURFACEANDNULLDATAHYPERSURFACE}
	\end{figure}
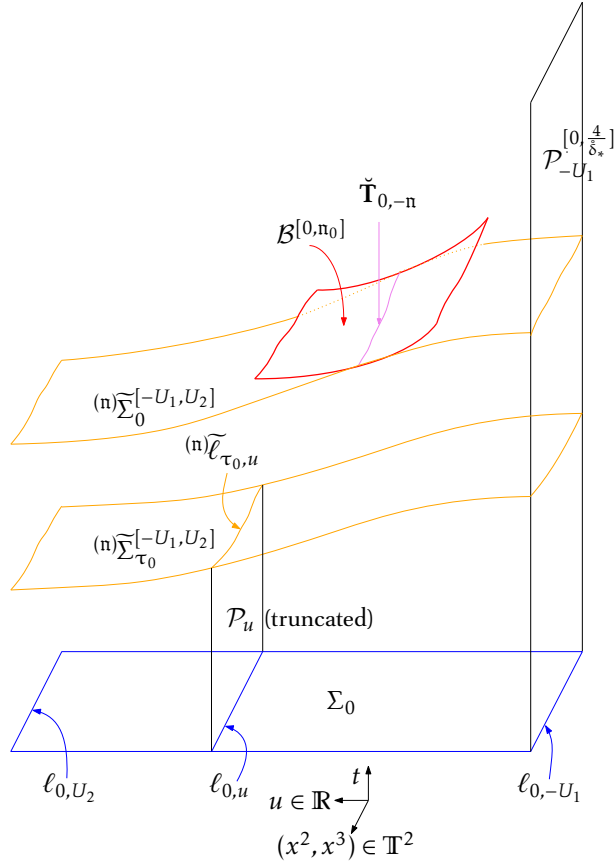
\end{center}

\subsubsection{The geometric div-curl-transport formulation of the flow}
\label{SSS:INTRONEWFORMULATION}
To derive estimates for the fluid variables, we fundamentally rely on 
a geometric reformulation of the flow, which was derived in
\cite{jS2019c} and which we restate below in Theorem~\ref{T:GEOMETRICWAVETRANSPORTSYSTEM}.
To aid our discussion in the remainder of the introduction, 
here we recall the formulation in schematic form:
\begin{subequations}
\begin{align} \label{E:WAVEINTRO}
 \square_{\gfour(\wavearray)} \Psi
	& 
	= 
	\nullform(\pmb{\partial} \wavearray, \pmb{\partial} \wavearray)
	+ 
	(\vortrenormalized,\GradEnt) \cdot \pmb{\partial} \wavearray 
	+ 
	(\VortVort, \DivGradEnt), 
		\\
	\Transport (\vortrenormalized,\GradEnt)
	& =  (\vortrenormalized,\GradEnt) \cdot \pmb{\partial} \wavearray, 
		\label{E:TRANSPORTINTRO} \\
\Transport (\VortVort,\DivGradEnt)
& = \nullform(\pmb{\partial} \wavearray,\partial \vortrenormalized)
		+
		\nullform(\pmb{\partial} \wavearray,\partial \GradEnt)
		+ 
		\GradEnt \cdot \nullform(\pmb{\partial} \wavearray, \pmb{\partial} \wavearray)
		+ 
		\GradEnt \cdot \GradEnt \cdot \pmb{\partial} \wavearray,
		\label{E:TRANSPORTFORMODIFIEDINTRO} \\
(\Flatdiv \vortrenormalized, \Flatcurl \GradEnt)
& = \partial \wavearray.
	\label{E:DIVCURLFORMODIFIEDINTRO}
\end{align}
\end{subequations}
In \eqref{E:WAVEINTRO}--\eqref{E:DIVCURLFORMODIFIEDINTRO},
$
\vortrenormalized = \frac{(\Flatcurl v)^i}{\exp (\LogDensity)}
$
is the specific vorticity
and
$
\GradEnt^i = \partial_i \Ent
$
is the entropy gradient vectorfield.
The modified fluid variables
$\VortVort$ and $\DivGradEnt$
are defined in Def.\,\ref{D:HIGHERORDERFLUIDVARIABLES}
and satisfy
$
\VortVort
\sim
\Flatcurl \vortrenormalized
$
and
$
\DivGradEnt
\sim
\Flatdiv \GradEnt
$,
where ``$\sim$'' denotes equality up to lower-order (in the sense of regularity) factors and terms.
Moreover,
$
\square_{\gfour(\wavearray)} f
$
is the covariant wave operator of $\gfour(\wavearray)$ acting on the scalar function $f$
(see Def.\,\ref{D:COVWAVEOP}),
$\pmb{\partial} f = (\partial_t f, \partial_1, \partial_2 f, \partial_3 f)$
is the array of Cartesian coordinate spacetime partial derivatives of $f$,
and
$\partial f = (\partial_1, \partial_2 f, \partial_3 f)$
is the array of Cartesian spatial partial derivatives of $f$.
The ``$\nullform$'' are \emph{null forms relative to} $\gfour$ (see Def.\,\ref{D:STANDARDNULLFORMS});
see Sect.\,\ref{SSS:GOODNULL} for a discussion of the crucial role they play in the proof.

We highlight that	when deriving
top-order $L^2$ estimates for the specific vorticity and entropy gradient, 
we crucially rely on the transport-div-curl equations
\eqref{E:TRANSPORTFORMODIFIEDINTRO}--\eqref{E:DIVCURLFORMODIFIEDINTRO}
to propagate sufficient regularity for the source terms $(\VortVort, \DivGradEnt)$
on RHS~\eqref{E:WAVEINTRO}; see Sect.\,\ref{S:TOPORDERELIPTICHYPERBOLICL2ESTIMATESFORSPECIFICVORTICITYANDENTROPYGRADIENT} for 
the details. 
That is, even though these terms formally satisfy
$
\VortVort
\sim
\Flatcurl \vortrenormalized
$
and
$
\DivGradEnt
\sim
\Flatdiv \GradEnt
$,
we cannot control $(\VortVort, \DivGradEnt)$
using only the transport equation
\eqref{E:TRANSPORTINTRO}. 
The reason is that \eqref{E:TRANSPORTINTRO} 
has a source term depending on
$\pmb{\partial} \wavearray$,
and thus, by estimating only the transport equation,
we could prove only that
$\Flatcurl \vortrenormalized, \, \Flatdiv \GradEnt$
are as regular as
$\pmb{\partial}^2 \wavearray$,
which would be insufficient regularity for treating
$\Flatcurl \vortrenormalized$
and
$\Flatdiv \GradEnt$
as source terms in the wave equation \eqref{E:TRANSPORTINTRO}.

\subsubsection{Three kinds of coordinate systems}
\label{SSS:INTROTHREECOORDINATES}
The proof of our main results is based on understanding the behavior of the solution
relative to three coordinate systems as well as understanding 
the degenerate (near the shock) transformation properties of the coordinate systems.

\begin{enumerate}
	\item (\textbf{The Cartesian coordinates} $(t,x^1,x^2,x^3)$).
		These are the fundamental coordinates relative to which the compressible
		Euler equations \eqref{E:INTROTRANSPORTVI}--\eqref{E:INTROBS} are posed.
		For the solution regime under study, the singularity coincides 
		with the blowup of $|\partial_1 \RRiemann|$, where $\RRiemann$ is the ``almost Riemann invariant'' introduced above.
	\item (\textbf{The geometric coordinates} $(t,u,x^2,x^3)$). As in the other works
		cited in Sect.\,\ref{SS:WORKSBUILDINGTOWARDSMAINTHEOREM},
		$u$ is an eikonal function, that is, a solution to the 
		eikonal equation $(\gfour^{-1})^{\alpha \beta}(\wavearray) \partial_{\alpha} u \partial_{\beta} u = 0$,
		where $\gfour$ is the acoustical metric. Its level sets $\nullhyparg{u}$ are characteristic for the
		wave operator $\square_{\gfour(\wavearray)}$ on LHS~\eqref{E:WAVEINTRO}.
		Our analytic framework is based on the fact that
		the solution remains rather smooth relative to the geometric coordinates,
		except that the high order geometric energies can blow up
		as the shock forms, which introduces severe technical
		difficulties into the analysis; see Sects.\,\ref{SSS:GEOMETRICENERGYESTIMATESONTHEROUGHFOLIATIONS}
		and \ref{SSS:INTROFLUIDVARIABLEELLITPICESTIMATES} for discussion
		of our geometric $L^2$ estimates.
	\item 
		\begin{itemize}
		\item
		(\textbf{The rough adapted coordinates} $(\timefunctionarg{\muxmulevelsetvalue},u,x^2,x^3)$). These are a new ingredient,
		fundamental for our main results. The $\timefunctionarg{\muxmulevelsetvalue}$ are a one-parameter family of 
		coordinate systems indexed by $\muxmulevelsetvalue \in [0,\muxmulevelsetvalue_0]$, where $\muxmulevelsetvalue_0 > 0$ is a real number
		depending on the solution regime under study. 
		They are constructed (see Sect.\,\ref{S:ROUGHTIMEFUNCTIONANDROUGHSUBSETS})
		so that their range is $[\timefunction_0,0]$ for some constant $\timefunction_0 < 0$,
		such that the initial data are given on the level set $\lbrace \timefunctionarg{\muxmulevelsetvalue} = \timefunction_0 \rbrace$
		(see Sect.\,\ref{SSS:PERTURBATIONSOFSIMPLEWAVESANDCAUCHYSTABILITY}),
		and such that the level set $\lbrace \timefunctionarg{\muxmulevelsetvalue} = 0 \rbrace$
		intersects the singular boundary precisely in one embedded, spacelike, two-dimensional torus.
		The union of the $\muxmulevelsetvalue$-indexed tori foliates the singular boundary;
		see Fig.\,\ref{F:INTROCARTESIANROUGHFOLIATIONS}, where the tori on the singular
		boundary are denoted by $\twoargmumuxtorus{0}{-\muxmulevelsetvalue}$.
		\item (\textbf{The rough adapted coordinates} $(\newtimefunction,u,x^2,x^3)$).
			With the help of the family of coordinate systems 
			$\lbrace (\timefunctionarg{\muxmulevelsetvalue},u,x^2,x^3) \rbrace_{\muxmulevelsetvalue \in [0,\muxmulevelsetvalue_0]}$,
			we construct another time function, denoted by $\newtimefunction$, 
			that foliates the entire region $\MInteresting$. $\newtimefunction$ is a $C^{1,1}$ function 
			-- \emph{and this is its optimal regularity regardless of the smoothness of the data} 
			(see Remark~\ref{R:NEWTIMEFUNCTIONISC11ANDNOTBETTERANDCONNECTIONTOCAUSALSTRUCTUREOFMUZEROLEVELSET}) -- 
			of the geometric coordinates $(t,u,x^2,x^3)$.
			In particular, the singular boundary portion under study is a submanifold-with-boundary
			that is contained in the level set $\lbrace \newtimefunction =  0 \rbrace$. 
			See Fig.\,\ref{F:INTERESTINGREGIONMAINRESULTS}, which displays two level sets of
			$\newtimefunction$, denoted by $\inthyp{0}{[- \rightu,\leftu]}$
			and
			$\inthyp{\timefunction_0}{[- \rightu,\leftu]}$.
\end{itemize}
		
\end{enumerate}

\subsubsection{The vanishing of $\upmu$ signifies the singularity formation}
\label{SSS:VANISHINGOFMUTIEDTOSINGULARITYFORMATION}
As in \cite{dC2007}, the formation of the shock singularity is precisely characterized by
the vanishing of the inverse foliation density $\upmu$, defined in \eqref{E:INTROINVERSEFOLIATION},
where in the present context, $\gfour$ is the acoustical metric.
We recall that the vanishing of $\upmu$ signifies the infinite density of the characteristics;
see Fig.\,\ref{F:INFININTEDENSITYOFCHARACTERISTICSONSINGULARBOUNDARY}, where $\upmu$
vanishes along $\mathcal{B}$ and the characteristics pile up there.
The blowup of some Cartesian partial derivative of the almost Riemann invariant $\RRiemann$
then follows from proving a bound of the schematic form 
$\upmu |X \RRiemann| \gtrsim 1$ in a neighborhood of the points where
$\upmu$ vanishes, where $X$ is a vectorfield that is transversal to the $\nullhyparg{u}$,
is $L^{\infty}$ close to the Cartesian partial derivative vectorfield $\partial_1$,
and has Euclidean length approximately equal to $1$.
This is essentially the same blowup-mechanism as in the irrotational and isentropic case,
as we described in Sect.\,\ref{SSS:CHRISTODOULOUBREAKTHROUGHONSHOCKFORMATION}.
Here we highlight that the vectorfield
$\muX \eqdef \upmu X$, which is a geometric replacement for the partial derivative vectorfield
$\geop{u}$ (in the geometric coordinate system), can be used to derive \emph{regular estimates} for the solution's transversal derivatives:
all fluid quantities and geometric tensors $Q$ 
that we use to study the solution satisfy $|\muX Q| \lesssim 1$ up to the singularity;
the $\upmu$-weight in the definition of $\muX$ precisely cancels out the singular
behavior $|X Q|$ as $\upmu \downarrow 0$;
see Sect.\,\ref{SSS:REGULARESTIMATESONROUGHFOLIATIONS} for further discussion.

\subsubsection{The causal structure of the singular boundary}
\label{SSS:STRUCTUREOFSINGULARITY}
The discussion in Sect.\,\ref{SSS:VANISHINGOFMUTIEDTOSINGULARITYFORMATION} suggests
that the singular boundary should be the entire hypersurface $\lbrace \upmu = 0 \rbrace$.
However, as we already discussed in Sect.\,\ref{SS:OVERVIEWOFDIFFICULTIES}, 
only a subset of $\lbrace \upmu = 0 \rbrace$ is relevant for the maximal development.
One reason is that formally, for the solutions under study,
$\lbrace \upmu = 0 \rbrace$ is the limit of hypersurfaces 
$\lbrace \upmu = \mulevelsetvalue \rbrace$ as $\mulevelsetvalue \downarrow 0$,
and for $\mulevelsetvalue > 0$, these hypersurfaces have a spacelike part, a null part, and timelike part.
For $\mulevelsetvalue > 0$, our analysis yields access to the portions of
$\lbrace \upmu = \mulevelsetvalue \rbrace$ that are $\gfour$-spacelike or $\gfour$-null, but not
necessarily the timelike portion. Hence, the portion of $\lbrace \upmu = 0 \rbrace$ that makes up the singular boundary
can be viewed as the limit as $\mulevelsetvalue \downarrow 0$, of the portion of
$\lbrace \upmu = \mulevelsetvalue \rbrace$ that is $\gfour$-spacelike or $\gfour$-null.
From a careful analysis of the causal structure of the hypersurfaces
$\lbrace \upmu = \mulevelsetvalue \rbrace$, which we provide in Lemma~\ref{L:IDENTITIESFORGRADIENTVECTORFIELDOFINVERSEFOLIATIONDENSITY}, 
we find that the accessible portion of
$\lbrace \upmu = 0 \rbrace$ that arises in the limit $\mulevelsetvalue \downarrow 0$ 
is $\lbrace \upmu = 0 \rbrace \cap \lbrace \muX \upmu \leq 0 \rbrace$,
where the vectorfield $\muX$ is the same one described in Sect.\,\ref{SSS:VANISHINGOFMUTIEDTOSINGULARITYFORMATION}.
More precisely, our results are local in spacetime, and our main theorem yields the structure of the singular
boundary portion
$\lbrace \upmu = 0 \rbrace \cap \lbrace - \muxmulevelsetvalue_0 \leq \muX \upmu \leq 0 \rbrace$,
which we denote by $\mathcal{B}^{[0,\muxmulevelsetvalue_0]}$ in Fig.\,\ref{F:INTERESTINGREGIONMAINRESULTS}.
We highlight the following key point: even if one were to construct the complementary surface portion
$\lbrace \upmu = 0 \rbrace \cap \lbrace \muX \upmu > 0 \rbrace$ in geometric coordinates,
the map from geometric to Cartesian coordinates would \emph{fail} to be one-to-one
on the full surface $\lbrace \upmu = 0 \rbrace$, though it \emph{is} one-to-one
on the subset $\mathcal{B}^{[0,\muxmulevelsetvalue_0]}$; see Theorem~\ref{T:ABBREVIATEDSTATEMENTOFMAINRESULTS}.
From a different perspective, in the context of Fig.\,\ref{F:MAXDEVELOPMENTINCARTESIAN},
the surface portion $\lbrace \upmu = 0 \rbrace \cap \lbrace \muX \upmu > 0 \rbrace$
is ``cut off'' by the Cauchy horizon $\underline{C}$, i.e., in the classical maximal globally hyperbolic development,
$\underline{C}$ develops ``before''
$\lbrace \upmu = 0 \rbrace \cap \lbrace \muX \upmu > 0 \rbrace$ has a chance to form.


In view of the above discussion, in our main theorem, we 
construct the singular boundary portion 
$
\mathcal{B}^{[0,\muxmulevelsetvalue_0]} 
\eqdef \bigcup_{\muxmulevelsetvalue \in [0,\muxmulevelsetvalue_0]}
\lbrace \upmu = 0 \rbrace 
\cap 
\lbrace \muX \upmu = - \muxmulevelsetvalue \rbrace
$;
see Fig.\,\ref{F:INTERESTINGREGIONMAINRESULTS}.
Our aforementioned \emph{transversal convexity} assumption on the initial data,
which we are able to propagate throughout the evolution
(see \eqref{E:MUTRANSVERSALCONVEXITY}),
ensures that
the hypersurfaces
$\lbrace \upmu = 0 \rbrace$
and
$\lbrace \muX \upmu = - \muxmulevelsetvalue \rbrace$
intersect transversally in embedded, two-dimensional $\gfour$-spacelike tori that we denote by
$\twoargmumuxtorus{0}{-\muxmulevelsetvalue}$. In particular,
we show that $\mathcal{B}^{[0,\muxmulevelsetvalue_0]}$ is an embedded three-dimensional
submanifold-with-boundary that is foliated by the tori
$\twoargmumuxtorus{0}{-\muxmulevelsetvalue}$.
Its past boundary, denoted by $\partial_-\mathcal{B}^{[0,\muxmulevelsetvalue_0]}$,
is the torus $\twoargmumuxtorus{0}{0}$, a set that we have also been referring to as the crease.

We next highlight that the crease plays a distinguished role in the 
shock development problem described in Sect.\,\ref{SSS:PROGRESSONSHOCKDEVELOPMENT}; 
once that problem is solved,
the hypersurface of discontinuity $\mathcal{K}$ will emanate
from the crease; see Fig.\,\ref{F:MAXDEVELOPMENTWITHSHOCKHYPERSURFACEINCARTESIAN} 
In particular, the crease is a crucial component of the ``data'' for the 
shock development problem.

\subsubsection{The one-parameter family of rough time functions $\lbrace \timefunctionarg{\muxmulevelsetvalue} \rbrace_{\muxmulevelsetvalue \in [0,\muxmulevelsetvalue_0]}$}
\label{SSS:INTROROUGHTIMEFUNCTIONS}
We follow the solution up to the singular boundary portion 
$\mathcal{B}^{[0,\muxmulevelsetvalue_0]} = \bigcup_{\muxmulevelsetvalue \in [0,\muxmulevelsetvalue_0]} \twoargmumuxtorus{0}{-\muxmulevelsetvalue}$
by constructing a one-parameter family of time functions
$\lbrace \timefunctionarg{\muxmulevelsetvalue} \rbrace_{\muxmulevelsetvalue \in [0,\muxmulevelsetvalue_0]}$
and controlling the solution on the level sets of each time function.
The $\timefunctionarg{\muxmulevelsetvalue}$ are constructed such that their range is $[\timefunction_0,0]$
for some small parameter $\timefunction_0 < 0$
and such that $\twoargmumuxtorus{0}{-\muxmulevelsetvalue} \subset \lbrace \timefunctionarg{\muxmulevelsetvalue} = 0 \rbrace$.
We construct $\timefunctionarg{\muxmulevelsetvalue}$ by setting it equal to $\upmu$ on the hypersurface
$\lbrace \muX \upmu = - \muxmulevelsetvalue \rbrace$ 
and then transporting it along the flow of a
well-constructed vectorfield that is transversal to $\lbrace \muX \upmu = - \muxmulevelsetvalue \rbrace$,
which ensures that the level sets of $\timefunctionarg{\muxmulevelsetvalue}$ are $\gfour$-spacelike
in the region of classical existence. 
In Fig.\,\ref{F:CARTESIANROUGHFOLIATIONFUTUREOFCREASE}, we denote the hypersurface
$\lbrace \muX \upmu = - \muxmulevelsetvalue \rbrace$ by 
$\datahypfortimefunctionarg{-\muxmulevelsetvalue}$
and we denote the transversal vectorfield by $\Wtransarg{\muxmulevelsetvalue}$.
See Sect.\,\ref{S:ROUGHTIMEFUNCTIONANDROUGHSUBSETS} 
for the details on the construction of $\timefunctionarg{\muxmulevelsetvalue}$.

A key feature of our construction is that for $\timefunction \in [\timefunction_0,0]$, we have:
\begin{align} \label{E:INTROMUMIN}
	\min_{\lbrace \timefunctionarg{\muxmulevelsetvalue} = \timefunction \rbrace} \upmu
	& = - \timefunction,
\end{align}
and the minimum value is achieved precisely along the set
$\twoargmumuxtorus{-\timefunction}{-\muxmulevelsetvalue} 
= 
\lbrace \timefunctionarg{\muxmulevelsetvalue} = \timefunction \rbrace
\cap
\lbrace \muX \upmu = - \muxmulevelsetvalue \rbrace
$.
Thus, for $\mulevelsetvalue \in [0,-\timefunction_0]$, we have
$
\twoargmumuxtorus{\mulevelsetvalue}{-\muxmulevelsetvalue} 
=
\lbrace \upmu = \mulevelsetvalue \rbrace
\cap
\lbrace \muX \upmu = - \muxmulevelsetvalue \rbrace
$.
Our high order energy estimates feature crucial factors of $1/\upmu$,
and by \eqref{E:INTROMUMIN}, 
in rough adapted coordinates $(\timefunctionarg{\muxmulevelsetvalue},u,x^2,x^3)$,
such factors are bounded in magnitude on
$\lbrace \timefunctionarg{\muxmulevelsetvalue} = \timefunction \rbrace$
by $\leq \frac{1}{|\timefunction|}$.
Thus, our energy estimates on the level sets of the rough time function 
feature difficult factors of $\frac{1}{|\timefunction|}$,
which leads to singular estimates at the high derivative levels.
While related difficulties were present in all the works that we cited above 
on shock formation without symmetry assumptions,
the estimate \eqref{E:INTROMUMIN} allows for a simplified approach
to handling the degeneracy; it directly connects the degeneracy to a coordinate function.
We refer to Sect.\,\ref{SSS:PROOFOFAPRIORIL2ESTIMATESWAVEVARIABLES} for our detailed
analysis of the factors of $\frac{1}{|\timefunction|}$ and the way they affect
our energy estimates, and to Sect.\,\ref{SSS:GEOMETRICENERGYESTIMATESONTHEROUGHFOLIATIONS} for an overview of these
estimates.

We emphasize the following key point:
\begin{quote}	
		The rough time functions $\timefunctionarg{\muxmulevelsetvalue}$ 
		are precisely adapted to the structure of the singularity and, 
		in particular, are not more regular than $\upmu$.
		Since it turns out that $\upmu$ is one degree less differentiable\footnote{In particular, in view of the transport
		equation \eqref{E:MUTRANSPORT}, we see that $\upmu$ cannot be more regular than the source term
		$\muX \wavearray$.
		This is optimal, as can already be seen in plane-symmetry
		via the explicit formula \eqref{AE:MUXMUPLANESYMMETRYWITHEXPLICITSOURCE}.} 
		than the fluid variables $\wavearray$,
		we cannot use commutation vectorfields adapted to the level sets of
		$\timefunctionarg{\muxmulevelsetvalue}$ to control the solution up to top-order;
		that approach would lead to derivative loss.
		For this reason, as in other works on shock formation,
		we use the eikonal function $u$
		to construct\footnote{It does not seem possible to adequately control the solution
		by commuting the equations with the standard Cartesian partial derivatives $\partial_{\alpha}$;
		the Cartesian vectorfields are not adapted to the singularity, and 
		that approach would lead to singular error terms that we have no obvious means to control.} 
		commutation vectorfields that allow us to control the solution.
\end{quote}

Finally, we point out that it is not too difficult to construct the $\timefunctionarg{\muxmulevelsetvalue}$;
the difficult part of the analysis is controlling the fluid solution and the acoustical geometry
all the way up the hypersurface $\lbrace \timefunctionarg{\muxmulevelsetvalue} = 0 \rbrace$.

\subsubsection{The bootstrap assumptions}
\label{SSS:BOOTSTRAPASSUMPTIONS}
To derive estimates, we fix $\muxmulevelsetvalue \in [0,\muxmulevelsetvalue_0]$ and make an elaborate set of bootstrap assumptions 
on an ``open-at-the-top'' slab $[\timefunction_0,\timefunctionboot) \times [- \rightu,\leftu] \times \mathbb{T}^2$
corresponding to the rough adapted coordinates $(\timefunctionarg{\muxmulevelsetvalue},u,x^2,x^3)$
(the corresponding region in geometric coordinate space is denoted by
$\twoargMrough{[\timefunction_0,\timefunction),[- \rightu,\leftu]}{\muxmulevelsetvalue}$ in the bulk of the paper).
More precisely, on the slab, our bootstrap assumptions describe the following:
\begin{itemize}
	\item The behavior of $\upmu$, including the 
	key property of transversal convexity;
	see Sect.\,\ref{SSS:BAFORINVERSEFOLIATIONDENSITY}.
	\item The rough time function $\timefunctionarg{\muxmulevelsetvalue}$;
		see Sect.\,\ref{SSS:BOOTSTRAPROUGHTIMEFUNCTION}.
	\item The properties of various change of variables maps; 
		see Sect.\,\ref{SSS:BACHOVMAPS}.
	\item The structure and locations of various embedded submanifolds;
		see Sect.\,\ref{SSS:BOOTSTRAPASSUMPTIONFORTORISTRUCTURE}.
	\item The size of the Cartesian coordinates $t$ and $x^1$ on the slab; 
		see Sect.\,\ref{SSS:BOOTSTRAPASSUMPTIONSSIZEOFCARTESIANTANDX1}.
	\item The ``soft'' regularity properties of various quantities on the closure of the slab; 
		see Sect.\,\ref{SSS:SOFTBACONCERNINGREGULARITY}.
	\item The $L^{\infty}$-size of the fluid variables and their $\nullhyparg{u}$-tangential derivatives up to mid-order;
		see Sect.\,\ref{SSS:FUNDAMENTALQUANTITATIVE}.
	\item The $L^{\infty}$-size of the transversal and mixed tangential-transversal derivatives of the fluid variables
		and the acoustic geometry;
		see Sect.\,\ref{SSS:AUXBOOTSTRAP}.
	\item The size of the energies and null-fluxes for the wave variables $\wavearray$ up to top-order;
		see Sect.\,\ref{SS:BOOTSTRAPASSUMPTIONSFORTHEWAVEENERGIES}.
\end{itemize}
Our primary analytic tasks in the paper are to derive strict improvements of the bootstrap assumptions
by making suitable assumptions on the data.
Once this has been accomplished, 
a standard continuity argument, 
provided by Prop.\,\ref{P:CONTINUATIONCRITERIA} and
the proof of Theorem~\ref{T:EXISTENCEUPTOTHESINGULARBOUNDARYATFIXEDKAPPA},
implies that the solution exists classically on the maximal slab 
$[\timefunction_0,0) \times [- \rightu,\leftu] \times \mathbb{T}^2$
and satisfies the bootstrap assumptions there. Hence, Theorem~\ref{T:EXISTENCEUPTOTHESINGULARBOUNDARYATFIXEDKAPPA}
yields the main results at fixed $\muxmulevelsetvalue$,
which in turn form the main ingredients in the proof of the central theorem on the singular boundary,
Theorem~\ref{T:DEVELOPMENTANDSTRUCTUREOFSINGULARBOUNDARY}.
In Sect.\,\ref{SS:SUMMARYOFIMPROVEMENTOFBOOTSTRAPASSUMPTIONS},
we provide a road map, indicating the spots in the article where we derive improvements
of the bootstrap assumptions.

\subsubsection{Regular estimates with respect to the geometric coordinates on the rough foliations}
\label{SSS:REGULARESTIMATESONROUGHFOLIATIONS}
A key idea of the proof, going back to \cite{dC2007},
is to prove that the solution remains rather smooth relative to the geometric
coordinates $(t,u,x^2,x^3)$. We implement this by following the approach of
\cites{dC2007,jSgHjLwW2016,jLjS2018,jLjS2021}
and using the eikonal function to construct a set of 
\emph{commutation vectorfields}:
\begin{align} \label{E:INTROCOMMUTATORS}
	\Fullset 
	& \eqdef \lbrace \Lunit, \muX, \Yvf{2}, \Yvf{3} \rbrace
\end{align}
that are adapted to the characteristics; see Fig.\,\ref{F:COMMUTATORVECTORFIELDSINCARTESIANCOORDINATES} for a picture of these
vectorfields (with one spatial dimension suppressed), 
and see Sect.\,\ref{S:ACOUSTICGEOMETRYANDCOMMUTATORVECTORFIELDS} for our detailed construction 
of the elements of $\Fullset$.
In \eqref{E:INTROCOMMUTATORS} and throughout,
$\Lunit$ is a $\gfour$-null vectorfield (i.e., $\gfour(\Lunit,\Lunit) = 0$) 
that is proportional to the (geodesic) gradient vectorfield $(\gfour^{-1})^{\alpha \beta} \partial_{\alpha} u \partial_{\beta}$,
normalized by $\Lunit t = 1$, and tangent to the characteristics $\nullhyparg{u}$.
$\Yvf{2}, \Yvf{3}$ are geometric replacements for $\partial_2$ and $\partial_3$
that span the tangent space of the \emph{smooth tori} $\ell_{t,u} \eqdef \nullhyparg{u} \cap \Sigma_t$.
The vectorfields $\lbrace \Lunit, \Yvf{2}, \Yvf{3} \rbrace$ span the tangent space of $\nullhyparg{u}$
at each of its points.
$\muX$ is tangent to the flat Cartesian hypersurface $\Sigma_t$, 
$\gfour$-orthogonal to $\ell_{t,u}$, and normalized\footnote{In this paper, this 
normalization condition is equivalent to $\gfour(\muX,\muX) = \upmu^2$.} 
by $\muX u = 1$. Note that $\muX$ is transversal to the characteristics. This is important because
all the degeneracies in the problem of shock formation are tied to the behavior of transversal derivatives.
We also note that, in accordance with the discussion in Sect.,\ref{SSS:VANISHINGOFMUTIEDTOSINGULARITYFORMATION},
we have $\muX = \upmu X$, where $\gfour(X,X) = 1$.

	\begin{center}
	\begin{figure}  
		\begin{overpic}[scale=.6, grid = false, tics=5, trim=-.5cm -1cm -1cm -.5cm, clip]{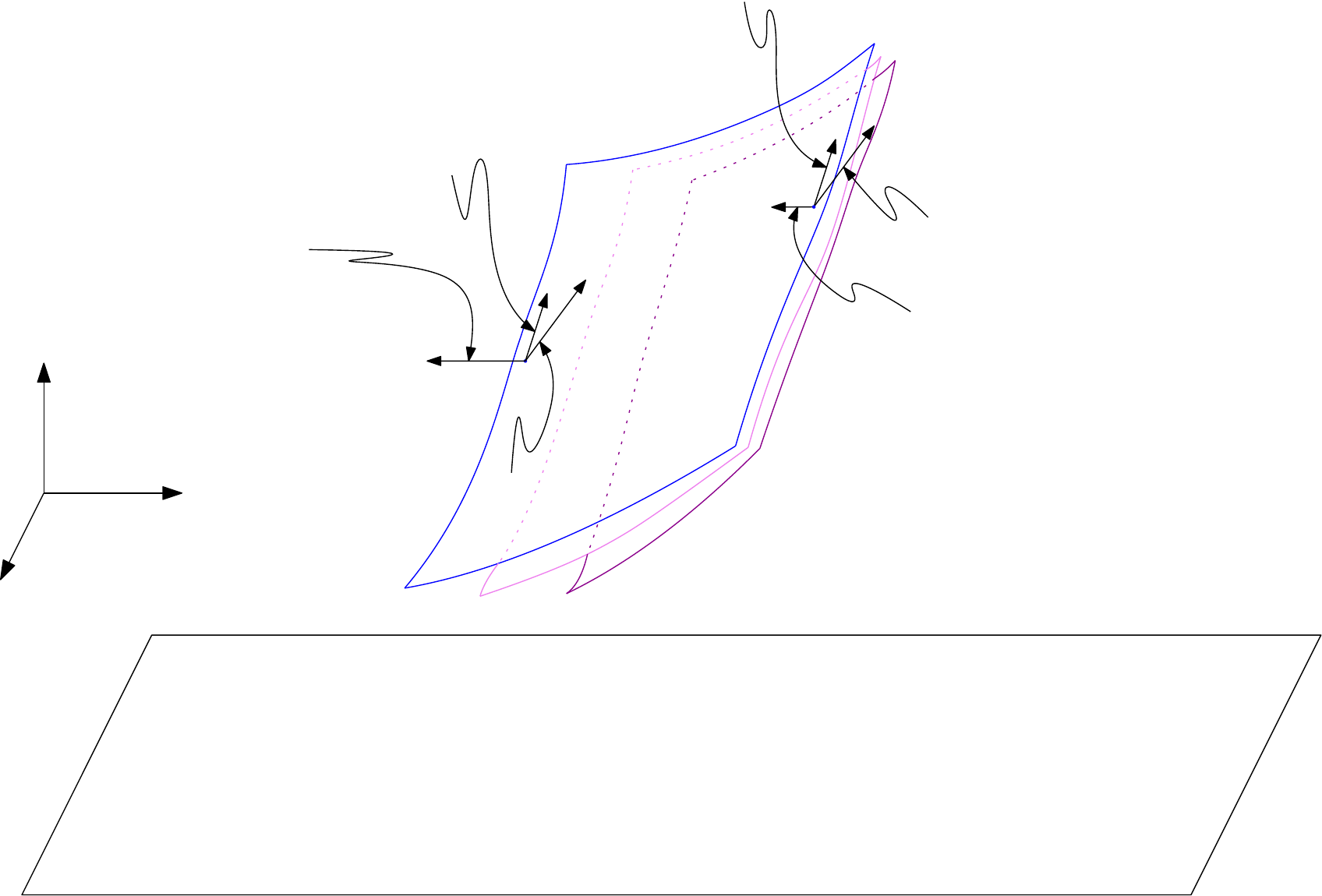}
			\put (50,15) {$\Sigma_0$}
			\put (13,35) {$(x^2,x^3) \in \mathbb{T}^2$}
			\put (4,26)	{$x^1 \in \R$}
			\put (3,38) {$t$}
			\put (21,50) {$\muX$}
			\put (30,57) {$\Yvf{A}$}
			\put (38,32) {$L$}
			\put (66,44) {$\muX$}
			\put (67,50) {$L$}
			\put (52,69) {$\Yvf{A}$}
		\end{overpic}
		\caption{The commutation vectorfields in Cartesian coordinates}
		\label{F:COMMUTATORVECTORFIELDSINCARTESIANCOORDINATES}
	\end{figure}
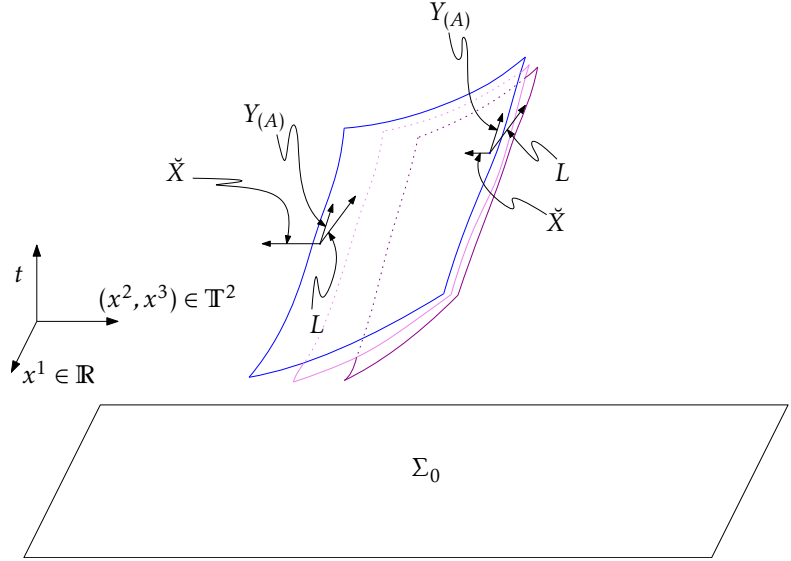
\end{center}	

As in \cite{jLjS2021}, we use the full set of commutators $\mathcal{Z}$ to derive
$L^{\infty}$ and H\"{o}lder estimates for solutions,
while for deriving energy estimates,
we only need to use the following $\nullhyparg{u}$-tangent subset:
\begin{align} \label{E:INTROTANGENTCOMMUTATORS}
	\Tanset 
	& \eqdef \lbrace \Lunit, \Yvf{2}, \Yvf{3} \rbrace.
\end{align}
In particular, our results show that, based on our assumptions on the data, 
the up-to-fourth order derivatives
of
$\wavearray$, $\vortrenormalized$, and $\GradEnt$
with respect to the elements of $\mathcal{Z}$ are bounded in $L^{\infty}$.
This is equivalent to the $L^{\infty}$-boundedness of the up-to-fourth order derivatives of these quantities
with respect to the geometric coordinate partial derivatives.
Similar results hold for the up-to-third order
derivatives of $\VortVort$ and $\DivGradEnt$,
and due to our assumptions on the data, 
all fluid variables enjoy additional regularity in directions tangent to 
the characteristics $\nullhyparg{u}$.
Similar results also hold for the rough time functions
$\timefunctionarg{\muxmulevelsetvalue}$
and various embeddings of various submanifolds into geometric or Cartesian coordinate space,
such as the hypersurfaces
$
\datahypfortimefunctionarg{-\muxmulevelsetvalue}
\eqdef
\lbrace
	\muX \upmu = - \muxmulevelsetvalue
\rbrace
$
and the \emph{$\upmu$-adapted tori}
$
\twoargmumuxtorus{\mulevelsetvalue}{-\muxmulevelsetvalue}
\eqdef
\lbrace
	\upmu = \mulevelsetvalue
\rbrace
\cap
\lbrace
	\muX \upmu = - \muxmulevelsetvalue
\rbrace
$.

To derive the $L^{\infty}$ and H\"{o}lder estimates, 
we make ``fundamental'' $L^{\infty}$ bootstrap assumptions for the fluid 
wave variables and their $\Tanset$ derivatives up to mid-order,
as we described in Sect.\,\ref{SSS:BOOTSTRAPASSUMPTIONS}.
By combining our fundamental bootstrap assumptions with some auxiliary ones,
we commute all relevant equations with the elements of
$\Fullset$ and treat them as transport equations with
derivative-losing source terms. This allows us to derive
$L^{\infty}$ and H\"{o}lder estimates for the solution's transversal
and mixed transversal-tangential derivatives
and to control the embeddings mentioned in the previous paragraph.
We carry out this analysis in
Sects.\ref{S:EMBEDDINGSANDFLOWMAPS}--\ref{S:LINFINITYFLUIDANDEIKONALANDIMPROVEMENTOFAUX}.
In Sect.\,\ref{S:SHARPCONTROLOFMUANDPROPERTIESOFCHOVGEOTOCARTESIAN},
we derive related, but much sharper, pointwise estimates for $\upmu$ 
at the low derivative levels, which in particular yield the crucial estimate
\eqref{E:INTROMUMIN}. 
Near the end of the paper, we use our energy estimates and Sobolev embedding
to improve the fundamental bootstrap assumptions; see
Sect.\,\ref{S:IMPROVEMENTSOFFUNDAMENTALQUANTITATIVEBOOTSTRAPASSUMPTIONS}.

Finally, we will briefly discuss the regularity of the Cartesian components of the elements of $\Fullset$.
They have regularity at the schematic level $\pmb{\partial} u$,
and it turns out that by using renormalizations (which we refer to a ``modified quantities'') 
and elliptic estimates on co-dimension $2$ hypersurfaces
(specifically, the \emph{rough tori}
$
\twoargroughtori{\timefunction,u}{\muxmulevelsetvalue}
\eqdef
\lbrace \timefunctionarg{\muxmulevelsetvalue} = \timefunction \rbrace
\cap
\nullhyparg{u}
$),
which are techniques that originated in \cites{dCsK1993,sKiR2003,dC2007},
we can show that 
the elements of $\Fullset$ have just enough regularity
to allow us to derive energy estimates up to top order.
We construct the modified quantities in Sect.\,\ref{S:CONSTRUCTIONOFMODIFIEDQUANTITIES},
we derive the elliptic estimates in
Sect.\,\ref{S:ELLIPTICESTIAMTESACOUSTICGEOMETRYONROUGHTORI},
and we derive the ``final'' $L^2$ estimates involving these quantities
Sect.\,\ref{S:WAVEANDACOUSTICGEOMETRYAPRIORIESTIMATES}.

\subsubsection{Smooth geometry versus rough geometry and elliptic estimates for $\upchi$}
\label{SSS:SMOOTHVSROUGH}
Although regularity considerations force us to commute the equations with 
elements of $\Fullset$ and $\Tanset$ (e.g., the rough adapted coordinate partial derivative vectorfields
do not have sufficient regularity for commuting the equations up to top-order),
in order to detect the precise structure of the singular boundary, 
we must derive estimates on the rough hypersurfaces
$
\lbrace \timefunctionarg{\muxmulevelsetvalue} = \timefunction \rbrace
$,
the characteristics $\nullhyparg{u}$,
and the rough tori 
$
\twoargroughtori{\timefunction,u}{\muxmulevelsetvalue}
\eqdef 
\lbrace \timefunctionarg{\muxmulevelsetvalue} = \timefunction \rbrace
\cap
\nullhyparg{u}
$.
For this reason, throughout the analysis, 
we have to control the geometry of these surfaces
(e.g., the Gauss curvature of $\twoargroughtori{\timefunction,u}{\muxmulevelsetvalue}$)
by quantitatively relating it to the elements of $\Fullset$ and $\Tanset$
and their derivatives.
A particularly noteworthy manifestation of this issue is the following:
we must control the top-order derivatives of the 
null second fundamental form 
$\upchi = \frac{1}{2} \angLie_{\Lunit} \gtorus$,
where $\gtorus$ is the Riemannian metric induced by $\gfour$ on $\ell_{t,u} = \Sigma_t \cap \nullhyparg{u}$.
Here, $\angLie_{\Lunit} \gtorus$ denotes Lie differentiation with respect to $\Lunit$ followed by
$\gfour$-orthogonal projection onto $\ell_{t,u}$.
Note that the  ``smooth torus'' $\ell_{t,u}$ is \emph{not} adapted (or necessarily fully contained in) the rough spacetime
regions under study. Even though $\upchi$ is $\ell_{t,u}$-tangent,
to avoid derivative loss in the top-order estimates on the rough spacetime regions,
we have to derive elliptic estimates for $\upchi$
on the \emph{rough tori}
$
\twoargroughtori{\timefunction,u}{\muxmulevelsetvalue}
$.
To close these elliptic estimates,
we must suitably relate $\upchi$ to a tensor on 
$
\twoargroughtori{\timefunction,u}{\muxmulevelsetvalue}
$
and control the Gauss curvature of
$
\twoargroughtori{\timefunction,u}{\muxmulevelsetvalue}
$
(viewed as a subset of spacetime equipped with the metric $\gfour$),
all while being mindful about crucial factors of $\upmu$,
whose vanishing signifies the shock singularity.
This delicate analysis is located in Sect.\,\ref{S:ELLIPTICESTIAMTESACOUSTICGEOMETRYONROUGHTORI},
where we construct two distinct frames on
$\nullhyparg{u}$,
one adapted to the smooth tori $\ell_{t,u}$ and one adapted to the rough tori 
$
\twoargroughtori{\timefunction,u}{\muxmulevelsetvalue}
$,
and we control the relationship between the two frames
as a key step in obtaining the elliptic estimates.

\subsubsection{The role of the good null structure}
\label{SSS:GOODNULL}
As in many of the aforementioned works on shock formation, in the present paper,
we crucially rely on the fact that all the derivative-quadratic inhomogeneous terms in the system 
\eqref{E:WAVEINTRO}--\eqref{E:DIVCURLFORMODIFIEDINTRO}
are null forms relative to $\gfour$.
From the point of view of analysis,
the crucial point is that when expanded relative to (say) the commutator frame \eqref{E:INTROCOMMUTATORS},
the $\upmu$-weighted\footnote{In our analysis, we introduce a $\upmu$ weight into various equations and estimates.} 
null forms satisfy, schematically,
$
\upmu \nullform(\pmb{\partial} \wavearray, \pmb{\partial} \wavearray)
= \muX \wavearray \cdot P \wavearray
+
\upmu P \wavearray \cdot P \wavearray
$,
where $P$ schematically denotes elements of the $\nullhyparg{u}$-tangent subset \eqref{E:INTROTANGENTCOMMUTATORS};
see Lemma~\ref{L:CRUCIALSTRUCTUREOFNULLFORMS} for a more detailed statement concerning the precise null forms
we encounter in our analysis.
In particular, in the expansion of $\upmu \nullform(\pmb{\partial} \wavearray, \pmb{\partial} \wavearray)$,
terms proportional to $\muX \wavearray \cdot \muX \wavearray$
are completely absent. This is crucial because, if present, signature considerations
would imply that such terms would have to be accompanied by a factor of $\upmu^{-1}$ 
(i.e., the term would be proportional to $\upmu^{-1} \muX \wavearray \cdot \muX \wavearray$),
and factor of $\upmu^{-1}$ (which blows up as $\upmu \to 0$) would obstruct our philosophy of
deriving regular estimates for the solution's $\mathcal{Z}$-derivatives.
Similar remarks apply to all the other null forms in \eqref{E:WAVEINTRO}--\eqref{E:DIVCURLFORMODIFIEDINTRO}.

The upshot is the following: all of the derivative-quadratic terms in the system
\eqref{E:WAVEINTRO}--\eqref{E:DIVCURLFORMODIFIEDINTRO}
that drive the formation of the shock are ``hidden'' in the definition of the covariant wave operator $\square_{\gfour(\wavearray)}$
on LHS~\eqref{E:WAVEINTRO} and become visible when one expands
$\square_{\gfour(\wavearray)} \Psi$ relative to the Cartesian coordinates (say, via the formula \eqref{E:WAVEOPERATORARBITRARYCOORDINATES}).
The virtue of $\square_{\gfour(\wavearray)}$ is that there is an advanced machinery
for deriving geometric commutator and multiplier energy estimates for such operators,
and we exploit this machinery throughout the paper.

\subsubsection{Geometric wave and transport energy estimates on the rough foliations, with singular high order behavior}
\label{SSS:GEOMETRICENERGYESTIMATESONTHEROUGHFOLIATIONS}
\begin{center}
	\begin{figure}  
		\begin{overpic}[scale=.4, grid = false, tics=5]{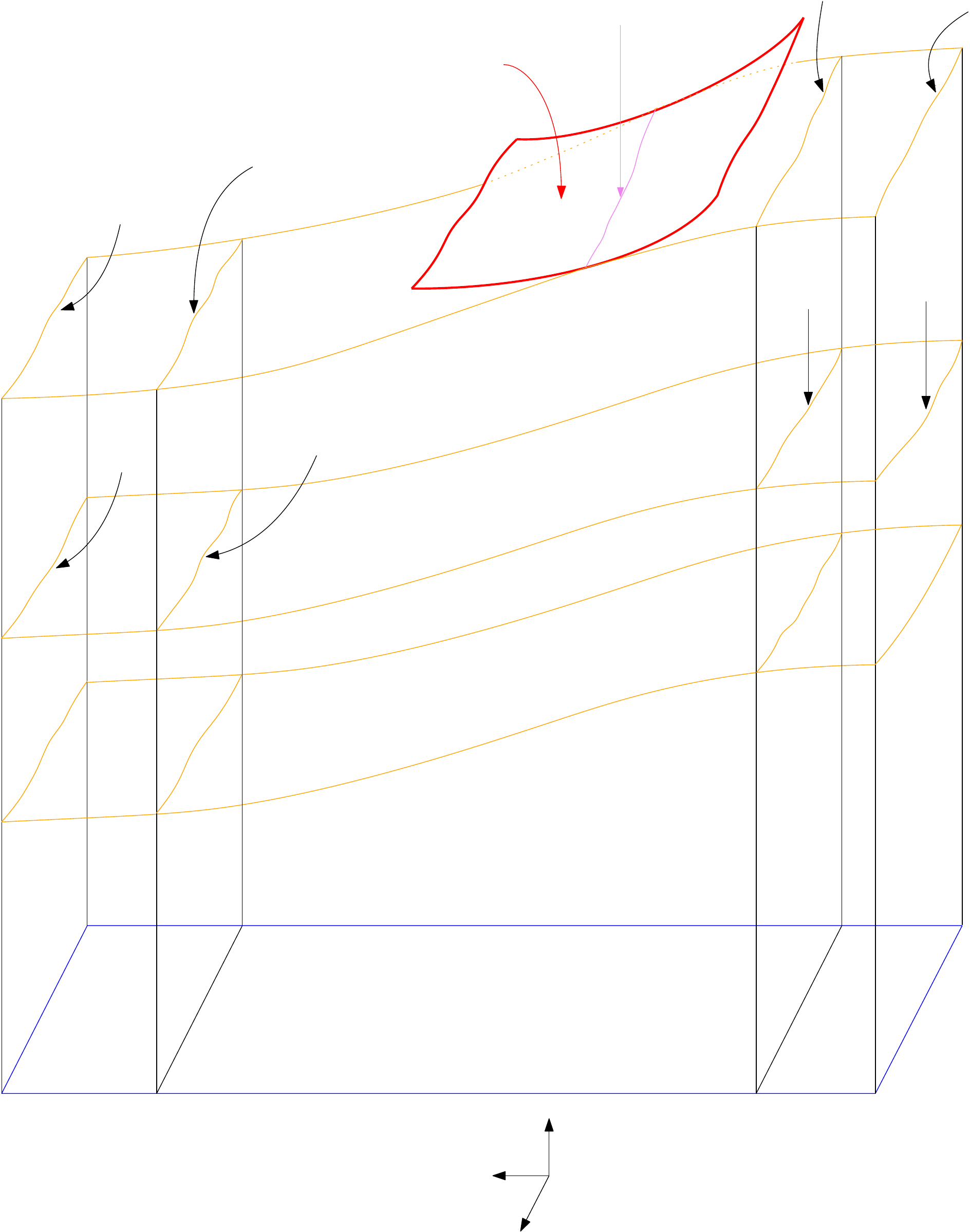}
			\put (48,99) {$\twoargmumuxtorus{0}{-\muxmulevelsetvalue}$}
			\put (37,97) {$\mathcal{B}^{[0,\muxmulevelsetvalue_0]}$}
			\put (40,17) {$\Sigma_0$}
			\put (37,-4) {$(x^2,x^3) \in \mathbb{T}^2$}
			\put (42,7) {$t$}
			\put (30.5,3.5) {$u \in\mathbb{R}$}
			\put (17,75) {$\hypthreearg{0}{[- \rightu,\leftu]}{\muxmulevelsetvalue}$}
			\put (33,59) {$\hypthreearg{\timefunction}{[- \rightu,\leftu]}{\muxmulevelsetvalue}$}
			\put (33,43) {$\hypthreearg{\timefunction_0}{[- \rightu,\leftu]}{\muxmulevelsetvalue}$}
			\put (3,64) {$\twoargroughtori{\timefunction,\leftu}{\muxmulevelsetvalue}$}
			\put (18,65) {$\twoargroughtori{\timefunction,\interestingu}{\muxmulevelsetvalue}$}
			\put (56,77) {$\twoargroughtori{\timefunction,-\interestingu}{\muxmulevelsetvalue}$}
			\put (70,77) {$\twoargroughtori{\timefunction,-\rightu}{\muxmulevelsetvalue}$}
			\put (3,84) {$\twoargroughtori{0,\leftu}{\muxmulevelsetvalue}$}
			\put (15,88) {$\twoargroughtori{0,\interestingu}{\muxmulevelsetvalue}$}
			\put (58,101) {$\twoargroughtori{0,-\interestingu}{\muxmulevelsetvalue}$}
			\put (74,100) {$\twoargroughtori{0,-\rightu}{\muxmulevelsetvalue}$}
			\put (30,47) {\rotatebox{30}{$\twoargMrough{[- \rightu,\leftu],[\timefunction_0,\timefunction]}{\muxmulevelsetvalue}$}}
			\put (-4,43) {$\nullhypthreearg{\muxmulevelsetvalue}{\leftu}{[\timefunction_0,\timefunction]}$}
			\put (10,44) {$\nullhypthreearg{\muxmulevelsetvalue}{\interestingu}{[\timefunction_0,\timefunction]}$}
			\put (57,57) {$\nullhypthreearg{\muxmulevelsetvalue}{-\interestingu}{[\timefunction_0,\timefunction]}$}
			\put (72,59) {$\nullhypthreearg{\muxmulevelsetvalue}{-\rightu}{[\timefunction_0,\timefunction]}$}
		\end{overpic}
		\caption{Regions in geometric coordinate space on which we derive estimates}
	\label{F:REGIONSWHEREWEDERIVEESTIMATES}
	\end{figure}
\end{center}

To derive $L^2$ estimates up to top order, we commute equations\footnote{More precisely, 
to avoid uncontrollable error terms,
we weight these equations by a factor of $\upmu$ before commuting.}
\eqref{E:WAVEINTRO}--\eqref{E:DIVCURLFORMODIFIEDINTRO}
with the elements of $\Tanset$ (see \eqref{E:INTROTANGENTCOMMUTATORS}) 
and derive energy and elliptic estimates on 
regions of the form $\twoargMrough{[\timefunction_0,\timefunction],[- \rightu,u]}{\muxmulevelsetvalue}$,
which are bounded by rough hypersurfaces on the top and bottom and null hypersurfaces on the sides;
see Fig\,\ref{F:REGIONSWHEREWEDERIVEESTIMATES}.

Our elliptic estimates, which are localized in spacetime, come in two flavors: those for the acoustic geometry,
as we mentioned in Sect.\,\ref{SSS:SMOOTHVSROUGH}, and elliptic estimates
for some of the fluid variables.
The localized elliptic estimates for the fluid variables 
are a difficult new feature of the present work,
and we discuss them in Sect.\,\ref{SSS:INTROFLUIDVARIABLEELLITPICESTIMATES};
here we discuss the energy estimates for the wave equations
\eqref{E:WAVEINTRO}, which in reality are coupled to the elliptic estimates.
More precisely, for solutions $\Psi$ to \eqref{E:WAVEINTRO},
we construct coercive energies $\mathbb{E}[\Psi](\timefunction,u)$
on rough hypersurface portions $\hypthreearg{\timefunction}{[- \rightu,u]}{\muxmulevelsetvalue}$
and null fluxes $\mathbb{F}[\Psi](\timefunction,u)$ on null hypersurface portions
$\nullhypthreearg{\muxmulevelsetvalue}{u}{[\timefunction_0,\timefunction)}$ 
(see Def.\,\ref{D:TRUNCATEDROUGHSUBSETS} for definitions of these portions),
and we derive suitable energy identities by using
a well-known framework based on the energy-momentum tensor for wave equations
and the well-chosen multiplier vectorfield
$\multipliervectorfield = (1 + 2 \upmu) \Lunit 
+ 
2 \muX$; see Prop.\,\ref{P:FUNDAMENTALENERGYNULLFLUXIDENTITIES}.
We note in passing that we derive similar energy estimates 
for the transport equations 
\eqref{E:TRANSPORTINTRO}--\eqref{E:TRANSPORTFORMODIFIEDINTRO},
and we will not discuss them in this introduction,
aside from mentioning that the $\gfour$-timelike
property $\gfour(\Transport,\Transport) = - 1$
(see Lemma~\ref{L:BASICPROPERTIESOFVECTORFIELDS}) is crucial for those estimates,
as it ensures that \emph{$\Transport$ is transversal to the characteristics}.

As we mentioned in Sect.\,\ref{SSS:REGULARESTIMATESONROUGHFOLIATIONS},
to control the top-order derivatives of 
the acoustic geometry without derivative loss,
we must rely on renormalizations, which we implement by constructing
appropriate ``modified quantities,''
originating in \cites{dCsK1993,sKiR2003,dC2007,jS2016b}; 
see Sect.\,\ref{S:CONSTRUCTIONOFMODIFIEDQUANTITIES}.
As in other works on shock waves without symmetry, this renormalization
leads to top-order energy identities for the wave variables $\Psi$ involving
singular terms. Because of the special properties of our rough foliations,
the singular terms are in fact factors of $\frac{1}{\timefunction}$,
and the resulting energy-flux inequalities can be caricatured as
follows, where the LHS is the sum of the energy, the null flux,
and the following coercive spacetime integral:
\begin{align} \label{E:INTROCOERCIVESPACETIME}
	\mathbb{K}[\tander^{\Ntop} \Psi](\timefunction,u)
	& \eqdef
	-
	\int_{(\timefunction',u',x^2,x^3) \in [\timefunction_0,\timefunction) \times [-\rightu,u] \times \mathbb{T}^2}
		\mathbf{1}_{\lbrace [-\interestingu,\interestingu] \rbrace}(u')
		(\argLrough{\muxmulevelsetvalue} \upmu) 
		|\angrmD \tander^{\Ntop} \Psi|_{\gtorus}^2 
	\, \mathrm{d} x^2 \mathrm{d} x^3 \mathrm{d} u' \mathrm{d} \timefunction'  
\end{align}
of the top-order derivatives of the wave variable $\Psi$ 
(see Def.\,\ref{D:WAVEANDTRANSPORTENERGIESANDNULLFLUXES} for the precise definitions):
\begin{align}  \label{E:INTROSCHEMATICTOPORDERENERGYID}
\begin{split}
	&
	\mathbb{E}[\tander^{\Ntop} \Psi](\timefunction,u)
	+
	\mathbb{F}[\tander^{\Ntop} \Psi](\timefunction,u)
	+
	\mathbb{K}[\tander^{\Ntop} \Psi](\timefunction,u)
		\\
	& 
	\leq 
	C \initialsmall^2
	+
	\boxed{A}
	\int_{\timefunction_0}^{\timefunction}
		\frac{1}{|\timefunction'|}
		\mathbb{E}[\tander^{\Ntop} \Psi](\timefunction',u)
	\, \mathrm{d} \timefunction' 
	+
	\cdots.
\end{split}
\end{align}
In \eqref{E:INTROSCHEMATICTOPORDERENERGYID}, 
$\initialsmall \geq 0$ is the small size of the perturbation
of the initial data away from a background plane-symmetric solution,
$\timefunction_0 < 0$ denotes the initial rough time (at which the data perturbation is assumed to be small),
and $\timefunction \in [\timefunction_0,0]$.
Moreover, $A > 0$ 
is a \emph{universal constant} 
(which we have placed in a box to highlight its importance)
that is independent of the equation of state,
and $\cdots$ denotes similar\footnote{In reality, in the energy estimates, there are other difficult error terms
that have to be handled with integration by parts with respect to $\frac{1}{\Lunit \timefunctionarg{\muxmulevelsetvalue}} \Lunit$, and these lead to
difficult singular boundary terms that we control in Lemma~\ref{L:DIFFICULTHYPERSURFACEIBPESTIMATE}.} 
or easier error terms or terms that can be handled with the 
help of the elliptic estimates described in Sect.\,\ref{SSS:INTROFLUIDVARIABLEELLITPICESTIMATES}.
In the spacetime integral $\mathbb{K}[\tander^{\Ntop} \Psi](\timefunction,u)$ on LHS~\eqref{E:INTROSCHEMATICTOPORDERENERGYID},
$|\angrmD \tander^{\Ntop} \Psi|_{\gtorus}^2$ denotes the square norm of the derivatives of
$\tander^{\Ntop} \Psi$ in directions tangent to the smooth tori 
$\ell_{t',u'} = \Sigma_{t'} \cap \nullhyparg{u'}$,
while
$\mathbf{1}_{\lbrace [-\interestingu,\interestingu] \rbrace}(u')$
is the characteristic function of a small interval $[-\interestingu,\interestingu]$ of $u$-values 
near the crease, and the factor $\argLrough{\muxmulevelsetvalue} \upmu$
is a null derivative of $\upmu$ that is \emph{quantitatively negative} in this interval.
Hence, in view of the minus sign in \eqref{E:INTROCOERCIVESPACETIME}, we see that
\emph{the overall sign of the spacetime integral on LHS~\eqref{E:INTROSCHEMATICTOPORDERENERGYID} is positive},
and this integral is crucially used in the proof to absorb some of the error terms in  ``$\cdots$'' 
on RHS~\eqref{E:INTROSCHEMATICTOPORDERENERGYID}.
We also note in passing that some of the error terms in ``$\cdots$'' are non-singular terms that can be handled
with the help of the null fluxes $\mathbb{F}[\tander^{\Ntop} \Psi](\timefunction,u)$,
and that some of these terms, when treated with Gr\"{o}nwall's inequality in $u$, allow for the possibility
of exponential growth in $u$, which is permissible within the scope of our approach ($u$ is confined to a compact set).
The key point is that by applying Gr\"{o}nwall's inequality to \eqref{E:INTROSCHEMATICTOPORDERENERGYID},
we deduce (ignoring the ``$\cdots$'' terms here) the following \emph{singular} (as $\timefunction \uparrow 0$) 
top-order energy estimate:
\begin{align} \label{E:INTROTOPORDERENERGYESTIMATE}
	\mathbb{E}[\tander^{\Ntop} \Psi](\timefunction,u)
	+
	\mathbb{F}[\tander^{\Ntop} \Psi](\timefunction,u)
	+
	\mathbb{K}[\tander^{\Ntop} \Psi](\timefunction,u)  
	& \lesssim \initialsmall^2 |\timefunction|^{-A}.
\end{align}
Thus, even though we have used geometrically defined commutators to try to turn
the problem of shock formation into a ``regular'' problem, a remnant of the singularity
survives\footnote{These singular top-order estimates are completely tied to the
energy estimates, and in particular, this difficulty does not arise in $1D$, where one can use transport estimates to close the problem.} at the top-derivative level, i.e, our ``unfolding'' of the characteristics
does not completely ``hide'' the singularity. 
This calls into question the basic philosophy of the approach.

However, two crucial structural features of the problem rescue the situation.
First, the universal constant $A$ -- and hence the singular factor $|\timefunction|^{-A}$ on RHS~\eqref{E:INTROTOPORDERENERGYESTIMATE} --
is \emph{independent} of $\Ntop$. This means that we can choose $\Ntop$ to be large without increasing the 
singularity strength. Second, below top-order, one can avoid using the renormalization procedure
to control the acoustic geometry. This leads to the loss of one derivative in the energy estimate hierarchy,
but the loss of one derivative is permissible below top-order. The price one pays is that
this procedure couples the below-top-order estimates to the singular top-order ones.
At one derivative below the top-order, the corresponding energy-null flux inequality
can be caricatured\footnote{In reality, the top-order and below-top-order energy estimates are highly coupled
and have to be derived simultaneously through an intricate version of Gr\"{o}nwall's lemma, 
which we provide in Sect.\,\ref{SSS:PROOFOFAPRIORIL2ESTIMATESWAVEVARIABLES}.} as follows:
\begin{align}  \label{E:INTROSCHEMATICBELOWTOPORDERENERGYID}
\begin{split}
	&
	\mathbb{E}[\tander^{\Ntop-1} \Psi](\timefunction,u)
	+
	\mathbb{F}[\tander^{\Ntop-1} \Psi](\timefunction,u)
	+
	\mathbb{K}[\tander^{\Ntop-1} \Psi](\timefunction,u)
	\\
& \leq 
	C \initialsmall^2
	+
	C
	\int_{\timefunction' = \timefunction_0}^{\timefunction} 
		\frac{1}{|\timefunction'|^{1/2}} \mathbb{E}^{1/2}[\tander^{\Ntop-1} \Psi](\timefunction',u) 
			\int_{\timefunction'' = \timefunction_0}^{\timefunction'} 
				\frac{1}{|\timefunction''|^{1/2}} 
				\mathbb{E}^{1/2}[\tander^{\Ntop} \Psi](\timefunction'',u) 
			\, \mathrm{d} \timefunction'' 
	\mathrm{d} \timefunction'
	+
	\cdots.
\end{split}
\end{align}
By applying Gr\"{o}nwall's inequality to \eqref{E:INTROSCHEMATICBELOWTOPORDERENERGYID}
and accounting for the singular top-order behavior stated in 
\eqref{E:INTROTOPORDERENERGYESTIMATE}, we deduce that:
\begin{align} \label{E:INTROBELOWTOPORDERENERGYESTIMATE}
	\mathbb{E}[\tander^{\Ntop-1} \Psi](\timefunction,u)
	+
	\mathbb{F}[\tander^{\Ntop-1} \Psi](\timefunction,u)
	+
	\mathbb{K}[\tander^{\Ntop-1} \Psi](\timefunction,u)
	& \lesssim 
	\initialsmall^2 |\timefunction|^{- (A-2)}.
\end{align}
\eqref{E:INTROBELOWTOPORDERENERGYESTIMATE} shows that if we descend one derivative level below top-order,
then the energy becomes less singular by a factor of $|\timefunction|^2$. 
As in \cites{dC2007,jS2016b,jLjS2018,jLjS2021},
one can continue the descent, proving that:
\begin{align} \label{E:INTROENERGYESTIMATEDESCENTTWOBELOWTOP}
\mathbb{E}[\tander^{\Ntop-2} \Psi](\timefunction,u)
	+
\mathbb{F}[\tander^{\Ntop-2} \Psi](\timefunction,u)
+
\mathbb{K}[\tander^{\Ntop-2} \Psi](\timefunction,u) 
&
\lesssim 
\initialsmall^2 |\timefunction|^{- (A-4)},
\end{align}
$\cdots$,
and finally arriving at the \emph{non-singular} estimate:
\begin{align} \label{E:INTRONONSINGULARENERGYESTIMATE}
\mathbb{E}[\tander^{\Ntop-\frac{A}{2}} \Psi](\timefunction,u)
	+
\mathbb{F}[\tander^{\Ntop-\frac{A}{2}} \Psi](\timefunction,u)
+
\mathbb{K}[\tander^{\Ntop-\frac{A}{2}} \Psi](\timefunction,u)
&
\lesssim 
\initialsmall^2.
\end{align}
The crucial non-singular energy estimate \eqref{E:INTRONONSINGULARENERGYESTIMATE} is what 
allows us to show, via Sobolev embedding,
that the solution is bounded with respect to the geometric coordinates $(t,u,x^2,x^3)$
at derivative levels $\approx \Ntop-\frac{A}{2}$ and below.

Finally, we note that a related energy estimate hierarchy holds for the other fluid variables
(including $\vortrenormalized$, $\GradEnt$, $\VortVort$, $\DivGradEnt$)
and the acoustic geometry (including quantities such as $\upmu$ and $\upchi$),
and that in practice, all these estimates are coupled
(though in certain spots in our bootstrap argument, we exploit weak coupling, which allows us to derive
some estimates before deriving others).
We refer to Sect.\,\ref{S:STATEMENTOFALLL2ESTIMATESANDBOOTSTRAPASSUMPTIONSFORWAVEENERGIES} for 
detailed statements of the full hierarchy of energy estimates.

\subsubsection{Localized elliptic fluid variable estimates via the characteristic current}
\label{SSS:INTROFLUIDVARIABLEELLITPICESTIMATES}
As we mentioned earlier, to close the top-order energy estimates
for the vorticity and entropy
on the rough domains
$
\twoargMrough{[\timefunction_0,\timefunction],[- \rightu,u']}{\muxmulevelsetvalue}
=
\lbrace \timefunction_0 \leq \timefunctionarg{\muxmulevelsetvalue} \leq \timefunction \rbrace
\cap
\lbrace - \rightu \leq u \leq u' \rbrace
$,
we cannot rely exclusively on the transport equation \eqref{E:TRANSPORTINTRO}; that approach
would lead to the loss of a derivative in our scheme.
Instead, we use the full structure of the transport-div-curl system
\eqref{E:TRANSPORTINTRO}--\eqref{E:DIVCURLFORMODIFIEDINTRO}.
We also have to show that our elliptic estimates, which are singular, 
are compatible with the
blowup-rates of the high order wave and transport energies, described in 
Sect.\,\ref{SSS:GEOMETRICENERGYESTIMATESONTHEROUGHFOLIATIONS}.
We used a version of this approach in \cite{jLjS2021}, where we followed the solution precisely to the constant-Cartesian-time
hypersurface of first blowup. The elliptic estimates in \cite{jLjS2021}
were much simpler because we derived them only on the flat hypersurfaces $\Sigma_t$, whose geometry is trivial,
and because we did not try to derive the localized structure of the singular boundary; this allowed us to close the proof
by deriving the elliptic estimates only across all of space, thereby 
(also exploiting the assumption of compactly supported data in \cite{jLjS2021}) avoiding the
difficult spatial boundary terms that we encounter in the present work.

To derive the desired estimates, in particular to handle the difficult boundary terms, 
we use the technology of \cite{lAjS2020},
which allows one to combine the Euclidean div-curl system
\eqref{E:TRANSPORTFORMODIFIEDINTRO}--\eqref{E:DIVCURLFORMODIFIEDINTRO}
with the transport equations \eqref{E:TRANSPORTINTRO}
to derive ``elliptic-hyperbolic identities'' that yield 
a sufficient amount of Sobolev regularity on \emph{any} region 
-- including $\twoargMrough{[\timefunction_0,\timefunction),[- \rightu,u']}{\muxmulevelsetvalue}$ --
that is globally hyperbolic with respect to $\gfour$.
See Prop.\,\ref{P:INTEGRALIDENTITYFORELLIPTICHYPERBOLICCURRENT} 
for the precise elliptic-hyperbolic identity that we use to close the top-order estimates.

There are three new ingredients in our elliptic estimates here compared to \cite{lAjS2020}:
\begin{itemize}
		\item The elliptic-hyperbolic identities depend, roughly, on certain
			curvature components of the boundaries of the domain.
			Some of these components become very singular as $\upmu \to 0$,
			and we need to ensure that the singularity strength
			is compatible with the blowup-rates of the high order wave and transport energies.
			In Prop.\,\ref{P:POINTWISEESTIMTAESFORELLIPTICHYPERBOLICIDENTITYERORTERMS},
			we provide pointwise estimates guaranteeing that indeed, 
			all of the error terms in the identities are controllable within the scope of our approach.
	\item We derive a new version of the identities from \cite{lAjS2020}
		based on applying the divergence theorem to a well-constructed 
		\emph{characteristic current} $\ehcurrent^{\alpha}$, 
		which is tangent to the characteristic hypersurfaces $\nullhyparg{u}$;
		see Def.\,\ref{D:PUTANGENTELLIPTICHYPERBOLICCURRENT}.
		Because of the $\nullhyparg{u}$-tangency, when we apply
		the divergence theorem to $\ehcurrent^{\alpha}$ on $\twoargMrough{[\timefunction_0,\timefunction],[- \rightu,u']}{\muxmulevelsetvalue}$,
		\emph{there are no boundary integrals along the lateral boundaries} $\nullhyparg{u}$. 
		This allows us to completely avoid
		error integrals on $\nullhyparg{u}$ that feature the top-order derivatives of $\upmu$;
		the point is that it is not possible to derive top-order $L^2$ estimates for $\upmu$ on $\nullhyparg{u}$
		because it satisfies a transport equation $\Lunit \upmu = \cdots$, where $\Lunit$ is \emph{tangent} to the $\nullhyparg{u}$.
		We also highlight that, unlike the currents in \cite{lAjS2020},
		our characteristic current here does not involve the future-directed normal $\hypunitnormalarg{\muxmulevelsetvalue}$
		to the constant-rough-time hypersurfaces $\hypthreearg{\timefunction}{[- \rightu,u']}{\muxmulevelsetvalue}$.
		Avoiding $\hypunitnormalarg{\muxmulevelsetvalue}$-dependent terms 
		is advantageous because the derivatives of
		$\hypunitnormalarg{\muxmulevelsetvalue}$ become very singular as $\upmu \to 0$,
		and it is not at all clear that such terms would have been compatible
		the blowup-rates of the high order wave and transport energies.
	\item Upon applying the divergence theorem to $\ehcurrent^{\alpha}$ on
		$\twoargMrough{[\timefunction_0,\timefunction),[- \rightu,u']}{\muxmulevelsetvalue}$,
		we generate boundary integrals along the top boundary portion
		$\hypthreearg{\timefunction}{[- \rightu,u']}{\muxmulevelsetvalue}$; this is a key difference from \cite{jLjS2021},
		where we integrated across all of space.
		To control some of these boundary
		integrals, we need to identify and integrate some perfect $\hypthreearg{\timefunction}{[- \rightu,u']}{\muxmulevelsetvalue}$-divergences,
		which leads to further co-dimension-two boundary integrals on the rough tori
		$\twoargroughtori{\timefunction,u'}{\muxmulevelsetvalue}$ and $\twoargroughtori{\timefunction,- \rightu}{\muxmulevelsetvalue}$.
		The key point is that, as in \cite{lAjS2020}, the rough tori boundary integrals either are controlled by the data or
		enter with a good sign,
		which is crucial for closing the elliptic estimates;
		see the tori integrals in the identity \eqref{E:INTEGRALIDENTITYFORELLIPTICHYPERBOLICCURRENT}.
		The analysis tied to these issues is very technical and is captured in divergence-form
		in Lemma~\ref{L:KEYIDPUTANGENTCURRENTCONTRACTEDAGAINSTVECTORFIELD}.
\end{itemize}

\subsection{Outline of the remainder of the paper}
\label{SS:PAPEROUTLINE}

\begin{itemize}
	\item In Sect.\,\ref{S:COMPRESSIBLEEULERFLOWANDITSGEOMETRIC}, we define the fluid variables that we use in our analysis and 
		recall the geometric formulation of the flow derived in \cite{jS2019c}.
	\item In Sects.\,\ref{S:ACOUSTICGEOMETRYANDCOMMUTATORVECTORFIELDS}--\ref{S:ROUGHACOUSTICALGEOMETRYANDCURVATURETENSORS},
		we derive basic properties and identities (not yet estimates) 
		tied to the eikonal function $u$, the rough time function $\timefunctionarg{\muxmulevelsetvalue}$,
		the corresponding geometries,
		and changes of variables between various coordinate systems and vectorfields.
	\item In Sect.\,\ref{S:ACOUSTICDOUBLENULLFRAMEANDITSRELATIONWITHTHEROUGHGEOMETRY}, we construct
		the ``ingoing'' $\gfour$-null vectorfield $\uLunit$, which is transversal to the characteristics $\nullhyparg{u}$ 
		and complements the $\gfour$-null vectorfield $\Lunit$, which generates the $\nullhyparg{u}$.
		We use $\uLunit$ in Sect.\,\ref{S:ELLIPTICHYPERBOLICIDENTITIES},
		when we derive the ``elliptic-hyperbolic'' identities
		that we will use to control the top-order derivatives of
		$\vortrenormalized$ and $\GradEnt$.
	\item In Sect.\,\ref{S:NORMSANDFORMS}, we define various norms, area forms, and volume forms.
		We also introduce notation for various strings of commutation vectorfields.
	\item In Sect.\,\ref{S:SCHEMATICSTRUCTUREANDIDENTITIES}, we introduce schematic notation and derive
		various identities in schematic form. These will be used throughout the remainder of the paper.
	\item In Sects.\,\ref{S:PARAMETERSANDSIZEASSUMPTIONSANDCONVENTIONSFORCONSTANTS}--\ref{S:ASSUMPTIONSONTHEDATA},
		we list various parameters corresponding to the solutions under study,
		state our assumptions on the initial data,
	\item In Sect.\,\ref{S:BOOTSTRAPEVERYTHINGEXCEPTENERGIES}, 
		we state all our bootstrap assumptions -- except for the energy bootstrap assumptions.
	\item In Sect.\,\ref{S:PRELIMINARYPOINTWISECOMMUTATORANDOPERATORCOMPARISON},
		we use the bootstrap assumptions to derive preliminary pointwise, commutator, and differential operator comparison estimates.
	\item In Sect.\,\ref{S:EMBEDDINGSANDFLOWMAPS},
		we use the bootstrap assumptions to analyze 
		the data hypersurface $\datahypfortimefunctiontwoarg{-\muxmulevelsetvalue}{[\timefunction_0,\timefunctionboot)}$
		for the rough time function $\timefunctionarg{\muxmulevelsetvalue}$,
		which solves the transport $\Wtransarg{\muxmulevelsetvalue} \timefunctionarg{\muxmulevelsetvalue} = 0$
		with data prescribed on $\datahypfortimefunctiontwoarg{-\muxmulevelsetvalue}{[\timefunction_0,\timefunctionboot)}$.
		We also derive basic properties of the flow map of $\Wtransarg{\muxmulevelsetvalue}$.
	\item In Sect.\,\ref{S:ESTIMATESFORROUGHTIMEFUNCTIONCONTINUOUSEXTENSIONSANDDIFFEOS},
		we use the bootstrap assumptions to derive estimates for
		$\timefunctionarg{\muxmulevelsetvalue}$. We also show that various quantities extend to the
		closure of the bootstrap region as elements of various H\"{o}lder spaces.
		Finally, we study the properties of the map
		$\CHOVroughtomumuxmu{\muxmulevelsetvalue}(\timefunction,u,x^2,x^3) = (\upmu,\muX \upmu,x^2,x^3)$,
		which is important for understanding the structure of the singular boundary.
	\item In Sect.\,\ref{S:ESTIMATESTIEDTOTHEFLOWMAPOFROUGHNULLVECTORFIELD},
		we use the bootstrap assumptions to control the flow map of the $\gfour$-null vectorfield
		$\argLrough{\muxmulevelsetvalue} = \frac{1}{\Lunit \timefunctionarg{\muxmulevelsetvalue}} \Lunit$,
		which is the principal operator in many of the transport equations that we will later study.
		construct $\timefunctionarg{\muxmulevelsetvalue}$,
		derive estimates for $\timefunctionarg{\muxmulevelsetvalue}$ and its geometry,
		and study the diffeomorphism properties of the change of variables map
		$(t,u,x^2,x^3) \rightarrow (\timefunctionarg{\muxmulevelsetvalue},u,x^2,x^3)$.
	\item In Sect.\,\ref{S:LINFINITYFLUIDANDEIKONALANDIMPROVEMENTOFAUX},
		we derive $L^{\infty}$ estimates that yield improvements of many of 
		our quantitative bootstrap assumptions.
	\item In Sect.\,\ref{S:SHARPCONTROLOFMUANDPROPERTIESOFCHOVGEOTOCARTESIAN}, we derive
		sharp estimates for $\upmu$. These are crucial for the energy estimates and for understanding
		the structure of the singular boundary.
		We also study the homeomorphism and diffeomorphism properties of the change of variables map
		$(t,u,x^2,x^3) \rightarrow (t,x^1,x^2,x^3)$.
	\item In Sect.\,\ref{S:CONSTRUCTIONOFMODIFIEDQUANTITIES}, we construct the modified quantities
		that we use to control the acoustic geometry in $L^2$ without derivative loss. 
	\item In Sect.\ref{S:BASICINGREDIENTSFORL2ANALYSIS},
		we construct some basic ingredients needed for the hyperbolic energy estimates.
		In particular, we define energies and null-fluxes, derive energy--null-flux identities for solutions
		to wave equations and transport equations, and exhibit the coerciveness of the energies
		and null fluxes.
	\item In Sect.\,\ref{S:ELLIPTICHYPERBOLICIDENTITIES},
		we derive the ``elliptic-hyperbolic'' identities
		that we will use to control the top-order derivatives of
		$\vortrenormalized$ and $\GradEnt$.
	\item In Sect.\,\ref{S:POINTWISESTIMATESFORWAVEEQUATIONS}, we commute the wave equations
		satisfied by the wave variables  $\wavearray = (\RRiemann,\LRiemann,v^2,v^3,\Ent)$ up to top-order and
		derive pointwise estimates for the inhomogeneous terms.
		These are a preliminary ingredient
		in our $L^2$ analysis of the inhomogeneous terms.
	\item In Sect.\,\ref{S:POINTWISEESTIMATESFORCONTROLLINGSPECIFICVORTICITYANDENTROPYGRADIENT},
		we provide an analog of Sect.\,\ref{S:POINTWISESTIMATESFORWAVEEQUATIONS} for the transport variables.
		That is, we commute the transport equations satisfied by
		$\vortrenormalized$, 
		$\GradEnt$, 
		$\VortVort$, and $\DivGradEnt$
		up to top-order and
		derive pointwise estimates for the inhomogeneous terms.
		We also derive pointwise estimates at the low derivative levels for the rough acoustical
		geometry. All the estimates in this section
		are preliminary ingredients for our derivation of
		energy and elliptic estimates.
	\item In Sect.\,\ref{S:STATEMENTOFALLL2ESTIMATESANDBOOTSTRAPASSUMPTIONSFORWAVEENERGIES},
		we state all of our a priori energy estimates.
		We also state bootstrap assumptions for the energies of the 
		``wave variables'' $\wavearray \eqdef (\RRiemann,\LRiemann,v^2,v^3,\Ent)$.
		The proof of the energy estimates occupies a substantial portion of the remainder of the paper, 
		all the way through Sect.\,\ref{S:WAVEANDACOUSTICGEOMETRYAPRIORIESTIMATES}.
	\item In Sect.\,\ref{S:PRELIMINARYL2ESTIMATESFORBELOWTOPORDERDERIVATIVESOFACOUSTICGEOMETRYANDDERIVATIVELOSING},
		we derive preliminary $L^2$ estimates for the below-top-order derivatives of the
		eikonal function quantities $\upmu$, $\Lunit^i$, $\upchi$, and $\mytr_{\gtorus}\upchi$.
		We also derive preliminary $L^2$ for $\wavearray$ that lose one derivative.
	\item In Sect.\,\ref{S:BELOWTOPORDERHYPERBOLICL2ESTIMATESFORSPECIFICVORTICITYANDENTROPYGRADIENT}, 
		we use the wave energy bootstrap assumptions to derive below-top-order energy estimates for 
		$\vortrenormalized$, 
		$\GradEnt$, 
		$\VortVort$, 
		and $\DivGradEnt$.
	\item In Sect.\,\ref{S:TOPORDERELIPTICHYPERBOLICL2ESTIMATESFORSPECIFICVORTICITYANDENTROPYGRADIENT},
		we use the wave energy bootstrap assumptions and the below-top-order energy estimates of
		Sect.\,\ref{S:BELOWTOPORDERHYPERBOLICL2ESTIMATESFORSPECIFICVORTICITYANDENTROPYGRADIENT}
		to derive the top-order ``elliptic-hyperbolic'' energy estimates
		$\vortrenormalized$, 
		$\GradEnt$, 
		$\VortVort$, 
		and $\DivGradEnt$.
	\item In Sect.\,\ref{S:ELLIPTICESTIAMTESACOUSTICGEOMETRYONROUGHTORI}, we derive general elliptic estimates on the rough
			tori, which we will later use to control the top-order derivatives of the acoustic geometry
			along the rough foliations.
	\item In Sect.\,\ref{S:WAVEANDACOUSTICGEOMETRYAPRIORIESTIMATES}, we derive up-to-top-order
		energy estimates for the wave variables and the acoustic geometry. 
		This completes the proof of the energy estimates
		stated in Sect.\,\ref{S:STATEMENTOFALLL2ESTIMATESANDBOOTSTRAPASSUMPTIONSFORWAVEENERGIES}
		and in particular yields a strict improvement of our wave energy bootstrap assumptions.
	\item In Sect.\,\ref{S:IMPROVEMENTSOFFUNDAMENTALQUANTITATIVEBOOTSTRAPASSUMPTIONS}, 
		we use the energy estimates and Sobolev embeddings to derive $L^{\infty}$ estimates
		that yield strict improvements of the remaining quantitative bootstrap assumptions.
		This closes the bootstrap argument and completes our proof of a priori estimates.
	\item In Sect.\,\ref{S:EXISTENCEUPTOSINGULARBOUNDARYATFIXEDKAPPA}, we use the a priori estimates
		and a continuation principle to show that we can extend the solution all the way up to the level 
		set $\lbrace \timefunctionarg{\muxmulevelsetvalue} = 0 \rbrace$, which contains 
		the two-dimensional torus $\twoargmumuxtorus{0}{-\muxmulevelsetvalue}$, which in turn is contained in
		the singular boundary. 
		We provide these results as Theorem~\ref{T:EXISTENCEUPTOTHESINGULARBOUNDARYATFIXEDKAPPA},
		which is the first main theorem of the paper. 
		This theorem provides a development of the data
		containing the portion of the singular boundary that is ``accessible'' via the foliation of spacetime
		by the level sets of $\timefunctionarg{\muxmulevelsetvalue}$.
	\item In Sect.\,\ref{S:DEVELOPMENTSOFTHEDATASINGULARBOUNDARYNEWTIMEFUNCTION}, we study the union
		of the developments as $\muxmulevelsetvalue$ varies, and we define an interesting sub-region,
		$\MInteresting$, which contains
		a portion of the singular boundary, namely $\mathcal{B}^{[0,\muxmulevelsetvalue_0]}$, and its past boundary 
		$\partial_- \mathcal{B}^{[0,\muxmulevelsetvalue_0]}$, the crease.
		We also construct a new rough time function
		$\newtimefunction$, whose level sets foliate $\MInteresting$.
		Finally, we derive various quantitative and qualitative properties of
		various geometric objects tied to $\MInteresting$
		and $\newtimefunction$.
	\item In Sect.\,\ref{S:HOMEOANDDIFFEOANDSINGULARBOUNDARYINCARTESIAN}, we study 
		the homeomorphism and diffeomorphism properties of the change of variables map
		$\Upsilon(t,u,x^2,x^3) = (t,x^1,x^2,x^3)$ on $\MInteresting$.
		We also exhibit the properties of $\Upsilon(\mathcal{B}^{[0,\muxmulevelsetvalue_0]})$,
		i.e., the embedding of the singular boundary in Cartesian coordinate space.
	\item In Sect.\,\ref{S:MAINRESULTS}, we state and prove Theorem~\ref{T:DEVELOPMENTANDSTRUCTUREOFSINGULARBOUNDARY},
		which is the main result of the paper. The theorem shows that
		$\MInteresting$ contains the portion $\mathcal{B}^{[0,\muxmulevelsetvalue_0]}$ 
		of the singular boundary and the crease, and it gives a detailed description of the solution
		in the different coordinate systems as well as the change of variables maps.
		The theorem is essentially a conglomeration of results derived 
		earlier in the paper.
\end{itemize}

 \section{Basic setup, compressible Euler flow, and its geometric reformulation} 
\label{S:COMPRESSIBLEEULERFLOWANDITSGEOMETRIC}
In this section, we first introduce some basic notational conventions and definitions.
We then provide a standard first-order quasilinear hyperbolic formulation of compressible Euler flow.
Next, we define a series of additional fluid variables and geometric tensors associated to the flow.
Finally, we recall the new formulation of the flow derived in \cite{jS2019c}. More precisely,
we use a slightly modified version of the formulation in \cite{jS2019c} that is adapted to the nearly
plane-symmetric solutions featured in our main results. The only difference with \cite{jS2019c}
is that here, we replace the density and  
velocity component $v^1$ with our ``almost Riemann invariants'' 
$\RRiemann$ and $\LRiemann$, which are useful
for capturing smallness in the regime under study.
In total, the new formulation comprises covariant wave equations, which govern the propagation of sound waves,
coupled to systems of transport-div-curl systems, which drive the evolution of the vorticity and entropy.

\subsection{Basic notation and conventions}
\label{SS:NOTATIONANDCONVENTIONS}
The precise definitions of some of the concepts referred to
here are provided later in the article.

\begin{itemize}
	\item  (\textbf{Cartesian coordinates}) 
	Our analysis takes place on subsets of the spacetime manifolds 
	$\mathbb{R} \times \Sigma$,
	where $\Sigma \eqdef \mathbb{R} \times \mathbb{T}^2$ is the spatial manifold. 
	We fix a standard Cartesian coordinate system $\lbrace x^{\alpha} \rbrace_{\alpha = 0,1,2,3}$
on $\mathbb{R} \times \Sigma$, where $t \eqdef x^0 \in \mathbb{R}$
is the time coordinate and $(x^1,x^2,x^3) \in \Sigma$
are the spatial coordinates (where $(x^2,x^3)$ are standard coordinates\footnote{While the coordinates $x^2,x^3$ on $\mathbb{T}^2$ are only locally defined, the corresponding partial derivative vectorfields $\partial_2,\partial_3$ can be extended
so as to form a global smooth frame on $\mathbb{T}^2$. Similar remarks apply to the one-forms $dx^2,dx^3$
These simple observations are relevant for this paper because when we derive estimates, the coordinate functions $x^2,x^3$
themselves are never directly relevant; what matters are estimates for the components of various tensorfields
with respect to the frame $\lbrace \partial_t, \partial_1, \partial_2, \partial_3 \rbrace$
and the basis dual co-frame $\lbrace dt, dx^1, dx^2, dx^3 \rbrace$, which are everywhere smooth. 
\label{FN:COORDINATESARENOTGLOBAL}.}  
on $\mathbb{T}^2$). By a \emph{plane-symmetric} solution,
we mean one whose fluid variables are independent of $(x^2,x^3)$ in this coordinate system. 
We sometimes refer to $t$ as the Cartesian time function.'' 
$\Sigma_{t'} \eqdef \lbrace (t,x^1,x^2,x^3) \in \R \times \R \times \T^2 | \ t \equiv t' \rbrace$
denotes the standard flat hypersurface of constant Cartesian time.
\item (\textbf{Cartesian coordinate partial derivatives}) 
We use the notation $\lbrace \partial_{\alpha} \rbrace_{\alpha = 0,1,2,3}$ 
(or $\partial_t \eqdef \partial_0$) to denote the Cartesian coordinate partial derivative
vectorfields. 
	\item (\textbf{Lowercase Greek index conventions}) Lowercase Greek spacetime indices 
	$\alpha$, $\beta$, etc.\ correspond to the Cartesian coordinate spacetime coordinates
	and vary over $0,1,2,3$. All lowercase Greek indices are lowered and raised with the acoustical metric
	$\gfour$ (see definition~\eqref{E:ACOUSTICALMETRIC}) 
	and its inverse $\gfour^{-1}$, and \emph{not with the Minkowski metric}.
	Throughout the article, if $\upxi$ is a type $\binom{m}{n}$ spacetime tensorfield,
	then unless we indicate otherwise, in our identities and estimates,
	\textbf{$\left\lbrace \upxi_{\beta_1 \cdots \beta_n}^{\alpha_1 \cdots \alpha_m} \right\rbrace_{
	\alpha_1, \cdots \alpha_m, \beta_1, \cdots, \beta_n = 0,1,2,3}$ denotes its 
	components with respect to the Cartesian coordinates}. 
	This is important because some of our identities
	and estimates hold only with respect to the Cartesian coordinates.
	\item (\textbf{Lowercase Latin index conventions})
	Lowercase Latin spatial indices $a$, $b$, etc.\ correspond to the Cartesian spatial coordinates and vary over $1,2,3$.
	Much like in the previous point, 
	if $\upxi$ is a type $\binom{m}{n}$ $\Sigma_t$-tangent tensorfield 
	(see Def.\,\ref{D:PROJECTIONTENSORFIELDSANDTANGENCYTOHYPERSURFACES}),
	then unless we indicate otherwise, in our identities and estimates,
	\textbf{$\left\lbrace \upxi_{b_1 \cdots b_n}^{a_1 \cdots a_m} \right\rbrace_{a_1,\cdots,a_m,b_1,\cdots,b_n = 1,2,3}$ 
	denotes its components with respect to the Cartesian spatial coordinates}. 
	\item (\textbf{Uppercase Latin index conventions})
	Uppercase Latin spatial indices $A$, $B$, etc.\ correspond to the coordinates 
	$(x^2,x^3)$ on $\mathbb{T}^2$ and vary over $2,3$. In particular, if $V$ is a vectorfield,
	and $A \in \lbrace 2,3 \rbrace$, then $V^A = V^{\alpha} \partial_{\alpha} x^A$, 
	where $(x^2,x^3)$ are the standard Cartesian coordinates on $\mathbb{T}^2$.
\item (\textbf{Einstein summation}) 
	We use Einstein's summation convention in that repeated indices are summed,
	e.g., $\Lunit^A X^A \eqdef \Lunit^2 X^2 + \Lunit^3 X^3$.
\item (\textbf{Use of ``$\cdot$''})  
		We sometimes use ``$\cdot$'' to denote the natural contraction between two tensors. 
		For example, if $\upxi$ is a spacetime one-form and $V$ is a 
		spacetime vectorfield,
		then $\upxi \cdot V \eqdef \upxi_{\alpha} V^{\alpha}$.
		At other times, we use ``$\cdot$'' to schematically denote products,
		e.g., $A_1 \cdot A_2 \cdot A_3$ is a trilinear form in $A_1, A_2, A_3$.
\item (\textbf{Tensor contractions})  
	If $V$ and $W$ are vectorfields, then $V_W \eqdef V^{\alpha} W_{\alpha} = \gfour_{\alpha \beta} V^{\alpha} W^{\beta}$.
	If $\upxi$ is a one-form and $V$ is a vectorfield, then $\upxi_V \eqdef \upxi_{\alpha} V^{\alpha}$.
	We use similar notation when contracting higher-order tensorfields against vectorfields.
	For example, if $\upxi$ is a type $\binom{0}{2}$ tensorfield and
	$V$ and $W$ are vectorfields, then $\upxi_{VW} \eqdef \upxi_{\alpha \beta} V^{\alpha} W^{\beta}$.
\item (\textbf{Commutator of operators})  If $Q_1$ and $Q_2$ are two operators, then $[Q_1,Q_2] \eqdef Q_1 Q_2 - Q_2 Q_1$ denotes their commutator. 
\item (\textbf{Constants}) We establish conventions for constants in Sect.\,\ref{SS:CONVENTIONSFORCONSTANTS}.
\end{itemize}

\subsection{Basic differential operators}
\label{SS:BASICDIFFERENTIALOPERATORS}
In our analysis, we will encounter many kinds of differential operators.
Here, we define some basic operators.

\begin{definition}[Gradient one-form of a scalar function]
	\label{D:GRADIENTONEFORMOFSCALARFUNCTION}
	If $f$ is a scalar function, then $\rmD f$
	denotes the gradient one-form associated to $f$,
	e.g., $(\rmD f)_{\alpha} \eqdef \rmD f \cdot \partial_{\alpha} = \partial_{\alpha} f$.
\end{definition}	

\begin{definition}[Vectorfield derivative of scalar functions]
	\label{D:VECTORFIELDDERIVATIVEOFSCALARFUNCTION}
	If $V$ is a vectorfield and $f$ is a scalar function, then 
	$V f \eqdef V^{\alpha} \partial_{\alpha} f = V \cdot \rmD f$
	denotes the derivative of $f$ in the direction $V$.
\end{definition}

\begin{definition}[Euclidean divergence and curl]
\label{D:EUCLIDEANDIVERGENCEANDCURL}
$\Flatdiv$ and $\Flatcurl$ respectively denote the Euclidean spatial divergence and curl operators. 
That is, given a $\Sigma_t$-tangent 
vectorfield $V = V^a \partial_a$, we define, 
relative to the Cartesian spatial coordinates, 
$\Flatdiv V$ and $\Flatcurl V$ to be the following
scalar function and $\Sigma_t$-tangent vectorfield:
\begin{align} \label{E:FLATDIVANDCURL}
		\Flatdiv V
		& \eqdef \partial_a V^a,
		&
		(\Flatcurl V)^i
		& \eqdef \upepsilon_{iab} \partial_a V^b,
	\end{align}
	where $\upepsilon_{iab}$ is the fully antisymmetric symbol normalized by $\upepsilon_{123} = 1$.
\end{definition}

\subsection{A first-order formulation involving the logarithmic density}
\label{SS:FIRSTORDERFORMULATIONWITHLOGDENSITY}

\subsubsection{The logarithmic density, assumptions on the equation of state, and normalizations}
\label{SSS:LOGDENSITYASSUMPTIONSONEOSANDNORMALIZATIONS}
We find it convenient to work with the logarithmic density featured in the next
definition, rather than the density.
In the rest of the paper,
	\begin{align} \label{E:BACKGROUNDDENSITY}
		\overline{\varrho} & > 0
	\end{align}
	denotes a fixed constant ``background density.''
	
\begin{definition}[Logarithmic density] 
\label{D:LOGDENS}
 We define the \emph{logarithmic density} $\LogDensity$ as follows:
\begin{align} \label{E:LOGDENS}
	\LogDensity 
	& \eqdef \ln\left(\frac{\varrho}{\overline\varrho}\right).
	\end{align}
\end{definition}

In the rest of the paper,
we view the speed of sound $\Speed$ (which is defined in \eqref{E:SOUNDSPEED}) 
to be a function of $(\LogDensity,\Ent)$.
Note that by \eqref{E:SOUNDSPEED} and the chain rule, 
we have $\Speed(\LogDensity,\Ent)
= \sqrt{(\overline{\varrho})^{-1} \exp(-\LogDensity) p;_{\LogDensity}}$,
where $p;_{\LogDensity} \eqdef \tfrac{\p p}{\p \LogDensity}$ denotes the partial derivative of the equation of state 
with respect to the logarithmic density at fixed $\Ent$. 

\begin{notation}[Partial differentiation with respect to state-space variables]
\label{N:PARTIALDIFFERENTIATIONWITHRESPECTTOSTATESPACEVARIABLES}
In accordance with the above notation, 
for any scalar function $f = f(\LogDensity,\Ent)$, we use the notation 
$f_{;\Ent} \eqdef \frac{\partial f}{\partial \LogDensity}$ to denote
the partial derivative of $f$ 
with respect to the logarithmic density at fixed $\Ent$.
Similarly, we denote the partial derivative of $f$ with respect to $\Ent$ at fixed $\LogDensity$
by $f_{;\Ent} \eqdef \frac{\partial f}{\partial \Ent}$.
We also write $ f_{;\LogDensity;\Ent} \eqdef \frac{\partial^2 f}{\partial \Ent \partial \LogDensity}$ and use similar notation for other higher partial derivatives of $f$ with respect to $\LogDensity, \Ent$.
\end{notation}

To ensure that shocks occur for solutions near static isentropic fluid states with constant density $\overline{\varrho} > 0$,
we assume the following non-degeneracy condition:
\begin{align}\label{E:NONDEGENCONDITION}
\overline{\Speed}^{-1} \overline{\Speed_{;\uprho}} + 1 
& \neq 0,
\end{align}
where LHS~\eqref{E:NONDEGENCONDITION} is defined to be the constant
obtained by evaluating $\Speed^{-1} \Speed_{;\LogDensity} + 1$ at $\LogDensity = \Ent \equiv 0$.
Equation \eqref{E:NONDEGENCONDITION} ensures that the null condition fails to hold for perturbations
of the background solution $\LogDensity = \Ent \equiv 0$; see Sect.\,\ref{SS:THEFACTORDRIVINGTHESHOCKFORMATION}.
Our main results hold for all equations of state except for that of the Chaplygin gas,
namely $p = C_0 - C_1 \exp(- \LogDensity)$, where $C_0 \geq 0$ and $C_1 > 0$ are constants.
This equation of state is degenerate in the following sense: $\Speed^{-1} \Speed_{;\LogDensity} + 1 \equiv 0$.

By rescaling Cartesian time if necessary, we can assume the following convenient normalization condition:
\begin{align} \label{E:BACKGROUNDSOUNDSPEEDISUNITY}
	\Speed(\LogDensity = 0, \Ent = 0) 
	& = 1.
\end{align}

\subsubsection{A first-order formulation involving the logarithmic density}
\label{SSS:FIRSTORDERFORMULATIONINVOLVINGLOGDENSITY}
From definition~\eqref{E:LOGDENS} and equations \eqref{E:INTROTRANSPORTVI}--\eqref{E:INTROBS}, 
it follows that
relative to the standard Cartesian coordinates on $\mathbb{R} \times \Sigma$,
the compressible Euler equations can be expressed as the following system
in $\LogDensity$, $v$, and $\Ent$:
 \begin{subequations}
\begin{align}
	\Transport v^i 
	& = 
	- \Speed^2 \updelta^{ia} \partial_a \LogDensity
	- \exp(-\LogDensity) \frac{p_{;\Ent}}{\overline{\varrho}} \updelta^{ia} \partial_a \Ent,
	\label{E:BVIEVOLUTION}
		\\
		 \label{E:BLOGDENSITYEVOLUTION}
	\Transport \LogDensity
	& = - \Flatdiv v,
		\\
	\Transport \Ent
	& = 0.
	\label{E:BENTROPYEVOLUTION}
\end{align}
\end{subequations}

\subsection{The almost Riemann invariants}
\label{SS:ALMOSTRIEMANNINVARIANTSANDLFLUIDVARIABLEARRAYS}
To study solutions close to simple isentropic plane-symmetric solutions, 
we find it convenient to replace
$\LogDensity$ and $v^1$ with a pair of ``almost Riemann invariants,''
denoted by $\RRiemann$ and $\LRiemann$.
In this paper, exact isentropic plane-symmetric simple solutions are such that 
$\RRiemann$ is a function of only $(t,x^1)$
and
$\LRiemann = \Ent = v^2 = v^3 \equiv 0$.


\begin{definition}[The almost Riemann invariants] \label{D:ALMOSTRIEMANNINVARIANTS}
We define the \emph{almost Riemann invariants\footnote{Compare $\almostRiemann_{(\pm)}$ with the true Riemann invariants for the plane-symmetric solutions given by \eqref{AE:RIEMANNINVARIANTS}.} away from symmetry} $\almostRiemann_{(\pm)}$ as follows:
\begin{align} \label{E:ALMOSTRIEMANNINVARIANTS}
\almostRiemann_{(\pm)} 
& \eqdef 
v^1 
\pm 
\almostRiemannfunction(\LogDensity,\Ent),
& 
\mbox{where } 
\almostRiemannfunction(\LogDensity,\Ent)
&
\eqdef 
\int_0^{\LogDensity} 
	\Speed(\LogDensity',\Ent) 
\, \rmD \LogDensity'.
\end{align}

\end{definition}

\begin{remark}[\textbf{Clarification on our approach to estimating $\LogDensity$ and $v^1$}]
	\label{R:HOWWEESTIMATEDENSITYANDV1}
	We have introduced $\almostRiemann_{(\pm)}$ because they are convenient for
	studying perturbations of simple isentropic plane-waves (for which only $\RRiemann$ is non-vanishing);
	$\almostRiemann_{(\pm)}$ allow us to capture various kinds of smallness of the
	perturbations. 
	It is well-known that for isentropic plane-symmetric solutions, 
	one can use $\lbrace \RRiemann, \LRiemann \rbrace$ as
	the unknowns in place of $\lbrace \LogDensity, v^1 \rbrace$;
	see Appendix~\ref{A:PS}.
	Away from symmetry, a similar remark also holds for our almost Riemann invariants,
	provided we take into account the entropy.
	Specifically, from \eqref{E:BACKGROUNDSOUNDSPEEDISUNITY} and definition~\eqref{E:ALMOSTRIEMANNINVARIANTS},
	it follows that $v^1 = \frac{1}{2}(\RRiemann - \LRiemann)$,
	and that when $\LogDensity$, $v^1$, and $\Ent$ are sufficiently small
	(as is captured by the smallness parameters $\mathring{\upalpha}$ and $\initialsmall$ that we
	introduce in Sect.\,\ref{S:PARAMETERSANDSIZEASSUMPTIONSANDCONVENTIONSFORCONSTANTS}),
	we have (via the implicit function theorem)
	$\LogDensity = 
	(\RRiemann - \LRiemann)
	\cdot
	\widetilde{F}(\RRiemann - \LRiemann,\Ent)$,
	where $\widetilde{F}$ is a smooth function.
	This allows us to control $\LogDensity$ and $v^1$ in terms of 
	$\RRiemann$, $\LRiemann$, and $\Ent$. 
	Throughout the article, we use this observation without explicitly pointing it out.
	In particular, even though many of the equations that we study explicitly involve
	$\LogDensity$ and $v^1$, it should be understood that we always
	estimate these quantities in terms of the ``wave variables''
	$\RRiemann$, $\LRiemann$, and $\Ent$,
	which are featured in the array \eqref{E:ARRAYOFWAVEVARIABLES}
	defined below.
\end{remark}


\subsection{The higher order fluid variables}
\label{SS:HIGHERORDERFLUIDVARIABLES}
The ``higher order'' fluid variables in the next definition appear in Theorem~\ref{T:GEOMETRICWAVETRANSPORTSYSTEM},
which provides the formulation of compressible Euler flow that we use throughout our analysis.

\begin{definition}[The higher order fluid variables] \label{D:HIGHERORDERANDMODFLUIDVARSDEF} \hfill
\label{D:HIGHERORDERFLUIDVARIABLES}
\begin{enumerate}
\item We define the \emph{specific vorticity} to be the $\Sigma_t$-tangent vectorfield whose Cartesian spatial components are:
\begin{align} \label{E:SPECIFICVORTICITYDEF}
	\vortrenormalized^i 
	& \eqdef 
	\frac{(\Flatcurl v)^i}{\exp(\LogDensity)} 
	= 
	\frac{\upepsilon_{ijk} \updelta^{jl}\p_lv^k}{\exp(\LogDensity)},
\end{align}
where $\updelta^{jl}$ is the Kronecker delta.
\item We define the \emph{entropy gradient} to be the $\Sigma_t$-tangent vectorfield whose Cartesian spatial components are: 
	\begin{align} \label{E:GRADENTDEF}
		\GradEnt^i 
		& \eqdef \updelta^{ia}\p_a \Ent = \p_i \Ent.
	\end{align}
\item We define the \emph{modified fluid variables} to be the $\Sigma_t$-tangent vectorfield $\VortVort$ and the scalar function 
$\DivGradEnt$ whose Cartesian spatial components are:
\begin{subequations}
	\begin{align} \label{E:MODIFIEDCURLOFVORTICITY}
		\VortVort^i
		& \eqdef
			\exp(-\LogDensity) (\Flatcurl \vortrenormalized)^i
			+
			\exp(-3\LogDensity) \Speed^{-2} \frac{p_{;\Ent}}{\overline{\varrho}} \GradEnt^a \partial_a v^i
			-
			\exp(-3\LogDensity) \Speed^{-2} \frac{p_{;\Ent}}{\overline{\varrho}} (\Flatdiv v) \GradEnt^i,
				\\
		\DivGradEnt 
		& \eqdef 
			\exp(-2 \LogDensity) \Flatdiv \GradEnt 
			-
			\exp(-2 \LogDensity) \GradEnt^a \partial_a \LogDensity.
			\label{E:MODIFIEDDIVERGENCEOFENTROPYGRADIENT}
\end{align}
\end{subequations}
\end{enumerate}
\end{definition}

\subsection{Arrays of fluid variables and array norm notation}
\label{SS:FLUIDVARIABLEARRAYSANDNORMNOTATION}
We provide the next definition for notational convenience.

\begin{definition}[The fluid variable array $\wavearray$ and the partial array $\wavearraypartial$]
\label{D:ARRAYSOFWAVEVARIABLES} 
We define the \emph{array of wave\footnote{These ``wave variables'' solve wave equations; 
see Theorem~\ref{T:GEOMETRICWAVETRANSPORTSYSTEM}.} variables} as follows:
\begin{subequations}
\begin{align} \label{E:ARRAYOFWAVEVARIABLES} 
\wavearray 
& \eqdef (\Psi_0,\Psi_1,\Psi_2,\Psi_3,\Psi_4) 
\eqdef (\RRiemann,\LRiemann,v^2,v^3,\Ent).
\end{align}

We define the \emph{partial array of wave variables} by:
\begin{align} \label{E:PARTIALWAVEARRAY}
\wavearraypartial 
& \eqdef (\LRiemann,v^2,v^3,\Ent).
\end{align}  
\end{subequations}

\end{definition} 
We view $\wavearray$ to be an array of scalar functions $\Psi_{\iota}$, where 
$\iota = 0,\cdots,4$. We will not attribute any tensorial structure to the labeling index $\iota$ besides simple contractions, denoted by $\diamond$, corresponding to the chain rule; see Def.\,\ref{D:DERIVATIVESOFARRAYS}.

In the next definition, we introduce notation for norms of arrays.
\begin{definition}[Conventions with variable arrays] \label{D:CONVENTIONESFORDIFFERENTIATION} \hfill
\begin{itemize}
\item Given the fluid variable array $\wavearray$ from Def.\,\ref{D:ARRAYSOFWAVEVARIABLES}, we define:
\begin{align} \label{E:ABSOLUTVALUEOFPSIARRAY}
	|\wavearray| 
	& \eqdef \max_{\iota \in \{0,\cdots,4\}} | \Psi_\iota|.
\end{align}
For any norm $\|\cdot\|$ on scalar functions that appears in the paper, we set:
\begin{align} \label{E:NORMOFPSIARRAY}
	\| \wavearray\| 
	& \eqdef \max_{\iota \in \{0,\cdots,4\}} \| \Psi_\iota \|.
\end{align}
We use a similar convention for $\vortrenormalized$:
$|\vortrenormalized| = \max_{a = 1,2,3} |\vortrenormalized^a|$, 
and similarly for $\wavearraypartial$, $\GradEnt$, $\VortVort$,
etc. 
\item We use the following convention when taking norms of more than one variable at a time:
\begin{align} \label{E:NORMSOFMORETHANONEVARIABLE}
 \| (\vortrenormalized,\GradEnt)\| 
	& \eqdef \max\{  \| \vortrenormalized\|,  \| \GradEnt \|\}.
\end{align}
\end{itemize}
\end{definition}

\subsection{The acoustical metric and related geometric objects}
In the following definition, we introduce the acoustical metric and its inverse. 
This Lorentzian\footnote{By ``Lorentzian,'' we mean that viewed as a quadratic form, the symmetric $4 \times 4$ matrix
$(\gfour_{\alpha \beta})_{\alpha,\beta=0,1,2,3}$ has signature $(-,+,+,+)$.} 
metric drives the propagation of sound waves and is necessary to reveal the full geometry of the singular boundary.

\begin{definition}[The acoustical metric] \label{D:ACOUSTICALMETRICDEF}
Relative to the Cartesian coordinates $(t,x^1,x^2,x^3)$, 
we define the acoustical metric $\gfour$ and the inverse acoustical metric $\gfour^{-1}$ as 
follows, where the material derivative vectorfield $\Transport$ is defined in \eqref{E:MATERIALDERIVATIVEVECOTRFIELD}
and the speed of sound $\Speed$ is defined in \eqref{E:SOUNDSPEED}:
\begin{subequations}
\begin{align}
\gfour & = - \rmD t \otimes \rmD t + \Speed^{-2} \sum_{a=1}^3 (\rmD x^a - v^a \rmD t)\otimes  (\rmD x^a - v^a \rmD t), 
\label{E:ACOUSTICALMETRIC} \\
\gfour^{-1} & = - \Transport \otimes \Transport + \Speed^2 \sum_{a=1}^3 \p_a\otimes \p_a.
\label{E:INVERSEACOUSTICALMETRIC}
\end{align}
\end{subequations}
\end{definition}

Straightforward calculations yield that $\gfour_{\alpha \gamma} (\gfour^{-1})^{\gamma \beta} = \updelta_{\alpha}^{\beta}$,
where $\updelta_{\alpha}^{\beta}$ is the Kronecker delta, i.e., $\gfour^{-1}$ is indeed the inverse of $\gfour$.
In the remainder of the article, we silently lower and raise lowercase Greek indices with $\gfour$ and $\gfour^{-1}$, 
e.g., $V^{\alpha} = (\gfour^{-1})^{\alpha \beta} V_{\beta}$.


In our forthcoming analysis, the undifferentiated quantities
$v^i$ and $\Speed - 1$ will be small, where 
we quantify their smallness via the parameters $\mathring{\upalpha}$ and $\initialsmall$,
which we introduce in Sect.\,\ref{S:PARAMETERSANDSIZEASSUMPTIONSANDCONVENTIONSFORCONSTANTS}.
Hence, in view of \eqref{E:ACOUSTICALMETRIC}, we find it convenient to introduce the following decomposition:
\begin{align} \label{E:SPLITMETRICINTOMINKOWSKIANDREMAINDERPART}
	\gfour_{\alpha \beta}(\wavearray) 
	& = m_{\alpha \beta} + \gfour_{\alpha \beta}^{(\textnormal{Small})}(\wavearray)
\end{align}
where $m_{\alpha \beta} = \diag(-1,1,1,1)$ is the Minkowski metric and $\gfour_{\alpha \beta}^{(\textnormal{Small})}(\wavearray)$ is a smooth function of $\wavearray$ satisfying: 
\begin{align} \label{E:METRICPERTURBATIONVANISHESATTRIVIALPSISOLUTION}
	\gfour_{\alpha \beta}^{(\textnormal{Small})}(\wavearray = 0) 
	& = 0.
\end{align}


The scalar functions $G^{\iota}_{\alpha \beta}$ in the following definition
will appear as coefficients in many of the equations that we study.

\begin{definition}[$\wavearray$-derivatives of $\gfour$] 
\label{D:DERIVATIVESOFMETRICWRTFLUID}
Viewing the Cartesian component functions
$\gfour_{\alpha \beta} = \gfour_{\alpha \beta}(\wavearray)$ 
as functions of the wave variables, for $\alpha,\beta = 0,1,2,3$ and $\iota = 0,1,2,3,4$, 
we define: 
\begin{subequations}
\begin{align} \label{E:DERIVATIVESOFMETRICWRTFLUID}
G^{\iota}_{\alpha \beta}(\wavearray) 
& \eqdef \frac{\p}{\p \Psi_{\iota}} \gfour_{\alpha \beta}(\wavearray), 
	\\
\vec{G}_{\alpha \beta} 
& = \vec{G}_{\alpha \beta}(\wavearray)  
\eqdef 
\left(G_{\alpha \beta}^0 (\wavearray), 
	G_{\alpha \beta}^1 (\wavearray), 
	G_{\alpha \beta}^2 (\wavearray), 
	G_{\alpha \beta}^3 (\wavearray), 
	G_{\alpha \beta}^4 (\wavearray)
\right).
 \label{E:ARRAYVERSIONDERIVATIVESOFMETRICWRTFLUID}
\end{align}
\end{subequations}
\end{definition}

For each fixed $\iota \in \{0,\cdots,4\}$, 
we view $\{G^{\iota}_{\alpha \beta}\}_{\alpha,\beta = 0,\cdots,3}$ to be the Cartesian components of the spacetime tensorfield
``$G^{\iota}$.''
Similarly, we view 
$\{ \vec{G}_{\alpha \beta}\}_{\alpha,\beta = 0,\cdots,3}$ to be the Cartesian components of 
the array-valued spacetime tensorfield $\vec{G}$. 

\begin{definition}[Operators involving $\wavearray$] \label{D:DERIVATIVESOFARRAYS} Let $V_1,V_2$ be vectorfields, 
and let $D$ be a differential operator. We define:
\begin{align} \label{D:DERIVATIVEOFARRAY}
D \wavearray 
& \eqdef (D \Psi_0, D \Psi_1, D \Psi_2, D \Psi_3, D \Psi_4),
&
\vec{G}_{V_1 V_2} \diamond D \wavearray 
& \eqdef \sum_{\iota = 0}^4 G_{\alpha \beta}^\iota V_1^{\alpha}V_2^\beta D \Psi_\iota.
\end{align}
\end{definition}

\subsection{Covariant wave operator and $\gfour$-null forms}
\label{SS:COVARIANTWAVEOPERATORANDGNULLFORMS}
In this section, we provide some definitions that we need
to state Theorem~\ref{T:GEOMETRICWAVETRANSPORTSYSTEM}, which provides
the geometric formulation of compressible Euler flow that we use throughout our analysis.

We start by recalling the standard definition of the covariant wave operator $\square_{\gfour}$.

\begin{definition}[Covariant wave operator of the acoustical metric]
\label{D:COVWAVEOP}
The covariant wave operator $\square_{\gfour}$ 
of the acoustical metric $\gfour = \gfour(\wavearray)$ 
acts on scalar-valued functions $\varphi$ as follows:\footnote{The formula \eqref{E:WAVEOPERATORARBITRARYCOORDINATES}
holds relative to arbitrary coordinates.}
\begin{align} \label{E:WAVEOPERATORARBITRARYCOORDINATES}
\square_{\gfour} \varphi
	:= \frac{1}{\sqrt{|\mbox{\upshape det} \gfour|}}
	\partial_{\alpha}
	\left\lbrace
			\sqrt{|\mbox{\upshape det} \gfour|} (\gfour^{-1})^{\alpha \beta}
			\partial_{\beta} \varphi
	\right\rbrace.
\end{align}
\end{definition}

We now recall the definition of a standard null form with respect to the acoustical metric
(``$\gfour$-null form'' for short).

\begin{definition}[Standard $\gfour$-Null forms]
\label{D:STANDARDNULLFORMS}
Let $\varphi$ and $\widetilde{\varphi}$ be scalar functions.
We define
$\nullform^{(g)}(\partial\varphi,\partial\widetilde{\varphi})$ to be the following 
derivative-quadratic term:
\begin{subequations}
\begin{align}\label{E:Q0DEF}
\nullform^{(\gfour)}(\partial \varphi,\partial \widetilde{\varphi})
& = 
(\gfour^{-1})^{\alpha \beta} 
\partial_{\alpha} \varphi \partial_{\beta} \widetilde{\varphi}.
\end{align}

For $0 \leq \alpha < \beta \leq 3$, 
we define $\nullform_{\alpha \beta}(\partial \varphi,\partial\widetilde{\varphi})$
to be the following derivative-quadratic term:
\begin{align}\label{E:QIJDEF}
\nullform_{\alpha \beta}(\partial\varphi,\partial\widetilde{\varphi})
	= 	\partial_{\alpha} \varphi \partial_{\beta} \widetilde{\varphi} 
			- 
			\partial_{\beta} \varphi \partial_{\alpha} \widetilde{\varphi}.
\end{align}
\end{subequations}
\end{definition}

In the rest of the paper, we use the terminology
\emph{null form relative to $\gfour$} 
or
\emph{$\gfour$-null form} 
to denote any linear combination
of the standard null forms
\eqref{E:Q0DEF}--\eqref{E:QIJDEF}
with (possibly solution dependent) coefficients that are controllable under the scope of our approach.


\subsection{The geometric wave-transport-divergence-curl formulation of the compressible Euler equations} 
\label{SS:GEOMETRICFORMULATIONOFFLOW}
Our main results fundamentally rely on the following formulation of the compressible Euler equations, first derived in \cite{jS2019c}.

\begin{theorem} [The geometric wave-transport-divergence-curl formulation of the compressible Euler equations]
	\label{T:GEOMETRICWAVETRANSPORTSYSTEM}
	Let $\overline{\varrho} > 0$ be any constant background density,\footnote{Recall
	that $\LogDensity$ depends on $\overline{\varrho}$; see Def.\,\ref{D:LOGDENS}.}
	and assume that $(\LogDensity,v^1,v^2,v^3,\Ent)$ is
	a solution
	to the compressible Euler equations 
	\eqref{E:BVIEVOLUTION}--\eqref{E:BENTROPYEVOLUTION}
	in three spatial dimensions under an arbitrary equation of state $p = p(\varrho,\Ent)$ with positive sound speed 
	$\Speed$ (see \eqref{E:SOUNDSPEED}).
	Let $\Transport$ be the material derivative vectorfield defined in \eqref{E:MATERIALDERIVATIVEVECOTRFIELD},
	let $\gfour$ be the acoustical metric from Def.\,\ref{D:ACOUSTICALMETRICDEF},
	let $\square_{\gfour}$ be the corresponding covariant wave operator from Def.\,\ref{D:COVWAVEOP},
	let $\almostRiemann_{(\pm)}$ be the almost Riemann invariants from Def.\,\ref{D:ALMOSTRIEMANNINVARIANTS}
	(see Remark~\ref{R:HOWWEESTIMATEDENSITYANDV1} concerning their significance for this paper),
	let $\almostRiemannfunction  = \almostRiemannfunction(\LogDensity,\Ent)$ be the function from
	Def.\,\ref{D:ALMOSTRIEMANNINVARIANTS},
	and let 
	$\vortrenormalized$,
	$\GradEnt$,
	$\VortVort$,
	and
	$\DivGradEnt$
	be the higher order variables from Def.\,\ref{D:HIGHERORDERFLUIDVARIABLES}.
	Then the scalar-valued functions
	$\LogDensity$,
	$v^i$, 
	$\almostRiemann_{(\pm)}$, 
	$\Ent$,
$\vortrenormalized^i$,
$\GradEnt^i$,
$\Flatdiv \vortrenormalized$,
$\VortVort^i$,
$\DivGradEnt$,
and
$(\Flatcurl \GradEnt)^i$,
	($i=1,2,3$),
	also solve the following equations,
	where $\upepsilon_{ijk}$ is the fully antisymmetric symbol normalized by $\upepsilon_{123}=1$, 
	and \textbf{the Cartesian component functions $v^i$ are 
	treated as scalar-valued functions
	under covariant differentiation on LHS~\eqref{E:VELOCITYWAVEEQUATION}}:

\medskip

\noindent \underline{\textbf{\upshape Covariant wave equations}}.
	\begin{subequations} \label{E:COVARIANTWAVEEQUATIONSWAVEVARIABLES}
	\begin{align}
		\square_{\gfour(\wavearray)} v^i
		& = 
			- 
			\Speed^2 \exp(2 \LogDensity) \VortVort^i
			+ 
			\nullform_{(v)}^i
			+ 
			\mathfrak{L}_{(v)}^i,
			\label{E:VELOCITYWAVEEQUATION}	\\
		\square_{\gfour(\wavearray)} \almostRiemann_{(\pm)} 
		& 
		= - \Speed^2 \exp(2 \LogDensity) \VortVort^1 
				\pm 
				\left\lbrace 
					F_{;\Ent} c^2 \exp(2\LogDensity)
					- 
					\Speed \exp(\LogDensity) \frac{p_{;\Ent}}{\overline{\varrho}} 
				\right\rbrace  
				\DivGradEnt 
				+ 
				\nullform_{(\pm)} 
				+ 
				\mathfrak{L}_{(\pm)}, \label{E:WAVEEQUATIONFORALMOSTRIEMANNINVARIANTS} \\
	\square_{\gfour(\wavearray)} \LogDensity
	& = 
		-
		\exp(\LogDensity) \frac{p_{;\Ent}}{\overline{\varrho}} \DivGradEnt
		+
		\nullform_{(\LogDensity)}
		+
		\mathfrak{L}_{(\LogDensity)},
			\label{E:RENORMALIZEDDENSITYWAVEEQUATION} 
				\\
	\square_{\gfour(\wavearray)} \Ent
	& = 
		\Speed^2 \exp(2 \LogDensity)  \DivGradEnt
		+
		\mathfrak{L}_{(\Ent)}.
	\label{E:ENTROPYWAVEEQUATION}
\end{align}
\end{subequations}

\medskip

\noindent \underline{\textbf{\upshape Transport equations}}.
\begin{subequations}
\begin{align}	\Transport \vortrenormalized^i
	& = \mathfrak{L}_{(\vortrenormalized)}^i,
		\label{E:RENORMALIZEDVORTICTITYTRANSPORTEQUATION}
		\\
	\Transport \Ent
	& = 0,	
		\label{E:ENTROPYTRANSPORTMAINSYSTEM}
			\\
	\Transport \GradEnt^i
	& = \mathfrak{L}_{(\GradEnt)}^i.
		\label{E:GRADENTROPYTRANSPORT}
\end{align}
\end{subequations}

\medskip	
	
\noindent \underline{\textbf{\upshape Transport-divergence-curl system for the specific vorticity}}.
\begin{subequations}
\begin{align} \label{E:FLATDIVOFRENORMALIZEDVORTICITY}
	\Flatdiv \vortrenormalized
	& = 
		\mathfrak{L}_{(\Flatdiv \vortrenormalized)},
		\\
\Transport \VortVort^i
& = 
	\mainnullform_{(\VortVort)}^i 
	+ 
	\nullform_{(\VortVort)}^i
	+	
	\mathfrak{L}_{(\VortVort)}^i. 
	\label{E:EVOLUTIONEQUATIONFLATCURLRENORMALIZEDVORTICITY} 
\end{align}	
\end{subequations}

\medskip

\noindent \underline{\textbf{\upshape Transport-divergence-curl system for the entropy gradient}}.
\begin{subequations}
\begin{align} 	
\Transport \DivGradEnt
	& =  \mainnullform_{(\DivGradEnt)} +
			\nullform_{(\DivGradEnt)},
 \label{E:TRANSPORTFLATDIVGRADENT}
			\\
	(\Flatcurl \GradEnt)^i & = 0.
	\label{E:CURLGRADENTVANISHES}
\end{align}
\end{subequations}
	Above, the main terms $\mainnullform_{(\VortVort)}^i$ and $\mainnullform_{(\DivGradEnt)}$ 
	in the transport equations for the modified fluid variables are 
	the \textbf{null forms relative to} $\gfour$ (see Def.\,\ref{D:STANDARDNULLFORMS})
	defined by:\footnote{Actually, the last the last term on
 RHS~\eqref{E:TRANSPORTDIVGRADENTMAINTERMS} is not a null form, but rather 
	a simpler harmless error term.}
	\begin{subequations}
	\begin{align}
	\begin{split} \label{E:TRANSPORTVORTVORTMAINTERMS}
		\mainnullform_{(\VortVort)}^i & \eqdef 
		- 
		2 \updelta_{jk} \upepsilon_{iab} \exp(-\LogDensity) (\partial_a v^j) \partial_b \vortrenormalized^k
		+
		\upepsilon_{ajk}
		\exp(-\LogDensity)
		(\partial_a v^i) 
		\partial_j \vortrenormalized^k
			\\
& \ \
		+ 
		\exp(-3 \LogDensity) \Speed^{-2} \frac{p_{;\Ent}}{\overline{\varrho}} 
		\left\lbrace
			(\Transport \GradEnt^a) \partial_a v^i
			-
			(\Transport v^i) \partial_a \GradEnt^a
		\right\rbrace
			\\
	& \ \
		+
		\exp(-3 \LogDensity) \Speed^{-2} \frac{p_{;\Ent}}{\overline{\varrho}}  
		\left\lbrace
			(\Transport v^a) \partial_a \GradEnt^i
			- 
			( \Transport \GradEnt^i) \partial_a v^a
		\right\rbrace,
	\end{split}
			\\
	\mainnullform_{(\DivGradEnt)} & = 	2 \exp(-2 \LogDensity) 
			\left\lbrace
				(\partial_a v^a) \partial_b \GradEnt^b
				-
				(\partial_a \GradEnt^b) \partial_b v^a
			\right\rbrace
			+
			\exp(-\LogDensity) \updelta_{ab} (\Flatcurl \vortrenormalized)^a \GradEnt^b. 
			\label{E:TRANSPORTDIVGRADENTMAINTERMS}  
	\end{align}
	\end{subequations} 
	Moreover, 
	$\nullform_{(v)}^i$,
	$\nullform_{(\pm)}$,
	$\nullform_{(\LogDensity)}$, 
	$\nullform_{(\VortVort)}^i$,
	and
	$\nullform_{(\DivGradEnt)}$
	are\footnote{The term 
	$\mainnullform_{(\VortVort)}^i$
	on RHS~\eqref{E:EVOLUTIONEQUATIONFLATCURLRENORMALIZEDVORTICITY} 
	and the term $\mainnullform_{(\DivGradEnt)}$
	on RHS~\eqref{E:TRANSPORTFLATDIVGRADENT} are also null forms relative to $\gfour$.
	We have isolated these two null forms with different notation because
	they 
	are more difficult to treat than 
	$\nullform_{(v)}^i$, $\nullform_{(\pm)}$
	$\nullform_{(\LogDensity)}$, 
	$\nullform_{(\VortVort)}^i$,
	and
	$\nullform_{(\DivGradEnt)}$;
	to bound the top-order derivatives of the ``$\mainnullform$'' terms,
	we rely on the delicate ``elliptic-hyperbolic'' identities
	that we derive in Sect.\,\ref{S:ELLIPTICHYPERBOLICIDENTITIES}.
	\label{FN:MORENULLFORMS}} 
	the null forms relative to $\gfour$ defined by:
	\begin{subequations}
		\begin{align}
		\nullform_{(v)}^i	
		& \eqdef 	-
					\left\lbrace
						1
						+
						\Speed^{-1} \Speed_{;\LogDensity}
					\right\rbrace
					(\gfour^{-1})^{\alpha \beta} (\partial_{\alpha} \LogDensity) \partial_{\beta} v^i,
			\label{E:VELOCITYNULLFORM} 
			\\
		\nullform_{(\pm)} 
			& 
			\eqdef \nullform_{(v)}^1 \mp 2 \Speed_{;\LogDensity} 
			(\gfour^{-1})^{\alpha \beta} \p_\alpha \LogDensity \p_\beta \LogDensity
			\pm 
			\Speed \left\lbrace (\p_a v^a) (\p_b v^b) - (\p_a v^b) \p_b v^a 
			\right\rbrace, 
				\label{E:RIEMANNINVARIANTWAVEEQUATIONNULLFORM}
		\\
		\nullform_{(\LogDensity)}
		& \eqdef 
		- 
		3 \Speed^{-1} \Speed_{;\LogDensity} 
		(\gfour^{-1})^{\alpha \beta} (\partial_{\alpha} \LogDensity) \partial_{\beta} \LogDensity
		+ 
		\left\lbrace
			(\partial_a v^a) \partial_b v^b
			-
			(\partial_a v^b) \partial_b v^a
		\right\rbrace,
			\label{E:DENSITYNULLFORM}
				\\
	\begin{split} \label{E:RENORMALIZEDVORTICITYCURLNULLFORM}
	\nullform_{(\VortVort)}^i
	& \eqdef
		\exp(-3 \LogDensity) \Speed^{-2} \frac{p_{;\Ent}}{\overline{\varrho}}  \GradEnt^i
		\left\lbrace
			(\partial_a v^b) \partial_b v^a
			-
			(\partial_a v^a) \partial_b v^b
		\right\rbrace
			 \\
& \ \
		+ 
		\exp(-3 \LogDensity) \Speed^{-2} \frac{p_{;\Ent}}{\overline{\varrho}}
		\GradEnt^b 
		\left\lbrace
			(\partial_a v^a) \partial_b v^i 
			- 
			(\partial_a v^i) \partial_b v^a 
		\right\rbrace
		\\
& \ \
		+ 
		2 \exp(-3 \LogDensity) \Speed^{-2} \frac{p_{;\Ent}}{\overline{\varrho}}
		\GradEnt^a
		\left\lbrace
			(\partial_a \LogDensity)  \Transport v^i
			 - 
		  (\partial_a v^i) \Transport \LogDensity
		\right\rbrace
			\\
	&  \ \
			+ 
			2 \exp(-3 \LogDensity) \Speed^{-3} \Speed_{;\LogDensity} \frac{p_{;\Ent}}{\overline{\varrho}}
			\GradEnt^a 
			\left\lbrace
				(\partial_a \LogDensity)  \Transport v^i
				- 
				(\partial_a v^i) \Transport \LogDensity
			\right\rbrace
				\\
	& \ \
		+ 
		\exp(-3 \LogDensity) \Speed^{-2} \frac{p_{;\Ent;\LogDensity}}{\overline{\varrho}}
		\GradEnt^a 
		\left\lbrace
			(\partial_a v^i)
			\Transport \LogDensity 
			- 
			(\partial_a \LogDensity) \Transport v^i
		\right\rbrace
			\\
	& \ \
		+
		\exp(-3 \LogDensity) \Speed^{-2} \frac{p_{;\Ent;\LogDensity}}{\overline{\varrho}} \GradEnt^i
		\left\lbrace
			(\Transport v^a) \partial_a \LogDensity
			-
			 (\Transport \LogDensity) \partial_a v^a
		\right\rbrace
				\\
	& \ \
		+
		2 \exp(-3 \LogDensity) \Speed^{-2} \frac{p_{;\Ent}}{\overline{\varrho}} \GradEnt^i
		\left\lbrace
			(\Transport \LogDensity) \partial_a v^a 
			- 
			(\Transport v^a) \partial_a \LogDensity
		\right\rbrace
				\\
	& \ \
		 	+
			2 \exp(-3 \LogDensity) \Speed^{-3} \Speed_{;\LogDensity} \frac{p_{;\Ent}}{\overline{\varrho}} \GradEnt^i
			\left\lbrace
				(\Transport \LogDensity) \partial_a v^a
		 		- 
		 		(\Transport v^a) \partial_a \LogDensity
		 	\right\rbrace,
\end{split} 
		 		\\
\label{E:DIVENTROPYGRADIENTNULLFORM}
\nullform_{(\DivGradEnt)} 
	& \eqdef
		2 \exp(-2 \LogDensity) 
		\GradEnt^a 
		\left\lbrace
			(\partial_a v^b) \partial_b \LogDensity
			-
			(\partial_a \LogDensity)
			\partial_b v^b 
		\right\rbrace.
\end{align}
\end{subequations}
	In addition, the terms
	$\mathfrak{L}_{(v)}^i$, $\mathfrak{L}_{(\pm)}$, 
	$\mathfrak{L}_{(\LogDensity)}$,
	$\mathfrak{L}_{(\Ent)}$,
	$\mathfrak{L}_{(\vortrenormalized)}^i$,
	$\mathfrak{L}_{(\GradEnt)}^i$,
		$\mathfrak{L}_{(\Flatdiv \vortrenormalized)}$,
		and
	$\mathfrak{L}_{(\VortVort)}^i$,
	which are at most linear in the derivatives of the unknowns, are defined as follows:
	\begin{subequations}
	\begin{align} 
	\begin{split} \label{E:VELOCITYILINEARORBETTER} 
		\mathfrak{L}_{(v)}^i
		& \eqdef 
		2 \exp(\LogDensity) \upepsilon_{iab} (\Transport v^a) \vortrenormalized^b
		-
		\frac{p_{;\Ent}}{\overline{\varrho}} \upepsilon_{iab} \vortrenormalized^a \GradEnt^b
			\\
	& \ \
		- 
		\frac{1}{2} \exp(-\LogDensity) \frac{p_{;\LogDensity;\Ent}}{\overline{\varrho}} \GradEnt^a \partial_a v^i
				\\
		& \ \
		- 
		2 \exp(-\LogDensity) \Speed^{-1} \Speed_{;\LogDensity} \frac{p_{;\Ent}}{\overline{\varrho}} 
		(\Transport \LogDensity) \GradEnt^i
		+
		\exp(-\LogDensity) \frac{p_{;\Ent;\LogDensity}}{\overline{\varrho}} (\Transport \LogDensity) \GradEnt^i,
		\end{split} 
			\\
	\begin{split} \label{E:RIEMANNLINEARORBETTER}  
	\mathfrak{L}_{(\pm)} 
		& 
		\eqdef 
		\mathfrak{L}^1_{(v)} \mp \frac{5}{2} \Speed \exp(-\LogDensity) \frac{p_{;\GradEnt^b\LogDensity}}{\overline{\varrho}} \GradEnt^a \p_a \LogDensity \pm 2 \Speed^2 c_{;\Ent} \GradEnt^a\p_a \LogDensity  
				\\
		& \ \  
		\mp \Speed \exp(-\LogDensity) \frac{p_{;\Ent;\Ent}}{\overline{\varrho}} \updelta_{ab}\GradEnt^a\GradEnt^b 
		\pm
		F_{;\Ent} \mathfrak{L}_{(\Ent)}, 
		\end{split} 
			\\
		\mathfrak{L}_{(\LogDensity)}  \label{E:DENSITYLINEARORBETTER}
		& \eqdef
		-
		\frac{5}{2} \exp(-\LogDensity) \frac{p_{;\Ent;\LogDensity}}{\overline{\varrho}} \GradEnt^a \partial_a \LogDensity 
		-
		\exp(-\LogDensity) \frac{p_{;\Ent;\Ent}}{\overline{\varrho}} \updelta_{ab} \GradEnt^a \GradEnt^b,
		  \\
		\mathfrak{L}_{(\Ent)}
		& \eqdef 
			\Speed^2  \GradEnt^a \partial_a \LogDensity
			- 
			\Speed \Speed_{;\LogDensity} \GradEnt^a \partial_a \LogDensity
			- 
			\Speed \Speed_{;\Ent} \updelta_{ab} \GradEnt^a \GradEnt^b,
			\label{E:ENTROPYLINEARORBETTER} 
				\\
		\mathfrak{L}_{(\vortrenormalized)}^i
		& \eqdef 
		\vortrenormalized^a \partial_a v^i
		-
		\exp(-2 \LogDensity) \Speed^{-2} \frac{p_{;\Ent}}{\overline{\varrho}} \upepsilon_{iab} (\Transport v^a) \GradEnt^b,
		\label{E:SPECIFICVORTICITYLINEARORBETTER}
			\\
		\mathfrak{L}_{(\GradEnt)}^i
		& \eqdef
			- 
			\GradEnt^a \partial_a v^i
			+ 
			\upepsilon_{iab} \exp(\LogDensity) \vortrenormalized^a \GradEnt^b,
			\label{E:ENTROPYGRADIENTLINEARORBETTER}
				\\
	\mathfrak{L}_{(\Flatdiv \vortrenormalized)}
		& \eqdef - \vortrenormalized^a \partial_a \LogDensity,
		\label{E:RENORMALIZEDVORTICITYDIVLINEARORBETTER} 
		\\
	\begin{split} \label{E:RENORMALIZEDVORTICITYCURLLINEARORBETTER} 
	\mathfrak{L}_{(\VortVort)}^i	 		
	& \eqdef
		 		2 \exp(-3 \LogDensity) \Speed^{-3} \Speed_{;\Ent} \frac{p_{;\Ent}}{\overline{\varrho}} 
				(\Transport v^i) \updelta_{ab} \GradEnt^a \GradEnt^b
				\\
		& \ \
			-
				2 \exp(-3 \LogDensity) \Speed^{-3} \Speed_{;\Ent} \frac{p_{;\Ent}}{\overline{\varrho}} 
				\updelta_{ab} \GradEnt^a (\Transport v^b) \GradEnt^i
					\\
		& \ \
			+ 
			\exp(-3 \LogDensity) \Speed^{-2} \frac{p_{;\Ent;\Ent}}{\overline{\varrho}} \updelta_{ab} (\Transport v^a) \GradEnt^b \GradEnt^i
				\\
		& \ \ 
			- 
			\exp(-3 \LogDensity) \Speed^{-2} \frac{p_{;\Ent;\Ent}}{\overline{\varrho}} (\Transport v^i) \updelta_{ab} \GradEnt^a \GradEnt^b.
		\end{split}
		\end{align}
	\end{subequations}
\end{theorem}

\begin{proof}[Discussion of the proof]
	Theorem~\ref{T:GEOMETRICWAVETRANSPORTSYSTEM} was essentially proved as \cite{jS2019c}*{Theorem 1},
	except that the wave equations \eqref{E:WAVEEQUATIONFORALMOSTRIEMANNINVARIANTS} for $\almostRiemann_{(\pm)}$ were not derived there.
	In \cite{jLjS2021}*{Theorem~5.1},
	those wave equations were derived 
	as a straightforward consequence of \cite{jS2019c}*{Theorem 1}.

\end{proof}


\section{The acoustic geometry and the controlling quantities $\controlvars$ and $\badcontrolvars$}
\label{S:ACOUSTICGEOMETRYANDCOMMUTATORVECTORFIELDS}
In this section, we construct the acoustic geometry,
reveal its basic properties,
and provide the evolution equations satisfied by various geometric tensors.
Our approach is based on the one pioneered by Christodoulou \cite{dC2007} in his study
of irrotational and isentropic solutions.
The fundamental object behind all the constructions is an acoustic eikonal function,
that is, a solution $u$ to the acoustic eikonal equation.
The eikonal function is fundamental for our approach because our proof shows that the fluid solution
remains rather smooth relative to the geometric coordinates $(t,u,x^2,x^3)$.
In particular, later on, we will use $u$ to construct suitable commutation and multiplier vectorfields
out of the geometric coordinates to control the solution and acoustic geometry up
to top-order.
This allows us, at least in some ways, to treat the problem of shock formation
as a long-time existence problem relative to the geometric coordinates.
We also introduce the solution variable arrays $\controlvars$ and $\badcontrolvars$,
which contain the wave variables and various components of the acoustic geometry.
We use these arrays throughout the paper to simplify the notation and to allow for convenient,
schematic expressions. We sometimes refer to 
$\controlvars$ and $\badcontrolvars$ as the ``controlling quantities'' because all of the
quantities that we analyze can, in principle, be constructed out of them.

\subsection{The eikonal function and inverse foliation density}
\label{SS:EIKONAFUNCTIONANDINVERSEFOLIATIONDENSITY}
In the following definition, we introduce the \emph{eikonal function} $u$ and the inverse foliation density $\upmu$. The level sets of $u$ are the characteristic for the wave operator $\square_{\gfour}$, 
while the inverse (i.e., reciprocal) 
of the inverse foliation density measures the density of these characteristics. 
In particular, the vanishing of $\upmu$ corresponds to the infinite density of the characteristics.
Our main results show that in the regime under study, 
the vanishing of $\upmu$ coincides with the blowup of
$|\partial_1 \RRiemann|$.

\begin{definition}[Eikonal function and inverse foliation density]
\label{D:EIKONALFUNCTIONANDMU}
The \emph{eikonal function} $u$ is the solution of the following fully nonlinear hyperbolic initial value problem,
where $\gfour$ is the acoustical metric defined in \eqref{E:ACOUSTICALMETRIC}
and the PDE is known as the \emph{acoustic eikonal equation}:
\begin{align} \label{E:EIKONALEQUATION}
\begin{cases}
(\gfour^{-1})^{\alpha \beta} \partial_{\alpha} u \partial_{\beta} u = 0, 
	\\
\p_t u > 0, 
	\\
u|_{\Sigma_0} = - x^1.
\end{cases}
\end{align}
We define the \emph{inverse foliation density} $\upmu$ by:
\begin{align} \label{E:MUDEF}
	\upmu 
	& 
	\eqdef 
	- 
	\frac{1}{(\gfour^{-1})^{\alpha \beta} \partial_{\alpha} t \partial_{\beta} u} > 0.
\end{align}
\end{definition}

\subsection{Acoustical subsets of spacetime}
\label{SS:ACOUSTICSUBSETSOFSPACETIME}

\begin{definition}[Acoustical subsets of spacetime] 
	\label{D:ACOUSTICSUBSETSOFSPACETIME}
We define the following ``acoustical subsets'' of spacetime:
\begin{subequations}
\begin{align}
\Sigma_{t'} 
& \eqdef 
\{(t,x^1,x^2,x^3) \in \R\times\R\times\T^2 \ | \ t = t'\}, 
	\label{E:SIGMAT}
	\\
\nullhyparg{u'} 
& \eqdef \{(t,x^1,x^2,x^3) \in \R\times\R\times\T^2 \ | \ u(t,x^1,x^2,x^3) = u'\}, 
	\label{E:NULLHYPERSURFACES}
	\\
\ell_{t',u'} 
& 
\eqdef 
\Sigma_{t'} \cap \nullhyparg{u'} 
= 
\{(t,x^1,x^2,x^3) \in \R\times\R\times\T^2 \ | \ t= t',\ u(t,x^1,x^2,x^3) = u'\}. 
\label{E:ELLTUSMOOTHTORI} 
\end{align}
\end{subequations}

Given real numbers $u_1 \leq u_2$ and $t_1 \leq t_2$, we define the following ``truncated'' subsets of spacetime:
\begin{subequations}
\begin{align}
\Sigma_{t'}^{[u_1,u_2]} 
& \eqdef \Sigma_{t'} \cap \{(t,x^1,x^2,x^3) \in \R\times\R\times\T^2 \ | \ u_1 \leq u(t,x^1,x^2,x^3) \leq u_2 \}, 
	\label{E:TRUNCATEDSIGMAT}
	\\
\nullhyparg{u'}^{[t_1,t_2]} 
& \eqdef \nullhyparg{u'} \cap \{(t,x^1,x^2,x^3) \in \R\times\R\times\T^2 \ | \ t_1 \leq t \leq t_2 \}.
	\label{E:TRUNCATEDNULLHYPERSURFACES}
\end{align}
\end{subequations}

\end{definition}

We refer to the $\Sigma_t$ as ``constant Cartesian-time hypersurfaces,'' 
the $\nullhyparg{u}$ as ``null hypersurfaces,''
``acoustic characteristics,''  
or ``characteristics,'' and the $\ell_{t,u}$ as ``smooth tori.'' 
We emphasize that with the exception of the appendices, 
in this paper,
we will \emph{not} derive estimates on $\Sigma_t$ 
or the $\ell_{t',u'}$. Instead, we will control the solution
on the rough hypersurfaces 
and the rough tori of Def.\,\ref{D:TRUNCATEDROUGHSUBSETS}.


\subsection{Projection tensorfields and related differential operators}
\label{SS:PROJECTIONSANDRELATEDDIFFERENTIALOPERATORS}

\begin{definition}[Projection tensorfields and tangency to hypersurfaces]  
\label{D:PROJECTIONTENSORFIELDSANDTANGENCYTOHYPERSURFACES}
\hfill
\begin{enumerate}
\item We define the type $\binom{1}{1}$
$\Sigma_t$-projection tensorfield $\Sigmatproject$ and the type $\binom{1}{1}$ $\ell_{t,u}$-projection tensorfield 
$\smoothtorusproject$ as follows, where $\updelta_{\beta}^{\alpha}$ denotes the Kronecker delta: 
\begin{subequations}
\begin{align}
\Sigmatproject_{\beta}^{\ \alpha}
& \eqdef 
	\updelta_{\beta}^{\alpha}
	+ 
	\Transport^{\alpha} 
	\Transport_\beta,  
	\label{E:CARTESIANSIGMAPROJECT} 
	\\
\smoothtorusproject^{\ \alpha}_\beta 
& 
\eqdef 
\updelta_{\beta}^{\alpha}
+  
\Transport^{\alpha} \Transport_\beta 
- 
X^{\alpha} X_\beta 
= 
\updelta_{\beta}^{\alpha}
- 
\Lunit^{\alpha} \updelta^0_\beta 
+
X^{\alpha} \Lunit_{\beta}.
 \label{E:SMOOTHTORUSPROJECT}
\end{align}
\end{subequations}
\item Given any type $\binom{m}{n}$ spacetime tensorfield $\upxi$, 
we respectively define its $\gfour$-orthogonal projection onto $\Sigma_t$, denoted by $\Sigmatproject \upxi$,
and its $\gfour$-orthogonal projection onto $\ell_{t,u}$, denoted by $\smoothtorusproject \upxi$,
as follows:
\begin{subequations}
\begin{align}
(\Sigmatproject \upxi)_{\beta_1 \cdots \beta_n}^{\alpha_1 \cdots \alpha_m}
& 
\eqdef 
\Sigmatproject_{\widetilde{\alpha}_1}^{\ \alpha_1}
\cdots 
\Sigmatproject_{\widetilde{\alpha}_m}^{\ \alpha_m}
\Sigmatproject_{\beta_1}^{\ \widetilde{\beta}_1}
\cdots 
\Sigmatproject_{ \beta_n}^{\ \widetilde{\beta}_n}
\upxi_{\widetilde{\beta}_1 \cdots \widetilde{\beta}_n}^{\widetilde{\alpha}_1 \cdots \widetilde{\alpha}_m}, 
	\label{E:PROJECTIONOFTENSORONTOCARTESIANSIGMAT} 
	\\
(\smoothtorusproject \upxi)_{\beta_1\cdots\beta_n}^{\alpha_1\cdots\alpha_m} 
& \eqdef 
\smoothtorusproject_{\widetilde{\alpha}_1}^{\ \alpha_1} 
\cdots 
\smoothtorusproject_{\widetilde{\alpha}_m}^{\ \alpha_m} 
\smoothtorusproject_{ \beta_1}^{\ \widetilde{\beta}_1}
\cdots 
\smoothtorusproject_{ \beta_n}^{\ \widetilde{\beta}_n}
\upxi_{\widetilde{\beta}_1 \cdots \widetilde{\beta}_n}^{\widetilde{\alpha}_1 \cdots \widetilde{\alpha}_m}.
	\label{E:PROJECTIONOFTENSORONTOFLATTORUS}
\end{align}
\end{subequations}
\item
We say that a spacetime tensorfield $\upxi$ is $\Sigma_t$-tangent if $\Sigmatproject \upxi = \upxi$. 
We say that a spacetime tensorfield $\upxi$ is $\ell_{t,u}$-tangent if $\smoothtorusproject \upxi = \upxi$. 
\item If $\upxi$ is a symmetric type $\binom{0}{2}$-spacetime tensor and $V$ is a vectorfield, 
then we define $\slashed{\upxi}_V \eqdef \smoothtorusproject(\upxi \cdot V)$, 
where $\upxi \cdot V$ is the one-form with components
$(\upxi \cdot V)_{\alpha} \eqdef \upxi_{\alpha \beta} V^{\beta}$.
\item If $\upxi$ is a spacetime tensor, then we define $\slashed{\upxi} = \smoothtorusproject \upxi$.
	From 3 above, it follows that $\xi$ is $\ell_{t,u}$-tangent if and only if $\slashed{\upxi} = \xi$.
\end{enumerate}
\end{definition}

It is straightforward to check that
$\Sigmatproject \Transport = 0$, while if $\SigmatTan$ is $\Sigma_t$-tangent, then 
$\Sigmatproject \SigmatTan = \SigmatTan$, i.e, in view of the properties of $\Transport$ from Lemma~\ref{L:BASICPROPERTIESOFVECTORFIELDS}, 
we see that $\Sigmatproject$ \emph{is} the $\gfour$-orthogonal projection onto $\Sigma_t$.
Similarly, $\smoothtorusproject \Lunit = \smoothtorusproject X = \smoothtorusproject \Transport = 0$,
while if $Y$ is $\ell_{t,u}$-tangent, then $\smoothtorusproject Y = Y$.


\subsection{First fundamental forms}
\label{SS:FIRSTFUDNAMENTALFORMS}

\begin{definition}[First fundamental forms] \hfill
\label{D:FIRSTFUNDAMENTALFORMS}
\begin{enumerate}
\item We define $g$, the first fundamental form of $\Sigma_t$ relative to $\gfour$, 
	to be the symmetric type $\binom{0}{2}$ tensorfield $\Sigmatproject \gfour$. Note that
	$g(Y,Z) = \gfour(Y,Z)$ for
	all pairs $(Y,Z)$ of $\Sigma_t$-tangent vectorfields.
	We define the corresponding inverse first fundamental form $g^{-1}$
	to be the symmetric type $\binom{2}{0}$ tensorfield that is $\gfour$-dual to $g$
	i.e., 
	$(g^{-1})^{\alpha \beta} 
	\eqdef 
	(\gfour^{-1})^{\alpha \widetilde{\alpha}} 
	(\gfour^{-1})^{\beta \widetilde{\beta}}
	g_{\widetilde{\alpha} \widetilde{\beta}}$.
	Note that the restriction of $g$ to $\Sigma_t$-tangent tensorfields is the Riemannian\footnote{It is Riemannian
	is because the $\gfour$-normal to $\Sigma_t$ is the $\gfour$-timelike vectorfield $\Transport$.} 
	metric on $\Sigma_t$ induced by $\gfour$. In particular, 
	relative to the Cartesian spatial coordinates, we have: 
	\begin{align} \label{E:CARTESIANCOMPONENTSOFFIRSTFUNDAMENTLAFORMOFSIGMAT}
		g_{ij} 
		& = \gfour_{ij}
		= \Speed^{-2} \updelta_{ij},
	\end{align}
	where $\updelta_{ij}$ denotes the Kronecker delta,
	and to obtain the last equality in \eqref{E:CARTESIANCOMPONENTSOFFIRSTFUNDAMENTLAFORMOFSIGMAT}, 
	we have used \eqref{E:ACOUSTICALMETRIC}.
 \item We define $\gtorus$, the first fundamental form of $\ell_{t,u}$ relative to $\gfour$, 
	to be the symmetric type $\binom{0}{2}$ tensorfield $\smoothtorusproject \gfour$. Note that
	$\gtorus(Y,Z) = \gfour(Y,Z)$ for
	all pairs $(Y,Z)$ of $\ell_{t,u}$-tangent vectorfields.
	Note that the restriction of $\gtorus$ to $\ell_{t,u}$-tangent tensorfields is the Riemannian\footnote{It is Riemannian
	is because $\ell_{t,u}$ is a submanifold of the spacelike hypersurface $\Sigma_t$.} 
	metric on $\ell_{t,u}$ induced by $\gfour$.
	We define the corresponding inverse first fundamental form $\gtorus^{-1}$
	to be the symmetric type $\binom{2}{0}$ tensorfield that is $\gfour$-dual to $\gtorus$
	i.e., 
	$(\gtorus^{-1})^{\alpha \beta} 
	\eqdef 
	(\gfour^{-1})^{\alpha \widetilde{\alpha}} 
	(\gfour^{-1})^{\beta \widetilde{\beta}}
	\gtorus_{\widetilde{\alpha} \widetilde{\beta}}$.
\end{enumerate}
\end{definition}

\subsection{Geometric coordinates and related vectorfields}
\label{SS:GEOMETRICCOORDIANTESANDPARTIALDERIVATIVEVECTORFIELDS}
\subsubsection{Geometric coordinates}
\label{SSS:GEOMETRICCOORDINATES}

\begin{definition}[The geometric coordinates and their corresponding partial derivative vectorfields] 
\label{D:GEOMETRICCOORDIANTESANDPARTIALDERIVATIVEVECTORFIELDS}
We define the \emph{geometric coordinate system} to be $(t,u,x^2,x^3)$. We define 
$\left\lbrace\geop{t},\geop{u},\geop{x^2},\geop{x^3} \right\rbrace$ to be the coordinate 
partial vectorfields in the geometric coordinate system.
\end{definition}

	\begin{remark}[Coordinate systems on $\ell_{t,u}$ and $\nullhyparg{u}$]
	Note that $(x^2,x^3)$ form a coordinate system on the smooth tori $\ell_{t,u}$ 
	and that $\left\lbrace \geop{x^2},\geop{x^3} \right\rbrace$ span the tangent space of $\ell_{t,u}$.
	Similarly, $(t,x^2,x^3)$ form a coordinate system on the null hypersurfaces $\nullhyparg{u}$ 
	and $\left\lbrace \geop{t}, \geop{x^2},\geop{x^3} \right\rbrace$ span the tangent space of 
	$\nullhyparg{u}$. We will silently use these basic facts throughout the rest of the article.
\end{remark}


\begin{notation}[Conventions used with $(x^2,x^3)$ and $\lbrace \geop{x^2}, \geop{x^3} \rbrace$]
	\label{N:NOTATIONINVOLVINGSMOOTHTORUSCOORDINATES}
	\hfill
\begin{enumerate}
\item If $V$ is a vectorfield, then for $A = 2,3$, $V^A \eqdef V x^A = V^{\alpha} \partial_{\alpha} x^A$.
	In particular, if $V$ is $\ell_{t,u}$-tangent, then $V = V^A \geop{x^A}$, and  
	$V^A$ are the components of $V$ with respect to geometric coordinates $(x^2,x^3)$ on $\ell_{t,u}$.
\item If $\upxi$ is a one-form, then we denote its contraction with $\geop{x^A}$ by using the abbreviated notation 
	$\upxi_A \eqdef \upxi\left( \geop{x^A} \right) = \upxi_{\alpha} (\geop{x^A})^{\alpha}$ for $A = 2,3$. 
\item We adopt a similar convention for contractions involving higher order tensorfields, 
	i.e., $\gtorus_{AB} = \gtorus(\geop{x^A},\geop{x^B})$ for $A,B = 2,3$. 
\item We sum repeated uppercase Latin indices over $A = 2,3$, e.g., 
$\upxi_{AA} \eqdef \upxi_{22} + \upxi_{33}$.
\end{enumerate}
\end{notation}

\subsubsection{The important acoustic vectorfields}
\label{SSS:IMPORTANTACOUSTICVECTORFIELDS}
The vectorfields in the next definition are fundamental for the rest of the paper.
We will use them to control the solution up to top-order.

\begin{definition}[The important acoustic vectorfields] \label{D:COMVECTORFIELDS} 
\hfill
\begin{enumerate}
\item We define the \emph{geodesic null vectorfield} by:
\begin{align} \label{E:LGEO}
\Lgeo^{\alpha}
& \eqdef - (\gfour^{-1})^{\alpha \beta}\p_\beta u,
\end{align}
and the \emph{rescaled null vectorfield} 
as follows, where $\upmu$ is the inverse foliation density defined in \eqref{E:MUDEF}:
\begin{align} \label{E:LUNIT}
\Lunit 
& \eqdef \upmu \Lgeo.
\end{align}
For $i = 1,2,3$, we define the scalar functions $\Lsmall^i$ as follows, 
where throughout the paper, $\Lunit^i$ denotes the Cartesian component $\Lunit x^i$:
\begin{align} \label{E:LSMALLDEF}
	\Lsmall^1 
	& \eqdef L^1 - 1, \qquad \Lsmall^2 \eqdef L^2,\qquad \Lsmall^3 \eqdef L^3.
\end{align}

\item We define $X$ to be the unique vectorfield that is $\Sigma_t$-tangent and $\gfour$-orthogonal to the 
smooth tori $\ell_{t,u}$, and normalized by: 
\begin{align} \label{E:X}
\gfour(\Lunit,X) 
& = -1,
\end{align}
and we define the rescaled vectorfield $\muX$ by: 
\begin{align} \label{E:MUX}
\muX 
& \eqdef \upmu X.
\end{align}
For $i = 1,2,3$, we define the scalar functions $\Xsmall^i$ as follows, 
where throughout the paper, $X^i$ denotes the Cartesian component $X x^i$:\footnote{For the solutions covered by our main results, the functions $\Lsmall^i$ and $\Xsmall^i$ will be $\ll 1$.}  
\begin{align} \label{E:XSMALL}
\Xsmall^1 
& \eqdef X^1 + 1, \qquad \Xsmall^2 \eqdef X^2, \qquad \Xsmall^3 \eqdef X^3.
\end{align}

\item We define $\Yvf{2}$, $\Yvf{3}$, to be the following $\ell_{t,u}$-tangent vectorfields:
\begin{align} \label{E:YCOMMUTATOR}
\Yvf{2} & \eqdef \p_2 - \gfour(\p_2,X) X,  & \Yvf{3} & \eqdef \p_3 - \gfour(\p_3,X) X.  
\end{align}
We also define $\Yvfsmall{2}$, $\Yvfsmall{3}$, to be the following vectorfields 
(which are not generally $\ell_{t,u}$-tangent):
\begin{align} \label{E:YSMALL}
\Yvfsmall{2} & \eqdef \Yvf{2} - \partial_2, & \Yvfsmall{3} & \eqdef \Yvf{3} - \partial_2. 
\end{align}
	We similarly define the Cartesian component functions $\Yvf{2}^i, \Yvfsmall{3}^i$, 
	analogously to \eqref{E:XSMALL}.
\item We define the \emph{commutation vectorfields} $\Fullset$, the 
$\nullhyparg{u}$-tangential subset $\Tanset$, and the $\ell_{t,u}$-tangential subset
$\Angularset$ as follows:
\begin{align} \label{E:COMMUTATIONVECTORFIELDS}
	\Fullset
	& \eqdef \lbrace \Lunit, \muX, \Yvf{2}, \Yvf{3} \rbrace,
	& 
	\Tanset
	& \eqdef \lbrace \Lunit, \Yvf{2}, \Yvf{3} \rbrace,
	&
	\Angularset
	& 
	\eqdef \lbrace \Yvf{2}, \Yvf{3} \rbrace.
\end{align}
\end{enumerate}
\end{definition}

Lemma~\ref{L:COMMUTATORSTOCOORDINATES} shows that
$\Fullset$ spans the tangent spaces of spacetime equipped with
the differential structure corresponding to the geometric coordinates $(t,u,x^2,x^3)$.
We sometimes refer to $\Fullset$ as the \emph{rescaled frame}
because the vectorfield $\muX = \upmu X$ degenerates with respect to the Cartesian differential
structure as $\upmu \downarrow 0$, i.e., $\muX^i = \upmu X^i$ tends to $0$.
Similarly, the lemma shows
that $\Tanset$ spans the tangent spaces of the characteristics $\nullhyparg{u}$
and that $\Angularset$ spans the tangent spaces of the $\ell_{t,u}$.
To derive $L^{\infty}$ and H\"{o}lder estimates, we commute various PDEs with 
elements of $\Fullset$. To derive energy estimates,
we will commute various PDEs with the elements of $\Tanset$.
For a handful of key estimates, we will refer to the set $\Angularset$.


Throughout the paper, we will often silently use the identities
featured in the following lemma.

\begin{lemma}[Basic properties of the vectorfields] 
\label{L:BASICPROPERTIESOFVECTORFIELDS} 
The follow results hold.

\begin{enumerate}
\item The vectorfield $\Lgeo$ is a geodesic and $\gfour$-null, i.e.,
 \begin{align} \label{E:LGEOISGEODESICANDNULL}
 \gfour(\Lgeo,\Lgeo) & = 0, 
	&  
\Dfour_{\Lgeo}\Lgeo & = 0.
\end{align}
The rescaled vectorfield $\Lunit$ is also $\gfour$-null:
\begin{align} \label{E:LUNITISNULL}
\gfour(\Lunit,\Lunit) 
& = 0,
\end{align}
and it satisfies the following identity,
where $\upmu$ is the inverse foliation density defined in \eqref{E:MUDEF}:
\begin{align} \label{E:COVARIANTLUNITDERIVATIVEOFLUNIT}
	\Dfour_{\Lunit} \Lunit
	& = \frac{(\Lunit \upmu)}{\upmu} \Lunit.
\end{align}

\item $\Lunit$ is $\gfour$-orthogonal to the characteristics $\nullhyparg{u}$, 
that is, for any vectorfield $\Singletan$ tangent to $\nullhyparg{u}$,
	we have:
	\begin{align} \label{E:LUNITISNORMALTOCHARACTERISTICS}
		\gfour(\Lunit,\Singletan)
		& = 0.
	\end{align}

\item The following identities hold:
\begin{align}
\Lunit u & = 0, & \Lunit t = \Lunit^0 & = 1, & \muX u & = 1, & \muX t = \muX^0 &= 0, \label{E:LULTMUXUMUXT} \\
\gfour(X,X) & = 1, & \gfour(\muX,\muX) & = \upmu^2, & \gfour(\Lunit,X) & = -1, & \gfour(\Lunit,\muX) &= -\upmu. \label{E:GOFLANDX}
\end{align}
\item The material vectorfield $\Transport$ is future directed, $\gfour$-orthogonal to $\Sigma_t$ (and hence also to $\ell_{t,u}$), 
and it is of $\gfour$-unit size: 
\begin{align} \label{E:TRANSPORTISUNITLENGTH}
\gfour(\Transport,\Transport) 
& = -1.
\end{align}
Moreover, we have:
\begin{align} \label{E:BISLPLUSX}
\Transport & = \Lunit + X, 
\end{align}
and relative to the Cartesian coordinates, we have:
\begin{align}
\Transport_{\alpha}
& = 
- 
\updelta_\alpha^0, \label{E:MATERIALDERIVATIVELOWEREDCARTESIANCOORDINATES}
\end{align}
where $\updelta_{\alpha}^{\beta}$ is the Kronecker delta. 

\item Finally, the following identities hold for $i=1,2,3$ and $A=2,3$:
\begin{subequations}
\begin{align}
	\Xsmall^i  & = - \Lsmall^i + v^i, 
	\label{E:XSMALLINTERMSOFLSMALLANDVELOCITY} 
		\\
	\Yvfsmall{A}^i
	& = 
			-\Speed^{-2} \Xsmall^A
			X^i
			=
			-\Speed^{-2} 
			(- \Lsmall^A + v^A)
			(-\Lunit^i + v^i).
			\label{E:YSMALLINTERMSOFLSMALLANDVELOCITY}
\end{align}
\end{subequations} 

\end{enumerate}
\end{lemma}


\begin{proof}
All aspects of the lemma except for 
\eqref{E:COVARIANTLUNITDERIVATIVEOFLUNIT}
and
\eqref{E:YSMALLINTERMSOFLSMALLANDVELOCITY} follow from
minor modifications of the proofs of \cite[(2.12), (2.13) and Lemma 2.1]{jSgHjLwW2016}.
The identity \eqref{E:COVARIANTLUNITDERIVATIVEOFLUNIT} follows from
definition \eqref{E:LUNIT},
\eqref{E:LGEOISGEODESICANDNULL},
and the Leibniz rule for the connection $\Dfour$.
The identity \eqref{E:YSMALLINTERMSOFLSMALLANDVELOCITY} follows from definitions
\eqref{E:XSMALL}--\eqref{E:YCOMMUTATOR}, the form \eqref{E:ACOUSTICALMETRIC} of $\gfour_{\alpha \beta}$,
and \eqref{E:XSMALLINTERMSOFLSMALLANDVELOCITY}.
\end{proof}


\subsection{Controlling quantities $\controlvars$ and $\badcontrolvars$} 
\label{SS:CONTROLVARS}
In the next definition, we introduce the solution variable arrays
$\controlvars$ and $\badcontrolvars$. These arrays allow us to provide
simple, schematic formulas in contexts where the precise details are not important
for the PDE analysis.

\begin{definition}[The controlling quantities] \label{D:CONTROLVARS} 
We define $\controlvars$ and $\badcontrolvars$ to be the following arrays of scalar functions:
\begin{subequations} 
\begin{align}
\controlvars 
& 
\eqdef (\wavearray, \Lsmall^1, \Lsmall^2, \Lsmall^3), 
	\label{E:CONTROLVARS} 
	\\
\badcontrolvars 
& \eqdef  
	(\wavearray, \upmu - 1, \Lsmall^1, \Lsmall^2, \Lsmall^3).
	\label{E:BADCONTROLVARS} 
\end{align}
\end{subequations}
\end{definition} 

For the solutions that we study in our main results, along the data hypersurface 
$\Sigma_0$, $\controlvars$ and $\badcontrolvars$ are small in $L^{\infty}$.

\subsection{Identities for the $\ell_{t,u}$-projection tensorfield and the first fundamental forms}
\label{SS:IDENTITIESFORFIRSTFUNDSANDSMOOTHTORUSPROJECTION}
The following lemma provides useful identities for
the $\ell_{t,u}$-projection tensorfield $\smoothtorusproject$,
the first fundamental form $\gtorus$ of $\ell_{t,u}$,
and the first fundamental form $g$ of $\Sigma_t$.


\begin{lemma}[Useful identities for the first fundamental forms]
\label{L:USEFULIDENTITIESFORFIRSTFUNDAMENTALFORM} 
Recall that $\smoothtorusproject$ is the $\ell_{t,u}$-projection tensorfield from Def.\,\ref{D:PROJECTIONTENSORFIELDSANDTANGENCYTOHYPERSURFACES},
that $\gtorus$ is the first fundamental form of $\ell_{t,u}$ from Def.\ref{D:FIRSTFUNDAMENTALFORMS},
and that $g$ is the first fundamental form of $\Sigma_t$ from Def.\ref{D:FIRSTFUNDAMENTALFORMS}.
Let $X$ be the vectorfield defined in Def.\,\ref{D:COMVECTORFIELDS}. 
Then the following identities hold relative to the geometric coordinates ($A,B = 2,3$):
\begin{subequations}
\begin{align}
\gtorus_{AB} & = \Speed^{-2} \updelta_{AB} +\Speed^{-2} \frac{X^AX^B}{(X^1)^2}, 
	\label{E:SMOOTHTORIGABEXPRESSION} 
		\\
(\gtorus^{-1})^{AB}
& 
= \Speed^2 \updelta^{AB} - X^A X^B, 
	\label{E:SMOOTHGINVERSEABEXPRESSION} 
		\\
\det \gtorus & = \frac{1}{\Speed^2 (X^1)^2}. 
\label{E:DETERMSMOOTHGTORUSRELTOGEOMETRICCOORDS}
\end{align}
\end{subequations}

Moreover, the following identities hold relative to arbitrary coordinates:
\begin{subequations}
	\begin{align} \label{E:FIRSTFUNDOFSIGMATINTERMSOFACOUSTICALMETRICANDTRANSPORT}
		g_{\alpha \beta}
		& =
		\gfour_{\alpha \beta} 
		+
		\Transport_{\alpha} \Transport_{\beta},
			\\
		(g^{-1})^{\alpha \beta}
		& 
		=
		(\gfour^{-1})^{\alpha \beta} 
		+
		\Transport^{\alpha} \Transport^{\beta}.
			\label{E:INVERSEFIRSTFUNDOFSIGMATINTERMSOFACOUSTICALMETRICANDTRANSPORT}
	\end{align}
\end{subequations}

Furthermore, the following identities hold relative to the Cartesian coordinates:
\begin{subequations}
	\begin{align}
		g
		& = \Speed^{-2} \sum_{a=1}^3 (\rmD x^a - v^a \rmD t)\otimes  (\rmD x^a - v^a \rmD t),
			\label{E:FIRSTFUNDOFSIGMATIDENTITY} 
				\\
		g^{-1}
		& = \Speed^2 \sum_{a=1}^3 \p_a \otimes \p_a.
		\label{E:INVERSEFIRSTFUNDOFSIGMATIDENTITY}
	\end{align}
\end{subequations}

In addition, the following identities hold relative to arbitrary coordinates:
\begin{subequations}
	\begin{align}
			\gtorus_{\alpha \beta}
			& = 
			g_{\alpha \beta} 
			-
			X_{\alpha} X_{\beta}
			=
			\gfour_{\alpha \beta} 
			+
			\Transport_{\alpha} \Transport_{\beta}
			-
			X_{\alpha} X_{\beta}
			= \gfour_{\alpha \beta} 
				+
				\Lunit_{\alpha} \Lunit_{\beta}
				+
				\Lunit_{\alpha} X_{\beta}
				+
				X_{\alpha} \Lunit_{\beta},
				\label{E:SMOOTHTORUSMETRICINTERMSOFSIGMATMETRICANDX} 
				\\
		(\gtorus^{-1})^{\alpha \beta}
		& = 
			(g^{-1})^{\alpha \beta} 
			-
			X^{\alpha} X^{\beta}
			=
			(\gfour^{-1})^{\alpha \beta} 
			+
			\Transport^{\alpha} \Transport^{\beta}
			-
			X^{\alpha} X^{\beta}
			=
				(\gfour^{-1})^{\alpha \beta}
				+
				\Lunit^{\alpha} \Lunit^{\beta}
				+
				\Lunit^{\alpha} X^{\beta}
				+
				X^{\alpha} \Lunit^{\beta},
			\label{E:SMOOTHTORUSINVERSEMETRICINTERMSOFINVERSESIGMATMETRICANDX}
				\\
		\smoothtorusproject^{\ \alpha}_\beta 
		& = 
		\gfour_{\beta \gamma}
		(\gtorus^{-1})^{\alpha \gamma}
		=
		\gtorus_{\beta \gamma}
		(\gtorus^{-1})^{\alpha \gamma}.
		\label{E:LTUSPROJECTIONISONEINDEXLOWERINGOFGTORUSINVERSE}
	\end{align}
\end{subequations}

Finally, relative to the Cartesian coordinates, the following identities hold for $\alpha, \beta = 0,1,2,3$:
\begin{subequations}
\begin{align} \label{E:GTORUSINVERSE0COMPONENTSVANISH}
		(\gtorus^{-1})^{0 \alpha}
		& 
		=
		(\gtorus^{-1})^{\alpha 0}
		= 0,
			\\
		(g^{-1})^{0 \alpha}
		& 
		=
		(g^{-1})^{\alpha 0}
		= 0,
		\label{E:FIRSTFUNDSIGMAT0COMPONENTSVANISH}
			\\
	\smoothtorusproject^{\ 0}_{\beta}
	& = 0.
	\label{E:SMOOTHTORUSUPPER0COMPONENTSVANISH}
\end{align}
\end{subequations}

\end{lemma}

\begin{proof}
The identities \eqref{E:SMOOTHTORIGABEXPRESSION}--\eqref{E:SMOOTHGINVERSEABEXPRESSION} were proved in \cite[Lemma 2.31]{jLjS2021}. 
The identity \eqref{E:DETERMSMOOTHGTORUSRELTOGEOMETRICCOORDS} follows from \eqref{E:SMOOTHTORIGABEXPRESSION} and the 
following identity:
\begin{align} \label{E:CARTESIANCOMPONENTSOFXSUMTOCSQUARED}
\Speed^2 
& = \sum_{i=1}^3 (X^i)^2,
\end{align}
which follows from \eqref{E:ACOUSTICALMETRIC} and $\gfour(X,X) = 1$ (see \eqref{E:GOFLANDX}). 

\eqref{E:FIRSTFUNDOFSIGMATINTERMSOFACOUSTICALMETRICANDTRANSPORT}
follows from definition \eqref{E:CARTESIANSIGMAPROJECT} 
and the fact that $g = \Sigmatproject \gfour$.
\eqref{E:INVERSEFIRSTFUNDOFSIGMATINTERMSOFACOUSTICALMETRICANDTRANSPORT}
follows from raising the indices in \eqref{E:FIRSTFUNDOFSIGMATINTERMSOFACOUSTICALMETRICANDTRANSPORT}
with $\gfour^{-1}$.

Since $\ell_{t,u} \subset \Sigma_t$, the first equality in \eqref{E:SMOOTHTORUSMETRICINTERMSOFSIGMATMETRICANDX}
follows from the fact that $X$ is $\Sigma_t$-tangent, $\gfour$-orthogonal to $\ell_{t,u}$, and normalized by
$\gfour(X,X) = g(X,X) = 1$. The second equality in \eqref{E:SMOOTHTORUSMETRICINTERMSOFSIGMATMETRICANDX}
follows from \eqref{E:FIRSTFUNDOFSIGMATINTERMSOFACOUSTICALMETRICANDTRANSPORT}.
The last equality in \eqref{E:SMOOTHTORUSMETRICINTERMSOFSIGMATMETRICANDX}
follows from \eqref{E:BISLPLUSX}.
\eqref{E:SMOOTHTORUSINVERSEMETRICINTERMSOFINVERSESIGMATMETRICANDX} follows from raising
the indices in \eqref{E:SMOOTHTORUSMETRICINTERMSOFSIGMATMETRICANDX} with $\gfour^{-1}$.

The first equality in \eqref{E:LTUSPROJECTIONISONEINDEXLOWERINGOFGTORUSINVERSE}
follows from definition \eqref{E:SMOOTHTORUSPROJECT} and the second equality in 
\eqref{E:SMOOTHTORUSINVERSEMETRICINTERMSOFINVERSESIGMATMETRICANDX}.
The second equality in \eqref{E:LTUSPROJECTIONISONEINDEXLOWERINGOFGTORUSINVERSE}
follows from the second equality in \eqref{E:SMOOTHTORUSMETRICINTERMSOFSIGMATMETRICANDX} 
and the fact that $\gtorus^{-1}$ vanishes when contracted against $\Transport$ or $X$.

\eqref{E:INVERSEFIRSTFUNDOFSIGMATIDENTITY}
follows from \eqref{E:INVERSEACOUSTICALMETRIC} 
and the last equality in \eqref{E:SMOOTHTORUSINVERSEMETRICINTERMSOFINVERSESIGMATMETRICANDX}.
\eqref{E:FIRSTFUNDOFSIGMATIDENTITY} follows from \eqref{E:INVERSEFIRSTFUNDOFSIGMATIDENTITY}
and \eqref{E:ACOUSTICALMETRIC}, which in particular implies that in Cartesian coordinates, 
the $\gfour$-dual of $\partial_a$
is $dx^a - v^a dt$.

\eqref{E:GTORUSINVERSE0COMPONENTSVANISH} follows 
from \eqref{E:INVERSEFIRSTFUNDOFSIGMATIDENTITY},
the first equality in \eqref{E:SMOOTHTORUSINVERSEMETRICINTERMSOFINVERSESIGMATMETRICANDX},
and the fact that $X^0 = 0$,
i.e., $X$ is $\Sigma_t$-tangent.
\eqref{E:FIRSTFUNDSIGMAT0COMPONENTSVANISH} follows from \eqref{E:INVERSEFIRSTFUNDOFSIGMATIDENTITY}.
\eqref{E:SMOOTHTORUSUPPER0COMPONENTSVANISH} follows from \eqref{E:LTUSPROJECTIONISONEINDEXLOWERINGOFGTORUSINVERSE}
and \eqref{E:GTORUSINVERSE0COMPONENTSVANISH}.
\end{proof}

\begin{definition}[Metric duality and musical notation]
	\label{D:METRICDUALITY}
	If $Y = Y^B \geop{x^B}$ is an $\ell_{t,u}$-tangent vectorfield, 
	then $Y_{\flat}$ denotes the $\ell_{t,u}$-tangent one-form that is $\gtorus$-dual to $Y$,
	i.e., for $A=2,3$, $(Y_{\flat})_A = \gtorusdoublearg{A}{B} Y^B$. 
	Similarly, if $\upxi$ is an $\ell_{t,u}$-tangent one-form, 
	then $\upxi^{\sharp}$ denotes the $\ell_{t,u}$-tangent vectorfield that is $\gtorus$-dual to $\upxi$,
	i.e., for $A=2,3$, $(\upxi^{\#})^A = (\gtorus^{-1})^{AB} \upxi_B$, where 
	$\upxi_B = \upxi \cdot \geop{x^B}$.
	Similarly, if $\upxi$ is a symmetric type $\binom{0}{2}$ $\ell_{t,u}$-tangent tensorfield, 
	then $\upxi^{\sharp}$ denotes the type $\binom{1}{1}$ $\ell_{t,u}$-tangent tensorfield
	obtained by raising one index of $\upxi$ with $\gtorus^{-1}$,
	while
	$\upxi^{\sharp \sharp}$ denotes the type $\binom{2}{0}$ $\ell_{t,u}$-tangent tensorfield
	obtained by raising both indices of $\upxi$ with $\gtorus^{-1}$.
\end{definition}

\begin{definition}[$\ell_{t,u}$-differential] 
\label{D:ANGULARDIFFERENTIAL}
If $\varphi$ is a scalar function, then we define $\angrmD \varphi$ 
to be the following $\ell_{t,u}$-tangent one-form:
\begin{align} \label{E:ANGULARDIFFERENTIAL}
\angrmD \varphi 
& \eqdef \smoothtorusproject \rmD \varphi.
\end{align}
\end{definition}
Note that $\argangrmD{A} \varphi = \angrmD \varphi \cdot \geop{x^A}  = \geop{x^A} \varphi$ for $A = 2,3$
and that $\argangrmD{\alpha} \varphi 
= 
\angrmD \varphi \cdot \partial_{\alpha} 
=  
\smoothtorusproject^{\ \beta}_{\alpha} \partial_{\beta} \varphi$ 
for $\alpha = 0,1,2,3$.

\begin{definition}[Levi-Civita connections and associated differential operators] \hfill
\label{D:CONNECTIONSANDDIFFERENTIALOPERATORS}
\begin{enumerate}
\item We denote the Levi-Civita connection of $\gfour$ by $\Dfour$. 
\item We denote the Levi-Civita connection of $\gtorus$ by $\newangD$. In particular,
	for $\ell_{t,u}$-tangent tensorfields $\upxi$, we have $\newangD \upxi = \smoothtorusproject \Dfour \upxi$.
\item If $\upxi$ is an $\ell_{t,u}$-tangent one-form, then we define its $\ell_{t,u}$-divergence 
	to be the scalar function
	$\angdiv \upxi \eqdef \gtorus^{-1}\cdot \newangD \upxi$. Similarly, if $Y$ is an $\ell_{t,u}$-tangent vectorfield, 
	then we define its $\ell_{t,u}$-divergence to be the scalar function $\angdiv Y \eqdef \gtorus^{-1} \cdot \newangD Y_{\flat}$.
	where $Y_{\flat}$ is the $\ell_{t,u}$-tangent one-form that is $\gfour$-dual to $Y$.
\item If $\upxi$ is a symmetric type $\binom{0}{2}$ $\ell_{t,u}$-tangent tensorfield, 
	then we define its $\ell_{t,u}$-divergence 
	to be the $\ell_{t,u}$-tangent one-form with the following $\ell_{t,u}$ components for $A=2,3$:
	$(\angdiv \upxi)_A \eqdef (\gtorus^{-1})^{BC} \cdot \newangDarg{B} \upxi_{CA}$. 
\item We denote the covariant wave operator of $\gfour$ by $\Box_{\gfour} \eqdef \gfour^{-1}\cdot\Dfour^2
	= (\gfour^{-1})^{\alpha \beta} \Dfour_{\alpha} \Dfour_{\beta}
$. 
\item We denote the $\ell_{t,u}$-Laplacian associated to $\gtorus$ by 
$\angLap \eqdef \gtorus^{-1} \cdot \newangD^2
=
(\gtorus^{-1})^{AB} \angDarg{A} \angDarg{B}
$.
\end{enumerate}
\end{definition}


\begin{definition}[Projected Lie derivatives] 
\label{D:PROJECTEDLIEDERIVATIVES}
Given a spacetime tensorfield $\upxi$
and a vectorfield $Z$,
we define $\SigmatLie_Z \upxi$ and $\angLie_Z \upxi$ to respectively be the following $\Sigma_t$-tangent and $\ell_{t,u}$-tangent tensorfields:
\begin{align} \label{E:PROJECTEDLIEDERIVATIVES}
\SigmatLie_Z \upxi 
	&	\eqdef \Sigmatproject \Lie_Z \upxi, 
&
\angLie_Z \upxi 
& \eqdef \smoothtorusproject \Lie_Z \upxi. 
	\end{align}
\end{definition}

We will use the following simple commutation lemma when deriving various equations.

\begin{lemma}[Angular differential $\angrmD$ commutes with $\angLie$] 
\label{L:ANGULARDIFFERENTIALCOMMUTESWITHANGLIE}
Let $f$ be a scalar function and let 
$Z \in \Fullset$ (see definition \eqref{E:COMMUTATIONVECTORFIELDS}). 
Then the following identity holds:
\begin{align} \label{E:ANGULARDIFFERENTIALCOMMUTESWITHANGLIE}
\angLie_Z \angrmD f 
& = \angrmD Z f.
\end{align}
\end{lemma}

\begin{proof}
The same proof of \cite[Lemma 2.10]{jSgHjLwW2016} holds.
\end{proof}


\subsection{Traces of tensorfields} 
\label{SS:TRACEOFTENSORS}
In our analysis, we will encounter various traces of tensorfields.

\begin{definition}[The trace of spacetime tensors and $\ell_{t,u}$-tangent tensorfields]
\label{D:TRACEOFTENSORS} \hfill
\begin{enumerate}
\item If $\upxi$ is a type $\binom{0}{2}$ spacetime tensorfield, 
then we define its $\gfour$-trace as follows:
\begin{subequations}
\begin{align} \label{E:SPACETIMETRACE}
	\mytr_{\gfour}\upxi 
	& \eqdef (\gfour^{-1})^{\alpha \beta} \upxi_{\alpha \beta}.
\end{align}
\item If $\upxi$ is a type $\binom{0}{2}$ spacetime tensorfield, then we define its $\gtorus$-trace as follows:
\begin{align} \label{E:SMOOTHTORUSTRACE}
\mytr_{\gtorus} \upxi 
& \eqdef 
(\gtorus^{-1})^{\alpha \beta} \upxi_{\alpha \beta}.
\end{align}
\end{subequations}
\end{enumerate}
\end{definition} 

\subsection{Pointwise norms and semi-norms of tensorfields}
\label{SS:POINTWISESEMINORMSOFTENSORS}
In the next definition, we define various pointwise norms and semi-norms
that we will use to measure the size of tensorfields.

\begin{definition}[Pointwise norms] \label{D:POINTWISESEMINORMS} \hfill
\begin{enumerate}
\item If $\upxi$ is a type $\binom{m}{n}$ spacetime tensorfield
such that $\gfour_{\alpha_1 \widetilde{\alpha}_1} \cdots \gfour_{\alpha_m \widetilde{\alpha}_m} (\gfour^{-1})^{\beta_1 \widetilde{\beta}_1} \cdots (\gfour^{-1})^{\beta_n \widetilde{\beta}_n} \upxi_{\beta_1\cdots\beta_n}^{\alpha_1\cdots\alpha_n} \upxi_{\widetilde{\beta}_1 \cdots \widetilde{\beta}_n}^{\widetilde{\alpha}_1\cdots\widetilde{\alpha}_m} \geq 0$,
 then we define $|\upxi|_{\gfour}\geq 0$ by:
\begin{subequations}
\begin{align} \label{E:SQUAREPOINTWISELORENTZIANNORMWITHRESPECTTOACOUSTICALMETRIC}
|\upxi|_{\gfour}^2 
& \eqdef \gfour_{\alpha_1 \widetilde{\alpha}_1} \cdots \gfour_{\alpha_m \widetilde{\alpha}_m} (\gfour^{-1})^{\beta_1 \widetilde{\beta}_1} \cdots (\gfour^{-1})^{\beta_n \widetilde{\beta}_n} \upxi_{\beta_1\cdots\beta_n}^{\alpha_1\cdots\alpha_n} \upxi_{\widetilde{\beta}_1 \cdots \widetilde{\beta}_n}^{\widetilde{\alpha}_1\cdots\widetilde{\alpha}_m}.
\end{align}
\item If $\upxi$ is a type $\binom{m}{n}$ tensorfield, then we define $|\upxi|_{\gtorus} \geq 0$ by:
\begin{align} \label{E:SQUAREPOINTWISESEMINORMWITHRESPECTTOFIRSTFUNDOFSMOOTHTORI}
|\upxi|_{\gtorus}^2 
& \eqdef \gtorus_{\alpha_1 \widetilde{\alpha}_1} \cdots \gtorus_{\alpha_m \widetilde{\alpha}_m} (\gtorus^{-1})^{\beta_1 \widetilde{\beta}_1} \cdots (\gtorus^{-1})^{\beta_n \widetilde{\beta}_n} \upxi_{\beta_1\cdots\beta_n}^{\alpha_1\cdots\alpha_n} \upxi_{\widetilde{\beta}_1 \cdots \widetilde{\beta}_n}^{\widetilde{\alpha}_1\cdots\widetilde{\alpha}_m}.
\end{align}
\item If $\upxi$ is a type $\binom{m}{n}$ tensorfield, then we define $|\upxi|_{g} \geq 0$ by:
\begin{align} \label{E:SQUAREPOINTWISESEMINORMWITHRESPECTTOFIRSTFUNDOFSIGMAT}
|\upxi|_{g}^2 \eqdef g_{\alpha_1 \widetilde{\alpha}_1} \cdots g_{\alpha_m \widetilde{\alpha}_m} (g^{-1})^{\beta_1 \widetilde{\beta}_1} \cdots (g^{-1})^{\beta_n \widetilde{\beta}_n} \upxi_{\beta_1\cdots\beta_n}^{\alpha_1\cdots\alpha_n} \upxi_{\widetilde{\beta}_1 \cdots \widetilde{\beta}_n}^{\widetilde{\alpha}_1\cdots\widetilde{\alpha}_m}.
\end{align}
\end{subequations}
\end{enumerate}
\end{definition}

\begin{remark}[Norms vs.\ semi-norms]
\label{R:NORMSANDSEMINORMS}
$|\cdot|_{\gfour}$ is a pointwise norm on the space of $\gfour$-spacelike tensorfields.
$|\cdot|_g$ is a pointwise norm on the space of $\Sigma_t$-tangent tensorfields and a pointwise seminorm on the space of all 
tensorfields.
$|\cdot|_{\gtorus}$ is a pointwise norm on the space of $\ell_{t,u}$-tangent tensorfields and a pointwise seminorm on the space of all 
tensorfields.

Similarly, the function $|\cdot|_{\gtorusroughfirstfund}$ from Def.\,\ref{D:POINTWISESEMINORMWITHRESPECTTOFIRSTFUNDOFROUGHTORI} below
is a pointwise norm on the space of $\twoargroughtori{\timefunction,u}{\muxmulevelsetvalue}$-tangent tensorfields
a pointwise semi-norm on the space of all tensorfields.
\end{remark}

\begin{remark}[Omitting the $0$ component in Cartesian coordinates]
\label{R:OMITTINGZEROCOMPONENTINCARTESIANCOORDINATES}
In view of 
\eqref{E:MATERIALDERIVATIVELOWEREDCARTESIANCOORDINATES},
Def.\,\ref{D:PROJECTIONTENSORFIELDSANDTANGENCYTOHYPERSURFACES},
definitions 
\eqref{E:SQUAREPOINTWISESEMINORMWITHRESPECTTOFIRSTFUNDOFSMOOTHTORI}--\eqref{E:SQUAREPOINTWISESEMINORMWITHRESPECTTOFIRSTFUNDOFSIGMAT}, 
and
\eqref{E:GTORUSINVERSE0COMPONENTSVANISH}--\eqref{E:FIRSTFUNDSIGMAT0COMPONENTSVANISH},
we see that if $\upxi$ is a type $\binom{m}{n}$ $\Sigma_t$-tangent tensorfield, 
then relative to the Cartesian coordinates, we have:
\begin{align} \label{E:NO0COMPONENTSQUAREPOINTWISESEMINORMWITHRESPECTTOFIRSTFUNDOFSMOOTHTORI}
|\upxi|_{\gtorus}^2 
& 
= 
\gtorus_{a_1 \widetilde{a}_1} 
\cdots 
\gtorus_{a_m \widetilde{a}_m} 
(\gtorus^{-1})^{b_1 \widetilde{b}_1} 
\cdots 
(\gtorus^{-1})^{b_n \widetilde{b}_n} 
\upxi_{b_1 \cdots b_n}^{a_1 \cdots a_n}
\upxi_{\widetilde{b}_1 \cdots \widetilde{b}_n}^{\widetilde{a}_1 \cdots \widetilde{a}_m},
	\\
|\upxi|_g^2 
& = 
g_{a_1 \widetilde{a}_1} 
\cdots 
g_{a_m \widetilde{a}_m} 
(g^{-1})^{b_1 \widetilde{b}_1} 
\cdots 
(g^{-1})^{b_n \widetilde{b}_n} 
\upxi_{b_1 \cdots b_n}^{a_1 \cdots a_n} 
\upxi_{\widetilde{b}_1 \cdots \widetilde{b}_n}^{\widetilde{a}_1 \cdots\widetilde{a}_m},
\label{E:NO0COMPONENTSQUAREPOINTWISESEMINORMWITHRESPECTTOFIRSTFUNDOFSIGMAT}
\end{align}
i.e., we can omit all ``$0$'' components on 
RHSs~\eqref{E:NO0COMPONENTSQUAREPOINTWISESEMINORMWITHRESPECTTOFIRSTFUNDOFSMOOTHTORI}--\eqref{E:NO0COMPONENTSQUAREPOINTWISESEMINORMWITHRESPECTTOFIRSTFUNDOFSIGMAT}.

Similarly, taking into account definition~\eqref{E:PROJECTIONOFTENSORONTOFLATTORUS} and \eqref{E:SMOOTHTORUSUPPER0COMPONENTSVANISH}, 
we see that if $\SigmatTan$ is a $\Sigma_t$-tangent vectorfield, 
then relative to the Cartesian coordinates, we have
$\smoothtorusproject_{\beta}^{\ \alpha} \partial_{\alpha} \SigmatTan^{\beta}
=
\smoothtorusproject_b^{\ a} \partial_a \SigmatTan^b
$.

In the rest of the paper, we will use these basic facts without always explicitly mentioning them.
\end{remark}


\subsection{Second fundamental forms and the torsion}
\label{SS:SECONDFUNDAMENTALFORMSANDTORSION}
In this section, we provide the standard definitions 
of the second fundamental form $k$ of $\Sigma_t$,
the null second fundamental form $\upchi$ of $\ell_{t,u}$,
and the one-form $\upzeta$.
These quantities will appear in various PDEs throughout the article.
It is well-known that there are many technical difficulties that have to be overcome to obtain
top-order energy estimates for
$\mytr_{\gtorus}\upchi$ and $\upchi$.
To achieve control, we will use the modified quantities defined in Sect.\,\ref{S:CONSTRUCTIONOFMODIFIEDQUANTITIES}
and elliptic estimates on the rough tori, which we derive in Sect.\,\ref{S:ELLIPTICESTIAMTESACOUSTICGEOMETRYONROUGHTORI}.

\begin{definition}[The second fundamental forms $k$ and $\upchi$, and the one-form $\upzeta$] 
\label{D:SECONDFUNDAMENTALFORMSANDZETAONEFORM}
\hfill
\begin{enumerate}
\item We define the \emph{second fundamental form} $k$ of $\Sigma_t$ as follows:
\begin{align} \label{E:SIGMATSECONDFUND}
	k 
	& 
	\eqdef \frac{1}{2} \SigmatLie_{\Transport} g.
\end{align}
\item We define the \emph{null second fundamental form} of $\ell_{t,u}$ as follows: 
\begin{align} \label{E:NULLSECONDFUND}
	\upchi  
	& 
	\eqdef \frac{1}{2} \angLie_\Lunit \gtorus.
\end{align}
\item We define $\upzeta$ to be the $\ell_{t,u}$-tangent one-form with the following components:
\begin{align} \label{E:TORISONTENSORFIELD}
	\upzeta_A 
	& 
	\eqdef 
	\gfour(\Dfour_A \Lunit, X).
\end{align}
\end{enumerate}
\end{definition}

\subsection{Transport equations for the eikonal function quantities}
\label{SS:TRANSPORTFORMUANDLUNITI}
To control the eikonal function quantities $\upmu$ and $\Lunit^i$, we will use the following transport equations.

\begin{lemma}[Transport equations satisfied by $\upmu$ and $\Lunit^i$] \label{L:TRANSPORTMUANDLUNITI}
$\upmu$ and $\Lunit^i$ satisfy the following transport equations:
\begin{align} \label{E:MUTRANSPORT}
\Lunit \upmu 
& 
= 
\frac{1}{2} \vec{G}_{\Lunit \Lunit} \diamond \muX \wavearray 
- 
\frac{1}{2} \upmu \vec{G}_{\Lunit \Lunit} \diamond \Lunit \wavearray 
- 
\upmu \vec{G}_{\Lunit X} \diamond \Lunit \wavearray, 
	\\
\Lunit \Lsmall^i 
& 
= \frac{1}{2} (\vec{G}_{\Lunit \Lunit}\diamond \Lunit \wavearray) X^i 
	- 
	(\angG_{\Lunit}^{\#} \diamond \Lunit \wavearray)\cdot \angrmD x^i 
	+ 
	\frac{1}{2} (\vec{G}_{\Lunit \Lunit} \diamond \angrmD^{\#}\wavearray)\cdot \angrmD x^i. 
	\label{E:LUNITITRANSPORT}
\end{align}
\end{lemma}
\begin{proof}
The same proof of \cite[Lemma 2.12]{jSgHjLwW2016} holds.
\end{proof}

\subsection{The factor driving the shock formation and formulas involving $G_{\Lunit \Lunit}$}
\label{SS:THEFACTORDRIVINGTHESHOCKFORMATION}
In the following lemma, we compute
an expression for the product $\frac{1}{2} \vec{G}_{\Lunit \Lunit} \diamond \muX \wavearray$
on the RHS of the evolution equation \eqref{E:MUTRANSPORT} for $\upmu$.
For every smooth equation of state besides that of a Chaplygin gas, 
there exist open sets of background densities $\overline{\varrho} > 0$ such that
the non-degeneracy condition \eqref{E:NONDEGENCONDITION} holds.
The identity \eqref{E:IDENTITYFORMAINTERMDRIVINGTHESHOCK} then shows
that for solutions that are close to the trivial solution $\wavearray \equiv 0$,
the expansion of $\frac{1}{2} \vec{G}_{\Lunit \Lunit} \diamond \muX \wavearray$
features a non-zero term proportional to $\muX \RRiemann$; 
the presence of this term is crucial for our main results, as it drives the formation of the shock, 
i.e., it drives $\upmu$ to $0$.
In contrast, for the equation of state $p = C_0 - C_1 \exp(- \LogDensity)$ of a Chaplygin gas, 
one can compute that $\Speed^{-1} \Speed_{;\LogDensity} + 1 \equiv 0$, and the non-degeneracy condition
\eqref{E:NONDEGENCONDITION} is therefore impossible. 
In this case, equation \eqref{E:IDENTITYFORMAINTERMDRIVINGTHESHOCK} shows that
the product $\frac{1}{2} \vec{G}_{\Lunit \Lunit} \diamond \muX \wavearray$ 
does not depend on the solution's $\muX$ derivative, and hence our main results do not apply. 
We note that one can show that for irrotational and isentropic solutions,
the equation $\Speed^{-1} \Speed_{;\LogDensity} + 1 = 0$ 
is equivalent to the statement that the quasilinear wave equation for a potential function
satisfies Klainerman's null condition \cite{sK1984}.

\begin{lemma}[Identity for the factor driving the shock formation] 
\label{L:IDENTIFYFORFACTORDRIVINGSHOCK}
For solutions to the compressible Euler equations
\eqref{E:BVIEVOLUTION}--\eqref{E:BENTROPYEVOLUTION},
the following identity holds,
where $\almostRiemannfunction(\LogDensity,\Ent)$ 
is the scalar function from \eqref{E:ALMOSTRIEMANNINVARIANTS}:
\begin{align} \label{E:IDENTITYFORMAINTERMDRIVINGTHESHOCK}
\begin{split}
\frac{1}{2} \vec{G}_{\Lunit \Lunit} \diamond \muX \wavearray 
& = - \frac{1}{2} \Speed^{-1}(\Speed^{-1} \Speed_{;\LogDensity} + 1)  
			\left\lbrace 
				\muX \RRiemann - \muX \LRiemann 
			\right\rbrace 
	\\
& - \frac{1}{2} \upmu \Speed^{-1} X^1 
\left\lbrace \Lunit \RRiemann + \Lunit \LRiemann\} - \upmu \Speed^{-2} \lbrace X^2 \Lunit v^2 + X^3 \Lunit v^3
\right\rbrace 
	\\
& - \upmu \Speed^{-1} \Speed_{;\Ent} X^a \GradEnt^a + \upmu \Speed^{-1} (\Speed^{-1} c_{;\LogDensity} + 1) F_{;\Ent} X^a \GradEnt^a.
\end{split}
\end{align}
\end{lemma}

\begin{proof}
The same proof of Lemma \cite[Lemma 2.37]{jLjS2021} holds.
\end{proof}

In the next lemma, we derive expressions for $G_{\Lunit \Lunit}^2$ and $G_{\Lunit \Lunit}^3$.
When deriving estimates, we will use the expressions to track smallness.

\begin{lemma}[Formulas for $G_{\Lunit \Lunit}^A$] 
\label{L:GLLAEXPRESSION}
The following identities hold for $A= 2,3$:
\begin{align}
G_{\Lunit \Lunit}^A 
& = 2 \Speed^{-2}(v^A - \Lunit^A) = 2\Speed^{-2} X^A. \label{E:GLLAEXPRESSION}
\end{align}
\end{lemma}

\begin{proof}
The identities follow from  
the expression \eqref{E:ACOUSTICALMETRIC} for the Cartesian component $\gfour_{\alpha \beta}$
viewed as a function of $(\RRiemann,\LRiemann,v^2,v^3,\Ent)$,
the identity $G_{\Lunit \Lunit}^A = (\frac{\p}{\p v^A} \gfour_{\alpha \beta}) \Lunit^{\alpha} \Lunit^{\beta}$,
and the identities $\Lunit^0 = 1$ and $\Lunit^A + X^A = v^A$,
which follow from Lemma~\ref{L:BASICPROPERTIESOFVECTORFIELDS} and \eqref{E:MATERIALDERIVATIVEVECOTRFIELD}.
\end{proof}

\subsection{Useful geometric decompositions}
\label{SS:USEFULGEOMETRIC}
In this section, we provide some geometric decompositions that we will use throughout the article.

We start with the following alternate expressions for $\upchi$ and $\angk$, which are useful for computations.

\begin{lemma}[Alternate expressions for $\upchi$ and $\angk$]
\label{L:ALTERNATEEXPRESSIONSFORSECONDFUNDAMENTALFORMS}
The following identities hold:
\begin{align} \label{E:SECONDFUNDSALTERNATE} 
\upchi_{AB} 
& 
= 
\gfour\left(\Dfour_A \Lunit, \geop{x^B}\right), 
&
\angk_{AB} 
& = \gfour\left(\Dfour_A \Transport, \geop{x^B}\right). 
\end{align}
\end{lemma}

\begin{proof}
The same proof of \cite[Lemma 3.61]{jS2016b} holds in this setting.
\end{proof}
We will use the following identities and decompositions when deriving estimates for 
$\upchi$, $\angk$, and $\upzeta$.
\begin{lemma}[Useful identities and decompositions for $\upchi$, $\angk$, and $\upzeta$] 
\label{L:USEFULIDENTITIESANDDECOMPOSITIONSFORSECONDFUNDAMENTALFORMSANDTORSION} 
The following\footnote{Here, $\otimesarray$ is defined by 
$\angG_{\Lunit} \otimesarray \angrmD \wavearray \eqdef \sum_{\iota = 0}^4 
\angG_{\Lunit}^\iota \otimes \angrmD \Psi_\iota$, 
and similarly for $\angrmD \Psi \otimesarray\angG_{\Lunit}$, $\angG_X\otimesarray \angrmD \wavearray$, etc.} identities hold:
\begin{subequations}
\begin{align}
	\upchi & = \gfour_{ab} \angrmD \Lunit^a \otimes \angrmD x^b + \frac{1}{2} \angG \diamond \Lunit \wavearray + \frac{1}{2} \angrmD \wavearray \otimesarray \angG_{\Lunit} - \frac{1}{2} \angG_{\Lunit} \otimesarray \angrmD \wavearray, 
		\label{E:CHIEXPRESSSIONINTERMSOFDERIVATIVESOFLUNITI}
	\\
	\mytr_{\gtorus}\upchi & = \gfour_{ab} \gtorus^{-1} \cdot \left\lbrace \angrmD \Lunit^a \otimes \angrmD x^b\right\rbrace + \frac{1}{2} \gtorus^{-1} \cdot \angG \diamond \Lunit \wavearray. \label{E:TRCHIEXPRESSSIONINTERMSOFDERIVATIVESOFLUNITI}
\end{align}
\end{subequations}

Moreover, we can decompose $\angk$ and $\upzeta$ into $\upmu^{-1}$-singular and $\upmu^{-1}$-regular pieces as follows:
\begin{subequations}
\begin{align} \label{E:CONNECTIONCOEFFICIENTDECOMPOSITIONS}
\upzeta 
& = \zetatan + \upmu^{-1} \zetatrans, & \angk & = \angktan + \upmu^{-1}\angktrans,
\end{align}
where:
\begin{align}
\angktan & \eqdef  
\frac{1}{2} \angG \diamond \Lunit \wavearray 
- 
\frac{1}{2} \angG_{\Lunit} \otimesarray \angrmD \wavearray 
- 
\frac{1}{2} \angrmD \wavearray \otimesarray \angG_{\Lunit} 
- 
\frac{1}{2} \angG_X \otimesarray \angrmD \wavearray 
- 
\frac{1}{2} \angrmD \wavearray \otimesarray \angG_X, \label{E:ANGKTAN}
	\\
\angktrans & \eqdef 
\frac{1}{2} \upmu^{-1} \angG\diamond \muX \wavearray, 
\label{E:ANGKTRANS} 
\\
\zetatan 
& \eqdef 
\frac{1}{2} \angG_X \diamond \Lunit \wavearray 
-
\frac{1}{2} \vec{G}_{\Lunit X} \diamond \angrmD \wavearray 
- 
\frac{1}{2} \vec{G}_{XX} \diamond \angrmD \wavearray, 
	\label{E:ZETATAN} 
	\\
\zetatrans 
& 
\eqdef 
- \frac{1}{2} \upmu^{-1} \angG_{\Lunit} \diamond \muX \wavearray. 
\label{E:ZETATRANS}
\end{align}
\end{subequations}
\end{lemma}

\begin{proof}
The same proofs of \cite[Lemmas 2.13, 2.15]{jLjS2018} holds with minor modifications accounting for the third spatial dimension.
\end{proof}


In the next lemma, we decompose the $\upmu$-weighted covariant 
wave operator relative to the rescaled frame $\{\Lunit,\muX, \Yvf{2}, \Yvf{3}\}$ 
and the second fundamental forms. 

\begin{lemma}[Frame decomposition of $\upmu \square_{\gfour(\wavearray)} f$] 
\label{L:FRAMEDCOMPOSITOINOFMUBOXG}
Let $f$ be a scalar function. Then relative to the 
rescaled frame $\{\Lunit,\muX, \Yvf{2}, \Yvf{3}\}$, 
the following identities hold:
\begin{subequations}
\begin{align}
\begin{split}
\upmu \square_{\gfour(\wavearray)} f 
& = - \Lunit (\upmu \Lunit f + 2 \muX f) 
	+ 
	\upmu \angLap f 
	- 
	(\mytr_{\gtorus} \upchi) \muX f
	\label{E:BOXDECOMPLOUTSIDE} 
		\\ 
& - 
\upmu \mytr_{\gtorus} \angk \Lunit f 
- 
2 \upmu \upzeta^{\#} \cdot \angrmD f,
\end{split}
	\\
\begin{split}
\upmu \square_{\gfour(\wavearray)} f 
& = 
- (\upmu \Lunit + 2 \muX)(\Lunit f) 
+ 
\upmu \angLap f 
- 
(\mytr_{\gtorus} \upchi)\muX f 
- 
(\Lunit \upmu) \Lunit f 
\label{E:BOXDECOMPMULBAROUTSIDE} 
	\\
& - \upmu \mytr_{\gtorus} \angk \Lunit f + 
2 \upmu \upzeta^{\#} \cdot \angrmD f 
+ 
2 (\angrmD^{\#} \upmu) \cdot \angrmD f. 
\end{split}
\end{align}
\end{subequations}
\end{lemma}

\begin{proof}
The same proof of \cite[Proposition 5.4]{jS2016b} holds with $\muX$ 
in the role of the vectorfield denoted by ``$\breve{R}$'' there. 
\end{proof}


\section{Rough time functions, rough adapted coordinates, and rough subsets}
\label{S:ROUGHTIMEFUNCTIONANDROUGHSUBSETS}
The geometric coordinates $(t,u,x^2,x^3)$ from Sect.\,\ref{S:ACOUSTICGEOMETRYANDCOMMUTATORVECTORFIELDS}
are fundamental for our construction of commutation and multiplier vectorfields.
However, these coordinates, in particular the Cartesian time function $t$, 
are not adapted to the shape of the singular boundary.
For this reason,
in this section, we construct a one-parameter family of rough time functions 
$\lbrace \timefunctionarg{\muxmulevelsetvalue} \rbrace_{\muxmulevelsetvalue \in [0,\muxmulevelsetvalue_0]}$
that \emph{are} adapted to structure of the singular boundary,
where $\muxmulevelsetvalue_0 > 0$ is a constant depending on the initial data on $\Sigma_0$. 
We refer to $(\timefunctionarg{\muxmulevelsetvalue},u,x^2,x^3)$ as \emph{rough adapted coordinates}.

In our forthcoming PDE analysis, we will derive estimates on
the level sets of the $\timefunctionarg{\muxmulevelsetvalue}$, 
which we will prove are $\gfour$-spacelike.
We construct $\timefunctionarg{\muxmulevelsetvalue}$ by
solving a well-chosen transport equation (see Def.\,\ref{D:ROUGHTIMEFUNCTION}) with
data equal to $- \upmu$ on the  ``initial hypersurface"
$\lbrace \muX \upmu = -\muxmulevelsetvalue \rbrace$.
Standard well-posedness and Cauchy stability results 
(see Appendices~\ref{A:PS} and \ref{A:OPENSETOFDATAEXISTS})
imply that for perturbations of simple isentropic plane-symmetric solutions,
our construction is well-defined
on short ``rough time'' intervals of the form $[\timefunction_0,\timefunctionboot)$,
where $\timefunction_0 < 0$ is a data-dependent constant (independent of $\muxmulevelsetvalue \in \muxmulevelsetvalue_0$)
and $\timefunction_0 < \timefunctionboot < 0$ is a ``bootstrap parameter.''
Our main results will show that each $\timefunctionarg{\muxmulevelsetvalue}$ exists
in a neighborhood of the singular boundary and has range $[\timefunction_0,0]$,

Finally, we remark that there are several subtleties tied to the analysis and regularity of
$\timefunctionarg{\muxmulevelsetvalue}$,
in particular since our main results crucially rely on our proofs that $\timefunctionarg{\muxmulevelsetvalue}$
is one degree more differentiable than the initial hypersurface
$\lbrace \muX \upmu = -\muxmulevelsetvalue \rbrace$
and that near the singular boundary, 
for $\mulevelsetvalue$ sufficiently small and positive,
$\lbrace \muX \upmu = -\muxmulevelsetvalue \rbrace
\cap
\lbrace \upmu = \mulevelsetvalue \rbrace
$
is an embedded two-dimensional, spacelike torus with sufficient regularity.
The proofs of these results and many supporting ones
are located in
Sects.\,\ref{S:EMBEDDINGSANDFLOWMAPS}--\ref{S:SHARPCONTROLOFMUANDPROPERTIESOFCHOVGEOTOCARTESIAN}.

\begin{quote}
	Starting now, we consider a fixed $\muxmulevelsetvalue \in [0,\muxmulevelsetvalue_0]$;
	see Sect.\,\ref{S:PARAMETERSANDSIZEASSUMPTIONSANDCONVENTIONSFORCONSTANTS} 
	for discussion of how $\muxmulevelsetvalue_0$ is tied to the initial data.
	$\muxmulevelsetvalue$ will remain fixed until Sect.\,\ref{S:MAINRESULTS},
	where we provide our main results.
	All of our estimates will involve constants that 
	can be chosen to be uniform with respect to $\muxmulevelsetvalue$
	over the interval $[0,\muxmulevelsetvalue_0]$.
\end{quote}

\subsection{Basic constructions} 
\label{SS:BASICCONSTRUCTIONSFORROUGHTIMEFUNCTION}
We now introduce some basic ingredients that we will use to construct $\timefunctionarg{\muxmulevelsetvalue}$.

\subsubsection{The constant $\interestingu > 0$, the cut-off function $\phi$, and the vectorfield $\Wtransarg{\muxmulevelsetvalue}$}
\label{SSS:CUTOFFANDWTRANS}
Our constructions rely on the vectorfield $\Wtransarg{\muxmulevelsetvalue}$
featured in the next definition, which is crucial for the rest of the paper.
In what follows, $\interestingu > 0$ is a positive constant;
see Sect.\,\ref{S:PARAMETERSANDSIZEASSUMPTIONSANDCONVENTIONSFORCONSTANTS} for further discussion of how 
the specific choice of $\interestingu$ that we make in our main theorem is tied to the initial data.

\begin{definition}[The cut-off $\phi$ and the vectorfield $\Wtransarg{\muxmulevelsetvalue}$] \label{D:WTRANSANDCUTOFF} 
Let $\psi \colon \R \to [0,1]$ be a fixed smooth cut-off function such that
$\psi(u) = 1$ when when $|u| \leq \frac{1}{2}$ and $\psi(u) = 0$ when $|u| \geq 1$.
Let $\phi$ be the cut-off function defined by
$\phi(u) \eqdef \psi\left(\frac{u}{\interestingu} \right)$.
In particular:
\begin{align} \label{E:CUTOFFFUNCTION}
\begin{cases}
\phi(u) = 1, & \mbox{when } |u| \leq  \frac{1}{2} \interestingu, 
	\\
0 \leq \phi(u) \leq 1, & \mbox{when } \frac{1}{2} \interestingu \leq |u| \leq \interestingu, 
	\\
\phi(u) = 0, & \mbox{when } |u| \geq \interestingu.
\end{cases} 
\end{align}

For fixed $\muxmulevelsetvalue \geq 0$, we define the \textbf{rough transversal vectorfield} $\Wtransarg{\muxmulevelsetvalue}$ as follows:
\begin{align}\label{E:WTRANSDEF}
\Wtransarg{\muxmulevelsetvalue} 
& 
\eqdef 
\muX 
+ 
\phi 
\frac{\muxmulevelsetvalue}{\Lunit \upmu} \Lunit.
\end{align}
\end{definition}

\subsubsection{$\upmu$-adapted subsets}
\label{SSS:MUADAPATEDSUBSETS}
To follow the solution up to the singular boundary, we will analyze it
on the following subsets (and others as well), 
which are adapted to the shape of the singular boundary.

\begin{definition}[Level sets of $\upmu$ and $\muX \upmu$ and the $\upmu$-adapted tori $\twoargmumuxtorus{\mulevelsetvalue}{-\muxmulevelsetvalue}$] 
\label{D:LEVELSETSOFMUANDXMUANDMUMUXTORI}
Recall that $\muX$ is defined in \eqref{E:MUX}.
Given real numbers $\mulevelsetvalue, \muxmulevelsetvalue \geq 0$, we define:
\begin{subequations}
\begin{align} \label{E:LEVELSETSOFMU}
		\mulevelsetarg{\mulevelsetvalue}
		& \eqdef
		\left\lbrace
			(t,u,x^2,x^3) 
			\in \R \times \R \times \T^2 
			\ | \
			\upmu(t,u,x^2,x^3) = \mulevelsetvalue
		\right\rbrace
		\cap
		\lbrace |u| \leq \interestingu \rbrace,
			\\
		\datahypfortimefunctionarg{-\muxmulevelsetvalue}
		& \eqdef
		\left\lbrace
			(t,u,x^2,x^3) 
			\in \R \times \R \times \T^2 
			\ | \
		\muX \upmu(t,u,x^2,x^3) = - \muxmulevelsetvalue
		\right\rbrace
		\cap
		\lbrace |u| \leq \interestingu \rbrace,
			\label{E:LEVELSETSOFMUXMU} \\
		\twoargmumuxtorus{\mulevelsetvalue}{-\muxmulevelsetvalue}
		& \eqdef 
		\mulevelsetarg{\mulevelsetvalue} 
		\cap
		\datahypfortimefunctionarg{-\muxmulevelsetvalue}.
		\label{E:MUXMUTORI}
\end{align}
\end{subequations}
\end{definition}

\begin{remark}[The values of $\mulevelsetvalue$ and $\muxmulevelsetvalue$ featured in our main results]
	\label{R:LAMBDAANDKAPPAINOURMAINRESULTS}
	Our main results concern solutions and values of $\mulevelsetvalue \in [0,\mupositive]$ and 
	$\muxmulevelsetvalue \in [0,\muxmulevelsetvalue_0]$ such that: 
	i) $\Wtransarg{\muxmulevelsetvalue}$ is transversal to 
	$\twoargmumuxtorus{\mulevelsetvalue}{-\muxmulevelsetvalue}$; 
	ii) $\twoargmumuxtorus{\mulevelsetvalue}{-\muxmulevelsetvalue}$ 
	is a torus, specifically, a
	$C^{1,1}$ graph over $\T^2$ in geometric coordinates.
\end{remark}

\begin{remark}[Differential structure with respect to the geometric coordinates] 
\label{R:DIFFERENTIALSTRUCTUREWITHRESPECTTOGEOMETRICCOORDINATES} 
In the rest of the paper, unless we explicitly mention otherwise, 
we implicitly consider $\mulevelsetarg{\mulevelsetvalue}$, $\datahypfortimefunctionarg{-\muxmulevelsetvalue}$,
and $\twoargmumuxtorus{\mulevelsetvalue}{-\muxmulevelsetvalue}$
to be subsets of spacetime with the differential structure 
\emph{induced by the geometric coordinates $(t,u,x^2,x^3)$};
this is already apparent from Def.\,\ref{D:LEVELSETSOFMUANDXMUANDMUMUXTORI}.
Similar remarks apply to the 
sets
$\hypthreearg{\timefunction}{I}{\muxmulevelsetvalue}$,
$\twoargroughtori{\timefunction,u}{\muxmulevelsetvalue}$,
$\nullhypthreearg{\muxmulevelsetvalue}{u}{I}$,
$\twoargMrough{I,J}{\muxmulevelsetvalue}$,
$\mulevelsettwoarg{\mulevelsetvalue}{I}$,
and $\datahypfortimefunctiontwoarg{-\muxmulevelsetvalue}{I}$
defined below. 
\end{remark}

\subsection{Rough time functions, the parameter $\timefunction_0$, and rough adapted coordinates}
\label{SS:ROUGHTIMEFUNCTIONANDROUGHADAPATEDCOORDINATES}
We now provide the transport equation initial value problem whose solution is the rough time function.

\subsubsection{Rough time functions and the parameter $\timefunction_0$}
\label{SSS:ROUGHTIMEFUNCTIONSANDINITIALROUGHSLICE}

\begin{definition}[The rough time function $\timefunctionarg{\muxmulevelsetvalue}$]  
\label{D:ROUGHTIMEFUNCTION}
Let $\mupositive > 0$ be a real number and let $\upmuboot \in [0,\mupositive]$
(see Sect.\,\ref{S:PARAMETERSANDSIZEASSUMPTIONSANDCONVENTIONSFORCONSTANTS} for further discussion of how 
the specific choice of $\mupositive$ that we make in our main theorem is tied to the initial data).
Let $\Wtransarg{\muxmulevelsetvalue} = \muX 
+ 
\phi 
\frac{\muxmulevelsetvalue}{\Lunit \upmu} \Lunit$ be the vectorfield defined in \eqref{E:WTRANSDEF},
and let $\twoargmumuxtorus{\mulevelsetvalue}{-\muxmulevelsetvalue}$ 
be the $\upmu$-adapted torus defined by \eqref{E:MUXMUTORI}.
For $\mulevelsetvalue \in [\upmuboot,\mupositive]$, 
we define the \textbf{rough time function} $\timefunctionarg{\muxmulevelsetvalue}$ 
to be the solution to the following transport equation 
initial value problem:\footnote{See Lemmas~\ref{L:FLOWMAPFORGENERATOROFROUGHTIMEFUNCTION} and 
\ref{L:ODESOLUTIONSTHATARESMOOTHERTHANTHEDATAHYPERSURFACE} for the well-posedness theory of this Cauchy problem.  
\label{FN:IVPFORTRANSPORT}}
\begin{subequations} \label{E:IVPFORROUGHTTIMEFUNCTION}
\begin{align} 	
	\Wtransarg{\muxmulevelsetvalue} \timefunctionarg{\muxmulevelsetvalue} 
	& = 0,  
	\label{E:TRANSPORTEQUATIONFORROUGHTIMEFUNCTION}
	\\
	\timefunctionarg{\muxmulevelsetvalue}|_{\twoargmumuxtorus{\mulevelsetvalue}{-\muxmulevelsetvalue}} 
	& =  
	- \mulevelsetvalue 
	= 
	- \upmu|_{\twoargmumuxtorus{\mulevelsetvalue}{-\muxmulevelsetvalue}}. 
	\label{E:INITIALCONDITIONFORROUGHTIMEFUNCTION}
\end{align} 
\end{subequations}
\end{definition}

We sometimes refer to $\twoargmumuxtorus{\mulevelsetvalue}{-\muxmulevelsetvalue}$ 
as a \emph{primal torus} for $\timefunctionarg{\muxmulevelsetvalue}$ because 
$\timefunctionarg{\muxmulevelsetvalue}$  is ``flowed out from'' it.

\begin{remark}[$\Wtransarg{\muxmulevelsetvalue}$ is tangent to the level sets of $\timefunctionarg{\muxmulevelsetvalue}$]
	\label{R:WTRANSISTANGENTTOLEVELSETSOFROUGHTIMEFUNCTION}
	Note that equation \eqref{E:TRANSPORTEQUATIONFORROUGHTIMEFUNCTION} implies that
	$\Wtransarg{\muxmulevelsetvalue}$ is tangent to the level sets of
	$\timefunctionarg{\muxmulevelsetvalue}$.
\end{remark}

\begin{definition}[The parameter $\timefunction_0$]
\label{D:RELATIONBETWEENINITIALROUGHTIMESLICEANDSMALLMUVALUE}
We define the parameter $\timefunction_0 < 0$ as follows:
\begin{align} \label{E:RELATIONBETWEENINITIALROUGHTIMESLICEANDSMALLMUVALUE}
	\timefunction_0 
	& \eqdef 
	- 
	\mupositive.
\end{align}
\end{definition}

In the bulk of the paper, portions of the rough hypersurface
$\lbrace \timefunctionarg{\muxmulevelsetvalue} = \timefunction_0 \rbrace$
will play the role of an ``initial data'' hypersurface near 
the singularity.
Note that by construction, we have $\timefunction_0 \leq \timefunctionarg{\muxmulevelsetvalue} \leq -\upmuboot \leq 0$.

\subsubsection{Rough adapted coordinates}
\label{SSS:ROUGHADAPTEDSUBSETS}
Having constructed the eikonal function $u$ and the
rough time function $\timefunctionarg{\muxmulevelsetvalue}$, 
we now define a system 
of coordinates adapted to them.

\begin{definition}[The rough adapted coordinates and their partial derivative vectorfields] 
\label{D:ROUGHCOORDINATESANDPARTIALDERIVATIVES}
We call $(\timefunctionarg{\muxmulevelsetvalue},u,x^2,x^3)$ the \emph{rough adapted coordinates}. We 
denote the corresponding rough adapted coordinate partial derivative vectorfields by 
$\left\lbrace\roughgeop{\timefunctionarg{\muxmulevelsetvalue}},\roughgeop{u},\roughgeop{x^2},\roughgeop{x^3} \right\rbrace$. 
\end{definition}

Our analysis will show that the map 
$(t,x^1,x^2,x^3) 
\rightarrow 
\left(\timefunctionarg{\muxmulevelsetvalue},u,x^2,x^3 \right)
$
is a homeomorphism -- all the way up to the singular boundary -- and that it
is a diffeomorphism away from the singular boundary.
Moreover, the map
$(t,u,x^2,x^3) 
\rightarrow 
\left(\timefunctionarg{\muxmulevelsetvalue},u,x^2,x^3 \right)
$
is a diffeomorphism all the way up to the singular boundary;
see Theorem~\ref{T:DEVELOPMENTANDSTRUCTUREOFSINGULARBOUNDARY}.

\begin{remark}[Suppressing the value of $\muxmulevelsetvalue$]	
	\label{R:SUPPRESSINGVALUEOFKAPPASOMETIMES}
	The notation $\left\lbrace\roughgeop{u},\roughgeop{x^2},\roughgeop{x^3} \right\rbrace$
	suppresses the dependence of these operators on $\muxmulevelsetvalue$.
	Moreover, we often write
	$\timefunction$ in place of $\timefunctionarg{\muxmulevelsetvalue}$ or
	$\roughgeop{\timefunction}$ in place of $\roughgeop{\timefunctionarg{\muxmulevelsetvalue}}$
	when there is no danger of confusion about the value of $\muxmulevelsetvalue$.
\end{remark}

\subsection{Rough subsets}
\label{SS:ROUGHSUBSETS}
In this section, we define various subsets of spacetime that are tied to
$u$ and $\timefunctionarg{\muxmulevelsetvalue}$. Most of the delicate PDE analysis in the bulk
of paper will take place on these subsets.
Our analysis will show that 
these subsets are well-adapted to the structure of the singular boundary.

\subsubsection{Truncated $\timefunctionarg{\muxmulevelsetvalue}$-adapted subsets}
\label{SSS:TRUNCATEDROUGHSUBSETS}

\begin{definition}[Truncated $\timefunctionarg{\muxmulevelsetvalue}$-adapted subsets]  \label{D:TRUNCATEDROUGHSUBSETS}
 Given intervals $I, J \subset \mathbb{R}$ and real numbers 
	$\timefunction, u \in \mathbb{R}$, we define:
	\begin{subequations}
	\begin{align} \label{E:ROUGHHYPERSURFACES}
		\hypthreearg{\timefunction}{J}{\muxmulevelsetvalue} 
		& \eqdef \{ (t,u,x^2,x^3) \in \R\times\R\times\T^2 \ | \ \timefunctionarg{\muxmulevelsetvalue}(t,u,x^2,x^3) 
		= 
		\timefunction, \ u \in J \},
			\\
	\twoargroughtori{\timefunction,u}{\muxmulevelsetvalue}
		& 
		\eqdef 
		\{ (t,u,x^2,x^3) \ | \ 
		(t,x^2,x^3) \in \mathbb{R} \times \mathbb{T}^2,
			\,
		\timefunctionarg{\muxmulevelsetvalue}(t,u,x^2,x^3) 
		= 
		\timefunction 
		\},
		 \label{E:ROUGHTORI}
			\\
	\nullhypthreearg{\muxmulevelsetvalue}{u}{I} 
	& 
	\eqdef 
 \bigcup_{\timefunction' \in I} \twoargroughtori{\timefunction',u}{\muxmulevelsetvalue},
	\label{E:NULLHYPERSURFACEROUGHTRUNCATED} 
		\\
\twoargMrough{I,J}{\muxmulevelsetvalue} 
& 
\eqdef \bigcup_{\timefunction' \in I} \hypthreearg{\timefunction'}{J}{\muxmulevelsetvalue} 
= 
\bigcup_{u' \in J} \nullhypthreearg{\muxmulevelsetvalue}{u'}{I}.
\label{E:TRUNCATEDMROUGH}
\end{align}
\end{subequations}
\end{definition}

We refer to the $\hypthreearg{\timefunction}{J}{\muxmulevelsetvalue}$ as
\emph{rough hypersurfaces}. 
We sometimes refer to $\hypthreearg{\timefunction_0}{J}{\muxmulevelsetvalue}$ as the \emph{initial rough hypersurface},
where $\timefunction_0 < 0$ is the parameter from Sect.\,\ref{SSS:ROUGHTIMEFUNCTIONSANDINITIALROUGHSLICE}.
We refer to the $\twoargroughtori{\timefunction,u}{\muxmulevelsetvalue}$ as
\emph{rough tori}.
We also note that $\nullhypthreearg{\muxmulevelsetvalue}{u}{I}$ is a portion of the $\gfour$-null surface $\nullhyparg{u}$.

From Defs.\,\ref{D:ROUGHCOORDINATESANDPARTIALDERIVATIVES} and \ref{D:TRUNCATEDROUGHSUBSETS},
it follows that $\left\lbrace\roughgeop{u},\roughgeop{x^2},\roughgeop{x^3} \right\rbrace$
spans the tangent space of $\hypthreearg{\timefunction}{I}{\muxmulevelsetvalue}$,
that $\left\lbrace\roughgeop{x^2},\roughgeop{x^3} \right\rbrace$
spans the tangent space of $\twoargroughtori{\timefunction,u}{\muxmulevelsetvalue}$,
that 
$\left\lbrace\roughgeop{t},\roughgeop{x^2},\roughgeop{x^3} \right\rbrace$
spans the tangent space of $\nullhypthreearg{\muxmulevelsetvalue}{u}{I}$,
and that $\left\lbrace\roughgeop{\timefunction}, \roughgeop{u},\roughgeop{x^2},\roughgeop{x^3} \right\rbrace$
spans the tangent space of $\twoargMrough{I,J}{\muxmulevelsetvalue}$.

\subsubsection{Truncated $\upmu$-adapted subsets}
\label{SSS:TRUNCATEDMUADAPATEDSUBSETS}
In our analysis, we will often derive estimates on truncated versions of the level sets of various functions,
which we now define.

\begin{definition}[Truncated level sets of $\upmu$ and $\muX \upmu$] 
\label{D:TRUNCATEDLEVELSETSOFMUANDXMUANDMUMUXTORI}
Let $I \subset [\timefunction_0,0]$ be an interval,
let $\mulevelsetvalue \in [0,\upmu_0$], and let $\mulevelsetarg{\mulevelsetvalue}$
and $\datahypfortimefunctionarg{-\muxmulevelsetvalue}$ be the sets from Def.\,\ref{D:LEVELSETSOFMUANDXMUANDMUMUXTORI}.
We define:
\begin{subequations}
\begin{align} \label{E:TRUNCATEDLEVELSETSOFMU}
		\mulevelsettwoarg{\mulevelsetvalue}{I}
		& 
		\eqdef
		\mulevelsetarg{\mulevelsetvalue}
		\cap
		\left\lbrace (t,u,x^2,x^3) \in \R \times \R \times \T^2 
			\ | \ 
			\timefunctionarg{\muxmulevelsetvalue}(t,u,x^2,x^3) \in I 
		\right\rbrace,
			\\
		\datahypfortimefunctiontwoarg{-\muxmulevelsetvalue}{I}
		& \eqdef
		\datahypfortimefunctionarg{-\muxmulevelsetvalue}
		\cap
		\left\lbrace (t,u,x^2,x^3) \in \R \times \R \times \T^2 
			\ | \ 
			\timefunctionarg{\muxmulevelsetvalue}(t,u,x^2,x^3) \in I 
		\right\rbrace.
			\label{E:TRUNCATEDLEVELSETSOFMUXMU} 
\end{align}
\end{subequations}
\end{definition}

Just as in Remark \ref{R:DIFFERENTIALSTRUCTUREWITHRESPECTTOGEOMETRICCOORDINATES}, 
we view 
$\hypthreearg{\timefunction}{J}{\muxmulevelsetvalue}$,
$\nullhypthreearg{\muxmulevelsetvalue}{u}{I}$,
$\twoargMrough{I,J}{\muxmulevelsetvalue}$, 
$\mulevelsettwoarg{\mulevelsetvalue}{I}$,
and
$\datahypfortimefunctiontwoarg{-\muxmulevelsetvalue}{I}$
as submanifolds of spacetime equipped with the differential structure induced by the geometric coordinates $(t,u,x^2,x^3)$.

\section{Coordinate transformations}
\label{S:COORDINATETRANSFORMATIONS}
To prove our main results, we will have to control the transformations between the Cartesian coordinates $(t,x^1,x^2,x^3)$,
the geometric coordinates $(t,u,x^2,x^3)$, 
the rough adapted coordinates $(\timefunctionarg{\muxmulevelsetvalue},u,x^2,x^3)$,
and a few other coordinate systems whose role will become clear later in the paper. 
In this section, we define the relevant change of variables maps and derive some basic
relationships between the partial derivative vectorfields in the different coordinate systems.

\subsection{Change of variables maps}
\label{SS:ALLTHECHOVMAPS}
In this short section, we define various change of variables maps that we use to prove our main results.

\begin{definition}[Change of variables maps]
\label{D:ALLTHECHOVMAPS}
	We define the change of variables map from geometric to Cartesian coordinates as follows:
	\begin{align} \label{E:CHOVGEOTOCARTESIAN}
		\Upsilon(t,u,x^2,x^3) 
		& 
		\eqdef (t,x^1,x^2,x^3).
	\end{align}

	We define the change of variables map from geometric coordinates 
	to rough adapted coordinates as follows:
	\begin{align} \label{E:CHOVGEOTOROUGH}
		\CHOVgeotorough{\muxmulevelsetvalue}(t,u,x^2,x^3) 
		& 
		\eqdef (\timefunctionarg{\muxmulevelsetvalue},u,x^2,x^3).
	\end{align}

	We define the map $\CHOVgeotomumuxmu$ 
	from geometric coordinates to ``$(\upmu,\muX \upmu,x^2,x^3)$-space''
	and its Jacobian 
	$\CHOVJacobiangeotomumuxmu$ as follows:
	\begin{subequations}
	\begin{align} \label{E:CHOVFROMGEOMETRICCOORDINATESTOMUWEGIGHTEDXMUCOORDINATES}
		\CHOVgeotomumuxmu(t,u,x^2,x^3)
		& \eqdef (\upmu,\muX \upmu,x^2,x^3),
			\\
		\CHOVJacobiangeotomumuxmu(t,u,x^2,x^3)
		& \eqdef 
		\frac{\partial \CHOVgeotomumuxmu(t,u,x^2,x^3)}{\partial(t,u,x^2,x^3)}
		=
		\frac{\partial (\upmu,\muX \upmu,x^2,x^3)}{\partial(t,u,x^2,x^3)}.
			\label{E:JACOBIANMATRIXFORCHOVFROMGEOMETRICCOORDINATESTOMUWEGIGHTEDXMUCOORDINATES}
	\end{align}
	\end{subequations}
	
	We define the map $\CHOVroughtomumuxmu{\muxmulevelsetvalue}$ 
	from rough adapted coordinates to ``$(\upmu,\muX \upmu,x^2,x^3)$-space''
	and its Jacobian 
	$\CHOVJacobianroughtomumuxmu{\muxmulevelsetvalue}$ as follows:
	\begin{subequations}
	\begin{align} \label{E:CHOVFROMROUGHCOORDINATESTOMUWEGIGHTEDXMUCOORDINATES}
		\CHOVroughtomumuxmu{\muxmulevelsetvalue}(\timefunctionarg{\muxmulevelsetvalue},u,x^2,x^3)
		& \eqdef (\upmu,\muX \upmu,x^2,x^3),
			\\
		\CHOVJacobianroughtomumuxmu{\muxmulevelsetvalue}(\timefunctionarg{\muxmulevelsetvalue},u,x^2,x^3)
		& \eqdef 
		\frac{\partial \CHOVroughtomumuxmu{\muxmulevelsetvalue}(\timefunction,u,x^2,x^3)}{\partial(\timefunctionarg{\muxmulevelsetvalue},u,x^2,x^3)}
		=
		\frac{\partial (\upmu,\muX \upmu,x^2,x^3)}{\partial(\timefunctionarg{\muxmulevelsetvalue},u,x^2,x^3)}.
			\label{E:JACOBIANMATRIXFORCHOVFROMROUGHCOORDINATESTOMUWEGIGHTEDXMUCOORDINATES}
	\end{align}
	\end{subequations}

\end{definition}
We note the following identity: 
\begin{align} \label{E:SIMPLERELATIONSHIPBETWEENCHOVMAPS}
	\CHOVgeotomumuxmu 
	& = \CHOVroughtomumuxmu{\muxmulevelsetvalue} \circ \CHOVgeotorough{\muxmulevelsetvalue}.
\end{align}

\begin{remark}[Invertibility of the change of variables maps]
\label{R:INVERTIBILITYOFCHOVMAPS}
We justify the invertibility of 
$\Upsilon$ in Prop.\,\ref{P:HOMEOMORPHICANDDIFFEOMORPHICEXTENSIONOFCARTESIANCOORDINATES}
and the invertibility of
$\CHOVgeotorough{\muxmulevelsetvalue}$ 
in Lemma~\ref{L:DIFFEOMORPHICEXTENSIONOFROUGHCOORDINATES}. 
We justify the \emph{local} invertibility of
$\CHOVroughtomumuxmu{\muxmulevelsetvalue}$
in Lemma~\ref{L:CHOVFROMROUGHCOORDINATESTOMUWEGIGHTEDXMUCOORDINATES}.
\end{remark}

\begin{remark}[Implicit functional dependence]
	\label{R:IMPLICITFUNCTIONALDEPENDENCE}
		In most of the paper, our convention is that functions and tensorfields
		should be viewed as depending on the geometric coordinates $(t,u,x^2,x^3)$, 
		\emph{unless we explicitly indicate otherwise}. 
		For example, it is understood that in
	\eqref{E:CHOVGEOTOCARTESIAN}--\eqref{E:JACOBIANMATRIXFORCHOVFROMGEOMETRICCOORDINATESTOMUWEGIGHTEDXMUCOORDINATES}, we are viewing the quantities on the RHSs as functions of $(t,u,x^2,x^3)$. 
	Whenever we make statements relative to the Cartesian coordinates and there 
	is the possibility of confusion, we explicitly indicate the presence of $\Upsilon$;
	see, for example, Prop.\,\ref{P:CHOVGEOMETRICTOCARTESIANISINJECTIVEONREGIONWECAREABOUT}.

	
	We now highlight some occasions when we
	abuse notation that involves functional dependence 
	on the rough adapted coordinates 
	$(\timefunctionarg{\muxmulevelsetvalue},u,x^2,x^3)$, 
	such as on 
	RHSs~\eqref{E:CHOVFROMROUGHCOORDINATESTOMUWEGIGHTEDXMUCOORDINATES}--\eqref{E:JACOBIANMATRIXFORCHOVFROMROUGHCOORDINATESTOMUWEGIGHTEDXMUCOORDINATES}.

	 \begin{itemize}
		 \item Whenever the wave variables $\wavearray$, the acoustic geometry variables,
		etc.\ are shown to depend on the rough adapted coordinates 
			(e.g., when we write $\upmu(\timefunctionarg{\muxmulevelsetvalue},u,x^2,x^3)$),
			it should be understood that we are implicitly composing 
		with $\CHOVgeotorough{\muxmulevelsetvalue}^{-1}$, e.g., 
		by writing $\upmu(\timefunction,u,x^2,x^3)$, we mean
		$\upmu \circ \CHOVgeotorough{\muxmulevelsetvalue}^{-1}(\timefunctionarg{\muxmulevelsetvalue},u,x^2,x^3)$. 
		Put differently, to avoid cluttering the notation, 
		we often avoid explicitly writing the composition with $\CHOVgeotorough{\muxmulevelsetvalue}^{-1}$.
	\item In view of the previous bullet point, in the language of differential geometry, 
	one can view Lemma~\ref{L:ROUGHPARTIALDERIVATIVESINTERMSOFGEOMETRICPARTIALDERIVATIVESANDVICEVERSA} 
	as describing the pushforward of $\{\roughgeop{\timefunction},\roughgeop{u},\roughgeop{x^2},\roughgeop{x^3}\}$ by 
	$\CHOVgeotorough{\muxmulevelsetvalue}^{-1}$ in terms of $\{\geop{t},\geop{u},\geop{x^2},\geop{x^3}\}$. 
	To avoid cluttering the discussion, we do not use the language of pushforwards in this article,
	except in Prop.\,\ref{P:DESCRIPTIONOFSINGULARBOUNDARYINCARTESIANSPACE}, where we use the notion of a pushforward
	to carefully address some degeneracies that occur along the singular boundary.
	Similar remarks apply to Lemma~\ref{L:GEOMETRICVECTORFIELDSINTERMSOFCARTESIANONES}
	(which could be described as identities involving pushfoward by $\Upsilon$).
	\item We highlight that we defined our area and volume forms on 
	the $\timefunctionarg{\muxmulevelsetvalue}$-adapted regions
	$\{\twoargroughtori{\timefunction,u}{\muxmulevelsetvalue}$, 
	$\nullhypthreearg{\muxmulevelsetvalue}{u}{I}$, 
	$\hypthreearg{\timefunction}{I}{\muxmulevelsetvalue}$, 
	and
	$\twoargMrough{I,J}{\muxmulevelsetvalue} \}$ 
	in terms of the rough adapted coordinates (see Def.\,\ref{D:ROUGHVOLFORMS}),
	and that we defined our $L^2$-type norms and energies 
	in terms of the rough adapted coordinates (see, for example \eqref{E:ROUGHTORUSINTEGRAL}--\eqref{E:SPACETIMEINTEGRAL}).
	Moreover, in
	Sects.\,\ref{S:STATEMENTOFALLL2ESTIMATESANDBOOTSTRAPASSUMPTIONSFORWAVEENERGIES}--\ref{S:WAVEANDACOUSTICGEOMETRYAPRIORIESTIMATES}, we derive our $L^2$ estimates in terms of the rough adapted coordinates.
 	\item In Lemmas~\ref{L:POINTWISEESTIMATESFORFULLYMODIFIEDQUANTITIY} and \ref{L:POINTWISEESTIMATEFORPARTIALLYMODIFIEDQUANT}, 
	we derive estimates for modified fluid variables $\fullymodquant{\tander^N}$ and $\partialmodquant{\tander^N}$ along integral curves of $\argLrough{\muxmulevelsetvalue}$ (see \eqref{E:LROUGH}) in rough adapted coordinates. 
	When stating and deriving these estimates, as is explained in the previous bullet points, 
	it should be understood that we are implicitly composing with 
	$\CHOVgeotorough{\muxmulevelsetvalue}^{-1}$. 
	Similar remarks apply during parts of the proof of Prop.\,\ref{P:IMPROVEMENTOFAUXILIARYBOOTSTRAP}.
	\end{itemize}

\end{remark}

\subsection{Coordinate partial derivative transformations}
\label{SS:COORDINATEPARTIALDERIVATIVETRANSFORMATIONS}

\begin{lemma}[Geometric coordinate vectorfields in terms of the Cartesian ones] The following identities hold,
where 
$\left\lbrace 
	\geop{t}, \geop{u}, \geop{x^2}, \geop{x^3} 
\right\rbrace
$ 
are the geometric coordinate partial derivative vectorfields
and $\lbrace \partial_t, \partial_1, \partial_2, \partial_3 \rbrace$ are the Cartesian coordinate partial derivative vectorfields:
\label{L:GEOMETRICVECTORFIELDSINTERMSOFCARTESIANONES}
\begin{align} \label{E:GEOMETRICVECTORFIELDSINTERMSOFCARTESIANONES}
\geop{t} 
& 
= \p_t 
+ 
\left\lbrace \frac{\Lunit^1 X^1 + \Lunit^2 X^2 + \Lunit^3 X^3}{X^1} \right\rbrace \p_1, 
& 
\geop{u} 
& 
= \frac{\upmu \Speed^2}{X^1} \p_1, 
& 
\geop{x^2} & = \p_2 - \frac{X^2}{X^1} \p_1, 
&
\geop{x^3} 
& = \p_3 - \frac{X^3}{X^1} \p_1.
\end{align}
\end{lemma}

\begin{proof}
The identities for $\geop{u}$, $\geop{x^2}$, and $\geop{x^3}$ were proved in \cite[Lemma 2.24]{jLjS2021}. 
To derive the identity \eqref{E:GEOMETRICVECTORFIELDSINTERMSOFCARTESIANONES} 
for $\geop{t}$,
we first equate the following two expressions for $\Lunit$: 
$\Lunit =\p_t +  \Lunit^1 \p_1 + \Lunit^A \p_A = \geop{t} + \Lunit^A \geop{x^A}$. 
We then solve for $\geop{t}$ 
to deduce that 
$\geop{t} = \p_t +  \Lunit^1 \p_1 + \Lunit^A \p_A - \Lunit^A \geop{x^A}$.
Finally, we use the identities for $\left\lbrace \geop{x^2}, \geop{x^3}\right\rbrace$ 
in \eqref{E:GEOMETRICVECTORFIELDSINTERMSOFCARTESIANONES}
to substitute for the factors of $\geop{x^A}$ in the expression $\Lunit^A \geop{x^A}$.
\end{proof}

The following lemma reveals the relationship between the geometric coordinate partial derivative 
vectorfields and the commutation vectorfields from Def.\,\ref{D:COMVECTORFIELDS}. 

\begin{lemma}[Relationship between $\left\lbrace\geop{t},\geop{u},\geop{x^2},\geop{x^3}\right\rbrace$ and $\{\Lunit,X,\Yvf{2},\Yvf{3}\}$] 
\label{L:COMMUTATORSTOCOORDINATES}
The following identities hold:
\begin{subequations} 
\begin{align}
\Lunit 
& = \geop{t} + \Lunit^A \geop{x^A},  
	\label{E:LUNITINTERMSOFGEOMETRICCOORDINATEVECTORFIELDS} 
	\\
\muX 
& = \geop{u} + \upmu X^A \geop{x^A}, 
	\label{E:MUXINTERMSOFGEOMETRICCOORDINATEVECTORFIELDS} 
	\\
\Yvf{2} 
& = \left\lbrace 1 - \Speed^{-2} (X^2)^2 \right\rbrace
		\geop{x^2}
		- 
		\Speed^{-2}X^2 X^3 \geop{x^3}, 
	\label{E:Y2INTERMSOFGEOMETRICCOORDINATEVECTORFIELDS}
		\\
\Yvf{3} 
& =  \left\lbrace
			1 - \Speed^{-2} (X^3)^2
			\right\rbrace
	\geop{x^3}
	- 
	\Speed^{-2}X^2X^3 \geop{x^2}. 
	\label{E:Y3INTERMSOFGEOMETRICCOORDINATEVECTORFIELDS}
\end{align}
\end{subequations}

Moreover, the following identities hold:
\begin{subequations}
\begin{align}
	\geop{t} 
	& = \Lunit - \Lunit^A \Yvf{A} 
	- 
	\left\lbrace
		\frac{\Lunit^2 X^2 + \Lunit^3 X^3}{(X^1)^2} 
	\right\rbrace	
	X^A \Yvf{A}, 
		\label{E:GEOPTOCOMMUTATORS} 
		\\
	\geop{u} 
	& = \muX - \frac{1}{(X^1)^2} \upmu \Speed^2 X^A \Yvf{A}, 
		\label{E:GEOPUTOCOMMUTATORS} 
		\\
	\geop{x^2} 
	& = 
	\left\lbrace
		\frac{(X^1)^2 + (X^2)^2}{(X^1)^2} \Yvf{2} 
	\right\rbrace
	+ 
	\frac{X^2X^3}{(X^1)^2} \Yvf{3}, 
		\label{E:GEOP2TOCOMMUTATORS} 
		\\
	\geop{x^3} 
	& = 
	\left\lbrace
		\frac{(X^1)^2 + (X^3)^2}{(X^1)^2} 
	\right\rbrace
	\Yvf{3} 
	+ 
	\frac{X^2 X^3}{(X^1)^2} \Yvf{2}.  
		\label{E:GEOP3TOCOMMUTATORS}
\end{align}
\end{subequations}
\end{lemma}

\begin{proof}
The identities \eqref{E:LUNITINTERMSOFGEOMETRICCOORDINATEVECTORFIELDS}--\eqref{E:Y3INTERMSOFGEOMETRICCOORDINATEVECTORFIELDS} 
were proved in \cite[Lemma 2.23]{jLjS2021}. The identities \eqref{E:GEOP2TOCOMMUTATORS}--\eqref{E:GEOP3TOCOMMUTATORS} then follow
follow from solving for $\geop{x^2}$ and $\geop{x^3}$ 
in \eqref{E:Y2INTERMSOFGEOMETRICCOORDINATEVECTORFIELDS}--\eqref{E:Y3INTERMSOFGEOMETRICCOORDINATEVECTORFIELDS} 
and using \eqref{E:CARTESIANCOMPONENTSOFXSUMTOCSQUARED}.
To derive \eqref{E:GEOPUTOCOMMUTATORS}, 
we solve for $\geop{u}$ in \eqref{E:MUXINTERMSOFGEOMETRICCOORDINATEVECTORFIELDS}
and use the expressions \eqref{E:Y2INTERMSOFGEOMETRICCOORDINATEVECTORFIELDS}--\eqref{E:Y3INTERMSOFGEOMETRICCOORDINATEVECTORFIELDS} 
as well as \eqref{E:CARTESIANCOMPONENTSOFXSUMTOCSQUARED}. To derive \eqref{E:GEOPTOCOMMUTATORS},
we solve for $\geop{t}$ in \eqref{E:LUNITINTERMSOFGEOMETRICCOORDINATEVECTORFIELDS} 
and use the expressions \eqref{E:Y2INTERMSOFGEOMETRICCOORDINATEVECTORFIELDS}--\eqref{E:Y3INTERMSOFGEOMETRICCOORDINATEVECTORFIELDS}. 
\end{proof}

The following lemma is an analog of Lemma~\ref{L:COMMUTATORSTOCOORDINATES}
with the Cartesian coordinate partial derivative vectorfields in place of the geometric ones.

\begin{lemma}[Relationship between $\{\p_t,\p_1,\p_2,\p_3\}$ and $\{\Lunit,X, \Yvf{2},\Yvf{3}\}$] \label{L:RELATIONSHIPBETWEENCARTESIANPARTIALDERIVATIVESANDSMOOTHGEOMETRICCOMMUTATORS} The following identities hold:
\begin{subequations}
\begin{align}
\p_t 
& = \Lunit - \frac{\Lunit^1X^1 + \Lunit^2X^2 + \Lunit^3X^3}{\Speed^2} X + \frac{\Lunit^1}{X^1} X^A \Yvf{A} - \Lunit^A \Yvf{A}, 
	\label{E:CARTESIANPARTIALTTOCOMMUTATORS} 
		\\
\p_1 
& = \frac{X^1}{\Speed^2} X - \frac{1}{X^1} X^A \Yvf{A}, 
	\label{E:CARTESIANPARTIAL1TOCOMMUTATORS}
	\\
\p_2 
& = \frac{X^2}{\Speed^2} X + \Yvf{2}, 
	\label{E:CARTESIANPARTIAL2TOCOMMUTATORS}
	\\
\p_3 
& = \frac{X^3}{\Speed^2} X + \Yvf{3}. 
	\label{E:CARTESIANPARTIAL3TOCOMMUTATORS}
\end{align}
\end{subequations}
\end{lemma}
\begin{proof}
To derive the identities \eqref{E:CARTESIANPARTIAL2TOCOMMUTATORS}--\eqref{E:CARTESIANPARTIAL3TOCOMMUTATORS},
we use the definition \eqref{E:YCOMMUTATOR} of the vectorfields $\Yvf{A}$
and the expression \eqref{E:ACOUSTICALMETRIC} for the acoustical metric in Cartesian coordinates.
To prove \eqref{E:CARTESIANPARTIAL1TOCOMMUTATORS}, 
we first use \eqref{E:GEOMETRICVECTORFIELDSINTERMSOFCARTESIANONES} 
and \eqref{E:GEOPUTOCOMMUTATORS} 
to deduce
$\p_1 = \frac{X^1}{\upmu \Speed^2} \geop{u} = \frac{X^1}{\Speed^2} X - \frac{X^1}{\Speed^2} X^A \geop{x^A}$.
We then use \eqref{E:GEOP2TOCOMMUTATORS}--\eqref{E:GEOP3TOCOMMUTATORS}
to deduce $\geop{x^A} = \Yvf{A} + \frac{X^A X^B}{(X^1)^2} \Yvf{B}$,
which we substitute into the RHS of the previous expression to obtain \eqref{E:CARTESIANPARTIAL1TOCOMMUTATORS}.
To prove \eqref{E:CARTESIANPARTIALTTOCOMMUTATORS}, we first use \eqref{E:GEOMETRICVECTORFIELDSINTERMSOFCARTESIANONES} and the already proven \eqref{E:CARTESIANPARTIAL1TOCOMMUTATORS} to express $\p_t$ in terms of $\geop{t}$ and the commutation vectorfields. To handle the $\geop{t}$ term, we first use \eqref{E:LUNITINTERMSOFGEOMETRICCOORDINATEVECTORFIELDS} to express
$\geop{t} = \Lunit - \Lunit^A \geop{x^A}$ and then use \eqref{E:GEOP2TOCOMMUTATORS}--\eqref{E:GEOP3TOCOMMUTATORS}
to substitute for the factors of $\geop{x^A}$.
\end{proof}

\begin{corollary}[Expressions for $\smoothtorusproject_{\beta}^{\ \alpha} \partial_{\alpha}$ in terms of $\lbrace \Yvf{2}, \Yvf{3} \rbrace$]
	\label{C:ELLTUPROJECTEDVERSIONOFCARTESIANPARTIALDERIVATIVES}
	The following identities hold relative to the Cartesian coordinates:
	\begin{align} \label{E:ELLTUPROJECTEDVERSIONOFCARTESIANPARTIAL0}
		\smoothtorusproject_0^{\ \alpha} \partial_{\alpha}
		& 
		= 
		\frac{\Lunit^1}{X^1} X^A \Yvf{A} - \Lunit^A \Yvf{A},
			\\
			\smoothtorusproject_1^{\ \alpha} \partial_{\alpha}
		& 
		= 
		- \frac{1}{X^1} X^A \Yvf{A},
			\label{E:ELLTUPROJECTEDVERSIONOFCARTESIANPARTIAL1} 
				\\
		\smoothtorusproject_2^{\ \alpha} \partial_{\alpha}
		& = \Yvf{2},
		&
		\smoothtorusproject_3^{\ \alpha} \partial_{\alpha}
		& = \Yvf{3}.
		\label{E:ELLTUPROJECTEDVERSIONOFCARTESIANPARTIAL2AND3}
	\end{align}
\end{corollary}

\begin{proof}
	To prove \eqref{E:ELLTUPROJECTEDVERSIONOFCARTESIANPARTIAL0},
	we first note that by \eqref{E:PROJECTIONOFTENSORONTOFLATTORUS},
	relative to the Cartesian coordinates, we have
	$(\smoothtorusproject \partial_t)^{\alpha} 
	= \smoothtorusproject_{\beta}^{\ \alpha} (\partial_t)^{\beta} = \smoothtorusproject_0^{\ \alpha}$.
	Hence, LHS~\eqref{E:ELLTUPROJECTEDVERSIONOFCARTESIANPARTIAL0} 
	is the projection of $\partial_t$ onto $\ell_{t,u}$.
	The identity \eqref{E:ELLTUPROJECTEDVERSIONOFCARTESIANPARTIAL0} therefore follows from \eqref{E:CARTESIANPARTIALTTOCOMMUTATORS}.
	The identities \eqref{E:ELLTUPROJECTEDVERSIONOFCARTESIANPARTIAL1}--\eqref{E:ELLTUPROJECTEDVERSIONOFCARTESIANPARTIAL2AND3}
	follow from similar arguments based on \eqref{E:CARTESIANPARTIAL1TOCOMMUTATORS}--\eqref{E:CARTESIANPARTIAL3TOCOMMUTATORS}.
\end{proof}

The next lemma reveals the relationship between 
the rough adapted coordinate partial derivative vectorfields
the geometric coordinate partial derivative vectorfields. 

\begin{lemma}[Relationship between $\{\geop{t},\geop{u},\geop{x^2},\geop{x^3}\}$ and $\{\roughgeop{\timefunction},\roughgeop{u},\roughgeop{x^2},\roughgeop{x^3}\}$] \label{L:ROUGHPARTIALDERIVATIVESINTERMSOFGEOMETRICPARTIALDERIVATIVESANDVICEVERSA} 
The following identities hold for $A = 2,3$, where $\timefunction$ denotes the rough time function
(with the superscript $\muxmulevelsetvalue$ suppressed):
\begin{subequations}
\begin{align}
\roughgeop{\timefunction} 
	& = \frac{1}{\geop{t}\timefunction} \geop{t}, 
	\label{E:ROUGHTIMEPARTIALDERIVATIVEINTERMSOFGEOMETRICTIMEPARTIALDERIVATIVE}
	\\
\roughgeop{u} 
& = 
\geop{u} 
- 
\frac{\geop{u}\timefunction}{\geop{t} \timefunctionarg{\muxmulevelsetvalue}} \geop{t}, 
\label{E:ROUGHGEOPUINTERMOFGEOMETRICCOORDINATEVECTORFIELDS}
	\\ 
\roughgeop{x^A} 
& 
= \geop{x^A} 
- 
\frac{\geop{x^A} \timefunctionarg{\muxmulevelsetvalue}}{\geop{t} \timefunctionarg{\muxmulevelsetvalue}} \geop{t} 
= \geop{x^A} 
+ 
\frac{\geop{x^A} \timefunction}{\geop{t} \timefunctionarg{\muxmulevelsetvalue}} \Lunit^B \geop{x^B} 
- 
\frac{\geop{x^A} \timefunction}{\geop{t} \timefunctionarg{\muxmulevelsetvalue}} \Lunit. 
\label{E:ROUGHANGULARPARTIALDERIVATIVESINTERMSOFGOODGEOMETRICPARTIALDERIVATIVES} 
\end{align}
\end{subequations}

Moreover, the following identity holds:
\begin{align} \label{E:SMOOTHANGULARVECTORFIELDSINTERMSOFROUGHONES}
\geop{x^A} 
& = 
\roughgeop{x^A} 
- 
\frac{\geop{x^A} \timefunction}{\Lunit \timefunctionarg{\muxmulevelsetvalue}} \Lunit^B \roughgeop{x^B} 
+ 
\frac{\geop{x^A} \timefunctionarg{\muxmulevelsetvalue}}{\Lunit \timefunctionarg{\muxmulevelsetvalue}} \Lunit.
\end{align}
\end{lemma}


\begin{proof}
The identity \eqref{E:ROUGHTIMEPARTIALDERIVATIVEINTERMSOFGEOMETRICTIMEPARTIALDERIVATIVE} follows from the chain rule identity 
$\geop{t} 
= 
(\geop{t} \timefunctionarg{\muxmulevelsetvalue}) \roughgeop{\timefunction} 
+ 
(\geop{t} u) \roughgeop{u} 
+ 
(\geop{t} x^A) \roughgeop{x^A}$ 
and the fact that $\geop{t} u = \geop{t} x^A = 0$. 
The identities 
\eqref{E:ROUGHGEOPUINTERMOFGEOMETRICCOORDINATEVECTORFIELDS}--\eqref{E:ROUGHANGULARPARTIALDERIVATIVESINTERMSOFGOODGEOMETRICPARTIALDERIVATIVES}follow from similar arguments and \eqref{E:LUNITINTERMSOFGEOMETRICCOORDINATEVECTORFIELDS}.

Finally, with the help of \eqref{E:LULTMUXUMUXT},
it is straightforward to confirm the identity \eqref{E:SMOOTHANGULARVECTORFIELDSINTERMSOFROUGHONES} 
by checking that both sides evaluate to the same values when acting on the rough adapted coordinate functions 
$\timefunction,u,x^2,x^3$. 
\end{proof}

\subsection{An identity for $\roughgeop{u}$}
\label{SS:ROUGHUDERIVATIVEINTERMSOFGEOMETRICVECTORFIELDSANDROUGHTIMEFUNCTION}
We will use the following simple identity in Sect.\,\ref{SS:PROPERTIESOFCHOVFROMGEOMETRICTOCARTESIAN},
when we study the homeomorphism and diffeomorphism properties of the change of variables map $\Upsilon$.

\begin{lemma}[An identity for $\roughgeop{u}$]
		\label{L:ROUGHUDERIVATIVEINTERMSOFGEOMETRICVECTORFIELDSANDROUGHTIMEFUNCTION}
	The following identity holds: 
	\begin{align} \label{E:ROUGHUDERIVATIVEINTERMSOFGEOMETRICVECTORFIELDSANDROUGHTIMEFUNCTION}
	\roughgeop{u} 
		& = \muX 
				+ 
				\phi
				\frac{\muxmulevelsetvalue}{\Lunit \upmu} 
				\Lunit
				-
				\left\lbrace
					\upmu X^A
					+
					\phi
					\frac{\muxmulevelsetvalue}{\Lunit \upmu} 
					\Lunit^A
				\right\rbrace
				\left\lbrace
					\geop{x^A} 
					+ 
					\frac{\geop{x^A} \timefunctionarg{\muxmulevelsetvalue}}{\geop{t} \timefunctionarg{\muxmulevelsetvalue}} \Lunit^B \geop{x^B} 
					- 
					\frac{\geop{x^A} \timefunctionarg{\muxmulevelsetvalue}}{\geop{t} \timefunctionarg{\muxmulevelsetvalue}} \Lunit
				\right\rbrace.
		\end{align}
\end{lemma}

\begin{proof}
	First, we use \eqref{E:WTRANSDEF}, \eqref{E:TRANSPORTEQUATIONFORROUGHTIMEFUNCTION}, 
	and Lemma~\ref{L:BASICPROPERTIESOFVECTORFIELDS} to deduce that
	$\Wtransarg{\muxmulevelsetvalue} u = 1$ and $\Wtransarg{\muxmulevelsetvalue} 
	\timefunctionarg{\muxmulevelsetvalue} = 0$.
	It follows that
	$\roughgeop{u} 
		= 
		\Wtransarg{\muxmulevelsetvalue}
		-
		\Wtranstwoarg{\muxmulevelsetvalue}{A} \roughgeop{x^A}
	$,
	where as usual, $\Wtranstwoarg{\muxmulevelsetvalue}{A} = \Wtransarg{\muxmulevelsetvalue} x^A$.
	From this identity, 
	\eqref{E:WTRANSDEF},
	\eqref{E:LUNITINTERMSOFGEOMETRICCOORDINATEVECTORFIELDS},
	and
	\eqref{E:ROUGHANGULARPARTIALDERIVATIVESINTERMSOFGOODGEOMETRICPARTIALDERIVATIVES},
	we conclude \eqref{E:ROUGHUDERIVATIVEINTERMSOFGEOMETRICVECTORFIELDSANDROUGHTIMEFUNCTION}.
\end{proof}

\section{The rough acoustical geometry and curvature tensors}
\label{S:ROUGHACOUSTICALGEOMETRYANDCURVATURETENSORS}
In this section, we set up the acoustical geometry of the rough foliations,
that is, the geometry associated to the rough hypersurfaces
$\twoargroughtori{\timefunction,u}{\muxmulevelsetvalue}$,
characteristics $\nullhyparg{u}$, and their
intersection $\twoargroughtori{\timefunction,u}{\muxmulevelsetvalue}$.
In particular, we define 
the first fundamental forms induced by the acoustical metric $\gfour$ on
the submanifolds $\twoargroughtori{\timefunction,u}{\muxmulevelsetvalue}$ and 
$\hypthreearg{\timefunction}{[- \rightu,\leftu]}{\muxmulevelsetvalue}$,
we define various geometric vectorfields tied to these submanifolds,
we exhibit various geometric decompositions,
and we introduce the Riemann and Ricci curvature of the acoustical metric $\gfour$ 
and various curvature tensors of the first fundamental form of
$\twoargroughtori{\timefunction,u}{\muxmulevelsetvalue}$.

\begin{remark}[Sometimes suppressing dependence on $\muxmulevelsetvalue$] 
\label{R:KAPPASUPPRESSION}
In view of \eqref{E:WTRANSDEF} and 
Def.\,\ref{D:ROUGHTIMEFUNCTION} for $\timefunctionarg{\muxmulevelsetvalue}$,
it follows that all our constructions in this section depend on the choice of $\muxmulevelsetvalue$. 
On the other hand, the constants ``$C$'' in our forthcoming estimates can be chosen to   
independent of $\muxmulevelsetvalue$ for $\muxmulevelsetvalue \in [0,\muxmulevelsetvalue_0]$. 
This is important because in Sects.\,\ref{S:DEVELOPMENTSOFTHEDATASINGULARBOUNDARYNEWTIMEFUNCTION}--\ref{S:MAINRESULTS}, 
we will vary $\muxmulevelsetvalue \in [0,\muxmulevelsetvalue_0]$ to obtain a continuum of time functions 
$\timefunctionarg{\muxmulevelsetvalue}$ and their respective foliations, and our results depend on the fact that the
``$C$'' can be chosen uniformly with respect to $\muxmulevelsetvalue$. 
However, to simplify the notation, 
there are many geometric objects for which we often suppress their dependence on 
$\muxmulevelsetvalue$ 
(for example, the tensorfield $\gtorusroughfirstfund$ introduced below depends on $\muxmulevelsetvalue$). 
Nonetheless, in some objects, such as $\hypthreearg{\timefunction}{[- \rightu,\leftu]}{\muxmulevelsetvalue}$,
we will retrain the explicit $\muxmulevelsetvalue$-dependence
to provide the reader mental reminders regarding our constructions. 
\end{remark}

\subsection{First fundamental forms of $\twoargroughtori{\timefunction,u}{\muxmulevelsetvalue}$ and 
$\hypthreearg{\timefunction}{[- \rightu,\leftu]}{\muxmulevelsetvalue}$} 
\label{SS:ROUGHFIRSTFUNDS}
We now introduce the first fundamental forms of the rough tori $\twoargroughtori{\timefunction,u}{\muxmulevelsetvalue}$ 
and the rough hypersurfaces $\hypthreearg{\timefunction}{[- \rightu,\leftu]}{\muxmulevelsetvalue}$
induced by the acoustical metric $\gfour$. 

\begin{definition}[The first fundamental forms of $\twoargroughtori{\timefunction,u}{\muxmulevelsetvalue}$ and 
$\hypthreearg{\timefunction}{[- \rightu,\leftu]}{\muxmulevelsetvalue}$] 
\label{D:ROUGHFIRSTFUNDS} \hfill

\noindent \underline{\textbf{First fundamental form of $\twoargroughtori{\timefunction,u}{\muxmulevelsetvalue}$}}:
\begin{itemize}
\item
We define the \emph{first fundamental form of} $\twoargroughtori{\timefunction,u}{\muxmulevelsetvalue}$ relative to $\gfour$ 
to be the symmetric type $\binom{0}{2}$ tensorfield 
$\gtorusroughfirstfund$ 
such that $\gtorusroughfirstfund(Y,Z) = \gfour(Y,Z)$
for all pairs $(Y,Z)$ of vectorfields tangent to $\twoargroughtori{\timefunction,u}{\muxmulevelsetvalue}$ 
and such that $\gtorusroughfirstfund(V,\cdot) = \gtorusroughfirstfund(\cdot,V) = 0$
if $V$ is $\gfour$-orthogonal to $\twoargroughtori{\timefunction,u}{\muxmulevelsetvalue}$.
\item
We define the \emph{inverse first fundamental form} $\gtorusroughinversefirstfund$
to be the dual of $\gtorusroughfirstfund$ relative to $\gfour$, 
i.e., relative to arbitrary coordinates,  
it is the symmetric type $\binom{2}{0}$ tensorfield with the following components:
\begin{align} \label{E:GTORUSROUGHINVERSEDEFINITION}
		(\gtorusroughinversefirstfund)^{\alpha \beta}
		& \eqdef
		(\gfour^{-1})^{\alpha \gamma}
		(\gfour^{-1})^{\beta \delta}
		\gtorusroughfirstfund_{\gamma \delta}.
\end{align}
\end{itemize}

\noindent \underline{\textbf{First fundamental form of $\hypthreearg{\timefunction}{[- \rightu,\leftu]}{\muxmulevelsetvalue}$}}:
\begin{itemize}
\item
We define the \emph{first fundamental form of} $\hypthreearg{\timefunction}{[- \rightu,\leftu]}{\muxmulevelsetvalue}$ relative to $\gfour$ 
to be the symmetric type $\binom{0}{2}$ tensorfield 
$\hypg$ 
such that $\hypg(Y,Z) = \gfour(Y,Z)$
for all pairs $(Y,Z)$ of vectorfields tangent to $\hypthreearg{\timefunction}{[- \rightu,\leftu]}{\muxmulevelsetvalue}$ 
and such that $\hypg(V,\cdot) = \hypg(\cdot,V) = 0$
if $V$ is $\gfour$-orthogonal to $\hypthreearg{\timefunction}{[- \rightu,\leftu]}{\muxmulevelsetvalue}$.
\item
We define the \emph{inverse first fundamental form} $\hypginverse$
to be the dual of $\hypg$ relative to $\gfour$, 
i.e., relative to arbitrary coordinates,  
it is the symmetric type $\binom{2}{0}$ tensorfield with the following components:
\begin{align} \label{E:ROUGHHYPERSURFACEINVERSEFIRSTFUNDDEFINITION}
		(\hypginverse)^{\alpha \beta}
		& \eqdef
		(\gfour^{-1})^{\alpha \gamma}
		(\gfour^{-1})^{\beta \delta}
		\hypg_{\gamma \delta}.
\end{align}
\end{itemize}

\end{definition}

\subsection{Geometric vectorfields associated to the rough foliations}
\label{SS:GEOMETRICVECTORFIELDSASSOCIATEDTOROUGHFOLIATIONS}
In this section, we define several vectorfields that play a crucial role in our analysis of the rough geometry,
and we reveal their basic properties.

We start by introducing the rough null vectorfield $\argLrough{\muxmulevelsetvalue}$,
which is the unique $\nullhyparg{u}$-tangent vectorfield normalized by $\argLrough{\muxmulevelsetvalue} \timefunctionarg{\muxmulevelsetvalue} = 1$,
i.e., $\argLrough{\muxmulevelsetvalue}$ is normalized relative to the rough time function. 
Roughly speaking, it plays a similar role that $\Lunit$ played in the works
\cites{dC2007,jSgHjLwW2016,jS2016b,jLjS2018},
which relied on the Cartesian time function.
However, $\argLrough{\muxmulevelsetvalue}$ enjoys less regularity than $\Lunit$,
so to obtain top-order estimates, we must commute the equations with $\Lunit$
(rather than $\argLrough{\muxmulevelsetvalue}$).

\begin{definition}[The rough null vectorfield] \label{D:LROUGH} 
 We define the \emph{rough null vectorfield} vectorfield $\argLrough{\muxmulevelsetvalue}$ as follows: 
\begin{align} \label{E:LROUGH}
\argLrough{\muxmulevelsetvalue} 
& \eqdef 
\frac{1}{\Lunit \timefunctionarg{\muxmulevelsetvalue}} \Lunit.
\end{align}
\end{definition}

We will often use the following basic property of $\argLrough{\muxmulevelsetvalue}$,
which follows immediately from \eqref{E:LROUGH}:
\begin{align} \label{E:LROUGHAPPLIEDTOROUGHTIMEFUNCTIONISUNITY}
\argLrough{\muxmulevelsetvalue} \timefunctionarg{\muxmulevelsetvalue} 
& = 1. 
\end{align}

\begin{definition}[The vectorfields $\Roughtoritangentvectorfieldarg{\muxmulevelsetvalue}$, $\Rtransarg{\muxmulevelsetvalue}$, $\Rtransunitarg{\muxmulevelsetvalue}$, $\hypnormalarg{\muxmulevelsetvalue}$, and $\hypunitnormalarg{\muxmulevelsetvalue}$] 
\label{D:GEOMETRICVECTORFIELDSADAPTEDTOROUGHFOLIATIONS} \hfill 
\begin{itemize}
\item We define $\Roughtoritangentvectorfieldarg{\muxmulevelsetvalue}$ to be the following $\twoargroughtori{\timefunction,u}{\muxmulevelsetvalue}$-tangent vectorfield:
\begin{align} \label{E:DEFROUGHTORITANGENTVECTORFIELD}
\Roughtoritangentvectorfieldarg{\muxmulevelsetvalue} 
& \eqdef 
	\gtorusroughinversefirstfund\left(dx^A,dx^B \right) 
	\frac{\geop{x^A} \timefunctionarg{\muxmulevelsetvalue}}{\geop{t} \timefunctionarg{\muxmulevelsetvalue}} \roughgeop{x^B}.
\end{align}
\item We define $\Rtransarg{\muxmulevelsetvalue}$ to be the following $\hypthreearg{\timefunction}{[- \rightu,\leftu]}{\muxmulevelsetvalue}$-tangent vectorfield:
\begin{align} \label{E:RTRANS} 
\Rtransarg{\muxmulevelsetvalue} 
& 
\eqdef 
\Wtransarg{\muxmulevelsetvalue} 
- 
\upmu \Roughtoritangentvectorfieldarg{\muxmulevelsetvalue} 
= \muX + \phi \frac{ \muxmulevelsetvalue}{\Lunit \upmu} \Lunit 
-  
\upmu \Roughtoritangentvectorfieldarg{\muxmulevelsetvalue},
\end{align}
where $\phi = \phi(u)$ is the cut-off function introduced in Def.\,\ref{D:WTRANSANDCUTOFF}. 

For the solutions featured in our main results,
$\Rtransarg{\muxmulevelsetvalue}$ will be $\gfour$-spacelike, i.e.,  
$|\Rtransarg{\muxmulevelsetvalue}|_{\gfour}^2 \eqdef \gfour(\Rtransarg{\muxmulevelsetvalue},\Rtransarg{\muxmulevelsetvalue}) > 0$. 
Moreover, we define
$\Rtransunitarg{\muxmulevelsetvalue}$ to be the $\gfour$-unit-length rescaling of $\Rtransarg{\muxmulevelsetvalue}$:
 \begin{align} \label{E:UNITLENGTHRTRANS} 
\Rtransunitarg{\muxmulevelsetvalue} 
&
\eqdef \frac{1}{|\Rtransarg{\muxmulevelsetvalue}|_{\gfour}} \Rtransarg{\muxmulevelsetvalue}.
\end{align}
\item We define $\hypnormalarg{\muxmulevelsetvalue}$ to be the following vectorfield: 
\begin{align} \label{E:HYPNORMALDEF}
	\hypnormalarg{\muxmulevelsetvalue} 
	& 
	\eqdef \Lunit 
	+ 
	\frac{\upmu}{|\Rtransarg{\muxmulevelsetvalue}|_{\gfour}^2} \Rtransarg{\muxmulevelsetvalue}.
\end{align}

For the solutions featured in our main results,
$\hypnormalarg{\muxmulevelsetvalue}$ will be $\gfour$-timelike, i.e., $\gfour(\hypnormalarg{\muxmulevelsetvalue},\hypnormalarg{\muxmulevelsetvalue}) < 0$.
Moreover, we define
$\hypunitnormalarg{\muxmulevelsetvalue}$ to be the $\gfour$-unit-length rescaling of $\hypnormalarg{\muxmulevelsetvalue}$:
\begin{align} \label{E:UNITHYPNORMALDEF}
\hypunitnormalarg{\muxmulevelsetvalue} 
& 
\eqdef 
\frac{1}{\sqrt{-\gfour(\hypnormalarg{\muxmulevelsetvalue},\hypnormalarg{\muxmulevelsetvalue})}} \hypnormalarg{\muxmulevelsetvalue}.
\end{align}
\end{itemize}
\end{definition}

\subsection{Identities involving the first fundamental forms and geometric vectorfields}
\label{SS:IDENTITIESINVOLVINGTHEFIRSTFUNDAMENTALFORMSANDGEOMETRICVECTORFIELDS}

\begin{lemma}[Identities for $\gtorusroughfirstfund$ and $\gtorusroughinversefirstfund$]
\label{L:IDENTITIESFORCOMPONENTSOFROUGHTORUSFIRSTFUNDAMENTALFORMANDITSINVERSE}
When restricted to the tangent space of $\twoargroughtori{\timefunction,u}{\muxmulevelsetvalue}$,
we have the following identity for $\gtorusroughfirstfund$ 
relative to the rough adapted coordinates:
\begin{align} \label{E:GTORUSFIRSTFUNDRESTRICTEDTOROUGHTORUS}
	\gtorusroughfirstfund
	& = 
	\gtorusroughfirstfund\left(\roughgeop{x^A},\roughgeop{x^B}\right)
	d x^A \otimes d x^B,
\end{align}
where with $\gtorus_{AB}$ as in \eqref{E:SMOOTHTORIGABEXPRESSION}, 
we have:
\begin{align} \label{E:GTORUSROUGHCOMPONENTS}
\gtorusroughfirstfund\left(\roughgeop{x^A},\roughgeop{x^B}\right)
& 
=  
\gtorus_{AB} 
+  
\frac{\geop{x^A} \timefunctionarg{\muxmulevelsetvalue}}{\geop{t} \timefunctionarg{\muxmulevelsetvalue}} \gtorus_{BC} \Lunit^C  
+   
\frac{\geop{x^B}\timefunction}{\geop{t} \timefunctionarg{\muxmulevelsetvalue}} \gtorus_{AC} \Lunit^C 
+ 
\frac{(\geop{x^A} \timefunction) \geop{x^B} \timefunction}{(\geop{t} \timefunctionarg{\muxmulevelsetvalue})^2} \gtorus_{CD} \Lunit^C \Lunit^D.
\end{align}

Moreover, relative to the rough adapted coordinates, the following identity holds:
\begin{align} \label{E:GTORUSINVERSEROUGHEXPRESSION}
\gtorusroughinversefirstfund
& = 
\gtorusroughinversefirstfund\left(dx^A,dx^B \right)
\roughgeop{x^A} \otimes \roughgeop{x^B},
\end{align}
where:
\begin{align} \label{E:ROUGHTORUSMETRICCOMPONENTSANDTHEINVERSECOMPONENTSRELATION}
	\gtorusroughinversefirstfund\left(dx^A,dx^C\right)
	\gtorusroughfirstfund\left(\roughgeop{x^C},\roughgeop{x^B}\right)
	&
	= \updelta_B^A,
\end{align}
and $\updelta_B^A$ is the Kronecker delta.

\end{lemma}	

\begin{proof}
\eqref{E:GTORUSFIRSTFUNDRESTRICTEDTOROUGHTORUS} is a simple consequence of the fact that
$(x^2,x^3)$ are coordinates on $\twoargroughtori{\timefunction,u}{\muxmulevelsetvalue}$.
The identity \eqref{E:GTORUSROUGHCOMPONENTS} follows from the fact that
$\gtorusroughfirstfund\left(\roughgeop{x^A},\roughgeop{x^B}\right)
=
\gfour\left(\roughgeop{x^A},\roughgeop{x^B}\right)$,
\eqref{E:ROUGHANGULARPARTIALDERIVATIVESINTERMSOFGOODGEOMETRICPARTIALDERIVATIVES},
and the fact that $\gfour(\Lunit,\Lunit) = \gfour(\Lunit,\roughgeop{x^A}) = 0$.
The identity \eqref{E:GTORUSINVERSEROUGHEXPRESSION} is a simple consequence of the fact that
the rough coordinate vectorfields 
$\left\lbrace \roughgeop{x^A} \right\rbrace_{A=2,3}$ 
span the tangent space of
$\twoargroughtori{\timefunction,u}{\muxmulevelsetvalue}$.
Next, we note that it is straightforward to check, using \eqref{E:GTORUSROUGHINVERSEDEFINITION}, 
that the type $\binom{1}{1}$ tensorfield with components
$(\gtorusroughinversefirstfund)^{\alpha \gamma} \gtorusroughfirstfund_{\gamma \beta}$
is the $\gfour$-orthogonal projection tensorfield onto $\twoargroughtori{\timefunction,u}{\muxmulevelsetvalue}$.
From this fact and the fact that
$(x^2,x^3)$ are coordinates on $\twoargroughtori{\timefunction,u}{\muxmulevelsetvalue}$,
the identity \eqref{E:ROUGHTORUSMETRICCOMPONENTSANDTHEINVERSECOMPONENTSRELATION} 
readily follows.
\end{proof}

\begin{lemma}[Relationship between components $\gtorus$ and $\gtorusroughfirstfund$] 
\label{L:RELATIONSHIPBETWEENSMOOTHTORIFIRSTFUNDANDROUGHTORIFIRSTFUND}
We define $\gtorusCOV$ to be the $2 \times 2$ matrix with the following entries:
\begin{align} \label{E:CHOVCOEFFICIENTSSMOOTHANGULARDERIVATIVESINTERMSOFROUGHONESANDL}
	\gtorusCOV_A^B 
	& \eqdef 
	\updelta_A^B
	- 
	\frac{\geop{x^A} \timefunction}{\Lunit \timefunctionarg{\muxmulevelsetvalue}} \Lunit^B,
\end{align}	
where $\updelta_A^B$ is the Kronecker delta.
Then the following identity holds:
\begin{align} \label{E:SMOOTHANGULARDERIVATIVESINTERMSOFROUGHONESANDL}
\geop{x^A} 
& 
= 
\gtorusCOV_A^B \roughgeop{x^B} 
+ 
\frac{\geop{x^A} \timefunctionarg{\muxmulevelsetvalue}}{\Lunit \timefunctionarg{\muxmulevelsetvalue}} \Lunit. 
\end{align}

Moreover, recalling that $\gtorus$ is the first fundamental form of $\ell_{t,u}$,
we have the following relationship between 
$\gtorus_{AB} \eqdef \gtorus(\geop{x^A},\geop{x^B})$
and 
$
\gtorusroughfirstfund\left(\roughgeop{x^A},\roughgeop{x^B} \right)
$:
\begin{align} \label{E:SMOOTHTORUSFIRSTFUNDCOMPONENTSINTERMSOFROUGHTORUSFIRSTFUNDCOMPONENTS}
	\gtorus_{AB}
	& 
	= \gtorusCOV_A^C \gtorusCOV_B^D \gtorusroughfirstfund\left(\roughgeop{x^C},\roughgeop{x^D} \right),
		\\
	(\gtorus^{-1})^{AB} 
	& 
	= (\gtorusCOV^{-1})_C^A (\gtorusCOV^{-1})_D^B \gtorusroughinversefirstfund\left(dx^C,dx^D \right).
		\label{E:SMOOTHTORUSINVERSEFIRSTFUNDCOMPONENTSINTERMSOFROUGHTORUSINVERSEFIRSTFUNDCOMPONENTS}
\end{align}

In addition, the inverse $(\gtorusCOV^{-1})_A^B$ of $\gtorusCOV_B^A$, defined by
$(\gtorusCOV^{-1})_B^C \gtorusCOV_C^A = \updelta_B^A$, can be expressed as follows:
\begin{align} \label{E:INVERSECHOVCOEFFICIENTSSMOOTHANGULARDERIVATIVESINTERMSOFROUGHONESANDL}
	(\gtorusCOV^{-1})_A^B 
	& 
	= 
	\updelta_A^B 
	+ 
	\frac{\geop{x^A} \timefunctionarg{\muxmulevelsetvalue}}{\geop{t} \timefunctionarg{\muxmulevelsetvalue}} \Lunit^B.
\end{align}	

Finally, we have the following identity:
\begin{align} \label{E:ROUGHTORUSCOMPONENTSINTERMSOFSMOOTHTORUSCOMPONENTS}
	\gtorusroughfirstfund\left(\roughgeop{x^A},\roughgeop{x^B} \right)
	& 
	= (\gtorusCOV^{-1})_A^{C} (\gtorusCOV^{-1})_B^D \gtorus_{CD}.
\end{align}
\end{lemma}

\begin{proof}
The identity \eqref{E:SMOOTHANGULARDERIVATIVESINTERMSOFROUGHONESANDL} 
is a restatement of
\eqref{E:SMOOTHANGULARVECTORFIELDSINTERMSOFROUGHONES}. 
To derive \eqref{E:INVERSECHOVCOEFFICIENTSSMOOTHANGULARDERIVATIVESINTERMSOFROUGHONESANDL},
we note that
\eqref{E:SMOOTHANGULARDERIVATIVESINTERMSOFROUGHONESANDL}
implies that
$
\roughgeop{x^A} 
= 
(\gtorusCOV^{-1})_A^B
\geop{x^B} 
+
f \Lunit
$
for some scalar function $f$.
We then note that the coefficients $(\gtorusCOV^{-1})_A^B$ are given by
\eqref{E:ROUGHANGULARPARTIALDERIVATIVESINTERMSOFGOODGEOMETRICPARTIALDERIVATIVES}.
Using \eqref{E:INVERSECHOVCOEFFICIENTSSMOOTHANGULARDERIVATIVESINTERMSOFROUGHONESANDL},
we see that
\eqref{E:ROUGHTORUSCOMPONENTSINTERMSOFSMOOTHTORUSCOMPONENTS} follows from \eqref{E:GTORUSROUGHCOMPONENTS}.
\eqref{E:SMOOTHTORUSFIRSTFUNDCOMPONENTSINTERMSOFROUGHTORUSFIRSTFUNDCOMPONENTS}
follows from applying two factors of $\gtorusCOV$ to each side
of \eqref{E:ROUGHTORUSCOMPONENTSINTERMSOFSMOOTHTORUSCOMPONENTS}.
\eqref{E:SMOOTHTORUSINVERSEFIRSTFUNDCOMPONENTSINTERMSOFROUGHTORUSINVERSEFIRSTFUNDCOMPONENTS}
follows from taking the inverse of \eqref{E:SMOOTHTORUSFIRSTFUNDCOMPONENTSINTERMSOFROUGHTORUSFIRSTFUNDCOMPONENTS}
and using \eqref{E:ROUGHTORUSMETRICCOMPONENTSANDTHEINVERSECOMPONENTSRELATION}.
\end{proof}

\begin{proposition}[Properties of the geometric vectorfields associated to the rough foliations] 
\label{P:BASICPROPERTIESOFROUGHVECTORFIELDS}
The vectorfield $\Rtransarg{\muxmulevelsetvalue}$ defined in \eqref{E:RTRANS}
is $\gfour$-orthogonal to $\twoargroughtori{\timefunction,u}{\muxmulevelsetvalue}$, 
i.e., 
$\gfour$-orthogonal to the elements of $\{\geop{x^C}\}_{C = 2,3}$. 
Moreover, its square norm as measured by $\gfour$ 
(which is the same as its square norm as measured by $\hypg$)
satisfies:
\begin{subequations}
\begin{align} \label{E:SIZEOFRTRANS}
|\Rtransarg{\muxmulevelsetvalue}|_{\gfour}^2  
&
=
|\Rtransarg{\muxmulevelsetvalue}|_{\hypg}^2 
= 
\upmu^2 (1 - \Rtransnormsmallfactorarg{\muxmulevelsetvalue}) 
- 
2 \muxmulevelsetvalue \phi \frac{\upmu}{\Lunit \upmu},
\end{align}
where $\phi = \phi(u)$ is the cut-off function from Def.\,\ref{D:WTRANSANDCUTOFF}
and $\Rtransnormsmallfactorarg{\muxmulevelsetvalue} \geq 0$ is defined by:
\begin{align} \label{E:RTRANSNORMSMALLFACTOR} 
\Rtransnormsmallfactorarg{\muxmulevelsetvalue} 
& 
\eqdef 
\frac{\gtorusroughinversefirstfund\left(dx^A,dx^B \right) (\geop{x^A} \timefunction) \geop{x^B} \timefunction}{(\geop{t}\timefunction)^2}.
\end{align}
In particular, if $\upmu > 0$, $\Rtransnormsmallfactorarg{\muxmulevelsetvalue} < 1$, and $\Lunit \upmu < 0$ on the support of $\phi$
(all of which are satisfied for the solutions featured in our main results),
then $\Rtransarg{\muxmulevelsetvalue}$ is $\gfour$-spacelike, and the vectorfield $\Rtransunitarg{\muxmulevelsetvalue}$ 
defined in \eqref{E:UNITLENGTHRTRANS} 
is the $\gfour$-unit normal to $\twoargroughtori{\timefunction,u}{\muxmulevelsetvalue}$ in 
$\hypthreearg{\timefunction}{[- \rightu,\leftu]}{\muxmulevelsetvalue}$. 

In addition, the vectorfield $\hypnormalarg{\muxmulevelsetvalue}$ defined in \eqref{E:HYPNORMALDEF}
is $\gfour$-normal to $\hypthreearg{\timefunction}{[- \rightu,\leftu]}{\muxmulevelsetvalue}$, i.e. 
$\gfour$-orthogonal to $\lbrace \Rtransarg{\muxmulevelsetvalue}, \roughgeop{x^2}, \roughgeop{x^3} \rbrace$. 
Moreover, its square size as measured by $\gfour$ satisfies:
\begin{align} \label{E:HYPNORMALSIZE}
\gfour(\hypnormalarg{\muxmulevelsetvalue},\hypnormalarg{\muxmulevelsetvalue}) 
& = - \frac{\upmu^2}{|\Rtransarg{\muxmulevelsetvalue}|_{\gfour}^2}.
\end{align}
In particular, in the solution regime under study, in which $\upmu > 0$ and $|\Rtransarg{\muxmulevelsetvalue}|_{\gfour} > 0$, 
$\hypnormalarg{\muxmulevelsetvalue}$ is $\gfour$-timelike and $\hypunitnormalarg{\muxmulevelsetvalue}$ is the future directed $\gfour$-unit normal to 
$\hypthreearg{\timefunction}{[- \rightu,\leftu]}{\muxmulevelsetvalue}$, which is consequently \textbf{$\gfour$-spacelike}.
 
Finally, the vectorfield $\hypunitnormalarg{\muxmulevelsetvalue}$ defined in \eqref{E:UNITHYPNORMALDEF} admits the following decomposition: 
\begin{align} \label{E:HYPUNITNORMALDECOMPOSITION}
\hypunitnormalarg{\muxmulevelsetvalue} 
& =  
\frac{|\Rtransarg{\muxmulevelsetvalue}|_{\gfour}}{\upmu} \Lunit 
+ 
\Rtransunitarg{\muxmulevelsetvalue}. 
\end{align}
\end{subequations}
\end{proposition}

\begin{proof}
We first prove the statements regarding $\Rtransarg{\muxmulevelsetvalue}$. 
To derive the orthogonality properties of $\Rtransarg{\muxmulevelsetvalue}$, 
we first use definition \eqref{E:RTRANS}, 
the second identity in \eqref{E:ROUGHANGULARPARTIALDERIVATIVESINTERMSOFGOODGEOMETRICPARTIALDERIVATIVES},
and the fact that $\gfour\left(\Lunit,\roughgeop{x^C} \right)=0$ to obtain:
\begin{align} \label{E:GINNERPRODUCTRTRANSANGEOPC}
\gfour\left(\Rtransarg{\muxmulevelsetvalue},\roughgeop{x^C} \right) 
 & = 
\gfour
\left(\muX, 
	\geop{x^C} 
	- 
	\frac{\geop{x^C}\timefunction}
		{\geop{t}\timefunction} 
	\Lunit 
	+ 
	\frac{\geop{x^C}\timefunction}
		{\geop{t} \timefunctionarg{\muxmulevelsetvalue}} 
	\Lunit^D \geop{x^D} \right) 
- 
\gfour
\left(\upmu 
	\gtorusroughinversefirstfund\left(dx^A,dx^B \right) 
	\frac{\geop{x^A} \timefunction}
		{\geop{t}\timefunction} 
	\roughgeop{x^B},
\roughgeop{x^C}
\right).
\end{align}
Since $\muX$ is $\gfour$-orthogonal to the elements of $\{\geop{x^C}\}_{C = 2,3}$, since $\gfour(\muX,\Lunit) = - \upmu$, 
and since \eqref{E:ROUGHTORUSMETRICCOMPONENTSANDTHEINVERSECOMPONENTSRELATION} implies that:
$$
\gtorusroughinversefirstfund\left(dx^A,dx^B \right)
\gfour\left(\roughgeop{x^B},\roughgeop{x^C}\right)
=
\gtorusroughinversefirstfund\left(dx^A,dx^B \right)
\gtorusroughfirstfund\left(\roughgeop{x^B},\roughgeop{x^C}\right)
=
\updelta_C^A,
$$
we conclude that $\mbox{RHS~\eqref{E:GINNERPRODUCTRTRANSANGEOPC}} = 0$, i.e., that
 $\Rtransarg{\muxmulevelsetvalue}$ 
is $\gfour$-orthogonal to $\twoargroughtori{\timefunction,u}{\muxmulevelsetvalue}$.

Next, using definition \eqref{E:RTRANS},
the relations $\gfour(\muX,\muX) = \upmu^2, \, \gfour(\Lunit,\Lunit) = 0$,
and the $\gfour$-orthogonality of $\Rtransarg{\muxmulevelsetvalue}$ to $\roughgeop{x^A}$,  
we conclude \eqref{E:SIZEOFRTRANS}--\eqref{E:RTRANSNORMSMALLFACTOR}.

Next, in view of definition \eqref{E:HYPNORMALDEF},
we note that the vectorfield $\hypnormalarg{\muxmulevelsetvalue}$ is $\gfour$-orthogonal to 
$\roughgeop{x^2}, \roughgeop{x^3}$ 
because both $\Lunit$ and $\Rtransarg{\muxmulevelsetvalue}$ are. 
Furthermore, we note that the relation $\gfour(\hypnormalarg{\muxmulevelsetvalue},\Rtransarg{\muxmulevelsetvalue})=0$
follows easily from definition \eqref{E:HYPNORMALDEF}
and the relations
$\gfour(\Lunit,\Lunit) = 0$ and $\gfour(\Lunit,\Rtransarg{\muxmulevelsetvalue}) = -\upmu$.
Hence, $\hypnormalarg{\muxmulevelsetvalue}$ is $\gfour$-orthogonal to the set $\lbrace \Rtransarg{\muxmulevelsetvalue}, \roughgeop{x^2}, \roughgeop{x^3} \rbrace$,
which spans the tangent space of $\hypthreearg{\timefunction}{[- \rightu,\leftu]}{\muxmulevelsetvalue}$
(see Remark~\ref{R:WTRANSISTANGENTTOLEVELSETSOFROUGHTIMEFUNCTION}). Therefore, 
$\hypnormalarg{\muxmulevelsetvalue}$ is $\gfour$-orthogonal to $\hypthreearg{\timefunction}{[- \rightu,\leftu]}{\muxmulevelsetvalue}$.

Similarly, \eqref{E:HYPNORMALSIZE} follows from definition \eqref{E:HYPNORMALDEF}
and the relations
$\gfour(\Lunit,\Lunit) = 0$ and $\gfour(\Lunit,\Rtransarg{\muxmulevelsetvalue}) = -\upmu$.

Finally, \eqref{E:HYPUNITNORMALDECOMPOSITION} follows from definition \eqref{E:HYPNORMALDEF}
and \eqref{E:HYPNORMALSIZE}.

\end{proof}

\subsubsection{Decompositions of $\gfour^{-1}$ and $\hypginverse$}
\label{SSS:DECOMPOSITIONSOFINVERSEACOUSTICALMETRICANDFIRSTFUND}

\begin{corollary}[Decompositions of $\gfour^{-1}$ and $\hypginverse$]
\label{C:DECOMPOSITIONSOFINVERSEACOUSTICALMETRICANDFIRSTFUND}
The inverse acoustical metric $\gfour^{-1}$ and the inverse first fundamental form $\hypginverse$ 
of $\hypthreearg{\timefunction}{[-\rightu,\leftu]}{\muxmulevelsetvalue}$
from Def.\,\ref{D:ROUGHFIRSTFUNDS} can be expressed as follows
relative to the vectorfields $\hypunitnormalarg{\muxmulevelsetvalue}$ and $\Rtransunitarg{\muxmulevelsetvalue}$ from Def.\,\ref{D:GEOMETRICVECTORFIELDSADAPTEDTOROUGHFOLIATIONS}
and the inverse first fundamental form $\gtorusroughinversefirstfund$ of $\twoargroughtori{\timefunction,u}{\muxmulevelsetvalue}$
from Def.\,\ref{D:ROUGHFIRSTFUNDS}:
\begin{subequations} 
\begin{align}
\gfour^{-1} 
& 
= 
- 
\hypunitnormalarg{\muxmulevelsetvalue} \otimes \hypunitnormalarg{\muxmulevelsetvalue} 
+  
\Rtransunitarg{\muxmulevelsetvalue} \otimes \Rtransunitarg{\muxmulevelsetvalue} 
+ 
\gtorusroughinversefirstfund, 
\label{E:INVERSEACOUSTICALMETRICINTEMRSOFNUNITANDRUNIT}
		\\
\hypginverse 
& 
= 
\Rtransunitarg{\muxmulevelsetvalue} \otimes \Rtransunitarg{\muxmulevelsetvalue} 
+ 
\gtorusroughinversefirstfund. 
\label{E:INVERSEFIRSTFUNDOFROUGHHYPERSURFACESINTERMSOFR}
\end{align}
\end{subequations}
\end{corollary}

\begin{proof}
	\eqref{E:INVERSEACOUSTICALMETRICINTEMRSOFNUNITANDRUNIT}--\eqref{E:INVERSEFIRSTFUNDOFROUGHHYPERSURFACESINTERMSOFR}
	are straightforward consequences of Prop.\,\ref{P:BASICPROPERTIESOFROUGHVECTORFIELDS}. 
\end{proof}

\subsection{The pointwise semi-norms of tensors with respect to $\gtorusroughfirstfund$ and the $\gtorusroughfirstfund$-trace}
\label{SS:POINTWISESEMINORMSOFTENSORSONROUGHTORI}

\begin{definition}[Pointwise norms] \label{D:POINTWISESEMINORMWITHRESPECTTOFIRSTFUNDOFROUGHTORI} 
If $\upxi$ is a type $\binom{m}{n}$ tensorfield, then we define $|\upxi|_{\gtorusroughfirstfund} \geq 0$ by:
\begin{align} \label{E:SQUAREPOINTWISESEMINORMWITHRESPECTTOFIRSTFUNDOFROUGHTORI}
|\upxi|_{\gtorusroughfirstfund}^2 
& \eqdef 
\gtorusroughfirstfund_{\alpha_1 \widetilde{\alpha}_1} \cdots 
\gtorusroughfirstfund_{\alpha_m \widetilde{\alpha}_m} 
(\gtorusroughinversefirstfund)^{\beta_1 \widetilde{\beta}_1} 
\cdots 
(\gtorusroughinversefirstfund)^{\beta_n \widetilde{\beta}_n} 
\upxi_{\beta_1\cdots\beta_n}^{\alpha_1\cdots\alpha_n} 
\upxi_{\widetilde{\beta}_1 \cdots \widetilde{\beta}_n}^{\widetilde{\alpha}_1\cdots\widetilde{\alpha}_m}.
\end{align}
\end{definition}

\begin{definition} [$\gtorusroughfirstfund$-trace]
\label{D:TRACEOFROUGHTORITANGENT02TENSORS}
If $\upxi$ is a type $\binom{0}{2}$ tensorfield, 
then we define its $\gtorusroughfirstfund$-trace $\mytr_{\gtorusroughfirstfund} \upxi$
as follows:
\begin{align} \label{E:TRACEOFROUGHTORITANGENT02TENSORS}
\mytr_{\gtorusroughfirstfund} \upxi
& 
\eqdef 
(\gtorusroughinversefirstfund)^{\alpha \beta}
\upxi_{\alpha \beta}.
\end{align}
\end{definition}

\subsection{$\gfour$-orthogonal projection onto the rough tori and $\roughangrmD$}
\label{SS:PROJECTIONONTOROUGHTORI}

\begin{definition}[$\gfour$-orthogonal projection 
onto the rough tori $\twoargroughtori{\timefunction,u}{\muxmulevelsetvalue}$ and $\twoargroughtori{\timefunction,u}{\muxmulevelsetvalue}$-tangency] 
\label{D:PROJECTIONONTOROUGHTORIANDROUGHTORITANGENCY}
\hfill
\begin{enumerate}
\item
We define the $\twoargroughtori{\timefunction,u}{\muxmulevelsetvalue}$-projection tensorfield $\roughtorusproject$ as follows, 
where $\hypunitnormalarg{\muxmulevelsetvalue}$ and $\Rtransunitarg{\muxmulevelsetvalue}$ are the vectorfields from Def.\,\ref{D:GEOMETRICVECTORFIELDSADAPTEDTOROUGHFOLIATIONS}
and $\updelta_{\beta}^{\alpha}$ denotes the Kronecker delta:
\begin{align} \label{E:ROUGHTORUSPROJECT}
\roughtorusproject_{\beta}^{\ \alpha} 
& 
\eqdef 
\updelta_{\beta}^{\alpha} 
+ 
\hypunitnormalarg{\muxmulevelsetvalue}_{\beta}
\hypunitnormalarg{\muxmulevelsetvalue}^{\alpha}  
- 
\Rtransunitarg{\muxmulevelsetvalue}_{\alpha}
\Rtransunitarg{\muxmulevelsetvalue}^{\beta}.
\end{align}
\item Given any type $\binom{m}{n}$ spacetime tensorfield $\upxi$, 
we define its $\gfour$-orthogonal projection onto $\twoargroughtori{\timefunction,u}{\muxmulevelsetvalue}$, denoted by $\roughtorusproject \upxi$,
as follows:
\begin{align} \label{E:PROJECTIONOFTENSORONTOROUGHTORI} 
(\roughtorusproject \upxi)_{\beta_1 \cdots \beta_n}^{\alpha_1 \cdots \alpha_m}
& 
\eqdef 
\roughtorusproject_{\widetilde{\alpha}_1}^{\ \alpha_1}
\cdots 
\roughtorusproject_{\widetilde{\alpha}_m}^{\ \alpha_m}
\roughtorusproject_{\beta_1}^{\ \widetilde{\beta}_1}
\cdots 
\roughtorusproject_{ \beta_n}^{\ \widetilde{\beta}_n}
\upxi_{\widetilde{\beta}_1 \cdots \widetilde{\beta}_n}^{\widetilde{\alpha}_1 \cdots \widetilde{\alpha}_m}.
\end{align}
\item
We say that a spacetime tensorfield $\upxi$ is $\twoargroughtori{\timefunction,u}{\muxmulevelsetvalue}$-tangent if 
$\roughtorusproject \upxi = \upxi$. 
\end{enumerate}

\end{definition}

With the help of Lemma~\ref{L:BASICPROPERTIESOFVECTORFIELDS}
and Prop.\,\ref{P:BASICPROPERTIESOFROUGHVECTORFIELDS},
it is straightforward to check that
$\roughtorusproject \hypunitnormalarg{\muxmulevelsetvalue} = \roughtorusproject \Rtransunitarg{\muxmulevelsetvalue} = \roughtorusproject \Lunit = 0$,
while if $Z$ is $\twoargroughtori{\timefunction,u}{\muxmulevelsetvalue}$-tangent, then $\roughtorusproject Z = Z$.
That is, $\roughtorusproject$ \emph{is} the $\gfour$-orthogonal projection onto  $\twoargroughtori{\timefunction,u}{\muxmulevelsetvalue}$.
Moreover, with the help of \eqref{E:INVERSEACOUSTICALMETRICINTEMRSOFNUNITANDRUNIT}, 
it is straightforward to check that
the first fundamental form $\gtorusroughfirstfund$ of $\twoargroughtori{\timefunction,u}{\muxmulevelsetvalue}$ 
from Def.\,\ref{D:ROUGHFIRSTFUNDS}
satisfies:
\begin{align} \label{E:FIRSTFUNDAMENTLFORMOFROUGHTORUSISEQUALTOTHEPROJECTIONOFTHEACOUSTICALMETRIC}
	\gtorusroughfirstfund
	& = 
	\roughtorusproject \gfour.
\end{align}

\begin{definition}[$\twoargroughtori{\timefunction,u}{\muxmulevelsetvalue}$-differential] 
\label{D:ROUGHTORUSDIFFERENTIAL}
Let $\varphi$ is a scalar function. We define $\roughangrmD \varphi$ 
to be the following $\twoargroughtori{\timefunction,u}{\muxmulevelsetvalue}$-tangent one-form:
\begin{align} \label{D:ROUGHTORUSDIFFERENTIAL}
\roughangrmD \varphi 
& \eqdef \roughtorusproject \rmD \varphi.
\end{align}
\end{definition}
Note that 
$[\roughangrmD \varphi]\left(\roughgeop{x^A} \right) 
= \roughgeop{x^A} \varphi$ 
for $A = 2,3$.

\subsection{The Levi-Civita connection $\roughangD$ of $\gtorusroughfirstfund$ and related differential operators}
\label{SS:CONNECTIONSANDDIFFERENTIALOPERATORSONROUGHTORI}

\begin{definition}[The Levi-Civita connection $\roughangD$ of $\gtorusroughfirstfund$ and related differential operators] 
\hfill
\label{D:CONNECTIONSANDDIFFERENTIALOPERATORSONROUGHTORI}
\begin{enumerate}
\item We denote the Levi-Civita connection of $\gtorusroughfirstfund$ by $\roughangD$.
	In particular,
	for $\twoargroughtori{\timefunction,u}{\muxmulevelsetvalue}$-tangent tensorfields $\upxi$, 
	we have $\roughangD \upxi = \roughtorusproject \Dfour \upxi$.
\item If $\upxi$ is an $\twoargroughtori{\timefunction,u}{\muxmulevelsetvalue}$-tangent one-form, 
	then we define its $\twoargroughtori{\timefunction,u}{\muxmulevelsetvalue}$-divergence 
	to be the scalar function
	$\roughangdiv \upxi \eqdef \gtorusroughinversefirstfund \cdot \roughangD \upxi$. 
	Similarly, if $V$ is an $\twoargroughtori{\timefunction,u}{\muxmulevelsetvalue}$-tangent vectorfield, 
	then we define its $\twoargroughtori{\timefunction,u}{\muxmulevelsetvalue}$-divergence to be the scalar function 
	$\roughangdiv V \eqdef \gtorus^{-1} \cdot \roughangD V_{\flat}$,
	where $V_{\flat}$ is the one-form that is $\gfour$-dual to $V$.
\item If $\upxi$ is a symmetric type $\binom{0}{2}$ $\twoargroughtori{\timefunction,u}{\muxmulevelsetvalue}$-tangent tensorfield, 
	then we define its $\twoargroughtori{\timefunction,u}{\muxmulevelsetvalue}$-divergence $\roughangdiv \upxi$ 
	to be the $\twoargroughtori{\timefunction,u}{\muxmulevelsetvalue}$-tangent one-form with the following 
	$\twoargroughtori{\timefunction,u}{\muxmulevelsetvalue}$ components for $A=2,3$:
	\begin{align} \label{E:ROUGHANGDIVOFTYPE02ROUGHTORITANGENTTENSORFIELD}
	[\roughangdiv \upxi]\left(\roughgeop{x^A} \right) 
	\eqdef 
	(\gtorusroughinversefirstfund)(dx^B,dx^C) 
	\left[\roughangDarg{\roughgeop{x^B}} \upxi \right]\left(\roughgeop{x^C},\roughgeop{x^A} \right).
	\end{align}
\end{enumerate}
\end{definition}


\subsection{Curvature tensors of $\gfour$ and $\gtorusroughfirstfund$}
\label{SS:CURVATURETENSORS}
The Riemann curvature tensors of $\gfour$ and $\gtorusroughfirstfund$ 
play a central role in the geometric analysis of the acoustical geometry. 


\begin{definition}[Curvature tensors of $\gfour$ and $\gtorusroughfirstfund$] \label{D:CURVATURETENSORS} The \emph{Riemann curvature tensor} $\Riemfour$ of the acoustical metric $\gfour$ is the type $\binom{0}{4}$ spacetime tensorfield 
defined by:
\begin{align} \label{E:ACOUSTICALCURVATURETENSOR}
\Riemfour(\mathbf{X},\mathbf{Y},\mathbf{Z},\mathbf{W}) 
& \eqdef \gfour(- \Dfour^2_{\mathbf{X} \mathbf{Y}} \mathbf{Z} 
+ 
\Dfour^2_{\mathbf{Y} \mathbf{X}} \mathbf{Z}, \mathbf{W}),
\end{align}
where $X$, $Y$, $Z$, and $W$ are arbitrary spacetime vectorfields, and
$\Dfour^2_{\mathbf{X} \mathbf{Y}} \mathbf{Z} 
\eqdef 
\mathbf{X}^{\alpha} \mathbf{Y}^\beta \Dfour_\alpha \Dfour_\beta \mathbf{Z}$. 

The \emph{Ricci curvature tensor} $\Ricfour$ of the acoustical metric $\gfour$ is the type $\binom{0}{2}$ spacetime tensor defined 
relative to arbitrary coordinates as follows: 
\begin{align} \label{E:SPACETIMERICCITENSOR}
\Ricfour_{\alpha \beta} 
& \eqdef 
(\gfour^{-1})^{\nu\sigma} \Riemfour_{\alpha \nu \beta\sigma}.
\end{align}

Similarly, the Riemann curvature tensor $\Riemtorus$ of the Riemannian metric $\gtorusroughfirstfund$ on 
$\twoargroughtori{\timefunction,u}{\muxmulevelsetvalue}$ is the type $\binom{0}{4}$ $\twoargroughtori{\timefunction,u}{\muxmulevelsetvalue}$ tensorfield defined 
as follows:
\begin{align} \label{E:RIEMANNCURVATURETENSOROFROUGHTORUS}
\Riemtorus(X,Y,Z,W) 
& \eqdef \gtorusroughfirstfund(- \roughangD^2_{XY}Z + \roughangD^2_{YX}Z,W),
\end{align}
where $X$, $Y$, $Z$, and $W$ are arbitrary $\twoargroughtori{\timefunction,u}{\muxmulevelsetvalue}$-tangent vectorfields 
and $\roughangD$ is the Levi-Civita connection of $\gtorusroughfirstfund$. 

The \emph{Ricci curvature tensor} $\Riemtorus$ of $\gtorusroughfirstfund$ is the type $\binom{0}{2}$ spacetime tensor defined 
relative to arbitrary coordinates as follows: 
\begin{align} \label{E:ROUGHTORUSRICCITENSOR}
\Rictorus_{\alpha \beta} 
& \eqdef 
(\gtorusroughinversefirstfund)^{\nu \sigma} \Riemtorus_{\alpha \nu \beta \sigma}.
\end{align}

The \emph{Scalar curvature} $\Scalartorus$ of $\gtorusroughfirstfund$ 
is the scalar function defined 
relative to arbitrary coordinates as follows: 
\begin{align} \label{E:ROUGHTORUSSCALARCURVATURE}
\Scalartorus
& \eqdef 
(\gtorusroughinversefirstfund)^{\alpha \beta} \Rictorus_{\alpha \beta}.
\end{align}
\end{definition}

It is well-known that because $\twoargroughtori{\timefunction,u}{\muxmulevelsetvalue}$ is two-dimensional,
the \emph{Gauss curvature} $\Gausstorus$ of $\gtorusroughfirstfund$ can be expressed as follows
in terms of it scalar curvature:
\begin{align} \label{E:2DGAUSSCURVATUREISTWICESCALARCURVATURE}
	\Gausstorus
	& = \frac{1}{2} \Scalartorus.
\end{align}


\section{The acoustic double-null frame and its relationship with the rough acoustical geometry}
\label{S:ACOUSTICDOUBLENULLFRAMEANDITSRELATIONWITHTHEROUGHGEOMETRY}
To control the top-order derivatives of $\vortrenormalized$ and $\GradEnt$,
we will rely on a family of ``elliptic-hyperbolic'' integral identities
that we derive in Sect.\,\ref{S:ELLIPTICHYPERBOLICIDENTITIES}.
In this section, we construct the acoustic double-null frame
that we use to derive the elliptic-hyperbolic identities.
Moreover, we provide various identities that relate the acoustic double-null frame
to the rough acoustical geometry constructed in Sect.\,\ref{S:ROUGHACOUSTICALGEOMETRYANDCURVATURETENSORS}.

\subsection{The acoustic double-null frame}
\label{SS:ACOUSTICDOUBLENULLFRAME}
The new ingredient in the acoustic double null frame is
the vectorfield $\uLunit$, which we now define.
In Lemma~\ref{L:BASICPROPERTIESOFULUNIT}, we will show that $\uLunit$ is 
$\gfour$-null and transversal to the characteristics $\nullhyparg{u}$. 

\begin{definition}[The vectorfield $\uLunit$] 
\label{D:ULUNIT}
We define $\uLunit$ to be the vectorfield whose Cartesian components are:  
\begin{align} \label{E:ULUNIT}
	\uLunit^{\alpha}
	& 
	\eqdef
	\Lunit^{\alpha}
	+
	2 X^{\alpha}.
\end{align}
\end{definition} 

\begin{lemma}[Basic properties of $\uLunit$]
	\label{L:BASICPROPERTIESOFULUNIT}
	The vectorfields $\uLunit$ and $\Transport$ satisfy the following identities:
	\begin{align} \label{E:SIMPLEIDENTITIESINVOLVINGULUNIT}
	\uLunit & = \Transport + X,
	&
	\Transport 
	& = \frac{1}{2}
			\left(
				\Lunit + \uLunit
			\right)
\end{align}

Moreover, $\uLunit$ is $\gfour$-null and transversal to the characteristics $\nullhyparg{u}$:
\begin{align} \label{E:ULUNITISNULLANDNORMALIZEDAGAINSTLUNIT}
	\gfour(\uLunit,\uLunit)
	& 
	=
	0,
		\\
\gfour(\uLunit,\Lunit) & = - 2. \label{E:INNERPRODUCTOFLANDULISMINUS2}
\end{align}

In addition, the acoustical metric $\gfour$ and its inverse $\gfour^{-1}$ satisfy the following identities,
where $\gtorus$ is the first fundamental form of $\ell_{t,u}$ from Def.\,\ref{D:FIRSTFUNDAMENTALFORMS}:
\begin{align} \label{E:ACOUSTICALMETRICINDOUBLENULLFRAME}
	\gfour_{\alpha \beta} 
	& = 
	- 
	\frac{1}{2} \Lunit_{\alpha} \uLunit_{\beta} 
	- 
	\frac{1}{2} \uLunit_{\alpha} \Lunit_{\beta} 
	+ 
	\gtorus_{\alpha \beta}, 
	& 
	(\gfour^{-1})^{\alpha \beta} 
	& =  
	- 
	\frac{1}{2} \Lunit^{\alpha} \uLunit^\beta 
	- 
	\frac{1}{2} \uLunit^{\alpha} \Lunit^{\beta} 
	+ 
	(\gtorus^{-1})^{\alpha \beta}. 
\end{align}

Finally, we have the following identity, where
$\updelta_{\beta}^{\alpha}$ is the Kronecker delta and $\smoothtorusproject$ is the 
$\ell_{t,u}$-projection tensorfield from Def.\,\ref{D:PROJECTIONTENSORFIELDSANDTANGENCYTOHYPERSURFACES}:
\begin{align} \label{E:ELLTUPROJECTIONINDOUBLENULLFRAME}
	\updelta_{\beta}^{\alpha}
	& =  
	- 
	\frac{1}{2} \Lunit^{\alpha} \uLunit_{\beta}
	- 
	\frac{1}{2} \uLunit^{\alpha} \Lunit_{\beta}{\beta} 
	+ 
	\smoothtorusproject_{\beta}^{\ \alpha}. 
\end{align}

\end{lemma}

\begin{proof}
Equations \eqref{E:SIMPLEIDENTITIESINVOLVINGULUNIT}--\eqref{E:INNERPRODUCTOFLANDULISMINUS2} are straightforward consequences of the definition of $\uLunit$ in \eqref{E:ULUNIT}, the identity $\Transport = \Lunit + X$ (see \eqref{E:BISLPLUSX}), 
and 
the identities
$\gfour(\Lunit,\Lunit) = 0$, 
$\gfour(\Lunit,X) = -1$,
and $\gfour(X,X) = 1$
from Lemma~\ref{L:BASICPROPERTIESOFVECTORFIELDS}.
The statement that $\uLunit$ is transversal to $\nullhyparg{u}$ follows from \eqref{E:INNERPRODUCTOFLANDULISMINUS2} 
and the fact that 
$\Lunit$ is $\gfour$-orthogonal to $\nullhyparg{u}$. 
The identity for $\gfour_{\alpha \beta}$
in \eqref{E:ACOUSTICALMETRICINDOUBLENULLFRAME} follows from contracting both sides of the identities against 
the frame $\left\lbrace\Lunit,\uLunit,\geop{x^2},\geop{x^3}\right\rbrace$ and computing that both sides are equal.
The identity for $(\gfour^{-1})^{\alpha \beta}$
in \eqref{E:ACOUSTICALMETRICINDOUBLENULLFRAME} follows from
raising the indices in the first identity with $\gfour^{-1}$.
The identity \eqref{E:ELLTUPROJECTIONINDOUBLENULLFRAME}
follows from using $\gfour$ to lower the $\beta$ index in the last identity stated in \eqref{E:ACOUSTICALMETRICINDOUBLENULLFRAME}
and using \eqref{E:LTUSPROJECTIONISONEINDEXLOWERINGOFGTORUSINVERSE}.
\end{proof}

\begin{definition}[The acoustic double-null frame]
	\label{D:ACOUSTICDOUBLENULLFRAME}
	We refer to $\left\lbrace\Lunit,\uLunit,\geop{x^2},\geop{x^3}\right\rbrace$ as the \emph{acoustic double-null frame}. 
\end{definition}

\subsection{Identities involving the acoustic double-null frame and the rough geometry}
\label{SS:IDENTITIESINVOLVINGACOUSTICDOUBLENULLFRAMEANDROUGHGEOMETRY}
The following lemma provides several identities involving the acoustic double-null frame 
and the rough acoustical geometry. 

\begin{lemma}[Identities involving the acoustic double-null frame and the rough acoustic geometry]
\label{L:IDENTITIESINVOLVINGROUGHGEOMETRY}
Let $\Rtransnormsmallfactorarg{\muxmulevelsetvalue}$ be as in \eqref{E:RTRANSNORMSMALLFACTOR}, and 
recall that $\angD$ is the Levi-Civita connection of $\gtorus$, the first
fundamental form of $\ell_{t,u} = \Sigma_t \cap \nullhyparg{u}$
relative to $\gfour$. 
Then the following identity holds:
\begin{align} \label{E:IDENTITYFORRTRANSNORMSMALLFACTORSQUARED}
\Rtransnormsmallfactorarg{\muxmulevelsetvalue} 
& = 
\frac{|\angD \timefunctionarg{\muxmulevelsetvalue}|_{\gtorus}^2}{(\Lunit \timefunctionarg{\muxmulevelsetvalue})^2}.
\end{align}

Moreover, the $\twoargroughtori{\timefunction,u}{\muxmulevelsetvalue}$-tangent vectorfield 
$\Roughtoritangentvectorfieldarg{\muxmulevelsetvalue}$ 
defined in \eqref{E:DEFROUGHTORITANGENTVECTORFIELD} 
admits the following decomposition: 
	\begin{align} \label{E:ROUGHTORITANGENTVECTORFIELDKEEPSAPPEARING}
		\Roughtoritangentvectorfieldarg{\muxmulevelsetvalue}
		& = 
			-
			\Rtransnormsmallfactorarg{\muxmulevelsetvalue} \Lunit
			+
			\frac{1}{\Lunit \timefunctionarg{\muxmulevelsetvalue}}
			\angDuparg{\#} \timefunctionarg{\muxmulevelsetvalue}.
	\end{align}
	
In addition, with $\Rtransarg{\muxmulevelsetvalue}$ the vectorfield defined in \eqref{E:RTRANS}, 
the vectorfield $\frac{1}{\upmu}\Rtransarg{\muxmulevelsetvalue}$
admits the following decomposition:
\begin{align} \label{E:RTRANSDIVIDEDBYMUIDENTITY}
	\frac{1}{\upmu} \Rtransarg{\muxmulevelsetvalue}
	& = 
	X
	+
	\left(
		\phi \frac{\muxmulevelsetvalue}{\upmu \Lunit \upmu}
		+
		\Rtransnormsmallfactorarg{\muxmulevelsetvalue} 
	\right)
	\Lunit
	-
	\frac{1}{\Lunit \timefunctionarg{\muxmulevelsetvalue}}
	\angDuparg{\#} \timefunctionarg{\muxmulevelsetvalue},
\end{align}
where we recall that the $\ell_{t,u}$-tangent vectorfield $\angDuparg{\#} \timefunctionarg{\muxmulevelsetvalue}$ is the dual of 
$\angD \timefunctionarg{\muxmulevelsetvalue}$ with respect to $\gfour$.

Furthermore, the following differentiation identities hold: 
		\begin{align}
		\Rtransarg{\muxmulevelsetvalue}
		\left[
			\frac{1}
			{\upmu - \phi \frac{\muxmulevelsetvalue}{\Lunit \upmu}}
		\right]
		& = - 
				\frac{\Rtransarg{\muxmulevelsetvalue} \upmu}{(\upmu - \phi \frac{\muxmulevelsetvalue}{\Lunit \upmu})^2}
				+ 
				\frac{\muxmulevelsetvalue \frac{\phi'}{\Lunit \upmu}}{(\upmu - \phi \frac{\muxmulevelsetvalue}{\Lunit \upmu})^2}
				-
				\frac{\phi \muxmulevelsetvalue \frac{\Rtransarg{\muxmulevelsetvalue} \Lunit \upmu}{(\Lunit \upmu)^2}}{(\upmu - \phi \frac{\muxmulevelsetvalue}{\Lunit \upmu})^2},
					\label{E:RTRANSDERIVATIVEOFSINGULARWEIGHT} 
					\\
		\Roughtoritangentvectorfieldarg{\muxmulevelsetvalue}
		\left[
			\frac{\upmu}
			{\upmu - \phi \frac{\muxmulevelsetvalue}{\Lunit \upmu}}
		\right]
		& = 
				-
				\frac{(\Roughtoritangentvectorfieldarg{\muxmulevelsetvalue} \upmu) \phi \frac{\muxmulevelsetvalue}{\Lunit \upmu}}
				{(\upmu - \phi \frac{\muxmulevelsetvalue}{\Lunit \upmu})^2}
				-
				\frac{\upmu \phi \muxmulevelsetvalue \frac{\Roughtoritangentvectorfieldarg{\muxmulevelsetvalue} \Lunit \upmu}
				{(\Lunit \upmu)^2}}{(\upmu - \phi \frac{\muxmulevelsetvalue}{\Lunit \upmu})^2}.
			\label{E:ROUGHTORIDERIVATIVEOFNONSINGULARWEIGHT}
	\end{align}	
	
In addition, we have the following decompositions for $\Lunit$ and $\uLunit$:
\begin{subequations}
	\begin{align}
		\Lunit
		& = \frac{\upmu}{\upmu - \phi \frac{\muxmulevelsetvalue}{\Lunit \upmu}}
				\Transport
				-
				\frac{1}{\upmu - \phi \frac{\muxmulevelsetvalue}{\Lunit \upmu}}
				\Rtransarg{\muxmulevelsetvalue}
				- 
				\frac{\upmu}{\upmu - \phi \frac{\muxmulevelsetvalue}{\Lunit \upmu}}
				\Roughtoritangentvectorfieldarg{\muxmulevelsetvalue},
				\label{E:LUNITINTERMSOFSIGMATILDETANGENTVECTORFIELDSANDTRANSPORT} \\
		\uLunit
		& = \frac{\upmu - \phi \frac{2 \muxmulevelsetvalue}{\Lunit \upmu}}{\upmu - \phi \frac{\muxmulevelsetvalue}{\Lunit \upmu}}
				\Transport
				+
				\frac{1}{\upmu - \phi \frac{\muxmulevelsetvalue}{\Lunit \upmu}}
				\Rtransarg{\muxmulevelsetvalue}
				+ 
				\frac{\upmu}{\upmu - \phi \frac{\muxmulevelsetvalue}{\Lunit \upmu}}
				\Roughtoritangentvectorfieldarg{\muxmulevelsetvalue}.
				\label{E:ULUNITINTERMSOFSIGMATILDETANGENTVECTORFIELDSANDTRANSPORT}
	\end{align}
	\end{subequations}
	
Moreover, the following identities hold:
\begin{subequations}
\begin{align}  \label{E:INNERPRODUCTOFRTRANSANDLUNITANDULUNIT}
		\Rtransarg{\muxmulevelsetvalue}^{\alpha} \Lunit_{\alpha}
		& = - \upmu,
		&
		\Rtransarg{\muxmulevelsetvalue}^{\alpha} \uLunit_{\alpha}
		& = \upmu(1 - 2 \Rtransnormsmallfactorarg{\muxmulevelsetvalue}) - 2 \phi \frac{\muxmulevelsetvalue}{\Lunit \upmu}, 
				\\
		\Rtransarg{\muxmulevelsetvalue}^{\alpha} \Transport_{\alpha}
		& = - \left(\upmu \Rtransnormsmallfactorarg{\muxmulevelsetvalue} + \phi \frac{\muxmulevelsetvalue}{\Lunit \upmu} \right).
		\label{E:INNERPRODUCTOFRTRANSANDULUNIT}
		&&
	\end{align}
	\end{subequations}
	
Finally, the following identities hold for any $\Sigma_t$-tangent vectorfield $\SigmatTan$,
where $g$ is the first fundamental form of $\Sigma_t$ relative to $\gfour$:
\begin{subequations}
\begin{align} \label{E:INNERPRODUCTOFRTRANSANDSIGMATTTANGENTVECTORFIELD} 
			\Rtransarg{\muxmulevelsetvalue}_{\alpha} \SigmatTan^{\alpha}
		& = 
				\left\lbrace
					\upmu(1 - \Rtransnormsmallfactorarg{\muxmulevelsetvalue}) 
					-
					\phi \frac{\muxmulevelsetvalue}{\Lunit \upmu}
				\right\rbrace
				X_a \SigmatTan^a
				-
				\frac{\upmu}{\Lunit \timefunctionarg{\muxmulevelsetvalue}}
				\angVarg{\alpha}
				\angDarg{\alpha}
				\timefunctionarg{\muxmulevelsetvalue},
					\\
				-
				\frac{1}{2}
				\Rtransarg{\muxmulevelsetvalue}_{\alpha} \SigmatTan^{\alpha} \SigmatTan^a X_a
				+
				\frac{1}{4}
				\Rtransarg{\muxmulevelsetvalue}_{\alpha} \uLunit^{\alpha}
				|\SigmatTan|_g^2
				& = 
					- \frac{1}{4} \upmu
							(\SigmatTan^a X_a)^2
						+
						\frac{1}{4}
						\left\lbrace
							\upmu(1 - 2 \Rtransnormsmallfactorarg{\muxmulevelsetvalue}) 
							- 
							2 \phi \frac{\muxmulevelsetvalue}{\Lunit \upmu}
						\right\rbrace
						|\angV|_{\gtorus}^2
					+
				\frac{1}{2}
				\frac{\upmu}{\Lunit \timefunctionarg{\muxmulevelsetvalue}}
				X_a \SigmatTan^a
				\angVarg{\alpha}
				\angDarg{\alpha}
				\timefunctionarg{\muxmulevelsetvalue}.
				\label{E:RTRANSCONTRACTIONIDENTITYINPROOFOFELLIPTICBOUNDARYTERMIDENTITY}
	\end{align}  	
	\end{subequations}
\end{lemma}

\begin{proof}
To prove \eqref{E:IDENTITYFORRTRANSNORMSMALLFACTORSQUARED}, 
we expand
$\frac{|\angD \timefunctionarg{\muxmulevelsetvalue}|_{\gtorus}^2}{(\Lunit \timefunctionarg{\muxmulevelsetvalue})^2} 
=
\frac{1}{(\Lunit \timefunctionarg{\muxmulevelsetvalue})^2} 
(\gtorus^{-1})^{AB} 
\left(\geop{x^A} \timefunctionarg{\muxmulevelsetvalue}\right)
\geop{x^B}\timefunctionarg{\muxmulevelsetvalue}
$ 
and use the identities
\eqref{E:SMOOTHTORUSINVERSEFIRSTFUNDCOMPONENTSINTERMSOFROUGHTORUSINVERSEFIRSTFUNDCOMPONENTS}--\eqref{E:INVERSECHOVCOEFFICIENTSSMOOTHANGULARDERIVATIVESINTERMSOFROUGHONESANDL}
and $\Lunit = \geop{t} + \Lunit^A \geop{x^A}$,
thereby confirming that 
$\frac{|\angD \timefunctionarg{\muxmulevelsetvalue}|_{\gtorus}^2}{(\Lunit \timefunctionarg{\muxmulevelsetvalue})^2} = \mbox{RHS~\eqref{E:RTRANSNORMSMALLFACTOR}}$ as desired.
\eqref{E:ROUGHTORITANGENTVECTORFIELDKEEPSAPPEARING} then follows from a similar argument based on
expanding 
$
-
\Rtransnormsmallfactorarg{\muxmulevelsetvalue} \Lunit 
+ 
\frac{1}{\Lunit \timefunctionarg{\muxmulevelsetvalue}} \angD^{\#}\timefunctionarg{\muxmulevelsetvalue} 
= 
- 
\frac{|\angD \timefunctionarg{\muxmulevelsetvalue}|^2}{(\Lunit \timefunctionarg{\muxmulevelsetvalue})^2} \Lunit 
+ 
\frac{1}{\Lunit \timefunctionarg{\muxmulevelsetvalue}} 
(\gtorus^{-1})^{AB} 
\left(
\geop{x^A} 
\timefunctionarg{\muxmulevelsetvalue}
\right) 
\geop{x^B}
$
and using
\eqref{E:CHOVCOEFFICIENTSSMOOTHANGULARDERIVATIVESINTERMSOFROUGHONESANDL}--\eqref{E:SMOOTHANGULARDERIVATIVESINTERMSOFROUGHONESANDL}
and
\eqref{E:SMOOTHTORUSINVERSEFIRSTFUNDCOMPONENTSINTERMSOFROUGHTORUSINVERSEFIRSTFUNDCOMPONENTS}--\eqref{E:INVERSECHOVCOEFFICIENTSSMOOTHANGULARDERIVATIVESINTERMSOFROUGHONESANDL}. 

The expression \eqref{E:RTRANSDIVIDEDBYMUIDENTITY} for $\frac{1}{\upmu} \Rtransarg{\muxmulevelsetvalue}$ follows from \eqref{E:RTRANS} and the already proved
identity \eqref{E:ROUGHTORITANGENTVECTORFIELDKEEPSAPPEARING}. 

The identities \eqref{E:RTRANSDERIVATIVEOFSINGULARWEIGHT}--\eqref{E:ROUGHTORIDERIVATIVEOFNONSINGULARWEIGHT} are straightforward consequences of the chain and Leibniz rules
and the fact that $\Rtransarg{\muxmulevelsetvalue} u = 1$, which follows from Lemma~\ref{L:BASICPROPERTIESOFVECTORFIELDS},
\eqref{E:RTRANS}, and \eqref{E:ROUGHTORITANGENTVECTORFIELDKEEPSAPPEARING}.

The decompositions of $\Lunit$ and $\uLunit$ stated in
\eqref{E:LUNITINTERMSOFSIGMATILDETANGENTVECTORFIELDSANDTRANSPORT}--\eqref{E:ULUNITINTERMSOFSIGMATILDETANGENTVECTORFIELDSANDTRANSPORT} follow from \eqref{E:RTRANS} and the identities $\Transport = \Lunit + X$, $\uLunit = \Transport + X$, and $\Transport  = \frac{1}{2}( \Lunit + \uLunit)$. 

Next, we note that if $\SigmatTan$ is $\Sigma_t$-tangent, then since 
$\Transport = \Lunit + X$, it follows from 
\eqref{E:MATERIALDERIVATIVELOWEREDCARTESIANCOORDINATES} 
that $\Lunit_{\alpha} \SigmatTan^{\alpha}= - X_a V^a$.
From this identity and \eqref{E:RTRANSDIVIDEDBYMUIDENTITY}, we deduce
\eqref{E:INNERPRODUCTOFRTRANSANDSIGMATTTANGENTVECTORFIELD}. 

The identities in \eqref{E:INNERPRODUCTOFRTRANSANDLUNITANDULUNIT} follow from 
Lemma~\ref{L:BASICPROPERTIESOFVECTORFIELDS},
\eqref{E:RTRANS}, 
definition~\eqref{E:ULUNIT},
and \eqref{E:ROUGHTORITANGENTVECTORFIELDKEEPSAPPEARING}.
Moreover, \eqref{E:INNERPRODUCTOFRTRANSANDULUNIT} follows from the same arguments.

Finally, the identity \eqref{E:RTRANSCONTRACTIONIDENTITYINPROOFOFELLIPTICBOUNDARYTERMIDENTITY} follows from
\eqref{E:INNERPRODUCTOFRTRANSANDLUNITANDULUNIT},
\eqref{E:INNERPRODUCTOFRTRANSANDSIGMATTTANGENTVECTORFIELD},
and the decomposition $g_{ab} = \gtorus_{ab} + X_a X_b$, which follows from \eqref{E:SMOOTHTORUSMETRICINTERMSOFSIGMATMETRICANDX}.
\end{proof}


\section{Norms, area and volume forms, and strings of commutation vectorfields}
\label{S:NORMSANDFORMS}
In this section, we define various norms on regions of spacetime that are tied to the rough geometry.
We also introduce the area and volume forms that we use in our $L^2$ analysis.
Finally, we introduce notation for repeated differentiation with respect to the commutation vectorfields.

\subsection{$L^{\infty}$-type Sobolev norms and H\"{o}lder norms} 
\label{SS:LINFTYTYPESOBOLEVNORMSANDHOLDERNORMS}

\subsubsection{Multi-index notation in various coordinate systems}
\label{SSS:MULTIINDEXCOORDINATENOTATION}

\begin{definition}[Multi-index notation in various coordinate systems]  \label{D:MULTIINDEXCOORDINATENOTATION}
Let $\alpha_1,\alpha_2,\alpha_3,\alpha_4 \in \mathbb{N}$, and let 
$\vec{\alpha} = (\alpha_1,\alpha_2,\alpha_3,\alpha_4)$ be the corresponding
multi-index of order $|\vec{\alpha}| \eqdef \sum_{i=1}^4 \alpha_i$.
We define the following order $|\vec{\alpha}|$ differential operator with respect to the geometric coordinates:
\begin{align} \label{E:GEOMETRICCOORDINATEMULTIINDEXDIFFERENTIALOPERATOR}
\frac{\p^{\vec{\alpha}}}{\p(t,u,x^2,x^3)} 
& \eqdef 
\left(\geop{t}\right)^{\alpha_1}
\left(\geop{u}\right)^{\alpha_2} 
\left(\geop{x^2}\right)^{\alpha_3} 
\left(\geop{x^3}\right)^{\alpha_4}.
\end{align} 
Similarly, we define the following order $|\vec{\alpha}|$ differential operator with respect to the rough adapted coordinates:
\begin{align} \label{E:ROUGHADAPATEDCOORDINATEMULTIINDEXDIFFERENTIALOPERATOR}
\frac{\widetilde{\p}^{\vec{\alpha}}}{\widetilde{\p}(\timefunctionarg{\muxmulevelsetvalue},u,x^2,x^3)} 
& \eqdef 
\left(\roughgeop{\timefunctionarg{\muxmulevelsetvalue}}\right)^{\alpha_1} 
\left(\roughgeop{u}\right)^{\alpha_2} 
\left( \roughgeop{x^2}\right)^{\alpha_3} 
\left(\roughgeop{x^3}\right)^{\alpha_4}.
\end{align}
\end{definition}

\subsubsection{Definitions of the essential sup-norm-type Sobolev norms and H\"{o}lder norms} 
\label{SSS:LINFTYTYPESOBOLEVNORMSANDHOLDERNORMS}

\begin{definition}[$L^{\infty}$-type Sobolev norms and H\"{o}lder norms] 
\label{D:ESSENTIALSUPNORMTYPENORMSANDHOLDERNORMS}
Let $f$ be a scalar function, and let $m \geq 0$ be an integer.
On the spacetime regions $\twoargMrough{I,J}{\muxmulevelsetvalue}$ 
defined in \eqref{E:TRUNCATEDMROUGH},
we define the following $L^{\infty}$-type Sobolev norms and H\"{o}lder norms of $f$
relative to the geometric coordinates: 
\begin{subequations}
\begin{align}
		\| f \|_{W_{\textnormal{geo}}^{m,\infty}\left(\twoargMrough{I,J}{\muxmulevelsetvalue}\right)} 
		& 
		\eqdef \sum_{|\vec{\alpha}| \leq m} \esssup_{p \in \twoargMrough{I,J}{\muxmulevelsetvalue}} 
		\left|
			\frac{\p^{\vec{\alpha}} f(p)}{\p(t,u,x^2,x^3)} 
		\right|, 
			\label{E:WKINFTYGEONORMS} \\
		\| f \|_{C^{m,1}_{\textnormal{geo}}\left(\twoargMrough{I,J}{\muxmulevelsetvalue}\right)} 
		& 
		\eqdef \sum_{|\alpha| \leq m} 
		\max_{p \in \twoargMrough{I,J}{\muxmulevelsetvalue}} 
		\left| 
			\frac{\p^{\vec{\alpha}} f(p)}{\p(t,u,x^2,x^3)} 
		\right| 
		+ 
		\sum_{|\alpha| = \kappa} \sup_{\substack{p_1,p_2 \in \twoargMrough{I,J}{\muxmulevelsetvalue} \\ p_1 \neq p_2}} 
			\frac{
			\left|
				\frac{\p^{\vec{\alpha}} f(p_1)}{\p(t,u,x^2,x^3)} 
				- 
				\frac{\p^{\vec{\alpha}}  f(p_2)}{\p(t,u,x^2,x^3)} 
			\right|}
				{\mbox{\upshape dist}_{\mbox{\upshape geo}}(p_1,p_2)},
		\label{E:CK1GEONORMS}
\end{align}
\end{subequations}
where $\mbox{\upshape dist}_{\mbox{\upshape geo}}(p_1,p_2)$ is the
standard Euclidean distance between $p_1$ and $p_2$ in the flat geometric coordinate space 
$\R_t \times \R_u \times \T^2$,
i.e., if
$p_i \eqdef (t_i,u_i,x_i^2,x_i^3)$,
$\Delta t \eqdef t_2 - t_1$,
and
$\Delta u \eqdef u_2 - u_1$,
then 
$\mbox{\upshape dist}_{\mbox{\upshape geo}}(p_1,p_2)
\eqdef
\sqrt{
|\Delta t|^2
+
|\Delta u|^2 
+ 
|\Delta x^2|_{\mathbb{T}}^2 
+ 
|\Delta x^3|_{\mathbb{T}}^2}
$,
where for $j=2,3$,
$|\Delta x^j|_{\mathbb{T}}$ is the Euclidean distance between $x_2^j$ and $x_1^j$ in the torus.

Similarly, for intervals $I,J \in \mathbb{R}$, 
we define the following norms in the rough adapted coordinate spacetime region 
$I \times J \times \T^2 \subset \R_{\timefunction} \times \R_u \times \T^2$: 
\begin{align} \label{E:CK1ROUGHNORMS}
	\| f \|_{C^{m,1}_{\textnormal{rough}}\left(I\times J \times \T^2 \right)}  
	& \eqdef \sum_{|\vec{\alpha}| \leq m}  
	\max_{q \in I \times J \times \T^2} 
	\left| 
		\frac{\widetilde{\p}^{\alpha}f(q)}{\widetilde{\p}(\timefunctionarg{\muxmulevelsetvalue},u,x^2,x^3)}  
	\right| 
	+ 
	\sum_{|\vec{\alpha}| = m} 
	\sup_{\substack{q_1,q_2 \in I \times J \times \T^2 \\ q_1 \neq q_2}} 
	\frac{ \left| \frac{\widetilde{\p}^{\vec{\alpha}}  f(q_1)}{\widetilde{\p}(\timefunctionarg{\muxmulevelsetvalue},u,x^2,x^3)}  
	- 
	\frac{\widetilde{\p}^{\vec{\alpha}}f(q_2)}{\widetilde{\p}(\timefunctionarg{\muxmulevelsetvalue},u,x^2,x^3)} 
	\right|}{\mbox{\upshape dist}_{\mbox{\upshape rough}}(q_1,q_2)},
\end{align}
where $\mbox{\upshape dist}_{\mbox{\upshape rough}}(q_1,q_2)$ is the
standard Euclidean distance between $q_1$ and $q_2$ in the flat rough adapted coordinate space 
$\R_{\timefunction} \times \R_u \times\T^2$,
i.e., if
$q_i \eqdef (\timefunction_i,u_i,x_i^2,x_i^3)$,
$\Delta \timefunction \eqdef \timefunction_2 - \timefunction_1$,
and
$\Delta u \eqdef u_2 - u_1$,
then 
$\mbox{\upshape dist}_{\mbox{\upshape rough}}(q_1,q_2)
\eqdef
\sqrt{
|\Delta \timefunction|^2
+
|\Delta u|^2 
+ 
|\Delta x^2|_{\mathbb{T}}^2 
+ 
|\Delta x^3|_{\mathbb{T}}^2}
$.

If $\vec{\varphi} = \lbrace \varphi_i \rbrace_{i=1,\cdots,M}$ is an array comprising $M$ scalar functions, 
then we extend the definitions of the above norms to $\vec{\varphi}$ by summing the scalar norms over $i$, e.g.,
$\| \vec{\varphi} \|_{C^{m,1}_{\textnormal{geo}}\left(\twoargMrough{I,J}{\muxmulevelsetvalue}\right)}
\eqdef \sum_{i=1}^M \| \varphi_i \|_{C^{m,1}_{\textnormal{geo}}\left(\twoargMrough{I,J}{\muxmulevelsetvalue}\right)}$.

Finally, if $\mathbf{J}$ is a matrix with scalar function entries defined on $\twoargMrough{I,J}{\muxmulevelsetvalue}$, 
then we extend the definitions of the above norms to $\mathbf{J}$ by $\| \mathbf{J}\|_{C^{m,1}_{\textnormal{geo}}\left(\twoargMrough{I,J}{\muxmulevelsetvalue}\right)} = \left\| | \mathbf{J} |_{\mbox{\upshape}Euc} \right\|_{C^{m,1}_{\textnormal{geo}}\left(\twoargMrough{I,J}{\muxmulevelsetvalue}\right)}$, where $|\cdot|_{\mbox{\upshape}Euc}$ is the standard Frobenius norm on matrices. 
Similar remarks apply to the $C_{\textnormal{rough}}^{m,1}$ norms of matrices defined on the rough coordinate
slab $I \times J \times \T^2$.
\end{definition}

\subsection{Area forms, volume forms, and corresponding $L^2$ norms} 
\label{SS:FORMSANDL2NORMS}
We now define the area and volume forms on the rough subsets (see Def.\,\ref{D:TRUNCATEDROUGHSUBSETS})
$\twoargroughtori{\timefunction,u}{\muxmulevelsetvalue}$,
 $\hypthreearg{\timefunction}{I}{\muxmulevelsetvalue}$, 
and
$\twoargMrough{I,J}{\muxmulevelsetvalue}$
that we will use in our analysis.
We also define corresponding $L^2$-type norms.
Our definitions are in terms of the rough adapted coordinates $(\timefunctionarg{\muxmulevelsetvalue},u,x^2,x^3)$
because those are the coordinates that we use
in our energy identities, where the forms arise.

\subsubsection{Geometric forms and related integrals}
\label{SSS:GEOMETRICFORMSANDINTEGRALS}

\begin{definition}[Geometric forms and related integrals]\label{D:ROUGHVOLFORMS}
\hfill
\begin{itemize}
\item Recall that $\gtorusroughfirstfund$ is the first fundamental form of the rough torus $\twoargroughtori{\timefunction,u}{\muxmulevelsetvalue}$.
	We define the canonical area form of $\twoargroughtori{\timefunction,u}{\muxmulevelsetvalue}$ 
	induced by $\gtorusroughfirstfund$
in the rough adapted coordinates $(\timefunction,u,x^2,x^3)$ by:
\begin{align} \label{E:AREAFORMROUGHTORUS}
	\volroughtorus 
	& = \volroughtorus(\timefunction,u',x^2,x^3) 
	\eqdef 
	\sqrt{\det \gtorusroughfirstfund(\timefunction,u',x^2,x^3)} \, \mathrm{d} x^2 \mathrm{d} x^3,
\end{align}
where $\det \gtorusroughfirstfund(\timefunction,u',x^2,x^3)$ is the determinant of the $2 \times 2$ matrix
$\left(\gtorusroughfirstfund(\timefunction,u',x^2,x^3) \left(\roughgeop{x^A},\roughgeop{x^B}\right) \right)_{A,B=2,3}$
(see RHS~\eqref{E:GTORUSFIRSTFUNDRESTRICTEDTOROUGHTORUS}).
\item We define the (non-canonical) area form $\volRoughHypersurface$ of $\hypthreearg{\timefunction}{I}{\muxmulevelsetvalue}$ 
in the rough adapted coordinates $(\timefunction,u,x^2,x^3)$ by:
\begin{align} \label{E:VOLUMEFORMROUGHHYPERSURFACE}
	\volRoughHypersurface 
	& = \volRoughHypersurface(\timefunction,u',x^2,x^3) 
	= \volroughtorus(\timefunction,u',x^2,x^3) \, \mathrm{d} u'.
\end{align}
\item We define the (non-canonical) area form $\volPuRoughCoordinates$ on $\nullhypthreearg{\muxmulevelsetvalue}{u}{J}$ 
in the rough adapted coordinates $(\timefunction,u,x^2,x^3)$ by:  
\begin{align} \label{E:VOLUMEFORMNULLHYPERSURFACEROUGHCOORDINATES}
\volPuRoughCoordinates 
= 
\volPuRoughCoordinates(\timefunction',u,x^2,x^3) 
& \eqdef \volroughtorus(\timefunction',u,x^2,x^3) \, \mathrm{d} \timefunction'.
\end{align} 
\item We define the (non-canonical) volume form $\volMRoughCoordinates$ of $\twoargMrough{I,J}{\muxmulevelsetvalue}$ 
in the rough adapted coordinates $(\timefunction,u,x^2,x^3)$ by:
\begin{align} \label{E:VOLUMEFORMSPACTEIMROUGHCOORDINATES}
	\volMRoughCoordinates 
	= \volMRoughCoordinates(\timefunction',u',x^2,x^3) 
	& 
	\eqdef 
	\volroughtorus(\timefunction',u',x^2,x^3) \, \mathrm{d} u' \, \mathrm{d} \timefunction'.
\end{align}
\end{itemize}

Unless we explicitly indicate otherwise, 
all integrals along 
$\twoargroughtori{\timefunction,u}{\muxmulevelsetvalue}$,
$\hypthreearg{\timefunction}{I}{\muxmulevelsetvalue}$,
and $\twoargMrough{I,J}{\muxmulevelsetvalue}$
are defined with respect to the above forms. 
Moreover, we will often suppress the variables with respect to which we integrate, 
e.g., we write:
\begin{subequations} 
\begin{align}
\int_{\twoargroughtori{\timefunction,u}{\muxmulevelsetvalue}} 
	f 
\, \volroughtorus 
& 
\eqdef 
\int_{(x^2,x^3) \in \T^2}  
	f(\timefunction,u,x^2,x^3) 
\, \volroughtorus(\timefunction,u,x^2,x^3),   
	\label{E:ROUGHTORUSINTEGRAL} \\
\int_{\nullhypthreearg{\muxmulevelsetvalue}{u}{J}} 
	f \, 
\volPuRoughCoordinates 
& 
\eqdef \int_{\timefunction' \in J}  
\int_{(x^2,x^3) \in \T^2} 
	f(\timefunction',u,x^2,x^3) \, \volroughtorus(\timefunction',u,x^2,x^3)
\, \mathrm{d} \timefunction', 
\label{E:NULLHYPINTEGRAL} \\
\int_{\hypthreearg{\timefunction}{I}{\muxmulevelsetvalue}} 
	f 
\, \volRoughHypersurface 
& \eqdef 
\int_{u' \in I} \int_{(x^2,x^3) \in \T^2} 
	f(\timefunction,u',x^2,x^3) 
\, \volroughtorus(\timefunction,u',x^2,x^3) \, \mathrm{d} u', 
	\label{E:ROUGHHYPERSURFACEINTEGRAL} \\
\int_{\twoargMrough{I,J}{\muxmulevelsetvalue}} 
	f 
\, \volMRoughCoordinates 
& 
\eqdef 
\int_{\timefunction' \in J} 
\int_{u' \in I} 
\int_{(x^2,x^3) \in \T^2} 
	f(\timefunction',u',x^2,x^3) 
\, \volroughtorus(\timefunction',u',x^2,x^3)
\, 
\mathrm{d} u' \, \mathrm{d}\timefunction'.
\label{E:SPACETIMEINTEGRAL}
\end{align}
\end{subequations}
\end{definition}

In a few of our calculations, we will also refer to the canonical volume forms of
$\hypg$ and $\gfour$ relative to the rough adapted coordinates, which we provide
in the following definition.

\begin{definition}[Canonical volume forms relative to the rough adapted coordinates] 
	\label{D:CANONICALVOLFORMSINROUGHADAPATEDCOORDINATES} 
	\hfill
\begin{itemize}
\item Recall that $\hypg$ is the first fundamental form of $\hypthreearg{\timefunction}{[- \rightu,\leftu]}{\muxmulevelsetvalue}$, 
		as in Def.\,\ref{D:ROUGHFIRSTFUNDS}.
		We define $\volcanonical_{\hypg} \eqdef \sqrt{\mbox{\upshape det} \hypg} \, \mathrm{d} x^2 \, \mathrm{d} x^3 \, \mathrm{d} u'$ 
		to be the canonical area form on $\hypthreearg{\timefunction}{[- \rightu,u]}{\muxmulevelsetvalue}$ induced by $\hypg$, 
		where $\mbox{\upshape det} \hypg$ is evaluated relative to the 
		rough adapted coordinates $(u',x^2,x^3)$ on $\hypthreearg{\timefunction}{[- \rightu,\leftu]}{\muxmulevelsetvalue}$.
\item Recall that $\gfour$ denotes the acoustical metric, defined relative to the Cartesian coordinates in \eqref{E:ACOUSTICALMETRIC}.
		We define 
		$\volcanonical_{\gfour} \eqdef \sqrt{|\det \gfour|}\, \mathrm{d} x^2 \, \mathrm{d} x^3\, \mathrm{d} u' \, \mathrm{d} \timefunction'$ 			to be the canonical area form on $\twoargMrough{[\timefunction_0,\timefunctionboot),[- \rightu,\leftu]}{\muxmulevelsetvalue}$ induced by 
			the acoustical metric $\gfour$, 
			where $\mbox{\upshape det} \gfour$ is evaluated relative to the rough adapted coordinates 
			$(\timefunction',u',x^2,x^3)$ on $\twoargMrough{[\timefunction_0,\timefunctionboot),[- \rightu,\leftu]}{\muxmulevelsetvalue}$.
\end{itemize}
\end{definition}

\subsubsection{Identities involving the forms}
\label{SSS:IDENTITIESINVOLVINGVOLUMEFORMSINROUGHADAPATEDCOORDINATES}
For future use, in the following lemma, we establish several identities involving 
$\volcanonical_{\hypg}, \volcanonical_{\gfour}$ 
and 
$\volRoughHypersurface, \volMRoughCoordinates$ relative to the rough adapted coordinates. 

\begin{lemma}[Identities involving $\volcanonical_{\hypg},\, \volcanonical_{\gfour}$ and $\volRoughHypersurface,\, \volMRoughCoordinates$]
\label{L:IDENTITIESINVOLVINGVOLUMEFORMSINROUGHADAPATEDCOORDINATES}
The following identities hold relative to the rough adapted coordinates $(\timefunction,u,x^2,x^3)$,
e.g., $\mbox{\upshape det} \gtorusroughfirstfund$
is the determinant of the $2 \times 2$ matrix
$\left(\gtorusroughfirstfund \left(\roughgeop{x^A},\roughgeop{x^B}\right) \right)_{A,B=2,3}$:
\begin{subequations}
\begin{align}
\mbox{\upshape det} \hypg 
& = 
|\Rtransarg{\muxmulevelsetvalue}|_{\gfour}^2 \det \gtorusroughfirstfund, \label{E:DETHYPGINROUGHADAPTEDCOORDS} 
	\\
\mbox{\upshape det} \gfour 
& = - \frac{\upmu^2}{(\Lunit \timefunctionarg{\muxmulevelsetvalue})^2} \mbox{\upshape det} \gtorusroughfirstfund. 
\label{E:DETACOUSTICALMETRICINROUGHADAPTEDCOORDS} 
\end{align}
\end{subequations}

Moreover, with $\volcanonical_{\hypg},
\volcanonical_{\gfour},
\volRoughHypersurface
$,
and
$\volroughtorus$ denoting the area and volume forms from
Defs.\,\ref{D:ROUGHVOLFORMS} and \ref{D:CANONICALVOLFORMSINROUGHADAPATEDCOORDINATES}, 
we have the following identities relative to the rough adapted coordinates:
\begin{subequations}
\begin{align}
\volcanonical_{\hypg} 
& = |\Rtransarg{\muxmulevelsetvalue}|_{\gfour} \volRoughHypersurface 
	= |\Rtransarg{\muxmulevelsetvalue}|_{\gfour} \volroughtorus \, \mathrm{d} u', 
	\label{E:VOLFORMCANONICALHYPGROUGHADAPTED} 
		\\
\volcanonical_{\gfour}
& = \frac{\upmu}{\Lunit \timefunctionarg{\muxmulevelsetvalue}} \volMRoughCoordinates 
= 
\frac{\upmu}{\Lunit \timefunctionarg{\muxmulevelsetvalue}} 
\volroughtorus
\, \mathrm{d} u' 
\, \mathrm{d} \timefunction. 
\label{E:VOLFORMACOUSTICALMETRICROUGHADAPTED}
\end{align}
\end{subequations}
\end{lemma}

\begin{proof}
Let $\Rtransunitarg{\muxmulevelsetvalue}$ be as in Def.\,\ref{D:GEOMETRICVECTORFIELDSADAPTEDTOROUGHFOLIATIONS}.
Recall (see Prop.\,\ref{P:BASICPROPERTIESOFROUGHVECTORFIELDS}) that
$\Rtransunitarg{\muxmulevelsetvalue}$ is tangent to $\hypthreearg{\timefunction}{[- \rightu,\leftu]}{\muxmulevelsetvalue}$ and
$\gfour$-orthogonal to the rough tori $\twoargroughtori{\timefunction,u}{\muxmulevelsetvalue}$. 
From 
\eqref{E:GEOP2TOCOMMUTATORS}--\eqref{E:GEOP3TOCOMMUTATORS},
\eqref{E:DEFROUGHTORITANGENTVECTORFIELD},
\eqref{E:RTRANS}, 
and Lemma~\ref{L:BASICPROPERTIESOFVECTORFIELDS}, it follows that
$\Rtransunitarg{\muxmulevelsetvalue} u = \frac{1}{|\Rtransarg{\muxmulevelsetvalue}|_{\gfour}}$. 
Also considering the identities 
\eqref{E:GTORUSINVERSEROUGHEXPRESSION}--\eqref{E:ROUGHTORUSMETRICCOMPONENTSANDTHEINVERSECOMPONENTSRELATION}
and
\eqref{E:INVERSEFIRSTFUNDOFROUGHHYPERSURFACESINTERMSOFR}
(where we view the components of RHS~\eqref{E:INVERSEFIRSTFUNDOFROUGHHYPERSURFACESINTERMSOFR} as entries
of a $3 \times 3$ matrix in rough adapted coordinates $(u,x^2,x^3)$
on $\hypthreearg{\timefunction}{[- \rightu,\leftu]}{\muxmulevelsetvalue}$),
we carry out straightforward calculations in the rough adapted coordinates 
to deduce that 
$\mbox{\upshape det} \hypginverse = |\Rtransarg{\muxmulevelsetvalue}|^{-2} \mbox{\upshape det} \gtorusroughinversefirstfund$,
where $\mbox{\upshape det} \gtorusroughinversefirstfund$
is the determinant of the $2 \times 2$ matrix
$\left(\gtorusroughinversefirstfund \left(d x^A,d x^B\right) \right)_{A,B=2,3}$
The desired result \eqref{E:DETHYPGINROUGHADAPTEDCOORDS} now readily follows. 
Similarly, recall that $\hypunitnormalarg{\muxmulevelsetvalue}$ denotes the $\gfour$-unit normal
to $\hypthreearg{\timefunction}{[- \rightu,\leftu]}{\muxmulevelsetvalue}$.
The identity \eqref{E:HYPUNITNORMALDECOMPOSITION} implies that 
$\hypunitnormalarg{\muxmulevelsetvalue} \timefunction = \frac{\Lunit \timefunctionarg{\muxmulevelsetvalue} |\Rtransarg{\muxmulevelsetvalue}|_{\gfour}}{\upmu}$ 
and $\hypunitnormalarg{\muxmulevelsetvalue} u = \frac{1}{|\Rtransarg{\muxmulevelsetvalue}|_{\gfour}}$.
Also considering the identities
\eqref{E:GTORUSINVERSEROUGHEXPRESSION}--\eqref{E:ROUGHTORUSMETRICCOMPONENTSANDTHEINVERSECOMPONENTSRELATION}
and
\eqref{E:INVERSEACOUSTICALMETRICINTEMRSOFNUNITANDRUNIT}
(where we view the components of RHS~\eqref{E:INVERSEACOUSTICALMETRICINTEMRSOFNUNITANDRUNIT} as entries
of a $4 \times 4$ matrix in rough adapted coordinates),
we carry out straightforward calculations in the rough adapted coordinates 
to deduce that
$\mbox{\upshape det} \gfour^{-1} = -\frac{(\Lunit \timefunctionarg{\muxmulevelsetvalue})^2}{\upmu^2} \det \gtorusroughinversefirstfund$,
from which \eqref{E:DETACOUSTICALMETRICINROUGHADAPTEDCOORDS} readily follows.

The identities \eqref{E:VOLFORMCANONICALHYPGROUGHADAPTED}--\eqref{E:VOLFORMACOUSTICALMETRICROUGHADAPTED} then follow from
\eqref{E:DETHYPGINROUGHADAPTEDCOORDS}--\eqref{E:DETACOUSTICALMETRICINROUGHADAPTEDCOORDS}
and Defs.\,\ref{D:ROUGHVOLFORMS} and \ref{D:CANONICALVOLFORMSINROUGHADAPATEDCOORDINATES}.

\end{proof}

\subsubsection{Geometric $L^2$ and $L^{\infty}$ norms}
\label{SSS:GEOMETRICL2NORMS}

\begin{definition}[Geometric $L^2$ norms] 
\label{D:GEOMETRICL2NORMS}
We define the $L^2$ norms with respect to the area and volume forms introduced 
in Def.\,\ref{D:ROUGHVOLFORMS}. 
Recall that we measure the norm of $\ell_{t,u}$-tangent tensorfields with $\gtorus$, i.e.,
if $\upxi$ is a type $\binom{0}{2}$ $\ell_{t,u}$-tangent tensorfield, then
$|\upxi|_{\gtorus}^2 \eqdef (\gtorus^{-1})^{\alpha \beta}(\gtorus^{-1})^{\sigma \delta} \upxi_{\alpha \sigma} \upxi_{\beta \delta}$. 
Then for scalar functions or $\ell_{t,u}$-tangent tensorfields $\upxi$, 
we define: 
\begin{subequations}
\begin{align}
\| \upxi \|_{L^2\left(\twoargroughtori{\timefunction,u}{\muxmulevelsetvalue}\right)}
& 
\eqdef 
\left( 
	\int_{\twoargroughtori{\timefunction,u}{\muxmulevelsetvalue}} |\upxi|_{\gtorus}^2 \, \volroughtorus
\right)^{1/2},
& 
\| \upxi \|_{L^2\left(\nullhypthreearg{\muxmulevelsetvalue}{u}{J}\right)} 
& \eqdef 
\left( 
	\int_{\nullhypthreearg{\muxmulevelsetvalue}{u}{J}} |\upxi|_{\gtorus}^2 \, \volPuRoughCoordinates 
\right)^{1/2}, 
	\label{E:GEOMETRICL2NORMSTORIANDNULLHYPERSURFACES} 
	\\
\| \upxi \|_{L^2\left(\hypthreearg{\timefunction}{I}{\muxmulevelsetvalue}\right)}
& 
\eqdef \left( 
	\int_{\hypthreearg{\timefunction}{I}{\muxmulevelsetvalue}} |\upxi|_{\gtorus}^2 
	\, \volRoughHypersurface 
\right)^{1/2}, 
& 
\| \upxi \|_{L^2\left(\twoargMrough{I,J}{\muxmulevelsetvalue}\right)} 
& \eqdef 
\left( 
	\int_{\twoargMrough{I,J}{\muxmulevelsetvalue}} |\upxi|_{\gtorus}^2 \, \volMRoughCoordinates 
\right)^{1/2}.
\label{E:GEOMETRICL2NORMSROUGHHYPERSURFACESANDSPACETIMEREGIONS}
\end{align}
\end{subequations}
\end{definition}

\begin{definition}[Geometric $L^{\infty}$ norms] 
\label{D:GEOLINFINITYROUGHTORUS}
 For scalar functions or $\ell_{t,u}$-tangent tensorfields $\upxi$, we define 
	the following $L^{\infty}$ norm on the rough tori $\twoargroughtori{\timefunction,u}{\muxmulevelsetvalue}$:
\begin{align} \label{E:GEOLINFINITYROUGHTORUS}
\| \upxi \|_{L^{\infty}\left(\twoargroughtori{\timefunction,u}{\muxmulevelsetvalue}\right)}
& 
\eqdef 
\mbox{ess sup}_{(x^2,x^3) \in \mathbb{T}^2}
|\upxi|_{\gtorus}(\timefunction,u,x^2,x^3),
\end{align}
where on RHS~\eqref{E:GEOLINFINITYROUGHTORUS}, we are viewing
$\upxi$ as a function of the rough adapted coordinates
$(\timefunction,u,x^2,x^3)$.
\end{definition}

\begin{remark}[Carefully note the role of $\gtorus$] \label{R:ROLEOFFIRSTFUNDOFSMOOTHTORI}
	We stress that $\gtorus$ is the Riemannian metric on the smooth tori,
	even though the integrals defining
	$\| \cdot \|_{L^2\left(\twoargroughtori{\timefunction,u}{\muxmulevelsetvalue}\right)}$,
	$\| \cdot \|_{L^2\left(\nullhypthreearg{\muxmulevelsetvalue}{u}{J}\right)}$,
	$\| \cdot \|_{L^2\left(\hypthreearg{\timefunction}{I}{\muxmulevelsetvalue}\right)}$,
	and
	$\| \cdot \|_{L^2\left(\twoargMrough{I,J}{\muxmulevelsetvalue}\right)}$ 
	are over regions and with respect to forms
	tied to the rough adapted coordinates.
	
	Similarly, on RHS~\eqref{E:GEOLINFINITYROUGHTORUS}, $|\upxi|_{\gtorus}$
	is the pointwise norm of $\upxi$ with respect to the 
	Riemannian metric $\gtorus$ on the smooth tori.
\end{remark}

\subsection{Strings of commutation vectorfields and vectorfield seminorms}
\label{SS:STRINGSOFCOMMUTATIONVECTORFIELDS}
To simplify the presentation of formulas and estimates, we will use
the notation from the following definition for strings of commutation vectorfields.

\begin{definition}[Strings of commutation vectorfields and vectorfield seminorms] 
\label{D:STRINGSOFCOMMUTATIONVECTORFIELDS}
Recall that 
$\Fullset = \{\Lunit, \muX, \Yvf{2},\Yvf{3}\}$, 
$\Tanset = \{\Lunit,\Yvf{2},\Yvf{3}\}$,
and 
$\Angularset = \{\Yvf{2},\Yvf{3}\}$ are the sets of commutation vectorfields 
defined in \eqref{E:COMMUTATIONVECTORFIELDS}.
We adopt the following definitions for differential operators constructed out of strings of commutation vectorfields.
In this section, $f$ denotes a scalar function.

\begin{itemize}
\item $\comder^{N;M} f$ denotes an arbitrary string of $N$ commutation vectorfields in $\comder$ applied to $f$, where the string contains \emph{at most} $M$ factors of $\muX$. We set $\comder^{0,0}f = f$. We often write
$\comder f$ instead of $\comder^{1;1} f$.
\item $\tander^N f$ denotes an arbitrary string of $N$ commutation vectorfields in $\tander$ applied to $f$. We set $\tander^0 f = f$. We often write $\tander f$ instead of $\tander^1 f$.
\item $\tanderY^N f$ denotes an arbitrary string of $N$ commutation vectorfields in $\mathscr{Y}$ applied to $f$. We set $\tanderY^0 f = f$.
\item $\comdersmall^{N;M} f$ denotes an arbitrary string of $N$ commutation vectorfields in $\comder$ applied to $f$, where the string contains \emph{at least} one factor of $\tander$ and \emph{at most} $M$ factors of $\muX$. \item $\tandersmall^N$ denotes an arbitrary string of $N$ commutation vectorfields in $\tander$ applied to $f$, where the string contains \emph{at least two factors} of $\Lunit$ or \emph{at least one factor of $\Yvf{2}, \Yvf{3}$}.
\item $\comderdoublesmall^{N;M} f$ denotes an arbitrary string of $N$ commutation vectorfields in $\comder$ applied to $f$, 
	where the string contains at least two factors of $\Lunit$ or at least one factor of $\Yvf{A}$ and \emph{at most} $M$ factors of $\muX$.
\item $\mathfrak{P}^{(N)}$ denotes the set of all differential operators of the form $\tander^N$.
\item $\mathfrak{Y}^{(N)}$ denotes the set of all differential operators of the form $\tanderY^N$.
\item  We define order $N$ strings of $\ell_{t,u}$-projected Lie derivatives
	such as $\angLie_{\tander}^N$ and $\angLie_{\comder}^{N;M}$ in an analogous fashion.  
	Such operators act on $\ell_{t,u}$-tangent tensorfields $\upxi$, e.g., $\angLie_{\tander}^N \upxi$.
\item $\angLie_{\mathfrak{P}}^{(N)}$ denotes the set of all differential operators of the form $\angLie_{\tander}^N$.
\item $\angLie_{\mathfrak{Y}}^{(N)}$ denotes the set of all differential operators of the form $\angLie_{\tanderY}^N$.
\item  $\comder^{\leq N;M} f$ denotes the array of all terms of the form
	$\comder^{N';M} f$, where $0 \leq N' \leq N$. 
\item If $N_1 < N_2$, then
	$\comder^{[N_1,N_2];M} f$ denotes the array of all terms of the form
	$\comder^{N';M} f$, where $N_1 \leq N' \leq N_2$.
\item We define arrays such as $\tanderY^{\leq 2} f$, 
	$\angLie_{\tanderY}^{[N_1,N_2]} \upxi$, etc.\
	in an analogous fashion.
\end{itemize}
We also define corresponding pointwise seminorms:
\begin{itemize}
\item $|\comder^{N;M} f|$ denotes the magnitude of $\comder^{N;M} f$ as defined above (there is no summation).
\item $|\comder^{\leq N;M} f|$ denotes the sum over all terms of the form $|\comder^{N';M} f|$ with $N' \leq N$. 
\item $| \comder^{[N_1,N_2];M}|$ is the sum over all terms of the form $|\comder^{N';M} f|$ with $N_1 \leq N' \leq N_2$. 
\item Terms such as 
	$|\tandersmall^{[N_1,N_2]}f|$, 
	$|\angLie_{\comder}^{\leq N;M} \upxi|_{\gtorus}$, 
	$|\tanderY^{\leq N} f|$,
	etc., are defined analogously,
	e.g., $|\angLie_{\comder}^{\leq N;M} \upxi|_{\gtorus}$
	is the sum over all terms of the form
	$|\angLie_{\comder}^{N';M} \upxi|_{\gtorus}$
	with $N' \leq N$.
\end{itemize}

\begin{itemize}
	\item We will freely combine the above definitions with Def.\,\ref{D:CONVENTIONESFORDIFFERENTIATION}, e.g.,
		\begin{align} \label{E:NOTATIONFORSTRINGDERIVATIVESOFMULTIPLEARRAYS}
			|\tander^N (\vortrenormalized,\GradEnt)| 
			& \eqdef 
				\max
				\left\lbrace 
					|\tander^N \vortrenormalized|, 
						\,
					|\tander^N \GradEnt|
				\right\rbrace.
		\end{align}
\end{itemize}
\end{definition}

\section{Schematic structure capturing structure and schematic identities}
\label{S:SCHEMATICSTRUCTUREANDIDENTITIES}
In this section, we introduce schematic notation that will help us succinctly exhibit
the important qualitative features of various equations.
We then provide a collection of identities expressed in schematic form;
they will be helpful when we derive estimates.

\subsection{Some schematic notation}
\label{SS:SOMESCHEMATICANOTATION}
\begin{notation}[Schematic functional dependence] 
\label{N:SCHEMATICFUNCTIONALDEPENDENCE}
We often use the notation $\smoothfunction(\upxi_{(1)},\cdots,\upxi_{(m)})$ to schematically depict an expression 
(often tensorial and involving contractions) that depends smoothly on the $\ell_{t,u}$-tangent
tensorfields $\upxi_{(1)},\cdots,\upxi_{(m)}$. 
Note that in general, $\smoothfunction(0) \neq 0$. 
\end{notation}

\begin{notation}[Schematic use of the symbol $\tander$] 
\label{N:SCHEMATICUSEOFTANGENTIALDIFFERENTATION}
Throughout the rest of the paper, 
$\tander$ 
schematically denotes a differential operator that is tangent to 
the characteristics $\nullhyparg{u}$. 
For example, $\tander f$ might denote $\angrmD f, \Lunit f$, or $\Yvf{2}f$. 
We use such notation when the details of $\tander$ are unimportant. 
\end{notation}

We use the notation $\vec{x}$ to denote the array of spatial Cartesian coordinates,
i.e.,
\begin{align} \label{E:ARRAYOFOFCARTESIANSPATIALCOORDIANTES}
	\vec{x}
	& \eqdef (x^1, x^2, x^3).
\end{align}
We use the same conventions from Def.\,\ref{D:DERIVATIVEOFARRAY} for
differential operators acting on $\vec{x}$, e.g.,
\begin{align} \label{E:ANDGULARDIFFERENTIALOFCARTESIANSPATIALCOORDIANTES}
	\angrmD \vec{x}
	& \eqdef (\angrmD x^1, \angrmD x^2, \angrmD x^3).
\end{align}

Finally, we recall that $\controlvars$ and $\badcontrolvars$ are the arrays from Def.\,\ref{D:CONTROLVARS}.

\subsection{Schematic structure of various tensorfields}
\label{SS:SCHEMATICSTRUCTUREOFVARIOUSTENSORFIELDS}

\begin{lemma}[Schematic structure of various tensorfields] 
\label{L:SCHEMATICSTRUCTUREOFVARIOUSTENSORSINTERMSOFCONTROLVARS} 
The following schematic relations hold for scalar functions
($\alpha,\beta = 0,1,2,3$, $\iota = 0,1,2,3,4$):
\begin{subequations}
\begin{align}
\gfour_{\alpha \beta}, 
	\,
(\gfour^{-1})^{\alpha \beta}, 
	\,
\gtorus_{\alpha \beta}, 
	\,
(\gtorus^{-1})^{\alpha \beta}, 
	\,
G_{\alpha \beta}^\iota, 
	\,
\smoothtorusproject_\beta^{\alpha}, 
	\,
\Lunit^{\alpha}, 
	\,
X^{\alpha}, 
	\,
\uLunit^{\alpha}, 
	\,
\Yvf{2}^{\alpha}, 
	\,
\Yvf{3}^{\alpha}, 
	\,
\Speed
& 
= \smoothfunction(\controlvars),  
	\label{E:SCHEMATICSTRUCTUREOFMETRICETC} 
	\\
\left(\geop{t}\right)^{\alpha},
	\,
\left(\geop{x^2}\right)^{\alpha},
	\,
\left(\geop{x^3}\right)^{\alpha}
& 
= \smoothfunction(\controlvars),
\label{E:SCHEMATICSTRUCTUREOFCARTESIANCOMPONENTSOFTANGENTIALGEOMETRICCOORDINATEVECTORFIELDS}
	\\
\left(\geop{u}\right)^{\alpha}
& 
= \smoothfunction(\badcontrolvars),
	\label{E:SCHEMATICSTRUCTUREOFCARTESIANCOMPONENTSOFTRANSVERSALGEOMETRICCOORDINATEVECTORFIELD} 
		\\
G_{\Lunit \Lunit}^\iota, 
	\,
G_{\Lunit X}^\iota, 
	\,
G_{XX}^\iota 
& = \smoothfunction(\controlvars), 
	\label{E:SCHEMATICSTRUCTUREOFGFRAMESCALARS} \\
\Xsmall^{\alpha}, 
	\, 
\Speed - 1 
& 
= \smoothfunction(\controlvars)\controlvars,
	\label{E:SCHEMATICSTRUCTUREOFXSMALL} \\
\muX^{\alpha}
& = \smoothfunction(\badcontrolvars).
	\label{E:SCHEMATICSTRUCTUREOFMUX}
\end{align}
\end{subequations}

Moreover, we have the following schematic relations for $\ell_{t,u}$-tangent tensorfields:
\begin{subequations}
\begin{align}
\gtorus, 
	\,
\angG_{\Lunit}, 
	\,
\angG_X,
	\,
\angG 
& = \smoothfunction(\controlvars, \angrmD \vec{x}\,), 
	\label{E:SCHEMATICSTRUCTUREOFGFRAMEELLTUTENSORS} 
		\\
\Yvf{2}, 
	\,
\Yvf{3} 
& = \smoothfunction(\controlvars, \gtorus^{-1}, \angrmD \vec{x}\,), 
	\label{E:SCHEMTAICSTRUCTUREOFSMOOTHANGULARCOMMUTATORVECTORFIELDS} 
		\\
\upchi 
& = \smoothfunction(\controlvars,  \angrmD \vec{x}\,) \tander \controlvars, 
	\label{E:SCHEMATICSTRUCTUREOFNULLSECONDFUNDAMENTALFORM} 
		\\
\mytr_{\gtorus} \upchi 
& 
= 
\smoothfunction(\controlvars, \gtorus^{-1}, \angrmD \vec{x}\,) \tander  \controlvars, 
	\label{E:SCHEMATICSTRUCTUREOFTRACEOFNULLSECONDFUNDAMENTALFORM}  
		\\
\zetatan, 
	\,
\angktan 
& 
= \smoothfunction(\controlvars, \angrmD \vec{x}\,) \tander \wavearray, 
	\label{E:SCHEMATICSTRUCTUREOFGOODPARTSOFCONNECTIONCOEFFICIENTS} 
		\\
\zetatrans, 
	\,
\angktrans 
& = \smoothfunction(\controlvars, \angrmD \vec{x}) \muX \wavearray. \label{E:SCHEMATICSTRUCTUREOFTRANSVERSALDIFFERENTIATIONPARTSOFCONNECTIONCOEFFICIENTS}
\end{align}
\end{subequations}

\end{lemma}

\begin{proof}
After one accounts for the third dimension,
the same proofs as in \cite[Lemma 2.19]{jSgHjLwW2016} hold for 
\eqref{E:SCHEMATICSTRUCTUREOFMETRICETC}--\eqref{E:SCHEMATICSTRUCTUREOFTRANSVERSALDIFFERENTIATIONPARTSOFCONNECTIONCOEFFICIENTS},
except \eqref{E:SCHEMATICSTRUCTUREOFMETRICETC} for
$\uLunit^{\alpha}$ was not stated there,
\eqref{E:SCHEMATICSTRUCTUREOFCARTESIANCOMPONENTSOFTANGENTIALGEOMETRICCOORDINATEVECTORFIELDS}--\eqref{E:SCHEMATICSTRUCTUREOFCARTESIANCOMPONENTSOFTRANSVERSALGEOMETRICCOORDINATEVECTORFIELD}  
were not stated there, and
\eqref{E:SCHEMATICSTRUCTUREOFMETRICETC} for
$\Speed - 1$ were not stated there. 
The identity $\uLunit^{\alpha} = \smoothfunction(\controlvars)$ stated in
\eqref{E:SCHEMATICSTRUCTUREOFMETRICETC} follows from definition \eqref{E:ULUNIT}
and the identity \eqref{E:SCHEMATICSTRUCTUREOFMETRICETC} for $\Lunit^{\alpha}$ and $X^{\alpha}$.
The identity
$\Speed - 1 = \smoothfunction(\controlvars)\controlvars$ 
stated in \eqref{E:SCHEMATICSTRUCTUREOFMETRICETC}
follows easily from \eqref{E:BACKGROUNDSOUNDSPEEDISUNITY}.
The identities
\eqref{E:SCHEMATICSTRUCTUREOFCARTESIANCOMPONENTSOFTANGENTIALGEOMETRICCOORDINATEVECTORFIELDS}--\eqref{E:SCHEMATICSTRUCTUREOFCARTESIANCOMPONENTSOFTRANSVERSALGEOMETRICCOORDINATEVECTORFIELD} 
follow from 
\eqref{E:GEOPTOCOMMUTATORS}--\eqref{E:GEOP3TOCOMMUTATORS}
and the remaining identities in the lemma.
\end{proof}

\subsection{Transversal derivatives in terms of  $\nullhyparg{u}$-tangential derivatives and structural properties of $\gfour$-null forms}
\label{SS:TRANSVERSALDERIVATIVESINTERMSOFTANGENTIALONES}
In this section, we use the transport equations from Theorem~\ref{T:GEOMETRICWAVETRANSPORTSYSTEM} to
derive expressions for the $\muX$ derivatives of 
$\vortrenormalized$, $\GradEnt$, $\VortVort$, and $\DivGradEnt$ in terms of $\nullhyparg{u}$-tangential derivatives.
We also exhibit some crucial structural properties of the 
inhomogeneous terms in the equations of Theorem~\ref{T:GEOMETRICWAVETRANSPORTSYSTEM},
including the $\gfour$-null forms.

\begin{lemma}[Expressions for $\muX \vortrenormalized$ and $\muX \GradEnt$ in terms of $\nullhyparg{u}$-tangential derivatives]
\label{L:TRANSVERSALDERIVATIVESOFTRANSPORTVARIABLESINTERMSOFTRANGENTIAL}
The following schematic identities hold for the $\muX$ derivatives of the Cartesian components of
$\vortrenormalized$ and $\GradEnt$:
	\begin{subequations}
	\begin{align} \label{E:VORTICITYTRANSVERSALTRANSPORTINTERMSOFTANGENTIAL}
		\muX \vortrenormalized^i
		& = 
		- 
		\upmu \Lunit \vortrenormalized 
		+
		\smoothfunction(\badcontrolvars, \GradEnt, \comder \wavearray) 
		\cdot
		(\vortrenormalized, \GradEnt).
			\\
		\muX \GradEnt^i
		& = 
		- 
		\upmu \Lunit \GradEnt 
		+
		\smoothfunction(\badcontrolvars, \GradEnt, \comder \wavearray) 
		\cdot
		(\vortrenormalized, \GradEnt).
		\label{E:ENTROPYGRADIENTTRANSVERSALTRANSPORTINTERMSOFTANGENTIAL}
	\end{align}
	\end{subequations}
\end{lemma}

\begin{proof}
The identity \eqref{E:VORTICITYTRANSVERSALTRANSPORTINTERMSOFTANGENTIAL} 
follows from multiplying the transport equation \eqref{E:RENORMALIZEDVORTICTITYTRANSPORTEQUATION} and 
by $\upmu$, using the identity $\muX  = - \upmu \Lunit + \upmu \Transport$ (see \eqref{E:BISLPLUSX}),
using Lemma~\ref{L:RELATIONSHIPBETWEENCARTESIANPARTIALDERIVATIVESANDSMOOTHGEOMETRICCOMMUTATORS} 
to write the Cartesian partial derivatives $\partial_{\alpha}$ in terms of the commutation vectorfields,
and using Lemma~\ref{L:SCHEMATICSTRUCTUREOFVARIOUSTENSORSINTERMSOFCONTROLVARS}.
\eqref{E:ENTROPYGRADIENTTRANSVERSALTRANSPORTINTERMSOFTANGENTIAL}
follows from a similar argument based on \eqref{E:GRADENTROPYTRANSPORT}
\end{proof}

\begin{lemma}[Crucial structural properties of $\gfour$-null forms]
	\label{L:CRUCIALSTRUCTUREOFNULLFORMS}
	The product of $\upmu$ and the terms defined in \eqref{E:TRANSPORTVORTVORTMAINTERMS}--\eqref{E:DIVENTROPYGRADIENTNULLFORM} 
	enjoy the following schematic structure:\footnote{All of these are $\gfour$-null forms, except the last term on
 RHS~\eqref{E:TRANSPORTDIVGRADENTMAINTERMS}, which turns out to be a harmless error term.}
\begin{subequations}
\begin{align} 
\upmu \mainnullform_{(\VortVort)}^i, 
	\,
\upmu \mainnullform_{(\DivGradEnt)} 
& = 
	\smoothfunction\left(\badcontrolvars, \GradEnt, \comder \wavearray \right) 
	\cdot
	(\tander^{\leq 1} \vortrenormalized, \tander^{\leq 1} \GradEnt),
		\label{E:NULLFORMSTRUCTUREMODIFIEDFLUIDVARIABLES}
		\\
\upmu \nullform_{(v)}^i,
	\, 
\upmu \nullform_{(\pm)}, 
	\, 
\upmu \nullform_{(\LogDensity)}
& 
= \smoothfunction(\badcontrolvars, \comder \wavearray) 
	\cdot
	\tander \wavearray,
		\label{E:NULLFORMSTRUCTUREWAVEVARIABLES} 
		\\
\upmu \nullform_{(\VortVort)}^i, 
	\, 
\upmu \nullform_{(\DivGradEnt)}
& 
= \smoothfunction(\badcontrolvars, \GradEnt, \comder \wavearray) 
	\cdot
	\GradEnt.
\label{E:NULLFORMSTRUCTURETRANSPORTVARIABLES}
\end{align}
\end{subequations}
\end{lemma}

\begin{proof}
	 All the results follow from \cite{jLjS2021}*{Lemma~8.2}, except the term stemming from last term on
	 RHS~\eqref{E:TRANSPORTDIVGRADENTMAINTERMS}, which is of the schematic form 
	 $\smoothfunction(\wavearray, \GradEnt) \cdot \upmu \partial_i \vortrenormalized$, 
	was not handled there. To handle this last term,
	we use Lemma~\ref{L:RELATIONSHIPBETWEENCARTESIANPARTIALDERIVATIVESANDSMOOTHGEOMETRICCOMMUTATORS} 
	to write the Cartesian spatial partial derivatives $\partial_i$ in terms of the commutation vectorfields,
	use the identity \eqref{E:VORTICITYTRANSVERSALTRANSPORTINTERMSOFTANGENTIAL}
	to substitute for the $\muX$ derivatives of $\vortrenormalized$,
	and use Lemma~\ref{L:SCHEMATICSTRUCTUREOFVARIOUSTENSORSINTERMSOFCONTROLVARS}.
\end{proof}

\begin{lemma}[Crucial structural properties of the linear inhomogeneous terms]
	\label{L:CRUCIALSTRUCTURELINEARINHOMOGENEOUS}
	The product of $\upmu$ and the terms $\VortVort, \DivGradEnt$-involving terms on
	RHSs \eqref{E:VELOCITYWAVEEQUATION}--\eqref{E:ENTROPYWAVEEQUATION},
	as well as the product of $\upmu$ and the terms defined in \eqref{E:VELOCITYILINEARORBETTER}--\eqref{E:RENORMALIZEDVORTICITYCURLLINEARORBETTER}
	enjoy the following schematic structure:
\begin{subequations}
\begin{align} 
\begin{split}  \label{E:STRUCTUREOFLINEARMODIFIEDFLUIDVARIABLETERMSONRHSWAVEEQUATIONS} 
	&
	\upmu \Speed^2 \exp(2 \LogDensity) \VortVort^i,
		\,
	\upmu
	\left\lbrace 
		F_{;\Ent} c^2 \exp(2\LogDensity)
		- 
		\Speed \exp(\LogDensity) \frac{p_{;\Ent}}{\overline{\varrho}} 
	\right\rbrace  
	\DivGradEnt,
		\\
	&
	\upmu
	\exp(\LogDensity) \frac{p_{;\Ent}}{\overline{\varrho}} \DivGradEnt,
		\,
	\upmu
	\Speed^2 \exp(2 \LogDensity) \DivGradEnt
		\\
	& = \upmu 
			\smoothfunction(\wavearray)
			\cdot
			(\VortVort,\DivGradEnt),
	\end{split}
		\\
	\begin{split} \label{E:ALLLINEARTERMSTRUCTURE} 
	&
	\upmu \mathfrak{L}_{(v)}^i, 
		\,
	\upmu \mathfrak{L}_{(\pm)},
		\,
	\upmu \mathfrak{L}_{(\LogDensity)},
		\,
	\upmu \mathfrak{L}_{(\Ent)},
		\,
	\upmu \mathfrak{L}_{(\vortrenormalized)}^i,
		\,
	\upmu \mathfrak{L}_{(\GradEnt)}^i,
		\,
	\upmu \mathfrak{L}_{(\Flatdiv \vortrenormalized)},
		\,
	\upmu \mathfrak{L}_{(\VortVort)}^i
		\\
	& 
	= 
	\smoothfunction(\badcontrolvars,\vortrenormalized,\GradEnt,\comder \wavearray) 
	\cdot
	(\vortrenormalized,\GradEnt).
\end{split}
\end{align}
\end{subequations}
		
\end{lemma}

\begin{proof}
	We use Lemma~\ref{L:RELATIONSHIPBETWEENCARTESIANPARTIALDERIVATIVESANDSMOOTHGEOMETRICCOMMUTATORS} 
	to write the Cartesian partial derivatives $\partial_{\alpha}$ in terms of the commutation vectorfields,
	and we use Lemma~\ref{L:SCHEMATICSTRUCTUREOFVARIOUSTENSORSINTERMSOFCONTROLVARS}.
\end{proof}

In the next lemma, we provide an analog of Lemma~\ref{L:TRANSVERSALDERIVATIVESOFTRANSPORTVARIABLESINTERMSOFTRANGENTIAL}
for the modified fluid variables.

\begin{lemma}[Expressions for the transversal derivatives of the modified fluid variables in terms of $\nullhyparg{u}$-tangential derivatives]
\label{L:TRANSVERSALDERIVATIVESOFMODIFIEDFLUIDVARIABLESINTERMSOFTRANGENTIAL}
Recall that the modified fluid variables $\VortVort$ and $\DivGradEnt$ are defined in Def.\,\ref{D:HIGHERORDERFLUIDVARIABLES}.
The following schematic identities hold for the $\muX$ derivatives of the Cartesian components of
$\VortVort$ and $\DivGradEnt$:
	\begin{align} \label{E:TRANSVERSALMODIFIEDINTERMSOFTANGENTIAL}
		(\muX \VortVort^i, \muX \DivGradEnt) 
		& 
		= 
		- 
		\upmu (\Lunit \VortVort,\Lunit \DivGradEnt) 
		+ 
		\smoothfunction(\badcontrolvars, \GradEnt, \comder \wavearray) 
		\cdot
		\tander \wavearray
		+ 
		\smoothfunction(\badcontrolvars, \GradEnt, \comder \wavearray) 
		\cdot
		(\tander^{\leq 1} \vortrenormalized, \tander^{\leq 1} \GradEnt).
	\end{align}

\end{lemma}

\begin{proof}
The identity \eqref{E:TRANSVERSALMODIFIEDINTERMSOFTANGENTIAL} follows from
an argument similar to the one we used to prove Lemma~\ref{L:TRANSVERSALDERIVATIVESOFTRANSPORTVARIABLESINTERMSOFTRANGENTIAL}, 
based on equations
\eqref{E:EVOLUTIONEQUATIONFLATCURLRENORMALIZEDVORTICITY} 
and
 \eqref{E:TRANSPORTFLATDIVGRADENT},
where we use
Lemma~\ref{L:CRUCIALSTRUCTUREOFNULLFORMS}
to handle the $\gfour$-null forms appearing on the RHSs of 
\eqref{E:EVOLUTIONEQUATIONFLATCURLRENORMALIZEDVORTICITY} 
and
\eqref{E:TRANSPORTFLATDIVGRADENT},
and we use 
\eqref{E:VORTICITYTRANSVERSALTRANSPORTINTERMSOFTANGENTIAL}--\eqref{E:ENTROPYGRADIENTTRANSVERSALTRANSPORTINTERMSOFTANGENTIAL} 
to substitute for 
$(\muX \vortrenormalized^i, \muX \GradEnt^i)$ whenever such terms arise.
\end{proof}

\begin{lemma}[Schematic identity for $\upmu \VortVort^i$ and $\upmu \DivGradEnt$]
	\label{L:SCHEMATICIDENTITYFORMUTIMESMODIFIEDFLUID}
	The $\upmu$-weighted Cartesian components of the modified fluid variables
	defined in \eqref{E:MODIFIEDCURLOFVORTICITY}--\eqref{E:MODIFIEDDIVERGENCEOFENTROPYGRADIENT}
	can be expressed as follows:
	\begin{align} \label{E:MUMODIFIEDFLUIDVARIABLESSCHEMATIC}
	(\upmu \VortVort^i, \upmu \DivGradEnt)
	& = 
			\smoothfunction(\badcontrolvars, \GradEnt, \comder \wavearray)
			\cdot
			\tander^{\le 1}(\vortrenormalized,\GradEnt).
	\end{align}
\end{lemma}

\begin{proof}
The identity \eqref{E:MUMODIFIEDFLUIDVARIABLESSCHEMATIC}
follows from definitions \eqref{E:MODIFIEDCURLOFVORTICITY}--\eqref{E:MODIFIEDDIVERGENCEOFENTROPYGRADIENT},
Lemma~\ref{L:SCHEMATICSTRUCTUREOFVARIOUSTENSORSINTERMSOFCONTROLVARS},
Lemma~\ref{L:RELATIONSHIPBETWEENCARTESIANPARTIALDERIVATIVESANDSMOOTHGEOMETRICCOMMUTATORS}
(which allows us to schematically express $\upmu \partial_{\alpha} 
= 
\smoothfunction(\controlvars) \muX + \upmu \smoothfunction(\controlvars) \tander$),
and
Lemma~\ref{L:TRANSVERSALDERIVATIVESOFTRANSPORTVARIABLESINTERMSOFTRANGENTIAL}
(which allows us to substitute for the $\muX$ derivatives of
$\vortrenormalized$ and $\GradEnt$).
\end{proof}

\subsection{Additional schematic identities involving differentiation}
\label{SS:ADDITIONALSCHEMATICIDENTITIESINVOLVINGDIFFERENTIATION}
For future use, in this section, we provide some additional schematic identities involving differentiation.

\begin{lemma}[Schematic identity for $\angLap \varphi$]
	\label{L:SCHEMATICEXPRESSIONFORANGULARLAPLACIAN}
	If $\varphi$ is a scalar function, then its angular Laplacian on $\ell_{t,u}$ can schematically
	be expressed as follows:
	\begin{align} \label{E:SCHEMATICEXPRESSIONFORANGULARLAPLACIAN}
		\angLap \varphi
		& = \smoothfunction(\controlvars) 
				\cdot
				\tanderY^2 \varphi
				+
				 \smoothfunction(\controlvars) 
				\cdot
				\tanderY \controlvars
				\cdot
				\tanderY \varphi.
	\end{align}
\end{lemma}

\begin{proof}
	Relative to the coordinates $(x^2,x^3)$ on $\ell_{t,u}$, we have
	$\angLap \varphi =
	\frac{1}{\sqrt{\mbox{\upshape det $\gtorus$}}} \geop{x^A} 
	\left\lbrace
		\sqrt{\mbox{\upshape det}\gtorus} 
		(\gtorus^{-1})^{AB}
		\geop{x^B} \varphi
	\right\rbrace
	$.
	From this identity, 
	\eqref{E:SMOOTHGINVERSEABEXPRESSION}--\eqref{E:DETERMSMOOTHGTORUSRELTOGEOMETRICCOORDS},
	\eqref{E:GEOP2TOCOMMUTATORS}--\eqref{E:GEOP3TOCOMMUTATORS},
	and
	Lemma~\ref{L:SCHEMATICSTRUCTUREOFVARIOUSTENSORSINTERMSOFCONTROLVARS},
	we conclude \eqref{E:SCHEMATICEXPRESSIONFORANGULARLAPLACIAN}.
\end{proof}

\begin{lemma}[Identity satisfied by $\muX \Lunit^i$] 
\label{L:MUXLISCHEMATICIDENTITY}
There exist smooth functions,
all schematically denoted by $\smoothfunction$,
such that the following identity holds:
\begin{align} \label{E:MUXLISCHEMATICIDENTITY}
	\muX \Lunit^i 
	& 
	= 
	\smoothfunction(\controlvars) \cdot \muX \wavearray \cdot(-\updelta_1^i + \Xsmall^i)
	+
	\smoothfunction(\controlvars) \cdot \muX \wavearraypartial
	+ 
	\upmu 
	\smoothfunction(\controlvars) \tander \wavearray
	+ 
	\smoothfunction(\controlvars) \tanderY \upmu.
\end{align}
\end{lemma}

\begin{proof}
	The same proof of \cite[Lemma 2.14]{jSgHjLwW2016} holds with minor modifications 
	that take into account 
	the expressions for $\gfour_{\alpha \beta}$ and $(\gfour^{-1})^{\alpha \beta}$
	given by
	\eqref{E:ACOUSTICALMETRIC} 
	and
	\eqref{E:INVERSEACOUSTICALMETRIC},
	the identity \eqref{E:BISLPLUSX},
	and the definition \eqref{E:XSMALL} of $\Xsmall^i$.
\end{proof}

\subsection{Deformation tensors}
\label{SS:DEFORMATIONTENSORS}
In our analysis, we encounter the deformation tensors of various vectorfields.

\begin{definition}[Deformation tensors] 
Let $Z$ be a spacetime vectorfield.
We define the \emph{deformation tensor} $\deform{Z}$ 
of $Z$ to be the following symmetric type $\binom{0}{2}$ tensorfield:
\begin{align} \label{E:DEFORMATIONTENSORDEF}
\deformarg{Z}{\alpha}{\beta} 
& 
\eqdef \Lie_Z \gfour_{\alpha \beta} 
= 
\Dfour_\alpha Z_\beta + \Dfour_\alpha Z_\beta,
\end{align}
where the final equality in \eqref{E:DEFORMATIONTENSORDEF} follows from the torsion-free property of $\Dfour$.
\end{definition}

\subsection{Rough-toroidal components of deformation tensors}
\label{S:ROUGHTORUSCOMPONENTSOFDEFORMATIONTENSORS}
The following lemma provides simple relationship between any deformation tensor 
$\deform{Z}$ contracted against $\roughgeop{x^A},\roughgeop{x^B}$ and $\Lie_Z \gtorusroughfirstfund$.

\begin{lemma}[Relating the $\twoargroughtori{\timefunction,u}{\muxmulevelsetvalue}$ components of $\Lie_Z \gtorusroughfirstfund$ and $\deform{Z}$]
\label{L:LIEGTORUSROUGH}
Let $Z$ be a spacetime vectorfield. 
Then the following identities hold for $A,B=2,3$:
\begin{align} \label{E:LIEGTORUSROUGH}
\deform{Z}\left(\roughgeop{x^A},\roughgeop{x^B}\right) 
& \eqdef
[\Lie_Z \gfour]\left(\roughgeop{x^A},\roughgeop{x^B}\right) 
=
[\Lie_Z \gtorusroughfirstfund]\left(\roughgeop{x^A},\roughgeop{x^B}\right).
\end{align}
\end{lemma}

\begin{proof}
Using \eqref{E:INVERSEACOUSTICALMETRICINTEMRSOFNUNITANDRUNIT},
we find that relative to arbitrary coordinates, we have
$\deformarg{Z}{\alpha}{\beta} 
\eqdef 
\Lie_{Z} \gfour_{\alpha \beta} = 
[\Lie_{Z}(-\hypunitnormalarg{\muxmulevelsetvalue}_{\flat} \otimes \hypunitnormalarg{\muxmulevelsetvalue}_{\flat}  
+
\Rtransunitarg{\muxmulevelsetvalue}_{\flat} \otimes \Rtransunitarg{\muxmulevelsetvalue}_{\flat} 
+ 
\gtorusroughfirstfund_{\alpha \beta}$,
where $\hypunitnormalarg{\muxmulevelsetvalue}_{\flat}$ denotes the one-form $\gfour$-dual to $\hypunitnormalarg{\muxmulevelsetvalue}$
and $\Rtransunitarg{\muxmulevelsetvalue}_{\flat}$ denotes the one-form $\gfour$-dual to $\Rtransunitarg{\muxmulevelsetvalue}$.
Since $\hypunitnormalarg{\muxmulevelsetvalue}$ and $\Rtransunitarg{\muxmulevelsetvalue}$ are $\gfour$-orthogonal to 
$\left\lbrace \roughgeop{x^2}, \roughgeop{x^3} \right\rbrace$, 
we conclude \eqref{E:LIEGTORUSROUGH}.
\end{proof}

\subsection{A schematic rewriting of the wave equations satisfied by $\wavearray$}
\label{SS:SCHEMATICREWRITINGOFGEOMETRICWAVEEQUATIONS}
The following lemma shows that $\muX \Psi$ obeys a transport equation with
source terms that are small but lose one derivative.
We will use it in Sect.\,\ref{S:LINFINITYFLUIDANDEIKONALANDIMPROVEMENTOFAUX}, 
when we derive improvements of the auxiliary bootstrap assumptions.

\begin{lemma}[A schematic rewriting of the wave equations satisfied by $\wavearray$]
\label{SS:SCHEMATICREWRITINGOFGEOMETRICWAVEEQUATIONS}
The covariant wave equations 
\eqref{E:COVARIANTWAVEEQUATIONSWAVEVARIABLES}
verified by $\Psi \in \{ \RRiemann,\LRiemann,v^2,v^3,\Ent\}$ 
can be expressed in the following schematic form:
\begin{align} \label{E:SCHEMATICREWRITINGOFWAVEEQUATIONSSATISFIEDBYWAVEVARIABLES}
\Lunit \muX \Psi 
& = 
\smoothfunction(\badcontrolvars) \tander^2 \wavearray 
+ 
\smoothfunction(\badcontrolvars,\comder \wavearray) \tander \controlvars  
+ 
\smoothfunction(\badcontrolvars, \GradEnt, \comder \wavearray)
\cdot
\tander^{\le 1}(\vortrenormalized,\GradEnt).
\end{align}

\end{lemma}

\begin{proof}
	We first decompose $\upmu \times \mbox{LHS~\eqref{E:COVARIANTWAVEEQUATIONSWAVEVARIABLES}}$
	using \eqref{E:BOXDECOMPLOUTSIDE}, 
	Lemma~\ref{L:SCHEMATICSTRUCTUREOFVARIOUSTENSORSINTERMSOFCONTROLVARS},
	and \eqref{E:SCHEMATICEXPRESSIONFORANGULARLAPLACIAN}.
	We then decompose $\upmu \times \mbox{RHS~\eqref{E:COVARIANTWAVEEQUATIONSWAVEVARIABLES}}$
	using Lemma~\ref{L:SCHEMATICSTRUCTUREOFVARIOUSTENSORSINTERMSOFCONTROLVARS},
	 \eqref{E:NULLFORMSTRUCTUREWAVEVARIABLES},
	and
	\eqref{E:MUMODIFIEDFLUIDVARIABLESSCHEMATIC}.
\end{proof}

\section{Parameters, their size assumptions, and conventions for constants}
\label{S:PARAMETERSANDSIZEASSUMPTIONSANDCONVENTIONSFORCONSTANTS}
In this section, we list and describe the ``size-parameters'' that appear throughout the paper.
These parameters will play a crucial role in Sect.\,\ref{S:ASSUMPTIONSONTHEDATA},
when we describe our assumptions on the data.
Then, in Sect.\,\ref{SS:CONVENTIONSFORCONSTANTS}, 
we state our conventions for how constants 
appearing in our analysis,
such as $C$ and $C_{\mydiam}$,
are allowed to depend on the parameters.

\subsection{Parameters} 
\label{SS:LISTOFPARAMETERS}

\subsubsection{Parameters of the background simple isentropic plane-symmetric solutions}
\label{SSS:PARAMETERSOFSIMPLEISENTROPICPLANESYMMMETRICSOLUTIONS}
First, we recall that in Appendix~\ref{A:PS},
we construct a large family of ``admissible''
simple isentropic plane-symmetric solutions,
where ``admissible'' means that it has properties such that
it falls under the scope of our main results; 
see Def.\,\ref{AD:ADMISSIBLEBACKGROUND} for the precise definition.
Each such admissible ``background'' solution has singularity-forming behavior that is described by 
the following (background solution-dependent) \textbf{positive} parameters, 
which we describe in detail in Appendix~\ref{A:PS}:
$
\farrightu,
	\,
\rightu, 
	\,
\interestingu, 
	\, 
\leftu, 
	\,
\timefunction_0,
	\,
\muxmulevelsetvalue_0, 
	\, 
\mupositive^{\text{PS}},
	\, 
\boringregionmupositive^{\text{PS}},
	\, 
\blowupdeltaPS,
	\, 
\transversalsizedeltaPS,
	\, 
\PSdataalpha,	 
$
and 
$\PSmutransversalHessiansize$.

\subsubsection{Parameters of the perturbed solutions}
\label{SSS:PARAMETERSOFPERTURBEDSOLUTIONS}
The positive parameters from Sect.\,\ref{SSS:PARAMETERSOFSIMPLEISENTROPICPLANESYMMMETRICSOLUTIONS} 
capture the behavior of various aspects of the background solution near its singular boundary.
In Appendix~\ref{A:OPENSETOFDATAEXISTS}, 
we use Cauchy stability arguments to show that there are open sets of initial data
on $\Sigma_0$, which are close to the data of one 
of the simple isentropic plane-symmetric solutions,
such that the state of the perturbed solution near its singular boundary 
(but still within the region of classical existence)
is described by the following parameters:
$
\farrightu,
	\,
\rightu, 
	\,
\interestingu, 
	\, 
\leftu, 
	\,
\timefunction_0,
	\,
\muxmulevelsetvalue_0, 
	\, 
\mupositive,
	\, 
\boringregionmupositive,
	\, 
\mathring{\updelta},
	\, 
\mathring{\updelta}_*,
	\, 
\mathring{\upalpha},	
	\, 
\secondtransversalderivativemulowerbound
$,
as well as our main new \emph{smallness parameter}: $\initialsmall$.
In Appendix~\ref{A:OPENSETOFDATAEXISTS}, we show that the parameters 
$
\farrightu,
	\,
\rightu, 
	\,
\interestingu, 
	\, 
\leftu,  
	\,
\timefunction_0,
	\,
\muxmulevelsetvalue_0,
	\,
\mupositive
$
for the perturbed solutions
can be chosen to be exactly the same as the ones for the background solutions. 
On the other hand, we will show that the remaining perturbed parameters
can be chosen to be close to the background ones in the following sense:
\begin{align} \label{E:PERTURBEDPARAMETERSARECLOSETOBACKGROUNDONES}
	\frac{1}{2}
	& \leq
	\frac{\boringregionmupositive}{\boringregionmupositive^{\text{PS}}},
		\,
	\frac{\mathring{\updelta}_*}{\transversalsizedeltaPS},
		\,
	\frac{\mathring{\updelta}}{\blowupdeltaPS},
		\,
	\frac{\mathring{\upalpha}}{\PSdataalpha},
		\,
	\frac{\secondtransversalderivativemulowerbound}{\secondtransversalderivativemulowerboundPS}
	\leq 2,
		\\
\initialsmall
& \geq 0
\mbox{\ can be chosen as small as we want.}
\label{E:DATAEPSILONCANBECHOSENASSMALLASWEWANT}
\end{align}

\subsubsection{Informal description of the parameters}
\label{SSS:INFORMALDESCRIPTIONOFTHEPARAMETERS}
Many of the parameters listed in Sect.\,\ref{SSS:PARAMETERSOFPERTURBEDSOLUTIONS}
will not appear until later in the paper, 
but to help guide the reader,
we now provide a summary of their role
and that of a few other parameters too.
We refer to Fig.\,\ref{F:REGIONSWHEREWEDERIVEESTIMATES} 
for an illustration that shows how some of the parameters
are tied to the location of various subsets.
\begin{itemize}
		\item The background density $\overline{\varrho} > 0$ is fixed throughout the article; see \eqref{E:BACKGROUNDDENSITY}.
		\item $\Ntop$ is an integer representing the maximum number of times we commute the equations
			when we derive energy estimates. In proving our main results, we assume that $\Ntop \geq 24$.
		\item The parameters $\timefunction_0$ and $\mupositive$ are related by $\timefunction_0 = - \mupositive$ 
			(see Def.\,\ref{D:RELATIONBETWEENINITIALROUGHTIMESLICEANDSMALLMUVALUE}).
			We view $\timefunction_0$ to be the ``initial rough time,'' i.e., the value of 
			$\timefunctionarg{\muxmulevelsetvalue}$ corresponding to the initial state of the solution near the singularity.
		\item The parameter $\mupositive$ is the minimum value of $\upmu$ along the initial rough hypersurface portion
			$\hypthreearg{\timefunction_0}{[-\interestingu,\interestingu]}{\muxmulevelsetvalue}$.
			For convenience, we will choose $\mupositive$ to be small. While the smallness of
			 $\mupositive$ is not essential, it allows us to
			focus on studying the solution only near the singularity and allows us to give short proofs of various estimates.
		\item The parameters $0 < \rightu$ and
			$0 < \interestingu < \leftu$ define the range of values for the eikonal function $u$ in the problem under study. 
			We will study the solution on various intervals of $u$ values, including:
			$[- \rightu, -\interestingu]$, 
			$[-\interestingu, \interestingu]$, 
			$[\interestingu,\leftu]$,
			and
			$[- \rightu,\leftu]$.
			The interesting analysis will take place in $[-\interestingu,\interestingu]$.
		\item The parameter $\mathring{\upalpha}$ measures the $L^{\infty}$-size of the amplitude of 
			$\RRiemann$ along the initial rough hypersurface $\hypthreearg{\timefunction_0}{[- \rightu,\leftu]}{\muxmulevelsetvalue}$.
		\item  The parameter $\boringregionmupositive$ quantifies the positivity of $\upmu$ away from the interesting region, 
			more precisely when $u \notin [-\interestingu,\interestingu]$;
			see \eqref{E:DATAMUISLARGEINBORINGREGION}.
		\item  The parameter $\mathring{\updelta}_*$ measures the size of 
			the crucial factor that drives the blowup; see definition \eqref{E:DELTASTARDEF}. 
		\item  The parameter $\mathring{\updelta}$ measures the $L^{\infty}$-size of the transversal derivatives of $\RRiemann$ along 
			the initial rough hypersurface $\hypthreearg{\timefunction_0}{[- \rightu,\leftu]}{\muxmulevelsetvalue}$. 
			It also controls the $L^{\infty}$ size of the transversal derivatives
			of various geometric quantities constructed out of the eikonal function.
			We make no smallness assumptions on $\mathring{\updelta}$. 
	\item The parameter $\initialsmall$ measures the extent to which the solution's data
		``break the simple isentropic plane-symmetry.'' 
	\item $\varepsilon$ is a small ``bootstrap parameter'' first appearing in 
		Sects.\,\ref{SS:MAINQUANTITATIVEBOOTSTRAPASSUMPTIONS}--\ref{SSS:AUXBOOTSTRAP}.
	\item The parameter $\secondtransversalderivativemulowerbound$ 
		quantifies the transversal convexity of $\upmu$ 
		(namely, the size of the second order $\nullhyparg{u}$-transversal derivatives 
		$\muX \muX \upmu$, $\Wtransarg{\muxmulevelsetvalue} \Wtransarg{\muxmulevelsetvalue} \upmu$, etc.) 
		in the interesting region $\hypthreearg{\timefunction_0}{[-\interestingu,\interestingu]}{\muxmulevelsetvalue}$; 
		see, for example, \eqref{E:DATATASSUMPTIONMUTRANSVERSALCONVEXITY}.
	\item The parameter $\muxmulevelsetvalue_0$ is such that the
		transversal convexity mentioned above holds
		for $\muxmulevelsetvalue \in [0,\muxmulevelsetvalue_0]$.
	\item The parameter $\farrightu$ 
		is defined by $\farrightu \eqdef \rightu + \frac{18}{\blowupdeltaPS}$
		(see \eqref{AE:FARRIGHTU})
		and plays a role only in Appendices~\ref{A:PS} and \ref{A:OPENSETOFDATAEXISTS}
		and the proof of Lemma~\ref{L:CONTROLOFDATAFORROUGHTORIENERGYESTIMATES},
		where we show that there exist open sets of initial data satisfying 
		all of our assumptions;
		see Sect.\,\ref{SS:PARAMETERSIZEASSUMPTIONS}.
\end{itemize}

\subsection{Parameter size assumptions}
\label{SS:PARAMETERSIZEASSUMPTIONS}
In this section, we state the size assumptions on the parameters that are sufficient
for our main results to hold, i.e., 
for Theorems~\ref{T:EXISTENCEUPTOTHESINGULARBOUNDARYATFIXEDKAPPA} 
and \ref{T:DEVELOPMENTANDSTRUCTUREOFSINGULARBOUNDARY} to hold.

\begin{quote}
For the remainder of the article, 
when we say that ``$A$ is small relative to $B$,''
we mean that $A \geq 0$, that $B > 0$, and that
there exists a continuous increasing function\footnote{Although we do not specify their form, 
the functions $f$ could always be chosen to be 
polynomials with positive coefficients or exponentials of such polynomials.}
$f :(0,\infty) \rightarrow (0,\infty)$ 
such that 
$
A < f(B)
$.
The functions $f$ are allowed to depend on the equation of state.
\end{quote}

\begin{assumption}[Size assumptions on the parameters]
	\label{ASSUMPTIONS:RELATIVESIZEASSUMPTIONS} \hfill
\begin{itemize}
	\item To close our estimates, we assume that the regularity-parameter $\Ntop$ is an integer satisfying:
		\begin{align} \label{E:NTOPLARGENESSASSUMPTION}
			\Ntop \geq 24.
		\end{align}
	\item We assume that the following parameters are positive, 
		but they do not have to be small or large:
		$\overline{\varrho}$, 
		$\rightu$,
		$\interestingu$,
		$\leftu$, 
		$\mathring{\updelta}$, 
		$\mathring{\updelta}_*$,
		and
		$\boringregionmupositive$.
	\item We assume that $0 < \secondtransversalderivativemulowerbound < 1$. We make no other assumptions
		on $\secondtransversalderivativemulowerbound$.
	\item We assume that $\mathring{\upalpha} > 0$, and that $\mathring{\upalpha}$ is small relative to $1$
		and small relative to the background density $\overline{\varrho}$.
	\item We assume that 
			$\mupositive$ and $|\timefunction_0|$
			are small relative to 
			$1$,
			$\mathring{\updelta}^{-1}$,
			$\mathring{\updelta}_*$,
			and $\secondtransversalderivativemulowerbound$.
			This is possible in view of Remark~\ref{R:INITIALROUGHHYPERSURFACEISCLOSETOSINGULARITY}.
	\item 
		We assume that:
		\begin{align} \label{E:KAPPA0ISSMALLERTHANM2RIGHTUOVER16}
			\muxmulevelsetvalue_0
			& \leq
			\frac{\secondtransversalderivativemulowerbound \interestingu}{16}.
		\end{align}
	\item We assume that $\initialsmall$ is small relative to 
		$\overline{\varrho}$, 
		$1$,
		$\mathring{\upalpha}$,
		$\rightu$,
		$\interestingu$,
		$\leftu$, 
		$\mathring{\updelta}^{-1}$, 
		$\mathring{\updelta}_*$,
		$\mupositive$, 
		and
		$\secondtransversalderivativemulowerbound$.
		\item 
			We assume that
			$0 < \interestingu \leq \leftu$,
			$\mupositive < \frac{\boringregionmupositive}{2}$,
			and
			$\timefunction_0 = - \mupositive$.
		\item Our main results will hold under the assumptions that
			$\fundbootsmall = C \initialsmall$ for some large constant $C$,
			where $\fundbootsmall$ is the bootstrap parameter first appearing in
			Sects.\,\ref{SS:MAINQUANTITATIVEBOOTSTRAPASSUMPTIONS}--\ref{SSS:AUXBOOTSTRAP}.
			This is consistent with the following parameter-size relations, which we assume in order to simplify our
			bootstrap argument:
			\begin{subequations}
			\begin{align} \label{E:DATAEPSILONISSMALLERTHANBOOTSTRAPEPSILONSMALLERTHANSQUAREOFDATAALPHA}
			\initialsmall 
			& \leq 
			\fundbootsmall
			\leq 
			\mathring{\upalpha}^2,
				\\
			\fundbootsmall^{3/2}
			& 
			\leq
			\initialsmall.
			\label{E:NONLINEARINEQUALITYRELATINGDATAEPSILONANDBOOTSTRAPEPSILON}
			\end{align}
			\end{subequations}
\end{itemize}
\end{assumption}

\subsection{Conventions for constants}
\label{SS:CONVENTIONSFORCONSTANTS}
In this section, we state our conventions for how the constants $C$, $\mathfrak{c}$, and $C_{\mydiam}$ 
appearing in our analysis are allowed to depend on the parameters introduced above.

\begin{itemize}
\item The constants $C$ and $\mathfrak{c}$ are free to vary from line to line,
	and we use them in a similar fashion; we mainly use ``$C$'' in our estimates, 
	introducing $\mathfrak{c}$ only in a few arguments in which multiple constants play a role.
These constants are allowed to depend on the nonlinearities (i.e., on the equation of state),
and they can continuously depend on the quantities
$\overline{\varrho}$, 
$\rightu$,
$\interestingu$,
$\leftu$, 
$\mathring{\updelta}^{-1}$, 
$\mathring{\updelta}_*$,
$\boringregionmupositive^{-1}$, 
and
$\secondtransversalderivativemulowerbound^{-1}$
from Sect.\,\ref{SSS:INFORMALDESCRIPTIONOFTHEPARAMETERS}. 
In particular, $C$ and $\mathfrak{c}$ are allowed, in principle, 
to \emph{increase} with respect to 
$\rightu$,
$\interestingu$,
$\leftu$, 
$\mathring{\updelta}^{-1}$, 
$\mathring{\updelta}_*$,
$\boringregionmupositive^{-1}$, 
and
$\secondtransversalderivativemulowerbound^{-1}$.
However, $C$ and $\mathfrak{c}$ can be chosen to be \textbf{independent of the parameters}
$\mathring{\upalpha}$, 
$\initialsmall$, 
$\fundbootsmall$, 
$\timefunction_0$,
$\mulevelsetvalue_0$,
and 
$\muxmulevelsetvalue_0$
under the smallness assumptions of Sect.\,\ref{SS:PARAMETERSIZEASSUMPTIONS}.
\item $A \lesssim B$ means that there exists a constant $C > 0$ (where $C$ has the properties described above) 
such that $A \leq C B$.
\item $A = \mathcal{O}(B)$ means that $|A| \lesssim |B|$. 
\item Constants $C_{\mydiam}$ are also allowed to vary from line to line.
\begin{quote}
However, unlike $C$ and $\mathfrak{c}$, 
the $C_{\mydiam}$ are \textbf{universal} 
in the sense that under the smallness assumptions of Sect.\,\ref{SS:PARAMETERSIZEASSUMPTIONS}, 
they can be chosen to be \textbf{independent} of
$
\overline{\varrho},
	\,
\farrightu,
	\,
\rightu, 
	\,
\interestingu, 
	\, 
\leftu, 
	\,
\timefunction_0,
	\,
\muxmulevelsetvalue_0, 
	\, 
\mupositive,
	\, 
\boringregionmupositive,
	\, 
\mathring{\updelta},
	\, 
\mathring{\updelta}_*,
	\, 
\mathring{\upalpha},	
	\, 
\secondtransversalderivativemulowerbound
$,
and $\initialsmall$,
and also independent of the equation of state. 
\end{quote}
\item $A = \mathcal{O}_{\mydiam}(B)$ means that there exists a constant $C_{\mydiam} > 0$
(where $C_{\mydiam}$ has the properties described above) such that 
$|A| \leq C_{\mydiam} |B|$.
\item As examples, we note that $\mathring{\updelta} \initialsmall \leq 1 \eqdef C_{\mydiam}$
and $\mathring{\updelta}^2 \initialsmall \leq 1 \eqdef C_{\mydiam}$
(because $\initialsmall$ is assumed to be small relative to $\mathring{\updelta}^{-1}$ and 
relative to increasing functions of $\mathring{\updelta}^{-1}$, such as $\mathring{\updelta}^{-2}$),
that $10 \mathring{\upalpha}^2 \leq C_{\mydiam} \mathring{\upalpha}$,
while we have only $\mathring{\updelta} \mathring{\upalpha} \leq C$
(i.e., our smallness assumptions on $\mathring{\upalpha}$ are not strong enough to ensure
that $\mathring{\updelta} \mathring{\upalpha}$ is small because $\mathring{\updelta}$ might be large).
\item As another example, in our estimates 
(e.g., the proof of \eqref{E:MUTRANSVERSALCONVEXITY}), 
by assuming that $|\timefunction_0|$ is small
and using that $C$ is independent of $|\timefunction_0|$
(in particular, $C$ does not implicitly contain any factors of $\frac{1}{|\timefunction_0|}$),
we can ensure that
$C |\timefunction_0| \leq \frac{1}{2 \secondtransversalderivativemulowerbound}$
and
$C |\timefunction_0| \leq \frac{\secondtransversalderivativemulowerbound}{2}$. 

\end{itemize}

\section{Assumptions on the data} 
\label{S:ASSUMPTIONSONTHEDATA}
In this section, we state our assumptions 
on the data in terms of the parameters listed in Sect.\,\ref{SS:LISTOFPARAMETERS}.
Moreover, in Appendix~\ref{A:OPENSETOFDATAEXISTS}, 
we show that our assumptions are satisfied by an open set of data that 
are close to the data of simple isentropic plane-wave solutions.
The analysis in Appendix~\ref{A:OPENSETOFDATAEXISTS} 
is based on the construction -- carried out in Appendix~\ref{A:PS} -- of simple isentropic plane-symmetric
solutions that satisfy the assumptions,
as well as (mostly) standard Cauchy stability arguments,
which we outline.

In the rest of the paper, $\Ntop$ denotes a fixed integer 
representing the maximum number of times
we need to commute the equations of Theorem~\ref{T:GEOMETRICWAVETRANSPORTSYSTEM}
for all estimates to close.
The proof of our main results relies on the following assumption:
\begin{align} \label{E:TOPORDERNUMBEROFDERIVATIES}
\Ntop 
& \ge 24.
\end{align}

\subsection{Background solutions $\RRiemannPS$ and bona fide initial data on $\Sigma_0$}
\label{SSS:BONAFIDEDATA}
Fix any of the ``admissible'' background (shock-forming) simple isentropic plane-symmetric 
solutions that we construct in Appendix~\ref{A:PS},
where we define ``admissible'' in Def.\,\ref{AD:ADMISSIBLEBACKGROUND}.
For such solutions, only a single Riemann invariant is non-vanishing;
we denote it by 
$\RRiemannPS$.
In the rest of Sect.\,\ref{SSS:BONAFIDEDATA}, 
we view $\RRiemannPS$ as a solution
in three spatial dimensions that is independent of the torus coordinates $(x^2,x^3)$.

We now discuss the initial data of a perturbation of one of the background solutions.
We consider the ``bona fide initial data'' of the perturbed solution to be:
\begin{align} \label{E:PERTURBEDBONAFIDEDATA}
\wavearray\big|_{\Sigma_0} 
= (\RRiemann,\LRiemann,v^2,v^3,\Ent) \big|_{\Sigma_0}
& \eqdef 
\left(\RRiemannpertinitial, \LRiemannpertinitial, \vtwopertinitial, \vthreepertinitial, \spertinitial \right),
\end{align}
where $\RRiemannpertinitial, \LRiemannpertinitial, \vtwopertinitial, \vthreepertinitial, \spertinitial: \R \times \T^2 \to\mathbb{R}$ 
are given scalar-valued functions. 
Let
$
\left(\vortrenormalizedpertinitial^i,
\GradEntpertinitial^i,
\VortVortpertinitial^i,
\DivGradEntpertinitial 
\right)_{i=1,2,3}
$
respectively denote the initial data on $\Sigma_0$
of 
$
\left(
\vortrenormalized^i,
\GradEnt^i,
\VortVort^i,
\DivGradEnt
\right)_{i=1,2,3}
$. Note that these data are determined by
$\left(\RRiemannpertinitial, \LRiemannpertinitial, \vtwopertinitial, \vthreepertinitial, \spertinitial \right)$,
the compressible Euler equations
\eqref{E:BVIEVOLUTION}--\eqref{E:BENTROPYEVOLUTION},
definition~\eqref{E:ALMOSTRIEMANNINVARIANTS},
and Def.\,\ref{D:HIGHERORDERFLUIDVARIABLES}.

\begin{remark}[The data of the eikonal function quantities on $\Sigma_0$ are determined]
\label{R:DATAOFEIKONAFUNCTIONQUANTITIES}
Recall that the initial condition of the eikonal function is 
$u|_{\Sigma_0} = -x^1$ (see \eqref{E:EIKONALEQUATION}).
It is straightforward to check that this initial condition and $\wavearray\big|_{\Sigma_0}$
together determine the data of all of the auxiliary quantities constructed out of $u$, 
such as $\upmu|_{\Sigma_0}$, $\Lunit^i|_{\Sigma_0}$, etc.
\end{remark}

We next note that relative to the Cartesian coordinates $(t,x^1,x^2,x^3)$,
we have (see definition \eqref{E:ELLTUSMOOTHTORI}): 
\begin{align} \label{E:SMOOTHTORIINCARTESIANCOORDINATESATTIME0}
	\ell_{0,u} 
	& = \lbrace (0,-u,x^2,x^3) \ | \ (x^2,x^3) \in \mathbb{T}^2 \rbrace.
\end{align}
Moreover, we recall (see \eqref{E:TRUNCATEDSIGMAT}) that for $u_1 \leq u_2$,
relative to the Cartesian coordinates,
we have:
\begin{align} \label{E:U1U2PORTIONOFTIME0CARTESIANHYPERSURFACE}
	\Sigma_0^{[u_1,u_2]} 
	=
	\left\lbrace 
		(0,-x^1,x^2,x^3) \ | \ 
		-x^1 \in [u_1,u_2], 
			\,
		(x^2,x^3) \in \mathbb{T}^2
	\right\rbrace.
\end{align}

In the next definition, we provide a family of norms of the data perturbations on $\Sigma_0$. 
Under suitable assumptions, smallness of the norms suffices for our main
results to hold.

\begin{definition}[Norms of the data perturbation on $\Sigma_0$]
\label{D:PERTURBATIONSMALLNESSINCARTESIANDIFFERENTIALSTRUCTURE}
Let $\Ntop$ be a fixed integer satisfying \eqref{E:TOPORDERNUMBEROFDERIVATIES},
and let $\RRiemannpertinitial^{\text{PS}} \eqdef \RRiemannPS|_{\Sigma_0}$
denote the initial data of $\RRiemannPS$ on $\Sigma_0$.
Given real numbers $u_1 < u_2$,
we define the $\mathring{\Delta}_{\Sigma_0^{[u_1,u_2]}}^{\Ntop+1}$ to be the following 
Sobolev norm
of the perturbation of the data from the data of the background solution
(that is, we subtract off $\RRiemannpertinitial^{\text{PS}}$ and then take the norm):
\begin{align} 
\begin{split} \label{E:PERTURBATIONSMALLNESSINCARTESIANDIFFERENTIALSTRUCTURE}
	\mathring{\Delta}_{\Sigma_0^{[u_1,u_2]}}^{\Ntop+1}
	&
	\eqdef
	\left\| 
		\left(\RRiemannpertinitial - \RRiemannpertinitial^{\text{PS}}, 
			\LRiemannpertinitial, \vtwopertinitial, \vthreepertinitial, \spertinitial 
		\right) 
	\right\|_{H_{\textnormal{Cartesian}}^{\Ntop+1}(\Sigma_0^{[u_1,u_2]})}
		\\
& \ \
	+
	\left\| 
		\left(\vortrenormalizedpertinitial^1,\vortrenormalizedpertinitial^2,\vortrenormalizedpertinitial^3,
		\GradEntpertinitial^1,\GradEntpertinitial^2,\GradEntpertinitial^3
		\right)
	\right\|_{H_{\textnormal{Cartesian}}^{\Ntop}(\Sigma_0^{[u_1,u_2]})}
	+
	\left\| 
		\left(\VortVortpertinitial^1,\VortVortpertinitial^2,\VortVortpertinitial^3,\DivGradEntpertinitial \right)
	\right\|_{H_{\textnormal{Cartesian}}^{\Ntop}(\Sigma_0^{[u_1,u_2]})}
			\\
	& \ \
	+
	\max_{u \in [u_1,u_2]}
	\sum_{|\vec{I}| \leq \Ntop}
	\left\| 
		\left(
			\partial_{\vec{I}} \vortrenormalizedpertinitial^1,
			\partial_{\vec{I}} \vortrenormalizedpertinitial^2,
			\partial_{\vec{I}} \vortrenormalizedpertinitial^3,
			\partial_{\vec{I}} \GradEntpertinitial^1,
			\partial_{\vec{I}} \GradEntpertinitial^2,
			\partial_{\vec{I}} \GradEntpertinitial^3
		\right)
	\right\|_{L_{\textnormal{Cartesian}}^2(\ell_{0,u})}.
\end{split}
\end{align}
In \eqref{E:PERTURBATIONSMALLNESSINCARTESIANDIFFERENTIALSTRUCTURE},
\begin{subequations}
\begin{align} \label{E:CARTESIANSOBOLEVNORMONSIGMA0}
\| f \|_{H_{\textnormal{Cartesian}}^N(\Sigma_0^{[u_1,u_2]} )} 
& 
\eqdef
\left\lbrace
\sum_{|\vec{I}| \leq N}
\int_{\Sigma_0^{[u_1,u_2]}}
	\left[\partial_{\vec{I}} f(t=0,x^1,x^2,x^3) \right]^2
\, \mathrm{d} x^1 \mathrm{d} x^2 \mathrm{d} x^3
\right\rbrace^{1/2},
	\\
\| f \|_{L_{\textnormal{Cartesian}}^2(\ell_{0,u})} 
& 
\eqdef
\left\lbrace
\int_{\mathbb{T}^2}
	\left[\partial_{\vec{I}} f(t=0,x^1 = - u,x^2,x^3) \right]^2
\, \mathrm{d} x^2 \mathrm{d} x^3
\right\rbrace^{1/2},
 \label{E:CARTESIANSOBOLEVNORMONELLTU}
\end{align}
\end{subequations}
where $\vec{I}$ denotes a multi-index of order $|\vec{I}|$ corresponding to
repeated partial differentiation with respect to the Cartesian \emph{spatial} coordinates,
i.e., repeated differentiation with respect to $\partial_1, \partial_2, \partial_3$.
In particular, $\| f \|_{H_{\textnormal{Cartesian}}^N(\Sigma_0^{[u_1,u_2]})}$
is the standard order $N$ Sobolev norm of $f$ along $\Sigma_0^{[u_1,u_2]}$,
while the 
$\| \cdot \|_{L_{\textnormal{Cartesian}}^2(\ell_{0,u})}$ norm sum on the last line 
of RHS~\eqref{E:PERTURBATIONSMALLNESSINCARTESIANDIFFERENTIALSTRUCTURE}
controls tangential \emph{and} transversal spatial derivatives of
$(\vortrenormalizedpertinitial,\GradEntpertinitial)$ along $\ell_{0,u}$.
\end{definition}

\subsection{The assumptions on the initial data}
\label{SS:ASSUMPTIONSONDATA}
In order for our main results to hold,
it suffices for $\mathring{\Delta}_{\Sigma_0^{[-\farrightu,\leftu]}}^{\Ntop+1}$
to be sufficiently small,
where $\blowupdeltaPS > 0$ and $\farrightu \eqdef \rightu + \frac{18}{\blowupdeltaPS} > 0$
are parameters associated to the background solution
(see Appendix~\ref{A:PS})
and $\Ntop$ is the fixed integer satisfying \eqref{E:TOPORDERNUMBEROFDERIVATIES}.

For the solutions under study, the analysis is difficult/interesting
only near the singularity. Moreover, the structures
we use to detect the singular boundary become evident only late
in the classical evolution, i.e., close to the singular boundary. 
For this reason, in Sect.\,\ref{SS:ASSUMPTIONSONDATA},
we find it convenient to describe the state of 
the fluid solution and acoustic geometry 
(e.g., $\upmu$, $\Lunit^i$, and $\upchi$)
on the ``late-time '' rough hypersurface portion 
$\hypthreearg{\timefunction_0}{[- \rightu,\leftu]}{\muxmulevelsetvalue}$,
the rough tori $\twoargroughtori{\timefunction_0,u}{\muxmulevelsetvalue}$,
as well as the null hypersurface portion
$\nullhyparg{- \rightu}^{4 \mathring{\updelta}_*^{-1}}$,
where $\mathring{\updelta}_* > 0$ is defined in \eqref{E:DELTASTARDEF}.
Our description is in terms of \emph{assumed bounds} for various norms of the solution on 
$\hypthreearg{\timefunction_0}{[- \rightu,\leftu]}{\muxmulevelsetvalue}$,
$\twoargroughtori{\timefunction_0,u}{\muxmulevelsetvalue}$,
and
$\nullhyparg{- \rightu}^{4 \mathring{\updelta}_*^{-1}}$
in terms of the parameters of 
Sect.\,\ref{SS:PARAMETERSIZEASSUMPTIONS}.
We refer to 
Fig.\,\ref{F:LATETIMEDATAHYPERSURFACEANDNULLDATAHYPERSURFACE} for an illustration of these data hypersurfaces.

In Appendix~\ref{A:OPENSETOFDATAEXISTS},
we use Cauchy stability-type arguments to sketch a proof 
that if $\mathring{\Delta}_{\Sigma_0^{[\farrightu,\leftu]}}^{\Ntop+1}$
is sufficiently small, then the assumptions we state in Sect.\,\ref{SS:ASSUMPTIONSONDATA} 
are satisfied,
where our main smallness parameter $\initialsmall$ (which vanishes for the background solutions)
satisfies $\initialsmall \lesssim \mathring{\Delta}_{\Sigma_0^{[\farrightu,\leftu]}}^{\Ntop+1}$
whenever $\mathring{\Delta}_{\Sigma_0^{[\farrightu,\leftu]}}^{\Ntop+1}$ is sufficiently small,
where the implicit constants depend on the background solution.
Since our main results apply whenever $\initialsmall$ is sufficiently small,
this in particular shows that there are open sets of data for which
our main theorem holds.

\begin{remark}[We don't need the background solutions to control the dynamics]
	\label{R:ROLEOFPLANESYMMETRICBACKGROUNDSOLUTIONS}
	We ``use'' the plane-symmetric background solutions only to show that there exist
	open sets of data that satisfy our assumptions on 
	$\hypthreearg{\timefunction_0}{[- \rightu,\leftu]}{\muxmulevelsetvalue}$,
	$\twoargroughtori{\timefunction_0,u}{\muxmulevelsetvalue}$,
	and
	$\nullhyparg{- \rightu}^{4 \mathring{\updelta}_*^{-1}}$.
	When studying the evolution to the future of $\hypthreearg{\timefunction_0}{[- \rightu,\leftu]}{\muxmulevelsetvalue}$,
	we never actually have to ``subtract off'' any background solution
	or even refer to one at all.  
\end{remark}

\subsubsection{Quantitative assumptions on the data of the fluid and eikonal function quantities along
$\hypthreearg{\timefunction_0}{[- \rightu,\leftu]}{\muxmulevelsetvalue}$,
	$\twoargroughtori{\timefunction_0,u}{\muxmulevelsetvalue}$,
	and
	$\nullhyparg{- \rightu}^{4 \mathring{\updelta}_*^{-1}}$} 
	\label{SSS:QUANTITATIVEASSUMPTIONSONDATAAWAYFROMSYMMETRY} 
In this section, we state quantitative assumptions on the data of the fluid variables and the eikonal
function quantities along $\hypthreearg{\timefunction_0}{[- \rightu,\leftu]}{\muxmulevelsetvalue}$,
	$\twoargroughtori{\timefunction_0,u}{\muxmulevelsetvalue}$,
	and
	$\nullhyparg{- \rightu}^{4 \mathring{\updelta}_*^{-1}}$.
We again emphasize that our assumptions hold for perturbations of 
the simple isentropic plane-symmetric solutions from 
Appendix~\ref{A:OPENSETOFDATAEXISTS}. 

We start by defining the data-parameter $\mathring{\updelta}_*$.

\begin{definition}[Key data-parameter tied to the Cartesian time of first blowup]
\label{D:DELTASTARDEF}
We define $\mathring{\updelta}_*$ by: 
\begin{align} \label{E:DELTASTARDEF}
	\mathring{\updelta}_* 
	& 
	\eqdef 
	\sup_{\hypthreearg{\timefunction_0}{[- \rightu,\leftu]}{\muxmulevelsetvalue}} 
	\frac{1}{2} 
	\left[ 
		(\Speed^{-1} \Speed;_{\LogDensity} + 1) \muX \RRiemann 
	\right]_+.
\end{align}
\end{definition}

We assume that:
\begin{align} \label{E:DELTASTARPOSITIVE}
	\mathring{\updelta}_* 
	& > 0.
\end{align}

\begin{remark}[Connection between $\mathring{\updelta}_*$ and the Cartesian time of first blowup]
\label{R:CONNECTOINBETWEENDELTASTARANDSHOCKTIME}
For simple isentropic plane-symmetric solutions, the Cartesian time of first shock formation is precisely $\frac{1}{\mathring{\updelta}_*}$; 
see Appendix~\ref{A:PS}.
Our main results show that for the perturbed solutions under study, 
the Cartesian time of first blowup, which we denote here by $t_{\mbox{\tiny First shock}}$,
satisfies $t_{\mbox{\tiny First shock}} = \lbrace 1 + \mathcal{O}(\initialsmall) \rbrace \frac{1}{\mathring{\updelta}_*}$.
\end{remark}

We refer to Sect.\,\ref{SS:STRINGSOFCOMMUTATIONVECTORFIELDS} for notation regarding
strings of commutation vectorfields and to Defs.\,\ref{D:GEOMETRICL2NORMS}
and \ref{D:GEOLINFINITYROUGHTORUS} for the definitions of our $L^2$ and $L^{\infty}$ norms.

\medskip

\noindent \underline{\textbf{$L^{\infty}$ assumptions on the wave variables and their pure transversal derivatives}}.
For $u \in [- \rightu,\leftu]$ and $M = 1,2,3,4$, we assume 
(recall that $\wavearray$ and $\wavearraypartial$ are defined in Def.\,\ref{D:ARRAYSOFWAVEVARIABLES}):
\begin{subequations}
\begin{align}  
	\left \| \RRiemann \right \|_{L^{\infty}\left(\twoargroughtori{\timefunction_0,u}{\muxmulevelsetvalue}\right)} 
	& \leq \mathring{\upalpha}, 
	\label{E:LINFINITYINITIALROUGHHYPERSURFACEBOUNDRRIEMANNAMPLITUDE}
		\\
	\left \| 
		\muX^M \RRiemann 
	\right\|_{L^{\infty}\left(\twoargroughtori{\timefunction_0,u}{\muxmulevelsetvalue}\right)} 
	& 
	\leq \mathring{\updelta},  
		\label{E:LINFINITYINITIALROUGHHYPERSURFACEBOUNDRRIEMANNTRANSVERSALDERIVATIVES}
		\\
	\left\| 
		\wavearraypartial
	\right\|_{L^{\infty}\left(\twoargroughtori{\timefunction_0,u}{\muxmulevelsetvalue}\right)},
		\,
	\left\| 
		\muX^M \wavearraypartial
	\right\|_{L^{\infty}\left(\twoargroughtori{\timefunction_0,u}{\muxmulevelsetvalue}\right)}
	& 
	\leq \initialsmall. 
	\label{E:LINFINITYINITIALROUGHHYPERSURFACEBOUNDSMALLPURETRANSVERSALDERIVATIVESSMALLWAVEVARIABLES} 
\end{align}
\end{subequations}

\medskip

\noindent \underline{\textbf{$L^{\infty}$ assumptions involving tangential derivatives of the wave variables}}.
For $u \in [- \rightu,\leftu]$, we assume:
	\begin{equation} 	\label{E:LINFINITYINITIALROUGHHYPERSURFACEBOUNDSWAVEVARIABLESTRANSVERSALANDTANDERIVATIVE} 
	\begin{split}
		&
		\left \| 
			\tander^{[1, \Ntop - 10]}\wavearray 
		\right \|_{L^{\infty}\left(\twoargroughtori{\timefunction_0,u}{\muxmulevelsetvalue}\right)}, 
		\, 
		\left\| 
			\comdersmall^{[1,\Ntop - 10;1]} \wavearray 
		\right \|_{L^{\infty}\left(\twoargroughtori{\timefunction_0,u}{\muxmulevelsetvalue}\right)}, 
		\, 
		\left\| 
			\comdersmall^{[1,6];2} \wavearray 
		\right \|_{L^{\infty}\left(\twoargroughtori{\timefunction_0,u}{\muxmulevelsetvalue}\right)}, 
		\\
		&
		\left\| 
			\comdersmall^{[1,5];3} \wavearray 
		\right \|_{L^{\infty}\left(\twoargroughtori{\timefunction_0,u}{\muxmulevelsetvalue}\right)},
			\,
		\left\| 
			\Lunit \muX \muX \muX \muX \wavearray 
		\right\|_{L^{\infty}\left(\twoargroughtori{\timefunction_0,u}{\muxmulevelsetvalue}\right)} 
		\leq \initialsmall. 
	\end{split}
	\end{equation}

\medskip

\noindent \underline{\textbf{$L^{\infty}$ assumptions involving tangential derivatives of the transport variables}}.
	For $u \in [- \rightu,\leftu]$, we assume:
	\begin{align} \label{E:LINFINITYINITIALROUGHHYPERSURFACEBOUNDTRANSPORTVARIABLES}
		\left\| 
			\tander^{\leq \Ntop - 11}(\vortrenormalized,\GradEnt) 
			\right\|_{L^{\infty}\left(\twoargroughtori{\timefunction_0,u}{\muxmulevelsetvalue}\right)}, 
			\,
		\left\| 
			\tander^{\leq \Ntop - 12}(\VortVort,\DivGradEnt) 
		\right\|_{L^{\infty}\left(\twoargroughtori{\timefunction_0,u}{\muxmulevelsetvalue}\right)}  
		& \leq \initialsmall. 
	\end{align}

\medskip

\noindent \underline{\textbf{$L^2$ assumptions along $\hypthreearg{\timefunction_0}{[- \rightu,\leftu]}{\muxmulevelsetvalue}$}}.
We assume:
\begin{subequations}
\begin{align}
\left\| \tander^{[1,\Ntop+1]} \wavearray\right \|_{L^2 \left(\hypthreearg{\timefunction_0}{[- \rightu,\leftu]}{\muxmulevelsetvalue} \right)} 
	& \leq  \initialsmall, 
	\label{E:TANGENTIALL2NORMSOFWAVEVARIABLESSMALLALONGINITIALROUGHHYPERSURFACE} 
		\\
\left\| 
	\muX \tander^{[1,\Ntop]} \wavearray 
\right\|_{L^2\left(\hypthreearg{\timefunction_0}{[- \rightu,\leftu]}{\muxmulevelsetvalue} \right)} 
& \leq  \initialsmall, 
	\label{E:TRANSVERSALDERIVATIVEOFTANGENTIALL2NORMSOFWAVEVARIABLESSMALLALONGINITIALROUGHHYPERSURFACE}  
		\\
\left\| 
	\tander^{\leq \Ntop} (\vortrenormalized,\GradEnt)
\right\|_{L^2 \left( \hypthreearg{\timefunction_0}{[- \rightu,\leftu]}{\muxmulevelsetvalue} \right)} & \leq  \initialsmall, 
	\label{E:TANGENTIALL2NORMSOFTRANSPORTVARIABLESSMALLALONGINITIALROUGHHYPERSURFACE} 
	\\
\left\| 
	\tander^{\leq \Ntop} (\VortVort,\DivGradEnt)
\right\|_{L^2 \left(\hypthreearg{\timefunction_0}{[- \rightu,\leftu]}{\muxmulevelsetvalue} \right)} 
& \leq  \initialsmall. 
	\label{E:TANGENTIALL2NORMSOFMODIFIEDFLUIDVARIABLESSMALLALONGINITIALROUGHHYPERSURFACE}
\end{align}
\end{subequations}
%

\medskip 

\noindent \underline{\textbf{$L^2$ assumptions along $\nullhyparg{- \rightu}^{[0,\frac{4}{\mathring{\updelta}_*}]}$}}. 
We assume:
\begin{subequations}
\begin{align}
\left\| \tander^{\leq \Ntop + 1} \wavearray \right\|_{L^2\left(\nullhyparg{- \rightu}^{[0,\frac{4}{\mathring{\updelta}_*}]}\right)} 
& \leq \initialsmall, 
\label{E:WAVESARESMALLONINITIALNULLHYERSURFACE}
\\
\left\| \tander^{\leq \Ntop } (\vortrenormalized,\GradEnt)\right\|_{L^2\left(\nullhyparg{- \rightu}^{[0,\frac{4}{\mathring{\updelta}_*}]} \right)} 
& \leq \initialsmall, 
\label{E:VORTICITYANDENTROPYGRADIENTARESMALLONINITIALNULLHYERSURFACE}
\\
\left\| \tander^{\leq \Ntop } (\VortVort,\DivGradEnt)\right\|_{L^2\left(\nullhyparg{- \rightu}^{[0,\frac{4}{\mathring{\updelta}_*}]}\right)} 
& \leq \initialsmall. 
\label{E:MODIFIEDFLUIDVARIABLESARESMALLONINITIALNULLHYERSURFACE}
\end{align}
\end{subequations}

\medskip

\noindent \underline{\textbf{$L^2$ assumptions along $\twoargroughtori{\timefunction_0,u}{\muxmulevelsetvalue}$}}.
For $u \in [- \rightu, \leftu]$, we assume:
\begin{subequations}
\begin{align}
	\left\| 
		\tander^{[1, \Ntop]} \wavearray 
	\right\|_{L^2(\twoargroughtori{\timefunction_0,u}{\muxmulevelsetvalue})}
	& \leq \initialsmall, 
		\label{E:SMALLDATAOFPSIONINITIALROUGHTORI}
			\\
	\left\| 
		\tander^{\leq \Ntop}(\vortrenormalized,\GradEnt) 
	\right\|_{L^2(\twoargroughtori{\timefunction_0,u}{\muxmulevelsetvalue})}
	& \leq \initialsmall,
	\label{E:SMALLDATAOFVORTICITYANDENTORPYGRADIENTONINITIALROUGHTORI}
		\\
	\left\| 
		\tander^{\leq \Ntop-1}(\VortVort,\DivGradEnt) 
	\right\|_{L^2(\twoargroughtori{\timefunction_0,u}{\muxmulevelsetvalue})}
	& \leq \initialsmall.
	\label{E:SMALLDATAOFMODIFIEDFLUIDVARIABLESONINITIALROUGHTORI}
\end{align}
\end{subequations}

\medskip

\noindent \underline{\textbf{$L^{\infty}$ assumptions tied to transversal derivatives of the eikonal function quantities}}.
For $u \in [- \rightu,\leftu]$ and $M = 0,1,2,3$, we assume:\footnote{Some of these ``assumptions'' 
can in fact be derived as a consequence of other assumptions,
up to constant factors that we absorb into the parameters
$\mathring{\updelta}$, 
$\mathring{\upalpha}$, 
and $\initialsmall$. 
For convenience, instead of deriving those ``assumptions,'' we just assume them.
The same is true for other assumptions on the eikonal function quantities stated below.
For example, even though the data of the $\muX$-derivatives of $\Lunit^i$
could be controlled via the identity \eqref{E:MUXLISCHEMATICIDENTITY},
we just assume \eqref{E:LINFINITYINITIALROUGHHYPERSURFACETRANSVERSALDERIVATIVESOFL1}.
We also refer to \eqref{E:MUTRANSPORT} and \eqref{E:IDENTITYFORMAINTERMDRIVINGTHESHOCK}
for intuition behind the form of 
RHS~\eqref{E:LINFINITYINITIALROUGHHYPERSURFACELDERIVATIVEOFMUANDTRANSVERSALDERIVATIVES},
to \eqref{E:PSMUINITIAL} for intuition behind the first term 
on RHS~\eqref{E:LINFINITYINITIALROUGHHYPERSURFACEMUANDTRANSVERSALDERIVATIVES}
and to Remark~\ref{R:CONNECTOINBETWEENDELTASTARANDSHOCKTIME} for intuition
behind the factor $\mathring{\updelta}_*^{-1}$
in the second term on RHS~\eqref{E:LINFINITYINITIALROUGHHYPERSURFACEMUANDTRANSVERSALDERIVATIVES}.
\label{FN:INTUITIONONMUDATAASSUMPTIONS}}
	\begin{subequations} 
	\begin{align}  \label{E:LINFINITYINITIALROUGHHYPERSURFACELDERIVATIVEOFMUANDTRANSVERSALDERIVATIVES} 
	\left\| \Lunit \muX^M \upmu \right\|_{L^{\infty}(\twoargroughtori{\timefunction_0,u}{\muxmulevelsetvalue})}
	& \leq 
			\frac{1}{2} 
			\left\|
				\muX^M
				\left\lbrace
				\Speed^{-1}(\Speed^{-1} \Speed_{;\LogDensity} + 1)  
				\muX \RRiemann 
				\right\rbrace
			\right\|_{L^{\infty}(\twoargroughtori{\timefunction_0,u}{\muxmulevelsetvalue})}
		+
		\initialsmall,
			\\
	\left\| \muX^M \upmu \right \|_{L^{\infty}(\twoargroughtori{\timefunction_0,u}{\muxmulevelsetvalue})}
	& \leq
		\left\|
			\muX^M
			\left\lbrace
				\Speed^{-1}
			\right\rbrace
		\right\|_{L^{\infty}(\twoargroughtori{\timefunction_0,u}{\muxmulevelsetvalue})}
		+
		\frac{1}{2}
		\mathring{\updelta}_*^{-1}
		\left\|
				\muX^M
				\left\lbrace
				\Speed^{-1}(\Speed^{-1} \Speed_{;\LogDensity} + 1)  
				\muX \RRiemann 
				\right\rbrace
			\right\|_{L^{\infty}(\twoargroughtori{\timefunction_0,u}{\muxmulevelsetvalue})}
		+
		\initialsmall.
		\label{E:LINFINITYINITIALROUGHHYPERSURFACEMUANDTRANSVERSALDERIVATIVES}  
\end{align}

Similarly, for $u \in [- \rightu,\leftu]$ and $M = 1,2,3,4$, we assume:
\begin{align}
	 \left\| \Lunit \muX^M \Lsmall^i \right\|_{L^{\infty}(\twoargroughtori{\timefunction_0,u}{\muxmulevelsetvalue})}
	& \leq 
			\initialsmall,
			\label{E:LINFINITYINITIALROUGHHYPERSURFACELDERIVATIVEOFTRANSVERSALDERIVATIVESOFLI}  
				\\
	\left\| \muX^M \Lsmall^1 \right \|_{L^{\infty}(\twoargroughtori{\timefunction_0,u}{\muxmulevelsetvalue})}
	& \leq
		\mathring{\updelta},
		\label{E:LINFINITYINITIALROUGHHYPERSURFACETRANSVERSALDERIVATIVESOFL1}
				\\
	\left\| \muX^M \Lsmall^A \right \|_{L^{\infty}(\twoargroughtori{\timefunction_0,u}{\muxmulevelsetvalue})}
		& \leq
		\initialsmall.
		\label{E:LINFINITYINITIALROUGHHYPERSURFACETRANSVERSALDERIVATIVESOFLA}
	\end{align}
	\end{subequations}

\medskip

\noindent \underline{\textbf{$L^{\infty}$ assumptions involving tangential derivatives of the eikonal function quantities}}.
For $u \in [- \rightu,\leftu]$, we assume:
	\begin{subequations}
	\begin{align}
	\left\| \tander_*^{[1,\Ntop-12]} \upmu \right\|_{L^{\infty}(\twoargroughtori{\timefunction_0,u}{\muxmulevelsetvalue})},
			\,
	\left\| \comderdoublesmall^{[1,5];1} \upmu \right\|_{L^{\infty}(\twoargroughtori{\timefunction_0,u}{\muxmulevelsetvalue})},
		\, 
	\left\| \comderdoublesmall^{[1,4];2} \upmu \right\|_{L^{\infty}(\twoargroughtori{\timefunction_0,u}{\muxmulevelsetvalue})}
	& \leq \initialsmall,
		\label{E:LINFINITYINITIALROUGHHYPERSURFACESMALLBOUNDSMUTRANSVERSALANDTANDERIVATIVE}  
			\\
	\left\|\Lsmall^1 \right\|_{L^{\infty}(\twoargroughtori{\timefunction_0,u}{\muxmulevelsetvalue})} 
	& 
	\leq \mathring{\upalpha},	
		\label{E:LINFINITYINITIALROUGHHYPERSURFACEBOUNDL1SMALLAMPLITUDE}
			\\
	\begin{split}
	\left\| 
		\Lsmall^A 
	\right\|_{L^{\infty}(\twoargroughtori{\timefunction_0,u}{\muxmulevelsetvalue})},
		\,
	\left\| 
		\tander^{\leq \Ntop-11} \Lsmall^i 
	\right\|_{L^{\infty}(\twoargroughtori{\timefunction_0,u}{\muxmulevelsetvalue})}, 
		\, 
	\left\| 
			\comdersmall^{[1,\Ntop-12];1} \Lsmall^i 
	\right\|_{L^{\infty}(\twoargroughtori{\timefunction_0,u}{\muxmulevelsetvalue})}, 
		\label{E:LINFINITYINITIALROUGHHYPERSURFACESMALLBOUNDLISMALLMIXEDTRANSVERSALTANGENTIAL}
			\\
	\left\| 
		\comdersmall^{[1,5];2} \Lsmall^i 
	\right\|_{L^{\infty}(\twoargroughtori{\timefunction_0,u}{\muxmulevelsetvalue})},
		\, 
	\left\| 
		\comdersmall^{[1,4];3} \Lsmall^i 
	\right\|_{L^{\infty}(\twoargroughtori{\timefunction_0,u}{\muxmulevelsetvalue})}
	& 
	\leq \initialsmall.
	\end{split}
	\end{align}
	\end{subequations}

\medskip

\noindent \underline{\textbf{$L^2$ assumptions on the eikonal function quantities}}.
We assume:
\begin{subequations}
\begin{align}
	\left\|  
		\tandersmall^{[1,\Ntop]} \upmu 
	\right \|_{L^2\left(\hypthreearg{\timefunction_0}{[- \rightu,\leftu]}{\muxmulevelsetvalue}\right)}, 
		\, 
	\left\| 
		\tander^{[1,\Ntop]} \Lsmall^i 
	\right \|_{L^2\left(\hypthreearg{\timefunction_0}{[- \rightu,\leftu]}{\muxmulevelsetvalue}\right)}, 
		\, 
	\left\| 
		\comdersmall^{[1,\Ntop];1} \Lsmall^i 
	\right\|_{L^2\left(\hypthreearg{\timefunction_0}{[- \rightu,\leftu]}{\muxmulevelsetvalue}\right)} 
	& \leq \initialsmall,
		\label{E:L2DATASMALLNESSFOREIKONALFUNCTIONQUANTITIES} 
			\\
	\left\| 
		\angLie_{\comder}^{\leq \Ntop-1;1} \upchi 
	\right\|_{L^2\left(\hypthreearg{\timefunction}{[- \rightu,\leftu]}{\muxmulevelsetvalue}\right)},
		\,
	\left\| 
		\angLie_{\tander}^{\leq \Ntop} \upchi 
	\right\|_{L^2\left(\hypthreearg{\timefunction}{[- \rightu,\leftu]}{\muxmulevelsetvalue}\right)}
	& \leq \initialsmall.
		\label{E:L2DATASMALLNESSFORCHIUPTOTOPORDER}
\end{align}
\end{subequations}

\subsubsection{Data assumptions for the size of $t$ and $x^1$ on $\hypthreearg{\timefunction_0}{[- \rightu,\leftu]}{\muxmulevelsetvalue}$}
\label{SSS:DATAASSUMPTIONSIZEOFCARTESIANTANDX1}
With $\mathring{\updelta}_*$ as in \eqref{E:DELTASTARDEF}, we assume:
\begin{subequations}
\begin{align} \label{E:DATAASSUMPTIONSIZEOFCARTESIANT}
		\frac{1}{2 \mathring{\updelta}_*}
		& 
		\leq
		\min_{\hypthreearg{\timefunction_0}{[- \rightu,\leftu]}{\muxmulevelsetvalue}} t
		\leq
		\sup_{\hypthreearg{\timefunction_0}{[- \rightu,\leftu]}{\muxmulevelsetvalue}} t
		\leq 
		\frac{2}{\mathring{\updelta}_*},
			\\
		- 
		\leftu
		+
		\frac{1}{2 \mathring{\updelta}_*}
		& 
		\leq
		\min_{\hypthreearg{\timefunction_0}{[- \rightu,\leftu]}{\muxmulevelsetvalue}} x^1
		\leq
		\sup_{\hypthreearg{\timefunction_0}{[- \rightu,\leftu]}{\muxmulevelsetvalue}} x^1
		\leq 
		\rightu
		+
		\frac{2}{\mathring{\updelta}_*}.
		\label{E:DATAASSUMPTIONSIZEOFCARTESIANX1}
	\end{align}
	\end{subequations}


\subsubsection{Localized assumptions on the data of $\upmu$ and its derivatives on the initial rough hypersurface}
\label{SSS:LOCALIZEDDATAASSUMPTIONSFORMUANDDERIVATIVES}
We now make localized quantitative and qualitative assumptions on the behavior of $\upmu$ along 
$\hypthreearg{\timefunction_0}{[- \rightu,\leftu]}{\muxmulevelsetvalue}$,
where we recall that $0 < \interestingu < \leftu$ (see Fig.\,\ref{F:REGIONSWHEREWEDERIVEESTIMATES}).

\medskip

\noindent  \underline{\textbf{Transversal convexity of $\upmu$.}}
We assume that there exists a constant 
$\secondtransversalderivativemulowerbound$
satisfying
$0 < \secondtransversalderivativemulowerbound < 1$ 
such that:\footnote{Most of the estimates assumed in \eqref{E:DATATASSUMPTIONMUTRANSVERSALCONVEXITY}
are redundant in the sense that they could derived, up to constant factors, as consequences of other assumptions;
it is only for convenience that we have chosen to make assumptions on all these quantities.}
\begin{equation}  \label{E:DATATASSUMPTIONMUTRANSVERSALCONVEXITY}
\begin{split}
			\secondtransversalderivativemulowerbound
			& \leq 
			\min_{\hypthreearg{\timefunction_0}{[-\interestingu,\interestingu]}{\muxmulevelsetvalue}}
				\left\lbrace
				\Wtransarg{\muxmulevelsetvalue} \Wtransarg{\muxmulevelsetvalue} \upmu,
					\,
				\Wtransarg{\muxmulevelsetvalue} \muX \upmu,
					\,
				\muX \muX \upmu,
					\,
				\muX \muX \upmu + \frac{\muxmulevelsetvalue \Lunit \muX \upmu}{\Lunit \upmu},
					\,
				\geop{u} \muX \upmu,
					\,
				\geop{u} \muX \upmu - \frac{(\geop{u} \upmu) \geop{t} \muX \upmu}{\geop{t} \upmu}, 
					\,
				\roughgeop{u} \muX \upmu
				\right\rbrace
					\\
			&
			\leq 
			\max_{\hypthreearg{\timefunction_0}{[-\interestingu,\interestingu]}{\muxmulevelsetvalue}} 
			\left\lbrace
				\Wtransarg{\muxmulevelsetvalue} \Wtransarg{\muxmulevelsetvalue} \upmu,
					\,
				\Wtransarg{\muxmulevelsetvalue} \muX \upmu,
					\,
				\muX \muX \upmu,
					\,
				\muX \muX \upmu + \frac{\muxmulevelsetvalue \Lunit \muX \upmu}{\Lunit \upmu},
					\,
				\geop{u} \muX \upmu,
					\,
				\geop{u} \muX \upmu - \frac{(\geop{u} \upmu) \geop{t} \muX \upmu}{\geop{t} \upmu}, 
					\,
				\roughgeop{u} \muX \upmu
				\right\rbrace
			\leq 
			\frac{1}{\secondtransversalderivativemulowerbound}.
\end{split}
\end{equation}

\medskip
	
\noindent  \underline{\textbf{$\datahypfortimefunctionarg{-\muxmulevelsetvalue}$
is located near $\lbrace u=0 \rbrace$.}}
With $\datahypfortimefunctionarg{-\muxmulevelsetvalue}$ as in \eqref{E:LEVELSETSOFMUXMU}, 
we assume that:
\begin{subequations}
\begin{align} \label{E:DATATASSUMPTIONMUXMUKAPPALEVELSETLOCATION}
			\datahypfortimefunctionarg{-\muxmulevelsetvalue} 
			& 
			\subset
			\hypthreearg{\timefunction_0}{[-\frac{1}{4}\interestingu, \frac{1}{4} \interestingu]}{\muxmulevelsetvalue},
				\\
			\min_{\hypthreearg{\timefunction_0}{[-\interestingu,\interestingu]}{\muxmulevelsetvalue}
			\backslash
			\hypthreearg{\timefunction_0}{[-\frac{1}{2} \interestingu, \frac{1}{2} \interestingu]}{\muxmulevelsetvalue}}
			|\muX \upmu + \muxmulevelsetvalue| 
			& \geq \frac{\secondtransversalderivativemulowerbound \interestingu}{4}.
			\label{E:DATATASSUMPTIONREGIONWHEREMUXMUKAPPALEVELSETISNOTLOCATED}
\end{align}
\end{subequations}

\medskip

\noindent \underline{\textbf{Quantitative positivity of $\upmu$ away from the interesting region.}}
We assume that there is a constant $\boringregionmupositive > 0$ such that:
\begin{align} \label{E:DATAMUISLARGEINBORINGREGION}
			\min_{\hypthreearg{\timefunction_0}{[- \rightu,\leftu]}{\muxmulevelsetvalue} 
			\backslash 
			\hypthreearg{\timefunction_0}{[-\interestingu,\interestingu]}{\muxmulevelsetvalue}} \upmu
			& \geq \boringregionmupositive.
\end{align}
We also assume (recall that $\mupositive = - \timefunction_0$):
		\begin{align} \label{E:MU1BIGGERTHANMU0}
			\frac{\boringregionmupositive}{2} > \mupositive.
		\end{align}

\medskip
	
\noindent \underline{\textbf{Quantitative negativity of $\Lunit \upmu$.}}
We assume that the following inequalities hold, 
where $\mathring{\updelta}_*$ is defined in \eqref{E:DELTASTARDEF}:
\begin{align} \label{E:DATAASIZEBOUNDSONLMUINTERESTINGREGION}
	- 
	\frac{17}{16}
	\mathring{\updelta}_*
	\leq
	\min_{\hypthreearg{\timefunction_0}{[-\interestingu,\interestingu]}{\muxmulevelsetvalue}} \Lunit \upmu
	\leq 
	\max_{\hypthreearg{\timefunction_0}{[-\interestingu,\interestingu]}{\muxmulevelsetvalue}} \Lunit \upmu
	\leq 
	- 
	\frac{15}{16}
	\mathring{\updelta}_*.
\end{align}

\medskip

\noindent \underline{\textbf{Quantitative bounds on the Jacobian $\CHOVJacobianroughtomumuxmu{\muxmulevelsetvalue}$.}}
We assume that the Jacobian matrix $\CHOVJacobianroughtomumuxmu{\muxmulevelsetvalue}(q)$
defined in \eqref{E:JACOBIANMATRIXFORCHOVFROMROUGHCOORDINATESTOMUWEGIGHTEDXMUCOORDINATES}
is invertible for every 
$q \in \{\timefunction_0\} \times [-\interestingu,\interestingu] \times \mathbb{T}^2$
and that:
\begin{align} \label{E:DATAJACOBIANDETERMINANTRATIOASSUMPTION}
	\sup_{q_1,q_2 \in \{\timefunction_0\} \times [-\interestingu,\interestingu] \times \mathbb{T}^2}
	\left|
		\InverseCHOVJacobianroughtomumuxmu{\muxmulevelsetvalue}(q_1) \CHOVJacobianroughtomumuxmu{\muxmulevelsetvalue}(q_2)
		-
		\mbox{\upshape ID} 
	\right|_{\mbox{\upshape}Euc}
	& \leq \frac{1}{3},
\end{align}
where $|\cdot|_{\mbox{\upshape}Euc}$ is the standard Frobenius norm on matrices
(equal to the square root of the sum of the squares of the matrix entries)
and $\mbox{\upshape ID}$ denotes the $4 \times 4$ identity matrix.

\medskip
\noindent \underline{\textbf{Quantitative bounds on the Jacobian $\CHOVJacobiangeotomumuxmu$.}}
We assume that the Jacobian matrix $\CHOVJacobiangeotomumuxmu(p)$
defined in \eqref{E:JACOBIANMATRIXFORCHOVFROMGEOMETRICCOORDINATESTOMUWEGIGHTEDXMUCOORDINATES}
is invertible for every 
$p \in \hypthreearg{\timefunction_0}{[-\interestingu,\interestingu]}{\muxmulevelsetvalue}$
and that the following holds, where 
$\twoargmumuxtorus{-\timefunction_0}{0} \subset \hypthreearg{\timefunction_0}{[-\interestingu,\interestingu]}{0}$
is the $\upmu$-adapted torus defined in \eqref{E:MUXMUTORI}:
\begin{align} \label{E:GEOTOMUMUXCOORDINATESDATAJACOBIANDETERMINANTRATIOASSUMPTION}
	\sup_{
		\substack{p_1 \in \twoargmumuxtorus{-\timefunction_0}{0}
				\\
			p_2 \in \hypthreearg{\timefunction_0}{[-\interestingu,\interestingu]}{\muxmulevelsetvalue}
			}}
	\left|
		\CHOVJacobiangeotomumuxmu(p_1)
		\InverseCHOVJacobiangeotomumuxmu(p_2) 
		-
		\mbox{\upshape ID} 
	\right|_{\mbox{\upshape}Euc}
	& \leq \frac{1}{3},
\end{align}
where $|\cdot|_{\mbox{\upshape}Euc}$ is the Frobenius norm.

\section{The bootstrap assumptions, except those concerning the wave energies}
\label{S:BOOTSTRAPEVERYTHINGEXCEPTENERGIES}
In proving our main results, we rely on a continuity argument
based on deriving improvements of a set of bootstrap assumptions
for the solution on a region of the form $\twoargMrough{[\timefunction_0,\timefunctionboot),[- \rightu,\leftu]}{\muxmulevelsetvalue}$.
In this section, we set up the bootstrap argument and 
state all the bootstrap assumptions, 
except for the ones concerning $L^2$-type energies, 
which we provide in Sect.\,\ref{SS:BOOTSTRAPASSUMPTIONSFORTHEWAVEENERGIES}.

\subsection{The start of the bootstrap argument: the bootstrap time interval $[\timefunction_0,\timefunctionboot)$ and $\upmuboot$}
\label{SS:BOOTSTRAPTIMEINTERVAL}
From now until Sect.\,\ref{S:EXISTENCEUPTOSINGULARBOUNDARYATFIXEDKAPPA}, 
we assume that there is a classical solution on an ``open-at-the-top'' region of the form
$\twoargMrough{[\timefunction_0,\timefunctionboot),[- \rightu,\leftu]}{\muxmulevelsetvalue}$,
and we will state our bootstrap assumptions on the same region. 
Here and throughout, 
\begin{align} \label{E:BOOTSTRAPTIME}
	\timefunctionboot & \in (3 \timefunction_0/4,0)
\end{align}
is the ``bootstrap rough-time,''
where we recall that the small parameter $\timefunction_0 < 0$ is the most negative value achieved by 
the rough time function $\timefunctionarg{\muxmulevelsetvalue}$. 
In Appendix~\ref{A:OPENSETOFDATAEXISTS}, we provide Cauchy stability arguments
guaranteeing the existence of a region of classical existence of the form
$\twoargMrough{[\timefunction_0,\timefunctionboot),[- \rightu,\leftu]}{\muxmulevelsetvalue}$
for some $\timefunctionboot$ satisfying \eqref{E:BOOTSTRAPTIME}.
Hence, at the start\footnote{In 
Lemma~\ref{L:DIFFEOMORPHICEXTENSIONOFROUGHCOORDINATES},
we will show that $\timefunctionarg{\muxmulevelsetvalue}$ can be suitably extended to have range $[\timefunction_0,\timefunctionboot]$,
and in Theorem~\ref{T:EXISTENCEUPTOTHESINGULARBOUNDARYATFIXEDKAPPA},
we will show that it can be suitably extended to have range $[\timefunction_0,0]$.} 
of our bootstrap argument,
$\timefunctionarg{\muxmulevelsetvalue}$ has range $[\timefunction_0,\timefunctionboot)$,
which contains $[\timefunction_0,3 \timefunction_0/4]$.
We also set:
\begin{align} \label{E:MUBOOTISNEGATIVETIMEFUNCTIONBOOT}
	\upmuboot 
	& \eqdef - \timefunctionboot
	> 0.
\end{align}

\subsection{Bootstrap assumptions tied to the fundamental scaffolding of the analysis}
\label{SS:BOOTSTRAPSCAFFOLDING}
The bootstrap assumptions in this section 
ensure that various fundamental aspects of our approach
(such as the change of variables maps from Sect.\,\ref{SS:ALLTHECHOVMAPS})
are well-defined and enjoy basic properties that
we use throughout the rest of the paper.

\subsubsection{Bootstrap assumptions for the inverse foliation density}
\label{SSS:BAFORINVERSEFOLIATIONDENSITY}

\begin{enumerate}
		\item  We assume that the following estimate holds on 
		$\twoargMrough{[\timefunction_0,\timefunctionboot),[- \rightu,\leftu]}{\muxmulevelsetvalue}$:
		\begin{align} \label{E:BAMUPOSITIVE}
			\upmu & > 0. \tag{\textbf{BA} $\upmu > 0$}
		\end{align}
	\item We assume that $\Lunit \upmu$ and $\geop{t} \upmu$ are quantitatively negative in 
	$\twoargMrough{[\timefunction_0,\timefunctionboot],[-\interestingu,\interestingu]}{\muxmulevelsetvalue}$,
	where $\mathring{\updelta}_*$ is defined in \eqref{E:DELTASTARDEF}:
	\begin{align} \label{E:BABOUNDSONLMUINTERESTINGREGION} \tag{\textbf{BA} $\Lunit \upmu$ neg}
	- 
	\frac{3}{4}
	\updelta_*
	& 
	\leq
	\min_{\twoargMrough{[\timefunction_0,\timefunctionboot),[-\interestingu,\interestingu]}{\muxmulevelsetvalue}} \Lunit \upmu
	\leq 
	\max_{\twoargMrough{[\timefunction_0,\timefunctionboot),[-\interestingu,\interestingu]}{\muxmulevelsetvalue}} \Lunit \upmu
	\leq 
	- 
	\frac{5}{4}
	\updelta_*,
		\\
	- 
	\frac{3}{4}
	\updelta_*
	& 
	\leq
	\min_{\twoargMrough{[\timefunction_0,\timefunctionboot),[-\interestingu,\interestingu]}{\muxmulevelsetvalue}} \geop{t} \upmu
	\leq 
	\max_{\twoargMrough{[\timefunction_0,\timefunctionboot),[-\interestingu,\interestingu]}{\muxmulevelsetvalue}} \geop{t} \upmu
	\leq 
	- 
	\frac{5}{4}
	\updelta_*.
		\label{E:BABOUNDSONGEOMETRICTDERIVATIVEMUINTERESTINGREGION}
		\tag{\textbf{BA} $\geop{t} \upmu$ neg}
\end{align}
		\item We assume that near $\datahypfortimefunctiontwoarg{-\muxmulevelsetvalue}{[\timefunction_0,\timefunctionboot)}$,
		$\upmu$ is quantitatively convex in directions transversal to 
		$\datahypfortimefunctiontwoarg{-\muxmulevelsetvalue}{[\timefunction_0,\timefunctionboot)}$. That is, we assume that
		the following estimates hold on $\twoargMrough{[\timefunction_0,\timefunctionboot),[-\interestingu,\interestingu]}{\muxmulevelsetvalue}$,
		where $0 < \secondtransversalderivativemulowerbound < 1$ is as in
		\eqref{E:DATATASSUMPTIONMUTRANSVERSALCONVEXITY}:
		\begin{align} \label{E:BAMUTRANSVERSALCONVEXITY} \tag{\textbf{BA} $\upmu$ cnvx}
			\begin{split}
			\frac{\secondtransversalderivativemulowerbound}{4}
			& \leq 
			\inf_{\twoargMrough{[\timefunction_0,\timefunctionboot),[-\interestingu,\interestingu]}{\muxmulevelsetvalue}}
				\left\lbrace
				\Wtransarg{\muxmulevelsetvalue} \Wtransarg{\muxmulevelsetvalue} \upmu,
					\,
				\Wtransarg{\muxmulevelsetvalue} \muX \upmu,
					\,
				\muX \muX \upmu,
					\,
				\muX \muX \upmu - \frac{(\muX \upmu) \Lunit \muX \upmu}{\Lunit \upmu},
					\,
				\geop{u} \muX \upmu,
					\,
				\geop{u} \muX \upmu - \frac{(\geop{u} \upmu) \geop{t} \muX \upmu}{\geop{t} \upmu}, 
					\,
				\roughgeop{u} \muX \upmu
				\right\rbrace
					\\
			&
			\leq 
			\sup_{\twoargMrough{[\timefunction_0,\timefunctionboot),[-\interestingu,\interestingu]}{\muxmulevelsetvalue}}
			\left\lbrace
				\Wtransarg{\muxmulevelsetvalue} \Wtransarg{\muxmulevelsetvalue} \upmu,
					\,
				\Wtransarg{\muxmulevelsetvalue} \muX \upmu,
					\,
				\muX \muX \upmu,
					\,
				\muX \muX \upmu - \frac{(\muX \upmu) \Lunit \muX \upmu}{\Lunit \upmu},
					\,
				\geop{u} \muX \upmu,
					\,
				\geop{u} \muX \upmu - \frac{(\geop{u} \upmu) \geop{t} \muX \upmu}{\geop{t} \upmu}, 
					\,
				\roughgeop{u} \muX \upmu
				\right\rbrace
			\leq 
			\frac{4}{\secondtransversalderivativemulowerbound}.
\end{split}
\end{align}
\end{enumerate}

\subsubsection{Bootstrap assumptions for the rough time function}
\label{SSS:BOOTSTRAPROUGHTIMEFUNCTION}
\begin{enumerate}
	\item We assume that
		$\geop t \timefunctionarg{\muxmulevelsetvalue} > 0$ in $\twoargMrough{[\timefunction_0,\timefunctionboot),[- \rightu,\leftu]}{\muxmulevelsetvalue}$.
	\item We assume that the following estimates hold, 
		where $\mathring{\updelta}_*$ is defined in \eqref{E:DELTASTARDEF}:
		\begin{align} \label{E:BALDERIVATIVEOFROUGHTIMEFUNCTIONISAPPROXIMATELYUNITY}
	\frac{3}{4} \updelta_*	
	&
	\leq
	\Lunit \timefunctionarg{\muxmulevelsetvalue}
	\leq 
	\frac{5}{4} \updelta_*,
	&&
	\mbox{on } \twoargMrough{[\timefunction_0,\timefunctionboot),[- \rightu,\leftu]}{\muxmulevelsetvalue}.
	\tag{ \textbf{BA} $\Lunit \timefunctionarg{\muxmulevelsetvalue}$}
\end{align}
\end{enumerate}

\subsubsection{Bootstrap assumptions for change of variables maps}
\label{SSS:BACHOVMAPS}
We make the following assumptions for various change of variables maps.

\begin{enumerate}
	\item The change of variables map 
			$\Upsilon(t,u,x^2,x^3) = (t,x^1,x^2,x^3)$
			defined in \eqref{E:CHOVGEOTOCARTESIAN}
			satisfies
			$\left\| 
				\Upsilon
			\right\|_{C^3_{\textnormal{geo}}\left(\twoargMrough{[\timefunction_0,\timefunctionboot),[- \rightu,\leftu]}{\muxmulevelsetvalue} \right)} 
			< \infty
			$
			and is a diffeomorphism 
			from $\twoargMrough{[\timefunction_0,\timefunctionboot),[- \rightu,\leftu]}{\muxmulevelsetvalue}$
			onto its image.
	\item The change of variables map $\CHOVgeotorough{\muxmulevelsetvalue}(t,u,x^2,x^3)  = (\timefunctionarg{\muxmulevelsetvalue},u,x^2,x^3)$
			defined in \eqref{E:CHOVGEOTOROUGH}
			satisfies
			$\left\| 
				\CHOVgeotorough{\muxmulevelsetvalue}
			\right\|_{C^2_{\textnormal{geo}}\left(\twoargMrough{[\timefunction_0,\timefunctionboot),[- \rightu,\leftu]}{\muxmulevelsetvalue} \right)} 
			< \infty
			$
			and is a diffeomorphism from 
			$\twoargMrough{[\timefunction_0,\timefunctionboot),[- \rightu,\leftu]}{\muxmulevelsetvalue}$
			onto its image, which is $[\timefunction_0,\timefunctionboot) \times [- \rightu,\leftu] \times \mathbb{T}^2$.
		\item The change of variables map $\CHOVroughtomumuxmu{\muxmulevelsetvalue}$ defined in 
		\eqref{E:CHOVFROMROUGHCOORDINATESTOMUWEGIGHTEDXMUCOORDINATES} satisfies 
		$\left\| 
			\CHOVroughtomumuxmu{\muxmulevelsetvalue} 
		\right\|_{C^{1,1}_{\textnormal{rough}}\left( [\timefunction_0,\timefunctionboot)\times[- \rightu,\leftu]\times\T^2\right)} 
		< 
		\infty$ 
		and is a diffeomorphism from 
		$[\timefunction_0,\timefunctionboot) \times [-\interestingu,\interestingu] \times \mathbb{T}^2$ 
		onto its image. 
		Furthermore, the Jacobian matrix $\CHOVJacobianroughtomumuxmu{\muxmulevelsetvalue}(q)$ defined in
		\eqref{E:JACOBIANMATRIXFORCHOVFROMROUGHCOORDINATESTOMUWEGIGHTEDXMUCOORDINATES}
		is invertible for every $q \in [\timefunction_0,\timefunctionboot) \times [-\interestingu,\interestingu]\times \T^2$
		and satisfies:
\begin{align} \label{E:BOOTSTRAPJACOBIANDETERMINANTRATIOASSUMPTION}
	\max_{q_1,q_2 \in  [\timefunction_0,\timefunctionboot) \times [-\interestingu,\interestingu]\times \T^2}
	\left|
		 \InverseCHOVJacobianroughtomumuxmu{\muxmulevelsetvalue}(q_1) \CHOVJacobianroughtomumuxmu{\muxmulevelsetvalue}(q_2)
		-
		\mbox{\upshape ID} 
	\right|_{\mbox{\upshape}Euc}
	& \leq \frac{2}{3}.
\end{align}
\end{enumerate}

\subsubsection{Bootstrap assumptions for the structure and location of $\twoargmumuxtorus{\mulevelsetvalue}{-\muxmulevelsetvalue}$
and $\datahypfortimefunctiontwoarg{-\muxmulevelsetvalue}{[\timefunction_0,\timefunctionboot)}$}
\label{SSS:BOOTSTRAPASSUMPTIONFORTORISTRUCTURE}
Recall that the $\upmu$-adapted tori $\twoargmumuxtorus{\mulevelsetvalue}{-\muxmulevelsetvalue}$ are defined in \eqref{E:MUXMUTORI}
and that the hypersurface portions $\datahypfortimefunctiontwoarg{-\muxmulevelsetvalue}{[\timefunction_0,\timefunctionboot)}$
are defined in \eqref{E:TRUNCATEDLEVELSETSOFMUXMU}.
Also recall that $0 < \upmuboot = - \timefunctionboot < \mupositive =  - \timefunction_0$.

\begin{enumerate}
\item
We assume that for each $\mulevelsetvalue \in (\upmuboot,\mupositive]$,
there exist scalar functions 
$\Cartesiantisafunctiononmumxtoriarg{\mulevelsetvalue}{-\muxmulevelsetvalue}, \, \Eikonalisafunctiononmumuxtoriarg{\mulevelsetvalue}{-\muxmulevelsetvalue} 
\in 
W^{2,\infty}(\mathbb{T}^2)$,
depending on $\mulevelsetvalue$ and $\muxmulevelsetvalue$, 
such that relative to the geometric coordinates $(t,u,x^2,x^3)$,
we have:
\begin{align} \label{E:BAMUTORI}
	\twoargmumuxtorus{\mulevelsetvalue}{-\muxmulevelsetvalue}^{[\timefunction_0,\timefunctionboot)}
	& 
	= \left\lbrace
			\left(\Cartesiantisafunctiononmumxtoriarg{\mulevelsetvalue}{-\muxmulevelsetvalue}(x^2,x^3), 
				\Eikonalisafunctiononmumuxtoriarg{\mulevelsetvalue}{-\muxmulevelsetvalue}(x^2,x^3),
				x^2,x^3 
			\right)
			\ | \
			(x^2,x^3) \in \mathbb{T}^2
		\right\rbrace.
		\tag{\textbf{BA} $\upmu-\textbf{TORI STRUCTURE}$}
\end{align}
In particular, in geometric coordinates, $\twoargmumuxtorus{\mulevelsetvalue}{-\muxmulevelsetvalue}$ is a $W^{2,\infty}$ graph over $\mathbb{T}^2$.
\item We assume that for each fixed $\timefunction \in [\timefunction_0,\timefunctionboot) = [-\mupositive,-\upmuboot)$,
		\begin{align} \label{E:BOOSTRAPTORILOCATION}
			\twoargmumuxtorus{\mulevelsetvalue}{-\muxmulevelsetvalue}
			& 
			\subset
			\hypthreearg{\timefunction}{[-\frac{3 \interestingu}{4},\frac{3 \interestingu}{4}]}{\muxmulevelsetvalue}. 
			\tag{\textbf{BA} $\twoargmumuxtorus{\mulevelsetvalue}{-\muxmulevelsetvalue}-\textbf{LOCATION}$}
		\end{align}
\item We assume that: 
			\begin{align} \label{E:BOOTSTRAPLEVELSETSTRUCTUREANDLOCATIONOFMUXEQUALSMINUSKAPPA}
			\datahypfortimefunctiontwoarg{-\muxmulevelsetvalue}{[\timefunction_0,\timefunctionboot)}
			& 
			\subset
			\twoargMrough{[\timefunction_0,\timefunctionboot),[-\frac{3 \interestingu}{4},\frac{3 \interestingu}{4}]}{\muxmulevelsetvalue}.
			\tag{\textbf{BA} $\datahypfortimefunctiontwoarg{-\muxmulevelsetvalue}{[\timefunction_0,\timefunctionboot)}-\textbf{LOCATION}$}
		\end{align}
\item We assume that the map 
$\embeddatahypersurfacearg{\muxmulevelsetvalue} : (\upmuboot,\mupositive] \times \mathbb{T}^2 
\rightarrow \twoargMrough{[\timefunction_0,\timefunctionboot),[-\frac{3 \interestingu}{4},\frac{3 \interestingu}{4}]}{\muxmulevelsetvalue}$ 
defined by:
	\begin{align} \label{E:EMBEDDATAHYPERSURFACE}
		\embeddatahypersurfacearg{\muxmulevelsetvalue}(\mulevelsetvalue,x^2,x^3)
		& \eqdef
			\left(\Cartesiantisafunctiononmumxtoriarg{\mulevelsetvalue}{-\muxmulevelsetvalue}(x^2,x^3), 
				\Eikonalisafunctiononmumuxtoriarg{\mulevelsetvalue}{-\muxmulevelsetvalue}(x^2,x^3),x^2,x^3
			\right)
	\end{align}
	is a diffeomorphism from $(\upmuboot,\mupositive] \times \mathbb{T}^2$ onto
	$\datahypfortimefunctiontwoarg{-\muxmulevelsetvalue}{[\timefunction_0,\timefunctionboot)}$
	such that 
	for every $\mulevelsetvalue' \in (\upmuboot,\mupositive)$,
	we have:
	\begin{align} \label{E:EMBEDDEDEDMUEQUALSMINUSKAPPAHYPERSURFACEFINITEC11NORMONCOMPACTSUBSETS}
		\embeddatahypersurfacearg{\muxmulevelsetvalue} 
		& \in C^{1,1}([\mulevelsetvalue',\mupositive] \times \mathbb{T}^2),
	\end{align}
	and such that: 
	\begin{align} \label{E:BAMUEQUALSMINUSKAPPAEMBEDDINGCARTESIANTIMEFUNCTIONNEGATIVEMUDERIVATIVE}
		- 
		\infty 
		& 
		<
		\inf_{(\mulevelsetvalue,x^2,x^3) \in (\upmuboot,\mupositive] \times \mathbb{T}^2} 
		\frac{\partial}{\partial \mulevelsetvalue} \Cartesiantisafunctiononmumxtoriarg{\mulevelsetvalue}{-\muxmulevelsetvalue}(x^2,x^3) 
		\leq
		\sup_{(\mulevelsetvalue,x^2,x^3) \in (\upmuboot,\mupositive] \times \mathbb{T}^2} 
		\frac{\partial}{\partial \mulevelsetvalue} \Cartesiantisafunctiononmumxtoriarg{\mulevelsetvalue}{-\muxmulevelsetvalue}(x^2,x^3) 
		< 
		0.
	\end{align}
	In particular,
	$\datahypfortimefunctiontwoarg{-\muxmulevelsetvalue}{[\timefunction_0,\timefunctionboot)}$ is a $C^{1,1}$
	embedded submanifold-with-boundary of 
	$\twoargMrough{[\timefunction_0,\timefunctionboot),[-\frac{3 \interestingu}{4},\frac{3 \interestingu}{4}]}{\muxmulevelsetvalue}$
	whose boundary is $\twoargmumuxtorus{\mupositive}{-\muxmulevelsetvalue}$,
	and:
	\begin{align} \label{E:BOOTSTRAPLEVELSETSTRUCTUREOFMUXEQUALSMINUSKAPPA}
			\datahypfortimefunctiontwoarg{-\muxmulevelsetvalue}{[\timefunction_0,\timefunctionboot)}
			& 
			=
			\bigcup_{\mulevelsetvalue \in (\upmuboot,\mupositive]}
			\twoargmumuxtorus{\mulevelsetvalue}{-\muxmulevelsetvalue}.
			\tag{\textbf{BA} $\datahypfortimefunctiontwoarg{-\muxmulevelsetvalue}{[\timefunction_0,\timefunctionboot)}-\textbf{FOLIATED}$}
		\end{align}
\end{enumerate}

\subsubsection{Bootstrap assumptions for the size of $t$ and $x^1$ on $\hypthreearg{\timefunction}{[- \rightu,\leftu]}{\muxmulevelsetvalue}$}
\label{SSS:BOOTSTRAPASSUMPTIONSSIZEOFCARTESIANTANDX1}
With $\mathring{\updelta}_*$ as in \eqref{E:DELTASTARDEF}, we assume
that the following inequalities hold for $\timefunction \in [\timefunction_0,\timefunctionboot)$:
\begin{subequations}
\begin{align} \label{E:BASIZEOFCARTESIANT} \tag{\textbf{BA} $t-\textbf{SIZE}$} 
		\frac{1}{4 \mathring{\updelta}_*}
		& 
		\leq
		\min_{\hypthreearg{\timefunction}{[- \rightu,\leftu]}{\muxmulevelsetvalue}} t
		\leq
		\max_{\hypthreearg{\timefunction}{[- \rightu,\leftu]}{\muxmulevelsetvalue}} t
		\leq 
		\frac{4}{\mathring{\updelta}_*},
				\\
		- \leftu
		+
		\frac{1}{4 \mathring{\updelta}_*}
		& 
		\leq
		\min_{\hypthreearg{\timefunction}{[- \rightu,\leftu]}{\muxmulevelsetvalue}} x^1
		\leq
		\max_{\hypthreearg{\timefunction}{[- \rightu,\leftu]}{\muxmulevelsetvalue}} x^1
		\leq 
		2 \interestingu
		+
		\frac{4}{\mathring{\updelta}_*}.
		\label{E:BASIZEOFCARTESIANX1} \tag{\textbf{BA} $x^1-\textbf{SIZE}$} 
	\end{align}
	\end{subequations}

\subsubsection{Soft bootstrap assumptions concerning regularity and Sobolev embedding}
\label{SSS:SOFTBACONCERNINGREGULARITY}
We assume that for every $\timefunction \in (\timefunction_0,\timefunctionboot)$, we have:
\begin{subequations}
\begin{align}
		\wavearray, \, \vortrenormalized^i, \, \GradEnt^i, \, \VortVort^i, \, \DivGradEnt 
			& \in C_{\textnormal{geo}}^{3,1}(\twoargMrough{[\timefunction_0,\timefunction],[- \rightu,\leftu]}{\muxmulevelsetvalue}),
				\tag{\textbf{BA} \textbf{Fluid Regularity}} 
					\label{E:FLUIDC31ONCOMPACTSUBSETS} \\
		\Upsilon & \in C_{\textnormal{geo}}^{3,1}(\twoargMrough{[\timefunction_0,\timefunction],[- \rightu,\leftu]}{\muxmulevelsetvalue}),	
			\tag{\textbf{BA} $\Upsilon$-\textbf{Regularity}} 
				\label{E:GEOTOCARTESIANCHOVC31ONCOMPACTSUBSETS} \\
		\Lunit^i, \, \upmu & \in C_{\textnormal{geo}}^{2,1}(\twoargMrough{[\timefunction_0,\timefunction],[- \rightu,\leftu]}{\muxmulevelsetvalue}).
			 \tag{\textbf{BA} \textbf{Geometry Regularity}}  
			\label{E:GEOMETRYC31ONCOMPACTSUBSETS}
\end{align}
\end{subequations}


\begin{remark}
	In \eqref{E:FLUIDC31ONCOMPACTSUBSETS}--\eqref{E:GEOMETRYC31ONCOMPACTSUBSETS},
	we are not making any quantitative assumptions on the size of the norms.
	That is, we are assuming only that the norms are finite,
	e.g., $\| \upmu \|_{C_{\textnormal{geo}}^{2,1}(\twoargMrough{[\timefunction_0,\timefunction],[- \rightu,\leftu]}{\muxmulevelsetvalue})} < \infty$.
	In Lemma~\ref{L:CONTINUOUSEXTNESION}, we will show that all of these norms are $\leq C$.
	
\end{remark}

\subsection{The main quantitative bootstrap assumptions} 
\label{SS:MAINQUANTITATIVEBOOTSTRAPASSUMPTIONS}
We now state our main quantitative bootstrap assumptions.
In the rest of the paper, 
$\fundbootsmall \geq 0$ denotes a small ``bootstrap'' parameter whose smallness we described in Sect.\,\ref{SS:PARAMETERSIZEASSUMPTIONS}.
Later on, we will close our bootstrap argument by setting 
$\fundbootsmall = C \initialsmall$ for some large constant $C$, 
where $\initialsmall$ is the data-size parameter from Sect.\,\ref{SSS:QUANTITATIVEASSUMPTIONSONDATAAWAYFROMSYMMETRY};
see, in particular, Prop.\,\ref{P:IMPROVEMENTOFFUNDAMENTALQUANTITATIVEBOOTSTRAPASSUMPTIONS}.

\subsubsection{Fundamental quantitative bootstrap assumptions}
\label{SSS:FUNDAMENTALQUANTITATIVE}
Our fundamental quantitative bootstrap assumptions for $\wavearray, \, \vortrenormalized, \, \GradEnt, \, \VortVort$, and $\DivGradEnt$ are that the following inequalities hold for $(\timefunction,u) \in [\timefunction_0,\timefunctionboot) \times [- \rightu,\leftu]$: 
\begin{align}
\left\| 
	\tander^{[1,\Ntop-10]} \wavearray
\right\|_{L^{\infty}\left(\twoargroughtori{\timefunction,u}{\muxmulevelsetvalue}\right)},  
	\,
\left\| 
	\tander^{\leq \Ntop-10} (\vortrenormalized,\GradEnt) 
\right\|_{L^{\infty}\left(\twoargroughtori{\timefunction,u}{\muxmulevelsetvalue}\right)} 
& \leq \fundbootsmall, \tag{\textbf{BA}$(\wavearray,\vortrenormalized,S)$ \textbf{FUND}}  
\label{E:FUNDAMENTALQUANTITATIVEBOOTWAVEANDTRANSPORT} 
	\\
\left\| 
	\tander^{\leq \Ntop-11}  (\VortVort, \DivGradEnt) 
\right\|_{L^{\infty}\left(\twoargroughtori{\timefunction,u}{\muxmulevelsetvalue}\right)} 
& \leq \fundbootsmall. \tag{\textbf{BA}$(\VortVort,\DivGradEnt)$ \textbf{FUND}}  
\label{E:FUNDAMENTALQUANTITATIVEBOOTMODIFIEDFLUIDVARS} 
\end{align}

\subsubsection{Auxiliary bootstrap assumptions} \label{SSS:AUXBOOTSTRAP}
To derive pointwise estimates, we find it convenient to make the following auxiliary bootstrap assumptions. 

\medskip

\noindent \underline{\textbf{Auxiliary bootstrap assumptions for small quantities}}.
 We assume that the following inequalities hold for 
$(\timefunction,u) \in [\timefunction_0,\timefunctionboot) \times [- \rightu,\leftu]$
(recall that $\wavearray$ and $\wavearraypartial$ are defined in Def.\,\ref{D:ARRAYSOFWAVEVARIABLES}):
\begin{align} \label{BA:AUXR+} \tag{\textbf{AUX $\RRiemann$ SMALL}} 
	\left\|\RRiemann \right\|_{L^{\infty}(\twoargroughtori{\timefunction,u}{\muxmulevelsetvalue})} 
	& \leq \mathring{\upalpha}^{1/2} + \auxbootsmall,  
		\\
	\left\| \wavearraypartial \right\|_{L^{\infty}(\twoargroughtori{\timefunction,u}{\muxmulevelsetvalue})} 
	& \leq \auxbootsmall,
		 \label{BA:AUXWAVEARRAYPARTIALLINFINITY} 
			\tag{\textbf{AUX $\wavearraypartial$ SMALL}} 
			\\
	 \label{BA:AUXWAVEARRAY} \tag{\textbf{AUX $\wavearray$ SMALL}} 
	\begin{split}  
	\left\|\Lunit \comder^{\leq \Ntop-11;1} \wavearray \right\|_{L^{\infty}(\twoargroughtori{\timefunction,u}{\muxmulevelsetvalue})},
		\,
	\left\|\comdersmall^{[1,\Ntop-11];1} \wavearray \right\|_{L^{\infty}(\twoargroughtori{\timefunction,u}{\muxmulevelsetvalue})},
			  \\
	\left\| \Lunit \comder^{[1,6];2} \wavearray \right\|_{L^{\infty}(\twoargroughtori{\timefunction,u}{\muxmulevelsetvalue})},
		\,
	\left\| \comdersmall^{[1,6];2} \wavearray \right\|_{L^{\infty}(\twoargroughtori{\timefunction,u}{\muxmulevelsetvalue})},
		&   \\
	\left\| \Lunit \comder^{[1,5];3} \wavearray \right\|_{L^{\infty}(\twoargroughtori{\timefunction,u}{\muxmulevelsetvalue})},
		\,
	\left\| \comdersmall^{[1,5];3} \wavearray \right\|_{L^{\infty}(\twoargroughtori{\timefunction,u}{\muxmulevelsetvalue})},
	&  \\
	\left\| \Lunit \muX \muX \muX \muX \wavearray \right\|_{L^{\infty}(\twoargroughtori{\timefunction,u}{\muxmulevelsetvalue})}
	& \leq 
	\auxbootsmall,
	\end{split} 
		\\
	\label{BA:AUXVORTGRADENT} \tag{\textbf{AUX $(\vortrenormalized,\GradEnt)$ SMALL}}  
	\begin{split}
	\left\| 
		\comder^{\leq \Ntop-11;1} (\vortrenormalized,\GradEnt) 
	\right\|_{L^{\infty}(\twoargroughtori{\timefunction,u}{\muxmulevelsetvalue})},
	&		
		\\
	\left\| 
		\comder^{\leq 6;2} (\vortrenormalized,\GradEnt) 
	\right\|_{L^{\infty}(\twoargroughtori{\timefunction,u}{\muxmulevelsetvalue})},	
	&		\\
	\left\| 
		\comder^{\leq 5;3} (\vortrenormalized,\GradEnt) 
	\right\|_{L^{\infty}(\twoargroughtori{\timefunction,u}{\muxmulevelsetvalue})},	
	&		
		\\
	\left\| 
		\muX \muX \muX \muX (\vortrenormalized,\GradEnt) 
	\right\|_{L^{\infty}(\twoargroughtori{\timefunction,u}{\muxmulevelsetvalue})}
	& \leq \auxbootsmall,
	\end{split}		
			\\
	\label{BA:AUXMODIFIEDFLUID} \tag{\textbf{AUX $(\VortVort,\DivGradEnt)$ SMALL}} 
	\begin{split}
	\left\| 
		\comder^{\leq \Ntop-12;1} (\VortVort,\DivGradEnt) 
	\right\|_{L^{\infty}(\twoargroughtori{\timefunction,u}{\muxmulevelsetvalue})},
	&		
		 \\
	\left\| 
		\comder^{\leq 6;2} (\VortVort,\DivGradEnt)
	\right\|_{L^{\infty}(\twoargroughtori{\timefunction,u}{\muxmulevelsetvalue})},	
	&  
	\\
	\left\| 
		\comder^{\leq 5;3} (\VortVort,\DivGradEnt)
	\right\|_{L^{\infty}(\twoargroughtori{\timefunction,u}{\muxmulevelsetvalue})},	
	& 
		\\
	\left\| 
		\muX \muX \muX \muX (\VortVort,\DivGradEnt)
	\right\|_{L^{\infty}(\twoargroughtori{\timefunction,u}{\muxmulevelsetvalue})}
	& \leq \auxbootsmall,
	\end{split}
			\\
	\label{BA:AUXMUSMALL} \tag{\textbf{AUX $\upmu$ SMALL}} 
	\begin{split}
	\left\| \Lunit \tander^{[1,\Ntop-12]} \upmu \right\|_{L^{\infty}(\twoargroughtori{\timefunction,u}{\muxmulevelsetvalue})}, 
		\,
	\left\| \tander_*^{[1,\Ntop-12]} \upmu \right\|_{L^{\infty}(\twoargroughtori{\timefunction,u}{\muxmulevelsetvalue})},
	& 
	 \\
	\left\| \Lunit \comdersmall^{[1,5];1}  \upmu \right\|_{L^{\infty}(\twoargroughtori{\timefunction,u}{\muxmulevelsetvalue})}, 
		\,
	\left\| \comderdoublesmall^{[1,5];1} \upmu \right\|_{L^{\infty}(\twoargroughtori{\timefunction,u}{\muxmulevelsetvalue})},
	&	
		\\
	\left\| \Lunit \comdersmall^{[1,4];2}  \upmu \right\|_{L^{\infty}(\twoargroughtori{\timefunction,u}{\muxmulevelsetvalue})}, 
		\,
	\left\| \comderdoublesmall^{[1,4];2} \upmu \right\|_{L^{\infty}(\twoargroughtori{\timefunction,u}{\muxmulevelsetvalue})}
	& \leq \auxbootsmall,
	\end{split}	
		\\
	\left\| \Lsmall^1 \right\|_{L^{\infty}(\twoargroughtori{\timefunction,u}{\muxmulevelsetvalue})} 
	& 
	\leq \mathring{\upalpha}^{1/2},
		\label{BA:AUXL1SMALL} \tag{\textbf{AUX $\Lsmall^1$ SMALL}}
		\\
	\left\|\Lsmall^A \right\|_{L^{\infty}(\twoargroughtori{\timefunction,u}{\muxmulevelsetvalue})}
	\label{BA:AUXLASMALL} \tag{\textbf{AUX $\Lsmall^A$ SMALL}} 
	& 
	\leq \auxbootsmall
				\\
	\label{BA:AUXTANGENTIALDERIVATIVESLISMALL} \tag{\textbf{AUX $\tander \Lsmall^i$ SMALL}} 
	\begin{split} 
	\left\| \Lunit \tander^{\leq \Ntop-11} \Lsmall^i \right\|_{L^{\infty}(\twoargroughtori{\timefunction,u}{\muxmulevelsetvalue})}, 
			\,
	\left\| \tander^{[1,\Ntop-11]} \Lsmall^i \right\|_{L^{\infty}(\twoargroughtori{\timefunction,u}{\muxmulevelsetvalue})}, 
			\\
	\left\| \Lunit \comder^{[1,\Ntop-12];1} \Lsmall^i \right\|_{L^{\infty}(\twoargroughtori{\timefunction,u}{\muxmulevelsetvalue})}, 
		\,
	\left\| \comdersmall^{[1,\Ntop-12];1} \Lsmall^i \right\|_{L^{\infty}(\twoargroughtori{\timefunction,u}{\muxmulevelsetvalue})}, 
			\\
	\left\| \Lunit \comder^{[1,5];2}\Lsmall^i \right\|_{L^{\infty}(\twoargroughtori{\timefunction,u}{\muxmulevelsetvalue})}, 
		\,
	\left\| \comdersmall^{[1,5];2} \Lsmall^i \right\|_{L^{\infty}(\twoargroughtori{\timefunction,u}{\muxmulevelsetvalue})},
		\\\
	\left\| \Lunit \comder^{[1,4];3}\Lsmall^i \right\|_{L^{\infty}(\twoargroughtori{\timefunction,u}{\muxmulevelsetvalue})}, 
		\,
	\left\| \comdersmall^{[1,4];3} \Lsmall^i \right\|_{L^{\infty}(\twoargroughtori{\timefunction,u}{\muxmulevelsetvalue})}
	& 
	\leq \auxbootsmall.
\end{split}
\end{align}

\medskip

\noindent \underline{\textbf{Auxiliary bootstrap assumptions tied to pure transversal derivatives}}. 
We assume that the following inequalities hold for 
$(\timefunction,u) \in [\timefunction_0,\timefunctionboot) \times [- \rightu,\leftu]$:
\begin{align}
	\left 
		\| \muX^M \RRiemann
	\right\|_{L^{\infty}(\twoargroughtori{\timefunction,u}{\muxmulevelsetvalue})} 
	& 
	\leq 
	\left\| 
		\muX^M \RRiemann 
	\right\|_{L^{\infty}(\twoargroughtori{\timefunction_0,u}{\muxmulevelsetvalue})} 
	+ 
	\auxbootsmall, 
	\label{BA:AUXTRANSVERSALPDERIVATIVESRRIEMANNLARGE} \tag{\textbf{AUX $\muX^M \RRiemann$ LARGE}} 
	&
	1 \leq M \leq 4,
		\\
	\left 
		\| \muX^M \wavearraypartial 
	\right\|_{L^{\infty}(\twoargroughtori{\timefunction,u}{\muxmulevelsetvalue})} 
	& 
	\leq 
	\auxbootsmall, \label{BA:AUXTRANSVERSALPDERIVATIVESPARTIALWAVEARRAYSMALL}   
	\tag{\textbf{AUX $\muX^M \wavearraypartial$ SMALL}} 
	&
	1 \leq M \leq 4,
		\\
	\left\| 
		\muX^M \Lsmall^1 
	\right\|_{L^{\infty}(\twoargroughtori{\timefunction,u}{\muxmulevelsetvalue})} 
	& 
	\leq	 
	\left\| 
		\muX^M \Lsmall^1 
	\right\|_{L^{\infty}(\twoargroughtori{\timefunction_0,u}{\muxmulevelsetvalue})}
	+ 
	\auxbootsmall, 
		\label{BA:AUXTRANSVERSALPDERIVATIVESL1LARGE}
		\tag{\textbf{AUX $\muX^{[1,3]}\Lsmall^1$ LARGE}} 
	&
	1 \leq M \leq 3,
		\\
	\left\| 
		\muX^M \Lsmall^A 
	\right\|_{L^{\infty}(\twoargroughtori{\timefunction,u}{\muxmulevelsetvalue})} 
	& \leq 
	\auxbootsmall, 
	\label{BA:AUXTRANSVERSALPDERIVATIVESLASMALL} \tag{\textbf{AUX $\muX^{[1,3]}\Lsmall^A$ SMALL}} 
	&
	1 \leq M \leq 3,
		\\
	\left\| 
		\muX^M \upmu 
	\right\|_{L^{\infty}(\twoargroughtori{\timefunction,u}{\muxmulevelsetvalue})} 
	& 
	\leq 
	\left\| 
		\muX^M \upmu 
	\right\|_{L^{\infty}(\twoargroughtori{\timefunction_0,u}{\muxmulevelsetvalue})} 
	+
	\left\| 
		\argLrough{\muxmulevelsetvalue} \muX^M \upmu 
	\right\|_{L^{\infty}(\twoargroughtori{\timefunction_0,u}{\muxmulevelsetvalue})} 
	+ 
	\auxbootsmall, 
		\label{BA:AUXMUSMALL} \tag{\textbf{AUX $\muX^{\le 3}\upmu$ LARGE}} 
	&
	0 \leq M \leq 3,
		\\
	\left\|
		\Lunit \muX^M \upmu 
	\right\|_{L^{\infty}(\twoargroughtori{\timefunction,u}{\muxmulevelsetvalue})}
	& 
	\leq 
	\left\| 
		\argLrough{\muxmulevelsetvalue} \muX^M \upmu 
	\right\|_{L^{\infty}(\twoargroughtori{\timefunction_0,u}{\muxmulevelsetvalue})} 
	+ 
	\auxbootsmall, \label{BA:AUXLMULARGE} 
	\tag{\textbf{AUX $\Lunit \muX^{\le 3}\upmu$ LARGE}}
	&
	0 \leq M \leq 3,
\end{align}
where $\argLrough{\muxmulevelsetvalue}$ is defined in \eqref{E:LROUGH}.

\subsection{Key running assumptions}
\label{SS:RUNNINGASSUMPTIONS}
From now until Sect.\,\ref{S:EXISTENCEUPTOSINGULARBOUNDARYATFIXEDKAPPA},
we will often silently use the parameter-size and initial data assumptions 
of Sects.\,\ref{SS:PARAMETERSIZEASSUMPTIONS}
and \ref{SS:ASSUMPTIONSONDATA}
and the bootstrap assumptions of 
Sects.\,\ref{SS:BOOTSTRAPSCAFFOLDING} and \ref{SS:MAINQUANTITATIVEBOOTSTRAPASSUMPTIONS}.
In particular, when we state lemmas, propositions, and corollaries,
we will not explicitly restate these assumptions.

\subsection{A summary of the forthcoming derivation of improvements of the bootstrap assumptions}
\label{SS:SUMMARYOFIMPROVEMENTOFBOOTSTRAPASSUMPTIONS}
In the subsequent sections, we will derive strict improvements of all the bootstrap assumptions
that we made throughout Sect.\,\ref{S:BOOTSTRAPEVERYTHINGEXCEPTENERGIES}.
For the reader's convenience, here we state all the forthcoming results that yield the desired
strict improvements. Here we clarify that by ``strict improvements,'' we mean one or more of the following
three things:
\begin{enumerate}
	\item (\textbf{Quantitative improvement}) By this, we mean 
		that some quantity $Q$ 
		was assumed to satisfy $A_1 \leq Q \leq A_2$ in the bootstrap assumptions
		(where $A_1,A_2$ are real numbers),
		and we derive the improved bound $B_1 \leq Q \leq B_2$, 
		where $A_1 < B_1 \leq B_2 < A_2$.
		\item (\textbf{From soft to quantitative}) By this, we mean that in the bootstrap assumptions,
		we assumed that some function $f$ belongs to some function space and has a finite norm in that space,
		and our improvement is a quantitative estimate for the norm of $f$.
	\item (\textbf{Extension to the closure})
		By this, we mean that our bootstrap assumptions involved an assumption on the ``open-at-the-top'' domain
			$\twoargMrough{[\timefunction_0,\timefunctionboot),[- \rightu,\leftu]}{\muxmulevelsetvalue}$,
			and we derive an improved result showing that the assumption holds
			on the closed domain $\twoargMrough{[\timefunction_0,\timefunctionboot],[- \rightu,\leftu]}{\muxmulevelsetvalue}$.
\end{enumerate}

Here are the precise spots in the article where we derive improvements of the bootstrap assumptions.
\begin{itemize}
	\item Regarding the bootstrap assumptions of Sect.\,\ref{SSS:BAFORINVERSEFOLIATIONDENSITY}:
		we derive improvements of \eqref{E:BAMUPOSITIVE} in \eqref{E:MINVALUEOFMUONFOLIATION},
		of \eqref{E:BABOUNDSONLMUINTERESTINGREGION} in \eqref{E:BOUNDSONLMUINTERESTINGREGION},
		of \eqref{E:BABOUNDSONGEOMETRICTDERIVATIVEMUINTERESTINGREGION} in \eqref{E:BOUNDSONGEOMETRICTDERIVATIVEMUINTERESTINGREGION}, 
		and of
		\eqref{E:BAMUTRANSVERSALCONVEXITY} in \eqref{E:MUTRANSVERSALCONVEXITY}.
	\item We derive improvements of the bootstrap assumptions of Sect.\,\ref{SSS:BOOTSTRAPROUGHTIMEFUNCTION}
		in Lemma~\ref{L:DIFFEOMORPHICEXTENSIONOFROUGHCOORDINATES}.
	\item We derive improvements of the bootstrap assumptions of Sect.\,\ref{SSS:BACHOVMAPS}
				in Lemmas~\ref{L:DIFFEOMORPHICEXTENSIONOFROUGHCOORDINATES}
				and
				\ref{L:CHOVFROMROUGHCOORDINATESTOMUWEGIGHTEDXMUCOORDINATES}
				and in Prop.\,\ref{P:HOMEOMORPHICANDDIFFEOMORPHICEXTENSIONOFCARTESIANCOORDINATES}.	
	\item Regarding the bootstrap assumptions of Sect.\,\ref{SSS:BOOTSTRAPASSUMPTIONFORTORISTRUCTURE}:
		we derive improvements of Items 1 and 4 of 
		in Cor.\,\ref{C:QUANTITATIVECONTROLOFEMBEDDINGSONCLOSURESOFTHEIRDOMAINS},
		of \eqref{E:BOOSTRAPTORILOCATION} in \eqref{E:IMPROVEDLEVELSETSTRUCTUREANDLOCATIONOFMIN},
		and of \eqref{E:BOOTSTRAPLEVELSETSTRUCTUREANDLOCATIONOFMUXEQUALSMINUSKAPPA} in \eqref{E:MUXMUKAPPALEVELSETLOCATION}.
	\item We derive improvements of the bootstrap assumptions of Sect.\,\ref{SSS:BOOTSTRAPASSUMPTIONSSIZEOFCARTESIANTANDX1}
		in Lemma~\ref{L:CONTROLOFCARTESIANTANDX1}.
	\item We derive improvements of the bootstrap assumptions of Sect.\,\ref{SSS:SOFTBACONCERNINGREGULARITY}
			in Lemma~\ref{L:CONTINUOUSEXTNESION}.		
	\item We derive improvements of the fundamental quantitative bootstrap assumptions of Sect.\,\ref{SSS:FUNDAMENTALQUANTITATIVE}
		in Prop.\,\ref{P:IMPROVEMENTOFFUNDAMENTALQUANTITATIVEBOOTSTRAPASSUMPTIONS}.
	\item We derive improvements of the auxiliary bootstrap assumptions of Sect.\,\ref{SSS:AUXBOOTSTRAP}
		in Prop.\,\ref{P:IMPROVEMENTOFAUXILIARYBOOTSTRAP}.
	\item In Sect.\,\ref{SS:BOOTSTRAPASSUMPTIONSFORTHEWAVEENERGIES}, 
		we state bootstrap assumptions for the 
		$L^2$-type energies for the wave variables $\wavearray$. 
		We derive improvements of them in Prop.\,\ref{P:APRIORIL2ESTIMATESWAVEVARIABLES}.
\end{itemize}

\section{Preliminary pointwise, commutator, and differential operator comparison estimates}
\label{S:PRELIMINARYPOINTWISECOMMUTATORANDOPERATORCOMPARISON}
In this section, we use the bootstrap assumptions to derive several preliminary estimates,
including pointwise estimates,
commutator identities and estimates, 
and differential operator estimates comparing $\angLie$ and $\newangD$. 

\subsection{The norm of the $\ell_{t,u}$-tangent commutation vectorfields and simple comparison estimates}
\label{SS:NORMOFSMOOTHTORITANGENTCOMMUTATORSANDSIMPLECOMPARISON}

\begin{lemma}[The norm of the $\ell_{t,u}$-tangent commutation vectorfields and simple comparison estimates]
\label{L:NORMOFSMOOTHTORITANGENTCOMMUTATORSANDSIMPLECOMPARISON}
The $\ell_{t,u}$-tangent commutation vectorfields $\Angularset = \lbrace \Yvf{2}, \Yvf{3} \rbrace$ 
satisfy the following pointwise estimates on $\twoargMrough{[\timefunction_0,\timefunctionboot),[- \rightu,\leftu]}{\muxmulevelsetvalue}$:
\begin{align} \label{E:POINTWISESEMINORMOFYVECTORFIELDS}
	|\Yvf{A}|_{\gtorus} 
	& 
	= 
	1 
	+ 
	\mathcal{O}_{\mydiam}(\mathring{\upalpha}^{1/2}).
\end{align}	

Moreover, for any type $\binom{m}{n}$ 
$\ell_{t,u}$-tangent tensorfield $\upxi_{\beta_1 \cdots \beta_n}^{\alpha_1 \cdots \alpha_m}$,
the following pointwise estimates hold on $\twoargMrough{[\timefunction_0,\timefunctionboot),[- \rightu,\leftu]}{\muxmulevelsetvalue}$,
where $\Angularset$ is defined in \eqref{E:COMMUTATIONVECTORFIELDS}:
\begin{align} \label{E:SMOOTHTORUSNORMCOMPARBLETOTANGENTIALCONTRACTIONS}
\begin{split}	
	|\upxi|_{\gtorus}^2
	& 
	=
	\left\lbrace
		1 + \mathcal{O}_{\mydiam}(\mathring{\upalpha}^{1/2})
	\right\rbrace
	\sum_{\substack{ {^{(1)}U}, \cdots, {^{(m)}U} \in \Angularset \\ {^{(1)}V}, \cdots, {^{(n)}V} \in \Angularset}} 
	\left|
		{^{(1)}U_{\alpha_1}} \cdots {^{(m)}U_{\alpha_m}} 
		{^{(1)}V^{\beta_1}} \cdots {^{(n)}V^{\beta_n}}  
		\upxi_{\beta_1 \cdots \beta_n}^{\alpha_1 \cdots \alpha_n}
	\right|^2
		\\
	&
	=
	\left\lbrace
		1 + \mathcal{O}_{\mydiam}(\mathring{\upalpha}^{1/2})
	\right\rbrace
	\sum_{\substack{A_1,\cdots,A_m =2,3 \\ B_1,\cdots,B_n =2,3}}
	\left|
		\upxi_{B_1 \cdots B_n}^{A_1 \cdots A_n}
	\right|^2.
\end{split}
\end{align}
\end{lemma}
\begin{proof}
We first prove \eqref{E:SMOOTHTORUSNORMCOMPARBLETOTANGENTIALCONTRACTIONS}. 
We only prove the result for one-forms $\upxi$ because
arbitrary type $\binom{m}{n}$ 
$\ell_{t,u}$-tangent tensorfields
can be handled via similar arguments.
Let $\upxi$ be an $\ell_{t,u}$-tangent one form. 
Then by \eqref{E:SMOOTHGINVERSEABEXPRESSION}, 
we have that $|\upxi|_{\gtorus}^2 = \Speed^2\{ (\upxi_2)^2 + (\upxi_3)^2\} - (X^A \upxi_A)^2$. 
Next, using the bootstrap assumptions and Lemma~\ref{L:SCHEMATICSTRUCTUREOFVARIOUSTENSORSINTERMSOFCONTROLVARS},
we deduce that 
$\Speed -1, \, X^2, \, X^3 = \mathcal{O}_{\mydiam}(\mathring{\upalpha}^{1/2})$
and thus 
$|\upxi|_{\gtorus}^2 = 
\left\lbrace 1 
	+ 
\mathcal{O}_{\mydiam}(\mathring{\upalpha}^{1/2})
\right\rbrace 
\left\lbrace 
	(\upxi_2)^2 + (\upxi_3)^2
\right\rbrace
$.
Using in addition \eqref{E:GEOP2TOCOMMUTATORS}--\eqref{E:GEOP3TOCOMMUTATORS} 
and the fact that $X^1 + 1 = \Xsmall^1 = \mathcal{O}_{\mydiam}(\mathring{\upalpha}^{1/2})$,
we have that 
$\upxi_A 
=
\upxi \cdot \Yvf{A} 
+ 
\mathcal{O}_{\mydiam}(\mathring{\upalpha}^{1/2}) 
\left\lbrace
	|\upxi \cdot \Yvf{2}|
	+
	|\upxi \cdot \Yvf{3}|
\right\rbrace
$. 
Combining these results, we conclude
\eqref{E:SMOOTHTORUSNORMCOMPARBLETOTANGENTIALCONTRACTIONS} for one-forms $\upxi$.

To prove \eqref{E:POINTWISESEMINORMOFYVECTORFIELDS}, 
we first use \eqref{E:SMOOTHTORUSNORMCOMPARBLETOTANGENTIALCONTRACTIONS}
to deduce that 
$|\Yvf{A}|_{\gtorus}^2 = 
\left\lbrace 1 
	+ 
\mathcal{O}_{\mydiam}(\mathring{\upalpha}^{1/2})
\right\rbrace 
\left\lbrace 
	(\Yvf{A}^2)^2 + (\Yvf{A}^3)^2
\right\rbrace
$.
Moreover, the identities
\eqref{E:Y2INTERMSOFGEOMETRICCOORDINATEVECTORFIELDS}--\eqref{E:Y3INTERMSOFGEOMETRICCOORDINATEVECTORFIELDS}
and the estimates noted in the previous paragraph yield
$
(\Yvf{A}^2)^2 + (\Yvf{A}^3)^2
=
1
+
\mathcal{O}_{\mydiam}(\mathring{\upalpha}^{1/2})
$.
Combining these estimates, we conclude \eqref{E:POINTWISESEMINORMOFYVECTORFIELDS}.

\end{proof}

\subsection{Basic facts that we use silently when deriving estimates}
\label{SS:SILENTFACTS}
In the rest of the paper, we silently use the following basic facts.
\begin{enumerate}
\item All of the estimates we derive hold on the bootstrap region 
	$\twoargMrough{[\timefunction_0,\timefunctionboot),[- \rightu,\leftu]}{\muxmulevelsetvalue}$. 
		Moreover, in deriving estimates, 
		we often rely on the parameter- and data-size assumptions
		of Sects.\,\ref{SS:PARAMETERSIZEASSUMPTIONS}
		and \ref{SS:ASSUMPTIONSONDATA}
		and the bootstrap assumptions of 
		Sects.\,\ref{SS:BOOTSTRAPSCAFFOLDING} and \ref{SS:MAINQUANTITATIVEBOOTSTRAPASSUMPTIONS}.
	\item All quantities that we estimate can be controlled in terms of the variables 
	$\badcontrolvars = \{\wavearray,\upmu-1,\Lsmall^1,\Lsmall^2,\Lsmall^3\}$, $\vortrenormalized$, and $\GradEnt$
	(though we also the modified fluid variables $\VortVort$ and $\DivGradEnt$ to help control $\vortrenormalized$ and $\GradEnt$).
\item We use the Leibniz rule for the operators $\angLie_Z$ and $\newangD$ when deriving pointwise estimates for the $\angLie_Z$ and $\newangD$ derivatives of tensor products of the schematic form $\prod_{i=1}^m \upxi_{(i)}$, where the $\upxi_{(i)}$ are scalar functions or $\ell_{t,u}$-tangent tensors. Our derivative counts are such that all the $\upxi_{(i)}$ except at most one are uniformly bounded in $L^{\infty}$ on $\twoargMrough{[\timefunction_0,\timefunctionboot),[- \rightu,\leftu]}{\muxmulevelsetvalue}$. Thus, our pointwise estimates often explicitly feature (on the right-hand sides) only one factor with many derivatives on it, multiplied by a constant that uniformly bounds the other factors. In some estimates, the right-hand sides also gain smallness factor such as $\fundbootsmall$, generated by the remaining 
$\upxi_{(i)}'s$. 
\item We use the conventions for constants $C, \mathfrak{c}$, and $C_{\mydiam}$
	stated in Sect.\,\ref{SS:CONVENTIONSFORCONSTANTS}.
\item We use the conventions for strings of commutation vectorfields stated in Sect.\,\ref{SS:STRINGSOFCOMMUTATIONVECTORFIELDS}.
\item We use the comparison estimates of Lemma~\ref{L:NORMOFSMOOTHTORITANGENTCOMMUTATORSANDSIMPLECOMPARISON}.
\end{enumerate}

\subsection{Pointwise estimates for Cartesian components of geometric vectorfields}
\label{SS:POINTWISEFORCARTESIANCOMPONENTSOFGEOMETRICVECTORFIELDS}
In this section, we provide simple pointwise estimates for the Cartesian components 
of the vectorfields $\{\Lunit,\muX, \Yvf{2},\Yvf{3}\}$ and their derivatives.

\begin{lemma}[Pointwise estimates for $x^i$ and the Cartesian components of the vectorfields 
$\{\Lunit,\muX, \Yvf{2},\Yvf{3}\}$] \label{L:POINTWISEESTIMATESFORCARTESIANCOMPONENTSOFVECTORFIELDS} 
Let $\Singletan \in \{\Lunit,\Yvf{2},\Yvf{3}\}$. For $i = 1,2,3$,
the following pointwise estimates hold on $\twoargMrough{[\timefunction_0,\timefunctionboot),[- \rightu,\leftu]}{\muxmulevelsetvalue}$:
\begin{subequations} 
\begin{align} 
|\Singletan^i| 
	& \lesssim 1 + |\controlvars|, 
	\label{E:TANGENTIALCOMMUTATORCARTESIANCOMPONENTPOINTWISEESTIMATE} 
	\\
|\tander^{[1,N]} \Singletan^i| 
& \lesssim  
| \tander^{[1,N]} \controlvars|, 
	\label{E:DERIVATIVESOFCARTESIANCOMPONENTSOFCOMMUTATORS} 
	\\
|\comdersmall^{[1,N];1} \Singletan^i| 
& \lesssim  
	|\comdersmall^{[1,N];1} \controlvars|,
	\label{E:COMDERSMALLPCARTESIANCOMPEST} 
	\\ 
|\comder^{[1,N];1} \Singletan^i| 
	& \lesssim  
	|\comder^{[1,N];1} \controlvars|, 
	\label{E:COMDERPCARTESIANCOMPEST} 
	\\ 
|\muX^i| 
& \lesssim 
1 
+ 
|\badcontrolvars|, 
	\label{E:MUXCARTESIANCOMPEST} 
	\\
|\tander^{[1,N]} \muX^i| 
	& \lesssim 
		|\tander^{[1,N]} \badcontrolvars|, 
		\label{E:DERMUXCARTESIANCOMPEST} 
	\\
|\comdersmall^{[1,N];1} \muX^i| 
& \lesssim 
	|\comdersmall^{[1,N];1} \badcontrolvars|, 
	\label{E:COMDERSMALLMUXCARTESIANCOMPEST}
	\\
|\comder^{[1,N];1} \muX^i| 
& \lesssim 
|\comder^{[1,N];1} \badcontrolvars|, 
	\label{E:COMDERMUXCARTESIANCOMPEST}
	\\
|\angrmD x^i|_{\gtorus} 
& \lesssim 
	1 
	+ 
	|\controlvars|, 
	\label{E:ANGDCARTESIANCOORDINATEIPOINTWISESTIMATE} 
	\\
|\angrmD \tander^{[1,N]} x^i |_{\gtorus} 
	& \lesssim 
		|\tander^{[1,N]} \controlvars|, 
		\label{E:ANGDTANDERXIEST} 
	\\
|\angrmD \comdersmall^{[1,N];1} x^i|_{\gtorus} 
	& \lesssim  
	|\comdersmall^{[1,N];1} \controlvars| 
	+ 
	|\tandersmall^{[1,N]} \badcontrolvars|, 
		\label{E:ANGDCOMDERSMALLXIEST} 
			\\
|\angrmD \comdersmall^{[1,N];1} x^i|_{\gtorus} 
	& \lesssim  
		| \tandersmall^{[1,N]} \controlvars  
		+ 
		|\tandersmall^{[1,N]} \badcontrolvars |. 
		\label{E:ANGDCOMDERxiEST}
\end{align}
\end{subequations}
\end{lemma}
 
\begin{proof} 
	The same proof of \cite[Lemma 8.4]{jLjS2018} holds with minor modifications accounting for the third dimension.
\end{proof}


\subsection{Pointwise estimates for various $\ell_{t,u}$-tangent tensorfields} 
In this section, we record several pointwise estimates of various $\ell_{t,u}$-tangent tensorfields. 

\begin{lemma}[Crude pointwise estimates for the Lie derivatives of $\gtorus$ and $\upchi$] 
\label{L:CRUDEPOINTWISEESTIMATESFORTENSORFIELDS} 
The following pointwise estimates hold on $\twoargMrough{[\timefunction_0,\timefunctionboot),[-\rightu,\leftu]}{\muxmulevelsetvalue}$:
\begin{subequations} 
\begin{align}  
	|\angLie_{\tander}^{N+1} \gtorus|_{\gtorus}, 
		\, 
	|\angLie_{\tander}^{N+1} \gtorus^{-1}|_{\gtorus},
		\, 
	|\angLie_{\tander}^N \upchi|_{\gtorus}, 
		\, 
	|\tander^N \mytr_{\gtorus}\upchi| 
	& \lesssim  
		|\tander^{[1,N+1]}\controlvars|, 
		\label{E:TANDERGANDCHIESTIMATE} 
			\\
	|\angLie_{\comdersmall}^{N+1;1} \gtorus|_{\gtorus},
		\,  
	|\angLie_{\comdersmall}^{N+1;1} \gtorus^{-1}|_{\gtorus},
		\,  
	|\angLie_{\comder}^{N;1} \upchi|_{\gtorus}, 
		\, 
	|\comder^{N;1} \mytr_{\gtorus}\upchi| 
	& \lesssim  
	|\comdersmall^{[1,N+1];1}\controlvars| 
	+ 
	|\tandersmall^{[1,N+1]} \badcontrolvars|, 
		\label{E:COMDERSMALLGANDCHIESTIMATE} 
		\\
	|\angLie_{\comder}^{N+1;1} \gtorus|_{\gtorus},
		\, 
	 |\angLie_{\comder}^{N+1;1} \gtorus^{-1}|_{\gtorus} 
	& \lesssim  
	|\comder^{[1,N+1];1}\controlvars| 
	+ 
	|\tandersmall^{[1,N+1]} \badcontrolvars|. 
		\label{E:COMDERGESTIMATE}
\end{align}
\end{subequations}
\end{lemma}
\begin{proof}
The same proof of \cite[Lemma 8.5]{jLjS2018} holds with minor modifications to account for the third spatial dimension.
\end{proof}


\subsection{Commutator identities and estimates}
\label{SS:COMMUTATORIDENTITIESANDESITMATES}
\begin{lemma}\cite[Simple commutator identities; Lemma 8.5]{jLjS2021}
	\label{L:SIMPLECOMMUTATORIDENTITY}
	For any $\nullhyparg{u}$-tangent vectorfields 
	$\Singletan, \Singletan_1, \Singletan_2 \in \{\Lunit,\Yvf{2},\Yvf{3}\}$, 
	the commutators $[\Singletan_1,\Singletan_2]$ and $[\muX,\Singletan]$ are $\ell_{t,u}$-tangent.
	Moreover, we have the following identities, where $\smoothfunction$
	schematically denotes a smooth function of its arguments:
	\begin{subequations}
	\begin{align} \label{E:TANGENTIALSIMPLECOMMUTATORIDENTITY}
		[\Singletan_1, \Singletan_2] 
		& = 
		\smoothfunction(\tander^{\le 1}\controlvars) \Yvf{2} 
		+
		\smoothfunction(\tander^{\le 1}\controlvars) \Yvf{3}.
	\end{align}
		
	For each $\Singletan \in \{\Lunit,\Yvf{2},\Yvf{3}\}$,
	there exist smooth functions, all schematically denoted by ``$\smoothfunction$,''
	such that the following identity holds:
	\begin{align} \label{E:TRANSVERSALTANGENTIALSIMPLECOMMUTATORIDENTITY}
	[\Singletan, \muX] 
	& = 
	\smoothfunction(\tander^{\le 1} \badcontrolvars,\muX \wavearray) \Yvf{2}
	+
	\smoothfunction(\tander^{\le 1} \badcontrolvars,\muX \wavearray) \Yvf{3}.
	\end{align}
	\end{subequations}
\end{lemma}

The following proposition provides various vectorfield commutator 
estimates that we use in our analysis.

\begin{proposition}[Pointwise commutator estimates] 
\label{P:COMMUTATORESTIMATES}
Let $1 \leq N \leq \Ntop$ be an integer, and let $\varphi$ be a scalar function. 
For any $\Singletan \in \Tanset$ (see definition \eqref{E:COMMUTATIONVECTORFIELDS}),
iterated commutators can be bounded pointwise as follows
on $\twoargMrough{[\timefunction_0,\timefunctionboot),[- \rightu,\leftu]}{\muxmulevelsetvalue}$:
\begin{subequations}
\begin{align} \label{E:COMMUTATOROFTANGENTIALANDTANGENTIALCOMMUTATORS}
\left|
	[\Singletan, \tander^N] \varphi
\right|
& \lesssim 
\auxbootsmall \left|\tander^{[1,N]} \varphi \right| 
+ 
\underbrace{\sum_{\substack {N_1+N_2 \leq N+1 \\ N_1,\,N_2\leq N}} 
\left|\tander^{[2,N_1]}\controlvars \right| \left|\tander^{[1,N_2]} \varphi \right|}_{\mbox{Absent if $N=1$}},
	\\
\begin{split}
\left|[\muX, \tander^N]\varphi \right|,
	\,
\left|
	[\Singletan, \comder^{N;1}] \varphi
\right|
 & \lesssim
 \left|\tander^{[1,N]} \varphi \right| 
	+ 
	\underbrace{\sum_{\substack {N_1+N_2 \leq N+1 \\ N_1,\,N_2\leq N}} 
	\left|\tander^{[2,N_1]} \badcontrolvars \right| \left|\tander^{[1,N_2]} \varphi \right|}_{\mbox{Absent if $N=1$}}
	\label{E:COMMUTATOROFMUXANDTANGENTIALCOMMUTATORS} 
	\\
 & \ \
	+ 
	\underbrace{\sum_{\substack {N_1+N_2 \leq N \\ N_1 \leq N-1}} 
	\left|\tander^{[1,N_1]} \muX \wavearray \right| \left|\tander^{[1,N_2]} \varphi \right|}_{\mbox{Absent if $N=1$}}.
\end{split}
\end{align}
\end{subequations}

In particular, 
for any $\Singletan \in \{\Lunit,\Yvf{2},\Yvf{3}\}$,
we have the following pointwise estimates
on $\twoargMrough{[\timefunction_0,\timefunctionboot),[- \rightu,\leftu]}{\muxmulevelsetvalue}$:
\begin{subequations}
\begin{align} \label{E:POINTWISEBOUNDCOMMUTATORSTANGENTIALANDTANGENTIALDERIVATIVESONSCALARFUNCTION}
\left|[\Singletan, \tander^N] \varphi \right| 
& 
\lesssim 
\auxbootsmall \left| \tander^{[1,N]}\varphi \right|,
&
\mbox{if } 1 \leq N \leq \Ntop - 11,
	\\
\left|[\muX, \tander^N]\varphi \right|,
	\,
\left|
	[\Singletan, \comder^{N;1}] \varphi
\right|
& 
\lesssim \left|\tander^{[1,N]} \varphi \right|,
& 
\mbox{if } 1 \leq N \leq \Ntop - 12.
\label{E:POINTWISEBOUNDCOMMUTATORSMUXANDTANGENTIALDERIVATIVESONSCALARFUNCTION}
\end{align}
\end{subequations}

Moreover, the following pointwise estimates hold:
\begin{subequations}
\begin{align} 
\left|
	[\Singletan, \comder^{N;2}] \varphi
\right|,
	\,
\left|
	[\muX, \comder^{N;1}] \varphi 
\right|
& 
\lesssim 
\left| \comdersmall^{[1,N];1} \varphi \right|,
&
\mbox{if } 1 \leq N \leq 5,
	 \label{E:POINTWISEBOUNDCOMMUTATORSUPTOTWOMUXANDTANGENTIALDERIVATIVESONSCALARFUNCTION}
	\\
\left|
	[\Singletan, \comder^{N;3}] \varphi
\right|,
	\,
\left|
	[\muX, \comder^{N;2}] \varphi 
\right|
& 
\lesssim 
\left| \comdersmall^{[1,N];2} \varphi \right|,
&
\mbox{if } 1 \leq N \leq 4.
 \label{E:POINTWISEBOUNDCOMMUTATORSUPTOTHREEMUXANDTANGENTIALDERIVATIVESONSCALARFUNCTION}
\end{align}
\end{subequations}

Finally, if $\upxi$ is an $\ell_{t,u}$-tangent one-form or an $\ell_{t,u}$-tangent type $\binom{0}{2}$-tensorfield 
and $\Singletan \in \{\Lunit,\Yvf{2},\Yvf{3}\}$,
then the following estimates hold:
\begin{align} \label{E:LANDANGLIECOMMUTATOR}
\left|
	[\angLie_{\Singletan}, \angLie_{\tander}^N] \upxi 
\right|_{\gtorus} 
&
\lesssim 
\auxbootsmall \left| \angLie_{\tander}^{\leq N} \upxi \right|_{\gtorus}  
+ 
\underbrace{\sum_{\substack{N_1+N_2 \leq N+1 \\ N_2\leq N}} 
\left|\tander^{[2,N_1]} \controlvars \right| 
\left|\angLie_{\tander}^{[1,N_2]}  \upxi \right|_{\gtorus}}_{\mbox{Absent if $N=1$}}.
\end{align}
\end{proposition}

\begin{proof}
We first prove \eqref{E:COMMUTATOROFTANGENTIALANDTANGENTIALCOMMUTATORS}.
Iterating \eqref{E:TANGENTIALSIMPLECOMMUTATORIDENTITY} and using 
Lemma~\ref{L:SCHEMATICSTRUCTUREOFVARIOUSTENSORSINTERMSOFCONTROLVARS}
and the bootstrap assumptions, we find that
$
\left|
	[\Singletan, \tander^N] \varphi
\right|
\lesssim
\sum_{\substack {N_1+N_2 \leq N+1 \\ N_1,\, N_2 \leq N}} 
\left|\tander^{[1,N_1]} \controlvars \right| \left|\tander^{[1,N_2]} \varphi \right|
$.
From this bound and the fact that the bootstrap assumptions imply that $|\tander \controlvars| \lesssim \auxbootsmall$,
we arrive at \eqref{E:COMMUTATOROFTANGENTIALANDTANGENTIALCOMMUTATORS}.

To prove \eqref{E:COMMUTATOROFMUXANDTANGENTIALCOMMUTATORS} for 
$\left|[\muX, \tander^N]\varphi \right|$,
we use a similar argument that also relies on
\eqref{E:TRANSVERSALTANGENTIALSIMPLECOMMUTATORIDENTITY} to deduce
$
\left|[\muX, \tander^N]\varphi \right|
\lesssim
 \sum_{\substack {N_1+N_2 \leq N+1 \\ N_1,\, N_2 \leq N}} 
	\left|\tander^{[1,N_1]} \badcontrolvars \right| \left|\tander^{[1,N_2]} \varphi \right|
+
	\sum_{\substack {N_1+N_2 \leq N \\ N_1 \leq N-1}} 
	\left|\tander^{\leq N_1} \muX \wavearray \right| \left|\tander^{[1,N_2]} \varphi \right|	
$.
From this bound and the fact that the bootstrap assumptions imply that $|\tander \badcontrolvars| \lesssim 1$
and $\left|\muX \wavearray \right| \lesssim 1$,
we conclude the desired bounds for $\left|[\muX, \tander^N]\varphi \right|$.
The estimate \eqref{E:COMMUTATOROFMUXANDTANGENTIALCOMMUTATORS} for
$
\left|
	[\Singletan, \comder^{N;1}] \varphi
\right|
$ can be proved through a similar argument, and we omit the details.

Except for \eqref{E:LANDANGLIECOMMUTATOR},
the remaining estimates in the proposition follow from similar arguments that take into
account the details of the auxiliary bootstrap assumptions of Sect.\,\ref{SSS:AUXBOOTSTRAP},
in particular the $L^{\infty}$ regularity of the solution variables with respect to $\muX$  
and $\nullhyparg{u}$-tangential differentiations.

We now prove \eqref{E:LANDANGLIECOMMUTATOR}.
We first consider the case in which $\upxi = \upxi_A \angrmD x^A$ is an $\ell_{t,u}$-tangent one-form,
where by our usual conventions, $\upxi_A \eqdef \upxi \cdot \geop{x^A}$.
By Lemma~\ref{L:ANGULARDIFFERENTIALCOMMUTESWITHANGLIE} and the Leibniz rule,
for any $\Singletan \in \Tanset$, we have:
\begin{align} \label{E:IDENTIFYFORANGLIEDIFFERENTIATEDMONEFORM}
	\angLie_{\Singletan} \upxi 
	& = 
	(\Singletan \upxi_A) \angrmD x^A 
	+ 
	\upxi_A \angrmD (\Singletan x^A).
\end{align}
Differentiating \eqref{E:IDENTIFYFORANGLIEDIFFERENTIATEDMONEFORM} with $\ell_{t,u}$-projected Lie derivatives,
using the Leibniz rule, using Lemma~\ref{L:ANGULARDIFFERENTIALCOMMUTESWITHANGLIE},
using the bootstrap assumptions, 
using 
\eqref{E:COMMUTATOROFTANGENTIALANDTANGENTIALCOMMUTATORS}
with the scalar functions $\upxi_A$ in the role of $\varphi$,
and using \eqref{E:ANGDCARTESIANCOORDINATEIPOINTWISESTIMATE}--\eqref{E:ANGDTANDERXIEST},
we see that 
$
\left|
[\angLie_{\Singletan}, \angLie_{\tander}^N] \upxi
\right|_{\gtorus}
\lesssim 
\sum_{A=2,3}
\auxbootsmall 
\left|\tander^{[1,N]} \upxi_A \right| 
+ 
\sum_{A=2,3}
\sum_{\substack {N_1+N_2 \leq N+1 \\ N_2 \leq N}} 
\left|\tander^{[2,N_1]}\controlvars \right| \left|\tander^{\leq N_2} \upxi_A \right|
$.
Moreover, by considering the coordinate components of the one-forms on each side of \eqref{E:IDENTIFYFORANGLIEDIFFERENTIATEDMONEFORM},
using the bootstrap assumptions,
and using \eqref{E:ANGDTANDERXIEST},
we find that 
$\sum_{A=2,3} |\Singletan \upxi_A| \lesssim |\angLie_{\Singletan} \upxi|_{\gtorus} + |\upxi|_{\gtorus}$.
Differentiating \eqref{E:IDENTIFYFORANGLIEDIFFERENTIATEDMONEFORM} and using similar arguments,
we use induction in $M$ to deduce that for $1 \leq M \leq \Ntop$, 
we have
$\sum_{A=2,3} |\tander^M \upxi_A| 
\lesssim 
\left|\angLie_{\tander}^{\leq M} \upxi \right|_{\gtorus}
+
\sum_{\substack {M_1+M_2 \leq M \\ M_2 \leq M-1}} 
\left|\tander^{[1,M_1]} \controlvars \right| 
\left|\angLie_{\tander}^{\leq M_2} \upxi \right|_{\gtorus}$.
Combining the estimates and using the bootstrap assumptions, 
we conclude \eqref{E:LANDANGLIECOMMUTATOR}
for $\ell_{t,u}$-tangent one-forms $\upxi$. 
For any $n \geq 2$, the estimate \eqref{E:LANDANGLIECOMMUTATOR}
for type $\binom{0}{n}$ $\ell_{t,u}$-tangent tensorfields $\upxi$ can be proved using similar arguments.
To prove \eqref{E:LANDANGLIECOMMUTATOR} for a general type $\binom{m}{n}$ 
$\ell_{t,u}$-tangent tensorfield $\upxi$, we use $\gtorus$-duality to 
express $\upxi$ as the dual of a type $\binom{0}{m+n}$ tensorfield $\upeta$,
use \eqref{E:LANDANGLIECOMMUTATOR} to deduce a commutator estimate for $\upeta$,
and then use the bootstrap assumptions and  
the estimates \eqref{E:TANDERGANDCHIESTIMATE} for $\ell_{t,u}$-projected Lie
derivatives of $\gtorus$ and $\gtorus^{-1}$ to obtain the estimates for 
$\upxi$ as a consequence of the ones for $\upeta$.

\end{proof}


\subsection{Differential operator estimates comparing $\angLie$ and $\newangD$}
In this section, we provide several pointwise estimates comparing different differential operators. 

\begin{lemma}[Differential operator pointwise comparison estimates]
Let $f$ be a scalar function on $\twoargMrough{[\timefunction_0,\timefunctionboot),[- \rightu,\leftu]}{\muxmulevelsetvalue}$. 
Then the following pointwise estimates hold:
\begin{subequations}
\begin{align}
	|\angrmD f|_{\gtorus}^2 
	& = \left\lbrace
				1 + \mathcal{O}_{\mydiam}(\mathring{\upalpha}^{1/2}) 
		\right\rbrace
		\sum_{A = 2}^3 |\Yvf{A} f|^2, 
		\label{E:ANGDFPOINTWISEBOUNDEDBYCOMMUTATORVECTORFIELDS} 
			\\
	|\angLap f|^2  
	& \leq 
	2(1 + C_{\mydiam} \mathring{\upalpha}^{1/2}) 
	\sum_{A = 2}^3 |\angrmD \Yvf{A} f|_{\gtorus}^2 
	+ 
	C \auxbootsmall |\angrmD f|_{\gtorus}^2. 
	\label{E:SMOOTHANGULARHESSIANOFFPOINTWISEBOUNDEDBYCOMMUTATORVECTORFIELDS} 
\end{align}
\end{subequations}
\end{lemma}

\begin{proof}
\eqref{E:ANGDFPOINTWISEBOUNDEDBYCOMMUTATORVECTORFIELDS}
follows from \eqref{E:SMOOTHTORUSNORMCOMPARBLETOTANGENTIALCONTRACTIONS}
with $\upxi \eqdef \angrmD f$.

We now prove \eqref{E:SMOOTHANGULARHESSIANOFFPOINTWISEBOUNDEDBYCOMMUTATORVECTORFIELDS}. 
We first note that \eqref{E:SMOOTHTORUSNORMCOMPARBLETOTANGENTIALCONTRACTIONS} and the Leibniz rule
yield the following estimate for any scalar function $f$:
\begin{align} \label{E:ANGHESSFESTIMATEINTERMEDIATE}
	|\newangDsquared f|_{\gtorus}^2
	& \leq 
			\left\lbrace
				1 + \mathcal{O}_{\mydiam}(\mathring{\upalpha}^{1/2}) 
			\right\rbrace
		\sum_{A = 2}^3 |\angrmD \Yvf{A} f|_{\gtorus}^2
		+
		C
		\sum_{A,B,C = 2}^3 
		|\gfour(\Dfour_{\Yvf{A}} \Yvf{B},\Yvf{C})|^2
		|\angrmD f|_{\gtorus}^2.
\end{align}
Next, we compute relative to the Cartesian coordinates and use 
Lemma~\ref{L:SCHEMATICSTRUCTUREOFVARIOUSTENSORSINTERMSOFCONTROLVARS} 
and the bootstrap assumptions
to deduce that schematically,
we have 
$|\gfour(\Dfour_{\Yvf{A}} \Yvf{B},\Yvf{C})| 
= 
|\smoothfunction(\controlvars) \cdot \Angularset \controlvars| 
\lesssim \auxbootsmall$. Inserting this bound into RHS~\eqref{E:ANGHESSFESTIMATEINTERMEDIATE} 
and using the inequality $|\angLap f|^2 \leq 2 |\newangD^2 f|_{\gtorus}^2$, 
we arrive at the desired estimate \eqref{E:SMOOTHANGULARHESSIANOFFPOINTWISEBOUNDEDBYCOMMUTATORVECTORFIELDS}.

\end{proof}


\subsection{Transport inequalities for the eikonal function quantities} \label{SS:TRANSPORTINEQUALITIES}
In this section, we provide transport inequalities satisfied by the eikonal function quantities 
$\upmu, \Lsmall^i, \upchi$, and $\mytr_{\gtorus}\upchi$.
 We also provide pointwise estimates for the differentiated quantity $\angLie_{\Lunit} \angLie_{\tander}^{N-1} \upchi$. 
These estimates \emph{involve a loss of one order of differentiability} relative to $\wavearray$ 
in the sense that the right-hand sides of the transport equations that we use to derive the inequalities
depend on the first-order derivatives of $\wavearray$. 
In Sect.\,\ref{S:PRELIMINARYL2ESTIMATESFORBELOWTOPORDERDERIVATIVESOFACOUSTICGEOMETRYANDDERIVATIVELOSING}, 
we use the transport inequalities to derive below-top-order energy estimates for the eikonal function quantities.

\begin{proposition}[Transport inequalities for the eikonal function quantities] 
\label{P:POINTWISETRANSPORTINEQUALITIESFOREIKFUNCTIONQUANTITIES} 
The following pointwise estimates hold on $\twoargMrough{[\timefunction_0,\timefunctionboot),[- \rightu,\leftu]}{\muxmulevelsetvalue}$:
\begin{subequations}
\begin{align}
	|\Lunit \upmu| 
	& 
	\lesssim 
	|\comder \wavearray|, 
	&&	\label{E:LMUPOINTWISE} 
		\\
	|\Lunit \tandersmall^N \upmu|, 
		\,
	|\tandersmall^N \Lunit \upmu| 
	& 
	\lesssim 
	|\comdersmall^{[1,N+1];1}\wavearray| 
	+
	|\tander^{[1,N]}\controlvars| 
	+ 
	\auxbootsmall 
	|\tandersmall^{[1,N]} \badcontrolvars|,
	&& \mbox{if } 1 \leq N \leq \Ntop,
	\label{E:LTANGENTIALMUPOINTWISE} 
		\\
	|\Lunit \tander^N \Lsmall^i|, 
		\,
	|\tander^N \Lunit \Lsmall^i| 
	& 
	\lesssim 
	|\tander^{[1,N+1]}\wavearray| 
	+ 
	\auxbootsmall 
	|\tander^{[1,N]}\controlvars|, 
	\label{E:LUNITTANGENTIALDERIVATIVESOFLUNITIPOINTWISE} 
	&& \mbox{if } 0 \leq N \leq \Ntop, \\
	|\Lunit \tander^{N-1} \mytr_{\gtorus} \upchi|, 
		\,
	|\tander^{N-1} \Lunit \mytr_{\gtorus}\upchi| 
	& 
	\lesssim 
	|\tander^{[1,N+1]}\wavearray| 
	+ 
	\auxbootsmall 
	|\tander^{[1,N]}\controlvars|,  
		\label{E:LTANGENTIALTRCHIPOINTWISE} 
	&& \mbox{if } 1 \leq N \leq \Ntop, \\
	|\angLie_{\Lunit} \angLie_{\tander}^{N-1} \upchi|_{\gtorus}, 
		\,
	|\angLie_{\tander}^{N-1}\angLie_{\Lunit} \upchi|_{\gtorus} 
	& 
	\lesssim |\tander^{[1,N+1]}\wavearray| 
	+ 
	\auxbootsmall |\tander^{[1,N]}\controlvars|,
	&& \mbox{if } 1 \leq N \leq \Ntop,
	\label{E:ANGLIELTANGENTIALCHIPOINTWISE} 
	\\
	|\Lunit \comder^{N;1} \Lsmall^i|, 
		\,
	|\comder^{N;1} \Lunit \Lsmall^i| 
	& 
	\lesssim 
	|\comdersmall^{[1,N+1];1} \wavearray| 
	+ 
	|\comdersmall^{[1,N];1} \controlvars|
	+ 
	\auxbootsmall
	|\tandersmall^{[1,N]} \badcontrolvars|, 
	&& \mbox{if } 1 \leq N \leq \Ntop,
		\label{E:LZLSMALLPOINTWISE} 
		\\
	|\Lunit \comder^{N-1;1} \mytr_{\gtorus} \upchi|, 
		\,
	|\comder^{N-1;1} \Lunit \mytr_{\gtorus}\upchi| 
	& 
	\lesssim 
	|\comdersmall^{[1,N+1];1} \wavearray| 
	+ 
	|\comdersmall^{[1,N];1} \controlvars| 
	+ 
	\auxbootsmall
	|\tandersmall^{[1,N]} \badcontrolvars|,  
	&& \mbox{if } 2 \leq N \leq \Ntop,
	\label{E:LZTRCHIPOINTWISE} 
		\\
	|\angLie_{\Lunit} \angLie_{\comder}^{N-1;1} \upchi|_{\gtorus}, 
		\,
	|\angLie_{\comder}^{N-1;1}\angLie_{\Lunit} \upchi|_{\gtorus} 
	& 
	\lesssim 
	|\comdersmall^{[1,N+1];1} \wavearray| 
	+ 
	| \comdersmall^{[1,N];1} \controlvars| 
	+ 
	\auxbootsmall 
	|\tandersmall^{[1,N]}\badcontrolvars|,
	&& \mbox{if } 2 \leq N \leq \Ntop.
	\label{E:ANGLIELZCHIPOINTWISE}
\end{align}
\end{subequations}
\end{proposition}

\begin{proof}
Thanks to the transport equations of Lemma~\ref{L:TRANSPORTMUANDLUNITI},
the identity \eqref{E:TRCHIEXPRESSSIONINTERMSOFDERIVATIVESOFLUNITI},
and the commutator estimates of Prop.\,\ref{P:COMMUTATORESTIMATES},
the estimates \eqref{E:LMUPOINTWISE}--\eqref{E:LTANGENTIALTRCHIPOINTWISE}, \eqref{E:LZLSMALLPOINTWISE},
and \eqref{E:LZTRCHIPOINTWISE}
follow from the same arguments given in \cite[Proposition 8.13]{jLjS2018},
and we omit the details. 
We now prove \eqref{E:ANGLIELTANGENTIALCHIPOINTWISE} for $\angLie_{\tander}^{N-1} \angLie_{\Lunit} \upchi$. Since angular Lie differentiation commutes with $\angrmD$ (see \eqref{E:ANGULARDIFFERENTIALCOMMUTESWITHANGLIE}), 
equation \eqref{E:CHIEXPRESSSIONINTERMSOFDERIVATIVESOFLUNITI} implies
\[ \angLie_{\Lunit} \upchi = (\vec{G}_{ab} \diamond \Lunit \wavearray)  \angrmD \Lunit^a \otimes \angrmD x^b + \gfour_{ab}\angrmD \Lunit^a \otimes \angrmD L^b + \gfour_{ab}\angrmD L \Lunit^a \otimes \angrmD x^b + \smoothfunction(\controlvars) PP\wavearray.\]
Substituting RHS~\eqref{E:LUNITITRANSPORT} for $\Lunit \Lunit^a$ and taking $\angLie_{\tander}^{N-1}$ derivatives of the resulting expression proves \eqref{E:CHIEXPRESSSIONINTERMSOFDERIVATIVESOFLUNITI} for $\angLie_{\tander}^{N-1}\angLie_{\Lunit} \upchi$. The estimate
\eqref{E:ANGLIELTANGENTIALCHIPOINTWISE}  
for $\angLie_{\Lunit} \angLie_{\tander}^{N-1}\upchi$ then follows from the commutator estimate 
\eqref{E:LANDANGLIECOMMUTATOR}, the crude estimate \eqref{E:TANDERGANDCHIESTIMATE}, 
and the bootstrap assumptions.
With the help of the crude estimate \eqref{E:COMDERSMALLGANDCHIESTIMATE},
the estimate \eqref{E:ANGLIELZCHIPOINTWISE} follows from similar arguments, and we omit the details.
\end{proof}

\subsection{Pointwise commutator estimates for $\upchi$ tied to a Codazzi-type identity}
\label{SS:CODAZZICOMMUTATOR} 

\begin{lemma}[Codazzi-type identity for $\upchi$]
	\label{L:CODAZZITYPEIDENTITY}
	There exist smooth functions, all schematically denoted by ``$\smoothfunction$,''
	such that the following identity holds:
	\begin{align} \label{E:CODAZZITYPEIDENTITY}
	\angdiv \upchi
	-
	\angrmD \mytr_{\gtorus} \upchi
	& 
	= 
	\smoothfunction(\tander^{\leq 1} \controlvars,\angrmD \vec{x}) \tander \controlvars
	+
	\smoothfunction(\controlvars,\angrmD \vec{x}) \tander^2 \wavearray.
\end{align}
\end{lemma}

\begin{proof}
We apply $\angdiv$ \eqref{E:CHIEXPRESSSIONINTERMSOFDERIVATIVESOFLUNITI}
and
$\angrmD$ to \eqref{E:TRCHIEXPRESSSIONINTERMSOFDERIVATIVESOFLUNITI}
and note that the RHSs of the resulting equations are
$\ell_{t,u}$-tangent one-forms such that the 
terms involving the second-order derivatives of $\Lunit$ 
have the components
$\gfour_{ab} (\gtorus^{-1})^{BC} [\newangDsquareddoublearg{A}{B} \Lunit^a] \argangrmD{C} x^b$,
i.e., the second-order-in-$\Lunit$ terms agree.
Also using
Lemmas~\ref{L:SCHEMATICSTRUCTUREOFVARIOUSTENSORSINTERMSOFCONTROLVARS}
and \ref{L:SCHEMATICEXPRESSIONFORANGULARLAPLACIAN},
we conclude \eqref{E:CODAZZITYPEIDENTITY}.
\end{proof}

\begin{lemma}[Pointwise commutator estimates for $\upchi$ tied to a Codazzi-type identity]
\label{L:CODAZZICOMMUTATORESTIMATES}
Let $1 \leq N \leq \Ntop$. Then the following pointwise estimates 
hold on $\twoargMrough{[\timefunction_0,\timefunctionboot),[- \rightu,\leftu]}{\muxmulevelsetvalue}$,
where on LHS~\eqref{E:CODAZZICOMMUTATORESTIMATES},
$\tander^{N-1}$ denotes the same order $N-1$ string of commutation vectorfields
in each of the two terms: 
\begin{align} \label{E:CODAZZICOMMUTATORESTIMATES} 
\left| 
	\angdiv \angLie_{\tander}^{N-1} \upchi 
	- 
	\angrmD \tander^{N-1} \mytr_{\gtorus} \upchi 
\right|_{\gtorus} 
& 
\lesssim 
\left| 
	\tander^{[1,N+1]} \wavearray 
\right| 
+ 
\left| 
	\tander^{[1,N]} \controlvars 
\right|.
\end{align}

\end{lemma}

\begin{proof}
		We apply $\angLie_{\tander}^{N-1}$ to \eqref{E:CODAZZITYPEIDENTITY}.
		By Lemma~\ref{L:ANGULARDIFFERENTIALCOMMUTESWITHANGLIE}, 
		the operator $\angLie_{\tander}^{N-1}$
		commutes under the operators $\angrmD$ on each side of
		\eqref{E:CODAZZITYPEIDENTITY}.
		This yields the main term $\angrmD \tander^{N-1} \mytr_{\gtorus} \upchi$
		on LHS~\eqref{E:CODAZZICOMMUTATORESTIMATES},
		while the bootstrap assumptions
		and Lemma~\ref{L:SCHEMATICEXPRESSIONFORANGULARLAPLACIAN}
		yield that
		$
		\left|
			\angLie_{\tander}^{N-1} \mbox{RHS~\eqref{E:CODAZZITYPEIDENTITY}}
		\right|_{\gtorus}
		\lesssim
		\mbox{RHS~\eqref{E:CODAZZICOMMUTATORESTIMATES}}
		$
		as desired.
		Hence, to complete the proof, we must show that
		$
		\left|
			\angLie_{\tander}^{N-1} \angdiv \upchi 
			-
			\angdiv \angLie_{\tander}^{N-1} \upchi 
		\right|_{\gtorus}
		\lesssim
		\mbox{RHS~\eqref{E:CODAZZICOMMUTATORESTIMATES}}
		$.
		To this end, we apply
		$\angLie_{\tander}^{N-1} \angdiv$ 
		and
		$\angdiv \angLie_{\tander}^{N-1}$
		to \eqref{E:CHIEXPRESSSIONINTERMSOFDERIVATIVESOFLUNITI},
		thereby obtaining expressions
		for
		$
		\angLie_{\tander}^{N-1} \angdiv \upchi 
		$
		and
		$	
		\angdiv \angLie_{\tander}^{N-1} \upchi 
		$
		respectively. 
		As above,
		Lemmas~\ref{L:ANGULARDIFFERENTIALCOMMUTESWITHANGLIE},
		\ref{L:COMMUTATORSTOCOORDINATES}
		and \ref{L:SCHEMATICEXPRESSIONFORANGULARLAPLACIAN}
		and the bootstrap assumptions
		imply that all terms except the principal ones, i.e., the ones involving the order $N+1$ derivatives of $\Lunit^i$,
		are bounded in the norm $| \cdot |_{\gtorus}$ by
		$
		\lesssim
		\mbox{RHS~\eqref{E:CODAZZICOMMUTATORESTIMATES}}$.
		The principal terms in 
		$
		\angLie_{\tander}^{N-1} \angdiv \upchi 
		$
		and
		$	
		\angdiv \angLie_{\tander}^{N-1} \upchi 
		$
		are respectively
		$\gfour_{ab} (\tander^{N-1} \angLap \Lunit^a) \angrmD x^b$
		and
		$\gfour_{ab} (\angLap \tander^{N-1} \Lunit^a) \angrmD x^b$.
		Hence, using
		Lemma~\ref{L:SCHEMATICEXPRESSIONFORANGULARLAPLACIAN},
		the commutator estimate \eqref{E:COMMUTATOROFTANGENTIALANDTANGENTIALCOMMUTATORS},
		and the bootstrap assumptions,
		we see that the principal terms in the difference 
		$
		\gfour_{ab} (\tander^{N-1} \angLap \Lunit^a) \angrmD x^b
		-
		\gfour_{ab} (\angLap \tander^{N-1} \Lunit^a) \angrmD x^b
		$
		cancel and that
		$
		\left|
			\angLie_{\tander}^{N-1} \angdiv \upchi 
			-
			\angdiv \angLie_{\tander}^{N-1} \upchi 
		\right|_{\gtorus}
		\lesssim
		\mbox{RHS~\eqref{E:CODAZZICOMMUTATORESTIMATES}}
		$
		as desired.
\end{proof}

\subsection{Pointwise estimates for the inhomogeneous terms in the commuted equations}
\label{SS:POINTWISEESTIMATESFORINHOMOGENEOUSTERMSINCOMMUTEDEQUATIONS}

\subsubsection{Pointwise estimates for the derivatives of the null forms}
\label{SSS:POINTWISEESTIMATESFORDERIVATIVESOFNULLFORMS}

\begin{lemma}[Pointwise estimates for the derivatives of the null forms]
\label{L:POINTWISEESTIMATESFORDERIVATIVESOFNULLFORMS}
Let $N \leq \Ntop$. The $\tander^N$-derivatives of the product of
$\upmu$ and the terms defined in \eqref{E:TRANSPORTVORTVORTMAINTERMS}--\eqref{E:DIVENTROPYGRADIENTNULLFORM} 
satisfy the following pointwise estimates
on $\twoargMrough{[\timefunction_0,\timefunctionboot),[- \rightu,\leftu]}{\muxmulevelsetvalue}$:
\begin{subequations}
\begin{align} 
\begin{split} \label{E:POINTWISEESIMATESFORDERIVATIVESOFNULLFORMSTRUCTUREMODIFIEDFLUIDVARIABLES}
\left|\tander^N (\upmu \mainnullform_{(\VortVort)}^i) \right|, 
	\,
\left|\tander^N (\upmu \mainnullform_{(\DivGradEnt)}) \right| 
	&
	\lesssim
	\left|
		\tander^{\leq N+1} (\vortrenormalized,\GradEnt)
	\right|
		\\
& \ \
	+
	\fundbootsmall
	\left|
		\muX \tander^{[1,N]} \wavearray
	\right| 
	+  
	\fundbootsmall
	\left|\tander^{[1,N+1]} \wavearray \right| 
	+ 
	\fundbootsmall
	\left|\tandersmall^{[1,N]} \badcontrolvars \right|,
\end{split}
		\\
\left|\tander^N (\upmu \nullform_{(v)}^i) \right|,
	\, 
\left|\tander^N (\upmu \nullform_{(\pm)}) \right|, 
	\, 
\left|\tander^N (\upmu \nullform_{(\LogDensity)}) \right| 
& 
\lesssim 
\fundbootsmall
\left|
	\muX \tander^{[1,N]} \wavearray
\right| 
+  
\left|
	\tander^{[1,N+1]} \wavearray 
\right| 
+ 
\fundbootsmall
\left|
	\tandersmall^{[1,N]} \badcontrolvars 
\right|,
		\label{E:POINTWISEESIMATESFORDERIVATIVESOFNULLFORMSTRUCTUREWAVEVARIABLES} 
		\\
\begin{split} \label{E:POINTWISEESIMATESFORDERIVATIVESOFNULLFORMSTRUCTURETRANSPORTVARIABLES}
\left| 
	\tander^N (\upmu \nullform_{(\VortVort)}^i) 
\right|, 
	\, 
\left|
	\tander^N (\upmu \nullform_{(\DivGradEnt)}) 
\right|
& 
\lesssim 
	\fundbootsmall
	\left|
		\tander^{\leq N} \vortrenormalized
	\right|
		\\
& \ \
	+
	\fundbootsmall
	\left|
		\muX \tander^{[1,N]} \wavearray
	\right| 
	+  
	\fundbootsmall
	\left|
		\tander^{[1,N+1]} \wavearray 
	\right| 
	+ 
	\fundbootsmall
	\left|
		\tandersmall^{[1,N]} \badcontrolvars 
	\right|.
\end{split}
\end{align}
\end{subequations}

\end{lemma}

\begin{proof}
	We differentiate the identities of Lemma~\ref{L:CRUCIALSTRUCTUREOFNULLFORMS} with $\tander^N$,
	use the bootstrap assumptions,
	and use the commutator estimate \eqref{E:COMMUTATOROFMUXANDTANGENTIALCOMMUTATORS}
	so that on the RHSs of the estimates, all terms featuring any $\muX$-differentiation of 
	$\wavearray$ are such that the $\muX$ operator acts last.
\end{proof}

\subsubsection{Pointwise estimates for the derivatives of the linear inhomogeneous terms}
\label{SSS:POINTWISEESTIMATESFORDERIVATIVESOFLINEARINHOMOGENEOUSTERMS}

\begin{lemma}[Pointwise estimates for the derivatives of the linear inhomogeneous terms]
\label{L:POINTWISEESTIMATESFORDERIVATIVESOFLINEARINHOMOGENEOUSTERMS}
Let $N \leq \Ntop$. 
Consider the product of $\upmu$ and the terms $\VortVort, \DivGradEnt$-involving terms on
RHSs \eqref{E:VELOCITYWAVEEQUATION}--\eqref{E:ENTROPYWAVEEQUATION},
as well as the product of $\upmu$ and the terms defined in \eqref{E:VELOCITYILINEARORBETTER}--\eqref{E:RENORMALIZEDVORTICITYCURLLINEARORBETTER}.
Then the $\tander^N$-derivatives of these terms
satisfy the following pointwise estimates
on $\twoargMrough{[\timefunction_0,\timefunctionboot),[- \rightu,\leftu]}{\muxmulevelsetvalue}$:
\begin{subequations}
\begin{align} 
\begin{split}  \label{E:POINTWISEESTIMATESFORLINEARMODIFIEDFLUIDVARIABLETERMSONRHSWAVEEQUATIONS} 
	&
	\left|
		\tander^N
		\left(
			\upmu \Speed^2 \exp(2 \LogDensity) \VortVort^i,
			\,
			\upmu
			\left\lbrace 
				F_{;\Ent} c^2 \exp(2\LogDensity)
				- 
				\Speed \exp(\LogDensity) \frac{p_{;\Ent}}{\overline{\varrho}} 
			\right\rbrace  
			\DivGradEnt
		\right)
		\right|,
			\\
	&
	\left|
		\tander^N
		\left(
			\upmu
			\exp(\LogDensity) \frac{p_{;\Ent}}{\overline{\varrho}} \DivGradEnt,
				\,
			\upmu
			\Speed^2 \exp(2 \LogDensity) \DivGradEnt
		\right)
	\right|
		\\
	& \lesssim 
		\upmu 
			\left|
				\tander^N (\VortVort,\DivGradEnt)
			\right|
			+
			\left|
				\tander^{\leq N-1} (\VortVort,\DivGradEnt)
			\right|
			+
			\varepsilon
			\left|
				\tander^{[1,N]} \wavearray 
			\right|,
	\end{split}
		\\
	\begin{split} \label{E:POINTWISEESTIMATESFORALLLINEARTERMS} 
	&
	\left|
	\tander^N
		\left(
		\upmu \mathfrak{L}_{(v)}^i, 
			\,
		\upmu \mathfrak{L}_{(\pm)},
			\,
		\upmu \mathfrak{L}_{(\LogDensity)},
			\,
		\upmu \mathfrak{L}_{(\Ent)},
			\,
		\upmu \mathfrak{L}_{(\vortrenormalized)}^i,
			\,
		\upmu \mathfrak{L}_{(\GradEnt)}^i,
			\,
		\upmu \mathfrak{L}_{(\Flatdiv \vortrenormalized)},
			\,
		\upmu \mathfrak{L}_{(\VortVort)}^i
		\right)
	\right|
		\\
	& 
\lesssim 
\left|
	\tander^{\leq N} (\vortrenormalized,\GradEnt)
\right|
+
\fundbootsmall
\left|
	\muX \tander^{[1,N]} \wavearray
\right| 
+  
\left|
	\tander^{[1,N+1]} \wavearray 
\right| 
+ 
\fundbootsmall
\left|
	\tandersmall^{[1,N]} \badcontrolvars 
\right|.
\end{split}
\end{align}
\end{subequations}
\end{lemma}

\begin{proof}
	We apply the same reasoning used in the proof of Lemma~\ref{L:POINTWISEESTIMATESFORDERIVATIVESOFNULLFORMS}
	to the identities provided by Lemma~\ref{L:CRUCIALSTRUCTURELINEARINHOMOGENEOUS}.
\end{proof}

\subsubsection{Pointwise estimates for the derivatives of the inhomogeneous terms in the commuted wave equations}
\label{SSS:POINTWISEESTIMATESFORALLINHOMOGENEOUSTERMSINCOMMUTEDWAVEEQUATIONS}

\begin{corollary}[Pointwise estimates for the derivatives of the inhomogeneous terms in the commuted wave equations]
\label{C:POINTWISEESTIMATESFORDERIVATIVESOFINHOMOGENEOUSTERMSINWAVEEQUATIONS}
Let 
$\wavearray \eqdef (\Psi_0,\Psi_1,\Psi_2,\Psi_3,\Psi_4) \eqdef (\RRiemann,\LRiemann,v^2,v^3,\Ent)$ be the solutions to the covariant wave equations \eqref{E:COVARIANTWAVEEQUATIONSWAVEVARIABLES}. 
We denote the product of $\upmu$ and the RHS of the covariant wave equation satisfied by $\Psi_{\iota}$ by $\mathfrak{G}_{\iota}$, 
i.e., 
$\upmu \Box_{\gfour} \Psi_{\iota} = \mathfrak{G}_{\iota}$. 
Let $1 \leq N \leq \Ntop$. 
Then the following pointwise estimates hold
on $\twoargMrough{[\timefunction_0,\timefunctionboot),[- \rightu,\leftu]}{\muxmulevelsetvalue}$:
\begin{align} \label{E:POINTWISEESTIMATESFORALLINHOMOGENEOUSTERMS}
\begin{split}
\left|
	\tander^N \mathfrak{G}_{\iota}
\right|
& \lesssim
	\upmu 
	\left|
		\tander^N (\VortVort,\DivGradEnt)
	\right|
	+
	\left|
		\tander^{\leq N-1} (\VortVort,\DivGradEnt)
	\right|
+
\left|
	\tander^{\leq N} (\vortrenormalized,\GradEnt)
\right|
+
\fundbootsmall
\left|
	\muX \tander^{[1,N]} \wavearray
\right| 
+  
\left|
	\tander^{[1,N+1]} \wavearray 
\right| 
+ 
\fundbootsmall
\left|
	\tandersmall^{[1,N]} \badcontrolvars 
\right|.
\end{split}
\end{align}
	
\end{corollary}

\begin{proof}
	The corollary is a direct consequence of the estimates
	\eqref{E:POINTWISEESIMATESFORDERIVATIVESOFNULLFORMSTRUCTUREWAVEVARIABLES}
	and
	\eqref{E:POINTWISEESTIMATESFORLINEARMODIFIEDFLUIDVARIABLETERMSONRHSWAVEEQUATIONS}--\eqref{E:POINTWISEESTIMATESFORALLLINEARTERMS}.
	
\end{proof}

\section{Embeddings of \texorpdfstring{$\datahypfortimefunctiontwoarg{-\muxmulevelsetvalue}{[\timefunction_0,\timefunctionboot)}$}{the data hypersurface for the rough time function} 
and the flow map of  
\texorpdfstring{$\Wtransarg{\muxmulevelsetvalue}$}{W}}
\label{S:EMBEDDINGSANDFLOWMAPS}
In this section, we derive quantitative control over how the level sets 
$\datahypfortimefunctiontwoarg{-\muxmulevelsetvalue}{[\timefunction_0,\timefunctionboot)}$ 
(see definition \eqref{E:TRUNCATEDLEVELSETSOFMUXMU})
are embedded in the spacetime region $\twoargMrough{[\timefunction_0,\timefunctionboot),[- \rightu,\leftu]}{\muxmulevelsetvalue}$. We also derive quantitative control on the flow map of $\Wtransarg{\muxmulevelsetvalue}$. 
We will use these results in Sect.\,\ref{S:ESTIMATESFORROUGHTIMEFUNCTIONCONTINUOUSEXTENSIONSANDDIFFEOS}
to demonstrate the viability of the ``transversal initial value" problem 
\eqref{E:TRANSPORTEQUATIONFORROUGHTIMEFUNCTION}--\eqref{E:INITIALCONDITIONFORROUGHTIMEFUNCTION}
that we used to construct the rough time function $\timefunctionarg{\muxmulevelsetvalue}$;
see Remark~\ref{R:SOLUTIONISMOREREGULARTHANHYPERSURFACE} concerning some subtleties
tied to the regularity theory in the construction.

\subsection{Embedded submanifolds and quantitative control over the embeddings}
\label{SS:EMBEDDEDSUBMANIFOLDSANDCONTROLOVERTHEEMBEDDINGS	}

\begin{lemma}[Embedded submanifolds and quantitative control over the embeddings]
	\label{L:XMUISMINUSCAPPISAGRAPH}
	Let $\Cartesiantisafunctiononmumxtoriarg{\mulevelsetvalue}{-\muxmulevelsetvalue}(x^2,x^3)$ 
	and 
	$\Eikonalisafunctiononmumuxtoriarg{\mulevelsetvalue}{-\muxmulevelsetvalue}(x^2,x^3)$ be the functions on $\mathbb{T}^2$
	from Sect.\,\ref{SSS:BOOTSTRAPASSUMPTIONFORTORISTRUCTURE}.
	For $\mulevelsetvalue \in (\upmuboot,\mupositive]$, we have:
	\begin{subequations}
	\begin{align} \label{E:W2INFTYBOUNDFORTORI}
		\| \Cartesiantisafunctiononmumxtoriarg{\mulevelsetvalue}{-\muxmulevelsetvalue} \|_{W^{2,\infty}(\mathbb{T}^2)},
			\,
		\| \Eikonalisafunctiononmumuxtoriarg{\mulevelsetvalue}{-\muxmulevelsetvalue} \|_{W^{2,\infty}(\mathbb{T}^2)}
		& \leq C,
			\\
		\left\| \left(\geop{x^2} \Cartesiantisafunctiononmumxtoriarg{\mulevelsetvalue}{-\muxmulevelsetvalue},
			\geop{x^3} \Cartesiantisafunctiononmumxtoriarg{\mulevelsetvalue}{-\muxmulevelsetvalue} \right) \right\|_{W^{1,\infty}(\mathbb{T}^2)},
			\,
		\left\| \left(\geop{x^2} \Eikonalisafunctiononmumuxtoriarg{\mulevelsetvalue}{-\muxmulevelsetvalue},
			\geop{x^3} \Eikonalisafunctiononmumuxtoriarg{\mulevelsetvalue}{-\muxmulevelsetvalue} \right) \right\|_{W^{1,\infty}(\mathbb{T}^2)}
		& \leq C \varepsilon^{1/2}.
		\label{E:W1INFTYBOUNDFORANGULARDERIVATIVESTORI}
	\end{align}
	\end{subequations}

	Moreover,
	\begin{align} \label{E:FORMULAFORMUDERIVATIVEOFCARTESIANTIMEFUNCTIONOFMUXISMINUSKAPPAEMBEEDDING}
		\frac{\partial}{\partial \mulevelsetvalue} \Cartesiantisafunctiononmumxtoriarg{\mulevelsetvalue}{-\muxmulevelsetvalue}
		& = \frac{1}{\geop{t} \upmu - \frac{(\geop{u} \upmu) \geop{t} \muX \upmu}{\geop{u} \muX \upmu}},
	\end{align}
	and there is a $C > 1$ such that 
	following estimate holds:
	\begin{align} \label{E:IMPROVEMENTMUEQUALSMINUSKAPPAEMBEDDINGCARTESIANTIMEFUNCTIONNEGATIVEMUDERIVATIVE}
		- 
		C
		& 
		<
		\inf_{(\mulevelsetvalue,x^2,x^3) \in (\upmuboot,\mupositive] \times \mathbb{T}^2} 
		\frac{\partial}{\partial \mulevelsetvalue} \Cartesiantisafunctiononmumxtoriarg{\mulevelsetvalue}{-\muxmulevelsetvalue}(x^2,x^3) 
		\leq
		\sup_{(\mulevelsetvalue,x^2,x^3) \in (\upmuboot,\mupositive] \times \mathbb{T}^2} 
		\frac{\partial}{\partial \mulevelsetvalue} \Cartesiantisafunctiononmumxtoriarg{\mulevelsetvalue}{-\muxmulevelsetvalue}(x^2,x^3) 
		< 
		- \frac{1}{C}.
	\end{align}
	
	\medskip
	
	\noindent \underline{\textbf{$\datahypfortimefunctiontwoarg{-\muxmulevelsetvalue}{[\timefunction_0,\timefunctionboot)}$ is a graph}}.
	Let
	\begin{align} \label{E:DOMAINOFEMBEDDINGFORXMUEQUALSMINUSKAPPAYSURFACE}
	\domainforembeddingdatahypfortimefunctiontwoarg{\muxmulevelsetvalue}{(\upmuboot,\mupositive]}
	&
	\eqdef 
	\left\lbrace 
		(t,x^2,x^3) \in \mathbb{R} \times \mathbb{T}^2 
		\ | \
		\Cartesiantisafunctiononmumxtoriarg{\mupositive}{-\muxmulevelsetvalue}(x^2,x^3) \leq t 
		< 
		\Cartesiantisafunctiononmumxtoriarg{\upmuboot}{-\muxmulevelsetvalue}(x^2,x^3)
	\right\rbrace.
	\end{align}
	Then $\domainforembeddingdatahypfortimefunctiontwoarg{\muxmulevelsetvalue}{(\upmuboot,\mupositive]}$ is precompact,
	and there exists an embedding
	$\embeddingdatahypfortimefunctionarg{\muxmulevelsetvalue}: \domainforembeddingdatahypfortimefunctiontwoarg{\muxmulevelsetvalue}{(\upmuboot,\mupositive]} 
	\rightarrow \twoargMrough{[\timefunction_0,\timefunctionboot),[-\frac{3}{4} \interestingu,\frac{3}{4} \interestingu]}{\muxmulevelsetvalue}$
	of the form
	$\embeddingdatahypfortimefunctionarg{\muxmulevelsetvalue}(t,x^2,x^3) 
	=\left(t,\scalarembeddingdatahypfortimefunctionarg{\muxmulevelsetvalue}(t,x^2,x^3),x^2,x^3 \right)$
	such that $\embeddingdatahypfortimefunctionarg{\muxmulevelsetvalue}\in W^{2,\infty}(\mbox{\upshape int}(\domainforembeddingdatahypfortimefunctiontwoarg{\muxmulevelsetvalue}{(\upmuboot,\mupositive]}))$ 
	and such that $\embeddingdatahypfortimefunctionarg{\muxmulevelsetvalue}$
	is a diffeomorphism from $\domainforembeddingdatahypfortimefunctiontwoarg{\muxmulevelsetvalue}{(\upmuboot,\mupositive]}$ onto 
	$\datahypfortimefunctiontwoarg{-\muxmulevelsetvalue}{[\timefunction_0,\timefunctionboot)}$,
	where 
	\begin{align} \label{E:INTERIORDOMAINOFEMBEDDINGFORXMUEQUALSMINUSKAPPAYSURFACE}
	\mbox{\upshape int}(\domainforembeddingdatahypfortimefunctiontwoarg{\muxmulevelsetvalue}{(\upmuboot,\mupositive]}) 
	&
	\eqdef
		\left\lbrace 
		(t,x^2,x^3) \in \mathbb{R} \times \mathbb{T}^2 
		\ | \
		\Cartesiantisafunctiononmumxtoriarg{\mupositive}{-\muxmulevelsetvalue}(x^2,x^3) 
		< 
		t 
		< \Cartesiantisafunctiononmumxtoriarg{\upmuboot}{-\muxmulevelsetvalue}(x^2,x^3)
	\right\rbrace
	\end{align}
	is the interior of $\domainforembeddingdatahypfortimefunctiontwoarg{\muxmulevelsetvalue}{(\upmuboot,\mupositive]}$.
	In particular, relative to the geometric coordinates $(t,u,x^2,x^3)$, we have
	\begin{align} \label{E:XMUISMINUSCAPPISAGRAPH}
		\datahypfortimefunctiontwoarg{-\muxmulevelsetvalue}{[\timefunction_0,\timefunctionboot)}
		& = 
		\left\lbrace 
			\left(t,\scalarembeddingdatahypfortimefunctionarg{\muxmulevelsetvalue}(t,x^2,x^3),x^2,x^3 \right) 
			\in \twoargMrough{[\timefunction_0,\timefunctionboot),[- \rightu,\leftu]}{\muxmulevelsetvalue}
			\ | \
			(t,x^2,x^3) 
			\subset 
			\domainforembeddingdatahypfortimefunctiontwoarg{\muxmulevelsetvalue}{(\upmuboot,\mupositive]}
		\right\rbrace.
	\end{align}
	Moreover, $\embeddingdatahypfortimefunctionarg{\muxmulevelsetvalue}$ is $C^{1,1}$ on every compact subset of 
	$\domainforembeddingdatahypfortimefunctiontwoarg{\muxmulevelsetvalue}{(\upmuboot,\mupositive]}$,
	and
	the following estimates hold:
	\begin{subequations}
	\begin{align} \label{E:BOUNDONEMBEDDINGOFXMUEQUALSMINUSKAPPA}
		\left\| 
			\embeddingdatahypfortimefunctionarg{\muxmulevelsetvalue} 
		\right\|_{W^{2,\infty}(\mbox{\upshape int}(\domainforembeddingdatahypfortimefunctiontwoarg{\muxmulevelsetvalue}{(\upmuboot,\mupositive]}))}
		& \leq C,
			\\
		\left\| 	
			\left(\geop{x^2} \scalarembeddingdatahypfortimefunctionarg{\muxmulevelsetvalue},
			\geop{x^3} \scalarembeddingdatahypfortimefunctionarg{\muxmulevelsetvalue}\right) 
		\right\|_{W^{1,\infty}(\mbox{\upshape int}(\domainforembeddingdatahypfortimefunctiontwoarg{\muxmulevelsetvalue}{(\upmuboot,\mupositive]}))}
		& \leq C \varepsilon^{1/2}.
		\label{E:SMALLNESSBOUNDONEMBEDDINGOFXMUEQUALSMINUSKAPPA}
	\end{align}
	Finally, on $\domainforembeddingdatahypfortimefunctiontwoarg{\muxmulevelsetvalue}{(\upmuboot,\mupositive]}$, the following estimate holds:
	\begin{align} \label{E:PARTIALTEMBEDDINGFUNCTIONKEYESTIMATE}
		\geop{t} \scalarembeddingdatahypfortimefunctionarg{\muxmulevelsetvalue}
		& = - \frac{\Lunit \muX \upmu \circ \embeddingdatahypfortimefunctionarg{\muxmulevelsetvalue}}{\muX \muX \upmu 
			\circ 
			\embeddingdatahypfortimefunctionarg{\muxmulevelsetvalue}} 
			+ 
			\mathcal{O}(\varepsilon^{1/2}).
	\end{align}
	\end{subequations}
\end{lemma}

\begin{proof}
	Throughout this proof, we silently use the soft bootstrap assumptions of
	Sect.\,\ref{SSS:SOFTBACONCERNINGREGULARITY}, which guarantee our needed
	qualitative regularity. When proving quantitative estimates, 
	we will use the bootstrap assumptions of Sects.\,\ref{SSS:FUNDAMENTALQUANTITATIVE} and \ref{SSS:AUXBOOTSTRAP}.
	
	To derive the existence of the embedding
	of the form 
	$\embeddingdatahypfortimefunctionarg{\muxmulevelsetvalue}(t,x^2,x^3) 
	= \left(t,\scalarembeddingdatahypfortimefunctionarg{\muxmulevelsetvalue}(t,x^2,x^3),x^2,x^3 \right)$,
	we first note that
	by 
	\eqref{E:EMBEDDEDEDMUEQUALSMINUSKAPPAHYPERSURFACEFINITEC11NORMONCOMPACTSUBSETS},
	\eqref{E:BAMUEQUALSMINUSKAPPAEMBEDDINGCARTESIANTIMEFUNCTIONNEGATIVEMUDERIVATIVE}, 
	and the
	inverse function theorem,
	the map
	$(\mulevelsetvalue,x^2,x^3) \rightarrow \left(\Cartesiantisafunctiononmumxtoriarg{\mulevelsetvalue}{-\muxmulevelsetvalue}(x^2,x^3),x^2,x^3 \right)$
	on $(\upmuboot,\mupositive] \times \mathbb{T}^2$ is a global diffeomorphism
	onto the set $\domainforembeddingdatahypfortimefunctiontwoarg{\muxmulevelsetvalue}{(\upmuboot,\mupositive]}$
	such that the inverse function of the map is of the form
	$(t,x^2,x^3) \rightarrow \left(I_t(x^2,x^3),x^2,x^3 \right)$,
	where for any $\mulevelsetvalue' \in (\upmuboot,\mupositive]$,
	the map
	$(t,x^2,x^3) \rightarrow I_t(x^2,x^3)$ is $C^{1,1}$ 
	on 
	$\left\lbrace 
		(t,x^2,x^3) \in \mathbb{R} \times \mathbb{T}^2 
		\ | \
		\Cartesiantisafunctiononmumxtoriarg{\mupositive}{-\muxmulevelsetvalue}(x^2,x^3) \leq t \leq \Cartesiantisafunctiononmumxtoriarg{\mulevelsetvalue'}{-\muxmulevelsetvalue}(x^2,x^3)
	\right\rbrace$.
	Thus, 
	the desired embedding $\embeddingdatahypfortimefunctionarg{\muxmulevelsetvalue}$ 
	(into $\twoargMrough{[\timefunction_0,\timefunctionboot),[-\frac{3}{4} \interestingu,\frac{3}{4} \interestingu]}{\muxmulevelsetvalue}$)
	is the composition of
	the diffeomorphism $(t,x^2,x^3) \rightarrow \left(I_t(x^2,x^3),x^2,x^3 \right)$
	with the embedding
	$
	\embeddatahypersurfacearg{\muxmulevelsetvalue}
	$
	defined in \eqref{E:EMBEDDATAHYPERSURFACE}.
	
	We now prove
	\eqref{E:BOUNDONEMBEDDINGOFXMUEQUALSMINUSKAPPA}--\eqref{E:SMALLNESSBOUNDONEMBEDDINGOFXMUEQUALSMINUSKAPPA}
	and \eqref{E:PARTIALTEMBEDDINGFUNCTIONKEYESTIMATE}.
	The estimate 
	$\| \scalarembeddingdatahypfortimefunctionarg{\muxmulevelsetvalue} \|_{L^{\infty}(\mbox{\upshape int}(\domainforembeddingdatahypfortimefunctiontwoarg{\muxmulevelsetvalue}{(\upmuboot,\mupositive]}))} 
	\leq C$
	follows trivially from the fact that $u \in [- \rightu,\leftu]$ in
	$\twoargMrough{[\timefunction_0,\timefunctionboot),[- \rightu,\leftu]}{\muxmulevelsetvalue}$.
	To control the derivatives of
	$\scalarembeddingdatahypfortimefunctionarg{\muxmulevelsetvalue}$,
	we differentiate the equation
	$[\muX \upmu]\left(t,\scalarembeddingdatahypfortimefunctionarg{\muxmulevelsetvalue}(t,x^2,x^3),x^2,x^3 \right) = - \muxmulevelsetvalue$
	(which holds since $\muX \upmu|_{\datahypfortimefunctiontwoarg{-\muxmulevelsetvalue}{[\timefunction_0,\timefunctionboot)}} = - \muxmulevelsetvalue$)
	with $\geop{t},\geop{x^2},\geop{x^3}$,
	use the chain rule to algebraically solve for the derivatives of $\scalarembeddingdatahypfortimefunctionarg{\muxmulevelsetvalue}$ relevant for the estimates
	\eqref{E:BOUNDONEMBEDDINGOFXMUEQUALSMINUSKAPPA}--\eqref{E:SMALLNESSBOUNDONEMBEDDINGOFXMUEQUALSMINUSKAPPA}
	and
	\eqref{E:PARTIALTEMBEDDINGFUNCTIONKEYESTIMATE},
	and then use the bootstrap assumptions.
	More precisely, we use the following consequences of 
	Lemmas~\ref{L:SCHEMATICSTRUCTUREOFVARIOUSTENSORSINTERMSOFCONTROLVARS}  and \ref{L:COMMUTATORSTOCOORDINATES}
	and the bootstrap assumptions:
	$\Lunit = \geop{t} + \mathcal{O}(\varepsilon^{1/2}) \geop{x^2} + \mathcal{O}(\varepsilon^{1/2}) \geop{x^3}$,
	$\muX = \geop{u} + \mathcal{O}(\varepsilon^{1/2}) \geop{x^2} + \mathcal{O}(\varepsilon^{1/2}) \geop{x^3}$,
	$\| \muX \upmu \|_{W_{\textnormal{geo}}^{2,\infty}(\twoargMrough{(\timefunction_0,\timefunctionboot),(- \rightu,\leftu)}{\muxmulevelsetvalue})} \leq C$,
	$\left\| 
		\left(
			\geop{x^2} \muX \upmu, \geop{x^3} \muX \upmu
		\right) 
	\right\|_{W_{\textnormal{geo}}^{1,\infty}(\twoargMrough{(\timefunction_0,\timefunctionboot),(- \rightu,\leftu)}{\muxmulevelsetvalue})} \leq C \varepsilon^{1/2}$,
	$\geop{t} \muX \upmu = \Lunit \muX \upmu + \mathcal{O}(\varepsilon^{1/2})$,
	$\geop{u} \muX \upmu = \muX \muX \upmu + \mathcal{O}(\varepsilon^{1/2})$,
	and $\geop{u} \muX \upmu \approx 1$ along $\datahypfortimefunctiontwoarg{-\muxmulevelsetvalue}{[\timefunction_0,\timefunctionboot)}$
	(see \eqref{E:BOOTSTRAPLEVELSETSTRUCTUREANDLOCATIONOFMUXEQUALSMINUSKAPPA} and \eqref{E:BAMUTRANSVERSALCONVEXITY}).
	We have therefore proved
	\eqref{E:BOUNDONEMBEDDINGOFXMUEQUALSMINUSKAPPA}--\eqref{E:SMALLNESSBOUNDONEMBEDDINGOFXMUEQUALSMINUSKAPPA}
	and \eqref{E:PARTIALTEMBEDDINGFUNCTIONKEYESTIMATE}.
	
	We now prove
	\eqref{E:W2INFTYBOUNDFORTORI}--\eqref{E:W1INFTYBOUNDFORANGULARDERIVATIVESTORI}.
	The estimate 
	$\| \Cartesiantisafunctiononmumxtoriarg{\mulevelsetvalue}{-\muxmulevelsetvalue} \|_{L^{\infty}(\mathbb{T}^2)} \leq C$
	follows trivially from \eqref{E:BASIZEOFCARTESIANT},
	as does the precompactness of $\domainforembeddingdatahypfortimefunctiontwoarg{\muxmulevelsetvalue}{(\upmuboot,\mupositive]}$.
	The estimate 
	$\| \Eikonalisafunctiononmumuxtoriarg{\mulevelsetvalue}{-\muxmulevelsetvalue} \|_{L^{\infty}(\mathbb{T}^2)} \leq C$
	follows trivially from the fact that $u \in [- \rightu,\leftu]$ in
	$\twoargMrough{[\timefunction_0,\timefunctionboot),[- \rightu,\leftu]}{\muxmulevelsetvalue}$.
	To control the derivatives of
	$\Cartesiantisafunctiononmumxtoriarg{\mulevelsetvalue}{-\muxmulevelsetvalue}$ and $\Eikonalisafunctiononmumuxtoriarg{\mulevelsetvalue}{-\muxmulevelsetvalue}$
	that are relevant for
	\eqref{E:W2INFTYBOUNDFORTORI}--\eqref{E:W1INFTYBOUNDFORANGULARDERIVATIVESTORI},
	we implicitly differentiate the identities 
	$(\upmu, \muX \upmu) 
	\circ 
	\left(\Cartesiantisafunctiononmumxtoriarg{\mulevelsetvalue}{-\muxmulevelsetvalue}(x^2,x^3),
	\Eikonalisafunctiononmumuxtoriarg{\mulevelsetvalue}{-\muxmulevelsetvalue}(x^2,x^3),
	x^2,x^3 \right) 
	= (\mulevelsetvalue,-\muxmulevelsetvalue)$
	with $\geop{x^2}$ and $\geop{x^3}$
	and argue as in the proof of 
	\eqref{E:BOUNDONEMBEDDINGOFXMUEQUALSMINUSKAPPA}--\eqref{E:SMALLNESSBOUNDONEMBEDDINGOFXMUEQUALSMINUSKAPPA},
	using in addition the fact that
	$\geop{t} \upmu \approx 1$ along $\datahypfortimefunctiontwoarg{-\muxmulevelsetvalue}{[\timefunction_0,\timefunctionboot)}$
	(see \eqref{E:BOOTSTRAPLEVELSETSTRUCTUREANDLOCATIONOFMUXEQUALSMINUSKAPPA} and \eqref{E:BABOUNDSONLMUINTERESTINGREGION}).
	We have therefore proved
	\eqref{E:W2INFTYBOUNDFORTORI}--\eqref{E:W1INFTYBOUNDFORANGULARDERIVATIVESTORI}.
	Similarly, to prove 
	\eqref{E:FORMULAFORMUDERIVATIVEOFCARTESIANTIMEFUNCTIONOFMUXISMINUSKAPPAEMBEEDDING}--\eqref{E:IMPROVEMENTMUEQUALSMINUSKAPPAEMBEDDINGCARTESIANTIMEFUNCTIONNEGATIVEMUDERIVATIVE},
	we implicitly differentiate
	the identities 
	$(\upmu, \muX \upmu) 
	\circ 
	\left(\Cartesiantisafunctiononmumxtoriarg{\mulevelsetvalue}{-\muxmulevelsetvalue}(x^2,x^3),
	\Eikonalisafunctiononmumuxtoriarg{\mulevelsetvalue}{-\muxmulevelsetvalue}(x^2,x^3),
	x^2,x^3 \right) 
	= (\mulevelsetvalue,-\muxmulevelsetvalue)$
	with $\frac{\partial}{\partial \mulevelsetvalue}$
	and use the bootstrap assumptions.
	
	To deduce that 
	$\embeddingdatahypfortimefunctionarg{\muxmulevelsetvalue}$ is $C^{1,1}$ 
	on every compact subset of $\domainforembeddingdatahypfortimefunctiontwoarg{\muxmulevelsetvalue}{(\upmuboot,\mupositive]}$,
	we first note that \eqref{E:BAMUEQUALSMINUSKAPPAEMBEDDINGCARTESIANTIMEFUNCTIONNEGATIVEMUDERIVATIVE} implies that for
	any compact subset $\mathfrak{K}$ of $\domainforembeddingdatahypfortimefunctiontwoarg{\muxmulevelsetvalue}{(\upmuboot,\mupositive]}$,
	there is a $\mulevelsetvalue' \in (\upmuboot,\mupositive]$
	and a compact set 
	$$
	\domainforembeddingdatahypfortimefunctiontwoarg{\muxmulevelsetvalue}{[\mulevelsetvalue',\mupositive]}
	\eqdef
	\left\lbrace 
		(t,x^2,x^3) \in \mathbb{R} \times \mathbb{T}^2 
		\ | \
		\Cartesiantisafunctiononmumxtoriarg{\mupositive}{-\muxmulevelsetvalue}(x^2,x^3) \leq t \leq \Cartesiantisafunctiononmumxtoriarg{\mulevelsetvalue'}{-		
		\muxmulevelsetvalue}(x^2,x^3)
	\right\rbrace
	$$
	such that
	$
	\mathfrak{K}
	\subset
	\domainforembeddingdatahypfortimefunctiontwoarg{\muxmulevelsetvalue}{[\mulevelsetvalue',\mupositive]}
	$.
	Since \eqref{E:EMBEDDEDEDMUEQUALSMINUSKAPPAHYPERSURFACEFINITEC11NORMONCOMPACTSUBSETS} implies that
	$\domainforembeddingdatahypfortimefunctiontwoarg{\muxmulevelsetvalue}{[\mulevelsetvalue',\mupositive]}$ has a $C^1$
	boundary,\footnote{In fact, the boundary is $C^{1,1}$, though we do not need this fact here.} 
	it is a standard Sobolev embedding result (see \cite{lE1998}*{Theorem~5 in Section~5.6})
	that 
	$
	 W_{\textnormal{geo}}^{1,\infty}(\domainforembeddingdatahypfortimefunctiontwoarg{\muxmulevelsetvalue}{(\mulevelsetvalue',\mupositive)})
	\hookrightarrow
	C^{0,1}(\domainforembeddingdatahypfortimefunctiontwoarg{\muxmulevelsetvalue}{[\mulevelsetvalue',\mupositive]})
	$,
	where by \eqref{E:BAMUEQUALSMINUSKAPPAEMBEDDINGCARTESIANTIMEFUNCTIONNEGATIVEMUDERIVATIVE},
	$
	\domainforembeddingdatahypfortimefunctiontwoarg{\muxmulevelsetvalue}{(\mulevelsetvalue',\mupositive)}
	\eqdef
	\left\lbrace 
		(t,x^2,x^3) \in \mathbb{R} \times \mathbb{T}^2 
		\ | \
		\Cartesiantisafunctiononmumxtoriarg{\mupositive}{-\muxmulevelsetvalue}(x^2,x^3) < t < \Cartesiantisafunctiononmumxtoriarg{\mulevelsetvalue'}{-\muxmulevelsetvalue}(x^2,x^3)
	\right\rbrace
	$
	is the interior of $\domainforembeddingdatahypfortimefunctiontwoarg{\muxmulevelsetvalue}{[\mulevelsetvalue',\mupositive]}$.
	Thus, in view of \eqref{E:BOUNDONEMBEDDINGOFXMUEQUALSMINUSKAPPA},
	we conclude that
	$\embeddingdatahypfortimefunctionarg{\muxmulevelsetvalue}$ is $C^{1,1}$ on 
	$\domainforembeddingdatahypfortimefunctiontwoarg{\muxmulevelsetvalue}{[\mulevelsetvalue',\mupositive]}$ 
	as desired.
\end{proof}

\subsection{Properties of the flow map of $\Wtransarg{\muxmulevelsetvalue}$ and 
	the viability of the data hypersurface $\datahypfortimefunctiontwoarg{-\muxmulevelsetvalue}{[\timefunction_0,\timefunctionboot)}$}
\label{SS:PROPERTIESOFFLOWMAPOFW}

\begin{lemma}[The flow map of $\Wtransarg{\muxmulevelsetvalue}$ and 
	the viability of the data hypersurface $\datahypfortimefunctiontwoarg{-\muxmulevelsetvalue}{[\timefunction_0,\timefunctionboot]}$]
	\label{L:FLOWMAPFORGENERATOROFROUGHTIMEFUNCTION} \hfill
	
	\noindent \underline{\textbf{Properties of the flow map of $\Wtransarg{\muxmulevelsetvalue}$}}.
	Let $\Wtransarg{\muxmulevelsetvalue}$ be the vectorfield defined in \eqref{E:WTRANSDEF}, and recall that
	$\Wtransarg{\muxmulevelsetvalue} \timefunction = 0$
	and
	$\Wtransarg{\muxmulevelsetvalue} u = 1$.
	Let $(\Delta u,t,u,x^2,x^3) \rightarrow \flowmapWtransargtwoarg{\muxmulevelsetvalue}{\Delta u}(t,u,x^2,x^3)$ 
	denote the flow map of $\Wtransarg{\muxmulevelsetvalue}$,
	i.e., at each fixed 
	$
	(t,u,x^2,x^3) \in \twoargMrough{[\timefunction_0,\timefunctionboot),[- \rightu,\leftu]}{\muxmulevelsetvalue}
	$,
	the components of $\flowmapWtransargtwoarg{\muxmulevelsetvalue}{\Delta u}(t,u,x^2,x^3)$ 
	solve the following ODE system initial value problem
	on the flow interval $\Delta u \in [- \rightu - u,\leftu-u]$:
	\begin{align} \label{E:FLOWMAPFORGENERATOROFROUGHTIMEFUNCTION}
		\frac{\partial}{\partial \Delta u} \flowmapWtransargtwoarg{\muxmulevelsetvalue}{\Delta u}(t,u,x^2,x^3) 
		& = \Wtransarg{\muxmulevelsetvalue} \circ \flowmapWtransargtwoarg{\muxmulevelsetvalue}{\Delta u}(t,u,x^2,x^3),
		&
		\flowmapWtransargtwoarg{\muxmulevelsetvalue}{0}(t,u,x^2,x^3)
		& = 
		(t,u,x^2,x^3).
	\end{align}
	
	Then for each fixed $\timefunction \in [\timefunction_0,\timefunctionboot)$
	and each pair $u_1, u_2 \in [- \rightu, \leftu]$,
	$\flowmapWtransargtwoarg{\muxmulevelsetvalue}{u_2 - u_1}$ is a diffeomorphism from 
	the rough torus $\twoargroughtori{\timefunction,u_1}{\muxmulevelsetvalue}$ onto the rough torus $\twoargroughtori{\timefunction,u_2}{\muxmulevelsetvalue}$.
	In particular, the integral curves of $\Wtransarg{\muxmulevelsetvalue}$ thread $\hypthreearg{\timefunction}{[- \rightu,\leftu]}{\muxmulevelsetvalue}$.
	Moreover, for every fixed $\timefunction \in [\timefunction_0,\timefunctionboot)$,
	each integral curve of $\Wtransarg{\muxmulevelsetvalue}$ passes through precisely one point on
	the $\upmu$-adapted torus
	$
	\twoargmumuxtorus{-\timefunction}{-\muxmulevelsetvalue}$
	defined in \eqref{E:MUXMUTORI}.

	Moreover, with $d_{\textnormal{geo}}$ denoting the differential with respect to the geometric coordinates,
	we have the following bounds, where the implicit constants in \eqref{E:BOUNDFORDETERMINANTWFLOWMAP}
	are independent of all $\Delta u$ such that $|\Delta u| \leq |2 \interestingu + \leftu|$:
	\begin{align} \label{E:W3INFTYBOUNDFORWFLOWMAP} 
		\sup_{|\Delta u| \leq |2 \interestingu + \leftu|}
		\| d_{\textnormal{geo}} \flowmapWtransargtwoarg{\muxmulevelsetvalue}{\Delta u} \|_{W_{\textnormal{geo}}^{2,\infty}
		(
		\twoargMrough{(\timefunction_0,\timefunctionboot),(- \rightu,\leftu)}{\muxmulevelsetvalue}
		\cap
		\twoargMrough{(\timefunction_0,\timefunctionboot),(- \rightu - \Delta u,\leftu - \Delta u)}{\muxmulevelsetvalue})
		}
		& \leq C,
	\end{align}
	\begin{align} \label{E:BOUNDFORDETERMINANTWFLOWMAP}
		\mbox{\upshape det } \left(d_{\textnormal{geo}} \flowmapWtransargtwoarg{\muxmulevelsetvalue}{\Delta u} \right)
		&
		\approx 1
		\mbox{ on }
		\twoargMrough{[\timefunction_0,\timefunctionboot),[- \rightu,\leftu]}{\muxmulevelsetvalue}
		\cap
		\twoargMrough{[\timefunction_0,\timefunctionboot),[- \rightu - \Delta u,\leftu - \Delta u]}{\muxmulevelsetvalue}.
	\end{align}
	
	\medskip
	\noindent \underline{\textbf{Estimates tied to the flow of
	$\datahypfortimefunctiontwoarg{-\muxmulevelsetvalue}{[\timefunction_0,\timefunctionboot)}$ by $\Wtransarg{\muxmulevelsetvalue}$}}.
	Let 
	$\embeddingdatahypfortimefunctionarg{\muxmulevelsetvalue}: \domainforembeddingdatahypfortimefunctiontwoarg{\muxmulevelsetvalue}{(\upmuboot,\mupositive]} 
	\rightarrow \twoargMrough{[\timefunction_0,\timefunctionboot),[-\frac{3}{4} \interestingu,\frac{3}{4} \interestingu]}{\muxmulevelsetvalue}$
	be the embedding of 
	$\datahypfortimefunctiontwoarg{-\muxmulevelsetvalue}{[\timefunction_0,\timefunctionboot)}$
	from Lemma~\ref{L:XMUISMINUSCAPPISAGRAPH},
	which is of the form $\embeddingdatahypfortimefunctionarg{\muxmulevelsetvalue}(t,x^2,x^3) 
	= 
	\left(t,\scalarembeddingdatahypfortimefunctionarg{\muxmulevelsetvalue}(t,x^2,x^3),x^2,x^3 \right)$,
	and let $\composedflowmapdiffeoarg{\muxmulevelsetvalue} : \domaincomposedflowmapdiffeoarg{\muxmulevelsetvalue} \rightarrow \twoargMrough{[\timefunction_0,\timefunctionboot),[- \rightu,\leftu]}{\muxmulevelsetvalue}$
	be the map and set defined by:
	\begin{align}
	\composedflowmapdiffeoarg{\muxmulevelsetvalue}(\Delta u,t,x^2,x^3)
		&
		\eqdef
		\flowmapWtransargtwoarg{\muxmulevelsetvalue}{\Delta u} 
		\circ \embeddingdatahypfortimefunctionarg{\muxmulevelsetvalue}(t,x^2,x^3),
	\end{align}
	\begin{align}
	\domaincomposedflowmapdiffeoarg{\muxmulevelsetvalue}
	& 
	\eqdef
		\left\lbrace 
			(\Delta u,t,x^2,x^3)
			\in \mathbb{R} \times \mathbb{R} \times \mathbb{T}^2 
			\ | \
			(t,x^2,x^3) \in \domainforembeddingdatahypfortimefunctiontwoarg{\muxmulevelsetvalue}{(\upmuboot,\mupositive]}
			\mbox{ and }
			\Delta u 
			\in 
			\left[- \rightu - \scalarembeddingdatahypfortimefunctionarg{\muxmulevelsetvalue}(t,x^2,x^3),
			\leftu - \scalarembeddingdatahypfortimefunctionarg{\muxmulevelsetvalue}(t,x^2,x^3) 
			\right]
		\right\rbrace.
	\end{align}
	Then $\composedflowmapdiffeoarg{\muxmulevelsetvalue}$ is a diffeomorphism from 
	$\domaincomposedflowmapdiffeoarg{\muxmulevelsetvalue}$ onto 
	$
	\twoargMrough{[\timefunction_0,\timefunctionboot),[- \rightu,\leftu]}{\muxmulevelsetvalue}
	$
	such that $\composedflowmapdiffeoarg{\muxmulevelsetvalue}$ and its inverse function 
	$(\composedflowmapdiffeoarg{\muxmulevelsetvalue})^{-1}$
	satisfy the following bounds:
	\begin{align} \label{E:W2INFTYBOUNDFORXMUISMINUSKAPPADIFFEOMORPHISM}
		\| \composedflowmapdiffeoarg{\muxmulevelsetvalue} \|_{W^{2,\infty}(\mbox{\upshape int}(\domaincomposedflowmapdiffeoarg{\muxmulevelsetvalue}))}
			& \leq C,
				\\
		\| (\composedflowmapdiffeoarg{\muxmulevelsetvalue})^{-1} \|_{W_{\textnormal{geo}}^{2,\infty}(\twoargMrough{(\timefunction_0,\timefunctionboot),(- \rightu,\leftu)}{\muxmulevelsetvalue})}
			& \leq C,
			\label{E:W2INFTYBOUNDFORXMUISMINUSKAPPAINVERSEDIFFEOMORPHISM}
	\end{align}
	where: 
	\begin{align} 
	\begin{split} \label{E:INTERIORDOMAINWTRANSFLOWMAP}
	(\mbox{\upshape int}(\domaincomposedflowmapdiffeoarg{\muxmulevelsetvalue})
	& =
	\Big\lbrace 
		(\Delta u,t,x^2,x^3)
		\in \mathbb{R} \times \mathbb{R} \times \mathbb{T}^2 
		\ | \
		(x^2,x^3) \in \mathbb{T}^2,
			\,
		\Cartesiantisafunctiononmumxtoriarg{\mupositive}{-\muxmulevelsetvalue}(x^2,x^3) < t < \Cartesiantisafunctiononmumxtoriarg{\upmuboot}{-\muxmulevelsetvalue}(x^2,x^3),
				\\
		& \ \ \ \ \ \ \
		\mbox{ and }
		\Delta u 
		\in 
		\left(- \rightu - \scalarembeddingdatahypfortimefunctionarg{\muxmulevelsetvalue}(t,x^2,x^3),
			\leftu-\scalarembeddingdatahypfortimefunctionarg{\muxmulevelsetvalue}(t,x^2,x^3)\right)
	\Big\rbrace
	\end{split}
	\end{align} 
	is the interior of $\domaincomposedflowmapdiffeoarg{\muxmulevelsetvalue}$,
	and $\Cartesiantisafunctiononmumxtoriarg{\mupositive}{-\muxmulevelsetvalue}$
	and $\Cartesiantisafunctiononmumxtoriarg{\upmuboot}{-\muxmulevelsetvalue}$
	are the functions appearing in \eqref{E:XMUISMINUSCAPPISAGRAPH}.
	Finally, $\composedflowmapdiffeoarg{\muxmulevelsetvalue}$ is $C^{1,1}$ on every
	compact subset of $\domaincomposedflowmapdiffeoarg{\muxmulevelsetvalue}$, and $(\composedflowmapdiffeoarg{\muxmulevelsetvalue})^{-1}$
	is $C^{1,1}$ on every compact subset of
	$
	\twoargMrough{[\timefunction_0,\timefunctionboot),[- \rightu,\leftu]}{\muxmulevelsetvalue}
	$.
	\end{lemma}

\begin{proof}
	Throughout this proof, we silently use the soft bootstrap assumptions of
	Sect.\,\ref{SSS:SOFTBACONCERNINGREGULARITY}, which guarantee sufficient
	qualitative regularity. For quantitative estimates, 
	we will use the bootstrap assumptions of Sects.\,\ref{SSS:FUNDAMENTALQUANTITATIVE} and \ref{SSS:AUXBOOTSTRAP}.
	
	From \eqref{E:BAMUTORI},
	\eqref{E:BOOSTRAPTORILOCATION},
	and the facts that $\Wtransarg{\muxmulevelsetvalue} \timefunctionarg{\muxmulevelsetvalue} = 0$
	and $\Wtransarg{\muxmulevelsetvalue} u = 1$,
	it follows that for each fixed $\timefunction \in [\timefunction_0,\timefunctionboot)$,
	every integral curve of
	$\Wtransarg{\muxmulevelsetvalue}$
	in $\hypthreearg{\timefunction}{[- \rightu,\leftu]}{\muxmulevelsetvalue}$
	must intersect 
	$\twoargmumuxtorus{-\timefunction}{-\muxmulevelsetvalue}$
	at one or more points in $\hypthreearg{\timefunction}{[-\frac{3}{4} \interestingu, \frac{3}{4} \interestingu]}{\muxmulevelsetvalue}$.
	Recalling that 
	$\upmu|_{\twoargmumuxtorus{-\timefunction}{-\muxmulevelsetvalue}} = - \timefunction$
	and $\Wtransarg{\muxmulevelsetvalue} \upmu|_{\twoargmumuxtorus{-\timefunction}{-\muxmulevelsetvalue}} = 0$,
	and using the transversal convexity bootstrap assumption \eqref{E:BAMUTRANSVERSALCONVEXITY} for
	$\Wtransarg{\muxmulevelsetvalue} \Wtransarg{\muxmulevelsetvalue} \upmu$,
	we see that the intersection occurs at a unique point.
	
	
	Next, we differentiate the evolution equation in \eqref{E:FLOWMAPFORGENERATOROFROUGHTIMEFUNCTION}
	and use the chain rule to deduce that
	$d_{\textnormal{geo}} \flowmapWtransargtwoarg{\muxmulevelsetvalue}{\Delta u}$
	satisfies the linear ODE system initial value problem:
	\begin{align} \label{E:ODEFORJACOBIANOFFLOWMAP}
	\frac{\partial}{\partial \Delta u} d_{\textnormal{geo}} \flowmapWtransargtwoarg{\muxmulevelsetvalue}{\Delta u}
	& 
	=
	(d_{\textnormal{geo}} \Wtransarg{\muxmulevelsetvalue}) 
	\circ 
	\flowmapWtransargtwoarg{\muxmulevelsetvalue}{\Delta u}(t,u,x^2,x^3)
	\cdot 
	d_{\textnormal{geo}} \flowmapWtransargtwoarg{\muxmulevelsetvalue}{\Delta u}(t,u,x^2,x^3),
		\\
	d_{\textnormal{geo}} \flowmapWtransargtwoarg{\muxmulevelsetvalue}{0} 
	&
	= \mbox{\upshape diag}(1,1,1,1).
	\label{E:DATAFORODEFORJACOBIANOFFLOWMAP}
	\end{align}
	Definition \eqref{E:WTRANSDEF},
	Lemma~\ref{L:COMMUTATORSTOCOORDINATES},
	and the bootstrap assumptions imply that:
	\begin{align} \label{E:BOUNDSFORGEOCOMPONENTSOFWTRANS}
	\| \left(
		\Wtransarg{\muxmulevelsetvalue} t, \Wtransarg{\muxmulevelsetvalue} u, \Wtranstwoarg{\muxmulevelsetvalue}{2}, \Wtranstwoarg{\muxmulevelsetvalue}{3} 
	\right)\|_{W_{\textnormal{geo}}^{3,\infty}(\twoargMrough{(\timefunction_0,\timefunctionboot),(- \rightu,\leftu)}{\muxmulevelsetvalue})}
	& \leq C
	\end{align}
	and thus:
	\begin{align} \label{E:W2INIFINTYBOUNDFORGEOMTRICPARTIALDERIVATIVESOFWTRANSCOMPONENTS}
		\| d_{\textnormal{geo}} \Wtransarg{\muxmulevelsetvalue} \|_{W_{\textnormal{geo}}^{2,\infty}(\twoargMrough{(\timefunction_0,\timefunctionboot),(- \rightu,\leftu)}{\muxmulevelsetvalue})}
	\leq C,
	\end{align}
	where we are viewing $d_{\textnormal{geo}} \Wtransarg{\muxmulevelsetvalue}$ to be the Jacobian matrix of the map
	$(t,u,x^2,x^3) 
	\rightarrow 
	\left(
		\Wtransarg{\muxmulevelsetvalue} t, \Wtransarg{\muxmulevelsetvalue} u, \Wtranstwoarg{\muxmulevelsetvalue}{2}, \Wtranstwoarg{\muxmulevelsetvalue}{3} 
	\right)$.
	From these facts,
	Gr\"{o}nwall's inequality, 
	and the fact that $|\Delta u| \leq \leftu + 2 \interestingu$ in the region under study,
	we conclude
	\eqref{E:W3INFTYBOUNDFORWFLOWMAP}.
	
	The estimate \eqref{E:BOUNDFORDETERMINANTWFLOWMAP}
	follows from a similar argument based on the linear ODE system:
	\begin{align} \label{E:ODEFORLOGOFDETERMINANTOFJACOBIANOFFLOWMAP}
	\frac{\partial}{\partial \Delta u} \ln \mbox{\upshape det } \left(\mathrm{d}_{\textnormal{geo}} \flowmapWtransargtwoarg{\muxmulevelsetvalue}{\Delta u} \right)
	& =
	\mbox{\upshape tr } (\mathrm{d}_{\textnormal{geo}} \Wtransarg{\muxmulevelsetvalue}) 
	\circ 
	\flowmapWtransargtwoarg{\muxmulevelsetvalue}{\Delta u}(t,u,x^2,x^3), \qquad 
	\det \left(\mathrm{d}_{\textnormal{geo}} \flowmapWtransargtwoarg{\muxmulevelsetvalue}{0} \right) = 1
	\end{align}
	and the estimate $\mytr(\mathrm{d}_{\textnormal{geo}} \Wtransarg{\muxmulevelsetvalue}) = \mathcal{O}(\auxbootsmall)$, 
	which is a consequence of definition \eqref{E:WTRANSDEF},
	Lemma~\ref{L:COMMUTATORSTOCOORDINATES},
	and the bootstrap assumptions.
	
	The estimate \eqref{E:W2INFTYBOUNDFORXMUISMINUSKAPPADIFFEOMORPHISM}
	follows from
	\eqref{E:BOUNDONEMBEDDINGOFXMUEQUALSMINUSKAPPA},
	\eqref{E:W3INFTYBOUNDFORWFLOWMAP},
	the chain rule,
	and finally from using the equation \eqref{E:FLOWMAPFORGENERATOROFROUGHTIMEFUNCTION}
	and the aforementioned estimate
	$
	\| \left(
		\Wtransarg{\muxmulevelsetvalue} t, \Wtransarg{\muxmulevelsetvalue} u, \Wtranstwoarg{\muxmulevelsetvalue}{2}, \Wtranstwoarg{\muxmulevelsetvalue}{3} 
	\right)\|_{W_{\textnormal{geo}}^{3,\infty}(\twoargMrough{(\timefunction_0,\timefunctionboot),(- \rightu,\leftu)}{\muxmulevelsetvalue})}
	\leq C
	$
	to control the partial derivatives of $\composedflowmapdiffeoarg{\muxmulevelsetvalue}$
	with respect to $\Delta u$.
	
	Next, we highlight that the map
	$
	\Delta u \rightarrow \flowmapWtransargtwoarg{\muxmulevelsetvalue}{\Delta u}(t,u,x^2,x^3)
	$
	is just a parameterization of the integral curve of $\Wtransarg{\muxmulevelsetvalue}$
	that passes through the point $(t,u,x^2,x^3)$.
	Moreover, we recall that we have already shown that
	every integral curve of
	$\Wtransarg{\muxmulevelsetvalue}$
	in $\hypthreearg{\timefunction}{[- \rightu,\leftu]}{\muxmulevelsetvalue}$
	must intersect 
	$
	\datahypfortimefunctiontwoarg{-\muxmulevelsetvalue}{[\timefunction_0,\timefunctionboot)}
	$
	at a unique point 
	(more precisely, a point in the torus 
	$
	\twoargmumuxtorus{-\timefunction}{-\muxmulevelsetvalue}
	$,
	which is contained in 
	$
	\datahypfortimefunctiontwoarg{-\muxmulevelsetvalue}{[\timefunction_0,\timefunctionboot)}
	$
	),
	where $\datahypfortimefunctiontwoarg{-\muxmulevelsetvalue}{[\timefunction_0,\timefunctionboot)}$
	is the image of the set $\domainforembeddingdatahypfortimefunctiontwoarg{\muxmulevelsetvalue}{(\upmuboot,\mupositive]}$ under
	$\embeddingdatahypfortimefunctionarg{\muxmulevelsetvalue}$ (see Lemma~\ref{L:XMUISMINUSCAPPISAGRAPH}).
	It follows that $\composedflowmapdiffeoarg{\muxmulevelsetvalue}$ is a bijection from 
	$\domaincomposedflowmapdiffeoarg{\muxmulevelsetvalue}$ onto $\twoargMrough{[\timefunction_0,\timefunctionboot),[- \rightu,\leftu]}{\muxmulevelsetvalue}$.
	Thus, to conclude that $\composedflowmapdiffeoarg{\muxmulevelsetvalue}$ is a diffeomorphism on $\domaincomposedflowmapdiffeoarg{\muxmulevelsetvalue}$, 
	it remains only for us to show that 
	its Jacobian matrix $d_{(\Delta u,t,x^2,x^3)} \composedflowmapdiffeoarg{\muxmulevelsetvalue}$ has non-vanishing determinant.
	To this end,
	we note that the standard theory of 
	flow maps, 
	the identity
	$\frac{\partial}{\partial \Delta u} \flowmapWtransargtwoarg{\muxmulevelsetvalue}{\Delta u}(t,u,x^2,x^3) 
	=
	\frac{\partial}{\partial u} \flowmapWtransargtwoarg{\muxmulevelsetvalue}{\Delta u}(t,u,x^2,x^3) 
	$
	(which relies on the fact that $\Wtransarg{\muxmulevelsetvalue} u = 1$),
	and the chain rule yield the identity
	$d_{(\Delta u,t,x^2,x^3)} \composedflowmapdiffeoarg{\muxmulevelsetvalue}(\Delta u,t,x^2,x^3)
	=
	 [(d_{\textnormal{geo}} \flowmapWtransargtwoarg{\muxmulevelsetvalue}{\Delta u}) \circ \embeddingdatahypfortimefunctionarg{\muxmulevelsetvalue}(t,x^2,x^3)]
	 \cdot
		 \flowmapWtransMatrixarg{\muxmulevelsetvalue}(t,x^2,x^3)
	$,
	where $\flowmapWtransMatrixarg{\muxmulevelsetvalue}$ is the $4 \times 4$ matrix-valued function on 
	$\domainforembeddingdatahypfortimefunctiontwoarg{\muxmulevelsetvalue}{(\upmuboot,\mupositive]}$
	whose first column is
	$(0,1,0,0)^{\top}$
	and whose last three columns are the Jacobian $d_{(t,x^2,x^3)} \embeddingdatahypfortimefunctionarg{\muxmulevelsetvalue}$.
	Since
	$\mbox{\upshape det } \left(d_{(\Delta u,t,x^2,x^3)} \composedflowmapdiffeoarg{\muxmulevelsetvalue}(\Delta u,t,x^2,x^3) \right)
	=
	 \mbox{\upshape det } \left((d_{\textnormal{geo}} \flowmapWtransargtwoarg{\muxmulevelsetvalue}{\Delta u}) \circ \embeddingdatahypfortimefunctionarg{\muxmulevelsetvalue}\right)
	 \cdot
		\mbox{\upshape det } \left(\flowmapWtransMatrixarg{\muxmulevelsetvalue}(t,x^2,x^3) \right)
	$,
	we see from \eqref{E:BOUNDFORDETERMINANTWFLOWMAP} that to prove:
	\begin{align} \label{E:BOUNDFORJACOBIANDETERMINANTOFCOMPOSITION}
		|\mbox{\upshape det } \left(d_{(\Delta u,t,x^2,x^3)} \composedflowmapdiffeoarg{\muxmulevelsetvalue} \right)|
		& \approx 1,
		\mbox{on the domain }
		\domaincomposedflowmapdiffeoarg{\muxmulevelsetvalue},
	\end{align}
	we need only to show that
	$
	|\mbox{\upshape det }  \flowmapWtransMatrixarg{\muxmulevelsetvalue}|
	\approx 1
	$
	on the domain $\domainforembeddingdatahypfortimefunctiontwoarg{\muxmulevelsetvalue}{(\upmuboot,\mupositive]}$.
	To this end,
	we use Lemma~\ref{L:XMUISMINUSCAPPISAGRAPH}
	(in particular \eqref{E:SMALLNESSBOUNDONEMBEDDINGOFXMUEQUALSMINUSKAPPA} and \eqref{E:PARTIALTEMBEDDINGFUNCTIONKEYESTIMATE})
	and the bootstrap assumptions to deduce that:
	\begin{align} \label{E:MATRIXARISINGINSTUDYOFFLOWMAPOFWTRANS}
	 \flowmapWtransMatrixarg{\muxmulevelsetvalue}
		& = \begin{pmatrix}
			0 & 1 & 0 & 0  \\
			1 & - \frac{\Lunit \muX \upmu \circ \embeddingdatahypfortimefunctionarg{\muxmulevelsetvalue}}{\muX \muX \upmu \circ \embeddingdatahypfortimefunctionarg{\muxmulevelsetvalue}}  + * & * & * \\
			0 & 0 & 1 & 0 \\
			0 & 0 & 0 & 1
		\end{pmatrix},
	\end{align}
	where
	``$*$'' denotes $\mathcal{O}(\varepsilon^{1/2})$ quantities.
	From this identity and the transversal convexity bootstrap assumption \eqref{E:BAMUTRANSVERSALCONVEXITY},
	we conclude that in fact,
	$
	\mbox{\upshape det } \flowmapWtransMatrixarg{\muxmulevelsetvalue}
	= -1
	$.
	
	\eqref{E:W2INFTYBOUNDFORXMUISMINUSKAPPAINVERSEDIFFEOMORPHISM} follows from
	differentiating the identity
	$
	\composedflowmapdiffeoarg{\muxmulevelsetvalue} \circ (\composedflowmapdiffeoarg{\muxmulevelsetvalue})^{-1}(t,u,x^2,x^3)
	= (t,u,x^2,x^3)
	$
	and using
	\eqref{E:W2INFTYBOUNDFORXMUISMINUSKAPPADIFFEOMORPHISM},
	\eqref{E:BOUNDFORJACOBIANDETERMINANTOFCOMPOSITION},
	and the chain rule.
	
	Finally,
	the facts that  
	$\composedflowmapdiffeoarg{\muxmulevelsetvalue}$ is $C^{1,1}$ on every
	compact subset of $\domaincomposedflowmapdiffeoarg{\muxmulevelsetvalue}$, and $(\composedflowmapdiffeoarg{\muxmulevelsetvalue})^{-1}$
	is $C^{1,1}$ on every compact subset of
	$
	\twoargMrough{[\timefunction_0,\timefunctionboot),[- \rightu,\leftu]}{\muxmulevelsetvalue}
	$
	follow from the bounds
	\eqref{E:W2INFTYBOUNDFORXMUISMINUSKAPPADIFFEOMORPHISM}--\eqref{E:W2INFTYBOUNDFORXMUISMINUSKAPPAINVERSEDIFFEOMORPHISM}
	and arguments similar to the ones
	given at the end of the proof of Lemma~\ref{L:XMUISMINUSCAPPISAGRAPH},
	so we omit the details.
	The main point in the proof for $\composedflowmapdiffeoarg{\muxmulevelsetvalue}$ is that
	compact subsets of $\domaincomposedflowmapdiffeoarg{\muxmulevelsetvalue}$
	are contained in a compact set of the form
	$
	\lbrace 
		(\Delta u,t,x^2,x^3)
		\in \mathbb{R} \times \mathbb{R} \times \mathbb{T}^2 
		\ | \
		(x^2,x^3) \in \mathbb{T}^2,
			\,
		\Cartesiantisafunctiononmumxtoriarg{\mupositive}{-\muxmulevelsetvalue}(x^2,x^3) \leq t \leq \Cartesiantisafunctiononmumxtoriarg{\upmu'}{-\muxmulevelsetvalue}(x^2,x^3),
			\,
		\mbox{ and }
		a \leq \Delta u + \scalarembeddingdatahypfortimefunctionarg{\muxmulevelsetvalue}(t,x^2,x^3) \leq b
	\rbrace
	$
	for some $\upmu' \in (\upmuboot,\mupositive]$ and $2\interestingu \leq a \leq b \leq \leftu$,
	and such sets have a $C^1$ boundary.
	Similarly, the main point in the proof for $(\composedflowmapdiffeoarg{\muxmulevelsetvalue})^{-1}$ 
	is that compact subsets of $\twoargMrough{[\timefunction_0,\timefunctionboot),[- \rightu,\leftu]}{\muxmulevelsetvalue}$
	are contained in a compact set of the form
	$
	\twoargMrough{[\timefunction_0,\timefunction],[- \rightu,\leftu]}{\muxmulevelsetvalue}
	$
	for some $\timefunction \in [\timefunction_0,\timefunctionboot)$,
	and such sets have a $C^1$ boundary.

\end{proof}

\section{Estimates for the rough time function, continuous extensions, and diffeomorphisms and homeomorphisms}
\label{S:ESTIMATESFORROUGHTIMEFUNCTIONCONTINUOUSEXTENSIONSANDDIFFEOS}
In this section, 
we use the results of Sect.\,\ref{S:EMBEDDINGSANDFLOWMAPS} 
to derive estimates for the rough time function
$\timefunctionarg{\muxmulevelsetvalue}$.
We then show that the map $\CHOVgeotorough{\muxmulevelsetvalue}(t,u,x^2,x^3) = (\timefunctionarg{\muxmulevelsetvalue},u,x^2,x^3)$ 
extends to a diffeomorphism on the closure of
$\twoargMrough{[\timefunction_0,\timefunctionboot),[- \rightu,\leftu]}{\muxmulevelsetvalue}$.
Finally, we establish related results for the map
$\CHOVroughtomumuxmu{\muxmulevelsetvalue}(\timefunctionarg{\muxmulevelsetvalue},u,x^2,x^3) = (\upmu,\muX \upmu,x^2,x^3)$	
from Def.\,\ref{D:ALLTHECHOVMAPS}.

\subsection{Transport equation solutions that are smoother than the data hypersurface}
\label{SS:TRANSPORTEQUATIONSOLUTIONSSMOOTHERTHANHYPERSURFACE}
We will use the following lemma to derive 
$W_{\textnormal{geo}}^{3,\infty}(\twoargMrough{(\timefunction_0,\timefunctionboot),(- \rightu,\leftu)}{\muxmulevelsetvalue})$ estimates for the rough time function. There is one nonstandard aspect that we carefully handle: 
the initial data are smoother than the initial hypersurface.

\begin{lemma}[Transport equation solutions that are smoother than the data hypersurface]
		\label{L:ODESOLUTIONSTHATARESMOOTHERTHANTHEDATAHYPERSURFACE}
		Consider the initial value problem for a scalar function $\phi$:
		\begin{align} \label{E:WTRANSPORTINITIALVALUEPROBLEM}
			\Wtransarg{\muxmulevelsetvalue} \phi
			& 
			= 0,
			&
			\phi|_{\datahypfortimefunctiontwoarg{-\muxmulevelsetvalue}{[\timefunction_0,\timefunctionboot)}}
			= 
			\ambient|_{\datahypfortimefunctiontwoarg{-\muxmulevelsetvalue}{[\timefunction_0,\timefunctionboot)}},
		\end{align}
		and assume that $\ambient \in W_{\textnormal{geo}}^{3,\infty}(\twoargMrough{(\timefunction_0,\timefunctionboot),(- \rightu,\leftu)}{\muxmulevelsetvalue})$ 
		is an ``ambient'' spacetime function of the geometric coordinates
		satisfying $(\Wtransarg{\muxmulevelsetvalue} \ambient)|_{\datahypfortimefunctiontwoarg{-\muxmulevelsetvalue}{[\timefunction_0,\timefunctionboot)}}
		= 0$,
		where the hypersurface $\datahypfortimefunctiontwoarg{-\muxmulevelsetvalue}{[\timefunction_0,\timefunctionboot)}$
		is $W^{2,\infty}$ (as was shown in Lemma~\ref{L:XMUISMINUSCAPPISAGRAPH} via the embedding $\embeddingdatahypfortimefunctionarg{\muxmulevelsetvalue}$).
		Then there exists a unique solution 
		$\phi \in W_{\textnormal{geo}}^{3,\infty}(\twoargMrough{(\timefunction_0,\timefunctionboot),(- \rightu,\leftu)}{\muxmulevelsetvalue})$
		satisfying the following estimates:
		\begin{subequations}
		\begin{align} \label{E:TRANSPORTEQUATIONSOLUTIONW3INFINITYESTIMATEWITHDATAONROUGHHYPERSURFACE}
			\| \phi \|_{W_{\textnormal{geo}}^{3,\infty}(\twoargMrough{(\timefunction_0,\timefunctionboot),(- \rightu,\leftu)}{\muxmulevelsetvalue})}
			& \lesssim 
				\| \ambient \|_{W_{\textnormal{geo}}^{3,\infty}(\twoargMrough{(\timefunction_0,\timefunctionboot),(- \rightu,\leftu)}{\muxmulevelsetvalue})},
					\\
			\left\| 
				\left(\geop{x^2} \phi,\geop{x^3} \phi \right)  
			\right\|_{W_{\textnormal{geo}}^{2,\infty}(\twoargMrough{(\timefunction_0,\timefunctionboot),(- \rightu,\leftu)}{\muxmulevelsetvalue})}
			& \lesssim 
				\left\| 
					\left(\geop{x^2} \ambient,\geop{x^3} \ambient \right)  
				\right\|_{W_{\textnormal{geo}}^{2,\infty}(\twoargMrough{(\timefunction_0,\timefunctionboot),(- \rightu,\leftu)}{\muxmulevelsetvalue})}.
				\label{E:TANGENTIALDERIVATIVESETIMATETRANSPORTEQUATIONSOLUTIONW3INFINITYESTIMATEWITHDATAONROUGHHYPERSURFACE}
		\end{align}
		\end{subequations}
	\end{lemma}
	
	\begin{remark}[The solution is more regular than the data hypersurface]
		\label{R:SOLUTIONISMOREREGULARTHANHYPERSURFACE}
		The main point of the lemma is that the solution has the same regularity as $\ambient$,
		even though the embedded data hypersurface $\datahypfortimefunctiontwoarg{-\muxmulevelsetvalue}{[\timefunction_0,\timefunctionboot)}$
		is one degree less regular,
		that is, the embedding $\embeddingdatahypfortimefunctionarg{\muxmulevelsetvalue}$ from Lemma~\ref{L:XMUISMINUSCAPPISAGRAPH}
		satisfies 
		$\embeddingdatahypfortimefunctionarg{\muxmulevelsetvalue}\in W^{2,\infty}(\mbox{\upshape int}
			(\domainforembeddingdatahypfortimefunctiontwoarg{\muxmulevelsetvalue}{(\upmuboot,\mupositive]}))$.
	\end{remark}
	
\begin{proof}
	We prove only \eqref{E:TRANSPORTEQUATIONSOLUTIONW3INFINITYESTIMATEWITHDATAONROUGHHYPERSURFACE}
	since \eqref{E:TANGENTIALDERIVATIVESETIMATETRANSPORTEQUATIONSOLUTIONW3INFINITYESTIMATEWITHDATAONROUGHHYPERSURFACE}
	can be proved using similar arguments.

	Throughout the proof, we silently use the fact that
	functions $f \in W_{\textnormal{geo}}^{1,\infty}(\twoargMrough{(\timefunction_0,\timefunctionboot),(- \rightu,\leftu)}{\muxmulevelsetvalue})$,
	are locally Lipschitz and thus have (by Rademacher's theorem) 
	a.e.\ differentiable locally Lipschitz traces along the $W^{2,\infty}$ hypersurface
	$\datahypfortimefunctiontwoarg{-\muxmulevelsetvalue}{[\timefunction_0,\timefunctionboot)}$
	such that 
	$
	\| f \|_{W^{1,\infty}(\mbox{\upshape int}(\datahypfortimefunctiontwoarg{-\muxmulevelsetvalue}{[\timefunction_0,\timefunctionboot)}))}
	\lesssim
	\| f \|_{W_{\textnormal{geo}}^{1,\infty}(\twoargMrough{(\timefunction_0,\timefunctionboot),(- \rightu,\leftu)}{\muxmulevelsetvalue})}$.
	
	By Lemma~\ref{L:FLOWMAPFORGENERATOROFROUGHTIMEFUNCTION},
	the evolution equation $\Wtransarg{\muxmulevelsetvalue} \phi = 0$
	is equivalent to the following ODE initial value problem:
	\begin{align} \label{E:ODEFORMULATIONOFWTRANSTRANSPORTEQUATION}
		\frac{\partial}{\partial \Delta u} \left(\phi \circ \composedflowmapdiffeoarg{\muxmulevelsetvalue}(\Delta u,t,x^2,x^3) \right)
		& = 
		0,
			\\
		\phi \circ \composedflowmapdiffeoarg{\muxmulevelsetvalue}(0,t,x^2,x^3)
		& = \ambient \circ \embeddingdatahypfortimefunctionarg{\muxmulevelsetvalue}(t,x^2,x^3),
			\label{E:DATAFORODEFORMULATIONOFWTRANSTRANSPORTEQUATION}
	\end{align}
	where $\composedflowmapdiffeoarg{\muxmulevelsetvalue}$ is the diffeomorphism
	from $\domaincomposedflowmapdiffeoarg{\muxmulevelsetvalue}$ onto 
	$\twoargMrough{[\timefunction_0,\timefunctionboot),[- \rightu,\leftu]}{\muxmulevelsetvalue}$
	from the lemma.
	The solution to 
	\eqref{E:ODEFORMULATIONOFWTRANSTRANSPORTEQUATION}--\eqref{E:DATAFORODEFORMULATIONOFWTRANSTRANSPORTEQUATION}
	is
	$
	\phi \circ \composedflowmapdiffeoarg{\muxmulevelsetvalue}(\Delta u,t,x^2,x^3)
	=
	\ambient \circ \embeddingdatahypfortimefunctionarg{\muxmulevelsetvalue}(t,x^2,x^3)
$.
From this formula,
the assumptions of the lemma,
\eqref{E:BOUNDONEMBEDDINGOFXMUEQUALSMINUSKAPPA},
and \eqref{E:W2INFTYBOUNDFORXMUISMINUSKAPPADIFFEOMORPHISM},
it immediately follows that 
$\phi \circ \composedflowmapdiffeoarg{\muxmulevelsetvalue} \in W^{2,\infty}(\mbox{\upshape int}(\domaincomposedflowmapdiffeoarg{\muxmulevelsetvalue}))$,
where we stress that the norm on $W^{2,\infty}(\mbox{\upshape int}(\domaincomposedflowmapdiffeoarg{\muxmulevelsetvalue}))$ is with respect to the
coordinates $(\Delta u,t,x^2,x^3)$.
Composing $\phi \circ \composedflowmapdiffeoarg{\muxmulevelsetvalue}$ with $(\composedflowmapdiffeoarg{\muxmulevelsetvalue})^{-1}$
(in particular, using that $\composedflowmapdiffeoarg{\muxmulevelsetvalue}$ and $(\composedflowmapdiffeoarg{\muxmulevelsetvalue})^{-1}$
are $W^{2,\infty}$ diffeomorphisms 
by \eqref{E:W2INFTYBOUNDFORXMUISMINUSKAPPADIFFEOMORPHISM} 
and
\eqref{E:W2INFTYBOUNDFORXMUISMINUSKAPPAINVERSEDIFFEOMORPHISM}),
we see that
$
\| \phi \|_{W_{\textnormal{geo}}^{2,\infty}(\twoargMrough{(\timefunction_0,\timefunctionboot),(- \rightu,\leftu)}{\muxmulevelsetvalue})}
	\lesssim
	\| \ambient \|_{W_{\textnormal{geo}}^{2,\infty}(\twoargMrough{(\timefunction_0,\timefunctionboot),(- \rightu,\leftu)}{\muxmulevelsetvalue})}
$.

To complete the proof, it suffices for us to show:
\begin{align} \label{E:PROOFSTEPODESOLUTIONSTHATARESMOOTHERTHANTHEDATAHYPERSURFACE}
\vec{\partial}_{\textnormal{geo}} \phi 
&
\in 
W_{\textnormal{geo}}^{2,\infty}(\twoargMrough{(\timefunction_0,\timefunctionboot),(- \rightu,\leftu)}{\muxmulevelsetvalue}),
&
\| \vec{\partial}_{\textnormal{geo}} \phi \|_{W_{\textnormal{geo}}^{2,\infty}(\twoargMrough{(\timefunction_0,\timefunctionboot),(- \rightu,\leftu)}{\muxmulevelsetvalue})},
&
\lesssim
\| \ambient \|_{W_{\textnormal{geo}}^{3,\infty}(\twoargMrough{(\timefunction_0,\timefunctionboot),(- \rightu,\leftu)}{\muxmulevelsetvalue})},
\end{align}
where
$\vec{\partial}_{\textnormal{geo}} \phi \eqdef \left(\geop{t} \phi, \geop{u} \phi,\geop{x^2} \phi,\geop{x^3} \phi \right)$
is the array of geometric coordinate partial derivatives of $\phi$.
Due to the limited regularity 
$\embeddingdatahypfortimefunctionarg{\muxmulevelsetvalue}\in W^{2,\infty}(\mbox{\upshape int}(\domainforembeddingdatahypfortimefunctiontwoarg{\muxmulevelsetvalue}{(\upmuboot,\mupositive]}))$,
the desired regularity of $\phi$ cannot be inferred directly from the formula 
$
	\phi \circ \composedflowmapdiffeoarg{\muxmulevelsetvalue}(\Delta u,t,x^2,x^3)
	=
	\ambient \circ \embeddingdatahypfortimefunctionarg{\muxmulevelsetvalue}(t,x^2,x^3)
$. 
Instead,
we commute \eqref{E:WTRANSPORTINITIALVALUEPROBLEM}
with the geometric coordinate partial derivative vectorfields to 
deduce that $\vec{\partial}_{\textnormal{geo}} \phi$ satisfies the following initial value problem:
\begin{align} \label{E:COMMUTEDWTRANSPORTINITIALVALUEPROBLEM}
	\Wtransarg{\muxmulevelsetvalue} \vec{\partial}_{\textnormal{geo}} \phi
	& = \vec{I} \cdot \vec{\partial}_{\textnormal{geo}} \phi,
	&
	\vec{\partial}_{\textnormal{geo}} \phi|_{\datahypfortimefunctiontwoarg{-\muxmulevelsetvalue}{[\timefunction_0,\timefunctionboot)}}
	 = 
	\vec{\partial}_{\textnormal{geo}} \ambient|_{\datahypfortimefunctiontwoarg{-\muxmulevelsetvalue}{[\timefunction_0,\timefunctionboot)}},
\end{align}
where the validity of the initial condition on RHS~\eqref{E:COMMUTEDWTRANSPORTINITIALVALUEPROBLEM}
relies on the compatibility condition assumption $(\Wtransarg{\muxmulevelsetvalue} \ambient) \circ \embeddingdatahypfortimefunctionarg{\muxmulevelsetvalue}= 0$,
and by \eqref{E:W2INIFINTYBOUNDFORGEOMTRICPARTIALDERIVATIVESOFWTRANSCOMPONENTS},
we have
$
\| \vec{I} \|_{W_{\textnormal{geo}}^{2,\infty}(\twoargMrough{(\timefunction_0,\timefunctionboot),(- \rightu,\leftu)}{\muxmulevelsetvalue})}
\lesssim
1
$.
As in \eqref{E:ODEFORMULATIONOFWTRANSTRANSPORTEQUATION}--\eqref{E:DATAFORODEFORMULATIONOFWTRANSTRANSPORTEQUATION}, 
we can rewrite \eqref{E:COMMUTEDWTRANSPORTINITIALVALUEPROBLEM} as the following linear ODE system
in the unknowns $(\vec{\partial}_{\textnormal{geo}} \phi) \circ \composedflowmapdiffeoarg{\muxmulevelsetvalue}(\Delta u,t,x^2,x^3)$:
\begin{align} \label{E:ALTERNATECOMMUTEDWTRANSPORTINITIALVALUEPROBLEM}
	\frac{\partial}{\partial \Delta u} 
	\left[
		(\vec{\partial}_{\textnormal{geo}} \phi) \circ \composedflowmapdiffeoarg{\muxmulevelsetvalue}(\Delta u,t,x^2,x^3)
	\right]
	& = 
		\left[
		\vec{I} \circ \composedflowmapdiffeoarg{\muxmulevelsetvalue}(\Delta u,t,x^2,x^3)
		\right]
		\cdot 
		(\vec{\partial}_{\textnormal{geo}} \phi) \circ \composedflowmapdiffeoarg{\muxmulevelsetvalue}(\Delta u,t,x^2,x^3),
			\\
	(\vec{\partial}_{\textnormal{geo}} \phi) \circ \composedflowmapdiffeoarg{\muxmulevelsetvalue}(0,t,x^2,x^3)
	 & = 
	 (\vec{\partial}_{\textnormal{geo}} \ambient)
	\circ \embeddingdatahypfortimefunctionarg{\muxmulevelsetvalue}(t,x^2,x^3),
	\label{E:DATAFORALTERNATECOMMUTEDWTRANSPORTINITIALVALUEPROBLEM}
\end{align}
where our assumptions on $\vec{\partial}_{\textnormal{geo}} \ambient$ and the regularity of $\embeddingdatahypfortimefunctionarg{\muxmulevelsetvalue}$ imply that
$(\vec{\partial}_{\textnormal{geo}} \ambient) \circ \embeddingdatahypfortimefunctionarg{\muxmulevelsetvalue}\in W^{2,\infty}(\mbox{\upshape int}(\domainforembeddingdatahypfortimefunctiontwoarg{\muxmulevelsetvalue}{(\upmuboot,\mupositive]}))$
(here $\domainforembeddingdatahypfortimefunctiontwoarg{\muxmulevelsetvalue}{(\upmuboot,\mupositive]}$ is the domain of $ \embeddingdatahypfortimefunctionarg{\muxmulevelsetvalue}$)
with
$
\| (\vec{\partial}_{\textnormal{geo}} \ambient) \circ \embeddingdatahypfortimefunctionarg{\muxmulevelsetvalue}\|_{W^{2,\infty}(\mbox{\upshape int}(\domainforembeddingdatahypfortimefunctiontwoarg{\muxmulevelsetvalue}{(\upmuboot,\mupositive]}))}
\lesssim
\| \ambient \|_{W_{\textnormal{geo}}^{3,\infty}(\twoargMrough{(\timefunction_0,\timefunctionboot),(- \rightu,\leftu)}{\muxmulevelsetvalue})}
$.
The standard theory of transport equations with $W^{2,\infty}(\mbox{\upshape int}(\domaincomposedflowmapdiffeoarg{\muxmulevelsetvalue}))$ coefficients
yields that 
\eqref{E:ALTERNATECOMMUTEDWTRANSPORTINITIALVALUEPROBLEM}--\eqref{E:DATAFORALTERNATECOMMUTEDWTRANSPORTINITIALVALUEPROBLEM} 
has a unique solution
(which must be $(\vec{\partial}_{\textnormal{geo}} \phi) \circ \composedflowmapdiffeoarg{\muxmulevelsetvalue}$, 
where $\phi$ is the solution to \eqref{E:WTRANSPORTINITIALVALUEPROBLEM})
satisfying $(\vec{\partial}_{\textnormal{geo}} \phi) \circ \composedflowmapdiffeoarg{\muxmulevelsetvalue} \in W^{2,\infty}(\mbox{\upshape int}(\domaincomposedflowmapdiffeoarg{\muxmulevelsetvalue}))$.
Moreover, with the help of the above bound for 
$
\| \vec{\partial}_{\textnormal{geo}} \ambient \|_{W^{2,\infty}(\mbox{\upshape int}(\domainforembeddingdatahypfortimefunctiontwoarg{\muxmulevelsetvalue}{(\upmuboot,\mupositive]}))}$, 
a standard argument
based on commuting \eqref{E:ALTERNATECOMMUTEDWTRANSPORTINITIALVALUEPROBLEM} 
up to two times with respect to the partial derivatives
in the coordinate system $(\Delta u,t,x^2,x^3)$,
integrating with respect to $\Delta u$,
applying Gr\"{o}nwall's inequality,
and using that $|\Delta u| \leq \leftu + \rightu$ in the region under study
yields the bound: 
\begin{align} \label{E:ALMOSTDONETRANSPORTEQUATIONSOLUTIONW3INFINITYESTIMATEWITHDATAONROUGHHYPERSURFACE}
\left\| 
	(\vec{\partial}_{\textnormal{geo}} \phi) \circ \composedflowmapdiffeoarg{\muxmulevelsetvalue}
\right\|_{W^{2,\infty}(\domaincomposedflowmapdiffeoarg{\muxmulevelsetvalue})}
\lesssim
\| \ambient \|_{W_{\textnormal{geo}}^{3,\infty}(\twoargMrough{(\timefunction_0,\timefunctionboot),(- \rightu,\leftu)}{\muxmulevelsetvalue})}
\left\lbrace
1
+
\| \vec{I} \circ \composedflowmapdiffeoarg{\muxmulevelsetvalue} \|_{W_{\textnormal{geo}}^{2,\infty}(\twoargMrough{(\timefunction_0,\timefunctionboot),(- \rightu,\leftu)}{\muxmulevelsetvalue})}
\right\rbrace.
\end{align}
Finally, composing 
$(\vec{\partial}_{\textnormal{geo}} \phi) \circ \composedflowmapdiffeoarg{\muxmulevelsetvalue}$
with $(\composedflowmapdiffeoarg{\muxmulevelsetvalue})^{-1}$,
and using \eqref{E:W2INFTYBOUNDFORXMUISMINUSKAPPADIFFEOMORPHISM}--\eqref{E:W2INFTYBOUNDFORXMUISMINUSKAPPAINVERSEDIFFEOMORPHISM}
as well as the bound
$\| \vec{I} \|_{W_{\textnormal{geo}}^{2,\infty}(\twoargMrough{(\timefunction_0,\timefunctionboot),(- \rightu,\leftu)}{\muxmulevelsetvalue})}
\lesssim
1
$
and the bound \eqref{E:ALMOSTDONETRANSPORTEQUATIONSOLUTIONW3INFINITYESTIMATEWITHDATAONROUGHHYPERSURFACE},
we conclude \eqref{E:TRANSPORTEQUATIONSOLUTIONW3INFINITYESTIMATEWITHDATAONROUGHHYPERSURFACE}.

\end{proof}

\subsection{Estimates for the rough time function and the change of variables map
from rough adapted coordinates to geometric coordinates}
Using Lemma~\ref{L:ODESOLUTIONSTHATARESMOOTHERTHANTHEDATAHYPERSURFACE},
we now derive estimates 
for the rough time function,
the change of variables map
from geometric coordinates to rough adapted coordinates,
and its inverse.

We start with the following simple lemma, which provides an identity for
$\Dfour \timefunctionarg{\muxmulevelsetvalue}|_{\datahypfortimefunctiontwoarg{-\muxmulevelsetvalue}{[\timefunction_0,\timefunctionboot)}}$.

\begin{lemma}[Identity for $\Dfour \timefunctionarg{\muxmulevelsetvalue}$ along 
$\datahypfortimefunctiontwoarg{-\muxmulevelsetvalue}{[\timefunction_0,\timefunctionboot)}$]
	\label{L:GRADIENTOFTIMEFUNCTIONAGREESWITHGRADIENTOFMUALONGMUXMUEQUALSMINUSKAPPHYPERSURFACE}
	Recall that $\datahypfortimefunctiontwoarg{-\muxmulevelsetvalue}{[\timefunction_0,\timefunctionboot)}$ 
	is the truncated $\muX \upmu$-level set defined in
	\eqref{E:TRUNCATEDLEVELSETSOFMUXMU}.
	With $\Dfour \varphi$ denoting the spacetime gradient one-form of the scalar function $\varphi$, we have the following
	identity:
	\begin{align} \label{E:GRADIENTOFTIMEFUNCTIONAGREESWITHGRADIENTOFMUALONGMUXMUEQUALSMINUSKAPPHYPERSURFACE}
		\Dfour \timefunctionarg{\muxmulevelsetvalue}|_{\datahypfortimefunctiontwoarg{-\muxmulevelsetvalue}{[\timefunction_0,\timefunctionboot)}}
		& = 
		- \Dfour \upmu|_{\datahypfortimefunctiontwoarg{-\muxmulevelsetvalue}{[\timefunction_0,\timefunctionboot)}}.
	\end{align}
	
\end{lemma}

\begin{proof}
	The bootstrap assumptions of Sect.\,\ref{SSS:BOOTSTRAPASSUMPTIONFORTORISTRUCTURE} imply that
	$\datahypfortimefunctiontwoarg{-\muxmulevelsetvalue}{[\timefunction_0,\timefunctionboot)}$ is a hypersurface portion foliated by 
	$\lbrace \twoargmumuxtorus{\mulevelsetvalue}{-\muxmulevelsetvalue} \rbrace_{\mulevelsetvalue \in (\upmuboot,\mupositive]}$. 
	From these facts and
	\eqref{E:INITIALCONDITIONFORROUGHTIMEFUNCTION},
	it follows that $\timefunctionarg{\muxmulevelsetvalue}$ and $-\upmu$ have the same derivatives
	in directions tangent to $\datahypfortimefunctiontwoarg{-\muxmulevelsetvalue}{[\timefunction_0,\timefunctionboot)}$.
	Moreover, since 
	definition \eqref{E:WTRANSDEF} implies that $\Wtransarg{\muxmulevelsetvalue} \upmu|_{\datahypfortimefunctionarg{-\muxmulevelsetvalue}} = 0$,
	we see from \eqref{E:TRANSPORTEQUATIONFORROUGHTIMEFUNCTION} that
	along $\datahypfortimefunctiontwoarg{-\muxmulevelsetvalue}{[\timefunction_0,\timefunctionboot)}$,
	$\timefunctionarg{\muxmulevelsetvalue}$ and $-\upmu$ have the same $\Wtransarg{\muxmulevelsetvalue}$-derivative.
	Since \eqref{E:BAMUTRANSVERSALCONVEXITY} and
	\eqref{E:BOOTSTRAPLEVELSETSTRUCTUREANDLOCATIONOFMUXEQUALSMINUSKAPPA} imply that
	$\Wtransarg{\muxmulevelsetvalue}$ is transversal to $\datahypfortimefunctiontwoarg{-\muxmulevelsetvalue}{[\timefunction_0,\timefunctionboot)}$,
	we conclude the desired identity \eqref{E:GRADIENTOFTIMEFUNCTIONAGREESWITHGRADIENTOFMUALONGMUXMUEQUALSMINUSKAPPHYPERSURFACE}.
\end{proof}

\begin{lemma}[Estimates for $\timefunctionarg{\muxmulevelsetvalue}$, $\CHOVgeotorough{\muxmulevelsetvalue}$, and $\InverseCHOVgeotorough{\muxmulevelsetvalue}$]
		\label{L:LINFTYESTIMATESFORROUGHTIMEFUNCTIONANDDERIVATIVES}
		The following estimates hold.
		
		\medskip
		
	\noindent \underline{\textbf{Estimates for $\timefunctionarg{\muxmulevelsetvalue}$}}:
		\begin{subequations}
		\begin{align}
			\left\| 
				\timefunctionarg{\muxmulevelsetvalue}
			\right\|_{W_{\textnormal{geo}}^{3,\infty}(\twoargMrough{(\timefunction_0,\timefunctionboot),(- \rightu,\leftu)}{\muxmulevelsetvalue})}
		& 
		\lesssim 1,
			\label{E:ALLDERIVATIVESLINFTYESTIMATESFORROUGHTIMEFUNCTIONANDDERIVATIVES}
			\\
			\left\| 
				\left(\geop{x^2} \timefunctionarg{\muxmulevelsetvalue},\geop{x^3} \timefunctionarg{\muxmulevelsetvalue} \right)  
			\right\|_{W_{\textnormal{geo}}^{2,\infty}(\twoargMrough{(\timefunction_0,\timefunctionboot),(- \rightu,\leftu)}{\muxmulevelsetvalue})}
			& 
			\lesssim \auxbootsmall.
			\label{E:SMALLDERIVATIVESLINFTYESTIMATESFORROUGHTIMEFUNCTIONANDDERIVATIVES}
		\end{align}
		\end{subequations}
	
		Moreover, 
		\begin{subequations}
		\begin{align} \label{E:PARTIALTIMEDERIVATIVEOFROUGHTIMEFUNCTIONISAPPROXIMATELYUNITY}
				\frac{1}{2} \updelta_*	
				&
				\leq
				\geop{t} \timefunctionarg{\muxmulevelsetvalue}
				\leq 
				\frac{3}{2} \updelta_*,
				&&
				\mbox{on } \twoargMrough{[\timefunction_0,\timefunctionboot),[- \rightu,\leftu]}{\muxmulevelsetvalue},	
					\\
				\frac{1}{2} \updelta_*	
				&
				\leq
				\Lunit \timefunctionarg{\muxmulevelsetvalue}
				\leq 
				\frac{3}{2} \updelta_*,
				&&
				\mbox{on } \twoargMrough{[\timefunction_0,\timefunctionboot),[- \rightu,\leftu]}{\muxmulevelsetvalue}.
			\label{E:LDERIVATIVEOFROUGHTIMEFUNCTIONISAPPROXIMATELYUNITY}
		\end{align}
		\end{subequations}

\medskip

\noindent \underline{\textbf{Estimates for $\CHOVgeotorough{\muxmulevelsetvalue}$ and $\InverseCHOVgeotorough{\muxmulevelsetvalue}$}}:
		The following estimate holds,
		where $\CHOVgeotorough{\muxmulevelsetvalue}$ is the change of variables map from geometric coordinates
		to rough adapted coordinates defined in \eqref{E:CHOVGEOTOROUGH}:
		\begin{align} \label{E:W3INFINITYBOUNDFORCHOVFROMGEOTOROUGHONOPEN}
			\left\| 
				\CHOVgeotorough{\muxmulevelsetvalue}
			\right\|_{W_{\textnormal{geo}}^{3,\infty}(\twoargMrough{(\timefunction_0,\timefunctionboot),(- \rightu,\leftu)}{\muxmulevelsetvalue})}
		& 
		\lesssim 1.
		\end{align}

		The following estimates hold, where
		$d_{\textnormal{geo}} \CHOVgeotorough{\muxmulevelsetvalue}$ 
		is the Jacobian matrix of $\CHOVgeotorough{\muxmulevelsetvalue}$:
		\begin{align} \label{E:KEYJACOBIANDETERMINANTESTIMATECHOVGEOTOROUGH}
				\mbox{\upshape det} \left(d_{\textnormal{geo}} \CHOVgeotorough{\muxmulevelsetvalue} \right)
				& \approx 1,
				&&
				\mbox{on } \twoargMrough{[\timefunction_0,\timefunctionboot),[- \rightu,\leftu]}{\muxmulevelsetvalue},
					\\
			\left\|
				[d_{\textnormal{geo}} \CHOVgeotorough{\muxmulevelsetvalue}]^{-1}
			\right\|_{W^{2,\infty} \left(\twoargMrough{[\timefunction_0,\timefunctionboot),[- \rightu,\leftu]}{\muxmulevelsetvalue} \right)}
			& \leq C.
			\label{E:INVERSEOFJACOBIANMATRIXFROMGEOTOROUGHW2INFINITYBOUND}
		\end{align}
		
		In addition, the following estimates hold, 
		where $\InverseCHOVgeotorough{\muxmulevelsetvalue}$ is the change of variables map from rough adapted coordinates to geometric coordinates
		and
		$d_{\textnormal{rough}} \InverseCHOVgeotorough{\muxmulevelsetvalue}$ 
		is its Jacobian matrix:
			\begin{align} \label{E:KEYJACOBIANDETERMINANTESTIMATECHOVROUGHTOGEO}
				\roughgeop{\timefunction} t 
				& =
				\mbox{\upshape det} \left(d_{\textnormal{rough}} \InverseCHOVgeotorough{\muxmulevelsetvalue} \right)
				\approx 1,
				&&
				\mbox{on } [\timefunction_0,\timefunctionboot) \times [- \rightu,\leftu] \times \mathbb{T}^2.
		\end{align}
		
		In addition, the following estimates hold:
		\begin{align} \label{E:W3INFINITYBOUNDFORINVERSECHOVGEOTOROUGH}
			\left\| 
				\InverseCHOVgeotorough{\muxmulevelsetvalue}
			\right\|_{W_{\textnormal{rough}}^{3,\infty}((\timefunction_0,\timefunctionboot) \times (- \rightu,\leftu) \times \mathbb{T}^2)}
			\leq C.
		\end{align}
		
		Finally, $\InverseCHOVgeotorough{\muxmulevelsetvalue}$ extends to 
		a $C^{2,1}([\timefunction_0,\timefunctionboot] \times [- \rightu,\leftu] \times \mathbb{T}^2)$
		map satisfying the estimate:
		\begin{align} \label{E:C21BOUNDFORINVERSECHOVGEOTOROUGH}
			\left\| 
				\InverseCHOVgeotorough{\muxmulevelsetvalue}
			\right\|_{C_{\textnormal{rough}}^{2,1}([\timefunction_0,\timefunctionboot] \times [- \rightu,\leftu] \times \mathbb{T}^2)}
			\leq C.
		\end{align}
		
\end{lemma}		
	
	\begin{proof}
	Recall from Def.\,\ref{D:ROUGHTIMEFUNCTION} that 
	$
	\Wtransarg{\muxmulevelsetvalue} \timefunctionarg{\muxmulevelsetvalue}
	= 0
	$
	and
	$
	\timefunctionarg{\muxmulevelsetvalue}|_{\datahypfortimefunctiontwoarg{-\muxmulevelsetvalue}{[\timefunction_0,\timefunctionboot)}}
	=
	- \upmu|_{\datahypfortimefunctiontwoarg{-\muxmulevelsetvalue}{[\timefunction_0,\timefunctionboot)}}
	$.
	Hence, \eqref{E:ALLDERIVATIVESLINFTYESTIMATESFORROUGHTIMEFUNCTIONANDDERIVATIVES}
	follows from applying Lemma~\ref{L:ODESOLUTIONSTHATARESMOOTHERTHANTHEDATAHYPERSURFACE}
	with $\ambient \eqdef \upmu$ and the bound
	$\| \upmu \|_{W_{\textnormal{geo}}^{3,\infty}(\twoargMrough{(\timefunction_0,\timefunctionboot),(- \rightu,\leftu)}{\muxmulevelsetvalue})} 
	\lesssim 1
	$,
	which follows from Lemma~\ref{L:COMMUTATORSTOCOORDINATES} and the bootstrap assumptions.
	
	Similarly, 
	\eqref{E:SMALLDERIVATIVESLINFTYESTIMATESFORROUGHTIMEFUNCTIONANDDERIVATIVES}
	follows from
	\eqref{E:TANGENTIALDERIVATIVESETIMATETRANSPORTEQUATIONSOLUTIONW3INFINITYESTIMATEWITHDATAONROUGHHYPERSURFACE}
	and the bound
	$
	\left\|
		\left(\geop{x^2} \upmu,\geop{x^3} \upmu \right)  
	\right\|_{W_{\textnormal{geo}}^{2,\infty}(\twoargMrough{(\timefunction_0,\timefunctionboot),(- \rightu,\leftu)}{\muxmulevelsetvalue})}
	\lesssim
	\auxbootsmall
$,
which follows from Lemma~\ref{L:COMMUTATORSTOCOORDINATES} and the bootstrap assumptions.

\eqref{E:W3INFINITYBOUNDFORCHOVFROMGEOTOROUGHONOPEN}
follows from
\eqref{E:ALLDERIVATIVESLINFTYESTIMATESFORROUGHTIMEFUNCTIONANDDERIVATIVES}
and the definition of $\CHOVgeotorough{\muxmulevelsetvalue}$.

To prove \eqref{E:PARTIALTIMEDERIVATIVEOFROUGHTIMEFUNCTIONISAPPROXIMATELYUNITY}
and \eqref{E:LDERIVATIVEOFROUGHTIMEFUNCTIONISAPPROXIMATELYUNITY},
we first commute the equation 
$
	\Wtransarg{\muxmulevelsetvalue} \timefunctionarg{\muxmulevelsetvalue}
	= 0
	$
with $\Lunit$ to obtain the transport equation
$
\Wtransarg{\muxmulevelsetvalue} \Lunit \timefunctionarg{\muxmulevelsetvalue}
= [\Wtransarg{\muxmulevelsetvalue},\Lunit] \timefunctionarg{\muxmulevelsetvalue}
$.
Using \eqref{E:WTRANSDEF},
Prop.\,\ref{P:COMMUTATORESTIMATES},
the bootstrap assumptions, and \eqref{E:SMALLDERIVATIVESLINFTYESTIMATESFORROUGHTIMEFUNCTIONANDDERIVATIVES},
we see that the source term can be pointwise 
bounded on $\twoargMrough{[\timefunction_0,\timefunctionboot),[- \rightu,\leftu]}{\muxmulevelsetvalue}$ as follows:
$
	\left|[\Wtransarg{\muxmulevelsetvalue},\Lunit] \timefunctionarg{\muxmulevelsetvalue} \right|
	\lesssim \auxbootsmall
$.
Moreover, since Lemma~\ref{L:GRADIENTOFTIMEFUNCTIONAGREESWITHGRADIENTOFMUALONGMUXMUEQUALSMINUSKAPPHYPERSURFACE} implies
that
$\Lunit \timefunctionarg{\muxmulevelsetvalue}|_{\datahypfortimefunctiontwoarg{-\muxmulevelsetvalue}{[\timefunction_0,\timefunctionboot)}} 
= 
-
\Lunit \upmu|_{\datahypfortimefunctiontwoarg{-\muxmulevelsetvalue}{[\timefunction_0,\timefunctionboot)}}
$,
we can use the bootstrap assumptions
\eqref{E:BABOUNDSONLMUINTERESTINGREGION}
and
\eqref{E:BOOTSTRAPLEVELSETSTRUCTUREANDLOCATIONOFMUXEQUALSMINUSKAPPA}
to deduce that
$
- 
	\frac{5}{4}
	\updelta_*
	\leq
	-
	\Lunit \timefunctionarg{\muxmulevelsetvalue}|_{\datahypfortimefunctiontwoarg{-\muxmulevelsetvalue}{[\timefunction_0,\timefunctionboot)}}
	\leq 
	- \frac{3}{4} \updelta_*
$.
Hence, recalling that $\Wtransarg{\muxmulevelsetvalue} u = 1$, 
we can integrate the transport equation for
$\Lunit \timefunctionarg{\muxmulevelsetvalue}$
starting from the $\Wtransarg{\muxmulevelsetvalue}$-transversal data hypersurface 
$
\datahypfortimefunctiontwoarg{-\muxmulevelsetvalue}{[\timefunction_0,\timefunctionboot)}
$
(see \eqref{E:BAMUTRANSVERSALCONVEXITY})
and use Gr\"{o}nwall's inequality and the fact that $|u| \lesssim 1$ 
on the region under study to conclude that:
\begin{align} \label{E:LUNITTIMEFUNCTIONUNIFORMLYPOSITIVE}
	- 
	\frac{5}{4} \updelta_*
	- 
	C \auxbootsmall 
	&
	\leq
	\min_{\twoargMrough{[\timefunction_0,\timefunctionboot),[- \rightu,\leftu]}{\muxmulevelsetvalue}} - \Lunit \timefunctionarg{\muxmulevelsetvalue}
	\leq 
	\max_{\twoargMrough{[\timefunction_0,\timefunctionboot),[- \rightu,\leftu]}{\muxmulevelsetvalue}} - \Lunit \timefunctionarg{\muxmulevelsetvalue}
	\leq 
	- 
	\frac{3}{4} \updelta_* 
	+ 
	C \auxbootsmall,
\end{align}
which yields \eqref{E:LDERIVATIVEOFROUGHTIMEFUNCTIONISAPPROXIMATELYUNITY}.
Finally, using \eqref{E:LUNITTIMEFUNCTIONUNIFORMLYPOSITIVE},
the identity
$
\geop{t} \timefunctionarg{\muxmulevelsetvalue}
=
\Lunit \timefunctionarg{\muxmulevelsetvalue}
- \Lunit^A \geop{x^A} \timefunctionarg{\muxmulevelsetvalue}
$,
and the bound $\left|\Lunit^A \geop{x^A} \timefunctionarg{\muxmulevelsetvalue} \right| \lesssim \auxbootsmall$
implied by the bootstrap assumptions and \eqref{E:SMALLDERIVATIVESLINFTYESTIMATESFORROUGHTIMEFUNCTIONANDDERIVATIVES},
we conclude \eqref{E:PARTIALTIMEDERIVATIVEOFROUGHTIMEFUNCTIONISAPPROXIMATELYUNITY}.

\eqref{E:KEYJACOBIANDETERMINANTESTIMATECHOVGEOTOROUGH} now follows from
\eqref{E:PARTIALTIMEDERIVATIVEOFROUGHTIMEFUNCTIONISAPPROXIMATELYUNITY} 
and the simple
identity
$
\geop{t} \timefunctionarg{\muxmulevelsetvalue} 
				=
				\mbox{\upshape det} \left(d_{\textnormal{geo}} \CHOVgeotorough{\muxmulevelsetvalue} \right)
$,
which
follows easily from the definition of $\CHOVgeotorough{\muxmulevelsetvalue}$.

\eqref{E:INVERSEOFJACOBIANMATRIXFROMGEOTOROUGHW2INFINITYBOUND}
now follows from \eqref{E:W3INFINITYBOUNDFORCHOVFROMGEOTOROUGHONOPEN}
and
\eqref{E:KEYJACOBIANDETERMINANTESTIMATECHOVGEOTOROUGH}.

The ``$=$'' in \eqref{E:KEYJACOBIANDETERMINANTESTIMATECHOVROUGHTOGEO}
follows easily from the definition of $\CHOVgeotorough{\muxmulevelsetvalue}$.
The ``$\approx$'' in \eqref{E:KEYJACOBIANDETERMINANTESTIMATECHOVROUGHTOGEO}
follows from the identity
$$
\mbox{\upshape det} \left([d_{\textnormal{rough}} \InverseCHOVgeotorough{\muxmulevelsetvalue}] \circ \CHOVgeotorough{\muxmulevelsetvalue} \right)
=
\left(\mbox{\upshape det} [d_{\textnormal{geo}} \CHOVgeotorough{\muxmulevelsetvalue}] \right)^{-1}
$$
and \eqref{E:INVERSEOFJACOBIANMATRIXFROMGEOTOROUGHW2INFINITYBOUND}.

\eqref{E:W3INFINITYBOUNDFORINVERSECHOVGEOTOROUGH}
follows from differentiating 
(up to two times, with the rough adapted coordinate partial derivative vectorfields) the identity
$
d_{\textnormal{rough}} \InverseCHOVgeotorough{\muxmulevelsetvalue}
=
[d_{\textnormal{geo}} \CHOVgeotorough{\muxmulevelsetvalue}]^{-1} \circ [\InverseCHOVgeotorough{\muxmulevelsetvalue}]
$
and using 
\eqref{E:BASIZEOFCARTESIANT},
\eqref{E:INVERSEOFJACOBIANMATRIXFROMGEOTOROUGHW2INFINITYBOUND},
and the chain rule.

\eqref{E:C21BOUNDFORINVERSECHOVGEOTOROUGH}
follows from
\eqref{E:W3INFINITYBOUNDFORINVERSECHOVGEOTOROUGH},
and the following Sobolev embedding result for scalar functions $f$
(see the proof of \cite{lE1998}*{Theorem~5 in Section~5.6}),
which relies on the convexity of the domain $(\timefunction_0,\timefunctionboot) \times (- \rightu,\leftu) \times \mathbb{T}^2$:
$
\| f \|_{C_{\textnormal{rough}}^{0,1}([\timefunction_0,\timefunctionboot] \times [- \rightu,\leftu] \times \mathbb{T}^2)}
\leq C 
\| f \|_{W_{\textnormal{rough}}^{1,\infty}((\timefunction_0,\timefunctionboot) \times (- \rightu,\leftu) \times \mathbb{T}^2)}
$.

\end{proof}

In the next lemma, we continue our analysis of the change of variables map
$\CHOVgeotorough{\muxmulevelsetvalue}$. More precisely, we exhibit its properties on the closure
of $\twoargMrough{[\timefunction_0,\timefunctionboot),[- \rightu,\leftu]}{\muxmulevelsetvalue}$.
We also derive the quasi-convexity of the closure of
$\twoargMrough{[\timefunction_0,\timefunctionboot),[- \rightu,\leftu]}{\muxmulevelsetvalue}$
and, as a consequence, prove a standard Sobolev embedding result.

\begin{lemma}[Properties of $\CHOVgeotorough{\muxmulevelsetvalue}$ on the closure of 
		$\twoargMrough{[\timefunction_0,\timefunctionboot),[- \rightu,\leftu]}{\muxmulevelsetvalue}$ and quasi-convexity]
		\label{L:DIFFEOMORPHICEXTENSIONOFROUGHCOORDINATES}
		The following results hold.
		
		\begin{enumerate}
			\item $\twoargMrough{[\timefunction_0,\timefunctionboot),[- \rightu,\leftu]}{\muxmulevelsetvalue}$ 
				is precompact in the topology of the geometric coordinates $(t,u,x^2,x^3)$.
			\item The change of variables map $\CHOVgeotorough{\muxmulevelsetvalue}(t,u,x^2,x^3) = (\timefunctionarg{\muxmulevelsetvalue},u,x^2,x^3)$
				extends to a $C^{2,1}$ diffeomorphism on the closure of
				$\twoargMrough{[\timefunction_0,\timefunctionboot),[- \rightu,\leftu]}{\muxmulevelsetvalue}$,
				which we denote by $\mbox{\upshape cl} \left(\twoargMrough{[\timefunction_0,\timefunctionboot),[- \rightu,\leftu]}{\muxmulevelsetvalue} \right)$.
				Moreover,
				 $\mbox{\upshape cl} \left(\twoargMrough{[\timefunction_0,\timefunctionboot),[- \rightu,\leftu]}{\muxmulevelsetvalue} \right)
				= \twoargMrough{[\timefunction_0,\timefunctionboot],[- \rightu,\leftu]}{\muxmulevelsetvalue}$,
				and
				$\CHOVgeotorough{\muxmulevelsetvalue}\left(\twoargMrough{[\timefunction_0,\timefunctionboot],[- \rightu,\leftu]}{\muxmulevelsetvalue} \right) 
			= [\timefunction_0,\timefunctionboot] \times [- \rightu,\leftu] \times \mathbb{T}^2$.
		\item The following estimates hold for the extended maps:
		\begin{subequations}
		\begin{align} \label{E:CLOSEDVERSIONKEYJACOBIANDETERMINANTESTIMATECHOVGEOTOROUGH}
				\frac{1}{2} \updelta_*	
				&
				\leq
				\geop{t} \timefunctionarg{\muxmulevelsetvalue}
				\leq 
				\frac{3}{2} \updelta_*,
				&&
				\mbox{on } \twoargMrough{[\timefunction_0,\timefunctionboot],[- \rightu,\leftu]}{\muxmulevelsetvalue},	
					\\
				\frac{1}{2} \updelta_*	
				&
				\leq
				\Lunit \timefunctionarg{\muxmulevelsetvalue}
				\leq 
				\frac{3}{2} \updelta_*,
				&&
				\mbox{on } \twoargMrough{[\timefunction_0,\timefunctionboot],[- \rightu,\leftu]}{\muxmulevelsetvalue}.
			\label{E:CLOSEDVERSIONLUNITROUGHTTIMEFUNCTION}
		\end{align}
		\end{subequations}
		\item 	\begin{align} \label{E:CLOSEDKEYJACOBIANDETERMINANTESTIMATECHOVROUGHTOGEO}
				\roughgeop{\timefunction} t 
				& =
				\mbox{\upshape det} \left(d_{\textnormal{rough}} \InverseCHOVgeotorough{\muxmulevelsetvalue} \right)
				\approx 1,
				&&
				\mbox{on } [\timefunction_0,\timefunctionboot] \times [- \rightu,\leftu] \times \mathbb{T}^2,
		\end{align}
		\item \begin{align}
							\left\| 
								\CHOVgeotorough{\muxmulevelsetvalue}
							\right\|_{C_{\textnormal{geo}}^{2,1}(\twoargMrough{[\timefunction_0,\timefunctionboot],[- \rightu,\leftu]}{\muxmulevelsetvalue})}
							& \leq C,
								\label{E:CLOSEDVERSIONC21BOUNDFORCHOVROUGHTOGEO} 
									\\
							\left\| 
								\CHOVgeotorough{\muxmulevelsetvalue}^{-1}
							\right\|_{C_{\textnormal{rough}}^{2,1}([\timefunction_0,\timefunctionboot] \times [- \rightu,\leftu] \times \mathbb{T}^2)}
							& \leq C,
							\label{E:CLOSEDVERSIONC21BOUNDFORINVERSECHOVGEOTOROUGH}
								\\
							\left\| 
								\left(\geop{x^2} \timefunctionarg{\muxmulevelsetvalue},\geop{x^3} \timefunctionarg{\muxmulevelsetvalue} \right)
							\right\|_{C_{\textnormal{geo}}^{1,1}(\twoargMrough{[\timefunction_0,\timefunctionboot],[- \rightu,\leftu]}{\muxmulevelsetvalue})}
							& \leq C \varepsilon^{1/2}.
							\label{E:SMALLC11ESTIMATESFORROUGHTIMEFUNCTION}
						\end{align}
				\item (\textbf{Quasi-convexity of} 
				$\twoargMrough{[\timefunction_0,\timefunctionboot],[- \rightu,\leftu]}{\muxmulevelsetvalue}$). 
					For every pair of points 
					$q_1,q_2 \in [\timefunction_0,\timefunctionboot] \times [- \rightu,\leftu] \times \mathbb{T}^2$,
					we have:
					\begin{align} \label{E:COMPARABLEDISTANCES}
						\mbox{\upshape dist}_{\mbox{\upshape flat}}
						\left(\CHOVgeotorough{\muxmulevelsetvalue}^{-1}(q_1),\CHOVgeotorough{\muxmulevelsetvalue}^{-1}(q_2) \right)
						& \approx
						\mbox{\upshape dist}_{\mbox{\upshape flat}}(q_1,q_2),
					\end{align}
					where $\mbox{\upshape dist}_{\mbox{\upshape flat}}(q_1,q_2)$ is the
					standard Euclidean distance between $q_1$ and $q_2$ in the flat space 
					$\mathbb{R}_{\timefunction} \times \mathbb{R}_u \times \mathbb{T}^2$.
					
					Moreover, $\twoargMrough{[\timefunction_0,\timefunctionboot],[- \rightu,\leftu]}{\muxmulevelsetvalue}$ is quasi-convex.\footnote{Here, 
					we are not just interested in a qualitative version of quasi-convexity, but rather in obtaining control over the constants ``$C$.''
					Similar remarks apply for the Sobolev embedding result \eqref{E:SOBOELVEMBEDDINGRELYINGONQUASICONVEXITY} and for other 
					quasi-convexity and Sobolev embedding results derived throughout the paper. \label{FN:QUANTITATIVEQUASICONVEXITY}}
					That is,
					every pair of points
					$p_1,p_2 \in \twoargMrough{[\timefunction_0,\timefunctionboot],[- \rightu,\leftu]}{\muxmulevelsetvalue}$
					are connected by a $C_{\textnormal{geo}}^1$ curve in $\twoargMrough{[\timefunction_0,\timefunctionboot],[- \rightu,\leftu]}{\muxmulevelsetvalue}$
					whose length with respect to the standard flat Euclidean metric on geometric coordinate space
					$\mathbb{R}_t \times \mathbb{R}_u \times \mathbb{T}^2$
					is $\leq C \mbox{\upshape dist}_{\mbox{\upshape flat}}(p_1,p_2)$.
				\item (\textbf{Sobolev embedding}).
					There is a constant $C > 0$ that is \textbf{independent of $\timefunctionboot$}
					such that the following Sobolev embedding result holds for scalar functions $f$ on 
					$\twoargMrough{(\timefunction_0,\timefunctionboot),(- \rightu,\leftu)}{\muxmulevelsetvalue}$:
					\begin{align} \label{E:SOBOELVEMBEDDINGRELYINGONQUASICONVEXITY}
						\| f \|_{C_{\textnormal{geo}}^{0,1}(\twoargMrough{[\timefunction_0,\timefunctionboot],[- \rightu,\leftu]}{\muxmulevelsetvalue})}
						& \leq
						C
						\| f \|_{W_{\textnormal{geo}}^{1,\infty}(\twoargMrough{(\timefunction_0,\timefunctionboot),(- \rightu,\leftu)}{\muxmulevelsetvalue})}.
					\end{align}
				\item ($\hypthreearg{\timefunction}{[- \rightu,\leftu]}{\muxmulevelsetvalue}$ \textbf{is a graph}).
					For $\timefunction \in [\timefunction_0,\timefunctionboot]$,
					there exists a function 
					$\Cartesiantisafunctiononlevelsetsofroughtimefunctionarg{\timefunction}{\muxmulevelsetvalue}: 
					[- \rightu,\leftu] \times \mathbb{T}^2 \rightarrow \mathbb{R}$,
					depending on $\timefunction$ and $\muxmulevelsetvalue$,
					such that:
					\begin{align} \label{E:C21BOUNDFORCONSTANTTIMEFUNCTIONGRAPH}
						\| \Cartesiantisafunctiononlevelsetsofroughtimefunctionarg{\timefunction}{\muxmulevelsetvalue} \|_{C^{2,1}([- \rightu,\leftu] \times \mathbb{T}^2)}
						& \leq C
					\end{align}
					and such that relative to the geometric coordinates, we have:
					\begin{align} \label{E:LEVELSETSOFTIMEFUNCTIONAREAGRAPH}
						\hypthreearg{\timefunction}{[- \rightu,\leftu]}{\muxmulevelsetvalue}
						& =
						\left\lbrace
							(t,u,x^2,x^3)
							\ | \
							t = \Cartesiantisafunctiononlevelsetsofroughtimefunctionarg{\timefunction}{\muxmulevelsetvalue}(u,x^2,x^3),
								\,
							(u,x^2,x^3) \in [- \rightu,\leftu] \times \mathbb{T}^2
						\right\rbrace.
					\end{align}
		\end{enumerate}

\end{lemma}

\begin{proof}
	In \eqref{E:C21BOUNDFORINVERSECHOVGEOTOROUGH} and just above it,
	we showed that $\InverseCHOVgeotorough{\muxmulevelsetvalue}$
	extends as a $C_{\textnormal{rough}}^{2,1}$ function to the compact, convex domain
	$[\timefunction_0,\timefunctionboot] \times [- \rightu,\leftu] \times \mathbb{T}^2$
	such that \eqref{E:KEYJACOBIANDETERMINANTESTIMATECHOVROUGHTOGEO} holds on
	$[\timefunction_0,\timefunctionboot] \times [- \rightu,\leftu] \times \mathbb{T}^2$.
	In particular, this yields \eqref{E:CLOSEDKEYJACOBIANDETERMINANTESTIMATECHOVROUGHTOGEO}.
	From these facts,
	\eqref{E:PARTIALTIMEDERIVATIVEOFROUGHTIMEFUNCTIONISAPPROXIMATELYUNITY},
	and the fact that $(t,u,x^2,x^3) = \InverseCHOVgeotorough{\muxmulevelsetvalue}(\timefunction,u,x^2,x^3)$,
	we conclude that the map $\InverseCHOVgeotorough{\muxmulevelsetvalue}$ 
	with domain $[\timefunction_0,\timefunctionboot] \times [- \rightu,\leftu] \times \mathbb{T}^2$
	has a global inverse,
	i.e., that $\CHOVgeotorough{\muxmulevelsetvalue}$ extends to the compact domain
	$\twoargMrough{[\timefunction_0,\timefunctionboot],[- \rightu,\leftu]}{\muxmulevelsetvalue}$
	as an invertible map such that $\timefunctionarg{\muxmulevelsetvalue}$ satisfies  
	\eqref{E:CLOSEDVERSIONKEYJACOBIANDETERMINANTESTIMATECHOVGEOTOROUGH}.

	We now prove \eqref{E:COMPARABLEDISTANCES}.
	For $i=1,2$, we set
	$q_i \eqdef (\timefunction_i,u_i,x_i^2,x_i^3)$
	and
	$p_i \eqdef \CHOVgeotorough{\muxmulevelsetvalue}^{-1}(q_i) \eqdef (t_i,u_i,x_i^2,x_i^3)$.
	We define 
	$\Delta \timefunction \eqdef \timefunction_2 - \timefunction_1$,
	$\Delta u \eqdef u_2 - u_1$,
	$\Delta x^2 \eqdef x_2^2 - x_1^2$,
	$\Delta x^3 \eqdef x_2^3 - x_1^3$, 
	and
	$
	\Delta q \eqdef (\Delta \timefunction,\Delta u,\Delta x^2,\Delta x^3) = q_2 - q_1
	$.
	Similarly, we define $\Delta t \eqdef t_2 - t_1$ and
	$
	\Delta p 
	\eqdef (\Delta t,\Delta u,\Delta x^2,\Delta x^3) = p_2 - p_1
	$.
	Without loss of generality, we assume 
	$\Delta \timefunction \geq 0$.
	Then by 
	\eqref{E:CLOSEDKEYJACOBIANDETERMINANTESTIMATECHOVROUGHTOGEO}
	and
	\eqref{E:CLOSEDVERSIONC21BOUNDFORINVERSECHOVGEOTOROUGH},
	there is a constant $C > 1$ such that
	$
	\frac{1}{C} \Delta \timefunction 
	-
	C(|\Delta u| + |\Delta x^2|_{\mathbb{T}} + |\Delta x^3|_{\mathbb{T}})
	\leq
	\Delta t 
	\leq C \Delta \timefunction + C(|\Delta u| + |\Delta x^2|_{\mathbb{T}} + |\Delta x^3|_{\mathbb{T}})
	$,
	where for $j=2,3$,
	$|\Delta x^j|_{\mathbb{T}}$ is the Euclidean distance between $x_2^j$ and $x_1^j$ in the torus.
	Hence, 
	with $|\Delta p|_{Taxi} \eqdef |\Delta t| + |\Delta u| + |\Delta x^2|_{\mathbb{T}} + |\Delta x^3|_{\mathbb{T}}$,
	we see that there exists a (different) constant $C > 1$ such that:
	\begin{align} \label{E:QUANTITATIVESTEPINPROOFOFCOMPARABLEDISTANCES}
	\left| 
	\frac{1}{C} \Delta \timefunction 
	-
	C\left(|\Delta u| + |\Delta x^2|_{\mathbb{T}} + |\Delta x^3|_{\mathbb{T}} \right)
	\right|
	+
	|\Delta u| + |\Delta x^2|_{\mathbb{T}} + |\Delta x^3|_{\mathbb{T}}
	& 
	\leq
	|\Delta p|_{Taxi}
	\leq
	C\left(|\Delta \timefunction| + |\Delta u| + |\Delta x^2|_{\mathbb{T}} + |\Delta x^3|_{\mathbb{T}} \right).
	\end{align}
	From \eqref{E:QUANTITATIVESTEPINPROOFOFCOMPARABLEDISTANCES}, 
	it follows that 
	$
	|\Delta p|_{Taxi}
	\approx
	|\Delta q|_{Taxi}
	\eqdef
	|\Delta \timefunction| + |\Delta u| + |\Delta x^2|_{\mathbb{T}} + |\Delta x^3|_{\mathbb{T}},
	$
	which implies \eqref{E:COMPARABLEDISTANCES}.
	
	We now prove the quasi-convexity of
	$\twoargMrough{[\timefunction_0,\timefunctionboot],[- \rightu,\leftu]}{\muxmulevelsetvalue}$.
	Let $p_1,p_2 \in \twoargMrough{[\timefunction_0,\timefunctionboot],[- \rightu,\leftu]}{\muxmulevelsetvalue}$,
	and let $q_1,q_2 \in [\timefunction_0,\timefunctionboot] \times [- \rightu,\leftu] \times \mathbb{T}^2$ be the unique points
	such that $p_i = \CHOVgeotorough{\muxmulevelsetvalue}^{-1}(q_i)$,
	as above. Let $\ell$ be a straight line in 
	$[\timefunction_0,\timefunctionboot] \times [- \rightu,\leftu] \times \mathbb{T}^2$
	whose flat length is equal to $\mbox{\upshape dist}_{\mbox{\upshape flat}}(q_1,q_2)$.
	From 
	\eqref{E:CLOSEDVERSIONC21BOUNDFORINVERSECHOVGEOTOROUGH},
	it follows that the image curve $\CHOVgeotorough{\muxmulevelsetvalue}^{-1}(\ell)$
	has a flat length that is $\lesssim \mbox{\upshape dist}_{\mbox{\upshape flat}}(q_1,q_2)$,
	which by \eqref{E:COMPARABLEDISTANCES} is 
	$\approx \mbox{\upshape dist}_{\mbox{\upshape flat}}(p_1,p_2)$
	as desired.
	
	\eqref{E:SOBOELVEMBEDDINGRELYINGONQUASICONVEXITY}
	is a standard Sobolev embedding result (see, for example, \cite[Theorem~7]{pHpKhT2008}),
	which relies on the
	quantitative quasi-convexity of
	$\twoargMrough{[\timefunction_0,\timefunctionboot],[- \rightu,\leftu]}{\muxmulevelsetvalue}$
	proved in the previous paragraph.
	
	\eqref{E:CLOSEDVERSIONC21BOUNDFORCHOVROUGHTOGEO} 
	follows from \eqref{E:W3INFINITYBOUNDFORCHOVFROMGEOTOROUGHONOPEN}
	and
	\eqref{E:SOBOELVEMBEDDINGRELYINGONQUASICONVEXITY}.
	
	\eqref{E:CLOSEDVERSIONLUNITROUGHTTIMEFUNCTION} follows from \eqref{E:LDERIVATIVEOFROUGHTIMEFUNCTIONISAPPROXIMATELYUNITY}
	and the fact $\timefunctionarg{\muxmulevelsetvalue}$ is an element of
	$C_{\textnormal{geo}}^{2,1}\left(\twoargMrough{[\timefunction_0,\timefunctionboot],[- \rightu,\leftu]}{\muxmulevelsetvalue} \right)$,
	as is shown by
	\eqref{E:CLOSEDVERSIONC21BOUNDFORCHOVROUGHTOGEO}.
	
	The existence of the function
	$\Cartesiantisafunctiononlevelsetsofroughtimefunctionarg{\timefunction}{\muxmulevelsetvalue}$
	such that \eqref{E:LEVELSETSOFTIMEFUNCTIONAREAGRAPH}
	holds, as well as the estimate \eqref{E:C21BOUNDFORCONSTANTTIMEFUNCTIONGRAPH},
	follow from the fact that $(t,u,x^2,x^3) = \InverseCHOVgeotorough{\muxmulevelsetvalue}(\timefunction,u,x^2,x^3)$
	and the estimate \eqref{E:CLOSEDVERSIONC21BOUNDFORINVERSECHOVGEOTOROUGH},
	i.e.,
	$\Cartesiantisafunctiononlevelsetsofroughtimefunctionarg{\timefunction}{\muxmulevelsetvalue}(u,x^2,x^3)$
	is the first component of $\InverseCHOVgeotorough{\muxmulevelsetvalue}(\timefunction,u,x^2,x^3)$.
\end{proof}

\subsection{H\"{o}lder-space extensions to the compact set $\twoargMrough{[\timefunction_0,\timefunctionboot],[- \rightu,\leftu]}{\muxmulevelsetvalue}$}
\label{SS:HOLDERSPACEEXTENSIONS}
With the help of the bootstrap assumptions and
Lemma~\ref{L:DIFFEOMORPHICEXTENSIONOFROUGHCOORDINATES}, we now show that various solution variables extend to
the compact set $\twoargMrough{[\timefunction_0,\timefunctionboot],[- \rightu,\leftu]}{\muxmulevelsetvalue}$
as functions with substantial H\"{o}lder regularity relative to the geometric coordinates.

\begin{lemma}[H\"{o}lder-space extensions to the compact set 
$\twoargMrough{[\timefunction_0,\timefunctionboot],[- \rightu,\leftu]}{\muxmulevelsetvalue}$]
	\label{L:CONTINUOUSEXTNESION}
	The following quantities extend to the compact set 
	$\twoargMrough{[\timefunction_0,\timefunctionboot],[- \rightu,\leftu]}{\muxmulevelsetvalue}$
	as elements of the following H\"{o}lder spaces, and their  
	corresponding spacetime H\"{o}lder norms on $\twoargMrough{[\timefunction_0,\timefunctionboot],[- \rightu,\leftu]}{\muxmulevelsetvalue}$
	are bounded by $\leq C$:
	\begin{itemize}
		\item $\wavearray, \, \vortrenormalized^i, \, \GradEnt^i, \, \VortVort^i, \, \DivGradEnt 
			\in C_{\textnormal{geo}}^{3,1}(\twoargMrough{[\timefunction_0,\timefunctionboot],[- \rightu,\leftu]}{\muxmulevelsetvalue})$
		\item $\Upsilon \in C_{\textnormal{geo}}^{3,1}(\twoargMrough{[\timefunction_0,\timefunctionboot],[- \rightu,\leftu]}{\muxmulevelsetvalue})$
		\item $\Lunit^i, \, \upmu \in C_{\textnormal{geo}}^{2,1}(\twoargMrough{[\timefunction_0,\timefunctionboot],[- \rightu,\leftu]}{\muxmulevelsetvalue})$
		\item $\timefunctionarg{\muxmulevelsetvalue} \in C_{\textnormal{geo}}^{2,1}(\twoargMrough{[\timefunction_0,\timefunctionboot],[- \rightu,\leftu]}{\muxmulevelsetvalue})$
	\end{itemize}
\end{lemma}

\begin{proof}
	The results for $\timefunctionarg{\muxmulevelsetvalue}$ were already proved as \eqref{E:CLOSEDVERSIONC21BOUNDFORCHOVROUGHTOGEO}.
	For the remaining results, we give the proof only for $\wavearray$ since the
	other solution variables can be handled using nearly identical arguments.
	To proceed, we note that
	Lemma~\ref{L:COMMUTATORSTOCOORDINATES},
	Lemma~\ref{L:SCHEMATICSTRUCTUREOFVARIOUSTENSORSINTERMSOFCONTROLVARS},
	and the bootstrap assumptions imply that
	$\| \wavearray \|_{W_{\textnormal{geo}}^{4,\infty}(\twoargMrough{(\timefunction_0,\timefunctionboot),(- \rightu,\leftu)}{\muxmulevelsetvalue})}
	\leq C
	$.
	From this bound and \eqref{E:SOBOELVEMBEDDINGRELYINGONQUASICONVEXITY},
	we conclude that
	$
	\| \wavearray \|_{C_{\textnormal{geo}}^{3,1}(\twoargMrough{[\timefunction_0,\timefunctionboot],[- \rightu,\leftu)}{\muxmulevelsetvalue}]}
	\leq C
	$
	as desired.
\end{proof}

\subsection{Properties of $\CHOVroughtomumuxmu{\muxmulevelsetvalue}$ and consequences}
\label{SS:CHOVFROMROUGHCOORDINATESTOMUWEGIGHTEDXMUCOORDINATES}
In this section, 
we provide a detailed analysis of the change of variables map 
$\CHOVroughtomumuxmu{\muxmulevelsetvalue} (\timefunctionarg{\muxmulevelsetvalue},u,x^2,x^3) = (\upmu, \muX \upmu,x^2,x^3)$ and its Jacobian matrix $\CHOVJacobianroughtomumuxmu{\muxmulevelsetvalue}$.

\begin{lemma}[Properties of $\CHOVroughtomumuxmu{\muxmulevelsetvalue}$]
	\label{L:CHOVFROMROUGHCOORDINATESTOMUWEGIGHTEDXMUCOORDINATES}
	The map $\CHOVroughtomumuxmu{\muxmulevelsetvalue}(\timefunction,u,x^2,x^3) = (\upmu,\muX \upmu,x^2,x^3)$
	from Def.\,\ref{D:ALLTHECHOVMAPS}
	extends to a $C_{\textnormal{rough}}^{1,1}$ map on $[\timefunction_0,\timefunctionboot] \times [- \rightu,\leftu] \times \mathbb{T}^2$
	satisfying:
	\begin{align} \label{E:MUXMUCHOVMAPC112BOUND}
		\left\| 	
			\CHOVroughtomumuxmu{\muxmulevelsetvalue} 
		\right\|_{C_{\textnormal{rough}}^{1,1}\left([\timefunction_0,\timefunctionboot] \times [- \rightu,\leftu] \times \mathbb{T}^2 \right)}
		& \leq C.
	\end{align}
	Moreover, the Jacobian matrix 
	$\CHOVJacobianroughtomumuxmu{\muxmulevelsetvalue}(\timefunction,u,x^2,x^3)
	= 
	\frac{\partial (\upmu,\muX \upmu,x^2,x^3)}{\partial(\timefunction,u,x^2,x^3)}$
	is invertible at every point $q \in [\timefunction_0,\timefunctionboot] \times [-\interestingu,\interestingu] \times \mathbb{T}^2$
	and	satisfies:
	\begin{align} \label{E:JACOBIANDETERMINANTRATIOBOUND}
	\max_{q_1,q_2 \in [\timefunction_0,\timefunctionboot] \times [-\interestingu,\interestingu] \times \mathbb{T}^2}
	\left|
		\InverseCHOVJacobianroughtomumuxmu{\muxmulevelsetvalue}(q_1) \CHOVJacobianroughtomumuxmu{\muxmulevelsetvalue}(q_2)
		-
		\mbox{\upshape ID} 
	\right|_{\mbox{\upshape}Euc}
	& \leq \frac{1}{2},
\end{align}
where $|\cdot|_{\mbox{\upshape}Euc}$ is the standard Frobenius norm on matrices
(equal to the square root of the sum of the squares of the matrix entries)
and $\mbox{\upshape ID}$ denotes the $4 \times 4$ identity matrix.

	Furthermore, $\CHOVroughtomumuxmu{\muxmulevelsetvalue}$ 
	is a diffeomorphism from the compact, convex set
	$[\timefunction_0,\timefunctionboot] \times [-\interestingu,\interestingu] \times \mathbb{T}^2$
	onto its (compact) image 
	$\CHOVroughtomumuxmu{\muxmulevelsetvalue}\left([\timefunction_0,\timefunctionboot] \times [-\interestingu,\interestingu] \times \mathbb{T}^2 \right)$,
	where $\CHOVroughtomumuxmu{\muxmulevelsetvalue}\left([\timefunction_0,\timefunctionboot] \times [-\interestingu,\interestingu] \times \mathbb{T}^2 \right)$ 
	enjoys the following properties,
	and we recall that $\upmuboot = -\timefunctionboot$ and $\mupositive = -\timefunction_0$:
	\begin{enumerate}
	\item It contains
	$
	[\upmuboot,\mupositive] \times \lbrace - \muxmulevelsetvalue \rbrace \times \mathbb{T}^2
	$.
	\item It contains
	$
	(\upmuboot,\mupositive) \times \lbrace - \muxmulevelsetvalue \rbrace \times \mathbb{T}^2
	$
	in its interior.
	\item It
	is quasi-convex in the following sense:
	every pair of points
	$r_1, r_2 \in \CHOVroughtomumuxmu{\muxmulevelsetvalue}\left([\timefunction_0,\timefunctionboot] \times [-\interestingu,\interestingu] \times \mathbb{T}^2 \right)$
	are connected by a $C^1$ curve in 
	$\CHOVroughtomumuxmu{\muxmulevelsetvalue}\left([\timefunction_0,\timefunctionboot] \times [-\interestingu,\interestingu] \times \mathbb{T}^2 \right)$
	whose length is $\lesssim \mbox{\upshape dist}_{\mbox{\upshape flat}}(r_1,r_2)$,
	where $\mbox{\upshape dist}_{\mbox{\upshape flat}}(r_1,r_2)$ is the
	standard Euclidean distance between $r_1$ and $r_2$ in the flat space 
	$\mathbb{R}_{\timefunction} \times \mathbb{R}_u \times \mathbb{T}^2$.
	\end{enumerate}
	
	Moreover, with $\InverseCHOVroughtomumuxmu{\muxmulevelsetvalue}$ denoting the inverse map
	and $\InverseCHOVgeotorough{\muxmulevelsetvalue}$ denoting the inverse of 
	the change of variables map $\CHOVgeotorough{\muxmulevelsetvalue}$ defined in \eqref{E:CHOVGEOTOROUGH},
	the following holds for $\mulevelsetvalue \in [\upmuboot,\mupositive]$:
	\begin{align} \label{E:PHIINVERSEIMAGEOFTORUSISTORUSCONTAINEDINROUGHTIMEFUNCTIONLEVELSET}
		\InverseCHOVgeotorough{\muxmulevelsetvalue}
		\circ
		\InverseCHOVroughtomumuxmu{\muxmulevelsetvalue}\left(\lbrace \mulevelsetvalue \rbrace \times \lbrace - \muxmulevelsetvalue \rbrace \times \mathbb{T}^2 \right)
		= 
		\twoargmumuxtorus{\mulevelsetvalue}{-\muxmulevelsetvalue} 
		\subset
		\hypthreearg{-\mulevelsetvalue}{[-\frac{\interestingu}{2},\frac{\interestingu}{2}]}{\muxmulevelsetvalue},
	\end{align}	
	and:
	\begin{align} \label{E:PHIINVERSEIMAGEOFTORUSCROSSMUINTERVALISTORUSCROSSINTERVALCONTAINEDININTERESTINGREGION}
		\InverseCHOVgeotorough{\muxmulevelsetvalue}
		\circ
		\InverseCHOVroughtomumuxmu{\muxmulevelsetvalue}
		\left([\upmuboot,\mupositive] \times 
		\lbrace 
			- 
			\muxmulevelsetvalue 
		\rbrace 
		\times \mathbb{T}^2 \right)
		= 
		\datahypfortimefunctiontwoarg{-\muxmulevelsetvalue}{[\timefunction_0,\timefunctionboot]}
		\subset
		\twoargMrough{[\timefunction_0,\timefunctionboot],[-\frac{\interestingu}{2},\frac{\interestingu}{2}]}{\muxmulevelsetvalue}.
	\end{align}	
	In addition, $\InverseCHOVroughtomumuxmu{\muxmulevelsetvalue}$ satisfies the following estimate:
	\begin{align} \label{E:INVERSEOFMUXMUCHOVMAPC112BOUND}
		\left\| 
			\InverseCHOVroughtomumuxmu{\muxmulevelsetvalue} \right 
		\|_{C^{1,1}\left( 
		\CHOVroughtomumuxmu{\muxmulevelsetvalue}\left([\timefunction_0,\timefunctionboot] \times [-\interestingu,\interestingu] \times \mathbb{T}^2 \right)
		\right)}
		& \leq C.
	\end{align}
	
	In addition, for each fixed
	$(\timefunction,x^2,x^3) \in [\timefunction_0,\timefunctionboot] \times \mathbb{T}^2$, 
	the map 
	$u \rightarrow \muX \upmu(\timefunction,u,x^2,x^3)$
	is strictly increasing on $[-\interestingu,\interestingu]$. 
	
	Furthermore,
	the map $(\timefunction,u,x^2,x^3) \rightarrow \left(\timefunction,\muX \upmu,x^2,x^3 \right)$
	is a $C_{\textnormal{rough}}^{1,1}$ diffeomorphism from 
	$[\timefunction_0,\timefunctionboot] \times [-\interestingu,\interestingu] \times \mathbb{T}^2$ onto its image,
	which contains 
	$[\timefunction_0,\timefunctionboot]
	\times 
	\left[-\muxmulevelsetvalue - \frac{\secondtransversalderivativemulowerbound \interestingu}{16}, 
	- \muxmulevelsetvalue + \frac{\secondtransversalderivativemulowerbound \interestingu}{16}\right] \times \mathbb{T}^2$.
	In particular, by \eqref{E:KAPPA0ISSMALLERTHANM2RIGHTUOVER16}, 
	the image set contains 
	$[\timefunction_0,\timefunctionboot] \times [-\muxmulevelsetvalue_0,0] \times \mathbb{T}^2$.

	Finally, there
	exists a family of functions
	$\lbrace \Eikonalisafunctiononmumuxtoriarg{\mulevelsetvalue}{-\muxmulevelsetvalue} \rbrace_{\mulevelsetvalue \in [\upmu_{Boot},\mupositive]}$
	on $\mathbb{T}^2$
	that, for $\mulevelsetvalue \in (\upmu_{Boot},\mupositive]$, are equal to the functions from
	from Sect.\,\ref{SSS:BOOTSTRAPASSUMPTIONFORTORISTRUCTURE},
	such that
	for each $\timefunction \in [\timefunction_0,\timefunctionboot]$,
	we have:
	\begin{align} \label{E:C11BOUNDFORUROUGHCOORDINATEGRAPHTORI}
	\sup_{\timefunction \in [\timefunction_0,\timefunctionboot]} 
	\| \Eikonalisafunctiononmumuxtoriarg{-\timefunction}{-\muxmulevelsetvalue} \|_{C^{1,1}(\mathbb{T}^2)}
	& \leq C,
		\\
	\CHOVgeotorough{\muxmulevelsetvalue}\left(\twoargmumuxtorus{-\timefunction}{-\muxmulevelsetvalue}\right)
	& 
	=
	\left\lbrace
		\left(\timefunction,\Eikonalisafunctiononmumuxtoriarg{-\timefunction}{-\muxmulevelsetvalue}(x^2,x^3),x^2,x^3 \right) 
		\ | \ 
		(x^2,x^3) \in \mathbb{T}^2
	\right\rbrace,
	\label{E:LASTSLICETORIAREGRAPHSABOVEFLATTORIINGEOMETRICCOORDINATES}
	\end{align}
	where $\CHOVgeotorough{\muxmulevelsetvalue}$ is the change of variables map defined in \eqref{E:CHOVGEOTOROUGH}.
	In particular,
	\begin{align} \label{E:LASTSLICELEVELSETSTRUCTUREANDLOCATIONOFMIN} 
		\twoargmumuxtorus{-\timefunctionboot}{-\muxmulevelsetvalue}
		& 
		\subset
		\hypthreearg{\timefunctionboot}{[-\frac{3}{4}\interestingu,\frac{3}{4}\interestingu]}{\muxmulevelsetvalue}.
	\end{align}
	
\end{lemma}

\begin{proof}
\noindent \textbf{Proof of \eqref{E:MUXMUCHOVMAPC112BOUND}}:
Lemmas~\ref{L:DIFFEOMORPHICEXTENSIONOFROUGHCOORDINATES} and \ref{L:CONTINUOUSEXTNESION} yield \eqref{E:MUXMUCHOVMAPC112BOUND}.

\medskip
\noindent \textbf{Proof of \eqref{E:JACOBIANDETERMINANTRATIOBOUND}}:
We first use Lemma~\ref{L:DIFFEOMORPHICEXTENSIONOFROUGHCOORDINATES} 
(in particular, \eqref{E:CLOSEDKEYJACOBIANDETERMINANTESTIMATECHOVROUGHTOGEO}), 
Lemma~\ref{L:CONTINUOUSEXTNESION},
\eqref{E:BABOUNDSONGEOMETRICTDERIVATIVEMUINTERESTINGREGION},
and \eqref{E:BAMUTRANSVERSALCONVEXITY} 
to deduce that on the compact, convex set
$\mathscr{D} \eqdef [\timefunction_0,\timefunctionboot] \times [-\interestingu,\interestingu] \times \mathbb{T}^2$, 
the Jacobian matrix 
$\CHOVJacobianroughtomumuxmu{\muxmulevelsetvalue} 
= 
\frac{\partial (\upmu,\muX \upmu,x^2,x^3)}{\partial(\timefunction,u,x^2,x^3)}$ 
satisfies 
$\| \CHOVJacobianroughtomumuxmu{\muxmulevelsetvalue} \|_{C_{\textnormal{rough}}^{0,1}(\mathscr{D})} \leq C$
and:
\begin{align} 
\begin{split} \label{E:DETERMINANTOFJACOBIANFROMROUGHCOORDINATESTOMUMUXMUCOORDINATES}
	\mbox{\upshape det} \CHOVJacobianroughtomumuxmu{\muxmulevelsetvalue} 
	& = \left(\roughgeop{\timefunction} \upmu \right) \roughgeop{u} \muX \upmu
				-
				\left(\roughgeop{u} \upmu \right) \roughgeop{\timefunction} \muX \upmu
			\\
			& =
			\mbox{\upshape det} \left(d_{\textnormal{rough}} \InverseCHOVgeotorough{\muxmulevelsetvalue} \right)
			\mbox{\upshape det} \frac{\partial (\upmu,\muX \upmu,x^2,x^3)}{\partial(t,u,x^2,x^3)}
				\\
	& =
			\mbox{\upshape det} \left(d_{\textnormal{rough}} \InverseCHOVgeotorough{\muxmulevelsetvalue} \right)
			\left\lbrace
				\left(\geop{t} \upmu \right) \geop{u} \muX \upmu
				-
				\left(\geop{u} \upmu \right) \geop{t} \muX \upmu
			\right\rbrace
				\\
	& < 
		- 1/C, \qquad (\mbox{on } \mathscr{D}).
\end{split}
\end{align}
From \eqref{E:DETERMINANTOFJACOBIANFROMROUGHCOORDINATESTOMUMUXMUCOORDINATES},
it follows that $\CHOVJacobianroughtomumuxmu{\muxmulevelsetvalue}$ is invertible on $\mathscr{D}$.
Also using the definition of Lipschitz continuity, we deduce the pointwise bound
$
\left|
\CHOVJacobianroughtomumuxmu{\muxmulevelsetvalue}(\timefunction,u,x^2,x^3)
-
\CHOVJacobianroughtomumuxmu{\muxmulevelsetvalue}(\timefunction_0,u,x^2,x^3)
\right|_{Euc}
\leq C |\timefunction - \timefunction_0|
\leq
C \mupositive
$.
From these bounds and the data assumption \eqref{E:DATAJACOBIANDETERMINANTRATIOASSUMPTION},
we conclude \eqref{E:JACOBIANDETERMINANTRATIOBOUND}
whenever $\mupositive$ is sufficiently small.

Next, using \eqref{E:MUXMUCHOVMAPC112BOUND},
the inverse function theorem
and \eqref{E:JACOBIANDETERMINANTRATIOBOUND},
we deduce that $\CHOVroughtomumuxmu{\muxmulevelsetvalue}$ is a $C_{\textnormal{rough}}^{1,1}$ diffeomorphism from 
$\mathscr{D}$ onto its image 
such that $\CHOVroughtomumuxmu{\muxmulevelsetvalue}(\mathscr{D})$ is quasi-convex.
We clarify that \eqref{E:MUXMUCHOVMAPC112BOUND}, \eqref{E:JACOBIANDETERMINANTRATIOBOUND}, and the convexity of
$\mathscr{D}$ together guarantee the injectivity of $\CHOVroughtomumuxmu{\muxmulevelsetvalue}$ on $\mathscr{D}$
and the quasi-convexity of $\CHOVroughtomumuxmu{\muxmulevelsetvalue}(\mathscr{D})$.

Next, we use
\eqref{E:BOOSTRAPTORILOCATION},
the fact that
$\CHOVroughtomumuxmu{\muxmulevelsetvalue}$
is a diffeomorphism on $\mathscr{D}$,
and the facts that
$\upmu|_{\twoargmumuxtorus{-\timefunction}{-\muxmulevelsetvalue}} = - \timefunction$
and $\Wtransarg{\muxmulevelsetvalue} \upmu|_{\twoargmumuxtorus{-\timefunction}{-\muxmulevelsetvalue}} = 0$
(i.e., that $\muX \upmu|_{\twoargmumuxtorus{-\timefunction}{-\muxmulevelsetvalue}} = - \muxmulevelsetvalue$)
to deduce that
$
\CHOVroughtomumuxmu{\muxmulevelsetvalue}(\mathscr{D})
$
contains 
$[\upmuboot,\mupositive] \times \lbrace - \muxmulevelsetvalue \rbrace \times \mathbb{T}^2$
and that
$
	(\upmuboot,\mupositive) \times \lbrace - \muxmulevelsetvalue \rbrace \times \mathbb{T}^2
$
is contained in the interior of 
$
\CHOVroughtomumuxmu{\muxmulevelsetvalue}(\mathscr{D})
$.

We now use \eqref{E:ROUGHGEOPUINTERMOFGEOMETRICCOORDINATEVECTORFIELDS}, 
\eqref{E:BAMUTRANSVERSALCONVEXITY}, and Lemma~\ref{L:CONTINUOUSEXTNESION}
to deduce that:
\begin{align} \label{E:SIMPLEESTIMATEFORROUGHUDERIVATIVEOFMUXMU}
			\frac{\secondtransversalderivativemulowerbound}{4}  
			&
			\leq 
			\min_{\twoargMrough{[\timefunction_0,\timefunctionboot],[-\interestingu,\interestingu]}{\muxmulevelsetvalue}}
			\roughgeop{u} \muX \upmu
			\leq \max_{\twoargMrough{[\timefunction_0,\timefunctionboot],[-\interestingu,\interestingu]}{\muxmulevelsetvalue}}
			\roughgeop{u} \muX \upmu
			\leq \frac{4}{\secondtransversalderivativemulowerbound}.
\end{align}
From \eqref{E:SIMPLEESTIMATEFORROUGHUDERIVATIVEOFMUXMU},
it follows that for each fixed 
$(\timefunction,x^2,x^3) \in [\timefunction_0,\timefunctionboot] \times \mathbb{T}^2$,
	the map 
	$u \rightarrow \muX \upmu(\timefunction,u,x^2,x^3)$
	is strictly increasing on $[-\interestingu,\interestingu]$.
	Hence the map $(\timefunction,u,x^2,x^3) \rightarrow \left(\timefunction,\muX \upmu,x^2,x^3 \right)$
	is a $C^{1,1}$ diffeomorphism from 
	$[\timefunction_0,\timefunctionboot] \times [-\interestingu,\interestingu] \times \mathbb{T}^2$ onto its image.
	Since \eqref{E:BOOSTRAPTORILOCATION}
	implies that there is a $u_* \in [-\frac{3}{4}\interestingu,\frac{3}{4}\interestingu]$
	such that the image of $(\timefunction,u_*,x^2,x^3)$ under the map is $(\timefunction,-\muxmulevelsetvalue,x^2,x^3)$,
	we further deduce from \eqref{E:SIMPLEESTIMATEFORROUGHUDERIVATIVEOFMUXMU} that
	the image of $[-\interestingu,\interestingu] \times \mathbb{T}^2$ under the map contains
	$[\timefunction_0,\timefunctionboot] \times [-\muxmulevelsetvalue - \frac{\secondtransversalderivativemulowerbound \interestingu}{16}, 
	- \muxmulevelsetvalue + \frac{\secondtransversalderivativemulowerbound \interestingu}{16}] \times \mathbb{T}^2$,
	as is desired.
	
Finally,
the properties of $\Eikonalisafunctiononmumuxtoriarg{\mulevelsetvalue}{-\muxmulevelsetvalue}$ 
and the estimates
\eqref{E:C11BOUNDFORUROUGHCOORDINATEGRAPHTORI}--\eqref{E:LASTSLICELEVELSETSTRUCTUREANDLOCATIONOFMIN} 
follow from the form  
\eqref{E:CHOVFROMROUGHCOORDINATESTOMUWEGIGHTEDXMUCOORDINATES}
of $\CHOVroughtomumuxmu{\muxmulevelsetvalue}$,
the fact that 
$\InverseCHOVroughtomumuxmu{\muxmulevelsetvalue}$ 
is $C^{1,1}$ on
the compact, quasi-convex set
$
\CHOVroughtomumuxmu{\muxmulevelsetvalue}(\mathscr{D})
$,
and the bootstrap assumption \eqref{E:BOOSTRAPTORILOCATION},
which implies that for $\mulevelsetvalue \in [\upmuboot,\mupositive]$
and $(x^2,x^3) \in \mathbb{T}^2$,
we have
$
\InverseCHOVroughtomumuxmu{\muxmulevelsetvalue}(\mulevelsetvalue,-\muxmulevelsetvalue,x^2,x^3) \in 
\hypthreearg{-\mulevelsetvalue}{[-\frac{3}{4}\interestingu,\frac{3}{4}\interestingu]}{\muxmulevelsetvalue}
$.
\end{proof}

\begin{corollary}[Quantitative control of the embeddings on the closures of their domains]
\label{C:QUANTITATIVECONTROLOFEMBEDDINGSONCLOSURESOFTHEIRDOMAINS}
\hfill

\noindent \underline{\textbf{Control over $\embeddatahypersurfacearg{\muxmulevelsetvalue}$}}.
	The map
	$\embeddatahypersurfacearg{\muxmulevelsetvalue}(\mulevelsetvalue,x^2,x^3)
	=
	\left(\Cartesiantisafunctiononmumxtoriarg{\mulevelsetvalue}{-\muxmulevelsetvalue}(x^2,x^3),\Eikonalisafunctiononmumuxtoriarg{\mulevelsetvalue}{-\muxmulevelsetvalue}(x^2,x^3),x^2,x^3 \right)
	\in
	\twoargmumuxtorus{\mulevelsetvalue}{-\muxmulevelsetvalue}
	$
	from \eqref{E:EMBEDDATAHYPERSURFACE} extends to a $C^{1,1}$ embedding
	from $[\upmuboot,\mupositive] \times \mathbb{T}^2$ 
	onto its image, which is $\datahypfortimefunctiontwoarg{-\muxmulevelsetvalue}{[\timefunction_0,\timefunctionboot]}$.
	In addition, there is a $C > 1$ such that the extended embedding satisfies:
	\begin{align} \label{E:EMBEDDINGOFDATAHYPERSURFACEC11EXTENDEDC11BOUND}
		\| \embeddatahypersurfacearg{\muxmulevelsetvalue} \|_{C^{1,1}([\upmuboot,\mupositive] \times \mathbb{T}^2)} 
		& \leq C
	\end{align}
	and: 
	\begin{align} \label{E:CLOSEDIMPROVEMENTMUEQUALSMINUSKAPPAEMBEDDINGCARTESIANTIMEFUNCTIONNEGATIVEMUDERIVATIVE}
		- 
		C
		& 
		<
		\min_{(\mulevelsetvalue,x^2,x^3) \in [\upmuboot,\mupositive] \times \mathbb{T}^2} 
		\frac{\partial}{\partial \mulevelsetvalue} \Cartesiantisafunctiononmumxtoriarg{\mulevelsetvalue}{-\muxmulevelsetvalue}(x^2,x^3) 
		\leq
		\max_{(\mulevelsetvalue,x^2,x^3) \in [\upmuboot,\mupositive] \times \mathbb{T}^2} 
		\frac{\partial}{\partial \mulevelsetvalue} \Cartesiantisafunctiononmumxtoriarg{\mulevelsetvalue}{-\muxmulevelsetvalue}(x^2,x^3) 
		< 
		- \frac{1}{C}.
	\end{align}
	
	Furthermore, for $\mulevelsetvalue \in [\upmuboot,\mupositive]$,
	we have:
	\begin{align} \label{E:GRAPHDESCRIPTIONOFMUXMUTORUS}
			\twoargmumuxtorus{\mulevelsetvalue}{-\muxmulevelsetvalue}
			& = 
			\left\lbrace
				\left(\Cartesiantisafunctiononmumxtoriarg{\mulevelsetvalue}{-\muxmulevelsetvalue}(x^2,x^3),\Eikonalisafunctiononmumuxtoriarg{\mulevelsetvalue}{-\muxmulevelsetvalue}(x^2,x^3),x^2,x^3 \right)
				 \ | \
				(x^2,x^3) \in \mathbb{T}^2
			\right\rbrace
	\end{align}
	and:
	\begin{align}
			\datahypfortimefunctiontwoarg{-\muxmulevelsetvalue}{[\timefunction_0,\timefunctionboot]}
			& 
			=
			\bigcup_{\mulevelsetvalue \in [\upmuboot,\mupositive]}
			\twoargmumuxtorus{\mulevelsetvalue}{-\muxmulevelsetvalue}.
			\label{E:CLOSEDIMPROVEMENTLEVELSETSTRUCTUREOFMUXEQUALSMINUSKAPPA}
	\end{align}
	
	\medskip
	
	\noindent \underline{\textbf{Control over $\embeddingdatahypfortimefunctionarg{\muxmulevelsetvalue}$}}.
	The map $\embeddingdatahypfortimefunctionarg{\muxmulevelsetvalue}(t,x^2,x^3) 
	= 
	\left(t,\scalarembeddingdatahypfortimefunctionarg{\muxmulevelsetvalue}(t,x^2,x^3),x^2,x^3 \right)$ from
	Lemma~\ref{L:XMUISMINUSCAPPISAGRAPH}
	extends to a $C^{1,1}$ embedding
	from
	$\domainforembeddingdatahypfortimefunctiontwoarg{\muxmulevelsetvalue}{[\upmuboot,\mupositive]}
	\eqdef 
	\lbrace 
		(t,x^2,x^3) \in \mathbb{R} \times \mathbb{T}^2 
		\ | \
		\Cartesiantisafunctiononmumxtoriarg{\mupositive}{-\muxmulevelsetvalue}(x^2,x^3) \leq t \leq \Cartesiantisafunctiononmumxtoriarg{\upmuboot}{-\muxmulevelsetvalue}(x^2,x^3)
	\rbrace
	$
	onto its image, which is $\datahypfortimefunctiontwoarg{-\muxmulevelsetvalue}{[\timefunction_0,\timefunctionboot]}$.
	Moreover,  the extended embedding satisfies the following estimate:
	\begin{align} \label{E:EXTENDEDEMBEDDINGOFXMUISMINUSKAPPADATAHYPERSURFACEC11BOUND}
	\| \embeddingdatahypfortimefunctionarg{\muxmulevelsetvalue} \|_{C^{1,1}\left(\domainforembeddingdatahypfortimefunctiontwoarg{\muxmulevelsetvalue}{[\upmuboot,\mupositive]} \right)} 
	& \leq C.
	\end{align}
	
\end{corollary}

\begin{proof}
	$\embeddatahypersurfacearg{\muxmulevelsetvalue}$ is the composition of the maps
	$(\mulevelsetvalue,x^2,x^3) \rightarrow \InverseCHOVroughtomumuxmu{\muxmulevelsetvalue}(\mulevelsetvalue,\muxmulevelsetvalue,x^2,x^3)$
	and 
	$\InverseCHOVgeotorough{\muxmulevelsetvalue}$
	(see Def.\,\ref{D:ALLTHECHOVMAPS})
	and thus
	all of the desired conclusions except for those concerning $\embeddingdatahypfortimefunctionarg{\muxmulevelsetvalue}$
	follow from
	Lemmas~\ref{L:LINFTYESTIMATESFORROUGHTIMEFUNCTIONANDDERIVATIVES}
	and 
	\ref{L:CHOVFROMROUGHCOORDINATESTOMUWEGIGHTEDXMUCOORDINATES}.
	
	To obtain the results for
	$\embeddingdatahypfortimefunctionarg{\muxmulevelsetvalue}$,
	we consider the map
	$
	(\mulevelsetvalue,x^2,x^3) 
	\rightarrow 
	\left(\Cartesiantisafunctiononmumxtoriarg{\mulevelsetvalue}{-\muxmulevelsetvalue}(x^2,x^3),x^2,x^3 \right)
	$
	on the domain $[\upmuboot,\mupositive] \times \mathbb{T}^2$.
	Using
	\eqref{E:EMBEDDINGOFDATAHYPERSURFACEC11EXTENDEDC11BOUND},
	\eqref{E:CLOSEDIMPROVEMENTMUEQUALSMINUSKAPPAEMBEDDINGCARTESIANTIMEFUNCTIONNEGATIVEMUDERIVATIVE},
	and the inverse function theorem, we solve for the global inverse of the map.
	Composing this inverse with $\embeddatahypersurfacearg{\muxmulevelsetvalue}$,
	we obtain precisely the map $\embeddingdatahypfortimefunctionarg{\muxmulevelsetvalue}$.
	The estimate \eqref{E:EXTENDEDEMBEDDINGOFXMUISMINUSKAPPADATAHYPERSURFACEC11BOUND}
	is a straightforward consequence of these arguments
	(including 
	\eqref{E:EMBEDDINGOFDATAHYPERSURFACEC11EXTENDEDC11BOUND}
	and
	\eqref{E:CLOSEDIMPROVEMENTMUEQUALSMINUSKAPPAEMBEDDINGCARTESIANTIMEFUNCTIONNEGATIVEMUDERIVATIVE}).
\end{proof}

\subsection{Estimates for $\InverseCHOVJacobiangeotomumuxmu$ tied to the inverse function theorem}
\label{SS:GEOTOMUMUXCOORDINATESJACOBIANDETERMINANTRATIOINEVOLUTIONBOUND}
Recall that 
$\CHOVgeotomumuxmu(t,u,x^2,x^3) = (\upmu,\muX \upmu,x^2,x^3)$
is the map defined in \eqref{E:CHOVFROMGEOMETRICCOORDINATESTOMUWEGIGHTEDXMUCOORDINATES}.
In the next lemma, we prove estimates for its Jacobian matrix in the
region $\lbrace |u| \leq \interestingu \rbrace$.
Near the end of the paper, in Prop.\,\ref{P:PROPERTIESOFMSINGULARANDCREASE},
we will use the estimates in our analysis of the invertibility properties
of $\CHOVgeotomumuxmu$, which ultimately will help us derive the structure of the
singular boundary.

\begin{lemma}[Estimates for $\CHOVJacobiangeotomumuxmu$ tied to the inverse function theorem]
	\label{L:GEOTOMUMUXCOORDINATESJACOBIANDETERMINANTRATIOINEVOLUTIONBOUND}
	The Jacobian matrix $\CHOVJacobiangeotomumuxmu$
	defined in \eqref{E:JACOBIANMATRIXFORCHOVFROMGEOMETRICCOORDINATESTOMUWEGIGHTEDXMUCOORDINATES}
	is invertible on $\twoargMrough{[\timefunction_0,\timefunctionboot],[-\interestingu,\interestingu]}{\muxmulevelsetvalue}$
	and satisfies the following bounds, where $C > 1$:
	\begin{align} 
	- C
	&
	\leq
	\min_{\twoargMrough{[\timefunction_0,\timefunctionboot],[-\interestingu,\interestingu]}{\muxmulevelsetvalue}} 
		\mbox{\upshape det} \CHOVJacobiangeotomumuxmu
\leq
	\max_{\twoargMrough{[\timefunction_0,\timefunctionboot],[-\interestingu,\interestingu]}{\muxmulevelsetvalue}} 
		\mbox{\upshape det} \CHOVJacobiangeotomumuxmu
	\leq 
		- \frac{1}{C},
		\label{E:JACOBIANDETBOUNDCHOVGEOTOMUXMUCOORDS} 
\end{align}
\begin{align}
	\sup_{
		\substack{p_1 \in \twoargmumuxtorus{-\timefunction_0}{0}
				\\
			p_2 \in \twoargMrough{[\timefunction_0,\timefunctionboot],[-\interestingu,\interestingu]}{\muxmulevelsetvalue}
			}}
	\left|
		\CHOVJacobiangeotomumuxmu(p_1)
		\InverseCHOVJacobiangeotomumuxmu(p_2) 
		-
		\mbox{\upshape ID} 
	\right|_{\mbox{\upshape}Euc}
	& \leq \frac{1}{2},
		\label{E:GEOTOMUMUXCOORDINATESJACOBIANDETERMINANTRATIOINEVOLUTIONBOUND}
\end{align}
where $|\cdot|_{\mbox{\upshape}Euc}$ is the standard Frobenius norm on matrices
(equal to the square root of the sum of the squares of the matrix entries)
and $\mbox{\upshape ID}$ denotes the $4 \times 4$ identity matrix.
\end{lemma}

\begin{proof}
	Thanks to the initial data assumption
	\eqref{E:GEOTOMUMUXCOORDINATESDATAJACOBIANDETERMINANTRATIOASSUMPTION}
	and the transversal convexity bootstrap assumption \eqref{E:BAMUTRANSVERSALCONVEXITY},
	the lemma can be proved using the same arguments we used to prove
	\eqref{E:JACOBIANDETERMINANTRATIOBOUND} 
	(see especially the estimates in \eqref{E:DETERMINANTOFJACOBIANFROMROUGHCOORDINATESTOMUMUXMUCOORDINATES}),
	which in particular relied on the smallness of $\mupositive$.
\end{proof}



\section{Control of the flow map of  $\argLrough{\muxmulevelsetvalue}$}
\label{S:ESTIMATESTIEDTOTHEFLOWMAPOFROUGHNULLVECTORFIELD}
In this section, we derive various estimates tied to the
flow map of the vectorfield $\argLrough{\muxmulevelsetvalue}$.
We then use these results to derive preliminary
estimates for solutions $f$ to transport equations of the form
$\argLrough{\muxmulevelsetvalue} f = F$.

\subsection{Basic properties of the flow map of $\argLrough{\muxmulevelsetvalue}$}
\label{SS:PROPERTIESOFFLOWMAPOFWIDETILDEL}

\begin{lemma}[Basic properties of the flow map of $\argLrough{\muxmulevelsetvalue}$]
	\label{L:PROPERTIESOFFLOWMAPOFWIDETILDEL}
	Let $\argLrough{\muxmulevelsetvalue}$ be the null vectorfield defined in \eqref{E:LROUGH},
	and let $\FlowmapLrougharg{\muxmulevelsetvalue}$ be the $\timefunction_0$-normalized flow map of
	$\argLrough{\muxmulevelsetvalue}$ with respect to the rough adapted coordinates
	$(\timefunctionarg{\muxmulevelsetvalue},u,x^2,x^3)$, i.e., the solution to
	the following initial value problem:
	\begin{align} \label{E:FLOWMAPOFLROUGHINROUGHCOORDINATES}
		\roughgeop{\timefunction}
		\FlowmapLrougharg{\muxmulevelsetvalue}(\timefunction,u,x^2,x^3)
		& = 
		\argLrough{\muxmulevelsetvalue} \circ \FlowmapLrougharg{\muxmulevelsetvalue}(\timefunction,u,x^2,x^3),
		&
		\FlowmapLrougharg{\muxmulevelsetvalue}(\timefunction_0,u,x^2,x^3)
		& = (\timefunction_0,u,x^2,x^3).
	\end{align}
	Then for $A = 2,3$ there exist functions 
	$\FlowmapLroughtwoarg{\muxmulevelsetvalue}{A}
	: [\timefunction_0,\timefunctionboot] \times [- \rightu,\leftu] \times \mathbb{T}^2
	\rightarrow \mathbb{T}$
	such that:
	\begin{align} \label{E:FORMOFROUGHNULLGENERATORFLOWMAP}
	\FlowmapLrougharg{\muxmulevelsetvalue}(\timefunction,u,x^2,x^3)
	& 
	=
	\left(
		\timefunction,u,
		\FlowmapLroughtwoarg{\muxmulevelsetvalue}{2}(\timefunction,u,x^2,x^3),
		\FlowmapLroughtwoarg{\muxmulevelsetvalue}{3}(\timefunction,u,x^2,x^3)
	\right).
	\end{align}
	Moreover, $\FlowmapLrougharg{\muxmulevelsetvalue}$
	is a $C^{1,1}$ diffeomorphism
	from $[\timefunction_0,\timefunctionboot] \times [- \rightu,\leftu] \times \mathbb{T}^2$
	onto $[\timefunction_0,\timefunctionboot] \times [- \rightu,\leftu] \times \mathbb{T}^2$
	satisfying:
	\begin{align} \label{E:C11ROUGHBOUNDFORFLOWMAPOFWIDETILDEL}
		\left\|
			\FlowmapLrougharg{\muxmulevelsetvalue}
			-
			\mbox{\upshape I}
		\right\|_{C_{\textnormal{rough}}^{1,1}([\timefunction_0,\timefunctionboot] \times [- \rightu, \leftu] \times \mathbb{T}^2)}
		& \lesssim \varepsilon^{1/2},
		& 
		\left\|
			\InverseFlowmapLrougharg{\muxmulevelsetvalue}
			-
			\mbox{\upshape I}
		\right\|_{C_{\textnormal{rough}}^{1,1}([\timefunction_0,\timefunctionboot] \times [- \rightu, \leftu] \times \mathbb{T}^2)}
		& \lesssim \varepsilon^{1/2},
	\end{align}
	where $\mbox{\upshape I}(\timefunction,u,x^2,x^3) \eqdef (\timefunction,u,x^2,x^3)$
	is the identity map and $\InverseFlowmapLrougharg{\muxmulevelsetvalue}$ 
	is the inverse function of $\FlowmapLrougharg{\muxmulevelsetvalue}$.
	In particular, for each fixed $(\timefunction,u) \in [\timefunction,\timefunctionboot] \times [- \rightu,\leftu]$,
	the map 
	$(x^2,x^3) 
	\mapsto 
	\left(\FlowmapLroughtwoarg{\muxmulevelsetvalue}{2}(\timefunction,u,x^2,x^3),\FlowmapLroughtwoarg{\muxmulevelsetvalue}{3}(\timefunction,u,x^2,x^3) \right)$
	is a $C^{1,1}$ diffeomorphism from $\twoargroughtori{\timefunction_0,u}{\muxmulevelsetvalue}$ onto 
	$\twoargroughtori{\timefunction,u}{\muxmulevelsetvalue}$.
	
\end{lemma}

\begin{proof}
	\eqref{E:FORMOFROUGHNULLGENERATORFLOWMAP} is a trivial consequence of \eqref{E:FLOWMAPOFLROUGHINROUGHCOORDINATES}
	the identities $\argLrough{\muxmulevelsetvalue} \timefunctionarg{\muxmulevelsetvalue} =  1$ and $\argLrough{\muxmulevelsetvalue} u = 0$.
	
	To prove \eqref{E:C11ROUGHBOUNDFORFLOWMAPOFWIDETILDEL},
	we first note that
	the functions $\left(\FlowmapLroughtwoarg{\muxmulevelsetvalue}{2} - x^2,\FlowmapLroughtwoarg{\muxmulevelsetvalue}{3} - x^3 \right)$
	solve the transport system:
	\begin{align} \label{E:LROUGHFLOWMAPODETORUSCOMPONENTS}
	\roughgeop{\timefunction}
	[\FlowmapLroughtwoarg{\muxmulevelsetvalue}{A} - x^A](\timefunction,u,x^2,x^3)
	& 
	= 
	\twoargLrough{\muxmulevelsetvalue}{A}\left(\timefunction,u,\FlowmapLroughtwoarg{\muxmulevelsetvalue}{2},\FlowmapLroughtwoarg{\muxmulevelsetvalue}{3}\right),
	&
	&(A=2,3)
	\end{align}
	with vanishing data at rough time $\timefunction_0$.
	Next,
	we use definition \eqref{E:LROUGH},
	Lemma~\ref{L:ROUGHPARTIALDERIVATIVESINTERMSOFGEOMETRICPARTIALDERIVATIVESANDVICEVERSA},
	the bootstrap assumptions,
	and Lemma~\ref{L:DIFFEOMORPHICEXTENSIONOFROUGHCOORDINATES}
	to deduce that:
	\begin{align} \label{E:L2ROUGHL3ROUGHC11BOUNDS}
	\left\| 
		\left(\twoargLrough{\muxmulevelsetvalue}{2},\twoargLrough{\muxmulevelsetvalue}{3} \right)
	\right\|_{C_{\textnormal{rough}}^{1,1}([\timefunction_0,\timefunctionboot] \times [- \rightu, \leftu] \times \mathbb{T}^2)}
	& 
	\lesssim \varepsilon^{1/2}.
	\end{align}
	Hence, commuting \eqref{E:LROUGHFLOWMAPODETORUSCOMPONENTS} 
	up to one time with the rough adapted coordinate partial derivatives,
	using \eqref{E:L2ROUGHL3ROUGHC11BOUNDS},
	and integrating with respect to rough time,
	we find that for $\timefunction \in [\timefunction_0,\timefunctionboot]$,
	we have:
	\begin{align} 
	\begin{split}\label{E:GRONWALLREADYC11ROUGHBOUNDFORFLOWMAPOFWIDETILDEL}
		\max_{A=2,3}
		\left\| 
			\FlowmapLroughtwoarg{\muxmulevelsetvalue}{A} - x^A
		\right\|_{C_{\textnormal{rough}}^{1,1}([\timefunction_0,\timefunction] \times [- \rightu, \leftu] \times \mathbb{T}^2)}
		& \leq C \varepsilon^{1/2}
			\\
		& \ \
			+
			C  \varepsilon^{1/2}
			\int_{\timefunction_0}^{\timefunction}
				\max_{A=2,3}
				\left\| 
					\FlowmapLroughtwoarg{\muxmulevelsetvalue}{A} - x^A
				\right\|_{C_{\textnormal{rough}}^{1,1}([\timefunction_0,\timefunction'] \times [- \rightu, \leftu] \times \mathbb{T}^2)}
			\, \mathrm{d} \timefunction'.
	\end{split}
	\end{align}
	From \eqref{E:GRONWALLREADYC11ROUGHBOUNDFORFLOWMAPOFWIDETILDEL} and Gr\"{o}nwall's inequality,
	we find that
	$
	\max_{A=2,3}
		\left\| 
			\FlowmapLroughtwoarg{\muxmulevelsetvalue}{A} - x^A
		\right\|_{C_{\textnormal{rough}}^{1,1}([\timefunction_0,\timefunction] \times [- \rightu, \leftu] \times \mathbb{T}^2)}
	\leq C \varepsilon^{1/2}
	$.
	From this bound and \eqref{E:FORMOFROUGHNULLGENERATORFLOWMAP},
	we conclude the first bound stated in \eqref{E:C11ROUGHBOUNDFORFLOWMAPOFWIDETILDEL}.
	The second bound in \eqref{E:C11ROUGHBOUNDFORFLOWMAPOFWIDETILDEL}
	follows from differentiating the identity
	$
	\InverseFlowmapLrougharg{\muxmulevelsetvalue} \circ \FlowmapLrougharg{\muxmulevelsetvalue} = \mbox{\upshape I}
	$
	and using the first bound.
\end{proof}


\subsection{Estimate for $\mbox{\upshape det} \gtorusroughfirstfund$ and Minkowski's integral inequality} 
\label{SS:ROUGHNULLVECTORFIELDMINKOWSKIINTEGRALINEQUALITY}
In the next lemma, we control the factor $\mbox{\upshape det} \gtorusroughfirstfund(\timefunction_2,u,x^2,x^3)$
featured in the area form $\volroughtorus$ (see \eqref{E:AREAFORMROUGHTORUS}) 
on the rough tori $\twoargroughtori{\timefunction,u}{\muxmulevelsetvalue}$.
We then use this estimates 
to prove a Minkowski's integral inequality-type estimate.

\begin{lemma}[Estimate for $\mbox{\upshape det} \gtorusroughfirstfund$ and Minkowski's integral inequality] 
\label{L:ROUGHNULLVECTORFIELDMINKOWSKIINTEGRALINEQUALITY}
Recall that $\gtorusroughfirstfund$ is the first fundamental form of $\twoargroughtori{\timefunction,u}{\muxmulevelsetvalue}$
(see Def.\,\ref{D:ROUGHFIRSTFUNDS}) and that
$\FlowmapLrougharg{\muxmulevelsetvalue}$ is the
$\timefunction_0$-normalized flow map of $\argLrough{\muxmulevelsetvalue}$ from Lemma~\ref{L:PROPERTIESOFFLOWMAPOFWIDETILDEL}.
Then for every $\timefunction_1, \timefunction_2 \in [\timefunction_0,\timefunctionboot]$
and every $(u,x^2,x^3) \in [- \rightu,\leftu] \times \mathbb{T}^2$,
the following estimates hold,
where $\mbox{\upshape det} \gtorusroughfirstfund$ is evaluated relative to the rough adapted coordinates
via the formula \eqref{E:GTORUSROUGHCOMPONENTS}:
\begin{subequations}
\begin{align} \label{E:VOLFORMESTIMATEROUGHTORI}
\mbox{\upshape det} \gtorusroughfirstfund(\timefunction_2,u,x^2,x^3)
& = 
\left\lbrace
	1 + \mathcal{O}(\auxbootsmall)
\right\rbrace
\mbox{\upshape det} \gtorusroughfirstfund(\timefunction_1,u,x^2,x^3)
=
1 + \mathcal{O}_{\mydiam}(\mathring{\upalpha}),
	\\
\mbox{\upshape det} \gtorusroughfirstfund \circ \FlowmapLrougharg{\muxmulevelsetvalue}(\timefunction_2,u,x^2,x^3)
& 
= 
\left\lbrace
	1 + \mathcal{O}(\auxbootsmall)
\right\rbrace
\mbox{\upshape det} \gtorusroughfirstfund(\timefunction_1,u,x^2,x^3)
=
1 + \mathcal{O}_{\mydiam}(\mathring{\upalpha}).
\label{E:NORMALIZEDROUGHNULLFLOWMAPDOESNOTSUBSTANTIALLYDISTORTVOLFORMESTIMATEROUGHTORI}
\end{align}
\end{subequations}

Moreover, for any scalar function $F$ on $\twoargMrough{[\timefunction_0,\timefunctionboot),[- \rightu,\leftu]}{\muxmulevelsetvalue}$,
the function $\varphi(\timefunction,u,x^2,x^3)$
defined by:
\begin{align} \label{E:ANTIROUGHTIMEDERIVATIVE}
\varphi(\timefunction,u,x^2,x^3) 
& \eqdef 
\int_{\timefunction' = \timefunction_0}^{\timefunction} 
	F(\timefunction',u,x^2,x^3) 
\, \mathrm{d} \timefunction'
\end{align}
satisfies the following estimate:
\begin{align} \label{E:MINKOWSKIFORINTEGRALS}
\| \varphi \|_{L^2\left(\hypthreearg{\timefunction}{[- \rightu,u]}{\muxmulevelsetvalue}\right)} 
& \leq 
\left\lbrace 
	1 + \mathcal{O}(\auxbootsmall)
\right\rbrace 
\int_{\timefunction' = \timefunction_0}^{\timefunction} 
	\| F \|_{L^2\left(\hypthreearg{\timefunction'}{[- \rightu,u]}{\muxmulevelsetvalue}\right)} 
\, \mathrm{d} \timefunction'.
\end{align}
\end{lemma}

\begin{proof}
First, using the expression of the rough metric components 
$\gtorusroughfirstfund\left(\roughgeop{x^A},\roughgeop{x^B}\right)$ 
in \eqref{E:GTORUSROUGHCOMPONENTS}
and the identity
$\mbox{\upshape det} \gtorus = \frac{1}{\Speed^2 (X^1)^2}$ 
stated in \eqref{E:DETERMSMOOTHGTORUSRELTOGEOMETRICCOORDS}, 
we compute that:
\begin{align} \label{E:IDENTITYNEEDEDFORVOLFORMESTIMATEROUGHTORI}
\mbox{\upshape det} \gtorusroughfirstfund 
& = 
\left(\frac{\Lunit \timefunctionarg{\muxmulevelsetvalue}}{\geop{t}\timefunction}\right)^2 
\frac{1}{\Speed^2(X^1)^2}.
\end{align}
Using \eqref{E:IDENTITYNEEDEDFORVOLFORMESTIMATEROUGHTORI}
and
our assumptions on the data from Sect.\,\ref{SSS:QUANTITATIVEASSUMPTIONSONDATAAWAYFROMSYMMETRY},
and recalling that
$
\Speed(\LogDensity = 0,\Ent = 0) = 1
$
and
$X^1 = - 1 + \Xsmall^1$,
we deduce that
$\mbox{\upshape det} \gtorusroughfirstfund(\timefunction_0,u,x^2,x^3)
=
1 + \mathcal{O}_{\mydiam}(\mathring{\upalpha})
$.
Moreover, using \eqref{E:IDENTITYNEEDEDFORVOLFORMESTIMATEROUGHTORI}, 
Lemma~\ref{L:COMMUTATORSTOCOORDINATES},
\eqref{E:ROUGHTIMEPARTIALDERIVATIVEINTERMSOFGEOMETRICTIMEPARTIALDERIVATIVE},
Lemma~\ref{L:DIFFEOMORPHICEXTENSIONOFROUGHCOORDINATES},
and the bootstrap assumptions, 
we deduce that
$
\roughgeop{\timefunction} \mbox{\upshape det} \gtorusroughfirstfund
=
\mathcal{O}(\auxbootsmall)
$
which, in view of the mean value theorem,
yields that for any 
$(\timefunction,u,x^2,x^3) \in [\timefunction_0,\timefunctionboot] \times [- \rightu,\leftu] \times \mathbb{T}^2$,
we have
$
\mbox{\upshape det} \gtorusroughfirstfund(\timefunction,u,x^2,x^3)
=
\mbox{\upshape det} \gtorusroughfirstfund(\timefunction_0,u,x^2,x^3)
+
\mathcal{O}(\auxbootsmall)
=
\left\lbrace
	1 + \mathcal{O}(\auxbootsmall)
\right\rbrace
\mbox{\upshape det} \gtorusroughfirstfund(\timefunction_0,u,x^2,x^3)
$.
In total, these estimates imply \eqref{E:VOLFORMESTIMATEROUGHTORI}.

To prove \eqref{E:NORMALIZEDROUGHNULLFLOWMAPDOESNOTSUBSTANTIALLYDISTORTVOLFORMESTIMATEROUGHTORI},
we first note that definition~\eqref{E:LROUGH} and the estimates cited above yield the following estimate
relative to the rough adapted coordinates:
$
\roughgeop{\timefunction} 
\left\lbrace
	\gtorusroughfirstfund \circ \FlowmapLrougharg{\muxmulevelsetvalue}(\timefunction,u,x^2,x^3)
\right\rbrace
=
[\argLrough{\muxmulevelsetvalue} \mbox{\upshape det} \gtorusroughfirstfund] \circ \FlowmapLrougharg{\muxmulevelsetvalue}(\timefunction,u,x^2,x^3)
=
\mathcal{O}(\auxbootsmall)
$.
Using this estimate, 
arguing as in the previous paragraph, and using the initial condition 
$\FlowmapLrougharg{\muxmulevelsetvalue}(\timefunction_0,u,x^2,x^3) = (\timefunction_0,u,x^2,x^3)$,
we deduce the following estimate
for every 
$(\timefunction,u,x^2,x^3) \in [\timefunction_0,\timefunctionboot) \times [- \rightu,\leftu] \times \mathbb{T}^2$:
$
\mbox{\upshape det}  \gtorusroughfirstfund \circ \FlowmapLrougharg{\muxmulevelsetvalue}(\timefunction,u,x^2,x^3)
=
\left\lbrace
	1 + \mathcal{O}(\auxbootsmall)
\right\rbrace
\mbox{\upshape det}  \gtorusroughfirstfund(\timefunction_0,u,x^2,x^3)
$.
Combining this estimate with \eqref{E:VOLFORMESTIMATEROUGHTORI},
we conclude \eqref{E:NORMALIZEDROUGHNULLFLOWMAPDOESNOTSUBSTANTIALLYDISTORTVOLFORMESTIMATEROUGHTORI}.

The inequality \eqref{E:MINKOWSKIFORINTEGRALS}
follows from the following estimates relative to the rough adapted coordinates,
which rely on \eqref{E:VOLFORMESTIMATEROUGHTORI}
and Minkowski's inequality for integrals:
\begin{align} 
\begin{split}\label{E:PROOFOFMINKOWSKIFORINTEGRALS}
& \| \varphi \|_{L^2\left(\hypthreearg{\timefunction}{[- \rightu,u]}{\muxmulevelsetvalue}\right)}^2 
	\\
& 
= 
\int_{u' = - \rightu}^u 
	\int_{\mathbb{T}^2} 
		\left( 
			\int_{\timefunction' = \timefunction_0}^{\timefunction} 
				F(\timefunction',u',x^2,x^3) 
			\, \mathrm{d} \timefunction' 
		\right)^2 
	\sqrt{\mbox{\upshape det} \gtorusroughfirstfund(\timefunction,u',x^2,x^3)}
	\, \mathrm{d} x^2 \mathrm{d} x^3 
\mathrm{d} u' 
	\\
& \leq 
	\left\lbrace 
		1 
		+ 
		\mathcal{O}(\auxbootsmall) 
	\right\rbrace 
	\int_{u' = - \rightu}^u 
		\int_{(x^2,x^3) \in \T^2} 
			\left( 
				\int_{\timefunction' = \timefunction_0}^{\timefunction} 
					F(\timefunction',u',x^2,x^3)
					[\mbox{\upshape det} \gtorusroughfirstfund(\timefunction',u',x^2,x^3)]^{1/4} 
				\, \mathrm{d} \timefunction'
			\right)^2 
			\, \mathrm{d}x^2 \mathrm{d}x^3 
		\mathrm{d} u' 
			\\
& \leq  
	\left\lbrace 
		1 
		+ 
		\mathcal{O}(\auxbootsmall) 
	\right\rbrace 
	\left\lbrace
		\int_{\timefunction' = \timefunction_0}^{\timefunction} 
			\left[
				\int_{u' = - \rightu}^u
					\int_{(x^2,x^3) \in \T^2} 
						\left(
							F(\timefunction',u',x^2,x^3)
						\right)^2 
					\, 
					\sqrt{\mbox{\upshape det} \gtorusroughfirstfund(\timefunction',u',x^2,x^3)}
					\, \mathrm{d}x^2 \mathrm{d} x^3 \mathrm{d} u' 
			\right]^{1/2} 
		\mathrm{d} \timefunction'			
	\right\rbrace^2 
		\\
& = 
\left\lbrace 
		1 
		+ 
		\mathcal{O}(\auxbootsmall) 
\right\rbrace 
\left\lbrace
	\int_{\timefunction' = \timefunction_0}^{\timefunction} 
		\| F \|_{L^2\left(\hypthreearg{\timefunction'}{[- \rightu,u]}{\muxmulevelsetvalue}\right)} 
	\, \mathrm{d} \timefunction'
\right\rbrace^2.
\end{split}
\end{align}
\end{proof}

\subsection{Estimates for solutions to $\argLrough{\muxmulevelsetvalue} f = F$}
\label{SS:ESTIAMTESFORTRANSPORTEQUATIONROUGHNULLVECTORFIELD}
The following lemma provides the simple transport equation estimates 
that we use to control solutions to $\argLrough{\muxmulevelsetvalue} f = F$. 

\begin{lemma}[Pointwise, $L^2$, and $L^{\infty}$ estimates tied to the integral curves of $\argLrough{\muxmulevelsetvalue}$]  
\label{L:TRANSPORTESTIMATESFORROUGHLFEQUALSSOURCE} 
Let $f$ be a function of the rough adapted coordinates on 
$[\timefunction_0,\timefunctionboot) \times [- \rightu,\leftu] \times \mathbb{T}^2$,
let $\FlowmapLrougharg{\muxmulevelsetvalue}$ be the $\timefunction_0$-normalized flow map of $\argLrough{\muxmulevelsetvalue}$ from 
Lemma~\ref{L:PROPERTIESOFFLOWMAPOFWIDETILDEL}.
Then relative to the rough adapted coordinates, the following identity holds for any 
$\timefunction_0 \le \timefunction_1 \le \timefunction_2 \leq \timefunctionboot \leq 0$
and any $(u,x^2,x^3) \in [- \rightu,\leftu] \times \mathbb{T}^2$:
\begin{align} \label{E:TRANSPORTIDENTITYALONGLROUGHINTEGRALCURVES}
f \circ \FlowmapLrougharg{\muxmulevelsetvalue}(\timefunction_2,u,x^2,x^3)
& = 
	f \circ \FlowmapLrougharg{\muxmulevelsetvalue}(\timefunction_1,u,x^2,x^3)
	+ 
	\int_{\timefunction' = \timefunction_1}^{\timefunction_2} (\argLrough{\muxmulevelsetvalue} f) 
	\circ 
	\FlowmapLrougharg{\muxmulevelsetvalue}(\timefunction',u,x^2,x^3) 
	\, \mathrm{d} \timefunction'.
\end{align}

Moreover, relative to the rough adapted coordinates, 
we have the following estimate \textbf{for the critically important factor $G_{\Lunit \Lunit}^0$}:
\begin{align} \label{E:CRUCIALGLL0TRANSPORTESTIMATE}
\left|
	G_{\Lunit \Lunit}^0 \circ \FlowmapLrougharg{\muxmulevelsetvalue}(\timefunction_2,u,x^2,x^3)
	- 
	G_{\Lunit \Lunit}^0 \circ \FlowmapLrougharg{\muxmulevelsetvalue}(\timefunction_1,u,x^2,x^3)
\right| 
& \lesssim \auxbootsmall |\timefunction_2 - \timefunction_1|.
\end{align}

In addition, we have the following $L^{\infty}$ and $L^2$ estimates:
\begin{subequations}
\begin{align} 
	\| f \|_{L^{\infty}\left(\twoargroughtori{\timefunction_2,u}{\muxmulevelsetvalue}\right)} 
	& \leq 
	\| f\|_{L^{\infty}\left(\twoargroughtori{\timefunction_1,u}{\muxmulevelsetvalue}\right)} 
	+
	\mupositive \sup_{\timefunction' \in [\timefunction_1,\timefunction_2]} 
	\| \argLrough{\muxmulevelsetvalue} f \|_{L^{\infty}\left(\twoargroughtori{\timefunction',u}{\muxmulevelsetvalue}\right)}, 
	\label{E:LINFINITYTRANSPORTROUGHLFESTIMATE} 
		\\
	\esssup_{(x^2,x^3) \in \mathbb{T}^2}
		\left|
			f \circ \FlowmapLrougharg{\muxmulevelsetvalue}(\timefunction_2,u,x^2,x^3)
			-
			f \circ \FlowmapLrougharg{\muxmulevelsetvalue}(\timefunction_1,u,x^2,x^3) 
		\right|
	& \leq 
	\mupositive \sup_{\timefunction' \in [\timefunction_1,\timefunction_2]} 
	\| \argLrough{\muxmulevelsetvalue} f \|_{L^{\infty}\left(\twoargroughtori{\timefunction',u}{\muxmulevelsetvalue}\right)}, 
	\label{E:MOREPRECISELINFINITYTRANSPORTROUGHLFESTIMATE} 
\end{align}
\end{subequations}

\begin{align} 	\label{E:L2ONROUGHCONSTANTTIMEHYPERSURFACESTRANSPORTROUGHLFESTIMATE}
	\| f \|_{L^2 \left(\hypthreearg{\timefunction_2}{[- \rightu,u]}{\muxmulevelsetvalue}\right)} 
	& \leq  
	\left\lbrace 
		1 
		+ 
		\mathcal{O}(\auxbootsmall)
	\right\rbrace 
	\| f \|_{L^2\left(\hypthreearg{\timefunction_1}{[- \rightu,u]}{\muxmulevelsetvalue}\right)} 
	+ 
	\left\lbrace 
		1 
		+ 
		\mathcal{O}(\auxbootsmall)
	\right\rbrace
	\int_{\timefunction' = \timefunction_1}^{\timefunction_2} 
		\| \argLrough{\muxmulevelsetvalue} f \|_{L^2\left(\hypthreearg{\timefunction'}{[- \rightu,u]}{\muxmulevelsetvalue}\right)} 
	\, \mathrm{d} \timefunction'. 
\end{align}

Furthermore, if $f$ is any function of the rough adapted coordinates on 
$[\timefunction_0,\timefunctionboot) \times [- \rightu,\leftu] \times \mathbb{T}^2$,
then for every $(\timefunction,u) \in [\timefunction_0,\timefunctionboot) \times [- \rightu,\leftu]$,
the following estimate holds:
\begin{align} \label{E:DISTORTINGWITHROUGHFLOWMAPDOESNOTCHANGEL2NORMSMUCH}
	\| f \circ \FlowmapLrougharg{\muxmulevelsetvalue} \|_{L^2\left(\hypthreearg{\timefunction}{[- \rightu,u]}{\muxmulevelsetvalue}\right)}
	&
	=
	\left\lbrace 
		1 
		+ 
		\mathcal{O}(\auxbootsmall)
	\right\rbrace 
	\| f \|_{L^2\left(\hypthreearg{\timefunction}{[- \rightu,u]}{\muxmulevelsetvalue}\right)}.
\end{align}

Finally, let $F$ be a scalar function of the rough adapted coordinates
on $[\timefunction_0,\timefunctionboot) \times [- \rightu,\leftu] \times \mathbb{T}^2$,
and let $\varphi$ be the function of the rough adapted coordinates defined by:
\begin{align} \label{E:AGAINANTIROUGHTIMEDERIVATIVE}
	\varphi(\timefunction,u,x^2,x^3) 
	& \eqdef 
	\int_{\timefunction' = \timefunction_0}^{\timefunction} 
		F \circ \FlowmapLrougharg{\muxmulevelsetvalue}(\timefunction',u,x^2,x^3) 
	\, \mathrm{d} \timefunction'.
\end{align}
Then for every $(\timefunction,u) \in [\timefunction_0,\timefunctionboot) \times [- \rightu,\leftu]$, 
the following estimate holds:
\begin{align} \label{E:L2ONROUGHCONSTANTTIMEHYPERSURFACEOFROUGHTIMEINTEGRALWITHFLOWMAPFACTORSBOUND}
	\| \varphi \|_{L^2 \left(\hypthreearg{\timefunction}{[- \rightu,u]}{\muxmulevelsetvalue}\right)} 
	& 
	\leq  
	\left\lbrace 
		1 
		+ 
		\mathcal{O}(\auxbootsmall)
	\right\rbrace
	\int_{\timefunction' = \timefunction_0}^{\timefunction} 
		\| F \|_{L^2\left(\hypthreearg{\timefunction'}{[- \rightu,u]}{\muxmulevelsetvalue}\right)} 
	\, \mathrm{d} \timefunction'. 
\end{align}

\end{lemma}

\begin{proof}
\eqref{E:TRANSPORTIDENTITYALONGLROUGHINTEGRALCURVES}
follows from
\eqref{E:FLOWMAPOFLROUGHINROUGHCOORDINATES}
and the fundamental theorem of calculus.

To prove \eqref{E:CRUCIALGLL0TRANSPORTESTIMATE},
we first note that 
\eqref{E:LROUGH},
\eqref{E:CLOSEDVERSIONLUNITROUGHTTIMEFUNCTION},
and the bootstrap assumptions imply that
$\left|\argLrough{\muxmulevelsetvalue} G_{\Lunit \Lunit}^0 \right|
\lesssim \auxbootsmall
$.
Hence, using this bound and
applying 
\eqref{E:TRANSPORTIDENTITYALONGLROUGHINTEGRALCURVES}
with $G_{\Lunit \Lunit}^0$ in the role of $f$,
we conclude \eqref{E:CRUCIALGLL0TRANSPORTESTIMATE}.

We now prove \eqref{E:LINFINITYTRANSPORTROUGHLFESTIMATE}. 
First, from \eqref{E:TRANSPORTIDENTITYALONGLROUGHINTEGRALCURVES},
we deduce:
\begin{align} \label{E:TRANSPORTINEEQUALITYALONGLROUGHINTEGRALCURVES}
	\left|
		f \circ \FlowmapLrougharg{\muxmulevelsetvalue}(\timefunction_2,u,x^2,x^3)
	\right|
	& 
	\leq
	\left|
		f \circ \FlowmapLrougharg{\muxmulevelsetvalue}(\timefunction_1,u,x^2,x^3)
	\right|
	+ 
	\int_{\timefunction' = \timefunction_1}^{\timefunction_2} 
	\left|
		(\argLrough{\muxmulevelsetvalue} f) 
		\circ 
		\FlowmapLrougharg{\muxmulevelsetvalue}(\timefunction',u,x^2,x^3) 
	\right|
	\, \mathrm{d} \timefunction'.
\end{align}
We now take the essential supremum norm of both sides of \eqref{E:TRANSPORTINEEQUALITYALONGLROUGHINTEGRALCURVES}
over $(x^2,x^3) \in \mathbb{T}^2$ and use Lemma~\ref{L:PROPERTIESOFFLOWMAPOFWIDETILDEL}
to deduce
$
\| f \|_{L^{\infty}\left(\twoargroughtori{\timefunction_2,u}{\muxmulevelsetvalue}\right)} 
\leq 
\| f \|_{L^{\infty}\left(\twoargroughtori{\timefunction_1,u}{\muxmulevelsetvalue}\right)} 
+
	\int_{\timefunction' = \timefunction_1}^{\timefunction_2} 
		\| \argLrough{\muxmulevelsetvalue} f \|_{L^{\infty}\left(\twoargroughtori{\timefunction',u}{\muxmulevelsetvalue}\right)} 
	\, \mathrm{d} \timefunction'
$.
From this bound and the fact that
$
|\timefunction_2 - \timefunction_1| \leq |\timefunction_0| = \mupositive
$,
we conclude \eqref{E:LINFINITYTRANSPORTROUGHLFESTIMATE}.
The estimate \eqref{E:MOREPRECISELINFINITYTRANSPORTROUGHLFESTIMATE} 
can be proved via a similar argument, and we omit the details.

To prove \eqref{E:L2ONROUGHCONSTANTTIMEHYPERSURFACESTRANSPORTROUGHLFESTIMATE},
we first take the 
$\| \cdot \|_{L^2\left(\hypthreearg{\timefunction_2}{[- \rightu,u]}{\muxmulevelsetvalue}\right)}$
norm of both sides of \eqref{E:TRANSPORTIDENTITYALONGLROUGHINTEGRALCURVES}.
We then use Lemma~\ref{L:PROPERTIESOFFLOWMAPOFWIDETILDEL}
(in particular \eqref{E:C11ROUGHBOUNDFORFLOWMAPOFWIDETILDEL}),
\eqref{E:NORMALIZEDROUGHNULLFLOWMAPDOESNOTSUBSTANTIALLYDISTORTVOLFORMESTIMATEROUGHTORI},
\eqref{E:MINKOWSKIFORINTEGRALS}, 
and the standard formula for changing variables in an integral
(these estimates allow us in particular to replace all terms
$
f \circ \FlowmapLrougharg{\muxmulevelsetvalue}(\timefunction,u,x^2,x^3)
$
under the $L^2$ norms with
$
f(\timefunction,u,x^2,x^3)
$,
up to $1 + \mathcal{O}(\auxbootsmall)$ multiplicative factors)
to deduce:
\begin{align} 	\label{E:FIRSTSTEPL2TRANSPORTROUGHLFESTIMATE}
	\| f \|_{L^2\left(\hypthreearg{\timefunction_2}{[- \rightu,u]}{\muxmulevelsetvalue}\right)} 
	& \leq  
	\left\lbrace 
		1 
		+ 
		\mathcal{O}(\auxbootsmall)
	\right\rbrace
	\| f(\timefunction_1,\cdot) \|_{L^2\left(\hypthreearg{\timefunction_2}{[- \rightu,u]}{\muxmulevelsetvalue}\right)}
	+ 
	\left\lbrace 
		1 
		+ 
		\mathcal{O}(\auxbootsmall)
	\right\rbrace
	\int_{\timefunction' = \timefunction_1}^{\timefunction_2} 
		\left\| 
			\argLrough{\muxmulevelsetvalue} f
		\right\|_{L^2\left(\hypthreearg{\timefunction'}{[- \rightu,u]}{\muxmulevelsetvalue}\right)} 
	\, \mathrm{d} \timefunction'.
\end{align}
We then use \eqref{E:VOLFORMESTIMATEROUGHTORI}
to deduce that the first term on RHS~\eqref{E:FIRSTSTEPL2TRANSPORTROUGHLFESTIMATE}
satisfies:
\begin{align} \label{E:BOUNDFORFIRSTTERMONRHSOFFIRSTSTEPL2TRANSPORTROUGHLFESTIMATE}
\| f(\timefunction_1,\cdot) \|_{L^2\left(\hypthreearg{\timefunction_2}{[- \rightu,u]}{\muxmulevelsetvalue}\right)}
&=
\left\lbrace 
		1 
		+ 
		\mathcal{O}(\auxbootsmall)
	\right\rbrace
\| f(\timefunction_1,\cdot) \|_{L^2\left(\hypthreearg{\timefunction_1}{[- \rightu,u]}{\muxmulevelsetvalue}\right)},
\end{align}
which in total yields \eqref{E:L2ONROUGHCONSTANTTIMEHYPERSURFACESTRANSPORTROUGHLFESTIMATE}.

The estimate \eqref{E:DISTORTINGWITHROUGHFLOWMAPDOESNOTCHANGEL2NORMSMUCH}
follows from \eqref{E:C11ROUGHBOUNDFORFLOWMAPOFWIDETILDEL},
\eqref{E:NORMALIZEDROUGHNULLFLOWMAPDOESNOTSUBSTANTIALLYDISTORTVOLFORMESTIMATEROUGHTORI},
and the standard formula for changing variables in an integral.

Finally, \eqref{E:L2ONROUGHCONSTANTTIMEHYPERSURFACEOFROUGHTIMEINTEGRALWITHFLOWMAPFACTORSBOUND}
follows from 
\eqref{E:MINKOWSKIFORINTEGRALS}
and
\eqref{E:DISTORTINGWITHROUGHFLOWMAPDOESNOTCHANGEL2NORMSMUCH}.
\end{proof}

\section{$L^{\infty}$ estimates and improvement of the auxiliary bootstrap assumptions}
\label{S:LINFINITYFLUIDANDEIKONALANDIMPROVEMENTOFAUX}
In this section, we derive $L^{\infty}$ estimates for the fluid variables and eikonal function quantities 
that in particular yield improvements of the auxiliary bootstrap assumptions stated in Sect.\,\ref{SSS:AUXBOOTSTRAP}.

\begin{proposition}[$L^{\infty}$ estimates and improvement of the auxiliary bootstrap assumptions] 
\label{P:IMPROVEMENTOFAUXILIARYBOOTSTRAP}
Under the parameter-size and initial data assumptions 
of Sects.\,\ref{SS:PARAMETERSIZEASSUMPTIONS}
and \ref{SS:ASSUMPTIONSONDATA}
and the bootstrap assumptions of 
Sects.\,\ref{SS:BOOTSTRAPSCAFFOLDING} and \ref{SS:MAINQUANTITATIVEBOOTSTRAPASSUMPTIONS},
the following estimates hold for
$(\timefunction,u) \in [\timefunction_0,\timefunctionboot) \times [- \rightu,\leftu]$
(where we recall that $\wavearray$ and $\wavearraypartial$ are defined in Def.\,\ref{D:ARRAYSOFWAVEVARIABLES},
and that in Sect.\,\ref{SS:STRINGSOFCOMMUTATIONVECTORFIELDS},
we introduced notation for strings of commutation vectorfields).
\medskip

\noindent \underline{\textbf{$L^{\infty}$ estimates for small quantities}}.

\begin{align}  
	\left\| \RRiemann \right\|_{L^{\infty}(\twoargroughtori{\timefunction,u}{\muxmulevelsetvalue})} 
	& \leq \mathring{\upalpha} 
		+ 
		C \fundbootsmall, 
		\label{E:LINFINITYIMPROVEMENTAUXR+} 
		\\
	\left\| \wavearraypartial \right\|_{L^{\infty}(\twoargroughtori{\timefunction,u}{\muxmulevelsetvalue})} 
	& \leq C \fundbootsmall,
		 \label{E:IMPROVEAUXWAVEARRAYPARTIALLINFINITY} 
					\\
	\left\| \argLrough{\muxmulevelsetvalue} \comder^{\leq \Ntop-11;1} \wavearray \right\|_{L^{\infty}(\twoargroughtori{\timefunction,u}{\muxmulevelsetvalue})},
		\,
	\left\| \comdersmall^{[1,\Ntop-11];1} \wavearray \right\|_{L^{\infty}(\twoargroughtori{\timefunction,u}{\muxmulevelsetvalue})},
			 \label{E:LINFINITYIMPROVEMENTAUXWAVEARRAY} 
		&		\\
	\left\| \argLrough{\muxmulevelsetvalue} \comder^{[1,6];2} \wavearray \right\|_{L^{\infty}(\twoargroughtori{\timefunction,u}{\muxmulevelsetvalue})},
		\,
	\left\| \comdersmall^{[1,6];2} \wavearray \right\|_{L^{\infty}(\twoargroughtori{\timefunction,u}{\muxmulevelsetvalue})},
		& \notag	\\
	\left\| \argLrough{\muxmulevelsetvalue} \comder^{[1,5];3} \wavearray \right\|_{L^{\infty}(\twoargroughtori{\timefunction,u}{\muxmulevelsetvalue})},
		\,
	\left\| \comdersmall^{[1,5];3} \wavearray \right\|_{L^{\infty}(\twoargroughtori{\timefunction,u}{\muxmulevelsetvalue})},
	& \notag	\\
	\left\| \argLrough{\muxmulevelsetvalue} \muX \muX \muX \muX \wavearray \right\|_{L^{\infty}(\twoargroughtori{\timefunction,u}{\muxmulevelsetvalue})}
	& \leq 
	C \fundbootsmall,
	\notag 
		\\
	\left\| 
		\comder^{\leq \Ntop-11;1} (\vortrenormalized,\GradEnt) 
	\right\|_{L^{\infty}(\twoargroughtori{\timefunction,u}{\muxmulevelsetvalue})},
	&		
		\label{E:LINFINITYIMPROVEMENTAUXVORTGRADENT} 
		\\
	\left\| 
		\comder^{\leq 6;2} (\vortrenormalized,\GradEnt) 
	\right\|_{L^{\infty}(\twoargroughtori{\timefunction,u}{\muxmulevelsetvalue})},	
	&		\notag \\
	\left\| 
		\comder^{\leq 5;3} (\vortrenormalized,\GradEnt) 
	\right\|_{L^{\infty}(\twoargroughtori{\timefunction,u}{\muxmulevelsetvalue})},	
	&		
		\notag \\
	\left\| 
		\muX \muX \muX \muX (\vortrenormalized,\GradEnt) 
	\right\|_{L^{\infty}(\twoargroughtori{\timefunction,u}{\muxmulevelsetvalue})}
	& \leq C \fundbootsmall,
			\notag \\
	\left\| 
		\comder^{\leq \Ntop-12;1} (\VortVort,\DivGradEnt) 
	\right\|_{L^{\infty}(\twoargroughtori{\timefunction,u}{\muxmulevelsetvalue})},
	&		
		\label{E:LINFINITYIMPROVEMENTAUXMODIFIEDFLUID}  
		\\
	\left\| 
		\comder^{\leq 6;2} (\VortVort,\DivGradEnt)
	\right\|_{L^{\infty}(\twoargroughtori{\timefunction,u}{\muxmulevelsetvalue})},	
	& \notag 
			\\
	\left\| 
		\comder^{\leq 5;3} (\VortVort,\DivGradEnt)
	\right\|_{L^{\infty}(\twoargroughtori{\timefunction,u}{\muxmulevelsetvalue})},	
	& \notag 
		\\
	\left\| 
		\muX \muX \muX \muX (\VortVort,\DivGradEnt)
	\right\|_{L^{\infty}(\twoargroughtori{\timefunction,u}{\muxmulevelsetvalue})}
	& \leq C \fundbootsmall,
			\notag 
			\\
	\left\| \argLrough{\muxmulevelsetvalue} \tander^{[1,\Ntop-12]} \upmu \right\|_{L^{\infty}(\twoargroughtori{\timefunction,u}{\muxmulevelsetvalue})}, 
		\,
	\left\| \tander_*^{[1,\Ntop-12]} \upmu \right\|_{L^{\infty}(\twoargroughtori{\timefunction,u}{\muxmulevelsetvalue})},
	& 
	 \label{E:LINFINITYIMPROVEMENTAUXMUSMALL}
		\\
	\left\| \argLrough{\muxmulevelsetvalue} \comdersmall^{[1,5];1}  \upmu \right\|_{L^{\infty}(\twoargroughtori{\timefunction,u}{\muxmulevelsetvalue})}, 
		\,
	\left\| \comderdoublesmall^{[1,5];1} \upmu \right\|_{L^{\infty}(\twoargroughtori{\timefunction,u}{\muxmulevelsetvalue})},
	&	\notag 
		\\
	\left\| \argLrough{\muxmulevelsetvalue} \comdersmall^{[1,4];2} \upmu \right\|_{L^{\infty}(\twoargroughtori{\timefunction,u}{\muxmulevelsetvalue})}, 
		\,
	\left\| \comderdoublesmall^{[1,4];2} \upmu \right\|_{L^{\infty}(\twoargroughtori{\timefunction,u}{\muxmulevelsetvalue})}
	& \leq C \fundbootsmall,
		\notag \\
	\left\| \Lsmall^1 \right\|_{L^{\infty}(\twoargroughtori{\timefunction,u}{\muxmulevelsetvalue})} 
	& 
	\leq \mathring{\upalpha}
			+
			C \fundbootsmall,
		\label{E:LINFINITYIMPROVEMENTAUXL1SMALL} 
		\\
	\left\| \Lsmall^A \right\|_{L^{\infty}(\twoargroughtori{\timefunction,u}{\muxmulevelsetvalue})}
	\label{E:LINFINITYIMPROVEMENTAUXLASMALL} 
	& 
	\leq C \fundbootsmall
				\\
	\left\| \argLrough{\muxmulevelsetvalue} \tander^{\leq \Ntop-11} \Lsmall^i \right\|_{L^{\infty}(\twoargroughtori{\timefunction,u}{\muxmulevelsetvalue})}, 
			\,
	\left\| \tander^{[1,\Ntop-11]} \Lsmall^i \right\|_{L^{\infty}(\twoargroughtori{\timefunction,u}{\muxmulevelsetvalue})}, 
			\label{E:LINFINITYIMPROVEMENTAUXMIXECTANGENTIALTRANSVERSALDERIVATIVESLISMALL}
			\\
	\left\| \argLrough{\muxmulevelsetvalue} \comder^{[1,\Ntop-12];1} \Lsmall^i \right\|_{L^{\infty}(\twoargroughtori{\timefunction,u}{\muxmulevelsetvalue})}, 
		\,
	\left\| \comdersmall^{[1,\Ntop-12];1} \Lsmall^i \right\|_{L^{\infty}(\twoargroughtori{\timefunction,u}{\muxmulevelsetvalue})}, 
		\notag 
			\\
	\left\| \argLrough{\muxmulevelsetvalue} \comder^{[1,5];2}\Lsmall^i \right\|_{L^{\infty}(\twoargroughtori{\timefunction,u}{\muxmulevelsetvalue})}, 
		\,
	\left\| \comdersmall^{[1,5];2} \Lsmall^i \right\|_{L^{\infty}(\twoargroughtori{\timefunction,u}{\muxmulevelsetvalue})},
		\notag \\
	\left\| \argLrough{\muxmulevelsetvalue} \comder^{[1,4];3}\Lsmall^i \right\|_{L^{\infty}(\twoargroughtori{\timefunction,u}{\muxmulevelsetvalue})}, 
		\,
	\left\| \comdersmall^{[1,4];3} \Lsmall^i \right\|_{L^{\infty}(\twoargroughtori{\timefunction,u}{\muxmulevelsetvalue})}
	& 
	\leq C \fundbootsmall.
		\notag
\end{align}

\medskip

\noindent \underline{\textbf{$L^{\infty}$ estimates tied to pure transversal derivatives}}. 
\begin{align}
	\left 
		\| \muX^M \RRiemann 
	\right\|_{L^{\infty}(\twoargroughtori{\timefunction,u}{\muxmulevelsetvalue})} 
	& 
	\leq 
	\left\| 
		\muX^M \RRiemann 
	\right\|_{L^{\infty}(\twoargroughtori{\timefunction_0,u}{\muxmulevelsetvalue})} 
	+ 
	C \fundbootsmall, 
	&
	1 \leq M \leq 4,
		\label{E:LINFINITYIMPROVEMENTAUXTRANSVERSALPDERIVATIVESRRIEMANNLARGE}  
			\\
	\left 
		\| \muX^M \wavearraypartial 
	\right\|_{L^{\infty}(\twoargroughtori{\timefunction,u}{\muxmulevelsetvalue})} 
	& 
	\leq 
	C \fundbootsmall, 
	&
	1 \leq M \leq 4,
		\label{E:LINFINITYIMPROVEMENTAUXTRANSVERSALPDERIVATIVESPARTIALWAVEARRAYSMALL}  
			\\
	\left\| 
		\muX^M \upmu 
	\right\|_{L^{\infty}(\twoargroughtori{\timefunction,u}{\muxmulevelsetvalue})} 
	& 
	\leq 
	\left\|
			\muX^M
			\left\lbrace
				\Speed^{-1}
			\right\rbrace
		\right\|_{L^{\infty}(\twoargroughtori{\timefunction_0,u}{\muxmulevelsetvalue})}
		+
		\frac{1}{2}
		\mathring{\updelta}_*^{-1}
		\left\|
				\muX^M
				\left\lbrace
				\Speed^{-1}(\Speed^{-1} \Speed_{;\LogDensity} + 1)  
				\muX \RRiemann 
				\right\rbrace
			\right\|_{L^{\infty}(\twoargroughtori{\timefunction_0,u}{\muxmulevelsetvalue})}
	+ 
	C \fundbootsmall, 
		\label{E:LINFINITYIMPROVEMENTAUXMUANDTRANSVERSALDERIVATIVES} 
	&
	0 \leq M \leq 3,
		\\
	\left\|
		\argLrough{\muxmulevelsetvalue} \muX^M \upmu 
	\right\|_{L^{\infty}(\twoargroughtori{\timefunction,u}{\muxmulevelsetvalue})}
	& 
	\leq 
	\frac{1}{\updelta_*} 
			\left\|
				\muX^M
				\left\lbrace
				\Speed^{-1}(\Speed^{-1} \Speed_{;\LogDensity} + 1)  
				\muX \RRiemann 
				\right\rbrace
			\right\|_{L^{\infty}(\twoargroughtori{\timefunction_0,u}{\muxmulevelsetvalue})}
	+ 
	C \fundbootsmall, 
		\label{E:LINFINITYIMPROVEMENTAUXNULLDERIVATIVEMUANDTRANSVERSALDERIVATIVES}  
	&
	0 \leq M \leq 3,
		\\
	\left\| 
		\muX^M \Lsmall^i
	\right\|_{L^{\infty}(\twoargroughtori{\timefunction,u}{\muxmulevelsetvalue})} 
	& 
	\leq	 
	C, 
		\label{E:LINFINITYESTIMATESFORMUXDERIVATIVESOFLI}
	&
	1 \leq M \leq 3,
\end{align}
where $\argLrough{\muxmulevelsetvalue}$ is defined in \eqref{E:LROUGH}.

\end{proposition}


\begin{proof} 
We refer to Sect.\,\ref{SS:SILENTFACTS} for various results that we will silently use throughout the analysis.
We will also silently use the assumptions on the initial data stated in 
Sect.\,\ref{SS:ASSUMPTIONSONDATA},
the parameter-relations \eqref{E:DATAEPSILONISSMALLERTHANBOOTSTRAPEPSILONSMALLERTHANSQUAREOFDATAALPHA},
and the estimate 
$
\frac{1}{\Lunit \timefunctionarg{\muxmulevelsetvalue}}
\approx 1
$
implied by \eqref{E:BALDERIVATIVEOFROUGHTIMEFUNCTIONISAPPROXIMATELYUNITY}.

\medskip

\noindent \textbf{Proof of \eqref{E:LINFINITYIMPROVEMENTAUXR+}:}
We use \eqref{E:LINFINITYTRANSPORTROUGHLFESTIMATE} with $\timefunction_1 \eqdef \timefunction_0$ and $\timefunction_2 \eqdef \timefunction$,
the bootstrap assumptions \eqref{E:FUNDAMENTALQUANTITATIVEBOOTWAVEANDTRANSPORT},
the data estimate \eqref{E:LINFINITYINITIALROUGHHYPERSURFACEBOUNDRRIEMANNAMPLITUDE},
and \eqref{E:DATAEPSILONISSMALLERTHANBOOTSTRAPEPSILONSMALLERTHANSQUAREOFDATAALPHA}
to conclude that:
\begin{align} \label{E:PROOFSTEPSIMPROVEAUXRRIEMANN}
\| \RRiemann \|_{L^{\infty}(\twoargroughtori{\timefunction,u}{\muxmulevelsetvalue})} 
& 
\leq
\| \RRiemann \|_{L^{\infty}(\twoargroughtori{\timefunction_0,u}{\muxmulevelsetvalue})} 	
+
\mupositive \sup_{\timefunction' \in [\timefunction_0,\timefunction]} 
\| \argLrough{\muxmulevelsetvalue} \RRiemann \|_{L^{\infty}\left(\twoargroughtori{\timefunction',u}{\muxmulevelsetvalue}\right)}
\leq 
\mathring{\upalpha} + C \fundbootsmall
\end{align}
as desired.

\medskip

\noindent \textbf{Proof of \eqref{E:LINFINITYIMPROVEMENTAUXL1SMALL}--\eqref{E:LINFINITYIMPROVEMENTAUXLASMALL}:}
We first use Lemma~\ref{L:SCHEMATICSTRUCTUREOFVARIOUSTENSORSINTERMSOFCONTROLVARS} 
to write the transport equation \eqref{E:LUNITITRANSPORT} for $\Lunit^i$ schematically as: 
\begin{align} \label{E:LUNITITRANSPORTSCHEMATIC}
	\Lunit \Lsmall^i
	& = \smoothfunction(\controlvars) 
			\cdot
			\tander \wavearray.
\end{align}
From \eqref{E:LUNITITRANSPORTSCHEMATIC} and the bootstrap assumptions, 
we find that
$\|\argLrough{\muxmulevelsetvalue} \Lunit^i \|_{L^{\infty}(\twoargroughtori{\timefunction,u}{\muxmulevelsetvalue})} 
= \mathcal{O}(\fundbootsmall)$. From this bound, 
the initial data assumptions, and the same arguments we used to prove \eqref{E:LINFINITYIMPROVEMENTAUXR+},
we conclude \eqref{E:LINFINITYIMPROVEMENTAUXL1SMALL}--\eqref{E:LINFINITYIMPROVEMENTAUXLASMALL}.

\medskip

\noindent \textbf{Proof of \eqref{E:LINFINITYIMPROVEMENTAUXTRANSVERSALPDERIVATIVESRRIEMANNLARGE} and
\eqref{E:LINFINITYIMPROVEMENTAUXTRANSVERSALPDERIVATIVESPARTIALWAVEARRAYSMALL} for $M = 1$}: 
Fix $\Psi \in \wavearray = (\RRiemann,\LRiemann,v^2,v^3,\Ent)$.
Multiplying \eqref{E:SCHEMATICREWRITINGOFWAVEEQUATIONSSATISFIEDBYWAVEVARIABLES}
by $\frac{1}{\Lunit \timefunctionarg{\muxmulevelsetvalue}} \approx 1$ and using the bootstrap assumptions, we see that:  
\begin{align} \label{E:LXPSISMALL}
 \left\| \argLrough{\muxmulevelsetvalue} \muX \Psi \right \|_{L^{\infty}\left(\twoargroughtori{\timefunction,u}{\muxmulevelsetvalue}\right)} 
& \leq 
C \fundbootsmall.
\end{align}
From \eqref{E:LXPSISMALL} and the same arguments we used to prove \eqref{E:LINFINITYIMPROVEMENTAUXR+},
we conclude 
\eqref{E:LINFINITYIMPROVEMENTAUXTRANSVERSALPDERIVATIVESRRIEMANNLARGE} and
\eqref{E:LINFINITYIMPROVEMENTAUXTRANSVERSALPDERIVATIVESPARTIALWAVEARRAYSMALL} 
for $M = 1$.

\medskip

\noindent \textbf{Proof of \eqref{E:LINFINITYIMPROVEMENTAUXWAVEARRAY} 
for $\left\| \argLrough{\muxmulevelsetvalue} \comder^{\leq \Ntop-11;1} \wavearray \right\|_{L^{\infty}(\twoargroughtori{\timefunction,u}{\muxmulevelsetvalue})}$
and $\left\| \comdersmall^{[1,\Ntop-11];1} \wavearray \right\|_{L^{\infty}(\twoargroughtori{\timefunction,u}{\muxmulevelsetvalue})}$}: 
We first commute \eqref{E:SCHEMATICREWRITINGOFWAVEEQUATIONSSATISFIEDBYWAVEVARIABLES} 
with $\tander^{[1,\Ntop-12]}$ and use the commutator estimate 
\eqref{E:POINTWISEBOUNDCOMMUTATORSTANGENTIALANDTANGENTIALDERIVATIVESONSCALARFUNCTION} 
and the bootstrap assumptions 
to deduce that
$ 
\left\| \Lunit \tander^{[1,\Ntop -12]} \muX \Psi \right\|_{L^{\infty}\left(\twoargroughtori{\timefunction,u}{\muxmulevelsetvalue}\right)}
\lesssim \fundbootsmall
$
and that the same bound holds with $\argLrough{\muxmulevelsetvalue}$ in place of $\Lunit$. 
In particular, this yields \eqref{E:LINFINITYIMPROVEMENTAUXWAVEARRAY} 
for the first term on the LHS.
Also using the same arguments we used to prove \eqref{E:LINFINITYIMPROVEMENTAUXR+}
(including our assumptions on the data), 
we find that
$ 
\left\| 
	\tander^{[1,\Ntop -12]} \muX \Psi \right\|_{L^{\infty}\left(\twoargroughtori{\timefunction,u}{\muxmulevelsetvalue}\right)}
\lesssim \fundbootsmall
$.
From this bound, the commutator estimate \eqref{E:POINTWISEBOUNDCOMMUTATORSMUXANDTANGENTIALDERIVATIVESONSCALARFUNCTION}, 
and the bootstrap assumptions, we deduce that
$\left\| \comdersmall^{[1,\Ntop-11];1} \wavearray \right\|_{L^{\infty}(\twoargroughtori{\timefunction,u}{\muxmulevelsetvalue})}
\lesssim 
\left\| \tander^{[1,\Ntop -12]} \muX \Psi \right\|_{L^{\infty}\left(\twoargroughtori{\timefunction,u}{\muxmulevelsetvalue}\right)}
+
\left\| \tander^{[1,\Ntop -12]} \Psi \right\|_{L^{\infty}\left(\twoargroughtori{\timefunction,u}{\muxmulevelsetvalue}\right)}
\lesssim 
\fundbootsmall
$,
which yields \eqref{E:LINFINITYIMPROVEMENTAUXWAVEARRAY}
for the second term on the LHS.

\medskip

\noindent \textbf{Proof of \eqref{E:LINFINITYIMPROVEMENTAUXMUSMALL} for 
$\left\| \argLrough{\muxmulevelsetvalue} \tander^{[1,\Ntop-12]} \upmu \right\|_{L^{\infty}(\twoargroughtori{\timefunction,u}{\muxmulevelsetvalue})}$
and $\left\| \tander_*^{[1,\Ntop-12]} \upmu \right\|_{L^{\infty}(\twoargroughtori{\timefunction,u}{\muxmulevelsetvalue})}$}: 
We begin by using Lemmas~\ref{L:IDENTIFYFORFACTORDRIVINGSHOCK} and
\ref{L:SCHEMATICSTRUCTUREOFVARIOUSTENSORSINTERMSOFCONTROLVARS} 
to write the transport equation \eqref{E:MUTRANSPORT} for $\upmu$ as: 
\begin{align} \label{E:MUTRANSPORTSCHEMATICVERSION}
	\Lunit \upmu
	& = 
	- \frac{1}{2} \Speed^{-1}(\Speed^{-1} \Speed_{;\LogDensity} + 1)  
		\muX \RRiemann
	+
	\smoothfunction(\controlvars) \cdot \muX \wavearraypartial
	+  
	\smoothfunction(\badcontrolvars) \cdot \tander \wavearray
	+ 
	\smoothfunction(\badcontrolvars) \cdot \GradEnt,
\end{align}
where the first product on RHS~\eqref{E:MUTRANSPORTSCHEMATICVERSION} is written exactly (for use later on)
and the remaining ones schematically.
Commuting \eqref{E:MUTRANSPORTSCHEMATICVERSION} with $\tander^N$ for $1 \leq N \leq \Ntop -12$ 
and using \eqref{E:POINTWISEBOUNDCOMMUTATORSTANGENTIALANDTANGENTIALDERIVATIVESONSCALARFUNCTION},
the bootstrap assumptions, 
and the already 
proved bound 
$
\left\| \tander^{[1,\Ntop -12]} \muX \Psi \right\|_{L^{\infty}\left(\twoargroughtori{\timefunction,u}{\muxmulevelsetvalue}\right)}
\lesssim
\fundbootsmall
$, 
we find that
$ 
\left\| \Lunit \tander^{[1,\Ntop-12]} \upmu \right \|_{L^{\infty}\left(\twoargroughtori{\timefunction,u}{\muxmulevelsetvalue}\right)}
\lesssim \fundbootsmall
$
and that the same bound holds with $\argLrough{\muxmulevelsetvalue}$ in place of $\Lunit$.  
In particular, this yields \eqref{E:LINFINITYIMPROVEMENTAUXMUSMALL} 
for the first term on the LHS.
In the case that $\tander^N = \tandersmall^N$, 
we can combine this bound with the initial data assumptions 
and the same arguments we used to prove \eqref{E:LINFINITYIMPROVEMENTAUXR+}
in order to conclude 
$\left\| \tandersmall^{[1,\Ntop-12]} \upmu \right\|_{L^{\infty}(\twoargroughtori{\timefunction,u}{\muxmulevelsetvalue})}
\lesssim \fundbootsmall
$,
which yields
\eqref{E:LINFINITYIMPROVEMENTAUXMUSMALL} 
for the second term on the LHS.

\medskip

\noindent \textbf{Proof of \eqref{E:LINFINITYIMPROVEMENTAUXMUANDTRANSVERSALDERIVATIVES} and
\eqref{E:LINFINITYIMPROVEMENTAUXNULLDERIVATIVEMUANDTRANSVERSALDERIVATIVES} in the case $M=0$}: 
The arguments given in the previous paragraph in particular yield that
$
\left\| \argLrough{\muxmulevelsetvalue} \Lunit \upmu \right\|_{L^{\infty}(\twoargroughtori{\timefunction,u}{\muxmulevelsetvalue})}
\lesssim \fundbootsmall
$.
From this bound and the same arguments we used to prove \eqref{E:LINFINITYIMPROVEMENTAUXR+},
we find that
$
\left\|
	\Lunit \upmu 
\right\|_{L^{\infty}(\twoargroughtori{\timefunction,u}{\muxmulevelsetvalue})}
\leq 
\left\| 
	\Lunit \upmu 
\right\|_{L^{\infty}(\twoargroughtori{\timefunction_0,u}{\muxmulevelsetvalue})} 
+ 
C \fundbootsmall
$.
From this bound, \eqref{E:MUTRANSPORTSCHEMATICVERSION}, 
the already proven bound
$
\left 
	\| \muX \wavearraypartial 
\right\|_{L^{\infty}(\twoargroughtori{\timefunction,u}{\muxmulevelsetvalue})} 
\lesssim 
\fundbootsmall
$,
and the bootstrap assumptions,
we further deduce that
$
\left\|
	\Lunit \upmu 
\right\|_{L^{\infty}(\twoargroughtori{\timefunction,u}{\muxmulevelsetvalue})}
\leq 
\frac{1}{2}
\left\| 
	\Speed^{-1}(\Speed^{-1} \Speed_{;\LogDensity} + 1)  
		\muX \RRiemann
\right\|_{L^{\infty}(\twoargroughtori{\timefunction_0,u}{\muxmulevelsetvalue})} 
+ 
C \fundbootsmall
$.
From this bound and \eqref{E:CLOSEDVERSIONLUNITROUGHTTIMEFUNCTION},
we find that
$
\left\|
	\argLrough{\muxmulevelsetvalue} \upmu 
\right\|_{L^{\infty}(\twoargroughtori{\timefunction,u}{\muxmulevelsetvalue})}
\leq
\frac{2}{\updelta_*} 
\left\| 
	 \Lunit \upmu 
\right\|_{L^{\infty}(\twoargroughtori{\timefunction_0,u}{\muxmulevelsetvalue})} 
+ 
C \fundbootsmall
$,
which yields \eqref{E:LINFINITYIMPROVEMENTAUXNULLDERIVATIVEMUANDTRANSVERSALDERIVATIVES}
in the case $M=0$. From this bound
and the same arguments we used to prove \eqref{E:LINFINITYIMPROVEMENTAUXR+},
we find that
$
\left\|
	\upmu 
\right\|_{L^{\infty}(\twoargroughtori{\timefunction,u}{\muxmulevelsetvalue})}
\leq
\left\|
	\upmu 
\right\|_{L^{\infty}(\twoargroughtori{\timefunction_0,u}{\muxmulevelsetvalue})}
+
\frac{2}{\updelta_*} 
\mupositive
\left\| 
	\Lunit \upmu 
\right\|_{L^{\infty}(\twoargroughtori{\timefunction_0,u}{\muxmulevelsetvalue})} 
+ 
C \fundbootsmall
$,
which (assuming $\mupositive \leq 1$) 
yields \eqref{E:LINFINITYIMPROVEMENTAUXMUANDTRANSVERSALDERIVATIVES} in the case $M=0$.

\medskip

\noindent \textbf{Proof of \eqref{E:LINFINITYIMPROVEMENTAUXMIXECTANGENTIALTRANSVERSALDERIVATIVESLISMALL}
for $\left\| \argLrough{\muxmulevelsetvalue} \tander^{\leq \Ntop-11} \Lsmall^i \right\|_{L^{\infty}(\twoargroughtori{\timefunction,u}{\muxmulevelsetvalue})}$
	and
	$\left\| \tander^{[1,\Ntop-11]} \Lsmall^i \right\|_{L^{\infty}(\twoargroughtori{\timefunction,u}{\muxmulevelsetvalue})}$}: 
We commute \eqref{E:LUNITITRANSPORTSCHEMATIC} with $\tander^N$ for $1 \leq N \leq \Ntop - 11$ to obtain: 
\begin{align} \label{E:COMMUTEDWITHTANGENTIALLUNITITRANSPORTSCHEMATIC}
	\Lunit \tander^N \Lunit^i
	& = 	
			[\Lunit, \tander^N] \Lunit^i
			+
			\tander^N
			\left\lbrace
			\smoothfunction(\controlvars) 
			\cdot
			\tander \wavearray
			\right\rbrace.
\end{align}
From \eqref{E:COMMUTEDWITHTANGENTIALLUNITITRANSPORTSCHEMATIC}, 
the commutator estimate \eqref{E:POINTWISEBOUNDCOMMUTATORSTANGENTIALANDTANGENTIALDERIVATIVESONSCALARFUNCTION},
and the bootstrap assumptions,
we find that
$\left\| \argLrough{\muxmulevelsetvalue} \tander^{\leq \Ntop-11} \Lsmall^i \right\|_{L^{\infty}(\twoargroughtori{\timefunction,u}{\muxmulevelsetvalue})}
\lesssim \fundbootsmall
$
and that the same bound holds with $\argLrough{\muxmulevelsetvalue}$ in place of $\Lunit$. 
The remainder of the proof relies on the same arguments used above.

\medskip

\noindent \textbf{Proof of \eqref{E:LINFINITYESTIMATESFORMUXDERIVATIVESOFLI} in the case $M=1$
and \eqref{E:LINFINITYIMPROVEMENTAUXMIXECTANGENTIALTRANSVERSALDERIVATIVESLISMALL}
for $\left\| \comdersmall^{[1,\Ntop-12];1} \Lsmall^i \right\|_{L^{\infty}(\twoargroughtori{\timefunction,u}{\muxmulevelsetvalue})}$}: 
The desired estimate \eqref{E:LINFINITYESTIMATESFORMUXDERIVATIVESOFLI} in the case $M=1$
follows from the identity \eqref{E:MUXLISCHEMATICIDENTITY}
and the bootstrap assumptions.

Next, we differentiate \eqref{E:MUXLISCHEMATICIDENTITY}
with $\tander^N$ for $1 \leq N \leq \Ntop - 13$
and use the bootstrap assumptions and the already
proven bounds
$\left\| \comdersmall^{[1,\Ntop-11];1} \wavearray \right\|_{L^{\infty}(\twoargroughtori{\timefunction,u}{\muxmulevelsetvalue})}
\lesssim 
\fundbootsmall
$,
$
\left\| \tandersmall^{[1,\Ntop-12]} \upmu \right\|_{L^{\infty}(\twoargroughtori{\timefunction,u}{\muxmulevelsetvalue})}
\lesssim \fundbootsmall
$,
and
$\left\|\tander^{[1,\Ntop-11]} \Lsmall^i \right\|_{L^{\infty}(\twoargroughtori{\timefunction,u}{\muxmulevelsetvalue})}
\lesssim \fundbootsmall
$
to deduce
$\left\| \tander^{[1,\Ntop-13]} \muX \Lsmall^i \right\|_{L^{\infty}(\twoargroughtori{\timefunction,u}{\muxmulevelsetvalue})}
\lesssim \fundbootsmall
$.
Also using the commutator estimate \eqref{E:POINTWISEBOUNDCOMMUTATORSMUXANDTANGENTIALDERIVATIVESONSCALARFUNCTION},
we conclude that
$\left\| \comdersmall^{[1,\Ntop-12];1} \Lsmall^i \right\|_{L^{\infty}(\twoargroughtori{\timefunction,u}{\muxmulevelsetvalue})}
\lesssim \fundbootsmall
$
as desired.

\medskip

\noindent \textbf{Proof of \eqref{E:LINFINITYIMPROVEMENTAUXMIXECTANGENTIALTRANSVERSALDERIVATIVESLISMALL}
for 
$\left\| \argLrough{\muxmulevelsetvalue} \comder^{[1,\Ntop-12];1} \Lsmall^i \right\|_{L^{\infty}(\twoargroughtori{\timefunction,u}{\muxmulevelsetvalue})}$}: 
We assume that $1 \leq N \leq \Ntop -12$ and that $\comder^{N;1}$ 
contains a factor of $\muX$. 
We commute \eqref{E:LUNITITRANSPORTSCHEMATIC} with $\comder^{N;1}$ to obtain:
\begin{align} \label{E:COMMUTEDWITHONETRANSVERSALANDTRANGENTIALLUNITITRANSPORTSCHEMATIC}
	\Lunit \comder^{N;1} \Lunit^i
	& = 	
			[\Lunit,\comder^{N;1}] \Lunit^i
			+
			\comder^{N;1}
			\left\lbrace
			\smoothfunction(\controlvars) 
			\cdot
			\tander \wavearray
			\right\rbrace.
\end{align}
Using the commutator estimate \eqref{E:POINTWISEBOUNDCOMMUTATORSMUXANDTANGENTIALDERIVATIVESONSCALARFUNCTION}
and the already proven bounds
$$
\left\|\tander^{[1,\Ntop-11]} \Lsmall^i \right\|_{L^{\infty}(\twoargroughtori{\timefunction,u}{\muxmulevelsetvalue})}
\lesssim \fundbootsmall,
	\qquad
\left\| \comdersmall^{[1,\Ntop-12];1} \Lsmall^i \right\|_{L^{\infty}(\twoargroughtori{\timefunction,u}{\muxmulevelsetvalue})}
\lesssim \fundbootsmall,
$$
we find that the first term on RHS~\eqref{E:COMMUTEDWITHONETRANSVERSALANDTRANGENTIALLUNITITRANSPORTSCHEMATIC}
is bounded in magnitude by $\lesssim \varepsilon$. From the bootstrap assumptions
and the already proven bound
$\left\| \comdersmall^{[1,\Ntop-11];1} \wavearray \right\|_{L^{\infty}(\twoargroughtori{\timefunction,u}{\muxmulevelsetvalue})}
\lesssim 
\fundbootsmall
$,
we deduce that the second term on RHS~\eqref{E:COMMUTEDWITHONETRANSVERSALANDTRANGENTIALLUNITITRANSPORTSCHEMATIC}
obeys the same bound. It follows that
$
\left\| \Lunit \comder^{[1,\Ntop - 12];1} \Lunit^i \right\|_{L^{\infty}\left(\twoargroughtori{\timefunction,u}{\muxmulevelsetvalue}\right)}
\lesssim \fundbootsmall
$
and that the same bound holds with $\argLrough{\muxmulevelsetvalue}$ in place of $\Lunit$, as desired.

\medskip
 
\noindent \textbf{Proof of \eqref{E:LINFINITYIMPROVEMENTAUXVORTGRADENT} for 
$\left\| 
		\comder^{\leq \Ntop-11;1} (\vortrenormalized,\GradEnt) 
	\right\|_{L^{\infty}(\twoargroughtori{\timefunction,u}{\muxmulevelsetvalue})}$
	and
	\eqref{E:LINFINITYIMPROVEMENTAUXMODIFIEDFLUID}  
	for
	$
	\left\| 
		\comder^{\leq \Ntop-12;1} (\VortVort,\DivGradEnt) 
	\right\|_{L^{\infty}(\twoargroughtori{\timefunction,u}{\muxmulevelsetvalue})}
	$}: 
We apply $\tander^{\leq \Ntop - 11}$ to 
\eqref{E:VORTICITYTRANSVERSALTRANSPORTINTERMSOFTANGENTIAL}--\eqref{E:ENTROPYGRADIENTTRANSVERSALTRANSPORTINTERMSOFTANGENTIAL} 
and use the bootstrap assumptions \eqref{E:FUNDAMENTALQUANTITATIVEBOOTWAVEANDTRANSPORT} 
to deduce that
$
\left\| 
	\tander^{\leq \Ntop - 11 }\muX(\vortrenormalized, \GradEnt) 
\right\|_{L^{\infty}(\twoargroughtori{\timefunction,u}{\muxmulevelsetvalue})}
\lesssim \fundbootsmall
$.
From this bound,
the commutator estimate \eqref{E:POINTWISEBOUNDCOMMUTATORSMUXANDTANGENTIALDERIVATIVESONSCALARFUNCTION},
and the bootstrap assumptions,
we conclude \eqref{E:LINFINITYIMPROVEMENTAUXVORTGRADENT} for the first term on the LHS.

The estimate \eqref{E:LINFINITYIMPROVEMENTAUXMODIFIEDFLUID}   
for 
$
\left\| 
		\comder^{\leq \Ntop-12;1} (\VortVort,\DivGradEnt) 
	\right\|_{L^{\infty}(\twoargroughtori{\timefunction,u}{\muxmulevelsetvalue})}
$
follows from a similar argument based on
\eqref{E:TRANSVERSALMODIFIEDINTERMSOFTANGENTIAL}.

\medskip

\noindent \textbf{Outline of the remainder of the proof}:
The remaining estimates in the proposition can be derived
by commuting the equations with one additional $\muX$ 
derivative followed by elements of $\Tanset$, 
using the above arguments and the estimates already proved,
and then repeating the process, adding one additional $\muX$ derivative each time. 
The estimates need to be derived in the same order
as above. The commutator terms can be handled with the help of
Prop.\,\ref{P:COMMUTATORESTIMATES}. 
The bootstrap assumptions guarantee that all terms that need to be controlled have 
sufficient $L^{\infty}(\twoargroughtori{\timefunction,u}{\muxmulevelsetvalue})$-regularity.
We omit the tedious, but straightforward details.

\end{proof}

The following corollary is an immediate consequence of the fact that we have improved the auxiliary bootstrap assumptions. 

\begin{corollary}[$\mathring{\upalpha}^{1/2}$ and $\auxbootsmall$ can be replaced by $C\mathring{\upalpha}$ and  $C \fundbootsmall$] 
\label{C:IMPROVEAUX}
All prior inequalities whose RHS feature an explicit factor of $\mathring{\upalpha}^{1/2}, \auxbootsmall$ remain true with $\mathring{\upalpha}^{1/2},\, \auxbootsmall$ respectively replaced by $C \mathring{\upalpha},\, C \fundbootsmall$.
\end{corollary}


\section{Sharp control of $\upmu$ and the properties of $\Upsilon$}
\label{S:SHARPCONTROLOFMUANDPROPERTIESOFCHOVGEOTOCARTESIAN}
In this section, we derive sharp control of $\upmu$ and its derivatives, as well as strict improvements 
of the bootstrap assumptions from Sects.\,\ref{SSS:BAFORINVERSEFOLIATIONDENSITY}--\ref{SSS:BOOTSTRAPASSUMPTIONFORTORISTRUCTURE}. 
Because our $L^2$-type energies feature $\upmu$ weights,
these sharp estimates will play a fundamental role in our proof of the energy estimates.
We also derive homeomorphism and diffeomorphism properties 
of the change of variables map from geometric coordinate to Cartesian coordinates,
which are crucial for understanding the structure of the singular boundary.

\subsection{Sharp control of $\upmu$ and its derivatives}
\label{SS:SHARPCONTROLOFMUANDDERIVATIVES}

\begin{proposition}[Sharp control of $\upmu$ and its derivatives]
\label{P:SHARPCONTROLOFMUANDDERIVATIVES}
The following estimates hold.

\medskip

\noindent \underline{\textbf{Minima of $\upmu$ occur precisely along $\twoargmumuxtorus{-\timefunction}{-\muxmulevelsetvalue}$}}.
For each fixed $\timefunction \in [\timefunction_0,\timefunctionboot] = [-\mupositive,-\upmuboot]$,
we have:
		\begin{align} \label{E:MINVALUEOFMUONFOLIATION}
			\min_{\hypthreearg{\timefunction}{[- \rightu,\leftu]}{\muxmulevelsetvalue}}
			\upmu
			& = 
			 - \timefunction,
		\end{align}
		and the minimum value of $- \timefunction$
		in \eqref{E:MINVALUEOFMUONFOLIATION} is achieved by $\upmu$ precisely on the $\upmu$-adapted torus 
		$
		\twoargmumuxtorus{-\timefunction}{-\muxmulevelsetvalue}
		$ defined in \eqref{E:MUXMUTORI}.
		In particular,
		$ 
		\upmuboot
		=
		\inf_{\twoargMrough{[\timefunction_0,\timefunctionboot),[- \rightu,\leftu]}{\muxmulevelsetvalue}} \upmu$,
		and the infimum is not achieved.

\medskip		
		
\noindent \underline{\textbf{$\upmu$ is large when $|u| \geq \interestingu$}}.
	The following lower bound holds, where $\boringregionmupositive > 0$ is the constant appearing in 
	\eqref{E:DATAMUISLARGEINBORINGREGION}:
	\begin{align} \label{E:MUISLARGEINBORINGREGION}
			\min_{\twoargMrough{[\timefunction_0,\timefunctionboot],[- \rightu,\leftu]}{\muxmulevelsetvalue}
			\backslash 
			\twoargMrough{[\timefunction_0,\timefunctionboot],[-\interestingu,\interestingu]}{\muxmulevelsetvalue}} \upmu
			& \geq \frac{\boringregionmupositive}{2}.
\end{align}
	
\medskip

\noindent \underline{\textbf{Location of $\datahypfortimefunctiontwoarg{-\muxmulevelsetvalue}{[\timefunction_0,\timefunctionboot]}$ 
and $\twoargmumuxtorus{-\timefunction}{-\muxmulevelsetvalue}$}}.
With
$\datahypfortimefunctiontwoarg{-\muxmulevelsetvalue}{[\timefunction_0,\timefunctionboot]}$
and
$\twoargmumuxtorus{-\timefunction}{-\muxmulevelsetvalue}$
denoting the sets defined in \eqref{E:LEVELSETSOFMUXMU}--\eqref{E:MUXMUTORI},
we have:
\begin{subequations}
\begin{align} \label{E:MUXMUKAPPALEVELSETLOCATION}
			\datahypfortimefunctiontwoarg{-\muxmulevelsetvalue}{[\timefunction_0,\timefunctionboot]}
			& 
			\subset
			\twoargMrough{[\timefunction_0,\timefunctionboot],[-\frac{1}{2}\interestingu,\frac{1}{2}\interestingu]}{\muxmulevelsetvalue},
				\\
		\text{for each $ \timefunction \in [\timefunction_0,\timefunctionboot]$,} \qquad
		\twoargmumuxtorus{-\timefunction}{-\muxmulevelsetvalue}, 
		& 
		\subset
		\hypthreearg{\timefunction}{[-\frac{1}{2}\interestingu,\frac{1}{2}\interestingu]}{\muxmulevelsetvalue}.
			\label{E:IMPROVEDLEVELSETSTRUCTUREANDLOCATIONOFMIN} 
\end{align}
\end{subequations}

Moreover, with $\InverseCHOVroughtomumuxmu{\muxmulevelsetvalue}$
denoting the inverse function of the function
$\CHOVroughtomumuxmu{\muxmulevelsetvalue}$ defined in \eqref{E:CHOVFROMROUGHCOORDINATESTOMUWEGIGHTEDXMUCOORDINATES},
we have:
\begin{align} \label{E:INTERMSOFPHICHOVMAPIMPROVEDLEVELSETSTRUCTUREANDLOCATIONOFMIN} 
\text{For each $\mulevelsetvalue \in [\mupositive,\upmuboot]$,} \qquad
	\InverseCHOVroughtomumuxmu{\muxmulevelsetvalue}
	\left(\lbrace \mulevelsetvalue \rbrace \times \lbrace - \muxmulevelsetvalue \rbrace \times \mathbb{T}^2 \right)
	\subset
	\lbrace - \mulevelsetvalue \rbrace \times \left[-\frac{\interestingu}{2},\frac{\interestingu}{2}\right] \times \mathbb{T}^2.
\end{align}

\medskip

\noindent \underline{\textbf{Transversal convexity of $\upmu$ and its consequences}}.
The following estimates hold,
where $\Wtransarg{\muxmulevelsetvalue}$ is the vectorfield from Def.\,\ref{D:WTRANSANDCUTOFF}:
\begin{align} 
	\begin{split} \label{E:MUTRANSVERSALCONVEXITY}
			\frac{\secondtransversalderivativemulowerbound}{2}
			& \leq 
			\min_{\twoargMrough{[\timefunction_0,\timefunctionboot],[-\interestingu,\interestingu]}{\muxmulevelsetvalue}}
				\left\lbrace
				\Wtransarg{\muxmulevelsetvalue} \Wtransarg{\muxmulevelsetvalue} \upmu,
					\,
				\Wtransarg{\muxmulevelsetvalue} \muX \upmu,
					\,
				\muX \muX \upmu,
					\,
				\muX \muX \upmu - \frac{(\muX \upmu) \Lunit \muX \upmu}{\Lunit \upmu},
					\,
				\geop{u} \muX \upmu,
					\,
				\geop{u} \muX \upmu - \frac{(\geop{u} \upmu) \geop{t} \muX \upmu}{\geop{t} \upmu}, 
					\,
				\roughgeop{u} \muX \upmu
				\right\rbrace
					\\
			&
			\leq 
			\max_{\twoargMrough{[\timefunction_0,\timefunctionboot],[-\interestingu,\interestingu]}{\muxmulevelsetvalue}}
			\left\lbrace
				\Wtransarg{\muxmulevelsetvalue} \Wtransarg{\muxmulevelsetvalue} \upmu,
					\,
				\Wtransarg{\muxmulevelsetvalue} \muX \upmu,
					\,
				\muX \muX \upmu,
					\,
				\muX \muX \upmu - \frac{(\muX \upmu) \Lunit \muX \upmu}{\Lunit \upmu},
					\,
				\geop{u} \muX \upmu,
					\,
				\geop{u} \muX \upmu - \frac{(\geop{u} \upmu) \geop{t} \muX \upmu}{\geop{t} \upmu}, 
					\,
				\roughgeop{u} \muX \upmu
				\right\rbrace
			\leq 
			\frac{2}{\secondtransversalderivativemulowerbound},
	\end{split}
	\end{align}

\begin{align} 
	\min_{\twoargMrough{[\timefunction_0,\timefunctionboot],[-\interestingu,\interestingu]}{\muxmulevelsetvalue}
			\backslash
			\mathcal{M}^{(\muxmulevelsetvalue)}_{[\timefunction_0,\timefunctionboot],[-\frac{\interestingu}{2},\frac{\interestingu}{2}]}} 
			|\muX \upmu + \muxmulevelsetvalue| 
			& \geq \frac{\secondtransversalderivativemulowerbound \interestingu}{8}.
			\label{E:REGIONWHEREMUXMUKAPPALEVELSETISNOTLOCATED}
\end{align}

Moreover, the following pointwise estimates hold on
$
\twoargMrough{[\timefunction_0,\timefunctionboot],[-\interestingu,\interestingu]}{\muxmulevelsetvalue}
$, where $\Rtransarg{\muxmulevelsetvalue}$ is the vectorfield defined in \eqref{E:RTRANS}:
\begin{align} \label{E:WMANDRTRANSMUUBOUNDEDBYSQRTMU}
	|\Wtransarg{\muxmulevelsetvalue} \upmu|,
			\,
	|\Rtransarg{\muxmulevelsetvalue} \upmu|
	& \leq 
	C \sqrt{\upmu}.
	\end{align}

\medskip

\noindent \underline{\textbf{Rough control of null and almost null derivatives of $\upmu$ in the interesting region}}.
The following estimates hold,
where $\argLrough{\muxmulevelsetvalue}$ is the vectorfield defined in \eqref{E:LROUGH}:
\begin{subequations}
\begin{align} \label{E:BOUNDSONLMUINTERESTINGREGION} 
	- 
	\frac{9}{8}
	\mathring{\updelta}_*
	\leq
	\min_{\twoargMrough{[\timefunction_0,\timefunctionboot],[-\interestingu,\interestingu]}{\muxmulevelsetvalue}} \Lunit \upmu
	\leq 
	\max_{\twoargMrough{[\timefunction_0,\timefunctionboot],[-\interestingu,\interestingu]}{\muxmulevelsetvalue}} \Lunit \upmu
	\leq 
	- 
	\frac{7}{8}
	\mathring{\updelta}_*,
		\\
	- 
	\frac{9}{8}
	\mathring{\updelta}_*
	\leq
	\min_{\twoargMrough{[\timefunction_0,\timefunctionboot],[-\interestingu,\interestingu]}{\muxmulevelsetvalue}} \geop{t} \upmu
	\leq 
	\max_{\twoargMrough{[\timefunction_0,\timefunctionboot],[-\interestingu,\interestingu]}{\muxmulevelsetvalue}} \geop{t} \upmu
	\leq 
	- 
	\frac{7}{8}
	\mathring{\updelta}_*,
		\label{E:BOUNDSONGEOMETRICTDERIVATIVEMUINTERESTINGREGION}
		\\
	- 
	\frac{3}{2}
	\leq
	\min_{\twoargMrough{[\timefunction_0,\timefunctionboot],[-\interestingu,\interestingu]}{\muxmulevelsetvalue}} \argLrough{\muxmulevelsetvalue} \upmu
	\leq 
	\max_{\twoargMrough{[\timefunction_0,\timefunctionboot],[-\interestingu,\interestingu]}{\muxmulevelsetvalue}} \argLrough{\muxmulevelsetvalue} \upmu
	\leq 
	- \frac{2}{3}.
	\label{E:SIZEBOUNDSONWIDETILDELMUINTERESTINGREGION}
	\end{align}
	\end{subequations}
	
	\noindent \underline{\textbf{Sharp control of $\Lunit \timefunctionarg{\muxmulevelsetvalue}$}}.
	The following estimates hold:
	\begin{align} \label{E:ROUGHTIMEFUNCTIONLDERIVATIVEBOUNDS}
	\frac{7}{9} 
	\mathring{\updelta}_*
	\leq
	\min_{\twoargMrough{[\timefunction_0,\timefunctionboot],[- \rightu,\leftu]}{\muxmulevelsetvalue}} 
		\Lunit \timefunctionarg{\muxmulevelsetvalue}
	\leq 
	\max_{\twoargMrough{[\timefunction_0,\timefunctionboot],[- \rightu,\leftu]}{\muxmulevelsetvalue}} 
		\Lunit \timefunctionarg{\muxmulevelsetvalue}
	\leq 
	\frac{9}{7} \mathring{\updelta}_*.
	\end{align}
	
	\medskip
	
	\noindent \underline{\textbf{Estimates for $\upmu$ along the flow map of $\FlowmapLrougharg{\muxmulevelsetvalue}$}}.
	Let $\FlowmapLrougharg{\muxmulevelsetvalue}$ denote the $\timefunction_0$-normalized
	flow map of $\argLrough{\muxmulevelsetvalue}$ from Lemma~\ref{L:PROPERTIESOFFLOWMAPOFWIDETILDEL}. 
	Then the following estimates hold:
	\begin{subequations}
	\begin{align} \label{E:MUMUSTDECREASEININTERESTINGREGION}
		\sup_{\substack{|u| \leq \interestingu 
					\\ 
					\timefunction_0 \leq \timefunction \leq \timefunction + \Delta \timefunction \leq \timefunctionboot 
					\\
					(x^2,x^3) \in \mathbb{T}^2}}
		\frac{\upmu \circ \FlowmapLrougharg{\muxmulevelsetvalue}(\timefunction + \Delta,u,x^2,x^3)}
			{\upmu \circ \FlowmapLrougharg{\muxmulevelsetvalue}(\timefunction,u,x^2,x^3)}
			& \leq 1,
			\\
		\sup_{\substack{u \in [- \rightu,\leftu] \backslash [-\interestingu,\interestingu] 
					\\ 
					\timefunction_0 \leq \timefunction \leq \timefunction + \Delta \timefunction \leq \timefunctionboot 
					\\
					(x^2,x^3) \in \mathbb{T}^2}}
		\frac{\upmu \circ \FlowmapLrougharg{\muxmulevelsetvalue}(\timefunction + \Delta,u,x^2,x^3)}
			{\upmu \circ \FlowmapLrougharg{\muxmulevelsetvalue}(\timefunction,u,x^2,x^3)}
		& \leq C.
		\label{E:MURATIOBOUNDEDINBORINGREGION}
	\end{align}
	\end{subequations}
		
	\medskip
	
	\noindent \underline{\textbf{A lower bound tied to the blowup in the interesting region}}.
	The following lower bounds holds:
	\begin{align} \label{E:LOWERBOUNDONMAGNITUDEOFXRPLUS}
	\min_{\twoargMrough{[\timefunction_0,\timefunctionboot],[-\interestingu,\interestingu]}{\muxmulevelsetvalue}}
	 \upmu |X \RRiemann|
	& \geq
		\frac{\mathring{\updelta}_*}{|\bar{\Speed}_{;\LogDensity} + 1|},
	\end{align}
	where $\bar{\Speed}_{;\LogDensity} \eqdef \Speed_{;\LogDensity}(\LogDensity = 0,\Ent=0)$
	is $\Speed_{;\LogDensity}$ evaluated at the trivial solution, 
	$\bar{\Speed}_{;\LogDensity} + 1$ is a \underline{non-zero} constant
	by assumption, and the vectorfield $X$ has Euclidean length satisfying
	$\sqrt{\sum_{a=1}^3 (X^a)^2} = 1 + \mathcal{O}(\mathring{\upalpha})$.
	
	\medskip
	
	\noindent \underline{\textbf{Especially sharp control in a small neighborhood}}.
	Recall that $\flowmapWtransargtwoarg{\muxmulevelsetvalue}{\Delta u}$ is the flow map of $\Wtransarg{\muxmulevelsetvalue}$
	(see Lemma~\ref{L:FLOWMAPFORGENERATOROFROUGHTIMEFUNCTION}).
	There exists a constant $\Delta U$ with $0 < \Delta U < \frac{\interestingu}{2}$,
	a neighborhood (in the geometric coordinate topology of $\twoargMrough{[\timefunction_0,\timefunctionboot],[- \rightu,\leftu]}{\muxmulevelsetvalue}$)
	$\smallneighborhoodofcreasetwoarg{[\timefunction_0,\timefunctionboot]}{\muxmulevelsetvalue}$ of 
	$\datahypfortimefunctiontwoarg{-\muxmulevelsetvalue}{[\timefunction_0,\timefunctionboot]}$
	of the form:
	\begin{align} \label{E:SMALLNEIGHBORHOOD}
		\smallneighborhoodofcreasetwoarg{[\timefunction_0,\timefunctionboot]}{\muxmulevelsetvalue}
		=
		\flowmapWtransargtwoarg{\muxmulevelsetvalue}{(-\Delta U,\Delta U)}
		\left(\datahypfortimefunctiontwoarg{-\muxmulevelsetvalue}{[\timefunction_0,\timefunctionboot]} \right)
		\eqdef 
		\bigcup_{\Delta u' \in (-\Delta U,\Delta U)}
		\flowmapWtransargtwoarg{\muxmulevelsetvalue}{\Delta u'}
		\left( \datahypfortimefunctiontwoarg{-\muxmulevelsetvalue}{[\timefunction_0,\timefunctionboot]} \right)
	\end{align}
	such that:
	\begin{align} \label{E:SMALLNEIGHBORHOODCONTAINEDININTERESTINGREGION}
	\smallneighborhoodofcreasetwoarg{[\timefunction_0,\timefunctionboot]}{\muxmulevelsetvalue} 
	&
	\subset 
	\twoargMrough{[\timefunction_0,\timefunctionboot],[-\interestingu,\interestingu]}{\muxmulevelsetvalue},
	\end{align}
	and a constant $\awayfromsmallneightborhoodnmupositive > 0$ defined by:
	\begin{align} \label{E:MU2DEF}
		\awayfromsmallneightborhoodnmupositive
		\eqdef \min \left\lbrace \frac{\secondtransversalderivativemulowerbound}{4} (\Delta U)^2, \frac{\boringregionmupositive}{2} \right\rbrace
	\end{align}
	(where $\boringregionmupositive$ is as in \eqref{E:DATAMUISLARGEINBORINGREGION})
	such that the following estimates hold:
	\begin{align} \label{E:WIDETILDELMUISALMOSTMINUSONEINSMALLNEIGHBORHOOD}
	- 1.01
	&
	\leq
	\min_{\smallneighborhoodofcreasetwoarg{[\timefunction_0,\timefunctionboot]}{\muxmulevelsetvalue}} \argLrough{\muxmulevelsetvalue} \upmu 
	\leq 
	\max_{\smallneighborhoodofcreasetwoarg{[\timefunction_0,\timefunctionboot]}{\muxmulevelsetvalue}} \argLrough{\muxmulevelsetvalue} \upmu 
	\leq 
	-.99,
		\\
	\awayfromsmallneightborhoodnmupositive
	& 
	\leq
	\min_{
	\twoargMrough{[\timefunction_0,\timefunctionboot],[- \rightu,\leftu]}{\muxmulevelsetvalue}
	\backslash
	\smallneighborhoodofcreasetwoarg{[\timefunction_0,\timefunctionboot]}{\muxmulevelsetvalue}} \upmu.
		\label{E:EASYREGIONLOWERBOUNDFORMU}
	\end{align}

\end{proposition}

\begin{proof}
	
\begin{remark}[Silent use of Lemma~\ref{L:CONTINUOUSEXTNESION}]
\label{R:SILENTUSEOFLEMMACONTINUOUSEXTNESION}
Throughout this proof, we sometimes silently use the continuous extension properties shown in Lemma~\ref{L:CONTINUOUSEXTNESION}.
In particular, these properties allow us to extend various results that
were proved prior to the lemma on the half-open rough time interval $[\timefunction_0,\timefunctionboot)$
to the closed rough time interval
$[\timefunction_0,\timefunctionboot]$.
\end{remark}
	
	\medskip
	
	\noindent \textbf{Proof of \eqref{E:MUTRANSVERSALCONVEXITY}}:
	To prove \eqref{E:MUTRANSVERSALCONVEXITY} for $\Wtransarg{\muxmulevelsetvalue} \Wtransarg{\muxmulevelsetvalue} \upmu$,
	we use first use \eqref{E:WTRANSDEF}, 
	\eqref{E:LROUGH},
	Lemmas~\ref{L:COMMUTATORSTOCOORDINATES} and \ref{L:CONTINUOUSEXTNESION},
	and the estimate
	$\frac{1}{\Lunit \timefunctionarg{\muxmulevelsetvalue}} \approx 1$ (see \eqref{E:LDERIVATIVEOFROUGHTIMEFUNCTIONISAPPROXIMATELYUNITY})
	to deduce that for $(\timefunction,u) \in [\timefunction_0,\timefunctionboot] \times [-\interestingu,\interestingu]$,
	we have
	$\| 
		\argLrough{\muxmulevelsetvalue} \Wtransarg{\muxmulevelsetvalue} \Wtransarg{\muxmulevelsetvalue} \upmu 
	\|_{L^{\infty}\left(\twoargroughtori{\timefunction,u}{\muxmulevelsetvalue}\right)} \leq C$.
	From this bound,
	\eqref{E:MOREPRECISELINFINITYTRANSPORTROUGHLFESTIMATE},
	and the data assumption
	\eqref{E:DATATASSUMPTIONMUTRANSVERSALCONVEXITY},
	we deduce that in the region $\twoargMrough{[\timefunction_0,\timefunctionboot],[-\interestingu,\interestingu]}{\muxmulevelsetvalue}$,
	we have 
	$
	\secondtransversalderivativemulowerbound
	-
	C \mupositive
	\leq
	\Wtransarg{\muxmulevelsetvalue} \Wtransarg{\muxmulevelsetvalue} \upmu
	\leq
	\frac{1}{\secondtransversalderivativemulowerbound}
	+
	C \mupositive
	$,
	which, for $\mupositive$ sufficiently small,
	implies \eqref{E:MUTRANSVERSALCONVEXITY}
	for $\Wtransarg{\muxmulevelsetvalue} \Wtransarg{\muxmulevelsetvalue} \upmu$
	(recall that, as we highlighted in Sect.\,\ref{SS:CONVENTIONSFORCONSTANTS}, the constants $C$
	can be chosen to be independent of $\mupositive$).
	The remaining estimates in \eqref{E:MUTRANSVERSALCONVEXITY}
	follow from a nearly identical argument,
	where we use the identity
	\eqref{E:ROUGHUDERIVATIVEINTERMSOFGEOMETRICVECTORFIELDSANDROUGHTIMEFUNCTION}
	when deriving the estimates involving $\roughgeop{u} \muX \upmu$.
	
	\medskip
	
	\noindent \textbf{Proofs of  
	\eqref{E:MUISLARGEINBORINGREGION}, 
	\eqref{E:REGIONWHEREMUXMUKAPPALEVELSETISNOTLOCATED},
	\eqref{E:BOUNDSONLMUINTERESTINGREGION},
	and \eqref{E:BOUNDSONGEOMETRICTDERIVATIVEMUINTERESTINGREGION}}:
	The arguments we used in the proof of \eqref{E:MUTRANSVERSALCONVEXITY} 
	also yield that $|\argLrough{\muxmulevelsetvalue} \Wtransarg{\muxmulevelsetvalue} \upmu| \leq C$.
	Hence, noting that $\Wtransarg{\muxmulevelsetvalue} \upmu = \muX \upmu + \muxmulevelsetvalue$,
	and using \eqref{E:MOREPRECISELINFINITYTRANSPORTROUGHLFESTIMATE} 
	and the data assumption \eqref{E:DATATASSUMPTIONREGIONWHEREMUXMUKAPPALEVELSETISNOTLOCATED},
	we find that
	in the spacetime region $\twoargMrough{[\timefunction_0,\timefunctionboot],[-\interestingu,\interestingu]}{\muxmulevelsetvalue}
			\backslash
			\mathcal{M}^{(\muxmulevelsetvalue)}_{[\timefunction_0,\timefunctionboot],[-\frac{\interestingu}{2},\frac{\interestingu}{2}]}$,
	we have
	$|\muX \upmu + \muxmulevelsetvalue| 
	\geq
	\frac{\secondtransversalderivativemulowerbound \interestingu}{4} - C \mupositive
	$
	which, for $\mupositive$ sufficiently small, implies \eqref{E:REGIONWHEREMUXMUKAPPALEVELSETISNOTLOCATED}.
	
	The estimates 
	\eqref{E:MUISLARGEINBORINGREGION} and \eqref{E:BOUNDSONLMUINTERESTINGREGION} 
	follow from similar arguments 
	based on the data assumptions 
	\eqref{E:DATAMUISLARGEINBORINGREGION} and \eqref{E:DATAASIZEBOUNDSONLMUINTERESTINGREGION}
	and the estimates 
	$|\argLrough{\muxmulevelsetvalue} \upmu| \leq C$ 
	and
	$|\argLrough{\muxmulevelsetvalue} \Lunit \upmu | \leq C \fundbootsmall$,
	which follow from \eqref{E:LDERIVATIVEOFROUGHTIMEFUNCTIONISAPPROXIMATELYUNITY}
	and the estimates of Prop.\,\ref{P:IMPROVEMENTOFAUXILIARYBOOTSTRAP}.
	We also conclude \eqref{E:BOUNDSONGEOMETRICTDERIVATIVEMUINTERESTINGREGION} 
	by combining these arguments with 
	the relation $\geop{t} \upmu = \Lunit \upmu + \mathcal{O}(\fundbootsmall)$,
	which follows from the identity $\geop{t} \upmu = \Lunit \upmu - \Lunit^A \geop{x^A} \upmu$,
	Lemma~\ref{L:COMMUTATORSTOCOORDINATES},
	Lemma~\ref{L:SCHEMATICSTRUCTUREOFVARIOUSTENSORSINTERMSOFCONTROLVARS},
	and the estimates of Prop.\,\ref{P:IMPROVEMENTOFAUXILIARYBOOTSTRAP}.

	\medskip
	
	\noindent \textbf{Proof of \eqref{E:ROUGHTIMEFUNCTIONLDERIVATIVEBOUNDS}}:
	This estimate follows from the proof of \eqref{E:LDERIVATIVEOFROUGHTIMEFUNCTIONISAPPROXIMATELYUNITY},
	except we now use the
	bound \eqref{E:BOUNDSONLMUINTERESTINGREGION}
	in place of the
	bootstrap assumption \eqref{E:BABOUNDSONLMUINTERESTINGREGION} used in the proof of
	\eqref{E:LDERIVATIVEOFROUGHTIMEFUNCTIONISAPPROXIMATELYUNITY}.
	
	\medskip
	\noindent \textbf{Proof of \eqref{E:SIZEBOUNDSONWIDETILDELMUINTERESTINGREGION} and \eqref{E:MUMUSTDECREASEININTERESTINGREGION}}:
	\eqref{E:SIZEBOUNDSONWIDETILDELMUINTERESTINGREGION} follows from definition \eqref{E:LROUGH}
	and the estimates \eqref{E:BOUNDSONLMUINTERESTINGREGION} and \eqref{E:ROUGHTIMEFUNCTIONLDERIVATIVEBOUNDS}.
	\eqref{E:MUMUSTDECREASEININTERESTINGREGION} then follows from \eqref{E:SIZEBOUNDSONWIDETILDELMUINTERESTINGREGION},
	which implies that $\upmu$ decreases along the integral curves of $\argLrough{\muxmulevelsetvalue}$ when $|u| \leq \interestingu$.
	
	\medskip
	
\noindent \textbf{Proof of \eqref{E:LOWERBOUNDONMAGNITUDEOFXRPLUS}}:
We use \eqref{E:MUTRANSPORT},
\eqref{E:IDENTITYFORMAINTERMDRIVINGTHESHOCK},
Lemma~\ref{L:SCHEMATICSTRUCTUREOFVARIOUSTENSORSINTERMSOFCONTROLVARS},
the estimates of Prop.\,\ref{P:IMPROVEMENTOFAUXILIARYBOOTSTRAP},
and
\eqref{E:BOUNDSONLMUINTERESTINGREGION}
to deduce that in $\twoargMrough{[\timefunction_0,\timefunctionboot],[-\interestingu,\interestingu]}{\muxmulevelsetvalue}$,
we have the estimate
$
\frac{1}{2}
\Speed^{-1}(\Speed^{-1} \Speed_{;\LogDensity} + 1)
|\muX \RRiemann| 
\geq 
\frac{7}{8} \mathring{\updelta}_* 
+ 
\mathcal{O}(\fundbootsmall) 
\geq 
\frac{3}{4} \mathring{\updelta}_*$,
and, also Taylor expanding
$\Speed^{-1}(\Speed^{-1} \Speed_{;\LogDensity} + 1)$
around the background solution $\wavearray = 0$
and using 
\eqref{E:NONDEGENCONDITION}
and
\eqref{E:BACKGROUNDSOUNDSPEEDISUNITY},
the estimate
$
\Speed^{-1}(\Speed^{-1} \Speed_{;\LogDensity} + 1)
|\muX \RRiemann|
=
\left\lbrace
	1 + \mathcal{O}(\mathring{\upalpha})
\right\rbrace
|\bar{\Speed}_{;\LogDensity} + 1|
|\muX \RRiemann|
$. 
Combining these two estimates, we conclude \eqref{E:LOWERBOUNDONMAGNITUDEOFXRPLUS}.

The fact that $\sqrt{\sum_{a=1}^3 (X^a)^2} = 1 + \mathcal{O}(\mathring{\upalpha})$
follows from
\eqref{E:XSMALL},
\eqref{E:SCHEMATICSTRUCTUREOFXSMALL},
and the estimates of Prop.\,\ref{P:IMPROVEMENTOFAUXILIARYBOOTSTRAP}.

	\medskip
	
	\noindent \textbf{Proof of \eqref{E:MUXMUKAPPALEVELSETLOCATION} and \eqref{E:IMPROVEDLEVELSETSTRUCTUREANDLOCATIONOFMIN}}:
	Since $|\muX \upmu + \muxmulevelsetvalue|  = 0$ along $\breve{\mathbb{X}}^{(-\muxmulevelsetvalue)}$,
	\eqref{E:MUXMUKAPPALEVELSETLOCATION} follows directly from \eqref{E:REGIONWHEREMUXMUKAPPALEVELSETISNOTLOCATED}.
	
	\eqref{E:IMPROVEDLEVELSETSTRUCTUREANDLOCATIONOFMIN} then follows from
	\eqref{E:MUXMUKAPPALEVELSETLOCATION} since
	$\twoargmumuxtorus{-\timefunction}{-\muxmulevelsetvalue} 
	\subset \datahypfortimefunctiontwoarg{-\muxmulevelsetvalue}{[\timefunction_0,\timefunctionboot]}$.
	
	\medskip
	
	\noindent \textbf{Proof of \eqref{E:MINVALUEOFMUONFOLIATION} and the fact
	that $\min_{\hypthreearg{\timefunction}{[- \rightu,\leftu]}{\muxmulevelsetvalue}} \upmu$
	occurs in
	$
		\twoargmumuxtorus{-\timefunction}{-\muxmulevelsetvalue}
	$}:
	Since the construction of the rough time function $\timefunctionarg{\muxmulevelsetvalue}$ is such that
	$\upmu|_{\datahypfortimefunctiontwoarg{-\muxmulevelsetvalue}{[\timefunction_0,\timefunctionboot]}} 
	= 
	- 
	\timefunctionarg{\muxmulevelsetvalue}|_{\datahypfortimefunctiontwoarg{-\muxmulevelsetvalue}{[\timefunction_0,\timefunctionboot]}}$,
	it follows from
	\eqref{E:MUISLARGEINBORINGREGION},
	\eqref{E:MUXMUKAPPALEVELSETLOCATION},
	and \eqref{E:MU1BIGGERTHANMU0}
	that for each fixed $\timefunction \in [\timefunction_0,\timefunctionboot] = [-\mupositive,-\upmuboot]$,
	$\min_{\hypthreearg{\timefunction}{[- \rightu,\leftu]}{\muxmulevelsetvalue}} \upmu$ is achieved only in
	the subset $\hypthreearg{\timefunction}{[-\frac{1}{2}\interestingu,\frac{1}{2}\interestingu]}{\muxmulevelsetvalue}$,
	which is interior to $\hypthreearg{\timefunction}{[- \rightu,\leftu]}{\muxmulevelsetvalue}$.
	Hence, since $\Wtransarg{\muxmulevelsetvalue}$ is tangent to 
	$\hypthreearg{\timefunction}{[- \rightu,\leftu]}{\muxmulevelsetvalue}$,
	it must be that $\Wtransarg{\muxmulevelsetvalue} \upmu = 0$ at the minima.
	Considering also the definition of $\datahypfortimefunctiontwoarg{-\muxmulevelsetvalue}{[\timefunction_0,\timefunctionboot]}$ 
	and the identity
	$\Wtransarg{\muxmulevelsetvalue} \upmu = \muX \upmu + \muxmulevelsetvalue$,
	we see that the minima of $\upmu$ in $\hypthreearg{\timefunction}{[- \rightu,\leftu]}{\muxmulevelsetvalue}$ must belong to 
	$\datahypfortimefunctiontwoarg{-\muxmulevelsetvalue}{[\timefunction_0,\timefunctionboot]} 
	\cap 
	\hypthreearg{\timefunction}{[-\frac{1}{2}\interestingu,\frac{1}{2}\interestingu]}{\muxmulevelsetvalue}$.
	Hence, we conclude that
	$\min_{\hypthreearg{\timefunction}{[- \rightu,\leftu]}{\muxmulevelsetvalue}} \upmu = - \timefunction$, which is 
	\eqref{E:MINVALUEOFMUONFOLIATION}.
	Moreover, also using
	Lemma~\ref{L:CHOVFROMROUGHCOORDINATESTOMUWEGIGHTEDXMUCOORDINATES},
	we conclude 
	that the torus
	$
	\twoargmumuxtorus{-\timefunction}{-\muxmulevelsetvalue}
	=
	\InverseCHOVroughtomumuxmu{\muxmulevelsetvalue}(\lbrace -\timefunction \rbrace \times \lbrace -\muxmulevelsetvalue \rbrace \times \mathbb{T}^2)
	$
	is exactly the set of points 
	within 
	$
	\hypthreearg{\timefunction}{[- \rightu,\leftu]}{\muxmulevelsetvalue}
	$
	where $\upmu$ achieves its minimum value of $-\timefunction$.

\medskip
	
\noindent \textbf{Proof of \eqref{E:WMANDRTRANSMUUBOUNDEDBYSQRTMU}}:	
First, we use
\eqref{E:WTRANSDEF},
\eqref{E:DEFROUGHTORITANGENTVECTORFIELD},
\eqref{E:RTRANS},
Lemma~\ref{L:COMMUTATORSTOCOORDINATES},
Lemma~\ref{L:ROUGHPARTIALDERIVATIVESINTERMSOFGEOMETRICPARTIALDERIVATIVESANDVICEVERSA},
Lemma~\ref{L:SCHEMATICSTRUCTUREOFVARIOUSTENSORSINTERMSOFCONTROLVARS},
and the estimates of Lemma~\ref{L:DIFFEOMORPHICEXTENSIONOFROUGHCOORDINATES},
Prop.\,\ref{P:IMPROVEMENTOFAUXILIARYBOOTSTRAP}, 
and Cor.\,\ref{C:IMPROVEAUX}
to deduce that
$|\Wtransarg{\muxmulevelsetvalue} \upmu| \leq C$
and
$\Rtransarg{\muxmulevelsetvalue} \upmu 
= 
\Wtrans \upmu + \upmu \mathcal{O}(\fundbootsmall)
= \Wtrans \upmu + \mathcal{O}(\fundbootsmall)$,
where we clarify that we will use the first equality in the next paragraph.
From these bounds 
and \eqref{E:MUISLARGEINBORINGREGION},
we see that the bounds stated in \eqref{E:WMANDRTRANSMUUBOUNDEDBYSQRTMU} 
hold in $\twoargMrough{[\timefunction_0,\timefunctionboot],[- \rightu,\leftu]}{\muxmulevelsetvalue}
			\backslash 
			\twoargMrough{[\timefunction_0,\timefunctionboot],[-\interestingu,\interestingu]}{\muxmulevelsetvalue}$.
			
It remains for us to prove the desired bounds in 
$\twoargMrough{[\timefunction_0,\timefunctionboot],[-\interestingu,\interestingu]}{\muxmulevelsetvalue}$. 
To this end, we assume that 
$\timefunction \in [\timefunction_0,\timefunctionboot]$ 
and
$q \in \hypthreearg{\timefunction}{[- \interestingu,\interestingu]}{\muxmulevelsetvalue}$.
By Lemma~\ref{L:FLOWMAPFORGENERATOROFROUGHTIMEFUNCTION},
there is a unique integral curve of $\Wtransarg{\muxmulevelsetvalue}$ 
that joins $q$ to a point $q_0 \in \twoargmumuxtorus{-\timefunction}{-\muxmulevelsetvalue}$ 
on the primal torus $\twoargmumuxtorus{-\timefunction}{-\muxmulevelsetvalue}$
(along which $\upmu \equiv - \timefunction$),
which, in view of the already proven result \eqref{E:IMPROVEDLEVELSETSTRUCTUREANDLOCATIONOFMIN}, 
is contained in $\hypthreearg{\timefunction}{[-\frac{1}{2}\interestingu,\frac{1}{2}\interestingu]}{\muxmulevelsetvalue}$.
Let $\iota = \iota(u')$ denote this integral curve, parameterized by the eikonal function,
and let $u_0$ and $u$ respectively denote the eikonal function values corresponding to $q_0$ and $q$.
In particular, $\iota(u) = q$, $\iota(u_0) = q_0$,
$\upmu \circ \iota(u_0) = - \timefunction$,
and by \eqref{E:WTRANSDEF}, $\Wtrans \upmu \circ \iota(u_0) = 0$. 
Using \eqref{E:MUTRANSVERSALCONVEXITY} and the mean value theorem, 
we see that $|\Wtrans \upmu \circ \iota(u)| \leq \frac{2}{\secondtransversalderivativemulowerbound} |u - u_0|$
and $\upmu \circ \iota(u) \geq - \timefunction + \frac{\secondtransversalderivativemulowerbound}{4} (u- u_0)^2 
\geq \frac{\secondtransversalderivativemulowerbound}{4} (u- u_0)^2$.
Combining the above results, we find that at $q$, we have 
$|\Wtrans \upmu| \leq \frac{2}{\secondtransversalderivativemulowerbound} 
\sqrt{\frac{4}{\secondtransversalderivativemulowerbound}} \sqrt{\upmu} 
\leq \frac{4}{\secondtransversalderivativemulowerbound^{3/2}} \sqrt{\upmu}$,
which yields \eqref{E:WMANDRTRANSMUUBOUNDEDBYSQRTMU} for the first term on the LHS.
From this estimate,
the bound
$\Rtransarg{\muxmulevelsetvalue} \upmu = \Wtrans \upmu + \upmu \mathcal{O}(\fundbootsmall)$ noted in the previous paragraph,
and the estimate $\upmu \leq C$ implied by \eqref{E:LINFINITYIMPROVEMENTAUXNULLDERIVATIVEMUANDTRANSVERSALDERIVATIVES},
we conclude the desired estimate \eqref{E:WMANDRTRANSMUUBOUNDEDBYSQRTMU}
for the second term on the LHS.

\medskip

\noindent \textbf{Proof of \eqref{E:WIDETILDELMUISALMOSTMINUSONEINSMALLNEIGHBORHOOD}}:	
Recalling that $\flowmapWtransargtwoarg{\muxmulevelsetvalue}{\Delta u}$ is the flow map of $\Wtransarg{\muxmulevelsetvalue}$,
for real numbers $\Delta u$ small and positive, we consider the set
$
\flowmapWtransargtwoarg{\muxmulevelsetvalue}{(-\Delta u,\Delta u)}
\left(\datahypfortimefunctiontwoarg{-\muxmulevelsetvalue}{[\timefunction_0,\timefunctionboot]} \right)
\eqdef
\bigcup_{\Delta u' \in (-\Delta u,\Delta u)}
\flowmapWtransargtwoarg{\muxmulevelsetvalue}{\Delta u'}
\left( \datahypfortimefunctiontwoarg{-\muxmulevelsetvalue}{[\timefunction_0,\timefunctionboot]} \right)
$,
which by Lemma~\ref{L:FLOWMAPFORGENERATOROFROUGHTIMEFUNCTION} 
and \eqref{E:MUXMUKAPPALEVELSETLOCATION}
is a neighborhood of $\datahypfortimefunctiontwoarg{-\muxmulevelsetvalue}{[\timefunction_0,\timefunctionboot]}$
in $\twoargMrough{[\timefunction_0,\timefunctionboot],[-\interestingu,\interestingu]}{\muxmulevelsetvalue}$.
Next, using \eqref{E:WTRANSDEF},
Lemma~\ref{L:COMMUTATORSTOCOORDINATES},
\eqref{E:LROUGH},
Lemma~\ref{L:SCHEMATICSTRUCTUREOFVARIOUSTENSORSINTERMSOFCONTROLVARS},
and the estimates of Lemma~\ref{L:DIFFEOMORPHICEXTENSIONOFROUGHCOORDINATES}
and
Prop.\,\ref{P:IMPROVEMENTOFAUXILIARYBOOTSTRAP}, 
we deduce that
$
|\Wtransarg{\muxmulevelsetvalue} \argLrough{\muxmulevelsetvalue} \upmu|
\leq C 
$.
Moreover, since Lemma~\ref{L:GRADIENTOFTIMEFUNCTIONAGREESWITHGRADIENTOFMUALONGMUXMUEQUALSMINUSKAPPHYPERSURFACE}
and definition~\eqref{E:LROUGH} imply that
$
\argLrough{\muxmulevelsetvalue} \upmu|_{\breve{\mathbb{X}}^{(-\muxmulevelsetvalue)}}
= -1
$,
we can use the mean value theorem to deduce that in 
$
\flowmapWtransargtwoarg{\muxmulevelsetvalue}{(-\Delta u,\Delta u)}
\left( \datahypfortimefunctiontwoarg{-\muxmulevelsetvalue}{[\timefunction_0,\timefunctionboot]} \right)
$,
the following estimates hold:
$
- 1 - C \Delta u
\leq
\argLrough{\muxmulevelsetvalue} \upmu \circ \iota(u) 
\leq -1 
+ C \Delta u$.
Choosing and fixing a value of $\Delta u$, which we denote by
$\Delta U$, to be sufficiently small such that $C \Delta U < .01$,
we arrive at \eqref{E:WIDETILDELMUISALMOSTMINUSONEINSMALLNEIGHBORHOOD}
with
$\smallneighborhoodofcreasetwoarg{[\timefunction_0,\timefunctionboot]}{\muxmulevelsetvalue}
\eqdef
\flowmapWtransargtwoarg{\muxmulevelsetvalue}{(-\Delta U,\Delta U)}
\left( \datahypfortimefunctiontwoarg{-\muxmulevelsetvalue}{[\timefunction_0,\timefunctionboot]} \right)
$.

\medskip

\noindent \textbf{Proof of \eqref{E:EASYREGIONLOWERBOUNDFORMU} and \eqref{E:MURATIOBOUNDEDINBORINGREGION}}:
We recall that
$
\datahypfortimefunctiontwoarg{-\muxmulevelsetvalue}{[\timefunction_0,\timefunctionboot]}
=
\bigcup_{\mulevelsetvalue \in [-\timefunctionboot,\mupositive]}
 \twoargmumuxtorus{\mulevelsetvalue}{-\muxmulevelsetvalue}
$
and that for $\timefunction \in [\timefunction_0,\timefunctionboot]$,
we have
$
\upmu|_{\twoargmumuxtorus{-\timefunction}{-\muxmulevelsetvalue}}
= -\timefunction
$.
Moreover, the arguments we used in the proof of \eqref{E:WMANDRTRANSMUUBOUNDEDBYSQRTMU}
imply that for $|\Delta U| \leq \frac{\interestingu}{2}$,
we have
$
	-\timefunction
	+
	\frac{\secondtransversalderivativemulowerbound}{4} (\Delta U)^2 
	\leq
	\min_{\flowmapWtransargtwoarg{\muxmulevelsetvalue}{(-\Delta U,\Delta U)}\left(\twoargmumuxtorus{-\timefunction}{-\muxmulevelsetvalue} \right)}
	\upmu
$.
Hence, in view of \eqref{E:MUISLARGEINBORINGREGION}
and definition~\eqref{E:MU2DEF},
we see that
\eqref{E:EASYREGIONLOWERBOUNDFORMU}
holds with $\Delta U$ defined to be the small constant fixed in the proof of
\eqref{E:WIDETILDELMUISALMOSTMINUSONEINSMALLNEIGHBORHOOD}.
\eqref{E:MURATIOBOUNDEDINBORINGREGION} then follows as a simple
consequence of \eqref{E:EASYREGIONLOWERBOUNDFORMU},
\eqref{E:FORMOFROUGHNULLGENERATORFLOWMAP},
and the bound 
$
\max_{\twoargMrough{[\timefunction_0,\timefunctionboot],[- \rightu,\leftu]}{\muxmulevelsetvalue}}
\upmu 
\lesssim 
1$
implied by the third item in
Lemma~\ref{L:CONTINUOUSEXTNESION}.

\end{proof}

\subsection{Homeomorphism and diffeomorphism properties of $\Upsilon$}
\label{SS:PROPERTIESOFCHOVFROMGEOMETRICTOCARTESIAN}
Our main goal in this section is to reveal the homeomorphism and diffeomorphism properties of
the change of variables map
$\Upsilon(t,u,x^2,x^3) = (t,x^1,x^2,x^3)$.
We start with the following monotonicity lemma, which plays an important role in
controlling $\Upsilon$.

\begin{lemma}[Monotonicity of $x^1$]
	\label{L:MONOTONICITYOFCARTESIANX1}
	The following identity holds, where $\phi = \phi(u)$ is the cut-off from Def.\,\ref{D:WTRANSANDCUTOFF}: 
	\begin{align} \label{E:ROUGHUDERIVATIVEOFCATERSIANX1}
	\begin{split}	
		\roughgeop{u} x^1
		& = 
			\upmu 
			\left\lbrace
				X^1
				+
				\frac{X^A X^A}{X^1}
				+ 
				\frac{X^A (\geop{x^A} \timefunctionarg{\muxmulevelsetvalue}) \Lunit^B X^B}{(\geop{t} \timefunctionarg{\muxmulevelsetvalue}) X^1} 
				+
				\frac{X^A \geop{x^A} \timefunctionarg{\muxmulevelsetvalue}}{\geop{t} \timefunctionarg{\muxmulevelsetvalue}}
				\Lunit^1
			\right\rbrace
				\\
		& \ \
			+ 
			\phi 
			\frac{\muxmulevelsetvalue}{\Lunit \upmu} 
			\left\lbrace
				\Lunit^1
				+
				\frac{\Lunit^A X^A}{X^1}
				+ 
				\frac{\Lunit^A (\geop{x^A} \timefunctionarg{\muxmulevelsetvalue}) \Lunit^B X^B}{(\geop{t} \timefunctionarg{\muxmulevelsetvalue}) X^1} 
				+
				\frac{\Lunit^A \geop{x^A} \timefunctionarg{\muxmulevelsetvalue}}{\geop{t} \timefunctionarg{\muxmulevelsetvalue}}
				\Lunit^1
			\right\rbrace.
	\end{split}
	\end{align}
	
	Moreover, the following estimate holds on 
	$\Mrough{[\timefunction_0,\timefunctionboot],[- \rightu,\leftu]}^{(\muxmulevelsetvalue)}$:
	\begin{align}  \label{E:X1MONOTONOCITYESTIMATE}
		\roughgeop{u} x^1
	& = - \upmu
				\left\lbrace
					1 + \mathcal{O}_{\mydiam}(\mathring{\upalpha})
				\right\rbrace
				+
				\phi 
				\frac{\muxmulevelsetvalue}{\Lunit \upmu}
				\left\lbrace
					1 + \mathcal{O}_{\mydiam}(\mathring{\upalpha})
				\right\rbrace.
	\end{align}
	
	Finally, 
	for every fixed
	$(\timefunction,x^2,x^3) \in [\timefunction_0,\timefunctionboot] \times \mathbb{T}^2$,
	the map $u \rightarrow x^1(\timefunction,u,x^2,x^3)$ is strictly decreasing 
	on $[- \rightu,\leftu]$.
	\end{lemma}

\begin{proof}
	\eqref{E:ROUGHUDERIVATIVEOFCATERSIANX1}
	follows from 
	\eqref{E:ROUGHUDERIVATIVEINTERMSOFGEOMETRICVECTORFIELDSANDROUGHTIMEFUNCTION}
	and
	\eqref{E:GEOMETRICVECTORFIELDSINTERMSOFCARTESIANONES}.
	
	\eqref{E:X1MONOTONOCITYESTIMATE} then follows from \eqref{E:ROUGHUDERIVATIVEOFCATERSIANX1},
	Lemma~\ref{L:SCHEMATICSTRUCTUREOFVARIOUSTENSORSINTERMSOFCONTROLVARS},
	the estimates of Lemma~\ref{L:LINFTYESTIMATESFORROUGHTIMEFUNCTIONANDDERIVATIVES}
	and Prop.\,\ref{P:IMPROVEMENTOFAUXILIARYBOOTSTRAP},
	and Cor.\,\ref{C:IMPROVEAUX},
	and \eqref{E:DATAEPSILONISSMALLERTHANBOOTSTRAPEPSILONSMALLERTHANSQUAREOFDATAALPHA},
	which in particular imply that
	$X^1 = -1 + \Xsmall^1 = - 1 + \mathcal{O}_{\mydiam}(\mathring{\upalpha})$,
	$\Lunit^1 = 1 + \Lsmall^1 = 1 + \mathcal{O}_{\mydiam}(\mathring{\upalpha})$,
	and 
	$X^A, \, \Lunit^A = \mathcal{O}(\varepsilon) = \mathcal{O}_{\mydiam}(\mathring{\upalpha})$.
	
	To prove the monotonicity of $x^1$, we note that
	\eqref{E:X1MONOTONOCITYESTIMATE} and Prop.\,\ref{P:SHARPCONTROLOFMUANDDERIVATIVES}
	imply that $\roughgeop{u} x^1 < 0$,
	except in the case $\timefunctionboot = \muxmulevelsetvalue = 0$,
	where $\roughgeop{u} x^1$ vanishes precisely along the torus
	$\twoargmumuxtorus{0}{0}$,
	which is contained
	in $\hypthreearg{0}{[-\frac{1}{2}\interestingu,\frac{1}{2}\interestingu]}{0}$.
	Moreover, \eqref{E:SIMPLEESTIMATEFORROUGHUDERIVATIVEOFMUXMU} 
	implies that
	$\roughgeop{u} \muX \upmu|_{\twoargmumuxtorus{0}{0}} > 0
	$.
	Thus, since \eqref{E:LASTSLICETORIAREGRAPHSABOVEFLATTORIINGEOMETRICCOORDINATES} implies that
	$\CHOVgeotorough{\muxmulevelsetvalue}(\twoargmumuxtorus{0}{0})$ 
	(i.e., the image of the crease in rough adapted coordinate space 
	$\mathbb{R}_{\timefunction} \times \mathbb{R}_u \times \mathbb{T}^2$)
	is a graph over $\mathbb{T}^2$,
	we see that every integral curve of
	$\roughgeop{u}$ in $\hypthreearg{0}{[- \rightu,\leftu]}{0}$
	intersects $\twoargmumuxtorus{0}{0}$ in precisely one point.
	In total, we have shown that for every fixed
	$(\timefunction,x^2,x^3) \in [\timefunction_0,\timefunctionboot] \times \mathbb{T}^2$,
	the map $u \rightarrow x^1(\timefunction,u,x^2,x^3)$ on the domain $[- \rightu,\leftu]$
	has a negative derivative, except at possibly a single point
	in $[-\frac{\interestingu}{2},\frac{\interestingu}{2}]$.
	From this fact, we conclude that the map is strictly decreasing
	as desired. 
	
\end{proof}

\begin{proposition}[Homeomorphism and diffeomorphism properties of $\Upsilon$ and
		 the embedded tori $\Upsilon\left(\twoargmumuxtorus{\mulevelsetvalue}{-\muxmulevelsetvalue} \right)$]
		\label{P:HOMEOMORPHICANDDIFFEOMORPHICEXTENSIONOFCARTESIANCOORDINATES}
		The change of variables map $\Upsilon(t,u,x^2,x^3) = (t,x^1,x^2,x^3)$
		is \textbf{injective} on the compact set
		$\Mrough{[\timefunction_0,\timefunctionboot],[- \rightu,\leftu]}^{(\muxmulevelsetvalue)}$
		and satisfies: 
		\begin{align} \label{E:C31BOUNDFORCHOVFROMGEOTOCARTESIAN}
		\| \Upsilon \|_{C_{\textnormal{geo}}^{3,1}(\Mrough{[\timefunction_0,\timefunctionboot],[- \rightu,\leftu]}^{(\muxmulevelsetvalue)})}
		\leq C.
		\end{align}
		In particular, $\Upsilon$ is a homeomorphism from $\Mrough{[\timefunction_0,\timefunctionboot],[- \rightu,\leftu]}^{(\muxmulevelsetvalue)}$
		onto its image.
		
		Moreover, with $d_{\textnormal{geo}} \Upsilon$ denoting the Jacobian matrix of $\Upsilon$,
		we have:
		\begin{align} \label{E:GEOTOCARTESIANJACOBIANDETERMINANTESTIMATE}
			\mbox{\upshape det} d_{\textnormal{geo}} \Upsilon
			& = 
				\upmu \frac{\Speed^2}{X^1}
				\approx
				- 
				\upmu.
		\end{align}
		
		In addition, if $\timefunctionboot < 0$, then $\Upsilon$ is a diffeomorphism from
		$\Mrough{[\timefunction_0,\timefunctionboot],[- \rightu,\leftu]}^{(\muxmulevelsetvalue)}$ onto its image.
		
		Finally, let $\twoargmumuxtorus{\mulevelsetvalue}{-\muxmulevelsetvalue}$ be the $\upmu$-adapted torus defined in
		\eqref{E:MUXMUTORI}, and let $\Upsilon\left(\twoargmumuxtorus{\mulevelsetvalue}{-\muxmulevelsetvalue} \right)$
		be the image in Cartesian coordinate space of
		$\twoargmumuxtorus{\mulevelsetvalue}{-\muxmulevelsetvalue}$ under $\Upsilon$.
		Then for $\mulevelsetvalue \in [\upmuboot,\mupositive] = [-\timefunctionboot,-\timefunction_0]$,
		$\Upsilon\left(\twoargmumuxtorus{\mulevelsetvalue}{-\muxmulevelsetvalue} \right)$ is $C^{1,1}$
		embedded submanifold of Cartesian
		space that is diffeomorphic to $\mathbb{T}^2$. More precisely,
		with 
		$\Cartesiantisafunctiononmumxtoriarg{\mulevelsetvalue}{-\muxmulevelsetvalue}(x^2,x^3)$
		and $\Eikonalisafunctiononmumuxtoriarg{\mulevelsetvalue}{-\muxmulevelsetvalue}$ denoting the functions on $\mathbb{T}^2$
		from \eqref{E:GRAPHDESCRIPTIONOFMUXMUTORUS}, the map $\embeddingofmuadapatedtorinCartesianspace{\muxmulevelsetvalue}$ defined by:
		\begin{align} \label{E:MUADAPATEDTORIINCARTESIANSPACE}
			\embeddingofmuadapatedtorinCartesianspace{\muxmulevelsetvalue}(x^2,x^3)
			\eqdef
			\Upsilon
			\circ
			\left(\Cartesiantisafunctiononmumxtoriarg{\mulevelsetvalue}{-\muxmulevelsetvalue}(x^2,x^3),
			\Eikonalisafunctiononmumuxtoriarg{\mulevelsetvalue}{-\muxmulevelsetvalue}(x^2,x^3),
			x^2,
			x^3 
			\right)
		\end{align}	
		is a $C^{1,1}$ embedding, i.e., a $C^{1,1}$
		diffeomorphism from $\mathbb{T}^2$ onto 
		$\Upsilon\left(\twoargmumuxtorus{\mulevelsetvalue}{-\muxmulevelsetvalue} \right)$.
\end{proposition}

\begin{proof}
	We already proved the bound \eqref{E:C31BOUNDFORCHOVFROMGEOTOCARTESIAN} in Lemma~\ref{L:CONTINUOUSEXTNESION}.
	
	Next, we use
	\eqref{E:GEOMETRICVECTORFIELDSINTERMSOFCARTESIANONES},
	Lemma~\ref{L:SCHEMATICSTRUCTUREOFVARIOUSTENSORSINTERMSOFCONTROLVARS},
	and
	Prop.\,\ref{P:IMPROVEMENTOFAUXILIARYBOOTSTRAP}
	to compute that
	$
	\frac{\partial \Upsilon(t,u,x^2,x^3)}{\partial (t,u,x^2,x^3)}
	= \begin{pmatrix}
			1 & 0 & 0 & 0  
				\\
			\Lunit^1 + * & \upmu \frac{\Speed^2}{X^1} & * & * 
				\\
			0 & 0 & 1 & 0 				
				\\
			0 & 0 & 0 & 1
		\end{pmatrix}
	$,
	where here and in the rest of the proof,
	``$*$'' denotes any quantity that is pointwise bounded in magnitude by $\mathcal{O}(\mathring{\upalpha})$.
	Thus, $\mbox{\upshape det} \frac{\partial \Upsilon(t,u,x^2,x^3)}{\partial (t,u,x^2,x^3)} = \upmu \frac{\Speed^2}{X^1}$,
	and therefore, using \eqref{E:XSMALLINTERMSOFLSMALLANDVELOCITY},
	Lemma~\ref{L:SCHEMATICSTRUCTUREOFVARIOUSTENSORSINTERMSOFCONTROLVARS},
	and
	Prop.\,\ref{P:IMPROVEMENTOFAUXILIARYBOOTSTRAP},
	we compute that
	$\mbox{\upshape det} \frac{\partial \Upsilon(t,u,x^2,x^3)}{\partial (t,u,x^2,x^3)} 
	=
	- \left\lbrace 1 + \mathcal{O}(\mathring{\upalpha}) \right\rbrace \upmu$,
	which yields \eqref{E:GEOTOCARTESIANJACOBIANDETERMINANTESTIMATE}.
	
	If $\timefunctionboot < 0$, then by Prop.\,\ref{P:SHARPCONTROLOFMUANDDERIVATIVES},
	$\upmu$ is uniformly positive on $\Mrough{[\timefunction_0,\timefunctionboot],[- \rightu,\leftu]}^{(\muxmulevelsetvalue)}$,
	and from \eqref{E:GEOTOCARTESIANJACOBIANDETERMINANTESTIMATE} and the inverse function
	theorem, we see that
	$\Upsilon$ is a local diffeomorphism on $\Mrough{[\timefunction_0,\timefunctionboot],[- \rightu,\leftu]}^{(\muxmulevelsetvalue)}$.
	Thus, to complete the proof, we need to show that $\Upsilon$ is injective on
	$\Mrough{[\timefunction_0,\timefunctionboot],[- \rightu,\leftu]}^{(\muxmulevelsetvalue)}$,
	even if $\timefunctionboot = 0$. We will achieve this by
	proving the injectivity of the map
	$(\timefunctionarg{\muxmulevelsetvalue},u,x^2,x^3) \rightarrow (\timefunctionarg{\muxmulevelsetvalue},x^1,x^2,x^3)$
	on the domain
	$[\timefunction_0,\timefunctionboot] \times [- \rightu,\leftu] \times \mathbb{T}^2$,
	and then the injectivity of the map 
	$(\timefunctionarg{\muxmulevelsetvalue},x^1,x^2,x^3) \rightarrow (t,x^1,x^2,x^3)$;
	since the composition of two injective functions is injective, this 
	would finish the proof.
	
	To proceed, we first note that the injectivity of the map 
	$(\timefunctionarg{\muxmulevelsetvalue},u,x^2,x^3) \rightarrow (\timefunctionarg{\muxmulevelsetvalue},x^1,x^2,x^3)$
	on
	$[\timefunction_0,\timefunctionboot] \times [- \rightu,\leftu] \times \mathbb{T}^2$
	follows from the monotonicity of the map $u \rightarrow x^1(\timefunction,u,x^2,x^3)$
	guaranteed by Lemma~\ref{L:MONOTONICITYOFCARTESIANX1}.
	For use below, we also note that by 
	\eqref{E:LSMALLDEF},
	\eqref{E:YSMALLINTERMSOFLSMALLANDVELOCITY},
	\eqref{E:ROUGHTIMEPARTIALDERIVATIVEINTERMSOFGEOMETRICTIMEPARTIALDERIVATIVE},
	Lemma~\ref{L:COMMUTATORSTOCOORDINATES},
	\eqref{E:PARTIALTIMEDERIVATIVEOFROUGHTIMEFUNCTIONISAPPROXIMATELYUNITY},
	and the estimates of Prop.\,\ref{P:IMPROVEMENTOFAUXILIARYBOOTSTRAP},
	we have:
	\begin{align} \label{E:ROUGHTIMEDERIVATIVEOFCARTESIANX1ISAPPROXIMATELYUNITY}
	\roughgeop{\timefunction} x^1
	\approx \geop{t} x^1 
	& = 
	\Lunit x^1 
	-
	\Lunit^A \geop{x^A} x^1
	= 1
	+ 
	\Lsmall
	-
	\Lunit^A \geop{x^A} x^1
	\approx 1.
	\end{align}
	
	Now for each fixed $(x^2,x^3) \in \mathbb{T}^2$,
	let ${^{(\muxmulevelsetvalue)}\mathcal{I}}_{x^2,x^3}$ denote the image of the set
	$[\timefunction_0,\timefunctionboot] \times [- \rightu,\leftu] \times \lbrace (x^2,x^3) \rbrace$
	under the map $(\timefunctionarg{\muxmulevelsetvalue},u,x^2,x^3) \rightarrow (\timefunctionarg{\muxmulevelsetvalue},x^1,x^2,x^3)$.
	The arguments given above, including the monotonicity guaranteed by \eqref{E:ROUGHTIMEDERIVATIVEOFCARTESIANX1ISAPPROXIMATELYUNITY}, 
	imply that for each fixed $(u,x^2,x^3) \in [- \rightu,\leftu] \times \mathbb{T}^2$,
	the map $\timefunction \rightarrow x^1(\timefunction,u,x^2,x^3)$
	is strictly increasing on $[\timefunction_0,\timefunctionboot]$, and that 
	for each fixed $(\timefunction,x^2,x^3) \in[\timefunction_0,\timefunctionboot] \times \mathbb{T}^2$,
	the map
	$u \rightarrow x^1(\timefunction,u,x^2,x^3)$ is strictly decreasing on $[- \rightu,\leftu]$.
	It follows that there exist scalar functions $\timefunction \rightarrow {^{(\muxmulevelsetvalue)}a}_{x^2,x^3}(\timefunction)$
	and
	$\timefunction \rightarrow {^{(\muxmulevelsetvalue)}b}_{x^2,x^3}(\timefunction)$
	on $[\timefunction_0,\timefunctionboot]$
	such that
	${^{(\muxmulevelsetvalue)}\mathcal{I}}_{x^2,x^3}
	=
	\lbrace (\timefunction,x^1,x^2,x^3) \ | \ 
	\timefunction \in [\timefunction_0,\timefunctionboot],
		\,
	{^{(\muxmulevelsetvalue)}a}_{x^2,x^3}(\timefunction) \leq x^1 \leq {^{(\muxmulevelsetvalue)}b}_{x^2,x^3}(\timefunction) \rbrace
	$,
	where ${^{(\muxmulevelsetvalue)}a}_{x^2,x^3}(\cdot)$ and ${^{(\muxmulevelsetvalue)}b}_{x^2,x^3}(\cdot)$ are $C^1$ functions of $\timefunction$
	such that ${^{(\muxmulevelsetvalue)}a}_{x^2,x^3}(\timefunction) < {^{(\muxmulevelsetvalue)}b}_{x^2,x^3}(\timefunction)$ and
	$
	\frac{d}{d \timefunction} {^{(\muxmulevelsetvalue)}a}_{x^2,x^3}, 
		\,
	\frac{d}{d \timefunction} {^{(\muxmulevelsetvalue)}b}_{x^2,x^3}
	\approx 1
	$.
	
	To complete the proof, it remains for us to show that the map
	$(\timefunctionarg{\muxmulevelsetvalue},x^1,x^2,x^3) \rightarrow (t,x^1,x^2,x^3)$
	is injective. Let ${^{(\muxmulevelsetvalue)}\mathcal{I}}_{x^2,x^3}$ be the set from the previous paragraph.
	It suffices to show that for each fixed $(x^2,x^3) \in \mathbb{T}^2$,
	any two distinct points in ${^{(\muxmulevelsetvalue)}\mathcal{I}}_{x^2,x^3}$ with the same $x^1$ coordinate
	must be mapped to distinct points under the map
	$(\timefunctionarg{\muxmulevelsetvalue},x^1,x^2,x^3) \rightarrow (t,x^1,x^2,x^3)$.
	The structure of ${^{(\muxmulevelsetvalue)}\mathcal{I}}_{x^2,x^3}$ revealed in the previous paragraph
	shows that for any two distinct points in ${^{(\muxmulevelsetvalue)}\mathcal{I}}_{x^2,x^3}$ with the same $x^1$ coordinate,
	the straight line segment joining them
	(along which $x^1,x^2,x^3$ are constant and $\newtimefunction$ varies)
	is contained in ${^{(\muxmulevelsetvalue)}\mathcal{I}}_{x^2,x^3}$
	(this can be thought of as the vertical convexity of ${^{(\muxmulevelsetvalue)}\mathcal{I}}_{x^2,x^3}$).
	Hence, to complete the proof, it suffices for us to show that the partial derivative of $t$
	with respect to $\timefunctionarg{\muxmulevelsetvalue}$ in the coordinate system
	$(\timefunctionarg{\muxmulevelsetvalue},x^1,x^2,x^3)$ 
	is positive, except possibly when $\timefunctionarg{\muxmulevelsetvalue} = 0$.
	To proceed, we note that the partial derivative of interest is equal to
	$\frac{1}{\partial_t \timefunctionarg{\muxmulevelsetvalue}}$,
	where $\partial_t$ is the Cartesian partial derivative.
	We now use 
	\eqref{E:LSMALLDEF},
	\eqref{E:XSMALL},
	\eqref{E:XSMALLINTERMSOFLSMALLANDVELOCITY},
	\eqref{E:WTRANSDEF},
	\eqref{E:IVPFORROUGHTTIMEFUNCTION},
	\eqref{E:CARTESIANPARTIALTTOCOMMUTATORS}, 
	Lemma~\ref{L:COMMUTATORSTOCOORDINATES},
	Lemma~\ref{L:SCHEMATICSTRUCTUREOFVARIOUSTENSORSINTERMSOFCONTROLVARS},
	and the estimates of 
	Lemma~\ref{L:DIFFEOMORPHICEXTENSIONOFROUGHCOORDINATES},
	Prop.\,\ref{P:IMPROVEMENTOFAUXILIARYBOOTSTRAP},
	and \eqref{E:BOUNDSONLMUINTERESTINGREGION}
	to compute that
	$\partial_t = \Lunit + (1 + *) X + * \Yvf{2} + * \Yvf{3}$
	and that
	$\frac{1}{\partial_t \timefunctionarg{\muxmulevelsetvalue}}
	\approx
	\frac{1}{1 + \frac{\muxmulevelsetvalue \phi}{\upmu}}
	=
	\frac{\upmu}{\upmu + \muxmulevelsetvalue \phi}
	$, where $\phi$ is the cut-off function from Def.\,\ref{D:WTRANSANDCUTOFF}.
	We now recall that by \eqref{E:MINVALUEOFMUONFOLIATION},
	$\upmu$ can vanish only when $\timefunctionarg{\muxmulevelsetvalue} = 0$.
	It
	follows that
	$
	\frac{1}{\partial_t \timefunctionarg{\muxmulevelsetvalue}} > 0
	$
	except possibly when $\timefunctionarg{\muxmulevelsetvalue} = 0$,
	which is the desired result.
	
	Finally, since the map $\embeddingofmuadapatedtorinCartesianspace{\muxmulevelsetvalue}$ from \eqref{E:MUADAPATEDTORIINCARTESIANSPACE}
	is the composition of the injective $C^{3,1}$ map $\Upsilon$ with the $C^{1,1}$ embedding
	of $\twoargmumuxtorus{\mulevelsetvalue}{-\muxmulevelsetvalue}$
	given by \eqref{E:GRAPHDESCRIPTIONOFMUXMUTORUS},
	it follows that $\embeddingofmuadapatedtorinCartesianspace{\muxmulevelsetvalue}$ is a $C^{1,1}$ injection from $\mathbb{T}^2$ into Cartesian coordinate space.
	Moreover, its differential $d_{(x^2,x^3)} \embeddingofmuadapatedtorinCartesianspace{\muxmulevelsetvalue}$ is a $4 \times 2$ matrix
	whose lower $2 \times 2$ block is the identity
	$
	\begin{pmatrix}
			1 & 0   
				\\
			0 & 1 
		\end{pmatrix}
	$.
	That is, $d_{(x^2,x^3)} \embeddingofmuadapatedtorinCartesianspace{\muxmulevelsetvalue}$ 
	is full rank, and therefore $\embeddingofmuadapatedtorinCartesianspace{\muxmulevelsetvalue}$ is a $C^{1,1}$ embedding.
	This concludes the proof of the proposition.
	
	\end{proof}

\subsection{Control of the size of $t$ and $x^1$}
\label{SS:CONTROLOFCARTESIANTANDX1}
In the next lemma, 
we derive improvements of the bootstrap assumptions of Sect.\,\ref{SSS:BOOTSTRAPASSUMPTIONSSIZEOFCARTESIANTANDX1}.

\begin{lemma}[Control of the size of $t$ and $x^1$ in $\twoargMrough{[\timefunction_0,\timefunctionboot],[- \rightu,\leftu]}{\muxmulevelsetvalue}$]
	\label{L:CONTROLOFCARTESIANTANDX1}
	The following estimates hold for $\timefunction \in [\timefunction_0,\timefunctionboot]$:
	\begin{subequations}
\begin{align} \label{E:SIZEOFCARTESIANT}
		\frac{1}{3 \mathring{\updelta}_*}
		& 
		\leq
		\min_{\hypthreearg{\timefunction}{[- \rightu,\leftu]}{\muxmulevelsetvalue}} t
		\leq
		\sup_{\hypthreearg{\timefunction}{[- \rightu,\leftu]}{\muxmulevelsetvalue}} t
		\leq 
		\frac{3}{\mathring{\updelta}_*},
			\\
		- \leftu
		+
		\frac{1}{3 \mathring{\updelta}_*}
		& 
		\leq
		\min_{\hypthreearg{\timefunction}{[- \rightu,\leftu]}{\muxmulevelsetvalue}} x^1
		\leq
		\sup_{\hypthreearg{\timefunction}{[- \rightu,\leftu]}{\muxmulevelsetvalue}} x^1
		\leq 
		\rightu
		+
		\frac{3}{\mathring{\updelta}_*}.
		\label{E:SIZEOFCARTESIANX1}
	\end{align}
	\end{subequations}
	
\end{lemma}

\begin{proof}
	We first prove \eqref{E:SIZEOFCARTESIANT}.
	From 
	\eqref{E:LULTMUXUMUXT},
	\eqref{E:LROUGH},
	and
	\eqref{E:BALDERIVATIVEOFROUGHTIMEFUNCTIONISAPPROXIMATELYUNITY},
	we have $\argLrough{\muxmulevelsetvalue} t = \frac{1}{\Lunit \timefunctionarg{\muxmulevelsetvalue}} \Lunit t 
	= 
	\frac{1}{\Lunit \timefunctionarg{\muxmulevelsetvalue}}
	\approx 1
	$.
	Hence, recalling that 
	$\argLrough{\muxmulevelsetvalue} u  = 0$
	and
	$\argLrough{\muxmulevelsetvalue} \timefunction = 1$, 
	we can integrate along the integral curves of
	$\argLrough{\muxmulevelsetvalue}$ starting from any point in
	$\hypthreearg{\timefunction_0}{[- \rightu,\leftu]}{\muxmulevelsetvalue}$
	and use the data assumption
	\eqref{E:DATAASSUMPTIONSIZEOFCARTESIANT}
	and the identity $|\timefunction_0| = \mupositive$
	to deduce that in $\twoargMrough{[\timefunction_0,\timefunctionboot],[- \rightu,\leftu]}{\muxmulevelsetvalue}$,
	we have 
	$
	\frac{1}{2 \mathring{\updelta}_*}
	+
	\frac{1}{C} \mupositive
	\leq
	t
	\leq
	\frac{2}{\mathring{\updelta}_*}
	+
	C \mupositive
	$,
	which, for $\mupositive$ sufficiently small,
	implies \eqref{E:SIZEOFCARTESIANT}.
	
	To prove \eqref{E:SIZEOFCARTESIANX1},
	we first use
	\eqref{E:LULTMUXUMUXT},
	\eqref{E:LROUGH},
	\eqref{E:BALDERIVATIVEOFROUGHTIMEFUNCTIONISAPPROXIMATELYUNITY},
	and \eqref{BA:AUXL1SMALL}
	to deduce that
	we have $\argLrough{\muxmulevelsetvalue} x^1 = \frac{1}{\Lunit \timefunctionarg{\muxmulevelsetvalue}} \Lunit x^1
	= 
	\frac{1}{\Lunit \timefunctionarg{\muxmulevelsetvalue}} (1 + \Lsmall^1)
	\approx 1
	$. We now argue 
	as in the proof of \eqref{E:SIZEOFCARTESIANT}
	using the data assumption \eqref{E:DATAASSUMPTIONSIZEOFCARTESIANX1},
	thereby arriving at \eqref{E:SIZEOFCARTESIANX1}.

\end{proof}


\section{Modified quantities for controlling the acoustic geometry} 
\label{S:CONSTRUCTIONOFMODIFIEDQUANTITIES}
In this section, we construct ``modified'' versions of the eikonal function quantity $\mytr_{\gtorus} \upchi$,
and we derive the transport equations that they satisfy.
There are two kinds of modified quantities: ``partially modified'' and ``fully modified.'' 
The partially modified quantities, when combined with integration by parts, will allow us to avoid 
uncontrollable error integrals in the wave equation energy identities.
The fully modified quantities,
when combined with elliptic estimates on the tori $\twoargroughtori{\timefunction,u}{\muxmulevelsetvalue}$, 
will allow us to control the top-order 
$\Yvf{A}$-derivatives of $\mytr_{\gtorus} \upchi$ without losing derivatives.

\subsection{Decompositions of $\Ricfour_{\Lunit \Lunit}$} 
\label{SS:RICLLDECOMPOSITIONS}
In the next lemma, we provide two key decompositions of $\Ricfour_{\Lunit \Lunit}$, where $\Ricfour$ is the Ricci curvature
of the acoustical metric $\gfour$. The decompositions are a crucial ingredient 
in our derivation of the transport equations satisfied by the modified quantities;
see the proofs of Props.\,\ref{P:TRANSPORTEQUATIONFORFULLYMODIFIEDQUANTITY} and \ref{P:TRANSPORTEQUATIONFORPARTIALMODIFIEDQUANTITY}.


\begin{lemma}[The key identities verified by $\Ricfour_{\Lunit \Lunit}$]
\label{L:RICLLDECOMPOSITIONS}
Assume that the entries of $\wavearray = (\RRiemann,\LRiemann, v^2,v^3,\Ent)$ solve the geometric wave equations 
\eqref{E:VELOCITYWAVEEQUATION}--\eqref{E:ENTROPYWAVEEQUATION}.
Then the following identity holds, where $\Ricfour$ is the Ricci curvature of $\gfour$:
\begin{align}  
\begin{split} \label{E:RICCICONTRACTEDLANDLFORFULLYMODIFIED}
\upmu \Ricfour_{\Lunit \Lunit} 
& = 
\Lunit 
\left\lbrace 
	- 
	\vec{G}_{\Lunit \Lunit} \diamond \muX \wavearray 
	- 
	\frac{1}{2} \upmu \mytr_{\gtorus} \angG \diamond \Lunit \wavearray 
	- 
	\frac{1}{2} \upmu \vec{G}_{\Lunit \Lunit} \diamond \Lunit \wavearray 
	+ 
	\upmu \angG_{\Lunit}^{\#} \diamond \cdot \angrmD \wavearray 
\right\rbrace 
	\\
& \ \ 
+ 
\mathfrak{A},
\end{split}
\end{align}
where $\mathfrak{A}$ has the following schematic structure:
\begin{align} \label{E:AINHOMRIC}
\mathfrak{A}  
& = 
\upmu 
\smoothfunction(\wavearray)
\cdot
(\VortVort,\DivGradEnt)
+			
\smoothfunction(\badcontrolvars,\comder \wavearray) 
\cdot
\tander \wavearray
+
\smoothfunction(\badcontrolvars,\vortrenormalized,\GradEnt,\comder \wavearray) 
\cdot
(\vortrenormalized,\GradEnt).
\end{align}

Moreover, \underline{without} the assumption that the geometric wave equations 
\eqref{E:VELOCITYWAVEEQUATION}--\eqref{E:ENTROPYWAVEEQUATION}
are satisfied, the following identity holds:
\begin{align} 
\begin{split} \label{E:RICCICONTRACTEDLANDLFORPARTIALLYMODIFIED}
\Ricfour_{\Lunit \Lunit} 
& 
= 
\frac{\Lunit \upmu}{\upmu} \mytr_{\gtorus} \upchi 
+ 
\Lunit 
\left\lbrace 
	-\frac{1}{2} \mytr_{\gtorus} \angG \diamond \Lunit \wavearray 
	+ 
	\angG_{\Lunit}^{\#} \diamond \cdot \angrmD \wavearray 
\right\rbrace 
	\\
& \ \
- 
\frac{1}{2} \vec{G}_{\Lunit \Lunit} \diamond \angLap \wavearray 
+ 
\mathfrak{B}, 
\end{split}
\end{align}
where $\mathfrak{B}$ has the following schematic structure:
\begin{align} \label{E:BINHOMRIC}
	\mathfrak{B} 
	& = 
	\smoothfunction(\controlvars)(\tander \wavearray) \cdot \tander \controlvars.
\end{align}

\end{lemma}

\begin{proof}[Sketch of a proof] 
In \cite[Lemma 6.1]{jLjS2018}, 
in the case of two spatial dimensions,
analogs of
\eqref{E:RICCICONTRACTEDLANDLFORFULLYMODIFIED} and \eqref{E:RICCICONTRACTEDLANDLFORPARTIALLYMODIFIED} 
were derived. 
Analogous identities were also derived in \cite[Corollary 11.4]{jS2016b}
in the case of quasilinear wave equations in three space dimensions.
The proof of \eqref{E:RICCICONTRACTEDLANDLFORPARTIALLYMODIFIED} 
given in \cite[Corollary 11.4]{jS2016b} is based 
on first writing $\Ricfour_{\Lunit \Lunit}$ relative to the Cartesian coordinates,
then writing all derivatives of $\wavearray$ in terms of derivatives with respect to
elements of $\lbrace \Lunit, X \rbrace$ and $\ell_{t,u}$-tangent differentiations,
and then finally expressing (with the help of Lemma~\ref{L:TRANSPORTMUANDLUNITI})
all of the principal terms (i.e., the terms that depend on the second derivatives of $\wavearray$)
as a perfect $\Lunit$ derivative up to lower-order terms, 
except for the term
$
- 
\frac{1}{2} \vec{G}_{\Lunit \Lunit} \diamond \angLap \wavearray 
$ 
on the last line of RHS~\eqref{E:RICCICONTRACTEDLANDLFORPARTIALLYMODIFIED}.
The detailed proof given in \cite[Corollary 11.4]{jS2016b} goes through nearly verbatim,
except we have used Lemma~\ref{L:SCHEMATICSTRUCTUREOFVARIOUSTENSORSINTERMSOFCONTROLVARS}
to simplify our schematic presentation of the term $\mathfrak{B}$.

Similarly, \eqref{E:RICCICONTRACTEDLANDLFORFULLYMODIFIED} can be proved using the same
arguments given in \cite[Corollary 11.4]{jS2016b}, but the new feature of
the present work is the structure of the terms on RHS~\eqref{E:RICCICONTRACTEDLANDLFORFULLYMODIFIED}.
To see how these terms arise, we explain how the proof of
\eqref{E:RICCICONTRACTEDLANDLFORFULLYMODIFIED} is connected to the identity \eqref{E:RICCICONTRACTEDLANDLFORPARTIALLYMODIFIED}.
To pass from \eqref{E:RICCICONTRACTEDLANDLFORPARTIALLYMODIFIED} to \eqref{E:RICCICONTRACTEDLANDLFORFULLYMODIFIED},
one uses the identity \eqref{E:BOXDECOMPLOUTSIDE}, 
the wave equations \eqref{E:VELOCITYWAVEEQUATION}--\eqref{E:ENTROPYWAVEEQUATION},
Lemma~\ref{L:TRANSPORTMUANDLUNITI}, and 
Lemma~\ref{L:SCHEMATICSTRUCTUREOFVARIOUSTENSORSINTERMSOFCONTROLVARS}
to express the product of $\upmu$ 
and the term 
$
-
\frac{1}{2} \vec{G}_{\Lunit \Lunit} \diamond \angLap \wavearray
$
on RHS~\eqref{E:RICCICONTRACTEDLANDLFORPARTIALLYMODIFIED} as follows:
\begin{align}
\begin{split} \label{E:DETAILEDINHOM}
-
\frac{1}{2} 
\upmu \vec{G}_{\Lunit \Lunit} \diamond \angLap \wavearray 
& = 
- 
\frac{1}{2} 
\Lunit 
\left\lbrace 
\vec{G}_{\Lunit \Lunit} 
\diamond 
(\upmu \Lunit\wavearray + 2 \muX \wavearray) 
\right\rbrace 
- 
\frac{1}{2} \mytr_{\gtorus} \upchi 
\vec{G}_{\Lunit \Lunit} 
\diamond 
\muX \wavearray 
	\\
& \ \
+ 
\smoothfunction(\controlvars) \cdot \text{Inhom} 
+  
\smoothfunction(\badcontrolvars,\comder \wavearray) 
\cdot
\tander \wavearray. 
\end{split}
\end{align}
In \eqref{E:DETAILEDINHOM}, ``Inhom'' denotes the inhomogeneous terms
$\upmu \times \mbox{RHS~\eqref{E:VELOCITYWAVEEQUATION}--\eqref{E:ENTROPYWAVEEQUATION}}$.
In a detailed proof (see \cite[Corollary 11.4]{jS2016b}), one finds that the term 
$-\frac{1}{2} \mytr_{\gtorus}\upchi \vec{G}_{\Lunit \Lunit} \diamond \muX \wavearray$
on RHS~\eqref{E:DETAILEDINHOM} is canceled
(and hence does not appear in \eqref{E:RICCICONTRACTEDLANDLFORFULLYMODIFIED}--\eqref{E:AINHOMRIC})
by part of the first product on RHS~\eqref{E:RICCICONTRACTEDLANDLFORPARTIALLYMODIFIED},
where one uses \eqref{E:MUTRANSPORT}
to substitute for the factor $\Lunit \upmu$ on RHS~\eqref{E:RICCICONTRACTEDLANDLFORPARTIALLYMODIFIED}.
Next, decomposing $\upmu \times \mbox{RHS~\eqref{E:VELOCITYWAVEEQUATION}--\eqref{E:ENTROPYWAVEEQUATION}}$
using \eqref{E:NULLFORMSTRUCTUREWAVEVARIABLES},
\eqref{E:STRUCTUREOFLINEARMODIFIEDFLUIDVARIABLETERMSONRHSWAVEEQUATIONS},
and
\eqref{E:ALLLINEARTERMSTRUCTURE},
we find that the term 
$\smoothfunction(\controlvars) \cdot \text{Inhom}$
on RHS~\eqref{E:DETAILEDINHOM} can be schematically expressed as follows:
$
\smoothfunction(\controlvars) \cdot \text{Inhom}
=
\upmu 
\smoothfunction(\wavearray)
\cdot
(\VortVort,\DivGradEnt)
+			
\smoothfunction(\badcontrolvars,\comder \wavearray) 
\cdot
\tander \wavearray
+
\smoothfunction(\badcontrolvars,\vortrenormalized,\GradEnt,\comder \wavearray) 
\cdot
(\vortrenormalized,\GradEnt)
$.
Placing these terms on RHS~\eqref{E:AINHOMRIC}, and also incorporating
the last term on RHS~\eqref{E:DETAILEDINHOM} into RHS~\eqref{E:AINHOMRIC},
we arrive at \eqref{E:RICCICONTRACTEDLANDLFORFULLYMODIFIED}--\eqref{E:AINHOMRIC}. 
\end{proof}

\subsection{Definition of the modified quantities} 
\label{SS:DEFSOFMODQUANTITIES}
We are now ready to define the modified quantities.
The definitions are motivated by the structure of the terms in Lemma~\ref{L:RICLLDECOMPOSITIONS};
this will become clear in the proofs of Props.\,\ref{P:TRANSPORTEQUATIONFORFULLYMODIFIEDQUANTITY}
and \ref{P:TRANSPORTEQUATIONFORPARTIALMODIFIEDQUANTITY}.

\begin{definition}[Modified versions of the $\nullhyparg{u}$-tangential derivatives of $\mytr_{\gtorus} \upchi$] 
\label{D:FULLYANDPARTIALLYMODIFIEDQUANTITIES} 
Let $N = \Ntop$, 
and let $\tander^N \in \mathfrak{P}^{(N)}$, where $\mathfrak{P}^{(N)}$ 
is the set of order $N$ $\nullhyparg{u}$-tangential commutator operators from
Sect.\,\ref{SS:STRINGSOFCOMMUTATIONVECTORFIELDS}.
We define the \emph{fully modified quantity} $\fullymodquant{\tander^N}$
as follows: 
\begin{subequations}
\begin{align}
\fullymodquant{\tander^N} 
& 
\eqdef 
\upmu \tander^N \mytr_{\gtorus} \upchi 
+ 
\tander^N \mathfrak{X},
	\label{E:FULLYMODIFIEDQUANTITY} 
		\\
\mathfrak{X} 
& 
\eqdef 
- 
\vec{G}_{\Lunit \Lunit} \diamond \muX \wavearray 
- 
\frac{1}{2} \upmu \mytr_{\gtorus} \angG \diamond \Lunit \wavearray 
- 
\frac{1}{2} \upmu \vec{G}_{\Lunit \Lunit} \diamond \Lunit \wavearray 
+ 
\upmu \angG_{\Lunit}^{\#} \diamond \cdot \angrmD \wavearray. 
\label{E:MODIFIEDQUANTITYINHOM}
\end{align}
\end{subequations}

Moreover, with $N = \Ntop - 1$ and $\tander^N \in \mathfrak{P}^{(N)}$,
we define the \emph{partially modified quantity} $\partialmodquant{\tander^N}$ 
as follows:
\begin{subequations}
\begin{align}
\partialmodquant{\tander^N} 
& 
\eqdef 
\tander^N \mytr_{\gtorus} \upchi 
+ 
\partialmodquantinhom{\tander^N}, 
\label{E:PARTIALMODIFIEDQUANTITY} 
	\\
\partialmodquantinhom{\tander^N} 
& \eqdef 
- 
\frac{1}{2} \mytr_{\gtorus} \angG \diamond \Lunit \tander^N \wavearray 
+ 
\angG_{\Lunit}^{\#} \diamond \cdot \angrmD \tander^N \wavearray. 
\label{E:PARTIALMODIFIEDQUANTITYINHOM}
\end{align}
\end{subequations}

Finally, we define the following
``$0^{\text{th}}$-order'' version of \eqref{E:PARTIALMODIFIEDQUANTITYINHOM}:
\begin{align} \label{E:PARTIALMODIFIEDQUANTITYINHOMZEROORDER}
\widetilde{\mathfrak{X}} 
& \eqdef 
- 
\frac{1}{2} \mytr_{\gtorus} \angG \diamond \Lunit \wavearray 
+ 
\angG_{\Lunit}^{\#} \diamond \cdot \angrmD \wavearray.
\end{align} 
\end{definition}

\subsection{Transport equations for the modified quantities}
\label{SS:TRANSPORTEQUATIONSFORMODQUANTITIES}
In this section, we derive transport equations for the fully modified quantities 
from Def.\,\ref{D:FULLYANDPARTIALLYMODIFIEDQUANTITIES}.
We start with the following lemma, which 
provides the transport equation satisfied by $\mytr_{\gtorus}\upchi$. 
This transport equation is an analog of
the well-known \emph{Raychaudhuri equation} in General Relativity \cite{aR1955}.


\begin{lemma}[Raychaudhuri-type transport equation for $\mytr_{\gtorus}\upchi$] 
\label{L:RAYCHAUDHURITRANSPORTCHI}
$\mytr_{\gtorus} \upchi$ obeys the following transport equation:
\begin{align}  \label{E:RAYCHAUDHURITRANSPORTCHI}
\upmu \Lunit \mytr_{\gtorus} \upchi 
& = 
(\Lunit \upmu) 
\mytr_{\gtorus} \upchi 
- 
\upmu \Ricfour_{\Lunit \Lunit} 
- 
\upmu |\upchi|_{\gtorus}^2. 
\end{align}
\end{lemma}

\begin{proof}
The same proof of \cite[(11.23)]{jS2016b} holds in the current setting.
\end{proof}



\begin{proposition}[Transport equation satisfied by $\fullymodquant{\tander^N}$] 
\label{P:TRANSPORTEQUATIONFORFULLYMODIFIEDQUANTITY}
Assume that $\wavearray = (\RRiemann,\LRiemann,v^2,v^3,\Ent)$ 
solve the geometric wave equations \eqref{E:VELOCITYWAVEEQUATION}--\eqref{E:ENTROPYWAVEEQUATION}. 
Let $N = \Ntop$, 
and let $\tander^N \in \mathfrak{P}^{(N)}$, where $\mathfrak{P}^{(N)}$
is the set of order $N$ $\nullhyparg{u}$-tangential commutator operators from
Sect.\,\ref{SS:STRINGSOFCOMMUTATIONVECTORFIELDS}.
Let $\fullymodquant{\tander^N}$ be the fully modified quantity defined in \eqref{E:FULLYMODIFIEDQUANTITY},
let $\mathfrak{X}$ be as defined in \eqref{E:MODIFIEDQUANTITYINHOM},
and let $\mathfrak{A}$ be the term on RHS~\eqref{E:RICCICONTRACTEDLANDLFORFULLYMODIFIED}.
Then $\fullymodquant{\tander^N}$ obeys the following transport equation,
where $\tander^N$ is the same differential operator every time it appears in \eqref{E:TRANSPORTEQUATIONFORFULLYMODIFIEDQUANTITY}:
\begin{align}
\begin{split} \label{E:TRANSPORTEQUATIONFORFULLYMODIFIEDQUANTITY}
\Lunit \fullymodquant{\tander^N} 
-  
\left(
2 \frac{\Lunit \upmu}{\upmu}
\right)
\fullymodquant{\tander^N} 
& = 
- 
\left(2 \frac{\Lunit \upmu}{\upmu}\right) \tander^N \mathfrak{X} 
+ 
\upmu [\Lunit,\tander^N] \mytr_{\gtorus} \upchi  
	\\
 & 
\ \ 
+ 
[\Lunit,\tander^N] \mathfrak{X} 
+ 
[\upmu,\tander^N] \Lunit \mytr_{\gtorus}\upchi 
+ 
[\tander^N,\Lunit \upmu] \mytr_{\gtorus}\upchi  
	\\ 
& \ \
- 
\tander^N \left(\upmu |\upchi|_{\gtorus}^2 \right)  
- 
\tander^N \mathfrak{A}.
\end{split}
\end{align}
\end{proposition}

\begin{proof}[Sketch of a proof]
	First, we use \eqref{E:RICCICONTRACTEDLANDLFORFULLYMODIFIED} to substitute for the product
	$\upmu \Ricfour_{\Lunit \Lunit}$ in \eqref{E:RAYCHAUDHURITRANSPORTCHI}.
	We then differentiate the resulting equation with $\tander^N$ and carry out tedious but straightforward
	commutations. Also taking into account definition \eqref{E:FULLYMODIFIEDQUANTITY},
	we conclude \eqref{E:TRANSPORTEQUATIONFORFULLYMODIFIEDQUANTITY}.
	We refer to the proof of \cite[Prop.~6.2]{jSgHjLwW2016} for more details.
\end{proof}

\begin{proposition}[Transport equation satisfied by $\partialmodquant{\tander^{N-1}}$] 
\label{P:TRANSPORTEQUATIONFORPARTIALMODIFIEDQUANTITY}
Let $N = \Ntop$, 
and let $\tander^{N-1} \in \mathfrak{P}^{(N-1)}$, where $\mathfrak{P}^{(N-1)}$ 
is the set of order $N-1$ $\nullhyparg{u}$-tangential commutator operators from
Sect.\,\ref{SS:STRINGSOFCOMMUTATIONVECTORFIELDS}.
Let $\partialmodquant{\tander^{N-1}}$ be the corresponding partially modified quantity defined in \eqref{E:PARTIALMODIFIEDQUANTITY},
let $\partialmodquantinhom{\tander^{N-1}}$ be the term defined in \eqref{E:PARTIALMODIFIEDQUANTITYINHOM},
let $\widetilde{\mathfrak{X}}$ be the term defined in \eqref{E:PARTIALMODIFIEDQUANTITYINHOMZEROORDER},
and let $\mathfrak{B}$ be the term on RHS~\eqref{E:RICCICONTRACTEDLANDLFORPARTIALLYMODIFIED}.
Then $\partialmodquant{\tander^{N-1}}$ obeys the following transport equation,
where $\tander^{N-1}$ is the same differential operator every time it appears in 
\eqref{E:TRANSPORTEQUATIONFORPARTIALMODIFIEDQUANTITY}:
\begin{align} \label{E:TRANSPORTEQUATIONFORPARTIALMODIFIEDQUANTITY}
\Lunit \partialmodquant{\tander^{N-1}} 
& =  
\frac{1}{2} \vec{G}_{\Lunit \Lunit}\diamond \angLap \tander^{N-1} \wavearray 
+ 
{^{(\tander^{N-1})}\mathfrak{B}},
\end{align}
where:
\begin{align}
\begin{split} \label{E:COMMUTEDPARTIALMODIFIEDQUANTITYINHOM}
{^{(\tander^{N-1})}\mathfrak{B}} 
& 
\eqdef 
- 
\tander^{N-1} \mathfrak{B} 
- 
\tander^{N-1} \left(|\upchi|_{\gtorus}^2 \right)
	\\
& \ \ 
+ 
\frac{1}{2} [\tander^{N-1},\vec{G}_{\Lunit \Lunit}] 
\diamond 
\angLap \wavearray 
+ 
\frac{1}{2} \vec{G}_{\Lunit \Lunit} 
\diamond 
[\tander^{N-1},\angLap]\wavearray  
	\\
& \ \ 
+ 
[\Lunit,\tander^{N-1}] \mytr_{\gtorus}\upchi 
+ 
[\Lunit,\tander^{N-1}] \widetilde{\mathfrak{X}} 
+ 
\Lunit 
\left\lbrace 
	\partialmodquantinhom{\tander^{N-1}} 
	- 
	\tander^{N-1} \widetilde{\mathfrak{X}}
\right\rbrace.
\end{split}
\end{align}

\end{proposition}

\begin{proof}
Substituting the identity \eqref{E:RICCICONTRACTEDLANDLFORPARTIALLYMODIFIED} into 
\eqref{E:RAYCHAUDHURITRANSPORTCHI}, dividing the resulting equation by $\upmu$,
and appealing to the definition
\eqref{E:PARTIALMODIFIEDQUANTITYINHOMZEROORDER}
of $\widetilde{\mathfrak{X}}$,
we deduce that:
\begin{align} \label{E:TRANSPORTEQUATIONFORPARTIALMODIFIEDQUANTITY0THORDER}
\Lunit 
\left(\mytr_{\gtorus} \upchi 
+ 
\widetilde{\mathfrak{X}}
\right) 
& 
= 
\frac{1}{2} \vec{G}_{\Lunit \Lunit} \diamond \angLap \wavearray 
- 
|\upchi|_{\gtorus}^2 
- 
\mathfrak{B}.
\end{align}
The transport equation \eqref{E:TRANSPORTEQUATIONFORPARTIALMODIFIEDQUANTITY} then follows 
from differentiating \eqref{E:TRANSPORTEQUATIONFORPARTIALMODIFIEDQUANTITY0THORDER} with $\tander^{N-1}$,
carrying out straightforward commutations,
and accounting for definitions 
\eqref{E:PARTIALMODIFIEDQUANTITY}--\eqref{E:PARTIALMODIFIEDQUANTITYINHOMZEROORDER}.

\end{proof}


\section{Basic ingredients in the $L^2$ analysis}
\label{S:BASICINGREDIENTSFORL2ANALYSIS}
In this section, we establish some preliminary ingredients that we will use when we derive energy estimates.
In Sect.\,\ref{S:ELLIPTICHYPERBOLICIDENTITIES}, we will derive one more crucial ingredient: 
elliptic-hyperbolic integral identities that we use to control the top-order derivatives of the specific vorticity and entropy gradient.
In Sect.\,\ref{SS:DIFFANDINTIDSONROUGHTORI}, we derive some differential and integral identities involving the rough tori 
$\twoargroughtori{\timefunction,u}{\muxmulevelsetvalue}$. 
In Sect.\,\ref{SS:SOBOLEVEMBEDDINGANDFTCESTIMATEONROUGHTORI}, we establish some basic Sobolev embedding estimates 
and fundamental theorem of calculus-type estimates on the rough tori.
In Sect.\,\ref{SS:IBPIDENTITIESNEEDFORTOPORDERENERGYESTIMATES}, 
we use the identities from Sect.\,\ref{SS:DIFFANDINTIDSONROUGHTORI} 
to prove various integration by parts identities 
that will play a key role in our energy estimates. 
Next, in Sect.\,\ref{SS:ENERGYIDENTITY},
we use the vectorfield multiplier method to construct 
the building block $L^2$-based energies and null-fluxes that we will use to control the wave variables $\wavearray$.
We also construct a companion set of building block energies and null-fluxes that we will use to control the transport variables
$\vortrenormalized$, $\GradEnt$, $\VortVort$, and $\DivGradEnt$.
Furthermore, in Prop.\,\ref{P:FUNDAMENTALENERGYNULLFLUXIDENTITIES}, 
we establish the fundamental energy-null-flux
integral identities that we exploit in our $L^2$ analysis.
In Sect.\,\ref{SS:FUNDAMENTALL2CONTROLLINGQUANTITIES}, we use the building block energies and null fluxes to 
define the quantities that we will use to control the solution in $L^2$.
Of particular interest are the spacetime integrals
 $\spacetimeintegralcontrolwave_N(\timefunction,u)$ 
and
$\spacetimeintegralcontrolwavepartial_N(\timefunction,u)$
defined in
\eqref{E:WAVESPACETIMEL2CONTROLLINGQUANTITY} and \eqref{E:PARTIALWAVESPACETIMEL2CONTROLLINGQUANTITY}
respectively. These spacetime integrals
appear on the left-hand side of our energy identity \eqref{E:FUNDAMENTALENERGYINTEGRALIDENTITCOVARIANTWAVES},
and they are fundamental for controlling error integrals that involve the quantities $\angrmD \tander^N \Psi$.
Finally, in Sect.\,\ref{SS:COERCIVENESSFUNDAMENTALL2CONTROLLINGQUANTITIES},
we exhibit the key coerciveness properties of our $L^2$-controlling quantities
with respect to the $L^2$ norms from Sect.\,\ref{SSS:GEOMETRICL2NORMS}.

\subsection{Differential and integral identities involving $\twoargroughtori{\timefunction,u}{\muxmulevelsetvalue}$} 
\label{SS:DIFFANDINTIDSONROUGHTORI}
The following lemma, though standard, plays an important role in our proof of the 
energy--null-flux identities for the wave variables
(see Prop.\,\ref{P:FUNDAMENTALENERGYNULLFLUXIDENTITIES}).
Moreover, the identity \eqref{E:RTRANSINTEGRALONROUGHTORI} plays a crucial role in 
our proof that one of the rough tori error integrals 
(specifically, the second term on LHS~\eqref{P:INTEGRALIDENTITYFORELLIPTICHYPERBOLICCURRENT})
in our elliptic-hyperbolic identities for the vorticity and entropy
has a \underline{favorable sign}.

\begin{lemma}[Differential and integral identities involving $\twoargroughtori{\timefunction,u}{\muxmulevelsetvalue}$] 
\label{L:DIFFERENTIALANDINTEGRALIDENTITIESONROUGHTORI}
Let $f$ be a scalar function,
and for 
$Z \in 
\lbrace \argLrough{\muxmulevelsetvalue}, \Rtransarg{\muxmulevelsetvalue} 
\rbrace$, let $\mytr_{\gtorusroughfirstfund} \deform{Z}$ be the $\gtorusroughfirstfund$-trace 
(see Def.\,\ref{D:TRACEOFROUGHTORITANGENT02TENSORS})
of the deformation tensor $\deform{Z}$ of $Z$.
Then the following integral identities hold:
\begin{subequations}
\begin{align}
	\roughgeop{\timefunction} 
	\left( 
		\int_{\twoargroughtori{\timefunction,u}{\muxmulevelsetvalue}} 
			f 
		\, \volroughtorus
	\right) 
	& = 
	\int_{\twoargroughtori{\timefunction,u}{\muxmulevelsetvalue}} 
	\left\lbrace 
		\argLrough{\muxmulevelsetvalue} f 
		+ 
		\frac{1}{2} f \mytr_{\gtorusroughfirstfund} \deform{\argLrough{\muxmulevelsetvalue}}
		\right\rbrace \, \volroughtorus,
	\label{E:IDENTITYROUGHTIMEDERIVATIVEOFROUGHTORUSINTEGRAL} 
	\\
	\roughgeop{u} \left(  \int_{\twoargroughtori{\timefunction,u}{\muxmulevelsetvalue}} f \, \volroughtorus\right) & = \int_{\twoargroughtori{\timefunction,u}{\muxmulevelsetvalue}} 
	\left\lbrace 
		\Rtransarg{\muxmulevelsetvalue} f 
		+ 
		\frac{1}{2} f \mytr_{\gtorusroughfirstfund} \deform{\Rtransarg{\muxmulevelsetvalue}}
	\right\rbrace \, \volroughtorus.
	\label{E:IDENTITYROUGHEIKONALDERIVATIVEOFROUGHTORUSINTEGRAL}
\end{align}
\end{subequations}
 
Moreover, for $u_1 \leq u_2$, we have:
\begin{align} \label{E:RTRANSINTEGRALONROUGHTORI}
\int_{\hypthreearg{\timefunction}{[u_1,u_2]}{\muxmulevelsetvalue}} 
	\left\lbrace
		\Rtransarg{\muxmulevelsetvalue} f 
		+ 
		\frac{1}{2} f \mytr_{\gtorusroughfirstfund} \deform{\Rtransarg{\muxmulevelsetvalue}}
	\right\rbrace
\, \volRoughHypersurface
& =
\int_{\twoargroughtori{\timefunction,u_2}{\muxmulevelsetvalue}} 
	f 
\, \volroughtorus
-
\int_{\twoargroughtori{\timefunction,u_1}{\muxmulevelsetvalue}} 
	f 
\, \volroughtorus.
\end{align}

\end{lemma}

\begin{proof}
We prove only 
\eqref{E:IDENTITYROUGHEIKONALDERIVATIVEOFROUGHTORUSINTEGRAL} and
\eqref{E:RTRANSINTEGRALONROUGHTORI} because 
\eqref{E:IDENTITYROUGHTIMEDERIVATIVEOFROUGHTORUSINTEGRAL} 
was proved in \cite[Lemma 6.3]{lAjS2020}\footnote{Note that since $\argLrough{\muxmulevelsetvalue}$ is $\gfour$-orthogonal to 
$\twoargroughtori{\timefunction,u}{\muxmulevelsetvalue}$ and $\argLrough{\muxmulevelsetvalue} \timefunction = 1$, the vectorfield $\argLrough{\muxmulevelsetvalue}$ agrees with the vectorfield denoted by ``$\breve{\underline{H}}$'' in \cite[Lemma 6.3]{lAjS2020}.} 
Let $\Phi_{(u')} = \Phi_{(u')}(\timefunction,u,x^2,x^3)$
be the flow map of $\Rtransarg{\muxmulevelsetvalue}$, normalized by $\Phi_{(0)}(\timefunction,u,x^2,x^3) = (\timefunction,u,x^2,x^3)$.
Since $\Rtransarg{\muxmulevelsetvalue} \timefunction = 0$ and $\Rtransarg{\muxmulevelsetvalue} u = 1$, 
it follows that if $\timefunction \in [\timefunction_0,\timefunctionboot)$,
$u \in (- \rightu,\leftu)$, and $|u'|$ is sufficiently small (depending on $u$), 
then $\Phi_{(u')}$ is a diffeomorphism from the rough torus 
$\twoargroughtori{\timefunction,u}{\muxmulevelsetvalue}$ onto the rough torus 
$\twoargroughtori{\timefunction,u + u'}{\muxmulevelsetvalue} 
\subset 
\twoargMrough{[\timefunction_0,\timefunctionboot),[- \rightu,\leftu]}{\muxmulevelsetvalue}$. 
Hence, in view of \eqref{E:AREAFORMROUGHTORUS}, 
\eqref{E:ROUGHTORUSINTEGRAL},
and the standard formula for change of variables in an integral, 
we deduce the following identity,
where $\Phi_{(u')}^*$ denotes pullback by $\Phi_{(u')}$
(in particular, $\Phi_{(u')}^* f = f \circ \Phi_{(u')}$),
and throughout this proof, determinants are taken relative to the coordinates $(x^2,x^3)$ on the rough tori:
\begin{align} \label{E:RTRANSDIFFINTINTEGRALONROUGHTORIINTERMEDIATE}
\int_{\twoargroughtori{\timefunction,u + u'}{\muxmulevelsetvalue}}
	f
\, \volroughtorus 
& = 
\int_{\twoargroughtori{\timefunction,u}{\muxmulevelsetvalue}} 
	[\Phi_{(u')}^* f] 
\volroughtorusarg{\Phi_{(u')}^* \gtorusroughfirstfund}
=
\int_{\mathbb{T}^2}
	f \circ \Phi_{(u')}(\timefunction,u,x^2,x^3)
	\sqrt{\mbox{\upshape det} [\Phi_{(u')}^* \gtorusroughfirstfund}]
\, \mathrm{d} x^2 \mathrm{d} x^3.
\end{align}
Next, using that $\Rtransarg{\muxmulevelsetvalue}$ is the infinitesimal generator of 
the flow map $\Phi_{(u')}$,
and using \eqref{E:LIEGTORUSROUGH},
we note the following differentiation identities:
$\frac{d}{du'}|_{u' = 0} \Phi_{(u')}^* f = \Rtransarg{\muxmulevelsetvalue} f$
and 
$
\frac{d}{du'}|_{u' = 0} \sqrt{\mbox{\upshape det} [\Phi_{(u')}^* \gtorusroughfirstfund}] 
=
\frac{1}{2} \sqrt{\mbox{\upshape det} \gtorusroughfirstfund}
\, 
\mytr_{\gtorusroughfirstfund} \deform{\Rtransarg{\muxmulevelsetvalue}}
$.
Using these identities,
we differentiate \eqref{E:RTRANSDIFFINTINTEGRALONROUGHTORIINTERMEDIATE} under the integral on the RHS
to obtain:
\begin{align} \label{E:RTRANSDIFFINTINTEGRALONROUGHTORISECONDINTERMEDIATE}
\roughgeop{u} 
\int_{\twoargroughtori{\timefunction,u}{\muxmulevelsetvalue}} 
	f
\, \volroughtorus 
= 
\left. \frac{d}{d u'}\right|_{u' = 0} 
\int_{\twoargroughtori{\timefunction,u + u'}{\muxmulevelsetvalue}}
	f
\, \volroughtorus
=  
\int_{\twoargroughtori{\timefunction,u}{\muxmulevelsetvalue}} 
	\left\lbrace 
		\Rtransarg{\muxmulevelsetvalue} f 
		+ 
		\frac{1}{2} f \mytr_{\gtorusroughfirstfund} \deform{\Rtransarg{\muxmulevelsetvalue}}
	\right\rbrace 
\, \volroughtorus.
\end{align}
We have therefore proved \eqref{E:IDENTITYROUGHEIKONALDERIVATIVEOFROUGHTORUSINTEGRAL}.

\eqref{E:RTRANSINTEGRALONROUGHTORI} then follows from integrating \eqref{E:IDENTITYROUGHEIKONALDERIVATIVEOFROUGHTORUSINTEGRAL}
with respect to $u$ and using the fundamental theorem of calculus,
\eqref{E:VOLUMEFORMROUGHHYPERSURFACE},
and \eqref{E:ROUGHHYPERSURFACEINTEGRAL}.
\end{proof}

\subsection{Sobolev embedding and fundamental theorem of calculus-type estimates on the rough tori}
\label{SS:SOBOLEVEMBEDDINGANDFTCESTIMATEONROUGHTORI}
In this section, 
on the rough tori $\twoargroughtori{\timefunction,u}{\muxmulevelsetvalue}$,
we derive $L^{\infty}$ Sobolev embedding estimates as well as simple $L^2$ estimates
that rely on the fundamental theorem of calculus.

\begin{lemma}[Sobolev embedding and fundamental theorem of calculus-type 
estimates on $\twoargroughtori{\timefunction,u}{\muxmulevelsetvalue}$] 
\label{L:SOBOLEVEMBEDDINGANDFTCL2ESTIMATESONROUGHTORI}
Let $f$ be a scalar function on $\twoargMrough{[\timefunction_0,\timefunctionboot),[- \rightu,\leftu]}{\muxmulevelsetvalue}$. 
Then the following estimates hold for $(\timefunction,u) \in [\timefunction_0,\timefunctionboot)\times [- \rightu,\leftu]$:
\begin{align}
	\| \tander^N f \|_{L^2\left(\twoargroughtori{\timefunction,u}{\muxmulevelsetvalue}\right)}^2 
	& 
	\lesssim 
	\| \tander^N f  \|_{L^2\left(\twoargroughtori{\timefunction_0,u}{\muxmulevelsetvalue}\right)}^2 
	+ 
	\int_{\nullhypthreearg{\muxmulevelsetvalue}{u}{[\timefunction_0,\timefunction]}} 
		\frac{1}{\Lunit \timefunctionarg{\muxmulevelsetvalue}} |\Lunit \tander^N f|^2 
	\, \volPuRoughCoordinates,
		\label{E:ROUGHTORUSINNULLHYPERSURFACEL2FUNDAMENTALTHEOREMOFCALCULUSESTIMATE} 
\end{align}


\begin{subequations}
\begin{align}
	\| f \|_{L^{\infty}\left(\twoargroughtori{\timefunction,u}{\muxmulevelsetvalue}\right)} 
	& \lesssim 
	\| \tander^{\leq 2} f \|_{L^2\left(\twoargroughtori{\timefunction,u}{\muxmulevelsetvalue}\right)},
		\label{E:H2LINFINITYSOBOLEVEMBEDDINGROUGHTORUS} 
			\\
	\| f\|_{L^{\infty}\left(\twoargroughtori{\timefunction,u}{\muxmulevelsetvalue}\right)}^2 
	& \lesssim 
	\| 
		\tander^{\leq 2} f  
	\|_{L^2\left(\twoargroughtori{\timefunction_0,u}{\muxmulevelsetvalue}\right)}^2 
		+ 
		\int_{\nullhypthreearg{\muxmulevelsetvalue}{u}{[\timefunction_0,\timefunction]}} \frac{1}{\Lunit \timefunctionarg{\muxmulevelsetvalue}} |\Lunit \tander^{\leq 2} f|^2 \, 		\volPuRoughCoordinates.  
	\label{E:H2LINFINITYFUNDAMENTALTHEOREMOFCALCULUSPLUSSOBOLEVEMBEDDINGONROUGHTORUS}
\end{align}
\end{subequations}
\end{lemma}

\begin{proof}
Using \eqref{E:IDENTITYROUGHTIMEDERIVATIVEOFROUGHTORUSINTEGRAL} with $f^2$ in place of $f$, 
\eqref{E:LROUGH},
the pointwise estimate \eqref{E:POINTWISEBOUNDFORROUGHTOROIDALTRACEOFDEFORMATIONTENSOROFROUGHNULLVECTORFIELD},
and Young's inequality, we deduce
$
\left|
\roughgeop{\timefunction} 
\| \tander^N f \|_{L^2(\twoargroughtori{\timefunction,u}{\muxmulevelsetvalue})}^2 
\right|
\leq 
\left\|
	\frac{1}{\sqrt{\Lunit \timefunctionarg{\muxmulevelsetvalue}}} \Lunit \tander^N f 
\right\|_{L^2(\twoargroughtori{\timefunction,u}{\muxmulevelsetvalue})}^2 
+ 
C \| \tander^N f\|_{L^2(\twoargroughtori{\timefunction,u}{\muxmulevelsetvalue})}^2
$. 
Integrating this inequality with respect to $\timefunction$, 
applying Gr\"{o}nwall's inequality, 
and also using the identity 
$\int_{\timefunction' = \timefunction_0}^{\timefunction} 
	\left\| 
		\frac{1}{\sqrt{\Lunit \timefunctionarg{\muxmulevelsetvalue}}} \tander^Nf
	\right\|_{L^2(\twoargroughtori{\timefunction',u}{\muxmulevelsetvalue})}^2 
\, \mathrm{d}\timefunction' 
= \left\| 
	\frac{1}{\sqrt{\Lunit \timefunctionarg{\muxmulevelsetvalue}}} \Lunit \tander^N f 
	\right \|_{L^2(\nullhypthreearg{\muxmulevelsetvalue}{u}{[\timefunction_0,\timefunction]})}^2
$, 
we conclude the inequality \eqref{E:ROUGHTORUSINNULLHYPERSURFACEL2FUNDAMENTALTHEOREMOFCALCULUSESTIMATE}.


We now prove 
\eqref{E:H2LINFINITYSOBOLEVEMBEDDINGROUGHTORUS}--\eqref{E:H2LINFINITYFUNDAMENTALTHEOREMOFCALCULUSPLUSSOBOLEVEMBEDDINGONROUGHTORUS}. 
We begin with the following estimate, which holds at fixed $(\timefunction,u)$ by virtue of
the standard Sobolev embedding result $H^2(\T^2) \hookrightarrow L^{\infty}(\T^2)$:
\begin{align} \label{E:STANDARDH2LINFINITYSOBOLEVEMBEDDINGONT2}
	\| f \|_{L^{\infty}\left(\twoargroughtori{\timefunction,u}{\muxmulevelsetvalue}\right)} 
	& \lesssim 
	\sum_{I + J \leq 2} 
		\left\lbrace 
			\int_{\T^2} 
				\left| 
					\left(\roughgeop{x^2}\right)^I \left(\roughgeop{x^3}\right)^J 
					f(\timefunction,u,x^2,x^3)
				\right|^2 
			\, 
			\mathrm{d} x^2 \mathrm{d} x^3
		\right\rbrace^{1/2}.
\end{align}
Using \eqref{E:ROUGHANGULARPARTIALDERIVATIVESINTERMSOFGOODGEOMETRICPARTIALDERIVATIVES},
Lemma~\ref{L:LINFTYESTIMATESFORROUGHTIMEFUNCTIONANDDERIVATIVES} to estimate derivatives of 
$\timefunction$ in \eqref{E:ROUGHANGULARPARTIALDERIVATIVESINTERMSOFGOODGEOMETRICPARTIALDERIVATIVES}, 
the identities \eqref{E:GEOP2TOCOMMUTATORS}--\eqref{E:GEOP3TOCOMMUTATORS},
Lemma~\ref{L:SCHEMATICSTRUCTUREOFVARIOUSTENSORSINTERMSOFCONTROLVARS},
the bootstrap assumptions,
and the area form element comparison estimate 
\eqref{E:VOLFORMESTIMATEROUGHTORI}, 
we deduce, in view of definitions~\eqref{E:AREAFORMROUGHTORUS} and \eqref{E:ROUGHTORUSINTEGRAL},
that $\mbox{RHS~\eqref{E:STANDARDH2LINFINITYSOBOLEVEMBEDDINGONT2}}
\lesssim 
\| \tander^{\leq 2} f \|_{L^2\left(\twoargroughtori{\timefunction,u}{\muxmulevelsetvalue}\right)}
$.
We have therefore proved \eqref{E:H2LINFINITYSOBOLEVEMBEDDINGROUGHTORUS}. 
The estimate 
\eqref{E:H2LINFINITYFUNDAMENTALTHEOREMOFCALCULUSPLUSSOBOLEVEMBEDDINGONROUGHTORUS} then
then follows from 
\eqref{E:H2LINFINITYSOBOLEVEMBEDDINGROUGHTORUS}
and
\eqref{E:ROUGHTORUSINNULLHYPERSURFACEL2FUNDAMENTALTHEOREMOFCALCULUSESTIMATE} with $N = 0,1,2$.  

\end{proof}

\subsection{Integration by parts identities} 
\label{SS:IBPIDENTITIESNEEDFORTOPORDERENERGYESTIMATES}
In this section, we establish several integration by parts identities that
we will exploit in our top-order energy estimates. 
We start with the following lemma, 
which provides a useful expression for the covariant divergence of a spacetime vectorfield.

\begin{lemma}[Covariant divergence identity for spacetime vectorfields] 
\label{L:COVARIANTDIVERGENCEOFSPACETIMEVECTORFIELDINTERMSOFRESCALEDFRAME}
Let $\mathscr{J}$ be a spacetime vectorfield.
Consider the decomposition $\upmu \mathscr{J} 
= 
- \upmu \mathscr{J}_{\Lunit} \Lunit 
- 
\mathscr{J}_{\muX} \Lunit 
- 
\mathscr{J}_{\Lunit} \muX 
+ 
\upmu \smoothtorusproject \mathscr{J}$ of $\upmu \mathscr{J}$ 
afforded by Lemma~\ref{L:BASICPROPERTIESOFVECTORFIELDS} ,
where $\mathscr{J}_{\Lunit} = \gfour(\mathscr{J},\Lunit)$, $\mathscr{J}_{\muX} = \gfour(\mathscr{J},\muX)$, 
and $\smoothtorusproject \mathscr{J}$ is the $\ell_{t,u}$-projection of $\mathscr{J}$ (see Def.\,\ref{D:PROJECTIONTENSORFIELDSANDTANGENCYTOHYPERSURFACES}).
Then the following identity holds:
\begin{align} \label{E:COVARIANTDIVERGENCEOFSPACETIMEVECTORFIELDINTERMSOFRESCALEDFRAME}
\upmu \Dfour_\alpha \mathscr{J}^{\alpha}= - L(\upmu \mathscr{J}_{\Lunit}) - L(\mathscr{J}_{\muX}) - \muX (\mathscr{J}_{\Lunit}) + \angdiv(\upmu \smoothtorusproject \mathscr{J}) - \upmu \mytr_{\gtorus} \angk \mathscr{J}_{\Lunit} - \mytr_{\gtorus} \upchi \mathscr{J}_{\muX}.
\end{align}
\end{lemma}

\begin{proof}
The same proof of \cite[Lemma 4.3]{jSgHjLwW2016} holds with minor modifications to account for the third space dimension.
\end{proof}

The following lemma provides some preliminary integration by parts identities.

\begin{lemma}[Preliminary integration by parts identities]
\label{L:PRELIMINARYINTEGRATIONBYPARTSIDENTITIES}
Let $(\timefunction,u) \in [\timefunction_0,\timefunctionboot)\times [- \rightu,\leftu]$,
and let $\upupsilon$ and $\upzeta$ be scalar functions on $\twoargMrough{[\timefunction_0,\timefunction),[- \rightu,u]}{\muxmulevelsetvalue}$. 
Then the following integration by parts identities hold,
where $A=2,3$ in \eqref{E:PRELIMINARYIBPINY}:
\begin{subequations}
\begin{align}
\begin{split} \label{E:PRELIMINARYIBPINLROUGH} 
	\int_{\twoargMrough{[\timefunction_0,\timefunction),[- \rightu,u]}{\muxmulevelsetvalue}} 
		(\argLrough{\muxmulevelsetvalue} \upupsilon) 
		\upzeta 
	\, \volMRoughCoordinates 
	& = 
	- 
	\int_{\twoargMrough{[\timefunction_0,\timefunction),[- \rightu,u]}{\muxmulevelsetvalue}} 
		\upupsilon 
		(\argLrough{\muxmulevelsetvalue} \upzeta) 
	\, \volMRoughCoordinates
		\\
& \ \
	- 
	\frac{1}{2} 
	\int_{\twoargMrough{[\timefunction_0,\timefunction),[- \rightu,u]}{\muxmulevelsetvalue}} 
		\mytr_{\gtorusroughfirstfund} \deform{\argLrough{\muxmulevelsetvalue}}
		\upupsilon \upzeta 
	\, \volMRoughCoordinates 
			\\
	& \ \ 
	+ 
	\int_{\hypthreearg{\timefunction}{[- \rightu,u]}{\muxmulevelsetvalue}} 
		\upupsilon \upzeta 
		\, \volRoughHypersurface 
	- 
	\int_{\hypthreearg{\timefunction_0}{u}{\muxmulevelsetvalue}} 
		\upupsilon \upzeta 
	\, \volRoughHypersurface, 
	\end{split}
		\\
\begin{split}	 \label{E:PRELIMINARYIBPINY}
	\int_{\twoargMrough{[\timefunction_0,\timefunction),[- \rightu,u]}{\muxmulevelsetvalue}} 
		\frac{1}{\Lunit \timefunctionarg{\muxmulevelsetvalue}} 
		(\Yvf{A} \upupsilon)
		\upzeta 
	\, \volMRoughCoordinates
	& = 
	-  
	\int_{\twoargMrough{[\timefunction_0,\timefunction),[- \rightu,u]}{\muxmulevelsetvalue}} 
		\frac{1}{\Lunit \timefunctionarg{\muxmulevelsetvalue}} 
		\upupsilon ( \Yvf{A} \upzeta ) 
	\, \volMRoughCoordinates 
		\\
	& \ \
	- 
	\frac{1}{2}  
	\int_{\twoargMrough{[\timefunction_0,\timefunction),[- \rightu,u]}{\muxmulevelsetvalue}} 
		\frac{1}{\Lunit \timefunctionarg{\muxmulevelsetvalue}} 
		\mytr_{\gtorus} \, \angdeform{\Yvf{A}}
		\upupsilon \upzeta 
	\, \volMRoughCoordinates 
		\\
	& \ \ 
	+ 
	\int_{\hypthreearg{\timefunction}{[- \rightu,u]}{\muxmulevelsetvalue}} 
		\frac{1}{\Lunit \timefunctionarg{\muxmulevelsetvalue}} 
		(\Yvf{A} \timefunction)  
		\upupsilon \upzeta 
	\, \volRoughHypersurface 
	- 
	\int_{\hypthreearg{\timefunction_0}{u}{\muxmulevelsetvalue}} 
		\frac{1}{\Lunit \timefunctionarg{\muxmulevelsetvalue}} 
		(\Yvf{A} \timefunction_0)  
		\upupsilon \upzeta 
	\, \volRoughHypersurface.  
\end{split}
\end{align}
\end{subequations}

%
\end{lemma}

\begin{proof}
The identity \eqref{E:PRELIMINARYIBPINLROUGH} follows from integrating \eqref{E:IDENTITYROUGHTIMEDERIVATIVEOFROUGHTORUSINTEGRAL} 
first with respect to $u$ and then with respect to $\timefunction$ with $f \eqdef \upupsilon \upzeta$.

We now prove \eqref{E:PRELIMINARYIBPINY}. 
We first expand $\Yvf{A}$ in terms of the rough adapted coordinate vectorfields as follows:
\begin{align} \label{E:YCOMMMUTATOREXPANDEDINROUGHADAPATEDCOORDIANTES}
	\Yvf{A} = (\Yvf{A} \timefunction) \roughgeop{\timefunction} + \Yvf{A}^B \roughgeop{x^B}.
\end{align}
Next, using \eqref{E:YCOMMMUTATOREXPANDEDINROUGHADAPATEDCOORDIANTES} 
and the fact that relative to the rough adapted coordinates, we have
$\sqrt{|\mbox{\upshape det} \gfour|} 
= 
\frac{\upmu}{\Lunit \timefunctionarg{\muxmulevelsetvalue}} \sqrt{\mbox{\upshape det} \gtorusroughfirstfund}$ 
(see \eqref{E:DETACOUSTICALMETRICINROUGHADAPTEDCOORDS}), 
we can expand the product of the covariant divergence of the vectorfield 
$\frac{1}{\upmu} \upupsilon \upzeta \Yvf{A}$ and $\sqrt{|\mbox{\upshape det} \gfour|}$ as follows
(again, relative to the rough adapted coordinates):
\begin{align} \label{E:IBPINYINTERMEDIATESTEP} 
\sqrt{|\mbox{\upshape det} \gfour|} 
\Dfour_\alpha \left( \frac{1}{\upmu} \upupsilon \upzeta \Yvf{A}\right)^{\alpha}
& = 
\roughgeop{\timefunction} 
\left( 
	\frac{1}{\Lunit \timefunctionarg{\muxmulevelsetvalue}} 
	(\Yvf{A} \timefunction) 
	\upupsilon 
	\upzeta \sqrt{\mbox{\upshape det} \gtorusroughfirstfund} \right) 
	+ 
	\roughgeop{x^B} 
	\left( \frac{1}{\Lunit \timefunctionarg{\muxmulevelsetvalue}} \Yvf{A}^B \upupsilon \upzeta \sqrt{\mbox{\upshape det} \gtorusroughfirstfund} \right).
\end{align}
Through a straightforward Leibniz rule-based expansion, we also have:
\begin{align} \label{E:IBPINYSECONDINTERMEDIATESTEP} 
\sqrt{|\mbox{\upshape det} \gfour|} 
\Dfour_\alpha \left( \frac{1}{\upmu} \upupsilon \upzeta \Yvf{A}\right)^{\alpha}
& 
= 
\left\lbrace 
	\frac{1}{\Lunit \timefunctionarg{\muxmulevelsetvalue}} 
	(\Yvf{A} \upupsilon) \upzeta  
	+ 
	\frac{1}{\Lunit \timefunctionarg{\muxmulevelsetvalue}} 
	\upupsilon 
	(\Yvf{A} \upzeta)  
	- 
	\frac{1}{\Lunit \timefunctionarg{\muxmulevelsetvalue}}
	\frac{\Yvf{A}\upmu}{\upmu} \upupsilon \upzeta  
	+ 
	\frac{1}{\Lunit \timefunctionarg{\muxmulevelsetvalue}} \upupsilon \upzeta 
	\left(\Dfour_\alpha \Yvf{A}^{\alpha}\right) 
\right\rbrace 
\sqrt{\mbox{\upshape det} \gtorusroughfirstfund}.
\end{align} 
Next we use \eqref{E:COVARIANTDIVERGENCEOFSPACETIMEVECTORFIELDINTERMSOFRESCALEDFRAME}, the fact that 
$\Yvf{A}$ is $\ell_{t,u}$-tangent, and Lemma~\ref{L:BASICPROPERTIESOFVECTORFIELDS}
to express the term $\Dfour_\alpha \Yvf{A}^{\alpha}$ on RHS~\eqref{E:IBPINYSECONDINTERMEDIATESTEP} 
as follows:
\begin{align} \label{E:SPACETIMECOVARIANTDIVERGENCEOFYCOMMUTATOR}
	\Dfour_\alpha \Yvf{A}^{\alpha}
	& 
	= 
	\frac{1}{\upmu} \Yvf{A} \upmu 
	+ 
	\angdiv \Yvf{A}
	= 
	\frac{1}{\upmu} \Yvf{A} \upmu 
	+ 
	\frac{1}{2} \mytr_{\gtorus} \angdeform{\Yvf{A}}.
\end{align}	
Using \eqref{E:SPACETIMECOVARIANTDIVERGENCEOFYCOMMUTATOR} to substitute for the
factor $\Dfour_\alpha \Yvf{A}^{\alpha}$ on RHS~\eqref{E:IBPINYSECONDINTERMEDIATESTEP}
and noting that the factor $\frac{1}{\upmu} \Yvf{A} \upmu$ on RHS~\eqref{E:SPACETIMECOVARIANTDIVERGENCEOFYCOMMUTATOR}
cancels the singular $\frac{\Yvf{A}\upmu}{\upmu} \upupsilon \upzeta$ term in \eqref{E:IBPINYINTERMEDIATESTEP},
we deduce that:
\begin{align} \label{E:IBPINYTHIRDINTERMEDIATESTEP} 
\sqrt{|\mbox{\upshape det} \gfour|} 
\Dfour_\alpha \left( \frac{1}{\upmu} \upupsilon \upzeta \Yvf{A}\right)^{\alpha}
& 
= 
\left\lbrace 
	\frac{1}{\Lunit \timefunctionarg{\muxmulevelsetvalue}} 
	(\Yvf{A} \upupsilon) \upzeta  
	+ 
	\frac{1}{\Lunit \timefunctionarg{\muxmulevelsetvalue}} \upupsilon 
	(\Yvf{A} \upzeta)  
	 + 
	\frac{1}{2} 
	\frac{1}{\Lunit \timefunctionarg{\muxmulevelsetvalue}} \upupsilon \upzeta 
	\mytr_{\gtorus} \angdeform{\Yvf{A}}
\right\rbrace 
\sqrt{\mbox{\upshape det} \gtorusroughfirstfund}.
\end{align} 
Next we integrate RHS~\eqref{E:IBPINYSECONDINTERMEDIATESTEP}
over $[\timefunction_0,\timefunction) \times [- \rightu,u] \times \mathbb{T}^2$
with respect to $\mathrm{d} x^2 \, \mathrm{d}x^3 \mathrm{d} u' \mathrm{d} \timefunction'$,
equate it to the integral of RHS~\eqref{E:IBPINYTHIRDINTERMEDIATESTEP},
and use Fubini's theorem. In view of definitions~\eqref{E:AREAFORMROUGHTORUS} and \eqref{E:VOLUMEFORMSPACTEIMROUGHCOORDINATES},
we see that the spacetime integrals of the terms on RHS~\eqref{E:IBPINYTHIRDINTERMEDIATESTEP}
appear in \eqref{E:PRELIMINARYIBPINY}.
The integral of the last term on RHS~\eqref{E:IBPINYSECONDINTERMEDIATESTEP} vanishes because $\mathbb{T}^2$
is a closed manifold. Finally, we note that the integral of the first term on RHS~\eqref{E:IBPINYSECONDINTERMEDIATESTEP}
yields, in view of \eqref{E:ROUGHHYPERSURFACEINTEGRAL} and the fundamental theorem of calculus, 
the two hypersurface integrals in \eqref{E:PRELIMINARYIBPINY}.

\end{proof}

We are now ready to establish the 
integration by parts identity
that we will use to control top-order wave equation error terms
that are tied to the partially modified quantities
defined in \eqref{E:PARTIALMODIFIEDQUANTITY}. 

\begin{lemma}[The key integration by parts identity tied to the partially modified quantities]
\label{L:KEYIBPIDENTIFYFORWAVEEQUATIONENERGYESTIMATESINVOLVINGPARTIALLYMODQUANT} 
Let $(\timefunction,u) \in [\timefunction_0,\timefunctionboot)\times [- \rightu,\leftu]$,
and let $\varphi$ and $\upeta$ be scalar functions on $\twoargMrough{[\timefunction_0,\timefunction),[- \rightu,u]}{\muxmulevelsetvalue}$.
Let $N \geq 1$.
Then the following integration by parts identity holds:
\begin{align}
\begin{split} \label{E:KEYIBPIDENTIFYFORWAVEEQUATIONENERGYESTIMATES}
	&
	\int_{\twoargMrough{[\timefunction_0,\timefunction),[- \rightu,u]}{\muxmulevelsetvalue}} 
		(1 + 2 \upmu) 
		(\muX \varphi) 
		(\argLrough{\muxmulevelsetvalue} \tander^N \varphi) 
		\Yvf{A} \upeta 
	\, \volMRoughCoordinates 
		\\
	& 
	= 
	\int_{\twoargMrough{[\timefunction_0,\timefunction),[- \rightu,u]}{\muxmulevelsetvalue}}  
		(1 + 2 \upmu) (\muX \varphi) 
		(\Yvf{A} \tander^N \varphi) 
		\argLrough{\muxmulevelsetvalue} \upeta 
	\, \volMRoughCoordinates 
		\\
	& \ \
		- \int_{\hypthreearg{\timefunction}{[- \rightu,u]}{\muxmulevelsetvalue} } 
				(1 + 2 \upmu) \muX \varphi (\Yvf{A}\tander^N \varphi)\upeta 
			\, \volRoughHypersurface 
		+ 
		\int_{\hypthreearg{\timefunction}{[- \rightu,u]}{\muxmulevelsetvalue}} 
			(\Yvf{A} \timefunction)
			(1 + 2 \upmu) 
			(\muX \varphi) 
			(\argLrough{\muxmulevelsetvalue} \tander^N \varphi) 
			\upeta
		\, \volRoughHypersurface 
		\\
	&  \ \ 
		+ 
		\int_{\hypthreearg{\timefunction_0}{u}{\muxmulevelsetvalue}}  
			(1 + 2 \upmu) \muX \varphi (\Yvf{A}\tander^N \varphi) \upeta 
		\, \volRoughHypersurface 
		- 
		\int_{\hypthreearg{\timefunction_0}{u}{\muxmulevelsetvalue}} 
			(\Yvf{A} \timefunction)
			(1 + 2 \upmu) 
			(\muX \varphi) 
			(\argLrough{\muxmulevelsetvalue} \tander^N \varphi) \upeta
		\, \volRoughHypersurface 
			\\
	& \ \ 
		+ 
		\int_{\twoargMrough{[\timefunction_0,\timefunction),[- \rightu,u]}{\muxmulevelsetvalue}} 
			\ErrorIBP[\tander^N\varphi;\upeta] 
		\, \volMRoughCoordinates,
\end{split}
\end{align} 
where:  
\begin{align} 
\begin{split} \label{E:ERRORTERM1KEYIBPIDENTIFYFORWAVEEQUATIONENERGYESTIMATES}
	\ErrorIBP[\tander^N \varphi;\upeta](\timefunction,u,x^2,x^3) 
	&  
	\eqdef 
	\frac{1}{\Lunit \timefunctionarg{\muxmulevelsetvalue}} 
	(1 + 2 \upmu) 
	(\muX \varphi) 
	([\Lunit, \Yvf{A}] \tander^N \varphi) 
	\upeta 
		\\
& \ \
	+ 
	2
	\frac{1}{\Lunit \timefunctionarg{\muxmulevelsetvalue}}
	(\Lunit \upmu)
	(\muX \varphi)
	(\Yvf{A} \tander^N \varphi)
	\upeta
	+
	\frac{1}{\Lunit \timefunctionarg{\muxmulevelsetvalue}} 
	(1 + 2 \upmu) 
	(\Lunit \muX \varphi) 
	(\Yvf{A} \tander^N \varphi)
	\upeta
		\\
& \ \
	-
	2
	\frac{1}{\Lunit \timefunctionarg{\muxmulevelsetvalue}}
	(\Yvf{A} \upmu)
	(\muX \varphi)
	(\Lunit \tander^N \varphi)
	\upeta
	-
	\frac{1}{\Lunit \timefunctionarg{\muxmulevelsetvalue}} 
	(1 + 2 \upmu) 
	(\Yvf{A} \muX \varphi) 
	(\Lunit \tander^N \varphi)
	\upeta.
\end{split} 
\end{align}

\end{lemma}
\begin{proof}
We now summarize the tedious but straightforward calculations that yield \eqref{E:KEYIBPIDENTIFYFORWAVEEQUATIONENERGYESTIMATES}. 
First, 
we use definition~\eqref{E:LROUGH}
and the integration by parts identity \eqref{E:PRELIMINARYIBPINY}
with $\upupsilon \eqdef \upeta$ and 
$\upzeta \eqdef (1 + 2 \upmu) 
		(\muX \varphi) 
		(\Lunit \tander^N \varphi)$
to remove the $\Yvf{A}$ operator from $\upeta$ on LHS~\eqref{E:KEYIBPIDENTIFYFORWAVEEQUATIONENERGYESTIMATES}. 
This produces the main integral 
$-
\int_{\twoargMrough{[\timefunction_0,\timefunction),[- \rightu,\leftu]}{\muxmulevelsetvalue}} 
	\frac{1}{\Lunit \timefunctionarg{\muxmulevelsetvalue}} 
	(1+2 \upmu) (\muX \varphi) 
	(\Yvf{A} \Lunit \tander^N \varphi) 
	\upeta 
\, \volMRoughCoordinates$ 
as well as many error integrals. 
We then commute $\Yvf{A}$ and $\Lunit$ and appeal to definition \eqref{E:LROUGH}
to rewrite the main integral as follows:
\begin{align} 
\begin{split} \label{E:PROOFSTEPKEYIBPIDENTIFYFORWAVEEQUATIONENERGYESTIMATES}
&
- \int_{\twoargMrough{[\timefunction_0,\timefunction),[- \rightu,\leftu]}{\muxmulevelsetvalue}} 
		(1+2 \upmu) 
		(\muX \varphi) 
		(\argLrough{\muxmulevelsetvalue} \Yvf{A} \tander^N \varphi) 
		\upeta
\, \volMRoughCoordinates
	\\
&
+
\int_{\twoargMrough{[\timefunction_0,\timefunction),[- \rightu,\leftu]}{\muxmulevelsetvalue}} 
		\frac{1}{\Lunit \timefunctionarg{\muxmulevelsetvalue}} 
		(1+2 \upmu) 
		(\muX \varphi)
		([\Lunit,\Yvf{A}] \tander^N \varphi)
		\upeta
\, \volMRoughCoordinates.
\end{split}
\end{align}
Finally, we use the integration by parts identity \eqref{E:PRELIMINARYIBPINLROUGH} 
 $\upupsilon \eqdef \Yvf{A} \tander^N \varphi$ 
and 
$\upzeta \eqdef 
(1+2 \upmu) 
(\muX \varphi) 
\upeta
$
to remove the factor of $\argLrough{\muxmulevelsetvalue}$ from the factor $\argLrough{\muxmulevelsetvalue} \Yvf{A}\tander^N \varphi$
in the first integral in \eqref{E:PROOFSTEPKEYIBPIDENTIFYFORWAVEEQUATIONENERGYESTIMATES}.


\end{proof}

\subsection{Fundamental energy identity} 
\label{SS:ENERGYIDENTITY}
In this section, 
we derive the fundamental energy--null-flux integral identities used to derive estimates. 

\subsubsection{Energy-momentum tensor, energy currents, and the multiplier vectorfield for the wave variables}
\label{SSS:ENERGYMOMENTUMTENSORANDCURRENTSFORWAVEEQUATIONS}
To derive energy identities for the wave variables $\wavearray$, which solve the quasilinear wave equations
\eqref{E:VELOCITYWAVEEQUATION}--\eqref{E:ENTROPYWAVEEQUATION},
we rely on the well-known multiplier method, which we now introduce.

Let $f$ be a scalar function (in our applications, $f$ will be some vectorfield derivative of one of the wave variables). 
We define the \emph{energy-momentum tensor} associated to $f$ to be the following symmetric type $\binom{0}{2}$ tensorfield,
where we recall that $\Dfour$ is the Levi-Civita connection of the acoustical metric $\gfour$:
\begin{align} \label{E:DEFOFENERGYMOMENTUM} 
\enmomem_{\alpha \beta} 
&	= 
\enmomem_{\alpha \beta}[f] 
\eqdef 
(\Dfour_\alpha f) 
\Dfour_\beta f 
-
\frac{1}{2} 
\gfour_{\alpha \beta} 
(\gfour^{-1})^{\kappa \lambda}
(\Dfour_{\kappa} f)
\Dfour_{\lambda} f.
\end{align}

Given any scalar function $f$ and any multiplier vectorfield $Z$, 
we define the corresponding \emph{energy current} vectorfield as follows:
\begin{align}  \label{E:WAVEEQUATIONENERGYCURRENT}
	\Jenarg{Z}{\alpha}[f] 
	& 
	\eqdef \enmomem^{\alpha \beta}[f] Z_{\beta}.
\end{align}
Recall that the deformation tensor of $Z$ is the following symmetric type $\binom{0}{2}$ tensorfield:
\begin{align} \label{E:AGAINDEFORMATIONTENSORDEF}
\deformarg{Z}{\alpha}{\beta} 
& \eqdef \Dfour_\alpha Z_\beta + \Dfour_\beta Z_\alpha.
\end{align}

The starting point for our derivation of our $L^2$-type integral identities 
for solutions of covariant wave equations is the following well-known identity,
which follows easily from the definitions:
\begin{align} \label{E:DIVERGENCEOFWAVEENERGYCURRENT}
 \Dfour_\alpha \Jenarg{Z}{\alpha}[f] 
& =  
(\square_{\gfour(\wavearray)} f) Z f 
+  
\frac{1}{2} \enmomem^{\alpha \beta} \deformarg{Z}{\alpha}{\beta}.
\end{align}

In order to obtain wave equation energy estimates that are sufficient to allow us
to track the solution up to the singular boundary, 
we use the multiplier vectorfield from the next definition.

\begin{definition}[The $\gfour$-timelike multiplier vectorfield $\multipliervectorfield $] 
We define $\multipliervectorfield$ to be the following vectorfield:
\begin{align} \label{E:MULTIPLIERVECTORFIELD}
\multipliervectorfield 
& \eqdef (1 + 2 \upmu) \Lunit 
+ 
2 \muX.
\end{align}
\end{definition}

Simple calculations based on Lemma~\ref{L:BASICPROPERTIESOFVECTORFIELDS} 
imply that $\gfour(\multipliervectorfield,\multipliervectorfield) = -4\upmu(1 + \upmu)$.
Hence, $\multipliervectorfield$ is $\gfour$-timelike whenever $\upmu > 0$. 
This property is important because it leads to coercive energy identities.

\subsubsection{Building block energies, null-fluxes, and spacetime integrals}
\label{SSS:BUILDINGBLOCKENERGIESANDNULLFLUXES}
We are now ready to define our building block energies and null-fluxes for the fluid variables.
We also define spacetime integrals that play a crucial role in our energy estimates;
see Def.\,\ref{D:COERCIVESPACETIMEWAVENERGYINTEGRAL}.

\begin{definition}[Energies and null-fluxes for the fluid variables] 
\label{D:WAVEANDTRANSPORTENERGIESANDNULLFLUXES} 
Let $f$ be a scalar function on $\twoargMrough{[\timefunction_0,\timefunctionboot),[- \rightu,\leftu]}{\muxmulevelsetvalue}$.
Recall that $\hypunitnormalarg{\muxmulevelsetvalue}$ denotes the future-directed $\gfour$-timelike unit normal of 
$\hypthreearg{\timefunction}{[- \rightu,u]}{\muxmulevelsetvalue}$,
that $\multipliervectorfield$ denotes the multiplier vectorfield defined in \eqref{E:MULTIPLIERVECTORFIELD},
and that the area forms $\volRoughHypersurface$ and $\volPuRoughCoordinates$ are defined in
Def.\,\ref{D:ROUGHVOLFORMS}.
For $(\timefunction,u) \in [\timefunction_0,\timefunctionboot)\times [- \rightu,\leftu]$,
we respectively define the \emph{wave energy} and the \emph{null-flux} associated to $f$ as follows:
\begin{subequations} 
\begin{align}
\mathbb{E}_{(\textnormal{Wave})}[f](\timefunction,u) 
& 
\eqdef 
\int_{\hypthreearg{\timefunction}{[- \rightu,u]}{\muxmulevelsetvalue}}  
	\enmomem[f](\multipliervectorfield,\hypunitnormalarg{\muxmulevelsetvalue})|\Rtransarg{\muxmulevelsetvalue}|_{\gfour} 
\, \volRoughHypersurface, 
	\label{E:WAVEENERGYDEF} 
		\\
\mathbb{F}_{(\textnormal{Wave})}[f] (\timefunction,u) 
& 
\eqdef \int_{\nullhypthreearg{\muxmulevelsetvalue}{u}{[\timefunction_0,\timefunction)}} 
	\frac{1}{\Lunit \timefunctionarg{\muxmulevelsetvalue}} \enmomem[f](\multipliervectorfield,\Lunit) 
\, \volPuRoughCoordinates.
\label{E:WAVENULLFLUXDEF}
\end{align}
\end{subequations}

If $f$ is a scalar function on $\twoargMrough{[\timefunction_0,\timefunctionboot),[- \rightu,\leftu]}{\muxmulevelsetvalue}$ 
and $(\timefunction,u) \in [\timefunction_0,\timefunctionboot)\times [- \rightu,\leftu]$,
then we respectively define the \emph{transport energy} and the {null-flux} associated to $f$ 
as follows, where $\phi = \phi(u)$ is the cut-off function introduced in Def.\,\ref{D:WTRANSANDCUTOFF}:
\begin{subequations}
\begin{align}
 \mathbb{E}_{(\textnormal{Transport})}[f](\timefunction,u) 
& 
\eqdef 
\int_{\hypthreearg{\timefunction}{[- \rightu,u]}{\muxmulevelsetvalue}} 
	\left(\upmu - \frac{\muxmulevelsetvalue \phi(u')}{\Lunit \upmu}\right) 
	f^2 
\, \volRoughHypersurface, \label{E:TRANSPORTENERGYDEF} 
	\\
 \mathbb{F}_{(\textnormal{Transport})}[f] (\timefunction,u) 
	& 
	\eqdef 	
	\int_{\nullhypthreearg{\muxmulevelsetvalue}{u}{[\timefunction_0,\timefunction)}} 
		\frac{1}{\Lunit \timefunctionarg{\muxmulevelsetvalue}} 
	f^2 
	\, \volPuRoughCoordinates. 
\label{E:TRANSPORTNULLFLUXDEF}
 \end{align}
\end{subequations}
\end{definition}

The integrals appearing in the next definition yield spacetime 
$L^2$ control of $\angrmD f$ \emph{without degenerate $\upmu$ weights}. 
These integrals arise from favorable terms in the wave equation energy identities (see Prop.\,\ref{P:FUNDAMENTALENERGYNULLFLUXIDENTITIES}).
They are of crucial importance for our energy estimates.

\begin{definition}[Key coercive spacetime integrals] 
\label{D:COERCIVESPACETIMEWAVENERGYINTEGRAL}
If $f$ is a scalar function on $\twoargMrough{[\timefunction_0,\timefunctionboot),[- \rightu,\leftu]}{\muxmulevelsetvalue}$
and $(\timefunction,u) \in [\timefunction_0,\timefunctionboot)\times [- \rightu,\leftu]$,
then we define the following spacetime integral,
where $\angrmD f$ is as in Def.\,\ref{D:ANGULARDIFFERENTIAL},
$\mathbf{1}_{[-\interestingu,\interestingu]} = \mathbf{1}_{[-\interestingu,\interestingu]}(u')$ 
denotes the characteristic function
of the interval $[-\interestingu,\interestingu]$,
$\phi = \phi(u')$ is the cut-off function from Def.\,\ref{D:WTRANSANDCUTOFF},
and 
$\volMRoughCoordinates 
= \volMRoughCoordinates(\timefunction',u',x^2,x^3)
$
is defined in \eqref{E:VOLUMEFORMSPACTEIMROUGHCOORDINATES}:
\begin{align} \label{E:COERCIVESPACETIMEWAVENERGYINTEGRAL}
\spacetimeintegralcontrolwave[f](\timefunction,u) 
&
\eqdef 
\int_{\twoargMrough{[\timefunction_0,\timefunction),[- \rightu,u]}{\muxmulevelsetvalue}} 
	\left\lbrace
		- 
		\frac{1}{2}
		\mathbf{1}_{[-\interestingu,\interestingu]}(u')
		(\argLrough{\muxmulevelsetvalue} \upmu) 
		+
		\frac{1}{\Lunit \timefunctionarg{\muxmulevelsetvalue}}
		\muxmulevelsetvalue \phi 
	\right\rbrace
	|\angrmD f|_{\gtorus}^2 
\, 
\volMRoughCoordinates.
\end{align} 

\end{definition}

\subsubsection{The fundamental energy--null-flux integral identities}
\label{SSS:ENERGYNULLFLUXINTEGRALIDENTITIES}
In the next proposition, 
we provide the fundamental energy--null-flux integral identities that form the foundation
of our hyperbolic energy estimates for the fluid variables. 

\begin{proposition}[Fundamental energy--null-flux identities] \label{P:FUNDAMENTALENERGYNULLFLUXIDENTITIES}
\hfill

\noindent \underline{\textbf{Wave equation energy--null-flux identity}}.
Suppose that on $\twoargMrough{[\timefunction_0,\timefunctionboot),[- \rightu,\leftu]}{\muxmulevelsetvalue}$,
the scalar function $f$ is a solution to the inhomogeneous covariant wave equation 
$\upmu \Box_{\gfour} f = \mathfrak{G}$. Then for $(\timefunction,u) \in [\timefunction_0,\timefunctionboot)\times [- \rightu,\leftu]$, 
the following identity holds,
where $\spacetimeintegralcontrolwave[f](\timefunction,u)$ is defined in 
\eqref{E:COERCIVESPACETIMEWAVENERGYINTEGRAL}:
\begin{align}
\begin{split} \label{E:FUNDAMENTALENERGYINTEGRALIDENTITCOVARIANTWAVES}
\mathbb{E}[f]_{(\textnormal{Wave})}(\timefunction,u) 
+ 
\mathbb{F}[f]_{(\textnormal{Wave})}(\timefunction,u) 
+ 
\spacetimeintegralcontrolwave[f](\timefunction,u) 
&
= 
\mathbb{E}[f]_{(\textnormal{Wave})}(\timefunction_0,u) 
+ \mathbb{F}[f]_{(\textnormal{Wave})}(\timefunction,- \rightu) 
	\\
&
- 
\int_{\twoargMrough{[\timefunction_0,\timefunctionboot),[- \rightu,\leftu]}{\muxmulevelsetvalue}} \frac{1}{\Lunit \timefunctionarg{\muxmulevelsetvalue}} \left\lbrace (1 + 2 \upmu)\Lunit f + 2\muX f\right\rbrace \mathfrak{G} \,  \volMRoughCoordinates  
	\\
& \ \
+ 
\int_{\twoargMrough{[\timefunction_0,\timefunctionboot),[- \rightu,\leftu]}{\muxmulevelsetvalue}} 
	\frac{1}{\Lunit \timefunctionarg{\muxmulevelsetvalue}} 
	{^{(\multipliervectorfield)}\mathfrak{B}}[f] 
\, \volMRoughCoordinates.
\end{split}
\end{align}
The error term ${^{(\multipliervectorfield)}\mathfrak{B}}[f]$ appearing 
in the last integral on RHS~\eqref{E:FUNDAMENTALENERGYINTEGRALIDENTITCOVARIANTWAVES}
can be decomposed as follows, 
where $\mathbf{1}_{[-\interestingu,\interestingu]^c} = \mathbf{1}_{[-\interestingu,\interestingu]^c}(u')$ 
denotes the characteristic function
of $[-\interestingu,\interestingu]^c = (-\infty,-\interestingu) \cup (\interestingu,\infty)$:
\begin{align}  \label{E:WAVEENERGYIDENTITYBULKERRORTERM}
{^{(\multipliervectorfield)}\mathfrak{B}}[f] 
& \eqdef 
	\frac{1}{2} \mathbf{1}_{[-\interestingu,\interestingu]^c} (\Lunit \upmu) |\angrmD f|_{\gtorus}^2 
	+ 
	\sum_{i=1}^6 {^{(\multipliervectorfield)}\mathfrak{B}_{(i)}[f]},
\end{align}
where: 
\begin{subequations}
\begin{align}
	{^{(\multipliervectorfield)}\mathfrak{B}_{(1)}[f]} 
	& 
	\eqdef (\Lunit f)^2 
	\left\lbrace 
		- \frac{1}{2} \Lunit \upmu 
		+ 
		\muX \upmu 
		- 
		\frac{1}{2} \upmu \mytr_{\gtorus} \upchi 
		+ 
		2 \upmu \mytr_{\gtorus} \angktan 
		+ 
		2 \mytr_{\gtorus} \angktrans
	\right\rbrace,  
	\label{E:WAVEENERGYIDENTITYBULKERRORTERM1} 
		\\
	 {^{(\multipliervectorfield)}\mathfrak{B}_{(2)}[f]} 
	& 
	\eqdef 
	- (\Lunit f) (\muX f) 
	\left\lbrace 
		\mytr_{\gtorus} \upchi 
		+ 
		2 \upmu \mytr_{\gtorus} \angktan 
		+ 
		2 \mytr_{\gtorus} \angktrans 
	\right\rbrace, 
	\label{E:WAVEENERGYIDENTITYBULKERRORTERM2} 
		\\
	 {^{(\multipliervectorfield)}\mathfrak{B}_{(3)}[f]} 
		& \eqdef  
		|\angrmD f|_{\gtorus}^2 
		\left\lbrace 
			\Rtransarg{\muxmulevelsetvalue} \upmu 
			+ 
			\upmu \Roughtoritangentvectorfieldarg{\muxmulevelsetvalue} \upmu 
			+ 
			2 \upmu  \Lunit \upmu 
			+ 
			\frac{1}{2}
			\upmu \mytr_{\gtorus}\upchi 
			+ 
			\upmu^2 \mytr_{\gtorus} \angktan 
			+ 
			\upmu \mytr_{\gtorus} \angktrans 
		\right\rbrace, 
			\label{E:WAVEENERGYIDENTITYBULKERRORTERM3}
			\\
	  {^{(\multipliervectorfield)}\mathfrak{B}_{(4)}[f]} 
		& 
		\eqdef 
		(\Lunit f) (\angrmD^{\#} f)
		\cdot
		\left\lbrace 
			(1 - 2 \upmu) \angrmD \upmu 
			+ 
			2 \upmu \zetatan 
			+ 
			2 \zetatrans 
		\right\rbrace, 
			\label{E:WAVEENERGYIDENTITYBULKERRORTERM4} 
				\\
	  {^{(\multipliervectorfield)}\mathfrak{B}_{(5)}[f]} 
		& \eqdef 
		-2 
		(\muX f) 
		(\angrmD^{\#} f)
		\cdot 
		\left\lbrace 
			\angrmD \upmu 
			+ 
			2 \upmu \zetatan 
			+ 
			2 \zetatrans 
		\right\rbrace, \label{E:WAVEENERGYIDENTITYBULKERRORTERM5} 
			\\
	  {^{(\multipliervectorfield)}\mathfrak{B}_{(6)}[f]} 
		& \eqdef 
		- 
		\upmu \angrmD^{\#} f \otimes \angrmD^{\#} f 
		\cdot 
		\left\lbrace 
			\upchi + 2 \upmu \angktan + 2 \angktrans 
		\right\rbrace. \label{E:WAVEENERGYIDENTITYBULKERRORTERM6}
 \end{align}
 \end{subequations}
	In \eqref{E:WAVEENERGYIDENTITYBULKERRORTERM1}--\eqref{E:WAVEENERGYIDENTITYBULKERRORTERM6}
	and below in \eqref{E:ENERGYIDENTITYFORBTRANSPORTEQUATIONS},
	the vectorfields $\Rtransarg{\muxmulevelsetvalue}$ and $\Roughtoritangentvectorfieldarg{\muxmulevelsetvalue}$ are as in Def.\,\ref{D:GEOMETRICVECTORFIELDSADAPTEDTOROUGHFOLIATIONS}
	and the $\ell_{t,u}$-tangent tensorfields $\upchi$,
	$\angktan$, $\angktrans$, $\zetatan$, and $\zetatrans$
	are as in Lemma~\ref{L:USEFULIDENTITIESANDDECOMPOSITIONSFORSECONDFUNDAMENTALFORMSANDTORSION}.

\medskip
\noindent \underline{\textbf{Transport equation energy--null-flux identity}}.
Suppose that on $\twoargMrough{[\timefunction_0,\timefunctionboot),[- \rightu,\leftu]}{\muxmulevelsetvalue}$,
$f$ is a solution to the inhomogeneous transport equation $\upmu \Transport f = \mathfrak{G}$. 
Then for $(\timefunction,u) \in [\timefunction_0,\timefunctionboot)\times [- \rightu,\leftu]$, 
the following identity holds:
\begin{align} 
\begin{split}
\label{E:ENERGYIDENTITYFORBTRANSPORTEQUATIONS}
\mathbb{E}[f]_{(\textnormal{Transport})}(\timefunction,u) 
+ 
\mathbb{F}[f]_{(\textnormal{Transport})}(\timefunction,u) 
= 
\mathbb{E}[f]_{(\textnormal{Transport})}(\timefunction_0,u) + \mathbb{F}[f]_{(\textnormal{Transport})}(\timefunction,-\rightu) 
	\\
+ 
2 
\int_{\twoargMrough{[\timefunction_0,\timefunctionboot),[- \rightu,\leftu]}{\muxmulevelsetvalue}} 
	\frac{1}{\Lunit \timefunctionarg{\muxmulevelsetvalue}} f \cdot G 
\, \volMRoughCoordinates 
+ 
\int_{\twoargMrough{[\timefunction_0,\timefunctionboot),[- \rightu,\leftu]}{\muxmulevelsetvalue}}\frac{1}{\Lunit \timefunctionarg{\muxmulevelsetvalue}} \left\lbrace 
	\Lunit \upmu + \upmu \mytr_{\gtorus} \angk 
\right\rbrace f^2 
\, \volMRoughCoordinates.
\end{split}
\end{align}

\end{proposition}

\begin{proof}

\noindent \underline{\textbf{Proof of \eqref{E:FUNDAMENTALENERGYINTEGRALIDENTITCOVARIANTWAVES}}}:

\hfill

We will apply the divergence identity \eqref{E:DIVERGENCEOFWAVEENERGYCURRENT} with $Z \eqdef \multipliervectorfield$,
where $\multipliervectorfield$ is defined in \eqref{E:MULTIPLIERVECTORFIELD}.
Throughout, we often use the abbreviated notation
$\mathbf{J} \eqdef \Jen{\multipliervectorfield}[f]$
to denote the energy current defined by \eqref{E:WAVEEQUATIONENERGYCURRENT}.
To proceed, we fist write $\mathbf{J}$
in terms of the rough adapted coordinate partial derivative vectorfields 
as follows:
\begin{align} \label{E:WAVEENERGYCURRENTINTERMSOFROUGHGEOMETRICCOORDINATES}
	\mathbf{J}
	& = 
	\mathbf{J}^{\timefunction} \roughgeop{\timefunction} 
	+ 
	\mathbf{J}^u \roughgeop{u} + \mathbf{J}^A \roughgeop{x^A}.
\end{align}

Next, we claim that the following identities hold, where $|\Rtransarg{\muxmulevelsetvalue}|_{\gfour}$
is as in \eqref{E:SIZEOFRTRANS}:
\begin{subequations}
\begin{align}
\mathbf{J}^{\timefunction} 
	& = 
	- \frac{(\Lunit \timefunctionarg{\muxmulevelsetvalue}) |\Rtransarg{\muxmulevelsetvalue}|_{\gfour}}{\upmu} \enmomem[f](\multipliervectorfield,\hypunitnormalarg{\muxmulevelsetvalue}), 
		\label{E:WAVECURRENTTIMEFUNCTIONCOMPONENT} 
	\\
\mathbf{J}^u 
	& = - \frac{1}{\upmu} \enmomem[f](\multipliervectorfield,\Lunit).
	\label{E:WAVECURRENTUCOMPONENT}
\end{align}
\end{subequations}
To prove \eqref{E:WAVECURRENTTIMEFUNCTIONCOMPONENT}, we first note that since
$\hypunitnormalarg{\muxmulevelsetvalue}$ is $\gfour$-orthogonal to $\hypthreearg{\timefunction}{[- \rightu,u]}{\muxmulevelsetvalue}$
while $\roughgeop{u}$ and $\roughgeop{x^A}$ are tangent to $\hypthreearg{\timefunction}{[- \rightu,u]}{\muxmulevelsetvalue}$,
we have, in view of \eqref{E:WAVEEQUATIONENERGYCURRENT} and \eqref{E:WAVEENERGYCURRENTINTERMSOFROUGHGEOMETRICCOORDINATES}:
\begin{align}  \label{E:ENMOMENTUMCOMPONENTMULTIPLIERHYPNORMAL}
	\enmomem[f](\multipliervectorfield,\hypunitnormalarg{\muxmulevelsetvalue}) 
	& = 
	\gfour\left(\Jen{\multipliervectorfield}[f],\hypunitnormalarg{\muxmulevelsetvalue}\right) 
	= 
	\mathbf{J}^{\timefunction} \gfour\left(\roughgeop{\timefunction}, \hypunitnormalarg{\muxmulevelsetvalue} \right).
\end{align}
Next, we use 
\eqref{E:LROUGH}--\eqref{E:LROUGHAPPLIEDTOROUGHTIMEFUNCTIONISUNITY}
and the fact that $\Lunit u = \roughgeop{\timefunction}  u = 0$
to deduce the identity
$\roughgeop{\timefunction} 
= 
\frac{1}{\Lunit \timefunctionarg{\muxmulevelsetvalue}} \Lunit 
- 
\frac{1}{\Lunit \timefunctionarg{\muxmulevelsetvalue}} \Lunit^C \roughgeop{x^C}$,
and we use Lemma~\ref{L:BASICPROPERTIESOFVECTORFIELDS},
\eqref{E:UNITLENGTHRTRANS},
\eqref{E:HYPUNITNORMALDECOMPOSITION},
and
\eqref{E:RTRANSDIVIDEDBYMUIDENTITY} to deduce the identity
$\gfour(\Lunit,\hypunitnormalarg{\muxmulevelsetvalue}) = - \frac{\upmu}{|\Rtransarg{\muxmulevelsetvalue}|_{\gfour}}$.
Combining these two identities with \eqref{E:ENMOMENTUMCOMPONENTMULTIPLIERHYPNORMAL}
and using that $\gfour(\roughgeop{x^C},\hypunitnormalarg{\muxmulevelsetvalue}) = 0$,
we find that:
\begin{align} \label{E:SECONDENMOMENTUMCOMPONENTMULTIPLIERHYPNORMAL}
\enmomem[f](\multipliervectorfield,\hypunitnormalarg{\muxmulevelsetvalue}) 
& 
=  
\frac{1}{\Lunit \timefunctionarg{\muxmulevelsetvalue}} 
\mathbf{J}^{\timefunction} \gfour (\Lunit,\hypunitnormalarg{\muxmulevelsetvalue}) 
= 
- \frac{\upmu}{(\Lunit \timefunctionarg{\muxmulevelsetvalue}) |\Rtransarg{\muxmulevelsetvalue}|_{\gfour}} \mathbf{J}^{\timefunction},
\end{align}
which yields \eqref{E:WAVECURRENTTIMEFUNCTIONCOMPONENT}. 
The identity \eqref{E:WAVECURRENTUCOMPONENT} can be proved using similar arguments 
based on taking the $\gfour$-inner product of $\Jen{\multipliervectorfield}[f]$ with $\Lunit$ and using that
$\Lunit$ is $\gfour$-orthogonal to $\nullhyparg{u}$ as well as the identity
$\roughgeop{u} = \Rtransarg{\muxmulevelsetvalue} - \Rtransarg{\muxmulevelsetvalue}^C \roughgeop{x^C}$, which follows from
Lemma~\ref{L:BASICPROPERTIESOFVECTORFIELDS} and
\eqref{E:DEFROUGHTORITANGENTVECTORFIELD}--\eqref{E:RTRANS}.

Next, we note the following formula, 
which follows by combining the standard identity for the divergence of a vectorfield expressed relative to the rough adapted coordinates 
with the identities
\eqref{E:DETACOUSTICALMETRICINROUGHADAPTEDCOORDS}
and
\eqref{E:WAVECURRENTTIMEFUNCTIONCOMPONENT}--\eqref{E:WAVECURRENTUCOMPONENT}:
\begin{align}
\begin{split}  \label{E:IDFORDETWEIGHTEDCOVARIANTDIVERGENCEOFWAVEENERGYCURRENT}
	\sqrt{|\mbox{\upshape det} \gfour|} \, \Dfour_\alpha \Jenarg{\multipliervectorfield}{\alpha}[f] 
	& = 
	\roughgeop{\timefunction} 
	\left( \sqrt{|\mbox{\upshape det} \gfour|} \mathbf{J}^{\timefunction} \right) 
	+ 
	\roughgeop{u} \left( \sqrt{|\mbox{\upshape det} \gfour|} \mathbf{J}^u\right) 
	+ 
	\roughgeop{x^A} \left( \sqrt{|\mbox{\upshape det} \gfour|} \mathbf{J}^A \right) 
		\\
	& = 
	- 
	\roughgeop{\timefunction} 
	\left( 
		\enmomem[f](\multipliervectorfield,\hypunitnormalarg{\muxmulevelsetvalue}) |\Rtransarg{\muxmulevelsetvalue}|_{\gfour}
		\sqrt{ \mbox{\upshape det} \gtorusroughfirstfund} 
	\right) 
	- 
	\roughgeop{u} 
	\left( 
		\frac{1}{\Lunit \timefunctionarg{\muxmulevelsetvalue}} 
		\enmomem[f](\multipliervectorfield,\Lunit) 
		\sqrt{ \mbox{\upshape det} \gtorusroughfirstfund} 
	\right) 
	+ 
	\roughgeop{x^A} 
	\left( 
		\sqrt{|\mbox{\upshape det} \gfour|} \mathbf{J}^A 
	\right).
\end{split}
\end{align}
Integrating \eqref{E:IDFORDETWEIGHTEDCOVARIANTDIVERGENCEOFWAVEENERGYCURRENT} over 
$\twoargMrough{[\timefunction_0,\timefunction),[- \rightu,u]}{\muxmulevelsetvalue}$
with respect to $\mathrm{d} x^2 \, \mathrm{d}x^3 \mathrm{d} u' \mathrm{d} \timefunction'$, 
using \eqref{E:VOLFORMACOUSTICALMETRICROUGHADAPTED} to relate the canonical volume form $\volcanonical_{\gfour}$ to $\volMRoughCoordinates$, using \eqref{E:DIVERGENCEOFWAVEENERGYCURRENT} to substitute for $\Dfour_\alpha \Jenarg{\multipliervectorfield}{\alpha}[f]$, 
applying Fubini's theorem,
and using that the integral of the last term on RHS~\eqref{E:IDFORDETWEIGHTEDCOVARIANTDIVERGENCEOFWAVEENERGYCURRENT} 
over $\T^2$ vanishes (since $\T^2$ is a closed manifold), we deduce,
in view of definitions \eqref{E:WAVEENERGYDEF}--\eqref{E:WAVENULLFLUXDEF} and 
\eqref{E:MULTIPLIERVECTORFIELD}, the following identity:
\begin{align} 
\begin{split} \label{E:FUNDAMENTALENERGYINTEGRALIDENTITCOVARIANTWAVESINTERMEDIATE}
\mathbb{E}[f]_{(\textnormal{Wave})}(\timefunction,u) 
+ 
\mathbb{F}[f]_{(\textnormal{Wave})}(\timefunction,u) 
& 
= 
\mathbb{E}[f]_{(\textnormal{Wave})}(\timefunction_0,u) 
+ 
\mathbb{F}[f]_{(\textnormal{Wave})}(\timefunction,- \rightu) 
	\\
& \ \
- 
\int_{\twoargMrough{[\timefunction_0,\timefunctionboot),[- \rightu,\leftu]}{\muxmulevelsetvalue}} 
	\frac{1}{\Lunit \timefunctionarg{\muxmulevelsetvalue}} 
	\left\lbrace 
		(1 + 2 \upmu) \Lunit f + 2 \muX f 
	\right\rbrace 
	\mathfrak{G} 
\,  \volMRoughCoordinates 
	\\
& \ \
- 
\frac{1}{2} 
\int_{\twoargMrough{[\timefunction_0,\timefunctionboot),[- \rightu,\leftu]}{\muxmulevelsetvalue}} 
	\frac{1}{\Lunit \timefunctionarg{\muxmulevelsetvalue}} 
	\upmu \enmomem^{\alpha \beta}[f] \deformarg{\multipliervectorfield}{\alpha}{\beta} 
\, \volMRoughCoordinates.
\end{split}
\end{align}

Next, we decompose the integrand in the last integral on RHS~\eqref{E:FUNDAMENTALENERGYINTEGRALIDENTITCOVARIANTWAVESINTERMEDIATE}
$ 
- 
\frac{1}{2} 
\frac{1}{\Lunit \timefunctionarg{\muxmulevelsetvalue}} 
\upmu  (\gfour^{-1})^{\alpha \delta} (\gfour^{-1})^{\beta \sigma} \enmomem_{\delta\sigma}[f]  
\deformarg{\multipliervectorfield}{\alpha}{\beta}
$ 
using 
definitions~\eqref{E:DEFOFENERGYMOMENTUM},
\eqref{E:AGAINDEFORMATIONTENSORDEF},
and \eqref{E:MULTIPLIERVECTORFIELD}
the identity 
$(\gfour^{-1})^{\alpha \delta} 
= 
- 
\Lunit^{\alpha} \Lunit^{\delta} 
- 
\Lunit^{\alpha} X^{\delta} 
- X^{\alpha}\Lunit^{\delta} 
+ 
(\gtorus^{-1})^{\alpha \delta}
$
which follows from \eqref{E:SMOOTHTORUSMETRICINTERMSOFSIGMATMETRICANDX},
and the analogous identity for $(\gfour^{-1})^{\beta \sigma}$.
Among the many terms that arise from the expansion 
are the following two:\footnote{The precise origin of these terms is $- \enmomem[f](\Lunit,X) \deformarg{\multipliervectorfield}{L}{\muX}$, which is found in the expansion of $-\frac{1}{2} \upmu \enmomem^{\alpha \beta}[f] \deformarg{\multipliervectorfield}{\alpha}{\beta}$.}
\begin{align} \label{E:FUNDAMENTALENERGYINTEGRALIDENTITCOVARIANTWAVESINTERMEDIATE2}
  &
	\frac{1}{2} 
	\frac{1}{\Lunit \timefunctionarg{\muxmulevelsetvalue}} 
	(\Lunit \upmu) 
	|\angrmD f|_{\gtorus}^2 
	+ 
	\frac{1}{\Lunit \timefunctionarg{\muxmulevelsetvalue}} 
	(\muX \upmu) |\angrmD f|_{\gtorus}^2.
 \end{align}
We now decompose the first product in \eqref{E:FUNDAMENTALENERGYINTEGRALIDENTITCOVARIANTWAVESINTERMEDIATE2} as follows,
where we use definition~\eqref{E:LROUGH}:
\begin{align} \label{E:DECOMPOSITIONOFLMUWAVENERGYTERM}
\frac{1}{2} 
\frac{1}{\Lunit \timefunctionarg{\muxmulevelsetvalue}} 
(\Lunit \upmu) 
|\angrmD f|_{\gtorus}^2 
& = 
\frac{1}{2}  
\mathbf{1}_{[-\interestingu,\interestingu]} 
(\argLrough{\muxmulevelsetvalue} \upmu) 
|\angrmD f|_{\gtorus}^2 
+ 
\frac{1}{2}  
\mathbf{1}_{[-\interestingu,\interestingu]^c}
\frac{1}{\Lunit \timefunctionarg{\muxmulevelsetvalue}} 
(\Lunit \upmu) 
|\angrmD f|_{\gtorus}^2.
\end{align}
Since the spacetime integral of 
$
\mathbf{1}_{[-\interestingu,\interestingu]} 
(\argLrough{\muxmulevelsetvalue} \upmu) 
|\angrmD f|_{\gtorus}^2 
$ 
is found in the \emph{negative} of the coercive spacetime integral defined in \eqref{E:COERCIVESPACETIMEWAVENERGYINTEGRAL}
(i.e., $- \spacetimeintegralcontrolwave[f](\timefunction,u)$), 
we bring this term to LHS~\eqref{E:FUNDAMENTALENERGYINTEGRALIDENTITCOVARIANTWAVESINTERMEDIATE} 
as part of $\spacetimeintegralcontrolwave[f](\timefunction,u)$.
The remaining term 
$
\frac{1}{2}  
\mathbf{1}_{[-\interestingu,\interestingu]}
\frac{1}{\Lunit \timefunctionarg{\muxmulevelsetvalue}}
(\Lunit \upmu) 
|\angrmD f|_{\gtorus}^2 
$ 
in the decomposition~\eqref{E:DECOMPOSITIONOFLMUWAVENERGYTERM}
is manifestly present on RHS~\eqref{E:FUNDAMENTALENERGYINTEGRALIDENTITCOVARIANTWAVES} 
as the first term in RHS~\eqref{E:WAVEENERGYIDENTITYBULKERRORTERM}.
Next, we examine the term
$
\frac{1}{\Lunit \timefunctionarg{\muxmulevelsetvalue}}
(\muX \upmu) |\angrmD f|_{\gtorus}^2$ present in \eqref{E:FUNDAMENTALENERGYINTEGRALIDENTITCOVARIANTWAVESINTERMEDIATE2}.  
Using \eqref{E:RTRANS}, we express this term as follows:
\begin{align} \label{E:MUXMUTIMESANGULARDERIVATIVEENERGYESTIMATEINTEGRANDFACTORDECOMPOSITION}
\frac{1}{\Lunit \timefunctionarg{\muxmulevelsetvalue}}
(\muX \upmu) 
|\angrmD f|_{\gtorus}^2 
& =
- 
\frac{1}{\Lunit \timefunctionarg{\muxmulevelsetvalue}}
\muxmulevelsetvalue \phi
|\angrmD f|_{\gtorus}^2
+
\frac{1}{\Lunit \timefunctionarg{\muxmulevelsetvalue}}
|\angrmD f|_{\gtorus}^2 
\left\lbrace 
	\Rtransarg{\muxmulevelsetvalue} \upmu 
	+ 
	\upmu \Roughtoritangentvectorfieldarg{\muxmulevelsetvalue} \upmu 
\right\rbrace
|\angrmD f|_{\gtorus}^2.
\end{align}
We bring the integral of the first product 
$- 
\frac{1}{\Lunit \timefunctionarg{\muxmulevelsetvalue}}
\muxmulevelsetvalue \phi |\angrmD f|_{\gtorus}^2$ 
on RHS~\eqref{E:MUXMUTIMESANGULARDERIVATIVEENERGYESTIMATEINTEGRANDFACTORDECOMPOSITION} 
over to LHS~\eqref{E:FUNDAMENTALENERGYINTEGRALIDENTITCOVARIANTWAVESINTERMEDIATE2} 
as the remaining part of $\spacetimeintegralcontrolwave[f](\timefunction,u)$
(see definition~\eqref{E:COERCIVESPACETIMEWAVENERGYINTEGRAL}).
These terms in the last product on RHS~\eqref{E:MUXMUTIMESANGULARDERIVATIVEENERGYESTIMATEINTEGRANDFACTORDECOMPOSITION}
are manifestly present in the term ${^{(\multipliervectorfield)}\mathfrak{B}_{(3)}[f]}$
defined in \eqref{E:WAVEENERGYIDENTITYBULKERRORTERM3}.
The remaining terms in the decomposition of
$ 
- 
\frac{1}{2} 
\frac{1}{\Lunit \timefunctionarg{\muxmulevelsetvalue}} 
\upmu  (\gfour^{-1})^{\alpha \delta} (\gfour^{-1})^{\beta \sigma} \enmomem_{\delta\sigma}[f]  
\deformarg{\multipliervectorfield}{\alpha}{\beta}
$ 
can be derived using the same arguments, 
based on straightforward but tedious calculations,
given in the proof of \cite[Lemma 3.3]{jSgHjLwW2016}.
We have therefore proved \eqref{E:FUNDAMENTALENERGYINTEGRALIDENTITCOVARIANTWAVES}.

\noindent \underline{\textbf{Proof of \eqref{E:FUNDAMENTALENERGYINTEGRALIDENTITCOVARIANTWAVES}}}:
We will apply the divergence theorem to the vectorfield $\mathbf{J} \eqdef f^2 \Transport$. 
We begin by expressing this vectorfield in terms of the rough geometric coordinate vectorfields as follows:
\begin{align} 
	\mathbf{J} 
	& = 
	\mathbf{J}^{\timefunction} \roughgeop{\timefunction} 
	+ 
	\mathbf{J}^u \roughgeop{u} 
	+ 
	\mathbf{J}^A \roughgeop{x^A}.
\end{align}
Next, we claim that the following identities hold:
\begin{subequations}
\begin{align} 
\mathbf{J}^{\timefunction} 
& 
= f^2 
\left(\Lunit \timefunctionarg{\muxmulevelsetvalue} 
	- 
	\frac{\muxmulevelsetvalue \phi}{\upmu \Lunit \upmu} 
	\Lunit \timefunctionarg{\muxmulevelsetvalue} 
\right) 
	\label{E:TRANSPORTENERGYCURRENTROUGHTIMECOMPONENT} 
		\\
\mathbf{J}^u 
& = \frac{1}{\upmu} f^2. 
\label{E:TRANSPORTENERGYCURRENTROUGHUCOMPONENT}
\end{align}
\end{subequations}
To prove \eqref{E:TRANSPORTENERGYCURRENTROUGHTIMECOMPONENT}, 
we note that 
$\mathbf{J}^{\timefunction} = \mathbf{J} \timefunctionarg{\muxmulevelsetvalue} = f^2 \Transport \timefunctionarg{\muxmulevelsetvalue}$
and use the equations
$0 = \Wtrans \timefunctionarg{\muxmulevelsetvalue} 
= 
\muX \timefunctionarg{\muxmulevelsetvalue} 
+ 
\frac{\muxmulevelsetvalue\phi}{\Lunit \upmu} 
\Lunit \timefunctionarg{\muxmulevelsetvalue}$
(see \eqref{E:WTRANSDEF} and \eqref{E:IVPFORROUGHTTIMEFUNCTION}) 
and 
$\Transport = \Lunit + X$ (see \eqref{E:BISLPLUSX}). 
Similarly, to prove \eqref{E:TRANSPORTENERGYCURRENTROUGHUCOMPONENT},
we note that
$\mathbf{J}^u = \mathbf{J} u = f^2 \Transport u$
and then use Lemma~\ref{L:BASICPROPERTIESOFVECTORFIELDS}.

Next, using the divergence theorem,
\eqref{E:DETACOUSTICALMETRICINROUGHADAPTEDCOORDS}, 
and
\eqref{E:TRANSPORTENERGYCURRENTROUGHTIMECOMPONENT}--\eqref{E:TRANSPORTENERGYCURRENTROUGHUCOMPONENT}, 
we deduce that:
\begin{align} 
\begin{split} \label{E:EXPRESSIONFORDETWEIGHTEDCOVARIANTDIVERGENCEOFTRANSPORTENERGYCURRENTINROUGHADATPEDCOORDINATES}
\sqrt{|\mbox{\upshape det} \gfour|} \Dfour_\alpha \mathbf{J}^{\alpha} 
&
=
\sqrt{|\mbox{\upshape det} \gfour|} \Dfour_\alpha (f^2 \Transport^{\alpha}) 
= 
\roughgeop{\timefunction} \left( \sqrt{|\mbox{\upshape det} \gfour|} \mathbf{J}^{\timefunction}\right) 
+ 
\roughgeop{u} 
\left( 
	\sqrt{|\mbox{\upshape det} \gfour|} \mathbf{J}^u
\right) 
+ 
\roughgeop{x^A} 
\left(
	\sqrt{|\mbox{\upshape det} \gfour|} \mathbf{J}^A
\right) 
	\\
& = \roughgeop{\timefunction} \left( f^2\left\lbrace \upmu - \frac{\muxmulevelsetvalue \phi}{\Lunit \upmu}\right\rbrace\sqrt{\mbox{\upshape det} \gtorusroughfirstfund} \right) + \roughgeop{u} \left( \frac{1}{\Lunit \timefunctionarg{\muxmulevelsetvalue}} f^2 \sqrt{\mbox{\upshape det} \gtorusroughfirstfund}\right) + \roughgeop{x^A} \left(\sqrt{|\mbox{\upshape det} \gfour|} \mathbf{J}^A\right).
\end{split}
\end{align}
Next, we note the following covariant divergence identity, 
which follows from the relation $\mathbf{J} = f^2 \Transport$, 
the Leibniz rule,
Lemma~\ref{L:BASICPROPERTIESOFVECTORFIELDS},
and
\eqref{E:COVARIANTDIVERGENCEOFSPACETIMEVECTORFIELDINTERMSOFRESCALEDFRAME}: 
\begin{align} \label{E:COVARIANTDIVERGENCEOFTRANSPORTENERGYCURRENT}
	\Dfour_\alpha \mathbf{J}^{\alpha} 
	& 
	= 
	\frac{1}{\upmu} (\Lunit \upmu) f^2 
	+ 
	\mytr_{\gtorus} \angk f^2 
	+ 
	2 f \Transport f.
\end{align}
Substituting \eqref{E:COVARIANTDIVERGENCEOFTRANSPORTENERGYCURRENT} into 
LHS~\eqref{E:EXPRESSIONFORDETWEIGHTEDCOVARIANTDIVERGENCEOFTRANSPORTENERGYCURRENTINROUGHADATPEDCOORDINATES}
and then integrating both sides of the resulting identity 
over $\twoargMrough{[\timefunction_0,\timefunction),[- \rightu,u]}{\muxmulevelsetvalue}$ 
with respect to $\mathrm{d} x^2 \, \mathrm{d}x^3 \mathrm{d} u' \mathrm{d} \timefunction'$,
using \eqref{E:VOLFORMACOUSTICALMETRICROUGHADAPTED} to relate the canonical volume form 
$\volcanonical_{\gfour}$ to $\volMRoughCoordinates$,  
applying Fubini's theorem,
and using that the integral of the last term on RHS~\eqref{E:IDFORDETWEIGHTEDCOVARIANTDIVERGENCEOFWAVEENERGYCURRENT} 
over $\T^2$ vanishes (since $\T^2$ is a closed manifold)
we conclude, 
in view of definitions 
\eqref{E:TRANSPORTENERGYDEF}--\eqref{E:TRANSPORTNULLFLUXDEF},
the desired identity \eqref{E:ENERGYIDENTITYFORBTRANSPORTEQUATIONS}.

\end{proof}

\subsection{The fundamental $L^2$-controlling quantities}
\label{SS:FUNDAMENTALL2CONTROLLINGQUANTITIES}
In this section, we use the building block quantities
from Sect.\,\ref{SSS:BUILDINGBLOCKENERGIESANDNULLFLUXES} to 
construct the ``fundamental $L^2$-controlling quantities'' that we use to control 
$\wavearray$, $\vortrenormalized$, $\GradEnt$, $\VortVort$, $\DivGradEnt$, 
and their derivatives in $L^2$ in various regions.

\begin{definition}[The fundamental $L^2$-controlling quantities] 
\label{D:MAINCOERCIVE} 
In terms of the energy-null-flux quantities of Def.\,\ref{D:WAVEANDTRANSPORTENERGIESANDNULLFLUXES} 
and the differentiated conventions of Def.\,\ref{D:CONVENTIONESFORDIFFERENTIATION}, 
we define the following $L^2$-controlling quantities:

\begin{itemize}[leftmargin=*]
\item \underline{\textbf{Total wave-controlling quantities}}.
\begin{subequations}
\begin{align}
	\hypersurfacecontrolwave_N(\timefunction,u) 
	& 
	\eqdef 
	\max_{\substack{\tander^N \in \mathfrak{P}^{(N)} \\ \Psi \in \{ \RRiemann, \LRiemann, v^2, v^3, \Ent\}}} 
	\sup_{\substack{\timefunction' \in [\timefunction_0,\timefunction] 
		\\
		u' \in [- \rightu,u]}} 
	\left\lbrace 
		\mathbb{E}_{(\textnormal{Wave})}[\tander^N \Psi](\timefunction',u') 
			+ 
		\mathbb{F}_{(\textnormal{Wave})}[\tander^N \Psi](\timefunction',u')
	\right\rbrace, \label{E:WAVESPACELIKEANDNULLHYPERSURFACEL2CONTROLLINGQUANTITY} 
\\
\spacetimeintegralcontrolwave_N(\timefunction,u) 
& 
\eqdef 
\max_{\substack{\tander^N \in \mathfrak{P}^{(N)} 
	\\ \Psi \in \{ \RRiemann, \LRiemann, v^2, v^3, \Ent\}}} 
\spacetimeintegralcontrolwave[\tander^N \Psi] (\timefunction,u), 
\label{E:WAVESPACETIMEL2CONTROLLINGQUANTITY} 
	\\
\totalcontrolwave_N(\timefunction,u) 
	& 
	\eqdef 
	\max
	\left\lbrace
	\hypersurfacecontrolwave_N(\timefunction,u),
		\,
	\spacetimeintegralcontrolwave_N(\timefunction,u)
	\right\rbrace,
\label{E:WAVETOTALL2CONTROLLINGQUANTITY}
\end{align}
\end{subequations}

\medskip

\item \underline{\textbf{Partial wave-controlling quantities}}.

\begin{subequations}
\begin{align}
\hypersurfacecontrolwavepartial_N(\timefunction,u) 
& 
\eqdef 
\max_{\substack{ \tander^N \in \mathfrak{P}^{(N)} 
	\\ \Psi \in \{\LRiemann, v^2, v^3, \Ent\}}}  
\sup_{\substack{\timefunction' \in [\timefunction_0,\timefunction] 
		\\
		u' \in [- \rightu,u]}}  
	\left\lbrace 
		\mathbb{E}_{(\textnormal{Wave})}[ \tander^N\Psi](\timefunction',u') 
			+ 
		\mathbb{F}_{(\textnormal{Wave})}[\tander^N \Psi](\timefunction',u')
	\right\rbrace, \label{E:WAVEPARTIALSPACELIKEANDNULLHYPERSURFACEL2CONTROLLINGQUANTITY}  
	 \\
\spacetimeintegralcontrolwavepartial_N(\timefunction,u) & \eqdef \max_{\substack{ \tander^N \in \mathfrak{P}^{(N)} \\ \Psi \in \{\LRiemann, v^2, v^3, \Ent\}}} \spacetimeintegralcontrolwave[\tander^N \Psi] (\timefunction,u), 
\label{E:PARTIALWAVESPACETIMEL2CONTROLLINGQUANTITY}  
	\\
\totalcontrolwavepartial_N(\timefunction,u) 
& 
\eqdef 
\max
\left\lbrace
	\hypersurfacecontrolwavepartial_N(\timefunction,u),
		\,
	\spacetimeintegralcontrolwavepartial_N(\timefunction,u)
\right\rbrace.
\label{E:WAVEPARTIALL2CONTROLLINGQUANTITY} 
\end{align}
\end{subequations}

\medskip

\item \underline{\textbf{Specific vorticity- and entropy-controlling quantities}}.
\begin{subequations}
\begin{align}
	\hypersurfacecontrolVort_N(\timefunction,u) 
	& \eqdef 
	\max_{\tander^N \in \mathfrak{P}^{(N)}} 
	\sup_{\substack{\timefunction' \in [\timefunction_0,\timefunction] 
					\\
				u' \in [- \rightu,u]}}  
	\left\lbrace 
		\mathbb{E}_{(\textnormal{Transport})}[\tander^N \vortrenormalized] (\timefunction',u') 
			+  
		\mathbb{F}_{(\textnormal{Transport})}[\tander^N \vortrenormalized] (\timefunction',u') 
		\right\rbrace, 
		\label{E:VORTICITYL2CONTROLLINGQUANTITY}
			\\
\hypersurfacecontrolGradEnt_N(\timefunction,u) 
& 
\eqdef 
\max_{\tander^N \in \mathfrak{P}^{(N)}} 
\sup_{\substack{\timefunction' \in [\timefunction_0,\timefunction] 
		\\
			u' \in [- \rightu,u]}} 
	\left\lbrace 
		\mathbb{E}_{(\textnormal{Transport})}[\tander^N \GradEnt] (\timefunction',u') 
			+ 
		\mathbb{F}_{(\textnormal{Transport})}[\tander^N \GradEnt] (\timefunction',u') 
	\right\rbrace, \label{E:ENTROPYGRADIENTL2CONTROLLINGQUANTITY}
\end{align}
\end{subequations}

\begin{subequations}
\begin{align}
\toricontrolVort_N(\timefunction,u) 
& 
\eqdef
\max_{\tander^N \in \mathfrak{P}^{(N)}} 
\left\| 
	\tander^N \vortrenormalized 
\right\|_{L^2(\twoargroughtori{\timefunction,u}{\muxmulevelsetvalue})}^2, 
	\label{E:TORIVORTICITYL2CONTROLLINGQUANTITY} 
\\
	\toricontrolGradEnt_N(\timefunction,u) 
	& 
	\eqdef 
	\max_{\tander^N \in \mathfrak{P}^{(N)}} 
	\left\| 
		\tander^N \GradEnt 
	\right\|_{L^2(\twoargroughtori{\timefunction,u}{\muxmulevelsetvalue})}^2.
		\label{E:TORIENTROPYGRADIENTL2CONTROLLINGQUANTITY} 
\end{align} 
\end{subequations}

\medskip

\item \underline{\textbf{Modified fluid variable-controlling quantities}}.
\begin{subequations}
\begin{align}
	\hypersurfacecontrolVortVort_N(\timefunction,u) 
	& 
	\eqdef 
	\max_{\tander^N \in \mathfrak{P}^{(N)}} 
	\sup_{\substack{\timefunction' \in [\timefunction_0,\timefunction] 
		\\
			u' \in [- \rightu,u]}} 
		\left\lbrace 
			\mathbb{E}_{(\textnormal{Transport})}[\tander^N \VortVort] (\timefunction',u') 
				+  
			\mathbb{F}_{(\textnormal{Transport})}[\tander^N \VortVort] (\timefunction',u') 
		\right\rbrace, \label{E:MODIFIEDVORTICITYVORTICITYL2CONTROLLINGQUANTITY} 
			\\
	\hypersurfacecontrolDivGradEnt_N(\timefunction,u) 
	& 
	\eqdef 
	\max_{\tander^N \in \mathfrak{P}^{(N)}} 
	\sup_{\substack{\timefunction' \in [\timefunction_0,\timefunction] 
		\\
			u' \in [- \rightu,u]}} 
	\left\lbrace 
		\mathbb{E}_{(\textnormal{Transport})}[\tander^N \DivGradEnt] (\timefunction',u') 
			+  
		\mathbb{F}_{(\textnormal{Transport})}[\tander^N \DivGradEnt] (\timefunction',u')
	\right\rbrace,
\label{E:MODIFIEDDIVGRADENTL2CONTROLLINGQUANTITY} 
\end{align}
\end{subequations}

\begin{subequations}
\begin{align}
\toricontrolVortVort_N(\timefunction,u) 
& 
\eqdef
\max_{\tander^N \in \mathfrak{P}^{(N)}} 
\left\| 
	\tander^N \VortVort 
\right\|_{L^2(\twoargroughtori{\timefunction,u}{\muxmulevelsetvalue})}^2, 
	\label{E:TORIMODIFIEDVORTICITYL2CONTROLLINGQUANTITY} 
			\\
	\toricontrolDivGradEnt_N(\timefunction,u) 
	& 
	\eqdef 
	\max_{\tander^N \in \mathfrak{P}^{(N)}} 
	\left\| 
		\tander^N \DivGradEnt 
	\right\|_{L^2(\twoargroughtori{\timefunction,u}{\muxmulevelsetvalue})}^2.
		\label{E:TORIDIVGRADENTL2CONTROLLINGQUANTITY} 
\end{align} 
\end{subequations}

\end{itemize}

\end{definition} 

\begin{remark}[Differences between $\hypersurfacecontrolwave_N$ and $\hypersurfacecontrolwave_N^{(Partial)}$] \label{R:WHYPARTIALENERGIES} Although $\hypersurfacecontrolwavepartial_N$ might seem to be a redundant quantity, it plays an important role in our energy estimates since in the top-order case, 
the partial energy $\hypersurfacecontrolwavepartial_N$ 
is only weakly influenced by the full energy $\hypersurfacecontrolwave_N$.
Similar remarks apply to
$\spacetimeintegralcontrolwavepartial_N(\timefunction,u)$
and
$\totalcontrolwavepartial_N(\timefunction,u)$.
\end{remark}

\begin{definition}[Summed $L^2$-controlling quantities] 
\label{D:SUMMEDL2CONTROLLINGQUANTITIES}
For positive integers $N_1 < N_2$ and non-negative integers $N$, 
we define the following summed $L^2$-controlling quantities:
\begin{align} \label{D:ENERGYSUMCONVENTIONS}
\hypersurfacecontrolwave_{[N_1,N_2]}(\timefunction,u) 
& 
\eqdef \sum_{M = N_1}^{N_2} \hypersurfacecontrolwave_M(\timefunction,u), 
&
\hypersurfacecontrolVort_{\leq N}(\timefunction,u) 
& = 
\sum_{M = 0}^N \hypersurfacecontrolVort_M(\timefunction,u),
\end{align}
and similarly for the other controlling quantities. When $N = 0$, 
we often omit the subscript, e.g., we write
$\hypersurfacecontrolVort(\timefunction,u)$
instead of
$\hypersurfacecontrolVort_0(\timefunction,u)$.
\end{definition}

\subsection{The coerciveness of the fundamental $L^2$-controlling quantities}
\label{SS:COERCIVENESSFUNDAMENTALL2CONTROLLINGQUANTITIES}
In this section, we exhibit the coerciveness properties of the $L^2$ controlling quantities
from Def.\,\ref{D:MAINCOERCIVE}.

\subsubsection{Decomposition of components of the energy-momentum tensor}
\label{SSS:DECOMPOSITIONOFENERGYMOMENTUMTENSORCOMPOENTNS}
We start with the following lemma, 
which yields identities for various components of the energy-momentum tensor. 
The components $\enmomem[f](\Transport,\Lunit)$ and $\enmomem[f](\hypunitnormalarg{\muxmulevelsetvalue},\multipliervectorfield)$
are of particular interest since they appear in our energy--null-flux identity
\eqref{E:FUNDAMENTALENERGYINTEGRALIDENTITCOVARIANTWAVES}.

\begin{lemma}[Decomposition of various components of the energy-momentum tensor]
\label{L:DECOMPOSITIONOFCOMPOENTNSOFENERGYMOMENTUMTENSOR} 
Let $f$ be a scalar function, and let $\enmomem$ be the corresponding energy-momentum tensor
defined in \eqref{E:DEFOFENERGYMOMENTUM}.
Then the following identities hold,
where $|\angD f|_{\gtorus}^2 = (\gtorus^{-1})^{AB} (\geop{x^A} f) \geop{x^B} f$:
\begin{subequations}
\begin{align}
	\enmomem[f](\Lunit,\Lunit) & = (\Lunit f)^2, 
		\label{E:QLL} 
		\\
	\enmomem[f](X,\Lunit) 
	& = 
	- \frac{1}{2} (\Lunit f)^2 + \frac{1}{2} | \angD f|_{\gtorus}^2, 
		\label{E:QLX} 
			\\
	\enmomem[f](\Transport,\Lunit)
	& 
	= \frac{1}{2} (\Lunit f)^2 + \frac{1}{2}|\angD f|_{\gtorus}^2
		\label{E:QTRANSPORTL}, 
			\\
	\enmomem[f](X,\Transport) 
		& = (\Lunit f) X f + (X f)^2 \label{E:QXTRANSPORT}, 
			\\
	\enmomem[f] \left(\Lunit,\roughgeop{x^A}\right) 
		& 
		= (\Lunit f) \roughgeop{x^A} f, \label{E:QLROUGHGEOPA}
			\\
\enmomem[f]\left(X,\roughgeop{x^A}\right) 
& = 
(X f) \roughgeop{x^A} f 
+ 
\frac{1}{2} \frac{\geop{x^A} \timefunctionarg{\muxmulevelsetvalue}}{\geop{t} \timefunctionarg{\muxmulevelsetvalue}} (\Lunit f)^2 
+  
\frac{\geop{x^A} \timefunctionarg{\muxmulevelsetvalue}}{\geop{t} \timefunctionarg{\muxmulevelsetvalue}} (\Lunit f) Xf 
- 
\frac{1}{2} \frac{\geop{x^A} \timefunctionarg{\muxmulevelsetvalue}}{\geop{t} \timefunctionarg{\muxmulevelsetvalue}} | \angD f|_{\gtorus}^2, 
	\label{E:QXROUGHGEOPA} 
		\\
\enmomem[f]\left(\Transport, \roughgeop{x^A}\right) 
& 
= (\Lunit f) \roughgeop{x^A} f 
+  
(Xf) \roughgeop{x^A} f 
+ 
\frac{1}{2} \frac{\geop{x^A} \timefunctionarg{\muxmulevelsetvalue}}{\geop{t} \timefunctionarg{\muxmulevelsetvalue}} (\Lunit f)^2 
+  
\frac{\geop{x^A} \timefunctionarg{\muxmulevelsetvalue}}{\geop{t} \timefunctionarg{\muxmulevelsetvalue}} (\Lunit f) Xf 
- 
\frac{1}{2}  
\frac{\geop{x^A} \timefunctionarg{\muxmulevelsetvalue}}{\geop{t} \timefunctionarg{\muxmulevelsetvalue}} 
| \angD f|_{\gtorus}^2. 
\label{E:QTRANSPORTROUGHGEOOPA}
\end{align}
\end{subequations}
 
Moreover, the following identities also hold, where $\multipliervectorfield$ is the multiplier vectorfield \eqref{E:MULTIPLIERVECTORFIELD}, 
$\phi$ is the cut-off function from Def.\,\ref{D:WTRANSANDCUTOFF},
$\Rtransarg{\muxmulevelsetvalue}$ is as in \eqref{E:RTRANS},
$\hypunitnormalarg{\muxmulevelsetvalue}$ is as in \eqref{E:UNITHYPNORMALDEF},
$|\Rtransarg{\muxmulevelsetvalue}|_{\gfour}$ is as in \eqref{E:SIZEOFRTRANS},
and
$\Rtransnormsmallfactorarg{\muxmulevelsetvalue}$ is defined by \eqref{E:RTRANSNORMSMALLFACTOR}: 
\begin{subequations} 
\begin{align}
	\enmomem[f](\multipliervectorfield,\Lunit) 
	& 
	= (1 + \upmu)(\Lunit f)^2 + \upmu |\angD f|_{\gtorus}^2, 
		\label{E:QTL} 
		\\
	\enmomem[f](\hypunitnormalarg{\muxmulevelsetvalue},\Lunit) 
	& 
	= 
	\frac{\upmu(1 - \Rtransnormsmallfactorarg{\muxmulevelsetvalue}) 
	- 
	\frac{\muxmulevelsetvalue\phi }{\Lunit \upmu}}{|\Rtransarg{\muxmulevelsetvalue}|_{\gfour}} \enmomem[f](\Lunit,\Lunit) 
	+ 
	\frac{\upmu}{|\Rtransarg{\muxmulevelsetvalue}|_{\gfour}}\enmomem[f](X,\Lunit) 
	- 
	\frac{\upmu}{|\Rtransarg{\muxmulevelsetvalue}|_{\gfour}} 
	\gtorusroughinversefirstfund\left(dx^A,dx^B \right)
	\frac{\geop{x^A} \timefunctionarg{\muxmulevelsetvalue}}{\geop{t} \timefunctionarg{\muxmulevelsetvalue}} 
	\enmomem[f](\roughgeop{x^B},\Lunit), \label{E:QNL} 
		\\ 
	\enmomem[f](\hypunitnormalarg{\muxmulevelsetvalue}, \Transport) 
	& 
	= 
	\frac{\upmu(1 - \Rtransnormsmallfactorarg{\muxmulevelsetvalue}) 
	- 
	\frac{\muxmulevelsetvalue\phi}{\Lunit \upmu}}{|\Rtransarg{\muxmulevelsetvalue}|_{\gfour}} \enmomem[f](\Transport,\Lunit) 
	+ 
	\frac{\upmu}{|\Rtransarg{\muxmulevelsetvalue}|_{\gfour}}\enmomem[f](\Transport,X) 
	- 
	\frac{\upmu}{|\Rtransarg{\muxmulevelsetvalue}|_{\gfour}}
	\gtorusroughinversefirstfund\left(dx^A,dx^B \right)
	\frac{\geop{x^A} \timefunctionarg{\muxmulevelsetvalue}}{\geop{t} \timefunctionarg{\muxmulevelsetvalue}} 
	\enmomem[f](\roughgeop{x^B},\Transport), 
		\label{E:QNB} 
		\\
	\begin{split} \label{E:QNMULTIPLIER}
	\enmomem[f](\hypunitnormalarg{\muxmulevelsetvalue},\multipliervectorfield) 
	& 
	=  
	\frac{\upmu^2(1 - \Rtransnormsmallfactorarg{\muxmulevelsetvalue}) - \upmu\frac{\muxmulevelsetvalue\phi}{\Lunit \upmu}}{|\Rtransarg{\muxmulevelsetvalue}|_{\gfour}} 
	\left\lbrace 
		(\Lunit f)^2 
		+ 
		|\angD f|_{\gtorus}^2 
	\right\rbrace 
	+ 
	\frac{2 \upmu^2}{|\Rtransarg{\muxmulevelsetvalue}|_{\gfour}} \left\lbrace (\Lunit f) Xf + (Xf)^2 \right\rbrace 
		\\
	& \ \ 
	+ 
	\frac{\upmu^2}{|\Rtransarg{\muxmulevelsetvalue}|_{\gfour}} 
	\Big\lbrace 
		- 2 (\Lunit f) 
		\frac{\gtorusroughinversefirstfund\left(dx^A,dx^B \right) \geop{x^A} \timefunctionarg{\muxmulevelsetvalue}}
		{ \geop{t} \timefunctionarg{\muxmulevelsetvalue}} \roughgeop{x^B} f 	
		- 
		2 (Xf)
		\gtorusroughinversefirstfund\left(dx^A,dx^B \right)
		\frac{\geop{x^A} \timefunctionarg{\muxmulevelsetvalue}}
		{ \geop{t} \timefunctionarg{\muxmulevelsetvalue}} \roughgeop{x^B} f 	
				\\
	& \ \ \ \ \ \ \ \ \ \ \ \ \ \ \ 
		-
		\Rtransnormsmallfactorarg{\muxmulevelsetvalue} (\Lunit f)^2 
		- 
		2 \Rtransnormsmallfactorarg{\muxmulevelsetvalue} (\Lunit f) Xf 
		+ 
		\Rtransnormsmallfactorarg{\muxmulevelsetvalue} |\angD f|_{\gtorus}^2 
	\Big\rbrace 
		\\
	& \ \ 
	+ 
	\frac{\upmu 
	\left( \frac{1}{2} - \Rtransnormsmallfactorarg{\muxmulevelsetvalue} \right) 
	- 
	\frac{\muxmulevelsetvalue \phi}{\Lunit \upmu}}{|\Rtransarg{\muxmulevelsetvalue}|_{\gfour}} (\Lunit f)^2 
	+ 
	\frac{\upmu}{2 |\Rtransarg{\muxmulevelsetvalue}|_{\gfour}} |\angD f|_{\gtorus}^2 
	- 
	\frac{\upmu}{|\Rtransarg{\muxmulevelsetvalue}|_{\gfour}} 
	\gtorusroughinversefirstfund\left(dx^A,dx^B \right) 
	\frac{\geop{x^A} \timefunctionarg{\muxmulevelsetvalue}}
	{\geop{t} \timefunctionarg{\muxmulevelsetvalue}} (\Lunit f) \roughgeop{x^B} f. 
\end{split} 
\end{align}
\end{subequations}

\end{lemma}

\begin{proof}
\eqref{E:QLL}--\eqref{E:QXTRANSPORT} are straightforward consequences of the definition  
\eqref{E:DEFOFENERGYMOMENTUM} of $\enmomem$,
Lemma~\ref{L:BASICPROPERTIESOFVECTORFIELDS},
and the decompositions of $\gfour^{-1}$ implied by \eqref{E:SMOOTHTORUSINVERSEMETRICINTERMSOFINVERSESIGMATMETRICANDX}. 

\eqref{E:QLROUGHGEOPA}--\eqref{E:QTRANSPORTROUGHGEOOPA} follow from combining similar
arguments with the identity
\eqref{E:ROUGHANGULARPARTIALDERIVATIVESINTERMSOFGOODGEOMETRICPARTIALDERIVATIVES}.

\eqref{E:QTL} follows from the identity
$\multipliervectorfield = \Lunit + 2 \upmu \Transport$ 
(which is a consequence of \eqref{E:BISLPLUSX} and \eqref{E:MULTIPLIERVECTORFIELD}), 
\eqref{E:QLL}, 
and \eqref{E:QTRANSPORTL}. 

To derive \eqref{E:QNL}--\eqref{E:QNB}, 
we first decompose $\hypunitnormalarg{\muxmulevelsetvalue}$ as follows using 
\eqref{E:DEFROUGHTORITANGENTVECTORFIELD},
\eqref{E:RTRANS}, 
\eqref{E:UNITLENGTHRTRANS},
\eqref{E:SIZEOFRTRANS},
\eqref{E:RTRANSNORMSMALLFACTOR}, 
and \eqref{E:HYPUNITNORMALDECOMPOSITION}:
\begin{align} 
\begin{split} \label{E:USEFULHYPUNITNORMALDECOMPOSITION}
\hypunitnormalarg{\muxmulevelsetvalue} 
& 
= \frac{\upmu(1 - \Rtransnormsmallfactorarg{\muxmulevelsetvalue}) - \frac{2\muxmulevelsetvalue \phi}{\Lunit \upmu}}{|\Rtransarg{\muxmulevelsetvalue}|_{\gfour}}\Lunit 
+ 
\frac{1}{|\Rtransarg{\muxmulevelsetvalue}|_{\gfour}} \muX 
+ 
\frac{\muxmulevelsetvalue \phi}{|\Rtransarg{\muxmulevelsetvalue}|_{\gfour} \Lunit \upmu} \Lunit 
- 
\frac{\upmu}{|\Rtransarg{\muxmulevelsetvalue}|_{\gfour}}
\gtorusroughinversefirstfund\left(dx^A,dx^B \right)
\frac{\geop{x^A} \timefunctionarg{\muxmulevelsetvalue} \roughgeop{x^B}}{\geop{t} \timefunctionarg{\muxmulevelsetvalue}} 
	\\
& =   
\frac{\upmu(1 - \Rtransnormsmallfactorarg{\muxmulevelsetvalue}) - \frac{\muxmulevelsetvalue \phi}{\Lunit \upmu}}{|\Rtransarg{\muxmulevelsetvalue}|_{\gfour}}\Lunit 
+ 
\frac{\upmu}{|\Rtransarg{\muxmulevelsetvalue}|_{\gfour}} X  
- 
\frac{\upmu}{|\Rtransarg{\muxmulevelsetvalue}|_{\gfour}}
\gtorusroughinversefirstfund\left(dx^A,dx^B \right) 
\frac{\geop{x^A} \timefunctionarg{\muxmulevelsetvalue} \roughgeop{x^B}}{\geop{t} \timefunctionarg{\muxmulevelsetvalue}}.
\end{split}
\end{align}
\eqref{E:QNL}--\eqref{E:QNB} now follow from
\eqref{E:USEFULHYPUNITNORMALDECOMPOSITION}
and
the linearity of the map $Z \rightarrow \enmomem[f](\hypunitnormalarg{\muxmulevelsetvalue},Z)$.

Finally, \eqref{E:QNMULTIPLIER} follows from 
\eqref{E:QLL}--\eqref{E:QTRANSPORTROUGHGEOOPA},
\eqref{E:QNL}--\eqref{E:QNB} 
and the following identity identity,
which follows from the identity $\multipliervectorfield = \Lunit + 2 \upmu \Transport$  mentioned above:
$\enmomem[f](\hypunitnormalarg{\muxmulevelsetvalue},\multipliervectorfield) = \enmomem[f](\hypunitnormalarg{\muxmulevelsetvalue},\Lunit) + 2 \upmu \enmomem[f](\hypunitnormalarg{\muxmulevelsetvalue},\Transport)$.
\end{proof}

\subsubsection{Coerciveness estimates on submanifolds}
\label{SSS:COERCIVENESSOFL2CONTROLLINGONSUBMANIFOLDS}
In the next lemma,
we exhibit the coerciveness of
the $L^2$-controlling quantities from Def.\,\ref{D:MAINCOERCIVE}.
In particular, we exhibit their $L^2$ coerciveness properties on
the submanifolds
$\hypthreearg{\timefunction}{[- \rightu,u]}{\muxmulevelsetvalue}$,
$\nullhypthreearg{\muxmulevelsetvalue}{u}{[\timefunction_0,\timefunction)}$,
and
$\twoargroughtori{\timefunction,u}{\muxmulevelsetvalue}$.
We reveal the coerciveness of the spacetime integrals
 $\spacetimeintegralcontrolwave_N(\timefunction,u)$ 
and
$\spacetimeintegralcontrolwavepartial_N(\timefunction,u)$
in a separate lemma, namely Lemma~\ref{L:COERCIVENESSOFSPACETIMEINTEGRAL}.

\begin{lemma}[The coerciveness on submanifolds of the fundamental $L^2$-controlling quantities] 
	\label{L:COERCIVENESSOFL2CONTROLLINGQUANITIES} 
Let $\hypersurfacecontrolwave_N(\timefunction,u)$,
$\cdots$,
$\hypersurfacecontrolDivGradEnt_N(\timefunction,u)$
be the $L^2$-controlling quantities from Def.\,\ref{D:MAINCOERCIVE},
and let $(\timefunction,u) \in [\timefunction_0,\timefunctionboot)\times [- \rightu,\leftu]$.
Then for $1 \leq N \leq \Ntop$,
the following lower bounds hold,
where $\mathfrak{P}^{(N)}$ is the set of order $N$ $\nullhyparg{u}$-tangential commutator operators from
Sect.\,\ref{SS:STRINGSOFCOMMUTATIONVECTORFIELDS}:
\begin{align}  
\begin{split} \label{E:COERCIVENESSOFHYPERSURFACECONTROLWAVE} 
	\hypersurfacecontrolwave_N(\timefunction,u) 
	& \geq 
	\max_{\substack{\tander^N \in \mathfrak{P}^{(N)} \\ \Psi \in \{ \RRiemann,\LRiemann,v^2,v^3,\Ent \}}}
	\left\lbrace 
		0.49 
		\left\| 
			\sqrt{\upmu - \frac{2 \muxmulevelsetvalue \phi}{\Lunit \upmu}} \Lunit \tander^N \Psi 
		\right\|_{L^2\left(\hypthreearg{\timefunction}{[- \rightu,u]}{\muxmulevelsetvalue}\right)}^2, 
			\,
		0.49 
		\left\| 
			\sqrt{\upmu} \left| \angrmD \tander^N \Psi \right|_{\gtorus} 
		\right\|_{L^2\left(\hypthreearg{\timefunction}{[- \rightu,u]}{\muxmulevelsetvalue}\right)}^2, 
		\right. 
		\\
	& 
	\left. 
	0.99 
	\left\| 
		\muX \tander^N \Psi 
	\right\|_{L^2\left(\hypthreearg{\timefunction}{[- \rightu,u]}{\muxmulevelsetvalue}\right)}^2, 
		\,
	\left\| 
		\frac{1}{\sqrt{\Lunit \timefunctionarg{\muxmulevelsetvalue}}} \Lunit \tander^N \Psi 
	\right\|_{L^2(\nullhypthreearg{\muxmulevelsetvalue}{u}{[\timefunction_0,\timefunction)})}^2, 
		\,
	\left\| 
		\frac{\sqrt{\upmu}}{\sqrt{\Lunit \timefunctionarg{\muxmulevelsetvalue}}}  
		\left| \angrmD\tander^N \Psi \right|_{\gtorus} 
	\right\|_{L^2(\nullhypthreearg{\muxmulevelsetvalue}{u}{[\timefunction_0,\timefunction)})}^2 
	\right\rbrace,
\end{split}		
		\\
\begin{split} 	\label{E:COERCIVENESSOFHYPERSURFACECONTROLWAVEPARTIAL}
\hypersurfacecontrolwavepartial_N(\timefunction,u) 
	& \geq 
	\max_{\substack{\tander^N \in \mathfrak{P}^{(N)} \\ \Psi \in \{ \LRiemann,v^2,v^3,\Ent \}}}
	\left\lbrace
		0.49 
		\left\| 
			\sqrt{\upmu - \frac{2 \muxmulevelsetvalue \phi}{\Lunit \upmu}} \Lunit \tander^N \Psi 
		\right\|_{L^2\left(\hypthreearg{\timefunction}{[- \rightu,u]}{\muxmulevelsetvalue}\right)}^2, 
			\,
		0.49 
		\left\| 
			\sqrt{\upmu} \left| \angrmD \tander^N \Psi \right|_{\gtorus} 
		\right\|_{L^2\left(\hypthreearg{\timefunction}{[- \rightu,u]}{\muxmulevelsetvalue}\right)}^2, 
	\right. 
			\\
	& \left. 
			0.99 
			\left\| 
				\muX \tander^N \Psi 
			\right\|_{L^2\left(\hypthreearg{\timefunction}{[- \rightu,u]}{\muxmulevelsetvalue}\right)}^2, 
				\,
			\left\| 
				\frac{1}{\sqrt{\Lunit \timefunctionarg{\muxmulevelsetvalue}}} \Lunit \tander^N \Psi 
			\right\|_{L^2(\nullhypthreearg{\muxmulevelsetvalue}{u}{[\timefunction_0,\timefunction)})}^2, 
				\,
			\left\| 
				\frac{\sqrt{\upmu}}{\sqrt{\Lunit \timefunctionarg{\muxmulevelsetvalue}}}  
				\left| \angrmD\tander^N \Psi\right|_{\gtorus} 
			\right\|_{L^2(\nullhypthreearg{\muxmulevelsetvalue}{u}{[\timefunction_0,\timefunction)})}^2 
		\right\rbrace.   
\end{split}
\end{align}

Moreover, for $N \leq \Ntop$, the following lower bounds hold:
\begin{subequations}
\begin{align}
\hypersurfacecontrolVort_N(\timefunction,u) 
& 
\geq 
\max 
\left\lbrace 
	\left\| 
		\sqrt{\upmu - \frac{\muxmulevelsetvalue \phi}{\Lunit \upmu}} \tander^N \vortrenormalized 
	\right\|_{L^2\left(\hypthreearg{\timefunction}{[- \rightu,u]}{\muxmulevelsetvalue}\right)}^2, 
		\,
	\left\| 
		\frac{1}{\sqrt{\Lunit \timefunctionarg{\muxmulevelsetvalue}}} \tander^N \vortrenormalized 
	\right\|_{L^2(\nullhypthreearg{\muxmulevelsetvalue}{u}{[\timefunction_0,\timefunction)})}^2, 
\right\rbrace \label{E:COERCIVENESSOFCONTROLVORT} 
		\\
\hypersurfacecontrolGradEnt_N(\timefunction,u) 
& \geq 
\max 
\left\lbrace 
	\left\| 
		\sqrt{\upmu - \frac{\muxmulevelsetvalue \phi}{\Lunit \upmu}} \tander^N \GradEnt 
	\right\|_{L^2\left(\hypthreearg{\timefunction}{[- \rightu,u]}{\muxmulevelsetvalue}\right)}^2, 
		\,
	\left\| 
		\frac{1}{\sqrt{\Lunit \timefunctionarg{\muxmulevelsetvalue}}} \tander^N \GradEnt 
	\right\|_{L^2(\nullhypthreearg{\muxmulevelsetvalue}{u}{[\timefunction_0,\timefunction)})}^2 
\right\rbrace, 
		\label{E:COERCIVENESSOFCONTROLGRADENT} 
\end{align}
\end{subequations}

\begin{subequations}
\begin{align}
\hypersurfacecontrolVortVort_N (\timefunction,u) 
& \geq 
\max 
\left\lbrace 
\left\| 
	\sqrt{\upmu - \frac{\muxmulevelsetvalue \phi}{\Lunit \upmu}} \tander^N \VortVort 
\right\|_{L^2\left(\hypthreearg{\timefunction}{[- \rightu,u]}{\muxmulevelsetvalue}\right)}^2,
	\,
 \left\| 
	\frac{1}{\sqrt{\Lunit \timefunctionarg{\muxmulevelsetvalue}}} \tander^N \VortVort 
\right\|_{L^2(\nullhypthreearg{\muxmulevelsetvalue}{u}{[\timefunction_0,\timefunction)})}^2 
\right\rbrace \label{E:COERCIVENESSHYPERSURFACECONTROLVORTVORT} 
		\\
\hypersurfacecontrolDivGradEnt_N (\timefunction,u) 
& \geq 
\max 
\left\lbrace 
	\left\| 
		\sqrt{\upmu - \frac{\muxmulevelsetvalue \phi}{\Lunit \upmu}} \tander^N \DivGradEnt 
\right\|_{L^2\left(\hypthreearg{\timefunction}{[- \rightu,u]}{\muxmulevelsetvalue}\right)}^2, 
	\,
\left\| 
	\frac{1}{\sqrt{\Lunit \timefunctionarg{\muxmulevelsetvalue}}} \tander^N \DivGradEnt 
\right\|_{L^2(\nullhypthreearg{\muxmulevelsetvalue}{u}{[\timefunction_0,\timefunction)})}^2 \right\rbrace. 
\label{E:COERCIVENESSHYPERSURFACEDIVGRADENT}
\end{align}
\end{subequations}

In addition, for $1 \leq N \leq \Ntop$ and $N' \leq \Ntop$, 
we have: 
	\begin{align} 
	\left\| \tander^N \Psi \right\|_{L^2(\twoargroughtori{\timefunction,u}{\muxmulevelsetvalue})}^2 
	& 
	\leq 
	C \initialsmall^2 
	+ 
	C \hypersurfacecontrolwave_N(\timefunction,u),  
		\label{E:LOSSOFONEDERIVATIVEL2ESTIMATESFORWAVEVARIABLESONROUGHTORIINTERMSOFDATAANDCONTROLLING} 
		\\
	\left\| \tander^{\leq N'}  \vortrenormalized \right\|_{L^2(\twoargroughtori{\timefunction,u}{\muxmulevelsetvalue})}^2 
	& \leq C \initialsmall^2 + C \hypersurfacecontrolVort_{\leq N'+1}(\timefunction,u), 
	\label{E:LOSSOFONEDERIVATIVEL2ESTIMATESFORTRANSPORTVARIABLESONROUGHTORIINTERMSOFDATAANDCONTROLLING}
		\\
	\left\| \tander^{\leq N'}  \GradEnt \right\|_{L^2(\twoargroughtori{\timefunction,u}{\muxmulevelsetvalue})}^2 
	& \leq C \initialsmall^2 
	+ 
	C \hypersurfacecontrolGradEnt_{\leq N'+1}(\timefunction,u), 
		\label{E:LOSSOFONEDERIVATIVETORIL2ESTIMATESFORGRADENTONROUGHTORIINTERMSOFDATAANDCONTROLLING}
			\\
\left\| \tander^{N'} \Psi \right\|_{L^2\left(\hypthreearg{\timefunction}{[- \rightu,u]}{\muxmulevelsetvalue}\right)}
& \leq 
	C \initialsmall
	+
	C
	\int_{\timefunction' = \timefunction_0}^{\timefunction} 
		\frac{\hypersurfacecontrolwave_{[1,N]}^{1/2}(\timefunction',u)}{|\timefunction'|^{1/2}}
	\, \mathrm{d} \timefunction'
	\leq
	C \initialsmall
	+
	C
	\hypersurfacecontrolwave_{[1,N]}^{1/2}(\timefunction,u).
	\label{E:L2ESTIMATESFORWAVEVARIABLESONROUGHHYPERSURFACELOSSOFONEDERIVATIVE}
	\end{align}

Finally, if $1 \leq N \leq \Ntop$,
$\Psi \in \{ \RRiemann, \LRiemann, v^2, v^3, \Ent\}$,
$\tander^N \in \mathfrak{P}^{(N)}$,
and $\fundbootsmall$ and $\mupositive$ are sufficiently small, 
then the following sharpened coerciveness estimates
hold whenever $0 \leq \upmu < \mupositive$:
\begin{align} \label{E:L2WAVECONTROLWITHBETTERCONSTANTIFMUSMALL}
\hypersurfacecontrolwave_N(\timefunction,u) 
& 
\geq
\mathbb{E}_{(\textnormal{Wave})}[\tander^N \Psi] (\timefunction,u) 
\geq 
1.99 
\left\| \muX \tander^N \Psi \right\|_{L^2\left(\hypthreearg{\timefunction}{[- \rightu,u]}{\muxmulevelsetvalue}\right)}^2.
\end{align}
\end{lemma}

\begin{proof} \hfill

\noindent \underline{\textbf{Proofs of \eqref{E:COERCIVENESSOFCONTROLVORT}--\eqref{E:COERCIVENESSHYPERSURFACEDIVGRADENT}}}:
These estimates follow directly from the definitions~\eqref{E:TRANSPORTENERGYDEF}--\eqref{E:TRANSPORTNULLFLUXDEF},
\eqref{E:VORTICITYL2CONTROLLINGQUANTITY}--\eqref{E:ENTROPYGRADIENTL2CONTROLLINGQUANTITY},
and \eqref{E:MODIFIEDVORTICITYVORTICITYL2CONTROLLINGQUANTITY}--\eqref{E:MODIFIEDDIVGRADENTL2CONTROLLINGQUANTITY}.

\medskip
\noindent 
\underline{\textbf{Proofs of \eqref{E:COERCIVENESSOFHYPERSURFACECONTROLWAVE}--\eqref{E:COERCIVENESSOFHYPERSURFACECONTROLWAVEPARTIAL}}}:
We fix any $(\timefunction,u) \in [\timefunction_0,\timefunctionboot)\times [- \rightu,\leftu]$,
$\Psi \in \{ \RRiemann,\LRiemann,v^2,v^3,\Ent \}$, 
and
$\tander^N \in \mathfrak{P}^{(N)}$.
First, using \eqref{E:WAVENULLFLUXDEF}
and \eqref{E:QTL} with $f \eqdef \tander^N \Psi$,
\eqref{E:GEOMETRICL2NORMSTORIANDNULLHYPERSURFACES},
and \eqref{E:WAVESPACELIKEANDNULLHYPERSURFACEL2CONTROLLINGQUANTITY},
we find that
$
\left\| 
	\frac{1}{\sqrt{\Lunit \timefunctionarg{\muxmulevelsetvalue}}} \Lunit \tander^N \Psi 
\right\|_{L^2(\nullhypthreearg{\muxmulevelsetvalue}{u}{[\timefunction_0,\timefunction)})}^2
\leq
\hypersurfacecontrolwave_N(\timefunction,u)
$
and
$
\left\| \frac{\sqrt{\upmu}}{\sqrt{\Lunit \timefunctionarg{\muxmulevelsetvalue}}}  
\left| \angrmD\tander^N \Psi\right|_{\gtorus} \right\|_{L^2(\nullhypthreearg{\muxmulevelsetvalue}{u}{[\timefunction_0,\timefunction)})}^2 
\leq
\hypersurfacecontrolwave_N(\timefunction,u)
$
as desired.
Next, we note that the product of $|\Rtransarg{\muxmulevelsetvalue}|_{\gfour}$ and the terms on first line of RHS~\eqref{E:QNMULTIPLIER} 
can be expressed as follows:
\begin{align} 
\begin{split} \label{E:COERCIVESTEP1} 
\left\lbrace 
	\upmu^2\left(\frac{0.001}{1.001} - \Rtransnormsmallfactorarg{\muxmulevelsetvalue}\right) 
	- 
	\upmu \frac{\muxmulevelsetvalue \phi}{\Lunit \upmu}
\right\rbrace 
(\Lunit \tander^N \Psi)^2 
+ 
\left\lbrace 
	\upmu^2 (1-\Rtransnormsmallfactorarg{\muxmulevelsetvalue}) 
	- 
	\upmu \frac{\muxmulevelsetvalue \phi}{\Lunit \upmu}
\right\rbrace |\angD \tander^N \Psi|_{\gtorus}^2  \\
+ 
\left( \frac{\upmu}{\sqrt{1.001}} \Lunit \tander^N \Psi + \sqrt{1.001} \muX \tander^N \Psi\right)^2 
+ 
0.999(\muX \Psi)^2.
\end{split}
\end{align}
Next, we multiply the terms on the second and third lines of RHS~\eqref{E:QNMULTIPLIER} by $|\Rtransarg{\muxmulevelsetvalue}|_{\gfour}$ 
and use 
\eqref{E:SMOOTHTORIGABEXPRESSION},
\eqref{E:GEOP2TOCOMMUTATORS}--\eqref{E:GEOP3TOCOMMUTATORS},
\eqref{E:ROUGHANGULARPARTIALDERIVATIVESINTERMSOFGOODGEOMETRICPARTIALDERIVATIVES},
\eqref{E:GTORUSROUGHCOMPONENTS}--\eqref{E:ROUGHTORUSMETRICCOMPONENTSANDTHEINVERSECOMPONENTSRELATION},
Lemma~\ref{L:SCHEMATICSTRUCTUREOFVARIOUSTENSORSINTERMSOFCONTROLVARS},
\eqref{E:ANGDFPOINTWISEBOUNDEDBYCOMMUTATORVECTORFIELDS},
\eqref{E:CLOSEDVERSIONKEYJACOBIANDETERMINANTESTIMATECHOVGEOTOROUGH},
\eqref{E:SMALLC11ESTIMATESFORROUGHTIMEFUNCTION},
the estimates of Prop.\,\ref{P:IMPROVEMENTOFAUXILIARYBOOTSTRAP},
Cor.\,\ref{C:IMPROVEAUX},
\eqref{E:POINTWISEBOUNDTANGENTIALANDTRANSVERSALDERIVATIVESOFRTRANSNORMSMALLFACTOR},
and Young's inequality to bound the magnitude of
the resulting terms as follows:
\begin{align} \label{E:ERRORTERMSQNT} 
& 
\leq
C \fundbootsmall \upmu^2 (\Lunit \tander^N \Psi)^2 
+
C \fundbootsmall \upmu^2 |\angD \tander^N \Psi|_{\gtorus}^2 
+
C \fundbootsmall (\muX \tander^N \Psi)^2.
\end{align}
The same reasoning, in conjunction with \eqref{E:BOUNDSONLMUINTERESTINGREGION}, 
yields that all the terms in \eqref{E:COERCIVESTEP1} are positive definite
and that the product of
$|\Rtransarg{\muxmulevelsetvalue}|_{\gfour}$ and the last term on RHS~\eqref{E:QNMULTIPLIER} can be absorbed by 
the product of $|\Rtransarg{\muxmulevelsetvalue}|_{\gfour}$ and the two terms
$
\frac{\upmu 
	\left( \frac{1}{2} - \Rtransnormsmallfactorarg{\muxmulevelsetvalue} \right) 
	- 
	\frac{\muxmulevelsetvalue \phi}{\Lunit \upmu}}{|\Rtransarg{\muxmulevelsetvalue}|_{\gfour}} (\Lunit f)^2 
	+ 
	\frac{\upmu}{2 |\Rtransarg{\muxmulevelsetvalue}|_{\gfour}} |\angD f|_{\gtorus}^2 
$
on the last line of \eqref{E:QNMULTIPLIER}.
In total, we see that if $\fundbootsmall$ is small enough, 
then in the product of $|\Rtransarg{\muxmulevelsetvalue}|_{\gfour}$ and \eqref{E:QNMULTIPLIER},
the overall coefficient of 
$(\muX \tander^N \Psi)^2$
can be bounded from below by $0.99$,
the overall coefficient of 
$\upmu (\Lunit \tander^N \Psi)^2$ can be bounded from below by $0.49 \left(\upmu - \frac{2 \muxmulevelsetvalue \phi}{\Lunit \upmu} \right)$,
and the overall coefficient of $\upmu |\angD \tander^N \Psi|_{\gtorus}^2$ 
can be bounded from below by $0.49$.
From these estimates,
\eqref{E:WAVEENERGYDEF} with $f \eqdef \tander^N \Psi$,
\eqref{E:GEOMETRICL2NORMSROUGHHYPERSURFACESANDSPACETIMEREGIONS},
and \eqref{E:WAVESPACELIKEANDNULLHYPERSURFACEL2CONTROLLINGQUANTITY},
we deduce that
$
0.99 
	\left\| 
		\muX \tander^N \Psi 
	\right\|_{L^2\left(\hypthreearg{\timefunction}{[- \rightu,u]}{\muxmulevelsetvalue}\right)}^2
\leq
\hypersurfacecontrolwave_N(\timefunction,u)
$,
$
0.49 
\left\| 
	\sqrt{\upmu - \frac{2 \muxmulevelsetvalue \phi}{\Lunit \upmu}} \Lunit \tander^N \Psi 
\right\|_{L^2\left(\hypthreearg{\timefunction}{[- \rightu,u]}{\muxmulevelsetvalue}\right)}^2
\leq
\hypersurfacecontrolwave_N(\timefunction,u)
 $
and
$
0.49 
		\left\| 
			\sqrt{\upmu} \left| \angrmD \tander^N \Psi \right|_{\gtorus} 
		\right\|_{L^2\left(\hypthreearg{\timefunction}{[- \rightu,u]}{\muxmulevelsetvalue}\right)}^2
\leq
\hypersurfacecontrolwave_N(\timefunction,u)
$.
Combining all five of the $L^2$ bounds along $\nullhypthreearg{\muxmulevelsetvalue}{u}{[\timefunction_0,\timefunction)}$
and $\hypthreearg{\timefunction}{[- \rightu,u]}{\muxmulevelsetvalue}$
that we derived in this paragraph, 
we conclude the desired lower bounds \eqref{E:COERCIVENESSOFHYPERSURFACECONTROLWAVE}.

Taking into account definition~\eqref{E:WAVEPARTIALSPACELIKEANDNULLHYPERSURFACEL2CONTROLLINGQUANTITY}, 
we can prove the lower bounds stated in \eqref{E:COERCIVENESSOFHYPERSURFACECONTROLWAVEPARTIAL}
via exactly the same argument.

\medskip
\noindent \underline{\textbf{Proof of \eqref{E:L2WAVECONTROLWITHBETTERCONSTANTIFMUSMALL}}}:
We again fix any $(\timefunction,u) \in [\timefunction_0,\timefunctionboot)\times [- \rightu,\leftu]$,
$\Psi \in \{ \RRiemann,\LRiemann,v^2,v^3,\Ent \}$, 
and
$\tander^N \in \mathfrak{P}^{(N)}$.
We consider the product of $|\Rtransarg{\muxmulevelsetvalue}|_{\gfour}$ and RHS~\eqref{E:QNMULTIPLIER}
with $f \eqdef \tander^N \Psi$.
The terms in the second braces on the first line of RHS~\eqref{E:QNMULTIPLIER} generate the terms
$2 \upmu (\Lunit \tander^N \Psi) \muX \tander^N \Psi + 2 (\muX \tander^N \Psi)^2$,
which by Young's inequality, can be pointwise bounded from below by
$1.999 (\muX \tander^N \Psi)^2 - 1000 \upmu^2 (\Lunit \tander^N \Psi)^2$.
If $\mupositive \leq 10^{-5}$,
then \eqref{E:MINVALUEOFMUONFOLIATION} implies
that on $\twoargMrough{[\timefunction_0,\timefunctionboot),[- \rightu,\leftu]}{\muxmulevelsetvalue}$,
we have the pointwise bound
$|1000 \upmu^2 (\Lunit \tander^N \Psi)^2| \leq \frac{1}{10} \upmu (\Lunit \tander^N \Psi)^2$.
Hence, if $\fundbootsmall$ is sufficiently small, 
then the same arguments we used in the proof of \eqref{E:COERCIVENESSOFHYPERSURFACECONTROLWAVE}
imply that we can absorb the term
$- 1000 \upmu^2 (\Lunit \tander^N \Psi)^2$ into
the terms 
$
\frac{\upmu 
	\left( \frac{1}{2} - \Rtransnormsmallfactorarg{\muxmulevelsetvalue} \right) 
	- 
	\frac{\muxmulevelsetvalue \phi}{\Lunit \upmu}}{|\Rtransarg{\muxmulevelsetvalue}|_{\gfour}} 
	(\Lunit \tander^N \Psi)^2 
$
on RHS~\eqref{E:QNMULTIPLIER}.
The arguments we used in the proof of \eqref{E:COERCIVENESSOFHYPERSURFACECONTROLWAVE}
also imply that all remaining terms on RHS~\eqref{E:QNMULTIPLIER}
are either positive definite or can be absorbed into the positive definite terms by
exploiting the smallness of $\fundbootsmall$.
In total, these arguments yield the pointwise estimate
$
|\Rtransarg{\muxmulevelsetvalue}|_{\gfour} \times \mbox{RHS~\eqref{E:QNMULTIPLIER}} \geq 1.99 (\muX \tander^N \Psi)^2
$.
From this estimate and the same arguments we used to prove \eqref{E:COERCIVENESSOFHYPERSURFACECONTROLWAVE},
we conclude the desired lower bound \eqref{E:L2WAVECONTROLWITHBETTERCONSTANTIFMUSMALL}.

\medskip
\noindent \underline{\textbf{Proofs of \eqref{E:LOSSOFONEDERIVATIVEL2ESTIMATESFORWAVEVARIABLESONROUGHTORIINTERMSOFDATAANDCONTROLLING}--\eqref{E:LOSSOFONEDERIVATIVETORIL2ESTIMATESFORGRADENTONROUGHTORIINTERMSOFDATAANDCONTROLLING}}}:
To prove \eqref{E:LOSSOFONEDERIVATIVEL2ESTIMATESFORWAVEVARIABLESONROUGHTORIINTERMSOFDATAANDCONTROLLING},
we first use \eqref{E:ROUGHTORUSINNULLHYPERSURFACEL2FUNDAMENTALTHEOREMOFCALCULUSESTIMATE} with 
$f \eqdef \tander^N \Psi$ 
as well as the already proved coerciveness result \eqref{E:COERCIVENESSOFHYPERSURFACECONTROLWAVE} 
to deduce that
$\left\| \tander^N \Psi \right\|_{L^2(\twoargroughtori{\timefunction,u}{\muxmulevelsetvalue})}^2  
\leq 
C
\left\| \tander^N \Psi \right\|_{L^2(\twoargroughtori{\timefunction_0,u}{\muxmulevelsetvalue})}^2 
+
C
\hypersurfacecontrolwave_{N}(\timefunction,u)$. 
From this estimate and the data estimate \eqref{E:SMALLDATAOFPSIONINITIALROUGHTORI},
we conclude \eqref{E:LOSSOFONEDERIVATIVEL2ESTIMATESFORWAVEVARIABLESONROUGHTORIINTERMSOFDATAANDCONTROLLING}.
The estimates 
\eqref{E:LOSSOFONEDERIVATIVETORIL2ESTIMATESFORGRADENTONROUGHTORIINTERMSOFDATAANDCONTROLLING}--\eqref{E:LOSSOFONEDERIVATIVETORIL2ESTIMATESFORGRADENTONROUGHTORIINTERMSOFDATAANDCONTROLLING} follow from a nearly identical argument based on the coerciveness results
\eqref{E:COERCIVENESSOFCONTROLVORT}--\eqref{E:COERCIVENESSOFCONTROLGRADENT}
and the data estimates \eqref{E:SMALLDATAOFVORTICITYANDENTORPYGRADIENTONINITIALROUGHTORI}.

\medskip
\noindent \underline{\textbf{Proof of \eqref{E:L2ESTIMATESFORWAVEVARIABLESONROUGHHYPERSURFACELOSSOFONEDERIVATIVE}}}:
We first use
\eqref{E:TANGENTIALL2NORMSOFWAVEVARIABLESSMALLALONGINITIALROUGHHYPERSURFACE}
and
\eqref{E:L2ONROUGHCONSTANTTIMEHYPERSURFACESTRANSPORTROUGHLFESTIMATE} 
to deduce that:
\begin{align} \label{E:FIRSTPROOFSTEPL2ESTIMATESFORWAVEVARIABLESONROUGHHYPERSURFACELOSSOFONEDERIVATIVE}
\left\| \tander^{N'} \Psi \right\|_{L^2\left(\hypthreearg{\timefunction}{[- \rightu,u]}{\muxmulevelsetvalue}\right)}
	&
	\lesssim 
	\initialsmall
	+
	\int_{\timefunction' = \timefunction_0}^{\timefunction} 
		\left\| \argLrough{\muxmulevelsetvalue} \tander^{N'} \Psi \right\|_{L^2\left(\hypthreearg{\timefunction'}{[- \rightu,u]}{\muxmulevelsetvalue}\right)}
	\, \mathrm{d} \timefunction.
\end{align}
Then, 
using \eqref{E:LROUGH},
\eqref{E:CLOSEDVERSIONLUNITROUGHTTIMEFUNCTION},
\eqref{E:MINVALUEOFMUONFOLIATION}, 
and
\eqref{E:COERCIVENESSOFHYPERSURFACECONTROLWAVE},
we bound the time integral on RHS~\eqref{E:FIRSTPROOFSTEPL2ESTIMATESFORWAVEVARIABLESONROUGHHYPERSURFACELOSSOFONEDERIVATIVE} by
$
\lesssim
\int_{\timefunction' = \timefunction_0}^{\timefunction} 
	\frac{\hypersurfacecontrolwave_{N'}^{1/2}(\timefunction',u)}{|\timefunction'|^{1/2}}
\, \mathrm{d} \timefunction
$,
which yields the first inequality stated in \eqref{E:L2ESTIMATESFORWAVEVARIABLESONROUGHHYPERSURFACELOSSOFONEDERIVATIVE}.
The second inequality stated in \eqref{E:L2ESTIMATESFORWAVEVARIABLESONROUGHHYPERSURFACELOSSOFONEDERIVATIVE}
follows from the fact that $\hypersurfacecontrolwave_N(\timefunction,u)$ is increasing in its arguments.

\end{proof}

\subsubsection{Coerciveness of the spacetime integrals $\spacetimeintegralcontrolwave$}
\label{SSS:COERCIVENESSOFSPACETIMEINTEGRALS}
The next lemma complements Lemma~\ref{L:COERCIVENESSOFL2CONTROLLINGQUANITIES} 
by exhibiting the coerciveness of the spacetime integrals 
$\spacetimeintegralcontrolwave_N(\timefunction,u)$ 
and
$\spacetimeintegralcontrolwavepartial_N(\timefunction,u)$  
defined in
\eqref{E:WAVESPACETIMEL2CONTROLLINGQUANTITY} and \eqref{E:PARTIALWAVESPACETIMEL2CONTROLLINGQUANTITY}
respectively; recall that these spacetime integrals
appear on the left-hand side of our energy identity \eqref{E:FUNDAMENTALENERGYINTEGRALIDENTITCOVARIANTWAVES}
for the wave variables.
The key point is that the integrands on
RHSs~\eqref{E:COERCIVENESSOFSPACETIMEINTEGRAL}--\eqref{E:COERCIVENESSOFPARTIALSPACETIMEINTEGRAL}
are quantitatively positive in the region $\lbrace |u| \leq \interestingu \rbrace$
where the shock can form
and in particular, these integrands \emph{do not contain any degenerate factor of $\upmu$}.
We fundamentally need the non-degenerate coerciveness 
guaranteed by \eqref{E:COERCIVENESSOFSPACETIMEINTEGRAL}--\eqref{E:COERCIVENESSOFPARTIALSPACETIMEINTEGRAL}
in order to control some of the error integrals that arise in our energy estimates.

\begin{lemma}[The coerciveness of the spacetime integrals $\spacetimeintegralcontrolwave$] 
\label{L:COERCIVENESSOFSPACETIMEINTEGRAL}
Let $1 \leq N \leq \Ntop$ be an integer,
let $\spacetimeintegralcontrolwave_N(\timefunction,u)$ 
and
$\spacetimeintegralcontrolwavepartial_N(\timefunction,u)$
be the spacetime integrals defined in 
\eqref{E:WAVESPACETIMEL2CONTROLLINGQUANTITY} and \eqref{E:PARTIALWAVESPACETIMEL2CONTROLLINGQUANTITY},
and let $\mathfrak{P}^{(N)}$ be the set of order $N$ $\nullhyparg{u}$-tangential commutator operators from
Sect.\,\ref{SS:STRINGSOFCOMMUTATIONVECTORFIELDS}.
Then the following lower bounds hold for 
$(\timefunction,u) \in [\timefunction_0,\timefunctionboot)\times [-\rightu,\leftu]$,
where $\mathbf{1}_{[-\interestingu,\interestingu]} = \mathbf{1}_{[-\interestingu,\interestingu]}(u')$ 
denotes the characteristic function of the interval $[-\interestingu,\interestingu]$
and $\phi$ is the cut-off from Def.\,\ref{D:WTRANSANDCUTOFF}:
\begin{subequations}
\begin{align} \label{E:COERCIVENESSOFSPACETIMEINTEGRAL}
\spacetimeintegralcontrolwave_N(\timefunction,u)  
&
\geq 
\max_{\substack{\tander^N \in \mathfrak{P}^{(N)} \\ \Psi \in \{ \RRiemann,\LRiemann,v^2,v^3,\Ent \}}}
	\frac{1}{4} 
	\int_{\twoargMrough{[\timefunction_0,\timefunction),[- \rightu,u]}{\muxmulevelsetvalue}} 
	\left\lbrace 
		\mathbf{1}_{[-\interestingu,\interestingu]}(u')
		+
		4 
		\frac{1}{\Lunit \timefunctionarg{\muxmulevelsetvalue}}
		\muxmulevelsetvalue \phi
	\right\rbrace
	\left|\angrmD \tander^N \Psi \right|_{\gtorus}^2 
	\, \volMRoughCoordinates,
		\\
\spacetimeintegralcontrolwavepartial_N(\timefunction,u)
&
\geq 
\max_{\substack{\tander^N \in \mathfrak{P}^{(N)} \\ \Psi \in \{\LRiemann,v^2,v^3,\Ent \}}}
	\frac{1}{4} 
	\int_{\twoargMrough{[\timefunction_0,\timefunction),[- \rightu,u]}{\muxmulevelsetvalue}} 
	 \left\lbrace 
		\mathbf{1}_{[-\interestingu,\interestingu]}(u')
		+
		4 
		\frac{1}{\Lunit \timefunctionarg{\muxmulevelsetvalue}}
		\muxmulevelsetvalue \phi
	\right\rbrace
	\left|\angrmD \tander^N \Psi \right|_{\gtorus}^2 
	\, \volMRoughCoordinates.
	\label{E:COERCIVENESSOFPARTIALSPACETIMEINTEGRAL}
\end{align}
\end{subequations}

\end{lemma}

\begin{proof}
Let $f$ be a scalar function, and let $\spacetimeintegralcontrolwave[f](\timefunction,u)$
be the corresponding spacetime integral defined in \eqref{E:COERCIVESPACETIMEWAVENERGYINTEGRAL}.
Using the estimate \eqref{E:SIZEBOUNDSONWIDETILDELMUINTERESTINGREGION}, we
deduce that:
\begin{align} \label{E:PROOFCOERCIVENESSOFSPACETIMEINTEGRAL}
\spacetimeintegralcontrolwave[f](\timefunction,u) 
&
\geq 
\frac{1}{4} 
	\int_{\twoargMrough{[\timefunction_0,\timefunction),[- \rightu,u]}{\muxmulevelsetvalue}} 
	\left\lbrace 
		\mathbf{1}_{[-\interestingu,\interestingu]}(u')
		+
		4 
		\frac{1}{\Lunit \timefunctionarg{\muxmulevelsetvalue}}
		\muxmulevelsetvalue \phi
	\right\rbrace
	\left|\angrmD f \right|_{\gtorus}^2
	\, \volMRoughCoordinates.
\end{align}
The desired bounds
\eqref{E:COERCIVENESSOFSPACETIMEINTEGRAL}--\eqref{E:COERCIVENESSOFPARTIALSPACETIMEINTEGRAL}
now follow from \eqref{E:PROOFCOERCIVENESSOFSPACETIMEINTEGRAL} and definitions
\eqref{E:WAVESPACETIMEL2CONTROLLINGQUANTITY} and \eqref{E:PARTIALWAVESPACETIMEL2CONTROLLINGQUANTITY}.

\end{proof}


\section{The elliptic-hyperbolic integral identities} \label{S:ELLIPTICHYPERBOLICIDENTITIES}
In this section, we set up the top-order
elliptic-hyperbolic regularity theory for the transport-div-curl system satisfied 
by $\vortrenormalized$ and $\GradEnt$, i.e., for the top-order derivatives of solutions to
equations
\eqref{E:RENORMALIZEDVORTICTITYTRANSPORTEQUATION}, \eqref{E:GRADENTROPYTRANSPORT},
\eqref{E:FLATDIVOFRENORMALIZEDVORTICITY}--\eqref{E:EVOLUTIONEQUATIONFLATCURLRENORMALIZEDVORTICITY},
and
\eqref{E:TRANSPORTFLATDIVGRADENT}--\eqref{E:CURLGRADENTVANISHES}.
More precisely, for any $\Sigma_t$-tangent vectorfield $\SigmatTan$,
we derive coercive integral identities --  featuring error terms --
that are localized to spacetime
regions of the form $\twoargMrough{[\timefunction_1,\timefunction_2),[u_1,u_2]}{\muxmulevelsetvalue}$;
see Prop.\,\ref{P:INTEGRALIDENTITYFORELLIPTICHYPERBOLICCURRENT} for the main identity,
which, in view of the coerciveness guaranteed by Lemma~\ref{L:COERCIVENESSOFELLITPICHYPERBOLICQUADRATICFORM},
yields $L^2$ spacetime control of 
$\pmb{\partial} \SigmatTan$ (see definition~\eqref{E:CARTESIANGRADIENTOFTENSORFIELD})
in terms of error terms.
In our forthcoming applications, we will apply the identity with 
$\tander^{\Ntop} \vortrenormalized$ and $\tander^{\Ntop} \GradEnt$ 
in the role of $\SigmatTan$, 
and we will use the special structure of the equations of
Theorem~\ref{T:GEOMETRICWAVETRANSPORTSYSTEM} and commutator estimates 
to control the error terms. Ultimately, this will yield (see Prop.\,\ref{P:ELLIPTICHYPERBOLICINTEGRALINEQUALITIES})
spacetime $L^2$ control over the top-order terms
$\pmb{\partial} \tander^{\Ntop} \vortrenormalized$ and $\pmb{\partial} \tander^{\Ntop} \GradEnt$.

The integral identities are adaptations of the framework we developed in \cite{lAjS2020} 
to handle the structure of the singular boundary of shock-forming solutions.
Unlike in \cite{lAjS2020}, our setup here avoids boundary integrals along
the characteristics $\nullhyparg{u}$. This allows us to avoid error integrals on $\nullhyparg{u}$
that involve the top-order derivatives of $\upmu$, which would have been uncontrollable.
Our setup also yields error terms whose singularity strength is controllable under the scope of our approach.
To achieve these goals, we rely on the following key ingredients:
\begin{itemize}
	\item New well-constructed \emph{characteristic currents} (see Sect.\,\ref{SS:CHARACTERISTICCURRENT}),
		the analysis of which incorporates both the elliptic and the hyperbolic sub-structures in the equations of 
		Theorem~\ref{T:GEOMETRICWAVETRANSPORTSYSTEM}.
	\item A delicate integration by parts identity to handle some difficult
		boundary integrals along $\hypthreearg{\timefunction}{[- \rightu,u]}{\muxmulevelsetvalue}$,
		which takes into account the rough acoustic geometry and 
		the precise structure of the equations of Theorem~\ref{T:GEOMETRICWAVETRANSPORTSYSTEM};
		see Lemma~\ref{L:KEYIDPUTANGENTCURRENTCONTRACTEDAGAINSTVECTORFIELD}
		for a differential version of the identity.
		Ultimately, this leads 
		to an integral identity (see Prop.\,\ref{P:INTEGRALIDENTITYFORELLIPTICHYPERBOLICCURRENT}) 
		that provides control (see Prop.\,\ref{P:ELLIPTICHYPERBOLICINTEGRALINEQUALITIES})
		of the spacetime integrals
		$
		\int_{\twoargMrough{[\timefunction_0,\timefunction),[-\rightu,u]}{\muxmulevelsetvalue}}
			|\pmb{\p}\tander^{\Ntop}\vortrenormalized|^2
		\, \volMRoughCoordinates
		$
		and
		$
		\int_{\twoargMrough{[\timefunction_0,\timefunction),[-\rightu,u]}{\muxmulevelsetvalue}}
			|\pmb{\p}\tander^{\Ntop} \GradEnt|^2
		\, \volMRoughCoordinates
		$
		\emph{as well as} the rough tori integrals
		$
		\int_{\twoargroughtori{\timefunction,u}{\muxmulevelsetvalue}}
			|\tander^{\Ntop}\vortrenormalized|^2
		\, \volroughtorus
		$
		and
		$
		\int_{\twoargroughtori{\timefunction,u}{\muxmulevelsetvalue}}
			|\tander^{\Ntop}\GradEnt|^2
		\, \volroughtorus
		$.
		We stress that the delicate integration by parts mentioned above yields rough tori integrals with \emph{favorable signs},
		and that our proof would not have closed if the integrals had the wrong signs.
		The availability of good signs for these terms is a key aspect of 
		the framework developed in \cite{lAjS2020}.
\end{itemize}

Throughout this section, we will use the observations provided by
Remark~\ref{R:OMITTINGZEROCOMPONENTINCARTESIANCOORDINATES}.

\subsection{Basic geometric constructions and definitions}
\label{SS:ELLIPTICHYPERBOLICIDENTITIESGEOMETRICCONSTRUCTIONSANDDEFINITIONS}
In this section, we define some basic geometric objects that play a role
in our derivation of the localized integral identities. 

\begin{definition}[Projection onto $\nullhyparg{u}$ and $\nullhyparg{u}$-tangency]  
\label{D:CHARACTERISTICHYPERSURFACEPROJECTIONTENSORFIELD}
Let $\uLunit$ be the $\gfour$-null vectorfield defined in \eqref{E:ULUNIT}
(which, in view of \eqref{E:INNERPRODUCTOFLANDULISMINUS2}, is transversal to the characteristics $\nullhyparg{u}$).
\begin{enumerate}
\item We define the type $\binom{1}{1}$ projection tensorfield $\Nullhypersurfaceproject$
onto the characteristic hypersurfaces $\nullhyparg{u}$ as follows,
where $\updelta_{\beta}^{\ \alpha}$ denotes the Kronecker delta: 
\begin{align} \label{E:NULLHYPERSURFACEPROJECTION}
	\Nullhypersurfaceproject_{\beta}^{\ \alpha}
	& 
	\eqdef 
	\updelta_{\beta}^{\ \alpha}
	+
	\frac{1}{2}
	\uLunit^{\alpha}
	\Lunit_{\beta}.
\end{align}
\item Given any type $\binom{m}{n}$ spacetime tensorfield $\upxi$, 
we define its $\nullhyparg{u}$-projection $\Nullhypersurfaceproject \upxi$ 
as follows:
\begin{align}
(\Nullhypersurfaceproject \upxi)_{\beta_1 \cdots \beta_n}^{\alpha_1 \cdots \alpha_m}
& 
\eqdef 
\Nullhypersurfaceproject_{\widetilde{\alpha}_1}^{\ \alpha_1}
\cdots 
\Nullhypersurfaceproject_{\widetilde{\alpha}_m}^{\ \alpha_m}
(\gfour^{-1})^{\widetilde{\beta}_1 \widetilde{\gamma}_1}
\gfour_{\beta_1 \gamma_1}
\Nullhypersurfaceproject_{\widetilde{\beta}_1}^{\ \gamma_1}
\cdots 
(\gfour^{-1})^{\widetilde{\beta}_n \widetilde{\gamma}_n}
\gfour_{\beta_n \gamma_n}
\Nullhypersurfaceproject_{\widetilde{\beta}_n}^{\ \gamma_n}
\upxi_{\widetilde{\gamma}_1 \cdots \widetilde{\gamma}_n}^{\widetilde{\alpha}_1 \cdots \widetilde{\alpha}_m}.
	\label{E:PROJECTIONOFTENSORONTONULLHYPERSURFACE} 
\end{align}
\item We say that a spacetime tensorfield $\upxi$ is $\nullhyparg{u}$-tangent
	if $\Nullhypersurfaceproject \upxi = \upxi$. 
\end{enumerate}
\end{definition}

\begin{remark}[Lack of symmetry]
	\label{R:LACKOFSYMMETRYNULLHYPORJECTION}
	The type $\binom{0}{2}$ tensorfield 
	$
	\gfour_{\alpha \widetilde{\alpha}}
	\Nullhypersurfaceproject_{\beta}^{\ \widetilde{\alpha}}
	$ is not symmetric. 
	This is the reason that in equation \eqref{E:PROJECTIONOFTENSORONTONULLHYPERSURFACE},
	there are factors of $\gfour^{-1}$ and $\gfour$ and we
	were careful about the placement of the indices on $\Nullhypersurfaceproject$
	that are contracted against the lower indices of $\upxi$,
	unlike in equations
	\eqref{E:PROJECTIONOFTENSORONTOCARTESIANSIGMAT},
	\eqref{E:PROJECTIONOFTENSORONTOFLATTORUS},
	and
	\eqref{E:PROJECTIONOFTENSORONTOROUGHTORI}.
\end{remark}

\begin{definition}[Additional geometric tensorfields used in the elliptic-hyperbolic identities]
\label{D:GEOMETRICTENSORFIELDSFORELLIPTICHYPERBOLIC}
We respectively define $\hfour$ and $\hfour^{-1}$ 
to be the type $\binom{0}{2}$ and $\binom{2}{0}$ tensorfields with the following
Cartesian components:
\begin{subequations}
\begin{align}
\hfour_{\alpha \beta} 
	& \eqdef 
	\gfour_{\alpha \beta} 
	+ 
	2 \Transport_{\alpha} \Transport_{\beta},
		\label{E:RIEMANNIANACOUSTICALMETRIC} 
		\\
	(\hfour^{-1})^{\alpha \beta} 
	& \eqdef (\gfour^{-1})^{\alpha \beta} 
	+ 
	2 \Transport^{\alpha} \Transport^{\beta}.
	\label{E:INVERSERIEMANNIANACOUSTICALMETRIC}
\end{align}
\end{subequations} 

By using $\hfour$ and $\hfour^{-1}$ to lower and raise the indices on the projection tensorfield $\Nullhypersurfaceproject$, 
we also define the type $\binom{0}{2}$ tensorfield $\Nullhypersurfacemetric$
and the type $\binom{2}{0}$ tensorfield $\Nullhypersurfaceinversemetric$
as follows:
\begin{align}
\Nullhypersurfacemetric_{\alpha \beta}
& \eqdef \hfour_{\alpha \sigma} \Nullhypersurfaceproject_{\beta}^{\ \sigma},		
	\label{E:NULLHYPERSURFACESRIEMANNIANMETRIC} 
		\\
\Nullhypersurfaceinversemetric^{\alpha \beta}	
& \eqdef (\hfour^{-1})^{\alpha \sigma} \Nullhypersurfaceproject_{\sigma}^{\ \beta}.  
\label{E:NULLHYPERSURFACESINVERSERIEMANNIANMETRIC}
\end{align}	
\end{definition}

In the next lemma, 
we exhibit some basic properties of the tensorfields from 
Defs.\,\ref{D:CHARACTERISTICHYPERSURFACEPROJECTIONTENSORFIELD} and \ref{D:GEOMETRICTENSORFIELDSFORELLIPTICHYPERBOLIC}.
Later on, we will use the positive definiteness of
$\hfour$ and $\hfour^{-1}$ (which are revealed by the lemma)
to exhibit the coerciveness properties of various energies

\begin{lemma}[Basic properties of the tensorfields from 
Defs.\,\ref{D:CHARACTERISTICHYPERSURFACEPROJECTIONTENSORFIELD} and \ref{D:GEOMETRICTENSORFIELDSFORELLIPTICHYPERBOLIC}]
\label{L:SIMPLEIDENTITIESFORELLIPTICHYPERBOLICTENSORFIELDS}
The tensorfield $\Nullhypersurfaceproject_{\beta}^{\ \alpha}$ is a projection onto $\nullhyparg{u}$ in the following sense: 
$\Nullhypersurfaceproject \uLunit =  0$,
while if $\Singletan$ is a $\nullhyparg{u}$-tangent vectorfield,
then $\Nullhypersurfaceproject \Singletan = \Singletan$.

Moreover, the tensorfields $\Nullhypersurfacemetric$ and $\Nullhypersurfaceinversemetric$ are 
symmetric and positive semi-definite, 
and $\Nullhypersurfacemetric$ restricts to a Riemannian metric on $\nullhyparg{u}$. 

In addition, for all pairs of $\Sigma_t$-tangent vectorfields $(X,Y)$,
we have $\hfour(X,Y) = \gfour(X,Y)$, and if $\SigmatTan$ is a $\Sigma_t$-tangent vectorfield, then
$\hfour(\SigmatTan,\Transport) = 0$.
Moreover, we have $\hfour(\Transport,\Transport) = 1$,
and $\hfour$ is a Riemannian metric on spacetime.
Furthermore, following identity holds:
\begin{align} \label{E:HINVERSEISTHEINVERSOFH}
	(\hfour^{-1})^{\alpha \gamma}
	\hfour_{\gamma \beta}
	& = \updelta_{\beta}^{\alpha},
\end{align}
where $\updelta_{\beta}^{\alpha}$ is the Kronecker delta. That is, $\hfour^{-1}$ is in fact
the inverse metric of $\hfour$.

In addition, the following identities hold, where
$g$ and $g^{-1}$ are respectively the first fundamental form of $\Sigma_t$ 
and inverse first fundamental form of $\Sigma_t$ from Def.\,\ref{D:FIRSTFUNDAMENTALFORMS}:
\begin{subequations}
\begin{align}
\hfour_{\alpha \beta} 
	& 
	=
	g_{\alpha \beta} 
	+ 
	\Transport_{\alpha} \Transport_{\beta},
		\label{E:RIEMANNIANACOUSTICALMETRICINTERMSOFFIRSTFUNDOFSIGMATANDTRANSPORT} 
		\\
	(\hfour^{-1})^{\alpha \beta} 
	& 
	=
	(g^{-1})^{\alpha \beta} 
	+ 
	\Transport^{\alpha} \Transport^{\beta}.
	\label{E:INVERSERIEMANNIANACOUSTICALMETRICINTERMSOFFIRSTFUNDOFSIGMATANDTRANSPORT} 
\end{align}
\end{subequations}

Finally, the following identities hold:
\begin{align}	
	\Nullhypersurfacemetric_{\alpha \beta}
	&  = 
		\hfour_{\alpha \beta}
		-
		\frac{1}{2}
		\Lunit_{\alpha} \Lunit_{\beta}
		= 
		\gtorus_{\alpha \beta}
		+
		\frac{1}{2}
		\uLunit_{\alpha} \uLunit_{\beta},
			\label{E:NULLHYPMETRICSINTERMSOFHFOUR}
				\\
	\Nullhypersurfaceinversemetric^{\alpha \beta}	
	& 
	= (\hfour^{-1})^{\alpha \beta}
		-
		\frac{1}{2}
		\uLunit^{\alpha} \uLunit^{\beta}
		= 
		(\gtorus^{-1})^{\alpha \beta}
		+
		\frac{1}{2}
		\Lunit^{\alpha} \Lunit^{\beta}, 
		\label{E:NULLHYPERSURFACEINVERSERIEMANNIANMETRICINTERMSOFNULLFRAME}
			\\
	\Nullhypersurfaceinversemetric^{\alpha \gamma}	
	\Nullhypersurfacemetric_{\gamma \beta}
	& = 
	\updelta_{\beta}^{\alpha}
	+ 
	\frac{1}{2}
	\uLunit^{\alpha}
	\Lunit_{\beta},
	\label{E:NULLHYPINVERSEMETRICCONTRACTEDWITHNULLHYPMETRIC}
		\\
	\Nullhypersurfacemetric_{\alpha \gamma}
	\Nullhypersurfaceinversemetric^{\gamma \beta}	
	& = 
	\updelta_{\alpha}^{\beta}
	+ 
	\frac{1}{2}
	\uLunit_{\alpha}
	\Lunit^{\beta},
	\label{E:NULLHYPMETRICCONTRACTEDWITHNULLHYPINVERSEMETRIC}
		\\
	\Nullhypersurfaceinversemetric^{\alpha \beta}	
	& 
	= (\hfour^{-1})^{\alpha \gamma}
		(\hfour^{-1})^{\beta \delta}
		\Nullhypersurfacemetric_{\gamma \delta},
		\label{E:NULLHYPSURFACEINVERSEMETRICFORMEDBYRAISINGINDICESWITHHINVERSEMETRIC}
	\end{align}
where $\gtorus$ and $\gtorus^{-1}$ are respectively the first fundamental form and inverse first fundamental form
of the \textbf{smooth tori} $\ell_{t,u} = \Sigma_t \cap \nullhyparg{u}$
from Def.\,\ref{D:FIRSTFUNDAMENTALFORMS}.
\end{lemma}

\begin{proof}
The facts that
$\Nullhypersurfaceproject \uLunit =  0$,
while if $\Singletan$ is $\nullhyparg{u}$-tangent, 
then $\Nullhypersurfaceproject \Singletan = \Singletan$
follow from \eqref{E:ULUNITISNULLANDNORMALIZEDAGAINSTLUNIT}--\eqref{E:INNERPRODUCTOFLANDULISMINUS2}
and the fact that $\Lunit$ is $\gfour$-orthogonal to $\nullhyparg{u}$. 

\eqref{E:HINVERSEISTHEINVERSOFH} follows from a straightforward computation based on
definitions \eqref{E:RIEMANNIANACOUSTICALMETRIC}--\eqref{E:INVERSERIEMANNIANACOUSTICALMETRIC}
and the identity \eqref{E:TRANSPORTISUNITLENGTH}.

\eqref{E:RIEMANNIANACOUSTICALMETRICINTERMSOFFIRSTFUNDOFSIGMATANDTRANSPORT} follows from
definition \eqref{E:RIEMANNIANACOUSTICALMETRIC} and the identity \eqref{E:FIRSTFUNDOFSIGMATINTERMSOFACOUSTICALMETRICANDTRANSPORT}.
Similarly,
\eqref{E:INVERSERIEMANNIANACOUSTICALMETRICINTERMSOFFIRSTFUNDOFSIGMATANDTRANSPORT} follows from
definition \eqref{E:INVERSERIEMANNIANACOUSTICALMETRIC} and
the identity \eqref{E:INVERSEFIRSTFUNDOFSIGMATINTERMSOFACOUSTICALMETRICANDTRANSPORT}.

The remaining properties of $\hfour$ stated in the lemma follow in a straightforward fashion from
\eqref{E:TRANSPORTISUNITLENGTH} and \eqref{E:MATERIALDERIVATIVELOWEREDCARTESIANCOORDINATES}.

To prove the first equality in \eqref{E:NULLHYPMETRICSINTERMSOFHFOUR}, 
we substitute the definition of 
\eqref{E:RIEMANNIANACOUSTICALMETRIC} 
$\hfour_{\alpha \sigma}$ and the definition \eqref{E:NULLHYPERSURFACEPROJECTION} 
$\Nullhypersurfaceproject_{\beta}^{\ \sigma}$ into 
RHS~\eqref{E:NULLHYPERSURFACESRIEMANNIANMETRIC}
and carry out straightforward algebraic computations
using Lemmas~\ref{L:BASICPROPERTIESOFVECTORFIELDS} and \ref{L:BASICPROPERTIESOFULUNIT}
(in particular, we use \eqref{E:SIMPLEIDENTITIESINVOLVINGULUNIT}).
The second equality in \eqref{E:NULLHYPMETRICSINTERMSOFHFOUR} follows from similar
arguments, where we in particular use 
\eqref{E:SIMPLEIDENTITIESINVOLVINGULUNIT}
and
\eqref{E:ACOUSTICALMETRICINDOUBLENULLFRAME}.
\eqref{E:NULLHYPMETRICSINTERMSOFHFOUR} also yields the symmetry of $\Nullhypersurfacemetric$.
The identities in \eqref{E:NULLHYPERSURFACEINVERSERIEMANNIANMETRICINTERMSOFNULLFRAME}
and
\eqref{E:NULLHYPINVERSEMETRICCONTRACTEDWITHNULLHYPMETRIC}--\eqref{E:NULLHYPMETRICCONTRACTEDWITHNULLHYPINVERSEMETRIC}
as well as the symmetry of $\Nullhypersurfaceinversemetric$
follow from similar arguments based on definitions 
\eqref{E:NULLHYPERSURFACESINVERSERIEMANNIANMETRIC}
and \eqref{E:NULLHYPERSURFACEPROJECTION}.

The fact that $\Nullhypersurfacemetric$ restricts to a Riemannian metric on $\nullhyparg{u}$
follows from \eqref{E:NULLHYPMETRICSINTERMSOFHFOUR},
Lemmas~\ref{L:BASICPROPERTIESOFVECTORFIELDS} and \ref{L:BASICPROPERTIESOFULUNIT},
the fact that the tangent space of $\nullhyparg{u}$ is the direct sum of 
the tangent space of $\ell_{t,u}$ and the span of $\Lunit$,
the fact that $\gtorus_{\alpha \beta}$ is 
positive definite on $\ell_{t,u}$-tangent vectorfields and satisfies $\gtorus(\Lunit,\cdot) = 0$.
The positive semi-definiteness of
$\Nullhypersurfacemetric$ and $\Nullhypersurfaceinversemetric$
follows from the second identities in
\eqref{E:NULLHYPMETRICSINTERMSOFHFOUR}--\eqref{E:NULLHYPERSURFACEINVERSERIEMANNIANMETRICINTERMSOFNULLFRAME}
and the positive semi-definiteness of $\gtorus_{\alpha \beta}$ and $(\gtorus^{-1})^{\alpha \beta}$. 

The identity \eqref{E:NULLHYPSURFACEINVERSEMETRICFORMEDBYRAISINGINDICESWITHHINVERSEMETRIC}
follows from definitions 
\eqref{E:NULLHYPERSURFACESRIEMANNIANMETRIC}--\eqref{E:NULLHYPERSURFACESINVERSERIEMANNIANMETRIC},
the identity \eqref{E:HINVERSEISTHEINVERSOFH},
and the symmetry of $\Nullhypersurfacemetric$ and $\Nullhypersurfaceinversemetric$.

\end{proof}

\subsection{Additional derivatives operators, pointwise norms, and comparison estimates}
\label{SS:ADDITIONALDERIVATIVEOPERATORSNORMSANDCOMPARISIONESTIMATE}
In this section, we define some additional derivative operators and pointwise norms,
and we establish some simple comparison estimates. Later, we will use them in our analysis of
the terms in the elliptic-hyperbolic integral identity provided by Prop.\,\ref{P:INTEGRALIDENTITYFORELLIPTICHYPERBOLICCURRENT}.

\subsubsection{The Cartesian gradient of $\upxi$ and $|Z \SigmatTan|_g$}
\label{SSS:DERIVATIVEOFATENSORFIELD}

\begin{definition}[The Cartesian gradient of $\upxi$ and $|Z \SigmatTan|_g$]
\label{D:DERIVATIVEOFATENSORFIELD}
\item 
	Given a type $\binom{m}{n}$ spacetime tensorfield $\upxi$, 
	we define its \emph{Cartesian gradient} 
	$\pmb{\partial} \upxi$ to be the type $\binom{m}{n+1}$ spacetime tensorfield with
	the following Cartesian components: 
	\begin{align} \label{E:CARTESIANGRADIENTOFTENSORFIELD}
	(\pmb{\partial} \upxi)_{\beta_1 \beta_2 \cdots \beta_{n+1}}^{\alpha_1 \cdots \alpha_m}
	& \eqdef
	\p_{\beta_1} \upxi_{\beta_2 \cdots \beta_{n+1}}^{\alpha_1 \cdots \alpha_m}.
	\end{align}
\item Let $\SigmatTan$ be a $\Sigma_t$-tangent vectorfield and let $Z$ be a spacetime vectorfield.
Relative to the Cartesian coordinates, 
we define (see Remark~\ref{R:OMITTINGZEROCOMPONENTINCARTESIANCOORDINATES}) $|Z \SigmatTan|_g \geq 0$
as follows:
\begin{align} \label{E:SIGMATPOINTWISENOROMOFZDERIVATIVEOFSIGMATTANGENTVECTORFIELD}
		 |Z \SigmatTan|_g^2 
		& \eqdef g_{ab} (Z \SigmatTan^a) Z \SigmatTan^b,
\end{align}
where as usual $Z \SigmatTan^a = Z^{\alpha} \partial_{\alpha} \SigmatTan^a$.
\end{definition}

Note that $|Z \SigmatTan|_g$ is the $|\cdot|_g$ norm of the $\Sigma_t$-tangent vectorfield with the Cartesian spatial components
$Z \SigmatTan^i$, $i=1,2,3$.

\subsubsection{The $\hfour$-norm of tensorfields}
\label{SSS:RIEMANNIANACOUSTICALMETRICSPACETIMENORMOFTENSORFIELDS}

\begin{definition}[The $\hfour$-norm of tensorfields]
\label{D:RIEMANNIANACOUSTICALMETRICSPACETIMENORMOFTENSORFIELDS}
	Recall that in Lemma~\ref{L:SIMPLEIDENTITIESFORELLIPTICHYPERBOLICTENSORFIELDS}, we showed that $\hfour$ is
	a Riemannian metric on spacetime.
	If $\upxi$ is a type $\binom{m}{n}$ spacetime tensorfield,
	then we define $|\upxi|_{\hfour}\geq 0$ by:
\begin{align} \label{E:SQUAREPOINTWISENORMWITHRESPECTTORIEMANNIANACOUSTICALMETRIC}
|\upxi|_{\hfour}^2 
& 
\eqdef 
\hfour_{\alpha_1 \widetilde{\alpha}_1} 
\cdots 
\hfour_{\alpha_m \widetilde{\alpha}_m} 
(\hfour^{-1})^{\beta_1 \widetilde{\beta}_1} 
\cdots 
(\hfour^{-1})^{\beta_n \widetilde{\beta}_n} 
\upxi_{\beta_1 \cdots \beta_n}^{\alpha_1 \cdots \alpha_n} 
\upxi_{\widetilde{\beta}_1 \cdots \widetilde{\beta}_n}^{\widetilde{\alpha}_1 \cdots \widetilde{\alpha}_m}.
\end{align}
\end{definition}

\subsubsection{Pointwise comparison results and the $\hfour$-size of
$\Nullhypersurfaceproject$, $\Nullhypersurfacemetric$, and $\Nullhypersurfaceinversemetric$}
\label{SSS:HRIEMANNIANMETRICOMPARISONANDNORMSOFNULLHYPERSURFACETENSORFIELDS}

\begin{lemma}[Pointwise comparison results and the $\hfour$-size of
$\Nullhypersurfaceproject$, $\Nullhypersurfacemetric$, and $\Nullhypersurfaceinversemetric$] 
\label{L:HRIEMANNIANMETRICOMPARISONANDNORMSOFNULLHYPERSURFACETENSORFIELDS}
For any type $\binom{m}{n}$ spacetime tensorfield
$\upxi^{\alpha_1 \cdots \alpha_m}_{\beta_1 \cdots \beta_n}$,
the following comparison estimates 
hold relative to the Cartesian coordinates on $\twoargMrough{[\timefunction_0,\timefunctionboot),[- \rightu,\leftu]}{\muxmulevelsetvalue}$:
\begin{subequations}
\begin{align} \label{E:HRIEMANNIANMETRICNORMCOMPARABLETOEUCLIDEANNORM}
	|\upxi|_{\hfour}
	& \approx
		\sum_{\substack{0 \leq \alpha_1, \cdots, \alpha_m \leq 3 \\ 0 \leq \beta_1, \cdots, \beta_n \leq 3}}
		\left|\upxi_{\beta_1 \cdots \beta_n}^{\alpha_1 \cdots \alpha_m} \right|,
			\\
	|\pmb{\partial} \upxi|_{\hfour}
	& \approx
		\sum_{\substack{0 \leq \alpha_1, \cdots, \alpha_m \leq 3 \\ 0 \leq \beta_1, \cdots, \beta_n \leq 3 \\ 0 \leq \gamma \leq 3}}
		\left|\partial_{\gamma} 
		\upxi_{\beta_1 \cdots \beta_n}^{\alpha_1 \cdots \alpha_m} \right|.
		\label{E:CARTERSIANDERIVATIVESHRIEMANNIANMETRICNORMCOMPARABLETOEUCLIDEANNORM}
\end{align}
\end{subequations}

Moreover, if $\SigmatTan$ is a $\Sigma_t$-tangent vectorfield and $Z$ is a spacetime vectorfield,
then the following comparison estimates 
hold relative to the Cartesian coordinates 
on $\twoargMrough{[\timefunction_0,\timefunctionboot),[- \rightu,\leftu]}{\muxmulevelsetvalue}$:
\begin{subequations}
\begin{align} 
|\SigmatTan|_g
	& \approx
	\sum_{a=1,2,3}
	\left| \SigmatTan^a \right|,
		\label{E:SIGMATTANGENTVECTORFIELDGNORMAPPROXIMATEDBYCARETSIANCOMPONENTNORMS} 
		\\
	|Z \SigmatTan|_g
	& \approx
		\sum_{a=1,2,3}
		\left| Z \SigmatTan^a \right|,
			\label{E:ZDERIVATIVEOFSIGMATTANGENTVECTORFIELDGNORMAPPROXIMATEDBYCARETSIANCOMPONENTNORMS}
				\\
	|\pmb{\partial} \SigmatTan|_{\hfour}
	& \approx
		\sum_{\substack{\alpha = 0,1,2,3 \\ a=1,2,3}}
		\left|\partial_{\alpha} \SigmatTan^a \right|
		\approx
		\sum_{\alpha = 0,1,2,3}
		|\partial_{\alpha} \SigmatTan|_g
		\approx
		|\Transport \SigmatTan|_g
		+
		\sum_{a=1,2,3}
		|\partial_a \SigmatTan|_g.
		\label{E:TRANSPORTANDCARTERSIANDERIVATIVESHRIEMANNIANMETRICNORMCOMPARABLETOEUCLIDEANNORM}
\end{align}
\end{subequations}

Finally, the following identities hold:
\begin{align}  \label{E:HRIEMANNIANNORMOFCHARACTERISTICTENSORFIELDSISROOT3}
	|\Nullhypersurfaceproject|_{\hfour} 
	&
	=
	|\Nullhypersurfacemetric|_{\hfour} 
	= 
	|\Nullhypersurfaceinversemetric|_{\hfour} 
	= 
	\sqrt{3}.
\end{align}
\end{lemma}

\begin{proof}
We prove \eqref{E:HRIEMANNIANMETRICNORMCOMPARABLETOEUCLIDEANNORM} when $\upxi$ is a type $\binom{1}{0}$ tensorfield;
the case of general type $\binom{m}{n}$ tensorfields can be handled through similar arguments.
To proceed, we first use
\eqref{E:MATERIALDERIVATIVELOWEREDCARTESIANCOORDINATES},
\eqref{E:FIRSTFUNDOFSIGMATIDENTITY},
and \eqref{E:RIEMANNIANACOUSTICALMETRICINTERMSOFFIRSTFUNDOFSIGMATANDTRANSPORT}
to derive the following identity relative to the Cartesian coordinates:
\begin{align} \label{E:HNORMSQUAREDOFVECTORFIELDEXPANDEDOUT}
|\upxi|_{\hfour}^2 
&
=
\hfour_{\alpha \beta} \upxi^{\alpha} \upxi^{\beta}
=
\Speed^{-2} \sum_{a=1,2,3} (\upxi^a)^2
+ 
\left\lbrace
	1
	+
	\Speed^{-2}
	\sum_{a=1,2,3} (v^a)^2
\right\rbrace
(\upxi^0)^2
-
2 \Speed^{-2}
\sum_{a=1,2,3}
v^a \upxi^a \upxi^0.
\end{align} 
From the bootstrap assumptions and \eqref{E:SCHEMATICSTRUCTUREOFXSMALL}, 
we deduce that 
$|v^a| \lesssim 1$ and $\Speed \approx 1$. From these estimates and \eqref{E:HNORMSQUAREDOFVECTORFIELDEXPANDEDOUT},
we conclude that $|\upxi|_{\hfour} \lesssim \sum_{\alpha = 0,1,2,3} |\upxi^{\alpha}|$.
To prove the reverse inequality, we note that
the Cauchy--Schwarz inequality and
Young's inequality imply that the cross term 
$
-
2 \Speed^{-2}
v^a \upxi^a \upxi^0$ on RHS~\eqref{E:HNORMSQUAREDOFVECTORFIELDEXPANDEDOUT}
is bounded in magnitude by
$
\leq 
\Speed^{-2} \sum_{a=1,2,3} (\upxi^a)^2
+
\Speed^{-2} 
\sum_{a=1,2,3} (v^a)^2
(\upxi^0)^2
$.
From this bound, \eqref{E:HNORMSQUAREDOFVECTORFIELDEXPANDEDOUT},
and the estimate $\Speed \approx 1$,
it follows that: 
\begin{align} \label{E:HNORMCONTROLSCARTESIAN0COMPONENTS}
	|\upxi|_{\hfour}^2 \geq (\upxi^0)^2.
\end{align}
Moreover, the Cauchy--Schwarz inequality and
Young's inequality imply that the cross term 
$
-
2 \Speed^{-2}
v^a \upxi^a \upxi^0$
is also bounded in magnitude by:
\[
\leq 
\left\lbrace
	\frac{\Speed^{-2} \sum_{d=1,2,3} (v^d)^2}{1 + \Speed^{-2} \sum_{b=1,2,3} (v^b)^2}
\right\rbrace
\times
\Speed^{-2} \sum_{a=1,2,3} (\upxi^a)^2
+
\left\lbrace
	1
	+
	\Speed^{-2}
	\sum_{a=1,2,3} (v^a)^2
\right\rbrace
(\upxi^0)^2.
\]
From this bound, \eqref{E:HNORMSQUAREDOFVECTORFIELDEXPANDEDOUT},
and the aforementioned estimates $|v^a| \lesssim 1$ and $\Speed \approx 1$,
it follows that: 
\begin{align} \label{E:HNORMCONTROLSCARTESIANSPATIALCOMPONENTS}
|\upxi|_{\hfour}^2 
&
\geq 
\left\lbrace
	\frac{1}{1 + \Speed^{-2} \sum_{b=1,2,3} (v^b)^2}
\right\rbrace
\times
\Speed^{-2}
\sum_{a=1,2,3} (\upxi^a)^2
\gtrsim 
\sum_{a=1,2,3} (\upxi^a)^2.
\end{align}
Combining \eqref{E:HNORMCONTROLSCARTESIAN0COMPONENTS} and \eqref{E:HNORMCONTROLSCARTESIANSPATIALCOMPONENTS},
we see that
$\sum_{\alpha = 0,1,2,3} |\upxi^{\alpha}| \lesssim |\upxi|_{\hfour}$, 
which completes the proof of \eqref{E:HRIEMANNIANMETRICNORMCOMPARABLETOEUCLIDEANNORM}.
The comparison estimates 
\eqref{E:SIGMATTANGENTVECTORFIELDGNORMAPPROXIMATEDBYCARETSIANCOMPONENTNORMS},
\eqref{E:ZDERIVATIVEOFSIGMATTANGENTVECTORFIELDGNORMAPPROXIMATEDBYCARETSIANCOMPONENTNORMS}, 
and
\eqref{E:TRANSPORTANDCARTERSIANDERIVATIVESHRIEMANNIANMETRICNORMCOMPARABLETOEUCLIDEANNORM}
follow from similar arguments, and we omit the straightforward details.

\eqref{E:CARTERSIANDERIVATIVESHRIEMANNIANMETRICNORMCOMPARABLETOEUCLIDEANNORM}
follows as a special case of \eqref{E:HRIEMANNIANMETRICNORMCOMPARABLETOEUCLIDEANNORM}
with $\pmb{\partial} \upxi$ in the role of $\upxi$.

We now prove \eqref{E:HRIEMANNIANNORMOFCHARACTERISTICTENSORFIELDSISROOT3}. 
First, using
\eqref{E:INNERPRODUCTOFLANDULISMINUS2},
definitions 
\eqref{E:NULLHYPERSURFACESRIEMANNIANMETRIC}--\eqref{E:NULLHYPERSURFACESINVERSERIEMANNIANMETRIC}, 
and \eqref{E:NULLHYPMETRICCONTRACTEDWITHNULLHYPINVERSEMETRIC},
we compute that
$|\Nullhypersurfaceproject|_{\hfour}^2 
= 
\Nullhypersurfaceinversemetric^{\alpha \beta}
\Nullhypersurfacemetric_{\beta \alpha} 
= 
\updelta_{\alpha}^{\alpha}
+
\frac{1}{2} \uLunit_{\alpha} \Lunit^{\alpha}
= 
4 - 1
= 3$,
as is desired. 
Next, we use 
\eqref{E:INNERPRODUCTOFLANDULISMINUS2},
\eqref{E:NULLHYPERSURFACEPROJECTION},
\eqref{E:NULLHYPERSURFACESRIEMANNIANMETRIC},
and
Lemma~\ref{L:SIMPLEIDENTITIESFORELLIPTICHYPERBOLICTENSORFIELDS}
to compute that:
\begin{align}
\begin{split} \label{E:INTERMEDIATESTEPINPROOFOFHRIEMANNIANNORMOFCHARACTERISTICTENSORFIELDS}
|\Nullhypersurfacemetric|_{\hfour}^2 
&
=
(\hfour^{-1})^{\alpha \beta} 
(\hfour^{-1})^{\gamma \delta} 
\Nullhypersurfacemetric_{\alpha \gamma} 
\Nullhypersurfacemetric_{\beta \delta}
=
\Nullhypersurfaceproject_{\gamma}^{\ \beta}
\Nullhypersurfaceproject_{\beta}^{\ \gamma}
=
\left\lbrace
	\updelta_{\gamma}^{\beta}
	+
	\frac{1}{2} \Lunit^{\beta} \uLunit_{\gamma}
\right\rbrace
\left\lbrace
	\updelta_{\beta}^{\gamma}
	+
	\frac{1}{2} \Lunit^{\gamma} \uLunit_{\beta}
\right\rbrace
	\\
&
= 4 
+ 
\Lunit^{\alpha} \uLunit_{\alpha}
+
\frac{1}{4} (\Lunit^{\alpha} \uLunit_{\alpha})^2
 = 3.
\end{split}
\end{align}
A similar calculation based on \eqref{E:NULLHYPERSURFACESINVERSERIEMANNIANMETRIC}
yields that $|\Nullhypersurfaceinversemetric|_{\hfour}^2 = 3$.
We have therefore proved \eqref{E:HRIEMANNIANNORMOFCHARACTERISTICTENSORFIELDSISROOT3}.

\end{proof}

\subsection{The coercive elliptic-hyperbolic quadratic form and its coerciveness}
\label{SS:NULLHYPERSURFACEADAPTEDCOERCIVEQUADRATICFORM}

\subsubsection{The coercive elliptic-hyperbolic quadratic form}
\label{SSS:DEFINITIONOFELLIPTICHYPERBOLICQUADRATICFORM}
In the next definition, we introduce the solution-adapted quadratic form $\ellipticCoerciveQuadratic$ 
that we will use to control the top-order derivatives of the specific vorticity and entropy gradient. 

\begin{definition}[The coercive elliptic-hyperbolic quadratic form]
\label{D:NULLHYPERSURFACEADAPTEDCOERCIVEQUADRATICFORM}
Let $\smoothtorusproject_b^{\ a}$ denote the $\Sigma_t$-components of the $\ell_{t,u}$ projection tensorfields defined in 
\eqref{E:SMOOTHTORUSPROJECT},
and let
and 
$\Nullhypersurfacemetric$
and
$\Nullhypersurfaceinversemetric$
be the tensorfields defined in
\eqref{E:NULLHYPERSURFACESRIEMANNIANMETRIC}--\eqref{E:NULLHYPERSURFACESINVERSERIEMANNIANMETRIC}
respectively.
Let $\SigmatTan$ be a $\Sigma_t$-tangent vectorfield. We define
$\ellipticCoerciveQuadratic[\pmb{\partial} \SigmatTan,\pmb{\partial} \SigmatTan]$
to be the following quadratic form associated to $\SigmatTan$:
\begin{align} 
\begin{split} \label{E:NULLHYPERSURFACEADAPTEDCOERCIVEQUADRATICFORM}
		\ellipticCoerciveQuadratic[\pmb{\partial} \SigmatTan,\pmb{\partial} \SigmatTan]
		&
		\eqdef
		\left(
			\Nullhypersurfaceinversemetric^{\alpha \beta}
			+
			4
			\Transport^{\alpha} \Transport^{\beta}
		\right)
		\left(
			\Nullhypersurfacemetric_{\gamma \delta}
			+
			4
			\Transport_{\gamma} \Transport_{\delta}
		\right)
		(\partial_{\alpha} \SigmatTan^{\gamma})
		\partial_{\beta} \SigmatTan^{\delta}
		 \\ 
		& \ \
		-
	\frac{1}{16}
	\left\lbrace
		- 
		3	
		(\Transport \SigmatTan^{\alpha})
		\uLunit_{\alpha}
		+
		(\Lunit \SigmatTan^{\alpha})
		\uLunit_{\alpha}
		+
		\smoothtorusproject_b^{\ a}
		\partial_a \SigmatTan^b
	\right\rbrace^2.
\end{split}
\end{align}
\end{definition}

\subsubsection{The coerciveness of the elliptic-hyperbolic quadratic form}
\label{SSS:COERCIVITYOFELLIPTICHYPERBOLICQUADRATICFORM}
In the next lemma, we exhibit the coerciveness of 
$\ellipticCoerciveQuadratic[\pmb{\partial} \SigmatTan,\pmb{\partial} \SigmatTan]$.

\begin{lemma}[Coerciveness of {$\ellipticCoerciveQuadratic[\pmb{\partial} \SigmatTan,\pmb{\partial} \SigmatTan]$}]
\label{L:COERCIVENESSOFELLITPICHYPERBOLICQUADRATICFORM}
	On $\twoargMrough{[\timefunction_0,\timefunctionboot),[- \rightu,\leftu]}{\muxmulevelsetvalue}$,
	the quadratic form
	$\ellipticCoerciveQuadratic$ from Def.\,\ref{D:NULLHYPERSURFACEADAPTEDCOERCIVEQUADRATICFORM}
	is quantitatively positive definite 
	on the space of Cartesian gradients of $\Sigma_t$-tangent vectorfields $\SigmatTan$
	in the following sense, 
	where
	$
	|\pmb{\partial} \SigmatTan|_{\hfour}^2
	\eqdef (\hfour^{-1})^{\alpha \beta} \hfour_{\gamma \delta} (\partial_{\alpha} \SigmatTan^{\gamma}) \partial_{\beta} \SigmatTan^{\delta}
	$:
	\begin{align} \label{E:COERCIVENESSOFELLIPTICHYPERBOLICQUADRATICFORM}
		\ellipticCoerciveQuadratic[\pmb{\partial} \SigmatTan,\pmb{\partial} \SigmatTan]
		& \approx 
			|\pmb{\partial} \SigmatTan|_{\hfour}^2
			\approx 
			\sum_{\alpha = 0}^3 |\partial_{\alpha} \SigmatTan|_g^2.
	\end{align}
\end{lemma}

\begin{proof}
	Throughout the proof, we will use the observations made in Remark~\ref{R:OMITTINGZEROCOMPONENTINCARTESIANCOORDINATES}.
	To start, we note that Young's inequality implies
	that the terms 
	$
		-
	\frac{1}{16}
	\left\lbrace
		- 
		3	
		(\Transport \SigmatTan^{\alpha})
		\uLunit_{\alpha}
		+
		(\Lunit \SigmatTan^{\alpha})
		\uLunit_{\alpha}
		+
		\smoothtorusproject_b^{\ a}
		\partial_a \SigmatTan^b
	\right\rbrace^2
	$
	on RHS~\eqref{E:NULLHYPERSURFACEADAPTEDCOERCIVEQUADRATICFORM}
	are bounded in magnitude 
	by:
	\begin{align}
	\begin{split} \label{E:PROOFSTEP1COERCIVENESSOFELLIPTICHYPERBOLICQUADRATICFORM}
	&
	\leq 
	\frac{4}{16}
	[3 (\Transport \SigmatTan^{\alpha}) \uLunit_{\alpha}]^2
	+
	\frac{4}{3}
	\times
	\frac{1}{16}
	\left\lbrace
		(\Lunit \SigmatTan^{\alpha})
		\uLunit_{\alpha}
		+
		\smoothtorusproject_b^{\ a}
		\partial_a \SigmatTan^b
	\right\rbrace^2 
		\\
	& 
	\leq
	\frac{9}{4}
	[(\Transport \SigmatTan^{\alpha}) \uLunit_{\alpha}]^2
	+
	\frac{1}{6}
	[(\Lunit \SigmatTan^{\alpha}) \uLunit_{\alpha}]^2
	+
	\frac{1}{6}
	(\smoothtorusproject_b^{\ a} \partial_a \SigmatTan^b)^2.
	\end{split}
	\end{align}
	Moreover, the Cauchy--Schwarz inequality and the fact that $|\smoothtorusproject|_{\gtorus}^2 = \mytr_{\gtorus} \gtorus = 2$
	together imply that:
	\begin{align} \label{E:PROOFSTEP2COERCIVENESSOFELLIPTICHYPERBOLICQUADRATICFORM}
	\frac{1}{6} (\smoothtorusproject_b^{\ a} \partial_a \SigmatTan^b)^2
	& \leq \frac{1}{3} (\gtorus^{-1})^{ab} \gtorus_{cd} (\partial_a \SigmatTan^c) \partial_b \SigmatTan^d
		\eqdef \frac{1}{3} |\partial \SigmatTan|_{\gtorus}^2.
	\end{align}
	Using 
	\eqref{E:NULLHYPMETRICSINTERMSOFHFOUR}--\eqref{E:NULLHYPERSURFACEINVERSERIEMANNIANMETRICINTERMSOFNULLFRAME},
	\eqref{E:NULLHYPERSURFACEADAPTEDCOERCIVEQUADRATICFORM},
	and
	\eqref{E:PROOFSTEP1COERCIVENESSOFELLIPTICHYPERBOLICQUADRATICFORM}--\eqref{E:PROOFSTEP2COERCIVENESSOFELLIPTICHYPERBOLICQUADRATICFORM}, 
	we compute that:
	\begin{align} \label{E:PROOFSTEP3COERCIVENESSOFELLIPTICHYPERBOLICQUADRATICFORM}
	\ellipticCoerciveQuadratic[\pmb{\partial} \SigmatTan,\pmb{\partial} \SigmatTan]
	&
	\approx
	\left(
			\Nullhypersurfaceinversemetric^{\alpha \beta}
			+
			4
			\Transport^{\alpha} \Transport^{\beta}
		\right)
		\left(
			\Nullhypersurfacemetric_{\gamma \delta}
			+
			4
			\Transport_{\gamma} \Transport_{\delta}
		\right)
		(\partial_{\alpha} \SigmatTan^{\gamma})
		\partial_{\beta} \SigmatTan^{\delta}.
	\end{align}
	From \eqref{E:PROOFSTEP3COERCIVENESSOFELLIPTICHYPERBOLICQUADRATICFORM},
	the fact that (in Cartesian coordinates)
	$
	\Transport_{\gamma} \partial_{\alpha} \SigmatTan^{\gamma} 
	=
	-
	\partial_{\alpha} \SigmatTan^0
	= 0$,
	the identities $\Lunit = \Transport - X$ and $\uLunit = \Transport + X$,
	the identities $g_{ab} = \gtorus_{ab} + X_a X_b$ and $(g^{-1})^{ab} = (\gtorus^{-1})^{ab} + X^a X^b$
	proved in 
	\eqref{E:SMOOTHTORUSMETRICINTERMSOFSIGMATMETRICANDX}--\eqref{E:SMOOTHTORUSINVERSEMETRICINTERMSOFINVERSESIGMATMETRICANDX},
	\eqref{E:GTORUSINVERSE0COMPONENTSVANISH},
	and the identities
	\eqref{E:RIEMANNIANACOUSTICALMETRICINTERMSOFFIRSTFUNDOFSIGMATANDTRANSPORT}--\eqref{E:INVERSERIEMANNIANACOUSTICALMETRICINTERMSOFFIRSTFUNDOFSIGMATANDTRANSPORT}
	and
	\eqref{E:NULLHYPMETRICSINTERMSOFHFOUR}--\eqref{E:NULLHYPERSURFACEINVERSERIEMANNIANMETRICINTERMSOFNULLFRAME},
	it follows that:
	\begin{align}
	\begin{split} \label{E:PROOFSTEP4COERCIVENESSOFELLIPTICHYPERBOLICQUADRATICFORM}
	\ellipticCoerciveQuadratic[\pmb{\partial} \SigmatTan,\pmb{\partial} \SigmatTan]
	&
	\approx
	|\Transport \SigmatTan|_g^2	
	+
	|\Lunit \SigmatTan|_g^2	
	+
	(\gtorus^{-1})^{ab} g_{cd} (\partial_a \SigmatTan^c) \partial_b \SigmatTan^d
	\approx
	|\Transport \SigmatTan|_g^2	
	+
	|X \SigmatTan|_g^2	
	+
	(\gtorus^{-1})^{ab} g_{cd} (\partial_a \SigmatTan^c) \partial_b \SigmatTan^d
		\\
	&
	\approx
	|\Transport \SigmatTan|_g^2
	+
	(g^{-1})^{ab} g_{cd} (\partial_a \SigmatTan^c) \partial_b \SigmatTan^d
	\approx
	|\pmb{\partial} \SigmatTan|_{\hfour}^2.
	\end{split}
	\end{align}
	\eqref{E:PROOFSTEP4COERCIVENESSOFELLIPTICHYPERBOLICQUADRATICFORM} implies
	the first ``$\approx$'' in \eqref{E:COERCIVENESSOFELLIPTICHYPERBOLICQUADRATICFORM}.
	The second ``$\approx$'' in \eqref{E:COERCIVENESSOFELLIPTICHYPERBOLICQUADRATICFORM}
	then follows from Lemma~\ref{L:HRIEMANNIANMETRICOMPARISONANDNORMSOFNULLHYPERSURFACETENSORFIELDS}.
\end{proof}

\subsection{The characteristic currents}
\label{SS:CHARACTERISTICCURRENT}
The $\nullhyparg{u}$-tangent vectorfields $\ehcurrent[\SigmatTan,\pmb{\p}\SigmatTan]$ in the next definition play a key role in our analysis. 
In our proof of Prop.\,\ref{P:INTEGRALIDENTITYFORELLIPTICHYPERBOLICCURRENT}, 
use them for bookkeeping when integrating by parts.
We sometimes refer to the $\ehcurrent[\SigmatTan,\pmb{\p}\SigmatTan]$ 
as ``characteristic currents'' since they are tangent to $\nullhyparg{u}$, 
or ``elliptic-hyperbolic currents'' since they are the basic ingredient for the elliptic-hyperbolic integral identities. 
In our prior work \cite{lAjS2020}, we used related -- but distinct --
$\hypthreearg{\timefunction}{[- \rightu,u]}{\muxmulevelsetvalue}$-tangent currents 
to derive elliptic-hyperbolic integral identities.
Compared to the currents in \cite{lAjS2020}, 
the ones featured in the next definition are better adapted to the structure of the singularity in the sense
that they do not generate any critical-strength error terms in our top-order $L^2$ estimates.
Moreover, when we integrate by parts over the spacetime region $\twoargMrough{[\timefunction_1,\timefunction_2),[u_1,u_2]}{\muxmulevelsetvalue}$,
the $\nullhyparg{u}$-tangency of the $\ehcurrent[\SigmatTan,\pmb{\p}\SigmatTan]$
allows us to avoid boundary integrals along $\nullhyparg{u}$. This is important because
some of the acoustic geometry error terms (such as the top-order derivatives of $\upmu$)
do not have sufficiently regularity to be controlled in $L^2$ along $\nullhyparg{u}$.

\subsubsection{The definition of the $\nullhyparg{u}$-tangent characteristic current}
\label{SSS:DEFOFCHARCURRENT}

\begin{definition}[The $\nullhyparg{u}$-tangent characteristic current]
\label{D:PUTANGENTELLIPTICHYPERBOLICCURRENT}
Let $\SigmatTan$ be a $\Sigma_t$-tangent vectorfield.
We define the characteristic current
to be the vectorfield $\ehcurrent = \ehcurrent[\SigmatTan,\pmb{\partial} \SigmatTan]$
with the following components, $(\alpha = 0,1,2,3)$:
\begin{align} \label{E:PUTANGENTELLIPTICHYPERBOLICCURRENT}
	\ehcurrent^{\alpha}[\SigmatTan,\pmb{\partial} \SigmatTan]
	& \eqdef 
	\SigmatTan^{\gamma} 
	\Nullhypersurfaceproject_{\gamma}^{\ \lambda} 
	\Nullhypersurfaceproject_{\kappa}^{\ \alpha} \partial_{\lambda} \SigmatTan^{\kappa}
	-
	\SigmatTan^{\gamma} \Nullhypersurfaceproject_{\gamma}^{\ \alpha} 
	\Nullhypersurfaceproject_{\lambda}^{\ \kappa}
	\partial_{\kappa} \SigmatTan^{\lambda}.
\end{align}
\end{definition}

\begin{remark}[{$\ehcurrent^{\alpha}[\SigmatTan,\pmb{\partial} \SigmatTan]$} is $\nullhyparg{u}$-tangent]
	\label{R:CHARACTERISTICCURRENTISPUTANGENT}
	Since for any vectorfield $Z$,
	the vectorfield $\Nullhypersurfaceproject_{\beta}^{\ \alpha} Z^{\beta}$
	is $\nullhyparg{u}$-tangent,
	it follows from \eqref{E:PUTANGENTELLIPTICHYPERBOLICCURRENT} 
	that indeed, $\ehcurrent^{\alpha}[\SigmatTan,\pmb{\partial} \SigmatTan]$ 
	is $\nullhyparg{u}$-tangent.
\end{remark}

\subsubsection{The covariant divergence identity satisfied by the elliptic-hyperbolic current}
\label{SSS:COVARIANTDIVERGENCEIDENTITYSATISFIEDBYCHARCURRENT}
In the next lemma, we provide the main covariant divergence identity
satisfied by the current $\ehcurrent^{\alpha}[\SigmatTan,\pmb{\partial} \SigmatTan]$ 
from Def.\,\ref{D:PUTANGENTELLIPTICHYPERBOLICCURRENT}.
The identity forms the starting point for the divergence-theorem based proof of Prop.\,\ref{P:INTEGRALIDENTITYFORELLIPTICHYPERBOLICCURRENT}.

\begin{lemma}[Covariant divergence identity for the elliptic-hyperbolic current]
	\label{L:COVARIANTDIVERGENCEIDENTITYFORELLIPTICHYPERBOLICCURRENT}
	Let 
	$\SigmatTan$ be a $\Sigma_t$-tangent vectorfield,
	let $\ellipticCoerciveQuadratic[\pmb{\partial} \SigmatTan,\pmb{\partial} \SigmatTan]$
	be the quadratic form defined by \eqref{E:NULLHYPERSURFACEADAPTEDCOERCIVEQUADRATICFORM},
	and let $\weight$ be a ``weight function.''
	Then the following identity holds relative to the Cartesian coordinates, 
	where $\Dfour$ is the Levi-Civita connection of $\gfour$:
	\begin{align}	
	\begin{split} \label{E:COVARIANTDIVERGENCEIDENTITYFORELLIPTICHYPERBOLICCURRENT}
		\weight \ellipticCoerciveQuadratic[\pmb{\partial} \SigmatTan,\pmb{\partial} \SigmatTan]
		& = \Dfour_{\alpha} (\weight \ehcurrent^{\alpha}[\SigmatTan,\pmb{\partial} \SigmatTan])
			+
			\weight
			\mathfrak{J}_{(\textnormal{Antisymmetric})}[\pmb{\partial} \SigmatTan,\pmb{\partial} \SigmatTan]
			+
			\weight
			\mathfrak{J}_{(\textnormal{Div})}[\pmb{\partial} \SigmatTan,\pmb{\partial} \SigmatTan]
				 \\
		& \ \
			+
			\mathfrak{J}_{(\pmb{\partial} \weight)}[\SigmatTan,\pmb{\partial} \SigmatTan]
			+
			\weight
			\mathfrak{J}_{(\textnormal{Absorb-1})}[\SigmatTan,\pmb{\partial} \SigmatTan]
			+
			\weight
			\mathfrak{J}_{(\textnormal{Absorb-2})}[\SigmatTan,\pmb{\partial} \SigmatTan]
				\\
		& \ \
			+
			\weight
			\mathfrak{J}_{(\textnormal{Material})}[\pmb{\partial} \SigmatTan,\pmb{\partial} \SigmatTan]
			+
			\weight
			\mathfrak{J}_{(\textnormal{Null Geometry})}[\SigmatTan,\pmb{\partial} \SigmatTan],
	\end{split}
	\end{align}
	where with 
	$(\mathrm{d} \SigmatTan_{\flat})_{\alpha \beta} 
	\eqdef \partial_{\alpha} \SigmatTan_{\beta} 
	- 
	\partial_{\beta} \SigmatTan_{\alpha}
	$,
	we have:
	\begin{subequations}
		\begin{align} \label{E:ANTISYMMETRICNULLCURRENTSPACETIMERRORTERM}
			\mathfrak{J}_{(\textnormal{Antisymmetric})}[\pmb{\partial} \SigmatTan,\pmb{\partial} \SigmatTan]
			& \eqdef
				\frac{1}{2}
				\Nullhypersurfaceinversemetric^{\alpha \beta}
				\Nullhypersurfaceinversemetric^{\gamma \delta} 
				(\mathrm{d} \SigmatTan_{\flat})_{\alpha \gamma}
				(\mathrm{d} \SigmatTan_{\flat})_{\beta \delta},
						\\
			\mathfrak{J}_{(\textnormal{Div})}[\pmb{\partial} \SigmatTan,\pmb{\partial} \SigmatTan]
			& \eqdef
				\frac{9}{16}
				(\partial_a \SigmatTan^a)^2,
				\label{E:DIVERGENCENULLCURRENTSPACETIMERRORTERM} 
				\\
			\mathfrak{J}_{(\pmb{\partial} \weight)}[\SigmatTan,\pmb{\partial} \SigmatTan]
			& \eqdef
			-
			\ehcurrent^{\alpha}[\SigmatTan,\pmb{\partial} \SigmatTan] \partial_{\alpha} \weight,
				\label{E:DERIVATIVEOFWEIGHTNULLCURRENTSPACETIMERRORTERM} 
				\\
			\begin{split} 		\label{E:ABSORBABLEANTISYMMETRICANDDIVERGENCENULLCURRENTSPACETIMERRORTERM} 
			\mathfrak{J}_{(\textnormal{Absorb-1})}[\SigmatTan,\pmb{\partial} \SigmatTan]
			& \eqdef
				-
			\Nullhypersurfaceinversemetric^{\alpha \beta}
			\Nullhypersurfaceinversemetric^{\gamma \delta}
			\SigmatTan^{\kappa} 
			(\partial_{\alpha} \gfour_{\delta \kappa})
			(\mathrm{d} \SigmatTan_{\flat})_{\beta \gamma}
				\\
		& \ \
			+
			\frac{3}{8}
			(\partial_a \SigmatTan^a)
			\left\lbrace
				- 
				3	
				(\Transport \SigmatTan^{\alpha})
				\uLunit_{\alpha}
				+
				(\Lunit \SigmatTan^{\alpha})
				\uLunit_{\alpha}
				+
				\smoothtorusproject_b^{\ a}
				\partial_a \SigmatTan^b
			\right\rbrace,
		\end{split}
					\\
			\mathfrak{J}_{(\textnormal{Absorb-2})}[\SigmatTan,\pmb{\partial} \SigmatTan]
			& \eqdef
				-
			\Chfour_{\alpha \ \beta}^{\ \alpha} 
			\ehcurrent^{\beta}[\SigmatTan,\pmb{\partial} \SigmatTan]
			-
			\Nullhypersurfaceinversemetric^{\alpha \beta}
			\Nullhypersurfaceproject_{\gamma}^{\ \delta}
			(\partial_{\beta} \gfour_{\delta \kappa}) 
			\SigmatTan^{\kappa}
			\partial_{\alpha} \SigmatTan^{\gamma}
			+
			\Nullhypersurfaceinversemetric^{\alpha \beta}
			\Nullhypersurfaceproject_{\gamma}^{\ \delta}
			\SigmatTan^{\kappa}
			(\partial_{\delta} \gfour_{\beta \kappa})
			\partial_{\alpha} \SigmatTan^{\gamma},
				\label{E:EASYABSORBABLENULLCURRENTSPACETIMERRORTERM} 
				\\
		\mathfrak{J}_{(\textnormal{Material})}[\pmb{\partial} \SigmatTan,\pmb{\partial} \SigmatTan] 
			& \eqdef
			4
			\Nullhypersurfacemetric_{\alpha \beta}
			(\Transport \SigmatTan^{\alpha})
			\Transport \SigmatTan^{\beta},
					\label{E:QUADRATICINMATERIALDERIVATIVESNULLCURRENTSPACETIMERRORTERM} 
						\\
		\mathfrak{J}_{(\textnormal{Null Geometry})}[\SigmatTan,\pmb{\partial} \SigmatTan]
		& \eqdef
			-
		\SigmatTan^{\gamma} 
		\left\lbrace
			\partial_{\alpha}
			\left(
				\Nullhypersurfaceproject_{\gamma}^{\ \lambda} 
				\Nullhypersurfaceproject_{\kappa}^{\ \alpha} 
			\right)
		\right\rbrace
		\partial_{\lambda} \SigmatTan^{\kappa}
		+
		\SigmatTan^{\gamma} 
		\left\lbrace
			\partial_{\alpha}
			\left(
				\Nullhypersurfaceproject_{\gamma}^{\ \alpha} 
				\Nullhypersurfaceproject_{\lambda}^{\ \kappa}
			\right)
		\right\rbrace
		\partial_{\kappa} 
		\SigmatTan^{\lambda},
		\label{E:DERIVATIVESOFNULLGEOMETRYNULLCURRENTSPACETIMERRORTERM}
		\end{align}
	\end{subequations}
	and on RHS~\eqref{E:EASYABSORBABLENULLCURRENTSPACETIMERRORTERM},
	$\Chfour_{\alpha \ \beta}^{\ \gamma}  
	\eqdef \frac{1}{2} (\gfour^{-1})^{\gamma \delta} 
	\left(\partial_{\alpha} \gfour_{\delta \beta} 
	+ 
	\partial_{\beta} \gfour_{\alpha \delta} 
	-
	\partial_{\delta} \gfour_{\alpha \beta}
	\right)$
	are the Cartesian Christoffel symbols of $\gfour$.

\end{lemma}

\begin{proof}
	Throughout the proof,
	we silently use the simple fact that $\Transport_{\alpha} \SigmatTan^{\alpha} = - \SigmatTan^0 = 0$
	and thus $\SigmatTan_{\alpha} \eqdef \gfour_{\alpha \beta} \SigmatTan^{\beta} = \hfour_{\alpha \beta} \SigmatTan^{\beta}$
	and 
	$(\gfour^{-1})^{\alpha \beta} \SigmatTan_{\beta} = (\hfour^{-1})^{\alpha \beta} \SigmatTan_{\beta}$.
	We also silently use the simple identity 
	$\partial_{\alpha} (\hfour^{-1})^{\beta \gamma} 
	= 
	- (\hfour^{-1})^{\beta \beta'} 
	(\hfour^{-1})^{\gamma \gamma'} 
	\partial_{\alpha} \hfour_{\beta' \gamma'}$
	and the fact that in Cartesian coordinates,
	$\partial_{\alpha} \gfour_{\beta \gamma} = \partial_{\alpha} \hfour_{\beta \gamma}$
	(this follows easily from \eqref{E:MATERIALDERIVATIVELOWEREDCARTESIANCOORDINATES} and \eqref{E:RIEMANNIANACOUSTICALMETRIC}). 
	Moreover, we frequently relabel indices from line to line whenever convenient.  
	We will also silently use the observations made in Remark~\ref{R:OMITTINGZEROCOMPONENTINCARTESIANCOORDINATES}
	and the symmetry of $\Nullhypersurfacemetric$ 
	and
	$\Nullhypersurfaceinversemetric$ 
	shown in Lemma~\ref{L:SIMPLEIDENTITIESFORELLIPTICHYPERBOLICTENSORFIELDS}.
	
	We start by using \eqref{E:PUTANGENTELLIPTICHYPERBOLICCURRENT}
	to compute that relative to the Cartesian coordinates, we have:
	\begin{align} 
	\begin{split} \label{E:FIRSTSTEPPROOFCOVARIANTDIVERGENCEIDENTITYFORELLIPTICHYPERBOLICCURRENT}
	\Dfour_{\alpha} \ehcurrent^{\alpha}[\SigmatTan,\pmb{\partial} \SigmatTan]
	& =
	\partial_{\alpha} \ehcurrent^{\alpha}[\SigmatTan,\pmb{\partial} \SigmatTan]
	+
	\Chfour_{\alpha \ \beta}^{\ \alpha} 
	\ehcurrent^{\beta}[\SigmatTan,\pmb{\partial} \SigmatTan]
		\\
	& 
	=
	\Nullhypersurfaceproject_{\kappa}^{\ \alpha}
	\Nullhypersurfaceproject_{\gamma}^{\ \lambda} 
	(\partial_{\alpha} \SigmatTan^{\gamma})
	\partial_{\lambda} \SigmatTan^{\kappa}
	-
	(\Nullhypersurfaceproject_{\beta}^{\ \alpha} 
	\partial_{\alpha} \SigmatTan^{\beta})^2
		\\
	& \ \
		+
		\Chfour_{\alpha \ \beta}^{\ \alpha} 
		\ehcurrent^{\beta}[\SigmatTan,\pmb{\partial} \SigmatTan]
		+
		\SigmatTan^{\gamma} 
		\left\lbrace
			\partial_{\alpha}
			\left(
				\Nullhypersurfaceproject_{\gamma}^{\ \lambda} 
				\Nullhypersurfaceproject_{\kappa}^{\ \alpha} 
			\right)
		\right\rbrace
		\partial_{\lambda} \SigmatTan^{\kappa}
		-
		\SigmatTan^{\gamma} 
		\left\lbrace
			\partial_{\alpha}
			\left(
				\Nullhypersurfaceproject_{\gamma}^{\ \alpha} 
				\Nullhypersurfaceproject_{\lambda}^{\ \kappa}
			\right)
		\right\rbrace
		\partial_{\kappa} 
		\SigmatTan^{\lambda},
	\end{split}
	\end{align}
	where we stress that due to cancellations, the second derivatives of $\SigmatTan$ are absent from 
	RHS~\eqref{E:FIRSTSTEPPROOFCOVARIANTDIVERGENCEIDENTITYFORELLIPTICHYPERBOLICCURRENT}.
	Next, we compute the following identity:
	\begin{align}
	\begin{split}  \label{E:SECONDSTEPPROOFCOVARIANTDIVERGENCEIDENTITYFORELLIPTICHYPERBOLICCURRENT}
		\partial_{\lambda} \SigmatTan^{\kappa}
		& = 
		\partial_{\lambda}[(\hfour^{-1})^{\kappa \sigma} \SigmatTan_{\sigma}]
		=
		(\hfour^{-1})^{\kappa \sigma}
		\partial_{\lambda} \SigmatTan_{\sigma}
		+
		[\partial_{\lambda}(\hfour^{-1})^{\kappa \sigma}]
		\SigmatTan_{\sigma}
			\\
	& =
		(\hfour^{-1})^{\kappa \sigma}
		\partial_{\sigma} \SigmatTan_{\lambda}
		+
		(\hfour^{-1})^{\kappa \sigma}
		(\partial_{\lambda} \SigmatTan_{\sigma} 
		-
		\partial_{\sigma} \SigmatTan_{\lambda}	
		)
		+
		[\partial_{\lambda}(\hfour^{-1})^{\kappa \sigma}]
		\SigmatTan_{\sigma}
				\\
	& 
	=
		(\hfour^{-1})^{\kappa \sigma}
		\partial_{\sigma} 
		(\hfour_{\lambda \lambda'} \SigmatTan^{\lambda'})
		+
		(\hfour^{-1})^{\kappa \sigma}
		(\partial_{\lambda} \SigmatTan_{\sigma} 
		-
		\partial_{\sigma} \SigmatTan_{\lambda}	
		)
		-
		(\hfour^{-1})^{\kappa \kappa'}
		\SigmatTan^{\sigma}
		(\partial_{\lambda} \hfour_{\kappa' \sigma})
				\\
	& =
		(\hfour^{-1})^{\kappa \sigma}
		\hfour_{\lambda \lambda'} 
		(\partial_{\sigma} \SigmatTan^{\lambda'})
	+
	(\hfour^{-1})^{\kappa \sigma}
	(\partial_{\lambda} \SigmatTan_{\sigma} 
		-
	\partial_{\sigma} \SigmatTan_{\lambda}	
	)
		+
		(\hfour^{-1})^{\kappa \sigma}
		(\partial_{\sigma} \hfour_{\lambda \lambda'}) 
		\SigmatTan^{\lambda'}
	-
	(\hfour^{-1})^{\kappa \kappa'}
	\SigmatTan^{\sigma}
	(\partial_{\lambda} \hfour_{\kappa' \sigma}).
	\end{split}
	\end{align}
	Next, taking into account the definition \eqref{E:NULLHYPERSURFACESINVERSERIEMANNIANMETRIC},
	we compute that the contribution of the second product 
	$
	(\hfour^{-1})^{\kappa \sigma}
	(\partial_{\lambda} \SigmatTan_{\sigma} 
		-
	\partial_{\sigma} \SigmatTan_{\lambda}	
	)
	$
	on
	RHS~\eqref{E:SECONDSTEPPROOFCOVARIANTDIVERGENCEIDENTITYFORELLIPTICHYPERBOLICCURRENT}
	to the first product on RHS~\eqref{E:FIRSTSTEPPROOFCOVARIANTDIVERGENCEIDENTITYFORELLIPTICHYPERBOLICCURRENT}
	is as follows:
	\begin{align} 
	\begin{split} \label{E:THIRDSTEPPROOFCOVARIANTDIVERGENCEIDENTITYFORELLIPTICHYPERBOLICCURRENT}
		&
		\Nullhypersurfaceproject_{\kappa}^{\ \alpha}
		\Nullhypersurfaceproject_{\gamma}^{\ \lambda} 
		(\partial_{\alpha} \SigmatTan^{\gamma})
		(\hfour^{-1})^{\kappa \sigma}
		(\partial_{\lambda} \SigmatTan_{\sigma} 
			-
		\partial_{\sigma} \SigmatTan_{\lambda}	
		)
			\\
		& = 
		\Nullhypersurfaceinversemetric^{\alpha \sigma}
		\Nullhypersurfaceproject_{\gamma}^{\ \lambda} 
		\left\lbrace
			\partial_{\alpha} [(\hfour^{-1})^{\gamma \gamma'} \SigmatTan_{\gamma'}]
		\right\rbrace
		(\mathrm{d} \SigmatTan_{\flat})_{\lambda \sigma}
				\\
		& = 
			\Nullhypersurfaceinversemetric^{\alpha \beta}
			\Nullhypersurfaceinversemetric^{\gamma \delta}
			(\partial_{\alpha} \SigmatTan_{\gamma})
			(\mathrm{d} \SigmatTan_{\flat})_{\delta \beta} 
			-
			\Nullhypersurfaceinversemetric^{\alpha \beta}
			\Nullhypersurfaceinversemetric^{\gamma \delta}
			\SigmatTan^{\kappa} 
			(\partial_{\alpha} \hfour_{\delta \kappa})
			(\mathrm{d} \SigmatTan_{\flat})_{\gamma \beta}
				\\
		& = 
			-
			\frac{1}{2}
			\Nullhypersurfaceinversemetric^{\alpha \beta}
			\Nullhypersurfaceinversemetric^{\gamma \delta}
			(\mathrm{d} \SigmatTan_{\flat})_{\alpha \gamma} 
			(\mathrm{d} \SigmatTan_{\flat})_{\beta \delta} 
			+
			\Nullhypersurfaceinversemetric^{\alpha \beta}
			\Nullhypersurfaceinversemetric^{\gamma \delta}
			\SigmatTan^{\kappa} 
			(\partial_{\alpha} \hfour_{\delta \kappa})
			(\mathrm{d} \SigmatTan_{\flat})_{\beta \gamma}.
	\end{split}
	\end{align}
	Next, using 
	\eqref{E:NULLHYPERSURFACESRIEMANNIANMETRIC},
	\eqref{E:NULLHYPERSURFACESINVERSERIEMANNIANMETRIC},
	\eqref{E:SECONDSTEPPROOFCOVARIANTDIVERGENCEIDENTITYFORELLIPTICHYPERBOLICCURRENT} 
	and \eqref{E:THIRDSTEPPROOFCOVARIANTDIVERGENCEIDENTITYFORELLIPTICHYPERBOLICCURRENT},
	we rewrite the first product on RHS~\eqref{E:FIRSTSTEPPROOFCOVARIANTDIVERGENCEIDENTITYFORELLIPTICHYPERBOLICCURRENT}
	as follows:
	\begin{align} 
	\begin{split} \label{E:FOURTHSTEPPROOFCOVARIANTDIVERGENCEIDENTITYFORELLIPTICHYPERBOLICCURRENT}
		\Nullhypersurfaceproject_{\kappa}^{\ \alpha}
		\Nullhypersurfaceproject_{\gamma}^{\ \lambda} 
		(\partial_{\alpha} \SigmatTan^{\gamma})
		\partial_{\lambda} \SigmatTan^{\kappa}
		& =
			\Nullhypersurfaceinversemetric^{\alpha \beta}
			\Nullhypersurfacemetric_{\gamma \delta}
			(\partial_{\alpha} \SigmatTan^{\gamma})
			\partial_{\beta} \SigmatTan^{\delta}
			-
			\frac{1}{2}
			\Nullhypersurfaceinversemetric^{\alpha \beta}
			\Nullhypersurfaceinversemetric^{\gamma \delta}
			(\mathrm{d} \SigmatTan_{\flat})_{\alpha \gamma} 
			(\mathrm{d} \SigmatTan_{\flat})_{\beta \delta}
				\\
		& \ \
			+
			\Nullhypersurfaceinversemetric^{\alpha \beta}
			\Nullhypersurfaceinversemetric^{\gamma \delta}
			\SigmatTan^{\kappa} 
			(\partial_{\alpha} \hfour_{\delta \kappa})
			(\mathrm{d} \SigmatTan_{\flat})_{\beta \gamma}
					\\
		& \ \
			+
			\Nullhypersurfaceinversemetric^{\alpha \beta}
			\Nullhypersurfaceproject_{\gamma}^{\ \delta}
			(\partial_{\beta} \hfour_{\delta \kappa}) 
			\SigmatTan^{\kappa}
			\partial_{\alpha} \SigmatTan^{\gamma}
			-
			\Nullhypersurfaceinversemetric^{\alpha \beta}
			\Nullhypersurfaceproject_{\gamma}^{\ \delta}
			\SigmatTan^{\kappa}
			(\partial_{\delta} \hfour_{\beta \kappa})
			\partial_{\alpha} \SigmatTan^{\gamma}.
	\end{split}
	\end{align}
	Next, using (see Lemma~\ref{L:BASICPROPERTIESOFULUNIT}) the identities
	$\Transport = \frac{1}{2} (\Lunit + \uLunit)$,
	$\Lunit = \Transport - X$,
	$\uLunit = \Transport + X$,
	$\Transport_{\alpha} \SigmatTan^{\alpha} = 0$,
	and 
	$\Transport_{\alpha} \partial_{\beta} \SigmatTan^{\alpha} = 0$,
	as well as the identity 
	$
	\partial_a \SigmatTan^a 
	= 
	X_a X \SigmatTan^a 
	+ 
	\smoothtorusproject_b^{\ a}
	\partial_a \SigmatTan^b
	$
	(see \eqref{E:SMOOTHTORUSPROJECT}),
	we compute that:
	\begin{align} \label{E:FIFTHSTEPPROOFCOVARIANTDIVERGENCEIDENTITYFORELLIPTICHYPERBOLICCURRENT}
		\Lunit_{\alpha}
		\uLunit \SigmatTan^{\alpha}
		& = 
		-
		\uLunit_{\alpha}
		\uLunit \SigmatTan^{\alpha}
		=
		\uLunit_{\alpha}
		\Lunit \SigmatTan^{\alpha}
		-
		2
		\uLunit_{\alpha}
		\Transport \SigmatTan^{\alpha},
			\\
		\Lunit_{\alpha}
		\uLunit \SigmatTan^{\alpha}
		& = 
		\Lunit_{\alpha}
		X \SigmatTan^{\alpha}
		+
		\Lunit_{\alpha}
		\Transport \SigmatTan^{\alpha}
		= 
		-
		X_a
		X \SigmatTan^a
		+
		\Lunit_{\alpha}
		\Transport \SigmatTan^{\alpha}
		=
		\smoothtorusproject_b^{\ a}
		\partial_a \SigmatTan^b
		-
		\partial_a \SigmatTan^a
		-
		\uLunit_{\alpha}
		\Transport \SigmatTan^{\alpha}.
		\label{E:SIXTHSTEPPROOFCOVARIANTDIVERGENCEIDENTITYFORELLIPTICHYPERBOLICCURRENT}
	\end{align}
	Next, using \eqref{E:NULLHYPERSURFACEPROJECTION},
	\eqref{E:FIFTHSTEPPROOFCOVARIANTDIVERGENCEIDENTITYFORELLIPTICHYPERBOLICCURRENT},
	and
	\eqref{E:SIXTHSTEPPROOFCOVARIANTDIVERGENCEIDENTITYFORELLIPTICHYPERBOLICCURRENT},
	we compute that:
	\begin{align} \label{E:SEVENTHSTEPPROOFCOVARIANTDIVERGENCEIDENTITYFORELLIPTICHYPERBOLICCURRENT}
		\Nullhypersurfaceproject_{\beta}^{\ \alpha} 
		\partial_{\alpha} \SigmatTan^{\beta}
		& = 
				\partial_a \SigmatTan^a
				+
				\frac{1}{4}
				\Lunit_{\alpha}
				\uLunit \SigmatTan^{\alpha}
				+
				\frac{1}{4}
				\Lunit_{\alpha}
				\uLunit \SigmatTan^{\alpha}
					\\
		& = 
			\frac{3}{4}
	\partial_a \SigmatTan^a
	-
	\frac{3}{4}
	(\Transport \SigmatTan^{\alpha})
	\uLunit_{\alpha}
	+
	\frac{1}{4}
	(\Lunit \SigmatTan^{\alpha})
	\uLunit_{\alpha}
	+
	\frac{1}{4}
	\smoothtorusproject_b^{\ a}
	\partial_a \SigmatTan^b.
		\label{E:EIGHTHSTEPPROOFCOVARIANTDIVERGENCEIDENTITYFORELLIPTICHYPERBOLICCURRENT}
	\end{align}
	Using 
	\eqref{E:SEVENTHSTEPPROOFCOVARIANTDIVERGENCEIDENTITYFORELLIPTICHYPERBOLICCURRENT}--\eqref{E:EIGHTHSTEPPROOFCOVARIANTDIVERGENCEIDENTITYFORELLIPTICHYPERBOLICCURRENT},
	we rewrite the second product on RHS~\eqref{E:FIRSTSTEPPROOFCOVARIANTDIVERGENCEIDENTITYFORELLIPTICHYPERBOLICCURRENT}
	as follows:
	\begin{align} \label{E:NINTHSTEPPROOFCOVARIANTDIVERGENCEIDENTITYFORELLIPTICHYPERBOLICCURRENT}
		(\Nullhypersurfaceproject_{\beta}^{\ \alpha} 
		\partial_{\alpha} \SigmatTan^{\beta})^2
		& = 
			\frac{1}{16}
			\left\lbrace
				3
				\partial_a \SigmatTan^a
				-	
				3
				(\Transport \SigmatTan^{\alpha})
				\uLunit_{\alpha}
				+
				(\Lunit \SigmatTan^{\alpha})
				\uLunit_{\alpha}
				+
				\smoothtorusproject_b^{\ a}
				\partial_a \SigmatTan^b
		\right\rbrace^2.
	\end{align}
	Combining
	\eqref{E:FIRSTSTEPPROOFCOVARIANTDIVERGENCEIDENTITYFORELLIPTICHYPERBOLICCURRENT}--\eqref{E:SIXTHSTEPPROOFCOVARIANTDIVERGENCEIDENTITYFORELLIPTICHYPERBOLICCURRENT}
	and \eqref{E:NINTHSTEPPROOFCOVARIANTDIVERGENCEIDENTITYFORELLIPTICHYPERBOLICCURRENT}, 
	we deduce the following identity:
	\begin{align} 
	\begin{split} \label{E:TENTHSTEPPROOFCOVARIANTDIVERGENCEIDENTITYFORELLIPTICHYPERBOLICCURRENT}
			\Nullhypersurfaceinversemetric^{\alpha \beta}
			\Nullhypersurfacemetric_{\gamma \delta}
			(\partial_{\alpha} \SigmatTan^{\gamma})
			\partial_{\beta} \SigmatTan^{\delta}
		& = 
			\Dfour_{\alpha} \ehcurrent^{\alpha}[\SigmatTan,\pmb{\partial} \SigmatTan]
				\\
	& \ \
			+
			\frac{1}{2}
			\Nullhypersurfaceinversemetric^{\alpha \beta}
			\Nullhypersurfaceinversemetric^{\gamma \delta}
			(\mathrm{d} \SigmatTan_{\flat})_{\alpha \gamma} 
			(\mathrm{d} \SigmatTan_{\flat})_{\beta \delta}
			+
			\frac{1}{16}
			\left\lbrace
				3
				\partial_a \SigmatTan^a
				-	
				3
				(\Transport \SigmatTan^{\alpha})
				\uLunit_{\alpha}
				+
				(\Lunit \SigmatTan^{\alpha})
				\uLunit_{\alpha}
				+
				\smoothtorusproject_b^{\ a}
				\partial_a \SigmatTan^b
		\right\rbrace^2
				 \\
		& \ \
		-
			\Nullhypersurfaceinversemetric^{\alpha \beta}
			\Nullhypersurfaceinversemetric^{\gamma \delta}
			\SigmatTan^{\kappa} 
			(\partial_{\alpha} \hfour_{\delta \kappa})
			(\mathrm{d} \SigmatTan_{\flat})_{\beta \gamma}
				\\
		& \ \
		-
		\SigmatTan^{\gamma} 
		\left\lbrace
			\partial_{\alpha}
			\left(
				\Nullhypersurfaceproject_{\gamma}^{\ \lambda} 
				\Nullhypersurfaceproject_{\kappa}^{\ \alpha} 
			\right)
		\right\rbrace
		\partial_{\lambda} \SigmatTan^{\kappa}
		+
		\SigmatTan^{\gamma} 
		\left\lbrace
			\partial_{\alpha}
			\left(
				\Nullhypersurfaceproject_{\gamma}^{\ \alpha} 
				\Nullhypersurfaceproject_{\lambda}^{\ \kappa}
			\right)
		\right\rbrace
		\partial_{\kappa} 
		\SigmatTan^{\lambda}
				\\
		& \ \
			-
			\Chfour_{\alpha \ \beta}^{\ \alpha} 
			\ehcurrent^{\beta}[\SigmatTan,\pmb{\partial} \SigmatTan]
			-
			\Nullhypersurfaceinversemetric^{\alpha \beta}
			\Nullhypersurfaceproject_{\gamma}^{\ \delta}
			(\partial_{\beta} \hfour_{\delta \kappa}) 
			\SigmatTan^{\kappa}
			\partial_{\alpha} \SigmatTan^{\gamma}
			+
			\Nullhypersurfaceinversemetric^{\alpha \beta}
			\Nullhypersurfaceproject_{\gamma}^{\ \delta}
			\SigmatTan^{\kappa}
			(\partial_{\delta} \hfour_{\beta \kappa})
			\partial_{\alpha} \SigmatTan^{\gamma}.
	\end{split}
	\end{align}
	Next, using \eqref{E:NULLHYPERSURFACEADAPTEDCOERCIVEQUADRATICFORM},
	we observe that:
	\begin{align} 
	\begin{split} \label{E:ELEVENTHSTEPPROOFCOVARIANTDIVERGENCEIDENTITYFORELLIPTICHYPERBOLICCURRENT}
			\ellipticCoerciveQuadratic[\pmb{\partial} \SigmatTan,\pmb{\partial} \SigmatTan]
			& =
			\Nullhypersurfaceinversemetric^{\alpha \beta}
			\Nullhypersurfacemetric_{\gamma \delta}
			(\partial_{\alpha} \SigmatTan^{\gamma})
			\partial_{\beta} \SigmatTan^{\delta}
			+
			4
			\Nullhypersurfacemetric_{\alpha \beta}
			(\Transport \SigmatTan^{\alpha})
			\Transport \SigmatTan^{\beta}
				\\
	& \ \
		+
		4
		\left(
			\Nullhypersurfaceinversemetric^{\alpha \beta}
			+
			4
			\Transport^{\alpha} \Transport^{\beta}
		\right)
		\Transport_{\gamma} \Transport_{\delta}
		(\partial_{\alpha} \SigmatTan^{\gamma})
		\partial_{\beta} \SigmatTan^{\delta}
			\\
	& \ \
			-
		\frac{1}{16}
		\left\lbrace
		- 
			3	
			(\Transport \SigmatTan^{\alpha})
			\uLunit_{\alpha}
			+
			(\Lunit \SigmatTan^{\alpha})
			\uLunit_{\alpha}
			+
			\smoothtorusproject_b^{\ a}
			\partial_a \SigmatTan^b
	\right\rbrace^2.
	\end{split}
	\end{align}
	Hence, noting that
	the terms on the second line of RHS~\eqref{E:ELEVENTHSTEPPROOFCOVARIANTDIVERGENCEIDENTITYFORELLIPTICHYPERBOLICCURRENT}
	vanish, we can add
	$$
	4
			\Nullhypersurfacemetric_{\alpha \beta}
			(\Transport \SigmatTan^{\alpha})
			\Transport \SigmatTan^{\beta}
			-
		\frac{1}{16}
		\left\lbrace
		- 
			3	
			(\Transport \SigmatTan^{\alpha})
			\uLunit_{\alpha}
			+
			(\Lunit \SigmatTan^{\alpha})
			\uLunit_{\alpha}
			+
			\smoothtorusproject_b^{\ a}
			\partial_a \SigmatTan^b
	\right\rbrace^2
	$$
	to each side of \eqref{E:TENTHSTEPPROOFCOVARIANTDIVERGENCEIDENTITYFORELLIPTICHYPERBOLICCURRENT}
	to obtain the following identity:
	\begin{align} 
	\begin{split} \label{E:TWELFTHSTEPPROOFCOVARIANTDIVERGENCEIDENTITYFORELLIPTICHYPERBOLICCURRENT}
			\ellipticCoerciveQuadratic[\pmb{\partial} \SigmatTan,\pmb{\partial} \SigmatTan]
		& = 
			\Dfour_{\alpha} \ehcurrent^{\alpha}[\SigmatTan,\pmb{\partial} \SigmatTan]
				\\
	& \ \
			+
			\frac{1}{2}
			\Nullhypersurfaceinversemetric^{\alpha \beta}
			\Nullhypersurfaceinversemetric^{\gamma \delta}
			(\mathrm{d} \SigmatTan_{\flat})_{\alpha \gamma} 
			(\mathrm{d} \SigmatTan_{\flat})_{\beta \delta}
			+
			\frac{9}{16}
			(\partial_a \SigmatTan^a)^2	
			 \\
	& \ \
			-
			\Nullhypersurfaceinversemetric^{\alpha \beta}
			\Nullhypersurfaceinversemetric^{\gamma \delta}
			\SigmatTan^{\kappa} 
			(\partial_{\alpha} \hfour_{\delta \kappa})
			(\mathrm{d} \SigmatTan_{\flat})_{\beta \gamma}
			+
			\frac{3}{8}
			(\partial_a \SigmatTan^a)
			\left\lbrace
				- 
				3	
				(\Transport \SigmatTan^{\alpha})
				\uLunit_{\alpha}
				+
				(\Lunit \SigmatTan^{\alpha})
				\uLunit_{\alpha}
				+
				\smoothtorusproject_b^{\ a}
				\partial_a \SigmatTan^b
			\right\rbrace
				 \\
		& \ \
			+
			4
			\Nullhypersurfacemetric_{\alpha \beta}
			(\Transport \SigmatTan^{\alpha})
			\Transport \SigmatTan^{\beta}
				\\
		& \ \
		-
		\SigmatTan^{\gamma} 
		\left\lbrace
			\partial_{\alpha}
			\left(
				\Nullhypersurfaceproject_{\gamma}^{\ \lambda} 
				\Nullhypersurfaceproject_{\kappa}^{\ \alpha} 
			\right)
		\right\rbrace
		\partial_{\lambda} \SigmatTan^{\kappa}
		+
		\SigmatTan^{\gamma} 
		\left\lbrace
			\partial_{\alpha}
			\left(
				\Nullhypersurfaceproject_{\gamma}^{\ \alpha} 
				\Nullhypersurfaceproject_{\lambda}^{\ \kappa}
			\right)
		\right\rbrace
		\partial_{\kappa} 
		\SigmatTan^{\lambda}
			 \\
		& \ \
			-
			\Chfour_{\alpha \ \beta}^{\ \alpha} 
			\ehcurrent^{\beta}[\SigmatTan,\pmb{\partial} \SigmatTan]
			-
			\Nullhypersurfaceinversemetric^{\alpha \beta}
			\Nullhypersurfaceproject_{\gamma}^{\ \delta}
			(\partial_{\beta} \hfour_{\delta \kappa}) 
			\SigmatTan^{\kappa}
			\partial_{\alpha} \SigmatTan^{\gamma}
			+
			\Nullhypersurfaceinversemetric^{\alpha \beta}
			\Nullhypersurfaceproject_{\gamma}^{\ \delta}
			\SigmatTan^{\kappa}
			(\partial_{\delta} \hfour_{\beta \kappa})
			\partial_{\alpha} \SigmatTan^{\gamma}.
	\end{split}
	\end{align}
	We have therefore proved \eqref{E:COVARIANTDIVERGENCEIDENTITYFORELLIPTICHYPERBOLICCURRENT}
	in the special case $\weight \eqdef 1$. To obtain \eqref{E:COVARIANTDIVERGENCEIDENTITYFORELLIPTICHYPERBOLICCURRENT}
	for a general weight $\weight$, we simply multiply this special case identity by $\weight$ 
	and use the commutation identity
	$
	\weight
	\Dfour_{\alpha}  \ehcurrent^{\alpha}[\SigmatTan,\pmb{\partial} \SigmatTan]
	=
	\Dfour_{\alpha} (\weight \ehcurrent^{\alpha}[\SigmatTan,\pmb{\partial} \SigmatTan])
	+
	\mathfrak{J}_{(\pmb{\partial} \weight)}[\SigmatTan,\pmb{\partial} \SigmatTan]
	$,
	where $\mathfrak{J}_{(\pmb{\partial} \weight)}[\SigmatTan,\pmb{\partial} \SigmatTan]$ is defined by 
	\eqref{E:DERIVATIVEOFWEIGHTNULLCURRENTSPACETIMERRORTERM}.
	
\end{proof}

\subsubsection{Key identity for the elliptic-hyperbolic boundary terms}
\label{SSS:KEYIDENTITYFORTHEELLIPTICHYPERBOLICOUNDARYTERMS}
To derive our main elliptic-hyperbolic integral identities, 
which we state as Prop.\,\ref{P:INTEGRALIDENTITYFORELLIPTICHYPERBOLICCURRENT},
we will start by integrating
\eqref{E:COVARIANTDIVERGENCEIDENTITYFORELLIPTICHYPERBOLICCURRENT} 
over the spacetime region $\twoargMrough{[\timefunction_1,\timefunction_2),[u_1,u_2]}{\muxmulevelsetvalue}$
and applying the divergence theorem.
This results in boundary integrals, including the integral of
$\hypunitnormalarg{\muxmulevelsetvalue}_{\alpha}
		\ehcurrent^{\alpha}[\SigmatTan,\pmb{\partial} \SigmatTan]$
over $\hypthreearg{\timefunction_2}{[u_1,u_2]}{\muxmulevelsetvalue}$,
where $\hypunitnormalarg{\muxmulevelsetvalue}$ is the future-directed $\gfour$-unit normal to 
$\hypthreearg{\timefunction_2}{[u_1,u_2]}{\muxmulevelsetvalue}$.
To avoid uncontrollable error terms in the boundary integral, 
we will integrate by parts over
$\hypthreearg{\timefunction_2}{[u_1,u_2]}{\muxmulevelsetvalue}$
with the help of the identity for
$
\hypunitnormalarg{\muxmulevelsetvalue}_{\alpha}
		\ehcurrent^{\alpha}[\SigmatTan,\pmb{\partial} \SigmatTan]
$
provided by the next lemma. Of crucial importance for 
our top-order $L^2$ estimates
is the sign of the first two terms
$
\Rtransarg{\muxmulevelsetvalue} 
			\CurrentboundaryerrorperfectRderivative[\SigmatTan,\SigmatTan]
			+
			\frac{1}{2} \CurrentboundaryerrorperfectRderivative[\SigmatTan,\SigmatTan]
			\mytr_{\gtorusroughfirstfund} \deform{\Rtransarg{\muxmulevelsetvalue}}
$
on RHS~\eqref{E:KEYIDPUTANGENTCURRENTCONTRACTEDAGAINSTVECTORFIELD};
in the proof of the integral identity \eqref{E:INTEGRALIDENTITYFORELLIPTICHYPERBOLICCURRENT}, 
we will integrate by parts and exploit the sign of these terms as well as 
the positive definiteness of the quadratic form $\CurrentboundaryerrorperfectRderivative[\SigmatTan,\SigmatTan]$
shown in \eqref{E:KEYCOERCIVITYPERFECTRDERIVAIVEERRORPUTANGENTCURRENTCONTRACTEDAGAINSTVECTORFIELD}.

\begin{lemma}[Key identity for the elliptic-hyperbolic boundary terms]
	\label{L:KEYIDPUTANGENTCURRENTCONTRACTEDAGAINSTVECTORFIELD}
	Let $\SigmatTan$ be a $\Sigma_t$-tangent vectorfield, 
	let $\angV$ be its $\gfour$-orthogonal projection onto the smooth tori $\ell_{t,u}$,
	let $(\mathrm{d} V_{\flat})_{\alpha \beta} \eqdef \partial_{\alpha} V_{\beta} - \partial_{\beta} V_{\alpha}$,
	and let
	$\ehcurrent^{\alpha}[\SigmatTan,\pmb{\partial} \SigmatTan]$
	be the corresponding characteristic current defined in
	\eqref{E:PUTANGENTELLIPTICHYPERBOLICCURRENT}.
	Then relative to the Cartesian coordinates,
	the following identity holds, where $\smoothtorusproject$ is the $\ell_{t,u}$-projection from Def.\,\ref{D:PROJECTIONTENSORFIELDSANDTANGENCYTOHYPERSURFACES},
	$\angD$ is as in Def.\,\ref{D:CONNECTIONSANDDIFFERENTIALOPERATORS},
	$|\cdot|_{\gtorus}$ and $|\cdot|_g$ are as in Def.\,\ref{D:POINTWISESEMINORMS},
	$\phi = \phi(u)$ is the cut-off function from Def.\,\ref{D:WTRANSANDCUTOFF},
	$\phi'$ is its derivative,
	$\Roughtoritangentvectorfieldarg{\muxmulevelsetvalue}$,
	$\hypunitnormalarg{\muxmulevelsetvalue}$
	and
	$\Rtransarg{\muxmulevelsetvalue}$
	are as in Def.\,\ref{D:GEOMETRICVECTORFIELDSADAPTEDTOROUGHFOLIATIONS},
	$|\Rtransarg{\muxmulevelsetvalue}|_{\hypg}$ is as in \eqref{E:SIZEOFRTRANS},
	$\Rtransnormsmallfactorarg{\muxmulevelsetvalue}$ is as in \eqref{E:RTRANSNORMSMALLFACTOR},
	$(\gtorusroughinversefirstfund)^{AB}$ is as in
	Lemma~\ref{L:IDENTITIESFORCOMPONENTSOFROUGHTORUSFIRSTFUNDAMENTALFORMANDITSINVERSE},
	$\roughangdiv$ is the $\twoargroughtori{\timefunctionarg{\muxmulevelsetvalue},u}{\muxmulevelsetvalue}$-divergence operator from
	Def.\,\ref{D:CONNECTIONSANDDIFFERENTIALOPERATORSONROUGHTORI},
	and $\mytr_{\gtorusroughfirstfund} \deform{\Rtransarg{\muxmulevelsetvalue}}$ is the $\gtorusroughfirstfund$-trace 
	(see Def.\,\ref{D:TRACEOFROUGHTORITANGENT02TENSORS})
	of the deformation tensor $\deform{\Rtransarg{\muxmulevelsetvalue}}$ of $\Rtransarg{\muxmulevelsetvalue}$:
	\begin{align} 
	\begin{split} \label{E:KEYIDPUTANGENTCURRENTCONTRACTEDAGAINSTVECTORFIELD}
		\frac{|\Rtransarg{\muxmulevelsetvalue}|_{\hypg}}
		{\upmu}
		\hypunitnormalarg{\muxmulevelsetvalue}_{\alpha}
		\ehcurrent^{\alpha}[\SigmatTan,\pmb{\partial} \SigmatTan] 
		& = 
			\Rtransarg{\muxmulevelsetvalue} 
			\CurrentboundaryerrorperfectRderivative[\SigmatTan,\SigmatTan]
			+
			\frac{1}{2} \CurrentboundaryerrorperfectRderivative[\SigmatTan,\SigmatTan]
			\mytr_{\gtorusroughfirstfund} \deform{\Rtransarg{\muxmulevelsetvalue}}
				\\
		& \ \
			- 
			\roughangdiv
			\left\lbrace
			\frac{\upmu}{\upmu - \phi \frac{\muxmulevelsetvalue}{\Lunit \upmu}}
			\left[
				-
				\frac{1}{2}
				\frac{1}{\upmu}
				\Rtransarg{\muxmulevelsetvalue}_{\alpha} \SigmatTan^{\alpha} \SigmatTan^a X_a
				+
				\frac{1}{4}
				\frac{1}{\upmu}
				\Rtransarg{\muxmulevelsetvalue}_{\alpha} \uLunit^{\alpha}
				|\SigmatTan|_g^2
		\right]
		\Roughtoritangentvectorfieldarg{\muxmulevelsetvalue}
		\right\rbrace
				 \\
		& \ \
			+
			\Currentboundaryerrorhavetocontrolprincipal[\SigmatTan,\pmb{\partial} \SigmatTan]
			+
			\Currentboundaryerrorhavetocontrollowerorder[\SigmatTan,\SigmatTan],
	\end{split}
	\end{align}
	where:
	\begin{align} \label{E:PERFECTRDERIVAIVEERRORPUTANGENTCURRENTCONTRACTEDAGAINSTVECTORFIELD}
		\CurrentboundaryerrorperfectRderivative[\SigmatTan,\SigmatTan]
		& \eqdef 
				\frac{1}{4}
				\left\lbrace
				\frac{1}
				{\upmu - \phi \frac{\muxmulevelsetvalue}{\Lunit \upmu}}
				\left[
					(X_a \SigmatTan^a)^2
					+
					(1 + 2 \Rtransnormsmallfactorarg{\muxmulevelsetvalue}) |\SigmatTan|_{\gtorus}^2
				-
				2
				\frac{1}{\Lunit \timefunctionarg{\muxmulevelsetvalue}}
				X_a \SigmatTan^a
				\angVarg{\alpha}
				\angDarg{\alpha}
				\timefunctionarg{\muxmulevelsetvalue}
			\right]
			\right\rbrace,
	\end{align}	
	
	\begin{align}
	\begin{split}  \label{E:PRINCIPALERRORTERMHAVETOCONTROLKEYIDPUTANGENTCURRENTCONTRACTEDAGAINSTVECTORFIELD}
		&
		\Currentboundaryerrorhavetocontrolprincipal[\SigmatTan,\pmb{\partial} \SigmatTan]
		\eqdef
				\frac{1}{\upmu}
				\Rtransarg{\muxmulevelsetvalue}_{\alpha}
				\Nullhypersurfaceproject_{\beta}^{\ \alpha}
				\SigmatTan^{\gamma}  
				(\gfour^{-1})^{\beta \delta} (\mathrm{d} V_{\flat})_{\gamma \delta}
				-
				\frac{1}{\upmu}
				\SigmatTan^{\beta} 
				\Lunit_{\beta} 
				\Lunit^{\gamma} 
				\Rtransarg{\muxmulevelsetvalue}^{\delta}
				(\mathrm{d} V_{\flat})_{\gamma \delta} 
				\\
		& \ \
			- 
			\frac{1}{\upmu}
			\Rtransarg{\muxmulevelsetvalue}_{\alpha}
			\Nullhypersurfaceproject_{\beta}^{\ \alpha} \SigmatTan^{\beta}
			\partial_a \SigmatTan^a
			+
			\frac{1}{\upmu}
			\Rtransarg{\muxmulevelsetvalue}_{\alpha}
			\SigmatTan^{\beta} \Lunit_{\beta} \Transport \SigmatTan^{\alpha}
			-
			\frac{1}{\upmu}
			\Rtransarg{\muxmulevelsetvalue}_{\alpha}
			\SigmatTan^{\alpha} \Lunit_{\beta} \Transport \SigmatTan^{\beta}
				\\
		& \ \
			-
			\frac{1}{2(\upmu - \phi \frac{\muxmulevelsetvalue}{\Lunit \upmu})}
			\Rtransarg{\muxmulevelsetvalue}_{\alpha} (\Transport \SigmatTan^{\alpha}) \SigmatTan^a X_a
			-
			\frac{1}{2(\upmu - \phi \frac{\muxmulevelsetvalue}{\Lunit \upmu})}
			\Rtransarg{\muxmulevelsetvalue}_{\alpha} \SigmatTan^{\alpha} (\Transport \SigmatTan^a) X_a
			+
			\frac{1}{2(\upmu - \phi \frac{\muxmulevelsetvalue}{\Lunit \upmu})}
			\Rtransarg{\muxmulevelsetvalue}_{\alpha} \uLunit^{\alpha}
			\SigmatTan_a \Transport \SigmatTan^a,
\end{split}
\end{align}
and:
\begin{align} 
\begin{split} \label{E:LOWERORDERERRORTERMHAVETOCONTROLKEYIDPUTANGENTCURRENTCONTRACTEDAGAINSTVECTORFIELD}
	& \Currentboundaryerrorhavetocontrollowerorder[\SigmatTan,\SigmatTan]
	 \eqdef
			-
			\frac{1}{2} 
			\CurrentboundaryerrorperfectRderivative[\SigmatTan,\SigmatTan]
			\mytr_{\gtorusroughfirstfund} \deform{\Rtransarg{\muxmulevelsetvalue}}
				\\
	& \ \
		+
		\frac{1}{2}
		\frac{(\Rtransarg{\muxmulevelsetvalue} \upmu)}{\upmu^2}
		|\angV|_{\gtorus}^2
		+
		\frac{1}{2}
		\frac{(\Rtransarg{\muxmulevelsetvalue} \upmu) \phi \frac{\muxmulevelsetvalue}{\Lunit \upmu}}{\upmu^2(\upmu - \phi \frac{\muxmulevelsetvalue}{\Lunit \upmu})}
		|\angV|_{\gtorus}^2
		+ 
		\frac{(\Rtransarg{\muxmulevelsetvalue} \upmu)}{\upmu (\upmu - \phi \frac{\muxmulevelsetvalue}{\Lunit \upmu})^2}
		\left\lbrace
				\frac{1}{2}
				\Rtransarg{\muxmulevelsetvalue}_{\alpha} \SigmatTan^{\alpha} \SigmatTan^a X_a
				-
				\frac{1}{4}
				\Rtransarg{\muxmulevelsetvalue}_{\alpha} \uLunit^{\alpha}
				|\SigmatTan|_g^2
		\right\rbrace
			\\
	& \ \
		+
		\frac{(\Lunit \upmu)}{\upmu^2}
		\left\lbrace
			-
				\frac{1}{2}
				\Rtransarg{\muxmulevelsetvalue}_{\alpha} \SigmatTan^{\alpha} \SigmatTan^a X_a
				+
				\frac{1}{4}
				\Rtransarg{\muxmulevelsetvalue}_{\alpha} \uLunit^{\alpha}
				|\SigmatTan|_g^2
		\right\rbrace
		+
		\frac{1}{2}
		\frac{(\Lunit \upmu) \phi \frac{\muxmulevelsetvalue}{\Lunit \upmu}}{\upmu^2}
		|\angV|_{\gtorus}^2
			\\
		& \ \
			-
			\left\lbrace
				\frac{(\Roughtoritangentvectorfieldarg{\muxmulevelsetvalue} \upmu) \phi \frac{\muxmulevelsetvalue}{\Lunit \upmu}}
				{(\upmu - \phi \frac{\muxmulevelsetvalue}{\Lunit \upmu})^2}
				+
				\frac{\upmu \phi \muxmulevelsetvalue \frac{\Roughtoritangentvectorfieldarg{\muxmulevelsetvalue} \Lunit \upmu}
				{(\Lunit \upmu)^2}}{(\upmu - \phi \frac{\muxmulevelsetvalue}{\Lunit \upmu})^2}
			\right\rbrace
			\left\lbrace
				-
				\frac{1}{2}
				\frac{1}{\upmu}
				\Rtransarg{\muxmulevelsetvalue}_{\alpha} \SigmatTan^{\alpha} \SigmatTan^a X_a
				+
				\frac{1}{4}
				\frac{1}{\upmu}
				\Rtransarg{\muxmulevelsetvalue}_{\alpha} \uLunit^{\alpha}
				|\SigmatTan|_g^2
		\right\rbrace
			+
			\frac{1}{2}
			\frac{(\Roughtoritangentvectorfieldarg{\muxmulevelsetvalue} \upmu) \phi \frac{\muxmulevelsetvalue}{\Lunit \upmu}}{\upmu(\upmu - \phi \frac{\muxmulevelsetvalue}{\Lunit \upmu})}
			|\angV|_{\gtorus}^2
					\\
		& \ \
			+
			\left\lbrace
			\frac{1}{\upmu - \phi \frac{\muxmulevelsetvalue}{\Lunit \upmu}}
			\left[
				-
				\frac{1}{2}
				\Rtransarg{\muxmulevelsetvalue}_{\alpha} \SigmatTan^{\alpha} \SigmatTan^a X_a
				+
				\frac{1}{4}
				\Rtransarg{\muxmulevelsetvalue}_{\alpha} \uLunit^{\alpha}
				|\SigmatTan|_g^2
		\right]
		\right\rbrace
		\roughangdiv \Roughtoritangentvectorfieldarg{\muxmulevelsetvalue}
		-
		\frac{\phi \muxmulevelsetvalue \frac{(\Rtransarg{\muxmulevelsetvalue} \Lunit \upmu)}{(\Lunit \upmu)^2}}{\upmu (\upmu - \phi \frac{\muxmulevelsetvalue}{\Lunit \upmu})^2}
		\left\lbrace
				-
				\frac{1}{2}
				\Rtransarg{\muxmulevelsetvalue}_{\alpha} \SigmatTan^{\alpha} \SigmatTan^a X_a
				+
				\frac{1}{4}
				\Rtransarg{\muxmulevelsetvalue}_{\alpha} \uLunit^{\alpha}
				|\SigmatTan|_g^2
		\right\rbrace
			\\
	& \ \
		+
		\frac{1}{2}
		\phi \frac{\muxmulevelsetvalue (\muX \Lunit \upmu)}{\upmu(\upmu - \phi \frac{\muxmulevelsetvalue}{\Lunit \upmu}) (\Lunit \upmu)^2}
		|\angV|_{\gtorus}^2
		+
		\frac{1}{2}
		\phi \frac{\muxmulevelsetvalue (\Lunit \Lunit \upmu)}{(\upmu - \phi \frac{\muxmulevelsetvalue}{\Lunit \upmu}) (\Lunit \upmu)^2}
		|\angV|_{\gtorus}^2
		 \\
	& \ \
			+
			\frac{1}{4 \upmu (\upmu - \phi \frac{\muxmulevelsetvalue}{\Lunit \upmu})}
			\Rtransarg{\muxmulevelsetvalue}_{\alpha}
			(\muX \uLunit^{\alpha})
			|\SigmatTan|_g^2
			+
			\frac{1}{\upmu}
			(\Rtransarg{\muxmulevelsetvalue} \Lunit_{\alpha})
			\Lunit_{\beta}
			\SigmatTan^{\alpha} 
			\SigmatTan^{\beta}
			-
			\frac{1}{2 \upmu (\upmu - \phi \frac{\muxmulevelsetvalue}{\Lunit \upmu})}
			\Rtransarg{\muxmulevelsetvalue}_{\alpha}
			\SigmatTan^{\alpha} \SigmatTan^a \muX X_a
				\\
	&  \ \
			+
			\frac{1}{\upmu - \phi \frac{\muxmulevelsetvalue}{\Lunit \upmu}}
			\left\lbrace
				-
				\frac{1}{2}
				(\muX X_{\alpha}) \SigmatTan^{\alpha} \SigmatTan^a X_a
				+
				\frac{1}{4}
				(\muX X_{\alpha}) \uLunit^{\alpha}
				|\SigmatTan|_g^2
			\right\rbrace
				 \\
	& \ \
	+
	\frac{1}{\upmu - \phi \frac{\muxmulevelsetvalue}{\Lunit \upmu}}
	\left(
		\phi \frac{\muxmulevelsetvalue}{\upmu \Lunit \upmu}
		+
		\Rtransnormsmallfactorarg{\muxmulevelsetvalue} 
	\right)
	\left\lbrace
				-
				\frac{1}{2}
				(\muX \Lunit_{\alpha}) \SigmatTan^{\alpha} \SigmatTan^a X_a
				+
				\frac{1}{4}
				(\muX \Lunit_{\alpha}) \uLunit^{\alpha}
				|\SigmatTan|_g^2
	\right\rbrace
	-
	\frac{1}{2}
	\frac{(\muX \Rtransnormsmallfactorarg{\muxmulevelsetvalue})}{\upmu - \phi \frac{\muxmulevelsetvalue}{\Lunit \upmu}} 
	|\angV|_{\gtorus}^2
			 \\
& \ \
	+
	\frac{1}{(\upmu - \phi \frac{\muxmulevelsetvalue}{\Lunit \upmu})}
	\frac{(\muX \Lunit \timefunctionarg{\muxmulevelsetvalue})}{(\Lunit \timefunctionarg{\muxmulevelsetvalue})^2}
	\left\lbrace
				-
				\frac{1}{2}
				(\smoothtorusproject_{\alpha}^{\ \beta} \partial_{\beta} \timefunctionarg{\muxmulevelsetvalue}) \SigmatTan^{\alpha} \SigmatTan^a X_a
				+
				\frac{1}{4}
				(\smoothtorusproject_{\alpha}^{\ \beta}
			\partial_{\beta} \timefunctionarg{\muxmulevelsetvalue}) 
			\uLunit^{\alpha}
			|\SigmatTan|_g^2
		\right\rbrace
			\\
& \ \
		+
			\frac{1}{(\upmu - \phi \frac{\muxmulevelsetvalue}{\Lunit \upmu})}
			\frac{1}{\Lunit \timefunctionarg{\muxmulevelsetvalue}}
			\left\lbrace
				\frac{1}{2}
				[\muX (\smoothtorusproject_{\alpha}^{\ \beta} \partial_{\beta} \timefunctionarg{\muxmulevelsetvalue})]
				\SigmatTan^{\alpha} \SigmatTan^a X_a
				-
				\frac{1}{4}
				[\muX (\smoothtorusproject_{\alpha}^{\ \beta} \partial_{\beta} \timefunctionarg{\muxmulevelsetvalue})]
				\uLunit^{\alpha}
				|\SigmatTan|_g^2
		\right\rbrace
			\\
		& \ \
		-	
		\frac{1}{2}
		\frac{1}{\upmu}
		(\Lunit \Rtransarg{\muxmulevelsetvalue}_{\alpha})
		\SigmatTan^{\alpha} 
		\Lunit_{\beta} \SigmatTan^{\beta}
		-
		\frac{1}{4}
		\frac{1}{\upmu}
		\left\lbrace
			\Lunit
			(\Rtransarg{\muxmulevelsetvalue}_{\alpha} \uLunit^{\alpha})
		\right\rbrace
		|\SigmatTan|_g^2
		-
		\frac{1}{2}
		\frac{1}{\upmu}
		\Rtransarg{\muxmulevelsetvalue}_{\alpha}
		\SigmatTan^{\alpha} 
		(\Lunit \Lunit_{\beta}) 
		\SigmatTan^{\beta}
		+
		\frac{1}{4(\upmu - \phi \frac{\muxmulevelsetvalue}{\Lunit \upmu})}
		\Rtransarg{\muxmulevelsetvalue}_{\alpha}
		(\Lunit \uLunit^{\alpha})
		|\SigmatTan|_g^2
				\\
	& \ \
				-
			\frac{1}{2(\upmu - \phi \frac{\muxmulevelsetvalue}{\Lunit \upmu})}
			\Rtransarg{\muxmulevelsetvalue}_{\alpha}
			\SigmatTan^{\alpha} \SigmatTan^a \Lunit X_a
			+
		\frac{\upmu}{\upmu - \phi \frac{\muxmulevelsetvalue}{\Lunit \upmu}}
		\left\lbrace
				-
				\frac{1}{2}
				(\Lunit X_{\alpha}) \SigmatTan^{\alpha} \SigmatTan^a X_a
				+
				\frac{1}{4}
				(\Lunit X_{\alpha}) \uLunit^{\alpha}
				|\SigmatTan|_g^2
		\right\rbrace
			 \\
	& \ \
	+
	\frac{\upmu}{\upmu - \phi \frac{\muxmulevelsetvalue}{\Lunit \upmu}}
	\left(
		\phi \frac{\muxmulevelsetvalue}{\upmu \Lunit \upmu}
		+
		\Rtransnormsmallfactorarg{\muxmulevelsetvalue} 
	\right)
	\left\lbrace
				-
				\frac{1}{2}
				(\Lunit \Lunit_{\alpha}) \SigmatTan^{\alpha} \SigmatTan^a X_a
				+
				\frac{1}{4}
				(\Lunit \Lunit_{\alpha}) \uLunit^{\alpha}
				|\SigmatTan|_g^2
	\right\rbrace
			\\
& \ \
	+
	\frac{\upmu}{(\upmu - \phi \frac{\muxmulevelsetvalue}{\Lunit \upmu})}
	\frac{(\Lunit \Lunit \timefunctionarg{\muxmulevelsetvalue})}{(\Lunit \timefunctionarg{\muxmulevelsetvalue})^2}
	\left\lbrace
				-
				\frac{1}{2}
				(\smoothtorusproject_{\alpha}^{\ \beta} \partial_{\beta} \timefunctionarg{\muxmulevelsetvalue}) \SigmatTan^{\alpha} \SigmatTan^a X_a
				+
				\frac{1}{4}
				(\smoothtorusproject_{\alpha}^{\ \beta}
			\partial_{\beta} \timefunctionarg{\muxmulevelsetvalue}) 
			\uLunit^{\alpha}
			|\SigmatTan|_g^2
		\right\rbrace
			\\
	& \ \
			+
			\frac{\upmu}{(\upmu - \phi \frac{\muxmulevelsetvalue}{\Lunit \upmu}}
			\frac{1}{\Lunit \timefunctionarg{\muxmulevelsetvalue}}
			\left\lbrace
				\frac{1}{2}
				\left[
					\Lunit (\smoothtorusproject_{\alpha}^{\ \beta} \partial_{\beta} \timefunctionarg{\muxmulevelsetvalue}) 
				\right]
				\SigmatTan^{\alpha} \SigmatTan^a X_a
				-
				\frac{1}{4}
				\left[
					\Lunit (\smoothtorusproject_{\alpha}^{\ \beta} \partial_{\beta} \timefunctionarg{\muxmulevelsetvalue}) 
				\right]
				\uLunit^{\alpha}
				|\SigmatTan|_g^2
		\right\rbrace
	-
	\frac{1}{2}
	\frac{\upmu (\Lunit \Rtransnormsmallfactorarg{\muxmulevelsetvalue})}{\upmu - \phi \frac{\muxmulevelsetvalue}{\Lunit \upmu}} 
	|\angV|_{\gtorus}^2
			\\
	& \ \
		-
		\frac{1}{\upmu}
		(\Rtransarg{\muxmulevelsetvalue} \gfour_{\alpha \beta})
		\Lunit^{\alpha}
		\Lunit_{\gamma} 
		\SigmatTan^{\gamma}
		\SigmatTan^{\beta}
		-
			\frac{1}{2}
			\frac{1}{\upmu}
			(\Rtransarg{\muxmulevelsetvalue} \gfour_{\alpha \beta})
			\SigmatTan^{\alpha} 
			\SigmatTan^{\beta} 
			+
			\frac{1}{4 \upmu (\upmu - \phi \frac{\muxmulevelsetvalue}{\Lunit \upmu})}
			\Rtransarg{\muxmulevelsetvalue}_{\alpha}
			\uLunit^{\alpha}
			(\muX g_{ab}) \SigmatTan^a \SigmatTan^b
			\\
	& \ \
		+
		\frac{1}{\upmu}
		\Rtransarg{\muxmulevelsetvalue}^{\alpha}
		(\Lunit \gfour_{\alpha \beta})
		\Lunit_{\gamma}
		\SigmatTan^{\gamma} 
		\SigmatTan^{\beta}
		-
			\frac{1}{4}
			\frac{1}{\upmu}
			\Rtransarg{\muxmulevelsetvalue}_{\alpha} \uLunit^{\alpha}
			(\Lunit \gfour_{\beta \gamma})
			\SigmatTan^{\beta} 
			\SigmatTan^{\gamma}
			+
			\frac{1}{4(\upmu - \phi \frac{\muxmulevelsetvalue}{\Lunit \upmu})}
			\Rtransarg{\muxmulevelsetvalue}_{\alpha}
			\uLunit^{\alpha}
			(\Lunit g_{ab}) \SigmatTan^a \SigmatTan^b
				\\
& \ \
		+
		\frac{1}{\upmu}
		\Rtransarg{\muxmulevelsetvalue}_{\alpha}
		\Nullhypersurfaceproject_{\beta}^{\ \alpha}
		(\gfour^{-1})^{\beta \gamma} 
		(\partial_{\gamma} \gfour_{\delta \kappa}) 
		\SigmatTan^{\delta}
		\SigmatTan^{\kappa}
		-
		\frac{1}{\upmu}
		\Rtransarg{\muxmulevelsetvalue}_{\alpha}
		\Nullhypersurfaceproject_{\beta}^{\ \alpha}
		(\gfour^{-1})^{\beta \gamma} 
		(\partial_{\delta} \gfour_{\gamma \kappa})
		\SigmatTan^{\delta}  
		\SigmatTan^{\kappa}	
			\\
	& \ \
			+
			\frac{\phi'\frac{\muxmulevelsetvalue}{\Lunit \upmu}}{\upmu (\upmu - \phi \frac{\muxmulevelsetvalue}{\Lunit \upmu})^2}
			\left\lbrace
				-
				\frac{1}{2}
				\Rtransarg{\muxmulevelsetvalue}_{\alpha} \SigmatTan^{\alpha} \SigmatTan^a X_a
				+
				\frac{1}{4}
				\Rtransarg{\muxmulevelsetvalue}_{\alpha} \uLunit^{\alpha}
				|\SigmatTan|_g^2
		\right\rbrace
			-
			\frac{1}{2}
			\frac{\phi' \frac{\muxmulevelsetvalue}{\Lunit \upmu}}{\upmu(\upmu - \phi \frac{\muxmulevelsetvalue}{\Lunit \upmu})}
			|\angV|_{\gtorus}^2.
\end{split}
\end{align}

	Moreover, on $\twoargMrough{[\timefunction_0,\timefunctionboot),[- \rightu,\leftu]}{\muxmulevelsetvalue}$, the term
	$\CurrentboundaryerrorperfectRderivative[\SigmatTan,\SigmatTan]$ defined in 
	\eqref{E:PERFECTRDERIVAIVEERRORPUTANGENTCURRENTCONTRACTEDAGAINSTVECTORFIELD} is quantitatively positive definite 
	in the following
	sense:
	\begin{align} \label{E:KEYCOERCIVITYPERFECTRDERIVAIVEERRORPUTANGENTCURRENTCONTRACTEDAGAINSTVECTORFIELD}
		\CurrentboundaryerrorperfectRderivative[\SigmatTan,\SigmatTan]
		& \approx 
			\frac{1}
			{\upmu - \phi \frac{\muxmulevelsetvalue}{\Lunit \upmu}}
			|\SigmatTan|_g^2.
	\end{align}

\end{lemma}

\begin{proof}
	First, using \eqref{E:NULLHYPERSURFACEPROJECTION} and \eqref{E:PUTANGENTELLIPTICHYPERBOLICCURRENT},
	we compute that:
	\begin{align} \label{E:FIRSTPROOFSTEPKEYIDPUTANGENTCURRENT}
		\ehcurrent^{\alpha}[\SigmatTan,\pmb{\partial} \SigmatTan]
		& = \SigmatTan^{\beta}  
				\Nullhypersurfaceproject_{\gamma}^{\ \alpha}
				\partial_{\beta}
				\SigmatTan^{\gamma}
				+
				\frac{1}{2}
				\SigmatTan^{\beta} \Lunit_{\beta} \uLunit \SigmatTan^{\alpha}
				-
				\frac{1}{2}
				\SigmatTan^{\alpha} \Lunit_{\beta} \uLunit \SigmatTan^{\beta}
				- 
				\SigmatTan^{\beta}
				\Nullhypersurfaceproject_{\beta}^{\ \alpha} 
				\partial_a \SigmatTan^a.
	\end{align}
	Next, we note the following commutation-type identity, where we recall
	that $(\mathrm{d} \SigmatTan_{\flat})_{\alpha \beta} \eqdef \partial_{\alpha} \SigmatTan_{\beta} - \partial_{\beta} \SigmatTan_{\alpha}$:
	\begin{align} \label{E:SECONDPROOFSTEPKEYIDPUTANGENTCURRENT}
		\partial_{\beta} \SigmatTan^{\gamma} 
		& = \gfour_{\beta \kappa} (\gfour^{-1})^{\gamma \lambda} \partial_{\lambda} \SigmatTan^{\kappa}
				+
				(\gfour^{-1})^{\gamma \kappa} (\mathrm{d} \SigmatTan_{\flat})_{\beta \kappa}
				+
				(\gfour^{-1})^{\gamma \kappa} (\partial_{\kappa} \gfour_{\beta \lambda}) \SigmatTan^{\lambda}
				-
				(\gfour^{-1})^{\gamma \kappa} 
				(\partial_{\beta} \gfour_{\kappa \lambda})
				\SigmatTan^{\lambda}.
	\end{align}
	Using \eqref{E:SECONDPROOFSTEPKEYIDPUTANGENTCURRENT}, 
	we rewrite the first product on RHS~\eqref{E:FIRSTPROOFSTEPKEYIDPUTANGENTCURRENT}
	as follows:
	\begin{align} 
	\begin{split} \label{E:THIRDPROOFSTEPKEYIDPUTANGENTCURRENT}
		\SigmatTan^{\beta}  
		\Nullhypersurfaceproject_{\gamma}^{\ \alpha}
		\partial_{\beta}
		\SigmatTan^{\gamma}
		& 
		= 
		\Nullhypersurfaceproject_{\beta}^{\ \alpha}
		(\gfour^{-1})^{\beta \gamma} 
		\SigmatTan_{\delta} 
		\partial_{\gamma} 
		\SigmatTan^{\delta}
		+
		\Nullhypersurfaceproject_{\beta}^{\ \alpha}
		\SigmatTan^{\gamma}  
		(\gfour^{-1})^{\beta \delta} (\mathrm{d} V_{\flat})_{\gamma \delta}
			\\
	& \ \
		+
		\Nullhypersurfaceproject_{\beta}^{\ \alpha}
		(\gfour^{-1})^{\beta \gamma} 
		(\partial_{\gamma} \gfour_{\delta \kappa}) 
		\SigmatTan^{\delta}
		\SigmatTan^{\kappa}
		-
		\Nullhypersurfaceproject_{\beta}^{\ \alpha}
		(\gfour^{-1})^{\beta \gamma} 
		(\partial_{\delta} \gfour_{\gamma \kappa})
		\SigmatTan^{\delta}  
		\SigmatTan^{\kappa}.
	\end{split}
	\end{align}
	Next, using the identity $\uLunit = 2 \Transport - \Lunit$ (see Lemma~\ref{L:BASICPROPERTIESOFULUNIT}),
	we rewrite the second and third products on RHS~\eqref{E:FIRSTPROOFSTEPKEYIDPUTANGENTCURRENT}
	as follows:
	\begin{align} \label{E:THIRDPROOFSTEPKEYIDPUTANGENTCURRENT}
		\frac{1}{2}
		\SigmatTan^{\beta} \Lunit_{\beta} \uLunit \SigmatTan^{\alpha}
		-
		\frac{1}{2}
		\SigmatTan^{\alpha} \Lunit_{\beta} \uLunit \SigmatTan^{\beta}
		&
		=
		-
		\frac{1}{2}
		\SigmatTan^{\beta} \Lunit_{\beta} \Lunit \SigmatTan^{\alpha}
		+
		\frac{1}{2}
		\SigmatTan^{\alpha} \Lunit_{\beta} \Lunit \SigmatTan^{\beta}
		+
		\SigmatTan^{\beta} \Lunit_{\beta} \Transport \SigmatTan^{\alpha}
		-
		\SigmatTan^{\alpha} \Lunit_{\beta} \Transport \SigmatTan^{\beta}.
	\end{align}
	Differentiating by parts in the second term on RHS~\eqref{E:THIRDPROOFSTEPKEYIDPUTANGENTCURRENT},
	we deduce that:
	\begin{align}
	\begin{split}  \label{E:FOURTHPROOFSTEPKEYIDPUTANGENTCURRENT}
		\frac{1}{2}
		\SigmatTan^{\beta} \Lunit_{\beta} \uLunit \SigmatTan^{\alpha}
		-
		\frac{1}{2}
		\SigmatTan^{\alpha} \Lunit_{\beta} \uLunit \SigmatTan^{\beta}
		& = 
		-
		\SigmatTan^{\beta} \Lunit_{\beta} \Lunit \SigmatTan^{\alpha}
		+
		\frac{1}{2}
		\Lunit (\SigmatTan^{\alpha} \Lunit_{\beta} \SigmatTan^{\beta})
			\\
	& \ \
		-
		\frac{1}{2}
		\SigmatTan^{\alpha} (\Lunit \Lunit_{\beta}) \SigmatTan^{\beta}
		+
		\SigmatTan^{\beta} \Lunit_{\beta} \Transport \SigmatTan^{\alpha}
		-
		\SigmatTan^{\alpha} \Lunit_{\beta} \Transport \SigmatTan^{\beta}.
	\end{split}
	\end{align}	
	Next, with the help of \eqref{E:SECONDPROOFSTEPKEYIDPUTANGENTCURRENT},
	we rewrite the first product on RHS~\eqref{E:FOURTHPROOFSTEPKEYIDPUTANGENTCURRENT}
	as follows:
	\begin{align} 
	\begin{split} \label{E:FIFTHPROOFSTEPKEYIDPUTANGENTCURRENT}
		\SigmatTan^{\beta} \Lunit_{\beta} \Lunit \SigmatTan^{\alpha}
		& = 
		\SigmatTan^{\beta} \Lunit_{\beta} (\gfour^{-1})^{\alpha \gamma} (\partial_{\gamma} \SigmatTan^{\delta}) \Lunit_{\delta}
		+
		\SigmatTan^{\beta} 
		\Lunit_{\beta} 
		\Lunit^{\gamma} 
		(\gfour^{-1})^{\alpha \delta}
		(\mathrm{d} V_{\flat})_{\gamma \delta}
			\\
	& \ \
		+
		(\gfour^{-1})^{\alpha \beta}
		(\partial_{\beta} \gfour_{\gamma \delta})
		\Lunit^{\gamma}
		\Lunit_{\kappa} 
		\SigmatTan^{\kappa}
		\SigmatTan^{\delta}
		-
		(\gfour^{-1})^{\alpha \beta}
		(\Lunit \gfour_{\beta \gamma})
		\Lunit_{\delta}
		\SigmatTan^{\delta} 
		\SigmatTan^{\gamma}.
	\end{split}
	\end{align}	
	Combining \eqref{E:FIRSTPROOFSTEPKEYIDPUTANGENTCURRENT}--\eqref{E:FIFTHPROOFSTEPKEYIDPUTANGENTCURRENT},
	we find that:
		\begin{align}
		\begin{split}  \label{E:SIXTHPROOFSTEPKEYIDPUTANGENTCURRENT}
		\ehcurrent^{\alpha}[\SigmatTan,\pmb{\partial} \SigmatTan]
		& = 
		\frac{1}{2}
		\Lunit (\SigmatTan^{\alpha} \Lunit_{\beta} \SigmatTan^{\beta})
		+
		\Nullhypersurfaceproject_{\beta}^{\ \alpha}
		(\gfour^{-1})^{\beta \gamma} 
		\SigmatTan_{\delta} 
		\partial_{\gamma} 
		\SigmatTan^{\delta}
		-
		\SigmatTan^{\beta} \Lunit_{\beta} 
		(\gfour^{-1})^{\alpha \gamma} 
		(\partial_{\gamma} \SigmatTan^{\delta}) 
		\Lunit_{\delta}
			\\
		& \ \
			+
		\Nullhypersurfaceproject_{\beta}^{\ \alpha}
		\SigmatTan^{\gamma}  
		(\gfour^{-1})^{\beta \delta} (\mathrm{d} V_{\flat})_{\gamma \delta}
		-
		\SigmatTan^{\beta} 
		\Lunit_{\beta} 
		\Lunit^{\gamma} 
		(\gfour^{-1})^{\alpha \delta}
		(\mathrm{d} V_{\flat})_{\gamma \delta}
			 \\
		& \ \
		+
		\Nullhypersurfaceproject_{\beta}^{\ \alpha}
		(\gfour^{-1})^{\beta \gamma} 
		(\partial_{\gamma} \gfour_{\delta \kappa}) 
		\SigmatTan^{\delta}
		\SigmatTan^{\kappa}
		-
		\Nullhypersurfaceproject_{\beta}^{\ \alpha}
		(\gfour^{-1})^{\beta \gamma} 
		(\partial_{\delta} \gfour_{\gamma \kappa})
		\SigmatTan^{\delta}  
		\SigmatTan^{\kappa}
				\\
	& \ \
		-
		(\gfour^{-1})^{\alpha \beta}
		(\partial_{\beta} \gfour_{\gamma \delta})
		\Lunit^{\gamma}
		\Lunit_{\kappa} 
		\SigmatTan^{\kappa}
		\SigmatTan^{\delta}
		+
		(\gfour^{-1})^{\alpha \beta}
		(\Lunit \gfour_{\beta \gamma})
		\Lunit_{\delta}
		\SigmatTan^{\delta} 
		\SigmatTan^{\gamma}
		-
		\frac{1}{2}
		\SigmatTan^{\alpha} (\Lunit \Lunit_{\beta}) \SigmatTan^{\beta}
				\\
		& \ \
		+
		\SigmatTan^{\beta} \Lunit_{\beta} \Transport \SigmatTan^{\alpha}
		-
		\SigmatTan^{\alpha} \Lunit_{\beta} \Transport \SigmatTan^{\beta}
		- 
		\Nullhypersurfaceproject_{\beta}^{\ \alpha} \SigmatTan^{\beta}
		\partial_a \SigmatTan^a.
	\end{split}
	\end{align}
	Next, taking into account definition \eqref{E:NULLHYPERSURFACEPROJECTION}
	and differentiating by parts,
	we express the second and third products on RHS~\eqref{E:SIXTHPROOFSTEPKEYIDPUTANGENTCURRENT}
	as follows:
	\begin{align} 
	\begin{split} \label{E:SEVENTHPROOFSTEPKEYIDPUTANGENTCURRENT}
		&
		\Nullhypersurfaceproject_{\beta}^{\ \alpha}
		(\gfour^{-1})^{\beta \gamma} 
		\SigmatTan_{\delta} 
		\partial_{\gamma} 
		\SigmatTan^{\delta}
		-
		\SigmatTan^{\beta} \Lunit_{\beta} 
		(\gfour^{-1})^{\alpha \gamma} 
		(\partial_{\gamma} \SigmatTan^{\delta}) \Lunit_{\delta}
			\\
		& =
			\frac{1}{2}
			(\gfour^{-1})^{\alpha \beta}
			\partial_{\beta}
			(|\SigmatTan|_g^2)
			+
			\frac{1}{4}
			\uLunit^{\alpha}
			\Lunit
			(|\SigmatTan|_g^2)
			-
			\frac{1}{2}
			(\gfour^{-1})^{\alpha \beta}
			\partial_{\beta}
			\left\lbrace
				(\SigmatTan^{\gamma} \Lunit_{\gamma})^2
			\right\rbrace
				 \\
			& \ \
				-
				\frac{1}{2}
				(\gfour^{-1})^{\alpha \beta}
				(\partial_{\beta} \gfour_{\gamma \delta})
				\SigmatTan^{\gamma} 
				\SigmatTan^{\delta} 
				-
			\frac{1}{4}
			\uLunit^{\alpha}
			(\Lunit \gfour_{\gamma \delta})
			\SigmatTan^{\gamma} 
			\SigmatTan^{\delta}
			+
			(\gfour^{-1})^{\alpha \beta}
			(\partial_{\beta} \Lunit_{\gamma})
			\SigmatTan^{\gamma} 
			\SigmatTan^{\delta} \Lunit_{\delta}.
	\end{split}
	\end{align}
	Using \eqref{E:SEVENTHPROOFSTEPKEYIDPUTANGENTCURRENT} to substitute for
	the second and third products on RHS~\eqref{E:SIXTHPROOFSTEPKEYIDPUTANGENTCURRENT},
	we deduce that:
	\begin{align} 
	\begin{split} \label{E:EIGHTHPROOFSTEPKEYIDPUTANGENTCURRENT}
		\ehcurrent^{\alpha}[\SigmatTan,\pmb{\partial} \SigmatTan]
		& = 
			\frac{1}{2}
			(\gfour^{-1})^{\alpha \beta}
			\partial_{\beta}
			(|\SigmatTan|_g^2)
			-
			\frac{1}{2}
			(\gfour^{-1})^{\alpha \beta}
			\partial_{\beta}
			\left\lbrace
				(\SigmatTan^{\gamma} \Lunit_{\gamma})^2
			\right\rbrace
			+
			\frac{1}{2}
			\Lunit (\SigmatTan^{\alpha} \Lunit_{\beta} \SigmatTan^{\beta})
			+
			\frac{1}{4}
			\uLunit^{\alpha}
			\Lunit
			(|\SigmatTan|_g^2)
			\\
		& \ \
		+
		\Nullhypersurfaceproject_{\beta}^{\ \alpha}
		\SigmatTan^{\gamma}  
		(\gfour^{-1})^{\beta \delta} (\mathrm{d} V_{\flat})_{\gamma \delta}
		-
		\SigmatTan^{\beta} 
		\Lunit_{\beta} 
		\Lunit^{\gamma} 
		(\gfour^{-1})^{\alpha \delta}
		(\mathrm{d} V_{\flat})_{\gamma \delta}
			 \\
		& \ \
		+
		\Nullhypersurfaceproject_{\beta}^{\ \alpha}
		(\gfour^{-1})^{\beta \gamma} 
		(\partial_{\gamma} \gfour_{\delta \kappa}) 
		\SigmatTan^{\delta}
		\SigmatTan^{\kappa}
		-
		\Nullhypersurfaceproject_{\beta}^{\ \alpha}
		(\gfour^{-1})^{\beta \gamma} 
		(\partial_{\delta} \gfour_{\gamma \kappa})
		\SigmatTan^{\delta}  
		\SigmatTan^{\kappa}
				\\
	& \ \
		-
		(\gfour^{-1})^{\alpha \beta}
		(\partial_{\beta} \gfour_{\gamma \delta})
		\Lunit^{\gamma}
		\Lunit_{\kappa} 
		\SigmatTan^{\kappa}
		\SigmatTan^{\delta}
		+
		(\gfour^{-1})^{\alpha \beta}
		(\Lunit \gfour_{\beta \gamma})
		\Lunit_{\delta}
		\SigmatTan^{\delta} 
		\SigmatTan^{\gamma}
		-
		\frac{1}{2}
		\SigmatTan^{\alpha} (\Lunit \Lunit_{\beta}) \SigmatTan^{\beta}
				\\
		& \ \
		+
		\SigmatTan^{\beta} \Lunit_{\beta} \Transport \SigmatTan^{\alpha}
		-
		\SigmatTan^{\alpha} \Lunit_{\beta} \Transport \SigmatTan^{\beta}
		- 
		\Nullhypersurfaceproject_{\beta}^{\ \alpha} \SigmatTan^{\beta}
		\partial_a \SigmatTan^a
					\\
		& \ \
			-
				\frac{1}{2}
				(\gfour^{-1})^{\alpha \beta}
				(\partial_{\beta} \gfour_{\gamma \delta})
				\SigmatTan^{\gamma} 
				\SigmatTan^{\delta} 
				-
			\frac{1}{4}
			\uLunit^{\alpha}
			(\Lunit \gfour_{\gamma \delta})
			\SigmatTan^{\gamma} 
			\SigmatTan^{\delta}
			+
			(\gfour^{-1})^{\alpha \beta}
			(\partial_{\beta} \Lunit_{\gamma})
			\SigmatTan^{\gamma} 
			\SigmatTan^{\delta} \Lunit_{\delta}.
	\end{split}
	\end{align}
	Next, we note that since $\ehcurrent[\SigmatTan,\pmb{\partial} \SigmatTan]$ is tangent to $\nullhyparg{u}$,
	we have $\Lunit_{\alpha} \ehcurrent^{\alpha}[\SigmatTan,\pmb{\partial} \SigmatTan] = 0$.
	From this fact and equations \eqref{E:UNITLENGTHRTRANS} and \eqref{E:HYPUNITNORMALDECOMPOSITION}, 
	it follows that
	$
	\frac{|\Rtransarg{\muxmulevelsetvalue}|_{\hypg}}
		{\upmu}
		\hypunitnormalarg{\muxmulevelsetvalue}_{\alpha}
		\ehcurrent^{\alpha}[\SigmatTan,\pmb{\partial} \SigmatTan] 
	=
	\frac{1}{\upmu} \Rtransarg{\muxmulevelsetvalue}_{\alpha} \ehcurrent^{\alpha}[\SigmatTan,\pmb{\partial} \SigmatTan]
	$.
	Using this identity,
	\eqref{E:EIGHTHPROOFSTEPKEYIDPUTANGENTCURRENT},
	the identity $\SigmatTan^{\gamma} \Lunit_{\gamma} = - \SigmatTan^a X_a$ implied by
	the fact that $\SigmatTan^{\gamma} \Transport_{\gamma} = 0$ and the identity $\Transport = \Lunit + X$ (see \eqref{E:BISLPLUSX}),
	and using that
	$
	|\SigmatTan|_g^2 
		-
	(\SigmatTan^a X_a)^2
	=
	|\SigmatTan|_{\gtorus}^2
	$
	(see \eqref{E:SMOOTHTORUSMETRICINTERMSOFSIGMATMETRICANDX})
	we deduce that:
	\begin{align} 
	\begin{split} \label{E:NINTHPROOFSTEPKEYIDPUTANGENTCURRENT}
		\frac{|\Rtransarg{\muxmulevelsetvalue}|_{\hypg}}
		{\upmu}
		\hypunitnormalarg{\muxmulevelsetvalue}_{\alpha}
		\ehcurrent^{\alpha}[\SigmatTan,\pmb{\partial} \SigmatTan] 
		& = 
		\frac{1}{2}
		\frac{1}{\upmu}
		\Rtransarg{\muxmulevelsetvalue}
		\left\lbrace
			|\angV|_{\gtorus}^2 
		\right\rbrace
		+
		\frac{1}{2}
		\frac{1}{\upmu}
		\Rtransarg{\muxmulevelsetvalue}_{\alpha} 
		\Lunit 
		(\SigmatTan^{\alpha} \Lunit_{\beta} \SigmatTan^{\beta})
		+
		\frac{1}{4}
		\frac{1}{\upmu}
		\Rtransarg{\muxmulevelsetvalue}_{\alpha}
		\uLunit^{\alpha}
		\Lunit
		(|\SigmatTan|_g^2)
				\\
		& \ \
			+
		\frac{1}{\upmu}
		\Rtransarg{\muxmulevelsetvalue}_{\alpha}
		\Nullhypersurfaceproject_{\beta}^{\ \alpha}
		\SigmatTan^{\gamma}  
		(\gfour^{-1})^{\beta \delta} (\mathrm{d} V_{\flat})_{\gamma \delta}
		-
		\frac{1}{\upmu}
		\SigmatTan^{\beta} 
		\Lunit_{\beta} 
		\Lunit^{\gamma} 
		\Rtransarg{\muxmulevelsetvalue}^{\delta}
		(\mathrm{d} V_{\flat})_{\gamma \delta}
			\\
		& \ \
		+
		\frac{1}{\upmu}
		\Rtransarg{\muxmulevelsetvalue}_{\alpha}
		\Nullhypersurfaceproject_{\beta}^{\ \alpha}
		(\gfour^{-1})^{\beta \gamma} 
		(\partial_{\gamma} \gfour_{\delta \kappa}) 
		\SigmatTan^{\delta}
		\SigmatTan^{\kappa}
		-
		\frac{1}{\upmu}
		\Rtransarg{\muxmulevelsetvalue}_{\alpha}
		\Nullhypersurfaceproject_{\beta}^{\ \alpha}
		(\gfour^{-1})^{\beta \gamma} 
		(\partial_{\delta} \gfour_{\gamma \kappa})
		\SigmatTan^{\delta}  
		\SigmatTan^{\kappa}
				\\
	& \ \
		-
		\frac{1}{\upmu}
		(\Rtransarg{\muxmulevelsetvalue} \gfour_{\alpha \beta})
		\Lunit^{\alpha}
		\Lunit_{\gamma} 
		\SigmatTan^{\gamma}
		\SigmatTan^{\beta}
		+
		\frac{1}{\upmu}
		\Rtransarg{\muxmulevelsetvalue}^{\alpha}
		(\Lunit \gfour_{\alpha \beta})
		\Lunit_{\gamma}
		\SigmatTan^{\gamma} 
		\SigmatTan^{\beta}
		-
		\frac{1}{2}
		\frac{1}{\upmu}
		\Rtransarg{\muxmulevelsetvalue}_{\alpha}
		\SigmatTan^{\alpha} 
		(\Lunit \Lunit_{\beta}) 
		\SigmatTan^{\beta}
				\\
		& \ \
			+
			\frac{1}{\upmu}
			\Rtransarg{\muxmulevelsetvalue}_{\alpha}
			\SigmatTan^{\beta} \Lunit_{\beta} \Transport \SigmatTan^{\alpha}
			-
			\frac{1}{\upmu}
			\Rtransarg{\muxmulevelsetvalue}_{\alpha}
			\SigmatTan^{\alpha} \Lunit_{\beta} \Transport \SigmatTan^{\beta}
			- 
			\frac{1}{\upmu}
			\Rtransarg{\muxmulevelsetvalue}_{\alpha}
			\Nullhypersurfaceproject_{\beta}^{\ \alpha} \SigmatTan^{\beta}
			\partial_a \SigmatTan^a
					\\
		& \ \
			-
			\frac{1}{2}
			\frac{1}{\upmu}
			(\Rtransarg{\muxmulevelsetvalue} \gfour_{\alpha \beta})
			\SigmatTan^{\alpha} 
			\SigmatTan^{\beta} 
			-
			\frac{1}{4}
			\frac{1}{\upmu}
			\Rtransarg{\muxmulevelsetvalue}_{\alpha} \uLunit^{\alpha}
			(\Lunit \gfour_{\beta \gamma})
			\SigmatTan^{\beta} 
			\SigmatTan^{\gamma}
			+
			\frac{1}{\upmu}
			(\Rtransarg{\muxmulevelsetvalue} \Lunit_{\alpha})
			\Lunit_{\beta}
			\SigmatTan^{\alpha} 
			\SigmatTan^{\beta}.
	\end{split}
	\end{align}
	Next, using 
	the identity $\Lunit_{\beta} \SigmatTan^{\beta} = - \SigmatTan^a X_a$ noted above
	and differentiating by parts,
	we rewrite the terms on the first line of RHS~\eqref{E:NINTHPROOFSTEPKEYIDPUTANGENTCURRENT}
	as follows:
	\begin{align} 
	\begin{split} \label{E:TENTHTHPROOFSTEPKEYIDPUTANGENTCURRENT}
		&
		\frac{1}{2}
		\frac{1}{\upmu}
		\Rtransarg{\muxmulevelsetvalue}
		\left\lbrace
			|\angV|_{\gtorus}^2 
		\right\rbrace
		+
		\frac{1}{2}
		\frac{1}{\upmu}
		\Rtransarg{\muxmulevelsetvalue}_{\alpha} 
		\Lunit 
		(\SigmatTan^{\alpha} \Lunit_{\beta} \SigmatTan^{\beta})
		+
		\frac{1}{4}
		\frac{1}{\upmu}
		\Rtransarg{\muxmulevelsetvalue}_{\alpha}
		\uLunit^{\alpha}
		\Lunit
		(|\SigmatTan|_g^2)
			\\
		& 
		=
		\frac{1}{2}
		\Rtransarg{\muxmulevelsetvalue}
		\left\lbrace
			\frac{1}{\upmu}
			|\angV|_{\gtorus}^2 
		\right\rbrace
		+
			\Lunit
			\left\lbrace
				-
				\frac{1}{2}
				\frac{1}{\upmu}
				\Rtransarg{\muxmulevelsetvalue}_{\alpha} \SigmatTan^{\alpha} \SigmatTan^a X_a
				+
				\frac{1}{4}
				\frac{1}{\upmu}
				\Rtransarg{\muxmulevelsetvalue}_{\alpha} \uLunit^{\alpha}
				|\SigmatTan|_g^2
		\right\rbrace
				\\
		& \ \
		+
		\frac{1}{2}
		\frac{(\Rtransarg{\muxmulevelsetvalue} \upmu)}{\upmu^2}
		|\angV|_{\gtorus}^2
		+
		\frac{(\Lunit \upmu)}{\upmu^2}
		\left\lbrace
			-
				\frac{1}{2}
				\Rtransarg{\muxmulevelsetvalue}_{\alpha} \SigmatTan^{\alpha} \SigmatTan^a X_a
				+
				\frac{1}{4}
				\Rtransarg{\muxmulevelsetvalue}_{\alpha} \uLunit^{\alpha}
				|\SigmatTan|_g^2
		\right\rbrace
				\\
		& \ \
		-	
		\frac{1}{2}
		\frac{1}{\upmu}
		(\Lunit \Rtransarg{\muxmulevelsetvalue}_{\alpha})
		\SigmatTan^{\alpha} 
		\Lunit_{\beta} \SigmatTan^{\beta}
		-
		\frac{1}{4}
		\frac{1}{\upmu}
		\left\lbrace
			\Lunit
			(\Rtransarg{\muxmulevelsetvalue}_{\alpha} \uLunit^{\alpha})
		\right\rbrace
		|\SigmatTan|_g^2.
	\end{split}
	\end{align}
	Next, we use the identities
	\eqref{E:LUNITINTERMSOFSIGMATILDETANGENTVECTORFIELDSANDTRANSPORT},
	\eqref{E:RTRANSCONTRACTIONIDENTITYINPROOFOFELLIPTICBOUNDARYTERMIDENTITY},
	\eqref{E:RTRANSDERIVATIVEOFSINGULARWEIGHT},
	and \eqref{E:ROUGHTORIDERIVATIVEOFNONSINGULARWEIGHT},
	and differentiation by parts to rewrite the 
	terms on the first line of RHS~\eqref{E:TENTHTHPROOFSTEPKEYIDPUTANGENTCURRENT}
	as follows, where $\roughangdiv$ is the divergence operator on
	$\twoargroughtori{\timefunctionarg{\muxmulevelsetvalue},u}{\muxmulevelsetvalue}$-tangent vectorfields
	from Def.\,\ref{D:CONNECTIONSANDDIFFERENTIALOPERATORSONROUGHTORI}:
	\begin{align} 
	\begin{split} \label{E:ELEVENTHPROOFSTEPKEYIDPUTANGENTCURRENT}
		&
		\frac{1}{2}
		\Rtransarg{\muxmulevelsetvalue}
		\left\lbrace
			\frac{1}{\upmu}
			|\angV|_{\gtorus}^2 
		\right\rbrace
		+
			\Lunit
			\left\lbrace
				-
				\frac{1}{2}
				\frac{1}{\upmu}
				\Rtransarg{\muxmulevelsetvalue}_{\alpha} \SigmatTan^{\alpha} \SigmatTan^a X_a
				+
				\frac{1}{4}
				\frac{1}{\upmu}
				\Rtransarg{\muxmulevelsetvalue}_{\alpha} \uLunit^{\alpha}
				|\SigmatTan|_g^2
		\right\rbrace
			\\
		& = 
		\frac{1}{4}
			\Rtransarg{\muxmulevelsetvalue}
			\left\lbrace
				\frac{1}
				{\upmu - \phi \frac{\muxmulevelsetvalue}{\Lunit \upmu}}
				\left[
					(X_a \SigmatTan^a)^2
					+
					(1 + 2 \Rtransnormsmallfactorarg{\muxmulevelsetvalue}) |\SigmatTan|_{\gtorus}^2
				-
				2
				\frac{1}{\Lunit \timefunctionarg{\muxmulevelsetvalue}}
				X_a \SigmatTan^a
				\angVarg{\alpha}
				\angDarg{\alpha}
				\timefunctionarg{\muxmulevelsetvalue}
			\right]
			\right\rbrace
						\\
		& \ \
			+
			\left\lbrace
				- 
				\frac{\Rtransarg{\muxmulevelsetvalue} \upmu}{(\upmu - \phi \frac{\muxmulevelsetvalue}{\Lunit \upmu})^2}
				+ 
				\frac{\muxmulevelsetvalue \frac{\phi'}{\Lunit \upmu}}{(\upmu - \phi \frac{\muxmulevelsetvalue}{\Lunit \upmu})^2}
				-
				\frac{\phi \muxmulevelsetvalue \frac{\Rtransarg{\muxmulevelsetvalue} \Lunit \upmu}{(\Lunit \upmu)^2}}{(\upmu - \phi \frac{\muxmulevelsetvalue}{\Lunit \upmu})^2}
			\right\rbrace
			\left\lbrace
				-
				\frac{1}{2}
				\frac{1}{\upmu}
				\Rtransarg{\muxmulevelsetvalue}_{\alpha} \SigmatTan^{\alpha} \SigmatTan^a X_a
				+
				\frac{1}{4}
				\frac{1}{\upmu}
				\Rtransarg{\muxmulevelsetvalue}_{\alpha} \uLunit^{\alpha}
				|\SigmatTan|_g^2
		\right\rbrace
				\\	
		& \ \
			- 
			\roughangdiv
			\left\lbrace
			\frac{\upmu}{\upmu - \phi \frac{\muxmulevelsetvalue}{\Lunit \upmu}}
			\left[
				-
				\frac{1}{2}
				\frac{1}{\upmu}
				\Rtransarg{\muxmulevelsetvalue}_{\alpha} \SigmatTan^{\alpha} \SigmatTan^a X_a
				+
				\frac{1}{4}
				\frac{1}{\upmu}
				\Rtransarg{\muxmulevelsetvalue}_{\alpha} \uLunit^{\alpha}
				|\SigmatTan|_g^2
		\right]
		\Roughtoritangentvectorfieldarg{\muxmulevelsetvalue}
		\right\rbrace
					\\
		& \ \
			+
			\left\lbrace
			\frac{1}{\upmu - \phi \frac{\muxmulevelsetvalue}{\Lunit \upmu}}
			\left[
				-
				\frac{1}{2}
				\Rtransarg{\muxmulevelsetvalue}_{\alpha} \SigmatTan^{\alpha} \SigmatTan^a X_a
				+
				\frac{1}{4}
				\Rtransarg{\muxmulevelsetvalue}_{\alpha} \uLunit^{\alpha}
				|\SigmatTan|_g^2
		\right]
		\right\rbrace
		\roughangdiv \Roughtoritangentvectorfieldarg{\muxmulevelsetvalue}
					\\
		& \ \
			-
			\left\lbrace
				\frac{(\Roughtoritangentvectorfieldarg{\muxmulevelsetvalue} \upmu) \phi \frac{\muxmulevelsetvalue}{\Lunit \upmu}}
				{(\upmu - \phi \frac{\muxmulevelsetvalue}{\Lunit \upmu})^2}
				+
				\frac{\upmu \phi \muxmulevelsetvalue \frac{\Roughtoritangentvectorfieldarg{\muxmulevelsetvalue} \Lunit \upmu}
				{(\Lunit \upmu)^2}}{(\upmu - \phi \frac{\muxmulevelsetvalue}{\Lunit \upmu})^2}
			\right\rbrace
			\left\lbrace
				-
				\frac{1}{2}
				\frac{1}{\upmu}
				\Rtransarg{\muxmulevelsetvalue}_{\alpha} \SigmatTan^{\alpha} \SigmatTan^a X_a
				+
				\frac{1}{4}
				\frac{1}{\upmu}
				\Rtransarg{\muxmulevelsetvalue}_{\alpha} \uLunit^{\alpha}
				|\SigmatTan|_g^2
		\right\rbrace
				\\
		& \ \
			+
			\frac{\upmu}{\upmu - \phi \frac{\muxmulevelsetvalue}{\Lunit \upmu}}
			\Transport
			\left\lbrace
				-
				\frac{1}{2}
				\frac{1}{\upmu}
				\Rtransarg{\muxmulevelsetvalue}_{\alpha} \SigmatTan^{\alpha} \SigmatTan^a X_a
				+
				\frac{1}{4}
				\frac{1}{\upmu}
				\Rtransarg{\muxmulevelsetvalue}_{\alpha} \uLunit^{\alpha}
				|\SigmatTan|_g^2
		\right\rbrace.
	\end{split}
	\end{align}
	Before proceeding, we note the following identity,
	which follows from substituting RHS~\eqref{E:RTRANSDIVIDEDBYMUIDENTITY} (after lowering indices on both sides via $\gfour$)
	for the factor of $\frac{\Rtransarg{\muxmulevelsetvalue}_{\alpha}}{\upmu}$ on LHS~\eqref{E:TWELFTHPROOFSTEPKEYIDPUTANGENTCURRENT},
	from using the Leibniz and chain rules,
	from using \eqref{E:LUNITINTERMSOFSIGMATILDETANGENTVECTORFIELDSANDTRANSPORT}
	to re-express the vectorfield differential operator
	$
		\frac{\upmu}{\upmu - \phi \frac{\muxmulevelsetvalue}{\Lunit \upmu}}
			\Transport
	$
	as differentiation with respect to
	$
	\Lunit
	+
	\frac{1}{\upmu - \phi \frac{\muxmulevelsetvalue}{\Lunit \upmu}}
	\Rtransarg{\muxmulevelsetvalue}
	+
	\frac{\upmu}{\upmu - \phi \frac{\muxmulevelsetvalue}{\Lunit \upmu}}
				\Roughtoritangentvectorfieldarg{\muxmulevelsetvalue}
	$
	when the operator falls on the factor of $\frac{1}{\upmu}$ on RHS~\eqref{E:RTRANSDIVIDEDBYMUIDENTITY},
	from using \eqref{E:BISLPLUSX}
	to re-express the vectorfield differential operator
	$
		\frac{\upmu}{\upmu - \phi \frac{\muxmulevelsetvalue}{\Lunit \upmu}}
			\Transport
	$
	as differentiation with respect to
	$
	\frac{1}{(\upmu - \phi \frac{\muxmulevelsetvalue}{\Lunit \upmu})}
	\muX
	+
	\frac{\upmu}{\upmu - \phi \frac{\muxmulevelsetvalue}{\Lunit \upmu}}
	\Lunit
	$
	when the operator falls on any other factor besides $\frac{1}{\upmu}$ on RHS~\eqref{E:RTRANSDIVIDEDBYMUIDENTITY},
	from using that 
	$\angDarg{\alpha} \timefunctionarg{\muxmulevelsetvalue} = \smoothtorusproject^{\ \beta}_{\alpha} \partial_{\beta} \timefunctionarg{\muxmulevelsetvalue}$,
	from using that $\muX \phi = \phi'$,
	and from using that $\Singletan \phi = 0$ for any $\nullhyparg{u}$-tangent vectorfield $\Singletan$:
	\begin{align} 
	\begin{split} \label{E:TWELFTHPROOFSTEPKEYIDPUTANGENTCURRENT}
		\frac{\upmu}{(\upmu - \phi \frac{\muxmulevelsetvalue}{\Lunit \upmu})}
		\Transport
		\left\lbrace
			\frac{\Rtransarg{\muxmulevelsetvalue}_{\alpha}}{\upmu}
		\right\rbrace
		& =
			-
			\frac{(\Rtransarg{\muxmulevelsetvalue} \upmu) \phi \frac{\muxmulevelsetvalue}{\Lunit \upmu}}{\upmu^2(\upmu - \phi \frac{\muxmulevelsetvalue}{\Lunit \upmu})}
			\Lunit_{\alpha}
			-
			\frac{(\Lunit \upmu) \phi \frac{\muxmulevelsetvalue}{\Lunit \upmu}}{\upmu^2}
			\Lunit_{\alpha}
			-
			\frac{(\Roughtoritangentvectorfieldarg{\muxmulevelsetvalue} \upmu) \phi \frac{\muxmulevelsetvalue}{\Lunit \upmu}}{\upmu(\upmu - \phi \frac{\muxmulevelsetvalue}{\Lunit \upmu})}
			\Lunit_{\alpha}
			+
			\frac{\phi' \frac{\muxmulevelsetvalue}{\Lunit \upmu}}{\upmu(\upmu - \phi \frac{\muxmulevelsetvalue}{\Lunit \upmu})}
			\Lunit_{\alpha}
				\\
	& \ \
		-
		\phi \frac{\muxmulevelsetvalue (\muX \Lunit \upmu)}{\upmu (\upmu - \phi \frac{\muxmulevelsetvalue}{\Lunit \upmu}) (\Lunit \upmu)^2}
		\Lunit_{\alpha}
		-
		\phi \frac{\muxmulevelsetvalue (\Lunit \Lunit \upmu)}{(\upmu - \phi \frac{\muxmulevelsetvalue}{\Lunit \upmu}) (\Lunit \upmu)^2}
		\Lunit_{\alpha}
			\\
	& \ \	
		+
		\frac{1}{\upmu - \phi \frac{\muxmulevelsetvalue}{\Lunit \upmu}}
		\muX X_{\alpha}
		+
		\frac{\upmu}{(\upmu - \phi \frac{\muxmulevelsetvalue}{\Lunit \upmu})}
		\Lunit X_{\alpha}
	+
	\frac{1}{\upmu - \phi \frac{\muxmulevelsetvalue}{\Lunit \upmu}}
	\left(
		\phi \frac{\muxmulevelsetvalue}{\upmu \Lunit \upmu}
		+
		\Rtransnormsmallfactorarg{\muxmulevelsetvalue} 
	\right)
	\muX \Lunit_{\alpha}
			\\
& \ \
	+
	\frac{\upmu}{\upmu - \phi \frac{\muxmulevelsetvalue}{\Lunit \upmu}}
	\left(
		\phi \frac{\muxmulevelsetvalue}{\upmu \Lunit \upmu}
		+
		\Rtransnormsmallfactorarg{\muxmulevelsetvalue} 
	\right)
	\Lunit \Lunit_{\alpha}
	+
	\frac{(\muX \Rtransnormsmallfactorarg{\muxmulevelsetvalue})}{\upmu - \phi \frac{\muxmulevelsetvalue}{\Lunit \upmu}}
	\Lunit_{\alpha}
	+
	\frac{\upmu (\Lunit \Rtransnormsmallfactorarg{\muxmulevelsetvalue})}{\upmu - \phi \frac{\muxmulevelsetvalue}{\Lunit \upmu}}
	\Lunit_{\alpha}
		 \\
& \ \
	+
	\frac{1}{(\upmu - \phi \frac{\muxmulevelsetvalue}{\Lunit \upmu})}
	\frac{(\muX \Lunit \timefunctionarg{\muxmulevelsetvalue})}{(\Lunit \timefunctionarg{\muxmulevelsetvalue})^2}
	(\smoothtorusproject_{\alpha}^{\ \beta}
	\partial_{\beta} \timefunctionarg{\muxmulevelsetvalue})
	+
	\frac{\upmu}{(\upmu - \phi \frac{\muxmulevelsetvalue}{\Lunit \upmu})}
	\frac{(\Lunit \Lunit \timefunctionarg{\muxmulevelsetvalue})}{(\Lunit \timefunctionarg{\muxmulevelsetvalue})^2}
	(\smoothtorusproject_{\alpha}^{\ \beta}
	\partial_{\beta} \timefunctionarg{\muxmulevelsetvalue})
			\\
	& \ \
			-
			\frac{1}{(\upmu - \phi \frac{\muxmulevelsetvalue}{\Lunit \upmu})}
			\frac{1}{\Lunit \timefunctionarg{\muxmulevelsetvalue}}
			\muX (\smoothtorusproject_{\alpha}^{\ \beta} \partial_{\beta} \timefunctionarg{\muxmulevelsetvalue})
			-
			\frac{\upmu}{(\upmu - \phi \frac{\muxmulevelsetvalue}{\Lunit \upmu})}
			\frac{1}{\Lunit \timefunctionarg{\muxmulevelsetvalue}}
			\Lunit (\smoothtorusproject_{\alpha}^{\ \beta} \partial_{\beta} \timefunctionarg{\muxmulevelsetvalue}).
	\end{split}
	\end{align}
	Next, we rewrite the terms on the last line of RHS~\eqref{E:ELEVENTHPROOFSTEPKEYIDPUTANGENTCURRENT}
	using the following strategy. 
	First, we use the Leibniz and chain rules to expand the differentiation with respect to
	$
		\frac{\upmu}{\upmu - \phi \frac{\muxmulevelsetvalue}{\Lunit \upmu}}
			\Transport
	$
	in Cartesian coordinates.
	Second, when this derivative operator falls on a Cartesian component of $\SigmatTan$,
	we make no further adjustments. Third,
	when the derivative operator
	$
		\frac{\upmu}{\upmu - \phi \frac{\muxmulevelsetvalue}{\Lunit \upmu}}
			\Transport
	$
	falls on either of the two factors of $\frac{1}{\upmu} \Rtransarg{\muxmulevelsetvalue}_{\alpha}$
	in braces on LHS~\eqref{E:THIRTEENTHPROOFSTEPKEYIDPUTANGENTCURRENT},
	we use \eqref{E:TWELFTHPROOFSTEPKEYIDPUTANGENTCURRENT} for substitution.
	Finally, 
	when the derivative operator
	$
		\frac{\upmu}{\upmu - \phi \frac{\muxmulevelsetvalue}{\Lunit \upmu}}
			\Transport
	$
	falls on any other quantity in braces on LHS~\eqref{E:THIRTEENTHPROOFSTEPKEYIDPUTANGENTCURRENT},
	we use \eqref{E:BISLPLUSX}
	to re-express the operator as differentiation with respect to
	$
	\frac{1}{\upmu - \phi \frac{\muxmulevelsetvalue}{\Lunit \upmu}}
	\muX
	+
	\frac{\upmu}{\upmu - \phi \frac{\muxmulevelsetvalue}{\Lunit \upmu}}
	\Lunit
	$. 
	These four steps allow us to deduce the following identity:
	\begin{align}  
	\begin{split} \label{E:THIRTEENTHPROOFSTEPKEYIDPUTANGENTCURRENT}
		&
		\frac{\upmu}{\upmu - \phi \frac{\muxmulevelsetvalue}{\Lunit \upmu}}
			\Transport
			\left\lbrace
				-
				\frac{1}{2}
				\frac{1}{\upmu}
				\Rtransarg{\muxmulevelsetvalue}_{\alpha} \SigmatTan^{\alpha} \SigmatTan^a X_a
				+
				\frac{1}{4}
				\frac{1}{\upmu}
				\Rtransarg{\muxmulevelsetvalue}_{\alpha} \uLunit^{\alpha}
				|\SigmatTan|_g^2
		\right\rbrace
			\\
	& =
				-
				\frac{1}{2(\upmu - \phi \frac{\muxmulevelsetvalue}{\Lunit \upmu})}
				\Rtransarg{\muxmulevelsetvalue}_{\alpha} (\Transport \SigmatTan^{\alpha}) \SigmatTan^a X_a
				-
				\frac{1}{2(\upmu - \phi \frac{\muxmulevelsetvalue}{\Lunit \upmu})}
				\Rtransarg{\muxmulevelsetvalue}_{\alpha} \SigmatTan^{\alpha} (\Transport \SigmatTan^a) X_a
				+
				\frac{1}{2(\upmu - \phi \frac{\muxmulevelsetvalue}{\Lunit \upmu})}
				\Rtransarg{\muxmulevelsetvalue}_{\alpha} \uLunit^{\alpha}
				\SigmatTan_a \Transport \SigmatTan^a 
					 \\
	& \ \
			-
			\frac{1}{2 \upmu (\upmu - \phi \frac{\muxmulevelsetvalue}{\Lunit \upmu})}
			\Rtransarg{\muxmulevelsetvalue}_{\alpha}
			\SigmatTan^{\alpha} \SigmatTan^a \muX X_a
			-
			\frac{1}{2(\upmu - \phi \frac{\muxmulevelsetvalue}{\Lunit \upmu})}
			\Rtransarg{\muxmulevelsetvalue}_{\alpha}
			\SigmatTan^{\alpha} \SigmatTan^a \Lunit X_a
				 \\
	& \ \
			+
			\frac{1}{4 \upmu (\upmu - \phi \frac{\muxmulevelsetvalue}{\Lunit \upmu})}
			\Rtransarg{\muxmulevelsetvalue}_{\alpha}
			(\muX \uLunit^{\alpha})
			|\SigmatTan|_g^2
			+
			\frac{1}{4(\upmu - \phi \frac{\muxmulevelsetvalue}{\Lunit \upmu})}
			\Rtransarg{\muxmulevelsetvalue}_{\alpha}
			(\Lunit \uLunit^{\alpha})
			|\SigmatTan|_g^2
						\\
	& \ \
			+
			\frac{1}{4 \upmu (\upmu - \phi \frac{\muxmulevelsetvalue}{\Lunit \upmu})}
			\Rtransarg{\muxmulevelsetvalue}_{\alpha}
			\uLunit^{\alpha}
			(\muX g_{ab}) \SigmatTan^a \SigmatTan^b
			+
			\frac{1}{4(\upmu - \phi \frac{\muxmulevelsetvalue}{\Lunit \upmu})}
			\Rtransarg{\muxmulevelsetvalue}_{\alpha}
			\uLunit^{\alpha}
			(\Lunit g_{ab}) \SigmatTan^a \SigmatTan^b
					\\
	& \ \			
			+
			\frac{(\Rtransarg{\muxmulevelsetvalue} \upmu) \phi \frac{\muxmulevelsetvalue}{\Lunit \upmu}}{\upmu^2(\upmu - \phi \frac{\muxmulevelsetvalue}{\Lunit \upmu})}
			\left\lbrace
				\frac{1}{2}
				\Lunit_{\alpha} \SigmatTan^{\alpha} \SigmatTan^a X_a
				-
				\frac{1}{4}
				\Lunit_{\alpha} \uLunit^{\alpha}
				|\SigmatTan|_g^2
			\right\rbrace
			+
			\frac{(\Lunit \upmu) \phi \frac{\muxmulevelsetvalue}{\Lunit \upmu}}{\upmu^2}
			\left\lbrace
				\frac{1}{2}
				\Lunit_{\alpha} \SigmatTan^{\alpha} \SigmatTan^a X_a
				-
				\frac{1}{4}
				\Lunit_{\alpha} \uLunit^{\alpha}
				|\SigmatTan|_g^2
			\right\rbrace
						\\
		& \ \
			+
			\frac{(\Roughtoritangentvectorfieldarg{\muxmulevelsetvalue} \upmu) \phi \frac{\muxmulevelsetvalue}{\Lunit \upmu}}{\upmu(\upmu - \phi \frac{\muxmulevelsetvalue}{\Lunit \upmu})}
			\left\lbrace
				\frac{1}{2}
				\Lunit_{\alpha} \SigmatTan^{\alpha} \SigmatTan^a X_a
				-
				\frac{1}{4}
				\Lunit_{\alpha} \uLunit^{\alpha}
				|\SigmatTan|_g^2
			\right\rbrace
			+
			\frac{\phi' \frac{\muxmulevelsetvalue}{\Lunit \upmu}}{\upmu(\upmu - \phi \frac{\muxmulevelsetvalue}{\Lunit \upmu})}
			\left\lbrace
				-
				\frac{1}{2}
				\Lunit_{\alpha} \SigmatTan^{\alpha} \SigmatTan^a X_a
				+
				\frac{1}{4}
				\Lunit_{\alpha} \uLunit^{\alpha}
				|\SigmatTan|_g^2
			\right\rbrace
				\\
	& \ \
		+
		\phi \frac{\muxmulevelsetvalue (\muX \Lunit \upmu)}{\upmu(\upmu - \phi \frac{\muxmulevelsetvalue}{\Lunit \upmu}) (\Lunit \upmu)^2}
		\left\lbrace
			\frac{1}{2}
			\Lunit_{\alpha} \SigmatTan^{\alpha} \SigmatTan^a X_a
			-
			\frac{1}{4}
				\Lunit_{\alpha} \uLunit^{\alpha}
				|\SigmatTan|_g^2
		\right\rbrace
			\\
	& \ \
		+
		\phi \frac{\muxmulevelsetvalue  (\Lunit \Lunit \upmu)}{(\upmu - \phi \frac{\muxmulevelsetvalue}{\Lunit \upmu}) (\Lunit \upmu)^2}
		\left\lbrace
			\frac{1}{2}
			\Lunit_{\alpha} \SigmatTan^{\alpha} \SigmatTan^a X_a
			-
			\frac{1}{4}
				\Lunit_{\alpha} \uLunit^{\alpha}
				|\SigmatTan|_g^2
		\right\rbrace
				\\
	& \ \	
		+
		\frac{1}{\upmu - \phi \frac{\muxmulevelsetvalue}{\Lunit \upmu}}
		\left\lbrace
				-
				\frac{1}{2}
				(\muX X_{\alpha}) \SigmatTan^{\alpha} \SigmatTan^a X_a
				+
				\frac{1}{4}
				(\muX X_{\alpha}) \uLunit^{\alpha}
				|\SigmatTan|_g^2
		\right\rbrace
				\\
		& \ \ 
		+
		\frac{\upmu}{\upmu - \phi \frac{\muxmulevelsetvalue}{\Lunit \upmu}}
		\left\lbrace
				-
				\frac{1}{2}
				(\Lunit X_{\alpha}) \SigmatTan^{\alpha} \SigmatTan^a X_a
				+
				\frac{1}{4}
				(\Lunit X_{\alpha}) \uLunit^{\alpha}
				|\SigmatTan|_g^2
		\right\rbrace
			\\
	& \ \
	+
	\frac{1}{\upmu - \phi \frac{\muxmulevelsetvalue}{\Lunit \upmu}}
	\left(
		\phi \frac{\muxmulevelsetvalue}{\upmu \Lunit \upmu}
		+
		\Rtransnormsmallfactorarg{\muxmulevelsetvalue} 
	\right)
	\left\lbrace
				-
				\frac{1}{2}
				(\muX \Lunit_{\alpha}) \SigmatTan^{\alpha} \SigmatTan^a X_a
				+
				\frac{1}{4}
				(\muX \Lunit_{\alpha}) \uLunit^{\alpha}
				|\SigmatTan|_g^2
	\right\rbrace
		\\
& \ \
	+
	\frac{\upmu}{\upmu - \phi \frac{\muxmulevelsetvalue}{\Lunit \upmu}}
	\left(
		\phi \frac{\muxmulevelsetvalue}{\upmu \Lunit \upmu}
		+
		\Rtransnormsmallfactorarg{\muxmulevelsetvalue} 
	\right)
	\left\lbrace
				-
				\frac{1}{2}
				(\Lunit \Lunit_{\alpha}) \SigmatTan^{\alpha} \SigmatTan^a X_a
				+
				\frac{1}{4}
				(\Lunit \Lunit_{\alpha}) \uLunit^{\alpha}
				|\SigmatTan|_g^2
	\right\rbrace 
		\\
& \ \
	+
	\frac{(\muX \Rtransnormsmallfactorarg{\muxmulevelsetvalue})}{\upmu - \phi \frac{\muxmulevelsetvalue}{\Lunit \upmu}} 
		\left\lbrace
				-
				\frac{1}{2}
				\Lunit_{\alpha} \SigmatTan^{\alpha} \SigmatTan^a X_a
				+
				\frac{1}{4}
				\Lunit_{\alpha} \uLunit^{\alpha}
			|\SigmatTan|_g^2
		\right\rbrace
			\\
& \ \
	+
	\frac{\upmu (\Lunit \Rtransnormsmallfactorarg{\muxmulevelsetvalue})}{\upmu - \phi \frac{\muxmulevelsetvalue}{\Lunit \upmu}} 
	\left\lbrace
				-
				\frac{1}{2}
				\Lunit_{\alpha} \SigmatTan^{\alpha} \SigmatTan^a X_a
				+
				\frac{1}{4}
				\Lunit_{\alpha} \uLunit^{\alpha}
			|\SigmatTan|_g^2
		\right\rbrace
			\\
& \ \
	+
	\frac{1}{(\upmu - \phi \frac{\muxmulevelsetvalue}{\Lunit \upmu})}
	\frac{(\muX \Lunit \timefunctionarg{\muxmulevelsetvalue})}{(\Lunit \timefunctionarg{\muxmulevelsetvalue})^2}
	\left\lbrace
				-
				\frac{1}{2}
				(\smoothtorusproject_{\alpha}^{\ \beta} \partial_{\beta} \timefunctionarg{\muxmulevelsetvalue}) \SigmatTan^{\alpha} \SigmatTan^a X_a
				+
				\frac{1}{4}
				(\smoothtorusproject_{\alpha}^{\ \beta}
			\partial_{\beta} \timefunctionarg{\muxmulevelsetvalue}) 
			\uLunit^{\alpha}
			|\SigmatTan|_g^2
		\right\rbrace 
			\\
& \ \
	+
	\frac{\upmu}{(\upmu - \phi \frac{\muxmulevelsetvalue}{\Lunit \upmu})}
	\frac{(\Lunit \Lunit \timefunctionarg{\muxmulevelsetvalue})}{(\Lunit \timefunctionarg{\muxmulevelsetvalue})^2}
	\left\lbrace
				-
				\frac{1}{2}
				(\smoothtorusproject_{\alpha}^{\ \beta} \partial_{\beta} \timefunctionarg{\muxmulevelsetvalue}) \SigmatTan^{\alpha} \SigmatTan^a X_a
				+
				\frac{1}{4}
				(\smoothtorusproject_{\alpha}^{\ \beta}
			\partial_{\beta} \timefunctionarg{\muxmulevelsetvalue}) 
			\uLunit^{\alpha}
			|\SigmatTan|_g^2
		\right\rbrace
				\\
	& \ \
			+
			\frac{1}{(\upmu - \phi \frac{\muxmulevelsetvalue}{\Lunit \upmu})}
			\frac{1}{\Lunit \timefunctionarg{\muxmulevelsetvalue}}
			\left\lbrace
				\frac{1}{2}
				[\muX (\smoothtorusproject_{\alpha}^{\ \beta} \partial_{\beta} \timefunctionarg{\muxmulevelsetvalue})]
				\SigmatTan^{\alpha} \SigmatTan^a X_a
				-
				\frac{1}{4}
				[\muX (\smoothtorusproject_{\alpha}^{\ \beta} \partial_{\beta} \timefunctionarg{\muxmulevelsetvalue})]
				\uLunit^{\alpha}
				|\SigmatTan|_g^2
		\right\rbrace 
			\\
		& \ \ +
			\frac{\upmu}{(\upmu - \phi \frac{\muxmulevelsetvalue}{\Lunit \upmu})}
			\frac{1}{\Lunit \timefunctionarg{\muxmulevelsetvalue}}
			\left\lbrace
				\frac{1}{2}
				[\Lunit (\smoothtorusproject_{\alpha}^{\ \beta} \partial_{\beta} \timefunctionarg{\muxmulevelsetvalue})]
				\SigmatTan^{\alpha} \SigmatTan^a X_a
				-
				\frac{1}{4} 
				[\Lunit (\smoothtorusproject_{\alpha}^{\ \beta} \partial_{\beta} \timefunctionarg{\muxmulevelsetvalue})]
				\uLunit^{\alpha}
				|\SigmatTan|_g^2
		\right\rbrace.
	\end{split}
	\end{align}
	Combining
	\eqref{E:NINTHPROOFSTEPKEYIDPUTANGENTCURRENT}--\eqref{E:THIRTEENTHPROOFSTEPKEYIDPUTANGENTCURRENT},
	using the identity 
	$
	\frac{1}{2}
			\Lunit_{\alpha} \SigmatTan^{\alpha} \SigmatTan^a X_a
			-
			\frac{1}{4}
				\Lunit_{\alpha} \uLunit^{\alpha}
				|\SigmatTan|_g^2
	= \frac{1}{2} |\angV|_{\gtorus}^2
	$
	(which follows from \eqref{E:INNERPRODUCTOFLANDULISMINUS2}
	and the identities 
	$\SigmatTan^{\gamma} \Lunit_{\gamma} = - \SigmatTan^a X_a$
	and
	$
	|\SigmatTan|_g^2 
		-
	(\SigmatTan^a X_a)^2
	=
	|\angV|_{\gtorus}^2
	$
	mentioned above)
	to substitute for the factors 
	$
	\left\lbrace
			\frac{1}{2}
			\Lunit_{\alpha} \SigmatTan^{\alpha} \SigmatTan^a X_a
			-
			\frac{1}{4}
				\Lunit_{\alpha} \uLunit^{\alpha}
				|\SigmatTan|_g^2
	\right\rbrace
	$
	on RHS~\eqref{E:THIRTEENTHPROOFSTEPKEYIDPUTANGENTCURRENT},
	and carrying out straightforward algebraic calculations,
	we deduce the desired identity \eqref{E:KEYIDPUTANGENTCURRENTCONTRACTEDAGAINSTVECTORFIELD}.
	We clarify that we explicitly placed the perfect $\Rtransarg{\muxmulevelsetvalue}$-derivative and $\roughangdiv$-derivative terms
	as the first and third terms on RHS~\eqref{E:KEYIDPUTANGENTCURRENTCONTRACTEDAGAINSTVECTORFIELD},
	that we placed the remaining terms involving a first derivative of $\SigmatTan$ 
	on RHS~\eqref{E:PRINCIPALERRORTERMHAVETOCONTROLKEYIDPUTANGENTCURRENTCONTRACTEDAGAINSTVECTORFIELD},
	that we placed all terms that are quadratic in $\SigmatTan$ 
	(without depending on the first derivatives of $\SigmatTan$)
	on RHS~\eqref{E:LOWERORDERERRORTERMHAVETOCONTROLKEYIDPUTANGENTCURRENTCONTRACTEDAGAINSTVECTORFIELD},
	and that we added the term
	$
	\frac{1}{2} \CurrentboundaryerrorperfectRderivative[\SigmatTan,\SigmatTan]
			\mytr_{\gtorusroughfirstfund} \deform{\Rtransarg{\muxmulevelsetvalue}}
	$
	to RHS~\eqref{E:KEYIDPUTANGENTCURRENTCONTRACTEDAGAINSTVECTORFIELD} (as the second term)
	and subtracted it on RHS~\eqref{E:LOWERORDERERRORTERMHAVETOCONTROLKEYIDPUTANGENTCURRENTCONTRACTEDAGAINSTVECTORFIELD}
	(as the first term).
	
	To prove \eqref{E:KEYCOERCIVITYPERFECTRDERIVAIVEERRORPUTANGENTCURRENTCONTRACTEDAGAINSTVECTORFIELD}, 
	we first use
	\eqref{E:SMOOTHGINVERSEABEXPRESSION},
	Lemma~\ref{L:COMMUTATORSTOCOORDINATES},
	and
	Lemma~\ref{L:SCHEMATICSTRUCTUREOFVARIOUSTENSORSINTERMSOFCONTROLVARS} 
to deduce that
$
|\angD \timefunctionarg{\muxmulevelsetvalue}|_{\gtorus}^2
= 
(\gtorus^{-1})^{AB} (\geop{x^A} \timefunctionarg{\muxmulevelsetvalue}) \geop{x^B} \timefunctionarg{\muxmulevelsetvalue}
= 
\smoothfunction(\controlvars) 
\cdot 
(\geop{x^2}\timefunctionarg{\muxmulevelsetvalue},\geop{x^3}\timefunctionarg{\muxmulevelsetvalue})
\cdot
(\geop{x^2}\timefunctionarg{\muxmulevelsetvalue},\geop{x^3}\timefunctionarg{\muxmulevelsetvalue})
$,
where the last expression is schematic.
From this expression,
\eqref{E:SMALLDERIVATIVESLINFTYESTIMATESFORROUGHTIMEFUNCTIONANDDERIVATIVES},
the bootstrap assumptions, 
and Cor.\,\ref{C:IMPROVEAUX},
we deduce that $|\angD \timefunctionarg{\muxmulevelsetvalue}|_{\gtorus} \lesssim \fundbootsmall$.
From this estimate,
\eqref{E:LDERIVATIVEOFROUGHTIMEFUNCTIONISAPPROXIMATELYUNITY}, 
the $\gtorus$ Cauchy--Schwarz inequality,
and Young's inequality,
we deduce the following bound for the last product on
RHS~\eqref{E:PERFECTRDERIVAIVEERRORPUTANGENTCURRENTCONTRACTEDAGAINSTVECTORFIELD}:
$
\left|
\frac{1}{\Lunit \timefunctionarg{\muxmulevelsetvalue}}
				X_a \SigmatTan^a
				\angVarg{\alpha}
				\angDarg{\alpha}
				\timefunctionarg{\muxmulevelsetvalue}
\right|
\lesssim
\fundbootsmall
(X_a \SigmatTan^a)^2
+
\fundbootsmall
|\angV|_{\gtorus}^2
$.
Using a similar argument that also takes into account
the identity \eqref{E:IDENTITYFORRTRANSNORMSMALLFACTORSQUARED}
and the estimate
\eqref{E:LDERIVATIVEOFROUGHTIMEFUNCTIONISAPPROXIMATELYUNITY},
we deduce that $|\Rtransnormsmallfactorarg{\muxmulevelsetvalue}| \lesssim \fundbootsmall$.
From these estimates,
definition \eqref{E:PERFECTRDERIVAIVEERRORPUTANGENTCURRENTCONTRACTEDAGAINSTVECTORFIELD},	
and the identity $
	|\SigmatTan|_g^2 
	=
	(\SigmatTan^a X_a)^2
	+
	|\SigmatTan|_{\gtorus}^2
	$
	noted just above \eqref{E:NINTHPROOFSTEPKEYIDPUTANGENTCURRENT},
we conclude \eqref{E:KEYCOERCIVITYPERFECTRDERIVAIVEERRORPUTANGENTCURRENTCONTRACTEDAGAINSTVECTORFIELD}.
		
\end{proof}

\subsection{The main elliptic-hyperbolic integral identity}
\label{SS:INTEGRALIDENTITYFORELLIPTICHYPERBOLICCURRENT}
We now state and prove the main elliptic-hyperbolic integral identity.

\begin{proposition}[The main elliptic-hyperbolic integral identity]
	\label{P:INTEGRALIDENTITYFORELLIPTICHYPERBOLICCURRENT}
	Let $\SigmatTan$ be a $\Sigma_t$-tangent vectorfield
	and let $\ellipticCoerciveQuadratic[\pmb{\partial} \SigmatTan,\pmb{\partial} \SigmatTan]$
	be the quadratic form from Def.\,\ref{D:NULLHYPERSURFACEADAPTEDCOERCIVEQUADRATICFORM}.
	Then for any $\timefunction_1 \leq \timefunction_2$ 
	and $u_1 \leq u_2$,
	the following integral identity holds:
	\begin{align}	
	\begin{split} \label{E:INTEGRALIDENTITYFORELLIPTICHYPERBOLICCURRENT}
		&
		\int_{\twoargMrough{[\timefunction_1,\timefunction_2),[u_1,u_2]}{\muxmulevelsetvalue}}
			\frac{1}{\Lunit \timefunctionarg{\muxmulevelsetvalue}}
			\ellipticCoerciveQuadratic[\pmb{\partial} \SigmatTan,\pmb{\partial} \SigmatTan]
		\, \volMRoughCoordinates
		+
		\int_{\twoargroughtori{\timefunction_2,u_2}{\muxmulevelsetvalue}}
			\CurrentboundaryerrorperfectRderivative[\SigmatTan,\SigmatTan]
		\, \volroughtorus
			\\
		& =	
			\int_{\twoargroughtori{\timefunction_2,u_1}{\muxmulevelsetvalue}}
				\CurrentboundaryerrorperfectRderivative[\SigmatTan,\SigmatTan]
			\, \volroughtorus
			+
			\int_{\twoargroughtori{\timefunction_1,u_2}{\muxmulevelsetvalue}}
				\CurrentboundaryerrorperfectRderivative[\SigmatTan,\SigmatTan]
			\, \volroughtorus
				-
			\int_{\twoargroughtori{\timefunction_1,u_1}{\muxmulevelsetvalue}}
				\CurrentboundaryerrorperfectRderivative[\SigmatTan,\SigmatTan]
			\, \volroughtorus
				   \\
		& \ \
			+
			\int_{\hypthreearg{\timefunction_1}{[u_1,u_2]}{\muxmulevelsetvalue}}
				\left\lbrace
						\Currentboundaryerrorhavetocontrolprincipal[\SigmatTan,\pmb{\partial} \SigmatTan]
						+
						\Currentboundaryerrorhavetocontrollowerorder[\SigmatTan,\SigmatTan]
					\right\rbrace
			\, \volRoughHypersurface
				 \\
		& \ \
			-
			\int_{\hypthreearg{\timefunction_2}{[u_1,u_2]}{\muxmulevelsetvalue}}
				\left\lbrace
						\Currentboundaryerrorhavetocontrolprincipal[\SigmatTan,\pmb{\partial} \SigmatTan]
						+
						\Currentboundaryerrorhavetocontrollowerorder[\SigmatTan,\SigmatTan]
					\right\rbrace
			\, \volRoughHypersurface
			+
			\int_{\twoargMrough{[\timefunction_1,\timefunction_2),[u_1,u_2]}{\muxmulevelsetvalue}}
				\EllipticHyperbolicCurrentIntegralIdentityTotalSpacetimeErrorTerm[\SigmatTan,\pmb{\partial} \SigmatTan]
			\, \volMRoughCoordinates,
	\end{split}
	\end{align}
	where $\CurrentboundaryerrorperfectRderivative[\SigmatTan,\SigmatTan]$
	is defined by \eqref{E:PERFECTRDERIVAIVEERRORPUTANGENTCURRENTCONTRACTEDAGAINSTVECTORFIELD}
	and is positive definite in the sense that the pointwise estimate \eqref{E:KEYCOERCIVITYPERFECTRDERIVAIVEERRORPUTANGENTCURRENTCONTRACTEDAGAINSTVECTORFIELD} holds,
	the error terms
	$\Currentboundaryerrorhavetocontrolprincipal[\SigmatTan,\pmb{\partial} \SigmatTan]$
	and
	$\Currentboundaryerrorhavetocontrollowerorder[\SigmatTan,\SigmatTan]$
	are defined in 
	\eqref{E:PRINCIPALERRORTERMHAVETOCONTROLKEYIDPUTANGENTCURRENTCONTRACTEDAGAINSTVECTORFIELD}--\eqref{E:LOWERORDERERRORTERMHAVETOCONTROLKEYIDPUTANGENTCURRENTCONTRACTEDAGAINSTVECTORFIELD},
	\begin{align}
	\begin{split}  \label{E:ELLIPTICHYPERBOLICINTEGRALIDENTITYBULKERRORTERM}
		\EllipticHyperbolicCurrentIntegralIdentityTotalSpacetimeErrorTerm[\SigmatTan,\pmb{\partial} \SigmatTan]
		& 
		\eqdef
		\frac{1}{\Lunit \timefunctionarg{\muxmulevelsetvalue}}
				\Big\lbrace
					\mathfrak{J}_{(\textnormal{Antisymmetric})}[\pmb{\partial} \SigmatTan,\pmb{\partial} \SigmatTan]
					+
					\mathfrak{J}_{(\textnormal{Div})}[\pmb{\partial} \SigmatTan,\pmb{\partial} \SigmatTan]
					+
					\upmu 
					\mathfrak{J}_{(\pmb{\partial} \frac{1}{\upmu})}[\SigmatTan,\pmb{\partial} \SigmatTan]
						\\
		& \ \ \ \ \ \ \ \ \ \ \ \ \ \ \
					+
					\mathfrak{J}_{(\textnormal{Absorb-1})}[\SigmatTan,\pmb{\partial} \SigmatTan]
					+
					\mathfrak{J}_{(\textnormal{Absorb-2})}[\SigmatTan,\pmb{\partial} \SigmatTan]
					+
					\mathfrak{J}_{(\textnormal{Material})}[\pmb{\partial} \SigmatTan,\pmb{\partial} \SigmatTan]
					+
					\mathfrak{J}_{(\textnormal{Null Geometry})}[\SigmatTan,\pmb{\partial} \SigmatTan]
			\Big\rbrace,
	\end{split}
	\end{align}
	and the error terms 
	$\mathfrak{J}_{(\textnormal{Antisymmetric})}[\pmb{\partial} \SigmatTan,\pmb{\partial} \SigmatTan], 
	\cdots, 
	\mathfrak{J}_{(\textnormal{Null Geometry})}[\SigmatTan,\pmb{\partial} \SigmatTan]$
	on RHS~\eqref{E:ELLIPTICHYPERBOLICINTEGRALIDENTITYBULKERRORTERM}
	are defined in
	\eqref{E:ANTISYMMETRICNULLCURRENTSPACETIMERRORTERM}--\eqref{E:DERIVATIVESOFNULLGEOMETRYNULLCURRENTSPACETIMERRORTERM}.
\end{proposition}

\begin{proof}
	We consider the divergence identity \eqref{E:COVARIANTDIVERGENCEIDENTITYFORELLIPTICHYPERBOLICCURRENT}
	with $\weight \eqdef \frac{1}{\upmu}$.
	We integrate the identity over $\twoargMrough{[\timefunction_1,\timefunction_2),[u_1,u_2]}{\muxmulevelsetvalue}$
	with respect to the canonical volume form \eqref{E:VOLFORMACOUSTICALMETRICROUGHADAPTED} of $\gfour$ in rough adapted coordinates,
	apply the divergence theorem,
	and take into account
	Def.\,\ref{D:ROUGHVOLFORMS}
	and
	the identity \eqref{E:VOLFORMCANONICALHYPGROUGHADAPTED},
	thereby obtaining the following identity,
	where $\hypunitnormalarg{\muxmulevelsetvalue}$ is the future-directed $\gfour$-unit normal to the $\gfour$-spacelike hypersurfaces
	$\hypthreearg{\timefunction}{[u_1,u_2]}{\muxmulevelsetvalue}$
	(see Prop.\,\ref{P:BASICPROPERTIESOFROUGHVECTORFIELDS}): 
		\begin{align}	
		\begin{split} \label{E:FIRSTPROOFSTEPINTEGRALIDENTITYFORELLIPTICHYPERBOLICCURRENT}
		&
		\int_{\twoargMrough{[\timefunction_1,\timefunction_2),[u_1,u_2]}{\muxmulevelsetvalue}}
			\frac{1}{\Lunit \timefunctionarg{\muxmulevelsetvalue}}
			\ellipticCoerciveQuadratic[\pmb{\partial} \SigmatTan,\pmb{\partial} \SigmatTan]
		\, \volMRoughCoordinates
			\\
		& =	
					-
					\int_{\hypthreearg{\timefunction_2}{[u_1,u_2]}{\muxmulevelsetvalue}}
						\frac{|\Rtransarg{\muxmulevelsetvalue}|_{\hypg}}
						{\upmu}
						\hypunitnormalarg{\muxmulevelsetvalue}_{\alpha}
						\ehcurrent^{\alpha}[\SigmatTan,\pmb{\partial} \SigmatTan] 
					\, \volRoughHypersurface
					+
					\int_{\hypthreearg{\timefunction_1}{[u_1,u_2]}{\muxmulevelsetvalue}}
						\frac{|\Rtransarg{\muxmulevelsetvalue}|_{\hypg}}
						{\upmu}
						\hypunitnormalarg{\muxmulevelsetvalue}_{\alpha}
						\ehcurrent^{\alpha}[\SigmatTan,\pmb{\partial} \SigmatTan] 
					\, \volRoughHypersurface
						\\
		& \ \
			+
			\int_{\twoargMrough{[\timefunction_1,\timefunction_2),[u_1,u_2]}{\muxmulevelsetvalue}}
				\EllipticHyperbolicCurrentIntegralIdentityTotalSpacetimeErrorTerm[\SigmatTan,\pmb{\partial} \SigmatTan]
			\, \volMRoughCoordinates.
	\end{split}
	\end{align}
	We emphasize that there are no null hypersurface
	boundary integrals in \eqref{E:FIRSTPROOFSTEPINTEGRALIDENTITYFORELLIPTICHYPERBOLICCURRENT} 
	because the current $\ehcurrent^{\alpha}[\SigmatTan,\pmb{\partial} \SigmatTan]$ defined in \eqref{E:PUTANGENTELLIPTICHYPERBOLICCURRENT}
	is tangent to the null hypersurface $\nullhyparg{u}$
	(and thus $\Lunit_{\alpha} \ehcurrent^{\alpha}[\SigmatTan,\pmb{\partial} \SigmatTan] = 0$),
	and the signs of the first two integrals on RHS~\eqref{E:FIRSTPROOFSTEPINTEGRALIDENTITYFORELLIPTICHYPERBOLICCURRENT}
	are tied to the Lorentzian nature of $\gfour$ and the fact that the $\gfour$-timelike vectorfield 
	$\hypunitnormalarg{\muxmulevelsetvalue}^{\alpha}$ is future-pointing
	(and thus outward-pointing to $\twoargMrough{[\timefunction_1,\timefunction_2),[u_1,u_2]}{\muxmulevelsetvalue}$ along 
	$\hypthreearg{\timefunction_2}{[u_1,u_2]}{\muxmulevelsetvalue}$
	and inward-pointing to $\twoargMrough{[\timefunction_1,\timefunction_2),[u_1,u_2]}{\muxmulevelsetvalue}$ along 
	$\hypthreearg{\timefunction_1}{[u_1,u_2]}{\muxmulevelsetvalue}$).
	Next, we use the identities \eqref{E:KEYIDPUTANGENTCURRENTCONTRACTEDAGAINSTVECTORFIELD}
	and \eqref{E:RTRANSINTEGRALONROUGHTORI}
	to re-express the integral $\int_{\hypthreearg{\timefunction_2}{[u_1,u_2]}{\muxmulevelsetvalue}} \cdots$
	on RHS~\eqref{E:FIRSTPROOFSTEPINTEGRALIDENTITYFORELLIPTICHYPERBOLICCURRENT}
	as follows, where we note that integral of the perfect divergence term
	$- \roughangdiv \cdots$ on RHS~\eqref{E:KEYIDPUTANGENTCURRENTCONTRACTEDAGAINSTVECTORFIELD} 
	vanishes when integrated over any rough torus $\twoargroughtori{\timefunction,u}{\muxmulevelsetvalue}$: 
	\begin{align} 
	\begin{split} \label{E:SECONDPROOFSTEPINTEGRALIDENTITYFORELLIPTICHYPERBOLICCURRENT}
			-
			\int_{\hypthreearg{\timefunction_2}{[u_1,u_2]}{\muxmulevelsetvalue}}
						\frac{|\Rtransarg{\muxmulevelsetvalue}|_{\hypg}}
						{\upmu}
						\hypunitnormalarg{\muxmulevelsetvalue}_{\alpha}
						\ehcurrent^{\alpha}[\SigmatTan,\pmb{\partial} \SigmatTan] 
					\, \volRoughHypersurface
		& = 
				-
				\int_{\twoargroughtori{\timefunction_2,u_2}{\muxmulevelsetvalue}}
					\CurrentboundaryerrorperfectRderivative[\SigmatTan,\SigmatTan]
				\, \volroughtorus
				+
				\int_{\twoargroughtori{\timefunction_2,u_1}{\muxmulevelsetvalue}}
					\CurrentboundaryerrorperfectRderivative[\SigmatTan,\SigmatTan]
				\, \volroughtorus
					\\
		& \ \
			-
				\int_{\hypthreearg{\timefunction_2}{[u_1,u_2]}{\muxmulevelsetvalue}}
				\left\lbrace
						\Currentboundaryerrorhavetocontrolprincipal[\SigmatTan,\pmb{\partial} \SigmatTan]
						+
						\Currentboundaryerrorhavetocontrollowerorder[\SigmatTan,\SigmatTan]
					\right\rbrace
			\, \volRoughHypersurface.
	\end{split}
	\end{align}
	Moreover, we note that \eqref{E:SECONDPROOFSTEPINTEGRALIDENTITYFORELLIPTICHYPERBOLICCURRENT}
	also holds with $\timefunction_1$ in place of $\timefunction_2$.
	Using these two identities to substitute for the first two integrals
	on RHS~\eqref{E:FIRSTPROOFSTEPINTEGRALIDENTITYFORELLIPTICHYPERBOLICCURRENT},
	we arrive at the identity \eqref{E:INTEGRALIDENTITYFORELLIPTICHYPERBOLICCURRENT}.
	
\end{proof}


\section{Pointwise estimates for the error terms in the commuted wave equations} 
\label{S:POINTWISESTIMATESFORWAVEEQUATIONS}
In this section, we derive pointwise estimates for the error terms that arise
when we commute the wave equations \eqref{E:VELOCITYWAVEEQUATION}--\eqref{E:ENTROPYWAVEEQUATION}
up to $\Ntop$ times, where we recall that $\Ntop$ is a fixed integer satisfying \eqref{E:NTOPLARGENESSASSUMPTION}.
More precisely,
for $\Psi \in \{\RRiemann,\LRiemann,v^2,v^3,\Ent\}$ 
and $1 \leq N \leq \Ntop$,
we derive pointwise estimates for the inhomogeneous term in the $\upmu$-weighted geometric wave equation
$\upmu \Box_{\gfour} \tander^N \Psi = \mathfrak{G}$ satisfied by $\tander^N \Psi$. 
These pointwise estimates are a preliminary ingredient for the $L^2$ estimates that we derive later on.
Many of the terms appearing in $\mathfrak{G}$ are harmless from the point of view of regularity 
and the strength of their singularity;  
the bulk of our effort goes towards the most difficult terms, 
which involve the top-order derivatives of the eikonal function,
which we handle by using the modified quantities from Def.\,\ref{D:FULLYANDPARTIALLYMODIFIEDQUANTITIES}.

\subsection{Identification of the most difficult error terms in the commuted wave equations}
\label{SS:MOSTDIFFICULTTERMSINCOMMUTEDWAVEEQUATIONS}
Most of the terms in the commuted wave equations are harmless from the point of view of regularity 
and the strength of their singularity. The next definition captures these ``harmless'' error terms.

\subsubsection{Harmless wave equation error terms}
\label{SSS:HARMLESSWAVEEQUATIONERRORTERMS}
\begin{definition}[Harmless wave equation error terms] 
\label{D:HARMLESSWAVE}
Let $1 \leq N \leq \Ntop$.
We define $\HarmlessWave{N}$ to be any term
that satisfies the following pointwise estimate 
on $\twoargMrough{[\timefunction_0,\timefunctionboot),[- \rightu,\leftu]}{\muxmulevelsetvalue}$:
\begin{align} \label{E:HARMLESSWAVE}
\left| \HarmlessWave{N} \right|
&
\lesssim 
\left| \comdersmall^{[1,N+1];1} \wavearray\right| 
+ 
\left|\comdersmall^{[1,N];1} \controlvars\right| 
+
\left|\tandersmall^{[1,N]} \badcontrolvars\right|.
\end{align}
\end{definition}

The following simple lemma shows that $\HarmlessWave{\Ntop - 12}$ terms are small in 
the norm $\| \cdot \|_{L^{\infty}(\twoargroughtori{\timefunction,u}{\muxmulevelsetvalue})}$.

\begin{lemma}[$L^{\infty}$ estimates for $\HarmlessWave{\Ntop - 12}$]  
\label{L:LINFINITYESTIMATESFORHARMLESSWAVEERRORTERMS}
Let $\HarmlessWave{\Ntop - 12}$ be as in Def.\,\ref{D:HARMLESSWAVE}.
Then the following estimate holds for $(\timefunction,u) \in [\timefunction_0,\timefunctionboot) \times [- \rightu,\leftu]$:
\begin{align} \label{E:LINFINITYESTIMATESFORHARMLESSWAVEERRORTERMS}
	\left\| 
		\HarmlessWave{\Ntop - 12} 
	\right\|_{L^{\infty}(\twoargroughtori{\timefunction,u}{\muxmulevelsetvalue})} 
	&
	\lesssim 
	\fundbootsmall.
\end{align}
\end{lemma}

\begin{proof}
The estimate \eqref{E:LINFINITYESTIMATESFORHARMLESSWAVEERRORTERMS} follows from definition~\eqref{E:HARMLESSWAVE} and
Prop.\,\ref{P:IMPROVEMENTOFAUXILIARYBOOTSTRAP}.
\end{proof}

\subsubsection{The most difficult error terms in the commuted wave equations}
\label{SSS:MOSTDIFFICULTWAVETERMS}
In the following proposition, we identify the most difficult error terms
in the commuted wave equations satisfied by the wave variables $\wavearray$.

\begin{proposition}[Identification of the most difficult error terms in the commuted wave equations] 
\label{P:MOSTDIFFICULTWAVETERMS} 
Let $\wavearray \eqdef (\Psi_0,\Psi_1,\Psi_2,\Psi_3,\Psi_4) \eqdef (\RRiemann,\LRiemann,v^2,v^3,\Ent)$ 
be solutions to the covariant wave equations \eqref{E:COVARIANTWAVEEQUATIONSWAVEVARIABLES},
and let $N \leq \Ntop-1$.
We denote the product of $\upmu$ and the RHS of the covariant wave equation satisfied by 
$\Psi_{\iota}$ by $\mathfrak{G}_{\iota}$, 
i.e., 
$\upmu \Box_{\gfour} \Psi_{\iota} = \mathfrak{G}_{\iota}$. 
Then the following wave equations hold for $\iota = 0,1,2,3,4$ and $A = 2,3$:
\begin{subequations}  
\begin{align} \label{E:TOPCOMMUTEDWAVELFIRSTTHENALLYS}
\upmu \Box_{\gfour} (\tanderY^{N-1} \Lunit \Psi_{\iota}) 
& 
= 
\angrmD^{\sharp} \Psi_{\iota} \cdot \upmu \angrmD \tanderY^{N-1} \mytr_{\gtorus}\upchi 
+ 
\tanderY^{N-1} \Lunit \mathfrak{G}_{\iota} 
+ 
\HarmlessWave{N},	
	\\
\begin{split} \label{E:TOPCOMMUTEDWAVEALLYS}
\upmu \Box_{\gfour} (\tanderY^{N-1}\Yvf{A} \Psi_{\iota}) 
& 
= 
(\muX \Psi_{\iota}) \tanderY^{N-1} \Yvf{A} \mytr_{\gtorus} \upchi 
+ 
(\Speed^{-2} X^A) \angrmD^{\sharp} \Psi_{\iota} \cdot \upmu \angrmD \tanderY^{N-1} \mytr_{\gtorus}\upchi 
\\
& \ \
+ 
\tanderY^{N-1}\Yvf{A} \mathfrak{G}_{\iota} 
+ 
\HarmlessWave{N}. 
\end{split}
\end{align}

Moreover, if $1 \leq N \leq \Ntop$ and $\tander^N$
denotes any order $N$ string of $\mathcal{P}_u$-tangent commutator
other than the ones appearing on LHSs~\eqref{E:TOPCOMMUTEDWAVELFIRSTTHENALLYS}--\eqref{E:TOPCOMMUTEDWAVEALLYS},
(i.e., if $\tander^N$ features at least two copies of $\Lunit$ or only a single $\Lunit$ but does not act first like it does in
\eqref{E:TOPCOMMUTEDWAVELFIRSTTHENALLYS}), 
then $\tander^N \Psi_{\iota}$ obeys the following wave equation:
\begin{align} \label{E:TOPCOMMUTEDWAVENOTDIFFICULT}
\upmu \Box_{\gfour} (\tander^N \Psi_{\iota}) 
	& = 
	\tander^N \mathfrak{G}_{\iota} 
	+ 
	\HarmlessWave{N}.
\end{align}
\end{subequations}
\end{proposition}

\begin{proof}
In two spatial dimensions, a detailed proof was provided in \cite[Proposition 13.2]{jLjS2018}, 
with the only difference being that there are no vorticity-involving terms in our definition of $\HarmlessWave{N}$ because we have soaked these terms up into our definition of $\mathfrak{G}_{\iota}$. 
The modification of the argument needed to account for the third space dimension is minimal, and we therefore omit it. 
\end{proof}

\subsection{Pointwise estimates for the difficult product $(\muX \Psi_{\iota}) \tanderY^N \mytr_{\gtorus} \upchi$} 
\label{SS:POINTWISEESTIMATESFORDIFFICULTPRODUCTINVOLVINGTOPORDERANGDERIVSOFCHI}
In this section, 
we derive pointwise estimates tied to the most difficult terms appearing in the commuted wave equations,
specifically the products
$(\muX \Psi_{\iota}) \tanderY^N \mytr_{\gtorus} \upchi$ 
on RHS~\eqref{E:TOPCOMMUTEDWAVEALLYS}
(note that $\tanderY^{N-1} \Yvf{A}$ can be expressed as $\tanderY^N$). 
Our analysis relies on the fully modified quantities 
$\fullymodquant{\tander^N}$ 
from Def.\,\ref{D:FULLYANDPARTIALLYMODIFIEDQUANTITIES}.

\subsubsection{Pointwise estimates for the inhomogeneous terms in the transport equations satisfied by the modified quantities} 
\label{SSS:POINTWISEESTIMATESFORINHOMOGENEOUSTERMSINTRANSPORTEQUATIONSFORMODIFIEDQUANTITIES}
We start with the following lemma, which provides pointwise estimates for 
the inhomogeneous terms in the transport equations satisfied by the fully modified quantities.
The lemma also provides,
for use in Sect.\,\ref{SS:POINTWISEESTIMATEFORPARTIALLYMODIFIEDQUANT}, 
pointwise estimates for the inhomogeneous terms in the transport equations satisfied by 
the partially modified quantities.

\begin{lemma}[Pointwise estimates for inhomogeneous terms tied to the modified quantities]
\label{L:POINTWISEESTIMATESFORINHOMOGENEOUSTERMSINTRANSPORTEQUATIONSFORMODIFIEDQUANTITIES}
Let $N = \Ntop$.

\medskip
\noindent \underline{\textbf{Estimates tied to the fully modified quantities}}.
Let $\tanderY^N \in \mathfrak{Y}^{(N)}$ where $\mathfrak{Y}^{(N)}$
is the set of order $N$ $\ell_{t,u}$-tangential commutator operators from Sect.\,\ref{SS:STRINGSOFCOMMUTATIONVECTORFIELDS}. 
Let $\mathfrak{X}$ be the term defined in \eqref{E:MODIFIEDQUANTITYINHOM}.
Then the following pointwise estimates hold on 
$\twoargMrough{[\timefunction_0,\timefunctionboot),[- \rightu,\leftu]}{\muxmulevelsetvalue}$:
\begin{subequations}
\begin{align}
	\left| 
		\tanderY^N \mathfrak{X} 
		+ 
		\vec{G}_{\Lunit \Lunit} 
		\diamond 
		\muX \tanderY^N \wavearray 
	\right|
	& 
	\lesssim \upmu \left|\tander^{[1,N+1]} \wavearray\right| 
	+ 
	\left|\comdersmall^{[1,N];1}\wavearray\right| 
	+ 
	\left|\tander^{[1,N]}\controlvars\right| 
	+ 
	\left|\tandersmall^{[1,N]} \badcontrolvars\right|, 
		\label{E:POINTWISESUMOFMODIFIEDQUANTITYINHOMANDGLL} 
			\\
	\left|\tander^N \mathfrak{X} \right| 
	& \lesssim 
	\left|\comdersmall^{[1,N+1];1}\wavearray\right| 
	+ 
	\left|\tander^{[1,N]}\controlvars\right| 
	+ 
	\left|\tandersmall^{[1,N]}\badcontrolvars \right|.
		\label{E:POINTWISEESTIMATETANGENTDERIVATIVESOFMODQUANTINHOM} 
\end{align}
\end{subequations}

\medskip

\noindent \underline{\textbf{Estimates tied to the partially modified quantities}}.
If $\tanderY^{N-1} \in \mathfrak{Y}^{(N-1)}$ 
and $\widetilde{\mathfrak{X}}$, 
$\partialmodquantinhom{\tander^{N-1}}$,
${^{(\tanderY^{N-1})}\mathfrak{B}}$
are as defined in \eqref{E:PARTIALMODIFIEDQUANTITYINHOMZEROORDER},
\eqref{E:PARTIALMODIFIEDQUANTITYINHOM} 
(with $\tanderY^{N-1}$ in the role of $\tander^N$), 
and \eqref{E:COMMUTEDPARTIALMODIFIEDQUANTITYINHOM} respectively
(with $\tanderY^{N-1}$ in the role of $\tander^{N-1}$),
then the following pointwise estimates hold on 
$\twoargMrough{[\timefunction_0,\timefunctionboot),[- \rightu,\leftu]}{\muxmulevelsetvalue}$:
\begin{subequations}
\begin{align}
\left|\partialmodquantinhom{\tanderY^{N-1}} \right|
& \lesssim 
\left| \tander^{[1,N]} \wavearray \right|, 
\label{E:POINTWISEBELOWTOPORDERPARTIALMODQUANTINHOM} 
		\\
\left| \Lunit \partialmodquantinhom{\tanderY^{N-1}} \right|, 
	\,
\left| \Yvf{A} \partialmodquantinhom{\tanderY^{N-1}}\right| 
& 
\lesssim 
\left| \tander^{[1,N+1]} \wavearray \right|, 
\label{E:POINTWISELANDYDERIVATIVESOFTOPORDERPARTIALMODQUANTINHOM} 
		\\
\left|{^{(\tanderY^{N-1})}\mathfrak{B}} \right| 
& 
\lesssim 
\left |\tander^{[1,N]} \controlvars\right|. 
	\label{E:POINTWISETANGENTDERIVATIVEOFBINHOM}
\end{align}
\end{subequations}

\end{lemma}

\begin{proof}
All estimates stated in the lemma are straightforward consequences of
Lemma~\ref{L:ANGULARDIFFERENTIALCOMMUTESWITHANGLIE},
Lemma~\ref{L:SCHEMATICSTRUCTUREOFVARIOUSTENSORSINTERMSOFCONTROLVARS},
Lemma~\ref{L:SCHEMATICEXPRESSIONFORANGULARLAPLACIAN},
the commutator estimate \eqref{E:POINTWISEBOUNDCOMMUTATORSTANGENTIALANDTANGENTIALDERIVATIVESONSCALARFUNCTION},
and the estimates of Prop.\,\ref{P:IMPROVEMENTOFAUXILIARYBOOTSTRAP}.
\end{proof}

\subsubsection{Preliminary pointwise estimates for $\fullymodquant{\tanderY^N}$} 
\label{SSS:PRELIMINARYRESULTPOINTWISEESTIMATESFORFULLYMODIFIED}
In the next lemma, 
we use the transport equation \eqref{E:TRANSPORTEQUATIONFORFULLYMODIFIEDQUANTITY}
to derive a preliminary pointwise estimate for the fully modified quantity $\fullymodquant{\tander^N}$
in the case $\tander^N = \tanderY^N$ with $N = \Ntop$, which in practice is the only case in which we need to use
the fully modified quantities. The estimates in the lemma are crucial ingredient
in our proof of Prop.\,\ref{P:MAINPOINTWISEESTIMATESFORFULLYMODIFIED}, 
in which we derive the main pointwise estimate for 
the difficult product $(\muX \RRiemann) \tanderY^N \mytr_{\gtorus} \upchi$.
The proof of the lemma is similar to the proof of \cite[Lemma 11.9]{jSgHjLwW2016}, but due to our 
reliance on the rough adapted coordinates, which are a new feature of the present paper,
we provide complete details here.

\begin{remark}[Boxed constants affect high order energy blowup-rates] 
\label{R:BOXEDCONSTANTS}
In Lemma~\ref{L:POINTWISEESTIMATESFORFULLYMODIFIEDQUANTITIY} and its proof, 
and also throughout the rest of the paper,
the important boxed constants such as $\boxed{2}$ 
affect the blowup-rate of our top-order energy estimates with respect 
to powers of $|\timefunction|^{-1}$; we therefore carefully track these boxed constants.
\end{remark}

\begin{lemma}[Pointwise estimates for $\fullymodquant{\tanderY^{\Ntop}}$] 
\label{L:POINTWISEESTIMATESFORFULLYMODIFIEDQUANTITIY}
Let $N = \Ntop$, and let $\mathfrak{Y}^{(N)}$ and $\angLie_{\mathfrak{Y}}^{(N)}$
be the sets of order $N$ $\ell_{t,u}$-tangential commutator operators from Sect.\,\ref{SS:STRINGSOFCOMMUTATIONVECTORFIELDS}.
Let $\tanderY^N \in \mathfrak{Y}^{(N)}$, and
let $\fullymodquant{\tanderY^N}$ be the corresponding fully modified quantity defined in
\eqref{E:FULLYMODIFIEDQUANTITY} (with $\tanderY^N$ in the role of $\tander^N$).
Let $\argLrough{\muxmulevelsetvalue}$ be the rough null vectorfield defined in \eqref{E:LROUGH},
and let $\FlowmapLrougharg{\muxmulevelsetvalue}$ be the $\timefunction_0$-normalized flow map of
$\argLrough{\muxmulevelsetvalue}$ with respect to the rough adapted coordinates
$(\timefunctionarg{\muxmulevelsetvalue},u,x^2,x^3)$ appearing in Lemma~\ref{L:PROPERTIESOFFLOWMAPOFWIDETILDEL}. 
Moreover, let $\smallneighborhoodofcreasetwoarg{[\timefunction_0,\timefunctionboot]}{\muxmulevelsetvalue}$ be the spacetime neighborhood constructed in 
Prop.\,\ref{P:SHARPCONTROLOFMUANDDERIVATIVES} (specifically, in \eqref{E:SMALLNEIGHBORHOOD}), 
on which we have derived especially sharp control of $\upmu$. 
In addition, if $\mathfrak{K}$ is a spacetime subset, 
let $\mathbf{1}_{\mathfrak{K}}$ denote the characteristic function $\mathfrak{K}$.
Then relative to the rough adapted coordinates $(\timefunction,u,x^2,x^3)$
(see Remark~\ref{R:IMPLICITFUNCTIONALDEPENDENCE}),
the following pointwise estimate holds on $\twoargMrough{[\timefunction_0,\timefunctionboot),[- \rightu,\leftu]}{\muxmulevelsetvalue}$,
where $\tanderY^N$ denotes the same operator in each term in \eqref{E:ESTIMATEFORFULLYMODIFIEDQUANTALONGLROUGH}
(except for the one term in which $\max_{\angLie_{\tanderY^N} \in \angLie_{\mathfrak{Y}}^{(N)}}$ is taken):
\begin{align} 
\begin{split} \label{E:ESTIMATEFORFULLYMODIFIEDQUANTALONGLROUGH}
&
\left| \fullymodquant{\tanderY^N} \right| 
\circ 
\FlowmapLrougharg{\muxmulevelsetvalue}(\timefunction,u,x^2,x^3)  
	\\
& 
\leq 
C 
\left| \fullymodquant{\tanderY^N} \right|(\timefunction_0,u,x^2,x^3) 
	\\
& \ \ 
+ 
\boxed{2}
\int_{\timefunction' = \timefunction_0}^{\timefunction}  
	\left\lbrace 
		\left| \frac{\argLrough{\muxmulevelsetvalue} \upmu}{\upmu} \right|
		\cdot 
		\left| \tanderY^N \mathfrak{X}\right| 
		\cdot 
		\mathbf{1}_{\{\hypthreearg{\timefunction'}{[- \rightu,u]}{\muxmulevelsetvalue} 
		\cap \smallneighborhoodofcreasetwoarg{[\timefunction_0,\timefunctionboot]}{\muxmulevelsetvalue} \}} 
	\right\rbrace  
	\circ 
	\FlowmapLrougharg{\muxmulevelsetvalue}(\timefunction',u,x^2,x^3)
\, \mathrm{d} \timefunction' 
	\\
	& \ \ 
	+ 
	C 
	\fundbootsmall
	\int_{\timefunction' = \timefunction_0}^{\timefunction} 
		\left\lbrace 
			\max_{\angLie_{\tanderY}^N \in \angLie_{\mathfrak{Y}}^{(N)}}
			|\upmu \angLie_{\tanderY}^N \upchi| 
		\right\rbrace 
		\circ 
		\FlowmapLrougharg{\muxmulevelsetvalue}(\timefunction',u,x^2,x^3)
	\, \mathrm{d} \timefunction'
		\\
	& \ \ 
	+ 
	C 
	\int_{\timefunction' = \timefunction_0}^{\timefunction} 
		\left\lbrace 
			\upmu 
			|\tanderY^N(\VortVort,\DivGradEnt)| 
			+ 
			|\tanderY^{\leq N-1}(\VortVort,\DivGradEnt)| 
		\right\rbrace 
		\circ 
		\FlowmapLrougharg{\muxmulevelsetvalue}(\timefunction',u,x^2,x^3)
	\, \mathrm{d} \timefunction'  
		\\
	& \ \ 
	+ 
	C 
	\int_{\timefunction' = \timefunction_0}^{\timefunction} 
		\left\lbrace 
			\upmu 
			|\tanderY^N(\Omega,\GradEnt)| 
			+ 
			|\tanderY^{\leq N-1}(\Omega,\GradEnt)| 
		\right\rbrace 
		\circ 
		\FlowmapLrougharg{\muxmulevelsetvalue}(\timefunction',u,x^2,x^3)
	\, \mathrm{d} \timefunction'  
		\\
	& \ \ 
	+ 
	C 
	\int_{\timefunction' = \timefunction_0}^{\timefunction}  
		\left\lbrace 
			| \comdersmall^{[1,N+1];1} \wavearray| 
			+ 
			|\tander^{[1,N]} \controlvars | 
			+ 
			|\tandersmall^{[1,N]}\badcontrolvars| 
		\right\rbrace 
		\circ \FlowmapLrougharg{\muxmulevelsetvalue}(\timefunction',u,x^2,x^3)
	\, \mathrm{d} \timefunction'. 
\end{split}
\end{align}

\end{lemma}

\begin{proof}
Our analysis relies on the transport equation \eqref{E:TRANSPORTEQUATIONFORFULLYMODIFIEDQUANTITY}.
To control solutions to this equation, 
we will use the integrating factor $I$, 
which we define as follows relative to the rough adapted coordinates:
\begin{align} 
\begin{split} \label{E:DELICATEPOINTWISEESTIMATEINTEGRATINGFACTOR}
I(\timefunction,u,x^2,x^3)  
& 
\eqdef 
\frac{\upmu^2 \circ \FlowmapLrougharg{\muxmulevelsetvalue}(\timefunction_0,u,x^2,x^3)}
{\upmu^2 \circ \FlowmapLrougharg{\muxmulevelsetvalue}(\timefunction,u,x^2,x^3)}
	\\
& = 
	\exp
	\left\lbrace
			-2
			\int_{\timefunction_0}^{\timefunction}
			\left(
				\frac{\argLrough{\muxmulevelsetvalue} \upmu}{\upmu}
			\right)
			\circ \FlowmapLrougharg{\muxmulevelsetvalue}(\timefunction',u,x^2,x^3)
		\, \mathrm{d} \timefunction'
	\right\rbrace,
\end{split}	
\end{align}
where to obtain the second equality in \eqref{E:DELICATEPOINTWISEESTIMATEINTEGRATINGFACTOR},
we have used \eqref{E:FLOWMAPOFLROUGHINROUGHCOORDINATES} and the fundamental theorem of calculus.
We also recall that by \eqref{E:FLOWMAPOFLROUGHINROUGHCOORDINATES}, 
we have $\FlowmapLrougharg{\muxmulevelsetvalue}(\timefunction_0,u,x^2,x^3) = (\timefunction_0,u,x^2,x^3)$
and thus $I(\timefunction_0,u,x^2,x^3) = 1$.

We now fix $\tanderY^N \in \mathfrak{Y}^{(N)}$.
Since $\argLrough{\muxmulevelsetvalue} = \frac{1}{\Lunit \timefunctionarg{\muxmulevelsetvalue}} \Lunit$ (see \eqref{E:LROUGH}),
we can multiply both sides of \eqref{E:TRANSPORTEQUATIONFORFULLYMODIFIEDQUANTITY} 
(with $\tanderY^N$ in the role of $\tander^N$)
by $\frac{1}{\Lunit \timefunctionarg{\muxmulevelsetvalue}}$,
evaluate at $\FlowmapLrougharg{\muxmulevelsetvalue}(\timefunction',u,x^2,x^3)$, 
then multiply both sides by $I(\timefunction',u,x^2,x^3)$,
integrate in rough time, use \eqref{E:FLOWMAPOFLROUGHINROUGHCOORDINATES} and the fundamental
theorem of calculus, and use that $I(\timefunction_0,u,x^2,x^3) = 1$
to deduce the following equation, valid
for $(\timefunction,u,x^2,x^3) \in [\timefunction_0,\timefunctionboot) \times [- \rightu,\leftu] \times \mathbb{T}^2$:
\begin{align} 
\begin{split} \label{E:POINTWISEFULLYMODIFIEDQUANTITYALONGLROUGHSTEP1}
\fullymodquant{\tanderY^N} \circ \FlowmapLrougharg{\muxmulevelsetvalue}(\timefunction,u,x^2,x^3)
& 
= 
\frac{1}{I(\timefunction,u,x^2,x^3)}
\fullymodquant{\tanderY^N}(\timefunction_0,u,x^2,x^3)
	\\
& \ \
+
\int_{\timefunction_0}^{\timefunction} 
	\frac{I(\timefunction',u,x^2,x^3)}{I(\timefunction,u,x^2,x^3)}
	\times 
	\left\lbrace
		\frac{1}{\Lunit \timefunctionarg{\muxmulevelsetvalue}}
		\times
		\text{(RHS~\eqref{E:TRANSPORTEQUATIONFORFULLYMODIFIEDQUANTITY})}
	\right\rbrace 
	\circ 
	\FlowmapLrougharg{\muxmulevelsetvalue}(\timefunction',u,x^2,x^3)
\, \mathrm{d}\timefunction'.
\end{split}
\end{align}

Next, we note the following estimate, which follows from 
\eqref{E:MUMUSTDECREASEININTERESTINGREGION}--\eqref{E:MURATIOBOUNDEDINBORINGREGION}:
\begin{align}\label{E:POINTWISEFULLYMODIFIEDQUANTITYALONGLROUGHSTEP2}
\sup_{\substack{\timefunction \in [\timefunction_0,\timefunctionboot]
	\\
	\timefunction' \in [\timefunction_0,\timefunction]
	\\ 
	u \in [- \rightu,\leftu] \\ (x^2,x^3) \in \mathbb{T}^2}} 
\frac{\upmu^2 \circ \FlowmapLrougharg{\muxmulevelsetvalue}(\timefunction,u,x^2,x^3)}
	{\upmu^2 \circ \FlowmapLrougharg{\muxmulevelsetvalue}(\timefunction',u,x^2,x^3)}  
& 
\leq C.
\end{align}
From \eqref{E:POINTWISEFULLYMODIFIEDQUANTITYALONGLROUGHSTEP2} and definition~\eqref{E:DELICATEPOINTWISEESTIMATEINTEGRATINGFACTOR},
we see that:
\begin{align} \label{E:INTEGRATINGFACTORRECIPROCALBOUNDED}
\sup_{(\timefunction,u,x^2,x^3) \in 
	[\timefunction_0,\timefunctionboot] \times [-\rightu,\leftu] \times \mathbb{T}^2}
\frac{1}{I(\timefunction,u,x^2,x^3)}
& 
\leq C.
\end{align}

We will now derive pointwise estimates for the terms on RHS~\eqref{E:POINTWISEFULLYMODIFIEDQUANTITYALONGLROUGHSTEP1}.
First, we use \eqref{E:INTEGRATINGFACTORRECIPROCALBOUNDED} to deduce that the first term
on RHS~\eqref{E:POINTWISEFULLYMODIFIEDQUANTITYALONGLROUGHSTEP1} is bounded in magnitude
by the first term on RHS~\eqref{E:ESTIMATEFORFULLYMODIFIEDQUANTALONGLROUGH} as desired.

We now bound the term on RHS~\eqref{E:POINTWISEFULLYMODIFIEDQUANTITYALONGLROUGHSTEP1}
generated by the first term on RHS~\eqref{E:TRANSPORTEQUATIONFORFULLYMODIFIEDQUANTITY}
(with $\tanderY^N$ in the role of $\tander^N$)
as follows,
where throughout the rest of the proof, we use that
$
\frac{I(\timefunction',u,x^2,x^3)}{I(\timefunction,u,x^2,x^3)}
=
\frac{\upmu^2 \circ \FlowmapLrougharg{\muxmulevelsetvalue}(\timefunction,u,x^2,x^3)}
{\upmu^2 \circ  \FlowmapLrougharg{\muxmulevelsetvalue}(\timefunction',u,x^2,x^3)} 
$:
\begin{align}
\begin{split} \label{E:TRANSPORTESTIMATEFORFULLYMODIFIEDQUANTSPLITUPOFSHARPCONTROLREGIONS}
&
\left|
- 2 
\int_{\timefunction_0}^{\timefunction} 
	\frac{\upmu^2 \circ \FlowmapLrougharg{\muxmulevelsetvalue}(\timefunction,u,x^2,x^3)}
		{\upmu^2 \circ \FlowmapLrougharg{\muxmulevelsetvalue}(\timefunction',u,x^2,x^3)} 
	\left\lbrace  
		\frac{\argLrough{\muxmulevelsetvalue} \upmu}{\upmu} \tanderY^N \mathfrak{X}
	\right\rbrace 
	\circ  
	\FlowmapLrougharg{\muxmulevelsetvalue}(\timefunction',u,x^2,x^3)
\, \mathrm{d}\timefunction' 
\right|
		\\
& 
\leq 
2
\int_{\timefunction_0}^{\timefunction} 
	\frac{\upmu^2 \circ \FlowmapLrougharg{\muxmulevelsetvalue}(\timefunction,u,x^2,x^3)}
		{\upmu^2 \circ \FlowmapLrougharg{\muxmulevelsetvalue}(\timefunction',u,x^2,x^3)}
	\left\lbrace  
		\left| \frac{\argLrough{\muxmulevelsetvalue} \upmu}{\upmu}\right| 
		\cdot 
		\left| \tanderY^N \mathfrak{X}\right| \cdot 
		\mathbf{1}_{\{\hypthreearg{\timefunction'}
		{[-\rightu,u]}{\muxmulevelsetvalue} \cap \smallneighborhoodofcreasetwoarg{[\timefunction_0,\timefunctionboot]}{\muxmulevelsetvalue} \}}   
	\right\rbrace 
	\circ 
	\FlowmapLrougharg{\muxmulevelsetvalue}(\timefunction',u,x^2,x^3)
\, \mathrm{d}\timefunction' 
	\\
& \ \ 
+  
2 
\int_{\timefunction_0}^{\timefunction} 
	\frac{\upmu^2 \circ \FlowmapLrougharg{\muxmulevelsetvalue}(\timefunction,u,x^2,x^3)}
		{\upmu^2 \circ \FlowmapLrougharg{\muxmulevelsetvalue}(\timefunction',u,x^2,x^3)}
	\left\lbrace  
		\left| \frac{\argLrough{\muxmulevelsetvalue} \upmu}{\upmu}\right| 
		\cdot 
		\left|  \tanderY^N \mathfrak{X}\right| 
		\cdot 
		\mathbf{1}_{\{\hypthreearg{\timefunction'}
		{[-\rightu,u]}{\muxmulevelsetvalue} \setminus \smallneighborhoodofcreasetwoarg{[\timefunction_0,\timefunctionboot]}{\muxmulevelsetvalue} \}} 
	\right\rbrace 
	\circ 
	\FlowmapLrougharg{\muxmulevelsetvalue}(\timefunction',u,x^2,x^3)
\, \mathrm{d}\timefunction'. 
\end{split}
\end{align}
To handle the first integral on RHS~\eqref{E:TRANSPORTESTIMATEFORFULLYMODIFIEDQUANTSPLITUPOFSHARPCONTROLREGIONS}, 
we use \eqref{E:MUMUSTDECREASEININTERESTINGREGION} to bound it by:
\begin{align} \label{E:POINTWISEFULLYMODIFIEDQUANTITYALONGLROUGHSTEP3}
& 
\leq
\boxed{2} 
\int_{\timefunction_0}^{\timefunction} 
	\left\lbrace  
		\left| \frac{\argLrough{\muxmulevelsetvalue} \upmu}{\upmu}\right| 
		\cdot 
		\left| \tanderY^N \mathfrak{X}\right| \cdot 
		\mathbf{1}_{\{\hypthreearg{\timefunction'}
		{[-\rightu,u]}{\muxmulevelsetvalue} \cap \smallneighborhoodofcreasetwoarg{[\timefunction_0,\timefunctionboot]}{\muxmulevelsetvalue} \}}   
	\right\rbrace 
	\circ 
	\FlowmapLrougharg{\muxmulevelsetvalue}(\timefunction',u,x^2,x^3)
\, {d}\timefunction',
\end{align}
which is $\leq \boxed{2} \int_{\timefunction' = \timefunction_0}^{\timefunction}  
	\left\lbrace 
		\left| \frac{\argLrough{\muxmulevelsetvalue} \upmu}{\upmu} \right|
		\cdot 
		\left| \tanderY^N \mathfrak{X} \right| 
		\cdot 
		\mathbf{1}_{\{\hypthreearg{\timefunction'}{[- \rightu,u]}{\muxmulevelsetvalue} 
		\cap \smallneighborhoodofcreasetwoarg{[\timefunction_0,\timefunctionboot]}{\muxmulevelsetvalue} \}} 
	\right\rbrace  
	\circ 
	\FlowmapLrougharg{\muxmulevelsetvalue}(\timefunction',u,x^2,x^3)
\, \mathrm{d} \timefunction'$ 
as desired.
To handle the second integral on RHS~\eqref{E:TRANSPORTESTIMATEFORFULLYMODIFIEDQUANTSPLITUPOFSHARPCONTROLREGIONS}, 
we use the crude bounds $|\argLrough{\muxmulevelsetvalue} \upmu| \lesssim 1$ 
and $|\upmu| \lesssim 1$
(which follow from the bootstrap assumptions), 
as well as \eqref{E:EASYREGIONLOWERBOUNDFORMU}
and \eqref{E:POINTWISEESTIMATETANGENTDERIVATIVESOFMODQUANTINHOM},
to bound it by the last integral on RHS~\eqref{E:ESTIMATEFORFULLYMODIFIEDQUANTALONGLROUGH} as desired.

Next, to handle the term on RHS~\eqref{E:POINTWISEFULLYMODIFIEDQUANTITYALONGLROUGHSTEP1}
generated by the second term on RHS~\eqref{E:TRANSPORTEQUATIONFORFULLYMODIFIEDQUANTITY},
we first use \eqref{E:POINTWISEFULLYMODIFIEDQUANTITYALONGLROUGHSTEP2}
to bound it in magnitude as follows:
\begin{align}  \label{E:POINTWISEFULLYMODIFIEDQUANTITYALONGLROUGHSTEP4}
&
\leq C
\int_{\timefunction' = \timefunction_0}^{\timefunction}  
	\left| 
		\frac{1}{\Lunit \timefunctionarg{\muxmulevelsetvalue}} 
		\times 
		\upmu
		[\Lunit,\tanderY^N] 
		\mytr_{\gtorus}\upchi 
	\right|
	\circ 
	\FlowmapLrougharg{\muxmulevelsetvalue}(\timefunction',u,x^2,x^3)
\, \mathrm{d}\timefunction'.
\end{align}
Next, using the commutator estimate \eqref{E:POINTWISEBOUNDCOMMUTATORSTANGENTIALANDTANGENTIALDERIVATIVESONSCALARFUNCTION},
Lemma~\ref{L:CRUDEPOINTWISEESTIMATESFORTENSORFIELDS},
\eqref{E:CLOSEDVERSIONLUNITROUGHTTIMEFUNCTION},
Prop.\,\ref{P:POINTWISETRANSPORTINEQUALITIESFOREIKFUNCTIONQUANTITIES},
Prop.\,\ref{P:IMPROVEMENTOFAUXILIARYBOOTSTRAP},
and
Cor.\,\ref{C:IMPROVEAUX},
we deduce the following bound:
\begin{align} 
\begin{split} \label{E:POINTWISEFULLYMODIFIEDQUANTITYALONGLROUGHSTEP4JARED}
|\mbox{RHS~\eqref{E:POINTWISEFULLYMODIFIEDQUANTITYALONGLROUGHSTEP4}}|
&
\leq 
C 
\fundbootsmall  
\int_{\timefunction' = \timefunction_0}^{\timefunction} 
	\max_{\angLie_{\tanderY}^N \in \angLie_{\mathfrak{Y}}^{(N)}}
	|\upmu \angLie_{\tanderY}^N \upchi| 
	\circ 
	\FlowmapLrougharg{\muxmulevelsetvalue}(\timefunction',u,x^2,x^3)
\, \mathrm{d} \timefunction' 
\\
& \ \ 
+ 
C 
\int_{\timefunction' = \timefunction_0}^{\timefunction}  
	\left\lbrace 
		| \comdersmall^{[1,N+1];1} \wavearray| 
		+ 
		|\tander^{[1,N]} \controlvars | 
		+ 
		|\tandersmall^{[1,N]}\badcontrolvars| 
	\right\rbrace 
	\circ 
	\FlowmapLrougharg{\muxmulevelsetvalue}(\timefunction',u,x^2,x^3)
\, \mathrm{d} \timefunction'.
\end{split}
\end{align}
The first term on RHS~\eqref{E:POINTWISEFULLYMODIFIEDQUANTITYALONGLROUGHSTEP4JARED}
is bounded by the $C \varepsilon$-multiplied term on RHS~\eqref{E:ESTIMATEFORFULLYMODIFIEDQUANTALONGLROUGH},
while the last term on RHS~\eqref{E:POINTWISEFULLYMODIFIEDQUANTITYALONGLROUGHSTEP4JARED}
is bounded by the last term on RHS~\eqref{E:ESTIMATEFORFULLYMODIFIEDQUANTALONGLROUGH}.

To handle the terms on RHS~\eqref{E:POINTWISEFULLYMODIFIEDQUANTITYALONGLROUGHSTEP1}
generated by the three commutator terms on the second line of RHS~\eqref{E:TRANSPORTEQUATIONFORFULLYMODIFIEDQUANTITY},
we can use the same arguments given above to bound them as follows:
\begin{align}  \label{E:POINTWISEFULLYMODIFIEDQUANTITYALONGLROUGHSTEP5}
	& 
	\leq
	C 
	\int_{\timefunction' = \timefunction_0}^{\timefunction}  
		\left\lbrace 
			| \comdersmall^{[1,N+1];1} \wavearray| 
			+ 
			|\tander^{[1,N]} \controlvars | 
			+ 
			|\tandersmall^{[1,N]}\badcontrolvars| 
		\right\rbrace 
		\circ 
		\FlowmapLrougharg{\muxmulevelsetvalue}(\timefunction',u,x^2,x^3)
	\, \mathrm{d} \timefunction',
\end{align}
which in turn is bounded by the last term on RHS~\eqref{E:ESTIMATEFORFULLYMODIFIEDQUANTALONGLROUGH}
as desired.

To handle the terms on RHS~\eqref{E:POINTWISEFULLYMODIFIEDQUANTITYALONGLROUGHSTEP1}
arising from the term $\tanderY^N(\upmu |\upchi|_{\gtorus}^2)$ on RHS~\eqref{E:TRANSPORTEQUATIONFORFULLYMODIFIEDQUANTITY},
we expand this term using the Leibniz rule for the operators $\angLie_{\Yvf{A}}$.
Then using the same arguments we used in proving \eqref{E:POINTWISEFULLYMODIFIEDQUANTITYALONGLROUGHSTEP4JARED}
(except no commutator estimates are needed), we find that:
\begin{align} \label{E:POINTWISEFULLYMODIFIEDQUANTITYALONGLROUGHSTEP6JARED}
\tanderY^N(\upmu |\upchi|_{\gtorus}^2)  
& 
\lesssim    
\fundbootsmall
\upmu 
\left|\upmu \angLie_{\tanderY}^N \upchi \right|_{\gtorus} 
+ 
\left|\tander^{[1,N]}\controlvars \right| 
+ 
\left|\tandersmall^{[1,N]}\badcontrolvars \right|.
\end{align}
Using \eqref{E:POINTWISEFULLYMODIFIEDQUANTITYALONGLROUGHSTEP6JARED} and 
\eqref{E:CLOSEDVERSIONLUNITROUGHTTIMEFUNCTION}, we see that the
time integral of the terms in \eqref{E:POINTWISEFULLYMODIFIEDQUANTITYALONGLROUGHSTEP1} 
generated by the term $\tanderY^N(\upmu |\upchi|_{\gtorus}^2)$ can be bounded as follows:
\begin{align} 
\begin{split} \label{E:POINTWISEFULLYMODIFIEDQUANTITYALONGLROUGHSTEP6} 
& 
\leq
C  
\fundbootsmall
\int_{\timefunction' = \timefunction_0}^{\timefunction}  
	\left\lbrace 
		|\upmu \angLie_{\tanderY}^N \upchi| 
	\right\rbrace  
	\circ  
	\FlowmapLrougharg{\muxmulevelsetvalue}(\timefunction',u,x^2,x^3)
\, \mathrm{d} \timefunction' 
	\\
& \ \
	+
	C \int_{\timefunction' = \timefunction_0}^{\timefunction}  
			\left\lbrace 
				| \comdersmall^{[1,N+1];1} \wavearray| 
					+ 
				|\tander^{[1,N]} \controlvars| 
					+ 
				|\tandersmall^{[1,N]}\badcontrolvars| 
			\right\rbrace 
			\circ 
			\FlowmapLrougharg{\muxmulevelsetvalue}(\timefunction',u,x^2,x^3)
		\, \mathrm{d} \timefunction',
\end{split}
\end{align}
which in turn is bounded by RHS~\eqref{E:ESTIMATEFORFULLYMODIFIEDQUANTALONGLROUGH}
as desired.

Finally, we handle the terms on RHS~\eqref{E:POINTWISEFULLYMODIFIEDQUANTITYALONGLROUGHSTEP1}
arising from the term $\tanderY^N \mathfrak{A}$ on RHS~\eqref{E:TRANSPORTEQUATIONFORFULLYMODIFIEDQUANTITY}, 
where $\mathfrak{A}$ has the schematic structure \eqref{E:AINHOMRIC}.
Using the estimates of Prop.\,\ref{P:IMPROVEMENTOFAUXILIARYBOOTSTRAP}
and \eqref{E:CLOSEDVERSIONLUNITROUGHTTIMEFUNCTION},
we bound the time integrals of 
the product of $\frac{1}{\Lunit \timefunctionarg{\muxmulevelsetvalue}}$ and
these terms by the sum of the last three integrals on
RHS~\eqref{E:ESTIMATEFORFULLYMODIFIEDQUANTALONGLROUGH} as desired.
We have therefore proved the lemma.
\end{proof}

\subsubsection{The main pointwise estimates for $(\muX \RRiemann) \tanderY^N \mytr_{\gtorus} \upchi$} 
\label{SSS:MAINPOINTWISEESTIMATESFORFULLYMODIFIED}
We are now ready to prove Prop.\,\ref{P:MAINPOINTWISEESTIMATESFORFULLYMODIFIED},
which provides the main pointwise estimate for the product
$(\muX \RRiemann) \tanderY^N \mytr_{\gtorus} \upchi$
in the case $N = \Ntop$.
The proof relies on the pointwise estimates for 
$\fullymodquant{\tanderY^{\Ntop}}$
provided by Lemma~\ref{L:POINTWISEESTIMATESFORFULLYMODIFIEDQUANTITIY}.

\begin{remark}[The role of the notation $C_*$]
\label{R:CSTAR} 
The constants $C_*$ on RHS~\eqref{E:MOSTDELICATEPOINTWISEESTIMATEFORRPLUS} 
have the same properties as the constants $C$ appearing throughout the paper. We 
have used the notation $C_*$ for some of the constants 
in \eqref{E:MOSTDELICATEPOINTWISEESTIMATEFORRPLUS} because this will aid our analysis of the 
coupling of the different wave energies, especially in
in Sect.\,\ref{SSS:PROOFOFAPRIORIL2ESTIMATESWAVEVARIABLES},
when we prove our Gr\"{o}nwall-type estimates.
Similar remarks apply for constants $C_*$ appearing throughout the rest of the paper. 
\end{remark}

\begin{proposition}[The key pointwise estimate for $(\muX \RRiemann) \tanderY^N \mytr_{\gtorus} \upchi$] 
\label{P:MAINPOINTWISEESTIMATESFORFULLYMODIFIED}
	Let $N = \Ntop$, and let $\tanderY^N \in \mathfrak{Y}^{(N)}$,
	where $\mathfrak{Y}^{(N)}$ is the set of order $N$ 
	$\ell_{t,u}$-tangential commutator operators from Sect.\,\ref{SS:STRINGSOFCOMMUTATIONVECTORFIELDS}.
	Let $\wavearraypartial = (\LRiemann,v^2,v^3,\Ent)$ be as in \eqref{E:PARTIALWAVEARRAY}.
	Then under the assumptions and notation of Lemma~\ref{L:POINTWISEESTIMATESFORFULLYMODIFIEDQUANTITIY}
	(see also Remark~\ref{R:IMPLICITFUNCTIONALDEPENDENCE}),
	the following pointwise estimate holds on 
	$\twoargMrough{[\timefunction_0,\timefunctionboot),[- \rightu,\leftu]}{\muxmulevelsetvalue}$,
	where $\tanderY^N$ is the same in all appearances in \eqref{E:MOSTDELICATEPOINTWISEESTIMATEFORRPLUS},
	(except for the one term in which $\max_{\angLie_{\tanderY^N} \in \angLie_{\mathfrak{Y}}^{(N)}}$ is taken):
	\begin{align}
	\begin{split} \label{E:MOSTDELICATEPOINTWISEESTIMATEFORRPLUS} 
	& 
	\left| 
		\frac{1}{\Lunit \timefunctionarg{\muxmulevelsetvalue}} 
		(\muX \RRiemann) 
		\tanderY^N \mytr_{\gtorus} \upchi 
	\right| 
	\circ 
	\FlowmapLrougharg{\muxmulevelsetvalue}(\timefunction,u,x^2,x^3) 
		\\
	 &  \leq 
	\boxed{2} 
	\left| 
		\frac{\argLrough{\muxmulevelsetvalue} \upmu}{\upmu} 
		\mathbf{1}_{\{\hypthreearg{\timefunction}{[-\rightu,u]}{\muxmulevelsetvalue} 
			\cap \smallneighborhoodofcreasetwoarg{[\timefunction_0,\timefunctionboot]}{\muxmulevelsetvalue}\}} 
	\right|  
	\circ \FlowmapLrougharg{\muxmulevelsetvalue}(\timefunction,u,x^2,x^3) 
	\cdot 
	|\muX \tanderY^N \RRiemann| \circ\FlowmapLrougharg{\muxmulevelsetvalue}(\timefunction,u,x^2,x^3)  
		\\
	& \ \ 
		+ 
		\frac{C_*}{|\timefunction|}   
		|\muX \tanderY^N \wavearraypartial| 
		\circ \FlowmapLrougharg{\muxmulevelsetvalue}(\timefunction,u,x^2,x^3)  
		\\
	& \ \ 
	+ 
	\boxed{4} 
	\left|
		\frac{\argLrough{\muxmulevelsetvalue} \upmu}{\upmu} 
		\mathbf{1}_{\{\hypthreearg{\timefunction}{[-\rightu,u]}{\muxmulevelsetvalue} 
		\cap \smallneighborhoodofcreasetwoarg{[\timefunction_0,\timefunctionboot]}{\muxmulevelsetvalue}\}} 
	\right| 
	\circ 
	\FlowmapLrougharg{\muxmulevelsetvalue}(\timefunction,u,x^2,x^3) 
		\\
	& \ \ \ \ \ \ \
		\times
		\int_{\timefunction' = \timefunction_0}^{\timefunction} 
		\left\lbrace 
			\left| \frac{\argLrough{\muxmulevelsetvalue} \upmu}{\upmu} 
				\mathbf{1}_{\{\hypthreearg{\timefunction'}{[-\rightu,u]}{\muxmulevelsetvalue}
					\cap \smallneighborhoodofcreasetwoarg{[\timefunction_0,\timefunctionboot]}{\muxmulevelsetvalue}\}} 
			\right|
			\left|
				\muX \tanderY^N \RRiemann 
			\right| 
		\right\rbrace 
		\circ \FlowmapLrougharg{\muxmulevelsetvalue}(\timefunction',u,x^2,x^3) 
	\, \mathrm{d} \timefunction' 
	\\
	& \ \ 
	+ 
	\frac{C_*}{|\timefunction|} 
	\int_{\timefunction' = \timefunction_0}^{\timefunction} 
		\frac{1}{|\timefunction'|}  
		|\muX \tanderY^N \wavearraypartial| \circ \FlowmapLrougharg{\muxmulevelsetvalue}(\timefunction',u,x^2,x^3) 
	\, \mathrm{d} \timefunction'  
		\\
	& \ \ 
	+ 
	\frac{C \fundbootsmall}{|\timefunction|} 
	\int_{\timefunction' = \timefunction_0}^{\timefunction}  
		\upmu
		\left|
			\max_{\tanderY^N \in \mathfrak{Y}^{(N)}}
				\angLie_{\tanderY}^N \upchi
		\right| 
		\circ 
		\FlowmapLrougharg{\muxmulevelsetvalue}(\timefunction',u,x^2,x^3)
	\, \mathrm{d} \timefunction' 
		\\
& \ \ 
	+   
	\frac{C}{|\timefunction|} 
	\int_{\timefunction' = \timefunction_0}^{\timefunction} 
		\left\lbrace 
			\upmu |\tanderY^N(\VortVort,\DivGradEnt,\Omega,\GradEnt)| 
			+ 
			|\tanderY^{\leq N-1}(\VortVort,\DivGradEnt,\Omega,\GradEnt)| 		
		\right\rbrace 
		\circ \FlowmapLrougharg{\muxmulevelsetvalue}(\timefunction',u,x^2,x^3)
	\, \mathrm{d} \timefunction' 
		\\
	& \ \ 
		+ 
		\Error \circ \FlowmapLrougharg{\muxmulevelsetvalue}(\timefunction,u,x^2,x^3), 
\end{split}
\end{align}
and the error term $\Error \circ \FlowmapLrougharg{\muxmulevelsetvalue}$ satisfies the following bound:
\begin{align}
\begin{split}  \label{E:ERRORTERMSINMOSTDELICATEPOINTWISEESTIMATEFORRPLUS}
\left| \Error \right|  
\circ 
\FlowmapLrougharg{\muxmulevelsetvalue}(\timefunction,u,x^2,x^3) 
& 
\lesssim 
\frac{1}{|\timefunction|} 
\left| 
	\fullymodquant{\tanderY^N}  
\right|(\timefunction_0,u,x^2,x^3)
			\\
& \ \ 
	+ 
	\frac{\fundbootsmall}{|\timefunction|} 
	\left| \muX \tanderY^N \wavearray \right| 
	\circ 
	\FlowmapLrougharg{\muxmulevelsetvalue}(\timefunction,u,x^2,x^3) 
	+ 
	\left| \comdersmall^{[1,N+1];1} \wavearray \right| 
	\circ 
	\FlowmapLrougharg{\muxmulevelsetvalue}(\timefunction,u,x^2,x^3) 
		\\
	& \ \ 
	+ 
	\frac{1}{|\timefunction|} 
	\left|\comdersmall^{[1,N];1}\wavearray \right| 
	\circ 
	\FlowmapLrougharg{\muxmulevelsetvalue}(\timefunction,u,x^2,x^3) 
	+ 
	\frac{1}{|\timefunction|} \left|
	\begin{pmatrix}
		\tander^{[1,N]}\controlvars  
		\\
	\tandersmall^{[1,N]} \badcontrolvars \end{pmatrix} \right| 
	\circ 
	\FlowmapLrougharg{\muxmulevelsetvalue}(\timefunction,u,x^2,x^3) 	
		\\
	& \ \ 
	+ 
	\frac{\fundbootsmall }{|\timefunction|} 
	\int_{\timefunction' = \timefunction_0}^{\timefunction} 
		\frac{1}{|\timefunction'|} 
		\left|  \comdersmall^{[1,N+1];1} \wavearray \right|  
		\circ 
		\FlowmapLrougharg{\muxmulevelsetvalue}(\timefunction',u,x^2,x^3)
	\, \mathrm{d} \timefunction' 
		\\
	& \ \ 
	+ 
	\frac{1}{|\timefunction|} 
	\int_{\timefunction' = \timefunction_0}^{\timefunction} 
		\left|  \comdersmall^{[1,N+1];1} \wavearray \right|  
		\circ 
		\FlowmapLrougharg{\muxmulevelsetvalue}(\timefunction',u,x^2,x^3)
	\, \mathrm{d} \timefunction' 
		\\
	& \ \ 
		+ 
		\frac{1}{|\timefunction|} 
		\int_{\timefunction' = \timefunction_0}^{\timefunction} 
			\frac{1}{|\timefunction'|} 
			\left\lbrace 
				\left| \comdersmall^{[1,N];1} \wavearray \right| 
				+ 
				\left| 
					\begin{pmatrix}
						\tander^{[1,N]}\controlvars  
						\\
						\tandersmall^{[1,N]} \badcontrolvars 
					\end{pmatrix} 
				\right|
			\right\rbrace  
			\circ 
			\FlowmapLrougharg{\muxmulevelsetvalue}(\timefunction',u,x^2,x^3)
\, \mathrm{d} \timefunction'.
\end{split}
\end{align}

Moreover, the following less precise pointwise estimate 
holds on $\twoargMrough{[\timefunction_0,\timefunctionboot),[- \rightu,\leftu]}{\muxmulevelsetvalue}$:
\begin{align}
\begin{split} \label{E:LESSPRECISEESTIMATETRACECHI}
\left|\upmu \tanderY^N \mytr_{\gtorus}\upchi \right| \circ \FlowmapLrougharg{\muxmulevelsetvalue}(\timefunction,u,x^2,x^3) 
& 
\lesssim 
\left| \fullymodquant{\tanderY^N}  \right|(\timefunction_0,u,x^2,x^3) 
	+
	\left\lbrace  \upmu \left| \tander^{[1,N+1]} \wavearray \right| 
	+ 
	\left| \muX \tander^{[1,N]} \wavearray \right| \right\rbrace \circ \FlowmapLrougharg{\muxmulevelsetvalue}(\timefunction,u,x^2,x^3) 
	\\
	& 
	\ \ 
	+ 
	\left\lbrace \left| \comdersmall^{[1,N];1} \wavearray \right| 
	+ 
	\left| \begin{pmatrix}
\tander^{[1,N]}\controlvars  
	\\
\tandersmall^{[1,N]} \badcontrolvars \end{pmatrix} \right| \right\rbrace \circ \FlowmapLrougharg{\muxmulevelsetvalue}(\timefunction,u,x^2,x^3)  
	\\
	& \ \ + \int_{\timefunction' = \timefunction_0}^{\timefunction} \frac{1}{|\timefunction'|} \left |\muX \tander^N \wavearray \right| \circ  \FlowmapLrougharg{\muxmulevelsetvalue}(\timefunction',u,x^2,x^3)\, \mathrm{d} \timefunction' 
		\\
	& \ \ +  \int_{\timefunction' =\timefunction_0}^{\timefunction} \left|  \comdersmall^{ N+1;1} \wavearray \right|  \circ \FlowmapLrougharg{\muxmulevelsetvalue}(\timefunction',u,x^2,x^3)\, \mathrm{d} \timefunction' 
		\\
	& \ \ 
	+ 
	\int_{\timefunction' = \timefunction_0}^{\timefunction} \frac{1}{|\timefunction'|} \left\lbrace \left|  \comdersmall^{[1,N];1} \wavearray \right| + \left| \begin{pmatrix}
\tander^{[1,N]}\controlvars 
	\\
\tandersmall^{[1,N]} \badcontrolvars \end{pmatrix} \right| \right\rbrace  \circ \FlowmapLrougharg{\muxmulevelsetvalue}(\timefunction',u,x^2,x^3)\, \mathrm{d} \timefunction' 
	\\
	& \ \ + \fundbootsmall \int_{\timefunction' = \timefunction_0}^{\timefunction}  
	\left|
		\upmu 
		\max_{\tanderY^N \in \mathfrak{Y}^{(N)}}
			\angLie_{\tanderY}^N \upchi
	\right| 
	\circ 
	\FlowmapLrougharg{\muxmulevelsetvalue}(\timefunction',u,x^2,x^3)\, \mathrm{d} \timefunction' 
		\\
	& \ \ 
	+  
	\int_{\timefunction' = \timefunction_0}^{\timefunction} \left\lbrace \upmu |\tanderY^N(\VortVort,\DivGradEnt,\Omega,\GradEnt)| + |\tanderY^{\leq N-1}(\VortVort,\DivGradEnt,\Omega,\GradEnt)| \right\rbrace \circ \FlowmapLrougharg{\muxmulevelsetvalue}(\timefunction',u,x^2,x^3)\, \mathrm{d} \timefunction'. 
\end{split}
\end{align}

\end{proposition}

\begin{proof}
Throughout this proof, 
$\Error$ denotes any term such that
$
\left| \Error \right|  
\circ 
\FlowmapLrougharg{\muxmulevelsetvalue}(\timefunction,u,x^2,x^3) 
$
satisfies \eqref{E:ERRORTERMSINMOSTDELICATEPOINTWISEESTIMATEFORRPLUS}.
We begin by using the definition \eqref{E:FULLYMODIFIEDQUANTITY} of $\fullymodquant{\tanderY^N}$
and the estimates \eqref{E:CLOSEDVERSIONLUNITROUGHTTIMEFUNCTION},
\eqref{E:LINFINITYIMPROVEMENTAUXTRANSVERSALPDERIVATIVESRRIEMANNLARGE},
\eqref{E:MINVALUEOFMUONFOLIATION},
and \eqref{E:POINTWISESUMOFMODIFIEDQUANTITYINHOMANDGLL} to deduce:
\begin{align} \label{E:PRECISEPOINTWISEESTIMATESTEP1}
\frac{1}{\Lunit \timefunctionarg{\muxmulevelsetvalue}}
(\muX \RRiemann) 
\tanderY^N \mytr_{\gtorus} \upchi 
& 
= 
\left(\frac{\muX\RRiemann}{\upmu \Lunit \timefunctionarg{\muxmulevelsetvalue}} \right)
\fullymodquant{\tanderY^N} 
+ 
\left(\frac{\muX \RRiemann}{\upmu \Lunit \timefunctionarg{\muxmulevelsetvalue}} \right)
\vec{G}_{\Lunit \Lunit}
\diamond 
\muX \tanderY^N \wavearray 
+ 
\Error,
\end{align}
where in \eqref{E:PRECISEPOINTWISEESTIMATESTEP1}, 
we view all terms as being evaluated at $\FlowmapLrougharg{\muxmulevelsetvalue}(\timefunction,u,x^2,x^3)$. 
Next, using the transport equation \eqref{E:MUTRANSPORT},
\eqref{E:LROUGH}, 
the fact that
$\vec{G}_{\Lunit \Lunit} 
\diamond 
\muX \wavearray 
\eqdef 
G_{\Lunit \Lunit}^0 
\muX \RRiemann 
+ 
G_{\Lunit \Lunit}^1 
\muX \LRiemann 
+ 
G_{\Lunit \Lunit}^2 \muX v^2 
+ 
G_{\Lunit \Lunit}^3 \muX v^3
+
G_{\Lunit \Lunit}^4 \muX \Ent$, 
and the identity 
$
1
=
\mathbf{1}_{\{\hypthreearg{\timefunction}{[-\rightu,u]}{\muxmulevelsetvalue} \cap \smallneighborhoodofcreasetwoarg{[\timefunction_0,\timefunctionboot]}{\muxmulevelsetvalue}\}}
+
\mathbf{1}_{\{\hypthreearg{\timefunction}{[-\rightu,u]}{\muxmulevelsetvalue} \setminus \smallneighborhoodofcreasetwoarg{[\timefunction_0,\timefunctionboot]}{\muxmulevelsetvalue}\}} 
$,
we deduce the following identity for the second product on RHS~\eqref{E:PRECISEPOINTWISEESTIMATESTEP1}: 
\begin{align} 
\begin{split} \label{E:PRECISEPOINTWISEESTIMATESTEP2}
&
\left(
\frac{\muX \RRiemann}{\upmu \Lunit \timefunctionarg{\muxmulevelsetvalue}}  
\right)
\vec{G}_{\Lunit \Lunit}
\diamond 
\muX \tanderY^N \wavearray 
	\\
& 
= 
2 \frac{\argLrough{\muxmulevelsetvalue} \upmu}{\upmu} 
\mathbf{1}_{\{\hypthreearg{\timefunction}{[-\rightu,u]}{\muxmulevelsetvalue} \cap \smallneighborhoodofcreasetwoarg{[\timefunction_0,\timefunctionboot]}{\muxmulevelsetvalue}\}} 
\muX \tanderY^N \RRiemann 
	\\
& \ \
+  
2 \frac{\argLrough{\muxmulevelsetvalue} \upmu}{\upmu} 
\mathbf{1}_{\{\hypthreearg{\timefunction}{[-\rightu,u]}{\muxmulevelsetvalue} 
	\setminus \smallneighborhoodofcreasetwoarg{[\timefunction_0,\timefunctionboot]}{\muxmulevelsetvalue}\}} 		
\muX \tanderY^N \RRiemann 
	\\
& \ \  
-  
\frac{1}{\upmu \Lunit \timefunctionarg{\muxmulevelsetvalue}} 
G_{\Lunit \Lunit}^1 
(\muX \LRiemann)
\muX \tanderY^N \RRiemann 
- 
\sum_{A = 2}^3 
\frac{1}{\upmu \Lunit \timefunctionarg{\muxmulevelsetvalue}} 
G_{\Lunit \Lunit}^A 
(\muX v^A)
\muX \tanderY^N \RRiemann 
-  
\frac{1}{\upmu \Lunit \timefunctionarg{\muxmulevelsetvalue}} 
G_{\Lunit \Lunit}^4 
(\muX \Ent)
\muX \tanderY^N \RRiemann 
	\\
& \ \ 
+ 
\left( 
	\vec{G}_{\Lunit \Lunit} 
	\diamond 
	\argLrough{\muxmulevelsetvalue} \wavearray 
\right) 
\muX \tanderY^N \RRiemann 
+ 
2
\left( 
	\vec{G}_{\Lunit X} \diamond \argLrough{\muxmulevelsetvalue} \wavearray 
\right) 
\muX \tanderY^N \RRiemann  
	\\
& \ \ 
+ 
\left(\frac{\muX \RRiemann}{\upmu \Lunit \timefunctionarg{\muxmulevelsetvalue}} \right)
G_{\Lunit \Lunit}^1 
\muX \tanderY^N \LRiemann 
+   
\left(\frac{\muX \RRiemann}{\upmu \Lunit \timefunctionarg{\muxmulevelsetvalue}} \right)
G_{\Lunit \Lunit}^A 
\muX \tanderY^N v^A
+ 
\left(\frac{\muX \RRiemann}{\upmu \Lunit \timefunctionarg{\muxmulevelsetvalue}} \right)
G_{\Lunit \Lunit}^4 
\muX \tanderY^N \Ent.
\end{split} 
\end{align}
Clearly, the first product
$
2 \frac{\argLrough{\muxmulevelsetvalue} \upmu}{\upmu} 
\mathbf{1}_{\{\hypthreearg{\timefunction}{[-\rightu,u]}{\muxmulevelsetvalue} \cap \smallneighborhoodofcreasetwoarg{[\timefunction_0,\timefunctionboot]}{\muxmulevelsetvalue}\}} 
\muX \tanderY^N \RRiemann 
$ 
on RHS~\eqref{E:PRECISEPOINTWISEESTIMATESTEP2} is bounded in magnitude
by the first product on RHS~\eqref{E:MOSTDELICATEPOINTWISEESTIMATEFORRPLUS} as desired.
Similarly, using 
Lemma~\ref{L:SCHEMATICSTRUCTUREOFVARIOUSTENSORSINTERMSOFCONTROLVARS},
\eqref{E:CLOSEDVERSIONLUNITROUGHTTIMEFUNCTION},
\eqref{E:MINVALUEOFMUONFOLIATION}, 
and the estimates of Prop.\,\ref{P:IMPROVEMENTOFAUXILIARYBOOTSTRAP},
we see that the products on the last line of RHS~\eqref{E:PRECISEPOINTWISEESTIMATESTEP2} 
are bounded by the $C_*$-multiplied product on the second line of RHS~\eqref{E:MOSTDELICATEPOINTWISEESTIMATEFORRPLUS}. 
Also using the estimate \eqref{E:EASYREGIONLOWERBOUNDFORMU},
we see that the product
$
2 \frac{\argLrough{\muxmulevelsetvalue} \upmu}{\upmu} 
\mathbf{1}_{\{\hypthreearg{\timefunction}{[-\rightu,u]}{\muxmulevelsetvalue} 
\setminus \smallneighborhoodofcreasetwoarg{[\timefunction_0,\timefunctionboot]}{\muxmulevelsetvalue}\}} 		
\muX \tanderY^N \RRiemann 
$
on RHS~\eqref{E:PRECISEPOINTWISEESTIMATESTEP2} is bounded by the second term 
$
\left| \comdersmall^{[1,N+1];1} \wavearray \right| 
\circ 
\FlowmapLrougharg{\muxmulevelsetvalue}(\timefunction,u,x^2,x^3) 
$
on the second line of 
RHS~\eqref{E:ERRORTERMSINMOSTDELICATEPOINTWISEESTIMATEFORRPLUS}. 
Finally, using
Lemma~\ref{L:SCHEMATICSTRUCTUREOFVARIOUSTENSORSINTERMSOFCONTROLVARS},
\eqref{E:CLOSEDVERSIONLUNITROUGHTTIMEFUNCTION},
\eqref{E:MINVALUEOFMUONFOLIATION}, 
and the estimates of Prop.\,\ref{P:IMPROVEMENTOFAUXILIARYBOOTSTRAP},
we see that the remaining products on RHS~\eqref{E:PRECISEPOINTWISEESTIMATESTEP2} 
are bounded in magnitude by 
$
\lesssim \frac{\fundbootsmall}{|\timefunction|} 
\left| \muX \tanderY^N \wavearray \right|
$
(and thus are of type $\Error$),
where we have used that all these remaining products gain an overall smallness factor of 
$\fundbootsmall$ from the factors 
$\muX \LRiemann$, $\muX v^A$, $\muX \Ent$, and $\argLrough{\muxmulevelsetvalue} \wavearray$.

We now bound the first product 
$
\left(\frac{\muX\RRiemann}{\upmu \Lunit \timefunctionarg{\muxmulevelsetvalue}} \right)
\fullymodquant{\tanderY^N}
$  
on RHS~\eqref{E:PRECISEPOINTWISEESTIMATESTEP1}. 
We start by multiplying the inequality \eqref{E:ESTIMATEFORFULLYMODIFIEDQUANTALONGLROUGH} by 
$
\left(
\frac{\muX\RRiemann}{\upmu \Lunit \timefunctionarg{\muxmulevelsetvalue}}
\right)
\circ 
\FlowmapLrougharg{\muxmulevelsetvalue}(\timefunction',u,x^2,x^3)
$. 
To bound the product corresponding to the term
$\boxed{2} \cdots$ on RHS~\eqref{E:ESTIMATEFORFULLYMODIFIEDQUANTALONGLROUGH}, 
we first use Lemma~\ref{L:SCHEMATICSTRUCTUREOFVARIOUSTENSORSINTERMSOFCONTROLVARS},
\eqref{E:CLOSEDVERSIONLUNITROUGHTTIMEFUNCTION},
\eqref{E:MINVALUEOFMUONFOLIATION},
the estimates of Prop.\,\ref{P:IMPROVEMENTOFAUXILIARYBOOTSTRAP},
and \eqref{E:POINTWISESUMOFMODIFIEDQUANTITYINHOMANDGLL} 
to express the product as follows:
\begin{align} 
\begin{split} \label{E:PRECISEPOINTWISEESTIMATESTEP3JARED}
&
\boxed{2} 
\left| \frac{\muX \RRiemann}{\upmu \Lunit \timefunctionarg{\muxmulevelsetvalue}} \right| 
\circ 
\FlowmapLrougharg{\muxmulevelsetvalue}(\timefunction,u,x^2,x^3) 
\int_{\timefunction' = \timefunction_0}^{\timefunction}  
	\left\lbrace 
		\left| 
			\frac{\argLrough{\muxmulevelsetvalue} \upmu}{\upmu} \right| \cdot \left| 
			\vec{G}_{\Lunit \Lunit} 
			\diamond 
			\muX \tanderY^N \wavearray  
		\right| 
		\cdot  
		\mathbf{1}_{\{\hypthreearg{\timefunction'}{[- \rightu,u]}{\muxmulevelsetvalue} 
			\cap \smallneighborhoodofcreasetwoarg{[\timefunction_0,\timefunctionboot]}{\muxmulevelsetvalue} \}} 
	\right\rbrace  
	\circ 
	\FlowmapLrougharg{\muxmulevelsetvalue}(\timefunction',u,x^2,x^3)
\, \mathrm{d} \timefunction'  
	\\
& \ \
+ 
\Error.
\end{split}
\end{align}
We now decompose the second factor in the integrand in \eqref{E:PRECISEPOINTWISEESTIMATESTEP3JARED} as follows:
\begin{align} \label{E:PRECISEPOINTWISEESTIMATESTEP4} 
\vec{G}_{\Lunit \Lunit} 
\diamond \muX 
\tanderY^N \wavearray 
& 
= 
G_{\Lunit \Lunit}^0 \muX \tanderY^N \RRiemann 
+ 
G_{\Lunit \Lunit}^1 \muX \tanderY^N \LRiemann 
+ 
\sum_{A=2}^3 G_{\Lunit \Lunit}^A \muX \tanderY^N v^A
+
G_{\Lunit \Lunit}^4 \muX \tanderY^N \Ent, 
\end{align}
where we view \eqref{E:PRECISEPOINTWISEESTIMATESTEP4} as being evaluated at $\FlowmapLrougharg{\muxmulevelsetvalue}(\timefunction',u,x^2,x^3)$.  
We now use \eqref{E:PRECISEPOINTWISEESTIMATESTEP4} to substitute for the integrand factor in \eqref{E:PRECISEPOINTWISEESTIMATESTEP3JARED}. 
Using 
Lemma~\ref{L:SCHEMATICSTRUCTUREOFVARIOUSTENSORSINTERMSOFCONTROLVARS},
\eqref{E:CLOSEDVERSIONLUNITROUGHTTIMEFUNCTION},
\eqref{E:MINVALUEOFMUONFOLIATION}, 
and the estimates of Prop.\,\ref{P:IMPROVEMENTOFAUXILIARYBOOTSTRAP},
we see that the
time integral corresponding to all the products on RHS~\eqref{E:PRECISEPOINTWISEESTIMATESTEP4}  
except the first one $G_{\Lunit \Lunit}^0 \muX \tanderY^N \RRiemann$ 
$G_{\Lunit \Lunit}^1 \muX \tanderY^N \LRiemann$ in \eqref{E:PRECISEPOINTWISEESTIMATESTEP4} 
are bounded in magnitude by the term
$
\frac{C_*}{|\timefunction|} 
	\int_{\timefunction' = \timefunction_0}^{\timefunction} 
		\frac{1}{|\timefunction'|}  
		|\muX \tanderY^N \wavearraypartial| \circ \FlowmapLrougharg{\muxmulevelsetvalue}(\timefunction',u,x^2,x^3) 
	\, \mathrm{d} \timefunction'  
$
on RHS~\eqref{E:MOSTDELICATEPOINTWISEESTIMATEFORRPLUS}. 
We now handle the one the remaining integral, which is:
\begin{align} \label{E:PRECISEPOINTWISEESTIMATESTEP5}
2 \left| \frac{ \muX \RRiemann}{\upmu \Lunit \timefunctionarg{\muxmulevelsetvalue}} \right| \circ \FlowmapLrougharg{\muxmulevelsetvalue}(\timefunction,u,x^2,x^3) \int_{\timefunction' = \timefunction_0}^{\timefunction}  
	\left\lbrace 
		\left| \frac{\argLrough{\muxmulevelsetvalue} \upmu}{\upmu} \right| 
		\cdot 
		\left| 
			G^0_{\Lunit \Lunit} \muX \tanderY^N \RRiemann 
		\right| 
		\cdot  
		\mathbf{1}_{\{\hypthreearg{\timefunction'}{[- \rightu,u]}{\muxmulevelsetvalue} \cap \mathcal{N}_{
			[\timefunction_0,\timefunctionboot]}^{(\muxmulevelsetvalue)} \}} 	
	\right\rbrace  
	\circ \FlowmapLrougharg{\muxmulevelsetvalue}(\timefunction',u,x^2,x^3)
\, \mathrm{d} \timefunction'.
\end{align}
To proceed, we first use \eqref{E:CRUCIALGLL0TRANSPORTESTIMATE} and Cor.\,\ref{C:IMPROVEAUX} to
obtain the following relation, valid for $\timefunction' \in [\timefunctionboot,\timefunction]$:
\begin{align}  
\begin{split} \label{E:PRECISEPOINTWISEESTIMATESTEP6JARED}
&
\left( G^0_{\Lunit \Lunit} \circ \FlowmapLrougharg{\muxmulevelsetvalue}(\timefunction',u,x^2,x^3)\right) 
\left( \muX \tanderY^N \RRiemann \circ \FlowmapLrougharg{\muxmulevelsetvalue}(\timefunction',u,x^2,x^3)\right) 
	\\
& = 
\left( G^0_{\Lunit \Lunit} \circ \FlowmapLrougharg{\muxmulevelsetvalue}(\timefunction,u,x^2,x^3) \right) 
\left( \muX \tanderY^N \RRiemann \circ \FlowmapLrougharg{\muxmulevelsetvalue}(\timefunction',u,x^2,x^3)\right) 
+ 
\mathcal{O}(\fundbootsmall) \muX \tanderY^N \RRiemann \circ \FlowmapLrougharg{\muxmulevelsetvalue}(\timefunction',u,x^2,x^3),
\end{split}
\end{align}
where we emphasize that the $G_{\Lunit \Lunit}^0$ factor on RHS~\eqref{E:PRECISEPOINTWISEESTIMATESTEP6JARED} 
is evaluated at rough-time $\timefunction$ (as opposed to $\timefunction'$). 
We now substitute \eqref{E:PRECISEPOINTWISEESTIMATESTEP6JARED} into the integral \eqref{E:PRECISEPOINTWISEESTIMATESTEP5}.
The integral corresponding to the last product on RHS~\eqref{E:PRECISEPOINTWISEESTIMATESTEP6JARED}
is bounded by the $C \fundbootsmall$-multiplied time integral on RHS~\eqref{E:ERRORTERMSINMOSTDELICATEPOINTWISEESTIMATEFORRPLUS} 
and thus is of type $\Error$. Next, consider the time integral corresponding to the first product on 
RHS~\eqref{E:PRECISEPOINTWISEESTIMATESTEP6JARED}. Since the factor
$G^0_{\Lunit \Lunit} \circ \FlowmapLrougharg{\muxmulevelsetvalue}(\timefunction,u,x^2,x^3)$
does not depend on the integration variable $\timefunction'$, 
we can pull it out of the integral to obtain the following integral: 
\begin{align} \label{E:PRECISEPOINTWISEESTIMATESTEP6}
2 
\left| 
	\frac{ \muX \RRiemann}{\upmu \Lunit \timefunctionarg{\muxmulevelsetvalue}} G_{\Lunit \Lunit}^0 
\right| 
\circ 
\FlowmapLrougharg{\muxmulevelsetvalue}(\timefunction,u,x^2,x^3) 
\int_{\timefunction' = \timefunction_0}^{\timefunction}  	
	\left\lbrace 
		\left| \frac{\argLrough{\muxmulevelsetvalue} \upmu}{\upmu} \right| 
		\cdot 
		\left|  \muX \tanderY^N \RRiemann \right|  
		\cdot  
		\mathbf{1}_{\{\hypthreearg{\timefunction'}{[-\rightu,u]}{\muxmulevelsetvalue} 
		\cap \smallneighborhoodofcreasetwoarg{[\timefunction_0,\timefunctionboot]}{\muxmulevelsetvalue} \}} 		
	\right\rbrace  
	\circ \FlowmapLrougharg{\muxmulevelsetvalue}(\timefunction',u,x^2,x^3)
\, \mathrm{d} \timefunction'.
\end{align}
Using the transport equation \eqref{E:MUTRANSPORT} satisfied by $\upmu$, 
Lemma~\ref{L:SCHEMATICSTRUCTUREOFVARIOUSTENSORSINTERMSOFCONTROLVARS},
\eqref{E:CLOSEDVERSIONLUNITROUGHTTIMEFUNCTION},
\eqref{E:MINVALUEOFMUONFOLIATION}, 
and the estimates of Prop.\,\ref{P:IMPROVEMENTOFAUXILIARYBOOTSTRAP},
we rewrite the product in \eqref{E:PRECISEPOINTWISEESTIMATESTEP6} that is outside of the integral 
as follows:
\begin{align}
\begin{split} \label{E:PRECISEPOINTWISEESTIMATESTEP7}
\frac{\muX \RRiemann}{\upmu \Lunit \timefunctionarg{\muxmulevelsetvalue}} G_{\Lunit \Lunit}^0 
& 
= 
2 \frac{\argLrough{\muxmulevelsetvalue} \upmu}{\upmu} 
- 
\frac{1}{\upmu \Lunit \timefunctionarg{\muxmulevelsetvalue}} G_{\Lunit \Lunit}^1 \muX \LRiemann 
- 
\sum_{A = 2}^3 \frac{1}{\upmu \Lunit \timefunctionarg{\muxmulevelsetvalue}} G_{\Lunit \Lunit}^A \muX v^A 
- 
\frac{1}{\upmu \Lunit \timefunctionarg{\muxmulevelsetvalue}} G_{\Lunit \Lunit}^4 \muX \Ent
	\\
& \ \
+ 
\vec{G}_{\Lunit \Lunit} \diamond \argLrough{\muxmulevelsetvalue} \wavearray 
+ 
2 \vec{G}_{\Lunit X} \diamond \argLrough{\muxmulevelsetvalue} \wavearray   
	\\
& = 
2 \frac{\argLrough{\muxmulevelsetvalue} \upmu}{\upmu} 
+ 
\mathcal{O}(\fundbootsmall) \frac{1}{|\timefunction|}. 
\end{split}
\end{align}
Substituting \eqref{E:PRECISEPOINTWISEESTIMATESTEP7} for the product outside the integral in \eqref{E:PRECISEPOINTWISEESTIMATESTEP6}, 
and using the bound
$
\left| \frac{\argLrough{\muxmulevelsetvalue} \upmu}{\upmu} \right| 
\circ 
\FlowmapLrougharg{\muxmulevelsetvalue}(\timefunction',u,x^2,x^3)
\lesssim \frac{1}{|\timefunction'|}
$
(which follows from 
\eqref{E:MINVALUEOFMUONFOLIATION}, 
and the estimates of Prop.\,\ref{P:IMPROVEMENTOFAUXILIARYBOOTSTRAP}),
we bound the resulting term as follows:
\begin{align} 
\begin{split} \label{E:PRECISEPOINTWISEESTIMATESTEP8}
&
\leq
4 
\left|\frac{\argLrough{\muxmulevelsetvalue} \upmu}{\upmu} \right| 
\circ 
\FlowmapLrougharg{\muxmulevelsetvalue}(\timefunction,u,x^2,x^3) 
\int_{\timefunction' = \timefunction_0}^{\timefunction}  
	\left\lbrace 
		\left| \frac{\argLrough{\muxmulevelsetvalue} \upmu}{\upmu} \right| 
		\cdot 
		\left| \muX \tanderY^N \RRiemann \right| 
		\cdot  
		\mathbf{1}_{\{\hypthreearg{\timefunction'}{[-\rightu,u]}{\muxmulevelsetvalue}
			\cap \smallneighborhoodofcreasetwoarg{[\timefunction_0,\timefunctionboot]}{\muxmulevelsetvalue} \}} 		
	\right\rbrace 
	\circ 
	\FlowmapLrougharg{\muxmulevelsetvalue}(\timefunction',u,x^2,x^3)
\, \mathrm{d} \timefunction' 
	\\
& \ \
+ 
\Error.
\end{split}
\end{align}
Next, we consider the simple estimate
$\left|\frac{\argLrough{\muxmulevelsetvalue} \upmu}{\upmu} \right| 
\leq 
\left|\frac{\argLrough{\muxmulevelsetvalue} \upmu}{\upmu} \mathbf{1}_{\{\hypthreearg{\timefunction}{[-\rightu,u]}{\muxmulevelsetvalue} \cap \smallneighborhoodofcreasetwoarg{[\timefunction_0,\timefunctionboot]}{\muxmulevelsetvalue}\}} \right|  
+  
\left|\frac{\argLrough{\muxmulevelsetvalue} \upmu}{\upmu} \mathbf{1}_{\{\hypthreearg{\timefunction}{[-\rightu,u]}{\muxmulevelsetvalue} \setminus \smallneighborhoodofcreasetwoarg{[\timefunction_0,\timefunctionboot]}{\muxmulevelsetvalue}\}} \right|$,
which we substitute into the factor outside the integral in \eqref{E:PRECISEPOINTWISEESTIMATESTEP8}.
The term corresponding to 
$
\left|\frac{\argLrough{\muxmulevelsetvalue} \upmu}{\upmu} \mathbf{1}_{\{\hypthreearg{\timefunction}{[-\rightu,u]}{\muxmulevelsetvalue} \cap \smallneighborhoodofcreasetwoarg{[\timefunction_0,\timefunctionboot]}{\muxmulevelsetvalue}\}} \right| 
$ 
is bounded by the term on the third line of RHS~\eqref{E:MOSTDELICATEPOINTWISEESTIMATEFORRPLUS}. 
Moreover, using \eqref{E:EASYREGIONLOWERBOUNDFORMU},
and the estimates of Prop.\,\ref{P:IMPROVEMENTOFAUXILIARYBOOTSTRAP},
we bound the term corresponding to
$
\left|\frac{\argLrough{\muxmulevelsetvalue} \upmu}{\upmu} \mathbf{1}_{\{\hypthreearg{\timefunction}{[-\rightu,u]}{\muxmulevelsetvalue} \setminus \smallneighborhoodofcreasetwoarg{[\timefunction_0,\timefunctionboot]}{\muxmulevelsetvalue}\}} \right|
$
by:
\begin{align} 
\begin{split} \label{E:SUBCRITICALTERMBOUNDINPROOFOFMOSTDELICATEPOINTWISEESTIMATEFORRPLUS}
& \leq
	C 
	\int_{\timefunction' = \timefunction_0}^{\timefunction} 
		\frac{1}{|\timefunction'|} 
		\left| \comdersmall^{[1,N+1];1} \wavearray \right| 
		\circ  
		\FlowmapLrougharg{\muxmulevelsetvalue}(\timefunction',u,x^2,x^3)
	\, \mathrm{d} \timefunction' 
		\\
	&
	\leq 
	C \frac{1}{|\timefunction|} 
	\int_{\timefunction' = \timefunction_0}^{\timefunction} 
		\left| \comdersmall^{[1,N+1];1} \wavearray \right|  
		\circ 
		\FlowmapLrougharg{\muxmulevelsetvalue}(\timefunction',u,x^2,x^3)
	\, \mathrm{d} \timefunction' 
	\leq 
	\Error
\end{split}
\end{align}
as desired.
Finally,
we use 
\eqref{E:CLOSEDVERSIONLUNITROUGHTTIMEFUNCTION},
\eqref{E:MINVALUEOFMUONFOLIATION},
and \eqref{E:LINFINITYIMPROVEMENTAUXTRANSVERSALPDERIVATIVESRRIEMANNLARGE}
to deduce 
$\left|
	\frac{\muX \RRiemann}{\upmu \Lunit \timefunctionarg{\muxmulevelsetvalue}}
	\circ 
\FlowmapLrougharg{\muxmulevelsetvalue}(\timefunction,u,x^2,x^3) 
\right|
\lesssim \frac{1}{|\timefunction|}
$,
and from this pointwise estimate, it is easy to see that the products of 
$\frac{\muX \RRiemann}{\upmu \Lunit \timefunctionarg{\muxmulevelsetvalue}}
\FlowmapLrougharg{\muxmulevelsetvalue}(\timefunction,u,x^2,x^3)
$ 
and the remaining terms on RHS~\eqref{E:ESTIMATEFORFULLYMODIFIEDQUANTALONGLROUGH} are 
bounded in magnitude by $\leq \mbox{RHS~\eqref{E:MOSTDELICATEPOINTWISEESTIMATEFORRPLUS}}$ 
as desired. We have therefore proved \eqref{E:MOSTDELICATEPOINTWISEESTIMATEFORRPLUS}.

The proof of \eqref{E:LESSPRECISEESTIMATETRACECHI} is similar, but much less delicate because it does not rely on careful decompositions like \eqref{E:PRECISEPOINTWISEESTIMATESTEP2} and \eqref{E:PRECISEPOINTWISEESTIMATESTEP7}; we omit the details.
\end{proof}


\subsection{Pointwise estimates for the partially modified quantities} 
\label{SS:POINTWISEESTIMATEFORPARTIALLYMODIFIEDQUANT}
In this section, we derive pointwise estimates for the partially modified quantities
from Def.\,\ref{D:FULLYANDPARTIALLYMODIFIEDQUANTITIES}.

\begin{lemma}[Pointwise estimates for partially modified quantities and their $\Lunit$-derivative] 
\label{L:POINTWISEESTIMATEFORPARTIALLYMODIFIEDQUANT}
Let $N = \Ntop$, and let $\tanderY^{N-1} \in \mathfrak{Y}^{(N-1)}$,
where $\mathfrak{Y}^{(N-1)}$ is the set of order $N-1$ 
$\ell_{t,u}$-tangential commutator operators from Sect.\,\ref{SS:STRINGSOFCOMMUTATIONVECTORFIELDS}. 
Let $\partialmodquant{\tanderY^{N-1}}$ be the corresponding partially modified quantity defined in
\eqref{E:PARTIALMODIFIEDQUANTITY} (with $N-1$ in the role of $N$).
Let $\wavearraypartial = (\LRiemann,v^2,v^3,\Ent)$ be as in \eqref{E:PARTIALWAVEARRAY}.
Let $\argLrough{\muxmulevelsetvalue}$ be the rough null vectorfield defined in \eqref{E:LROUGH},
and let $\FlowmapLrougharg{\muxmulevelsetvalue}$ be the $\timefunction_0$-normalized flow map of
$\argLrough{\muxmulevelsetvalue}$ with respect to the rough adapted coordinates
$(\timefunctionarg{\muxmulevelsetvalue},u,x^2,x^3)$ appearing in Lemma~\ref{L:PROPERTIESOFFLOWMAPOFWIDETILDEL}.
Then there exist constants $C > 0$ and $C_* > 0$ such that the following pointwise estimate holds 
relative to the rough adapted coordinates on 
$\twoargMrough{[\timefunction_0,\timefunctionboot),[- \rightu,\leftu]}{\muxmulevelsetvalue}$:
\begin{subequations}
\begin{align}
\begin{split} \label{E:POINTWISEESTIMATELROUGHPARTIALMODIFIEDQUANT}
	 \left| \argLrough{\muxmulevelsetvalue} \partialmodquant{\tanderY^{N-1}} \right| 
	\circ 
	\FlowmapLrougharg{\muxmulevelsetvalue}(\timefunction,u,x^2,x^3) 
	& 
	\leq 
	\frac{1}{2} 
	\left| \frac{1}{\Lunit \timefunctionarg{\muxmulevelsetvalue}} 
	(G_{\Lunit \Lunit}^0) \angLap \tanderY^{N-1} \RRiemann \right|   
	\circ 
	\FlowmapLrougharg{\muxmulevelsetvalue}(\timefunction,u,x^2,x^3) 
		\\
	& \ \ 
	+ 
	C_* 
	\left| \angLap \tanderY^{N-1} \wavearraypartial \right|  
	\circ 
	\FlowmapLrougharg{\muxmulevelsetvalue}(\timefunction,u,x^2,x^3) 
		\\
	& \ \ 
	+ 
	C 
	\left|\tander^{[1,N]} \controlvars \right|  
	\circ 
	\FlowmapLrougharg{\muxmulevelsetvalue}(\timefunction,u,x^2,x^3).  
\end{split}
\end{align}

Moreover, 
the following pointwise estimate holds 
relative to the rough adapted coordinates on 
$\twoargMrough{[\timefunction_0,\timefunctionboot),[- \rightu,\leftu]}{\muxmulevelsetvalue}$:
\begin{align}
\begin{split} \label{E:POINTWISEESTIMATEPARTIALMODIFIEDQUANTITY}
	 \left| \partialmodquant{\tanderY^{N-1}} \right| 
	\circ 
	\FlowmapLrougharg{\muxmulevelsetvalue}(\timefunction,u,x^2,x^3) 
	& \leq 
	\left| \partialmodquant{\tanderY^{N-1}} \right|(\timefunction_0,u,x^2,x^3) 
		\\
	& \ \ 
	+ 
	\frac{1}{2} 
	\left\lbrace 
		\left|\frac{1}{\Lunit \timefunctionarg{\muxmulevelsetvalue}} G_{\Lunit \Lunit}^0 \right| 
		\circ 
		\FlowmapLrougharg{\muxmulevelsetvalue}(\timefunction,u,x^2,x^3) 
	\right\rbrace 
	\int_{\timefunction' = \timefunction_0}^{\timefunction}  
		\left| \angLap \tanderY^{N-1} \RRiemann\right|  
		\circ 
		\FlowmapLrougharg{\muxmulevelsetvalue}(\timefunction',u,x^2,x^3)
	\, \mathrm{d}\timefunction' 
		\\
	& \ \ 
	+ 
	C_* 
	\int_{\timefunction' = \timefunction_0}^{\timefunction}  
		\left|\angLap \tanderY^{N-1} \wavearraypartial \right|   
		\circ 
		\FlowmapLrougharg{\muxmulevelsetvalue}(\timefunction',u,x^2,x^3)
	\, \mathrm{d} \timefunction' 
		\\
	& \ \ 
	+ 
	C  
	\int_{\timefunction' = \timefunction_0}^{\timefunction} 
		\left\lbrace 
			\fundbootsmall  
			\left| \tander^{[1,N+1]}\wavearray \right| 
			+  
			\left| \tander^{[1,N]} \controlvars  \right| 
		\right\rbrace 
		\circ 
		\FlowmapLrougharg{\muxmulevelsetvalue}(\timefunction',u,x^2,x^3)
	\, \mathrm{d} \timefunction'.
\end{split}
\end{align}
\end{subequations}
	
\end{lemma}

\begin{proof}
To prove \eqref{E:POINTWISEESTIMATELROUGHPARTIALMODIFIEDQUANT}, we first expand
the first term on RHS~\eqref{E:TRANSPORTEQUATIONFORPARTIALMODIFIEDQUANTITY} as follows:
\begin{align} \label{E:POINTWISEESTIMATELROUGHPARTIALMODIFIEDQUANTINTERMEDIATE}
		\frac{1}{2} G_{\Lunit \Lunit}^0 \angLap \tanderY^{N-1} \RRiemann 
		+ 
		\frac{1}{2} G_{\Lunit \Lunit}^1 \angLap \tanderY^{N-1} \LRiemann 
		+ 
		\sum_{A=2,3}
		\frac{1}{2} G_{\Lunit \Lunit}^2 \angLap \tanderY^{N-1} v^A 
		+ 
		\frac{1}{2} G_{\Lunit \Lunit}^4 \angLap \tanderY^{N-1} \Ent.
\end{align}
Multiplying \eqref{E:TRANSPORTEQUATIONFORPARTIALMODIFIEDQUANTITY}
by $\frac{1}{\Lunit \timefunctionarg{\muxmulevelsetvalue}}$ (in view of \eqref{E:LROUGH})
and evaluating at $\FlowmapLrougharg{\muxmulevelsetvalue}(\timefunction,u,x^2,x^3)$,
we see that the first product $\frac{1}{2} G_{\Lunit \Lunit}^0 \angLap \tanderY^{N-1} \RRiemann$
in \eqref{E:POINTWISEESTIMATELROUGHPARTIALMODIFIEDQUANTINTERMEDIATE} 
yields precisely the first term on RHS~\eqref{E:POINTWISEESTIMATELROUGHPARTIALMODIFIEDQUANT}. 
Next, using the bound $\sum_{\iota = 1}^4 |G_{\Lunit \Lunit}^{\iota}| \lesssim 1$ 
(which follows from Lemma~\ref{L:SCHEMATICSTRUCTUREOFVARIOUSTENSORSINTERMSOFCONTROLVARS} and Prop.\,\ref{P:IMPROVEMENTOFAUXILIARYBOOTSTRAP}), 
as well as the estimate $\frac{1}{|\Lunit \timefunctionarg{\muxmulevelsetvalue}|} \approx 1$ (see \eqref{E:CLOSEDVERSIONLUNITROUGHTTIMEFUNCTION}),
we see that the remaining products in \eqref{E:POINTWISEESTIMATELROUGHPARTIALMODIFIEDQUANTINTERMEDIATE} 
are bounded in magnitude by the $C_*$-multiplied term on RHS~\eqref{E:POINTWISEESTIMATELROUGHPARTIALMODIFIEDQUANT}. 
Also using \eqref{E:POINTWISETANGENTDERIVATIVEOFBINHOM}, 
we conclude \eqref{E:POINTWISEESTIMATELROUGHPARTIALMODIFIEDQUANT}.

To prove \eqref{E:POINTWISEESTIMATEPARTIALMODIFIEDQUANTITY},
we start by integrating \eqref{E:POINTWISEESTIMATELROUGHPARTIALMODIFIEDQUANT} in rough time
and using \eqref{E:TRANSPORTIDENTITYALONGLROUGHINTEGRALCURVES}. 
The integrals of all products on RHS~\eqref{E:POINTWISEESTIMATELROUGHPARTIALMODIFIEDQUANT} except the first one
are clearly bounded in magnitude by RHS~\eqref{E:POINTWISEESTIMATEPARTIALMODIFIEDQUANTITY}. 
To handle the remaining the integral of the remaining product
$
\frac{1}{2} 
	\left| \frac{1}{\Lunit \timefunctionarg{\muxmulevelsetvalue}} 
	(G_{\Lunit \Lunit}^0) \angLap \tanderY^{N-1} \RRiemann \right|   
	\circ 
	\FlowmapLrougharg{\muxmulevelsetvalue}
$,
we use \eqref{E:CRUCIALGLL0TRANSPORTESTIMATE} and Cor.\,\ref{C:IMPROVEAUX} 
to replace the integrand factor $|G_{\Lunit \Lunit}^0 \circ \FlowmapLrougharg{\muxmulevelsetvalue}(\timefunction',u,x^2,x^3)|$ 
with $|G_{\Lunit \Lunit}^0 \circ \FlowmapLrougharg{\muxmulevelsetvalue}(\timefunction,u,x^2,x^3)|$
(which we can pull out of the integral, as is indicated in the first product on RHS~\eqref{E:POINTWISEESTIMATEPARTIALMODIFIEDQUANTITY})
factor from the $\timefunction$-integral, 
at the expense of error terms featuring a small $\fundbootsmall$ factor. 
Using the comparison estimates 
\eqref{E:ANGDFPOINTWISEBOUNDEDBYCOMMUTATORVECTORFIELDS}--\eqref{E:SMOOTHANGULARHESSIANOFFPOINTWISEBOUNDEDBYCOMMUTATORVECTORFIELDS},
we bound these error terms by the terms on the last line of RHS~\eqref{E:POINTWISEESTIMATEPARTIALMODIFIEDQUANTITY}
as desired.
\end{proof}

\section{Pointwise estimates for controlling the specific vorticity, entropy gradient, and rough acoustic geometry}
\label{S:POINTWISEESTIMATESFORCONTROLLINGSPECIFICVORTICITYANDENTROPYGRADIENT}
In this section, we derive the pointwise estimates that we will use
in Sects.\,\ref{S:BELOWTOPORDERHYPERBOLICL2ESTIMATESFORSPECIFICVORTICITYANDENTROPYGRADIENT}--\eqref{S:TOPORDERELIPTICHYPERBOLICL2ESTIMATESFORSPECIFICVORTICITYANDENTROPYGRADIENT}, when we derive $L^2$ estimates for 
$\vortrenormalized$ and $\GradEnt$.
A key ingredient, which we will use in various spots throughout the paper,
is pointwise estimates for various geometric quantities that are tied to the rough acoustic geometry;
see Lemma~\ref{L:POINTWISEESTIMATESINVOLVINGROUGHACOUSTICGEOMETRY}.

\subsection{Pointwise estimates for $\vortrenormalized$, $\GradEnt$, $\VortVort$, $\DivGradEnt$, and their derivatives}
\label{SS:POINTWISEESTIMATESFORTRANSPORTVARIABLESANDTHEIRDERIVATIVES}
In this section, we derive
pointwise estimates for $\vortrenormalized$, $\GradEnt$, $\VortVort$, $\DivGradEnt$, and various derivatives
of these quantities. We provide the main estimates in Prop.\,\ref{P:POINTWISESTIAMTESFORALLTHETRANSPORTVARIABLES}.

\subsubsection{A simple identity for $\mathrm{d} \SigmatTan_{\flat}$}
\label{SSS:SIMPLEIDENTITYEXTERIORDERIVATIVEOFSIGMATTANGENTONEFORM}
In our proof of Prop.\,\ref{P:POINTWISESTIAMTESFORALLTHETRANSPORTVARIABLES},
we will use the following lemma, which provides an identity for $\mathrm{d} \SigmatTan_{\flat}$ when $\SigmatTan$
is $\Sigma_t$-tangent.

\begin{lemma}\label{L:MAINDECOMPOSITIONOFANTISYMMETRICPARTOFGRADIENTOFSIGMATTANGENTONEFORM}  
Let $\SigmatTan$ be a $\Sigma_t$-tangent vectorfield, and let $\mathrm{d} \SigmatTan_{\flat}$
be the two-form with the following components:
$(\mathrm{d} \SigmatTan_{\flat})_{\alpha \beta} 
	\eqdef \partial_{\alpha} \SigmatTan_{\beta} 
	- 
	\partial_{\beta} \SigmatTan_{\alpha}
	$.
Then relative to the
Cartesian coordinates, the following identity holds,
where $\Speed = \Speed(\LogDensity,\Ent)$ is
the speed of sound:
	\begin{align} 
	\begin{split} \label{E:MAINDECOMPOSITIONOFANTISYMMETRICPARTOFGRADIENTOFSIGMATTANGENTONEFORM}
		(\mathrm{d} \SigmatTan_{\flat})_{\alpha \beta}
		&
		=
		\partial_{\alpha} \SigmatTan_{\beta}
		-
		\partial_{\beta} \SigmatTan_{\alpha}
			= 
				2 (\partial_{\beta} \ln \Speed) \SigmatTan_{\alpha} 
				- 
				2 (\partial_{\alpha} \ln \Speed) \SigmatTan_{\beta}
				+
				\updelta_{\alpha}^0 \SigmatTan_a \partial_{\beta} v^a
				-
				\updelta_{\beta}^0 \SigmatTan_a \partial_{\alpha} v^a
				\\
				& \ \
				+
				\left\lbrace
				\updelta_{\alpha}^0 
				\gfour_{\beta \gamma}
				-
				\updelta_{\beta}^0 
				\gfour_{\alpha \gamma}
				\right\rbrace
				\Transport \SigmatTan^{\gamma}
			+
			\Speed^{-2}
			\upepsilon_{\alpha \beta \gamma \delta}
			\Transport^{\gamma}
			(\Flatcurl \SigmatTan)^{\delta}.
	\end{split}
	\end{align}
\end{lemma}	

\begin{proof}
The same proof of \cite[Lemma 5.6]{lAjS2020} holds.
\end{proof}

\subsubsection{Simple commutator estimates involving $\vortrenormalized$ and $\GradEnt$}
\label{SSS:SIMPLECOMMUTATORLEMMAINVOLVINGTRANSPORTVARIABLE}
We will use the following simple commutator estimates in our proof of Prop.\,\ref{P:POINTWISESTIAMTESFORALLTHETRANSPORTVARIABLES}.

\begin{lemma}[Commuting $\upmu$-weighted Cartesian derivatives with the geometric vectorfields] 
\label{L:POINTWISEESTIMATESFORMUWEIGHTEDCARTESIANCOMMUTATOR} 
Let $1 \leq N \leq \Ntop$.
Then the following commutator estimates hold on $\twoargMrough{[\timefunction_0,\timefunctionboot),[- \rightu,\leftu]}{\muxmulevelsetvalue}$:
\begin{align}
\begin{split} \label{E:POINTWISEESTIMATESFORMUWEIGHTEDCARTESIANCOMMUTATOR}
	& 
	\left| [\upmu \partial_i, \tander^N] (\vortrenormalized,\GradEnt) \right|, 
		\, 
	\left| [\upmu \Flatcurl, \tander^N] (\vortrenormalized,\GradEnt) \right|, 
		\, 
	\left| [\upmu \Flatdiv, \tander^N] (\vortrenormalized,\GradEnt) \right|  
		\\
	&  
	\lesssim 
	\left| \tander^{\leq N} (\vortrenormalized, \GradEnt) \right| 
	+ 
	\left| \muX \tander^{\leq N-1}(\vortrenormalized,\GradEnt)\right| 
	+  
	\fundbootsmall |\muX \tander^{[1,N-1]} \wavearray| 
	+
	\fundbootsmall \upmu |\tander^{[1,N]} \wavearray| 
	+ 
	\fundbootsmall |\tander^{[2,N]} \badcontrolvars|. 
\end{split}
\end{align}
\end{lemma}

\begin{proof}
To derive the estimates for $|[\upmu \partial_i, \tander^N] (\vortrenormalized,\GradEnt)|$,
we first use Lemma~\ref{L:RELATIONSHIPBETWEENCARTESIANPARTIALDERIVATIVESANDSMOOTHGEOMETRICCOMMUTATORS} to express the $\upmu$-weighted Cartesian partial derivatives in terms of the geometric commutation vectorfields
$\lbrace \Lunit, \muX, \Yvf{2}, \Yvf{3} \rbrace$.
Also using Lemma~\ref{L:SCHEMATICSTRUCTUREOFVARIOUSTENSORSINTERMSOFCONTROLVARS},
the commutator estimates 
\eqref{E:POINTWISEBOUNDCOMMUTATORSTANGENTIALANDTANGENTIALDERIVATIVESONSCALARFUNCTION}--\eqref{E:POINTWISEBOUNDCOMMUTATORSMUXANDTANGENTIALDERIVATIVESONSCALARFUNCTION},
and the estimates of Prop.\,\ref{P:IMPROVEMENTOFAUXILIARYBOOTSTRAP},
we conclude the desired estimate \eqref{E:POINTWISEESTIMATESFORMUWEIGHTEDCARTESIANCOMMUTATOR} for
$|[\upmu \partial_i, \tander^N] (\vortrenormalized,\GradEnt)|$.
The desired estimates for $| [\upmu \Flatcurl, \tander^N] (\vortrenormalized,\GradEnt)|$ 
and 
$| [\upmu \Flatdiv, \tander^N] (\vortrenormalized,\GradEnt) |$ follow immediately from the estimate for 
$| [\upmu \partial_i, \tander^N] (\vortrenormalized,\GradEnt)|$.
\end{proof}

\subsection{The main pointwise estimates}
\label{SSS:POINTWISEESTIMATESFORTRANSPORTVARIABLESANDTHEIRDERIVATIVES}

\begin{proposition}[Pointwise estimates for $\vortrenormalized,\GradEnt,\VortVort,\DivGradEnt$, and their derivatives]
	\label{P:POINTWISESTIAMTESFORALLTHETRANSPORTVARIABLES}
	The following pointwise estimates hold on $\twoargMrough{[\timefunction_0,\timefunctionboot),[- \rightu,\leftu]}{\muxmulevelsetvalue}$:
	
	\medskip
	
	\noindent \underline{\textbf{Transport estimates}}.
	For $0 \leq N \leq \Ntop$, we have:	
	\begin{align}
		|\upmu \Transport \tander^N (\vortrenormalized,\GradEnt)|
		& \lesssim 
			|\tander^{\leq N} (\vortrenormalized,\GradEnt)|
			+
			\fundbootsmall |\muX \tander^{[1,N]} \wavearray|
			+
			\fundbootsmall \upmu |\tander^{N+1} \wavearray|
			+
			\fundbootsmall |\tander^{[2,N]} \badcontrolvars|,
				\label{E:COMMUTEDTRANSPORTPOINTWISEESTIMATESFORSPECIFICVORTICITYANDENTROPYGRADIENT} 
	\end{align}
	
	\begin{subequations}
	\begin{align}
		|\upmu \Transport \tander^N \VortVort|
		& \lesssim 
			|\tander^{\leq N} \VortVort|
			+
			|\tander^{\leq N+1} (\vortrenormalized,\GradEnt)|
			+
			\varepsilon |\muX \tander^{[1,N]} \wavearray|
			+
			\varepsilon |\tander^{N+1} \wavearray|
			+
			\varepsilon |\tander^{[2,N]} \badcontrolvars|,
			\label{E:COMMUTEDTRANSPORTPOINTWISEESTIMATESFORMODIFIEDCURLOFVORT}
				\\
		|\upmu \Transport \tander^N \DivGradEnt|
		& \lesssim 
			|\tander^{\leq N} \DivGradEnt|
			+
			|\tander^{\leq N+1} (\vortrenormalized,\GradEnt)|
			+
			\varepsilon |\muX \tander^{[1,N]} \wavearray|
			+
			\varepsilon |\tander^{N+1} \wavearray|
			+
			\varepsilon |\tander^{[2,N]} \badcontrolvars|.
			\label{E:COMMUTEDTRANSPORTPOINTWISEESTIMATESFORMODIFIEDDIVERGENCEOFENTROPYGRADIENT}
	\end{align}
	\end{subequations}

	\medskip
	
	\noindent \underline{\textbf{Algebraic estimates for transversal derivatives in terms of tangential derivatives}}.
	For $0 \leq N \leq \Ntop$, we have:	
	\begin{subequations}
	\begin{align}
			|\muX \tander^N \vortrenormalized|,
				\,
			|\tander^N \muX \vortrenormalized| 
			& \lesssim 
			\upmu |\Lunit \tander^N \vortrenormalized|
			+
			|\tander^{\leq N} (\vortrenormalized,\GradEnt)|
			+
			\varepsilon |\muX \tander^{[1,N]} \wavearray|
			+
			\varepsilon \upmu |\tander^{N+1} \wavearray|
			+
			\varepsilon |\tander^{[2,N]} \badcontrolvars|,
				\label{E:COMMUTEDTRANSVERSALDERVIATVIESOFVORTICITYINTERMSOFTANGENTIAL} 
					\\
			|\muX \tander^N \GradEnt| 
				\,
			|\tander^N \muX \GradEnt|
			& \lesssim 
			\upmu |\Lunit \tander^N \GradEnt|
			+
			|\tander^{\leq N} (\vortrenormalized,\GradEnt)|
			+
			\varepsilon |\muX \tander^{[1,N]} \wavearray|
			+
			\varepsilon \upmu |\tander^{N+1} \wavearray|
			+
			\varepsilon |\tander^{[2,N]} \badcontrolvars|.
				\label{E:COMMUTEDTRANSVERSALDERVIATVIESOFENTROPYGRADIENTINTERMSOFTANGENTIAL} 
		\end{align}
		\end{subequations}

	\medskip
	
	\noindent \underline{\textbf{Algebraic estimates for $(\Flatdiv \vortrenormalized, \Flatdiv \GradEnt)$ and 
	 $(\Flatcurl \vortrenormalized, \Flatcurl \GradEnt)$ in terms of $(\VortVort,\DivGradEnt)$}}.
	For $0 \leq N \leq \Ntop$, we have:	
	\begin{subequations}
	\begin{align}
	\begin{split} 	\label{E:COMMUTEDPOINTWISEEUCLIDEANCURLOFVORTICITY}
		\left|
			\Flatcurl \tander^N \vortrenormalized
		\right|
		& \lesssim 
			\left|
				\tander^N \VortVort 
			\right|
			+
			\frac{1}{\upmu}
			\left|
				\tander^{\leq N-1} \VortVort 
			\right|
			+
			\frac{1}{\upmu}
			\left|
				\tander^{\leq N} (\vortrenormalized,\GradEnt) 
			\right|
				\\
		& \ \
			+
			\frac{\varepsilon}{\upmu} 
			\left|	
				\muX \tander^{[1,N]} \wavearray 
			\right|
			+
			\varepsilon 
			\left|
				\tander^{N+1} \wavearray 
			\right|
			+
			\frac{\varepsilon}{\upmu}
			\left|
				\tander^{[2,N]} \badcontrolvars 
			\right|,
			\end{split} 
				\\
		\begin{split} 	\label{E:COMMUTEDPOINTWISEEUCLIDEANEUCLIDEANDIVERGENCEOFENTROPYGRADIEENT} 
		\left|
			\Flatdiv \tander^N \GradEnt
		\right|
		& \lesssim	
			\left|
				\tander^N \DivGradEnt 
			\right|
			+
			\frac{1}{\upmu}
			\left|
				\tander^{\leq N-1} \DivGradEnt 
			\right|
			+
			\frac{1}{\upmu}
			\left|
				\tander^{\leq N} (\vortrenormalized,\GradEnt) 
			\right|
				\\
		& \ \
			+
			\frac{\varepsilon}{\upmu}
			\left|	
				\muX \tander^{[1,N]} \wavearray 
			\right|
			+
			\varepsilon 
			\left|
				\tander^{N+1} \wavearray 
			\right|
			+
			\frac{\varepsilon}{\upmu}
			\left|
				\tander^{[2,N]} \badcontrolvars 
			\right|,
		\end{split}
				\\
		\left|
			\Flatdiv \tander^N \vortrenormalized
		\right|,
			\,
		\left|
			\Flatcurl \tander^N \GradEnt
		\right|
		& \lesssim 
			\frac{1}{\upmu}
			|\tander^{\leq N} (\vortrenormalized,\GradEnt)|
			+
			\frac{\varepsilon}{\upmu}
			|\muX \tander^{[1,N]} \wavearray|
			+
			\varepsilon |\tander^{N+1} \wavearray|
			+
			\frac{\varepsilon}{\upmu}
			|\tander^{[2,N]} \badcontrolvars|.
				\label{E:COMMUTEDPOINTWISEEUCLIDEANDIVERGENCEOFVORTICITYANDEUCLIDEANCURLOFENTROPYGRADIENT}
	\end{align}
	\end{subequations}
	
	\medskip
	
	\noindent \underline{\textbf{Estimates for the exterior derivative of $\tander^{\Ntop} \vortrenormalized_{\flat}$ and 
	$\tander^{\Ntop} \GradEnt_{\flat}$}}.
	The following estimates hold:
	\begin{subequations}
	\begin{align}
	\begin{split} \label{E:POINTWISEBOUNDEXTERIORDERIVATIVEOFTOPORDERDERIVATIVESOFVORTICITY} 
		|(\mathrm{d} \tander^{\Ntop} \vortrenormalized)_{\flat}|_{\hfour}
		& \lesssim 
			\left|
				\tander^{\Ntop} \VortVort 
			\right|
			+
			\frac{1}{\upmu}
			\left|
				\tander^{\leq \Ntop-1} \VortVort 
			\right|
			+
			\frac{1}{\upmu}
			\left|
				\tander^{\leq \Ntop} (\vortrenormalized,\GradEnt) 
			\right|
				\\
		& \ \
			+
			\frac{\varepsilon}{\upmu}
			\left|	
				\muX \tander^{[1,\Ntop]} \wavearray 
			\right|
			+
			\varepsilon 
			\left|
				\tander^{\Ntop+1} \wavearray 
			\right|
			+
			\frac{\varepsilon}{\upmu}
			\left|
				\tander^{[2,\Ntop]} \badcontrolvars 
			\right|,
		\end{split}		
					\\
		|(\mathrm{d} \tander^{\Ntop} \GradEnt)_{\flat}|_{\hfour}
		& \lesssim 
			\frac{1}{\upmu}
			\left|
				\tander^{\leq \Ntop} (\vortrenormalized,\GradEnt) 
			\right|
			+
			\frac{\varepsilon}{\upmu}
			\left|	
				\muX \tander^{[1,\Ntop]} \wavearray 
			\right|
			+
			\varepsilon 
			\left|
				\tander^{\Ntop+1} \wavearray 
			\right|
			+
			\frac{\varepsilon}{\upmu}
			\left|
				\tander^{[2,\Ntop]} \badcontrolvars 
			\right|.
				\label{E:POINTWISEBOUNDEXTERIORDERIVATIVEOFTOPORDERDERIVATIVESOFENTROPYGRADIENT} 
			\end{align}
		\end{subequations}
	
\end{proposition}

\begin{proof}
\noindent \textbf{Proof of \eqref{E:COMMUTEDTRANSPORTPOINTWISEESTIMATESFORSPECIFICVORTICITYANDENTROPYGRADIENT}}:
We first prove the estimate
\eqref{E:COMMUTEDTRANSPORTPOINTWISEESTIMATESFORSPECIFICVORTICITYANDENTROPYGRADIENT} for $\vortrenormalized^i$. 
We start by multiplying the transport equation \eqref{E:RENORMALIZEDVORTICTITYTRANSPORTEQUATION} 
by $\upmu$ and commuting with $\tander^N$ to deduce that
$|\upmu \Transport \tander^N \vortrenormalized^i| 
\lesssim 
|\tander^N (\upmu \mathfrak{L}_{(\vortrenormalized)}^i) | 
+ 
|[\upmu \Transport, \tander^N] \vortrenormalized^i|$. 
Next, using Lemma~\ref{L:POINTWISEESTIMATESFORDERIVATIVESOFLINEARINHOMOGENEOUSTERMS},
we see that 
$|\tander^N (\upmu \mathfrak{L}_{(\vortrenormalized)}^i)| 
\lesssim \mbox{RHS~\eqref{E:COMMUTEDTRANSPORTPOINTWISEESTIMATESFORSPECIFICVORTICITYANDENTROPYGRADIENT}}$. 
Next, 
using the relation $\upmu \Transport = \upmu \Lunit + \muX$ (see \eqref{E:BISLPLUSX}), 
we derive the commutator identity 
$[\upmu \Transport, \tander^N] = \upmu[\Lunit,\tander^N] + [\upmu,\tander^N] \Lunit + [\muX,\tander^N]$. Using the Leibniz rule, 
the commutator estimates 
\eqref{E:COMMUTATOROFTANGENTIALANDTANGENTIALCOMMUTATORS}--\eqref{E:COMMUTATOROFMUXANDTANGENTIALCOMMUTATORS},
and the estimates of Prop.\,\ref{P:IMPROVEMENTOFAUXILIARYBOOTSTRAP}, 
we find that 
$|[\upmu \Transport, \tander^N] \vortrenormalized^i| 
\lesssim \mbox{RHS~\eqref{E:COMMUTEDTRANSPORTPOINTWISEESTIMATESFORSPECIFICVORTICITYANDENTROPYGRADIENT}}$. 
Combining these estimates, 
we conclude the desired estimate 
\eqref{E:COMMUTEDTRANSPORTPOINTWISEESTIMATESFORSPECIFICVORTICITYANDENTROPYGRADIENT} for $\vortrenormalized^i$.
Using similar arguments, based on the transport equation
\eqref{E:GRADENTROPYTRANSPORT},
we also conclude 
\eqref{E:COMMUTEDTRANSPORTPOINTWISEESTIMATESFORSPECIFICVORTICITYANDENTROPYGRADIENT} for $\GradEnt^i$.

\medskip
\noindent \textbf{Proof of \eqref{E:COMMUTEDTRANSVERSALDERVIATVIESOFVORTICITYINTERMSOFTANGENTIAL} 
and 
\eqref{E:COMMUTEDTRANSVERSALDERVIATVIESOFENTROPYGRADIENTINTERMSOFTANGENTIAL}}:	
To deduce \eqref{E:COMMUTEDTRANSVERSALDERVIATVIESOFVORTICITYINTERMSOFTANGENTIAL} for 
$\tander^N \muX \vortrenormalized$ 
we differentiate the identity \eqref{E:VORTICITYTRANSVERSALTRANSPORTINTERMSOFTANGENTIAL}
with $\tander^N$
and use the estimates of Prop.\,\ref{P:IMPROVEMENTOFAUXILIARYBOOTSTRAP}
and the commutator estimate 
\eqref{E:COMMUTATOROFTANGENTIALANDTANGENTIALCOMMUTATORS} 
(to commute $\tander^N$ under the factor of $\Lunit$ in the first term on RHS~\eqref{E:VORTICITYTRANSVERSALTRANSPORTINTERMSOFTANGENTIAL}).
From this estimate for $\tander^N \muX \vortrenormalized$,
the commutator estimate \eqref{E:COMMUTATOROFMUXANDTANGENTIALCOMMUTATORS},
and the estimates of Prop.\,\ref{P:IMPROVEMENTOFAUXILIARYBOOTSTRAP},
we also conclude the desired bound \eqref{E:COMMUTEDTRANSVERSALDERVIATVIESOFVORTICITYINTERMSOFTANGENTIAL}
for $\muX \tander^N \vortrenormalized$.
The estimates stated in \eqref{E:COMMUTEDTRANSVERSALDERVIATVIESOFENTROPYGRADIENTINTERMSOFTANGENTIAL} 
follow from a nearly identical argument based
on the identity \eqref{E:ENTROPYGRADIENTTRANSVERSALTRANSPORTINTERMSOFTANGENTIAL}.

\medskip
\noindent \textbf{Proof of \eqref{E:COMMUTEDTRANSPORTPOINTWISEESTIMATESFORMODIFIEDCURLOFVORT}
and 
\eqref{E:COMMUTEDTRANSPORTPOINTWISEESTIMATESFORMODIFIEDDIVERGENCEOFENTROPYGRADIENT}}:	
To prove \eqref{E:COMMUTEDTRANSPORTPOINTWISEESTIMATESFORMODIFIEDCURLOFVORT},
we first multiply the transport equation \eqref{E:EVOLUTIONEQUATIONFLATCURLRENORMALIZEDVORTICITY} 
by $\upmu$ and commute with $\tander^N$ to deduce that
$|\upmu \Transport \tander^N \VortVort^i| 
\lesssim 
\left|\tander^N \left\lbrace\upmu \mainnullform_{(\VortVort)}^i 
+ 
\upmu \nullform_{(\VortVort)}^i 
+ 
\upmu \mathfrak{L}_{(\VortVort)}^i \right\rbrace
\right| 
+ 
|[\upmu \Transport, \tander^N]\VortVort|$. 
Next, using
Lemmas~\ref{L:POINTWISEESTIMATESFORDERIVATIVESOFNULLFORMS}
and \ref{L:POINTWISEESTIMATESFORDERIVATIVESOFLINEARINHOMOGENEOUSTERMS},
we see that
$
\left|\tander^N \left\lbrace\upmu \mainnullform_{(\VortVort)}^i 
+ 
\upmu \nullform_{(\VortVort)}^i 
+ 
\upmu \mathfrak{L}_{(\VortVort)}^i \right\rbrace
\right|  
\lesssim 
\mbox{RHS~\eqref{E:COMMUTEDTRANSPORTPOINTWISEESTIMATESFORMODIFIEDCURLOFVORT}}
$.
Next, using 
\eqref{E:POINTWISEESTIMATESFORALLLINEARTERMS},
we find that 
$|\tander^N (\upmu \mathfrak{L}_{(\VortVort)}^i)| 
\lesssim \mbox{RHS~\eqref{E:COMMUTEDTRANSPORTPOINTWISEESTIMATESFORMODIFIEDCURLOFVORT}}$.
To show that
$|[\upmu \Transport, \tander^N] \VortVort| \lesssim \mbox{RHS~\eqref{E:COMMUTEDTRANSPORTPOINTWISEESTIMATESFORMODIFIEDCURLOFVORT}}$,
we can use the same argument that we used in the proof of
\eqref{E:COMMUTEDTRANSPORTPOINTWISEESTIMATESFORSPECIFICVORTICITYANDENTROPYGRADIENT}.
We have therefore proved the estimate \eqref{E:COMMUTEDTRANSPORTPOINTWISEESTIMATESFORMODIFIEDCURLOFVORT}.
The estimate \eqref{E:COMMUTEDTRANSPORTPOINTWISEESTIMATESFORMODIFIEDDIVERGENCEOFENTROPYGRADIENT} 
can be proved to applying nearly identical arguments based on the transport equation \eqref{E:TRANSPORTFLATDIVGRADENT}.

\medskip

\noindent 
\textbf{Proofs of \eqref{E:COMMUTEDPOINTWISEEUCLIDEANCURLOFVORTICITY}--\eqref{E:COMMUTEDPOINTWISEEUCLIDEANDIVERGENCEOFVORTICITYANDEUCLIDEANCURLOFENTROPYGRADIENT}}:	
	To prove \eqref{E:COMMUTEDPOINTWISEEUCLIDEANCURLOFVORTICITY},
	we first multiply \eqref{E:MODIFIEDCURLOFVORTICITY} by $\upmu$ and 
	use Lemmas~\ref{L:RELATIONSHIPBETWEENCARTESIANPARTIALDERIVATIVESANDSMOOTHGEOMETRICCOMMUTATORS}
	and
	\ref{L:SCHEMATICSTRUCTUREOFVARIOUSTENSORSINTERMSOFCONTROLVARS}
	to write the resulting equation in the schematic form
	$	
		\upmu (\Flatcurl \vortrenormalized)^i
		\eqdef
		F
		=
		\upmu \smoothfunction(\wavearray) \VortVort^i
		+
		\smoothfunction(\controlvars) \cdot \GradEnt^a \cdot \muX \wavearray
		+
		\upmu \smoothfunction(\controlvars) \cdot \GradEnt^a \cdot \tander \wavearray
	$.
	Hence, if $0 \leq N \leq \Ntop$, we can commute this equation with 
	$\tander^N$ and then divide by $\upmu$ to deduce
	$
	\Flatcurl \tander^N \vortrenormalized^i
	=
	\frac{1}{\upmu}
	[\upmu \Flatcurl, \tander^N ] \vortrenormalized^i
	+
	\frac{1}{\upmu}
	\tander^N F
	$.
	The bootstrap assumptions and the estimates of Prop.\,\ref{P:IMPROVEMENTOFAUXILIARYBOOTSTRAP}
	imply that 
	$
	\frac{1}{\upmu}
	|\tander^N F|
	\lesssim 
	\left|
		\tander^N \VortVort 
	\right|
			+
			\frac{1}{\upmu}
			\left|
				\tander^{\leq N-1} \VortVort 
			\right|
			+
			\frac{1}{\upmu}
			\left|
				\tander^{\leq N} (\vortrenormalized,\GradEnt) 
			\right|
			+
			\frac{\varepsilon}{\upmu} 
			\left|	
				\muX \tander^{[1,N]} \wavearray 
			\right|
			+
			\varepsilon 
			\left|
				\tander^{N+1} \wavearray 
			\right|
			+
			\frac{\varepsilon}{\upmu}
			\left|
				\tander^{[2,N]} \badcontrolvars 
			\right|
	\lesssim
	\mbox{RHS~\eqref{E:COMMUTEDPOINTWISEEUCLIDEANCURLOFVORTICITY}}
	$
	as desired.
	To show that
	$
	\frac{1}{\upmu}
	\left|
		[\upmu \Flatcurl, \tander^N] \vortrenormalized^i
	\right|
	\lesssim
	\mbox{RHS~\eqref{E:COMMUTEDPOINTWISEEUCLIDEANCURLOFVORTICITY}}
	$, 
	we use \eqref{E:POINTWISEESTIMATESFORMUWEIGHTEDCARTESIANCOMMUTATOR}.
	
	To prove \eqref{E:COMMUTEDPOINTWISEEUCLIDEANEUCLIDEANDIVERGENCEOFENTROPYGRADIEENT},
	we first multiply \eqref{E:MODIFIEDDIVERGENCEOFENTROPYGRADIENT} by $\upmu$ 
	and use Lemmas~\ref{L:RELATIONSHIPBETWEENCARTESIANPARTIALDERIVATIVESANDSMOOTHGEOMETRICCOMMUTATORS}
	and
	\ref{L:SCHEMATICSTRUCTUREOFVARIOUSTENSORSINTERMSOFCONTROLVARS}
	to deduce the schematic equation
	$\upmu \Flatdiv \GradEnt = \upmu f(\wavearray) \DivGradEnt + \smoothfunction(\controlvars) \cdot \GradEnt^a \cdot \muX \wavearray +  \upmu \smoothfunction(\controlvars) \cdot \GradEnt^a \cdot \tander \wavearray$. 
We now argue as in the proof of \eqref{E:COMMUTEDPOINTWISEEUCLIDEANCURLOFVORTICITY},
where we use \eqref{E:POINTWISEESTIMATESFORMUWEIGHTEDCARTESIANCOMMUTATOR} to bound the commutator term
$
\frac{1}{\upmu}
	\left|
		[\upmu \Flatdiv, \tander^N] \GradEnt^i
	\right|
$,
thereby concluding \eqref{E:COMMUTEDPOINTWISEEUCLIDEANEUCLIDEANDIVERGENCEOFENTROPYGRADIEENT}.

We now prove \eqref{E:COMMUTEDPOINTWISEEUCLIDEANDIVERGENCEOFVORTICITYANDEUCLIDEANCURLOFENTROPYGRADIENT} 
for $\Flatdiv \tander^N \vortrenormalized$. 
Multiplying \eqref{E:FLATDIVOFRENORMALIZEDVORTICITY} by $\upmu$, 
commuting with $\tander^N$, 
and then dividing by $\upmu$,
we find that
$
|\Flatdiv \tander^N \vortrenormalized| 
\lesssim
\frac{1}{\upmu}
|\tander^N
(\upmu \mathfrak{L}_{(\Flatdiv \vortrenormalized)})|
+
\frac{1}{\upmu}
|[\upmu \Flatdiv, \tander^N] \vortrenormalized|
$.
The desired estimate now follows 
from Lemma~\ref{L:POINTWISEESTIMATESFORDERIVATIVESOFLINEARINHOMOGENEOUSTERMS} 
and the commutator estimate \eqref{E:POINTWISEESTIMATESFORMUWEIGHTEDCARTESIANCOMMUTATOR}.
The proof of \eqref{E:COMMUTEDPOINTWISEEUCLIDEANDIVERGENCEOFVORTICITYANDEUCLIDEANCURLOFENTROPYGRADIENT} 
for $\Flatcurl \tander^N \GradEnt$
follows from a similar argument based on equation \eqref{E:CURLGRADENTVANISHES}
and the commutator estimate \eqref{E:POINTWISEESTIMATESFORMUWEIGHTEDCARTESIANCOMMUTATOR}.

\medskip
\noindent 
\textbf{Proofs of \eqref{E:POINTWISEBOUNDEXTERIORDERIVATIVEOFTOPORDERDERIVATIVESOFVORTICITY}--\eqref{E:POINTWISEBOUNDEXTERIORDERIVATIVEOFTOPORDERDERIVATIVESOFENTROPYGRADIENT}}:	
To prove \eqref{E:POINTWISEBOUNDEXTERIORDERIVATIVEOFTOPORDERDERIVATIVESOFVORTICITY},
we first use the identity \eqref{E:MAINDECOMPOSITIONOFANTISYMMETRICPARTOFGRADIENTOFSIGMATTANGENTONEFORM}
with $\tander^{\Ntop} \vortrenormalized$ in the role of $\SigmatTan$,
Lemmas~\ref{L:RELATIONSHIPBETWEENCARTESIANPARTIALDERIVATIVESANDSMOOTHGEOMETRICCOMMUTATORS}
and \ref{L:SCHEMATICSTRUCTUREOFVARIOUSTENSORSINTERMSOFCONTROLVARS},
the estimates of Prop.\,\ref{P:IMPROVEMENTOFAUXILIARYBOOTSTRAP},
and \eqref{E:HRIEMANNIANMETRICNORMCOMPARABLETOEUCLIDEANNORM} to deduce that
$
|(\mathrm{d} \tander^{\Ntop} \vortrenormalized)_{\flat}|_{\hfour}
\lesssim
\frac{1}{\upmu}
|\tander^{\Ntop} \vortrenormalized|
+
|\Transport \tander^{\Ntop} \vortrenormalized|
+
|\Flatcurl \tander^{\Ntop} \vortrenormalized|
$.
From this estimate and the pointwise estimates
\eqref{E:COMMUTEDTRANSPORTPOINTWISEESTIMATESFORSPECIFICVORTICITYANDENTROPYGRADIENT}
and
\eqref{E:COMMUTEDPOINTWISEEUCLIDEANCURLOFVORTICITY}
with $N \eqdef \Ntop$, we conclude the desired estimate \eqref{E:POINTWISEBOUNDEXTERIORDERIVATIVEOFTOPORDERDERIVATIVESOFVORTICITY}.
The estimate
\eqref{E:POINTWISEBOUNDEXTERIORDERIVATIVEOFTOPORDERDERIVATIVESOFENTROPYGRADIENT}
follows from a similar argument that relies on the pointwise estimate \eqref{E:COMMUTEDTRANSPORTPOINTWISEESTIMATESFORSPECIFICVORTICITYANDENTROPYGRADIENT}
for $|\Transport \tander^{\Ntop} \GradEnt|$ and the pointwise estimate
\eqref{E:COMMUTEDPOINTWISEEUCLIDEANDIVERGENCEOFVORTICITYANDEUCLIDEANCURLOFENTROPYGRADIENT}
	for
	$
	\left|
			\Flatcurl \tander^{\Ntop} \GradEnt
	\right|
	$.
	
\end{proof}

\subsection{Pointwise estimates tied to the rough acoustic geometry}
\label{SS:POINTWISEESTIMATESINVOLVINGROUGHACOUSTICGEOMETRY}
Recall that \eqref{E:INTEGRALIDENTITYFORELLIPTICHYPERBOLICCURRENT} is
the main elliptic-hyperbolic integral identity that we will use to control
the top-order derivatives of $\vortrenormalized$ and $\GradEnt$.
In Prop.\,\ref{P:POINTWISEESTIMTAESFORELLIPTICHYPERBOLICIDENTITYERORTERMS}, 
we derive pointwise estimates for the error integrands on RHS~\eqref{E:INTEGRALIDENTITYFORELLIPTICHYPERBOLICCURRENT}.
Some of these error integrands
(e.g., the term $\Currentboundaryerrorhavetocontrollowerorder[\SigmatTan,\SigmatTan]$ defined in 
\eqref{E:LOWERORDERERRORTERMHAVETOCONTROLKEYIDPUTANGENTCURRENTCONTRACTEDAGAINSTVECTORFIELD})
depend on geometric quantities that are tied to the rough acoustic geometry. 
In the next lemma, we derive pointwise estimates for these quantities.
This serves as a preliminary step for our proof of Prop.\,\ref{P:POINTWISEESTIMTAESFORELLIPTICHYPERBOLICIDENTITYERORTERMS}.

\begin{lemma}[Pointwise estimates involving the rough acoustic geometry]
\label{L:POINTWISEESTIMATESINVOLVINGROUGHACOUSTICGEOMETRY}
The following pointwise estimates hold on $\twoargMrough{[\timefunction_0,\timefunctionboot),[- \rightu,\leftu]}{\muxmulevelsetvalue}$,
where $\mathbf{1}_{\lbrace \upmu < - \frac{\muxmulevelsetvalue }{\Lunit \upmu} \rbrace}$
denotes the characteristic function of the set
$\lbrace (t,u,x^2,x^3) \ | \ \upmu(t,u,x^2,x^3) <  - \frac{\muxmulevelsetvalue }{\Lunit \upmu(t,u,x^2,x^3)}\rbrace$ 
and similarly for $\mathbf{1}_{\lbrace \upmu \geq - \frac{\muxmulevelsetvalue }{\Lunit \upmu} \rbrace}$,
$\deform{\Rtransarg{\muxmulevelsetvalue}}$ denotes the deformation tensor of $\Rtransarg{\muxmulevelsetvalue}$,
$\mytr_{\gtorusroughfirstfund} \deform{\Rtransarg{\muxmulevelsetvalue}}$ denotes its 
trace with respect to $\gtorusroughfirstfund$,
and $\roughangdiv \Roughtoritangentvectorfieldarg{\muxmulevelsetvalue}$ denotes the
$\twoargroughtori{\timefunctionarg{\muxmulevelsetvalue},u}{\muxmulevelsetvalue}$-divergence 
(see Def.\,\ref{D:CONNECTIONSANDDIFFERENTIALOPERATORSONROUGHTORI})
of the vectorfield $\Roughtoritangentvectorfieldarg{\muxmulevelsetvalue}$:
\label{L:POINTWISEESTIMATESINVOLVINGROUGHACOUSTICGEOMETRY}
\begin{align} 
	\sum_{\alpha=0,1,2,3}
	|\Rtransarg{\muxmulevelsetvalue}^{\alpha}|
	& \lesssim 
			\mathbf{1}_{\lbrace \upmu < - \frac{\muxmulevelsetvalue }{\Lunit \upmu} \rbrace}
			\muxmulevelsetvalue
			+
			\mathbf{1}_{\lbrace \upmu \geq - \frac{\muxmulevelsetvalue }{\Lunit \upmu} \rbrace}
			\upmu,     
			\label{E:RTRANSCARTESIANPOINTWISE}
				\\
	\sum_{\alpha=0,1,2,3}
	|\Roughtoritangentvectorfieldarg{\muxmulevelsetvalue}^{\alpha}|
	& \lesssim \fundbootsmall,  
	\label{E:ROUGHTORITANGENTCARTESIANPOINTWISE}
		\\
	|\tander^{\leq 1} \Roughtoritangentvectorfieldarg{\muxmulevelsetvalue} \upmu|
	& \lesssim \fundbootsmall^2,  
	\label{E:TANGENTIALDERIVATIVESOFROUGHTORIVECTORFIELDAPPLIEDTOMUPOINTWISE}
		\\
	|\tander^{\leq 1} \Rtransnormsmallfactorarg{\muxmulevelsetvalue}|,
		\, 
	|\muX \Rtransnormsmallfactorarg{\muxmulevelsetvalue}|
	& \lesssim \fundbootsmall^2, 
		\label{E:POINTWISEBOUNDTANGENTIALANDTRANSVERSALDERIVATIVESOFRTRANSNORMSMALLFACTOR} 
		\\
	|\Rtransarg{\muxmulevelsetvalue} \Lunit \upmu|
	& \lesssim 1, 
	\label{E:RTRANSAPPLIEDTOLUNITMUPOINTWISE}
		\\
	|\Roughtoritangentvectorfieldarg{\muxmulevelsetvalue} \Lunit \upmu|
	& \lesssim \fundbootsmall^2, 
	\label{E:ROUGHTORIVECTORFIELDAPPLIEDTOLUNITMUPOINTWISE}
		\\
	|\roughangdiv \Roughtoritangentvectorfieldarg{\muxmulevelsetvalue}|
	& \lesssim \fundbootsmall,
	\label{E:DIVERGENCEOFROUGHTORITANGENTVECTORFIELD}
		\\
	|\mytr_{\gtorusroughfirstfund} \deform{\argLrough{\muxmulevelsetvalue}}|
	& \lesssim 1,
	\label{E:POINTWISEBOUNDFORROUGHTOROIDALTRACEOFDEFORMATIONTENSOROFROUGHNULLVECTORFIELD}
		\\
	|\mytr_{\gtorusroughfirstfund} \deform{\Rtransarg{\muxmulevelsetvalue}}|
	& \lesssim 1.
	\label{E:POINTWISEBOUNDFORROUGHTOROIDALTRACEOFDEFORMATIONTENSOROFRTRANS}
\end{align}

\end{lemma}

\begin{proof}
We refer to Sect.\,\ref{SS:SILENTFACTS} for some results that we will silently use throughout the analysis.

\medskip

\noindent \textbf{Proofs of \eqref{E:ROUGHTORITANGENTCARTESIANPOINTWISE}--\eqref{E:POINTWISEBOUNDTANGENTIALANDTRANSVERSALDERIVATIVESOFRTRANSNORMSMALLFACTOR}:}
To prove \eqref{E:POINTWISEBOUNDTANGENTIALANDTRANSVERSALDERIVATIVESOFRTRANSNORMSMALLFACTOR},
we first use 
\eqref{E:IDENTITYFORRTRANSNORMSMALLFACTORSQUARED},
\eqref{E:SMOOTHGINVERSEABEXPRESSION},
Lemma~\ref{L:COMMUTATORSTOCOORDINATES},
and
Lemma~\ref{L:SCHEMATICSTRUCTUREOFVARIOUSTENSORSINTERMSOFCONTROLVARS} 
to deduce that:
\begin{align} \label{E:SCHEMATICEXPRESSIONFORRTRANSSMALL}
\Rtransnormsmallfactorarg{\muxmulevelsetvalue}
&
= 
\frac{1}{(\Lunit \timefunctionarg{\muxmulevelsetvalue})^2}
(\gtorus^{-1})^{AB} (\geop{x^A} \timefunctionarg{\muxmulevelsetvalue}) \geop{x^B} \timefunctionarg{\muxmulevelsetvalue}
= 
\frac{1}{(\Lunit \timefunctionarg{\muxmulevelsetvalue})^2} 
\smoothfunction(\controlvars) 
\cdot 
(\geop{x^2}\timefunctionarg{\muxmulevelsetvalue},\geop{x^3}\timefunctionarg{\muxmulevelsetvalue})
\cdot
(\geop{x^2}\timefunctionarg{\muxmulevelsetvalue},\geop{x^3}\timefunctionarg{\muxmulevelsetvalue}),
\end{align}
where the expression to the right of the last equality is schematic.
From \eqref{E:SCHEMATICEXPRESSIONFORRTRANSSMALL},
\eqref{E:SMALLDERIVATIVESLINFTYESTIMATESFORROUGHTIMEFUNCTIONANDDERIVATIVES},
\eqref{E:LDERIVATIVEOFROUGHTIMEFUNCTIONISAPPROXIMATELYUNITY},
and the bootstrap assumptions, 
and Cor.\,\ref{C:IMPROVEAUX},
we deduce that $|\Rtransnormsmallfactorarg{\muxmulevelsetvalue}| \lesssim \fundbootsmall^2$ as desired.
To prove that $|\tander \Rtransnormsmallfactorarg{\muxmulevelsetvalue}|, \, |\muX \Rtransnormsmallfactorarg{\muxmulevelsetvalue}| \lesssim \fundbootsmall^2$,
we differentiate the expression for
$
\Rtransnormsmallfactorarg{\muxmulevelsetvalue}
$
with elements $\tander \in \lbrace \Lunit, \Yvf{2}, \Yvf{3} \rbrace$ and $\muX$ and apply a similar argument,
where we use Lemma~\ref{L:COMMUTATORSTOCOORDINATES}
to express the vectorfield derivatives in terms of geometric coordinate partial derivatives
when they fall on $\timefunctionarg{\muxmulevelsetvalue}$.

The estimates \eqref{E:ROUGHTORITANGENTCARTESIANPOINTWISE} 
and
\eqref{E:TANGENTIALDERIVATIVESOFROUGHTORIVECTORFIELDAPPLIEDTOMUPOINTWISE}
follow from similar arguments based on the identity
\eqref{E:ROUGHTORITANGENTVECTORFIELDKEEPSAPPEARING}.

\medskip
\noindent \textbf{Proof of \eqref{E:RTRANSCARTESIANPOINTWISE}:}
The estimate \eqref{E:RTRANSCARTESIANPOINTWISE} follows from the decomposition 
$\Rtransarg{\muxmulevelsetvalue}^{\alpha}= \upmu X^{\alpha}+ \frac{\muxmulevelsetvalue \phi}{\Lunit \upmu} \Lunit^{\alpha} - \upmu \Roughtoritangentvectorfieldarg{\muxmulevelsetvalue}^{\alpha}$
(see \eqref{E:RTRANS}),
Lemma~\ref{L:SCHEMATICSTRUCTUREOFVARIOUSTENSORSINTERMSOFCONTROLVARS},
the bootstrap assumptions (in particular \eqref{E:BABOUNDSONLMUINTERESTINGREGION}),
and \eqref{E:ROUGHTORITANGENTCARTESIANPOINTWISE}.

\medskip
\noindent \textbf{Proofs of \eqref{E:RTRANSAPPLIEDTOLUNITMUPOINTWISE}--\eqref{E:ROUGHTORIVECTORFIELDAPPLIEDTOLUNITMUPOINTWISE}:}
These estimates follow from the decompositions 
\eqref{E:RTRANS} and
\eqref{E:ROUGHTORITANGENTVECTORFIELDKEEPSAPPEARING}
and the arguments given in the previous paragraphs.

\medskip
\noindent \textbf{Proof of \eqref{E:DIVERGENCEOFROUGHTORITANGENTVECTORFIELD}:}
First, computing relative to 
the coordinates $(x^2,x^3)$ on the rough tori $\twoargroughtori{\timefunctionarg{\muxmulevelsetvalue},u}{\muxmulevelsetvalue}$,
we deduce that:
\begin{align} \label{E:EXPANSIONOFROUGHDIVERGENCEOFROUGHTORITANGENTVECTORFIELD}
\begin{split}
\roughangdiv \Roughtoritangentvectorfieldarg{\muxmulevelsetvalue} 
& 
= 
\roughgeop{x^A} \Roughtoritangentvectorfieldarg{\muxmulevelsetvalue}^A
	\\
& \ \
+ 
\frac{1}{2} \Roughtoritangentvectorfieldarg{\muxmulevelsetvalue}^A 
(\gtorusroughinversefirstfund)(dx^B,dx^C)
\left\lbrace
	\roughgeop{x^B} \gtorusroughfirstfund\left(\roughgeop{x^A},\roughgeop{x^C}\right) 
	+ 
	\roughgeop{x^C} \gtorusroughfirstfund\left(\roughgeop{x^A},\roughgeop{x^B}\right) 
	- 
	\roughgeop{x^A} \gtorusroughfirstfund\left(\roughgeop{x^B},\roughgeop{x^C}\right) 
\right\rbrace.
\end{split}
\end{align}
Next, we use 
\eqref{E:DEFROUGHTORITANGENTVECTORFIELD} to write 
$\Roughtoritangentvectorfieldarg{\muxmulevelsetvalue}^A  
= (\gtorusroughinversefirstfund)(dx^A,dx^B)
\frac{\geop{x^B}\timefunctionarg{\muxmulevelsetvalue}}{\geop{t} \timefunctionarg{\muxmulevelsetvalue}}$. 
We then use \eqref{E:GEOP2TOCOMMUTATORS}--\eqref{E:GEOP3TOCOMMUTATORS},
\eqref{E:ROUGHANGULARPARTIALDERIVATIVESINTERMSOFGOODGEOMETRICPARTIALDERIVATIVES},
and Lemma~\ref{L:SCHEMATICSTRUCTUREOFVARIOUSTENSORSINTERMSOFCONTROLVARS}
to schematically express 
$\roughgeop{x^A}
= 
\smoothfunction
\left(\controlvars,
	\frac{1}{\geop{t} \timefunctionarg{\muxmulevelsetvalue}}, 
	\tanderY \timefunctionarg{\muxmulevelsetvalue} 
\right)
\Yvf{2}
+
\smoothfunction
\left(\controlvars,
	\frac{1}{\geop{t} \timefunctionarg{\muxmulevelsetvalue}}, 
	\tanderY \timefunctionarg{\muxmulevelsetvalue} 
\right)
\Yvf{3}
+
\smoothfunction
\left(\controlvars,
	\frac{1}{\geop{t} \timefunctionarg{\muxmulevelsetvalue}}, 
\right)
\cdot
(\tanderY \timefunctionarg{\muxmulevelsetvalue})
\Lunit
$,
and we also use
\eqref{E:SMOOTHTORIGABEXPRESSION},
\eqref{E:SMOOTHGINVERSEABEXPRESSION}
\eqref{E:GTORUSROUGHCOMPONENTS},
\eqref{E:CHOVCOEFFICIENTSSMOOTHANGULARDERIVATIVESINTERMSOFROUGHONESANDL},
\eqref{E:SMOOTHTORUSINVERSEFIRSTFUNDCOMPONENTSINTERMSOFROUGHTORUSINVERSEFIRSTFUNDCOMPONENTS},
and Lemma~\ref{L:SCHEMATICSTRUCTUREOFVARIOUSTENSORSINTERMSOFCONTROLVARS}
to deduce the schematic identities
$\gtorusroughfirstfund\left(\roughgeop{x^A},\roughgeop{x^B}\right)
=
\smoothfunction
\left(\controlvars,
	\frac{1}{\geop{t} \timefunctionarg{\muxmulevelsetvalue}}, 
	\tanderY \timefunctionarg{\muxmulevelsetvalue} 
\right)
$
and
$(\gtorusroughinversefirstfund)(dx^A,dx^B) 
= 
\smoothfunction
\left(\controlvars,
	\frac{1}{\Lunit \timefunctionarg{\muxmulevelsetvalue}}, 
	\tanderY \timefunctionarg{\muxmulevelsetvalue} 
\right)
$.
Also using \eqref{E:GEOPTOCOMMUTATORS}, we arrive at the following schematic identity: 
\begin{align} \label{E:SCHEMATICIDENTITYFORROUGHDIVERGENCEOFROUGHTORITANGENTVECTORFIELD}
\roughangdiv \Roughtoritangentvectorfieldarg{\muxmulevelsetvalue} 
& 
= 
\smoothfunction
\left(\tander^{\leq 1} \controlvars,
	\frac{1}{\Lunit \timefunctionarg{\muxmulevelsetvalue}}, 
	\frac{1}{\geop{t} \timefunctionarg{\muxmulevelsetvalue}}, 
	\tander^{\leq 2} \timefunctionarg{\muxmulevelsetvalue} 
\right)
\cdot
\tanderY^{[1,2]} \timefunctionarg{\muxmulevelsetvalue}.
\end{align}
From \eqref{E:SCHEMATICIDENTITYFORROUGHDIVERGENCEOFROUGHTORITANGENTVECTORFIELD},
Lemma~\ref{L:SCHEMATICSTRUCTUREOFVARIOUSTENSORSINTERMSOFCONTROLVARS},
the estimates 
\eqref{E:ALLDERIVATIVESLINFTYESTIMATESFORROUGHTIMEFUNCTIONANDDERIVATIVES}--\eqref{E:SMALLDERIVATIVESLINFTYESTIMATESFORROUGHTIMEFUNCTIONANDDERIVATIVES}
and
\eqref{E:PARTIALTIMEDERIVATIVEOFROUGHTIMEFUNCTIONISAPPROXIMATELYUNITY}--\eqref{E:LDERIVATIVEOFROUGHTIMEFUNCTIONISAPPROXIMATELYUNITY}
for the rough time function, the bootstrap assumptions,
and Cor.\,\ref{C:IMPROVEAUX},
we arrive at the desired estimate
\eqref{E:DIVERGENCEOFROUGHTORITANGENTVECTORFIELD}.

\medskip
\noindent \textbf{Proof of \eqref{E:POINTWISEBOUNDFORROUGHTOROIDALTRACEOFDEFORMATIONTENSOROFRTRANS}:}
First, using \eqref{E:LIEGTORUSROUGH}, we compute that relative to the coordinates $(\timefunctionarg{\muxmulevelsetvalue},u,x^2,x^3)$, we have
$
\mytr_{\gtorusroughfirstfund} \deform{\Rtransarg{\muxmulevelsetvalue}}
= 
\gtorusroughinversefirstfund\left(dx^A,dx^B\right)
\Rtransarg{\muxmulevelsetvalue} \gtorusroughfirstfund \left(\roughgeop{x^A},\roughgeop{x^B}\right) 
+
2 \roughgeop{x^A} \Rtransarg{\muxmulevelsetvalue}^A
$.
Considering this identity,
arguing as in the proof of \eqref{E:SCHEMATICIDENTITYFORROUGHDIVERGENCEOFROUGHTORITANGENTVECTORFIELD},
and also using \eqref{E:MUXINTERMSOFGEOMETRICCOORDINATEVECTORFIELDS} and \eqref{E:RTRANS},
we deduce the following schematic identity:
\begin{align} \label{E:SCHEMATICIDENTITYFORROUGHTORUSTRACEOFDEFORMATIONTENSOROFRTRANS}
\mytr_{\gtorusroughfirstfund} \deform{\Rtransarg{\muxmulevelsetvalue}}
& 
= 
\smoothfunction
\left(\tander^{\leq 1} \badcontrolvars,
	\comder \controlvars,
	\frac{1}{\Lunit \timefunctionarg{\muxmulevelsetvalue}}, 
	\frac{1}{\geop{t} \timefunctionarg{\muxmulevelsetvalue}}, 
	\tander^{\leq 2} \timefunctionarg{\muxmulevelsetvalue},
	\muX \tander \timefunctionarg{\muxmulevelsetvalue}
\right).
\end{align}
From \eqref{E:SCHEMATICIDENTITYFORROUGHTORUSTRACEOFDEFORMATIONTENSOROFRTRANS},
Lemmas~\ref{L:COMMUTATORSTOCOORDINATES} and \ref{L:SCHEMATICSTRUCTUREOFVARIOUSTENSORSINTERMSOFCONTROLVARS},
the estimates 
\eqref{E:ALLDERIVATIVESLINFTYESTIMATESFORROUGHTIMEFUNCTIONANDDERIVATIVES}
and
\eqref{E:PARTIALTIMEDERIVATIVEOFROUGHTIMEFUNCTIONISAPPROXIMATELYUNITY}--\eqref{E:LDERIVATIVEOFROUGHTIMEFUNCTIONISAPPROXIMATELYUNITY}
for the rough time function, and the bootstrap assumptions,
we arrive at the desired estimate \eqref{E:POINTWISEBOUNDFORROUGHTOROIDALTRACEOFDEFORMATIONTENSOROFRTRANS}.

\medskip
\noindent \textbf{Proof of \eqref{E:POINTWISEBOUNDFORROUGHTOROIDALTRACEOFDEFORMATIONTENSOROFROUGHNULLVECTORFIELD}:}
Considering definition~\eqref{E:LROUGH} and arguing as in the proof of
\eqref{E:SCHEMATICIDENTITYFORROUGHTORUSTRACEOFDEFORMATIONTENSOROFRTRANS},
we find that:
\begin{align} \label{E:SCHEMATICIDENTITYFORROUGHTORUSTRACEOFDEFORMATIONTENSOROFROUGHNULLVECTORFIELD}
\mytr_{\gtorusroughfirstfund} \deform{\argLrough{\muxmulevelsetvalue}}
& 
= 
\smoothfunction
\left(\tander^{\leq 1} \controlvars,
	\frac{1}{\Lunit \timefunctionarg{\muxmulevelsetvalue}}, 
	\frac{1}{\geop{t} \timefunctionarg{\muxmulevelsetvalue}}, 
	\tander^{\leq 2} \timefunctionarg{\muxmulevelsetvalue}
\right).
\end{align}
From \eqref{E:SCHEMATICIDENTITYFORROUGHTORUSTRACEOFDEFORMATIONTENSOROFROUGHNULLVECTORFIELD},
Lemmas~\ref{L:COMMUTATORSTOCOORDINATES} and \ref{L:SCHEMATICSTRUCTUREOFVARIOUSTENSORSINTERMSOFCONTROLVARS},
the estimates 
\eqref{E:ALLDERIVATIVESLINFTYESTIMATESFORROUGHTIMEFUNCTIONANDDERIVATIVES}
and
\eqref{E:PARTIALTIMEDERIVATIVEOFROUGHTIMEFUNCTIONISAPPROXIMATELYUNITY}--\eqref{E:LDERIVATIVEOFROUGHTIMEFUNCTIONISAPPROXIMATELYUNITY}
for the rough time function, and the bootstrap assumptions,
we arrive at the desired estimate \eqref{E:POINTWISEBOUNDFORROUGHTOROIDALTRACEOFDEFORMATIONTENSOROFROUGHNULLVECTORFIELD}.

\end{proof}

\subsection{Pointwise estimates for the elliptic-hyperbolic integral identity error terms}
\label{SS:POINTWISEESTIMTAESFORELLIPTICHYPERBOLICIDENTITYERORTERMS} 
Recall that to prove our top-order $L^2$ estimates for the specific vorticity and entropy gradient,
we will rely on the elliptic-hyperbolic integral identity \eqref{E:INTEGRALIDENTITYFORELLIPTICHYPERBOLICCURRENT}
with $\tander^{\Ntop} \vortrenormalized$ and $\tander^{\Ntop} \GradEnt$
in the role of $\SigmatTan$.
In the next proposition, we derive pointwise estimates for the error terms appearing in the identity.

\begin{proposition}[Pointwise estimates for the elliptic-hyperbolic integral identity error terms]
	\label{P:POINTWISEESTIMTAESFORELLIPTICHYPERBOLICIDENTITYERORTERMS} 
	Let $\varsigma \in (0,1]$.
	Then the error terms appearing in the elliptic-hyperbolic integral identity \eqref{E:INTEGRALIDENTITYFORELLIPTICHYPERBOLICCURRENT}
	(with $\tander^{\Ntop} \vortrenormalized$ and $\tander^{\Ntop} \GradEnt$ in the role of $\SigmatTan$)
	satisfy the following pointwise estimates 
	on $\twoargMrough{[\timefunction_0,\timefunctionboot),[- \rightu,\leftu]}{\muxmulevelsetvalue}$,
	where the implicit constants are independent of $\varsigma$.
	
	\medskip
	\noindent \underline{\textbf{Estimates for controlling spacetime error integrals}}.
		\begin{subequations}
		\begin{align}
		\begin{split} \label{E:ANTISYMMETRICVORTICITYELLIPTICHYPERBOLICIDENTITYERORTERMPOINTWISE}	
		\left|
			\mathfrak{J}_{(\textnormal{Antisymmetric})}
			[\pmb{\partial} \tander^{\Ntop} \vortrenormalized,\pmb{\partial} \tander^{\Ntop} \vortrenormalized]
		\right|
		&
		\lesssim
			\left|
				\tander^{\Ntop} \VortVort 
			\right|^2
			+
			\frac{1}{\upmu^2}
			\left|
				\tander^{\leq \Ntop-1} \VortVort 
			\right|^2				
			+
			\frac{1}{\upmu^2}
			\left|
				\tander^{\leq \Ntop} (\vortrenormalized,\GradEnt) 
			\right|^2
				\\
		& \ \
			+
			\frac{\varepsilon^2}{\upmu^2}
			\left|	
				\muX \tander^{[1,\Ntop]} \wavearray 
			\right|^2
			+
			\varepsilon^2 
			\left|
				\tander^{\Ntop+1} \wavearray 
			\right|^2
			+
			\frac{\varepsilon^2}{\upmu^2}
			\left|
				\tander^{[2,\Ntop]} \badcontrolvars 
			\right|^2,
		\end{split}
			\\
		\begin{split} \label{E:ANTISYMMETRICENTROPYGRADIENTELLIPTICHYPERBOLICIDENTITYERORTERMPOINTWISE} 
		\left|
			\mathfrak{J}_{(\textnormal{Antisymmetric})}[\pmb{\partial} \tander^{\Ntop}\GradEnt,\pmb{\partial} \tander^{\Ntop} \GradEnt]
		\right|
		& \lesssim
			\frac{1}{\upmu^2}
			\left|
				\tander^{\leq \Ntop} (\vortrenormalized,\GradEnt) 
			\right|^2
				\\
		&  
			\ \
			+
			\frac{\varepsilon^2}{\upmu^2}
			\left|	
				\muX \tander^{[1,\Ntop]} \wavearray 
			\right|^2
			+
			\varepsilon^2 
			\left|
				\tander^{\Ntop+1} \wavearray 
			\right|^2
			+
			\frac{\varepsilon^2}{\upmu^2}
			\left|
				\tander^{[2,\Ntop]} \badcontrolvars 
			\right|^2,				
	\end{split}
	\end{align}
\end{subequations}	
	
\begin{subequations}
	\begin{align}
	\begin{split} \label{E:DIVVORTICITYELLIPTICHYPERBOLICIDENTITYERORTERMPOINTWISE}	 
	\left|
			\mathfrak{J}_{(\textnormal{Div})}[\pmb{\partial} \tander^{\Ntop} \vortrenormalized,\pmb{\partial} \tander^{\Ntop}\vortrenormalized]
		\right| 			
	& \lesssim
			\frac{1}{\upmu^2}
			|\tander^{\leq \Ntop} (\vortrenormalized,\GradEnt)|^2
				\\
	& \ \
			+
			\frac{\varepsilon^2}{\upmu^2}
			|\muX \tander^{[1,\Ntop]} \wavearray|^2
			+
			\varepsilon^2 |\tander^{\Ntop+1} \wavearray|^2
			+
			\frac{\varepsilon^2}{\upmu^2}
			|\tander^{[2,\Ntop]} \badcontrolvars|^2, 
	\end{split}
					\\
	\begin{split} \label{E:DIVENTROPYGRADIENTELLIPTICHYPERBOLICIDENTITYERORTERMPOINTWISE} 
		\left|
			\mathfrak{J}_{(\textnormal{Div})}[\pmb{\partial} \tander^{\Ntop} \GradEnt,\pmb{\partial} \tander^{\Ntop}\GradEnt]
		\right| 			
		& \lesssim
			\left|
				\tander^{\Ntop} \DivGradEnt 
			\right|^2
			+
			\frac{1}{\upmu^2}
			\left|
				\tander^{\leq \Ntop-1} \DivGradEnt 
			\right|^2
			+
			\frac{1}{\upmu^2}
			|\tander^{\leq \Ntop} (\vortrenormalized,\GradEnt)|^2
				\\
		& \ \
			+
			\frac{\varepsilon^2}{\upmu^2}
			|\muX \tander^{[1,\Ntop]} \wavearray|^2
			+
			\varepsilon^2 |\tander^{\Ntop+1} \wavearray|^2
			+
			\frac{\varepsilon^2}{\upmu^2}
			|\tander^{[2,\Ntop]} \badcontrolvars|^2,
	\end{split}
	\end{align}
\end{subequations}
	
\begin{subequations}
\begin{align}
		\upmu
		\left|
			\mathfrak{J}_{(\pmb{\partial} \frac{1}{\upmu})}
				[\tander^{\Ntop} \vortrenormalized,\pmb{\partial} \tander^{\Ntop} \vortrenormalized]
		\right|
		& \lesssim
			\varsigma 
			\frac{1}{\Lunit \timefunctionarg{\muxmulevelsetvalue}}
			\ellipticCoerciveQuadratic[\pmb{\partial}\tander^{\Ntop} \vortrenormalized,\pmb{\partial} \tander^{\Ntop}			
			\vortrenormalized]
			+
			\frac{1}{\varsigma}
			\frac{1}{\upmu^2}
			|\tander^{\leq \Ntop} \vortrenormalized|^2,
				\label{E:DERIVATIVEOFONEOVERMUWEIGHTVORTICITYELLIPTICHYPERBOLICIDENTITYERORTERMPOINTWISE} 
					\\
		\upmu
		\left|
			\mathfrak{J}_{(\pmb{\partial} \frac{1}{\upmu})}[\tander^{\Ntop} \GradEnt,\pmb{\partial} \tander^{\Ntop} \GradEnt]
		\right|
		& \lesssim
			\varsigma 
			\frac{1}{\Lunit \timefunctionarg{\muxmulevelsetvalue}}
			\ellipticCoerciveQuadratic[\pmb{\partial}\tander^{\Ntop} \GradEnt,\pmb{\partial} \tander^{\Ntop}			
			\GradEnt]
			+
			\frac{1}{\varsigma} 
			\frac{1}{\upmu^2}
			|\tander^{\leq \Ntop} \GradEnt|^2,	
				\label{E:DERIVATIVEOFONEOVERMUWEIGHTENTROPYGRADIENTELLIPTICHYPERBOLICIDENTITYERORTERMPOINTWISE} 
	\end{align}
	\end{subequations}
	
\begin{subequations}
	\begin{align}
	\begin{split} \label{E:FIRSTABSORBTERMVORTICITYELLIPTICHYPERBOLICIDENTITYERORTERMPOINTWISE}  
	\left|
			\mathfrak{J}_{(\textnormal{Absorb-1})}[\tander^{\Ntop} \vortrenormalized,\pmb{\partial} \tander^{\Ntop}\vortrenormalized]
		\right| 
		& 
		\lesssim
			\varsigma 
			\frac{1}{\Lunit \timefunctionarg{\muxmulevelsetvalue}}
			\ellipticCoerciveQuadratic[\pmb{\partial}\tander^{\Ntop} \vortrenormalized,\pmb{\partial} \tander^{\Ntop}\vortrenormalized]
					\\
			& \ \
			+
			\left(1 + \frac{1}{\varsigma} \right)
			\frac{1}{\upmu^2}
			|\tander^{\leq \Ntop} (\vortrenormalized,\GradEnt)|^2
				\\
		& \ \
			+
			\left(1 + \frac{1}{\varsigma} \right)
			\frac{\varepsilon^2}{\upmu^2}
			|\muX \tander^{[1,\Ntop]} \wavearray|^2
			+
			\left(1 + \frac{1}{\varsigma} \right)
			\varepsilon^2 |\tander^{\Ntop+1} \wavearray|^2
				\\
		& \ \
			+
			\left(1 + \frac{1}{\varsigma} \right)
			\frac{\varepsilon^2}{\upmu^2}
			|\tander^{[2,\Ntop]} \badcontrolvars|^2,
		\end{split}
			\\
		\begin{split} \label{E:FIRSTABSORBTERMENTROPYGRADIENTELLIPTICHYPERBOLICIDENTITYERORTERMPOINTWISE} 
		\left|
			\mathfrak{J}_{(\textnormal{Absorb-1})}[\tander^{\Ntop} \GradEnt,\pmb{\partial} \tander^{\Ntop} \GradEnt]
		\right|
		 &
		\lesssim
			\varsigma 
			\frac{1}{\Lunit \timefunctionarg{\muxmulevelsetvalue}}
			\ellipticCoerciveQuadratic[\pmb{\partial}\tander^{\Ntop} \GradEnt,\pmb{\partial} \tander^{\Ntop} \GradEnt]
				\\
		& \ \
			+
			\left(1 + \frac{1}{\varsigma} \right)
			\left|
				\tander^{\Ntop} \DivGradEnt 
			\right|^2
			+
			\left(1 + \frac{1}{\varsigma} \right)
			\frac{1}{\upmu^2}
			\left|
				\tander^{\leq \Ntop-1} \DivGradEnt 
			\right|^2
					\\
			& \ \
				+
			\left(1 + \frac{1}{\varsigma} \right)
			\frac{1}{\upmu^2}
			|\tander^{\leq \Ntop} (\vortrenormalized,\GradEnt)|^2
				\\
		& \ \
			+
			\left(1 + \frac{1}{\varsigma} \right)
			\frac{\varepsilon^2}{\upmu^2}
			|\muX \tander^{[1,\Ntop]} \wavearray|^2
			+
			\left(1 + \frac{1}{\varsigma} \right)
			\varepsilon^2 |\tander^{\Ntop+1} \wavearray|^2
				\\
		& \ \
			+
			\left(1 + \frac{1}{\varsigma} \right)
			\frac{\varepsilon^2}{\upmu^2}
			|\tander^{[2,\Ntop]} \badcontrolvars|^2,
	\end{split}
	\end{align}
\end{subequations}

	\begin{subequations}
		\begin{align}
		\left|
			\mathfrak{J}_{(\textnormal{Absorb-2})}[\tander^{\Ntop} \vortrenormalized,\pmb{\partial} \tander^{\Ntop}\vortrenormalized]
		\right|
		& \lesssim
			\varsigma 
			\frac{1}{\Lunit \timefunctionarg{\muxmulevelsetvalue}}
			\ellipticCoerciveQuadratic[\pmb{\partial}\tander^{\Ntop} \vortrenormalized,\pmb{\partial} \tander^{\Ntop}\vortrenormalized]
			+
			\frac{1}{\varsigma} 
			\frac{1}{\upmu^2}
			|\tander^{\Ntop} \vortrenormalized|^2,
			 \label{E:SECONDABSORBTERMVORTICITYELLIPTICHYPERBOLICIDENTITYERORTERMPOINTWISE} 
				\\
		\left|
			\mathfrak{J}_{(\textnormal{Absorb-2})}[\tander^{\Ntop} \GradEnt,\pmb{\partial} \tander^{\Ntop} \GradEnt]
		\right|
		& \lesssim
			\varsigma 
			\frac{1}{\Lunit \timefunctionarg{\muxmulevelsetvalue}}
			\ellipticCoerciveQuadratic[\pmb{\partial}\tander^{\Ntop} \GradEnt,\pmb{\partial} \tander^{\Ntop} \GradEnt]
			+
			\frac{1}{\varsigma} 
			\frac{1}{\upmu^2}
			|\tander^{\Ntop} \GradEnt|^2,
			\label{E:SECONDABSORBTERMENTROPYGRADIENTELLIPTICHYPERBOLICIDENTITYERORTERMPOINTWISE} 
	\end{align}
	\end{subequations}
	
	\begin{subequations}
	\begin{align}
	\begin{split} \label{E:QUADRATICINMATERIALDERIVATIVETERMVORTICITYELLIPTICHYPERBOLICIDENTITYERORTERMPOINTWISE} 
		\left|
			\mathfrak{J}_{(\textnormal{Material})}[\tander^{\Ntop} \vortrenormalized,\pmb{\partial} \tander^{\Ntop}\vortrenormalized]
		\right| 
		 &
		\lesssim
			\frac{1}{\upmu^2}
			|\tander^{\leq \Ntop} (\vortrenormalized,\GradEnt)|^2
					\\
			& \ \
				+
			\frac{\varepsilon^2}{\upmu^2}
			|\muX \tander^{[1,\Ntop]} \wavearray|^2
			+
			\varepsilon^2 |\tander^{\Ntop+1} \wavearray|^2
			+
			\frac{\varepsilon^2}{\upmu^2}
			|\tander^{[2,\Ntop]} \badcontrolvars|^2,
		\end{split} 	
						\\
		\begin{split} \label{E:QUADRATICINMATERIALDERIVATIVETERMENTROPYGRADIENTELLIPTICHYPERBOLICIDENTITYERORTERMPOINTWISE} 
		\left|
			\mathfrak{J}_{(\textnormal{Material})}[\tander^{\Ntop} \GradEnt,\pmb{\partial} \tander^{\Ntop} \GradEnt]
		\right|
		 &
		\lesssim
		\frac{1}{\upmu^2}
			|\tander^{\leq \Ntop} (\vortrenormalized,\GradEnt)|^2
				\\
		& \ \
			+
			\frac{\varepsilon^2}{\upmu^2}
			|\muX \tander^{[1,\Ntop]} \wavearray|^2
			+
			\varepsilon^2 |\tander^{\Ntop+1} \wavearray|^2
			+
			\frac{\varepsilon^2}{\upmu^2}
			|\tander^{[2,\Ntop]} \badcontrolvars|^2,
	\end{split}
	\end{align}
\end{subequations}
	
	\begin{subequations}
	\begin{align}
		\left|
			\mathfrak{J}_{(\textnormal{Null Geometry})}[\tander^{\Ntop} \vortrenormalized,\pmb{\partial} \tander^{\Ntop} \vortrenormalized]
		\right|
		& \lesssim
			\varsigma 
			\frac{1}{\Lunit \timefunctionarg{\muxmulevelsetvalue}}
			\ellipticCoerciveQuadratic[\pmb{\partial}\tander^{\Ntop} \vortrenormalized,\pmb{\partial} \tander^{\Ntop}\vortrenormalized]
			+
			\frac{1}{\varsigma} 
			\frac{1}{\upmu^2}
			|\tander^{\Ntop} \vortrenormalized|^2,
			 \label{E:NULLGEOMETRYVORTICITYELLIPTICHYPERBOLICIDENTITYERORTERMPOINTWISE}  
				\\
		\left|
			\mathfrak{J}_{(\textnormal{Null Geometry})}[\tander^{\Ntop} \GradEnt,\pmb{\partial} \tander^{\Ntop} \GradEnt]
		\right|
		& \lesssim
			\varsigma 
			\frac{1}{\Lunit \timefunctionarg{\muxmulevelsetvalue}}
			\ellipticCoerciveQuadratic[\pmb{\partial}\tander^{\Ntop} \GradEnt,\pmb{\partial} \tander^{\Ntop} \GradEnt]
			+
			\frac{1}{\varsigma} 
			\frac{1}{\upmu^2}
			|\tander^{\Ntop} \GradEnt|^2.
			\label{E:NULLGEOMETRYENTROPYGRADIENTELLIPTICHYPERBOLICIDENTITYERORTERMPOINTWISE}
	\end{align}
	\end{subequations}
	
	\medskip
	\noindent \underline{\textbf{Estimates for controlling spatial error integrals}}.
	
	\begin{subequations}
	\begin{align}
	\begin{split} \label{E:PRINCIPALSIGAMTILDEELLIPTICHYPERBOLICVORTICITYPOINTWISE}	
		\left|
			\Currentboundaryerrorhavetocontrolprincipal[\tander^{\Ntop} \vortrenormalized,\pmb{\partial} \tander^{\Ntop} \vortrenormalized]
		\right|
		& 
		\lesssim
			\left\lbrace
				\upmu
				-
				\phi \frac{\muxmulevelsetvalue}{\Lunit \upmu}
			\right\rbrace
			\left|
				\tander^{\Ntop} \VortVort 
			\right|^2
			+
			\frac{1}{\upmu^{3/2}}
			\left\lbrace
				\upmu
				-
				\phi \frac{\muxmulevelsetvalue}{\Lunit \upmu}
			\right\rbrace
			\left|
				\tander^{\leq \Ntop-1} \VortVort 
			\right|^2
					\\
		& \ \
			+
			\frac{1}{\upmu^{5/2}}
			\left\lbrace
				\upmu
				-
				\phi \frac{\muxmulevelsetvalue}{\Lunit \upmu}
			\right\rbrace
			\left|
				\tander^{\leq \Ntop} (\vortrenormalized,\GradEnt) 
			\right|^2
					\\
		& \ \
			+
			\frac{\varepsilon^2}{\upmu^{3/2}}
			\left|	
				\muX \tander^{[1,\Ntop]} \wavearray 
			\right|^2
			+
			\varepsilon^2 
			\left\lbrace
				\upmu
				-
				\phi \frac{\muxmulevelsetvalue}{\Lunit \upmu}
			\right\rbrace
			\left|
				\tander^{\Ntop+1} \wavearray 
			\right|^2
				\\
		& \ \
			+
			\frac{\varepsilon^2}{\upmu^{3/2}}
			\left|
				\tander^{[2,\Ntop]} \badcontrolvars 
			\right|^2,
	\end{split}
				\\
	\begin{split} \label{E:PRINCIPALSIGAMTILDEELLIPTICHYPERBOLICENTROPYGRADIENTPOINTWISE}
		\left|
			\Currentboundaryerrorhavetocontrolprincipal[\tander^{\Ntop} \GradEnt,\pmb{\partial} \tander^{\Ntop} \GradEnt]
		\right|
		& 
		\lesssim
			\left\lbrace
				\upmu
				-
				\phi \frac{\muxmulevelsetvalue}{\Lunit \upmu}
			\right\rbrace
			\left|
				\tander^{\Ntop} \DivGradEnt
			\right|^2
			+
			\frac{1}{\upmu^{3/2}}
			\left\lbrace
				\upmu
				-
				\phi \frac{\muxmulevelsetvalue}{\Lunit \upmu}
			\right\rbrace
			\left|
				\tander^{\leq \Ntop-1} \DivGradEnt 
			\right|^2
						\\
		& \ \
			+
			\frac{1}{\upmu^{5/2}}
			\left\lbrace
				\upmu
				-
				\phi \frac{\muxmulevelsetvalue}{\Lunit \upmu}
			\right\rbrace
			\left|
				\tander^{\leq \Ntop} (\vortrenormalized,\GradEnt) 
			\right|^2
				 \\
		& \ \
			+
			\frac{\varepsilon^2}{\upmu^{3/2}}
			\left|	
				\muX \tander^{[1,\Ntop]} \wavearray 
			\right|^2
			+
			\varepsilon^2 
			\left\lbrace
				\upmu
				-
				\phi \frac{\muxmulevelsetvalue}{\Lunit \upmu}
			\right\rbrace
			\left|
				\tander^{\Ntop+1} \wavearray 
			\right|^2
				\\
		& \ \
			+
			\frac{\varepsilon^2}{\upmu^{3/2}}
			\left|
				\tander^{[2,\Ntop]} \badcontrolvars 
			\right|^2,
	\end{split} 
	\end{align}
	\end{subequations}
	
	\begin{subequations}
	\begin{align}
		\left|
			\Currentboundaryerrorhavetocontrollowerorder[\tander^{\Ntop} \vortrenormalized,\tander^{\Ntop} \vortrenormalized]
		\right|
		& \lesssim
			\frac{1}{\upmu^{5/2}}
			\left\lbrace
				\upmu
				-
				\phi \frac{\muxmulevelsetvalue}{\Lunit \upmu}
			\right\rbrace
			\left|
				\tander^{\Ntop} \vortrenormalized
			\right|^2,
			\label{E:LOWERORDERSIGAMTILDEELLIPTICHYPERBOLICVORTICITYPOINTWISE}	
				\\
		\left|
			\Currentboundaryerrorhavetocontrollowerorder[\tander^{\Ntop} \GradEnt,\tander^{\Ntop} \GradEnt]
		\right|
		& \lesssim
			\frac{1}{\upmu^{5/2}}
			\left\lbrace
				\upmu
				-
				\phi \frac{\muxmulevelsetvalue}{\Lunit \upmu}
			\right\rbrace
			\left|
				\tander^{\Ntop} \GradEnt
			\right|^2.
			\label{E:LOWERORDERSIGAMTILDEELLIPTICHYPERBOLIENTROPYGRADIENTPOINTWISE}	
	\end{align}
	\end{subequations}
	
\end{proposition}

\begin{proof}
\hfill

\noindent \textbf{Proofs of \eqref{E:ANTISYMMETRICVORTICITYELLIPTICHYPERBOLICIDENTITYERORTERMPOINTWISE}--\eqref{E:ANTISYMMETRICENTROPYGRADIENTELLIPTICHYPERBOLICIDENTITYERORTERMPOINTWISE}
and \eqref{E:DIVVORTICITYELLIPTICHYPERBOLICIDENTITYERORTERMPOINTWISE}--\eqref{E:DIVENTROPYGRADIENTELLIPTICHYPERBOLICIDENTITYERORTERMPOINTWISE}:}
	First, using \eqref{E:HRIEMANNIANNORMOFCHARACTERISTICTENSORFIELDSISROOT3}
	and \eqref{E:ANTISYMMETRICNULLCURRENTSPACETIMERRORTERM},
	we deduce that
	$|\mathfrak{J}_{(\textnormal{Antisymmetric})}[\pmb{\partial} \SigmatTan,\pmb{\partial} \SigmatTan]|
	\lesssim 
	|\mathrm{d} \SigmatTan_{\flat}|_{\hfour}^2
	$.
From this bound with 
$\tander^{\Ntop} \vortrenormalized$
and
$\tander^{\Ntop} \GradEnt$
in the role of $\SigmatTan$
and
the pointwise estimates
\eqref{E:POINTWISEBOUNDEXTERIORDERIVATIVEOFTOPORDERDERIVATIVESOFVORTICITY}--\eqref{E:POINTWISEBOUNDEXTERIORDERIVATIVEOFTOPORDERDERIVATIVESOFENTROPYGRADIENT},
we conclude
the desired bounds
\eqref{E:ANTISYMMETRICVORTICITYELLIPTICHYPERBOLICIDENTITYERORTERMPOINTWISE}--\eqref{E:ANTISYMMETRICENTROPYGRADIENTELLIPTICHYPERBOLICIDENTITYERORTERMPOINTWISE}.

\eqref{E:DIVVORTICITYELLIPTICHYPERBOLICIDENTITYERORTERMPOINTWISE}--\eqref{E:DIVENTROPYGRADIENTELLIPTICHYPERBOLICIDENTITYERORTERMPOINTWISE}	
follow from similar arguments
based on \eqref{E:DIVERGENCENULLCURRENTSPACETIMERRORTERM}
and the pointwise estimates \eqref{E:COMMUTEDPOINTWISEEUCLIDEANEUCLIDEANDIVERGENCEOFENTROPYGRADIEENT}--\eqref{E:COMMUTEDPOINTWISEEUCLIDEANDIVERGENCEOFVORTICITYANDEUCLIDEANCURLOFENTROPYGRADIENT}.

\medskip
\noindent \textbf{Proofs of \eqref{E:DERIVATIVEOFONEOVERMUWEIGHTVORTICITYELLIPTICHYPERBOLICIDENTITYERORTERMPOINTWISE}--\eqref{E:DERIVATIVEOFONEOVERMUWEIGHTENTROPYGRADIENTELLIPTICHYPERBOLICIDENTITYERORTERMPOINTWISE}:}
Let $\SigmatTan$ be any $\Sigma_t$-tangent vectorfield. 
We first note that since the elliptic hyperbolic current
$
\ehcurrent^{\alpha}[\SigmatTan,\pmb{\partial} \SigmatTan]
$
defined by \eqref{E:PUTANGENTELLIPTICHYPERBOLICCURRENT}
is tangent to $\nullhyparg{u}$,
we can use the identity \eqref{E:ELLTUPROJECTIONINDOUBLENULLFRAME} to obtain the following
decomposition, where $\smoothtorusproject$ is the $\ell_{t,u}$-projection
tensorfield from Def.\,\ref{D:PROJECTIONTENSORFIELDSANDTANGENCYTOHYPERSURFACES}:
$
\ehcurrent^{\alpha}[\SigmatTan,\pmb{\partial} \SigmatTan]
=
- \frac{1}{2}
	\uLunit_{\beta}
	\ehcurrent^{\beta}[\SigmatTan,\pmb{\partial} \SigmatTan]
	\Lunit^{\alpha}
	+
	\smoothtorusproject_{\beta}^{\ \alpha} 
	\ehcurrent^{\beta}[\SigmatTan,\pmb{\partial} \SigmatTan]
$.
From this decomposition
and definition~\eqref{E:DERIVATIVEOFWEIGHTNULLCURRENTSPACETIMERRORTERM},
we deduce the following pointwise bound, where
$
\smoothtorusproject \ehcurrent[\SigmatTan,\pmb{\partial} \SigmatTan]
$
is the $\ell_{t,u}$-projection of $\ehcurrent[\SigmatTan,\pmb{\partial} \SigmatTan]$:
\begin{align} \label{E:FIRSTPOINTWISEBOUNDFORCHARCURRENTDERIVATIVEOF1OVERMUWEIGHTERRORTERM}
\upmu
|\mathfrak{J}_{(\pmb{\partial} \frac{1}{\upmu})}[\SigmatTan,\pmb{\partial} \SigmatTan]|
\lesssim
\frac{1}{\upmu}
|\Lunit \upmu|
|\uLunit_{\beta} \ehcurrent^{\beta}[\SigmatTan,\pmb{\partial} \SigmatTan]|+
|\angD \upmu|_{\gtorus}
|\smoothtorusproject \ehcurrent[\SigmatTan,\pmb{\partial} \SigmatTan]|_{\gtorus}.
\end{align}
Next, using \eqref{E:FIRSTPOINTWISEBOUNDFORCHARCURRENTDERIVATIVEOF1OVERMUWEIGHTERRORTERM},
definition \eqref{E:PUTANGENTELLIPTICHYPERBOLICCURRENT},
Lemma~\ref{L:SCHEMATICSTRUCTUREOFVARIOUSTENSORSINTERMSOFCONTROLVARS},
the estimates of Prop.\,\ref{P:IMPROVEMENTOFAUXILIARYBOOTSTRAP},
and the pointwise comparison estimates provided by Lemma~\ref{L:HRIEMANNIANMETRICOMPARISONANDNORMSOFNULLHYPERSURFACETENSORFIELDS},
we deduce the following pointwise bound:
\begin{align} \label{E:SECONDPOINTWISEBOUNDFORCHARCURRENTDERIVATIVEOF1OVERMUWEIGHTERRORTERM}
\upmu
|\mathfrak{J}_{(\pmb{\partial} \frac{1}{\upmu})}[\SigmatTan,\pmb{\partial} \SigmatTan]|
& 
\lesssim
\frac{1}{\upmu}
|\SigmatTan|_g
|\pmb{\partial} \SigmatTan|_{\hfour}.
\end{align}
Using \eqref{E:SECONDPOINTWISEBOUNDFORCHARCURRENTDERIVATIVEOF1OVERMUWEIGHTERRORTERM} with 
$\tander^{\Ntop} \vortrenormalized$
and
$\tander^{\Ntop} \GradEnt$
in the role of $\SigmatTan$,
\eqref{E:COERCIVENESSOFELLIPTICHYPERBOLICQUADRATICFORM},
the bound $\frac{1}{\Lunit \timefunctionarg{\muxmulevelsetvalue}} \approx 1$
(see \eqref{E:ROUGHTIMEFUNCTIONLDERIVATIVEBOUNDS}),
\eqref{E:SIGMATTANGENTVECTORFIELDGNORMAPPROXIMATEDBYCARETSIANCOMPONENTNORMS},
and Young's inequality (where we multiply and divide by powers of $\varsigma$), 
we deduce the desired bounds
\eqref{E:DERIVATIVEOFONEOVERMUWEIGHTVORTICITYELLIPTICHYPERBOLICIDENTITYERORTERMPOINTWISE}--\eqref{E:DERIVATIVEOFONEOVERMUWEIGHTENTROPYGRADIENTELLIPTICHYPERBOLICIDENTITYERORTERMPOINTWISE}.

\medskip
\noindent \textbf{Proofs of \eqref{E:FIRSTABSORBTERMVORTICITYELLIPTICHYPERBOLICIDENTITYERORTERMPOINTWISE}--\eqref{E:FIRSTABSORBTERMENTROPYGRADIENTELLIPTICHYPERBOLICIDENTITYERORTERMPOINTWISE}:}
Let $\SigmatTan$ be any $\Sigma_t$-tangent vectorfield. 
We first use definition \eqref{E:ABSORBABLEANTISYMMETRICANDDIVERGENCENULLCURRENTSPACETIMERRORTERM},
\eqref{E:NULLHYPERSURFACEINVERSERIEMANNIANMETRICINTERMSOFNULLFRAME},
Lemma~\ref{L:SCHEMATICSTRUCTUREOFVARIOUSTENSORSINTERMSOFCONTROLVARS},
the estimates of Prop.\,\ref{P:IMPROVEMENTOFAUXILIARYBOOTSTRAP},
the bound $|\partial_{\alpha} \gfour_{\beta \gamma}| \lesssim |\pmb{\partial} \wavearray| \lesssim \frac{1}{\upmu}$
(which follows from Lemma~\ref{L:RELATIONSHIPBETWEENCARTESIANPARTIALDERIVATIVESANDSMOOTHGEOMETRICCOMMUTATORS}
Lemma~\ref{L:SCHEMATICSTRUCTUREOFVARIOUSTENSORSINTERMSOFCONTROLVARS}, and Prop.\,\ref{P:IMPROVEMENTOFAUXILIARYBOOTSTRAP}),
and the pointwise comparison estimates provided by Lemma~\ref{L:HRIEMANNIANMETRICOMPARISONANDNORMSOFNULLHYPERSURFACETENSORFIELDS}
to deduce that:
\begin{align} \label{E:FIRSTPOINTWISEBOUNDFIRSTABSORBELLIPTICHYPERBOLICCURRENTERRORTERM}
\left|
	\mathfrak{J}_{(\textnormal{Absorb-1})}[\SigmatTan,\pmb{\partial} \SigmatTan]
\right|
&
\lesssim
\frac{1}{\upmu}
|\SigmatTan|_g
|\pmb{\partial} \SigmatTan|_{\hfour}
+
|\partial_a \SigmatTan^a|
|\pmb{\partial} \SigmatTan|_{\hfour}.
\end{align}
Using \eqref{E:FIRSTPOINTWISEBOUNDFIRSTABSORBELLIPTICHYPERBOLICCURRENTERRORTERM}
with
$\tander^{\Ntop} \vortrenormalized$
and
$\tander^{\Ntop} \GradEnt$
in the role of $\SigmatTan$,
\eqref{E:COERCIVENESSOFELLIPTICHYPERBOLICQUADRATICFORM},
the bound $\frac{1}{\Lunit \timefunctionarg{\muxmulevelsetvalue}} \approx 1$
(see \eqref{E:ROUGHTIMEFUNCTIONLDERIVATIVEBOUNDS}),
the pointwise estimates
\eqref{E:COMMUTEDPOINTWISEEUCLIDEANEUCLIDEANDIVERGENCEOFENTROPYGRADIEENT}--\eqref{E:COMMUTEDPOINTWISEEUCLIDEANDIVERGENCEOFVORTICITYANDEUCLIDEANCURLOFENTROPYGRADIENT},
\eqref{E:SIGMATTANGENTVECTORFIELDGNORMAPPROXIMATEDBYCARETSIANCOMPONENTNORMS},
and Young's inequality (where we multiply and divide by powers of $\upmu$ and $\varsigma$ as needed), 
we deduce the desired bounds
\eqref{E:FIRSTABSORBTERMVORTICITYELLIPTICHYPERBOLICIDENTITYERORTERMPOINTWISE}--\eqref{E:FIRSTABSORBTERMENTROPYGRADIENTELLIPTICHYPERBOLICIDENTITYERORTERMPOINTWISE}.

\medskip
\noindent \textbf{Proofs of \eqref{E:SECONDABSORBTERMVORTICITYELLIPTICHYPERBOLICIDENTITYERORTERMPOINTWISE}--\eqref{E:SECONDABSORBTERMENTROPYGRADIENTELLIPTICHYPERBOLICIDENTITYERORTERMPOINTWISE}:}
Let $\SigmatTan$ be any $\Sigma_t$-tangent vectorfield.
We first use definitions \eqref{E:EASYABSORBABLENULLCURRENTSPACETIMERRORTERM}
and
\eqref{E:PUTANGENTELLIPTICHYPERBOLICCURRENT}
and the same arguments used in the proof of \eqref{E:FIRSTPOINTWISEBOUNDFIRSTABSORBELLIPTICHYPERBOLICCURRENTERRORTERM}
to deduce the pointwise bound
$
\left|
	\mathfrak{J}_{(\textnormal{Absorb-2})}[\SigmatTan,\pmb{\partial} \SigmatTan]
\right|
\lesssim
\frac{1}{\upmu}
|\SigmatTan|_g
|\pmb{\partial} \SigmatTan|_{\hfour}
$.
Then using the same arguments given just below 
\eqref{E:SECONDPOINTWISEBOUNDFORCHARCURRENTDERIVATIVEOF1OVERMUWEIGHTERRORTERM},
we conclude \eqref{E:SECONDABSORBTERMVORTICITYELLIPTICHYPERBOLICIDENTITYERORTERMPOINTWISE}--\eqref{E:SECONDABSORBTERMENTROPYGRADIENTELLIPTICHYPERBOLICIDENTITYERORTERMPOINTWISE}.

\medskip
\noindent \textbf{Proofs of \eqref{E:QUADRATICINMATERIALDERIVATIVETERMVORTICITYELLIPTICHYPERBOLICIDENTITYERORTERMPOINTWISE}--\eqref{E:QUADRATICINMATERIALDERIVATIVETERMENTROPYGRADIENTELLIPTICHYPERBOLICIDENTITYERORTERMPOINTWISE}:}
Let $\SigmatTan$ be any $\Sigma_t$-tangent vectorfield.
We first use definition \eqref{E:QUADRATICINMATERIALDERIVATIVESNULLCURRENTSPACETIMERRORTERM},
\eqref{E:HRIEMANNIANNORMOFCHARACTERISTICTENSORFIELDSISROOT3},
and the pointwise comparison estimates provided by Lemma~\ref{L:HRIEMANNIANMETRICOMPARISONANDNORMSOFNULLHYPERSURFACETENSORFIELDS}
to deduce the pointwise bound
$
\left|
	\mathfrak{J}_{(\textnormal{Material})}[\pmb{\partial} \SigmatTan,\pmb{\partial} \SigmatTan]
\right|
\lesssim
|\Transport \SigmatTan|_{\hfour}^2
\lesssim
\sum_{a=1,2,3}
|\Transport \SigmatTan^a|^2
$.
From this bound with
$\tander^{\Ntop} \vortrenormalized$
and
$\tander^{\Ntop} \GradEnt$
in the role of $\SigmatTan$
and the pointwise estimates \eqref{E:COMMUTEDTRANSPORTPOINTWISEESTIMATESFORSPECIFICVORTICITYANDENTROPYGRADIENT},
we conclude
\eqref{E:QUADRATICINMATERIALDERIVATIVETERMVORTICITYELLIPTICHYPERBOLICIDENTITYERORTERMPOINTWISE}--\eqref{E:QUADRATICINMATERIALDERIVATIVETERMENTROPYGRADIENTELLIPTICHYPERBOLICIDENTITYERORTERMPOINTWISE}.

\medskip
\noindent \textbf{Proofs of \eqref{E:NULLGEOMETRYVORTICITYELLIPTICHYPERBOLICIDENTITYERORTERMPOINTWISE}--\eqref{E:NULLGEOMETRYENTROPYGRADIENTELLIPTICHYPERBOLICIDENTITYERORTERMPOINTWISE}:}
Let $\SigmatTan$ be any $\Sigma_t$-tangent vectorfield.
We first use definitions
\eqref{E:DERIVATIVESOFNULLGEOMETRYNULLCURRENTSPACETIMERRORTERM}
and
\eqref{E:NULLHYPERSURFACEPROJECTION},
Lemma~\ref{L:RELATIONSHIPBETWEENCARTESIANPARTIALDERIVATIVESANDSMOOTHGEOMETRICCOMMUTATORS},
Lemma~\ref{L:SCHEMATICSTRUCTUREOFVARIOUSTENSORSINTERMSOFCONTROLVARS},
and the estimates of Prop.\,\ref{P:IMPROVEMENTOFAUXILIARYBOOTSTRAP}
to deduce the pointwise bound
$
\left|
	\mathfrak{J}_{(\textnormal{Null Geometry})}[\SigmatTan,\pmb{\partial} \SigmatTan]
\right|
\lesssim
\frac{1}{\upmu}
|\SigmatTan|_g
|\pmb{\partial} \SigmatTan|_{\hfour}
$.
Then using the same arguments given just below 
\eqref{E:SECONDPOINTWISEBOUNDFORCHARCURRENTDERIVATIVEOF1OVERMUWEIGHTERRORTERM},
we conclude \eqref{E:NULLGEOMETRYVORTICITYELLIPTICHYPERBOLICIDENTITYERORTERMPOINTWISE}--\eqref{E:NULLGEOMETRYENTROPYGRADIENTELLIPTICHYPERBOLICIDENTITYERORTERMPOINTWISE}.

\medskip
\noindent \textbf{Proofs of \eqref{E:PRINCIPALSIGAMTILDEELLIPTICHYPERBOLICVORTICITYPOINTWISE}--\eqref{E:PRINCIPALSIGAMTILDEELLIPTICHYPERBOLICENTROPYGRADIENTPOINTWISE}:}
Let $\SigmatTan$ be any $\Sigma_t$-tangent vectorfield. 
We first use
definitions~\eqref{E:PRINCIPALERRORTERMHAVETOCONTROLKEYIDPUTANGENTCURRENTCONTRACTEDAGAINSTVECTORFIELD}
and \eqref{E:NULLHYPERSURFACEPROJECTION},
Lemma~\ref{L:SCHEMATICSTRUCTUREOFVARIOUSTENSORSINTERMSOFCONTROLVARS},
the estimates of Prop.\,\ref{P:IMPROVEMENTOFAUXILIARYBOOTSTRAP},
\eqref{E:RTRANSCARTESIANPOINTWISE},
\eqref{E:BOUNDSONLMUINTERESTINGREGION},
and the pointwise comparison estimates provided by Lemma~\ref{L:HRIEMANNIANMETRICOMPARISONANDNORMSOFNULLHYPERSURFACETENSORFIELDS}
to deduce the following pointwise estimate, where $\phi$ is the cut-off function from Def.\,\ref{D:WTRANSANDCUTOFF}:
	\begin{align} \label{E:FIRSTPOINTWISEPROOFSTEPPRINCIPALSIGAMTILDEELLIPTICHYPERBOLICVORTICITYPOINTWISE}
	\left|
		\Currentboundaryerrorhavetocontrolprincipal[\SigmatTan,\pmb{\partial} \SigmatTan]
	\right|
	&
	\lesssim 
	\frac{1}{\upmu}
			\left\lbrace
				\upmu
				-
				\phi \frac{\muxmulevelsetvalue}{\Lunit \upmu}
			\right\rbrace
			|\SigmatTan|_g
			\left\lbrace
				|\Transport \SigmatTan|_g
				+
				|\mathrm{d} \SigmatTan_{\flat}|_{\hfour}
				+
				|\partial_a \SigmatTan^a|
			\right\rbrace.
\end{align}
Using \eqref{E:FIRSTPOINTWISEPROOFSTEPPRINCIPALSIGAMTILDEELLIPTICHYPERBOLICVORTICITYPOINTWISE} with 
$\tander^{\Ntop} \vortrenormalized$
and
$\tander^{\Ntop} \GradEnt$
in the role of $\SigmatTan$,
the pointwise estimates
\eqref{E:COMMUTEDTRANSPORTPOINTWISEESTIMATESFORSPECIFICVORTICITYANDENTROPYGRADIENT},
\eqref{E:COMMUTEDPOINTWISEEUCLIDEANCURLOFVORTICITY}--\eqref{E:COMMUTEDPOINTWISEEUCLIDEANDIVERGENCEOFVORTICITYANDEUCLIDEANCURLOFENTROPYGRADIENT},
and
\eqref{E:POINTWISEBOUNDEXTERIORDERIVATIVEOFTOPORDERDERIVATIVESOFVORTICITY}--\eqref{E:POINTWISEBOUNDEXTERIORDERIVATIVEOFTOPORDERDERIVATIVESOFENTROPYGRADIENT} 
with $\Ntop$ in the role of $N$,
the pointwise comparison estimates
\eqref{E:SIGMATTANGENTVECTORFIELDGNORMAPPROXIMATEDBYCARETSIANCOMPONENTNORMS}--\eqref{E:ZDERIVATIVEOFSIGMATTANGENTVECTORFIELDGNORMAPPROXIMATEDBYCARETSIANCOMPONENTNORMS},
the estimate 
$\left|
				\upmu
				-
				\phi \frac{\muxmulevelsetvalue}{\Lunit \upmu}
			\right| \lesssim 1$
(which follows from Prop.\,\ref{P:IMPROVEMENTOFAUXILIARYBOOTSTRAP} and \eqref{E:BOUNDSONLMUINTERESTINGREGION})
and Young's inequality (where we multiply and divide by powers of $\upmu$ and $\varsigma$ as needed), 
we conclude
that the desired estimates 
\eqref{E:PRINCIPALSIGAMTILDEELLIPTICHYPERBOLICVORTICITYPOINTWISE}--\eqref{E:PRINCIPALSIGAMTILDEELLIPTICHYPERBOLICENTROPYGRADIENTPOINTWISE}
hold for any $\varsigma \in (0,1]$.

\medskip
\noindent \textbf{Proofs of \eqref{E:LOWERORDERSIGAMTILDEELLIPTICHYPERBOLICVORTICITYPOINTWISE}--\eqref{E:LOWERORDERSIGAMTILDEELLIPTICHYPERBOLIENTROPYGRADIENTPOINTWISE}:}
Let $\SigmatTan$ be any $\Sigma_t$-tangent vectorfield.
We first use definition~\eqref{E:LOWERORDERERRORTERMHAVETOCONTROLKEYIDPUTANGENTCURRENTCONTRACTEDAGAINSTVECTORFIELD},
Lemma~\ref{L:COMMUTATORSTOCOORDINATES},
Cor.\,\ref{C:ELLTUPROJECTEDVERSIONOFCARTESIANPARTIALDERIVATIVES},
Lemma~\ref{L:SCHEMATICSTRUCTUREOFVARIOUSTENSORSINTERMSOFCONTROLVARS},
Lemma~\ref{L:DIFFEOMORPHICEXTENSIONOFROUGHCOORDINATES},
Prop.\,\ref{P:IMPROVEMENTOFAUXILIARYBOOTSTRAP},
Lemma~\ref{L:POINTWISEESTIMATESINVOLVINGROUGHACOUSTICGEOMETRY},
\eqref{E:WMANDRTRANSMUUBOUNDEDBYSQRTMU},
\eqref{E:BOUNDSONLMUINTERESTINGREGION},
\eqref{E:KEYCOERCIVITYPERFECTRDERIVAIVEERRORPUTANGENTCURRENTCONTRACTEDAGAINSTVECTORFIELD},
and the pointwise comparison estimates provided by Lemma~\ref{L:HRIEMANNIANMETRICOMPARISONANDNORMSOFNULLHYPERSURFACETENSORFIELDS}
to deduce the following pointwise bound, where $\phi$ is the cut-off function from Def.\,\ref{D:WTRANSANDCUTOFF}:
$
	\left|
		\Currentboundaryerrorhavetocontrollowerorder[\SigmatTan,\SigmatTan]
	\right|
	\lesssim
	\frac{1}{\upmu^{5/2}}
			\left\lbrace
				\upmu
				-
				\phi \frac{\muxmulevelsetvalue}{\Lunit \upmu}
			\right\rbrace
			\left|
				\SigmatTan
			\right|_g^2
	$.
	Using this bound with 
$\tander^{\Ntop} \vortrenormalized$
and
$\tander^{\Ntop} \GradEnt$
in the role of $\SigmatTan$
and the comparison estimate \eqref{E:SIGMATTANGENTVECTORFIELDGNORMAPPROXIMATEDBYCARETSIANCOMPONENTNORMS},
we conclude the desired bounds
\eqref{E:LOWERORDERSIGAMTILDEELLIPTICHYPERBOLICVORTICITYPOINTWISE}--\eqref{E:LOWERORDERSIGAMTILDEELLIPTICHYPERBOLIENTROPYGRADIENTPOINTWISE}.
	
\end{proof}

\section{Statement of the a priori $L^2$ estimates, 
data estimates for the $L^2$-controlling quantities, and 
bootstrap assumptions for the wave variable energies}
\label{S:STATEMENTOFALLL2ESTIMATESANDBOOTSTRAPASSUMPTIONSFORWAVEENERGIES}
In Sect.\,\ref{SS:STATEMENTOFAPRIORIL2ESTIMATES}, we state all of the a priori energy estimates
for the fluid variables and the acoustic geometry.
In particular, we state the main estimates for the $L^2$-controlling quantities from
Sect.\,\ref{SS:FUNDAMENTALL2CONTROLLINGQUANTITIES}.
The proofs take considerable effort and form the focus of the paper
through Sect.\,\ref{S:WAVEANDACOUSTICGEOMETRYAPRIORIESTIMATES}.
In Sect.\,\ref{SS:DATAESTIMATESFORL2CONTROLLINGQUANTITIES}, 
we show that along the data hypersurfaces 
$\hypthreearg{\timefunction_0}{[- \rightu,\leftu]}{\muxmulevelsetvalue}$
and
$\nullhypthreearg{\muxmulevelsetvalue}{-\rightu}{[\timefunction_0,\timefunctionboot)}$,
the $L^2$-controlling quantities are bounded by $\lesssim \initialsmall^2$,
i.e., the $L^2$-controlling quantities have small data.
Finally, in Sect.\,\ref{SS:BOOTSTRAPASSUMPTIONSFORTHEWAVEENERGIES},
we state bootstrap assumptions for the $\totalcontrolwave_N$,
i.e., for the $L^2$-controlling quantities of the wave variables.

\subsection{Statement of the a priori $L^2$ estimates}
\label{SS:STATEMENTOFAPRIORIL2ESTIMATES}

\subsubsection{Statement of the a priori $L^2$ estimates for the wave variables}
\label{SSS:STATEMENTOFAPRIORIL2ESTIMATESFORWAVEVARIABLES}
In the following proposition, we state our main a priori energy estimates
for the wave variables. Its proof is located in Sect.\,\ref{SSS:PROOFOFAPRIORIL2ESTIMATESWAVEVARIABLES}.

\begin{proposition}[The main a priori estimates for $\totalcontrolwave_{[1,\Ntop]}$] \label{P:APRIORIL2ESTIMATESWAVEVARIABLES} 
	Let $\totalcontrolwave_{N}(\timefunction,u)$ 
	be the $L^2$-controlling quantity for the wave variables 
	$\wavearray$, as defined in \eqref{E:WAVETOTALL2CONTROLLINGQUANTITY}.
	Under the data assumptions of Sect.\,\ref{S:ASSUMPTIONSONTHEDATA},
	the parameter size-assumptions of Sect.\,\ref{SS:PARAMETERSIZEASSUMPTIONS},
	and the bootstrap assumptions of Sect.\,\ref{S:BOOTSTRAPEVERYTHINGEXCEPTENERGIES},
	there exists a constant $C > 0$
	such that the following estimates hold for 
	$(\timefunction,u) \in [\timefunction_0,\timefunctionboot) \times [- \rightu,\leftu]$: 
\begin{subequations} 
\begin{align}
\totalcontrolwave_{\Ntop - K}(\timefunction,u)  
& \leq C \initialsmall^2 |\timefunction|^{-15.6 + 2K} 
&& 
\mbox{if } 0 \leq K \leq 7, 
	\label{E:MAINWAVEENERGYESTIMATESBLOWUP}
			\\
\totalcontrolwave_{N}(\timefunction,u) 
& \leq C \initialsmall^2
&& \mbox{if } 1 \leq N \leq \Ntop - 8.
	\label{E:MAINWAVEENERGYESTIMATESREGULAR}
\end{align}
\end{subequations}

\end{proposition}

\subsubsection{Statement of the a priori $L^2$ estimates for the transport variables}
\label{SSS:STATEMENTOFAPRIORIL2ESTIMATESFORTRANSPORTVARIABLES}
The following proposition provides an analog of
Prop.\,\ref{P:APRIORIL2ESTIMATESWAVEVARIABLES} for the transport variables.
Its proof is located in Sects.\,\ref{SS:PROOFOFBELOWTOPORDERTRANSPORTENERGYESTIMATES}
and \ref{SS:PROOFOFMAINVORTVORTDIVGRADENTTOPORDERBLOWUP}.

\begin{proposition}[The main a priori $L^2$ estimates for the transport variables on hypersurfaces] 
	\label{P:MAINHYPERSURFACEENERGYESTIMATESFORTRANSPORTVARIABLES} 
	Let
	$\hypersurfacecontrolVort_{N}(\timefunction,u)$
	and 
	$\hypersurfacecontrolGradEnt_{N}(\timefunction,u)$ 
	be the $L^2$-controlling quantities for
	$\vortrenormalized$ and $\GradEnt$,
	as defined in \eqref{E:VORTICITYL2CONTROLLINGQUANTITY}--\eqref{E:ENTROPYGRADIENTL2CONTROLLINGQUANTITY},
	and let $\hypersurfacecontrolVortVort_{N}(\timefunction,u)$
	and
		$\hypersurfacecontrolDivGradEnt_{N}(\timefunction,u)$
	be the $L^2$-controlling quantities for
	the modified fluid variables $\VortVort$ and $\DivGradEnt$,
	as defined in \eqref{E:MODIFIEDVORTICITYVORTICITYL2CONTROLLINGQUANTITY}--\eqref{E:MODIFIEDDIVGRADENTL2CONTROLLINGQUANTITY}.
	Under the data assumptions of Sect.\,\ref{S:ASSUMPTIONSONTHEDATA},
	the parameter size-assumptions of Sect.\,\ref{SS:PARAMETERSIZEASSUMPTIONS},
	and the bootstrap assumptions of Sect.\,\ref{S:BOOTSTRAPEVERYTHINGEXCEPTENERGIES},
	there exists a constant $C > 0$
	such that the following estimates hold for 
	$(\timefunction,u) \in [\timefunction_0,\timefunctionboot) \times [- \rightu,\leftu]$: 
\begin{subequations} 
\begin{align}
	\hypersurfacecontrolVort_{\Ntop - K}(\timefunction,u), 
		\,
	\hypersurfacecontrolGradEnt_{\Ntop - K}(\timefunction,u) 
	& \leq C \initialsmall^2 |\timefunction|^{-14.6 + 2K }, & & \text{for } 0 \leq K \leq 7,
	\label{E:MAINL2ESTIMATESVORTANDENTROPYGRADIENTBLOWUP}
		\\
	\hypersurfacecontrolVortVort_{\Ntop}(\timefunction,u), 
		\,
	\hypersurfacecontrolDivGradEnt_{\Ntop}(\timefunction,u) 
	& \leq C \initialsmall^2 |\timefunction|^{-17.1}, 
		\label{E:MAINTOPORDERENERGYESTIMATESMODIFIEDFLUIDVARIABLESBLOWUP} 
			\\
	\hypersurfacecontrolVortVort_{\Ntop - 1 -K}(\timefunction,u), 
		\,
	\hypersurfacecontrolDivGradEnt_{\Ntop - 1 - K}(\timefunction,u) 
	& \leq C \initialsmall^2 |\timefunction|^{-14.6 + 2K},  & & \text{for } 0 \leq K \leq 7.
		\label{E:MAINL2BELOWTOPORDERESTIMATESMODIFIEDFLUIDBLOWUP} 
\end{align}

Moreover, for $ 0 \leq N_1 \leq \Ntop - 8$  and $0 \leq N_2 \leq \Ntop - 9$, we have:
\begin{align}
	\hypersurfacecontrolVort_{N_1}(\timefunction,u), 
		\,
	\hypersurfacecontrolGradEnt_{N_1}(\timefunction,u) 
	& \leq C \initialsmall^2, 
		\label{E:MAINL2BELOWTOPORDERESTIMATESVORTANDENTROPYGRADIENTREGULAR} 
			\\
\hypersurfacecontrolVortVort_{N_2}(\timefunction,u), 
	\,
\hypersurfacecontrolDivGradEnt_{N_2}(\timefunction,u) & \leq C \initialsmall^2.
	\label{E:MAINL2BELOWTOPORDERESTIMATESMODIFIEDFLUIDVARIABLESREGULAR} 
\end{align}
\end{subequations}
\end{proposition}

In addition to the $L^2$ estimates of Prop.\,\ref{P:MAINHYPERSURFACEENERGYESTIMATESFORTRANSPORTVARIABLES}, 
which are estimates for the transport variables on constant-rough-time 
hypersurfaces and null hypersurfaces,
we also derive $L^2$ estimates for the transport variables
on the rough tori $\twoargroughtori{\timefunction,- \rightu}{\muxmulevelsetvalue}$.
We need these estimates because rough tori integrals
are featured in the main elliptic-hyperbolic integral identity
(see Prop.\,\ref{P:INTEGRALIDENTITYFORELLIPTICHYPERBOLICCURRENT})
that we use to control various top-order spacetime $L^2$ norms
of the $\vortrenormalized$ and $\GradEnt$.
The rough tori estimates of interest are provided by the next proposition.
Its proof is located in Sects.\,\ref{SS:PROOFOFBELOWTOPORDERROUGHTORIENERGYESTIMATES}
and \ref{SS:PROOFOFMAINTORIVORTGRADENTTOPORDERBLOWUP}.

\begin{proposition}[The main a priori $L^2$ estimates for the transport variables on the rough tori] 
	\label{P:ROUGHTORIENERGYESTIMATES} 
	Let
	$\toricontrolVort_{N}(\timefunction,u)$,
	$\toricontrolGradEnt_{N}(\timefunction,u)$,
	$\toricontrolVortVort_N$,
	and
	$\toricontrolDivGradEnt_N(\timefunction,u)$
	be the tori-$L^2$-controlling quantities for
	$\vortrenormalized$,
	$\GradEnt$,
	$\VortVort$,
	and $\DivGradEnt$
	as defined in \eqref{E:TORIVORTICITYL2CONTROLLINGQUANTITY}--\eqref{E:TORIENTROPYGRADIENTL2CONTROLLINGQUANTITY}
	and
	\eqref{E:TORIMODIFIEDVORTICITYL2CONTROLLINGQUANTITY}--\eqref{E:TORIDIVGRADENTL2CONTROLLINGQUANTITY}.
	Under the data assumptions of Sect.\,\ref{S:ASSUMPTIONSONTHEDATA},
	the parameter size-assumptions of Sect.\,\ref{SS:PARAMETERSIZEASSUMPTIONS},
	and the bootstrap assumptions of Sect.\,\ref{S:BOOTSTRAPEVERYTHINGEXCEPTENERGIES},
	there exists a constant $C > 0$
	such that the following estimates hold for 
	$(\timefunction,u) \in [\timefunction_0,\timefunctionboot) \times [- \rightu,\leftu]$: 
\begin{subequations}
\begin{align}
	\toricontrolVort_{\Ntop}(\timefunction,u), 
		\,
	\toricontrolGradEnt_{\Ntop} 
	& 
	\leq C \initialsmall^2 |\timefunction|^{-17.1}, 
		\label{E:MAINTORIL2VORTGRADENTTOPORDERBLOWUP} 
			\\
	\toricontrolVort_{\Ntop - 1 - K}(\timefunction,u), 
		\,
	\toricontrolGradEnt_{\Ntop - 1 - K} 
	& \leq C \initialsmall^2 |\timefunction|^{-14.6  + 2 K}, & & \text{for } 0 \leq K \leq 7, 
		\label{E:MAINTORIL2VORTGRADENTBELOWTOPORDERBLOWUP} 
		\\
\toricontrolVortVort_{\Ntop-1}(\timefunction,u),
	\,
\toricontrolDivGradEnt_{\Ntop-1}(\timefunction,u)
& \leq C \initialsmall^2 |\timefunction|^{-17.1}, 
	\label{E:MAINTORIL2MODIFIEDFLUIDVARIABLESTOPORDERBLOWUP}
		\\
\toricontrolVortVort_{\Ntop-2-K}(\timefunction,u),
	\,
\toricontrolDivGradEnt_{\Ntop-2-K}(\timefunction,u)
& \leq C \initialsmall^2 |\timefunction|^{-14.6 + 2K},  & & \text{for } 0 \leq K \leq 7,
	\label{E:MAINTORIL2MODIFIEDFLUIDVARIABLESBELOWTOPORDERBLOWUP}
\end{align}
\end{subequations}

\begin{subequations}
\begin{align}
\toricontrolVort_N(\timefunction,u), 
	\,
\toricontrolGradEnt_N(\timefunction,u) 
& \leq C \initialsmall^2,  
& & \text{for } 0 \leq N \leq \Ntop - 9, 
\label{E:MAINTORIL2VORTGRADENTBELOWTOPORDERREGULAR}
	\\
\toricontrolVortVort_N(\timefunction,u),
	\,
\toricontrolDivGradEnt_N(\timefunction,u) 
& \leq C \initialsmall^2,  
& & \text{for } 0 \leq N \leq \Ntop - 10.
\label{E:MAINTORIL2MODIFIEDFLUIDVARIABLESBELOWTOPORDERREGULAR}
\end{align}
\end{subequations}

\end{proposition}

\subsubsection{Statement of the a priori $L^2$ estimates for the acoustic geometry}
\label{SSS:STATEMENTOFAPRIORIL2ESTIMATESFORACOUSTICGEOMETRY}
The following proposition provides our main a priori energy estimates
for the acoustic geometry. Its proof is located in

\begin{proposition}[The main a priori estimates for the acoustic geometry along the rough foliations] 
\label{P:APRIORIL2ESTIMATESACOUSTICGEOMETRY}
	Under the data assumptions of Sect.\,\ref{S:ASSUMPTIONSONTHEDATA},
	the parameter size-assumptions of Sect.\,\ref{SS:PARAMETERSIZEASSUMPTIONS},
	and the bootstrap assumptions of Sect.\,\ref{S:BOOTSTRAPEVERYTHINGEXCEPTENERGIES},
	there exists a constant $C > 0$
	such that the following estimates hold for 
	$(\timefunction,u) \in [\timefunction_0,\timefunctionboot) \times [- \rightu,\leftu]$: 
\begin{subequations}
\begin{align}
	\left\| 
		\tander^{\Ntop} \mytr_{\gtorus} \upchi, 
			\, 
		\angLie_{\tander}^{\Ntop} \upchi  
	\right\|_{L^2\left( \hypthreearg{\timefunction}{[- \rightu,\leftu]}{\muxmulevelsetvalue}\right)} 
	& \leq C \initialsmall |\timefunction|^{-8.8},
		\label{E:MAINL2CHITOPORDERBLOWUP} \\
	\left\| 
		\tander^{N - 1} \mytr_{\gtorus} \upchi,
			\, 
		\angLie_{\tander}^{N - 1} \upchi, 
			\,
		\tandersmall^N\upmu, 
			\, 
		\tander^N \Lsmall^i  
	\right\|_{L^2\left( \hypthreearg{\timefunction}{[- \rightu,\leftu]}{\muxmulevelsetvalue}\right)} 
	& \leq C \initialsmall |\timefunction|^{- 7.3 + \Ntop - N}, \qquad \Ntop - 7 \leq N \leq \Ntop, 
		\label{E:MAINACOUSTGEOMETRYBELOWTOPORDERBLOWUP} 
		\\
 	\left\| 
		\tander^{\leq \Ntop - 9}  \mytr_{\gtorus} \upchi, 
			\, 
		\angLie_{\tander}^{\Ntop - 9} \upchi, 
			\,
		\tandersmall^{[1,\Ntop - 8]}\upmu,
			\, 
		\tander^{\leq \Ntop - 8} \Lsmall^i  
	\right\|_{L^2\left( \hypthreearg{\timefunction}{[- \rightu,\leftu]}{\muxmulevelsetvalue}\right)} 
	& \leq C \initialsmall.  
	\label{E:MAINACOUSTGEOMETRYBELOWTOPORDERREGULAR}
	\end{align}
	\end{subequations}
\end{proposition}

\subsection{Data estimates for the $L^2$-controlling quantities} 
\label{SS:DATAESTIMATESFORL2CONTROLLINGQUANTITIES}
In this section, we prove initial data estimates for 
the $L^2$-controlling quantities
$\totalcontrolwave_{[1,\Ntop]}$, $\hypersurfacecontrolVort_{\leq \Ntop }$, etc. 
We will use these data estimates in our proofs of
Props.\,\ref{P:APRIORIL2ESTIMATESWAVEVARIABLES}--\ref{P:APRIORIL2ESTIMATESACOUSTICGEOMETRY}.


\begin{lemma}[The $L^2$-controlling quantities are initially small] 
\label{L:ALLL2CONTROLLINGQUANTITIESINITIALLYSMALL}
The following estimates hold for 
$\timefunction \in [\timefunction_0,\timefunctionboot)$ and $u \in [-\rightu,\leftu]$:  
\begin{align} \label{E:WAVEL2CONTROLLINGINITIALLYSMALL}
\totalcontrolwave_{[1,\Ntop]}(\timefunction,-\rightu)
& \leq  
C \initialsmall^2, 
& 
\totalcontrolwave_{[1,\Ntop]}(\timefunction_0,u) 
& \leq C \initialsmall^2,
\end{align}

\begin{subequations}
\begin{align}
\hypersurfacecontrolVort_{\leq \Ntop }(\timefunction_0,u) & \leq  C \initialsmall^2, 
& \hypersurfacecontrolVort_{\leq \Ntop}(\timefunction,- \rightu) & \leq  C \initialsmall^2, 
	\label{E:VORTICITYL2CONTROLLINGINITIALLYSMALL}  
		\\
\hypersurfacecontrolGradEnt_{\leq \Ntop}(\timefunction_0,u) & \leq  C \initialsmall^2, 
& \hypersurfacecontrolGradEnt_{\leq \Ntop}(\timefunction,- \rightu) & \leq  C \initialsmall^2, 
	\label{E:GRADENTL2CONTROLLINGINITIALLYSMALL} 
\end{align}
\end{subequations}

\begin{subequations}
\begin{align} 
\hypersurfacecontrolVortVort_{\leq \Ntop}(\timefunction_0,u) & \leq  C \initialsmall^2,  
& \hypersurfacecontrolVortVort_{\leq \Ntop}(\timefunction,- \rightu) & \leq  C \initialsmall^2, 
	\label{E:MODIVIEDVORTL2CONTROLLINGINITIALLYSMALL}  
	\\
\hypersurfacecontrolDivGradEnt_{\leq \Ntop}(\timefunction_0,u)& \leq  C \initialsmall^2, 
& \hypersurfacecontrolDivGradEnt_{\leq \Ntop}(\timefunction,- \rightu) & \leq  C \initialsmall^2,
\label{E:MODIFIEDDIVGRADENTL2CONTROLLINGINITIALLYSMALL} 
\end{align}
\end{subequations}

\begin{subequations}
\begin{align} 
\toricontrolVort_{\leq \Ntop }(\timefunction_0,u) & \leq  C \initialsmall^2, 
	\label{E:VORTICITYTORIL2CONTROLLINGINITIALLYSMALL}  
		\\
\toricontrolGradEnt_{\leq \Ntop}(\timefunction_0,u)
& \leq C \initialsmall^2,
\label{E:GRADENTTORIL2CONTROLLINGINITIALLYSMALL} 
\end{align}
\end{subequations}

\begin{subequations}
\begin{align} 
\toricontrolVortVort_{\leq \Ntop-1}(\timefunction_0,u) & \leq  C \initialsmall^2, 
	\label{E:MODIFIEDVORTICITYTORIL2CONTROLLINGINITIALLYSMALL}  
		\\
\toricontrolDivGradEnt_{\leq \Ntop-1}(\timefunction_0,u)
& \leq C \initialsmall^2.
\label{E:DIVGRADENTTORIL2CONTROLLINGINITIALLYSMALL} 
\end{align}
\end{subequations}

\end{lemma}

\begin{proof}
	The estimates stated in \eqref{E:WAVEL2CONTROLLINGINITIALLYSMALL} are straightforward consequences of
	the data assumptions
	\eqref{E:TANGENTIALL2NORMSOFWAVEVARIABLESSMALLALONGINITIALROUGHHYPERSURFACE},
	\eqref{E:TRANSVERSALDERIVATIVEOFTANGENTIALL2NORMSOFWAVEVARIABLESSMALLALONGINITIALROUGHHYPERSURFACE},
	and \eqref{E:WAVESARESMALLONINITIALNULLHYERSURFACE},
	definitions 
	\eqref{E:WAVEENERGYDEF},
	\eqref{E:WAVENULLFLUXDEF},
	\eqref{E:COERCIVESPACETIMEWAVENERGYINTEGRAL},
	\eqref{E:WAVETOTALL2CONTROLLINGQUANTITY}, 
	and \eqref{E:WAVEPARTIALL2CONTROLLINGQUANTITY},
	the identities 
	\eqref{E:QLL}--\eqref{E:QNMULTIPLIER},
	\eqref{E:SMOOTHTORIGABEXPRESSION},
	\eqref{E:GEOP2TOCOMMUTATORS}--\eqref{E:GEOP3TOCOMMUTATORS},
	\eqref{E:ROUGHANGULARPARTIALDERIVATIVESINTERMSOFGOODGEOMETRICPARTIALDERIVATIVES},
	\eqref{E:SIZEOFRTRANS}, 
	\eqref{E:IDENTITYFORRTRANSNORMSMALLFACTORSQUARED},
	and
	\eqref{E:GTORUSROUGHCOMPONENTS}--\eqref{E:ROUGHTORUSMETRICCOMPONENTSANDTHEINVERSECOMPONENTSRELATION},
	Lemma~\ref{L:SCHEMATICSTRUCTUREOFVARIOUSTENSORSINTERMSOFCONTROLVARS},
	the bootstrap assumptions,
	and the estimates 
	\eqref{E:LDERIVATIVEOFROUGHTIMEFUNCTIONISAPPROXIMATELYUNITY}
	and
	\eqref{E:CLOSEDVERSIONC21BOUNDFORCHOVROUGHTOGEO}.
	
	The estimates \eqref{E:VORTICITYL2CONTROLLINGINITIALLYSMALL}--\eqref{E:GRADENTL2CONTROLLINGINITIALLYSMALL} 
	follow from applying similar reasoning 
	based on definitions 
	\eqref{E:TRANSPORTENERGYDEF}--\eqref{E:TRANSPORTNULLFLUXDEF}
	and
	\eqref{E:VORTICITYL2CONTROLLINGQUANTITY}--\eqref{E:ENTROPYGRADIENTL2CONTROLLINGQUANTITY}
	and
	the data assumptions
	\eqref{E:TANGENTIALL2NORMSOFTRANSPORTVARIABLESSMALLALONGINITIALROUGHHYPERSURFACE} 
	and \eqref{E:VORTICITYANDENTROPYGRADIENTARESMALLONINITIALNULLHYERSURFACE}.
	
	The estimates 
	\eqref{E:MODIVIEDVORTL2CONTROLLINGINITIALLYSMALL}--\eqref{E:MODIFIEDDIVGRADENTL2CONTROLLINGINITIALLYSMALL} 
	follow from applying similar reasoning 
	based on definitions 
	\eqref{E:TRANSPORTENERGYDEF}--\eqref{E:TRANSPORTNULLFLUXDEF}
	and
	\eqref{E:MODIFIEDVORTICITYVORTICITYL2CONTROLLINGQUANTITY}--\eqref{E:MODIFIEDDIVGRADENTL2CONTROLLINGQUANTITY}
	and
	the data assumptions
	\eqref{E:TANGENTIALL2NORMSOFMODIFIEDFLUIDVARIABLESSMALLALONGINITIALROUGHHYPERSURFACE} 
	and \eqref{E:MODIFIEDFLUIDVARIABLESARESMALLONINITIALNULLHYERSURFACE}.
		
	The estimates 
	\eqref{E:VORTICITYTORIL2CONTROLLINGINITIALLYSMALL}--\eqref{E:GRADENTTORIL2CONTROLLINGINITIALLYSMALL}  
	follow from
	definitions~\eqref{E:TORIVORTICITYL2CONTROLLINGQUANTITY}--\eqref{E:TORIENTROPYGRADIENTL2CONTROLLINGQUANTITY}	
	and the data assumptions
	\eqref{E:SMALLDATAOFVORTICITYANDENTORPYGRADIENTONINITIALROUGHTORI}.
	
	The estimates
	\eqref{E:MODIFIEDVORTICITYTORIL2CONTROLLINGINITIALLYSMALL}--\eqref{E:DIVGRADENTTORIL2CONTROLLINGINITIALLYSMALL} 
	follow from
	definitions~\eqref{E:TORIMODIFIEDVORTICITYL2CONTROLLINGQUANTITY}--\eqref{E:TORIDIVGRADENTL2CONTROLLINGQUANTITY}
	and the data assumptions 
	\eqref{E:SMALLDATAOFMODIFIEDFLUIDVARIABLESONINITIALROUGHTORI}.
\end{proof}

\subsection{Bootstrap assumptions for the $\totalcontrolwave_N$}
\label{SS:BOOTSTRAPASSUMPTIONSFORTHEWAVEENERGIES}
In proving Props.\,\ref{P:APRIORIL2ESTIMATESWAVEVARIABLES}--\ref{P:APRIORIL2ESTIMATESACOUSTICGEOMETRY},
we find it convenient to make bootstrap assumptions for the $L^2$-controlling quantities
for the wave variables.
Specifically, with
$\totalcontrolwave_{N}(\timefunction,u)$ 
(see definition \eqref{E:WAVETOTALL2CONTROLLINGQUANTITY})
denoting the $L^2$-controlling quantity for the wave variables $\wavearray$,
we assume that the following bootstrap assumptions hold for 
$(\timefunction,u) \in [\timefunction_0,\timefunctionboot) \times [- \rightu,\leftu]$, 
where $\fundbootsmall$ is the bootstrap parameter from Sect.\,\ref{SSS:FUNDAMENTALQUANTITATIVE}:
\begin{subequations} 
\begin{align}
\totalcontrolwave_{\Ntop - K}(\timefunction,u) 
	&  
	\leq \fundbootsmall |\timefunction|^{-15.6 + 2K}, 
\qquad 0 \leq K \leq 7 
\label{E:MAINWAVEENERGYBOOTSTRAPBLOWUP} 
	\\
\totalcontrolwave_{[1,\Ntop - 8]}(\timefunction,u) 
& 
\leq \fundbootsmall. 
	\label{E:MAINWAVEENERGYBOOTSTRAPREGULAR} 
\end{align}
\end{subequations}

\begin{remark}[The wave variable energy estimates improve the bootstrap assumptions]
	\label{R:WAVEENERGYESTIMATESYIELDSTRICTIMPROVEMENTSOVERWAVEENERGYBOOTSTRAPASSUMPTIONS}
	Note that when $\initialsmall$ is sufficiently small,
	the estimates of Prop.\,\ref{P:APRIORIL2ESTIMATESWAVEVARIABLES} 
	yield strict improvements over the bootstrap assumptions
	\eqref{E:MAINWAVEENERGYBOOTSTRAPBLOWUP}--\eqref{E:MAINWAVEENERGYBOOTSTRAPREGULAR}.
\end{remark}

\section{Preliminary below-top-order $L^2$ estimates for the acoustic geometry and a derivative-losing estimate} 
\label{S:PRELIMINARYL2ESTIMATESFORBELOWTOPORDERDERIVATIVESOFACOUSTICGEOMETRYANDDERIVATIVELOSING}
In this short section, 
we derive preliminary $L^2$ estimates for the below-top-order derivatives of the
eikonal function quantities $\upmu$, $\Lunit^i$, $\upchi$, and $\mytr_{\gtorus}\upchi$. 
We state the bounds in terms of the wave $L^2$-controlling quantities $\hypersurfacecontrolwave_{[1,N]}(\timefunction,u)$
and the initial data size parameter $\initialsmall$.
We also derive related $L^2$ estimates for top-order derivatives of $\upchi$ and $\mytr_{\gtorus} \upchi$ in the case that one 
$\Lunit$-differentiation is involved. We provide the main estimates in Lemma~\ref{L:EIKONALFUNCTIONBELOWTOPORDERL2ESTIMATES}.
The estimates are rather straightforward 
consequences of the transport inequalities provided by
Prop.\,\ref{P:POINTWISETRANSPORTINEQUALITIESFOREIKFUNCTIONQUANTITIES}
and
Lemma~\ref{L:TRANSPORTESTIMATESFORROUGHLFEQUALSSOURCE}. 
In Cor.\,\ref{C:NONSINGULARL2ESTIMATESFORWAVEVARIABLESTHATLOSEONEDERIVATIVE},
we derive $L^2$ estimates for $\wavearray$ that do not involve
any singular powers of $|\timefunction|^{-1}$. 

Most of the $L^2$ estimates we derive in this section lose one derivative.
In Prop.\,\ref{P:TOPORDERL2ESTIMATEMUCHI}, we will derive complementary 
$L^2$ estimates for the top-order derivatives of $\upchi$ and $\mytr_{\gtorus}\upchi$.
Those estimates are much harder to prove because we cannot afford to lose 
any derivatives, which forces us to rely on the modified quantities from Sect.\,\ref{S:CONSTRUCTIONOFMODIFIEDQUANTITIES}
and elliptic estimates for the top-order derivatives of $\upchi$.

\subsection{Preliminary below-top-order $L^2$ estimates for the eikonal function quantities}
\label{SS:PRELIMINARYL2ESTIMATESFORBELOWTOPORDERDERIVATIVESOFACOUSTICGEOMETRY}

\begin{lemma}[Preliminary below-top-order $L^2$ estimates for the eikonal function quantities] 
\label{L:EIKONALFUNCTIONBELOWTOPORDERL2ESTIMATES}
Let $N \leq \Ntop$. Then the following estimates hold for 
$(\timefunction,u) \in [\timefunction_0,\timefunctionboot) \times [-\rightu,\leftu]$:
\begin{subequations} 
\begin{align} 
\begin{pmatrix} 
\left\| \Lunit \tandersmall^{[1,N]} \upmu \right \|_{L^2\left(\hypthreearg{\timefunction}{[-\rightu,u]}{\muxmulevelsetvalue}\right)} 
	\\
\left\| \Lunit \tander^{\leq N} \Lsmall^i \right \|_{L^2\left(\hypthreearg{\timefunction}{[-\rightu,u]}{\muxmulevelsetvalue}\right)} 
	\\
\left\| \Lunit \tander^{\leq N-1} \mytr_{\gtorus} \upchi \right \|_{L^2\left(\hypthreearg{\timefunction}{[\rightu,u]}{\muxmulevelsetvalue}\right)} 
	\\
\left\| \angLie_{\Lunit} \angLie_{ \tander}^{\leq N-1} \upchi \right \|_{L^2\left(\hypthreearg{\timefunction}{[-\rightu,u]}{\muxmulevelsetvalue}\right)} \\
\left\| \Lunit \comder^{\leq N;1} \Lsmall^i  \right \|_{L^2\left(\hypthreearg{\timefunction}{[-\rightu,u]}{\muxmulevelsetvalue}\right)} 
	\\
\left\| \Lunit \comder^{\leq N-1;1} \mytr_{\gtorus}\upchi  \right \|_{L^2\left(\hypthreearg{\timefunction}{[-\rightu,u]}{\muxmulevelsetvalue}\right)} 
	\\
\left\| \angLie_{\Lunit} \angLie_{\comder}^{\leq N-1;1} \upchi \right \|_{L^2\left(\hypthreearg{\timefunction}{[-\rightu,u]}{\muxmulevelsetvalue}\right)} 
\\
\end{pmatrix} 
& \lesssim \initialsmall 
+ 
\frac{\hypersurfacecontrolwave_{[1,N]}^{1/2} (\timefunction,u)}{|\timefunction|^{1/2}},
	\label{E:PRELIMINARYEIKONALWITHL} 
	\\
\begin{pmatrix}
\left\|  \tandersmall^{[1,N]} \upmu \right \|_{L^2\left(\hypthreearg{\timefunction}{[-\rightu,u]}{\muxmulevelsetvalue}\right)} 
	\\
\left\| \tander^{[1,N]} \Lsmall^i\right \|_{L^2\left(\hypthreearg{\timefunction}{[-\rightu,u]}{\muxmulevelsetvalue}\right)} 	
	\\
\left\| \tander^{\leq N-1} \mytr_{\gtorus} \upchi \right \|_{L^2\left(\hypthreearg{\timefunction}{[-\rightu,u]}{\muxmulevelsetvalue}\right)} 
	\\
\left\|  \angLie_{\tander}^{\leq N-1} \upchi \right \|_{L^2\left(\hypthreearg{\timefunction}{[-\rightu,u]}{\muxmulevelsetvalue}\right)} 
	\\
\left\| \comdersmall^{[1,N];1} \Lsmall^i \right\|_{L^2\left(\hypthreearg{\timefunction}{[-\rightu,u]}{\muxmulevelsetvalue}\right)} 
	\\
\left\| \comder^{\leq N-1;1} \mytr_{\gtorus}\upchi \right\|_{L^2\left(\hypthreearg{\timefunction}{[-\rightu,u]}{\muxmulevelsetvalue}\right)} 
	\\
\left\| \angLie_{\comder}^{\leq N-1;1} \upchi \right\|_{L^2\left(\hypthreearg{\timefunction}{[-\rightu,u]}{\muxmulevelsetvalue}\right)}
\end{pmatrix} 
& \lesssim 
\initialsmall 
+ 
\int_{\timefunction' = \timefunction_0}^{\timefunction} 
	\frac{\hypersurfacecontrolwave_{[1,N]}^{1/2} (\timefunction',u)}{|\timefunction'|^{1/2}}  
\, \mathrm{d} \timefunction'. 
\label{E:PRELIMINARYEIKONALWITHOUTL}
\end{align}
\end{subequations}
\end{lemma}

\begin{proof}
We fix $u \in [-\rightu,\leftu]$ and define the following
functions for $\timefunction \in [\timefunction_0,\timefunctionboot)$:
\begin{align}
	q_N(\timefunction) 
	& 
	\eqdef 
	\sum_{i = 1}^3 
	\left\| \tander^{[1,N]} \Lsmall^i  \right\|_{L^2\left(\hypthreearg{\timefunction}{[-\rightu,u]}{\muxmulevelsetvalue}\right)} 
	+ 
	\left\| \tander^{\leq N-1} \mytr_{\gtorus}\upchi \right\|_{L^2\left(\hypthreearg{\timefunction}{[-\rightu,u]}{\muxmulevelsetvalue}\right)} 
	+ 
	\left\| \angLie_{\tander}^{\leq N-1}\upchi \right\|_{L^2\left(\hypthreearg{\timefunction}{[-\rightu,u]}{\muxmulevelsetvalue}\right)}, 
		\label{E:QNDEFPROOFOFPRELIMINARYEIKONALL2}
		\\
	\begin{split} \label{E:PNDEFPROOFOFPRELIMINARYEIKONALL2}
	p_N(\timefunction) 
	& \eqdef 
	\sum_{i = 1}^3 \left\| \comdersmall^{[1,N];1} \Lsmall^i \right\|_{L^2\left(\hypthreearg{\timefunction}{[-\rightu,u]}{\muxmulevelsetvalue}\right)} 
	+ 
	\left\| \comder^{\leq N-1;1} \mytr_{\gtorus}\upchi \right\|_{L^2\left(\hypthreearg{\timefunction}{[-\rightu,u]}{\muxmulevelsetvalue}\right)} 
		 \\
	&
	+ 
	\left\| \angLie_{\comder}^{\leq N-1;1}\upchi \right\|_{L^2\left(\hypthreearg{\timefunction}{[-\rightu,u]}{\muxmulevelsetvalue}\right)} 
	+ 
	\left\| \tandersmall^{[1,N]}\upmu \right\|_{L^2\left(\hypthreearg{\timefunction}{[-\rightu,u]}{\muxmulevelsetvalue}\right)}.
	\end{split}
\end{align}
Next, we record the following pointwise estimate: 
$\left|\Lunit\left(|\angLie_{\tander}^{N-1}\upchi|_{\gtorus} \right) \right| 
\lesssim 
|\angLie_{\Lunit} \angLie_{\tander}^{N-1} \upchi|_{\gtorus}  
+ 
\fundbootsmall |\angLie_{\tander}^{N-1} \upchi|_{\gtorus} $, 
which follows from the Leibniz rule for $\ell_{t,u}$-projected Lie derivatives,
the identity $\angLie_{\Lunit} \gtorus^{-1} = - 2 \upchi^{\# \#}$,
and the estimate 
$|\upchi|_{\gtorus} \lesssim \fundbootsmall$, 
which follows from \eqref{E:CHIEXPRESSSIONINTERMSOFDERIVATIVESOFLUNITI}, 
Lemma~\ref{L:SCHEMATICSTRUCTUREOFVARIOUSTENSORSINTERMSOFCONTROLVARS},
the bootstrap assumptions,
and Prop.\,\ref{P:IMPROVEMENTOFAUXILIARYBOOTSTRAP}.
Multiplying 
\eqref{E:LUNITTANGENTIALDERIVATIVESOFLUNITIPOINTWISE}--\eqref{E:ANGLIELTANGENTIALCHIPOINTWISE} by 
$\frac{1}{\Lunit \timefunctionarg{\muxmulevelsetvalue}}$ 
and using \eqref{E:L2ONROUGHCONSTANTTIMEHYPERSURFACESTRANSPORTROUGHLFESTIMATE},
Cor.\,\ref{C:IMPROVEAUX},
as well as the estimates 
\eqref{E:CLOSEDVERSIONLUNITROUGHTTIMEFUNCTION},
\eqref{E:MINVALUEOFMUONFOLIATION}, 
\eqref{E:COERCIVENESSOFHYPERSURFACECONTROLWAVE}, 
\eqref{E:L2ESTIMATESFORWAVEVARIABLESONROUGHHYPERSURFACELOSSOFONEDERIVATIVE},
and the pointwise estimate recorded above,
we deduce:
\begin{align} \label{E:QNGRONWALLREADYPROOFOFPRELIMINARYEIKONALL2}
q_N(\timefunction) 
&
\lesssim 
q_N(\timefunction_0) 
+
\initialsmall
+ 
\fundbootsmall
\int_{\timefunction' = \timefunction_0}^{\timefunction} 
	q_N(\timefunction') 
\, \mathrm{d} \timefunction' 
+ 
\int_{\timefunction' = \timefunction_0}^{\timefunction} 
	\frac{1}{|\timefunction'|^{1/2}} 
	\hypersurfacecontrolwave_{[1,N]}^{1/2}(\timefunction',u) 
\, \mathrm{d} \timefunction'.
\end{align}
The $L^2$ assumptions on the data stated in Sect.\,\ref{SSS:QUANTITATIVEASSUMPTIONSONDATAAWAYFROMSYMMETRY} imply that
$q_N(\timefunction_0) \lesssim \initialsmall$. Inserting this estimate into RHS~\eqref{E:QNGRONWALLREADYPROOFOFPRELIMINARYEIKONALL2}
and applying Gr\"{o}nwall's inequality, we find that 
$
q_N(\timefunction) 
\lesssim
\initialsmall
+
\int_{\timefunction' = \timefunction_0}^{\timefunction} 
	\frac{1}{|\timefunction'|^{1/2}} 
	\hypersurfacecontrolwave_{[1,N]}^{1/2}(\timefunction',u) 
\, \mathrm{d} \timefunction'
$,
which yields
\eqref{E:PRELIMINARYEIKONALWITHOUTL} for 
$\tander^{[1,N]}\Lsmall^i$, $\tander^{\leq N-1} \mytr_{\gtorus}\upchi$, and $\angLie_{\tander}^{\leq N-1}\upchi$. 

We now prove \eqref{E:PRELIMINARYEIKONALWITHOUTL} for 
$\comdersmall^{[1,N];1} \Lsmall^i$, 
$\comder^{\leq N-1;1} \mytr_{\gtorus} \upchi$, 
$\angLie_{\comder}^{\leq N-1;1} \upchi$, 
and $\tandersmall^{[1,N]}\upmu$. 
We begin by examining the first term on RHS~\eqref{E:LTANGENTIALMUPOINTWISE}. 
Using the commutator estimates 
\eqref{E:COMMUTATOROFTANGENTIALANDTANGENTIALCOMMUTATORS}--\eqref{E:COMMUTATOROFMUXANDTANGENTIALCOMMUTATORS},
the bootstrap assumptions,
and the estimates of Prop.\,\ref{P:IMPROVEMENTOFAUXILIARYBOOTSTRAP}, 
we deduce:
\begin{align} \label{E:PRELIMINARYEIKONALWITHOUTLINTERMEDIATE}
		|\comdersmall^{[1,N+1];1} \wavearray| 
		& 
		\lesssim 
		|\muX \tander^{[1,N]} \wavearray| 
		+ 
		|\tander^{[1,N+1]} \wavearray| 
		+ 
		\fundbootsmall 
		|\comdersmall^{[1,N];1}\controlvars| 
		+ 
		\fundbootsmall |\tandersmall^{[1,N]} \badcontrolvars|.
\end{align}
Arguing as in the proof of \eqref{E:QNGRONWALLREADYPROOFOFPRELIMINARYEIKONALL2},
but using \eqref{E:LTANGENTIALMUPOINTWISE} and \eqref{E:LZLSMALLPOINTWISE}--\eqref{E:ANGLIELZCHIPOINTWISE}
in place of \eqref{E:LUNITTANGENTIALDERIVATIVESOFLUNITIPOINTWISE}--\eqref{E:ANGLIELTANGENTIALCHIPOINTWISE},
and using \eqref{E:PRELIMINARYEIKONALWITHOUTLINTERMEDIATE},
\eqref{E:COERCIVENESSOFHYPERSURFACECONTROLWAVE},
and \eqref{E:L2ESTIMATESFORWAVEVARIABLESONROUGHHYPERSURFACELOSSOFONEDERIVATIVE},
we find that:
\begin{align} \label{E:PRELIMINARYEIKONALWITHOUTLINTERMEDIATE2}
	p_N(\timefunction) 
	& \lesssim 
	p_N(\timefunction_0) 
	+
	\initialsmall
	+ 
	\int_{\timefunction' = \timefunction_0}^{\timefunction}  p_N(\timefunction') \, \mathrm{d} \timefunction' 
	+ 
	\int_{\timefunction' = \timefunction_0}^{\timefunction} 
		\frac{1}{|\timefunction'|^{1/2}} \hypersurfacecontrolwave_{[1,N]}^{1/2}(\timefunction',u) \, 
	\mathrm{d} \timefunction'.
	\end{align}
As above,
the $L^2$ assumptions on the data stated in Sect.\,\ref{SSS:QUANTITATIVEASSUMPTIONSONDATAAWAYFROMSYMMETRY} imply that
$p(\timefunction_0) \lesssim \initialsmall$,
and we can use Gr\"{o}nwall's inequality to deduce
$
q_N(\timefunction) 
\lesssim
\initialsmall
+
\int_{\timefunction' = \timefunction_0}^{\timefunction} 
	\frac{1}{|\timefunction'|^{1/2}} 
	\hypersurfacecontrolwave_{[1,N]}^{1/2}(\timefunction',u) 
\, \mathrm{d} \timefunction'
$,
thereby concluding
\eqref{E:PRELIMINARYEIKONALWITHOUTL} for 
$\comdersmall^{[1,N];1} \Lsmall^i$, 
$\comder^{\leq N-1;1} \mytr_{\gtorus} \upchi$,
$\angLie_{\comder}^{\leq N-1;1} \upchi$, 
and $\tandersmall^{[1,N]}\upmu$.

We now prove \eqref{E:PRELIMINARYEIKONALWITHL}. 
Taking the norm $\| \cdot \|_{L^2\left(\hypthreearg{\timefunction}{[-\rightu,u]}{\muxmulevelsetvalue}\right)}$ of inequalities 
\eqref{E:LTANGENTIALMUPOINTWISE}--\eqref{E:ANGLIELZCHIPOINTWISE}, 
and arguing as above using the already proven estimate \eqref{E:PRELIMINARYEIKONALWITHOUTL}, 
we obtain the desired result. We clarify that this argument generates the time integrals
$	
\int_{\timefunction' = \timefunction_0}^{\timefunction} 
	\frac{1}{|\timefunction'|^{1/2}} \hypersurfacecontrolwave_{[1,N]}^{1/2}(\timefunction',u)
	\, \mathrm{d} \timefunction'$,
which we bound by 
$\lesssim \hypersurfacecontrolwave_{[1,N]}^{1/2}(\timefunction,u)
\lesssim \mbox{RHS~\eqref{E:PRELIMINARYEIKONALWITHL}}
$
by exploiting the monotonicity of $\hypersurfacecontrolwave_{[1,N]}(\timefunction,u)$ with respect to
its arguments.
\end{proof}

\subsection{$L^2$ estimates for $\wavearray$ that lose one derivative}
\label{SS:WAVEVARIABLEL2ESTIAMTESTHATLOSEONEDERIVATIVE}

\begin{corollary}[Non-singular $L^2$ estimates for $\wavearray$ that lose one derivative]
\label{C:NONSINGULARL2ESTIMATESFORWAVEVARIABLESTHATLOSEONEDERIVATIVE}
Let $1 \leq N \leq \Ntop$. The following estimates hold for 
$(\timefunction,u) \in [\timefunction_0,\timefunctionboot) \times [-\rightu,\leftu]$:
\begin{align} \label{E:NONDEGENERATEBUTDERIVATIVELOSINGL2ESTIMATEPSI}
	\left\| \comdersmall^{N;1} \wavearray \right\|_{L^2\left(\hypthreearg{\timefunction}{[-\rightu,u]}{\muxmulevelsetvalue}\right)} 
	& 
	\lesssim \initialsmall 
	+ 
	\hypersurfacecontrolwave_{[1,N]}^{1/2}(\timefunction,u).
\end{align}
\end{corollary}

\begin{proof}
Using \eqref{E:PRELIMINARYEIKONALWITHOUTLINTERMEDIATE} with $N+1$ replaced by $N$, we deduce the pointwise estimate
\begin{align} \label{E:MAINPOINTWISESEMINORMNEEDEDTOPROVENONDEGENERATEBUTDERIVATIVELOSINGL2ESTIMATEPSI}
\left|\comdersmall^{N;1} \wavearray \right| 
&
\lesssim 
\left|\muX \tander^{[1,N-1]} \wavearray \right| 
+ 
\left|\tander^{[1,N]} \wavearray \right| 
+ 
\fundbootsmall 
\left|\comdersmall^{[1,N-1];1}\controlvars \right| 
+ 
\fundbootsmall 
\left|\tandersmall^{[1,N-1]} \badcontrolvars \right|,
\end{align} 
where when $N = 1$, we must have $\comdersmall^{1;1}  = \tander$ on 
LHS~\eqref{E:MAINPOINTWISESEMINORMNEEDEDTOPROVENONDEGENERATEBUTDERIVATIVELOSINGL2ESTIMATEPSI}
and then only the second term on RHS~\eqref{E:MAINPOINTWISESEMINORMNEEDEDTOPROVENONDEGENERATEBUTDERIVATIVELOSINGL2ESTIMATEPSI} is present.
Taking the norm 
$\| \cdot \|_{L^2\left(\hypthreearg{\timefunction}{[-\rightu,u]}{\muxmulevelsetvalue}\right)}$ of 
\eqref{E:MAINPOINTWISESEMINORMNEEDEDTOPROVENONDEGENERATEBUTDERIVATIVELOSINGL2ESTIMATEPSI}
and using
\eqref{E:COERCIVENESSOFHYPERSURFACECONTROLWAVE},
\eqref{E:L2ESTIMATESFORWAVEVARIABLESONROUGHHYPERSURFACELOSSOFONEDERIVATIVE},
\eqref{E:PRELIMINARYEIKONALWITHOUTL},
and the estimate
$	
\int_{\timefunction' = \timefunction_0}^{\timefunction} 
	\frac{1}{|\timefunction'|^{1/2}} \hypersurfacecontrolwave_{[1,N]}^{1/2}(\timefunction',u)
	\, \mathrm{d} \timefunction'
\lesssim \hypersurfacecontrolwave_{[1,N]}^{1/2}(\timefunction,u)$,
which follows from the fact that
$\hypersurfacecontrolwave_{[1,N]}(\timefunction,u)$ is increasing in its arguments,
we conclude \eqref{E:NONDEGENERATEBUTDERIVATIVELOSINGL2ESTIMATEPSI}.

\end{proof}

\section{Below-top-order hyperbolic $L^2$ estimates for the specific vorticity and entropy gradient}
\label{S:BELOWTOPORDERHYPERBOLICL2ESTIMATESFORSPECIFICVORTICITYANDENTROPYGRADIENT}
In this short section, 
we prove the below-top-order $L^2$ estimates for the specific vorticity and entropy gradient. 
Specifically, we prove
\eqref{E:MAINL2ESTIMATESVORTANDENTROPYGRADIENTBLOWUP},
\eqref{E:MAINL2BELOWTOPORDERESTIMATESMODIFIEDFLUIDBLOWUP},
\eqref{E:MAINL2BELOWTOPORDERESTIMATESVORTANDENTROPYGRADIENTREGULAR}--\eqref{E:MAINL2BELOWTOPORDERESTIMATESMODIFIEDFLUIDVARIABLESREGULAR},
\eqref{E:MAINTORIL2VORTGRADENTBELOWTOPORDERBLOWUP}, 
\eqref{E:MAINTORIL2MODIFIEDFLUIDVARIABLESBELOWTOPORDERBLOWUP},
and
\eqref{E:MAINTORIL2VORTGRADENTBELOWTOPORDERREGULAR}--\eqref{E:MAINTORIL2MODIFIEDFLUIDVARIABLESBELOWTOPORDERREGULAR}.
We derive a preliminary energy integral inequality 
in Sect.\,\ref{SS:BELOWTOPORDERHYPERBOLICL2ESTIMATESFORSPECIFICVORTICITYANDENTROPYGRADIENT},
and we prove the final estimates in 
Sects.\,\ref{SS:PROOFOFBELOWTOPORDERTRANSPORTENERGYESTIMATES}--\ref{SS:PROOFOFBELOWTOPORDERROUGHTORIENERGYESTIMATES}.
These estimates are relatively straightforward consequences of
the transport energy identity \eqref{E:ENERGYIDENTITYFORBTRANSPORTEQUATIONS}
and various pointwise estimates we have already derived,
including the ones provided by Prop.\,\ref{P:POINTWISESTIAMTESFORALLTHETRANSPORTVARIABLES}.
In Sect.\,\ref{S:TOPORDERELIPTICHYPERBOLICL2ESTIMATESFORSPECIFICVORTICITYANDENTROPYGRADIENT},
we will prove the top-order estimate \eqref{E:MAINTOPORDERENERGYESTIMATESMODIFIEDFLUIDVARIABLESBLOWUP}
and the related estimates \eqref{E:MAINTORIL2VORTGRADENTTOPORDERBLOWUP} and \eqref{E:MAINTORIL2MODIFIEDFLUIDVARIABLESTOPORDERBLOWUP}.
The proofs of these estimates are much more difficult because they rely on
the intricate elliptic-hyperbolic integral identity \eqref{E:INTEGRALIDENTITYFORELLIPTICHYPERBOLICCURRENT}.

In Sect.\,\ref{S:WAVEANDACOUSTICGEOMETRYAPRIORIESTIMATES}, we will use the estimates for 
$\vortrenormalized$,
$\GradEnt$,
$\VortVort$, and $\DivGradEnt$
that we derive in this section in our proof of the wave a priori estimates, 
which we stated as Prop.\,\ref{P:APRIORIL2ESTIMATESWAVEVARIABLES}.
Hence, we highlight that for the logic of the paper, it is important that 
\textbf{the estimates we derive in this section do not rely on the wave estimates of Prop.\,\ref{P:APRIORIL2ESTIMATESWAVEVARIABLES}};
our proofs of 
\eqref{E:MAINL2ESTIMATESVORTANDENTROPYGRADIENTBLOWUP},
\eqref{E:MAINL2BELOWTOPORDERESTIMATESMODIFIEDFLUIDBLOWUP}--\eqref{E:MAINL2BELOWTOPORDERESTIMATESMODIFIEDFLUIDVARIABLESREGULAR}
and
\eqref{E:MAINTORIL2VORTGRADENTBELOWTOPORDERBLOWUP}--\eqref{E:MAINTORIL2VORTGRADENTBELOWTOPORDERREGULAR}
instead rely on the bootstrap assumptions
\eqref{E:MAINWAVEENERGYBOOTSTRAPBLOWUP}--\eqref{E:MAINWAVEENERGYBOOTSTRAPREGULAR},
for the wave energies, which are \emph{weaker} than the estimates that we derive in Prop.\,\ref{P:APRIORIL2ESTIMATESWAVEVARIABLES}.

\subsection{Integral inequalities for the below-top-order vorticity- and entropy gradient-controlling quantities}
\label{SS:BELOWTOPORDERHYPERBOLICL2ESTIMATESFORSPECIFICVORTICITYANDENTROPYGRADIENT}
We begin with the following preliminary lemma,
which provides integral inequalities
for the below-top-order vorticity- and entropy gradient-controlling quantities.

\begin{lemma}[Integral inequalities for the below-top-order vorticity- and entropy gradient-controlling quantities] 
\label{L:BELOWTOPORDERENERGYINTEGRALINEQUALITIESSPECIFICVORTICITYANDENTROPYGRADIENT}
	Let $0 \leq N \leq \Ntop$. 
	Then the following integral inequalities hold for
	$(\timefunction,u) \in [\timefunction_0,\timefunctionboot) \times [-\rightu,\leftu]$:
	\begin{align}   \label{E:BELOWTOPORDERENERGYINTEGRALINEQUALITIESSPECIFICVORTICITYANDENTROPYGRADIENT}
		\hypersurfacecontrolVort_{\leq N}(\timefunction,u)
		+
		\hypersurfacecontrolGradEnt_{\leq N}(\timefunction,u)
		& 
		\lesssim
		\initialsmall^2
		+
		\int_{u' = - \rightu}^u
			\left\lbrace
				\hypersurfacecontrolVort_{\leq N}(\timefunction,u')
				+
				\hypersurfacecontrolGradEnt_{\leq N}(\timefunction,u')
			\right\rbrace
		\, \mathrm{d} u'
		+
		\varepsilon^2 
		\int_{\timefunction' = \timefunction_0}^{\timefunction} 
			\hypersurfacecontrolwave_{[1,N]}(\timefunction',u) 
		\, \mathrm{d} \timefunction'.
	\end{align}
	
	Moreover, the following estimates hold for $0 \leq N \leq \Ntop - 1$:
	\begin{align} 
	\begin{split} \label{E:BELOWTOPORDERENERGYINTEGRALINEQUALITIESFORMODIFIEDFLUID}
		\hypersurfacecontrolVortVort_{\leq N}(\timefunction,u)
		+
		\hypersurfacecontrolDivGradEnt_{\leq N}(\timefunction,u)
		& \lesssim
		\initialsmall^2
		+
		\int_{u' = - \rightu}^u
			\left\lbrace
				\hypersurfacecontrolVortVort_{\leq N}(\timefunction,u')
				+
				\hypersurfacecontrolDivGradEnt_{\leq N}(\timefunction,u)
			\right\rbrace
		\, \mathrm{d} u'
			\\
	& \ \
		+
		\int_{u' = - \rightu}^u 
			\left\lbrace
				\hypersurfacecontrolVort_{\leq N+1}(\timefunction,u')
				+
				\hypersurfacecontrolGradEnt_{\leq N+1}(\timefunction,u')
			\right\rbrace
		\, \mathrm{d} u'	
			\\
	& \ \
		+
		\varepsilon^2 
		\int_{\timefunction' = \timefunction_0}^{\timefunction} 
			\hypersurfacecontrolwave_{[1,N]}(\timefunction',u) 
		\, \mathrm{d} \timefunction'
		+
		\varepsilon^2 
		\int_{u' = - \rightu}^u 
			\hypersurfacecontrolwave_{[1,N]}(\timefunction,u')
		\, \mathrm{d} u' 
		+ 
		\fundbootsmall^2 \spacetimeintegralcontrolwave_{[1,N]}(\timefunction,u).
	\end{split}
	\end{align}
	
\end{lemma}

\begin{proof}
	We first prove \eqref{E:BELOWTOPORDERENERGYINTEGRALINEQUALITIESSPECIFICVORTICITYANDENTROPYGRADIENT}.
	For $0 \leq N \leq \Ntop$, we consider the transport energy identity \eqref{E:ENERGYIDENTITYFORBTRANSPORTEQUATIONS}
	with $(\tander^N \vortrenormalized,\tander^N \GradEnt)$ in the role of $f$.
	We use the bootstrap assumptions and Prop.\,\ref{P:IMPROVEMENTOFAUXILIARYBOOTSTRAP} to deduce the bound
	$|\Lunit \upmu + \upmu \mytr_{\gtorus} \angk| \lesssim 1$ for the integrand factors in the last
	integral on RHS~\eqref{E:ENERGYIDENTITYFORBTRANSPORTEQUATIONS},
	and we use \eqref{E:ROUGHTIMEFUNCTIONLDERIVATIVEBOUNDS} to deduce that
	the integrand factors
	$\frac{1}{\Lunit \timefunctionarg{\muxmulevelsetvalue}}$ 
	in \eqref{E:ENERGYIDENTITYFORBTRANSPORTEQUATIONS}
	verify $\frac{1}{\Lunit \timefunctionarg{\muxmulevelsetvalue}} \approx 1$.
	We also use
	Lemma~\ref{L:ALLL2CONTROLLINGQUANTITIESINITIALLYSMALL}
	to bound the data-dependent terms 
	$\mathbb{E}[f]_{(\textnormal{Transport})}(\timefunction_0,u) + \mathbb{F}[f]_{(\textnormal{Transport})}(\timefunction,-\rightu)$
	on RHS~\eqref{E:ENERGYIDENTITYFORBTRANSPORTEQUATIONS} by $\lesssim \initialsmall^2$.
	Considering also Def.\,\ref{D:MAINCOERCIVE}, Lemma~\ref{L:COERCIVENESSOFL2CONTROLLINGQUANITIES},
	\eqref{E:ROUGHTIMEFUNCTIONLDERIVATIVEBOUNDS},
	and using Young's inequality,
	we deduce:
		\begin{align}  \label{E:FIRSTSTEPPROOFBELOWTOPORDERENERGYINTEGRALINEQUALITIESSPECIFICVORTICITYANDENTROPYGRADIENT}
		\begin{split}
		\hypersurfacecontrolVort_{\leq N}(\timefunction,u)
		+
		\hypersurfacecontrolGradEnt_{\leq N}(\timefunction,u)
		& 
		\lesssim
		\initialsmall^2
		+
		\int_{u' = - \rightu}^u
			\left\lbrace
				\hypersurfacecontrolVort_{\leq N}(\timefunction,u')
				+
				\hypersurfacecontrolGradEnt_{\leq N}(\timefunction,u')
			\right\rbrace
		\, \mathrm{d} u'
			\\
	& \ \
		+
		\int_{\twoargMrough{[\timefunction_0,\timefunction),[- \rightu,u]}{\muxmulevelsetvalue}}
			|\upmu \Transport \tander^{\leq N} (\vortrenormalized,\GradEnt)|^2
		\, \volMRoughCoordinates.
	\end{split}
	\end{align}
	Next, we use the pointwise estimate \eqref{E:COMMUTEDTRANSPORTPOINTWISEESTIMATESFORSPECIFICVORTICITYANDENTROPYGRADIENT}
	to bound the integrand factors of $|\upmu \Transport \tander^{\leq N} (\vortrenormalized,\GradEnt)|$
	on RHS~\eqref{E:FIRSTSTEPPROOFBELOWTOPORDERENERGYINTEGRALINEQUALITIESSPECIFICVORTICITYANDENTROPYGRADIENT}.
	Again appealing to Def.\,\ref{D:MAINCOERCIVE} and Lemma~\ref{L:COERCIVENESSOFL2CONTROLLINGQUANITIES},
	and also using \eqref{E:MINVALUEOFMUONFOLIATION},
	\eqref{E:ROUGHTIMEFUNCTIONLDERIVATIVEBOUNDS},
	and \eqref{E:PRELIMINARYEIKONALWITHOUTL},
	we conclude \eqref{E:BELOWTOPORDERENERGYINTEGRALINEQUALITIESSPECIFICVORTICITYANDENTROPYGRADIENT},
	but with the additional double integral
	$
	\varepsilon^2 
		\int_{\timefunction' = \timefunction_0}^{\timefunction} 
			\left\lbrace
			\int_{\timefunction'' = \timefunction_0}^{\timefunction'}
				\frac{\hypersurfacecontrolwave_{[1,N]}^{1/2}(\timefunction'',u)}{|\timefunction''|^{1/2}}
			\, \mathrm{d} \timefunction''
			\right\rbrace^2
		\, \mathrm{d} \timefunction'
	$
	on the RHS arising from the estimate \eqref{E:PRELIMINARYEIKONALWITHOUTL},
	which we use to handle the terms $\fundbootsmall |\tander^{[2,N]} \badcontrolvars|$
	on RHS~\eqref{E:COMMUTEDTRANSPORTPOINTWISEESTIMATESFORSPECIFICVORTICITYANDENTROPYGRADIENT}.
	By using that $\hypersurfacecontrolwave_{[1,N]}(\timefunction,u)$ is increasing in its arguments,
	we can bound this double integral by
	$
	\lesssim
	\varepsilon^2 
		\int_{\timefunction' = \timefunction_0}^{\timefunction} 
			\hypersurfacecontrolwave_{[1,N]}(\timefunction',u) 
		\, \mathrm{d} \timefunction'
	$,
	which in turn is bounded by RHS~\eqref{E:BELOWTOPORDERENERGYINTEGRALINEQUALITIESSPECIFICVORTICITYANDENTROPYGRADIENT} as desired.
	We have therefore proved \eqref{E:BELOWTOPORDERENERGYINTEGRALINEQUALITIESSPECIFICVORTICITYANDENTROPYGRADIENT}.
	
	The estimate \eqref{E:BELOWTOPORDERENERGYINTEGRALINEQUALITIESFORMODIFIEDFLUID}
	can be proved using similar arguments
	based on the pointwise estimates 
	\eqref{E:COMMUTEDTRANSPORTPOINTWISEESTIMATESFORMODIFIEDCURLOFVORT}--\eqref{E:COMMUTEDTRANSPORTPOINTWISEESTIMATESFORMODIFIEDDIVERGENCEOFENTROPYGRADIENT}
	for $0 \leq N \leq \Ntop-1$. 
	However, we need two additional ingredients, which we now point out:
	to bound the spacetime integrals
	$
		\varepsilon^2 
		\int_{\twoargMrough{[\timefunction_0,\timefunction),[- \rightu,u]}{\muxmulevelsetvalue}}
			|\tander^{N+1} \wavearray|^2
		\, \volMRoughCoordinates
	$
	generated by the terms $\varepsilon |\tander^{N+1} \wavearray|$
	on	
	RHSs~\eqref{E:COMMUTEDTRANSPORTPOINTWISEESTIMATESFORMODIFIEDCURLOFVORT}--\eqref{E:COMMUTEDTRANSPORTPOINTWISEESTIMATESFORMODIFIEDDIVERGENCEOFENTROPYGRADIENT},
	we also need to use 
	\textbf{i)} 
	the $L^2(\nullhypthreearg{\muxmulevelsetvalue}{u}{[\timefunction_0,\timefunction)})$-control 
	guaranteed by 
	\eqref{E:ROUGHTIMEFUNCTIONLDERIVATIVEBOUNDS}
	and
	\eqref{E:COERCIVENESSOFHYPERSURFACECONTROLWAVE},
	which leads to the presence of
	the term 
	$
		\varepsilon^2 
		\int_{u' = - \rightu}^u 
			\hypersurfacecontrolwave_{[1,N]}(\timefunction,u')
		\, \mathrm{d} u' 
	$
	on RHS~\eqref{E:BELOWTOPORDERENERGYINTEGRALINEQUALITIESFORMODIFIEDFLUID};
	and \textbf{ii)} the spacetime integral coerciveness estimate
	\eqref{E:COERCIVENESSOFSPACETIMEINTEGRAL},
	which leads the presence of the term
	$\fundbootsmall^2 \spacetimeintegralcontrolwave_{[1,N]}(\timefunction,u)$
	on RHS~\eqref{E:BELOWTOPORDERENERGYINTEGRALINEQUALITIESFORMODIFIEDFLUID}.
\end{proof}

\subsection{Proof of the estimates \eqref{E:MAINL2ESTIMATESVORTANDENTROPYGRADIENTBLOWUP} and
\eqref{E:MAINL2BELOWTOPORDERESTIMATESMODIFIEDFLUIDBLOWUP}--\eqref{E:MAINL2BELOWTOPORDERESTIMATESMODIFIEDFLUIDVARIABLESREGULAR}}
\label{SS:PROOFOFBELOWTOPORDERTRANSPORTENERGYESTIMATES}
We first prove \eqref{E:MAINL2BELOWTOPORDERESTIMATESVORTANDENTROPYGRADIENTREGULAR}. We set
$\mathbb{T}(\timefunction,u)  
\eqdef 	
	\hypersurfacecontrolVort_{\leq \Ntop-8}(\timefunction,u) 
	+ 
	\hypersurfacecontrolGradEnt_{\leq \Ntop-8}(\timefunction,u)$.
From
\eqref{E:BELOWTOPORDERENERGYINTEGRALINEQUALITIESSPECIFICVORTICITYANDENTROPYGRADIENT},
the wave energy bootstrap assumptions \eqref{E:MAINWAVEENERGYBOOTSTRAPREGULAR},
and \eqref{E:NONLINEARINEQUALITYRELATINGDATAEPSILONANDBOOTSTRAPEPSILON},
we find that
$
\mathbb{T}(\timefunction,u)
\leq 
C \initialsmall^2
+
C
\int_{u' = - \rightu}^u 
	\mathbb{T}(\timefunction,u')
\, \mathrm{d} u'
$.
Applying Gr\"{o}nwall's inequality, we conclude that $\mathbb{T}(\timefunction,u)
\leq 
C \initialsmall^2$, which yields \eqref{E:MAINL2BELOWTOPORDERESTIMATESVORTANDENTROPYGRADIENTREGULAR}.

Similarly, to prove \eqref{E:MAINL2BELOWTOPORDERESTIMATESMODIFIEDFLUIDVARIABLESREGULAR},
we set
$\mathbb{T}(\timefunction,u)  
\eqdef 	
	\hypersurfacecontrolVortVort_{\leq \Ntop-9}(\timefunction,u) 
+ 
\hypersurfacecontrolDivGradEnt_{\leq \Ntop-9}(\timefunction,u)$.
From
\eqref{E:BELOWTOPORDERENERGYINTEGRALINEQUALITIESFORMODIFIEDFLUID},
the wave energy bootstrap assumptions \eqref{E:MAINWAVEENERGYBOOTSTRAPREGULAR},
\eqref{E:NONLINEARINEQUALITYRELATINGDATAEPSILONANDBOOTSTRAPEPSILON},
and the already proved estimates \eqref{E:MAINL2BELOWTOPORDERESTIMATESVORTANDENTROPYGRADIENTREGULAR},
we find that
$
\mathbb{T}(\timefunction,u)
\leq 
C \initialsmall^2
+
C
\int_{u' = - \rightu}^u 
	\mathbb{T}(\timefunction,u')
\, \mathrm{d} u'
$.
Applying Gr\"{o}nwall's inequality, we conclude that $\mathbb{T}(\timefunction,u)
\leq 
C \initialsmall^2$, which yields \eqref{E:MAINL2BELOWTOPORDERESTIMATESMODIFIEDFLUIDVARIABLESREGULAR}.

The estimates 
\eqref{E:MAINL2ESTIMATESVORTANDENTROPYGRADIENTBLOWUP}
and
\eqref{E:MAINL2BELOWTOPORDERESTIMATESMODIFIEDFLUIDBLOWUP}
can be proved by combining similar arguments with the wave energy bootstrap assumptions 
\eqref{E:MAINWAVEENERGYBOOTSTRAPBLOWUP}--\eqref{E:MAINWAVEENERGYBOOTSTRAPREGULAR}.

\hfill $\qed$

\subsection{Proof of the estimates 
\eqref{E:MAINTORIL2VORTGRADENTBELOWTOPORDERBLOWUP}, 
\eqref{E:MAINTORIL2MODIFIEDFLUIDVARIABLESBELOWTOPORDERBLOWUP},
and
\eqref{E:MAINTORIL2VORTGRADENTBELOWTOPORDERREGULAR}--\eqref{E:MAINTORIL2MODIFIEDFLUIDVARIABLESBELOWTOPORDERREGULAR}}
\label{SS:PROOFOFBELOWTOPORDERROUGHTORIENERGYESTIMATES}
Fix any integer $N$ with $0 \leq N \leq \Ntop-1$.
Using \eqref{E:ROUGHTORUSINNULLHYPERSURFACEL2FUNDAMENTALTHEOREMOFCALCULUSESTIMATE},
\eqref{E:CLOSEDVERSIONLUNITROUGHTTIMEFUNCTION},
\eqref{E:COERCIVENESSOFCONTROLVORT}--\eqref{E:COERCIVENESSOFCONTROLGRADENT},
and the data estimates
\eqref{E:VORTICITYTORIL2CONTROLLINGINITIALLYSMALL}--\eqref{E:GRADENTTORIL2CONTROLLINGINITIALLYSMALL},
we deduce that:
\begin{align} 
	\begin{split} \label{E:FIRSTSTEPPROOFOFBELOWTOPORDERESTIMATESP:ROUGHTORIENERGYESTIMATES}
	\left\| 
		\tander^N (\Omega,\GradEnt) 
	\right\|_{L^2\left(\twoargroughtori{\timefunction,u}{\muxmulevelsetvalue}\right)}^2 
	& 
	\lesssim 
	\left\| 
		\tander^N (\Omega,\GradEnt)  
	\right\|_{L^2\left(\twoargroughtori{\timefunction_0,u}{\muxmulevelsetvalue}\right)}^2 
	+ 
	\int_{\nullhypthreearg{\muxmulevelsetvalue}{u}{[\timefunction_0,\timefunction]}} 
		\frac{1}{\Lunit \timefunctionarg{\muxmulevelsetvalue}} \left|\Lunit \tander^N (\Omega,\GradEnt) \right|^2 
	\, \volPuRoughCoordinates
		\\
	& 
	\lesssim 
	\initialsmall^2 
	+ 
	\left\| 
		\frac{1}{\sqrt{\Lunit \timefunctionarg{\muxmulevelsetvalue}}}
		\tander^{\leq N}(\vortrenormalized,\GradEnt)
	\right\|_{L^2(\nullhypthreearg{\muxmulevelsetvalue}{u}{[\timefunction_0,\timefunction]})}^2
	\lesssim
	\initialsmall^2
	+
	\hypersurfacecontrolVort_{N+1}(\timefunction,u)
	+
	\hypersurfacecontrolGradEnt_{N+1}(\timefunction,u).
	\end{split}
	\end{align}
	From \eqref{E:FIRSTSTEPPROOFOFBELOWTOPORDERESTIMATESP:ROUGHTORIENERGYESTIMATES}
	and the already proven estimates
	\eqref{E:MAINL2ESTIMATESVORTANDENTROPYGRADIENTBLOWUP} and \eqref{E:MAINL2BELOWTOPORDERESTIMATESVORTANDENTROPYGRADIENTREGULAR}
	for
	$
	\hypersurfacecontrolVort_{N+1}(\timefunction,u),
	$
	and
	$
	\hypersurfacecontrolGradEnt_{N+1}(\timefunction,u)
	$,
	we conclude,
	in view of definitions~\eqref{E:TORIVORTICITYL2CONTROLLINGQUANTITY}--\eqref{E:TORIENTROPYGRADIENTL2CONTROLLINGQUANTITY},
	the desired estimates \eqref{E:MAINTORIL2VORTGRADENTBELOWTOPORDERBLOWUP} and \eqref{E:MAINTORIL2VORTGRADENTBELOWTOPORDERREGULAR}. 
	
	The estimates
	\eqref{E:MAINTORIL2MODIFIEDFLUIDVARIABLESBELOWTOPORDERBLOWUP}
	and
	\eqref{E:MAINTORIL2MODIFIEDFLUIDVARIABLESBELOWTOPORDERREGULAR}
	can be proved via similar arguments based on the data estimates
	\eqref{E:MODIFIEDVORTICITYTORIL2CONTROLLINGINITIALLYSMALL}--\eqref{E:DIVGRADENTTORIL2CONTROLLINGINITIALLYSMALL},
	\eqref{E:COERCIVENESSHYPERSURFACECONTROLVORTVORT}--\eqref{E:COERCIVENESSHYPERSURFACEDIVGRADENT},
	the already proven estimates
	\eqref{E:MAINL2BELOWTOPORDERESTIMATESMODIFIEDFLUIDBLOWUP} and \eqref{E:MAINL2BELOWTOPORDERESTIMATESMODIFIEDFLUIDVARIABLESREGULAR},
	and definitions~\eqref{E:TORIMODIFIEDVORTICITYL2CONTROLLINGQUANTITY}--\eqref{E:TORIDIVGRADENTL2CONTROLLINGQUANTITY}.
	
\hfill $\qed$

\section{Top-order elliptic-hyperbolic $L^2$ estimates for the specific vorticity and entropy gradient}
\label{S:TOPORDERELIPTICHYPERBOLICL2ESTIMATESFORSPECIFICVORTICITYANDENTROPYGRADIENT}
Our main goal in this section is to derive the top-order $L^2$ estimate 
\eqref{E:MAINTOPORDERENERGYESTIMATESMODIFIEDFLUIDVARIABLESBLOWUP}
for the modified fluid variables $\VortVort$ and $\DivGradEnt$.
It turns out that by virtue of the elliptic-hyperbolic 
integral identity \eqref{E:INTEGRALIDENTITYFORELLIPTICHYPERBOLICCURRENT},
the proof of \eqref{E:MAINTOPORDERENERGYESTIMATESMODIFIEDFLUIDVARIABLESBLOWUP}
is coupled to the proof of the top-order energy estimates 
\eqref{E:MAINTORIL2VORTGRADENTTOPORDERBLOWUP}
for $\vortrenormalized$ and $\GradEnt$
along the rough tori. Hence, we also prove \eqref{E:MAINTORIL2VORTGRADENTTOPORDERBLOWUP}
in this section. Finally, as a simple consequence of \eqref{E:MAINTOPORDERENERGYESTIMATESMODIFIEDFLUIDVARIABLESBLOWUP}, 
we will also prove the top-order
rough tori energy estimate \eqref{E:MAINTORIL2MODIFIEDFLUIDVARIABLESTOPORDERBLOWUP} for
$\VortVort$ and $\DivGradEnt$.
In Sect.\,\ref{S:WAVEANDACOUSTICGEOMETRYAPRIORIESTIMATES}, we will use the estimates for $\VortVort$ and $\DivGradEnt$
that we derive in this section in our proof of the wave a priori estimates, 
which we stated as Prop.\,\ref{P:APRIORIL2ESTIMATESWAVEVARIABLES}.
Hence, we highlight that for the logic of the paper, it is important that 
\textbf{the estimates we derive in this section do not rely on the wave estimates of Prop.\,\ref{P:APRIORIL2ESTIMATESWAVEVARIABLES}};
our proofs of \eqref{E:MAINTOPORDERENERGYESTIMATESMODIFIEDFLUIDVARIABLESBLOWUP} 
and \eqref{E:MAINTORIL2VORTGRADENTTOPORDERBLOWUP}
instead rely on the bootstrap assumptions
\eqref{E:MAINWAVEENERGYBOOTSTRAPBLOWUP}--\eqref{E:MAINWAVEENERGYBOOTSTRAPREGULAR},
for the wave energies, which are \emph{weaker} than the estimates that we derive in Prop.\,\ref{P:APRIORIL2ESTIMATESWAVEVARIABLES}.

To explain the main challenges in the analysis,
we recall that the transport equations \eqref{E:EVOLUTIONEQUATIONFLATCURLRENORMALIZEDVORTICITY} 
and
\eqref{E:TRANSPORTFLATDIVGRADENT}
satisfied by
$\VortVort$ and $\DivGradEnt$
feature some difficult source terms,
denoted by $\mainnullform_{(\VortVort)}^i$ and $\mainnullform_{(\DivGradEnt)}$,
that depend on the general first-order derivatives
of $\vortrenormalized$ and $\GradEnt$. 
These source terms have the potential to cause the loss of a derivative
at the top-order because they cannot be bounded using pure transport estimates. 
In the below-top-order estimates of 
Sect.\,\ref{S:BELOWTOPORDERHYPERBOLICL2ESTIMATESFORSPECIFICVORTICITYANDENTROPYGRADIENT}, we allowed
the loss of a derivative, as is signified by the terms
$
\int_{u' = - \rightu}^u 
			\left\lbrace
				\hypersurfacecontrolVort_{\leq N+1}(\timefunction,u')
				+
				\hypersurfacecontrolGradEnt_{\leq N+1}(\timefunction,u')
			\right\rbrace
		\, \mathrm{d} u'	
$
on RHS~\eqref{E:BELOWTOPORDERENERGYINTEGRALINEQUALITIESFORMODIFIEDFLUID}.
To avoid the loss at the top-order, we handle the difficult source terms in a different way,
one that is based on combining the elliptic-hyperbolic integral identity provided by 
Prop.\,\ref{P:INTEGRALIDENTITYFORELLIPTICHYPERBOLICCURRENT}
with pointwise estimates that take into account the special structure of the equations 
of Theorem~\ref{T:GEOMETRICWAVETRANSPORTSYSTEM}, sharp estimates for the acoustic geometry 
and the rough time function, and the below-top-order estimates that we already derived in
Sect.\,\ref{S:BELOWTOPORDERHYPERBOLICL2ESTIMATESFORSPECIFICVORTICITYANDENTROPYGRADIENT}.

We organize this section this as follows: 
\begin{itemize}
\item
	In Sect.\,\ref{SS:PRELIMINARYENERGYINTEGRALIDENTITIESFORTOPDERIVATIVESOFMODIFIEDFLUID},
	we derive some preliminary energy integral inequalities for $\VortVort$ and $\DivGradEnt$, 
	which feature the difficult source terms; this part of the proof is not more difficult than the
	proof of the below-top-order energy integral inequalities we derived in 
	Lemma~\ref{L:BELOWTOPORDERENERGYINTEGRALINEQUALITIESSPECIFICVORTICITYANDENTROPYGRADIENT}.
\item In Sect.\,\ref{SS:MAINENERGYINTEGRALINEQUALITIESFORTOPORDERMODIFIEDFLUIDCONDITIONALONSOURCETERMPROPOSITION},
	we use the preliminary energy integral inequalities to derive the main energy integral inequalities for the top-order
	derivatives of $\VortVort$ and $\DivGradEnt$. These main energy integral inequalities
	are conditional on having $L^2$ estimates for the difficult source terms, which we derive independently as
	Prop.\,\ref{P:ELLIPTICHYPERBOLICINTEGRALINEQUALITIES}
	in Sect.\,\ref{SS:MAINELLIPTICHYPERBOLICINTEGRALINEQUALITIES}.
\item In Sect.\,\ref{SS:PROOFOFMAINVORTVORTDIVGRADENTTOPORDERBLOWUP},
	we use the main energy integral inequalities
	to prove the top-order $L^2$ estimate \eqref{E:MAINTOPORDERENERGYESTIMATESMODIFIEDFLUIDVARIABLESBLOWUP}.
\item In Sect.\,\ref{SS:CONTROLOFDATAFORROUGHTORIENERGYESTIMATES},
	to initiate the proof of Prop.\,\ref{P:ELLIPTICHYPERBOLICINTEGRALINEQUALITIES},
	we derive estimates for the rough tori error integrals
	$\left\| 
		\tander^{\leq \Ntop}(\vortrenormalized,\GradEnt)
	\right\|_{L^2(\roughtori{\timefunction,- \rightu})}$ for $\timefunction \in [\timefunction_0,\timefunctionboot]$,
	which appear on the right-hand side of the elliptic-hyperbolic integral identity
	\eqref{E:INTEGRALIDENTITYFORELLIPTICHYPERBOLICCURRENT}
	(with $\tander^{\leq \Ntop}(\vortrenormalized,\GradEnt)$ in the role of $\SigmatTan$)
	when we use the identity in our proof of Prop.\,\ref{P:ELLIPTICHYPERBOLICINTEGRALINEQUALITIES}.
	It might be tempting to think of these rough tori integrals as ``data terms'' since
	the rough tori $\roughtori{\timefunction,- \rightu}$ of interest
	are contained in the 
	``data null hypersurface'' portion 
	$\nullhypthreearg{\muxmulevelsetvalue}{-\rightu}{[\timefunction_0,\timefunctionboot]}$.
	In particular, if the data on $\Sigma_0$ are compactly supported in 
	$\Sigma_0 \cap \lbrace - \rightu \leq u \leq \leftu \rbrace$, then standard domain of dependence considerations
	imply that $(\vortrenormalized,\GradEnt)$ vanish along $\roughtori{\timefunction,- \rightu}$.
	However, in general, 
	the integrals
	$\left\| 
		\tander^{\leq \Ntop}(\vortrenormalized,\GradEnt)
	\right\|_{L^2(\roughtori{\timefunction,- \rightu})}$
	are not true ``data terms'' because their size depends
	on various norms of the rough time function $\timefunctionarg{\muxmulevelsetvalue}$, 
	which in turn depends on the behavior of the fluid near the singular boundary.
	Hence, to bound the integrals 
	$\left\| 
		\tander^{\leq \Ntop}(\vortrenormalized,\GradEnt)
	\right\|_{L^2(\roughtori{\timefunction,- \rightu})}$,
	we combine slight extensions of the Cauchy-stability results 
	that we derive in Appendix~\ref{A:OPENSETOFDATAEXISTS}
	with suitable $C_{\textnormal{geo}}^{2,1}$ estimates for $\timefunctionarg{\muxmulevelsetvalue}$.
 \item In Sect.\,\ref{SS:ELLIPTICHYPERBOLICIDENTITYL2ESTIMATESFORERRORINTEGRALS}, 
	in service of the proof of Prop.\,\ref{P:ELLIPTICHYPERBOLICINTEGRALINEQUALITIES},
	we derive estimates for the error integrals appearing in
	the elliptic-hyperbolic integral identity
	\eqref{E:INTEGRALIDENTITYFORELLIPTICHYPERBOLICCURRENT}.
\item In Sect.\,\ref{SS:MAINELLIPTICHYPERBOLICINTEGRALINEQUALITIES},
	we prove Prop.\,\ref{P:ELLIPTICHYPERBOLICINTEGRALINEQUALITIES}.
\item Finally, in Sect.\,\ref{SS:PROOFOFMAINTORIVORTGRADENTTOPORDERBLOWUP}, 
	we prove the top-order estimate \eqref{E:MAINTORIL2VORTGRADENTTOPORDERBLOWUP}
	for $\vortrenormalized$ and $\GradEnt$ along the rough tori
	as well as the top-order estimate \eqref{E:MAINTORIL2MODIFIEDFLUIDVARIABLESTOPORDERBLOWUP}
	for $\VortVort$ and $\DivGradEnt$ along the rough tori.
\end{itemize}

\subsection{Preliminary energy integral inequalities for the top-order derivatives of $\VortVort$ and $\DivGradEnt$}
\label{SS:PRELIMINARYENERGYINTEGRALIDENTITIESFORTOPDERIVATIVESOFMODIFIEDFLUID}

\begin{lemma}[Preliminary energy integral inequalities for the top-order derivatives of $\VortVort$ and $\DivGradEnt$]
\label{L:PRELIMINARYINTEGRALINEQUALITIESFORTOPDERIVATIVESOFMODIFIEDFLUID}
	For any $\varsigma \in (0,1]$, 
	the following integral inequalities hold for
	$(\timefunction,u) \in [\timefunction_0,\timefunctionboot) \times [-\rightu,\leftu]$,
	where the pointwise norm $|\cdot|_{\hfour}$ is defined in \eqref{E:SQUAREPOINTWISENORMWITHRESPECTTORIEMANNIANACOUSTICALMETRIC}
	and the implicit constants are independent of $\varsigma$:
	\begin{align} 
	\begin{split} \label{E:PRELIMINARYINTEGRALINEQUALITIESFORTOPDERIVATIVESOFMODIFIEDFLUID}
		\hypersurfacecontrolVortVort_{\Ntop}(\timefunction,u)
		+
		\hypersurfacecontrolDivGradEnt_{\Ntop}(\timefunction,u)
		&
		\lesssim
		\initialsmall^2
		+
		\varsigma
		\int_{\twoargMrough{[\timefunction_0,\timefunction),[- \rightu,u]}{\muxmulevelsetvalue}}
			\left\lbrace
				\left|
					\pmb{\partial} \tander^{\Ntop} \vortrenormalized 
				\right|_{\hfour}^2
				+
				\left|
					\pmb{\partial} \tander^{\Ntop} \GradEnt
				\right|_{\hfour}^2
			\right\rbrace
		\, \volMRoughCoordinates
			\\
	& \ \
		+
		\left(1 + \frac{1}{\varsigma} \right)
		\int_{u' = - \rightu}^u
			\left\lbrace
				\hypersurfacecontrolVortVort_{\Ntop}(\timefunction,u')
				+
				\hypersurfacecontrolDivGradEnt_{\Ntop}(\timefunction,u')
			\right\rbrace
		\, \mathrm{d} u'
			\\
	& \ \
		+
		\int_{u' = - \rightu}^u
			\left\lbrace \hypersurfacecontrolVortVort_{\leq \Ntop-1}(\timefunction,u') + \hypersurfacecontrolDivGradEnt_{\leq \Ntop-1}(\timefunction,u')\right\rbrace
		\, \mathrm{d} u'
			\\
	& \ \
		+
		\int_{u' = - \rightu}^u 
			\left\lbrace
				\hypersurfacecontrolVort_{\leq \Ntop}(\timefunction,u')
				+
				\hypersurfacecontrolGradEnt_{\leq \Ntop}(\timefunction,u')
			\right\rbrace
		\, \mathrm{d} u'	
			\\
	& \ \
		+
		\varepsilon^2 
		\int_{\timefunction' = \timefunction_0}^{\timefunction} 
			\hypersurfacecontrolwave_{[1,\Ntop]}(\timefunction',u) 
		\, \mathrm{d} \timefunction'
		+
		\varepsilon^2 
		\int_{u' = - \rightu}^u 
			\hypersurfacecontrolwave_{[1,\Ntop]}(\timefunction,u')
		\, \mathrm{d} u'.
	\end{split}
	\end{align}
	
\end{lemma}

\begin{proof}
	The proof is almost identical to the proof of \eqref{E:BELOWTOPORDERENERGYINTEGRALINEQUALITIESFORMODIFIEDFLUID}
	with $\Ntop$ in the role of $N$,
	except we separate the spacetime error integrals 
	generated by the top-order derivatives of $(\vortrenormalized,\GradEnt)$.
	More precisely, with $N \eqdef \Ntop$ in the pointwise estimates
	\eqref{E:COMMUTEDTRANSPORTPOINTWISEESTIMATESFORMODIFIEDCURLOFVORT}--\eqref{E:COMMUTEDTRANSPORTPOINTWISEESTIMATESFORMODIFIEDDIVERGENCEOFENTROPYGRADIENT}
for $|\upmu \Transport (\tander^{\Ntop} \VortVort,\tander^{\Ntop} \DivGradEnt)|$,	
we isolate the contribution of the error terms	
	$|\tander^{\Ntop + 1} (\vortrenormalized,\GradEnt)|$ on the RHSs.
	Since (schematically)
	$\tander^{\Ntop + 1} (\vortrenormalized^i,\GradEnt^i) = 
	\Singletan^{\alpha} \partial_{\alpha} \tander^{\Ntop} (\vortrenormalized^i,\GradEnt^i)
	$ for some $\Singletan \in \{\Lunit,\Yvf{2},\Yvf{3}\}$,
	we can use the simple Cartesian component bound $|\Singletan^{\alpha}| \lesssim 1$ 
	(which follows from Lemma~\ref{L:SCHEMATICSTRUCTUREOFVARIOUSTENSORSINTERMSOFCONTROLVARS} and the bootstrap assumptions)
	and \eqref{E:CARTERSIANDERIVATIVESHRIEMANNIANMETRICNORMCOMPARABLETOEUCLIDEANNORM}
	to pointwise bound these error terms in magnitude by 
	$\lesssim 
	\left|
		\pmb{\partial} \tander^{\Ntop} \vortrenormalized
	\right|_{\hfour}
	+
	\left|
		\pmb{\partial} \tander^{\Ntop} \GradEnt
	\right|_{\hfour}$.
	Thus, in the energy identity 
	(that is, \eqref{E:ENERGYIDENTITYFORBTRANSPORTEQUATIONS} with 
	$f \eqdef (\tander^{\Ntop}\VortVort,\tander^{\Ntop}\DivGradEnt)$), 
	the spacetime error integral corresponding to these terms is bounded by:
	\begin{align} \label{E:MAINERRORINTEGRALINPROOFOFPRELIMINARYINTEGRALINEQUALITIESFORTOPDERIVATIVESOFMODIFIEDFLUID}
	&
	\lesssim
	\int_{\twoargMrough{[\timefunction_0,\timefunction),[- \rightu,u]}{\muxmulevelsetvalue}}
			\left\lbrace
				\left|
					\tander^{\Ntop} \VortVort
				\right|_{\hfour}
				+
				\left|
					\tander^{\Ntop} \DivGradEnt
				\right|_{\hfour}
			\right\rbrace
			\left\lbrace
				\left|
					\pmb{\partial} \tander^{\Ntop} \vortrenormalized 
				\right|_{\hfour}
				+
				\left|
					\pmb{\partial} \tander^{\Ntop} \GradEnt
				\right|_{\hfour}
			\right\rbrace
		\, \volMRoughCoordinates.
	\end{align}
	Using the estimate \eqref{E:CLOSEDVERSIONLUNITROUGHTTIMEFUNCTION},
	\eqref{E:COERCIVENESSHYPERSURFACECONTROLVORTVORT}--\eqref{E:COERCIVENESSHYPERSURFACEDIVGRADENT},
	and Young's inequality, 
	for any $\varsigma \in (0,1]$,
	we find that
	$
	\mbox{RHS~\eqref{E:MAINERRORINTEGRALINPROOFOFPRELIMINARYINTEGRALINEQUALITIESFORTOPDERIVATIVESOFMODIFIEDFLUID}}
	\lesssim
		\frac{1}{\varsigma} 
		\int_{u' = - \rightu}^u
			\left\lbrace
				\hypersurfacecontrolVortVort_{\Ntop}(\timefunction,u')
				+
				\hypersurfacecontrolDivGradEnt_{\Ntop}(\timefunction,u')
			\right\rbrace
		\, \mathrm{d} u'
	+
	\varsigma
		\int_{\twoargMrough{[\timefunction_0,\timefunction),[- \rightu,u]}{\muxmulevelsetvalue}}
			\left\lbrace
				\left|
					\pmb{\partial} \tander^{\Ntop} \vortrenormalized 
				\right|_{\hfour}^2
				+
				\left|
					\pmb{\partial} \tander^{\Ntop} \GradEnt
				\right|_{\hfour}^2
			\right\rbrace
		\, \volMRoughCoordinates
	$,
	which is bounded by RHS~\eqref{E:PRELIMINARYINTEGRALINEQUALITIESFORTOPDERIVATIVESOFMODIFIEDFLUID}
	as desired.
\end{proof}

\subsection{The main energy integral inequalities for the top-order derivatives of $\VortVort$ and $\DivGradEnt$, conditional on
Prop.\,\ref{P:ELLIPTICHYPERBOLICINTEGRALINEQUALITIES}}
\label{SS:MAINENERGYINTEGRALINEQUALITIESFORTOPORDERMODIFIEDFLUIDCONDITIONALONSOURCETERMPROPOSITION}
Most of our effort in Sect.\,\ref{S:TOPORDERELIPTICHYPERBOLICL2ESTIMATESFORSPECIFICVORTICITYANDENTROPYGRADIENT} 
is dedicated towards bounding the
spacetime integrals:
$$
\varsigma
\int_{\twoargMrough{[\timefunction_0,\timefunctionboot),[- \rightu,\leftu]}{\muxmulevelsetvalue}}
			\left\lbrace
				\left|
					\pmb{\partial} \tander^{\Ntop} \vortrenormalized 
				\right|_{\hfour}^2
				+
				\left|
					\pmb{\partial} \tander^{\Ntop} \GradEnt
				\right|_{\hfour}^2
			\right\rbrace
		\, \volMRoughCoordinates,
	$$which appear on RHS~\eqref{E:PRELIMINARYINTEGRALINEQUALITIESFORTOPDERIVATIVESOFMODIFIEDFLUID}.
	We derive the needed estimates in Prop.\,\ref{P:ELLIPTICHYPERBOLICINTEGRALINEQUALITIES}.
	Given Prop.\,\ref{P:ELLIPTICHYPERBOLICINTEGRALINEQUALITIES} 
	and Lemma~\ref{L:PRELIMINARYINTEGRALINEQUALITIESFORTOPDERIVATIVESOFMODIFIEDFLUID}, 
	it is easy to derive energy integral inequalities
	that can be used to obtain the desired top-order energy estimates
	for $\VortVort$ and $\DivGradEnt$; we derive these integral inequalities in the next lemma.
	
\begin{lemma}[The main energy integral inequalities for the top-order derivatives of $\VortVort$ and $\DivGradEnt$]
\label{L:MAININTEGRALINEQUALITIESFORTOPDERIVATIVESOFMODIFIEDFLUID}
	For any $\varsigma \in (0,1]$, 
	the following integral inequalities hold for
	$(\timefunction,u) \in [\timefunction_0,\timefunctionboot) \times [-\rightu,\leftu]$,
	where the implicit constants are independent of $\varsigma$:
	\begin{align}
	\begin{split}  \label{E:MAININTEGRALINEQUALITIESFORTOPDERIVATIVESOFMODIFIEDFLUID}
		\hypersurfacecontrolVortVort_{\Ntop}(\timefunction,u)
		+
		\hypersurfacecontrolDivGradEnt_{\Ntop}(\timefunction,u)
		& \lesssim
		\varsigma 
		\left\lbrace
			\hypersurfacecontrolVortVort_{\Ntop}(\timefunction,u)
			+
			\hypersurfacecontrolDivGradEnt_{\Ntop}(\timefunction,u)
		\right\rbrace
		+
		\frac{\initialsmall^2}{|\timefunction|}
				\\
		& \ \
		+
		\left(1 + \frac{1}{\varsigma} \right)
		\int_{u' = - \rightu}^u
			\left\lbrace
				\hypersurfacecontrolVortVort_{\Ntop}(\timefunction,u')
				+
				\hypersurfacecontrolDivGradEnt_{\Ntop}(\timefunction,u')
			\right\rbrace
		\, \mathrm{d} u'
		+
		\varepsilon^2 
		\frac{1}{|\timefunction|^{3/2}}
		\hypersurfacecontrolwave_{[1,\Ntop]}(\timefunction,u)
			\\
		&  \ \
		+
		\frac{1}{|\timefunction|^2}
		\left\lbrace
			\hypersurfacecontrolVortVort_{\leq \Ntop-1}(\timefunction,u)
			+
			\hypersurfacecontrolDivGradEnt_{\leq \Ntop-1}(\timefunction,u)
		\right\rbrace
		+
		\frac{1}{|\timefunction|^{5/2}}
		\left\lbrace
			\hypersurfacecontrolVort_{\leq \Ntop}(\timefunction,u) 
			+
			\hypersurfacecontrolGradEnt_{\leq \Ntop}(\timefunction,u)
		\right\rbrace.
	\end{split}
	\end{align}

\end{lemma}

\begin{proof}
	We start with inequality \eqref{E:PRELIMINARYINTEGRALINEQUALITIESFORTOPDERIVATIVESOFMODIFIEDFLUID}.
	Using Lemma~\ref{L:COERCIVENESSOFELLITPICHYPERBOLICQUADRATICFORM}
	and the estimate
	$\frac{1}{\Lunit \timefunctionarg{\muxmulevelsetvalue}} \approx 1$
	implied by \eqref{E:ROUGHTIMEFUNCTIONLDERIVATIVEBOUNDS},
	we deduce the pointwise bounds
	$
	\left|
				\pmb{\partial} \tander^{\Ntop} \vortrenormalized 
	\right|_{\hfour}^2
	\lesssim
	\frac{1}{\Lunit \timefunctionarg{\muxmulevelsetvalue}}
			\ellipticCoerciveQuadratic[\pmb{\partial} \tander^{\Ntop} \vortrenormalized, \pmb{\partial} \tander^{\Ntop} \vortrenormalized]
	$
	and
	$
	\left|
					\pmb{\partial} \tander^{\Ntop} \GradEnt
				\right|_{\hfour}^2
	\lesssim
	\frac{1}{\Lunit \timefunctionarg{\muxmulevelsetvalue}}
			\ellipticCoerciveQuadratic[\pmb{\partial} \tander^{\Ntop} \GradEnt, \pmb{\partial} \tander^{\Ntop} \GradEnt]
	$
	for the error integrands on the first line of RHS~\eqref{E:PRELIMINARYINTEGRALINEQUALITIESFORTOPDERIVATIVESOFMODIFIEDFLUID}.
	Hence, thanks to the estimates
	\eqref{E:ELLIPTICHYPERBOLICINTEGRALINEQUALITYFORSPECIFICVORTICITY}--\eqref{E:ELLIPTICHYPERBOLICINTEGRALINEQUALITYFORENTROPYGRADIENT},
	which we prove independently in Sect.\,\ref{SS:MAINELLIPTICHYPERBOLICINTEGRALINEQUALITIES},
	the spacetime integral
	$
\varsigma
\int_{\twoargMrough{[\timefunction_0,\timefunctionboot),[- \rightu,\leftu]}{\muxmulevelsetvalue}}
			\left\lbrace
				\left|
					\pmb{\partial} \tander^{\Ntop} \vortrenormalized 
				\right|_{\hfour}^2
				+
				\left|
					\pmb{\partial} \tander^{\Ntop} \GradEnt
				\right|_{\hfour}^2
			\right\rbrace
		\, \volMRoughCoordinates
	$
	on the first line of RHS~\eqref{E:PRELIMINARYINTEGRALINEQUALITIESFORTOPDERIVATIVESOFMODIFIEDFLUID}
	is bounded by
\[	\varsigma
	\left\lbrace
		\mbox{RHS~\eqref{E:ELLIPTICHYPERBOLICINTEGRALINEQUALITYFORSPECIFICVORTICITY}}
		+
		\mbox{RHS~\eqref{E:ELLIPTICHYPERBOLICINTEGRALINEQUALITYFORENTROPYGRADIENT}}
	\right\rbrace.\]
The RHS of the resulting inequality features the
error integrals
$\left(1 + \frac{1}{\varsigma} \right)
		\int_{u' = - \rightu}^u
			\left\lbrace
				\hypersurfacecontrolVortVort_{\Ntop}(\timefunction,u')
				+
				\hypersurfacecontrolDivGradEnt_{\Ntop}(\timefunction,u')
			\right\rbrace
		\, \mathrm{d} u'
$,
which we place directly on RHS~\eqref{E:MAININTEGRALINEQUALITIESFORTOPDERIVATIVESOFMODIFIEDFLUID}.
Finally, we further bound the remaining error integrals on the RHS
by using the monotonicity of the controlling quantities 
(with respect to their arguments $\timefunction, u$) to pull the controlling quantities out of the integrals and to gain a power of $|\timefunction|$ upon
integration with respect to $\mathrm{d} \timefunction'$,
e.g.,
$
\int_{\timefunction' = \timefunction_0}^{\timefunction} 
			\frac{1}{|\timefunction'|^3}
			\hypersurfacecontrolVort_{\leq \Ntop}(\timefunction',u) 
		\, \mathrm{d} \timefunction'
\lesssim
\frac{1}{|\timefunction|^2}
\hypersurfacecontrolVort_{\leq \Ntop}(\timefunction,u) 
$
and
$
\int_{u' = - \rightu}^u 
			\hypersurfacecontrolDivGradEnt_{\Ntop}(\timefunction,u') 
		\, \mathrm{d} u'
\lesssim 
\hypersurfacecontrolDivGradEnt_{\Ntop}(\timefunction,u)
$.
In total, these arguments yield \eqref{E:MAININTEGRALINEQUALITIESFORTOPDERIVATIVESOFMODIFIEDFLUID}.
	
\end{proof}

\subsection{Proof of the main top-order energy estimate \eqref{E:MAINTOPORDERENERGYESTIMATESMODIFIEDFLUIDVARIABLESBLOWUP}}
\label{SS:PROOFOFMAINVORTVORTDIVGRADENTTOPORDERBLOWUP}
Given Lemma~\ref{L:MAININTEGRALINEQUALITIESFORTOPDERIVATIVESOFMODIFIEDFLUID},
we are now ready to prove our main top-order a priori energy estimates \eqref{E:MAINTOPORDERENERGYESTIMATESMODIFIEDFLUIDVARIABLESBLOWUP}
for the modified fluid variables.
We again emphasize that Lemma~\ref{L:MAININTEGRALINEQUALITIESFORTOPDERIVATIVESOFMODIFIEDFLUID}
is conditional on the estimates of Prop.\,\ref{P:ELLIPTICHYPERBOLICINTEGRALINEQUALITIES},
which we prove independently below.

To proceed, we choose and fix $\varsigma > 0$ sufficiently small such that the term
$
\varsigma 
		\left\lbrace
			\hypersurfacecontrolVortVort_{\Ntop}(\timefunction,u)
			+
			\hypersurfacecontrolDivGradEnt_{\Ntop}(\timefunction,u)
		\right\rbrace
$
on RHS~\eqref{E:MAININTEGRALINEQUALITIESFORTOPDERIVATIVESOFMODIFIEDFLUID}
can be absorbed back into the LHS at the expense of increasing the implicit
constants. Next, we 
use the already proven estimates
\eqref{E:MAINL2ESTIMATESVORTANDENTROPYGRADIENTBLOWUP} and
\eqref{E:MAINL2BELOWTOPORDERESTIMATESMODIFIEDFLUIDBLOWUP}--\eqref{E:MAINL2BELOWTOPORDERESTIMATESMODIFIEDFLUIDVARIABLESREGULAR},
the bootstrap assumptions \eqref{E:MAINWAVEENERGYBOOTSTRAPBLOWUP}--\eqref{E:MAINWAVEENERGYBOOTSTRAPREGULAR},
and \eqref{E:NONLINEARINEQUALITYRELATINGDATAEPSILONANDBOOTSTRAPEPSILON}
to bound all terms on RHS~\eqref{E:MAININTEGRALINEQUALITIESFORTOPDERIVATIVESOFMODIFIEDFLUID}
except for the integral
$
\int_{u' = - \rightu}^u
			\left\lbrace
				\hypersurfacecontrolVortVort_{\Ntop}(\timefunction,u')
				+
				\hypersurfacecontrolDivGradEnt_{\Ntop}(\timefunction,u')
			\right\rbrace
		\, \mathrm{d} u'
$.
In total, this leads to the following inequality:
\begin{align} \label{E:GRONWALLREADYMAINTOPORDERENERGYESTIMATESMODIFIEDFLUIDVARIABLESBLOWUP}
\hypersurfacecontrolVortVort_{\Ntop}(\timefunction,u)
+
\hypersurfacecontrolDivGradEnt_{\Ntop}(\timefunction,u)
& 
	\leq 
	C
	\initialsmall^2 |\timefunction|^{-17.1}
		+
	C 
		\int_{u' = - \rightu}^u
			\left\lbrace
				\hypersurfacecontrolVortVort_{\Ntop}(\timefunction,u')
				+
				\hypersurfacecontrolDivGradEnt_{\Ntop}(\timefunction,u')
			\right\rbrace
		\, \mathrm{d} u'.
\end{align}
From \eqref{E:GRONWALLREADYMAINTOPORDERENERGYESTIMATESMODIFIEDFLUIDVARIABLESBLOWUP} and
Gr\"{o}nwall's inequality,
we conclude
that 
$
\hypersurfacecontrolVortVort_{\Ntop}(\timefunction,u)
+
\hypersurfacecontrolDivGradEnt_{\Ntop}(\timefunction,u)
\leq 
	C
	\initialsmall^2 |\timefunction|^{-17.1}
$,
which yields the desired bound \eqref{E:MAINTOPORDERENERGYESTIMATESMODIFIEDFLUIDVARIABLESBLOWUP}.

\hfill $\qed$


\subsection{Control of $\left\| 
	\tander^{\leq \Ntop}(\vortrenormalized,\GradEnt)
\right\|_{L^2(\roughtori{\timefunction,- \rightu})}$
for $\timefunction \in [\timefunction_0,\timefunctionboot]$}
\label{SS:CONTROLOFDATAFORROUGHTORIENERGYESTIMATES}
Recall that our proof of the top-order $L^2$ estimates
\eqref{E:MAINTOPORDERENERGYESTIMATESMODIFIEDFLUIDVARIABLESBLOWUP}
and 
\eqref{E:MAINTORIL2VORTGRADENTTOPORDERBLOWUP}
relies on Prop.\,\ref{P:ELLIPTICHYPERBOLICINTEGRALINEQUALITIES},
whose proof relies on 
the integral identity 
\eqref{E:INTEGRALIDENTITYFORELLIPTICHYPERBOLICCURRENT}.
In order to exploit the identity \eqref{E:INTEGRALIDENTITYFORELLIPTICHYPERBOLICCURRENT}
with $u_1 = - \rightu$,
we in particular have to first control the rough tori error integrals
$\left\| 
	\tander^{\leq \Ntop}(\vortrenormalized,\GradEnt)
\right\|_{L^2(\roughtori{\timefunction,- \rightu})}$ for $\timefunction \in [\timefunction_0,\timefunctionboot]$;
in view of \eqref{E:MUISLARGEINBORINGREGION}
and 
\eqref{E:KEYCOERCIVITYPERFECTRDERIVAIVEERRORPUTANGENTCURRENTCONTRACTEDAGAINSTVECTORFIELD}
with $\tander^{\leq \Ntop}(\vortrenormalized,\GradEnt)$ in the role of $\SigmatTan$, 
we see that control of 
$\left\| 
	\tander^{\leq \Ntop}(\vortrenormalized,\GradEnt)
\right\|_{L^2(\roughtori{\timefunction,- \rightu})}
$ is sufficient to bound the
first integral $\int_{\twoargroughtori{\timefunction_2,u_1}{\muxmulevelsetvalue}} \cdots$ on
RHS~\eqref{E:INTEGRALIDENTITYFORELLIPTICHYPERBOLICCURRENT}.
As we explained at the beginning of Sect.\,\ref{S:TOPORDERELIPTICHYPERBOLICL2ESTIMATESFORSPECIFICVORTICITYANDENTROPYGRADIENT},
these integrals are not pure ``data terms'' because their size depends
on various norms of the rough time function, which in turn 
depends on the behavior of the fluid near the singular boundary.
Hence, in the next lemma,
we derive estimates for
$\left\| 
	\tander^{\leq \Ntop}(\vortrenormalized,\GradEnt)
\right\|_{L^2(\roughtori{\timefunction,- \rightu})}$.
The proof relies on the Cauchy stability-type estimates proved in
Appendix~\ref{A:OPENSETOFDATAEXISTS},
which also rely on various applications of the integral identity
\eqref{E:INTEGRALIDENTITYFORELLIPTICHYPERBOLICCURRENT}.
We emphasize that our bounds for
$\left\| 
	\tander^{\leq \Ntop}(\vortrenormalized,\GradEnt)
\right\|_{L^2(\roughtori{\timefunction,- \rightu})}$
cannot be proved by combining
the data estimates \eqref{E:VORTICITYANDENTROPYGRADIENTARESMALLONINITIALNULLHYERSURFACE}
along the null hypersurface $\nullhyparg{- \rightu}$
with trace estimates because of the usual loss of differentiability incurred by trace estimates.
However, if we had assumed that the data are one degree more differentiable, more precisely that
\eqref{E:VORTICITYANDENTROPYGRADIENTARESMALLONINITIALNULLHYERSURFACE} holds 
with $\Ntop$ replaced by $\Ntop+1$, then we could have used trace estimates to give a simpler proof 
of the desired estimates for
$\left\| 
	\tander^{\leq \Ntop}(\vortrenormalized,\GradEnt)
\right\|_{L^2(\roughtori{\timefunction,- \rightu})}$.
Though simpler, we avoided that approach because it would have led to estimates such that
the solution is less differentiable than the data.

\begin{lemma}[Control of $\left\| 
	\tander^{\leq \Ntop}(\vortrenormalized,\GradEnt)
\right\|_{L^2(\roughtori{\timefunction,- \rightu})}$
for {$\timefunction \in [\timefunction_0,\timefunctionboot]$}]
	\label{L:CONTROLOFDATAFORROUGHTORIENERGYESTIMATES}
Let 
$
\mathring{\Delta}_{\Sigma_0^{[-\farrightu,\leftu]}}^{\Ntop+1}
$
be the norm of the perturbation of the data away from the background solution,
as defined in \eqref{E:PERTURBATIONSMALLNESSINCARTESIANDIFFERENTIALSTRUCTURE}.
If $
\mathring{\Delta}_{\Sigma_0^{[-\farrightu,\leftu]}}^{\Ntop+1}
$
is sufficiently small, then the following estimates hold for 
$\timefunction \in [\timefunction_0,\timefunctionboot]$:
\begin{align} \label{E:SMALLDATAOFVORTICITYANDENTORPYGRADIENTONOUTERROUGHTORI}
\left\| 
	\tander^{\leq \Ntop}(\vortrenormalized,\GradEnt)
\right\|_{L^2(\roughtori{\timefunction,- \rightu})} 
& 
\leq \initialsmall,
\end{align}
where 
$\initialsmall = \mathcal{O}(\mathring{\Delta}_{\Sigma_0^{[-\farrightu,\leftu]}}^{\Ntop+1})$,
and the implicit constants depend on the background solution.
\end{lemma}

\begin{proof}
	The bootstrap assumption \eqref{E:BASIZEOFCARTESIANT} 
	implies that for $\timefunction \in [\timefunction_0,\timefunctionboot]$,
	the rough tori $\roughtori{\timefunction,- \rightu}$ are contained in the 
	``data null hypersurface'' $\nullhyparg{- \rightu}^{[0,4 \mathring{\updelta}_*]}$.
	For this reason,
	the proof of \eqref{E:SMALLDATAOFVORTICITYANDENTORPYGRADIENTONOUTERROUGHTORI}
	relies on the smallness results we derived for the solution
	on $\nullhyparg{- \rightu}^{[0,4 \mathring{\updelta}_*]}$ in
	Appendix~\ref{A:OPENSETOFDATAEXISTS}.
	In particular, in our proof here, we will refer to various steps in the
	proof of Prop.\,\ref{P:CAUCHYSTABILITYANDEXISTENCEOFOPENSETS}.
	We clarify that although the proof of Prop.\,\ref{P:CAUCHYSTABILITYANDEXISTENCEOFOPENSETS}
	relies on ideas from the bulk of the paper, 
	its proof is independent of the results of Lemma~\ref{L:CONTROLOFDATAFORROUGHTORIENERGYESTIMATES}.
	In the rest of the proof, we will
	silently assume that $\mathring{\Delta}_{\Sigma_0^{[-\farrightu,\leftu]}}^{\Ntop+1}$
	is sufficiently small.
	
	There are two broad steps in the proof of \eqref{E:SMALLDATAOFVORTICITYANDENTORPYGRADIENTONOUTERROUGHTORI}: 
	\textbf{I)} extend 
	the rough time function $\timefunctionarg{\muxmulevelsetvalue}$ into a subset of the ``smallness region'' $\CSregion_{Small}^{[0,5 \blowuptimePS]}$
	from Prop.\,\ref{P:CAUCHYSTABILITYANDEXISTENCEOFOPENSETS} (see Fig.\,\ref{F:CAUCHYSTABILITYREGION})
	and derive standard $C_{\textnormal{geo}}^{2,1}$ estimates for $\timefunctionarg{\muxmulevelsetvalue}$ in the extended region;
	and \textbf{II)} combine a version of the integral identity 
	\eqref{E:INTEGRALIDENTITYFORELLIPTICHYPERBOLICCURRENT} in a well-chosen subset of
	$\CSregion_{Small}^{[0,5 \blowuptimePS]}$
	with the estimates for $\timefunctionarg{\muxmulevelsetvalue}$ from step \textbf{I)},
	and use some results from the proof of Prop.\,\ref{P:CAUCHYSTABILITYANDEXISTENCEOFOPENSETS}
	to conclude \eqref{E:SMALLDATAOFVORTICITYANDENTORPYGRADIENTONOUTERROUGHTORI}.
	
	\medskip
	
	\noindent \textbf{Step I: Extending and controlling $\timefunctionarg{\muxmulevelsetvalue}$}.
	In Steps 1 and 2 of the proof of Prop.\,\ref{P:CAUCHYSTABILITYANDEXISTENCEOFOPENSETS}, 
	we show that up to the top-order derivative level
	(i.e., the derivative level corresponding to the energy estimates of
	Props.\,\ref{P:APRIORIL2ESTIMATESWAVEVARIABLES},
	\ref{P:MAINHYPERSURFACEENERGYESTIMATESFORTRANSPORTVARIABLES},
	\ref{P:ROUGHTORIENERGYESTIMATES},
	and
	\ref{P:APRIORIL2ESTIMATESACOUSTICGEOMETRY}),
	the fluid variables, $\upmu - 1$, $\Lsmall^i$, and $\upchi$ 
	are bounded by
	$\lesssim \mathring{\Delta}_{\Sigma_0^{[-\farrightu,\leftu]}}^{\Ntop+1}$
	on the region
	$\CSregion_{Small}^{[0,5 \blowuptimePS]}$ defined in \eqref{E:CAUCHYSTABILITYSOLUTIONISSMALLREGION},
	which, by \eqref{E:PERTURBEDBLOWUPDELTAISCLOSETOBACKGROUNDONE}, 
	contains $\nullhyparg{- \rightu}^{[0,4 \mathring{\updelta}_*]}$.
	That is, in $\CSregion_{Small}^{[0,5 \blowuptimePS]}$,
	the solution is close to the trivial fluid solution with Euclidean acoustic geometry.
	Unlike in 
	Props.\,\ref{P:APRIORIL2ESTIMATESWAVEVARIABLES},
	\ref{P:MAINHYPERSURFACEENERGYESTIMATESFORTRANSPORTVARIABLES},
	\ref{P:ROUGHTORIENERGYESTIMATES},
	and
	\ref{P:APRIORIL2ESTIMATESACOUSTICGEOMETRY},
	in Steps 1 and 2 of the proof of Prop.\,\ref{P:CAUCHYSTABILITYANDEXISTENCEOFOPENSETS}
	we derive the smallness 
	in $\CSregion_{Small}^{[0,5 \blowuptimePS]}$ 
	with respect to foliations by portions
	of \emph{Cartesian} time slices $\Sigma_t$,
	null hypersurfaces $\nullhyparg{u}$,
	and \emph{smooth tori} $\ell_{t,u}$.
	We will now explain how we can combine this smallness with 
	the transport equation \eqref{E:IVPFORROUGHTTIMEFUNCTION}
	and the data of $\timefunctionarg{\muxmulevelsetvalue}$ on
	$
	\nullhypthreearg{\muxmulevelsetvalue}{-\rightu}{[\timefunction_0,\timefunctionboot]}
	$
	(which is already ``known'' to us in view of the last item in Lemma~\ref{L:CONTINUOUSEXTNESION})
	to extend $\timefunctionarg{\muxmulevelsetvalue}$
	to a larger domain with the following two properties:
	\textbf{a)} $\timefunctionarg{\muxmulevelsetvalue} = \timefunctionarg{\muxmulevelsetvalue}(t,u,x^2,x^3)$ is defined on
	$
		\twoargMrough{[\timefunction_0,\timefunctionboot],[-U_*,\leftu]}{\muxmulevelsetvalue}
	$,
	where $U_* > 0$ is the number in \eqref{E:MINUSUSTARTLOCATION},
	and \textbf{b)}
		$
		\twoargMrough{[\timefunction_0,\timefunctionboot],[-U_*,-\rightu]}{\muxmulevelsetvalue}
		\subset 
		\twoargMrough{[\timefunction_0,\timefunctionboot],[-\rightu,-\leftu]}{\muxmulevelsetvalue} 
		\cup
		\CSregion_{Small}^{[0,5 \blowuptimePS]}
	$.
	To carry out this extension, we first note that
	\eqref{AE:PSBLOWUPTIME},
	\eqref{E:BASIZEOFCARTESIANT}, 
	\eqref{E:PERTURBEDBLOWUPDELTAISCLOSETOBACKGROUNDONE}, 
	and \eqref{E:CAUCHYSTABILITYSOLUTIONISSMALLREGION}
	imply that
	the data null hypersurface
	$
	\nullhypthreearg{\muxmulevelsetvalue}{- \rightu}{[\timefunction_0,\timefunctionboot]}
	$
	is contained in $\CSregion_{Small}^{[0,5 \blowuptimePS]}$
	and that the distance (with respect to the standard flat metric on geometric coordinate space)
	between
	$
	\nullhypthreearg{\muxmulevelsetvalue}{- \rightu}{[\timefunction_0,\timefunctionboot]}
	$
	and the top boundary of $\CSregion_{Small}^{[0,5 \blowuptimePS]}$
	(which we denote by ``$\Sigma_{5 \blowuptimePS}^{[15 \blowuptimePS - \farrightu,- \rightu]}$'' 
	in Fig.\,\ref{F:CAUCHYSTABILITYREGION}) is at least $\frac{1}{2} \blowuptimePS$.
	We next note that Def.\,\ref{D:WTRANSANDCUTOFF} implies that in $\CSregion_{Small}^{[0,5 \blowuptimePS]}$
	(a region in which $u < \interestingu$),
	the transport equation \eqref{E:IVPFORROUGHTTIMEFUNCTION}
	takes the form $\muX \timefunctionarg{\muxmulevelsetvalue} = 0$,
	where we recall that
	$\geop{u} = \muX - \muX^A \geop{x^A}$.
	The smallness provided by Steps 1 and 2 
	of the proof of Prop.\,\ref{P:CAUCHYSTABILITYANDEXISTENCEOFOPENSETS}
	implies that
	$
	\sum_{A=2,3}
	\| \muX^A \|_{C_{\textnormal{geo}}^{2,1}\left(\CSregion_{Small}^{[0,5 \blowuptimePS]} \right)}
	\lesssim 
	\mathring{\Delta}_{\Sigma_0^{[-\farrightu,\leftu]}}^{\Ntop+1}
	$.
	Using this smallness
	and deriving standard $C_{\textnormal{geo}}^{2,1}$-estimates for solutions
	to $\muX \timefunctionarg{\muxmulevelsetvalue} = 0$
	starting from the data of $\timefunctionarg{\muxmulevelsetvalue}$ on 
	the data null hypersurface portion
	$
	\nullhypthreearg{\muxmulevelsetvalue}{- \rightu}{[\timefunction_0,\timefunctionboot]}
	$
	(which is contained in $\twoargMrough{[\timefunction_0,\timefunctionboot],[- \rightu,\leftu]}{\muxmulevelsetvalue}$),
	we find that:
	\begin{align} \label{E:C21GEOBOUNDFORROUGHTIMEFUNCTIONINEXTENDEDREGION}
	\| \timefunctionarg{\muxmulevelsetvalue} \|_{C_{\textnormal{geo}}^{2,1}\left(\twoargMrough{[\timefunction_0,\timefunctionboot],[-U_*,-\rightu]}{\muxmulevelsetvalue} \right)}
	&
	\leq 
	\| \timefunctionarg{\muxmulevelsetvalue} \|_{C_{\textnormal{geo}}^{2,1}\left(\twoargMrough{[\timefunction_0,\timefunctionboot],[- \rightu,\leftu]}{\muxmulevelsetvalue} \right)}
	+
	\mathcal{O}\left(\mathring{\Delta}_{\Sigma_0^{[-\farrightu,\leftu]}}^{\Ntop+1} \right).
	\end{align}
	From \eqref{E:C21GEOBOUNDFORROUGHTIMEFUNCTIONINEXTENDEDREGION}
	and \eqref{E:CLOSEDVERSIONC21BOUNDFORCHOVROUGHTOGEO},
	we deduce the following bound, which we will use below:
	\begin{align} \label{E:C21GEOESTIMATEFORROUGHTIMEFUNCTIONINSMALLCACUHYSTABILITYREGION}
		\| \timefunctionarg{\muxmulevelsetvalue} \|_{C_{\textnormal{geo}}^{2,1}\left(\twoargMrough{[\timefunction_0,\timefunctionboot],[-U_*,\leftu]}{\muxmulevelsetvalue} \right)}
	\leq C.
	\end{align}
	
	\noindent \textbf{Step II: Using energy estimates and applications of the integral identity 
	\eqref{E:INTEGRALIDENTITYFORELLIPTICHYPERBOLICCURRENT} to finish the proof}.
	In Step 5 of the proof of Prop.\,\ref{P:CAUCHYSTABILITYANDEXISTENCEOFOPENSETS},
	we derived -- independently of Lemma~\ref{L:CONTROLOFDATAFORROUGHTORIENERGYESTIMATES} 
	(see Remark~\ref{R:ROUGHTORIESTIMATESONPUSTARAREINDEPENDENTOFLEMMACONTROLOFDATAFORROUGHTORIENERGYESTIMATES}) --
	geometric energy estimates in subsets
	$
	\twoargMrough{[\timefunction_*,\timefunction_0/2],[-U_*,\leftu]}{\muxmulevelsetvalue}
	$
	of $\CSregion_{Small}^{[0,5 \blowuptimePS]}$,
	where $\timefunction_* \in [2 \timefunction_0, (3/2) \timefunction_0]$ is the number from \eqref{E:LOCATIONOFTIMEFUNCTIONSTAR},
	$U_* > 0$ is the number from \eqref{E:MINUSUSTARTLOCATION},
	and $\timefunction_2$ and $u_2$ are any numbers satisfying
	$\timefunction_2 \in [\timefunction_*,\timefunction_0/2]$
	and
	$u_2 \in [-U_*,\leftu]$.
	These estimates in particular yield the rough tori $L^2$ estimates stated in
	\eqref{E:INCSSMALLREGIONENERGIESAREBOUNDEDONALLTORIINTEGRALSALONDATANULLHYPERSURFACEUPTOBOOSTRAPTIME}.
	The proof of \eqref{E:INCSSMALLREGIONENERGIESAREBOUNDEDONALLTORIINTEGRALSALONDATANULLHYPERSURFACEUPTOBOOSTRAPTIME}
	relies on applying the integral identity \eqref{E:INTEGRALIDENTITYFORELLIPTICHYPERBOLICCURRENT} on the region 
	$\twoargMrough{[\timefunction_*,\timefunction_2],[-U_*,u_2]}{\muxmulevelsetvalue}$,
	i.e., with $\timefunction_1 = \timefunction_*$
	and
	$u_1 = -U_*$ in \eqref{E:INTEGRALIDENTITYFORELLIPTICHYPERBOLICCURRENT}.
	We highlight two key ingredients that are needed 
	to control the error terms in that application of \eqref{E:INTEGRALIDENTITYFORELLIPTICHYPERBOLICCURRENT}:
	\textbf{i)} bounds for the rough tori integrals 
	arising from RHS~\eqref{E:INTEGRALIDENTITYFORELLIPTICHYPERBOLICCURRENT}
	(see Remark~\ref{R:ROUGHTORIESTIMATESONPUSTARAREINDEPENDENTOFLEMMACONTROLOFDATAFORROUGHTORIENERGYESTIMATES}),
	i.e., the integrals 
	$\int_{\twoargroughtori{\timefunction_2,-U_*}{\muxmulevelsetvalue}} \cdots$,
			$
			\int_{\twoargroughtori{\timefunction_*,u_2}{\muxmulevelsetvalue}}
				\cdots
			$,
			$
			\int_{\twoargroughtori{\timefunction_*,-U_*}{\muxmulevelsetvalue}}
				\cdots
			$,
	which we bound in the arguments given above 
	\eqref{E:INCSSMALLREGIONENERGIESAREBOUNDEDONALLTORIINTEGRALSALONDATANULLHYPERSURFACEUPTOBOOSTRAPTIME};
	and \textbf{ii)} the bound
	$
	\| \timefunctionarg{\muxmulevelsetvalue} \|_{C_{\textnormal{geo}}^{2,1}\left(\twoargMrough{[\timefunction_*,\timefunction_0/2],[-U_*,\leftu]}{\muxmulevelsetvalue} \right)}
	\leq C
	$,
	which is needed to bound various error terms on RHS~\eqref{E:INTEGRALIDENTITYFORELLIPTICHYPERBOLICCURRENT}
	that arise when we apply it on the region 
	$\twoargMrough{[\timefunction_*,\timefunction_0/2],[-U_*,\leftu]}{\muxmulevelsetvalue}$.
	The key point is that the estimate \eqref{E:C21GEOESTIMATEFORROUGHTIMEFUNCTIONINSMALLCACUHYSTABILITYREGION}
	shows that the same bound for $\timefunctionarg{\muxmulevelsetvalue}$ holds on the larger 
	subset $\twoargMrough{[\timefunction_0,\timefunctionboot],[-U_*,-\rightu]}{\muxmulevelsetvalue}$
	of $\CSregion_{Small}^{[0,5 \blowuptimePS]}$.
	Hence, the same arguments that yield
	\eqref{E:INCSSMALLREGIONENERGIESAREBOUNDEDONALLTORIINTEGRALSALONDATANULLHYPERSURFACEUPTOBOOSTRAPTIME} 
	can be used to show that:
	\begin{align} \label{E:ENERGIESAREBOUNDEDONALLTORIINTEGRALSALONDATANULLHYPERSURFACEUPTOBOOSTRAPTIME}
		\max_{\timefunction \in [\timefunction_*,\timefunctionboot]}
		\max_{u \in [-U_*,\leftu]}
		\int_{\roughtori{\timefunction,u}}
			|\tander^{\leq \Ntop} (\Omega,\GradEnt)|^2
		\, \volroughtorus
		&
		\lesssim
		\left(\mathring{\Delta}_{\Sigma_0^{[-\farrightu,\leftu]}}^{\Ntop+1} \right)^2.
\end{align}
Finally, we note that since 
$[\timefunction_0,\timefunctionboot] \subset [\timefunction_*,\timefunctionboot]$
and $- \rightu \in [-U_*,\leftu]$,
\eqref{E:ENERGIESAREBOUNDEDONALLTORIINTEGRALSALONDATANULLHYPERSURFACEUPTOBOOSTRAPTIME}
implies \eqref{E:SMALLDATAOFVORTICITYANDENTORPYGRADIENTONOUTERROUGHTORI}
with $\initialsmall = \mathcal{O}(\mathring{\Delta}_{\Sigma_0^{[-\farrightu,\leftu]}}^{\Ntop+1})$.	
	
\end{proof}

\subsection{$L^2$ estimates for the error terms}
\label{SS:ELLIPTICHYPERBOLICIDENTITYL2ESTIMATESFORERRORINTEGRALS}
In the next lemma, we derive $L^2$ estimates for the error integrals
in the elliptic-hyperbolic identities provided by Prop.\,\ref{P:INTEGRALIDENTITYFORELLIPTICHYPERBOLICCURRENT}.

\begin{lemma}[$L^2$ estimates for the error terms in the elliptic-hyperbolic identities]
	\label{L:L2ESTIMTAESFORELLIPTICHYPERBOLICIDENTITYERORTERMS}
	Let $\varsigma \in [0,1)$,
	let $(\timefunction,u) \in [\timefunction_0,\timefunctionboot) \times [-\rightu,\leftu]$,
	and let
	$
	\EllipticHyperbolicCurrentIntegralIdentityTotalSpacetimeErrorTerm[\pmb{\partial} \tander^{\Ntop} \vortrenormalized, \pmb{\partial} 			\tander^{\Ntop} \vortrenormalized]
	$
	and
	$
	\EllipticHyperbolicCurrentIntegralIdentityTotalSpacetimeErrorTerm[\pmb{\partial} \tander^{\Ntop} \GradEnt, 
	\pmb{\partial} \tander^{\Ntop} \GradEnt]
	$
	be the error terms defined by \eqref{E:ELLIPTICHYPERBOLICINTEGRALIDENTITYBULKERRORTERM}.
	Then the following spacetime integral estimates hold,
	where the implicit constants are independent of $\varsigma$:
	\begin{subequations}
	\begin{align} 
	\begin{split} \label{E:SPACETIMEERRORTERMESTIMATEFORSPECIFICVORTICITYELLIPTICHYPERBOLICIDENTITY}
		\left|
		\int_{\twoargMrough{[\timefunction_0,\timefunction),[- \rightu,u]}{\muxmulevelsetvalue}}
			\EllipticHyperbolicCurrentIntegralIdentityTotalSpacetimeErrorTerm[ \tander^{\Ntop} \vortrenormalized, \pmb{\partial} 			
			\tander^{\Ntop} \vortrenormalized]
		\, \volMRoughCoordinates
		\right|
		& \lesssim 
			\varsigma
			\int_{\twoargMrough{[\timefunction_0,\timefunction),[- \rightu,\leftu]}{\muxmulevelsetvalue}}
				\frac{1}{\Lunit \timefunctionarg{\muxmulevelsetvalue}}
				\ellipticCoerciveQuadratic[\pmb{\partial} \tander^{\Ntop} \vortrenormalized, \pmb{\partial} \tander^{\Ntop} \vortrenormalized]
			\, \volMRoughCoordinates
				\\
	& \ \
		+
			\int_{u' = - \rightu}^u 
				\hypersurfacecontrolVortVort_{\Ntop}(\timefunction,u') 
			\, \mathrm{d} u'
		+
		\int_{\timefunction' = \timefunction_0}^{\timefunction} 
			\frac{1}{|\timefunction'|^{3}}
			\hypersurfacecontrolVortVort_{\leq \Ntop-1}(\timefunction',u) 
		\, \mathrm{d} \timefunction'
			 \\
	&  \ \
		+
		\left(1 + \frac{1}{\varsigma} \right)
		\int_{\timefunction' = \timefunction_0}^{\timefunction} 
			\frac{1}{|\timefunction'|^{3}}
			\left\lbrace
				\hypersurfacecontrolVort_{\leq \Ntop}(\timefunction',u) 
				+	
				\hypersurfacecontrolGradEnt_{\leq \Ntop}(\timefunction',u) 
			\right\rbrace
		\, \mathrm{d} \timefunction'
			\\
	& \ \
		+
		\varepsilon^2 
		\int_{\timefunction' = \timefunction_0}^{\timefunction} 
			\frac{1}{|\timefunction'|^{2}}
			\hypersurfacecontrolwave_{[1,\Ntop]}(\timefunction',u) 
		\, \mathrm{d} \timefunction',
	\end{split}
			\\
	\begin{split}  \label{E:SPACETIMEERRORTERMESTIMATEFORENTROPYGRADIENTELLIPTICHYPERBOLICIDENTITY}
		\left|
		\int_{\twoargMrough{[\timefunction_0,\timefunction),[- \rightu,u]}{\muxmulevelsetvalue}}
			\EllipticHyperbolicCurrentIntegralIdentityTotalSpacetimeErrorTerm[ \tander^{\Ntop} \GradEnt, \pmb{\partial} \tander^{\Ntop} \GradEnt]
		\, \volMRoughCoordinates
		\right|
		& \lesssim
			\varsigma
			\int_{\twoargMrough{[\timefunction_0,\timefunction),[- \rightu,\leftu]}{\muxmulevelsetvalue}}
				\frac{1}{\Lunit \timefunctionarg{\muxmulevelsetvalue}}
				\ellipticCoerciveQuadratic[\pmb{\partial} \tander^{\Ntop} \GradEnt, \pmb{\partial} \tander^{\Ntop} \GradEnt]
			\, \volMRoughCoordinates
			 \\
		& \ \
			+
			\int_{u' = - \rightu}^u 
			\hypersurfacecontrolDivGradEnt_{\Ntop}(\timefunction,u') 
		\, \mathrm{d} u'
		+
		\int_{\timefunction' = \timefunction_0}^{\timefunction} 
			\frac{1}{|\timefunction'|^{3}}
			\hypersurfacecontrolDivGradEnt_{\leq \Ntop-1}(\timefunction',u) 
		\, \mathrm{d} \timefunction'
			 \\
	&  \ \
		+
		\left(1 + \frac{1}{\varsigma} \right)
		\int_{\timefunction' = \timefunction_0}^{\timefunction} 
			\frac{1}{|\timefunction'|^{3}}
			\left\lbrace
				\hypersurfacecontrolVort_{\leq \Ntop}(\timefunction',u) 
				+
				\hypersurfacecontrolGradEnt_{\leq \Ntop}(\timefunction',u) 
			\right\rbrace
		\, \mathrm{d} \timefunction'
			 \\
	& \ \
		+
		\varepsilon^2 
		\int_{\timefunction' = \timefunction_0}^{\timefunction} 
			\frac{1}{|\timefunction'|^{2}}
			\hypersurfacecontrolwave_{[1,\Ntop]}(\timefunction',u) 
		\, \mathrm{d} \timefunction'.
	\end{split}
	\end{align}
	\end{subequations}
	
	Moreover, the error terms
	$
	\Currentboundaryerrorhavetocontrolprincipal[\tander^{\Ntop} \vortrenormalized, \pmb{\partial} \tander^{\Ntop} \vortrenormalized],
	\cdots,
	\Currentboundaryerrorhavetocontrollowerorder[\tander^{\Ntop} \GradEnt, \tander^{\Ntop} \GradEnt]
	$
	defined by
	\eqref{E:PRINCIPALERRORTERMHAVETOCONTROLKEYIDPUTANGENTCURRENTCONTRACTEDAGAINSTVECTORFIELD}--\eqref{E:LOWERORDERERRORTERMHAVETOCONTROLKEYIDPUTANGENTCURRENTCONTRACTEDAGAINSTVECTORFIELD}
	verify the following rough hypersurface integral estimates:
	\begin{subequations}
	\begin{align} 
	\begin{split} \label{E:PRINCIPALSPATIALERRORTERMESTIMATEFORSPECIFICVORTICITYELLIPTICHYPERBOLICIDENTITY}
		& \left|
			\int_{\hypthreearg{\timefunction}{[- \rightu,u]}{\muxmulevelsetvalue}}
				\Currentboundaryerrorhavetocontrolprincipal
					[\tander^{\Ntop} \vortrenormalized, \pmb{\partial} \tander^{\Ntop} \vortrenormalized]
			\, \volRoughHypersurface
		\right| 
			\\
		& \ \ 
		 \lesssim
			\initialsmall^2
			+
			\hypersurfacecontrolVortVort_{\Ntop}(\timefunction,u)
			+
			\frac{1}{|\timefunction|^{3/2}}
			\hypersurfacecontrolVortVort_{\leq \Ntop-1}(\timefunction,u)
			+
			\frac{1}{|\timefunction|^{3/2}}
			 \hypersurfacecontrolDivGradEnt_{\leq \Ntop-1}(\timefunction,u) 
				\\
		& \ \
			+
			\varepsilon^2 
			\frac{1}{|\timefunction|^{3/2}}
			\hypersurfacecontrolwave_{[1,\Ntop]}(\timefunction,u)
			+
			\frac{1}{|\timefunction|^{5/2}}
			\hypersurfacecontrolVort_{\leq \Ntop}(\timefunction,u) 
			+
			\frac{1}{|\timefunction|^{5/2}}
			\hypersurfacecontrolGradEnt_{\leq \Ntop}(\timefunction,u),
		\end{split}
					 \\
		\begin{split} \label{E:PRINCIPALSPATIALERRORTERMESTIMATEFORENTROPYGRADIENTELLIPTICHYPERBOLICIDENTITY} 
		& \left|
			\int_{\hypthreearg{\timefunction}{[- \rightu,u]}{\muxmulevelsetvalue}}
				\Currentboundaryerrorhavetocontrolprincipal[\tander^{\Ntop} \GradEnt, \pmb{\partial} \tander^{\Ntop} \GradEnt]
			\, \volRoughHypersurface
		\right| 				
			\\
		& 
		 \ \ \lesssim
			\initialsmall^2
			+
			\hypersurfacecontrolDivGradEnt_{\Ntop}(\timefunction,u)
			+
			\frac{1}{|\timefunction|^{3/2}}
			\hypersurfacecontrolVortVort_{\leq 
			\Ntop-1}(\timefunction,u)
			+
			\frac{1}{|\timefunction|^{3/2}}
			 \hypersurfacecontrolDivGradEnt_{\leq \Ntop-1}(\timefunction,u)
					\\
		& \ \
			+
			\varepsilon^2 
			\frac{1}{|\timefunction|^{3/2}}
			\hypersurfacecontrolwave_{[1,\Ntop]}(\timefunction,u)
			+
			\frac{1}{|\timefunction|^{5/2}}
			\hypersurfacecontrolVort_{\leq \Ntop}(\timefunction,u) 
			+
			\frac{1}{|\timefunction|^{5/2}}
			\hypersurfacecontrolGradEnt_{\leq \Ntop}(\timefunction,u),
	\end{split}
	\end{align}
	\end{subequations}
	
	\begin{subequations}
	\begin{align} \label{E:LOWERORDERSPATIALERRORTERMESTIMATEFORSPECIFICVORTICITYELLIPTICHYPERBOLICIDENTITY}
		\left|
			\int_{\hypthreearg{\timefunction}{[- \rightu,u]}{\muxmulevelsetvalue}}
				\Currentboundaryerrorhavetocontrollowerorder[\tander^{\Ntop} \vortrenormalized, \tander^{\Ntop} \vortrenormalized]
			\, \volRoughHypersurface
		\right|
		& \lesssim
			\frac{1}{|\timefunction|^{5/2}}
			\hypersurfacecontrolVort_{\Ntop}(\timefunction,u),
			\\
	\left|
		\int_{\hypthreearg{\timefunction}{[- \rightu,u]}{\muxmulevelsetvalue}}
			\Currentboundaryerrorhavetocontrollowerorder[\tander^{\Ntop} \GradEnt, \tander^{\Ntop} \GradEnt]
		\, \volRoughHypersurface
	\right|
	& \lesssim
			\frac{1}{|\timefunction|^{5/2}}
			\hypersurfacecontrolGradEnt_{\Ntop}(\timefunction,u).
			\label{E:LOWERORDERSPATIALERRORTERMESTIMATEFORENTROPYGRADIENTELLIPTICHYPERBOLICIDENTITY}
	\end{align}
	\end{subequations}
	
	Finally, the error terms
	$\CurrentboundaryerrorperfectRderivative[\tander^{\Ntop} \vortrenormalized, \tander^{\Ntop} \vortrenormalized],
	\cdots,
	\CurrentboundaryerrorperfectRderivative[\tander^{\Ntop} \GradEnt, \tander^{\Ntop} \GradEnt]
	$
	defined by \eqref{E:PERFECTRDERIVAIVEERRORPUTANGENTCURRENTCONTRACTEDAGAINSTVECTORFIELD}
	verify the following rough tori integral estimates:
	\begin{subequations}
	\begin{align} 
	 \label{E:TOPORDERSMALLNESSOFVORTICITYANDENTROPYGRADIENTDATAALONGNULLHYPTORI}
		\int_{\twoargroughtori{\timefunction,- \rightu}{\muxmulevelsetvalue}}
			\CurrentboundaryerrorperfectRderivative[\tander^{\Ntop} \vortrenormalized, \tander^{\Ntop} \vortrenormalized]
		\, \volroughtorus,
			\,
		\int_{\twoargroughtori{\timefunction,- \rightu}{\muxmulevelsetvalue}}
			\CurrentboundaryerrorperfectRderivative[\tander^{\Ntop} \GradEnt, \tander^{\Ntop} \GradEnt]
		\, \volroughtorus
		& \lesssim
			\initialsmall^2,
				\\
		\int_{\twoargroughtori{\timefunction_0,u}{\muxmulevelsetvalue}}
			\CurrentboundaryerrorperfectRderivative[\tander^{\Ntop} \vortrenormalized, \tander^{\Ntop} \vortrenormalized]
		\, \volroughtorus,
			\,
		\int_{\twoargroughtori{\timefunction_0,u}{\muxmulevelsetvalue}}
			\CurrentboundaryerrorperfectRderivative[\tander^{\Ntop} \GradEnt, \tander^{\Ntop} \GradEnt]
		\, \volroughtorus
		& \lesssim
			\frac{\initialsmall^2}{|\timefunction_0|}.
			\label{E:TOPORDERSMALLNESSOFVORTICITYANDENTROPYGRADIENTDATAALONGROUGHHYPTORI}
\end{align}
\end{subequations}	
	
\end{lemma}

\begin{proof}
We first prove \eqref{E:TOPORDERSMALLNESSOFVORTICITYANDENTROPYGRADIENTDATAALONGNULLHYPTORI}.
We provide the details only for the $\tander^{\Ntop} \vortrenormalized$-dependent integral
on LHS~\eqref{E:TOPORDERSMALLNESSOFVORTICITYANDENTROPYGRADIENTDATAALONGNULLHYPTORI} 
since the $\tander^{\Ntop} \GradEnt$-dependent integral can be treated using identical arguments.
To proceed, we first use \eqref{E:MUISLARGEINBORINGREGION} to deduce that $\upmu$ is bounded from below by $\gtrsim 1$ on
$\twoargroughtori{\timefunction,- \rightu}{\muxmulevelsetvalue}$.
Hence,
using \eqref{E:KEYCOERCIVITYPERFECTRDERIVAIVEERRORPUTANGENTCURRENTCONTRACTEDAGAINSTVECTORFIELD},
\eqref{E:SMALLDATAOFVORTICITYANDENTORPYGRADIENTONOUTERROUGHTORI},
and \eqref{E:SIGMATTANGENTVECTORFIELDGNORMAPPROXIMATEDBYCARETSIANCOMPONENTNORMS},
we conclude that
$
\int_{\twoargroughtori{\timefunction,- \rightu}{\muxmulevelsetvalue}}
	\CurrentboundaryerrorperfectRderivative[\tander^{\Ntop} \vortrenormalized, \tander^{\Ntop} \vortrenormalized]
\, \volroughtorus
\lesssim
\int_{\twoargroughtori{\timefunction,- \rightu}{\muxmulevelsetvalue}}
	|\tander^{\Ntop} \vortrenormalized|_g^2
\, \volroughtorus
\lesssim
\initialsmall^2 
$
as desired.
Similarly, to prove \eqref{E:TOPORDERSMALLNESSOFVORTICITYANDENTROPYGRADIENTDATAALONGROUGHHYPTORI}
for the $\tander^{\Ntop} \vortrenormalized$-dependent integral on the LHS,
we use \eqref{E:MINVALUEOFMUONFOLIATION} with $\timefunction \eqdef \timefunction_0$,
\eqref{E:KEYCOERCIVITYPERFECTRDERIVAIVEERRORPUTANGENTCURRENTCONTRACTEDAGAINSTVECTORFIELD},
\eqref{E:SMALLDATAOFVORTICITYANDENTORPYGRADIENTONINITIALROUGHTORI},
and \eqref{E:SIGMATTANGENTVECTORFIELDGNORMAPPROXIMATEDBYCARETSIANCOMPONENTNORMS}
to conclude that
$
\int_{\twoargroughtori{\timefunction_0,u}{\muxmulevelsetvalue}}
	\CurrentboundaryerrorperfectRderivative[\tander^{\Ntop} \vortrenormalized, \tander^{\Ntop} \vortrenormalized]
\, \volroughtorus
\lesssim
\frac{1}{|\timefunction_0|}
\int_{\twoargroughtori{\timefunction_0,u}{\muxmulevelsetvalue}}
	|\tander^{\Ntop} \vortrenormalized|_g^2
\, \volroughtorus
\lesssim
\frac{1}{|\timefunction|_0}
\initialsmall^2 
$
as desired. The $\tander^{\Ntop} \GradEnt$-dependent integral 
on LHS~\eqref{E:TOPORDERSMALLNESSOFVORTICITYANDENTROPYGRADIENTDATAALONGROUGHHYPTORI}
can be bounded using identical arguments.

	The remaining estimates 
	\eqref{E:SPACETIMEERRORTERMESTIMATEFORSPECIFICVORTICITYELLIPTICHYPERBOLICIDENTITY}--\eqref{E:LOWERORDERSPATIALERRORTERMESTIMATEFORENTROPYGRADIENTELLIPTICHYPERBOLICIDENTITY}
	are straightforward consequences of 
	the pointwise estimates provided by Prop.\,\ref{P:POINTWISEESTIMTAESFORELLIPTICHYPERBOLICIDENTITYERORTERMS},
	the pointwise estimates
	$
	|\upmu| \lesssim 1
	$
	and
	$
	\left|
		\phi \frac{\muxmulevelsetvalue}{\Lunit \upmu}
	\right|
	\lesssim 1
	$
	(which follow from the bootstrap assumptions),
	Def.\,\ref{D:MAINCOERCIVE} of the $L^2$-controlling quantities,
	the coerciveness guaranteed by Lemma~\ref{L:COERCIVENESSOFL2CONTROLLINGQUANITIES}
	together with
	the estimates \eqref{E:CLOSEDVERSIONLUNITROUGHTTIMEFUNCTION} and
	\eqref{E:MINVALUEOFMUONFOLIATION},
	the already proven $L^2$ estimates \eqref{E:PRELIMINARYEIKONALWITHOUTL},
	and the fact that the $L^2$-controlling quantities
	$\hypersurfacecontrolwave_M(\timefunction,u)$, $\hypersurfacecontrolVortVort_M(\timefunction,u)$, etc.\
	are increasing in their arguments.
\end{proof}

\subsection{The main elliptic-hyperbolic integral inequalities}
\label{SS:MAINELLIPTICHYPERBOLICINTEGRALINEQUALITIES}
Thanks to the availability of 
Lemma~\ref{L:L2ESTIMTAESFORELLIPTICHYPERBOLICIDENTITYERORTERMS},
we are now ready to prove Prop.\,\ref{P:ELLIPTICHYPERBOLICINTEGRALINEQUALITIES}.

\begin{proposition}[The main elliptic-hyperbolic integral inequalities]
		\label{P:ELLIPTICHYPERBOLICINTEGRALINEQUALITIES}
	Let $\ellipticCoerciveQuadratic[\pmb{\partial} \SigmatTan,\pmb{\partial} \SigmatTan]$
	be the quadratic form from Def.\,\ref{D:NULLHYPERSURFACEADAPTEDCOERCIVEQUADRATICFORM},
	and let
	$\CurrentboundaryerrorperfectRderivative[\SigmatTan,\SigmatTan]$
	be the quadratic form defined by \eqref{E:PERFECTRDERIVAIVEERRORPUTANGENTCURRENTCONTRACTEDAGAINSTVECTORFIELD}.
	Then the following spacetime integral estimates hold for 
$(\timefunction,u) \in [\timefunction_0,\timefunctionboot) \times [-\rightu,\leftu]$:
	\begin{subequations}
	\begin{align}	
	\begin{split} \label{E:ELLIPTICHYPERBOLICINTEGRALINEQUALITYFORSPECIFICVORTICITY}
		&
		\int_{\twoargMrough{[\timefunction_0,\timefunction),[- \rightu,u]}{\muxmulevelsetvalue}}
			\frac{1}{\Lunit \timefunctionarg{\muxmulevelsetvalue}}
			\ellipticCoerciveQuadratic[\pmb{\partial} \tander^{\Ntop} \vortrenormalized, \pmb{\partial} \tander^{\Ntop} \vortrenormalized]
		\, \volMRoughCoordinates
		+
		\int_{\twoargroughtori{\timefunction,u}{\muxmulevelsetvalue}}
			\CurrentboundaryerrorperfectRderivative[\tander^{\Ntop} \vortrenormalized,\tander^{\Ntop} \vortrenormalized]
		\, \volroughtorus
			\\
		& \lesssim
			\frac{\initialsmall^2}{|\timefunction|}
			+
			\hypersurfacecontrolVortVort_{\Ntop}(\timefunction,u)
			+
			\frac{1}{|\timefunction|^{2}}
			\hypersurfacecontrolVortVort_{\leq \Ntop-1}(\timefunction,u)
			+
			\frac{1}{|\timefunction|^{2}}
			 \hypersurfacecontrolDivGradEnt_{\leq \Ntop-1}(\timefunction,u)
			 \\
		& \ \
			+
			\varepsilon^2 
			\frac{1}{|\timefunction|^{3/2}}
			\hypersurfacecontrolwave_{[1,\Ntop]}(\timefunction,u)
			+
			\frac{1}{|\timefunction|^{5/2}}
			\hypersurfacecontrolVort_{\leq \Ntop}(\timefunction,u) 
			+
			\frac{1}{|\timefunction|^{5/2}}
			\hypersurfacecontrolGradEnt_{\leq \Ntop}(\timefunction,u) 
					 \\
	& \ \
		+
		\int_{u' = - \rightu}^u 
			\hypersurfacecontrolVortVort_{\Ntop}(\timefunction,u') 
		\, \mathrm{d} u'
		+
		\int_{\timefunction' = \timefunction_0}^{\timefunction} 
			\frac{1}{|\timefunction'|^{3}}
			\hypersurfacecontrolVortVort_{\leq \Ntop-1}(\timefunction',u) 
		\, \mathrm{d} \timefunction'
			 \\
	&  \ \
		+
		\int_{\timefunction' = \timefunction_0}^{\timefunction} 
			\frac{1}{|\timefunction'|^{3}}
			\hypersurfacecontrolVort_{\leq \Ntop}(\timefunction',u) 
		\, \mathrm{d} \timefunction'
		+
		\int_{\timefunction' = \timefunction_0}^{\timefunction} 
			\frac{1}{|\timefunction'|^{3}}
			\hypersurfacecontrolGradEnt_{\leq \Ntop}(\timefunction',u) 
		\, \mathrm{d} \timefunction'
		+
		\varepsilon^2 
		\int_{\timefunction' = \timefunction_0}^{\timefunction} 
			\frac{1}{|\timefunction'|^{2}}
			\hypersurfacecontrolwave_{[1,\Ntop]}(\timefunction',u) 
		\, \mathrm{d} \timefunction',
	\end{split}	
			\\
	\begin{split} \label{E:ELLIPTICHYPERBOLICINTEGRALINEQUALITYFORENTROPYGRADIENT} 
	&
		\int_{\twoargMrough{[\timefunction_0,\timefunction),[- \rightu,u]}{\muxmulevelsetvalue}}
			\frac{1}{\Lunit \timefunctionarg{\muxmulevelsetvalue}}
			\ellipticCoerciveQuadratic[\pmb{\partial} \tander^{\Ntop} \GradEnt, \pmb{\partial} \tander^{\Ntop} \GradEnt]
		\, \volMRoughCoordinates
		+
		\int_{\twoargroughtori{\timefunction,u}{\muxmulevelsetvalue}}
			\CurrentboundaryerrorperfectRderivative[\tander^{\Ntop} \GradEnt,\tander^{\Ntop} \GradEnt]
		\, \volroughtorus
				\\
		& \lesssim
			\frac{\initialsmall^2}{|\timefunction|}
			+
			\hypersurfacecontrolDivGradEnt_{\Ntop}(\timefunction,u)
			+
			\frac{1}{|\timefunction|^{3/2}}
			\hypersurfacecontrolVortVort_{\leq \Ntop-1}(\timefunction,u)
			+
			\frac{1}{|\timefunction|^{3/2}}
			 \hypersurfacecontrolDivGradEnt_{\leq \Ntop-1}(\timefunction,u)
				 \\
		& \ \
			+
			\varepsilon^2 
			\frac{1}{|\timefunction|^{3/2}}
			\hypersurfacecontrolwave_{[1,\Ntop]}(\timefunction,u)
			+
			\frac{1}{|\timefunction|^{5/2}}
			\hypersurfacecontrolVort_{\leq \Ntop}(\timefunction,u) 
			+
			\frac{1}{|\timefunction|^{5/2}}
			\hypersurfacecontrolGradEnt_{\leq \Ntop}(\timefunction,u)
					\\
	& \ \
		+
		\int_{u' = - \rightu}^u 
			\hypersurfacecontrolDivGradEnt_{\Ntop}(\timefunction,u') 
		\, \mathrm{d} u'
		+
		\int_{\timefunction' = \timefunction_0}^{\timefunction} 
			\frac{1}{|\timefunction'|^{3}}
			\hypersurfacecontrolDivGradEnt_{\leq \Ntop-1}(\timefunction',u) 
		\, \mathrm{d} \timefunction'
				\\
	&  \ \
		+
		\int_{\timefunction' = \timefunction_0}^{\timefunction} 
			\frac{1}{|\timefunction'|^{3}}
			\hypersurfacecontrolVort_{\leq \Ntop}(\timefunction',u) 
		\, \mathrm{d} \timefunction'
		+
		\int_{\timefunction' = \timefunction_0}^{\timefunction} 
			\frac{1}{|\timefunction'|^{3}}
			\hypersurfacecontrolGradEnt_{\leq \Ntop}(\timefunction',u) 
		\, \mathrm{d} \timefunction'
		+
		\varepsilon^2 
		\int_{\timefunction' = \timefunction_0}^{\timefunction} 
			\frac{1}{|\timefunction'|^{2}}
			\hypersurfacecontrolwave_{[1,\Ntop]}(\timefunction',u) 
		\, \mathrm{d} \timefunction'.
\end{split}
\end{align}
\end{subequations}
\end{proposition}

\begin{proof}
	We consider the integral identity
	\eqref{E:INTEGRALIDENTITYFORELLIPTICHYPERBOLICCURRENT}
	with $\tander^{\Ntop} \vortrenormalized$ and $\tander^{\Ntop} \GradEnt$ in the role of $\SigmatTan$.
	Using Lemma~\ref{L:L2ESTIMTAESFORELLIPTICHYPERBOLICIDENTITYERORTERMS} with 
	$\timefunction_1 \eqdef \timefunction_0$,
	$\timefunction_2 \eqdef \timefunction$,
	$u_1 \eqdef - \rightu$,
	and
	$u_2 \eqdef u$,
	we bound the integrals on RHS~\eqref{E:INTEGRALIDENTITYFORELLIPTICHYPERBOLICCURRENT},
	where we can discard the integrals
	$
		-
		\int_{\twoargroughtori{\timefunction_1,u_1}{\muxmulevelsetvalue}}
				\CurrentboundaryerrorperfectRderivative[\SigmatTan,\SigmatTan]
	\, \volroughtorus
	$
	because they are non-positive in view of
	\eqref{E:KEYCOERCIVITYPERFECTRDERIVAIVEERRORPUTANGENTCURRENTCONTRACTEDAGAINSTVECTORFIELD}.
	Finally, by choosing and fixing $\varsigma > 0$ to be sufficiently small,
we can absorb the first term
$
\varsigma
			\int_{\twoargMrough{[\timefunction_0,\timefunction),[- \rightu,\leftu]}{\muxmulevelsetvalue}}
				\frac{1}{\Lunit \timefunctionarg{\muxmulevelsetvalue}}
				\ellipticCoerciveQuadratic[\pmb{\partial} \tander^{\Ntop} \vortrenormalized, \pmb{\partial} \tander^{\Ntop} \vortrenormalized]
			\, \volMRoughCoordinates
$
on RHS~\eqref{E:SPACETIMEERRORTERMESTIMATEFORSPECIFICVORTICITYELLIPTICHYPERBOLICIDENTITY}
and 
the first term 
$
\varsigma
			\int_{\twoargMrough{[\timefunction_0,\timefunction),[- \rightu,\leftu]}{\muxmulevelsetvalue}}
				\frac{1}{\Lunit \timefunctionarg{\muxmulevelsetvalue}}
				\ellipticCoerciveQuadratic[\pmb{\partial} \tander^{\Ntop} \GradEnt, \pmb{\partial} \tander^{\Ntop} \GradEnt]
			\, \volMRoughCoordinates
$
on RHS~\eqref{E:SPACETIMEERRORTERMESTIMATEFORENTROPYGRADIENTELLIPTICHYPERBOLICIDENTITY}
into LHS~\eqref{E:ELLIPTICHYPERBOLICINTEGRALINEQUALITYFORSPECIFICVORTICITY}
and
LHS~\eqref{E:ELLIPTICHYPERBOLICINTEGRALINEQUALITYFORENTROPYGRADIENT} respectively.
This yields the desired estimates
\eqref{E:ELLIPTICHYPERBOLICINTEGRALINEQUALITYFORSPECIFICVORTICITY}--\eqref{E:ELLIPTICHYPERBOLICINTEGRALINEQUALITYFORENTROPYGRADIENT}.

		
\end{proof}

\subsection{Proof of the top-order rough tori energy estimates \eqref{E:MAINTORIL2VORTGRADENTTOPORDERBLOWUP} and
\eqref{E:MAINTORIL2MODIFIEDFLUIDVARIABLESTOPORDERBLOWUP}}
\label{SS:PROOFOFMAINTORIVORTGRADENTTOPORDERBLOWUP}
We first prove the estimate \eqref{E:MAINTORIL2VORTGRADENTTOPORDERBLOWUP}.
We consider the estimates
\eqref{E:ELLIPTICHYPERBOLICINTEGRALINEQUALITYFORSPECIFICVORTICITY}--\eqref{E:ELLIPTICHYPERBOLICINTEGRALINEQUALITYFORENTROPYGRADIENT}.
In view of 
definitions \eqref{E:TORIVORTICITYL2CONTROLLINGQUANTITY}--\eqref{E:TORIENTROPYGRADIENTL2CONTROLLINGQUANTITY},
the quantitative positive definiteness estimate
\eqref{E:KEYCOERCIVITYPERFECTRDERIVAIVEERRORPUTANGENTCURRENTCONTRACTEDAGAINSTVECTORFIELD},
the estimate \eqref{E:BOUNDSONLMUINTERESTINGREGION},
the fact that the cutoff function $\phi$ on 
RHS~\eqref{E:KEYCOERCIVITYPERFECTRDERIVAIVEERRORPUTANGENTCURRENTCONTRACTEDAGAINSTVECTORFIELD}
is supported in the $u$-interval $[-\interestingu,\interestingu]$ (see Def.\,\ref{D:WTRANSANDCUTOFF}),
and the estimate $|\upmu| \lesssim 1$ (which follows from the bootstrap assumptions),
we see that up to $\mathcal{O}(1)$ factors, the rough tori integrals on   
LHSs~\eqref{E:ELLIPTICHYPERBOLICINTEGRALINEQUALITYFORSPECIFICVORTICITY}--\eqref{E:ELLIPTICHYPERBOLICINTEGRALINEQUALITYFORENTROPYGRADIENT} 
bound the terms on LHS~\eqref{E:MAINTORIL2VORTGRADENTTOPORDERBLOWUP} from above. 
Moreover, using the already proven estimates
\eqref{E:MAINL2ESTIMATESVORTANDENTROPYGRADIENTBLOWUP}--\eqref{E:MAINL2BELOWTOPORDERESTIMATESMODIFIEDFLUIDVARIABLESREGULAR},
the wave energy bootstrap assumptions \eqref{E:MAINWAVEENERGYBOOTSTRAPBLOWUP}--\eqref{E:MAINWAVEENERGYBOOTSTRAPREGULAR},
and \eqref{E:NONLINEARINEQUALITYRELATINGDATAEPSILONANDBOOTSTRAPEPSILON},
we see that all terms on 
RHSs~\eqref{E:ELLIPTICHYPERBOLICINTEGRALINEQUALITYFORSPECIFICVORTICITY}--\eqref{E:ELLIPTICHYPERBOLICINTEGRALINEQUALITYFORENTROPYGRADIENT}
are bounded by $\lesssim \initialsmall^2 |\timefunction|^{-17.1}$,
which yields the desired result.

The estimate \eqref{E:MAINTORIL2MODIFIEDFLUIDVARIABLESTOPORDERBLOWUP} follows from combining the same
arguments we used to prove \eqref{E:FIRSTSTEPPROOFOFBELOWTOPORDERESTIMATESP:ROUGHTORIENERGYESTIMATES}
with the data estimates
\eqref{E:MODIFIEDVORTICITYTORIL2CONTROLLINGINITIALLYSMALL}--\eqref{E:DIVGRADENTTORIL2CONTROLLINGINITIALLYSMALL},
the already proven estimate \eqref{E:MAINTOPORDERENERGYESTIMATESMODIFIEDFLUIDVARIABLESBLOWUP},
\eqref{E:COERCIVENESSHYPERSURFACECONTROLVORTVORT}--\eqref{E:COERCIVENESSHYPERSURFACEDIVGRADENT},
and definitions~\eqref{E:TORIMODIFIEDVORTICITYL2CONTROLLINGQUANTITY}--\eqref{E:TORIDIVGRADENTL2CONTROLLINGQUANTITY}.

\hfill $\qed$

\section{Elliptic estimates for the acoustic geometry on the rough tori \texorpdfstring{$\twoargroughtori{\timefunction,u}{\muxmulevelsetvalue}$}{ell}} 
\label{S:ELLIPTICESTIAMTESACOUSTICGEOMETRYONROUGHTORI}
In this section, we derive elliptic 
$L^2$ estimates for symmetric $\binom{0}{2}$-type tensorfields that are tangent to the smooth tori $\ell_{t,u}$. 
As we will explain, our analysis fundamentally relies on the \emph{rough} tori 
$\twoargroughtori{\timefunction,u}{\muxmulevelsetvalue}$, which, unlike the smooth tori, 
are adapted to our foliations by level sets of $\timefunctionarg{\muxmulevelsetvalue}$.
We provide the main estimate in Prop.\,\ref{P:ELLIPTICESTIMATESONROUGHTORI}.
In Sect.\,\ref{SS:TOPORDERL2ESTIMATSFORCHI}, 
we combine the elliptic estimates of Prop.\,\ref{P:ELLIPTICESTIMATESONROUGHTORI}
with hyperbolic $L^2$ estimates for the
fully modified quantities defined in Sect.\,\ref{S:CONSTRUCTIONOFMODIFIEDQUANTITIES}
to obtain top-order $L^2$ estimates for the null second fundamental form $\upchi$,
which is tangent to $\ell_{t,u}$.
We fundamentally need these top-order estimates for $\upchi$ to avoid the loss of a derivative 
in the top-order commuted wave equations.
Specifically, in the top-order case $N = \Ntop$,
we need these $L^2$ estimates to handle the terms
on RHSs \eqref{E:TOPCOMMUTEDWAVELFIRSTTHENALLYS}--\eqref{E:TOPCOMMUTEDWAVEALLYS}
that explicitly depend on the order $N$ derivatives $\mytr_{\gtorus}\upchi$.
The point is that when $N = \Ntop$, $\tanderY^N \mytr_{\gtorus}\upchi$ cannot be controlled in $L^2$ through
pure transport estimates; using only transport estimates at the top order would result in the loss of one derivative 
due to the presence of the source term 
$|\upchi|_{\gtorus}^2$ on the RHS of the transport equation \eqref{E:RAYCHAUDHURITRANSPORTCHI}
satisfied by $\mytr_{\gtorus}\upchi$,
which depends not only on
$\mytr_{\gtorus}\upchi$, but also on its trace-free part $\hat{\upchi}$.
The strategy of avoiding derivative loss in $\upchi$ via a combination of elliptic $L^2$ estimates
on co-dimension $2$ surfaces and hyperbolic $L^2$ estimates 
was originally employed in the context of Einstein's equations in \cite{dCsK1993}.
Later, this strategy was used in many other works on wave and wave-like equations,
for example, in the context of low regularity local well-posedness for quasilinear wave equations in \cite{sKiR2003},
in the context of irrotational shock formation in \cite{dC2007},
and in the context of shock formation in $3D$ with vorticity and entropy in \cite{jLjS2021}.
What is new here compared to these works is our reliance 
on the rough tori to obtain the needed top-order estimates, even though the operators
$\tanderY^N$ and $\upchi$ are adapted to the \emph{smooth tori} $\ell_{t,u} = \Sigma_t \cap \nullhyparg{u}$.

To obtain the desired elliptic estimates, we decompose
symmetric type $\binom{0}{2}$ $\ell_{t,u}$-tangent tensorfields $\upxi$ into 
a main piece that is tangent to $\twoargroughtori{\timefunction,u}{\muxmulevelsetvalue}$,
which we control with standard elliptic estimates on the rough tori (see Lemma~\ref{L:STANDARDELLIPTICONROUGHTORI}),
and error terms, which we must control with separate (easier) arguments.
Our primary application will be to apply the main elliptic estimate \eqref{E:ELLIPTICESTIMATECHI}
with $\angLie_{\tander}^{\Ntop-1} \upchi$ in the role of $\upxi$,
which will yield $L^2$ control of $\angLie_{\tander}^{\Ntop} \upchi$;
see the proof of \eqref{E:PRELIMINARYTOPORDERLIECHI}.
There are many equivalent ways we could have carried out the decompositions and analysis
of this section.
We have chosen to use orthonormal frames on the smooth tori and the rough tori
and to quantitatively control the relationship between the two frames;
see Sect.\,\ref{SS:ORTHONORMALFRAMESONSMOOTHANDROUGHTORI}.


\subsection{Statement of the main elliptic estimates}
\label{SS:ELLIPTICESTIMATESFORACOUSTICGEOMETRY}
In this section, we state the proposition that yields the main elliptic estimates of interest.
Its proof is located in Sect.\,\ref{SS:PROOFOFP:ELLIPTICESTIMATESONROUGHTORI}.

\begin{proposition}[The main elliptic estimates for symmetric type-$\binom{0}{2}$ $\ell_{t,u}$-tangent tensorfields] 
\label{P:ELLIPTICESTIMATESONROUGHTORI}
Let $\upxi$ be a symmetric type-$\binom{0}{2}$ 
$\ell_{t,u}$-tangent tensorfield. 
Then the following estimate holds
for $(\timefunction,u) \in [\timefunction_0,\timefunctionboot) \times [-\rightu,\leftu]$:
\begin{align}
\begin{split} \label{E:ELLIPTICESTIMATECHI}
\sum_{\Singletan \in \{\Lunit,\Yvf{2},\Yvf{3}\}} 
\int_{\hypthreearg{\timefunction}{[-\rightu,u]}{\muxmulevelsetvalue}} 
	\upmu^2 |\angLie_{\Singletan} \upxi|_{\gtorus}^2 
\, \volRoughHypersurface 
& 
\leq 
C \int_{\hypthreearg{\timefunction}{[-\rightu,u]}{\muxmulevelsetvalue}}
	\upmu^2 |\angLie_{\Lunit} \upxi|_{\gtorus}^2 
\, \volRoughHypersurface 
+ 
C 
\int_{\hypthreearg{\timefunction}{[-\rightu,u]}{\muxmulevelsetvalue}}
	\upmu^2 |\angdiv \upxi|_{\gtorus}^2 
\, \volRoughHypersurface 
	\\
& \ \
+ 
C 
\sum_{A=2,3} 
\int_{\hypthreearg{\timefunction}{[-\rightu,u]}{\muxmulevelsetvalue}} 
	\upmu^2 (\Yvf{A} \mytr_{\gtorus}\upxi)^2 
\, \volRoughHypersurface 
+ 
C 
\fundbootsmall 
\int_{\hypthreearg{\timefunction}{[-\rightu,u]}{\muxmulevelsetvalue}} 
	|\upxi|_{\gtorus}^2 
\, \volRoughHypersurface. 
\end{split}
\end{align}
\end{proposition}

\subsection{Orthonormal frames on the smooth tori and the rough tori}
\label{SS:ORTHONORMALFRAMESONSMOOTHANDROUGHTORI}
For use throughout Sect.\,\ref{S:ELLIPTICESTIAMTESACOUSTICGEOMETRYONROUGHTORI}, we recall that 
$\gtorus$ denotes the first fundamental form of the smooth tori $\ell_{t,u}$
and $\gtorusroughfirstfund$ 
denotes the first fundamental form of the rough tori $\twoargroughtori{\timefunction,u}{\muxmulevelsetvalue}$.
In our ensuing analysis, we will use the pairs of orthonormal frames featured in the next definition.

\begin{definition}[The frames $\{e_A\}_{A = 2,3}$ and $\{f_A\}_{A = 2,3}$]
 \label{D:ORTHONORMALFRAMESONSMOOTHANDROUGHTORI}
$\{e_A\}_{A = 2,3}$ is defined to be the orthonormal frame on the smooth torus $\ell_{t,u}$ 
obtained from applying the Gram--Schmidt process to the geometric coordinate partial derivative vectorfields
$\left\{\geop{x^A} \right\}_{A = 2,3}$ with respect to $\gtorus$, starting with 
$e_2 \eqdef \frac{1}{\sqrt{\gtorus(\geop{x^2},\geop{x^2})}}\geop{x^2}$.
Similarly, $\left\{f_A \right\}_{A = 2,3}$ is defined to be 
the orthonormal frame on the rough torus $\twoargroughtori{\timefunction,u}{\muxmulevelsetvalue}$ 
obtained from applying the Gram--Schmidt process to the rough adapted coordinate partial derivative vectorfields
$\left\{\roughgeop{x^A} \right\}_{A = 2,3}$ with respect to $\gtorusroughfirstfund$,
starting with 
$f_2 \eqdef \frac{1}{\sqrt{\gtorusroughfirstfund(\roughgeop{x^2},\roughgeop{x^2})}}\roughgeop{x^2}$.
\end{definition}

In the next lemma, we provide standard expressions for
$\gtorus^{-1}$ and $\gtorusroughinversefirstfund$ relative to the orthonormal frames.

\begin{lemma}[Expressions for $\gtorus^{-1}$ and $\gtorusroughinversefirstfund$ relative to orthonormal frames]
\label{L:EXPRESSIONSFORFIRSTFUNDFORMSOFTORIRELATIVETOORTHONORMALFRAMES}
Let $\gtorus^{-1}$ be the inverse first fundamental form of $\ell_{t,u}$ from Def.\,\ref{D:FIRSTFUNDAMENTALFORMS},
let $\gtorusroughinversefirstfund$ be the inverse first fundamental form of $\twoargroughtori{\timefunction,u}{\muxmulevelsetvalue}$
from Def.\,\ref{D:ROUGHFIRSTFUNDS},
and let $\{e_A \}_{A = 2,3}$ and $\{f_A\}_{A = 2,3}$ be the orthonormal frames from 
Def.\,\ref{D:ORTHONORMALFRAMESONSMOOTHANDROUGHTORI}.
Then the following identities hold, where $\updelta^{AB}$ is the Kronecker delta:
\begin{subequations}
\begin{align}
	\gtorus^{-1}
	& = \updelta^{AB} e_A \otimes e_B,
		\label{E:INVERSEFIRSTFUNDOFSMOOTHTORIRELATIVETOORTHONORMALFRAME} 
			\\
	\gtorusroughinversefirstfund
	& = \updelta^{AB} f_A \otimes f_B.
		\label{E:INVERSEFIRSTFUNDOFROUGHTORIRELATIVETOORTHONORMALFRAME} 
\end{align}
\end{subequations}
\end{lemma}

\begin{proof}
	\eqref{E:INVERSEFIRSTFUNDOFSMOOTHTORIRELATIVETOORTHONORMALFRAME}--\eqref{E:INVERSEFIRSTFUNDOFROUGHTORIRELATIVETOORTHONORMALFRAME}
	are standard identities for inverse metrics
	relative to orthonormal frames.
\end{proof}

In the next lemma, we exhibit the relationships between 
the two frames $\{e_A\}_{A = 2,3}$ and $\{f_A\}_{A = 2,3}$.
\begin{lemma}[Relationship between $\{e_A\}_{A = 2,3}$ and $\{f_A\}_{A = 2,3}$] 
\label{L:CHANGEOFORTHONORMALFRAMES}
On $\twoargMrough{[\timefunction_0,\timefunctionboot),[- \rightu,\leftu]}{\muxmulevelsetvalue}$,
there exists a $2 \times 2$ orthogonal-matrix-valued function $\COVframe$ with components $\{\COVframe_{AB}\}_{A,B = 2,3}$ 
and scalar functions $\{\COVL_A \}_{A = 2,3}$ such that: 
\begin{subequations}
\begin{align}
e_A & = 
\COVframe_{AB} f_B 
+ 
\COVL_A \Lunit 
= \COVframe_{AB}
	\left\lbrace
		f_B + \COVframe_{CB} \COVL_C \Lunit
	\right\rbrace, 
	\label{E:COVGEOMETRICTORIFRAMETOROUGHTORIFRAME} 
			\\
f_A & = (\COVframe^{-1})_{AB} 
				\left\lbrace
					e_B - \COVL_B \Lunit
				\right\rbrace
				=
				\COVframe_{BA}
				\left\lbrace
					e_B - \COVL_B \Lunit
				\right\rbrace. \label{E:COVROUGHTORIFRAMETOGEOMETRICTORIFRAME}
\end{align}
\end{subequations}

Moreover, $\gtorusroughinversefirstfund$ and $\gtorus^{-1}$ are related
through the following identity:
\begin{align} \label{E:RELATIONSHIPBETWEENINVERSEFIRSTFUNDSOFSMOOTHANDROUGHTORI}
	\gtorusroughinversefirstfund
	& = 
	\gtorus^{-1}
	-
	\COVL_A e_A \otimes \Lunit
	-
	\COVL_A \Lunit \otimes e_A
	+
	\COVL_A 
	\COVL_A
	 \Lunit \otimes \Lunit.
\end{align}

Finally, the following estimates 
hold on $\twoargMrough{[\timefunction_0,\timefunctionboot),[- \rightu,\leftu]}{\muxmulevelsetvalue}$:
\begin{subequations} 
\begin{align}
|\COVframe_{AB}|
& \leq 1, \label{E:COVFRAMEABSIZE}
	\\
|\COVL_A| & \leq C \fundbootsmall.
	\label{E:COVFRAMELSIZE}
\end{align}
\end{subequations}
\end{lemma}

\begin{proof}
The existence of a matrix $\COVframe$ and scalar functions $\{\COVL_A\}_{A = 2,3}$ satisfying 
\eqref{E:COVGEOMETRICTORIFRAMETOROUGHTORIFRAME} follows 
from the fact that the frames $\{e_2,e_3,\Lunit\}$ and $\{f_2,f_3,\Lunit\}$ both span the tangent space 
of $\nullhyparg{u}$ at any of its points. 
To see that $\COVframe$ is orthogonal,
we use the fact that $\Lunit$ is $\gfour$-orthogonal to the tangent space of $\nullhyparg{u}$
and the fact that both frames are $\gfour$-orthonormal
to deduce, with $\updelta_{AB}$ denoting the Kronecker delta, that:
\begin{align} \label{E:CHANGEOFFRAMEMATRIXISORTHOGONAL}
\updelta_{AB} 
& 
= 
\gtorus(e_A,e_B) = \gfour(e_A,e_B) = 
\COVframe_{AC} \COVframe_{BD} 
\gfour(f_C,f_D) =  
\COVframe_{AC}\COVframe_{BD} \gtorusroughfirstfund (f_C,f_D) 
= 
\COVframe_{AC}\COVframe_{BD} \updelta_{CD}
=
\COVframe_{AC} \COVframe_{BC}.
\end{align}
From \eqref{E:CHANGEOFFRAMEMATRIXISORTHOGONAL}, 
we see that
$(\COVframe^{-1})_{AB} = \COVframe_{BA}$ (i.e., $\COVframe$ is an orthogonal matrix)
and thus \eqref{E:COVROUGHTORIFRAMETOGEOMETRICTORIFRAME} follows. 

\eqref{E:RELATIONSHIPBETWEENINVERSEFIRSTFUNDSOFSMOOTHANDROUGHTORI}
follows from using
\eqref{E:COVROUGHTORIFRAMETOGEOMETRICTORIFRAME}
to substitute for the frame vectorfields $\lbrace f_A \rbrace_{A=2,3}$ in
\eqref{E:INVERSEFIRSTFUNDOFSMOOTHTORIRELATIVETOORTHONORMALFRAME},
using the orthogonality of the matrix $\COVframe_{AB}$, 
and comparing with \eqref{E:INVERSEFIRSTFUNDOFROUGHTORIRELATIVETOORTHONORMALFRAME}.

The estimate \eqref{E:COVFRAMEABSIZE} follows trivially since any orthogonal matrix has Euclidean-orthonormal rows 
and thus its entries are $\leq 1$ in magnitude.

To derive the estimate \eqref{E:COVFRAMELSIZE}, 
we first use \eqref{E:COVGEOMETRICTORIFRAMETOROUGHTORIFRAME} and
Lemma~\ref{L:BASICPROPERTIESOFVECTORFIELDS}
to deduce that
$
0 = \gfour(e_A,X)
=
\COVframe_{AB} \gfour(f_B,X)
- 
\COVL_A
$. From this identity and \eqref{E:COVFRAMEABSIZE},
we see that $|\COVL_A| \lesssim |\gfour(f_B,X)|$. 
Hence, 
\eqref{E:COVFRAMELSIZE} will follow once we show that:
\begin{align} \label{E:INNERPRODUCTOFFBANDXISSMALL}
	|\gfour(f_B,X)| \lesssim \fundbootsmall.
\end{align}
To prove \eqref{E:INNERPRODUCTOFFBANDXISSMALL}, 
we first use 
\eqref{E:SMOOTHTORIGABEXPRESSION},
Lemma~\ref{L:SCHEMATICSTRUCTUREOFVARIOUSTENSORSINTERMSOFCONTROLVARS},
the bootstrap assumptions, and Cor.\,\ref{C:IMPROVEAUX}
to deduce that
$\gtorus_{AB} 
= 
\Speed^{-2} \updelta_{AB} + \mathcal{O}(\fundbootsmall)
= 
\left\lbrace
	1 + \mathcal{O}(\mathring{\upalpha})
\right\rbrace
\updelta_{AB}
+ 
\mathcal{O}(\fundbootsmall)
$,
where $\updelta_{AB}$ is the Kronecker delta.
Also using \eqref{E:GTORUSROUGHCOMPONENTS}
and the estimates of Lemma~\ref{L:DIFFEOMORPHICEXTENSIONOFROUGHCOORDINATES}
for $\timefunctionarg{\muxmulevelsetvalue}$,
we find that 
$\gtorusroughfirstfund\left(\roughgeop{x^A},\roughgeop{x^B}\right)
=
\left\lbrace
	1 + \mathcal{O}(\mathring{\upalpha})
\right\rbrace
\updelta_{AB} 
+ 
\mathcal{O}(\fundbootsmall)$.
Next, since $\left\lbrace \roughgeop{x^2}, \roughgeop{x^3} \right\rbrace$ 
spans the tangent space of
the rough tori $\twoargroughtori{\timefunction,u}{\muxmulevelsetvalue}$,
for $A=2,3$,
there exist scalar functions $\upalpha_A$ and $\upbeta_A$
such that
$f_A = \upalpha_A \roughgeop{x^2} + \upbeta_A \roughgeop{x^3}$.
Since 
$1 = \gtorusroughfirstfund(f_A,f_A)$ 
(with no summation over $A$) by assumption,
it follows from the estimate 
$\gtorusroughfirstfund\left(\roughgeop{x^A},\roughgeop{x^B}\right)
=
\left\lbrace
	1 + \mathcal{O}(\mathring{\upalpha})
\right\rbrace
\updelta_{AB}
+ 
\mathcal{O}(\fundbootsmall)
$
that 
$|\upalpha_A|, \, |\upbeta_A| \lesssim 1$,
i.e., $f_A = \mathcal{O}(1) \roughgeop{x^2} +\mathcal{O}(1) \roughgeop{x^3}$.
From this relation,
\eqref{E:ROUGHANGULARPARTIALDERIVATIVESINTERMSOFGOODGEOMETRICPARTIALDERIVATIVES},
Lemma~\ref{L:SCHEMATICSTRUCTUREOFVARIOUSTENSORSINTERMSOFCONTROLVARS},
the bootstrap assumptions, 
Lemma~\ref{L:DIFFEOMORPHICEXTENSIONOFROUGHCOORDINATES},
and Cor.\,\ref{C:IMPROVEAUX},
we deduce that
$f_A = \mathcal{O}(1) \geop{x^2} + \mathcal{O}(1) \geop{x^3}
+
\mathcal{O}\left(\frac{\geop{x^A} \timefunction}{\geop{t} \timefunctionarg{\muxmulevelsetvalue}} \right) \Lunit
=
\mathcal{O}(1) \geop{x^2} + \mathcal{O}(1) \geop{x^3}
+ 
\mathcal{O}(\fundbootsmall) \Lunit
$.
From this identity and Lemma~\ref{L:BASICPROPERTIESOFVECTORFIELDS}, we conclude
\eqref{E:INNERPRODUCTOFFBANDXISSMALL}.

\end{proof}

\subsection{An alternate representation of $\gfour$-orthogonal projection onto the rough tori}
\label{SS:PROJECTIONOFSMOOTHTORITENSORONTOROUGHTORI}

\begin{definition}[$\gfour$-orthogonal projection onto the rough tori $\twoargroughtori{\timefunction,u}{\muxmulevelsetvalue}$] 
\label{D:PROJECTIONOFSMOOTHTORITENSORONTOROUGHTORI}
Let $\upxi$ be a symmetric type-$\binom{0}{2}$ $\ell_{t,u}$-tangent tensorfield.
We define $\roughangxi$ to be the symmetric type $\binom{0}{2}$ $\twoargroughtori{\timefunction,u}{\muxmulevelsetvalue}$-tangent tensorfield
whose type $\binom{2}{0}$ $\gfour$-dual, which we denote by $\roughangxitwouparg{\#}{\#}$, 
has the following components relative to the orthonormal frame $\lbrace f_2, f_3 \rbrace$
on $\twoargroughtori{\timefunction,u}{\muxmulevelsetvalue}$
from Def.\,\ref{D:ORTHONORMALFRAMESONSMOOTHANDROUGHTORI}:
\begin{align} \label{E:ROUGHTORICOMPONENTSOFFLATTORITENSORS}
\roughangxitwouparg{\#}{\#}
& 
\eqdef \upxi(f_A,f_B) f_A \otimes f_B. 
\end{align}
\end{definition}

\begin{remark}[\eqref{E:ROUGHTORICOMPONENTSOFFLATTORITENSORS} is
	$\gfour$-orthogonal projection onto $\twoargroughtori{\timefunction,u}{\muxmulevelsetvalue}$]
	\label{R:WIDETILDEXIISACTUALLYGORTHOGONALPROJECTIONONTOROUGHTORI}
	It is straightforward to check that the 
	tensorfield $\roughangxi$ defined by \eqref{E:ROUGHTORICOMPONENTSOFFLATTORITENSORS}
	is the $\gfour$-orthogonal projection of $\upxi$ onto
	$\twoargroughtori{\timefunction,u}{\muxmulevelsetvalue}$,
	i.e.,
	$\roughangxi = \roughtorusproject \upxi$,
	where $\roughtorusproject \upxi$ is defined by \eqref{E:PROJECTIONOFTENSORONTOROUGHTORI}.
\end{remark}

\subsection{Identities involving symmetric type $\binom{0}{2}$ tensorfields}
\label{SS:IDENTITIESINVOLVINGXIANDITSPROJECTIONONTOTHEROUGHTORI}
In the remainder of Sect.\,\ref{S:ELLIPTICESTIAMTESACOUSTICGEOMETRYONROUGHTORI}, 
we will work with $\nullhyparg{u}$-tangent tensorfields, as defined in
Def.\,\ref{D:CHARACTERISTICHYPERSURFACEPROJECTIONTENSORFIELD}.
It is straightforward to check that $\upeta$ is	
$\nullhyparg{u}$-tangent if and only if
any contraction of it with $\Lunit$ (which is $\gfour$-orthogonal to $\nullhyparg{u}$)
vanishes.

\begin{lemma}[Identities involving symmetric type $\binom{0}{2}$ tensorfields] 
\label{L:IDENTITIESINVOLVINGXIANDITSPROJECTIONONTOTHEROUGHTORI}
Let $\upxi$ be a symmetric type $\binom{0}{2}$ $\ell_{t,u}$-tangent tensorfield,
and let $\roughangxi$ be the corresponding $\twoargroughtori{\timefunction,u}{\muxmulevelsetvalue}$-tangent tensorfield from
Def.\,\ref{D:PROJECTIONOFSMOOTHTORITENSORONTOROUGHTORI}.
Then the following identity holds:
\begin{align}
\begin{split} \label{E:ROUGHTORICOMPONENTSINTERMSOFSMOOTHTORITANGENTTENSORANDLERRORS}
	\widetilde{\upxi}^{\# \#}
	& 
	=
	\upxi^{\# \#}
	-
	\COVL_A 
	\upxi(e_A,e_B) \Lunit 
	\otimes 
	e_B 
	-
	\COVL_B 
	\upxi(e_A,e_B) 
	e_A 
	\otimes 
	\Lunit
	+
	\COVL_A
	\COVL_B
	\upxi(e_A,e_B) 
	\Lunit 
	\otimes 
	\Lunit.
\end{split}
\end{align}

Moreover, if $\upeta$ is any symmetric type $\binom{0}{2}$ $\nullhyparg{u}$-tangent tensorfield,
then the following identities hold, 
where $\Dfour$ is the Levi-Civita connection of $\gfour$
and $\angLie_{\Lunit}$ is $\ell_{t,u}$-projected Lie derivative operator from
Def.\,\ref{D:PROJECTEDLIEDERIVATIVES}:
\begin{subequations}
\begin{align}
		[\Dfour \upeta](\Lunit,\Lunit)
	& = 0,
			\label{E:LLCOMPONENTOFFULLCOVARIANTDERIVATIVEOFLNORMALTENSORFIELDVANISHES} 
				\\
	[\Dfour_{\Lunit} \upeta](\Lunit,\cdot)
	& = 0,
			\label{E:LCOMPONENTOFLCOVARIANTDERIVATIVEOFLNORMALTENSORFIELDVANISHES} 
				\\
	[\Dfour_{e_A} \upeta](\Lunit,e_B)
	& = 
	- 
	\upeta(e_C,e_B) \upchi(e_A,e_C),
		\label{E:LCOMPONENTOFSMOOTHTORUSCOMPONENTCOVARIANTDERIVATIVEOFLNORMALTENSORFIELDEXPRESSIBLEINTERMSOFNULLSECONDFUND} 
	\end{align}
	\end{subequations}

	\begin{align} \label{E:SMOOTHTORUSCOMPONENTOFLCOVARIANTDERIVATIVEOFLNORMALTENSORFIELDEXPRESSIBLEINTERMSOFNULLSECONDFUND} 
		[\Dfour_{\Lunit} \upeta](e_A,e_B)
		& = 
		[\angLie_{\Lunit} \upeta](e_A,e_B)
		-
		\upeta[e_C,e_B] \upchi(e_A,e_C)
		-
		\upeta[e_A,e_C] \upchi(e_B,e_C).
	\end{align}

\end{lemma}

\begin{proof}
The identity \eqref{E:ROUGHTORICOMPONENTSINTERMSOFSMOOTHTORITANGENTTENSORANDLERRORS}
follows from using \eqref{E:COVROUGHTORIFRAMETOGEOMETRICTORIFRAME}
to substitute for the frame vectorfields $\lbrace f_A \rbrace_{A=2,3}$ in \eqref{E:ROUGHTORICOMPONENTSOFFLATTORITENSORS}
and using that $\upxi(\Lunit,\cdot) = 0$.

To prove \eqref{E:LLCOMPONENTOFFULLCOVARIANTDERIVATIVEOFLNORMALTENSORFIELDVANISHES},
we differentiate the identity $\upeta(\Lunit,\Lunit) = 0$ with $\Dfour$, use the Leibniz rule,
and use that $\upeta(\Lunit,\cdot) = 0$.

To prove \eqref{E:LCOMPONENTOFLCOVARIANTDERIVATIVEOFLNORMALTENSORFIELDVANISHES},
we differentiate the identity $\upeta(\Lunit,\cdot) = 0$ with $\Dfour_{\Lunit}$, 
use the Leibniz rule and the identity \eqref{E:COVARIANTLUNITDERIVATIVEOFLUNIT},
and use that $\upeta(\Lunit,\cdot) = 0$.

To prove \eqref{E:LCOMPONENTOFSMOOTHTORUSCOMPONENTCOVARIANTDERIVATIVEOFLNORMALTENSORFIELDEXPRESSIBLEINTERMSOFNULLSECONDFUND},
we differentiate the identity $\upeta(\Lunit,\cdot) = 0$ with $\Dfour_{e_A}$, 
use the Leibniz rule and the identity \eqref{E:SECONDFUNDSALTERNATE},
and use that $\upeta(\Lunit,\cdot) = 0$.

To prove \eqref{E:SMOOTHTORUSCOMPONENTOFLCOVARIANTDERIVATIVEOFLNORMALTENSORFIELDEXPRESSIBLEINTERMSOFNULLSECONDFUND},
we contract the Lie differentiation identity
$
	\Lie_{\Lunit} \upxi_{\alpha \beta}
	= 
	\Dfour_{\Lunit} \upxi_{\alpha \beta}
	+
	\upxi_{\kappa \beta} \Dfour_{\alpha} \Lunit^{\kappa}
	+
	\upxi_{\alpha \kappa} \Dfour_{\beta} \Lunit^{\kappa}
$
against $e_A^{\alpha} e_B^{\beta}$,
use that $\upeta(\Lunit,\cdot) = 0$,
and use the identity $\Dfour_{e_A} \Lunit = \upchi(e_A,e_C) e_C - \upzeta_A \Lunit$,
which follows from 
Lemma~\ref{L:BASICPROPERTIESOFVECTORFIELDS},
\eqref{E:TORISONTENSORFIELD},
and
\eqref{E:SECONDFUNDSALTERNATE}.

\end{proof}

\subsection{Identities connecting $\gtorus$, $\gtorusroughfirstfund$, and their corresponding differential operators}
\label{SS:IDENTITIESCONNECTINGTWOFIRSTFUNDFORMSANDTHEIRDIFFOPPERATORS}

\begin{lemma}[Identities involving contractions with $\gtorus$ and $\gtorusroughfirstfund$]
\label{L:SMOOTHANDROUGHFIRSTFUNDFORMSAGREEONPUTANGENTENSORS}
Let $\upeta$ be a $\nullhyparg{u}$-tangent tensorfield.
Recall that $|\cdot|_{\gtorus}$ and $|\cdot|_{\gtorusroughfirstfund}$
are defined in \eqref{E:SQUAREPOINTWISESEMINORMWITHRESPECTTOFIRSTFUNDOFSMOOTHTORI} and
\eqref{E:SQUAREPOINTWISESEMINORMWITHRESPECTTOFIRSTFUNDOFROUGHTORI}
respectively. The following identity holds:
\begin{align} \label{E:SEMINNORMSOFPUTANGENTTENSORSINDEPENDENTOFMETRIC}
|\upeta|_{\gtorus} 
& 
= |\upeta|_{\gtorusroughfirstfund}. 
\end{align}

Moreover, if $\upxi$ and $\roughangxi$ are as in Def.\,\ref{D:PROJECTIONOFSMOOTHTORITENSORONTOROUGHTORI}, then:
\begin{align}
|\upxi|_{\gtorus} & = |\roughangxi|_{\gtorusroughfirstfund}. 
\label{E:NORMSOFROUGHTORICOMPONENTSANDFLATTORIAREEQUAL}
\end{align}

In addition, if $\upeta$ is a symmetric type $\binom{0}{2}$ $\nullhyparg{u}$-tangent tensorfield, then
with $\mytr_{\gtorus}$ and $\mytr_{\gtorusroughfirstfund}$ as defined in
\eqref{E:SMOOTHTORUSTRACE}
and 
\eqref{E:TRACEOFROUGHTORITANGENT02TENSORS} respectively, the following identity holds:
\begin{align} \label{E:TRACESOFPUTANGENTTENSORINDEPENDENTOFMETRIC} 
\mytr_{\gtorus} \upeta 
& = 
\mytr_{\gtorusroughfirstfund} \upeta. 
\end{align}

Finally, if $\upxi$ and $\roughangxi$ are as in Def.\,\ref{D:PROJECTIONOFSMOOTHTORITENSORONTOROUGHTORI}, then:
\begin{align} \label{E:TRACESOFROUGHTORICOMPONENTSANDFLATTORIAREEQUAL} 
\mytr_{\gtorus} \upxi 
& = 
\mytr_{\gtorusroughfirstfund} \roughangxi. 
\end{align}
	
\end{lemma}

\begin{proof}
	The lemma follows easily from the fact that 
	any $\nullhyparg{u}$-tangent tensorfield
	(including $\gtorus$
	and
	$\gtorusroughfirstfund$)
	vanishes upon any contraction with $\Lunit$,
	and the identities
	\eqref{E:RELATIONSHIPBETWEENINVERSEFIRSTFUNDSOFSMOOTHANDROUGHTORI}
	and
	\eqref{E:ROUGHTORICOMPONENTSINTERMSOFSMOOTHTORITANGENTTENSORANDLERRORS}.
\end{proof}

\subsection{Differential operator pointwise comparison estimates needed for the elliptic estimates}
\label{SS:DIFFOPERATORSPOINTWISECOMPARISONNEEDEDFORELLIPTICESTIMATESONROUGHTORI}

\begin{lemma}[Differential operator pointwise comparison estimates needed for the elliptic estimates]
\label{L:DIFFOPERATORSPOINTWISECOMPARISONNEEDEDFORELLIPTICESTIMATESONROUGHTORI}
Let $\varphi$ be a scalar function, 
let $\angrmD \varphi$ be the $\ell_{t,u}$-tangent one-form from Def.\,\ref{D:ANGULARDIFFERENTIAL},
and let $\roughangrmD \varphi$ be the $\twoargroughtori{\timefunction,u}{\muxmulevelsetvalue}$-tangent one-form
from Def.\,\ref{D:ROUGHTORUSDIFFERENTIAL}.
Let $\upxi$ be a symmetric type $\binom{0}{2}$ $\ell_{t,u}$-tangent tensorfield,
and let $\roughangxi$ 
be the corresponding symmetric type $\binom{0}{2}$ $\twoargroughtori{\timefunction,u}{\muxmulevelsetvalue}$-tangent tensorfield
from Def.\,\ref{D:PROJECTIONOFSMOOTHTORITENSORONTOROUGHTORI}.
Then the following pointwise estimates hold
on $\twoargMrough{[\timefunction_0,\timefunctionboot),[- \rightu,\leftu]}{\muxmulevelsetvalue}$:
\begin{subequations}
\begin{align}
	|\roughangrmD \varphi|_{\gtorusroughfirstfund} 
	&
	=
	\left\lbrace
		1 + \mathcal{O}(\fundbootsmall)
	\right\rbrace
	|\angrmD \varphi|_{\gtorus}
	+ 
	\mathcal{O}(\fundbootsmall)
	|\Lunit \varphi|,
		\label{E:ROUGHTORUSDIFFERENTIALAPPROXIMATEDBYSMOOTHTORUSDIFFERENTIAL} 
			\\
	|\roughangdiv \roughangxi|_{\gtorusroughfirstfund}
	& =
	\left\lbrace
		1 + \mathcal{O}(\fundbootsmall)
	\right\rbrace
	|\angdiv \upxi|_{\gtorus}
	+
	\mathcal{O}(\fundbootsmall)
	|\angLie_{\Lunit} \upxi|_{\gtorus}
	+
	\mathcal{O}(\fundbootsmall)
	|\upxi|_{\gtorus},
		\label{E:ROUGHTORUSDIVERGENCEAPPROXIMATEDBYSMOOTHTORUSDIVERGENCE} 
			\\
	|\roughangD \roughangxi|_{\gtorusroughfirstfund}
	& =
	\left\lbrace
		1 + \mathcal{O}(\fundbootsmall)
	\right\rbrace
	|\angD \upxi|_{\gtorus}
	+
	\mathcal{O}(\fundbootsmall)
	|\angLie_{\Lunit} \upxi|_{\gtorus}
	+
	\mathcal{O}(\fundbootsmall)
	|\upxi|_{\gtorus}.
	\label{E:ROUGHTORUSCOVARIANTDERIVATIVEAPPROXIMATEDBYSMOOTHTORUSCOVARIANTDERIVATIVE}
\end{align}
\end{subequations}

\end{lemma}

\begin{proof}
We prove only \eqref{E:ROUGHTORUSDIVERGENCEAPPROXIMATEDBYSMOOTHTORUSDIVERGENCE}
since \eqref{E:ROUGHTORUSDIFFERENTIALAPPROXIMATEDBYSMOOTHTORUSDIFFERENTIAL} 
and \eqref{E:ROUGHTORUSCOVARIANTDERIVATIVEAPPROXIMATEDBYSMOOTHTORUSCOVARIANTDERIVATIVE}
can be proved by similar arguments.

We start by noting the following identities,
where $\COVframe_{AB}$ is the orthogonal matrix from Lemma~\ref{L:CHANGEOFORTHONORMALFRAMES}:
\begin{align} \label{E:ANGDIVSMOOTHTORITANGENTTENSORNORMSQUAREDINTERMSOFSMOOTHORTHONORMALFRAME}
		|\angdiv \upxi|_{\gtorus}^2 
		& 
		= 
		[\Dfour_{e_A}\upxi](e_A,e_C)[\Dfour_{e_B}\upxi](e_B,e_C),
			\\
	|\roughangdiv \roughangxi|_{\gtorusroughfirstfund}^2  
		& 
		= 
		[\Dfour_{f_A} \roughangxi](f_A,f_C)
		[\Dfour_{f_B} \roughangxi](f_B,f_C)
		= 
		\left\lbrace
			\COVframe_{CD}
			[\Dfour_{f_A} \roughangxi](f_A,f_D)
		\right\rbrace
		\left\lbrace
		\COVframe_{CE}
		[\Dfour_{f_B} \roughangxi](f_B,f_E)
		\right\rbrace.
		 \label{E:ANGDIVROUGHTORITANGENTTENSORNORMSQUAREDINTERMSOFROUGHORTHONORMALFRAME}
\end{align}
	The identity \eqref{E:ANGDIVSMOOTHTORITANGENTTENSORNORMSQUAREDINTERMSOFSMOOTHORTHONORMALFRAME}
	follows from \eqref{E:INVERSEFIRSTFUNDOFSMOOTHTORIRELATIVETOORTHONORMALFRAME} 
	and the fact that $\angD \upxi = \smoothtorusproject \Dfour \upxi$,
	where $\smoothtorusproject$ is the $\gfour$-orthogonal projection onto $\ell_{t,u}$
	from Def.\,\ref{D:PROJECTIONTENSORFIELDSANDTANGENCYTOHYPERSURFACES}.
	Similarly, the first equality in \eqref{E:ANGDIVROUGHTORITANGENTTENSORNORMSQUAREDINTERMSOFROUGHORTHONORMALFRAME}
	follows from \eqref{E:INVERSEFIRSTFUNDOFROUGHTORIRELATIVETOORTHONORMALFRAME}
	and the fact that  $\roughangD \roughangxi = \roughtorusproject \Dfour \roughangxi$,
	where $\roughtorusproject$ is the $\gfour$-orthogonal projection onto $\twoargroughtori{\timefunction,u}{\muxmulevelsetvalue}$
	from Def.\,\ref{D:PROJECTIONONTOROUGHTORIANDROUGHTORITANGENCY}.
	The second equality in \eqref{E:ANGDIVROUGHTORITANGENTTENSORNORMSQUAREDINTERMSOFROUGHORTHONORMALFRAME}
	follows from the orthogonality of $\COVframe$.

We now note the following identity, which we derive just below:
\begin{align}
	\begin{split} \label{E:FDERIVATIVESOFROUGHTORITANGENTTENSORINTERMSOFEDERIVATIVESOFSMOOTHTORITANGENTTENSOR}
		\COVframe_{CD}
		[\Dfour_{f_A} \roughangxi](f_A,f_D)
		& 
		= 
		[\Dfour_{e_A} \upxi](e_A,e_C)
		-
		\COVL_A
		[\angLie_{\Lunit} \upxi](e_A,e_C)
			\\
	& 	
		\	\
		-
		\COVL_D 
		\upxi(e_D,e_C) 
		\mytr_{\gtorus} \upchi
		+
		\COVL_C 
		\upxi(e_A,e_D)\upchi(e_A,e_D)
		+
		2
		\COVL_A
		\upxi(e_C,e_D)\upchi(e_A,e_D).
	\end{split}
	\end{align}	
From \eqref{E:FDERIVATIVESOFROUGHTORITANGENTTENSORINTERMSOFEDERIVATIVESOFSMOOTHTORITANGENTTENSOR},
the identities 
\eqref{E:INVERSEFIRSTFUNDOFSMOOTHTORIRELATIVETOORTHONORMALFRAME}--\eqref{E:INVERSEFIRSTFUNDOFROUGHTORIRELATIVETOORTHONORMALFRAME}
and
\eqref{E:ANGDIVSMOOTHTORITANGENTTENSORNORMSQUAREDINTERMSOFSMOOTHORTHONORMALFRAME}--\eqref{E:ANGDIVROUGHTORITANGENTTENSORNORMSQUAREDINTERMSOFROUGHORTHONORMALFRAME},
the orthogonality of the matrix $\COVframe$,
the Cauchy--Schwarz inequality with respect to $\gtorus$,
Lemma~\ref{L:SMOOTHANDROUGHFIRSTFUNDFORMSAGREEONPUTANGENTENSORS},
the estimate \eqref{E:COVFRAMELSIZE},
and the pointwise estimate $|\upchi|_{\gtorus} \lesssim \fundbootsmall$
noted below \eqref{E:PNDEFPROOFOFPRELIMINARYEIKONALL2},
we conclude \eqref{E:ROUGHTORUSDIVERGENCEAPPROXIMATEDBYSMOOTHTORUSDIVERGENCE}.

It remains for us to prove \eqref{E:FDERIVATIVESOFROUGHTORITANGENTTENSORINTERMSOFEDERIVATIVESOFSMOOTHTORITANGENTTENSOR}.
To proceed, we first use
\eqref{E:COVARIANTLUNITDERIVATIVEOFLUNIT},
\eqref{E:COVGEOMETRICTORIFRAMETOROUGHTORIFRAME}--\eqref{E:COVROUGHTORIFRAMETOGEOMETRICTORIFRAME},
\eqref{E:ROUGHTORICOMPONENTSINTERMSOFSMOOTHTORITANGENTTENSORANDLERRORS},
\eqref{E:LLCOMPONENTOFFULLCOVARIANTDERIVATIVEOFLNORMALTENSORFIELDVANISHES}--\eqref{E:LCOMPONENTOFSMOOTHTORUSCOMPONENTCOVARIANTDERIVATIVEOFLNORMALTENSORFIELDEXPRESSIBLEINTERMSOFNULLSECONDFUND},
\eqref{E:SMOOTHTORUSCOMPONENTOFLCOVARIANTDERIVATIVEOFLNORMALTENSORFIELDEXPRESSIBLEINTERMSOFNULLSECONDFUND},
the orthogonality of the matrix $\COVframe_{AB}$,
the fact that $\Lunit$ is null and $\gfour$-orthogonal to $\nullhyparg{u}$,
the fact that $\roughangxi$ and $\upxi$ are symmetric and satisfy
$\roughangxi(\Lunit,\cdot) = \upxi(\Lunit,\cdot) = 0$,
and the fact that $\angLie_{\Lunit} \Lunit = 0$
to compute that:
	\begin{align}
	\begin{split} \label{E:FIRSTPROOFSTEPROUGHTORUSDIVERGENCEAPPROXIMATEDBYSMOOTHTORUSDIVERGENCE}
		[\Dfour_{e_A} \roughangxi](e_A,e_C)
		& = 
			\COVframe_{CD}
			[\Dfour_{f_A} \roughangxi](f_A,f_D)
				\\
		& \ \
			+
			\COVL_C 
			[\Dfour_{e_A} \roughangxi](e_A,\Lunit)
			+
			\COVL_A
			[\Dfour_{e_A} \roughangxi](\Lunit,e_C)
			-
			\COVL_C
			\COVL_A
			[\Dfour_{e_A} \roughangxi](\Lunit,\Lunit)
					\\
		& \ \
			+
			\COVL_A
			[\Dfour_{\Lunit} \roughangxi](e_A,e_C)
			-
			\COVL_A
			\COVL_A
			[\Dfour_{\Lunit} \roughangxi](\Lunit,e_C)
			-
			\COVL_C
			\COVL_A
			[\Dfour_{\Lunit} \roughangxi](e_A,\Lunit)
			+
			\COVL_C
			\COVL_A
			\COVL_A
			[\Dfour_{\Lunit} \roughangxi](\Lunit,\Lunit)
					\\
		& = 
			\COVframe_{CD}
			[\Dfour_{f_A} \roughangxi](f_A,f_D)
			+
			\COVL_A
			[\angLie_{\Lunit} \upxi](e_A,e_C)	
				\\
		& \ \
			-
			\COVL_C 
			\upxi(e_A,e_D) \upchi(e_A,e_D)
			-
			2
			\COVL_A
			\upxi(e_C,e_D)\upchi(e_A,e_D)
			-
			\COVL_A
			\upxi(e_A,e_D)\upchi(e_C,e_D).
		\end{split}
	\end{align}
	Next, using 
	\eqref{E:ROUGHTORICOMPONENTSINTERMSOFSMOOTHTORITANGENTTENSORANDLERRORS},
	the fact that $\Lunit$ is $\gfour$-orthogonal to $\nullhyparg{u}$,
	the fact that $\roughangxi(\Lunit,\cdot) = 0$, 
	the identity $\gfour(\Dfour_{e_A}\Lunit,e_B) = \upchi(e_A,e_B)$
	(which follows from \eqref{E:SECONDFUNDSALTERNATE}),
	and the identity $\upchi(e_A,e_A) = \mytr_{\gtorus} \upchi$,
	we deduce:
	\begin{align} \label{E:SECONDPROOFSTEPROUGHTORUSDIVERGENCEAPPROXIMATEDBYSMOOTHTORUSDIVERGENCE}
	[\Dfour_{e_A} \roughangxi](e_A,e_C)
	& =
	[\Dfour_{e_A} \upxi](e_A,e_C)
		-
	\COVL_D 
	\upxi(e_D,e_C) 
	\mytr_{\gtorus} \upchi
	-
	\COVL_D 
	\upxi(e_A,e_D) 
	\upchi(e_A,e_C).
	\end{align}
	Using \eqref{E:SECONDPROOFSTEPROUGHTORUSDIVERGENCEAPPROXIMATEDBYSMOOTHTORUSDIVERGENCE}
	to substitute for LHS~\eqref{E:FIRSTPROOFSTEPROUGHTORUSDIVERGENCEAPPROXIMATEDBYSMOOTHTORUSDIVERGENCE}
	and rearranging terms,
	we conclude \eqref{E:FDERIVATIVESOFROUGHTORITANGENTTENSORINTERMSOFEDERIVATIVESOFSMOOTHTORITANGENTTENSOR}.

\end{proof}

\subsection{Standard elliptic estimates for symmetric type $\binom{0}{2}$ tensorfields on the rough tori}
\label{SS:STANDARDELLIPTICONROUGHTORI}
In the next lemma, we provide standard elliptic estimates for
symmetric type $\binom{0}{2}$ $\twoargroughtori{\timefunction,u}{\muxmulevelsetvalue}$-tangent tensorfields.
Its proof is located in Sect.\,\ref{SSS:PROOFOFL:STANDARDELLIPTICONROUGHTORI}.

\begin{lemma}[Standard elliptic estimates for symmetric type $\binom{0}{2}$ tensorfields on the rough tori]
\label{L:STANDARDELLIPTICONROUGHTORI}
Let $\Upxi$ be a symmetric type $\binom{0}{2}$ $\twoargroughtori{\timefunction,u}{\muxmulevelsetvalue}$-tangent tensorfield.
Then the following estimate
holds for $(\timefunction,u) \in [\timefunction_0,\timefunctionboot) \times [- \rightu,\leftu]$:
\begin{align} \label{E:STANDARDELLIPTICONROUGHTORI}
\int_{\twoargroughtori{\timefunction,u}{\muxmulevelsetvalue}} 
	\upmu^2 |\roughangD \Upxi|_{\gtorusroughfirstfund}^2 
\, \volroughtorus
& \leq 
6 
\int_{\twoargroughtori{\timefunction,u}{\muxmulevelsetvalue}} 
	\upmu^2 |\roughangdiv \Upxi|_{\gtorusroughfirstfund}^2 
\, \volroughtorus
+ 
3 
\int_{\twoargroughtori{\timefunction,u}{\muxmulevelsetvalue}} 
	\upmu^2 |\roughangD \mytr_{\gtorusroughfirstfund} \Upxi|_{\gtorusroughfirstfund}^2
\, \volroughtorus
+ 
C 
\fundbootsmall 
\int_{\twoargroughtori{\timefunction,u}{\muxmulevelsetvalue}} 
	|\Upxi|_{\gtorusroughfirstfund}^2
\, \volroughtorus.
\end{align}

\end{lemma}

\subsubsection{The Gauss curvature of $\gtorusroughfirstfund$}
\label{SSS:GAUSSCURVATUREOFROUGHTORI}
Our proof of Lemma~\ref{L:STANDARDELLIPTICONROUGHTORI} relies on the next lemma, 
which provides $L^{\infty}$ estimates for the Gauss curvature 
$\Gausstorus$ of $(\twoargroughtori{\timefunction,u}{\muxmulevelsetvalue},\gtorusroughfirstfund)$.

\begin{lemma}[$L^{\infty}$ estimates for the Gauss curvature of $\gtorusroughfirstfund$] 
\label{L:GAUSSCURVATUREOFROUGHTORILINFINITYESTIMATE}
Recall that $\Gausstorus$ denotes the Gauss curvature of $(\twoargroughtori{\timefunction,u}{\muxmulevelsetvalue},\gtorusroughfirstfund)$
(see Sect.\,\ref{SS:CURVATURETENSORS}).
Then the following estimate holds for
$(\timefunction,u) \in [\timefunction_0,\timefunctionboot] \times [-\rightu,\leftu]$:
\begin{align} \label{E:LINFINITYBOUNDFORROUGHGAUSSCURVATURE}
\left\| \Gausstorus \right \|_{L^{\infty}\left(\twoargroughtori{\timefunction,u}{\muxmulevelsetvalue} \right)} 
& \leq C \fundbootsmall.  
\end{align}
\end{lemma}

\begin{proof}
Recall that $\Gausstorus$ is equal to half the scalar curvature of $\gtorusroughfirstfund$
(see \eqref{E:2DGAUSSCURVATUREISTWICESCALARCURVATURE}).
Hence, at fixed $(\timefunction,u)$, relative to the coordinates $(x^2,x^3)$ on the rough tori $\twoargroughtori{\timefunction,u}{\muxmulevelsetvalue}$,
using the standard expression for curvature in terms of the 
coordinate components of $\gtorusroughfirstfund$ and their partial derivatives,
we can schematically express $\Gausstorus$ as follows, where 
$\gtorusroughinversefirstfund$ schematically denotes the component functions
$\gtorusroughinversefirstfund(\mathrm{d} x^A, \mathrm{d} x^B)$
and
$\gtorusroughfirstfund$
schematically denotes the component functions
$\gtorusroughfirstfund\left(\roughgeop{x^A},\roughgeop{x^B}\right)$:
\begin{align} \label{E:ROUGHTORIGAUSSCURVATURESCHMEATIC}
 \Gausstorus 
	& 
	= 
	\gtorusroughinversefirstfund \cdot \gtorusroughinversefirstfund \cdot \roughgeop{x^A} \roughgeop{x^B} \gtorusroughfirstfund 
	+ 
	\gtorusroughinversefirstfund \cdot \gtorusroughinversefirstfund \cdot \gtorusroughinversefirstfund 
	\cdot \roughgeop{x^A} \gtorusroughfirstfund 
	\cdot \roughgeop{x^B} \gtorusroughfirstfund.
 \end{align}
From \eqref{E:ROUGHTORIGAUSSCURVATURESCHMEATIC} and arguments similar to the ones we used to prove
\eqref{E:SCHEMATICIDENTITYFORROUGHDIVERGENCEOFROUGHTORITANGENTVECTORFIELD},
we further deduce that there is a smooth function $\smoothfunction$ such that
schematically, we have:
\begin{align} \label{E:MOREPRECISEROUGHTORIGAUSSCURVATURESCHMEATIC}
 \Gausstorus 
	& =
	\smoothfunction
	\left(\tander^{\leq 2} \controlvars,
	\frac{1}{\Lunit \timefunctionarg{\muxmulevelsetvalue}}, 
	\frac{1}{\geop{t} \timefunctionarg{\muxmulevelsetvalue}}, 
	\tander^{[1,3]} \timefunction 
\right)
\cdot
\left(
	\tander^{[1,2]} \controlvars,
		\,
	\tander^{\leq 2} \tanderY \timefunction,
\right)
\end{align}
From \eqref{E:MOREPRECISEROUGHTORIGAUSSCURVATURESCHMEATIC},
the results of Lemma~\ref{L:DIFFEOMORPHICEXTENSIONOFROUGHCOORDINATES},
including
the estimates 
\eqref{E:CLOSEDVERSIONKEYJACOBIANDETERMINANTESTIMATECHOVGEOTOROUGH},
\eqref{E:CLOSEDVERSIONLUNITROUGHTTIMEFUNCTION},
\eqref{E:CLOSEDVERSIONC21BOUNDFORCHOVROUGHTOGEO},
and 
\eqref{E:SMALLC11ESTIMATESFORROUGHTIMEFUNCTION},
Rademacher's theorem,
the bootstrap assumptions,
and Cor.\,\ref{C:IMPROVEAUX},
we arrive at the desired estimate
\eqref{E:LINFINITYBOUNDFORROUGHGAUSSCURVATURE}.

\end{proof}

\subsubsection{Proof of Lemma~\ref{L:STANDARDELLIPTICONROUGHTORI}}
\label{SSS:PROOFOFL:STANDARDELLIPTICONROUGHTORI}
We now prove Lemma~\ref{L:STANDARDELLIPTICONROUGHTORI}.
In this proof only, we will use capital Latin indices to denote
the components of $\twoargroughtori{\timefunction,u}{\muxmulevelsetvalue}$-tangent tensorfields
with respect to the frame $\left\lbrace \roughgeop{x^A} \right\rbrace_{A=2,3}$
and co-frame $\lbrace \roughangrmD x^A \rbrace_{A=2,3}$ on $\twoargroughtori{\timefunction,u}{\muxmulevelsetvalue}$,
and we raise and lower indices with $\gtorusroughinversefirstfund$ and $\gtorusroughfirstfund$.
In particular, 
$\Upxi 
= 
\Upxi\left(\roughgeop{x^A},\roughgeop{x^B} \right)
\roughangrmD x^A \otimes \roughangrmD x^B
$.
We start by defining $\widetilde{I} = \widetilde{I}^A \roughgeop{x^A}$ to be the
$\twoargroughtori{\timefunction,u}{\muxmulevelsetvalue}$-tangent vectorfield
with the following components relative to the coordinates $(x^2,x^3)$ on
$\twoargroughtori{\timefunction,u}{\muxmulevelsetvalue}$:
\begin{align} \label{E:CURRENTFORELLIPTICIDENTITYONROUGHTORIMANIFOLD}
	\widetilde{I}^A
	& 
	\eqdef  
	\upmu^2 \Upxi_{\widetilde B \widetilde C} \roughangDuparg{B} \Upxi^{AC} 
	- 
	\upmu^2 \Upxi^{AB} (\roughangdiv \Upxi)_B.
\end{align}
Next, with $\Gausstorus$ denoting the Gauss curvature of the (two-dimensional) rough tori
$\twoargroughtori{\timefunction,u}{\muxmulevelsetvalue}$,
we note the following standard identity,
which follows from the symmetry of $\Upxi$
(see \cite[Lemma 18.9]{jS2016b} for the main ideas of the proof,
where we note that only trace-free tensorfields were handled in
\cite[Lemma 18.9]{jS2016b} and thus RHSs~\eqref{E:STANDARDELLIPTICONROUGHTORI} and \eqref{E:ELLIPTICIDENTITYONROUGHTORIMANIFOLD} feature
additional $\mytr_{\gtorusroughfirstfund} \Upxi$-dependent terms compared to \cite[Lemma 18.9]{jS2016b}):
\begin{align}
\begin{split} \label{E:ELLIPTICIDENTITYONROUGHTORIMANIFOLD}
	\upmu^2 |\roughangD \Upxi|_{\gtorusroughfirstfund}^2 
	+ 
	2 \upmu^2 \Gausstorus |\Upxi|_{\gtorusroughfirstfund}^2
	& 
	= 
	2 \upmu^2 |\roughangdiv \Upxi|_{\gtorusroughfirstfund}^2 
	+ 
	\upmu^2 \Gausstorus (\mytr_{\gtorusroughfirstfund} \Upxi)^2 
	+
	\upmu^2 |\widetilde{\angrmD} \mytr_{\gtorusroughfirstfund} \Upxi|_{\gtorusroughfirstfund}^2  
		\\
& \ \ 
		- 
		2 \upmu^2 
		\gtorusroughinversefirstfund(\roughangdiv \Upxi,\widetilde{\angrmD} 
		\mytr_{\gtorusroughfirstfund} \Upxi) 
		+ 
		2 
		\upmu  
		\Upxi^{AB} 
		\left(\roughgeop{x^A} \upmu \right)
		(\roughangdiv \Upxi)_{\widetilde B} 
	\\
& \ \
- 
2 \upmu 
\Upxi_{\widetilde{B} \widetilde{C}} 
\left(\roughgeop{x^A} \upmu\right) 
\roughangDuparg{B} \Upxi^{AC} 
+ 
\roughangdiv \widetilde{I}.
\end{split}
\end{align}
We then integrate \eqref{E:ELLIPTICIDENTITYONROUGHTORIMANIFOLD} over $\twoargroughtori{\timefunction,u}{\muxmulevelsetvalue}$
with respect to the area form $\volroughtorus$ defined in \eqref{E:AREAFORMROUGHTORUS}
and note that the integral of the perfect divergence term $\roughangdiv \widetilde{I}$ vanishes.
Next, we use the $\gtorusroughfirstfund$-Cauchy--Schwarz inequality and
Young's inequality to pointwise bound the three cross terms on RHS~\eqref{E:ELLIPTICIDENTITYONROUGHTORIMANIFOLD}
as follows:
\begin{align}
2 \left| 
	\upmu^2 \gtorusroughinversefirstfund(\roughangdiv \Upxi, \roughangD \mytr_{\gtorusroughfirstfund} \Upxi) 
\right| 
& \leq 
	\upmu^2 |\roughangdiv \Upxi|_{\gtorusroughfirstfund}^2 
+ 
\upmu^2 
|\roughangD \mytr_{\gtorusroughfirstfund} \Upxi|_{\gtorusroughfirstfund}^2, 
		\label{E:ROUGHTORIELLIPTICELEMENTARYINEQUALITY1} 
			\\
2 \left| 
	\upmu \Upxi^{AB} \left(\roughgeop{x^A} \upmu \right)(\roughangdiv \Upxi)_{\widetilde B} 
\right| 
& \leq 
\upmu^2 
|\roughangdiv \Upxi|_{\gtorusroughfirstfund}^2 
+ 
|\roughangD \upmu|_{\gtorusroughfirstfund}^2 
|\Upxi|_{\gtorusroughfirstfund}^2, 
	\label{E:ROUGHTORIELLIPTICELEMENTARYINEQUALITY2}  
		\\
2 \left| \upmu \Upxi_{\widetilde{B} \widetilde{C}}  
\left(\roughgeop{x^A} \upmu \right) \roughangD^B \Upxi^{AC}\right| 
& 
\leq 
\frac{1}{3} 
\upmu^2 
|\roughangD \Upxi|_{\gtorusroughfirstfund}^2 
+ 
3 
|\widetilde{\angrmD} \upmu|_{\gtorusroughfirstfund}^2 
|\Upxi|_{\gtorusroughfirstfund}^2.
 \label{E:ROUGHTORIELLIPTICELEMENTARYINEQUALITY3} 
\end{align}
Just below, we will show that:
\begin{align} \label{E:POINTWISEBOUNDROUGHGRADIENTOFMU}
|\roughangD \upmu|_{\gtorusroughfirstfund}
&
\lesssim 
|\angD \upmu|_{\gtorus}
+
\fundbootsmall
|\Lunit \upmu|
\lesssim
|\tanderY \upmu|
+
\fundbootsmall
|\argLrough{\muxmulevelsetvalue} \upmu|
\lesssim
\fundbootsmall.
\end{align}
Using \eqref{E:POINTWISEBOUNDROUGHGRADIENTOFMU} to control the relevant factors on
RHSs~\eqref{E:ROUGHTORIELLIPTICELEMENTARYINEQUALITY2}--\eqref{E:ROUGHTORIELLIPTICELEMENTARYINEQUALITY3},
using the Gauss curvature estimate
\eqref{E:LINFINITYBOUNDFORROUGHGAUSSCURVATURE}
to control the factors of
$\Gausstorus$ in \eqref{E:ELLIPTICIDENTITYONROUGHTORIMANIFOLD},
using the elementary inequality 
$|\mytr_{\gtorusroughfirstfund} \Upxi| 
\lesssim 
|\Upxi|_{\gtorusroughfirstfund}$,
and using the estimate $|\upmu| \lesssim 1$ (which follows from the bootstrap assumptions),
we conclude \eqref{E:STANDARDELLIPTICONROUGHTORI}.

To prove \eqref{E:POINTWISEBOUNDROUGHGRADIENTOFMU},
we use \eqref{E:ROUGHTORUSDIFFERENTIALAPPROXIMATEDBYSMOOTHTORUSDIFFERENTIAL},
the bootstrap assumptions,
and Cor.\,\ref{C:IMPROVEAUX}
to deduce that 
$
|\roughangD \upmu|_{\gtorusroughfirstfund}
\lesssim 
|\angD \upmu|_{\gtorus}
+
\fundbootsmall
|\Lunit \upmu|
\lesssim
|\tanderY \upmu|
+
\fundbootsmall
|\argLrough{\muxmulevelsetvalue} \upmu|
\lesssim
\fundbootsmall
$
as desired.

\hfill $\qed$

\subsection{Proof of Prop.\,\ref{P:ELLIPTICESTIMATESONROUGHTORI}}
\label{SS:PROOFOFP:ELLIPTICESTIMATESONROUGHTORI}
We now prove Prop.\,\ref{P:ELLIPTICESTIMATESONROUGHTORI}.
Let $\upxi$ be a symmetric type-$\binom{0}{2}$ 
$\ell_{t,u}$-tangent tensorfield,
and let $\roughangxi$ 
be the corresponding symmetric type $\binom{0}{2}$ $\twoargroughtori{\timefunction,u}{\muxmulevelsetvalue}$-tangent tensorfield
from Def.\,\ref{D:PROJECTIONOFSMOOTHTORITENSORONTOROUGHTORI}.
Since the term
$
\int_{\hypthreearg{\timefunction}{[-\rightu,u]}{\muxmulevelsetvalue}} 
	\upmu^2 |\angLie_{\Lunit} \upxi|_{\gtorus}^2 
\, \volRoughHypersurface 
$
on LHS~\eqref{E:ELLIPTICESTIMATECHI} 
is manifestly bounded by RHS~\eqref{E:ELLIPTICESTIMATECHI}, 
we only have to show that for $A=2,3$,
the term
$
\int_{\hypthreearg{\timefunction}{[-\rightu,u]}{\muxmulevelsetvalue}} 
	\upmu^2 |\angLie_{\Yvf{A}} \upxi|_{\gtorus}^2 
\, \volRoughHypersurface 
$
on LHS~\eqref{E:ELLIPTICESTIMATECHI} is $\leq \mbox{RHS~\eqref{E:ELLIPTICESTIMATECHI}}$. 
To proceed, we consider the inequality \eqref{E:STANDARDELLIPTICONROUGHTORI} 
with $\roughangxi$ in the role of $\Upxi$.
Integrating the inequality with respect to $u'$ and 
using Lemmas~\ref{L:SMOOTHANDROUGHFIRSTFUNDFORMSAGREEONPUTANGENTENSORS}
and \ref{L:DIFFOPERATORSPOINTWISECOMPARISONNEEDEDFORELLIPTICESTIMATESONROUGHTORI},
we deduce, in view of definition \eqref{E:VOLUMEFORMROUGHHYPERSURFACE}, that:
\begin{align}
\begin{split} \label{E:FIRSTSTEPPROOFELLIPTICESTIMATECHI}
\int_{\hypthreearg{\timefunction}{[-\rightu,u]}{\muxmulevelsetvalue}}
	\upmu^2 |\angD \upxi|_{\gtorus}^2 
\, \volRoughHypersurface 
& \lesssim
\int_{\hypthreearg{\timefunction}{[-\rightu,u]}{\muxmulevelsetvalue}} 
	\upmu^2 |\angdiv \upxi|_{\gtorus}^2 
\, \volRoughHypersurface
+ 
\int_{\hypthreearg{\timefunction}{[-\rightu,u]}{\muxmulevelsetvalue}} 
	\upmu^2 |\angD \mytr_{\gtorus} \upxi|_{\gtorus}^2
\, \volRoughHypersurface
	\\
& \ \
+ 
\fundbootsmall
\int_{\hypthreearg{\timefunction}{[-\rightu,u]}{\muxmulevelsetvalue}} 
	\upmu^2 |\Lunit \mytr_{\gtorus} \upxi|_{\gtorus}^2
\, \volRoughHypersurface
+ 
\fundbootsmall
\int_{\hypthreearg{\timefunction}{[-\rightu,u]}{\muxmulevelsetvalue}} 
	\upmu^2 |\angLie_{\Lunit} \upxi|_{\gtorus}^2
\, \volRoughHypersurface
+ 
\fundbootsmall 
\int_{\hypthreearg{\timefunction}{[-\rightu,u]}{\muxmulevelsetvalue}}
	|\upxi|_{\gtorus}^2
\, \volRoughHypersurface.
\end{split}
\end{align}
Next, using the Leibniz rule, 
\eqref{E:TANDERGANDCHIESTIMATE},
and the bootstrap assumptions, we find that
$|\Lunit \mytr_{\gtorus} \upxi|_{\gtorus}
\lesssim 
|\angLie_{\Lunit} \upxi|_{\gtorus}
+
|\angLie_{\Lunit} \gtorus^{-1}|_{\gtorus}
|\upxi|_{\gtorus}
\lesssim 
|\angLie_{\Lunit} \upxi|_{\gtorus}
+
|\upxi|_{\gtorus}
$
and thus the third integral
$
\fundbootsmall
\int_{\hypthreearg{\timefunction}{[-\rightu,u]}{\muxmulevelsetvalue}} 
	\upmu^2 |\Lunit \mytr_{\gtorus} \upxi|_{\gtorus}^2
\, \volRoughHypersurface
$
on RHS~\eqref{E:FIRSTSTEPPROOFELLIPTICESTIMATECHI} 
is
bounded by the last two integrals on RHS~\eqref{E:FIRSTSTEPPROOFELLIPTICESTIMATECHI}.
Next, we use the torsion-free property of the connection $\angD$ 
and the $\gtorus$-Cauchy--Schwarz inequality to deduce the pointwise estimate
$
|\angLie_{\Yvf{A}} \upxi|_{\gtorus}
\leq 
|\angDarg{\Yvf{A}} \upxi|_{\gtorus}
+
2 |\upxi|_{\gtorus}|\angD \Yvf{A}|_{\gtorus}
$.
Also using 
\eqref{E:POINTWISESEMINORMOFYVECTORFIELDS},
\eqref{E:SMOOTHTORUSNORMCOMPARBLETOTANGENTIALCONTRACTIONS},
the estimate for $|\gfour(\Dfour_{\Yvf{A}} \Yvf{B},\Yvf{C})|$
given in the proof of \eqref{E:SMOOTHANGULARHESSIANOFFPOINTWISEBOUNDEDBYCOMMUTATORVECTORFIELDS},
and Cor.\,\ref{C:IMPROVEAUX},
we find that
$
|\angLie_{\Yvf{A}} \upxi|_{\gtorus}
\lesssim
|\angD \upxi|_{\gtorus}
+
\fundbootsmall
|\upxi|_{\gtorus}
$.
From this bound,
\eqref{E:FIRSTSTEPPROOFELLIPTICESTIMATECHI},
the estimates proved above,
and the pointwise bound 
$|\angD \mytr_{\gtorus} \upxi|_{\gtorus} 
\lesssim
\sum_{A=2,3} 
|\Yvf{A} \mytr_{\gtorus} \upxi|_{\gtorus}$,
which follows from \eqref{E:SMOOTHTORUSNORMCOMPARBLETOTANGENTIALCONTRACTIONS},
we conclude that
$
\int_{\hypthreearg{\timefunction}{[-\rightu,u]}{\muxmulevelsetvalue}}
	\upmu^2 |\angLie_{\Yvf{A}} \upxi|_{\gtorus}
\, \volRoughHypersurface 
\lesssim
\mbox{RHS~\eqref{E:ELLIPTICESTIMATECHI}}
$.
We have therefore proved \eqref{E:ELLIPTICESTIMATECHI},
which completes the proof of Prop.\,\ref{P:ELLIPTICESTIMATESONROUGHTORI}.

\hfill $\qed$


\section{Proof of the $L^2$ estimates for the wave variables and the acoustic geometry}  
\label{S:WAVEANDACOUSTICGEOMETRYAPRIORIESTIMATES}
In this section, we prove Props.\,\ref{P:APRIORIL2ESTIMATESWAVEVARIABLES} and \ref{P:APRIORIL2ESTIMATESACOUSTICGEOMETRY}, 
which provide the main a priori energy estimates for the wave variables and the acoustic geometry
along the rough foliations. We accomplish this via a bootstrap argument that relies on 
the energy estimates of Prop.\,\ref{P:MAINHYPERSURFACEENERGYESTIMATESFORTRANSPORTVARIABLES}
for the transport variables, which we already proved in
Sects.\,\ref{SS:PROOFOFBELOWTOPORDERTRANSPORTENERGYESTIMATES} and \ref{SS:PROOFOFMAINVORTVORTDIVGRADENTTOPORDERBLOWUP}.

We recall that the fundamental $L^2$-controlling
quantities, such as $\totalcontrolwave_{[1,N]}$ and $\totalcontrolwavepartial_{[1,N]}$, 
are defined in Sect.\,\ref{SS:FUNDAMENTALL2CONTROLLINGQUANTITIES}
(see in particular Def.\,\ref{D:SUMMEDL2CONTROLLINGQUANTITIES}).

\subsection{Statement of the integral inequalities used in proving a priori $L^2$ estimates for the wave variables} 
\label{SS:INTEGRALINEQUALITIIESFORWAVEENERGIES}
In this section, we state Prop.\,\ref{P:MAINWAVEENERGYINTEGRALINEQUALITIES},
which provides a coupled system of integral inequalities 
for the wave energies $\totalcontrolwave_N$ and the partial wave energies $\totalcontrolwavepartial_N$. 
As we will see in Sect.\,\ref{SSS:PROOFOFAPRIORIL2ESTIMATESWAVEVARIABLES}, 
these integral inequalities are the main ingredients in our proof of the $L^2$ 
a priori estimates of Prop.\,\ref{P:APRIORIL2ESTIMATESWAVEVARIABLES}.
Most of our effort in Sect.\,\ref{S:WAVEANDACOUSTICGEOMETRYAPRIORIESTIMATES} is dedicated towards
proving preliminary estimates that we will use in proving Prop.\,\ref{P:MAINWAVEENERGYINTEGRALINEQUALITIES}.

We now state the proposition. Its proof is located in Sect.\,\ref{SS:PROOFOF:MAINWAVEENERGYINTEGRALINEQUALITIES}.

\begin{proposition}[The main integral inequalities for the $\totalcontrolwave_N$] 
\label{P:MAINWAVEENERGYINTEGRALINEQUALITIES} 
Let $\varsigma \in (0,1]$. Let $\timefunction \in [\timefunction_0,\timefunctionboot]$ and let 
$\smallneighborhoodofcreasetwoarg{[\timefunction_0,\timefunctionboot]}{\muxmulevelsetvalue}$ denote the set 
from \eqref{E:SMALLNEIGHBORHOOD}. 
For $\wavearray = (\Psi_0,\Psi_1,\Psi_2,\Psi_3,\Psi_4) = ( \RRiemann,\LRiemann,v^2,v^3,\Ent)$, 
let $\vec{\mathfrak{G}} \eqdef (\mathfrak{G}_0,\cdots,\mathfrak{G}_4)$ be the vector array of the inhomogeneous terms 
in the covariant wave equations $\upmu \Box_{\gfour}\Psi_{\iota} = \mathfrak{G}_{\iota}$, 
i.e., $\upmu \mathfrak{G}$ is equal to $\upmu \times \mbox{RHS~\eqref{E:COVARIANTWAVEEQUATIONSWAVEVARIABLES}}$. 
Similarly, we define $\vec{\mathfrak{G}}_{(Partial)} \eqdef (\mathfrak{G}_1,\cdots,\mathfrak{G}_4)$ to be $\upmu$-weighted
inhomogeneous terms in the covariant wave equations satisfied by
$\wavearraypartial = (\LRiemann,v^2,v^3,\Ent)$.
There exist constants
$C > 0$ and $C_* > 0$ 
(see Remark~\ref{R:CSTAR})
that are independent of $\varsigma$ 
such that the following estimate for $\totalcontrolwave_N$ holds
for $(\timefunction,u) \in [\timefunction_0,\timefunctionboot) \times [- \rightu,\leftu]$: 

\medskip

\noindent \underline{\textbf{Top-order integral inequalities for $\wavearray$}.}
In the case $N = \Ntop$, 
we have the following estimates for the $L^2$-controlling quantity $\totalcontrolwavepartial_N$ 
defined in \eqref{E:WAVETOTALL2CONTROLLINGQUANTITY}:
\begin{align}
\begin{split}  \label{E:TOPORDERWAVEL2CONTROLLINGINTEGRALINEQUALITY}
	\totalcontrolwave_N(\timefunction,u) 
	& \lesssim 
	\boxed{\left\lbrace \frac{4 \times 1.01}{1.99} + 4.13 \right\rbrace} 
	\int_{\timefunction' = \timefunction_0}^{\timefunction} 
		\frac{1}{|\timefunction'|} \hypersurfacecontrolwave_N (\timefunction',u) 
	\, \mathrm{d} \timefunction' 
	 \\
	& \ \
	+ 
	\boxed{\frac{8\times(1.01)^2}{1.99}}  
	\int_{\timefunction' = \timefunction_0}^{\timefunction}
		\frac{1}{|\timefunction'|} \hypersurfacecontrolwave_{N}^{1/2}(\timefunction',u) 		
		\int_{\timefunction'' = \timefunction_0}^{\timefunction'} 
			\frac{1}{|\timefunction''|} \hypersurfacecontrolwave_N^{1/2}(\timefunction'',u) 		
		\, \mathrm{d}\timefunction'' 
	\, \mathrm{d} \timefunction' 
		 \\
	& \ \
	+ 
	\boxed{4.13}  
	\frac{1}{|\timefunction|^{1/2}} 
	\hypersurfacecontrolwave_N^{1/2}(\timefunction,u) 
	\int_{\timefunction' = \timefunction_0}^{\timefunction} 	
		\frac{1}{|\timefunction'|^{1/2}} \hypersurfacecontrolwave_N^{1/2}(\timefunction',u) 
	\, \mathrm{d} \timefunction' 
	 \\
	& \ \ 
	+ 
	C_* 
	\int_{\timefunction' = \timefunction_0}^{\timefunction} 
		\frac{1}{|\timefunction'|} 
		\hypersurfacecontrolwave_N^{1/2}(\timefunction',u) 
		\left(\hypersurfacecontrolwavepartial_N\right)^{1/2}(\timefunction',u) 
	\, \mathrm{d} \timefunction' 
	 \\
	& \ \
	+ 
	C_* 
	\int_{\timefunction' = \timefunction_0}^{\timefunction} 
		\frac{1}{|\timefunction'|} \hypersurfacecontrolwave_N^{1/2}(\timefunction',u)  
			\int_{\timefunction'' = \timefunction_0}^{\timefunction'} 
				\frac{1}{|\timefunction''|} \left(\hypersurfacecontrolwavepartial_N\right)^{1/2}(\timefunction'',u) 
				\, \mathrm{d} \timefunction'' 
			\, \mathrm{d}\timefunction' 
			 \\
	& \ \
	+ 
	C_*  
	\frac{1}{|\timefunction|^{1/2}} 
	\hypersurfacecontrolwave_N^{1/2}(\timefunction,u) 
	\int_{\timefunction' = \timefunction_0}^{\timefunction} 
		\frac{1}{|\timefunction'|^{1/2}} \left(\hypersurfacecontrolwavepartial_N \right)^{1/2}(\timefunction',u) 
	\, \mathrm{d} \timefunction' 
	 \\
	& \ \ 
	+ 
	\Errortoparg{N}(\timefunction,u),
\end{split}
\end{align} 
where for any $\varsigma \in (0,1]$,
$\Errortoparg{N}(\timefunction,u)$ satisfies the following estimate, 
where the implicit constants are independent of $\varsigma$:
\begin{align}
\begin{split} \label{E:ERRORTOPORDERWAVEESTIMATES}
|\Errortoparg{N}|(\timefunction,u) 
&
\lesssim 
\left(1 + \varsigma^{-1} \right)\initialsmall^2 \frac{1}{|\timefunction|^{3/2}} 
	\\
& \ \ 
		+
		\int_{\timefunction' = \timefunction_0}^{\timefunction} 
		\hypersurfacecontrolwave_N^{1/2}(\timefunction',u) 
		\left\lbrace 
			\hypersurfacecontrolVortVort_N^{1/2}
			+ 
			\hypersurfacecontrolDivGradEnt_N^{1/2}
		\right\rbrace
		(\timefunction',u) 
	\, \mathrm{d} \timefunction'
		\\
& \ \
		+
		\int_{\timefunction' = \timefunction_0}^{\timefunction} 
			\frac{1}{|\timefunction'|^{4/3}} 
			\left\lbrace
				\int_{\timefunction'' = \timefunction_0}^{\timefunction'}
					\left[
						\hypersurfacecontrolVortVort_N^{1/2} + \hypersurfacecontrolDivGradEnt_N^{1/2}
					\right]
					(\timefunction'',u) 
				\, \mathrm{d} \timefunction'' 
			\right\rbrace^2 
		\, \mathrm{d}\timefunction' 
		 \\
	& \ \
	+ 
	\int_{\timefunction' = \timefunction_0}^{\timefunction} 
		\frac{1}{|\timefunction'|^{4/3}} 
		\left\lbrace \int_{\timefunction'' = \timefunction_0}^{\timefunction'} 
			\frac{1}{|\timefunction''|^{1/2}} 
			\left[\hypersurfacecontrolVortVort_{\leq N-1}^{1/2} 
			+ 
			\hypersurfacecontrolDivGradEnt_{\leq N-1}^{1/2}
			\right]
			(\timefunction'',u) 
		\, \mathrm{d} \timefunction'' \right\rbrace^2 
	\, \mathrm{d}\timefunction' 
	 \\
& \ \
	+
	\int_{u' = -\rightu}^u
		\left\lbrace 
			\hypersurfacecontrolVortVort_{\leq N-1}
			+ 
			\hypersurfacecontrolDivGradEnt_{\leq N-1}
		\right\rbrace
		(\timefunction,u') 
	\, \mathrm{d} u'
		\\
& \ \
+ 
\int_{\timefunction' = \timefunction_0}^{\timefunction} 
	\frac{1}{|\timefunction'|^{3/2}} 
	\left\lbrace \int_{\timefunction'' = \timefunction_0}^{\timefunction'} 
		\left[
			\hypersurfacecontrolVort_N^{1/2} + \hypersurfacecontrolGradEnt_N^{1/2}
		\right]
			(\timefunction'',u) 
	\, \mathrm{d} \timefunction'' \right\rbrace^2 
\, \mathrm{d}\timefunction' 
	\\
	& \ \
	+ 
	\int_{\timefunction' = \timefunction_0}^{\timefunction} 
		\frac{1}{|\timefunction'|^{3/2}} 
		\left\lbrace 
			\int_{\timefunction'' = \timefunction_0}^{\timefunction'} 
				\frac{1}{|\timefunction''|^{1/2}} 
				\left[
					\hypersurfacecontrolVort_{\leq N-1}^{1/2} 
					+ 
					\hypersurfacecontrolGradEnt_{\leq N-1}^{1/2}
				\right]
				(\timefunction'',u) 
			\, \mathrm{d} \timefunction'' 
		\right\rbrace^2 
		\, \mathrm{d}\timefunction'
		\\
	& \ \
	+
	\int_{u' = -\rightu}^u
		\left\lbrace 
			\hypersurfacecontrolVort_{\leq N}
			+ 
			\hypersurfacecontrolGradEnt_{\leq N}
		\right\rbrace
		(\timefunction,u') 
	\, \mathrm{d} u'
		\\
	& \ \
	+ 
	\fundbootsmall 
	\hypersurfacecontrolwave_N(\timefunction,u) 
	+ 
	\varsigma \hypersurfacecontrolwave_N(\timefunction,u) 
	+ 
	\varsigma \totalcontrolwave_N(\timefunction,u) 
		\\ 
	& \ \
	+ 
	\fundbootsmall 
	\frac{1}{|\timefunction|^{1/2}} 
	\hypersurfacecontrolwave_N^{1/2}(\timefunction,u)  
	\int_{\timefunction' = \timefunction_0}^{\timefunction} 
		\frac{1}{|\timefunction'|^{1/2}} \hypersurfacecontrolwave_{[1,N]}^{1/2}(\timefunction',u)  
	\, \mathrm{d} \timefunction' 
			\\
	& \ \
	+
	\frac{1}{|\timefunction|^{1/2}} 
	\hypersurfacecontrolwave_N^{1/2}(\timefunction,u)  
	\int_{\timefunction' = \timefunction_0}^{\timefunction} 
			\hypersurfacecontrolwave_{[1,N]}^{1/2}(\timefunction',u) 
	\, \mathrm{d} \timefunction'
		\\
	& \ \
	+ 
	\hypersurfacecontrolwave_N^{1/2}(\timefunction,u)  
	\int_{\timefunction' = \timefunction_0}^{\timefunction} 
		\frac{1}{|\timefunction'|^{1/2}} \hypersurfacecontrolwave_{[1,N]}^{1/2}(\timefunction',u)  
	\, \mathrm{d} \timefunction' 
			\\
	& \ \
	+ 
	\fundbootsmall  
	\int_{\timefunction' = \timefunction_0}^{\timefunction} 
		\frac{1}{|\timefunction'|} \hypersurfacecontrolwave_N(\timefunction',u) 
	\, \mathrm{d} \timefunction' 
		\\
	& \ \		
		+ 
	\left(1 + \varsigma^{-1} \right)  
	\int_{\timefunction' = \timefunction_0}^{\timefunction} 
		\frac{1}{|\timefunction'|^{2/3}} \hypersurfacecontrolwave_{[1,N]}(\timefunction',u) 
	\, \mathrm{d} \timefunction'   
		\\
	& \ \
		+ 
	\left(1 + \varsigma^{-1} \right) 
	\int_{u' = -\rightu}^u 
		\hypersurfacecontrolwave_N(\timefunction,u') 
	\, \mathrm{d} u' 
		\\
	& \ \
	+ 
	\fundbootsmall 
	\int_{\timefunction' = \timefunction_0}^{\timefunction} 
		\hypersurfacecontrolwave_{ N}^{1/2}(\timefunction',u) \frac{1}{|\timefunction'|}  
		\int_{\timefunction'' = \timefunction_0}^{\timefunction'}  
			\frac{1}{|\timefunction''|} \hypersurfacecontrolwave_{[1,N]}^{1/2}(\timefunction'',u) 
		\, \mathrm{d} \timefunction'' 
	\, \mathrm{d}\timefunction' 
		\\
	& \ \
	+ 
	\int_{\timefunction' = \timefunction_0}^{\timefunction} 
		\hypersurfacecontrolwave_{ N}^{1/2}(\timefunction',u) \frac{1}{|\timefunction'|}  
		\int_{\timefunction'' = \timefunction_0}^{\timefunction'}  
			\frac{1}{|\timefunction''|^{1/2}} \hypersurfacecontrolwave_{[1,N]}^{1/2}(\timefunction'',u) 
		\, \mathrm{d} \timefunction'' 
	\, \mathrm{d}\timefunction'
	\\
	& \ \
	+ 
	\int_{\timefunction' = \timefunction_0}^{\timefunction} 
		\hypersurfacecontrolwave_{ N}^{1/2}(\timefunction',u) \frac{1}{|\timefunction'|}  
		\int_{\timefunction'' = \timefunction_0}^{\timefunction'}  
			\frac{1}{|\timefunction''|} 
			\int_{\timefunction''' = \timefunction_0}^{\timefunction''} 
				\frac{1}{|\timefunction'''|^{1/2}} 
				\hypersurfacecontrolwave_{[1,N]}^{1/2}(\timefunction''',u) 
			\, \mathrm{d} \timefunction''' 
		\, \mathrm{d} \timefunction'' 
	\, \mathrm{d} \timefunction' 
		\\
	& \ \
	+ 
	\int_{\timefunction' = \timefunction_0}^{\timefunction} 
		\frac{1}{|\timefunction'|^{5/2}} \hypersurfacecontrolwave_{[1, N-1]}(\timefunction',u) 
	\, \mathrm{d} \timefunction'.
\end{split}
\end{align} 

\medskip

\noindent \underline{\textbf{Top-order integral inequalities for $\wavearraypartial$}.}
In the case $N = \Ntop$, 
we have the following estimates for the $L^2$-controlling quantity $\totalcontrolwavepartial_N$
defined in \eqref{E:WAVEPARTIALL2CONTROLLINGQUANTITY}:
\begin{align}
\begin{split} \label{E:MAINWAVEPARTIALINTEGRALINEQUALITIES}
	\totalcontrolwavepartial_N(\timefunction,u) 
	& 
	\leq 
	\Errortop(\timefunction,u), 
\end{split}
\end{align} 
where $\Errortop(\timefunction,u)$ satisfies \eqref{E:ERRORTOPORDERWAVEESTIMATES}.

\medskip

\noindent \underline{\textbf{Below-top-order integral inequalities for $\wavearray$}.}
Finally, if $2 \leq N \leq \Ntop$, 
then we have the following estimates for the $L^2$-controlling quantity $\totalcontrolwave_{[1,N-1]}$
defined by \eqref{E:WAVETOTALL2CONTROLLINGQUANTITY} and Def.\,\ref{D:SUMMEDL2CONTROLLINGQUANTITIES}:
\begin{align}
\begin{split} \label{E:MAINWAVEBELOWTOPINTEGRALINEQUALITIES}
	\totalcontrolwave_{[1,N-1]} (\timefunction,u)  
	& \leq
		C
		\int_{\timefunction' = \timefunction_0}^{\timefunction} 
			\frac{1}{|\timefunction'|^{1/2}} \hypersurfacecontrolwave_{[1,N-1]}^{1/2}(\timefunction',u) 
			\int_{\timefunction'' = \timefunction_0}^{\timefunction'} 
				\frac{1}{|\timefunction''|^{1/2}} 
				\hypersurfacecontrolwave_N^{1/2}(\timefunction'',u) 
			\, \mathrm{d} \timefunction'' 
		\, \mathrm{d} \timefunction' 
			\\
& \ \
		+
		\Errorsubcriticalarg{N-1}(\timefunction,u), 
\end{split}
\end{align}
where for any integer $M \geq 1$,
$\Errorsubcriticalarg{M}(\timefunction,u)$ is defined to be any term that 
satisfies the following estimate, where the implicit constants depend only on the background solution (and hence are independent of $\varsigma$):
\begin{align}
\begin{split} \label{E:ERRORBELOWTOPORDERWAVEESTIMATES}
\left|
	\Errorsubcriticalarg{M}
\right|(\timefunction,u)  
& 
\lesssim 
\initialsmall^2 
+ 
\left(1 + \varsigma^{-1} \right) 
\int_{\timefunction' = \timefunction_0}^{\timefunction} 
	\frac{1}{|\timefunction'|^{1/2}} 
	\hypersurfacecontrolwave_{[1,M]}(\timefunction',u) 
\, \mathrm{d} \timefunction' 
	\\
& \ \ 
		+
		\int_{\timefunction' = \timefunction_0}^{\timefunction} 
			\hypersurfacecontrolwave_M^{1/2}(\timefunction',u) 
			\left\lbrace 
				\hypersurfacecontrolVortVort_M^{1/2}
				+ 
				\hypersurfacecontrolDivGradEnt_M^{1/2}
			\right\rbrace
			(\timefunction',u) 
		\, \mathrm{d} \timefunction'
		\\
& \ \
	+
	\int_{u' = -\rightu}^u
		\left\lbrace 
			\hypersurfacecontrolVortVort_{\leq M-1}
			+ 
			\hypersurfacecontrolDivGradEnt_{\leq M-1}
		\right\rbrace
		(\timefunction,u') 
	\, \mathrm{d} u'
		\\
& \ \
	+
	\int_{u' = -\rightu}^u
		\left\lbrace 
			\hypersurfacecontrolVort_{\leq M}
			+ 
			\hypersurfacecontrolGradEnt_{\leq M}
		\right\rbrace
		(\timefunction,u') 
	\, \mathrm{d} u'
		\\
& \ \ 
	+  
	\left(1 + \varsigma^{-1}\right) 
	\int_{u' = - \rightu}^u 
		\hypersurfacecontrolwave_{[1,M]}(\timefunction,u') 
	\, \mathrm{d} u' 
		\\
	& \ \
	+ 
	\varsigma 
	\spacetimeintegralcontrolwave_{[1,M]}(\timefunction,u).
\end{split}
\end{align}
\end{proposition}

\subsection{Estimates for the easiest error integrals}
\label{SS:ESTIMATESFOREASIESTERRORINTEGRALS}

\subsubsection{Estimates for the error integrals generated by the error terms $\HarmlessWave{N}$}
\label{SSS:ESTIMATESFOREERRORINTEGRALSGENERATEDBYHARMLESSWAVETERMS}
In the following lemma, we derive bounds for all the wave equation error integrals 
that involve the $\HarmlessWave{N}$ terms defined in \eqref{E:HARMLESSWAVE}. 

\begin{lemma}[Bounds for error integrals involving $\HarmlessWave{N}$ terms] 
\label{L:HARMLESSWAVETERMSERRORINTEGRALBOUNDS} 
Let $1 \leq N \leq \Ntop$, let $\Psi \in \wavearray = \{\RRiemann,\LRiemann,v^2,v^3,\Ent\}$,
and let $\varsigma \in (0,1]$. Recall that
terms of type $\HarmlessWave{N}$ are defined in Def.\,\ref{D:HARMLESSWAVE}.
The following estimates hold for 
$(\timefunction,u) \in [\timefunction_0,\timefunctionboot) \times [-\rightu,\leftu]$, 
where the implicit constants are independent of $\varsigma$:
\begin{align}
\begin{split} \label{E:ENERGYESTIMATEHARMLESSWAVEERRORTERMS}
	 \int_{\twoargMrough{[\timefunction_0,\timefunction],[-\rightu,u]}{\muxmulevelsetvalue}}   
		\frac{1}{\Lunit \timefunctionarg{\muxmulevelsetvalue}} 
		\left| 
			\begin{pmatrix} 
				(1 + 2 \upmu) \Lunit \tander^N \Psi 
				\\ 
				2 \muX \tander^N \Psi 
			\end{pmatrix} 
		\right|  
		\left| \HarmlessWave{N} \right| 
		\, \volMRoughCoordinates 
		& 
		\lesssim 
		\left( 1+ \varsigma^{-1} \right) 
		\int_{\timefunction' = \timefunction_0}^{\timefunction} 
			\frac{1}{|\timefunction'|^{1/2}} \hypersurfacecontrolwave_{[1,N]}(\timefunction',u) 
		\, \mathrm{d} \timefunction' 
		\\
	& \ \
	+ 
	\left( 1+ \varsigma^{-1}\right) 
	\int_{u' = -\rightu}^u 
		\hypersurfacecontrolwave_{[1,N]} (\timefunction,u') 
	\, \mathrm{d} u' 
		\\
	& \ \
	+
	\varsigma \spacetimeintegralcontrolwave_{[1,N]}(\timefunction,u) 
	+ 
	\initialsmall^2.
\end{split}
\end{align}
In particular,
RHS~\eqref{E:ENERGYESTIMATEHARMLESSWAVEERRORTERMS} is of type
$\Errorsubcriticalarg{N}(\timefunction,u)$,
where $\Errorsubcriticalarg{N}(\timefunction,u)$
satisfies \eqref{E:ERRORBELOWTOPORDERWAVEESTIMATES} 
(hence, in the case $N = \Ntop$, 
it is also is of type $\Errortoparg{N}(\timefunction,u)$, 
i.e., it satisfies the weaker estimate \eqref{E:ERRORTOPORDERWAVEESTIMATES}). 
\end{lemma}

\begin{proof}
We will give a detailed proof for terms of type
$\int_{\twoargMrough{[\timefunction_0,\timefunction],[-\rightu,u]}{\muxmulevelsetvalue}} 
	\left|\Lunit \tander^N \Psi \right| 
	\cdot 
	\left|\Yvf{A} \tander^{\leq N} \Psi \right|
\, \volMRoughCoordinates$, 
which are the most difficult terms generated by the LHS~\eqref{E:ENERGYESTIMATEHARMLESSWAVEERRORTERMS}.  
By using the pointwise bound 
$\left|\Yvf{A} \tander^N \Psi \right| 
\leq
|\Yvf{A}|_{\gtorus}
\left|\angrmD \tander^N \Psi \right|_{\gtorus}
\lesssim 
\left|\angrmD \tander^N \Psi \right|_{\gtorus}$ 
(see \eqref{E:POINTWISESEMINORMOFYVECTORFIELDS}),
the estimates
\eqref{E:ROUGHTIMEFUNCTIONLDERIVATIVEBOUNDS},
\eqref{E:COERCIVENESSOFHYPERSURFACECONTROLWAVE},
\eqref{E:L2ESTIMATESFORWAVEVARIABLESONROUGHHYPERSURFACELOSSOFONEDERIVATIVE}, 
and \eqref{E:COERCIVENESSOFSPACETIMEINTEGRAL},
and Young's inequality, 
and noting that \eqref{E:MUISLARGEINBORINGREGION} implies that
$
1
=
\mathbf{1}_{[-\interestingu,\interestingu]}(u') 
+
\mathbf{1}_{[-\interestingu,\interestingu]^c}(u')
\leq
\mathbf{1}_{[-\interestingu,\interestingu]}(u')
+
C \upmu
$,
we deduce the following estimate for any $\varsigma \in (0,1]$, 
where the implicit constants are independent of $\varsigma$:
	\begin{align}
	\begin{split} \label{E:REPTERMENERGYESTIMATEHARMLESSWAVEERRORTERMS}
	 &
		\int_{\twoargMrough{[\timefunction_0,\timefunction],[-\rightu,u]}{\muxmulevelsetvalue}} \frac{1}{\Lunit \timefunctionarg{\muxmulevelsetvalue}} 
				\left|\Lunit \tander^N \Psi \right| 
				\cdot 
				\left| \Yvf{A} \tander^{\leq N} \Psi \right| 
		\, \volMRoughCoordinates 
				\\
		& 
		\lesssim 
		\left( 1+ \varsigma^{-1}\right)  
		\int_{u' = -\rightu}^u 
			\int_{\nullhypthreearg{\muxmulevelsetvalue}{u'}{[\timefunction_0,\timefunction]}} 
				\frac{1}{\Lunit \timefunctionarg{\muxmulevelsetvalue}} 
				\left| \Lunit \tander^N \Psi \right|^2 
			\, \volPuRoughCoordinates 
		\, \mathrm{d} u'  
			\\
	& \ \
		+  
		\int_{u' = -\rightu}^u 
			\int_{\nullhypthreearg{\muxmulevelsetvalue}{u'}{[\timefunction_0,\timefunction]}} 
				\frac{1}{\Lunit \timefunctionarg{\muxmulevelsetvalue}} 
				\upmu 
				\left| \angrmD \tander^N \Psi \right|_{\gtorus}^2 
			\, \volPuRoughCoordinates 
		\, \mathrm{d} u' 
			\\
		& 
		\ \
		+ 
		\varsigma 
		\int_{\twoargMrough{[\timefunction_0,\timefunction],[-\rightu,u]}{\muxmulevelsetvalue}}
			\mathbf{1}_{[-\interestingu,\interestingu]}(u') \frac{1}{\Lunit \timefunctionarg{\muxmulevelsetvalue}}  
			\left| \angrmD \tander^N \Psi \right|_{\gtorus}^2 
		\, \volMRoughCoordinates
				\\
		& \ \
		+
		\int_{\timefunction' = \timefunction_0}^{\timefunction} 
			\int_{\hypthreearg{\timefunction'}{[-\rightu,u]}{\muxmulevelsetvalue}}
				\left| \tander^{[1,N]} \Psi \right|^2
			\, \volRoughHypersurface 
		\, \mathrm{d} \timefunction'
				\\
		& \lesssim 
			\left( 1+ \varsigma^{-1}\right) 
			\int_{u' = - \rightu}^u 
				\hypersurfacecontrolwave_{[1,N]}(\timefunction,u) 
			\, \mathrm{d} u' 
		+  
		\varsigma \spacetimeintegralcontrolwave_{[1,N]}(\timefunction,u)
		+ 
		\initialsmall^2
		+
		\int_{\timefunction' = \timefunction_0}^{\timefunction} 
			\hypersurfacecontrolwave_{[1,N]}(\timefunction',u) 
		\, \mathrm{d} \timefunction',
	\end{split}
	\end{align}
which is $\lesssim \mbox{RHS~\eqref{E:ENERGYESTIMATEHARMLESSWAVEERRORTERMS}}$ as desired.
The remaining terms on LHS~\eqref{E:ENERGYESTIMATEHARMLESSWAVEERRORTERMS} can be bounded 
by combining similar arguments with 
the estimates of Prop.\,\ref{P:SHARPCONTROLOFMUANDDERIVATIVES} (especially \eqref{E:MINVALUEOFMUONFOLIATION}),
the coerciveness estimates of Lemmas~\ref{L:COERCIVENESSOFL2CONTROLLINGQUANITIES} and \ref{L:COERCIVENESSOFSPACETIMEINTEGRAL},
the estimates \eqref{E:PRELIMINARYEIKONALWITHOUTL} and \eqref{E:NONDEGENERATEBUTDERIVATIVELOSINGL2ESTIMATEPSI},
the Cauchy--Schwarz inequality for integrals, and Young's inequalities;
we refer to the proof of \cite[Lemma 14.12]{jLjS2018} for further details.  
\end{proof}

\subsubsection{Estimates for the error integrals generated by the inhomogeneous terms in the covariant wave equations}
\label{SSS:ESTIMATESFORERRORINTEGRALSGENERATEDBYTHEINHOMOGENEOUSTERMS}

\begin{lemma}[Estimates for the error integrals generated by the inhomogeneous terms in the covariant wave equations] 
\label{L:ESTIMATESFORERRORINTEGRALSGENERATEDBYTHEINHOMOGENEOUSTERMS}
Let $\varsigma \in (0,1]$. For 
$\Psi \in \wavearray = (\Psi_0,\Psi_1,\Psi_2,\Psi_3,\Psi_4) = (\RRiemann,\LRiemann, v^2,v^3,\Ent)$, 
let $\vec{\mathfrak{G}} = (\mathfrak{G}_0,\cdots,\mathfrak{G}_4)$ 
be the array of the inhomogeneous terms in the covariant wave equations 
$\upmu \Box_{\gfour}\Psi_{\iota} = \mathfrak{G}_{\iota}$
(see \eqref{E:COVARIANTWAVEEQUATIONSWAVEVARIABLES}). 
Recall that $\multipliervectorfield$ is the multiplier vectorfield defined in \eqref{E:MULTIPLIERVECTORFIELD}.
Then the following estimates hold for 
$(\timefunction,u) \in [\timefunction_0,\timefunctionboot) \times [-\rightu,\leftu]$.
\medskip

\noindent \underline{\textbf{Top-order estimates}}.
If $N = \Ntop$, then we have: 
\begin{align}  \label{E:ESTIMATESFORTOPORDERERRORINTEGRALSINCOVARIANTWAVEEQ} 
	\int_{\twoargMrough{[\timefunction_0,\timefunction],[-\farrightu,u]}{\muxmulevelsetvalue}} 
		\frac{1}{\Lunit \timefunctionarg{\muxmulevelsetvalue}}
		\left|\multipliervectorfield \tander^N \wavearray \right| 
		\left|\tander^N \vec{\mathfrak{G}} \right| 
	\, \volMRoughCoordinates 
	& 
	\leq
	\Errortoparg{N}(\timefunction,u), 
\end{align}
where $\Errortoparg{N}(\timefunction,u)$ satisfies \eqref{E:ERRORTOPORDERWAVEESTIMATES}.

\medskip

\noindent \underline{\textbf{Below-top-order estimates}}.
If $2 \leq N \leq \Ntop$, then we have:
\begin{align} \label{E:ESTIMATESFORBELOWTOPORDERERRORINTEGRALSINCOVARIANTWAVEEQ} 
	\int_{\twoargMrough{[\timefunction_0,\timefunction],[-\farrightu,u]}{\muxmulevelsetvalue}} 
		\frac{1}{\Lunit \timefunctionarg{\muxmulevelsetvalue}}
		\left|\multipliervectorfield \tander^{N-1} \wavearray \right| 
		\left|\tander^{N-1} \vec{\mathfrak{G}} \right| 
	\, \volMRoughCoordinates 
	& 
	\leq
	\Errorsubcriticalarg{N-1}(\timefunction,u), 
\end{align}
where
$\Errorsubcriticalarg{N-1}(\timefunction,u)$ satisfies \eqref{E:ERRORBELOWTOPORDERWAVEESTIMATES} 
(with $N-1$ in the role of $M$ in \eqref{E:ERRORBELOWTOPORDERWAVEESTIMATES}).

\end{lemma}

\begin{proof}
	Throughout the proof, we will silently use 
	the estimate
	$\frac{1}{\Lunit \timefunctionarg{\muxmulevelsetvalue}} \approx 1$
	implied by \eqref{E:ROUGHTIMEFUNCTIONLDERIVATIVEBOUNDS}.
	
We first prove \eqref{E:ESTIMATESFORTOPORDERERRORINTEGRALSINCOVARIANTWAVEEQ}, 
i.e, we handle the case $N = \Ntop$.
We pointwise bound the term $|\tander^N \vec{\mathfrak{G}}|$
on LHS~\eqref{E:ESTIMATESFORTOPORDERERRORINTEGRALSINCOVARIANTWAVEEQ}.
Specifically, using the pointwise estimate
\eqref{E:POINTWISEESTIMATESFORALLINHOMOGENEOUSTERMS},
the bootstrap assumptions, 
\eqref{E:POINTWISESEMINORMOFYVECTORFIELDS},
Young's inequality,
and noting (as in the proof of \eqref{E:REPTERMENERGYESTIMATEHARMLESSWAVEERRORTERMS})
that \eqref{E:MUISLARGEINBORINGREGION} implies that
$
1
=
\mathbf{1}_{[-\interestingu,\interestingu]}(u') 
+
\mathbf{1}_{[-\interestingu,\interestingu]^c}(u')
\leq
\mathbf{1}_{[-\interestingu,\interestingu]}(u')
+
C \upmu
$,
we deduce the following pointwise estimate 
for the integrand on LHS~\eqref{E:ESTIMATESFORTOPORDERERRORINTEGRALSINCOVARIANTWAVEEQ},
valid for any $\varsigma \in (0,1]$, 
with implicit constants that are independent of $\varsigma$:
\begin{align}
\begin{split} \label{E:PROOFSTEPESTIMATESFORTOPORDERERRORINTEGRALSINCOVARIANTWAVEEQ}
	\frac{1}{\Lunit \timefunctionarg{\muxmulevelsetvalue}}
		\left|\multipliervectorfield \tander^N \wavearray \right| 
		\left|\tander^N \vec{\mathfrak{G}} \right| 
	& \lesssim 
	 \left|
		\muX \tander^N \wavearray
	\right|
	\cdot
	\left|
		\upmu^{1/2} \tander^N (\VortVort,\DivGradEnt)
	\right|
	+
	 \left|
		\upmu^{1/2} \Lunit \tander^N \wavearray
	\right|
	\cdot
	\left|
		\upmu^{1/2} \tander^N (\VortVort,\DivGradEnt)
	\right|
		\\
& \ \
	+
	\left|
		\tander^{\leq N-1} (\VortVort,\DivGradEnt)
	\right|^2
		\\
& \ \
+
\left(1 + \varsigma^{-1}\right)
\left|
	\muX \tander^{[1,N]} \wavearray
\right|^2
+
\left(1 + \varsigma^{-1}\right)
\left|
	\Lunit \tander^{[1,N]} \wavearray
\right|^2
	\\
& \ \
+
\varsigma
\mathbf{1}_{[-\interestingu,\interestingu]}
\left|
	\angrmD \tander^{[1,N]} \wavearray
\right|_{\gtorus}^2
+
\left|
	\upmu^{1/2} \angrmD \tander^{[1,N]} \wavearray
\right|_{\gtorus}^2
	\\
& \ \
+
\left|
	\tander^{\leq N} (\vortrenormalized,\GradEnt)
\right|^2
+ 
\fundbootsmall
\left|
	\tandersmall^{[1,N]} \badcontrolvars 
\right|^2.
\end{split}
\end{align}
We now integrate RHS~\eqref{E:PROOFSTEPESTIMATESFORTOPORDERERRORINTEGRALSINCOVARIANTWAVEEQ}
over
$
\twoargMrough{[\timefunction_0,\timefunction],[-\farrightu,u]}{\muxmulevelsetvalue}
$.
Using \eqref{E:COERCIVENESSOFHYPERSURFACECONTROLWAVE},
\eqref{E:COERCIVENESSHYPERSURFACECONTROLVORTVORT}--\eqref{E:COERCIVENESSHYPERSURFACEDIVGRADENT},
and the Cauchy--Schwarz inequality, we find that the integrals of the first
two terms on RHS~\eqref{E:PROOFSTEPESTIMATESFORTOPORDERERRORINTEGRALSINCOVARIANTWAVEEQ}
are 
$\lesssim
\int_{\timefunction' = \timefunction_0}^{\timefunction} 
		\hypersurfacecontrolwave_N^{1/2}(\timefunction',u) 
		\left\lbrace 
			\hypersurfacecontrolVortVort_N^{1/2}
			+ 
			\hypersurfacecontrolDivGradEnt_N^{1/2}
		\right\rbrace
		(\timefunction',u) 
	\, \mathrm{d} \timefunction'
$,
which in turn is manifestly bounded by RHS~\eqref{E:ERRORTOPORDERWAVEESTIMATES} as desired.
Moreover, using \eqref{E:COERCIVENESSHYPERSURFACECONTROLVORTVORT}--\eqref{E:COERCIVENESSHYPERSURFACEDIVGRADENT},
we see that 
the integral of 
$\left|
	\tander^{\leq N-1} (\VortVort,\DivGradEnt)
\right|^2
$
is
$\lesssim
\int_{u' = -\rightu}^u
		\left\lbrace 
			\hypersurfacecontrolVortVort_{\leq N-1}
			+ 
			\hypersurfacecontrolDivGradEnt_{\leq N-1}
		\right\rbrace
		(\timefunction,u') 
	\, \mathrm{d} u'
$, 
which is manifestly bounded by RHS~\eqref{E:ERRORTOPORDERWAVEESTIMATES} as desired.
Finally, using \eqref{E:MINVALUEOFMUONFOLIATION},
Lemma~\ref{L:COERCIVENESSOFL2CONTROLLINGQUANITIES},
Lemma~\ref{L:COERCIVENESSOFSPACETIMEINTEGRAL},
the estimate \eqref{E:PRELIMINARYEIKONALWITHOUTL},
Young's inequality,
and the fact that the $L^2$-controlling quantities
$\hypersurfacecontrolwave_M(\timefunction,u)$, $\hypersurfacecontrolVortVort_M(\timefunction,u)$, etc.\
are increasing in their arguments,
it is straightforward to check that the integrals of the 
remaining terms on RHS~\eqref{E:PROOFSTEPESTIMATESFORTOPORDERERRORINTEGRALSINCOVARIANTWAVEEQ}
are of type $\Errortoparg{N}(\timefunction,u)$, i.e., they satisfy \eqref{E:ERRORTOPORDERWAVEESTIMATES}.

To prove \eqref{E:ESTIMATESFORBELOWTOPORDERERRORINTEGRALSINCOVARIANTWAVEEQ},
we repeat the above arguments with $N-1$ in the role of $N$ and simply note
that the same arguments imply something stronger than what was claimed above: they imply that
the error integrals generated by
LHS~\eqref{E:ESTIMATESFORBELOWTOPORDERERRORINTEGRALSINCOVARIANTWAVEEQ}
are of type $\Errorsubcriticalarg{N-1}(\timefunction,u)$,
i.e., they are all bounded by
RHS~\eqref{E:ERRORBELOWTOPORDERWAVEESTIMATES} 
(with $N-1$ in the role of $M$ in \eqref{E:ERRORBELOWTOPORDERWAVEESTIMATES}).
\end{proof}

\subsubsection{Estimates for the error integrals generated by the multiplier vectorfield}
\label{SSS:ESTIMATESFORERRORTERMSGENERATEDBYMULTIPLIERVECTORFIELD}
In the next lemma, we bound the error integrals that are generated by the multiplier vectorfield. 

\begin{lemma}[Estimates for the error integrals generated by the multiplier vectorfield]
\label{L:ESTIMATESFORERRORTERMSGENERATEDBYMULTIPLIERVECTORFIELD}
Assume that $1 \leq N \leq \Ntop$,
$\Psi \in \wavearray = \{\RRiemann,\LRiemann, v^2,v^2,\Ent\}$,
and $\varsigma \in (0,1]$.
Recall that the multiplier vectorfield $\multipliervectorfield$ is defined in \eqref{E:MULTIPLIERVECTORFIELD}
and that ${^{(\multipliervectorfield)}\mathfrak{B}}[\tander^N \Psi]$ 
is the error term defined by 
\eqref{E:WAVEENERGYIDENTITYBULKERRORTERM} and \eqref{E:WAVEENERGYIDENTITYBULKERRORTERM1}--\eqref{E:WAVEENERGYIDENTITYBULKERRORTERM6}
and appearing on RHS~\eqref{E:FUNDAMENTALENERGYINTEGRALIDENTITCOVARIANTWAVES}.
Then the following estimates hold for 
$(\timefunction,u) \in [\timefunction_0,\timefunctionboot) \times [-\rightu,\leftu]$, 
where the implicit constants are independent of $\varsigma$:
\begin{align}
\begin{split} \label{E:WAVEENERGYESTIMATEMULTIPLIERERRORTERMS}
	\int_{\twoargMrough{[\timefunction_0,\timefunction],[-\rightu,u]}{\muxmulevelsetvalue}} 
		\left|
			\frac{1}{\Lunit \timefunctionarg{\muxmulevelsetvalue}} 
			{^{(\multipliervectorfield)}\mathfrak{B}}[\tander^N \Psi] 
		\right|
	\, \volMRoughCoordinates	
& \lesssim
\left( 1 + \varsigma^{-1}\right) 
\int_{\timefunction' = \timefunction_0}^{\timefunction} 
	\frac{1}{|\timefunction'|^{1/2}}
	\hypersurfacecontrolwave_N(\timefunction',u) 
\, \mathrm{d} \timefunction'  
	\\
& \ \ 
+ 
\left( 1 + \varsigma^{-1}\right) 
\int_{u' = -\rightu}^u 
	\hypersurfacecontrolwave_N(\timefunction,u') 
\, \mathrm{d} u' 
		\\
& \ \ 
+ 
\varsigma \spacetimeintegralcontrolwave_N(\timefunction,u). 
\end{split}
\end{align}
In particular,
RHS~\eqref{E:WAVEENERGYESTIMATEMULTIPLIERERRORTERMS} is of type
$\Errorsubcriticalarg{N}(\timefunction,u)$,
where $\Errorsubcriticalarg{N}(\timefunction,u)$
satisfies \eqref{E:ERRORBELOWTOPORDERWAVEESTIMATES} 
(hence, in the case $N = \Ntop$, 
it is also is of type $\Errortoparg{N}(\timefunction,u)$, 
i.e., it satisfies the weaker estimate \eqref{E:ERRORTOPORDERWAVEESTIMATES}). 
\end{lemma} 

\begin{proof}
Throughout this proof, we will silently use the estimate \eqref{E:CLOSEDVERSIONLUNITROUGHTTIMEFUNCTION},
which implies that the factor
$
\frac{1}{\Lunit \timefunctionarg{\muxmulevelsetvalue}} 
$
on LHS~\eqref{E:WAVEENERGYESTIMATEMULTIPLIERERRORTERMS} is $\approx 1$.

We first bound the error integral generated by the first term on RHS~\eqref{E:WAVEENERGYIDENTITYBULKERRORTERM}, 
i.e., in view of the factor 
$
\frac{1}{\Lunit \timefunctionarg{\muxmulevelsetvalue}}
$
on LHS~\eqref{E:WAVEENERGYESTIMATEMULTIPLIERERRORTERMS} and definition~\eqref{E:LROUGH},
the integral of 
$\frac{1}{2} \mathbf{1}_{[-\interestingu,\interestingu]^c} 
(\argLrough{\muxmulevelsetvalue} \upmu) 
\left|\angrmD \tander^N \Psi \right|_{\gtorus}^2$. 
Using \eqref{E:EASYREGIONLOWERBOUNDFORMU}, 
we find that the integrand satisfies the pointwise bound
$
	\frac{1}{\Lunit \timefunctionarg{\muxmulevelsetvalue}}
		\mathbf{1}_{[-\interestingu,\interestingu]^c} 
		|\Lunit \upmu| \left|\angrmD \tander^N \Psi \right|_{\gtorus}^2 
\lesssim	
\upmu 
 \left|\angrmD \tander^N \Psi \right|_{\gtorus}^2 
$.
Hence, by \eqref{E:COERCIVENESSOFHYPERSURFACECONTROLWAVE}, 
we deduce:
\begin{align} \label{E:PROOFSTEPWAVEENERGYESTIMATEMULTIPLIERERRORTERMS}
	\frac{1}{2} 
	\int_{\twoargMrough{[\timefunction_0,\timefunction],[-\rightu,u]}{\muxmulevelsetvalue}} 
		\left|
			\frac{1}{\Lunit \timefunctionarg{\muxmulevelsetvalue}}
			\mathbf{1}_{[-\interestingu,\interestingu]^c} 
			(\Lunit \upmu) \left|\angrmD \tander^N \Psi \right|_{\gtorus}^2 
		\right|
	\, \volMRoughCoordinates 
	& 
	\lesssim 
		\int_{u' = -\rightu}^u 
			\hypersurfacecontrolwave_N(\timefunction,u') 
	\, \mathrm{d} u',
\end{align}
which is $\lesssim \mbox{RHS~\eqref{E:WAVEENERGYESTIMATEMULTIPLIERERRORTERMS}}$ as desired.

Next, we bound the integral of the term ${^{(\multipliervectorfield)}\mathfrak{B}_{(3)}}[\tander^N \Psi]$ 
defined in \eqref{E:WAVEENERGYIDENTITYBULKERRORTERM3}.
First, using the crucial pointwise bound \eqref{E:WMANDRTRANSMUUBOUNDEDBYSQRTMU}, 
we bound the term $(\Rtransarg{\muxmulevelsetvalue} \upmu) \left|\angrmD \tander^N \Psi \right|_{\gtorus}^2$ 
from \eqref{E:WAVEENERGYIDENTITYBULKERRORTERM3} 
in magnitude by $\lesssim \sqrt{\upmu} \left|\angrmD \tander^N \Psi \right|_{\gtorus}^2$.
Hence, using the coerciveness estimate
\eqref{E:COERCIVENESSOFHYPERSURFACECONTROLWAVE}
and \eqref{E:MINVALUEOFMUONFOLIATION},
we can bound the integral of this term over the region
$\twoargMrough{[\timefunction_0,\timefunction],[-\rightu,u]}{\muxmulevelsetvalue}$
by
$
\lesssim 
	\int_{\timefunction' = \timefunction_0}^{\timefunction}  
		\frac{1}{|\timefunction'|^{1/2}} \hypersurfacecontrolwave_N(\timefunction',u) 
	\mathrm{d}\timefunction'
$
as desired.

We now handle the remaining terms in the definition \eqref{E:WAVEENERGYIDENTITYBULKERRORTERM3}
of ${^{(\multipliervectorfield)}\mathfrak{B}_{(3)}}$
as well as the
remaining bulk terms ${^{(\multipliervectorfield)}\mathfrak{B}_{(i)}[\tander^N \Psi]}$ with $i \in\{1,2,4,5,6\}$,
i.e., the terms defined in 
\eqref{E:WAVEENERGYIDENTITYBULKERRORTERM1},
\eqref{E:WAVEENERGYIDENTITYBULKERRORTERM2},
\eqref{E:WAVEENERGYIDENTITYBULKERRORTERM4},
\eqref{E:WAVEENERGYIDENTITYBULKERRORTERM5},
and
\eqref{E:WAVEENERGYIDENTITYBULKERRORTERM6}.
To this end, we first use Lemma~\ref{L:SCHEMATICSTRUCTUREOFVARIOUSTENSORSINTERMSOFCONTROLVARS}, 
\eqref{E:TANGENTIALDERIVATIVESOFROUGHTORIVECTORFIELDAPPLIEDTOMUPOINTWISE},
the bootstrap assumptions, 
and Young's inequality, 
and we split
$
1
=
\mathbf{1}_{[-\interestingu,\interestingu]}(u) 
+
\mathbf{1}_{[-\interestingu,\interestingu]^c}(u)
$,
to deduce the following pointwise estimates,
valid for any $\varsigma \in (0,1]$ with implicit constants that are independent of $\varsigma$:
\begin{align}
\begin{split}  \label{E:PROOFSTEP3WAVEENERGYESTIMATEMULTIPLIERERRORTERMS}
&
\left| 
			\upmu \Roughtoritangentvectorfieldarg{\muxmulevelsetvalue} \upmu 
			+ 
			2 \upmu  \Lunit \upmu 
			+ 
			\frac{1}{2}
			\upmu \mytr_{\gtorus}\upchi 
			+ 
			\upmu^2 \mytr_{\gtorus} \angktan 
			+ 
			\upmu \mytr_{\gtorus} \angktrans 
		\right|
		\left|\angrmD \tander^N \Psi \right|_{\gtorus}^2,
			\\
& \max_{i \in\{1,2,4,5,6\}}
\left|{^{(\multipliervectorfield)}\mathfrak{B}_{(i)}[\tander^N \Psi]} \right|
	\\
&
\lesssim 
\left( 1 + \varsigma^{-1}\right) (\Lunit \tander^N\Psi)^2 
+ 
\left( 1 + \varsigma^{-1}\right) 
(\muX \tander^N \Psi)^2 
	\\
& \ \
+
\upmu \left|\angrmD \tander^N \Psi \right|_{\gtorus}^2
+ 
\mathbf{1}_{[-\interestingu,\interestingu]^c} \left|\angrmD \tander^N \Psi \right|_{\gtorus}^2 
+ 
\varsigma \mathbf{1}_{[-\interestingu,\interestingu]} \left|\angrmD \tander^N \Psi \right|_{\gtorus}^2. 
\end{split}
\end{align}
Next, we use \eqref{E:MUISLARGEINBORINGREGION}
to deduce the following pointwise estimate for the next-to-last term on RHS~\eqref{E:PROOFSTEP3WAVEENERGYESTIMATEMULTIPLIERERRORTERMS}:
$\mathbf{1}_{[-\interestingu,\interestingu]^c} \left|\angrmD \tander^N \Psi \right|_{\gtorus}^2
\lesssim \upmu \left|\angrmD \tander^N \Psi \right|_{\gtorus}^2
$.
Combining this estimate with the coerciveness guaranteed by 
\eqref{E:COERCIVENESSOFHYPERSURFACECONTROLWAVE} and \eqref{E:COERCIVENESSOFSPACETIMEINTEGRAL}, 
we find that the integral of
RHS~\eqref{E:PROOFSTEP3WAVEENERGYESTIMATEMULTIPLIERERRORTERMS} 
over the spacetime region $\twoargMrough{[\timefunction_0,\timefunction],[-\rightu,u]}{\muxmulevelsetvalue}$ 
is $\lesssim \mbox{RHS~\eqref{E:WAVEENERGYESTIMATEMULTIPLIERERRORTERMS}}$ as desired.

\end{proof}

\subsection{Top-order $L^2$ estimates for $\upchi$}
\label{SS:TOPORDERL2ESTIMATSFORCHI}
In this section, we derive $L^2$ estimates for the top-order terms 
$\upmu \tander^{\Ntop} \mytr_{\gtorus} \upchi$ and $\upmu \angLie_\tander^{\Ntop} \upchi$ 
in terms of the $L^2$-controlling quantities. 
The proofs rely on the pointwise estimate \eqref{E:LESSPRECISEESTIMATETRACECHI} satisfied 
by the fully-modified quantity $\fullymodquant{\tander^N}$ 
as well as the elliptic estimate \eqref{E:ELLIPTICESTIMATECHI}. 

\subsubsection{Statement of the top-order $L^2$ estimates for $\upchi$}
\label{SSS:STATEMENTOFTOPORDERL2ESTIMATSFORCHI}
In the next proposition, we state the estimates. Its proof is located in Sect.\,\ref{SSS:PROOFOFTOPORDERL2ESTIMATSFORCHI}.

\begin{proposition}[Top-order $L^2$ estimates for $\upchi$] \label{P:TOPORDERL2ESTIMATEMUCHI}
Let $N = \Ntop$, and let $\mathfrak{P}^{(N)}$ and $\angLie_{\mathfrak{P}}^{(N)}$ denote the sets
of order $N$ $\nullhyparg{u}$-tangential operators from Sect.\,\ref{SS:STRINGSOFCOMMUTATIONVECTORFIELDS}.
Then the following estimates hold for 
$(\timefunction,u) \in [\timefunction_0,\timefunctionboot) \times [-\rightu,\leftu]$:
\begin{align}
\begin{split} \label{E:TOPORDERL2ESTIMATEMUCHI}
	& 
	\max_{\tander^N \in \mathfrak{P}^{(N)}}
	\left\| \upmu \tander^N \mytr_{\gtorus} \upchi \right\|_{L^2\left(\hypthreearg{\timefunction}{[-\rightu,u]}{\muxmulevelsetvalue}\right)}, 
		\,
	\max_{\angLie_{\tander}^N \in \angLie_{\mathfrak{P}}^{(N)}}
	\left\| \upmu \angLie_{\tander}^N \upchi \right\|_{L^2\left(\hypthreearg{\timefunction}{[-\rightu,u]}{\muxmulevelsetvalue}\right)}  
		\\
	&
	\lesssim 
	\initialsmall \ln \left( |\timefunction|^{-1} \right)
	+
	\hypersurfacecontrolwave_{[1,N]}^{1/2}(\timefunction,u) 
	+ 
	\int_{\timefunction' = \timefunction_0}^{\timefunction} 
		\frac{1}{|\timefunction'|} \hypersurfacecontrolwave_{[1,N]}^{1/2}(\timefunction',u) 
	\, \mathrm{d} \timefunction' 
		\\
	& \ \
	+ 
	\int_{\timefunction' = \timefunction_0}^{\timefunction} 
		\left\lbrace 
			\hypersurfacecontrolVortVort_N^{1/2} 
			+ 
			\hypersurfacecontrolDivGradEnt_N^{1/2} 
			+ 
			\hypersurfacecontrolVort_N^{1/2} 
			+ 
			\hypersurfacecontrolGradEnt_N^{1/2}
		\right\rbrace
		(\timefunction',u) 
	\, \mathrm{d} \timefunction'  
	\\
	& \ \
		+ 
		\int_{\timefunction' = \timefunction_0}^{\timefunction} 
			\frac{1}{|\timefunction'|^{1/2}} 
			\left\lbrace 
				\hypersurfacecontrolVortVort_{\leq N-1}^{1/2} 
					+ 
				\hypersurfacecontrolDivGradEnt_{\leq N-1}^{1/2}  
					+ 
				\hypersurfacecontrolVort_{\leq N-1}^{1/2} 
					+ 
				\hypersurfacecontrolGradEnt_{\leq N-1}^{1/2} 
			\right\rbrace(\timefunction',u) 
		\, \mathrm{d} \timefunction'. 
	\end{split}
	\end{align}
\end{proposition}

\subsubsection{Preliminary estimates}
\label{SSS:PRELIMINARYESTIAMTSFORTOPORDERL2ESTIMATSFORCHI}
In the following lemma, we derive preliminary $L^2$
estimates that we will use in the proof of Prop.\,\ref{P:TOPORDERL2ESTIMATEMUCHI}. 

\begin{lemma}[Preliminary top-order $L^2$ estimates for $\upchi$] 
Let $N = \Ntop$, and let
$\mathfrak{P}^{(N)}$, $\angLie_{\mathfrak{P}}^{(N)}$, and $\angLie_{\mathfrak{Y}}^{(N)}$
be the sets of order $N$ $\nullhyparg{u}$-tangential commutator operators from
Sect.\,\ref{SS:STRINGSOFCOMMUTATIONVECTORFIELDS}.
Then the following estimates hold for
$(\timefunction,u) \in [\timefunction_0,\timefunctionboot) \times [-\rightu,\leftu]$:
\begin{subequations} 
\begin{align}
\begin{split} \label{E:PRELIMINARYTOPORDERTRACECHI} 
\max_{\tander^N \in \mathfrak{P}^{(N)}}
\left\| \upmu \tander^N \mytr_{\gtorus} \upchi \right\|_{L^2\left(\hypthreearg{\timefunction}{[-\rightu,u]}{\muxmulevelsetvalue}\right)} 
&   
\lesssim 
\initialsmall \ln \left( |\timefunction|^{-1} \right)
+
\hypersurfacecontrolwave_{[1,N]}^{1/2}(\timefunction,u) 
	+ 
\int_{\timefunction' = \timefunction_0}^{\timefunction} 
	\frac{1}{|\timefunction'|} 
	\hypersurfacecontrolwave_{[1,N]}^{1/2}(\timefunction',u) 
\, \mathrm{d} \timefunction' 
	 \\
& \ \ 
+ 
\fundbootsmall 
\int_{\timefunction' = \timefunction_0}^{\timefunction} 
	\left\| 
		\upmu 
		\max_{\angLie_{\tanderY}^N \in \angLie_{\mathfrak{Y}}^{(N)}}
		\angLie_{\tanderY}^N \upchi 
	\right\|_{L^2\left(\hypthreearg{\timefunction'}{[-\rightu,u]}{\muxmulevelsetvalue}\right)} 
\, \mathrm{d} \timefunction'  
	\\
& \ \
+ 
\int_{\timefunction' = \timefunction_0}^{\timefunction} 
	\left\lbrace 
		\hypersurfacecontrolVortVort_N^{1/2} 
		+ 
		\hypersurfacecontrolDivGradEnt_N^{1/2} 
		+ 
		\hypersurfacecontrolVort_N^{1/2} 
		+ 
		\hypersurfacecontrolGradEnt_N^{1/2}\right\rbrace(\timefunction',u) 
	\, \mathrm{d} \timefunction' 
		\\
	& 
	\ \ 
	+ 
	\int_{\timefunction' = \timefunction_0}^{\timefunction} 
		\frac{1}{|\timefunction'|^{1/2}} 
		\left\lbrace 
			\hypersurfacecontrolVortVort_{\leq N-1}^{1/2} 
			+ 
			\hypersurfacecontrolDivGradEnt_{\leq N-1}^{1/2} 
			+ 
			\hypersurfacecontrolVort_{\leq N-1}^{1/2} 
			+ 
			\hypersurfacecontrolGradEnt_{\leq N-1}^{1/2} 
		\right\rbrace  (\timefunction',u) 
		\, \mathrm{d} \timefunction',
	\end{split}	
		\\
\max_{\angLie_{\tander}^N \in \angLie_{\mathfrak{P}}^{(N)}}
\left\| 
	\upmu 
	\angLie_{\tander}^N \upchi 
\right\|_{L^2\left(\hypthreearg{\timefunction}{[-\rightu,u]}{\muxmulevelsetvalue}\right)} 
& 
\lesssim  
\initialsmall
+ 
\hypersurfacecontrolwave_{[1,N]}^{1/2}(\timefunction,u) 
+ 
\max_{\tander^N \in \mathfrak{P}^{(N)}}
\left\|
	\upmu \tander^N \mytr_{\gtorus}\upchi 
\right\|_{L^2\left(\hypthreearg{\timefunction}{[-\rightu,u]}{\muxmulevelsetvalue}\right)}.
	\label{E:PRELIMINARYTOPORDERLIECHI}
\end{align}
\end{subequations}
\end{lemma}

\begin{proof}
\noindent \textbf{Proof of \eqref{E:PRELIMINARYTOPORDERTRACECHI}}:
It suffices to show that for any $\tander^N \in \mathfrak{P}^{(N)}$, we have:
$\left\| \upmu \tander^N \mytr_{\gtorus} \upchi \right\|_{L^2\left(\hypthreearg{\timefunction}{[-\rightu,u]}{\muxmulevelsetvalue}\right)}
\lesssim \mbox{RHS~\eqref{E:PRELIMINARYTOPORDERTRACECHI}}$.
We first consider the case 
$\tander^N = \tanderY^N \in \mathfrak{Y}^{(N)}$.
We take the 
$\| \cdot \|_{L^2\left(\hypthreearg{\timefunction}{[-\rightu,u]}{\muxmulevelsetvalue}\right)}$
norm of the imprecise estimate \eqref{E:LESSPRECISEESTIMATETRACECHI}.
In the rest of the proof, we sometimes silently use
\eqref{E:DISTORTINGWITHROUGHFLOWMAPDOESNOTCHANGEL2NORMSMUCH}
and Cor.\,\ref{C:IMPROVEAUX},
which together imply that the flow map factors $\FlowmapLrougharg{\muxmulevelsetvalue}$ in \eqref{E:LESSPRECISEESTIMATETRACECHI}
distort $\| \cdot \|_{L^2\left(\hypthreearg{\timefunction}{[-\rightu,u]}{\muxmulevelsetvalue}\right)}$ norms
only by overall factors of $1 + \mathcal{O}(\fundbootsmall)$.
For this reason, in this proof, we often suppress the factors of $\FlowmapLrougharg{\muxmulevelsetvalue}$
to simplify the notation.

To proceed, we use
\eqref{E:CLOSEDVERSIONLUNITROUGHTTIMEFUNCTION},
the estimate $|\upmu| \lesssim 1$ (which follows from the bootstrap assumptions),
\eqref{E:MINKOWSKIFORINTEGRALS},
and Lemma~\ref{L:COERCIVENESSOFL2CONTROLLINGQUANITIES}
to deduce that the norm
$\| \cdot \|_{L^2\left(\hypthreearg{\timefunction}{[-\rightu,u]}{\muxmulevelsetvalue}\right)}$ 
of the terms
$\upmu |\tander^{[1,N+1]}\wavearray|, 
	\, 
|\muX \tander^{[1,N]} \wavearray|$, 
and
$ \int_{\timefunction' = \timefunction_0}^{\timefunction}  
	\frac{1}{|\timefunction'|} 
	|\muX \tander^N \wavearray| 
\, \mathrm{d} \timefunction'$ 
on RHS~\eqref{E:LESSPRECISEESTIMATETRACECHI}
are $\lesssim$ 
the sum of the first three terms on RHS~\eqref{E:PRELIMINARYTOPORDERTRACECHI}. 
To bound the norm of the term $|\comdersmall^{N;1} \wavearray|$ on RHS~\eqref{E:LESSPRECISEESTIMATETRACECHI},
we simply use the already proven estimate \eqref{E:NONDEGENERATEBUTDERIVATIVELOSINGL2ESTIMATEPSI},
while to handle the terms
$
\left| \begin{pmatrix}
\tander^{[1,N]}\controlvars  
	\\
\tandersmall^{[1,N]} \badcontrolvars \end{pmatrix}
\right|
$
on RHS~\eqref{E:LESSPRECISEESTIMATETRACECHI}, 
we use \eqref{E:NONDEGENERATEBUTDERIVATIVELOSINGL2ESTIMATEPSI}
and the already proven estimate \eqref{E:PRELIMINARYEIKONALWITHOUTL}.
To bound the norms of the terms
$
\int_{\timefunction' = \timefunction_0}^{\timefunction} 
	\frac{1}{|\timefunction'|} 
	\left| \comdersmall^{[1,N];1} \wavearray \right|
\, \mathrm{d} \timefunction'
$
and
$
\int_{\timefunction' = \timefunction_0}^{\timefunction} 
	\frac{1}{|\timefunction'|} 
	\left| 
		\begin{pmatrix}
			\tander^{[1,N]}\controlvars 
				\\
			\tandersmall^{[1,N]} \badcontrolvars \end{pmatrix} \right| 
\, \mathrm{d} \timefunction'
$
on RHS~\eqref{E:LESSPRECISEESTIMATETRACECHI},
we use similar arguments, where we note that
\eqref{E:NONDEGENERATEBUTDERIVATIVELOSINGL2ESTIMATEPSI}
and \eqref{E:PRELIMINARYEIKONALWITHOUTL}
generate the error term
$
C
\int_{\timefunction_0}^{\timefunction}\frac{\initialsmall}{|\timefunction'|}  \, \mathrm{d}\timefunction'$,
which is $\lesssim$ the term
$\initialsmall \ln \left(|\timefunction|^{-1} \right)$ on RHS~\eqref{E:PRELIMINARYTOPORDERTRACECHI}. 
To bound the norm $\| \cdot \|_{L^2\left(\hypthreearg{\timefunction}{[-\rightu,u]}{\muxmulevelsetvalue}\right)}$ 
of the first term $\left| \fullymodquant{\tanderY^N}  \right|(\timefunction_0,u,x^2,x^3)$
on RHS~\eqref{E:LESSPRECISEESTIMATETRACECHI},
we first use \eqref{E:VOLFORMESTIMATEROUGHTORI} 
to deduce that:
$\left\| 
	\fullymodquant{\tanderY^N}(\timefunction_0,\cdot) 
\right\|_{L^2\left(\hypthreearg{\timefunction}{[-\rightu,u]}{\muxmulevelsetvalue}\right)} 
\lesssim 
\left\| 
	\fullymodquant{\tanderY^N} 
\right\|_{L^2\left(\hypthreearg{\timefunction_0}{[-\rightu,u]}{\muxmulevelsetvalue}\right)}
$. 
Next, using Def.\,\ref{D:FULLYANDPARTIALLYMODIFIEDQUANTITIES} and the data assumptions 
stated in Sect.\,\ref{SSS:QUANTITATIVEASSUMPTIONSONDATAAWAYFROMSYMMETRY}, 
we find that 
$
\left\| 
	\fullymodquant{\tanderY^N} 
\right\|_{L^2\left(\hypthreearg{\timefunction_0}{[-\rightu,u]}{\muxmulevelsetvalue}\right)}
\lesssim \initialsmall$,
which is $\lesssim \mbox{RHS~\eqref{E:PRELIMINARYTOPORDERTRACECHI}}$ as desired.
The $L^2$ norm of the remaining time integrals on RHS~\eqref{E:LESSPRECISEESTIMATETRACECHI}
can be bounded using similar arguments, 
the estimate $|\upmu| \lesssim 1$ (which follows from the bootstrap assumptions),
and Lemma~\ref{L:COERCIVENESSOFL2CONTROLLINGQUANITIES},
which we use to control the
$(\VortVort,\DivGradEnt,\Omega,\GradEnt)$-involving terms.
We have therefore proved \eqref{E:PRELIMINARYTOPORDERTRACECHI} in the case 
$\tander^N =  \tanderY^N$. 

We now prove \eqref{E:PRELIMINARYTOPORDERTRACECHI} in the case that the operator
$\tander^N$ on the LHS is not of type $\tanderY^N$, i.e., the case in which 
$\tander^N$ contains at least one factor of $\Lunit$.
In this case, we can use \eqref{E:COMMUTATOROFTANGENTIALANDTANGENTIALCOMMUTATORS}
and the bootstrap assumptions
to commute the factor of $\Lunit$ so that it acts last
and then use the pointwise estimate \eqref{E:LTANGENTIALTRCHIPOINTWISE} 
to deduce:
$
|\tander^N \mytr_{\gtorus} \upchi|
\lesssim
|\Lunit \tander^{N-1} \mytr_{\gtorus} \upchi|
+
|\tander^{\leq N-1} \mytr_{\gtorus} \upchi|
+
\left|\tander^{[1,N-1]}\controlvars \right|
\lesssim 
|\tander^{[1,N+1]}\wavearray| 
	+ 
|\tander^{[1,N]}\controlvars|
$.
Multiplying this estimate by $\upmu$
and using the arguments given above,
including \eqref{E:PRELIMINARYEIKONALWITHOUTL}
and the fact that $\hypersurfacecontrolwave_{[1,N]}(\timefunction,u)$ is increasing in its arguments,
we find that
$\| 
\upmu \tander^N \mytr_{\gtorus} \upchi 
\|_{L^2\left(\hypthreearg{\timefunction}{[-\rightu,u]}{\muxmulevelsetvalue}\right)}
\lesssim
\initialsmall 
+ 
\hypersurfacecontrolwave_{[1,N]}^{1/2}(\timefunction,u) 
$, 
which is $\lesssim \mbox{RHS~\eqref{E:PRELIMINARYTOPORDERTRACECHI}}$ as desired.
We have therefore proved \eqref{E:PRELIMINARYTOPORDERTRACECHI}.

\medskip

\noindent \textbf{Proof of \eqref{E:PRELIMINARYTOPORDERLIECHI}}:
It suffices to show that for any $\angLie_{\tander}^N \in \mathfrak{\mathfrak{P}}^{(N)}$, we have:
$
\left\| 
	\upmu 
	\angLie_{\tander}^N \upchi 
\right\|_{L^2\left(\hypthreearg{\timefunction}{[-\rightu,u]}{\muxmulevelsetvalue}\right)} 
\lesssim 
\mbox{RHS~\eqref{E:PRELIMINARYTOPORDERLIECHI}}
$.
To this end, we first apply the elliptic estimate \eqref{E:ELLIPTICESTIMATECHI} with 
$\upxi \eqdef \angLie_{\tander}^{N-1}\upchi$ to deduce:
\begin{align}  
\begin{split}	\label{E:PRELIMINARYTOPORDERLIECHIINTERMEDIATESTEP1}  
	\int_{\hypthreearg{\timefunction}{[-\rightu,u]}{\muxmulevelsetvalue}}  
		\upmu^2 |\angLie_{\tander}^N \upchi |_{\gtorus}^2 
	\, \volRoughHypersurface 
	& 
	\lesssim  
	\int_{\hypthreearg{\timefunction}{[-\rightu,u]}{\muxmulevelsetvalue}}
		\upmu^2 |\angLie_{\Lunit} \angLie_{\tander}^{N-1}\upchi|_{\gtorus}^2 
	\, \volRoughHypersurface 
	+  
	\int_{\hypthreearg{\timefunction}{[-\rightu,u]}{\muxmulevelsetvalue}}
		\upmu^2 |\angdiv \angLie_{\tander}^{N-1}\upchi|_{\gtorus}^2 
	\, \volRoughHypersurface 
		\\
	& \ \ 
	+  
	\int_{\hypthreearg{\timefunction}{[-\rightu,u]}{\muxmulevelsetvalue}} 
		\upmu^2 (\tander^{\le 1} \mytr_{\gtorus}\angLie_{\tander}^{N-1}\upchi)^2 
	\, \volRoughHypersurface  
	+ 
	\int_{\hypthreearg{\timefunction}{[-\rightu,u]}{\muxmulevelsetvalue}} 
		|\angLie_{\tander}^{N-1}\upchi|_{\gtorus}^2 
	\, \volRoughHypersurface.  
\end{split}
\end{align}
The same arguments given in the previous paragraph, starting
from the pointwise estimate \eqref{E:LTANGENTIALTRCHIPOINTWISE}, 
imply that
$$\| 
	\upmu \angLie_{\Lunit} \angLie_{\tander}^{N-1} \upchi 
\|_{L^2\left(\hypthreearg{\timefunction}{[-\rightu,u]}{\muxmulevelsetvalue}\right)}
\lesssim
\initialsmall 
+ 
\hypersurfacecontrolwave_{[1,N]}^{1/2}(\timefunction,u),
$$which yields the desired bound for the first term 
on RHS~\eqref{E:PRELIMINARYTOPORDERLIECHIINTERMEDIATESTEP1}.
Those arguments also imply, based on \eqref{E:PRELIMINARYEIKONALWITHOUTL},
that 
$\| 
	\angLie_{\tander}^{N-1} \upchi
\|_{L^2\left(\hypthreearg{\timefunction}{[-\rightu,u]}{\muxmulevelsetvalue}\right)}
\lesssim
\initialsmall 
+ 
\hypersurfacecontrolwave_{[1,N]}^{1/2}(\timefunction,u) 
$,
which yields the desired bound for the last term 
on RHS~\eqref{E:PRELIMINARYTOPORDERLIECHIINTERMEDIATESTEP1}.
To handle the third term on RHS~\eqref{E:PRELIMINARYTOPORDERLIECHIINTERMEDIATESTEP1},
we start with the following pointwise triangle inequality estimate:
\begin{align} \label{E:PRELIMINARYTOPORDERLIECHIINTERMEDIATESTEP2}
\upmu 
\left|
	\tander^{\leq 1} \mytr_{\gtorus} \angLie_{\tander}^{N-1} \upchi
\right|
& \lesssim
\upmu \left|\tander^{\leq 1} \tander^{N-1} \mytr_{\gtorus} \upchi\right|
+
\upmu
\left|
	\tander^{\leq 1} 
	\left(
		\mytr_{\gtorus} \angLie_{\tander}^{N-1} \upchi 
		- 
		\tander^{N-1} \mytr_{\gtorus} \upchi 
	\right)
\right|.
\end{align}
The norm $\| \cdot \|_{L^2\left(\hypthreearg{\timefunction}{[-\rightu,u]}{\muxmulevelsetvalue}\right)}$
of the first term $\upmu |\tander^{\le 1} \tander^{N-1} \mytr_{\gtorus} \angLie_{\tander}^{N-1}\upchi|$
on RHS~\eqref{E:PRELIMINARYTOPORDERLIECHIINTERMEDIATESTEP2}
is $\lesssim$ the last term on RHS~\eqref{E:PRELIMINARYTOPORDERLIECHI}.
To handle the second term on RHS~\eqref{E:PRELIMINARYTOPORDERLIECHIINTERMEDIATESTEP2},
we first note the following pointwise commutator estimate,
which follows easily from the Leibniz rule, 
the bootstrap assumptions,
and \eqref{E:TANDERGANDCHIESTIMATE}:
\begin{align} \label{E:LIEANDTRACECOMMUTATOR}
	\left|
		\tander^{\le 1} \left(\mytr_{\gtorus} \angLie_{\tander}^{N-1}\upchi 
		- 
		\tander^{N-1} \mytr_{\gtorus} \upchi\right)
	\right| 
	& \lesssim 	
	\left| 
		\tander^{[1,N+1]}\wavearray  
	\right| 
	+ 
	\left| \tander^{[1,N]} \controlvars \right|. 
\end{align}
Multiplying \eqref{E:LIEANDTRACECOMMUTATOR} by $\upmu$
and arguing as in the proof of \eqref{E:PRELIMINARYTOPORDERTRACECHI},
using in particular 
Lemma~\ref{L:COERCIVENESSOFL2CONTROLLINGQUANITIES} and \eqref{E:PRELIMINARYEIKONALWITHOUTL},
we bound the norm $\| \cdot \|_{L^2\left(\hypthreearg{\timefunction}{[-\rightu,u]}{\muxmulevelsetvalue}\right)}$
of the resulting RHS by
$
\lesssim
\initialsmall 
+ 
\hypersurfacecontrolwave_{[1,N]}^{1/2}(\timefunction,u) 
$,
which is $\lesssim \mbox{RHS~\eqref{E:PRELIMINARYTOPORDERLIECHI}}$ as desired.
It remains for us to bound the $\angdiv \angLie_{\tander}^{N-1} \upchi$-involving
integral on RHS~\eqref{E:PRELIMINARYTOPORDERLIECHIINTERMEDIATESTEP1}.
We start with the following pointwise triangle inequality estimate:
\begin{align} \label{E:PRELIMINARYTOPORDERLIECHIINTERMEDIATESTEP3}
\upmu
\left|
	\angdiv \angLie_{\tander}^{N-1} \upchi
\right|
& \lesssim
\upmu
\left|
	\angrmD \tander^{N-1} \mytr_{\gtorus} \upchi
\right|
+
\upmu
\left| 
	\angdiv \angLie_{\tander}^{N-1} \upchi 
	- 
	\angrmD \tander^{N-1} \mytr_{\gtorus} \upchi 
\right|.
\end{align}
The norm $\| \cdot \|_{L^2\left(\hypthreearg{\timefunction}{[-\rightu,u]}{\muxmulevelsetvalue}\right)}$
of the first term 
$
\upmu
\left|
	\angrmD \tander^{N-1} \mytr_{\gtorus} \upchi
\right|$
on RHS~\eqref{E:PRELIMINARYTOPORDERLIECHIINTERMEDIATESTEP3}
is $\lesssim$ the last term on RHS~\eqref{E:PRELIMINARYTOPORDERLIECHI}.
To handle the second term on RHS~\eqref{E:PRELIMINARYTOPORDERLIECHIINTERMEDIATESTEP3},
we simply multiply the pointwise estimate \eqref{E:CODAZZICOMMUTATORESTIMATES} by $\upmu$,
take the norm $\| \cdot \|_{L^2\left(\hypthreearg{\timefunction}{[-\rightu,u]}{\muxmulevelsetvalue}\right)}$
of the resulting inequality,
and use the same arguments we used to control
$\upmu \times \mbox{RHS~\eqref{E:LIEANDTRACECOMMUTATOR}}$.
We have therefore proved \eqref{E:PRELIMINARYTOPORDERLIECHI}, which finishes the proof of the lemma.

\end{proof}

\subsubsection{Proof of Prop.\,\ref{P:TOPORDERL2ESTIMATEMUCHI}}
\label{SSS:PROOFOFTOPORDERL2ESTIMATSFORCHI}
The estimate \eqref{E:TOPORDERL2ESTIMATEMUCHI}
for 
$
\max_{\tander^N \in \mathfrak{P}^{(N)}}
\left\|\upmu \tander^N \mytr_{\gtorus} \upchi \right\|_{L^2\left(\hypthreearg{\timefunction}{[-\rightu,u]}{\muxmulevelsetvalue}\right)}$ follows from inserting the estimate \eqref{E:PRELIMINARYTOPORDERLIECHI} for 
$
\max_{\angLie_{\tander}^N \in \angLie_{\mathfrak{P}}^{(N)}}
\left\| \upmu \angLie_{\tander}^N \upchi \right\|_{L^2\left(\hypthreearg{\timefunction}{[-\rightu,u]}{\muxmulevelsetvalue}\right)}$ 
into RHS~\eqref{E:PRELIMINARYTOPORDERTRACECHI} and applying Gr\"{o}nwall's inequality. 
We then insert the already proved estimate \eqref{E:TOPORDERL2ESTIMATEMUCHI} for 
$
\max_{\tander^N \in \mathfrak{P}^{(N)}}
\left\|\upmu \tander^N \mytr_{\gtorus} \upchi \right\|_{L^2\left(\hypthreearg{\timefunction}{[-\rightu,u]}{\muxmulevelsetvalue}\right)}$ 
into RHS~\eqref{E:PRELIMINARYTOPORDERLIECHI}, 
thereby obtaining the desired estimate for
$
\max_{\angLie_{\tander}^N \in \angLie_{\mathfrak{P}}^{(N)}}
\left\| \upmu \angLie_{\tander}^N \upchi \right\|_{L^2\left(\hypthreearg{\timefunction}{[-\rightu,u]}{\muxmulevelsetvalue}\right)}$
and completing the proof of the proposition.

\hfill $\qed$

\subsection{Estimates for the easy top-order eikonal function-involving error integrals}
\label{SS:ESTIMATESFOREASYTOPORDEREIKONALFUNCTIONERRORINTEGRALS}
In Prop.\,\ref{P:TOPORDERL2ESTIMATEMUCHI}, we derived preliminary top-order $L^2$ estimates for $\upchi$. 
With the help of these estimates, we are now ready to control
the wave equation error integrals that depend on these terms. 
More precisely, in the next lemma, we use these preliminary estimates to control ``easy'' error integrals,
which are generated by the first product $\angrmD^{\sharp} \Psi_{\iota} \cdot \upmu \angrmD \tanderY^{N-1} \mytr_{\gtorus}\upchi $
on RHS~\eqref{E:TOPCOMMUTEDWAVELFIRSTTHENALLYS}
and the second product 
$(\Speed^{-2} X^A) \angrmD^{\sharp} \Psi_{\iota} \cdot \upmu \angrmD \tanderY^{N-1}\mytr_{\gtorus}\upchi$
on RHS~\eqref{E:TOPCOMMUTEDWAVEALLYS}.
The corresponding error integrals are easy in the sense that the integrands contain a helpful factor of $\upmu$.
In Sect.\,\ref{SS:ESTIMATESFORMOSTDIFFICULTEIKONALFUNCTIONERRORINTEGRALS}, we
will control the analogous -- but much more difficult -- error integral generated by the first product
$
(\muX \Psi_{\iota}) \tanderY^{N-1} \Yvf{A} \mytr_{\gtorus} \upchi 
$
on RHS~\eqref{E:TOPCOMMUTEDWAVEALLYS},
which lacks the factor of $\upmu$.

\begin{lemma}[Estimates for the easy top-order eikonal function-involving error integrals]
\label{L:ESTIMATESFOREASYTOPORDEREIKONALFUNCTIONERRORINTEGRALS}
Let $N = \Ntop$ and $\Psi \in \wavearray = \{\RRiemann,\LRiemann,v^2,v^3,\Ent\}$.
Then the following estimates hold for 
$(\timefunction,u) \in [\timefunction_0,\timefunctionboot) \times [-\rightu,\leftu]$:
\begin{align}
\begin{split}  \label{E:ESTIMATESFOREASYTOPORDEREIKONALFUNCTIONERRORINTEGRALS}
	& 
	 \int_{\twoargMrough{[\timefunction_0,\timefunction],[-\rightu,u]}{\muxmulevelsetvalue}}   
		\frac{1}{\Lunit \timefunctionarg{\muxmulevelsetvalue}} 
		\left| 
			\begin{pmatrix} 
				(1 + 2 \upmu) \Lunit \tander^N \Psi
				\\ 
				2 \muX \tander^N \Psi 
			\end{pmatrix} 
		\right|  
		\left| 
			\begin{pmatrix} 
				(\angrmD^{\sharp} \Psi) \cdot \upmu \angrmD \tanderY^{N-1} \mytr_{\gtorus} \upchi
				\\ 
				(\Speed^{-2} X^A) \angrmD^{\sharp} \Psi_{\iota} \cdot \upmu \angrmD \tanderY^{N-1}\mytr_{\gtorus}\upchi 
			\end{pmatrix} 
		\right|  
	\, \volMRoughCoordinates 
		\\
	& =
		\Errortoparg{N}(\timefunction,u),
\end{split}
\end{align}
where $\Errortoparg{N}(\timefunction,u)$ satisfies the estimate \eqref{E:ERRORTOPORDERWAVEESTIMATES}. 
\end{lemma}

\begin{proof}
	The bootstrap assumptions and \eqref{E:MINVALUEOFMUONFOLIATION}
	imply that the integrand on LHS~\eqref{E:ESTIMATESFOREASYTOPORDEREIKONALFUNCTIONERRORINTEGRALS} is
	pointwise bounded by
	$
	\lesssim
	\left(
		\frac{1}{|\timefunction'|^{1/2}}
		\frac{1}{\sqrt{\Lunit \timefunctionarg{\muxmulevelsetvalue}}} 
		\upmu^{1/2}
		|\Lunit \tander^N \Psi|
		+
		|\muX \tander^N \Psi|
	\right)
	\cdot
	\upmu |\tanderY^N \mytr_{\gtorus} \upchi|
	$.
	Hence, integrating over $\twoargMrough{[\timefunction_0,\timefunction],[-\rightu,u]}{\muxmulevelsetvalue}$
	and using \eqref{E:COERCIVENESSOFHYPERSURFACECONTROLWAVE}
	and Young's inequality,
	we bound LHS~\eqref{E:ESTIMATESFOREASYTOPORDEREIKONALFUNCTIONERRORINTEGRALS} by:
	\begin{align} 
	\begin{split} \label{E:INTERMEDIATEPROOFSTEPESTIMATESFOREASYTOPORDEREIKONALFUNCTIONERRORINTEGRALS}
	& 
	\lesssim 
	\int_{\timefunction' = \timefunction_0}^{\timefunction} 
		\frac{1}{|\timefunction'|^{1/2}}
		\hypersurfacecontrolwave_N^{1/2}(\timefunction',u) 
		\left\| 
			\upmu \tanderY^N \mytr_{\gtorus} \upchi
		\right\|_{L^2\left(\hypthreearg{\timefunction'}{[-\rightu,u]}{\muxmulevelsetvalue}\right)}
	\, \mathrm{d} \timefunction'
		\\
	&
	\lesssim 
	\int_{\timefunction' = \timefunction_0}^{\timefunction} 
		\frac{1}{|\timefunction'|^{1/2}}
		\hypersurfacecontrolwave_N(\timefunction',u) 
	\, \mathrm{d} \timefunction'	
	+
	\int_{\timefunction' = \timefunction_0}^{\timefunction} 
		\frac{1}{|\timefunction'|^{1/2}}
		\left\| 
			\upmu \tanderY^N \mytr_{\gtorus} \upchi
		\right\|_{L^2\left(\hypthreearg{\timefunction'}{[-\rightu,u]}{\muxmulevelsetvalue}\right)}^2
	\, \mathrm{d} \timefunction'.
	\end{split}
	\end{align}
	Substituting the estimate \eqref{E:TOPORDERL2ESTIMATEMUCHI} into the last time integral
	on RHS~\eqref{E:INTERMEDIATEPROOFSTEPESTIMATESFOREASYTOPORDEREIKONALFUNCTIONERRORINTEGRALS},
	we find that LHS~\eqref{E:ESTIMATESFOREASYTOPORDEREIKONALFUNCTIONERRORINTEGRALS}
	is bounded by:
	\begin{align} \label{E:PROOFSTEPESTIMATESFOREASYTOPORDEREIKONALFUNCTIONERRORINTEGRALS}
	\begin{split}
		& 
		\lesssim 
		\initialsmall^2
		+
		\int_{\timefunction' = \timefunction_0}^{\timefunction} 
			\frac{1}{|\timefunction'|^{1/2}}
			\hypersurfacecontrolwave_{[1,N]}(\timefunction',u) 
		\, \mathrm{d} \timefunction'
		\\
	& \ \
			+
			\int_{\timefunction' = \timefunction_0}^{\timefunction} 
				\left\lbrace 
					\int_{\timefunction'' = \timefunction_0}^{\timefunction'} 
						\frac{1}{|\timefunction''|} 
						\hypersurfacecontrolwave_{[1,N]}^{1/2}(\timefunction'',u) 
					\, \mathrm{d} \timefunction'' 
				\right\rbrace^2 
			\, \mathrm{d}\timefunction'
				\\
	& \  \
			+
			\int_{\timefunction' = \timefunction_0}^{\timefunction} 
			\frac{1}{|\timefunction'|^{1/2}} 
			\left\lbrace 
				\int_{\timefunction'' = \timefunction_0}^{\timefunction'} 
					\left(
						\hypersurfacecontrolVortVort_N^{1/2} 
							+ 
						\hypersurfacecontrolDivGradEnt_N^{1/2} 
						+ 
						\hypersurfacecontrolVort_N^{1/2} 
						+ 
						\hypersurfacecontrolGradEnt_N^{1/2}
					\right)
					(\timefunction'',u) 
				\, \mathrm{d} \timefunction'' \right\rbrace^2 
		\, \mathrm{d}\timefunction'
			\\
		& \ \
			+
			\int_{\timefunction' = \timefunction_0}^{\timefunction} 
			\frac{1}{|\timefunction'|^{1/2}} 
			\left\lbrace 
				\int_{\timefunction'' = \timefunction_0}^{\timefunction'} 
					\frac{1}{|\timefunction''|^{1/2}} 
					\left(
						\hypersurfacecontrolVortVort_{\leq N-1}^{1/2} 
						+ 
						\hypersurfacecontrolDivGradEnt_{\leq N-1}^{1/2}  
						+ 
						\hypersurfacecontrolVort_{\leq N-1}^{1/2} 
						+ 
						\hypersurfacecontrolGradEnt_{\leq N-1}^{1/2}
					\right)
					(\timefunction'',u) 
				\, \mathrm{d} \timefunction'' \right\rbrace^2 
		\, \mathrm{d}\timefunction'.
	\end{split}
	\end{align}
	Also using that the $L^2$-controlling quantities
	$\hypersurfacecontrolwave_M(\timefunction,u)$, $\hypersurfacecontrolVortVort_M(\timefunction,u)$, etc.\
	are increasing in their arguments, we conclude that
	$\mbox{RHS~\eqref{E:PROOFSTEPESTIMATESFOREASYTOPORDEREIKONALFUNCTIONERRORINTEGRALS}} 
	\lesssim \mbox{RHS~\eqref{E:ERRORTOPORDERWAVEESTIMATES}}$ as desired.
\end{proof}

\subsection{Estimates for the most difficult top-order eikonal function-involving error integrals}
\label{SS:ESTIMATESFORMOSTDIFFICULTEIKONALFUNCTIONERRORINTEGRALS}
In this section, we control the most difficult eikonal function-involving 
terms appearing in the commuted wave equations of Prop.\,\ref{P:MOSTDIFFICULTWAVETERMS}. 
More precisely, the most difficult product is the first one
$(\muX \Psi_{\iota}) \tanderY^{N-1} \Yvf{A} \mytr_{\gtorus} \upchi$
on RHS~\eqref{E:TOPCOMMUTEDWAVEALLYS}.
When we derive energy estimates using the fundamental energy--null-flux identity \eqref{E:FUNDAMENTALENERGYINTEGRALIDENTITCOVARIANTWAVES}, these difficult terms are multiplied by $\multipliervectorfield \tanderY^N \Psi_{\iota}$,
where the multiplier vectorfield
$\multipliervectorfield$ is defined in \eqref{E:MULTIPLIERVECTORFIELD}.
This leads to the following difficult error integrals:
\begin{subequations}
\begin{align}
\int_{\twoargMrough{[\timefunction_0,\timefunction],[-\farrightu,\leftu]}{\muxmulevelsetvalue}} 
	\frac{1}{\Lunit \timefunctionarg{\muxmulevelsetvalue}} 
	\left\lbrace 
		2 \muX \tanderY^N \Psi_{\iota}
	\right\rbrace 
	\left\lbrace 
		(\muX \Psi_{\iota}) \tanderY^N \mytr_{\gtorus} \upchi 
	\right\rbrace 
	\, \volMRoughCoordinates, 
		\label{E:OVERVIEWMUXPSIDIFFICULTERRORINTEGAL} 
			\\
	\int_{\twoargMrough{[\timefunction_0,\timefunction],[-\farrightu,\leftu]}{\muxmulevelsetvalue}} 
		\frac{1}{\Lunit \timefunctionarg{\muxmulevelsetvalue}} 
		\left\lbrace 
			(1 + 2 \upmu) \Lunit \tanderY^N \Psi_{\iota}
		\right\rbrace 
		\left\lbrace 
			(\muX \Psi_{\iota}) \tanderY^N \mytr_{\gtorus}\upchi 
		\right\rbrace 
		\, \volMRoughCoordinates.
		\label{E:OVERVIEWLPSIDIFFICULTERRORINTEGAL}
\end{align}
\end{subequations}
We bound the integral \eqref{E:OVERVIEWMUXPSIDIFFICULTERRORINTEGAL} 
in Sect.\,\ref{SSS:EIKONALTOPORDERENERGYESTIMATESWITHOUTIBP} 
and the integral \eqref{E:OVERVIEWLPSIDIFFICULTERRORINTEGAL},
which we control via a further integration by parts with respect to $\Lunit$,
in Sect.\,{\ref{SSS:EIKONALTOPORDERENERGYESTIMATEWITHIBP}. 
It turns out that these two integrals are the main ones driving
top-order wave energy blowup-rate, i.e., RHS~\eqref{E:MAINWAVEENERGYESTIMATESBLOWUP} with $K=0$.
We provide the main estimates for these two integrals in Lemmas~\ref{L:BOUNDSFORMOSTDIFFICULTWAVEERRORNTEGRALS}
and \ref{L:ESTIMATESFORDIFFICULTSPACETIMEERRORINTEGRALSINVOLVINGIBPWRTL}.

\subsubsection{Estimates that do not involve integration by parts}  
\label{SSS:EIKONALTOPORDERENERGYESTIMATESWITHOUTIBP}
We begin our analysis of the integral \eqref{E:OVERVIEWMUXPSIDIFFICULTERRORINTEGAL} 
with the following lemma, which provides $L^2$ estimates for the difficult product 
$(\muX \Psi_{\iota}) \tanderY^N \mytr_{\gtorus} \upchi$
on RHS~\eqref{E:TOPCOMMUTEDWAVEALLYS}.
More precisely, in the lemma, we handle the most difficult case,
which is $\Psi_{\iota} = \RRiemann$. The products
corresponding to remaining wave variables $\lbrace \LRiemann,v^2,v^3,\Ent \rbrace$
are much easier to handle in the energy estimates because in these cases,
we gain a smallness factor of $\varepsilon$ from
the factor $\muX \Psi_{\iota}$; see \eqref{E:LINFINITYIMPROVEMENTAUXTRANSVERSALPDERIVATIVESPARTIALWAVEARRAYSMALL}.

\begin{lemma}[$L^2$ estimates for the most difficult product] 
\label{L:L2ESTIMATESFORMOSTDIFFICULTPRODUCT}
Let $N = \Ntop$.
The following estimates hold for 
$(\timefunction,u) \in [\timefunction_0,\timefunctionboot) \times [-\rightu,\leftu]$:
\begin{align}
 \begin{split} \label{E:L2WAVEESTIMATEMOSTDIFFICULTPRODUCT}
	& \left\| 
		\frac{1}{\Lunit \timefunctionarg{\muxmulevelsetvalue}} 
		(\muX \RRiemann) \tanderY^N \mytr_{\gtorus}\upchi 
	\right\|_{L^2\left(\hypthreearg{\timefunction}{[-\rightu,u]}{\muxmulevelsetvalue}\right)}   
		\\ 
	& \leq \boxed{\frac{2 \times 1.01}{\sqrt{1.99}}} 
	\frac{1}{|\timefunction|} \hypersurfacecontrolwave_N^{1/2}(\timefunction,u)  
		\\
	& \ \
	+ 
	\boxed{\frac{4 \times(1.01)^2}{\sqrt{1.99}}}  
	\frac{1}{|\timefunction|} 
	\int_{\timefunction' = \timefunction_0}^{\timefunction} 
		\frac{1}{|\timefunction'|} 
		\hypersurfacecontrolwave_N^{1/2}(\timefunction',u) 
	\, \mathrm{d}\timefunction'   
	\\
	& \ \ 
	+ 
	\frac{C_*}{|\timefunction|} \left(\hypersurfacecontrolwavepartial_N\right)^{1/2}(\timefunction,u) 
	+  
	\frac{C_*}{|\timefunction|} 
	\int_{\timefunction' = \timefunction_0}^{\timefunction}\frac{1}{|\timefunction'|} 
		\left(\hypersurfacecontrolwavepartial_N \right)^{1/2}(\timefunction',u) 
	\, \mathrm{d} \timefunction'  
	\\
	& \ \ 
	+ 
	\frac{C \fundbootsmall}{|\timefunction|} 
	\int_{\timefunction' = \timefunction_0}^{\timefunction} 
		\hypersurfacecontrolwave_{[1,N]}^{1/2}(\timefunction',u) 
	\, \mathrm{d} \timefunction'   
	\\
	& \ \ 
	+  
	\frac{C \fundbootsmall}{|\timefunction|} 
	\int_{\timefunction' = \timefunction_0}^{\timefunction}
		\int_{\timefunction'' = \timefunction_0}^{\timefunction'} 
			\left\lbrace 
				\frac{1}{|\timefunction''|} \hypersurfacecontrolwave_{[1,N]}^{1/2} 
					+ 
				\hypersurfacecontrolVortVort_N^{1/2} 
				+ 
				\hypersurfacecontrolDivGradEnt_N^{1/2} 
				+ 
				\hypersurfacecontrolVort_N^{1/2} 
				+ 
				\hypersurfacecontrolGradEnt_N^{1/2}
			\right\rbrace 
			\circ (\timefunction'',u)
		\, \mathrm{d} \timefunction'' 
	\mathrm{d} \timefunction'   
	  \\
	& \ \ 
	+  
	\frac{C \fundbootsmall}{|\timefunction|} 
	\int_{\timefunction' = \timefunction_0}^{\timefunction}
		\int_{\timefunction'' = \timefunction_0}^{\timefunction'} 
			\frac{1}{|\timefunction''|^{1/2}}
			\left\lbrace 
				\hypersurfacecontrolVortVort_{\leq N-1}^{1/2} 
				+ 
				\hypersurfacecontrolDivGradEnt_{\leq N-1}^{1/2}
				+
				\hypersurfacecontrolVort_{\leq N-1}^{1/2} 
					+ 
				\hypersurfacecontrolGradEnt_{\leq N-1}^{1/2}
			 \right\rbrace 
			\circ (\timefunction'',u)
		\, \mathrm{d} \timefunction'' 
	\mathrm{d} \timefunction'   
	  \\
	& \ \ 
	+ 
	\frac{C}{|\timefunction|} 
	\int_{\timefunction' = \timefunction_0}^{\timefunction} 
		\left\lbrace 
			\hypersurfacecontrolVortVort_N^{1/2} 
				+ 
			\hypersurfacecontrolDivGradEnt_N^{1/2} 
				+ 
			\hypersurfacecontrolVort_N^{1/2} 
				+ 
			\hypersurfacecontrolGradEnt_N^{1/2}
		\right\rbrace \circ (\timefunction',u) 
	\, \mathrm{d} \timefunction'   
	\\
	& \ \ 
	+ 
	\frac{C}{|\timefunction|} 
	\int_{\timefunction' = \timefunction_0}^{\timefunction} 
		\frac{1}{|\timefunction'|^{1/2}}
		\left\lbrace 
			\hypersurfacecontrolVortVort_{\leq N-1}^{1/2}
				+ 
			\hypersurfacecontrolDivGradEnt_{\leq N-1}^{1/2}
				+
			\hypersurfacecontrolVort_{\leq N-1}^{1/2} 
				+ 
			\hypersurfacecontrolGradEnt_{\leq N-1}^{1/2}
		\right\rbrace \circ (\timefunction',u) 
	\, \mathrm{d} \timefunction'   
	\\
	& \ \
	+ 
	\frac{C \fundbootsmall}{|\timefunction|} 
	\hypersurfacecontrolwave_{ N}^{1/2}(\timefunction,u) 
	+ 
	\frac{C}{|\timefunction|^{1/2}} 
	\hypersurfacecontrolwave_{[1,N]}^{1/2}(\timefunction,u)  
	+ 
	\frac{C}{|\timefunction|^{3/2}} 
	\hypersurfacecontrolwave_{[1, N-1]}^{1/2}(\timefunction,u)  
	 \\
	& \ \ 
	+ 
	\frac{C}{|\timefunction|}  
	\int_{\timefunction' = \timefunction_0}^{\timefunction} 
		\frac{1}{|\timefunction'|^{1/2}} 
		\hypersurfacecontrolwave_{[1,N]}^{1/2}(\timefunction',u) 
	\, \mathrm{d}\timefunction'   
	\\
	& \ \
	+  
	\frac{C}{|\timefunction|} 
	\int_{\timefunction' = \timefunction_0}^{\timefunction} 
		\frac{1}{|\timefunction'|} 
		\int_{\timefunction'' = \timefunction_0}^{\timefunction'} 
			\frac{1}{|\timefunction''|^{1/2}} \hypersurfacecontrolwave_{[1,N]}^{1/2} (\timefunction'',u) 
		\, \mathrm{d} \timefunction'' \mathrm{d} \timefunction'    
	\\
	& \ \
	+ 
	\frac{C \fundbootsmall }{|\timefunction|} 
	\int_{\timefunction' = \timefunction_0}^{\timefunction} 
		\frac{1}{|\timefunction'|} 
		\hypersurfacecontrolwave_{[1,N]}^{1/2}(\timefunction',u) 
	\, \mathrm{d}\timefunction'
	+ 
	\frac{C \initialsmall}{|\timefunction|^{3/2}}. 
\end{split}
\end{align}

\end{lemma}

\begin{proof}	
We consider the pointwise estimate \eqref{E:MOSTDELICATEPOINTWISEESTIMATEFORRPLUS}. 
In this proof, we sometimes silently use
\eqref{E:DISTORTINGWITHROUGHFLOWMAPDOESNOTCHANGEL2NORMSMUCH}
and Cor.\,\ref{C:IMPROVEAUX},
which together imply that the flow map factors $\FlowmapLrougharg{\muxmulevelsetvalue}$ in \eqref{E:MOSTDELICATEPOINTWISEESTIMATEFORRPLUS}
distort $\| \cdot \|_{L^2\left(\hypthreearg{\timefunction}{[-\rightu,u]}{\muxmulevelsetvalue}\right)}$ norms
only by overall factors of $1 + \mathcal{O}(\fundbootsmall)$;
the $\mathcal{O}(\fundbootsmall)$ factors lead to small error terms on RHS~\eqref{E:L2WAVEESTIMATEMOSTDIFFICULTPRODUCT}.

We now use \eqref{E:WIDETILDELMUISALMOSTMINUSONEINSMALLNEIGHBORHOOD} to bound 
$\left| 
\frac{\argLrough{\muxmulevelsetvalue} \upmu}{\upmu} 
\mathbf{1}_{\{\hypthreearg{\timefunction}{[-\rightu,u]}{\muxmulevelsetvalue} \cap \smallneighborhoodofcreasetwoarg{[\timefunction_0,\timefunctionboot]}{\muxmulevelsetvalue}\}} \circ \FlowmapLrougharg{\muxmulevelsetvalue}(\timefunction,u,x^2,x^3) 
 \right| \leq \frac{1.01}{|\timefunction|}$ 
everywhere it appears on RHS~\eqref{E:MOSTDELICATEPOINTWISEESTIMATEFORRPLUS}. 
In particular, we bound
the first and third terms on RHS~\eqref{E:MOSTDELICATEPOINTWISEESTIMATEFORRPLUS}
(which are multiplied by boxed constants)
by 
$\frac{2(1.01)}{|\timefunction|} 
\left|\muX \tanderY^N \RRiemann \right| 
\circ 
\FlowmapLrougharg{\muxmulevelsetvalue}(\timefunction,u,x^2,x^3)$
and 
$\frac{4(1.01)^2}{|\timefunction|} 
\int_{\timefunction' = \timefunction_0}^{\timefunction}
	\frac{1}{|\timefunction'|} 
\left|\muX \tanderY^N \RRiemann \right|
\circ \FlowmapLrougharg{\muxmulevelsetvalue}(\timefunction',u,x^2,x^3) 
\, \mathrm{d}\timefunction'$ respectively. 
We now take the norm
$\| \cdot \|_{L^2\left(\hypthreearg{\timefunction}{[-\rightu,u]}{\muxmulevelsetvalue}\right)}$
of the resulting pointwise inequality and use \eqref{E:L2ONROUGHCONSTANTTIMEHYPERSURFACEOFROUGHTIMEINTEGRALWITHFLOWMAPFACTORSBOUND}
and Cor.\,\ref{C:IMPROVEAUX}.
Also using the sharpened coerciveness estimate \eqref{E:L2WAVECONTROLWITHBETTERCONSTANTIFMUSMALL},
we see that these two terms lead, respectively, 
to the presence of the two boxed-constant-involving products
$\boxed{\frac{2\times 1.01}{\sqrt{1.99}}}\cdots, \boxed{\frac{4\times(1.01)^2}{\sqrt{1.99}}}\cdots$ 
on RHS~\eqref{E:L2WAVEESTIMATEMOSTDIFFICULTPRODUCT}
plus some error terms with $C \fundbootsmall$ factors. 
Similarly, since the $C_*$-multiplied terms  
on RHS~\eqref{E:MOSTDELICATEPOINTWISEESTIMATEFORRPLUS} 
involve 
$\muX \tanderY^N \wavearraypartial$,
we can use 
Lemma~\ref{L:COERCIVENESSOFL2CONTROLLINGQUANITIES} (specifically \eqref{E:COERCIVENESSOFHYPERSURFACECONTROLWAVEPARTIAL})
to bound their $\| \cdot \|_{L^2\left(\hypthreearg{\timefunction}{[-\rightu,u]}{\muxmulevelsetvalue}\right)}$ norms
by the $C_*$-multiplied terms on the third line of RHS~\eqref{E:L2WAVEESTIMATEMOSTDIFFICULTPRODUCT}.
	
Next, we use 
\eqref{E:L2ONROUGHCONSTANTTIMEHYPERSURFACEOFROUGHTIMEINTEGRALWITHFLOWMAPFACTORSBOUND},
\eqref{E:TOPORDERL2ESTIMATEMUCHI}, 
and Lemma~\ref{L:COERCIVENESSOFL2CONTROLLINGQUANITIES} to bound the 
$\| \cdot \|_{L^2\left(\hypthreearg{\timefunction}{[-\rightu,u]}{\muxmulevelsetvalue}\right)}$ norm of the term
$\frac{C \fundbootsmall}{|\timefunction|} 
\int_{\timefunction' = \timefunction_0}^{\timefunction} 
	\upmu 
	|\angLie_\tanderY^N \upchi| 
	\circ \FlowmapLrougharg{\muxmulevelsetvalue}(\timefunction',u,x^2,x^3) 
\, \mathrm{d} \timefunction'$ on RHS~\eqref{E:MOSTDELICATEPOINTWISEESTIMATEFORRPLUS}. 
We find that the this term is bounded by the sum of the
double time integrals on the fifth and sixth lines of RHS~\eqref{E:L2WAVEESTIMATEMOSTDIFFICULTPRODUCT},
plus as a few other terms on RHS~\eqref{E:L2WAVEESTIMATEMOSTDIFFICULTPRODUCT}.
	
Next, using 
\eqref{E:L2ONROUGHCONSTANTTIMEHYPERSURFACEOFROUGHTIMEINTEGRALWITHFLOWMAPFACTORSBOUND},
\eqref{E:LINFINITYIMPROVEMENTAUXMUANDTRANSVERSALDERIVATIVES},
\eqref{E:MINVALUEOFMUONFOLIATION},
and Lemma~\ref{L:COERCIVENESSOFL2CONTROLLINGQUANITIES},
we see that the $\| \cdot \|_{L^2\left(\hypthreearg{\timefunction}{[-\rightu,u]}{\muxmulevelsetvalue}\right)}$ norms of the 
time-integrals of
$\upmu |\tanderY^N(\VortVort,\DivGradEnt,\Omega,\GradEnt)| \circ \FlowmapLrougharg{\muxmulevelsetvalue}$ 
and
$|\tanderY^{\leq N-1}(\VortVort,\DivGradEnt,\Omega,\GradEnt)| \circ \FlowmapLrougharg{\muxmulevelsetvalue}$ 
on RHS~\eqref{E:MOSTDELICATEPOINTWISEESTIMATEFORRPLUS} 
are bounded by the seventh and eighth lines of RHS~\eqref{E:L2WAVEESTIMATEMOSTDIFFICULTPRODUCT}. 
	
It remains for us to bound the 
$\| \cdot \|_{L^2\left(\hypthreearg{\timefunction}{[-\rightu,u]}{\muxmulevelsetvalue}\right)}$ norm of the terms 
$\Error \circ \FlowmapLrougharg{\muxmulevelsetvalue}(\timefunction,u,x^2,x^3)$ on RHS~\eqref{E:MOSTDELICATEPOINTWISEESTIMATEFORRPLUS},
which satisfy the pointwise bound \eqref{E:ERRORTERMSINMOSTDELICATEPOINTWISEESTIMATEFORRPLUS}.
To handle the first term 
$
\frac{1}{|\timefunction|} \left|\fullymodquant{\tanderY^N}  \right|(\timefunction_0,u,x^2,x^3)
$
on RHS~\eqref{E:ERRORTERMSINMOSTDELICATEPOINTWISEESTIMATEFORRPLUS},
we note that in our proof of \eqref{E:PRELIMINARYTOPORDERTRACECHI},
we showed that
$\left\| \fullymodquant{\tanderY^N} \right\|_{L^2\left(\hypthreearg{\timefunction}{[-\rightu,u]}{\muxmulevelsetvalue}\right)} 
\lesssim 
\left\| \fullymodquant{\tanderY^N} \right\|_{L^2\left(\hypthreearg{\timefunction_0}{[-\rightu,u]}{\muxmulevelsetvalue}\right)}
\lesssim \initialsmall$.
Hence, the first term on RHS~\eqref{E:ERRORTERMSINMOSTDELICATEPOINTWISEESTIMATEFORRPLUS}
is 
$
\lesssim 
\frac{\initialsmall}{|\timefunction|}
\lesssim
\mbox{RHS~\eqref{E:L2WAVEESTIMATEMOSTDIFFICULTPRODUCT}}
$
as desired.
With the help of \eqref{E:MINVALUEOFMUONFOLIATION},
Lemma~\ref{L:COERCIVENESSOFL2CONTROLLINGQUANITIES},
and the estimate \eqref{E:PRELIMINARYEIKONALWITHOUTL},
we can bound
the $\| \cdot \|_{L^2\left(\hypthreearg{\timefunction}{[-\rightu,u]}{\muxmulevelsetvalue}\right)}$ norm of the remaining terms on
RHS~\eqref{E:MOSTDELICATEPOINTWISEESTIMATEFORRPLUS} by
$\lesssim \mbox{RHS~\eqref{E:L2WAVEESTIMATEMOSTDIFFICULTPRODUCT}}$
by using a subset of the ideas we used above; 
we refer to the proofs of \cite[Lemma 14.8]{jSgHjLwW2016} and \cite[Lemma 14.14]{jLjS2018} for more details.
\end{proof}

With the help of Lemma~\ref{L:L2ESTIMATESFORMOSTDIFFICULTPRODUCT},
we now establish the following lemma,
which is the main result of Sect.\,\ref{SSS:EIKONALTOPORDERENERGYESTIMATESWITHOUTIBP}.

\begin{lemma}[Bounds for the most difficult error integrals in the wave equation energy estimates] 
\label{L:BOUNDSFORMOSTDIFFICULTWAVEERRORNTEGRALS}
Let $N = \Ntop$.
The following estimates hold for 
$(\timefunction,u) \in [\timefunction_0,\timefunctionboot) \times [-\rightu,\leftu]$:
\begin{align} 
\begin{split} \label{E:SPACETIMEBOUNDSMOSTDIFFICULTWAVEPRODUCT}
	& 2 
	\left| 
		\int_{\twoargMrough{[\timefunction_0,\timefunction],[-\rightu,u]}{\muxmulevelsetvalue}} 
			\frac{1}{\Lunit \timefunctionarg{\muxmulevelsetvalue}} 
			(\muX \tanderY^N \RRiemann) 
			(\muX \RRiemann) 
			\tanderY^N \mytr_{\gtorus} \upchi 
		\, \volMRoughCoordinates  
	\right|  
	\\
	& \leq 
	\boxed{\frac{4\times 1.01}{1.99}} 
	\int_{\timefunction' = \timefunction_0}^{\timefunction} 
		\frac{1}{|\timefunction'|} 
		\hypersurfacecontrolwave_N(\timefunction',u) 
	\, \mathrm{d} \timefunction'  
	\\
	&  \ \
	+ 
	\boxed{\frac{8\times (1.01)^2}{1.99}}  
	\int_{\timefunction' = \timefunction_0}^{\timefunction} 
		\frac{1}{|\timefunction'|} \hypersurfacecontrolwave_{ N}^{1/2}(\timefunction',u) 
		\int_{\timefunction'' = \timefunction_0}^{\timefunction'} 
			\frac{1}{|\timefunction''|} 
		\hypersurfacecontrolwave_N^{1/2}(\timefunction'',u) 
		\, \mathrm{d}\timefunction''  
	\mathrm{d} \timefunction'  
	\\
	& \ \ 
		+ 
		C_* 
		\int_{\timefunction' = \timefunction_0}^{\timefunction} 
			\frac{1}{|\timefunction'|} 
			\hypersurfacecontrolwave_N^{1/2}(\timefunction',u) 
			\left(\hypersurfacecontrolwavepartial_N \right)^{1/2}(\timefunction',u) 
		\, \mathrm{d} \timefunction'  
		\\
	& \ \
		+ 
		C_* 
		\int_{\timefunction' = \timefunction_0}^{\timefunction} 
			\frac{1}{|\timefunction'|} 
			\hypersurfacecontrolwave_N^{1/2}(\timefunction',u)  
			\int_{\timefunction'' = \timefunction_0}^{\timefunction'} 
				\frac{1}{|\timefunction''|} 
				\left(\hypersurfacecontrolwavepartial_N\right)^{1/2}(\timefunction'',u) 
			\, \mathrm{d} \timefunction'' 
		\, \mathrm{d}\timefunction'  
			\\
	& \ \
		+ \Errortop(\timefunction,u), 
\end{split}
\end{align}
where $\Errortop$ satisfies the estimate \eqref{E:ERRORTOPORDERWAVEESTIMATES}. 

Moreover, for $\Psi \in \wavearraypartial = \{ \LRiemann,v^2,v^3,\Ent\}$, 
we have the following less degenerate estimates:
\begin{align}
\begin{split}  \label{E:SPACETIMEBOUNDSMOSTDIFFICULTWAVEPRODUCTPARTIAL}
	2 \left| 
				\int_{\twoargMrough{[\timefunction_0,\timefunction],[-\rightu,u]}{\muxmulevelsetvalue}} 
					\frac{1}{\Lunit \timefunctionarg{\muxmulevelsetvalue}} 
					(\muX \tanderY^N \Psi) (\muX \Psi) 
					\tanderY^N \mytr_{\gtorus} \upchi 
				\, \volMRoughCoordinates  
			\right| 
	& 
	\lesssim
	\Errortop(\timefunction,u),
\end{split}
\end{align}
where $\Errortop$ satisfies the estimate \eqref{E:ERRORTOPORDERWAVEESTIMATES}. 
\end{lemma}

\begin{proof}
We first prove \eqref{E:SPACETIMEBOUNDSMOSTDIFFICULTWAVEPRODUCT}. 
By H\"{o}lder's inequality, we have: 
\begin{align} 
\begin{split} \label{E:FIRSTSTEPSPACETIMEBOUNDSMOSTDIFFICULTWAVEPRODUCT}
	&
	2 \left| 
			\int_{\twoargMrough{[\timefunction_0,\timefunction],[-\rightu,u]}{\muxmulevelsetvalue}}
				\frac{1}{\Lunit \timefunctionarg{\muxmulevelsetvalue}} (\muX \tanderY^N \RRiemann) (\muX \RRiemann) \tanderY^N \mytr_{\gtorus} \upchi 
			\, \volMRoughCoordinates  
		\right| 
		\\
	& \leq 
	2 
	\int_{\timefunction' = \timefunction_0}^{\timefunction} 
		\left\| 
			\muX \tanderY^N \RRiemann 
		\right\|_{L^2\left( \hypthreearg{\timefunction'}{[-\rightu,u]}{\muxmulevelsetvalue}\right)} 
		\left\| 
			\frac{1}{\Lunit \timefunctionarg{\muxmulevelsetvalue}} (\muX \RRiemann) \tanderY^N \mytr_{\gtorus}\upchi 
		\right\|_{L^2\left( \hypthreearg{\timefunction'}{[-\rightu,u]}{\muxmulevelsetvalue}\right)} 
	\, \mathrm{d}\timefunction'.
\end{split}
\end{align}
Using the sharpened coerciveness bound 
\eqref{E:L2WAVECONTROLWITHBETTERCONSTANTIFMUSMALL}, we find that: 
\begin{align} \label{E:SECONDSTEPSPACETIMEBOUNDSMOSTDIFFICULTWAVEPRODUCT}
\mbox{RHS~\eqref{E:FIRSTSTEPSPACETIMEBOUNDSMOSTDIFFICULTWAVEPRODUCT}}
& \leq
\frac{2}{\sqrt{1.99}} 
\int_{\timefunction' = \timefunction_0}^{\timefunction}  
	\hypersurfacecontrolwave_N^{1/2}(\timefunction',u) 
	\left\| 
		\frac{1}{\Lunit \timefunctionarg{\muxmulevelsetvalue}} (\muX \RRiemann) \tanderY^N \mytr_{\gtorus}\upchi 
	\right\|_{L^2\left(\hypthreearg{\timefunction'}{[-\rightu,u]}{\muxmulevelsetvalue}\right)} 
\, \mathrm{d}\timefunction'.
\end{align}
We now insert the estimate \eqref{E:L2WAVEESTIMATEMOSTDIFFICULTPRODUCT} 
into RHS~\eqref{E:SECONDSTEPSPACETIMEBOUNDSMOSTDIFFICULTWAVEPRODUCT}.
The desired estimate \eqref{E:SPACETIMEBOUNDSMOSTDIFFICULTWAVEPRODUCT}
then follows from a series of standard applications of Young's inequality,
as we now explain.
We will control several representative terms in detail and leave
the remaining details to the reader.
First, the $\timefunction$-integrals of the product of 
$\frac{2}{\sqrt{1.99}} \hypersurfacecontrolwave_N^{1/2}(\timefunction',u)$ 
and the first four terms on RHS~\eqref{E:L2WAVEESTIMATEMOSTDIFFICULTPRODUCT} are clearly bounded 
by the first four terms on RHS~\eqref{E:SPACETIMEBOUNDSMOSTDIFFICULTWAVEPRODUCT} as desired. 
Next, we observe that the $\timefunction$-integrals of the product of 
$\frac{2}{\sqrt{1.99}} \hypersurfacecontrolwave_N^{1/2}(\timefunction',u)$ 
and the single $\timefunction$-integrals on the seventh line of RHS~\eqref{E:L2WAVEESTIMATEMOSTDIFFICULTPRODUCT} are
bounded in magnitude by:
\begin{align} \label{E:SPACETIMEBOUNDSMOSTDIFFICULTWAVEPRODUCTINTERMEDIATESTEP2}
\begin{split}
& 
\lesssim
\int_{\timefunction' = \timefunction_0}^{\timefunction}  
	\left\lbrace 
		\frac{1}{|\timefunction'|^{1/3}} 
		\hypersurfacecontrolwave_N^{1/2}(\timefunction',u) 
	\right\rbrace 
	\left\lbrace 
		\frac{1}{|\timefunction'|^{2/3}}  
	\int_{\timefunction'' = \timefunction_0}^{\timefunction'} 
		\left[ 
				\hypersurfacecontrolVortVort_N^{1/2} 
				+ 
				\hypersurfacecontrolDivGradEnt_N^{1/2} 
				+
				\hypersurfacecontrolVort_N^{1/2} 
				+ 
				\hypersurfacecontrolGradEnt_N^{1/2}
		\right]
		(\timefunction'',u) 
	\, \mathrm{d}\timefunction'' 
	\right\rbrace 
\, \mathrm{d} \timefunction' 
	\\
& \lesssim 
\int_{\timefunction' = \timefunction_0}^{\timefunction} 
	\frac{1}{|\timefunction'|^{2/3}}  
	\hypersurfacecontrolwave_N(\timefunction',u) 
\, \mathrm{d} \timefunction' 
+  
\int_{\timefunction' = \timefunction_0}^{\timefunction} 
	\frac{1}{|\timefunction'|^{4/3}} 
	\left\lbrace 
		\int_{\timefunction'' = \timefunction_0}^{\timefunction'} 
			\left[
				\hypersurfacecontrolVortVort_N^{1/2} 
				+ 
				\hypersurfacecontrolDivGradEnt_N^{1/2} 
				+
				\hypersurfacecontrolVort_N^{1/2} 
				+ 
				\hypersurfacecontrolGradEnt_N^{1/2}
				\right]
			(\timefunction'',u) 
		\, \mathrm{d} \timefunction'' 
		\right\rbrace^2 
\, \mathrm{d}\timefunction'.
\end{split}
\end{align}
Accounting for the term $\Errortop(\timefunction,u)$ on RHS~\eqref{E:SPACETIMEBOUNDSMOSTDIFFICULTWAVEPRODUCT},
we conclude that 
$\mbox{RHS~\eqref{E:SPACETIMEBOUNDSMOSTDIFFICULTWAVEPRODUCTINTERMEDIATESTEP2}}
\lesssim 
\mbox{RHS~\eqref{E:SPACETIMEBOUNDSMOSTDIFFICULTWAVEPRODUCT}}
$
as desired.
Moreover, since 
$
\hypersurfacecontrolwave_{[1,N]}^{1/2}(\timefunction,u)$,
$\hypersurfacecontrolVortVort_N^{1/2}(\timefunction,u)$,
$\hypersurfacecontrolDivGradEnt_N^{1/2}(\timefunction,u)$,
$\hypersurfacecontrolVort_N^{1/2}(\timefunction,u)$,
and
$\hypersurfacecontrolGradEnt_N^{1/2}(\timefunction,u)$
are increasing in their arguments, the terms involving double $\timefunction$-integrals on the fifth line of 
RHS~\eqref{E:L2WAVEESTIMATEMOSTDIFFICULTPRODUCT} are bounded by 
the single $\timefunction$-integrals on the fifth line of RHS~\eqref{E:L2WAVEESTIMATEMOSTDIFFICULTPRODUCT}
plus 
$\frac{C}{|\timefunction|} 
	\int_{\timefunction' = \timefunction_0}^{\timefunction} 
		\ln\left( |\timefunction'|^{-1}\right)
		\hypersurfacecontrolwave_{[1,N]}^{1/2}(\timefunction',u) 
	\, \mathrm{d} \timefunction' 
$.
Hence, the integral of the product of 
$\frac{2}{\sqrt{1.99}} \hypersurfacecontrolwave_N^{1/2}(\timefunction',u)$ 
and the double $\timefunction$-integrals on the fifth line of 
RHS~\eqref{E:L2WAVEESTIMATEMOSTDIFFICULTPRODUCT} are bounded by 
RHS~\eqref{E:SPACETIMEBOUNDSMOSTDIFFICULTWAVEPRODUCTINTERMEDIATESTEP2} plus:
\begin{align} \label{E:ADDITIONALENERGYESTIMATEINTEGRALERRORTERM}
\int_{\timefunction' = \timefunction_0}^{\timefunction} 
		\hypersurfacecontrolwave_N^{1/2}(\timefunction',u) 
		\frac{1}{|\timefunction'|}  
		\int_{\timefunction'' = \timefunction_0}^{\timefunction'}  
			\ln\left( |\timefunction''|^{-1}\right) 
			\hypersurfacecontrolwave_{[1,N]}(\timefunction'',u) 
		\, \mathrm{d} \timefunction'' 
	\, \mathrm{d}\timefunction'.
\end{align}
Using the trivial bound $\ln\left( |\timefunction''|^{-1}\right) \lesssim \frac{1}{|\timefunction''|^{1/2}}$,
we bound \eqref{E:ADDITIONALENERGYESTIMATEINTEGRALERRORTERM} by
$
\lesssim
\int_{\timefunction' = \timefunction_0}^{\timefunction} 
		\hypersurfacecontrolwave_N^{1/2}(\timefunction',u) \frac{1}{|\timefunction'|}  
		\int_{\timefunction'' = \timefunction_0}^{\timefunction'}  
			\frac{1}{|\timefunction''|^{1/2}} 
			\hypersurfacecontrolwave_{[1,N]}^{1/2}\timefunction'',u) 
		\, \mathrm{d} \timefunction'' 
	\, \mathrm{d}\timefunction'
$,
which in turn is bounded by the term $\Errortop(\timefunction,u)$ on RHS~\eqref{E:SPACETIMEBOUNDSMOSTDIFFICULTWAVEPRODUCT}
(more precisely, by the third-from-last term on RHS~\eqref{E:ERRORTOPORDERWAVEESTIMATES}).
As our last representative term, we note that the integral of the product of 
$\frac{2}{\sqrt{1.99}} \hypersurfacecontrolwave_N^{1/2}(\timefunction',u)$ 
and the term 
$\frac{C\initialsmall}{|\timefunction|^{3/2}}$ term on RHS~\eqref{E:L2WAVEESTIMATEMOSTDIFFICULTPRODUCT} 
is:
\begin{align}
\begin{split}
&
\lesssim 
\int_{\timefunction' = \timefunction_0}^{\timefunction} 
	\left\lbrace 
		\frac{\initialsmall}{|\timefunction'|^{5/4}}
	\right\rbrace 
	\left\lbrace 
		\frac{1}{|\timefunction'|^{1/4}} 
		\hypersurfacecontrolwave_N^{1/2}(\timefunction',u) 
	\right\rbrace  
\, \mathrm{d} \timefunction' 
	\\
& 
\lesssim \initialsmall^2 
\int_{\timefunction' = \timefunction_0}^{\timefunction} 
	\frac{1}{|\timefunction'|^{5/2}} 
\, \mathrm{d}\timefunction' 
+  
\int_{\timefunction' = \timefunction_0}^{\timefunction} 
	\frac{1}{|\timefunction'|^{1/2}}\hypersurfacecontrolwave_N 
\, \mathrm{d}\timefunction' 
	\\
& 
\lesssim \frac{\initialsmall^2}{|\timefunction|^{3/2}} 
	+  
\int_{\timefunction' = \timefunction_0}^{\timefunction} 
	\frac{1}{|\timefunction'|^{2/3}}\hypersurfacecontrolwave_N 
\, \mathrm{d}\timefunction' 
\lesssim 
\Errortop(\timefunction,u),
\end{split}
\end{align} 
as desired.

To prove \eqref{E:SPACETIMEBOUNDSMOSTDIFFICULTWAVEPRODUCTPARTIAL} 
for $\Psi \in \wavearraypartial = \{ \LRiemann,v^2,v^3,\Ent\}$,
we first use \eqref{E:LINFINITYIMPROVEMENTAUXTRANSVERSALPDERIVATIVESPARTIALWAVEARRAYSMALL}
to deduce that the magnitude of the integrand on the LHS of \eqref{E:SPACETIMEBOUNDSMOSTDIFFICULTWAVEPRODUCTPARTIAL} is $\lesssim  \fundbootsmall \left|\tanderY^N \Psi \right| \left|\tanderY^N \mytr_{\gtorus} \upchi \right|$. 
We can now argue as in the proof of \eqref{E:SPACETIMEBOUNDSMOSTDIFFICULTWAVEPRODUCT}, 
except that due to the smallness factor $\fundbootsmall$, we do not have to keep careful track
of any boxed constant-involving terms
(we can relegate such terms to the $\fundbootsmall$-multiplied terms in
$\Errortop(\timefunction,u)$ on RHS~\eqref{E:SPACETIMEBOUNDSMOSTDIFFICULTWAVEPRODUCTPARTIAL}),
and on RHS~\eqref{E:SPACETIMEBOUNDSMOSTDIFFICULTWAVEPRODUCTPARTIAL},
we can bound all wave error terms in terms of the full wave energies 
$\hypersurfacecontrolwave$, i.e., without reference to the partial wave energies
$\hypersurfacecontrolwavepartial$.

\end{proof}

\subsubsection{Estimates involving integration by parts with respect to $\argLrough{\muxmulevelsetvalue}$}
\label{SSS:EIKONALTOPORDERENERGYESTIMATEWITHIBP}
In this section, we bound the difficult top-order integrals highlighted in \eqref{E:OVERVIEWLPSIDIFFICULTERRORINTEGAL}.
The proof involves several rather technical steps, and we therefore provide some preliminary lemmas
before proving the main estimates in Lemma~\ref{L:MAINESTIMATESFORDIFFICULTSPACETIMEERRORINTEGRALSINVOLVINGIBPWRTL}.
%
%
%
We will control the error integrals \eqref{E:OVERVIEWLPSIDIFFICULTERRORINTEGAL} 
by integrating by parts with respect to the rough null vectorfield 
$\argLrough{\muxmulevelsetvalue}$ defined in \eqref{E:LROUGH}.
As we will see in the proof of Lemma~\ref{L:MAINESTIMATESFORDIFFICULTSPACETIMEERRORINTEGRALSINVOLVINGIBPWRTL},
the analysis fundamentally relies on the decomposition
$
\tanderY^N \mytr_{\gtorus}\upchi 
= 
\Yvf{A} \tanderY^{N-1} \mytr_{\gtorus} \upchi 
= 
\Yvf{A} 
\partialmodquant{\tanderY^{N-1}} 
- 
\Yvf{A}
\partialmodquantinhom{\tanderY^{N-1}}
$,
where
$\partialmodquant{\tanderY^{N-1}}$
is the partially modified quantity from Def.\,\ref{D:FULLYANDPARTIALLYMODIFIEDQUANTITIES}.
The integrals involving $\partialmodquant{\tanderY^{N-1}}$ are the most difficult to estimate
and are the ones that we treat via integration by parts with respect to $\argLrough{\muxmulevelsetvalue}$,
specifically by invoking the identity \eqref{E:KEYIBPIDENTIFYFORWAVEEQUATIONENERGYESTIMATES}.

Before bounding the error integrals on RHS~\eqref{E:KEYIBPIDENTIFYFORWAVEEQUATIONENERGYESTIMATES},
we first establish a preliminary lemma in which we handle the most difficult part of the analysis.

\begin{lemma}[Difficult top-order hypersurface $L^2$ estimates related to integration by parts with respect to $\argLrough{\muxmulevelsetvalue}$] 
\label{L:DIFFICULTHYPERSURFACEIBPESTIMATE}
Let $N = \Ntop$, and let $\partialmodquant{\tanderY^{N-1}}$ be the partially modified quantity defined by
\eqref{E:PARTIALMODIFIEDQUANTITY}.
The following estimates hold for 
$(\timefunction,u) \in [\timefunction_0,\timefunctionboot) \times [-\rightu,\leftu]$:
\begin{subequations}
\begin{align}
\begin{split} \label{E:PARTIALMODIFIEDQUANTITYPRELIMINARYSPACETIMEL2ESTIMATE} 
	 \left\| 
		\frac{1}{\sqrt\upmu} 
		(\muX \RRiemann) 
		\argLrough{\muxmulevelsetvalue} 
		\partialmodquant{\tanderY^{N-1}} 
	\right\|_{L^2\left(\hypthreearg{\timefunction}{[-\rightu,u]}{\muxmulevelsetvalue}\right)} 
	&  
	\leq 
	\boxed{2.89} 
	\frac{1}{|\timefunction|} 
	\hypersurfacecontrolwave_N^{1/2}(\timefunction,u) 
		 \\
	&  \ \
	+ 
	\frac{C_*}{|\timefunction|} (\hypersurfacecontrolwavepartial_N)^{1/2}(\timefunction,u) 
	+ 
	\frac{C \fundbootsmall}{|\timefunction|} \hypersurfacecontrolwave_{[1,N]}^{1/2}(\timefunction,u) 
	+  
	\frac{C }{|\timefunction|^{1/2}} \hypersurfacecontrolwave_N^{1/2}(\timefunction,u) 
		\\
	& \ \
	+ 
	\frac{C}{|\timefunction|} \hypersurfacecontrolwave_{[1,N-1]}^{1/2}(\timefunction,u) 
		+ 
	\frac{C\initialsmall}{|\timefunction|^{1/2}}, 
\end{split}
			\\
\begin{split} 	\label{E:PARTIALMODIFIEDQUANTITYPRELIMINARYHYPERSURFACEL2ESTIMATE}
	\left\| \frac{1}{\sqrt\upmu} 
		(\muX \RRiemann) 
		\partialmodquant{\tanderY^{N-1}} 
	\right\|_{L^2\left(\hypthreearg{\timefunction}{[-\rightu,u]}{\muxmulevelsetvalue}\right)} 
	& \leq 
	\boxed{2.89} 
	\frac{1}{|\timefunction|^{1/2}} 
	\int_{\timefunction' = \timefunction_0}^{\timefunction} 
		\frac{1}{|\timefunction'|^{1/2}} 
		\hypersurfacecontrolwave_N^{1/2}(\timefunction',u) 
	\, \mathrm{d} \timefunction' 
		\\
	& \ \
	+ 
	\frac{C_*}{|\timefunction|^{1/2}}  
	\int_{\timefunction' = \timefunction_0}^{\timefunction} 
		\frac{1}{|\timefunction'|^{1/2}} 
		(\hypersurfacecontrolwavepartial_{[1,N]})^{1/2}(\timefunction',u) 
	\, \mathrm{d} \timefunction'  
		\\
	& \ \
	+ 
	\frac{C \fundbootsmall}{|\timefunction|^{1/2}}  
	\int_{\timefunction' = \timefunction_0}^{\timefunction} 
		\frac{1}{|\timefunction'|^{1/2}} 
		\hypersurfacecontrolwave_{[1,N]}^{1/2}(\timefunction',u) 
	\, \mathrm{d} \timefunction'  
		\\
	& \ \
	 + 
	C 
	\frac{1}{|\timefunction|^{1/2}} 
	\int_{\timefunction' = \timefunction_0}^{\timefunction} 
			\hypersurfacecontrolwave_{[1,N]}^{1/2}(\timefunction',u) 
	\, \mathrm{d} \timefunction'
		\\
	& \ \
	 + 
	C \int_{\timefunction' = \timefunction_0}^{\timefunction} 
		\frac{1}{|\timefunction'|^{1/2}} 
		\hypersurfacecontrolwave_{[1,N]}^{1/2}(\timefunction',u) 
	\, \mathrm{d} \timefunction'  
		\\
	& \ \
	 + 
	\frac{C}{|\timefunction|^{1/2}}  
	\int_{\timefunction' = \timefunction_0}^{\timefunction} 
		\frac{1}{|\timefunction'|^{1/2}} 
		\hypersurfacecontrolwave_{[1,N-1]}^{1/2}(\timefunction',u) 
	\, \mathrm{d} \timefunction' 
	+ 
	\frac{C \initialsmall}{|\timefunction|^{1/2}}.
\end{split}
\end{align}
\end{subequations}

Moreover, we have the following less precise estimates:
\begin{subequations}
\begin{align} \label{E:LDERIVATIVEOFPARTIALMODIFIEDQUANTITYPRELIMINARYHYPERSURFACEL2ESTIMATELESSPRECISE} 
\left\| 
	\argLrough{\muxmulevelsetvalue} \partialmodquant{\tanderY^{N-1}}
\right\|_{L^2\left(\hypthreearg{\timefunction}{[-\rightu,u]}{\muxmulevelsetvalue}\right)}
& \lesssim \frac{1}{|\timefunction|^{1/2}} \hypersurfacecontrolwave_{[1,N]}^{1/2}(\timefunction,u) 
+ 
\initialsmall, 
	\\
\left\| 
	\partialmodquant{\tanderY^{N-1}} 
\right\|_{L^2\left(\hypthreearg{\timefunction}{[-\rightu,u]}{\muxmulevelsetvalue}\right)} 
& 
\lesssim 
\int_{\timefunction' = \timefunction_0}^{\timefunction} 
	\frac{1}{|\timefunction'|^{1/2}} 
	\hypersurfacecontrolwave_{[1,N]}^{1/2}(\timefunction',u) 
\, \mathrm{d} \timefunction' 
+ 
\initialsmall
\lesssim
\hypersurfacecontrolwave_{[1,N]}^{1/2}(\timefunction,u)
+
\initialsmall.
\label{E:PARTIALMODIFIEDQUANTITYPRELIMINARYHYPERSURFACEL2ESTIMATELESSPRECISE}
\end{align}
\end{subequations}
\end{lemma}

\begin{proof}
We first prove \eqref{E:PARTIALMODIFIEDQUANTITYPRELIMINARYSPACETIMEL2ESTIMATE}. 
In the proof, we sometimes silently use
\eqref{E:DISTORTINGWITHROUGHFLOWMAPDOESNOTCHANGEL2NORMSMUCH}
and Cor.\,\ref{C:IMPROVEAUX},
which together imply that the flow map factors $\FlowmapLrougharg{\muxmulevelsetvalue}$ in \eqref{E:POINTWISEESTIMATELROUGHPARTIALMODIFIEDQUANT}
distort $\| \cdot \|_{L^2\left(\hypthreearg{\timefunction}{[-\rightu,u]}{\muxmulevelsetvalue}\right)}$ norms
only by overall factors of $1 + \mathcal{O}(\fundbootsmall)$;
the $\mathcal{O}(\fundbootsmall)$ factors lead to small error terms on the RHS
of the estimates. In particular, to simplify the notation, in the following discussion,
we will sometimes suppress the factors of $\FlowmapLrougharg{\muxmulevelsetvalue}$
in \eqref{E:POINTWISEESTIMATELROUGHPARTIALMODIFIEDQUANT}.
We start by multiplying \eqref{E:POINTWISEESTIMATELROUGHPARTIALMODIFIEDQUANT} by $\frac{1}{\sqrt{\upmu}} \muX \RRiemann$. 
We now consider the product generated by the first term on RHS~\eqref{E:POINTWISEESTIMATELROUGHPARTIALMODIFIEDQUANT}. 
Multiplying the first equality in \eqref{E:PRECISEPOINTWISEESTIMATESTEP7} by $\frac{\sqrt{\upmu}}{2}$
and using definition 
\eqref{E:LROUGH},
the bootstrap assumptions, 
Cor.\,\ref{C:IMPROVEAUX}, 
and \eqref{E:MINVALUEOFMUONFOLIATION},
we deduce: 
	\begin{align} \label{E:PARTIALMODIFIEDQUANTITYPRELIMINARYSPACETIMEL2ESTIMATEINTERMEDIATE} 
		\frac{1}{2 \Lunit \timefunctionarg{\muxmulevelsetvalue}} 
		\frac{1}{ \sqrt{\upmu}} G_{\Lunit \Lunit}^0 \muX \RRiemann 
		= 
		\frac{\argLrough{\muxmulevelsetvalue} \upmu}{\sqrt{\upmu}} 
		+ 
		\mathcal{O}(\fundbootsmall) \frac{1}{|\timefunction|^{1/2}}. 
	\end{align}
Using \eqref{E:PARTIALMODIFIEDQUANTITYPRELIMINARYSPACETIMEL2ESTIMATEINTERMEDIATE} 
to substitute for the product generated by $\frac{1}{\sqrt{\upmu}} \muX \RRiemann$ times the first term 
on RHS~\eqref{E:POINTWISEESTIMATELROUGHPARTIALMODIFIEDQUANT}, 
we pointwise bound this difficult product as follows,
where $\smallneighborhoodofcreasetwoarg{[\timefunction_0,\timefunctionboot]}{\muxmulevelsetvalue}$
is the set from \eqref{E:SMALLNEIGHBORHOOD}:
\begin{align} 
\begin{split} \label{E:PARTIALMODIFIEDQUANTITYPRELIMINARYSPACETIMEL2ESTIMATEINTERMEDIATE2}
		&
		\lesssim
		\left| 
			\frac{\argLrough{\muxmulevelsetvalue} \upmu}{\upmu} 
			\mathbf{1}_{\{\hypthreearg{\timefunction}{[-\rightu,u]}{\muxmulevelsetvalue} 
				\cap \smallneighborhoodofcreasetwoarg{[\timefunction_0,\timefunctionboot]}{\muxmulevelsetvalue}\}}
		\right| 
			\cdot 
			\left| 
				\sqrt{\upmu} \angLap \tanderY^{N-1} \RRiemann 
			\right|  
			\\
	& \ \ + 
		\left| 
			\frac{\argLrough{\muxmulevelsetvalue} \upmu}{\upmu} 
			\mathbf{1}_{\{\hypthreearg{\timefunction}{[-\rightu,u]}{\muxmulevelsetvalue} 
				\setminus \smallneighborhoodofcreasetwoarg{[\timefunction_0,\timefunctionboot]}{\muxmulevelsetvalue}\}}\right| 
			\cdot 
			\left| 
				\sqrt{\upmu} \angLap \tanderY^{N-1} \RRiemann 
			\right|  
		+ 
		\mathcal{O}(\fundbootsmall) 
		\frac{1}{|\timefunction|^{1/2} \sqrt{\upmu}} 
		\left| \sqrt{\upmu} \angLap \tanderY^{N-1} \RRiemann \right|. 
	\end{split}
	\end{align}
Using \eqref{E:EASYREGIONLOWERBOUNDFORMU}
and the crude estimate $|\argLrough{\muxmulevelsetvalue} \upmu| \lesssim 1$
(see \eqref{E:LINFINITYIMPROVEMENTAUXNULLDERIVATIVEMUANDTRANSVERSALDERIVATIVES}), 
we deduce the pointwise bound
$\left| \frac{\argLrough{\muxmulevelsetvalue} \upmu}{\upmu} \mathbf{1}_{\{\hypthreearg{\timefunction}{[-\rightu,u]}{\muxmulevelsetvalue} \setminus \smallneighborhoodofcreasetwoarg{[\timefunction_0,\timefunctionboot]}{\muxmulevelsetvalue}\}}\right| \lesssim 1$. 
Hence, using the pointwise comparison estimate 
\eqref{E:SMOOTHANGULARHESSIANOFFPOINTWISEBOUNDEDBYCOMMUTATORVECTORFIELDS}, 
we see that the second product on RHS~\eqref{E:PARTIALMODIFIEDQUANTITYPRELIMINARYSPACETIMEL2ESTIMATEINTERMEDIATE2} 
is pointwise bounded by
$ 
\lesssim \left|\sqrt{\upmu} \angrmD \tanderY^{[N-1,N]} \RRiemann \right|
$. 
Again using \eqref{E:MINVALUEOFMUONFOLIATION}, we find that
the third product on RHS~\eqref{E:PARTIALMODIFIEDQUANTITYPRELIMINARYSPACETIMEL2ESTIMATEINTERMEDIATE2} is pointwise bounded by
$\lesssim \frac{\fundbootsmall}{|\timefunction|} \left|\sqrt{\upmu} \angrmD \tanderY^{[N-1,N]} \RRiemann \right|$. 
From these pointwise bounds and \eqref{E:COERCIVENESSOFHYPERSURFACECONTROLWAVE}, we find that the
$\| \cdot \|_{L^2\left(\hypthreearg{\timefunction}{[-\rightu,u]}{\muxmulevelsetvalue}\right)}$ 
norms of the second and third terms in 
\eqref{E:PARTIALMODIFIEDQUANTITYPRELIMINARYSPACETIMEL2ESTIMATEINTERMEDIATE2} are $\leq$ the sum of the
third, fourth, and fifth terms on RHS~\eqref{E:PARTIALMODIFIEDQUANTITYPRELIMINARYSPACETIMEL2ESTIMATE}, 
as desired.
Next, we use 
\eqref{E:MINVALUEOFMUONFOLIATION},
\eqref{E:WIDETILDELMUISALMOSTMINUSONEINSMALLNEIGHBORHOOD}, 
\eqref{E:SMOOTHANGULARHESSIANOFFPOINTWISEBOUNDEDBYCOMMUTATORVECTORFIELDS},
and Cor.\,\ref{C:IMPROVEAUX} to pointwise bound the first term in 
\eqref{E:PARTIALMODIFIEDQUANTITYPRELIMINARYSPACETIMEL2ESTIMATEINTERMEDIATE2} as follows:
\begin{align} 
\begin{split} \label{E:SHARPPOINTWISEBOUNDFIRSTTERMINPARTIALMODIFIEDQUANTITYPRELIMINARYSPACETIMEL2ESTIMATEINTERMEDIATE2}
\left| 
	\frac{\argLrough{\muxmulevelsetvalue} \upmu}{\upmu} \mathbf{1}_{\{\hypthreearg{\timefunction}{[-\rightu,u]}{\muxmulevelsetvalue} 
	\cap \smallneighborhoodofcreasetwoarg{[\timefunction_0,\timefunctionboot]}{\muxmulevelsetvalue}\}}
\right| 
\cdot 
\left| 
	\sqrt{\upmu} \angLap \tanderY^{N-1} \RRiemann 
\right|  
& 
\leq 
\sqrt{2}
\left\lbrace
	1 + \mathcal{O}_{\mydiam}(\mathring{\upalpha})
\right\rbrace
\frac{1.01}{|\timefunction|}
\sqrt{
\sum_{A = 2}^3 
\left|\sqrt{\upmu} \angrmD \Yvf{A} \tanderY^{N-1} \RRiemann \right|_{\gtorus}^2}
	\\
& \ \
+
\frac{C}{|\timefunction|}
\left|\sqrt{\upmu} \angrmD \tanderY^{N-1} \RRiemann \right|_{\gtorus}.
\end{split}
\end{align}
From \eqref{E:SHARPPOINTWISEBOUNDFIRSTTERMINPARTIALMODIFIEDQUANTITYPRELIMINARYSPACETIMEL2ESTIMATEINTERMEDIATE2}, 
\eqref{E:COERCIVENESSOFHYPERSURFACECONTROLWAVE}
and in particular its implication:
\begin{align} \label{E:PARTIALMODIFIEDQUANTITYPRELIMINARYSPACETIMEL2ESTIMATEINTERMEDIATE3}
		\left\| 
			\sqrt{\sum_{A = 2}^3 
				\left|\sqrt{\upmu} \angrmD \Yvf{A} \tanderY^{N-1} \RRiemann \right|_{\gtorus}^2}
		\right\|_{L^2\left(\hypthreearg{\timefunction}{[-\rightu,u]}{\muxmulevelsetvalue}\right)} 
		& \leq  
		\sqrt{\frac{2}{0.49}} \hypersurfacecontrolwave_N^{1/2}(\timefunction,u),
	\end{align}
our assumption that $\mathring{\upalpha}$ 
and $\fundbootsmall$ are sufficiently small,
and the inequality $\sqrt{2} \times 1.01 \times \sqrt{\frac{2}{0.49}} < 2.886$,
we deduce that the
$\| \cdot \|_{L^2\left(\hypthreearg{\timefunction}{[-\rightu,u]}{\muxmulevelsetvalue}\right)}$ norm of the first term in 
\eqref{E:PARTIALMODIFIEDQUANTITYPRELIMINARYSPACETIMEL2ESTIMATEINTERMEDIATE2} 
is 
$\leq$ the sum of the $\boxed{2.89}$-multiplied 
term on RHS~\eqref{E:PARTIALMODIFIEDQUANTITYPRELIMINARYSPACETIMEL2ESTIMATE} 
and the $\frac{C}{|\timefunction|} \hypersurfacecontrolwave_{[1,N-1]}^{1/2}(\timefunction,u)$ term. 
We have therefore obtained the desired estimates for the product of
$\frac{1}{\sqrt{\upmu}} \muX \RRiemann$ and the first term 
on RHS~\eqref{E:POINTWISEESTIMATELROUGHPARTIALMODIFIEDQUANT}.
Combining similar but simpler arguments with \eqref{E:L2ESTIMATESFORWAVEVARIABLESONROUGHHYPERSURFACELOSSOFONEDERIVATIVE}
and the pointwise bound $\frac{1}{\upmu^{1/2}}|\muX \RRiemann| \lesssim \frac{1}{|\timefunction|^{1/2}}$
implied by \eqref{E:LINFINITYIMPROVEMENTAUXTRANSVERSALPDERIVATIVESRRIEMANNLARGE} and \eqref{E:MINVALUEOFMUONFOLIATION},
we find that the $\| \cdot \|_{L^2\left(\hypthreearg{\timefunction}{[-\rightu,u]}{\muxmulevelsetvalue}\right)}$ norms 
of the product of $\frac{1}{\sqrt{\upmu}} \muX \RRiemann$ and the terms
$C_* \left|\angLap \tanderY^{N-1} \wavearraypartial \right|$ and $C \fundbootsmall \left| \tander^{[1,N+1]}\wavearray \right|$ 
on RHS~\eqref{E:POINTWISEESTIMATELROUGHPARTIALMODIFIEDQUANT}
are bounded by the sum of the last five terms on
RHS~\eqref{E:PARTIALMODIFIEDQUANTITYPRELIMINARYSPACETIMEL2ESTIMATE}.
Finally, using the pointwise bound $\frac{1}{\upmu^{1/2}}|\muX \RRiemann| \lesssim \frac{1}{|\timefunction|^{1/2}}$
noted above,
\eqref{E:L2ESTIMATESFORWAVEVARIABLESONROUGHHYPERSURFACELOSSOFONEDERIVATIVE},
\eqref{E:PRELIMINARYEIKONALWITHOUTL},
\eqref{E:MINVALUEOFMUONFOLIATION}, 
and the fact that the $\hypersurfacecontrolwave_M(\timefunction,u)$ are increasing in their arguments, 
we bound the $\| \cdot \|_{L^2\left(\hypthreearg{\timefunction}{[-\rightu,u]}{\muxmulevelsetvalue}\right)}$ norm
of the product of $\frac{1}{\sqrt{\upmu}} \muX \RRiemann$ 
and the last term $C \left|\Tanset^{[1,N]} \controlvars \right|$ 
on RHS~\eqref{E:POINTWISEESTIMATELROUGHPARTIALMODIFIEDQUANT},
by 
$
\lesssim
\initialsmall 
\frac{1}{|\timefunction|^{1/2}} 
+ 
\frac{1}{|\timefunction|^{1/2}} 
\int_{\timefunction' = \timefunction_0}^{\timefunction} 
	\frac{1}{|\timefunction'|^{1/2}} 
	\hypersurfacecontrolwave_{[1,N]}^{1/2}(\timefunction',u) 
\, \mathrm{d}\timefunction' 
\lesssim 
\initialsmall\frac{1}{|\timefunction|^{1/2}} 
+ 
\frac{1}{|\timefunction|^{1/2}} 
\hypersurfacecontrolwave_{[1,N]}^{1/2}(\timefunction,u)
$. We have therefore proved \eqref{E:PARTIALMODIFIEDQUANTITYPRELIMINARYSPACETIMEL2ESTIMATE}.

We now prove \eqref{E:PARTIALMODIFIEDQUANTITYPRELIMINARYHYPERSURFACEL2ESTIMATE}. 
We start by multiplying \eqref{E:POINTWISEESTIMATEPARTIALMODIFIEDQUANTITY} by 
$\frac{1}{\sqrt{\upmu}} \muX \RRiemann$. 
We now focus on the most difficult product, which is generated by the second term 
on RHS~\eqref{E:POINTWISEESTIMATEPARTIALMODIFIEDQUANTITY}, i.e.,
\begin{align} \label{E:PARTIALMODIFIEDQUANTITYPRELIMINARYHYPERSURFACEL2ESTIMATEINTERMEDIATE}
		\frac{1}{2} 
		\left\lbrace 
			\left|
				\frac{1}{\Lunit \timefunctionarg{\muxmulevelsetvalue}}
				\frac{1}{\sqrt{\upmu}} 
				G_{\Lunit \Lunit}^0 
				\muX \RRiemann 
			\right| 
		\circ \FlowmapLrougharg{\muxmulevelsetvalue}(\timefunction,u,x^2,x^3) 
		\right\rbrace 
		\int_{\timefunction' = \timefunction_0}^{\timefunction}  
			\left| \angLap \tanderY^{N-1} \RRiemann \right|   
			\circ 
			\FlowmapLrougharg{\muxmulevelsetvalue}(\timefunction',u,x^2,x^3) 
		\, \mathrm{d}\timefunction'.
	\end{align}
Using \eqref{E:PARTIALMODIFIEDQUANTITYPRELIMINARYSPACETIMEL2ESTIMATEINTERMEDIATE}, 
we substitute 
$\frac{\argLrough{\muxmulevelsetvalue} \upmu}{\sqrt{\upmu}} 
+ 
\mathcal{O}(\fundbootsmall) \frac{1}{|\timefunction|^{1/2}}
$ 
for the first product 
$\left| \frac{1}{2\Lunit \timefunctionarg{\muxmulevelsetvalue}}\frac{1}{\sqrt{\upmu}}  G_{\Lunit \Lunit}^0 \muX \RRiemann\right|$ 
in \eqref{E:PARTIALMODIFIEDQUANTITYPRELIMINARYHYPERSURFACEL2ESTIMATEINTERMEDIATE}. 
Next, we bound the $\| \cdot \|_{L^2\left(\hypthreearg{\timefunction}{[-\rightu,u]}{\muxmulevelsetvalue}\right)}$ norm of 
the product generated by 
the factor $\mathcal{O}(\fundbootsmall) \frac{1}{|\timefunction|^{1/2}}$ 
by 
using  
\eqref{E:SMOOTHANGULARHESSIANOFFPOINTWISEBOUNDEDBYCOMMUTATORVECTORFIELDS}, 
\eqref{E:L2ONROUGHCONSTANTTIMEHYPERSURFACEOFROUGHTIMEINTEGRALWITHFLOWMAPFACTORSBOUND}, 
\eqref{E:MINVALUEOFMUONFOLIATION}, 
and \eqref{E:COERCIVENESSOFHYPERSURFACECONTROLWAVE}. 
We find that these error terms are $\leq$ 
the $C \fundbootsmall$-multiplied 
term on the third line of RHS~\eqref{E:PARTIALMODIFIEDQUANTITYPRELIMINARYHYPERSURFACEL2ESTIMATE}. 
Next, we use the triangle inequality,
\eqref{E:MINVALUEOFMUONFOLIATION},
\eqref{E:WIDETILDELMUISALMOSTMINUSONEINSMALLNEIGHBORHOOD},
\eqref{E:EASYREGIONLOWERBOUNDFORMU},
and the crude estimate
$|\argLrough{\muxmulevelsetvalue} \upmu| \lesssim 1$ noted earlier
to deduce the following pointwise bound:
\begin{align} 
\begin{split} \label{E:ROUGHLMUOVERSQTMUKEYPOINTWISEBOUND}
\left|\frac{\argLrough{\muxmulevelsetvalue} \upmu}{\sqrt{\upmu}} \right| 
&
\leq 
\left|
	\frac{\argLrough{\muxmulevelsetvalue} \upmu}{\sqrt{\upmu}} \mathbf{1}_{\{\hypthreearg{\timefunction}{[-\rightu,u]}{\muxmulevelsetvalue} 
	\cap \smallneighborhoodofcreasetwoarg{[\timefunction_0,\timefunctionboot]}{\muxmulevelsetvalue}\}} 
\right|  
+  
\left|
	\frac{\argLrough{\muxmulevelsetvalue} \upmu}{\sqrt{\upmu}} 
	\mathbf{1}_{\{\hypthreearg{\timefunction}{[-\rightu,u]}{\muxmulevelsetvalue} 
	\setminus \smallneighborhoodofcreasetwoarg{[\timefunction_0,\timefunctionboot]}{\muxmulevelsetvalue}\}} 
\right|
	\\
& \leq 
\frac{1.01}{|\timefunction|^{1/2}}
+
C.
\end{split}
\end{align}
From \eqref{E:ROUGHLMUOVERSQTMUKEYPOINTWISEBOUND}
and
\eqref{E:MINVALUEOFMUONFOLIATION},
it follows that
$\left|\frac{\argLrough{\muxmulevelsetvalue} \upmu}{\sqrt{\upmu}} \right|  
\times 
\mbox{RHS~\eqref{E:PARTIALMODIFIEDQUANTITYPRELIMINARYHYPERSURFACEL2ESTIMATEINTERMEDIATE}}$
is pointwise bounded by:
\begin{align} 
\begin{split} \label{E:ALTERNATEPARTIALMODIFIEDQUANTITYPRELIMINARYHYPERSURFACEL2ESTIMATEINTERMEDIATE}
		&
		\leq
		\frac{1.01}{|\timefunction|^{1/2}}
		\int_{\timefunction' = \timefunction_0}^{\timefunction}  
			\frac{1}{|\timefunction'|^{1/2}}
			\left| \sqrt{\upmu} \angLap \tanderY^{N-1} \RRiemann \right|   
			\circ 
			\FlowmapLrougharg{\muxmulevelsetvalue}(\timefunction',u,x^2,x^3) 
		\, \mathrm{d}\timefunction'
			\\
	& \ \
		+
		C
		\int_{\timefunction' = \timefunction_0}^{\timefunction}  
			\frac{1}{|\timefunction'|^{1/2}}
			\left| \sqrt{\upmu} \angLap \tanderY^{N-1} \RRiemann \right|   
			\circ 
			\FlowmapLrougharg{\muxmulevelsetvalue}(\timefunction',u,x^2,x^3) 
		\, \mathrm{d}\timefunction'.
	\end{split}
	\end{align}
Using
\eqref{E:SMOOTHANGULARHESSIANOFFPOINTWISEBOUNDEDBYCOMMUTATORVECTORFIELDS}, 
\eqref{E:L2ONROUGHCONSTANTTIMEHYPERSURFACEOFROUGHTIMEINTEGRALWITHFLOWMAPFACTORSBOUND},
and
\eqref{E:COERCIVENESSOFHYPERSURFACECONTROLWAVE},
we find that the 
$\| \cdot \|_{L^2\left(\hypthreearg{\timefunction}{[-\rightu,u]}{\muxmulevelsetvalue}\right)}$ 
norm of the last term on RHS~\eqref{E:ALTERNATEPARTIALMODIFIEDQUANTITYPRELIMINARYHYPERSURFACEL2ESTIMATEINTERMEDIATE}
is bounded by the terms on the last and penultimate lines of 
RHS~\eqref{E:PARTIALMODIFIEDQUANTITYPRELIMINARYHYPERSURFACEL2ESTIMATE}.
Next, we use
\eqref{E:SMOOTHANGULARHESSIANOFFPOINTWISEBOUNDEDBYCOMMUTATORVECTORFIELDS}
and Cor.\,\ref{C:IMPROVEAUX} to
deduce the following pointwise bound for the 
first term on RHS~\eqref{E:ALTERNATEPARTIALMODIFIEDQUANTITYPRELIMINARYHYPERSURFACEL2ESTIMATEINTERMEDIATE}:
\begin{align}  
\begin{split} \label{E:PARTIALMODIFIEDQUANTITYPRELIMINARYHYPERSURFACEL2ESTIMATEINTERMEDIATE2}
		&
		\frac{1.01}{|\timefunction|^{1/2}} 
		\int_{\timefunction' = \timefunction_0}^{\timefunction}  
			\frac{1}{|\timefunction'|^{1/2}}
			\left| \sqrt{\upmu} \angLap \tanderY^{N-1} \RRiemann \right| 
			\circ 
			\FlowmapLrougharg{\muxmulevelsetvalue}(\timefunction',u,x^2,x^3)
		\, \mathrm{d}\timefunction'
			\\
		&
		\leq
		\frac{1.01}{|\timefunction|^{1/2}} 
		\int_{\timefunction' = \timefunction_0}^{\timefunction}  
			\frac{1}{|\timefunction'|^{1/2}}
			\left\lbrace  
				\sqrt{
				2
				\left[
					1 + \mathcal{O}_{\mydiam}(\mathring{\upalpha})
				\right]
				\sum_{A = 2}^3 
				\left|\sqrt{\upmu} \angrmD \Yvf{A} \tanderY^{N-1} \RRiemann \right|_{\gtorus}^2}
				\right\rbrace 
			\circ 
			\FlowmapLrougharg{\muxmulevelsetvalue}(\timefunction',u,x^2,x^3)
		\, \mathrm{d}\timefunction'
			\\
		& \ \
			+
		\frac{C}{|\timefunction|^{1/2}} 
		\int_{\timefunction' = \timefunction_0}^{\timefunction}  
			\frac{1}{|\timefunction'|^{1/2}}
			\left|
				\sqrt{\upmu} \angrmD \tanderY^{N-1} \RRiemann 
			\right|_{\gtorus}
			\circ 
			\FlowmapLrougharg{\muxmulevelsetvalue}(\timefunction',u,x^2,x^3)
		\, \mathrm{d}\timefunction'.
\end{split}
\end{align}
By combining the same arguments we used to bound
the 
$\| \cdot \|_{L^2\left(\hypthreearg{\timefunction}{[-\rightu,u]}{\muxmulevelsetvalue}\right)}$ 
norm of RHS~\eqref{E:SHARPPOINTWISEBOUNDFIRSTTERMINPARTIALMODIFIEDQUANTITYPRELIMINARYSPACETIMEL2ESTIMATEINTERMEDIATE2}
with 
\eqref{E:L2ONROUGHCONSTANTTIMEHYPERSURFACEOFROUGHTIMEINTEGRALWITHFLOWMAPFACTORSBOUND}
and Cor.\,\ref{C:IMPROVEAUX},
we find that when
$\mathring{\upalpha}$ 
and $\fundbootsmall$ are sufficiently small,
the $\| \cdot \|_{L^2\left(\hypthreearg{\timefunction}{[-\rightu,u]}{\muxmulevelsetvalue}\right)}$ 
norm of RHS~\eqref{E:PARTIALMODIFIEDQUANTITYPRELIMINARYHYPERSURFACEL2ESTIMATEINTERMEDIATE2}
is $\leq$ the sum of the
$\boxed{2.89}$-multiplied 
and $C \fundbootsmall$-multiplied
terms on RHS~\eqref{E:PARTIALMODIFIEDQUANTITYPRELIMINARYHYPERSURFACEL2ESTIMATE}
and the next-to-last term on RHS~\eqref{E:PARTIALMODIFIEDQUANTITYPRELIMINARYHYPERSURFACEL2ESTIMATE}.
We have therefore obtained the desired bound for the 
$\| \cdot \|_{L^2\left(\hypthreearg{\timefunction}{[-\rightu,u]}{\muxmulevelsetvalue}\right)}$ norm of the product of
$\frac{1}{\sqrt{\upmu}} \muX \RRiemann$ and the first term on RHS~\eqref{E:POINTWISEESTIMATEPARTIALMODIFIEDQUANTITY}.
Similarly, by combining the pointwise bound $\frac{1}{\upmu^{1/2}}|\muX \RRiemann| \lesssim \frac{1}{|\timefunction|^{1/2}}$
noted above with
\eqref{E:L2ONROUGHCONSTANTTIMEHYPERSURFACEOFROUGHTIMEINTEGRALWITHFLOWMAPFACTORSBOUND},
\eqref{E:MINVALUEOFMUONFOLIATION},
and \eqref{E:COERCIVENESSOFHYPERSURFACECONTROLWAVEPARTIAL},
we bound the
$\| \cdot \|_{L^2\left(\hypthreearg{\timefunction}{[-\rightu,u]}{\muxmulevelsetvalue}\right)}$ 
norm of the product of 
$\frac{1}{\upmu^{1/2}} \muX \RRiemann$ 
and the term
$
C_* \int_{\timefunction' = \timefunction_0}^{\timefunction}  \left| \angLap \tanderY^{N-1} \wavearraypartial \right|   \circ \FlowmapLrougharg{\muxmulevelsetvalue}(\timefunction',u,x^2,x^3)\, \mathrm{d} \timefunction'
$
from RHS~\eqref{E:POINTWISEESTIMATEPARTIALMODIFIEDQUANTITY} 
by
$\leq$ the sum of the
$C_*$-multiplied term on RHS~\eqref{E:PARTIALMODIFIEDQUANTITYPRELIMINARYHYPERSURFACEL2ESTIMATE}
and the next-to-last term on RHS~\eqref{E:PARTIALMODIFIEDQUANTITYPRELIMINARYHYPERSURFACEL2ESTIMATE}.
Finally, using the bound $\frac{1}{\upmu^{1/2}}|\muX \RRiemann| \lesssim \frac{1}{|\timefunction|^{1/2}}$ noted above,
\eqref{E:L2ONROUGHCONSTANTTIMEHYPERSURFACEOFROUGHTIMEINTEGRALWITHFLOWMAPFACTORSBOUND},
and
\eqref{E:MINVALUEOFMUONFOLIATION},
we bound the $\| \cdot \|_{L^2\left(\hypthreearg{\timefunction}{[-\rightu,u]}{\muxmulevelsetvalue}\right)}$ 
norms of the product of 
$\frac{1}{\upmu^{1/2}} \muX \RRiemann$ and the terms on the last line of RHS~\eqref{E:POINTWISEESTIMATEPARTIALMODIFIEDQUANTITY} 
by
$
\frac{C \fundbootsmall}{|\timefunction|^{1/2}}
\int_{\timefunction' = \timefunction_0}^{\timefunction} 
	\frac{1}{|\timefunction'|^{1/2}}
	\left \| \sqrt{\upmu} \tander^{N+1} \wavearray \right \|_{L^2\left(\hypthreearg{\timefunction}{[-\rightu,u]}{\muxmulevelsetvalue}\right)}
\, \mathrm{d} \timefunction'
$
$
+
\frac{C}{|\timefunction|^{1/2}}
\int_{\timefunction' = \timefunction_0}^{\timefunction} 
	\left\| \tander^{[1,N]} \controlvars  
\right \|_{L^2\left(\hypthreearg{\timefunction}{[-\rightu,u]}{\muxmulevelsetvalue}\right)}
\, \mathrm{d} \timefunction'
$,
and using
\eqref{E:COERCIVENESSOFHYPERSURFACECONTROLWAVE},
\eqref{E:L2ESTIMATESFORWAVEVARIABLESONROUGHHYPERSURFACELOSSOFONEDERIVATIVE},
\eqref{E:PRELIMINARYEIKONALWITHOUTL},
and the fact that the $\hypersurfacecontrolwave_M(\timefunction,u)$ are increasing in their arguments,
we conclude that these time integrals are $\leq$ the sum of the non-boxed-constant-multiplied terms on
RHS~\eqref{E:PARTIALMODIFIEDQUANTITYPRELIMINARYHYPERSURFACEL2ESTIMATE} as desired. 

The estimates 
\eqref{E:LDERIVATIVEOFPARTIALMODIFIEDQUANTITYPRELIMINARYHYPERSURFACEL2ESTIMATELESSPRECISE}--\eqref{E:PARTIALMODIFIEDQUANTITYPRELIMINARYHYPERSURFACEL2ESTIMATELESSPRECISE} can be proved using only a subset of the arguments we gave above; 
we omit the details, which are simpler 
since they \emph{do not} involve sharp constants or delicate
decompositions as in \eqref{E:PARTIALMODIFIEDQUANTITYPRELIMINARYSPACETIMEL2ESTIMATEINTERMEDIATE},
and the LHSs of the estimates are less degenerate by a factor of $\sqrt{\upmu}$ compared to 
\eqref{E:PARTIALMODIFIEDQUANTITYPRELIMINARYSPACETIMEL2ESTIMATE}--\eqref{E:PARTIALMODIFIEDQUANTITYPRELIMINARYHYPERSURFACEL2ESTIMATE}. 
\end{proof}

With the help of the preliminary estimates provided by Lemma~\ref{L:DIFFICULTHYPERSURFACEIBPESTIMATE},
we are now ready to bound the most difficult error integral integrals
that arise when we integrate by parts with respect to $\argLrough{\muxmulevelsetvalue}$
using the identity
\eqref{E:KEYIBPIDENTIFYFORWAVEEQUATIONENERGYESTIMATES}.
Specifically, we bound the first two error integrals on
RHS~\eqref{E:KEYIBPIDENTIFYFORWAVEEQUATIONENERGYESTIMATES}.

\begin{lemma}[Estimates for difficult top-order error integrals related to integration by parts
 with respect to $\argLrough{\muxmulevelsetvalue}$] 
\label{L:ESTIMATESFORDIFFICULTSPACETIMEERRORINTEGRALSINVOLVINGIBPWRTL}
Assume that $N = \Ntop$.
Then the following estimates hold for 
$(\timefunction,u) \in [\timefunction_0,\timefunctionboot) \times [-\rightu,\leftu]$:
\begin{align}
\begin{split} \label{E:MAINSPACETIMEWAVEESTIMATEIBP}  
	& \left| \int_{\twoargMrough{[\timefunction_0,\timefunction],[-\rightu,u]}{\muxmulevelsetvalue}} 
		(1 + 2 \upmu) 
		(\Yvf{A}\tanderY^N \RRiemann) 
		(\muX \RRiemann) 
		\argLrough{\muxmulevelsetvalue} 
	\partialmodquant{\tanderY^{N-1}} 
	\,\volMRoughCoordinates \right| 
		\\
	& \leq
	\boxed{4.13} 
		\int_{\timefunction' = \timefunction_0}^{\timefunction} \frac{1}{|\timefunction'|}  
		\hypersurfacecontrolwave_N (\timefunction',u) 
	\, \mathrm{d} \timefunction' 
		\\
	& \ \
	+ 
	C_* \int_{\timefunction' = \timefunction_0}^{\timefunction}  \frac{1}{|\timefunction'|}
	\hypersurfacecontrolwave^{1/2}_N (\timefunction',u) 
	\left(\hypersurfacecontrolwavepartial_N\right)^{1/2}(\timefunction',u) \, \mathrm{d} \timefunction' 
		\\
	 & \ \ 
		+ \Errortop(\timefunction,u), 
	\end{split}
	\\
	\begin{split}  \label{E:MAINHYPERSURFACEESTIMATEIBP}
	 &  
	\left| 
		\int_{\hypthreearg{\timefunction}{[-\rightu,u]}{\muxmulevelsetvalue}} 
			(1 + 2 \upmu) (\Yvf{A}\tanderY^N \RRiemann) 
			(\muX \RRiemann)  
			\partialmodquant{\tanderY^{N-1}} 
		\, \volRoughHypersurface 
	\right| 
		\\
	 & \leq 
	\boxed{4.13}  
	\frac{1}{|\timefunction|^{1/2}} 
	\hypersurfacecontrolwave_N^{1/2}(\timefunction,u) 
	\int_{\timefunction' = \timefunction_0}^{\timefunction}
		\frac{1}{|\timefunction'|^{1/2}} 
		\hypersurfacecontrolwave_N^{1/2}(\timefunction',u) 
	\, \mathrm{d} \timefunction' 
		\\
	& \ \
	+ C_* \frac{1}{|\timefunction|^{1/2}} 
	\hypersurfacecontrolwave_N^{1/2}(\timefunction,u) 
	\int_{\timefunction' = \timefunction_0}^{\timefunction} 
		\frac{1}{|\timefunction'|^{1/2}} 
		\left(\hypersurfacecontrolwavepartial_N\right)^{1/2}(\timefunction',u) 
	\, \mathrm{d} \timefunction'
		\\
	& \ \ 
		+ \Errortop(\timefunction,u),
\end{split}
\end{align}
where $\Errortop(\timefunction,u)$ satisfies \eqref{E:ERRORTOPORDERWAVEESTIMATES}.

Moreover, for every $\Psi \in \wavearraypartial = \{ \LRiemann,v^2,v^3,\Ent\}$, 
we have the following less degenerate estimates:
 \begin{align} 	\label{E:MAINSPACETIMEWAVEESTIMATEIBPPARTIAL}
	 \left| 
		\int_{\twoargMrough{[\timefunction_0,\timefunction],[-\rightu,u]}{\muxmulevelsetvalue}} 
			(1 + 2 \upmu) 
			(\Yvf{A} \tanderY^N \Psi) 
			(\muX \Psi) 
			\argLrough{\muxmulevelsetvalue} 
			\partialmodquant{\tanderY^{N-1}} 
		\, \volMRoughCoordinates 
	\right| 
	& \lesssim \Errortop(\timefunction,u), 
		\\
	 \left| 
		\int_{\hypthreearg{\timefunction}{[-\rightu,u]}{\muxmulevelsetvalue}} 
			(1 + 2 \upmu) 
			(\Yvf{A}\tanderY^N \Psi) 
			(\muX \Psi) 
			\partialmodquant{\tanderY^{N-1}} 
		\, \volRoughHypersurface 
	\right| 
	& \lesssim \Errortop(\timefunction,u). \label{E:MAINHYPERSURFACEESTIMATEIBPPARTIAL}
 \end{align}

\end{lemma}

\begin{proof}
We first prove \eqref{E:MAINSPACETIMEWAVEESTIMATEIBP}. 
We start by noting the following estimates, 
which follow from the bootstrap assumptions:
$|\muX \RRiemann| \lesssim 1$, $|\upmu| \lesssim 1$.
Also using \eqref{E:POINTWISESEMINORMOFYVECTORFIELDS},
Cor.\,\ref{C:IMPROVEAUX},
\eqref{E:MINVALUEOFMUONFOLIATION}, 
and the Cauchy--Schwarz inequality for integrals,
we bound LHS~\eqref{E:MAINSPACETIMEWAVEESTIMATEIBP} by:
	\begin{align}
	\begin{split} \label{E:PROOFMAINSPACETIMEWAVEESTIMATEIBP}
		& \leq 
			(1 + C_{\mydiam} \mr\upalpha) 
			\int_{\timefunction' = \timefunction_0}^{\timefunction} 
				\left\| 
					\sqrt{\upmu} \angrmD \tanderY^N \RRiemann 
				\right\|_{L^2\left( \hypthreearg{\timefunction'}{[-\rightu,u]}{\muxmulevelsetvalue}\right)} 			\left\| \frac{1}{\sqrt\upmu} (\muX \RRiemann) 
				\argLrough{\muxmulevelsetvalue} \partialmodquant{\tanderY^{N-1}} 
			\right\|_{L^2\left(\hypthreearg{\timefunction'}{[-\rightu,u]}{\muxmulevelsetvalue}\right)} 
			\, \mathrm{d} \timefunction' 
		\\
		& \ \
			+  
			C  
			\int_{\timefunction' = \timefunction_0}^{\timefunction}  
			\left\| 
				\sqrt{\upmu} \angrmD \tanderY^N \RRiemann 
			\right\|_{L^2\left( \hypthreearg{\timefunction'}{[-\rightu,u]}{\muxmulevelsetvalue}\right)} 		
			\left\| 
				\argLrough{\muxmulevelsetvalue} \partialmodquant{\tanderY^{N-1}} 
			\right\|_{L^2\left(\hypthreearg{\timefunction'}{[-\rightu,u]}{\muxmulevelsetvalue}\right)} 
			\, \mathrm{d} \timefunction'.
	\end{split}
	\end{align}
The desired estimate \eqref{E:MAINSPACETIMEWAVEESTIMATEIBP} now
follows from inserting \eqref{E:PARTIALMODIFIEDQUANTITYPRELIMINARYSPACETIMEL2ESTIMATE} and  
\eqref{E:LDERIVATIVEOFPARTIALMODIFIEDQUANTITYPRELIMINARYHYPERSURFACEL2ESTIMATELESSPRECISE} into the relevant factors 
in the integrals on RHS~\eqref{E:PROOFMAINSPACETIMEWAVEESTIMATEIBP}
and using Young's inequality in the form $ab \leq \frac{1}{f} a^2 + f b^2$ (for appropriately chosen $f$)
as well as the coerciveness estimate 
$\left\| 
	\sqrt{\upmu} \angrmD \tanderY^N \RRiemann
\right\|_{L^2\hypthreearg{\timefunction'}{[-\rightu,u]}{\muxmulevelsetvalue})} 
\leq 
\frac{1}{\sqrt{0.49}} \hypersurfacecontrolwave_N^{1/2}(\timefunction',u)
$
(see \eqref{E:COERCIVENESSOFHYPERSURFACECONTROLWAVE}).  
We clarify that the factor $\boxed{4.13}$
is a consequence of the factor $\boxed{2.89}$ 
in the first term on RHS~\eqref{E:PARTIALMODIFIEDQUANTITYPRELIMINARYSPACETIMEL2ESTIMATE}
and our assumed smallness of $\mathring{\upalpha}$.
We further clarify that the integral 
$C 
\int_{\timefunction' = \timefunction_0}^{\timefunction} 
	\initialsmall 
	\frac{1}{|\timefunction'|^{1/2}}\hypersurfacecontrolwave_N^{1/2}(\timefunction',u) 
\, \mathrm{d}\timefunction'$,
which is generated by the last terms on RHSs~\eqref{E:PARTIALMODIFIEDQUANTITYPRELIMINARYSPACETIMEL2ESTIMATE} 
and \eqref{E:PARTIALMODIFIEDQUANTITYPRELIMINARYHYPERSURFACEL2ESTIMATE}, 
is 
$\lesssim 
\int_{\timefunction' = \timefunction_0}^{\timefunction} 
	\left(\frac{\initialsmall^2}{|\timefunction'|^{1/3}} 
	+ 
	\frac{1}{|\timefunction|^{2/3}} \hypersurfacecontrolwave_N(\timefunction',u)\right) 
\, \mathrm{d} \timefunction' 
\lesssim 
\initialsmall^2 
+ 
\Errortop(\timefunction,u) \lesssim \Errortop(\timefunction,u)
$.

The estimate \eqref{E:MAINHYPERSURFACEESTIMATEIBP} 
follows from arguments similar to the ones we used in 
proving \eqref{E:MAINSPACETIMEWAVEESTIMATEIBP}, but we now rely on 
\eqref{E:PARTIALMODIFIEDQUANTITYPRELIMINARYHYPERSURFACEL2ESTIMATE} and 
\eqref{E:PARTIALMODIFIEDQUANTITYPRELIMINARYHYPERSURFACEL2ESTIMATELESSPRECISE} in place of  
\eqref{E:PARTIALMODIFIEDQUANTITYPRELIMINARYSPACETIMEL2ESTIMATE} and  
\eqref{E:LDERIVATIVEOFPARTIALMODIFIEDQUANTITYPRELIMINARYHYPERSURFACEL2ESTIMATELESSPRECISE}. 
We clarify that the last term on RHS~\eqref{E:PARTIALMODIFIEDQUANTITYPRELIMINARYHYPERSURFACEL2ESTIMATE} 
generates the error term $C \initialsmall \frac{1}{|\timefunction|^{1/2}} \hypersurfacecontrolwave_{[1,N]}^{1/2}(\timefunction,u)$, 
which, for any $\varsigma \in (0,1]$,
by Young's inequality, we can bound by 
$ \leq C \varsigma^{-1} 
	\initialsmall^2 
	\frac{1}{|\timefunction|} 
	+ 
	C \varsigma \hypersurfacecontrolwave_{[1,N]}(\timefunction,u) 
	\lesssim 
	\Errortop(\timefunction,u)
	$. 
	We omit the remaining details. 

The estimates \eqref{E:MAINSPACETIMEWAVEESTIMATEIBPPARTIAL}--\eqref{E:MAINHYPERSURFACEESTIMATEIBPPARTIAL} 
follow from applying similar arguments that rely on the bounds 
\eqref{E:LDERIVATIVEOFPARTIALMODIFIEDQUANTITYPRELIMINARYHYPERSURFACEL2ESTIMATELESSPRECISE}--\eqref{E:PARTIALMODIFIEDQUANTITYPRELIMINARYHYPERSURFACEL2ESTIMATELESSPRECISE}. 
The desired estimates are in fact much simpler to deduce since there is an overall gain in smallness
stemming from the bound
$|\muX \Psi| \lesssim \fundbootsmall$ for $\Psi \in \wavearraypartial$, which we proved in
\eqref{E:LINFINITYIMPROVEMENTAUXTRANSVERSALPDERIVATIVESPARTIALWAVEARRAYSMALL}; 
we omit the details.
\end{proof}

Before proving the main estimates of this section, we first bound the remaining
error integrals on the right-hand side of 
the integration by parts identity \eqref{E:KEYIBPIDENTIFYFORWAVEEQUATIONENERGYESTIMATES}.
The estimates are much easier to derive compared to the ones we established in 
Lemma~\ref{L:ESTIMATESFORDIFFICULTSPACETIMEERRORINTEGRALSINVOLVINGIBPWRTL}.
We split the analysis into two lemmas.
In the next lemma, we bound the error integral of the term $\ErrorIBP$ 
on RHS~\eqref{E:KEYIBPIDENTIFYFORWAVEEQUATIONENERGYESTIMATES}.

\begin{lemma}[Estimates for easy error integrals that arise during the integration by parts with respect to $\argLrough{\muxmulevelsetvalue}$]
	\label{L:IBPEASYSPACETIMEERRRORINTEGRALSESITMATES}
	Assume that $N = \Ntop$. 
	Let $\Psi \in \{ \RRiemann,\LRiemann,v^2,v^3,\Ent\}$,
	let $\partialmodquant{\tanderY^{N-1}}$ be the partially modified quantity defined by
	\eqref{E:PARTIALMODIFIEDQUANTITY},
	and let
	$\ErrorIBP[\tanderY^N \Psi;\partialmodquant{\tanderY^{N-1}}]$
	be the error term defined in
	\eqref{E:ERRORTERM1KEYIBPIDENTIFYFORWAVEEQUATIONENERGYESTIMATES}.
	Then the following estimate holds for 
	$(\timefunction,u) \in [\timefunction_0,\timefunctionboot) \times [-\rightu,\leftu]$:
	\begin{align} \label{E:IBPFIRSTEASYSPATIALERRRORINTEGRALESITMATE}
		\int_{\twoargMrough{[\timefunction_0,\timefunction),[- \rightu,u]}{\muxmulevelsetvalue}} 
			\left|
				\ErrorIBP[\tanderY^N \Psi;\partialmodquant{\tanderY^{N-1}}]
			\right|
		\, \volMRoughCoordinates 
		& 
		\lesssim
		\int_{\timefunction' = \timefunction_0}^{\timefunction}
			\frac{1}{|\timefunction'|^{1/2}} 
			\hypersurfacecontrolwave_{[1,N]}(\timefunction',u) 
		\, \mathrm{d} \timefunction'
		+
		\initialsmall^2.
	\end{align}

In particular, the error integrals 
on LHS~\eqref{E:IBPFIRSTEASYSPATIALERRRORINTEGRALESITMATE}
is of type $\Errortoparg{N}$, i.e., they satisfy the bound \eqref{E:ERRORTOPORDERWAVEESTIMATES}.
\end{lemma}

\begin{proof}
First, using the bootstrap assumptions
and \eqref{E:MINVALUEOFMUONFOLIATION}, 
we deduce the following pointwise bound
for the error term defined in \eqref{E:ERRORTERM1KEYIBPIDENTIFYFORWAVEEQUATIONENERGYESTIMATES}:
\begin{align} \label{E:POINTWISEBOUNDFORIBPFIRSTEASYSPATIALERRRORINTEGRALESITMATE}
	\left|
		\ErrorIBP[\tanderY^N \Psi;\partialmodquant{\tanderY^{N-1}}]
	\right|
	& 
	\lesssim
	\frac{1}{|\timefunction|^{1/2}}
	\left|
		\sqrt{\upmu} \tander^{N+1} \Psi
	\right|
	\left|
		\partialmodquant{\tanderY^{N-1}}
	\right|.
\end{align}
From \eqref{E:POINTWISEBOUNDFORIBPFIRSTEASYSPATIALERRRORINTEGRALESITMATE}
and the Cauchy--Schwarz inequality,
we deduce that:
\begin{align}  \label{E:FIRSTINTEGRALBOUNDFORIBPFIRSTEASYSPATIALERRRORINTEGRALESITMATE}
	&
	\int_{\twoargMrough{[\timefunction_0,\timefunction),[- \rightu,u]}{\muxmulevelsetvalue}} 
		\left|
			\ErrorIBP[\tanderY^N \Psi;\partialmodquant{\tanderY^{N-1}}]
		\right|
	\, \volMRoughCoordinates 
		\\
	& 
	\lesssim
	\int_{\timefunction' = \timefunction_0}^{\timefunction}
		\frac{1}{|\timefunction'|^{1/2}}
		\left\| 
			\sqrt{\upmu} \tander^{N+1} \Psi
		\right\|_{L^2\left(\hypthreearg{\timefunction'}{[-\rightu,u]}{\muxmulevelsetvalue}\right)}
		\left\| 
			\partialmodquant{\tanderY^{N-1}}
		\right\|_{L^2\left(\hypthreearg{\timefunction'}{[-\rightu,u]}{\muxmulevelsetvalue}\right)}
	\, \mathrm{d} \timefunction'.
\end{align}
From 
\eqref{E:COERCIVENESSOFHYPERSURFACECONTROLWAVE}, 
\eqref{E:PARTIALMODIFIEDQUANTITYPRELIMINARYHYPERSURFACEL2ESTIMATELESSPRECISE},
and Young's inequality,
we conclude that
$\mbox{RHS~\eqref{E:FIRSTINTEGRALBOUNDFORIBPFIRSTEASYSPATIALERRRORINTEGRALESITMATE}}
\lesssim 
\mbox{RHS~\eqref{E:IBPFIRSTEASYSPATIALERRRORINTEGRALESITMATE}}$
as desired.
\end{proof}

The hypersurface error integrals that we treat in the next lemma 
appear on RHS~\eqref{E:KEYIBPIDENTIFYFORWAVEEQUATIONENERGYESTIMATES}.
The integrals involve the $\Yvf{A}$ derivatives of the rough time function
and are therefore new compared to earlier works on shocks, 
such as \cites{dC2007,jSgHjLwW2016,jLjS2021}.

\begin{lemma}[Estimates for additional hypersurface error terms related to integration by parts with respect to $\argLrough{\muxmulevelsetvalue}$] 
\label{L:L2ESTIMATESFORNEWHYPERSURFACEERRORTERMS}
Assume that $N = \Ntop$, and let $\varsigma \in (0,1]$. 
Let $\Psi \in \{ \RRiemann,\LRiemann,v^2,v^3,\Ent\}$,
and let $\partialmodquant{\tanderY^{N-1}}$ be the partially modified quantity defined by
\eqref{E:PARTIALMODIFIEDQUANTITY}. 
Then the following estimates hold for 
$(\timefunction,u) \in [\timefunction_0,\timefunctionboot) \times [-\rightu,\leftu]$:
\begin{align} 
\begin{split} 	\label{E:YDERIVATIVETIMEFUNCTIONERRORINTEGRALIBP}
	\int_{\hypthreearg{\timefunction}{[-\rightu,u]}{\muxmulevelsetvalue}} 
		\left| 
			(\Yvf{A} \timefunctionarg{\muxmulevelsetvalue}) 
			(1 + 2 \upmu) \muX \Psi 
			(\argLrough{\muxmulevelsetvalue} \tanderY^N \Psi) 
			\partialmodquant{\tanderY^{N-1}}
		\right| 
	\, \volRoughHypersurface
	& 
	\lesssim 
	\fundbootsmall 
	\frac{1}{|\timefunction|^{1/2}} 
	\hypersurfacecontrolwave_N^{1/2}(\timefunction,u)  
	\int_{\timefunction' = \timefunction_0}^{\timefunction} 
		\frac{1}{|\timefunction'|^{1/2}} 
		\hypersurfacecontrolwave_{[1,N]}(\timefunction',u) 
	\, \mathrm{d} \timefunction'
		\\
	& \ \
		+ 
		\varepsilon \hypersurfacecontrolwave_N(\timefunction,u)
		+
		\initialsmall^2 \frac{1}{|\timefunction|},
\end{split} 
	\\
\int_{\hypthreearg{\timefunction_0}{[-\rightu,u]}{\muxmulevelsetvalue}} 
	\left| 
	(\Yvf{A} \timefunctionarg{\muxmulevelsetvalue}_0) 
	(1 + 2 \upmu) \muX \Psi 
	(\argLrough{\muxmulevelsetvalue} \tanderY^N \wavearray) 
	\partialmodquant{\tanderY^{N-1}}
	\right| 
\volRoughHypersurface 
\label{E:YDERIVATIVETIMEFUNCTION0ERRORINTEGRALIBP} 
& 
\lesssim 
\initialsmall^2.
\end{align}
In particular, the error integrals 
on LHSs~\eqref{E:YDERIVATIVETIMEFUNCTIONERRORINTEGRALIBP}--\eqref{E:YDERIVATIVETIMEFUNCTION0ERRORINTEGRALIBP} 
are of type $\Errortoparg{N}$, i.e., they satisfy the bound \eqref{E:ERRORTOPORDERWAVEESTIMATES}. 
\end{lemma}

\begin{proof}
We first prove \eqref{E:YDERIVATIVETIMEFUNCTIONERRORINTEGRALIBP}.
we first use
the bootstrap assumptions,
Lemma~\ref{L:NORMOFSMOOTHTORITANGENTCOMMUTATORSANDSIMPLECOMPARISON},
\eqref{E:SMALLDERIVATIVESLINFTYESTIMATESFORROUGHTIMEFUNCTIONANDDERIVATIVES},
Cor.\,\ref{C:IMPROVEAUX},
and \eqref{E:MINVALUEOFMUONFOLIATION}
to pointwise bound the integrand on LHS~\eqref{E:YDERIVATIVETIMEFUNCTIONERRORINTEGRALIBP} by 
$\leq 
\varepsilon
\frac{1}{|\timefunction|^{1/2}}
\left| \sqrt{\upmu} \Lunit \tander^N \RRiemann \right| 
\left|\partialmodquant{\tanderY^{N-1}} \right|
$.
From this estimate, 
the Cauchy--Schwarz inequality, 
\eqref{E:COERCIVENESSOFHYPERSURFACECONTROLWAVE}, 
the first inequality in \eqref{E:PARTIALMODIFIEDQUANTITYPRELIMINARYHYPERSURFACEL2ESTIMATELESSPRECISE},
and Young's inequality,
it follows that 
$
|\mbox{LHS~\eqref{E:YDERIVATIVETIMEFUNCTIONERRORINTEGRALIBP}}|
\lesssim
\initialsmall^2
\frac{1}{|\timefunction|} 
+
\fundbootsmall 
\hypersurfacecontrolwave_N(\timefunction,u)  
+
\fundbootsmall 
\frac{1}{|\timefunction|^{1/2}} 
\hypersurfacecontrolwave_N^{1/2}(\timefunction,u)  
\int_{\timefunction' = \timefunction_0}^{\timefunction} 
	\frac{1}{|\timefunction'|^{1/2}} \hypersurfacecontrolwave_{[1,N]}^{1/2}(\timefunction',u)  
\, \mathrm{d} \timefunction' 
$,
which is in turn bounded by RHS~\eqref{E:ERRORTOPORDERWAVEESTIMATES}, i.e.,
this term is of type $\Errortop(\timefunction,u)$ as desired.

To prove \eqref{E:YDERIVATIVETIMEFUNCTION0ERRORINTEGRALIBP}, we note that the integral
on LHS~\eqref{E:YDERIVATIVETIMEFUNCTION0ERRORINTEGRALIBP}  
is a data integral that, by virtue of the arguments we used to prove \eqref{E:YDERIVATIVETIMEFUNCTIONERRORINTEGRALIBP},
but now with $\timefunction_0$ in the role of $\timefunction$,
can be bounded by 
$
\lesssim
\initialsmall^2
+
\hypersurfacecontrolwave_{[1,N]}(\timefunction_0,u)  
$.
Using \eqref{E:WAVEL2CONTROLLINGINITIALLYSMALL}, we see that the RHS of the previous
expression is $\lesssim \initialsmall^2$, which is
bounded by RHS~\eqref{E:ERRORTOPORDERWAVEESTIMATES} as desired.

\end{proof}

We are now ready to combine the results of the lemmas established above
to obtain the main estimates for the top-order error integrals highlighted in \eqref{E:OVERVIEWLPSIDIFFICULTERRORINTEGAL}.

\begin{lemma}[The main estimates for difficult top-order spacetime error integrals requiring integration by parts in $\argLrough{\muxmulevelsetvalue}$] 
\label{L:MAINESTIMATESFORDIFFICULTSPACETIMEERRORINTEGRALSINVOLVINGIBPWRTL}
Assume that $N = \Ntop$.
Then the following estimates hold for 
$(\timefunction,u) \in [\timefunction_0,\timefunctionboot) \times [-\rightu,\leftu]$:
\begin{align}
\begin{split} \label{E:FINALMAINSPACETIMEWAVEESTIMATEIBP}  
	& 
	\left| 
	\int_{\twoargMrough{[\timefunction_0,\timefunction],[-\rightu,u]}{\muxmulevelsetvalue}} 
		\frac{1}{\Lunit \timefunctionarg{\muxmulevelsetvalue}} 
			(1 + 2 \upmu) 
			(\Lunit \tanderY^N \RRiemann)
			(\muX \RRiemann) 
			\tanderY^N \mytr_{\gtorus} \upchi 
	\, \volMRoughCoordinates
	\right| 
		\\
	& 
	\leq
	\boxed{4.13} 
		\int_{\timefunction' = \timefunction_0}^{\timefunction} \frac{1}{|\timefunction'|}  
		\hypersurfacecontrolwave_N (\timefunction',u) 
	\, \mathrm{d} \timefunction' 
		\\
	& \ \
	+
	\boxed{4.13}  
	\frac{1}{|\timefunction|^{1/2}} 
	\hypersurfacecontrolwave_N^{1/2}(\timefunction,u) 
	\int_{\timefunction' = \timefunction_0}^{\timefunction}
		\frac{1}{|\timefunction'|^{1/2}} 
		\hypersurfacecontrolwave_N^{1/2}(\timefunction',u) 
	\, \mathrm{d} \timefunction' 
		\\
	& \ \
	+ 
	C_* \int_{\timefunction' = \timefunction_0}^{\timefunction}  \frac{1}{|\timefunction'|}
	\hypersurfacecontrolwave^{1/2}_N (\timefunction',u) 
	\left(\hypersurfacecontrolwavepartial_N\right)^{1/2}(\timefunction',u) \, \mathrm{d} \timefunction' 
		\\
	& \ \
	+ C_* \frac{1}{|\timefunction|^{1/2}} 
	\hypersurfacecontrolwave_N^{1/2}(\timefunction,u) 
	\int_{\timefunction' = \timefunction_0}^{\timefunction} 
		\frac{1}{|\timefunction'|^{1/2}} 
		\left(\hypersurfacecontrolwavepartial_N\right)^{1/2}(\timefunction',u) 
	\, \mathrm{d} \timefunction' 
		\\
	& \ \ 
		+ \Errortoparg{N}(\timefunction,u),
	\end{split}
\end{align}
where $\Errortoparg{N}(\timefunction,u)$ satisfies \eqref{E:ERRORTOPORDERWAVEESTIMATES}.
	
Moreover, for every $\Psi \in \wavearraypartial = \{\LRiemann,v^2,v^3,\Ent\}$, 
we have the following less degenerate estimates:
 \begin{align} 	\label{E:FINALMAINSPACETIMEWAVEESTIMATEIBPPARTIAL}
	 \left| 
	\int_{\twoargMrough{[\timefunction_0,\timefunction],[-\rightu,u]}{\muxmulevelsetvalue}} 
		\frac{1}{\Lunit \timefunctionarg{\muxmulevelsetvalue}} 
			(1 + 2 \upmu) 
			(\Lunit \tanderY^N \Psi)
			(\muX \Psi) 
			\tanderY^N \mytr_{\gtorus} \upchi 
	\, \volMRoughCoordinates
	\right| 
	& 
	\lesssim \Errortoparg{N}(\timefunction,u),
 \end{align}
	where $\Errortoparg{N}(\timefunction,u)$ satisfies \eqref{E:ERRORTOPORDERWAVEESTIMATES}.

\end{lemma}

\begin{proof}
We first prove \eqref{E:FINALMAINSPACETIMEWAVEESTIMATEIBP}. 
The operator $\tanderY^N$ on LHS~\eqref{E:FINALMAINSPACETIMEWAVEESTIMATEIBP} is
of the form $\tanderY^N = \Yvf{A} \tanderY^{N-1}$ for some $A \in \lbrace 2,3 \rbrace$.
We now use \eqref{E:PARTIALMODIFIEDQUANTITY} 
to decompose the factor $\tanderY^N \mytr_{\gtorus}\upchi$ on LHS~\eqref{E:FINALMAINSPACETIMEWAVEESTIMATEIBP}
as follows:
$\tanderY^N \mytr_{\gtorus}\upchi 
= 
\Yvf{A} 
\tanderY^{N-1} \mytr_{\gtorus} \upchi 
= 
\Yvf{A} \partialmodquant{\tanderY^{N-1}} 
- 
\Yvf{A} \partialmodquantinhom{\tanderY^{N-1}}$.
We insert this decomposition into LHS~\eqref{E:FINALMAINSPACETIMEWAVEESTIMATEIBP}
and will handle each of the two integrals separately, 
starting with the one generated by the piece $\Yvf{A} \partialmodquantinhom{\tanderY^{N-1}}$,
which is easier. Specifically, 
the pointwise estimate \eqref{E:POINTWISELANDYDERIVATIVESOFTOPORDERPARTIALMODQUANTINHOM},
the estimate $|\muX \RRiemann| \lesssim 1$ 
(see \eqref{E:LINFINITYIMPROVEMENTAUXMUANDTRANSVERSALDERIVATIVES}),
and Def.\,\ref{D:HARMLESSWAVE}
imply that the product
$
(\muX \RRiemann)
\Yvf{A} \partialmodquantinhom{\tanderY^{N-1}}
$ is of type
$\HarmlessWave{N}$.
Hence, the estimate \eqref{E:ENERGYESTIMATEHARMLESSWAVEERRORTERMS}
implies that the corresponding integral
$-\int_{\twoargMrough{[\timefunction_0,\timefunction],[-\rightu,u]}{\muxmulevelsetvalue}} 
	\frac{1}{\Lunit \timefunctionarg{\muxmulevelsetvalue}} 
	\left\lbrace 
		(1 + 2 \upmu) 
		\Lunit \tanderY^N \RRiemann 
	\right\rbrace 
	\left\lbrace 
		(\muX \RRiemann)  
		\Yvf{A} \partialmodquantinhom{\tanderY^{N-1}} 
	\right\rbrace 
\, \volMRoughCoordinates
$
is of type $\Errortoparg{N}(\timefunction,u)$ as desired.

To complete the proof of \eqref{E:FINALMAINSPACETIMEWAVEESTIMATEIBP},
it remains for us to bound the following spacetime integral:
\begin{align} \label{E:PROOFOFMAINL2ESTIMATESANNOYINGSPACETIMEINTEGRALINVOLVINGIBP}
\int_{\twoargMrough{[\timefunction_0,\timefunction],[-\rightu,u]}{\muxmulevelsetvalue}} 
	\frac{1}{\Lunit \timefunctionarg{\muxmulevelsetvalue}} 
	\left\lbrace 
		(1 + 2 \upmu) \Lunit \tanderY^N \RRiemann 
	\right\rbrace 
	\left\lbrace 
		(\muX \RRiemann) 
		\Yvf{A} \partialmodquant{\tanderY^{N-1}} 
	\right\rbrace 
\, \volMRoughCoordinates.
\end{align}
To bound \eqref{E:PROOFOFMAINL2ESTIMATESANNOYINGSPACETIMEINTEGRALINVOLVINGIBP},
we first integrate by parts using the identity \eqref{E:KEYIBPIDENTIFYFORWAVEEQUATIONENERGYESTIMATES}
with $\RRiemann$ in the role of $\varphi$ and $\partialmodquant{\tanderY^{N-1}}$ in the role of $\upeta$.
The first two integrals on RHS~\eqref{E:KEYIBPIDENTIFYFORWAVEEQUATIONENERGYESTIMATES},
namely
\begin{align*}
&
\int_{\twoargMrough{[\timefunction_0,\timefunction),[- \rightu,u]}{\muxmulevelsetvalue}}  
		(1 + 2 \upmu)
		(\muX \RRiemann) 
		(\Yvf{A} \tander^N \RRiemann) 
		\argLrough{\muxmulevelsetvalue} 
		\partialmodquant{\tanderY^{N-1}} 
	\, \volMRoughCoordinates,
		\\
&
- 
\int_{\hypthreearg{\timefunction}{[- \rightu,u]}{\muxmulevelsetvalue} } 
			(1 + 2 \upmu)
			(\muX \RRiemann)
			(\Yvf{A} \tander^N \RRiemann)
			\partialmodquant{\tanderY^{N-1}} 
	\, \volRoughHypersurface,
\end{align*}
are the main ones,
and in \eqref{E:MAINSPACETIMEWAVEESTIMATEIBP}--\eqref{E:MAINHYPERSURFACEESTIMATEIBP},
we showed that they are bounded in magnitude by $\leq \mbox{RHS~\eqref{E:FINALMAINSPACETIMEWAVEESTIMATEIBP}}$ as desired.
The third and fifth integrals on RHS~\eqref{E:KEYIBPIDENTIFYFORWAVEEQUATIONENERGYESTIMATES}, namely
$ 
	\int_{\hypthreearg{\timefunction}{[- \rightu,u]}{\muxmulevelsetvalue}} 
		(\Yvf{A} \timefunctionarg{\muxmulevelsetvalue})
		(1 + 2 \upmu) 
		(\muX \RRiemann) 
		(\argLrough{\muxmulevelsetvalue} \tander^N \RRiemann) 
		 \partialmodquant{\tanderY^{N-1}} 
	\, \volRoughHypersurface 
$
and
$ 
	-
	\int_{\hypthreearg{\timefunction_0}{[- \rightu,u]}{\muxmulevelsetvalue}} 
		(\Yvf{A} \timefunctionarg{\muxmulevelsetvalue})
		(1 + 2 \upmu) 
		(\muX \RRiemann) 
		(\argLrough{\muxmulevelsetvalue} \tander^N \RRiemann) 
		 \partialmodquant{\tanderY^{N-1}} 
	\, \volRoughHypersurface 
$,
were shown to be of type $\Errortop(\timefunction,u)$
in Lemma~\ref{L:L2ESTIMATESFORNEWHYPERSURFACEERRORTERMS}.
The remaining integrals on RHS~\eqref{E:KEYIBPIDENTIFYFORWAVEEQUATIONENERGYESTIMATES},
which have integrands equal to
$\ErrorIBP[\tanderY^N \Psi;\partialmodquant{\tanderY^{N-1}}]$,
were shown to be of type $\Errortop(\timefunction,u)$
in Lemma~\ref{L:IBPEASYSPACETIMEERRRORINTEGRALSESITMATES}.

The estimate \eqref{E:FINALMAINSPACETIMEWAVEESTIMATEIBPPARTIAL}
can be proved using nearly identical arguments.
The difference compared to \eqref{E:FINALMAINSPACETIMEWAVEESTIMATEIBP} is that the integrand factor
$\muX \Psi$ on LHS~\eqref{E:FINALMAINSPACETIMEWAVEESTIMATEIBPPARTIAL} satisfies the
pointwise bound $|\muX \Psi| \lesssim \fundbootsmall$ (see \eqref{E:LINFINITYIMPROVEMENTAUXTRANSVERSALPDERIVATIVESPARTIALWAVEARRAYSMALL});
this provides a smallness factor that allows us to relegate all the terms
to the error term
$\Errortoparg{N}(\timefunction,u)$
on RHS~\eqref{E:FINALMAINSPACETIMEWAVEESTIMATEIBPPARTIAL}.

\end{proof}

\subsection{Estimates for the less-degenerate top-order eikonal function error integrals} 
\label{SS:ESTIMATESFORLESSDEGENERATETOPORDEREIKONALFUNCTIONERRORINTEGRALS}
In the next lemma, we bound error integrals that feature an additional factor of $\upmu$ compared to those 
on LHS~\eqref{E:SPACETIMEBOUNDSMOSTDIFFICULTWAVEPRODUCT}. The integrals are much easier to estimate
and lead to much less degenerate estimates compared to \eqref{E:SPACETIMEBOUNDSMOSTDIFFICULTWAVEPRODUCT}.

\begin{lemma}[Estimates for less-degenerate top-order integrals] 
\label{L:ESTIMATESFORLESSDEGENERATETOPORDEREIKONALFUNCTIONERRORINTEGRALS}
Assume that $N = \Ntop$,
$\Psi \in \wavearray = \{\RRiemann,\LRiemann, v^2,v^2,\Ent\}$,
and recall that the multiplier vectorfield $\multipliervectorfield$ is defined in \eqref{E:MULTIPLIERVECTORFIELD}.
Then for $A = 2,3$, the following estimates hold for 
$(\timefunction,u) \in [\timefunction_0,\timefunctionboot) \times [-\rightu,\leftu]$,
\begin{align} 
\begin{split} \label{E:LESSDEGENERATETOPORDERWAVEERRORTERMS}
	& 
	\left| 
		\int_{\twoargMrough{[\timefunction_0,\timefunction],[-\rightu,u]}{\muxmulevelsetvalue}} 
			\frac{1}{\Lunit \timefunctionarg{\muxmulevelsetvalue}} 
			(\multipliervectorfield \tander^N \Psi) 
			(\angrmD^{\#} \Psi) 
			\cdot 
			\upmu \angrmD \tander^{N-1} \mytr_{\gtorus} \upchi 
		\, \volMRoughCoordinates 	
	\right|, 
		\\
	& \left| 
			\int_{\twoargMrough{[\timefunction_0,\timefunction],[-\rightu,u]}{\muxmulevelsetvalue}} 
				\frac{1}{\Lunit \timefunctionarg{\muxmulevelsetvalue}} 
				(\multipliervectorfield \tander^N \Psi) (\Speed^{-2} X^A) 
				(\angrmD^{\#} \Psi) 
				\cdot 
				\upmu \angrmD \tander^{N-1} \mytr_{\gtorus} \upchi 
			\, \volMRoughCoordinates 
	\right|
		\\
	&
	\lesssim 
	\Errortop(\timefunction,u), 
\end{split}
\end{align}
where $\Errortop$ satisfies the estimate \eqref{E:ERRORTOPORDERWAVEESTIMATES}. 
\end{lemma}

\begin{proof}
We provide the proof for the first spacetime integral on LHS~\eqref{E:LESSDEGENERATETOPORDERWAVEERRORTERMS}; the proof for the second spacetime integral follows from the same arguments and the pointwise estimate 
$|\Speed^{-2} X^A| \lesssim 1$, which follows from Lemma~\ref{L:SCHEMATICSTRUCTUREOFVARIOUSTENSORSINTERMSOFCONTROLVARS} 
and the bootstrap assumptions. 
To proceed, we start by considering the terms containing an $\Lunit$-derivative in the expansion
$\multipliervectorfield \tander^N \Psi 
= 
(1 + 2 \upmu) \Lunit \tander^N \Psi 
+
2 \muX \tander^N \Psi$. 
Using Young's inequality, 
\eqref{E:ANGDFPOINTWISEBOUNDEDBYCOMMUTATORVECTORFIELDS},  
\eqref{E:CLOSEDVERSIONLUNITROUGHTTIMEFUNCTION},
\eqref{E:COERCIVENESSOFHYPERSURFACECONTROLWAVE},
and the bootstrap assumptions,
we bound the corresponding error integrals as follows:
\begin{align}
\begin{split} \label{E:PROOFLESSDEGENERATETOPORDERWAVEERRORTERMS}
	& 
	\left| 
		\int_{\twoargMrough{[\timefunction_0,\timefunction],[-\rightu,u]}{\muxmulevelsetvalue}} 
			\frac{1}{\Lunit \timefunctionarg{\muxmulevelsetvalue}} 
			(1 + 2 \upmu) 
			(\Lunit \tander^N \Psi) 
			(\angrmD^{\#} \Psi) 
			\cdot 
			\upmu \angrmD \tander^{N-1} \mytr_{\gtorus} \upchi 
		\, \volMRoughCoordinates 
	\right| 
		\\
	&  
	\lesssim 
		\int_{\twoargMrough{[\timefunction_0,\timefunction],[-\rightu,u]}{\muxmulevelsetvalue}} 
			\frac{1}{\Lunit \timefunctionarg{\muxmulevelsetvalue}}
			\left|
				\Lunit \tander^N \Psi 
			\right|^2 
		\, \volMRoughCoordinates 
		+ 
		\int_{\twoargMrough{[\timefunction_0,\timefunction],[-\rightu,u]}{\muxmulevelsetvalue}} 
			\left|\upmu \angrmD \tander^{N-1} \mytr_{\gtorus} \upchi\right|^2 
		\, \volMRoughCoordinates 
		\\
	& \lesssim 
	\int_{u' = -\rightu}^u 
		\hypersurfacecontrolwave_{[1,N]}(\timefunction,u') 
	\, \mathrm{d} u' 
	+ 
	\int_{\timefunction' = \timefunction_0}^{\timefunction} 
		\left\| 
			\upmu \tander^N \mytr_{\gtorus} \upchi 
		\right\|^2_{L^2\left(\hypthreearg{\timefunction'}{[-\rightu,u]}{\muxmulevelsetvalue}\right)} 
	\, \mathrm{d}\timefunction'.
\end{split}
\end{align}
The term 
$
\int_{u' = -\rightu}^u 
		\hypersurfacecontrolwave_{[1,N]}(\timefunction,u') 
	\, \mathrm{d} u' 
$
is manifestly bounded by RHS~\eqref{E:LESSDEGENERATETOPORDERWAVEERRORTERMS}
(see RHS~\eqref{E:ERRORTOPORDERWAVEESTIMATES}).
To bound the last integral on RHS~\eqref{E:PROOFLESSDEGENERATETOPORDERWAVEERRORTERMS}
by RHS~\eqref{E:LESSDEGENERATETOPORDERWAVEERRORTERMS},
we insert the estimate \eqref{E:TOPORDERL2ESTIMATEMUCHI} into the integrand
and use arguments, based on Young's inequality, similar to the ones we used
below \eqref{E:SECONDSTEPSPACETIMEBOUNDSMOSTDIFFICULTWAVEPRODUCT}.

To bound the error integrals on LHS~\eqref{E:LESSDEGENERATETOPORDERWAVEERRORTERMS} 
that are generated by the term $2 \muX \tander^N \Psi$ in the decomposition
$\multipliervectorfield \tander^N \Psi 
= 
(1 + 2 \upmu) \Lunit \tander^N \Psi 
+
2 \muX \tander^N \Psi$,
we first argue in the proof of \eqref{E:PROOFLESSDEGENERATETOPORDERWAVEERRORTERMS}
to deduce:
\begin{align}
\begin{split} \label{E:SECONDPROOFLESSDEGENERATETOPORDERWAVEERRORTERMS}
	& 
	\left| 
		2
		\int_{\twoargMrough{[\timefunction_0,\timefunction],[-\rightu,u]}{\muxmulevelsetvalue}} 
			\frac{1}{\Lunit \timefunctionarg{\muxmulevelsetvalue}} 
			(\muX \tander^N \Psi) 
			(\angrmD^{\#} \Psi) 
			\cdot 
			\upmu \angrmD \tander^{N-1} \mytr_{\gtorus} \upchi 
		\, \volMRoughCoordinates 
	\right| 
		\\
	&  
	\lesssim 
		\int_{\twoargMrough{[\timefunction_0,\timefunction],[-\rightu,u]}{\muxmulevelsetvalue}} 
			\left|
				\muX \tander^N \Psi 
			\right|^2 
		\, \volMRoughCoordinates 
		+ 
		\int_{\twoargMrough{[\timefunction_0,\timefunction],[-\rightu,u]}{\muxmulevelsetvalue}} 
			\left|\upmu \angrmD \tander^{N-1} \mytr_{\gtorus} \upchi\right|^2 
		\, \volMRoughCoordinates 
		\\
	& \lesssim 
	\int_{\timefunction' = \timefunction_0}^{\timefunction} 
		\hypersurfacecontrolwave_{[1,N]}(\timefunction',u) 
	\, \mathrm{d} \timefunction' 
	+ 
	\int_{\timefunction' = \timefunction_0}^{\timefunction} 
		\left\| 
			\upmu \tander^N \mytr_{\gtorus} \upchi 
		\right\|^2_{L^2\left(\hypthreearg{\timefunction'}{[-\rightu,u]}{\muxmulevelsetvalue}\right)} 
	\, \mathrm{d}\timefunction'.
\end{split}
\end{align}
Arguing as in the end of the previous paragraph, we conclude that
$
\mbox{RHS~\eqref{E:SECONDPROOFLESSDEGENERATETOPORDERWAVEERRORTERMS}}
\lesssim
\mbox{RHS~\eqref{E:LESSDEGENERATETOPORDERWAVEERRORTERMS}}
$
as desired.

\end{proof}

\subsection{Estimates for error integrals involving a loss of one derivative}
\label{SS:ESTIMATESFORWAVEERRORINTEGRALSTHATLOSEADERIVATIVE}
Prop.\,\ref{P:APRIORIL2ESTIMATESWAVEVARIABLES} states that the wave energies
become less singular by a factor of $|\timefunction|^2$
at every level of descent below top-order, until one reaches a level
where the energies are bounded. Our proof of this ``descent scheme'' relies on bounding
the difficult error integrals generated by the $\tanderY^N \mytr_{\gtorus}\upchi$-involving products on
RHSs~\eqref{E:TOPCOMMUTEDWAVELFIRSTTHENALLYS}--\eqref{E:TOPCOMMUTEDWAVEALLYS}
by controlling the $\tanderY^N \mytr_{\gtorus}\upchi$ factors
via transport equation estimates that \emph{lose one derivative};
this is very different compared to, for example,
the top-order error integrals we bounded in
Lemmas~\ref{L:BOUNDSFORMOSTDIFFICULTWAVEERRORNTEGRALS} and
\ref{L:MAINESTIMATESFORDIFFICULTSPACETIMEERRORINTEGRALSINVOLVINGIBPWRTL},
where we could not afford any derivative loss.
While losing one derivative is permissible below top-order,
this approach couples the below top-order energy estimates to the top-order ones.
The main merit of this approach is that it leads to estimates that are less singular 
with respect to powers of $\timefunction$, which ultimately allows us to implement
the energy estimate descent scheme.
In the next lemma, we prove the main estimates for 
below-top-order error integrals that lose one derivative.

\begin{lemma}[Estimates for wave equation error integrals involving a loss of one derivative]
\label{L:ESTIMATESFORWAVEERRORINTEGRALSTHATLOSEADERIVATIVE}
Assume that $2 \leq N \leq \Ntop$ and
$\Psi \in \wavearray = \{\RRiemann,\LRiemann, v^2,v^2,\Ent\}$,
and recall that $\multipliervectorfield$ is the multiplier vectorfield defined in \eqref{E:MULTIPLIERVECTORFIELD}.
Then the following estimates hold for 
$(\timefunction,u) \in [\timefunction_0,\timefunctionboot) \times [-\rightu,\leftu]$:
\begin{align}
\begin{split} \label{E:ESTIMATESFORWAVEERRORINTEGRALSTHATLOSEADERIVATIVE}
	& \int_{\twoargMrough{[\timefunction_0,\timefunction],[-\rightu,u]}{\muxmulevelsetvalue}} 
			\frac{1}{\Lunit \timefunctionarg{\muxmulevelsetvalue}} 
			\left| 
				\multipliervectorfield \tander^{N-1} \Psi 
			\right| 
			\left|
				\begin{pmatrix}
					(\muX \Psi) \tander^{N-1} \mytr_{\gtorus}\upchi 
						\\
					(\angrmD^{\#} \Psi) \cdot \upmu \angrmD \tander^{N-2} \mytr_{\gtorus} \upchi
						\\
					\Speed^{-2} X^A (\angrmD^{\#} \Psi) \cdot \upmu \angrmD \tander^{N-2} \mytr_{\gtorus}\upchi
				\end{pmatrix}
		  \right|_{\gtorus} 
		\, \volMRoughCoordinates  
			\\
	&  \lesssim 
			\int_{\timefunction' = \timefunction_0}^{\timefunction} 
				\frac{1}{|\timefunction'|^{1/2} } 
				\hypersurfacecontrolwave_{[1,N-1]}^{1/2}(\timefunction',u) 
				\int_{\timefunction'' = \timefunction_0}^{\timefunction'} 
					\frac{1}{|\timefunction''|^{1/2}} 
					\hypersurfacecontrolwave_N^{1/2}(\timefunction'',u) 
				\, \mathrm{d} \timefunction'' 
			\mathrm{d} \timefunction'   
			\\
	& \ \ 
		+ 
		\Errorsubcriticalarg{N-1},
\end{split}
\end{align}
where $\Errorsubcriticalarg{N-1}$ satisfies \eqref{E:ERRORBELOWTOPORDERWAVEESTIMATES} with $N-1$ in the role of $M$.

\end{lemma}

\begin{proof}
We prove \eqref{E:ESTIMATESFORWAVEERRORINTEGRALSTHATLOSEADERIVATIVE} only
for the term generated by the first entry $(\muX\Psi) \tander^{N-1} \mytr_{\gtorus} \upchi$ on 
LHS~\eqref{E:ESTIMATESFORWAVEERRORINTEGRALSTHATLOSEADERIVATIVE}; 
the remaining terms on the LHS are easier to estimate
because they enjoy an additional power of $\upmu$, and we omit the details. 
Since $\multipliervectorfield \tander^{N-1} \Psi = (1 + 2 \upmu) \Lunit \tander^{N-1} \Psi + 2\muX \tander^{N-1} \Psi$, 
we can use 
\eqref{E:CLOSEDVERSIONLUNITROUGHTTIMEFUNCTION},
\eqref{E:COERCIVENESSOFHYPERSURFACECONTROLWAVE} 
\eqref{E:PRELIMINARYEIKONALWITHOUTL},
the bootstrap assumptions,
Young's inequality, H\"{o}lder's inequality,
and the fact that the
$\hypersurfacecontrolwave_M(\timefunction,u)$ are increasing in their arguments 
to deduce:
\begin{align}
\begin{split} \label{E:PROOFENERGYWAVEESTIMATESBELOWTOPORDERTERMS}
	& 
	\int_{\twoargMrough{[\timefunction_0,\timefunction],[-\rightu,u]}{\muxmulevelsetvalue}} 
		\frac{1}{\Lunit \timefunctionarg{\muxmulevelsetvalue}} 
		\left| \multipliervectorfield \tander^{N-1} \Psi \right| 
		\left| (\muX \Psi) \tander^{N-1} \mytr_{\gtorus}\upchi \right| 
	\, \volMRoughCoordinates 
		\\
	& \lesssim 
		\int_{\timefunction' = \timefunction_0}^{\timefunction}
				\left\lbrace
					\left\| 
						\Lunit \tander^{N-1} \Psi
					\right\|_{L^2\left(\hypthreearg{\timefunction'}{[-\rightu,u]}{\muxmulevelsetvalue}\right)}
					+
					\left\| 
						\muX \tander^{N-1} \Psi
					\right\|_{L^2\left(\hypthreearg{\timefunction'}{[-\rightu,u]}{\muxmulevelsetvalue}\right)}
				\right\rbrace
				\left\| 
					\tander^{N-1} \mytr_{\gtorus}\upchi
				\right\|_{L^2\left(\hypthreearg{\timefunction'}{[-\rightu,u]}{\muxmulevelsetvalue}\right)}
		\, \mathrm{d} \timefunction' 
		\\
	& \lesssim 
	\int_{\timefunction' = \timefunction_0}^{\timefunction} 
		\frac{1}{|\timefunction'|^{1/2} } \hypersurfacecontrolwave_{[1,N-1]}^{1/2}(\timefunction',u)
		\left\lbrace  
			\initialsmall 
			+ 
			\int_{\timefunction'' = \timefunction_0}^{\timefunction'} 
				\frac{\hypersurfacecontrolwave_{[1,N]}^{1/2} (\timefunction,u)}{|\timefunction''|^{1/2}}  
			\, \mathrm{d} \timefunction'' 
		\right\rbrace 
	\, \mathrm{d} \timefunction' 
		\\
	&  \lesssim 
			\int_{\timefunction' = \timefunction_0}^{\timefunction} 
				\frac{1}{|\timefunction'|^{1/2} } 
				\hypersurfacecontrolwave_{[1,N-1]}^{1/2}(\timefunction',u) 
				\int_{\timefunction'' = \timefunction_0}^{\timefunction'} 
					\frac{1}{|\timefunction''|^{1/2}} \hypersurfacecontrolwave_N^{1/2}(\timefunction'',u) 
				\, \mathrm{d} \timefunction'' 
			\, \mathrm{d} \timefunction' 
				\\
	& \ \ 
		+ 
		\int_{\timefunction' = \timefunction_0}^{\timefunction} 
			\frac{1}{|\timefunction'|^{1/2}} \hypersurfacecontrolwave_{[1,N-1]}(\timefunction',u)
		\, \mathrm{d} \timefunction' 
		+ 
		\initialsmall^2 
		\int_{\timefunction' = \timefunction_0}^{\timefunction} 
			\frac{1}{|\timefunction'|^{1/2}} 
		\, \mathrm{d} \timefunction'.
\end{split}
\end{align}
Noting that the last term on RHS~\eqref{E:PROOFENERGYWAVEESTIMATESBELOWTOPORDERTERMS} is
$\lesssim \initialsmall^2$, 
we conclude that 
$
\mbox{RHS~\eqref{E:PROOFENERGYWAVEESTIMATESBELOWTOPORDERTERMS}}
\lesssim 
\mbox{RHS~\eqref{E:ESTIMATESFORWAVEERRORINTEGRALSTHATLOSEADERIVATIVE}}
$
as desired.

\end{proof}

\subsection{Proof of Prop.\,\ref{P:MAINWAVEENERGYINTEGRALINEQUALITIES}}
\label{SS:PROOFOF:MAINWAVEENERGYINTEGRALINEQUALITIES}
We now prove the main integral inequalities for the wave variables, i.e., Prop.\,\ref{P:MAINWAVEENERGYINTEGRALINEQUALITIES}. 

\medskip

\noindent \textbf{Proof of \eqref{E:TOPORDERWAVEL2CONTROLLINGINTEGRALINEQUALITY}:}
We first prove the top-order estimate \eqref{E:TOPORDERWAVEL2CONTROLLINGINTEGRALINEQUALITY}.
Let $N = \Ntop$ and $\Psi \in \{\RRiemann,\LRiemann,v^2,v^3,\Ent\}$, and let $\mathfrak{G}$ denote
$\upmu$ times the inhomogeneous term in the geometric wave equation \eqref{E:COVARIANTWAVEEQUATIONSWAVEVARIABLES} satisfied by 
$\Psi$, i.e. $\upmu \Box_{\gfour} \Psi = \mathfrak{G}$. 
Fix any $\tander^N \in \mathfrak{P}^{(N)}$, where  
$\mathfrak{P}^{(N)}$ is defined in Def.\,\ref{D:STRINGSOFCOMMUTATIONVECTORFIELDS},
and let $(\timefunction,u) \in [\timefunction_0,\timefunctionboot)\times [- \rightu,\leftu]$.
The starting point of the proof is the fundamental energy--null-flux identity 
\eqref{E:FUNDAMENTALENERGYINTEGRALIDENTITCOVARIANTWAVES} with $f \eqdef \tander^N \Psi$:
\begin{align}  
\begin{split} \label{E:MAINWAVEINTEGRALIDENTITIESINTERMEDIATESTEP1}
	& 
	\mathbb{E}[\tander^N \Psi]_{(\textnormal{Wave})}(\timefunction,u)  
	+ 
	\mathbb{F}[\tander^N \Psi]_{(\textnormal{Wave})}(\timefunction,u) 
	+ 
	\spacetimeintegralcontrolwave[\tander^N \Psi](\timefunction,u) 
		\\
	& 
	= 
	\mathbb{E}[\tander^N \Psi]_{(\textnormal{Wave})}(\timefunction_0,u) 
	+ 
	\mathbb{F}[\tander^N \Psi]_{(\textnormal{Wave})}(\timefunction,-\rightu)  
		\\
	&  \ \ 
	+ 
	\int_{\twoargMrough{[\timefunction_0,\timefunction],[-\rightu,u]}{\muxmulevelsetvalue}} 
		\frac{1}{\Lunit \timefunctionarg{\muxmulevelsetvalue}} 
		{^{(\multipliervectorfield)}\mathfrak{B}}[\tander^N \Psi] 
	\, \volMRoughCoordinates  
	\\
	& \ \ 
		- 
		\int_{\twoargMrough{[\timefunction_0,\timefunction],[-\rightu,u]}{\muxmulevelsetvalue}} \frac{1}{\Lunit \timefunctionarg{\muxmulevelsetvalue}} 
			\left\lbrace (1 + 2 \upmu)\Lunit \tander^N \Psi + 2\muX \tander^N \Psi\right\rbrace 
			\upmu \Box_{\gfour} \tander^N \Psi 
		\,  \volMRoughCoordinates.
	\end{split}
\end{align}
We will show that
$
\mbox{RHS~\eqref{E:MAINWAVEINTEGRALIDENTITIESINTERMEDIATESTEP1}} 
\leq 
\mbox{RHS~\eqref{E:TOPORDERWAVEL2CONTROLLINGINTEGRALINEQUALITY}}$,
which is the difficult step. After that, we
can take the supremum of the resulting estimate over the relevant values of $\timefunction$ and $u$, 
then take the maximum over $\Psi \in \{\RRiemann,\LRiemann,v^2,v^3,\Ent\}$ and over all $\Tanset^N \in \mathfrak{P}^{(N)}$,
and appeal to definitions~\eqref{E:WAVESPACELIKEANDNULLHYPERSURFACEL2CONTROLLINGQUANTITY}--\eqref{E:WAVETOTALL2CONTROLLINGQUANTITY},
finally concluding the desired top-order estimate \eqref{E:TOPORDERWAVEL2CONTROLLINGINTEGRALINEQUALITY}.

To complete the proof of \eqref{E:TOPORDERWAVEL2CONTROLLINGINTEGRALINEQUALITY}, 
it remains for us to show that
$
\mbox{RHS~\eqref{E:MAINWAVEINTEGRALIDENTITIESINTERMEDIATESTEP1}} 
\leq 
\mbox{RHS~\eqref{E:TOPORDERWAVEL2CONTROLLINGINTEGRALINEQUALITY}}$.
In the rest of the proof, $\Errortoparg{N}(\timefunction,u)$ denotes a term of type 
$\Errortoparg{N}(\timefunction,u)$ on RHS~\eqref{E:TOPORDERWAVEL2CONTROLLINGINTEGRALINEQUALITY},
i.e., any term that satisfies \eqref{E:ERRORTOPORDERWAVEESTIMATES}.
First, using \eqref{E:WAVEL2CONTROLLINGINITIALLYSMALL}, 
we see that the initial data energy and null-flux terms on RHS~\eqref{E:MAINWAVEINTEGRALIDENTITIESINTERMEDIATESTEP1}
satisfy 
$\mathbb{E}[\tander^N \Psi]_{(\textnormal{Wave})}(\timefunction_0,u) 
+ 
\mathbb{F}[\tander^N \Psi]_{(\textnormal{Wave})}(\timefunction,-\rightu) 
\lesssim \initialsmall^2 
= 
\Errortoparg{N}(\timefunction,u)$ as desired.  

Next, using Lemma~\ref{L:ESTIMATESFORERRORTERMSGENERATEDBYMULTIPLIERVECTORFIELD}, 
we see that the
${^{(\multipliervectorfield)}\mathfrak{B}}[\tander^N \Psi]$-involving integral
on RHS~\eqref{E:MAINWAVEINTEGRALIDENTITIESINTERMEDIATESTEP1} 
is type $\Errorsubcriticalarg{N}(\timefunction,u)$
(and hence of type $\Errortoparg{N}(\timefunction,u)$) as desired.  

\medskip
It remains for us to bound the last integral 
$
\int_{\twoargMrough{[\timefunction_0,\timefunction],[-\rightu,u]}{\muxmulevelsetvalue}} \frac{1}{\Lunit \timefunctionarg{\muxmulevelsetvalue}} 
			\left\lbrace (1 + 2 \upmu)\Lunit \tander^N \Psi + 2\muX \tander^N \Psi\right\rbrace 
			\upmu \Box_{\gfour} \tander^N \Psi 
		\,  \volMRoughCoordinates 
$
on RHS~\eqref{E:MAINWAVEINTEGRALIDENTITIESINTERMEDIATESTEP1}. We split the argument into Steps 1, 2, 3A, and 3B.
We stress that
\emph{none of the error integrals we have treated thus far
and none of the error integrals that we treat in Steps 1, 2, or 3B
generate the boxed-constant-multiplied integrals or $C_*$-multiplied integrals on
$RHS~\eqref{E:TOPORDERWAVEL2CONTROLLINGINTEGRALINEQUALITY}$; they
are generated only in Step 3A.}

\medskip

\noindent \emph{Step 1: the case $\tander^N \notin \left\lbrace \tanderY^{N-1} \Lunit, \, \tanderY^N \right\rbrace$}.
If $\tander^N$ is any string of tangential commutation vectorfields other than $\tanderY^{N-1} \Lunit$ or $\tanderY^N$, 
then by substituting \eqref{E:TOPCOMMUTEDWAVENOTDIFFICULT} 
for $\upmu \Box_{\gfour} \tander^N \Psi$ on 
RHS~\eqref{E:MAINWAVEINTEGRALIDENTITIESINTERMEDIATESTEP1}
and using Lemma~\ref{L:HARMLESSWAVETERMSERRORINTEGRALBOUNDS} 
to handle the $\HarmlessWave{\Ntop}$ terms
and Lemma~\ref{L:ESTIMATESFORERRORINTEGRALSGENERATEDBYTHEINHOMOGENEOUSTERMS} 
to handle the terms generated by $\mathfrak{G}$,
we find that the corresponding error integrals
are of type $\Errortoparg{N}(\timefunction,u)$ 
In total, we have shown that, except for the cases of 
$\tander^N = \tanderY^{N-1} \Lunit$ or $\tander^N = \tanderY^N$,
all error integrals on RHS~\eqref{E:MAINWAVEINTEGRALIDENTITIESINTERMEDIATESTEP1} 
are of type $\Errortoparg{N}(\timefunction,u)$.
In particular,
\emph{none of the integrals handled thus far generate the boxed-constant-multiplied integrals or $C_*$-multiplied integrals on
$RHS~\eqref{E:TOPORDERWAVEL2CONTROLLINGINTEGRALINEQUALITY}$.}

\medskip

\noindent \emph{Step 2: the case $\tander^N = \tanderY^{N-1} \Lunit$}.
We now consider the last integral on RHS~\eqref{E:MAINWAVEINTEGRALIDENTITIESINTERMEDIATESTEP1}
in the case $\tander^N = \tanderY^{N-1} \Lunit$.
We use \eqref{E:TOPCOMMUTEDWAVELFIRSTTHENALLYS} to substitute for $\upmu \Box_{\gfour} \tander^N \Psi$ on 
RHS~\eqref{E:MAINWAVEINTEGRALIDENTITIESINTERMEDIATESTEP1}. All error integrals
except for the one generated by the first product on RHS~\eqref{E:TOPCOMMUTEDWAVELFIRSTTHENALLYS} 
can be handled using the same arguments given in
Step 1.
The error integral generated by the first product on RHS~\eqref{E:TOPCOMMUTEDWAVELFIRSTTHENALLYS}
is: 
\begin{align} \label{E:EASIERTOPORDERERRORINTEGRALLFIRSTPROOFOFWAVEL2ESTIMATES}
	&
	\int_{\twoargMrough{[\timefunction_0,\timefunction],[-\rightu,u]}{\muxmulevelsetvalue}} 
		\frac{1}{\Lunit \timefunctionarg{\muxmulevelsetvalue}} 
		\left\lbrace
			(1 + 2 \upmu) \Lunit \tanderY^N \Psi 
			+ 
			2 \muX \tanderY^N \Psi 
		\right\rbrace 
		(\angrmD^{\sharp} \Psi) 
		\cdot \upmu \angrmD \tanderY^{N-1} \mytr_{\gtorus}\upchi 
	\, \volMRoughCoordinates,
\end{align}
and the estimate \eqref{E:ESTIMATESFOREASYTOPORDEREIKONALFUNCTIONERRORINTEGRALS}
implies that the integral is also of type $\Errortoparg{N}(\timefunction,u)$.

\medskip

\noindent \emph{Step 3A: the case $\tander^N = \tanderY^N$ and $\Psi = \RRiemann$}.
We now consider the most difficult case, $\tander^N = \tanderY^N$ and
$\Psi = \RRiemann$. We substitute RHS~\eqref{E:TOPCOMMUTEDWAVEALLYS} 
for $\upmu \Box_{\gfour} \tanderY^N \Psi$ on RHS~\eqref{E:MAINWAVEINTEGRALIDENTITIESINTERMEDIATESTEP1}.
All error integrals
except for the ones generated by the first two products on
RHS~\eqref{E:TOPCOMMUTEDWAVEALLYS} can be handled using the same arguments given in
Step 1.
The error integral generated by the second product on
RHS~\eqref{E:TOPCOMMUTEDWAVELFIRSTTHENALLYS} is: 
\begin{align} \label{E:EASIERTOPORDERERRORINTEGRALALLYSPROOFOFWAVEL2ESTIMATES}
	&
	\int_{\twoargMrough{[\timefunction_0,\timefunction],[-\rightu,u]}{\muxmulevelsetvalue}} 
		\frac{1}{\Lunit \timefunctionarg{\muxmulevelsetvalue}} 
		\left\lbrace
			(1 + 2 \upmu) \Lunit \tanderY^N \RRiemann 
			+ 
			2 \muX \tanderY^N \RRiemann 
		\right\rbrace 
		(\Speed^{-2} X^A) 
		(\angrmD^{\sharp} \RRiemann) 
		\cdot 
		\upmu 
		\angrmD \tanderY^{N-1} \mytr_{\gtorus} \upchi
	\, \volMRoughCoordinates,
\end{align} 
and the estimate \eqref{E:ESTIMATESFOREASYTOPORDEREIKONALFUNCTIONERRORINTEGRALS}
implies that the integral is also of type $\Errortoparg{N}(\timefunction,u)$.

The first product on
RHS~\eqref{E:TOPCOMMUTEDWAVELFIRSTTHENALLYS} generates the following two difficult error integrals:
\begin{align} \label{E:DIFFICULTTOPORDERERRORINTEGRAL1ALLYSPROOFOFWAVEL2ESTIMATES}
	&
	\int_{\twoargMrough{[\timefunction_0,\timefunction],[-\rightu,u]}{\muxmulevelsetvalue}} 
		\frac{1}{\Lunit \timefunctionarg{\muxmulevelsetvalue}} 
		\left\lbrace
			(1 + 2 \upmu) \Lunit \tanderY^N \RRiemann 
		\right\rbrace 
		(\muX \RRiemann) 
		\tanderY^{N-1} \Yvf{A} \mytr_{\gtorus} \upchi
	\, \volMRoughCoordinates,
			\\
		& 2
		\int_{\twoargMrough{[\timefunction_0,\timefunction],[-\rightu,u]}{\muxmulevelsetvalue}} 
		\frac{1}{\Lunit \timefunctionarg{\muxmulevelsetvalue}} 
		\left\lbrace
			 \muX \tanderY^N \RRiemann 
		\right\rbrace 
		(\muX \RRiemann) 
		\tanderY^{N-1} \Yvf{A} \mytr_{\gtorus} \upchi
	\, \volMRoughCoordinates,
	\label{E:DIFFICULTTOPORDERERRORINTEGRAL2ALLYSPROOFOFWAVEL2ESTIMATES}
\end{align} 
and in
\eqref{E:FINALMAINSPACETIMEWAVEESTIMATEIBP}
and
\eqref{E:SPACETIMEBOUNDSMOSTDIFFICULTWAVEPRODUCT}
respectively, we showed
that the integrals
\eqref{E:DIFFICULTTOPORDERERRORINTEGRAL1ALLYSPROOFOFWAVEL2ESTIMATES}--\eqref{E:DIFFICULTTOPORDERERRORINTEGRAL2ALLYSPROOFOFWAVEL2ESTIMATES}
are bounded in magnitude by $\leq$ 
RHS~\eqref{E:TOPORDERWAVEL2CONTROLLINGINTEGRALINEQUALITY} as desired.
\emph{It is precisely this step that
generates all the boxed-constant-multiplied integrals and $C_*$-multiplied integrals on
$RHS~\eqref{E:TOPORDERWAVEL2CONTROLLINGINTEGRALINEQUALITY}$.}

\medskip

\noindent \emph{Step 3B: the case $\tander^N = \tanderY^N$ and $\Psi \in \{\LRiemann,v^2,v^3, \Ent\}$}.
This case can be handled as in Step 3A, but the analogs of the error integrals
\eqref{E:DIFFICULTTOPORDERERRORINTEGRAL1ALLYSPROOFOFWAVEL2ESTIMATES}--\eqref{E:DIFFICULTTOPORDERERRORINTEGRAL2ALLYSPROOFOFWAVEL2ESTIMATES},
specifically the following integrals:
\begin{align} \label{E:LESSDEGENERATETOPORDERERRORINTEGRAL1ALLYSPROOFOFWAVEL2ESTIMATES}
	&
	\int_{\twoargMrough{[\timefunction_0,\timefunction],[-\rightu,u]}{\muxmulevelsetvalue}} 
		\frac{1}{\Lunit \timefunctionarg{\muxmulevelsetvalue}} 
		\left\lbrace
			(1 + 2 \upmu) \Lunit \tanderY^N \Psi
		\right\rbrace 
		(\muX \Psi) 
		\tanderY^{N-1} \Yvf{A} \mytr_{\gtorus} \upchi
	\, \volMRoughCoordinates,
			\\
		& 2
		\int_{\twoargMrough{[\timefunction_0,\timefunction],[-\rightu,u]}{\muxmulevelsetvalue}} 
		\frac{1}{\Lunit \timefunctionarg{\muxmulevelsetvalue}} 
		\left\lbrace
			 \muX \tanderY^N \Psi
		\right\rbrace 
		(\muX \Psi) 
		\tanderY^{N-1} \Yvf{A} \mytr_{\gtorus} \upchi
	\, \volMRoughCoordinates,
	\label{E:LESSDEGENERATETOPORDERERRORINTEGRAL2ALLYSPROOFOFWAVEL2ESTIMATES}
\end{align} 
where $\Psi \in \{\LRiemann,v^2,v^3, \Ent\}$,
can be bounded in magnitude via the less degenerate estimates
\eqref{E:FINALMAINSPACETIMEWAVEESTIMATEIBPPARTIAL}
and
\eqref{E:SPACETIMEBOUNDSMOSTDIFFICULTWAVEPRODUCTPARTIAL}.
In particular, all error integrals
that we encounter in Step 3B are of type $\Errortoparg{N}(\timefunction,u)$.

We have therefore proved \eqref{E:TOPORDERWAVEL2CONTROLLINGINTEGRALINEQUALITY}.

\medskip

\noindent \textbf{Proof of \eqref{E:MAINWAVEPARTIALINTEGRALINEQUALITIES}:}
We now prove the top-order estimate \eqref{E:MAINWAVEPARTIALINTEGRALINEQUALITIES}.
The proof mirrors the proof of \eqref{E:TOPORDERWAVEL2CONTROLLINGINTEGRALINEQUALITY},
except that, in view of definitions
\eqref{E:WAVEPARTIALSPACELIKEANDNULLHYPERSURFACEL2CONTROLLINGQUANTITY}--\eqref{E:WAVEPARTIALL2CONTROLLINGQUANTITY},
we do not have to derive energy estimates for
$\RRiemann$. Consequently, the proof of \eqref{E:MAINWAVEPARTIALINTEGRALINEQUALITIES} does not involve
the difficult error integrals 
\eqref{E:DIFFICULTTOPORDERERRORINTEGRAL1ALLYSPROOFOFWAVEL2ESTIMATES}--\eqref{E:DIFFICULTTOPORDERERRORINTEGRAL2ALLYSPROOFOFWAVEL2ESTIMATES}
from Step 3A, which are the only error integrals that generate the difficult
boxed-constant-involving terms. This explains why there are no such 
boxed-constant-involving terms on RHS~\eqref{E:MAINWAVEPARTIALINTEGRALINEQUALITIES}.
In particular, all error integrals 
that we encounter in the proof of \eqref{E:MAINWAVEPARTIALINTEGRALINEQUALITIES} are of type $\Errortoparg{N}(\timefunction,u)$.

\medskip

\noindent \textbf{Proof of \eqref{E:MAINWAVEBELOWTOPINTEGRALINEQUALITIES}:}
We now prove the below-top-order estimates \eqref{E:MAINWAVEBELOWTOPINTEGRALINEQUALITIES}.
We fix any $\Psi \in \{\RRiemann,\LRiemann,v^2,v^3,\Ent\}$, 
and we consider the integral identity \eqref{E:MAINWAVEINTEGRALIDENTITIESINTERMEDIATESTEP1}
with $N'$ in the role of $N$, i.e.,
\begin{align}  
\begin{split} \label{E:MAINWAVEINTEGRALIDENTITIESBELOWTOPORDERINTERMEDIATESTEP1}
	& 
	\mathbb{E}[\tander^{N'} \Psi]_{(\textnormal{Wave})}(\timefunction,u)  
	+ 
	\mathbb{F}[\tander^{N'} \Psi]_{(\textnormal{Wave})}(\timefunction,u) 
	+ 
	\spacetimeintegralcontrolwave[\tander^{N'} \Psi](\timefunction,u) 
		\\
	& 
	= 
	\mathbb{E}[\tander^{N'} \Psi]_{(\textnormal{Wave})}(\timefunction_0,u) 
	+ 
	\mathbb{F}[\tander^{N'} \Psi]_{(\textnormal{Wave})}(\timefunction,-\rightu)  
		\\
	&  \ \ 
	+ 
	\int_{\twoargMrough{[\timefunction_0,\timefunction],[-\rightu,u]}{\muxmulevelsetvalue}} 
		\frac{1}{\Lunit \timefunctionarg{\muxmulevelsetvalue}} 
		{^{(\multipliervectorfield)}\mathfrak{B}}[\tander^{N'} \Psi] 
	\, \volMRoughCoordinates  
	\\
	& \ \ 
		- 
		\int_{\twoargMrough{[\timefunction_0,\timefunction],[-\rightu,u]}{\muxmulevelsetvalue}} \frac{1}{\Lunit \timefunctionarg{\muxmulevelsetvalue}} 
			\left\lbrace (1 + 2 \upmu)\Lunit \tander^{N'} \Psi + 2\muX \tander^{N'} \Psi\right\rbrace 
			\upmu \Box_{\gfour} \tander^{N'} \Psi 
		\,  \volMRoughCoordinates.
	\end{split}
\end{align}
We will show that if $2 \leq N \leq \Ntop$ and
$1 \leq N' \leq N-1$, 
then 
$
\mbox{RHS~\eqref{E:MAINWAVEINTEGRALIDENTITIESBELOWTOPORDERINTERMEDIATESTEP1}} 
\leq 
\mbox{RHS~\eqref{E:MAINWAVEBELOWTOPINTEGRALINEQUALITIES}}$. 
Then for the same reasons given just below \eqref{E:MAINWAVEINTEGRALIDENTITIESINTERMEDIATESTEP1},
we see, taking into account the Def.\,\ref{D:SUMMEDL2CONTROLLINGQUANTITIES}
of $\totalcontrolwave_{[1,N-1]} (\timefunction,u)$,
that this implies the desired estimate \eqref{E:MAINWAVEBELOWTOPINTEGRALINEQUALITIES}.
The analogs of all the estimates through Step 1 above can be carried out
just as before; the arguments we have given show
that since $1 \leq N' \leq N-1$,
all the corresponding error integrals are of type $\Errorsubcriticalarg{N-1}$,
i.e., that they are $\lesssim \mbox{RHS~\eqref{E:ERRORBELOWTOPORDERWAVEESTIMATES}}$
with $N-1$ in the role of $M$ in \eqref{E:ERRORBELOWTOPORDERWAVEESTIMATES}.
The big difference occurs in Steps 2 and 3, where we now use a different
method to control the error integrals generated by the terms on 
RHSs~\eqref{E:TOPCOMMUTEDWAVELFIRSTTHENALLYS}--\eqref{E:TOPCOMMUTEDWAVEALLYS}
that explicitly depend on $\mytr_{\gtorus} \upchi$.
More precisely, we control them all using the derivative-losing estimate
\eqref{E:ESTIMATESFORWAVEERRORINTEGRALSTHATLOSEADERIVATIVE},
which shows that they are bounded in magnitude by:
\begin{align} \label{E:LOSSOFONEDERIVATIVEINTEGRALSINPROOFOFMAINWAVEBELOWTOPINTEGRALINEQUALITIES}
& \lesssim 
			\int_{\timefunction' = \timefunction_0}^{\timefunction} 
				\frac{1}{|\timefunction'|^{1/2} } 
				\hypersurfacecontrolwave_{[1,N']}^{1/2}(\timefunction',u) 
				\int_{\timefunction'' = \timefunction_0}^{\timefunction'} 
					\frac{1}{|\timefunction''|^{1/2}} 
					\hypersurfacecontrolwave_{N'+1}^{1/2}(\timefunction'',u) 
				\, \mathrm{d} \timefunction'' 
			\, \mathrm{d} \timefunction'   
		+ 
		\Errorsubcriticalarg{N'}.
\end{align}
Since $1 \leq N' \leq N-1$, we have $\Errorsubcriticalarg{N'} \lesssim \Errorsubcriticalarg{N-1}$,
which is $\leq \mbox{RHS~\eqref{E:MAINWAVEBELOWTOPINTEGRALINEQUALITIES}}$ as desired.
We handle the remaining double integral term in 
\eqref{E:LOSSOFONEDERIVATIVEINTEGRALSINPROOFOFMAINWAVEBELOWTOPINTEGRALINEQUALITIES}
by splitting it into two cases, the first being $N' = N-1$.
Then the double integral is bounded by the first term on RHS~\eqref{E:MAINWAVEBELOWTOPINTEGRALINEQUALITIES}.
In the second case, which is $1 \leq N' \leq N-2$, using the fact that the
$\hypersurfacecontrolwave_M(\timefunction,u)$ are increasing in their arguments,
we see that
the double time integral in \eqref{E:LOSSOFONEDERIVATIVEINTEGRALSINPROOFOFMAINWAVEBELOWTOPINTEGRALINEQUALITIES}
is
$
\lesssim 
\int_{\timefunction' = \timefunction_0}^{\timefunction} 
		\frac{1}{|\timefunction'|^{1/2} } 
		\hypersurfacecontrolwave_{[1,N-2]}^{1/2}(\timefunction',u) 
		\int_{\timefunction'' = \timefunction_0}^{\timefunction'} 
			\frac{1}{|\timefunction''|^{1/2}} 
			\hypersurfacecontrolwave_{[1,N-1]}^{1/2}(\timefunction'',u) 
		\, \mathrm{d} \timefunction'' 
\, \mathrm{d} \timefunction'
\lesssim
		\int_{\timefunction' = \timefunction_0}^{\timefunction} 
			\frac{1}{|\timefunction'|^{1/2} } 
			\hypersurfacecontrolwave_{[1,N-1]}(\timefunction',u) 
		\, \mathrm{d} \timefunction'
$,
which in turn is of type $\Errorsubcriticalarg{N-1}$, i.e., it is
$\lesssim \mbox{RHS~\eqref{E:ERRORBELOWTOPORDERWAVEESTIMATES}}$
with $N-1$ in the role of $M$ in \eqref{E:ERRORBELOWTOPORDERWAVEESTIMATES}, 
as desired.
We have therefore proved \eqref{E:MAINWAVEBELOWTOPINTEGRALINEQUALITIES}, 
which finishes the proof of the proposition.

\hfill $\qed$

\subsubsection{The proof of Prop.\,\ref{P:APRIORIL2ESTIMATESWAVEVARIABLES}}
\label{SSS:PROOFOFAPRIORIL2ESTIMATESWAVEVARIABLES}
In this section, we prove Prop.\,\ref{P:APRIORIL2ESTIMATESWAVEVARIABLES}, 
which provides the main a priori estimates for the wave variables.
Throughout the proof, we will refer to the $L^2$-controlling
quantities defined in Defs.\,\ref{D:MAINCOERCIVE} and \ref{D:SUMMEDL2CONTROLLINGQUANTITIES}.
Our proof relies on the a priori estimates for the transport variables 
from Prop.\,\ref{P:MAINHYPERSURFACEENERGYESTIMATESFORTRANSPORTVARIABLES},
which we already proved in
Sects.\,\ref{SS:PROOFOFBELOWTOPORDERTRANSPORTENERGYESTIMATES}
and \ref{SS:PROOFOFMAINVORTVORTDIVGRADENTTOPORDERBLOWUP}.
To help the reader follow the global structure of the paper,
we recall that our proof of Prop.\,\ref{P:MAINHYPERSURFACEENERGYESTIMATESFORTRANSPORTVARIABLES}
relied on the wave energy bootstrap assumptions \eqref{E:MAINWAVEENERGYBOOTSTRAPBLOWUP}--\eqref{E:MAINWAVEENERGYBOOTSTRAPREGULAR},
and that the conclusions of Prop.\,\ref{P:APRIORIL2ESTIMATESWAVEVARIABLES}
yield strict improvements of those bootstrap assumptions.

\medskip

\noindent \textbf{Estimates for $ \totalcontrolwave_{[1,\Ntop]}, \, \totalcontrolwavepartial_{[1,\Ntop]}$, and $\totalcontrolwave_{[1,\Ntop-1]}$.} 

\noindent \emph{The setup}: 
The proof of the a priori estimates for $\totalcontrolwave_{\Ntop}$ 
is coupled to the ones for $\totalcontrolwavepartial_{\Ntop}$ and $\totalcontrolwave_{[1,\Ntop-1]}$. 
Hence, we prove the a priori estimates for all three energies simultaneously via a coupled Gr\"{o}nwall argument. 
To start, we set:
\begin{subequations} 
\begin{align}
	\apriorimain(\timefunction,u) 
	& 
	\eqdef \sup_{(\timefunction',u') \in [\timefunction_0,\timefunction]\times[-\rightu,u]} 
		\iota_{\apriorimain}^{-1} (\timefunction',u') 
		\totalcontrolwave_{[1,\Ntop]}(\timefunction',u') 
			\label{E:RENORMALIZEDTOPWAVEENERGY} 
			\\
	\aprioripartial(\timefunction,u) 
	& 
	\eqdef 
	\sup_{(\timefunction',u') \in [\timefunction_0,\timefunction]\times[-\rightu,u]} 
		\iota_{\aprioripartial}^{-1} (\timefunction',u') 
		\totalcontrolwavepartial_{[1,\Ntop]}(\timefunction',u'), 
		 \label{E:RENORMALIZEDPARTIALWAVEENERGY} 
		\\
	\apriorilower(\timefunction,u) 
	& 
	\eqdef 
	\sup_{(\timefunction',u') \in [\timefunction_0,\timefunction]\times[-\rightu,u]} 
		\iota_{\apriorilower}^{-1} (\timefunction',u') \totalcontrolwave_{[1,\Ntop-1]}(\timefunction',u'), 
		\label{E:RENORMALIZEDBELOWTOPWAVEENERGY}
\end{align}
\end{subequations}
where:
\begin{subequations}
\begin{align}
	\mathfrak{l}(\timefunction) 
	& 
	\eqdef 
	\exp 
	\left(
		\int_{\timefunction' = \timefunction_0}^{\timefunction}
			\frac{1}{|\timefunction'|^{9/10}} 
		\, \mathrm{d} \timefunction'
	\right), 
		\label{E:TIMEONLYGRONWALLMULTIPLICATIONFACTORTOPFULLANDPARTIAL} 
			\\
	\iota_{\apriorimain}(\timefunction,u) 
	= 
	\iota_{\aprioripartial} (\timefunction,u) 
	& \eqdef 
	|\timefunction|^{-15.6} \mathfrak{l}^{\mathfrak{c}}(\timefunction) e^{\mathfrak{c} u}, 
		\label{E:GRONWALLMULTIPLICATIONFACTORTOPFULLANDPARTIAL}
			\\
	\iota_{\apriorilower}(\timefunction,u) 
		& \eqdef |\timefunction|^{-13.6}  \mathfrak{l}^{\mathfrak{c}}(\timefunction) e^{\mathfrak{c} u}, 
			\label{E:GRONWALLMULTIPLICATIONFACTORBELOWTOP}
\end{align}
\end{subequations}
and $\mathfrak{c}$ is a sufficiently large positive constant that we choose below. 
For future use, 
we note that when $\mathfrak{c}$ is fixed, 
the functions $\mathfrak{l}^{\mathfrak{c}}(\timefunction) $ and $e^{\mathfrak{c}u}$ are uniformly bounded from above by 
a positive ($\mathfrak{c}$-dependent) constant
for $(\timefunction,u) \in [\timefunction_0,\timefunctionboot) \times [- \rightu,\leftu]$. 
We will also silently use the basic facts that the functions $\timefunction \rightarrow \mathfrak{l}^{\mathfrak{c}}(\timefunction)$ 
and $u \rightarrow e^{\mathfrak{c} u}$ are increasing. 
Finally, we will also use the following estimates, 
whose straightforward proofs we omit: 
\begin{align}
	\int_{\timefunction' = \timefunction_0}^{\timefunction}
		\frac{\mathfrak{l}^{\mathfrak{c}}(\timefunction')}{|\timefunction'|^{9/10}} 
	\, \mathrm{d} \timefunction' 
	& 
	\leq 
	\frac{1}{\mathfrak{c}} 
	\mathfrak{l}^{\mathfrak{c}}(\timefunction), 
	&
	\int_{u' = -\rightu}^u  
		e^{\mathfrak{c} u'} 
	\,\mathrm{d} u' 
	& 
	\leq \frac{1}{\mathfrak{c}} e^{\mathfrak{c} u}. 
		\label{E:INTEGRALOFINTEGRATINGFACTORS}
\end{align}

From the above discussion, it follows that
the top-order estimate 
(i.e., \eqref{E:MAINWAVEENERGYESTIMATESBLOWUP} with $K \eqdef 0$)
and the just-below-top-order estimate
(i.e., \eqref{E:MAINWAVEENERGYESTIMATESBLOWUP} with $K \eqdef 1$)
follow once we prove that there is 
a uniform $C > 0$ (independent of all $\mathfrak{c} \geq 1$ and all sufficiently small $\varsigma \in (0,1]$)
and a $\mathfrak{c} \gg 1$
such that the following estimates hold for $(\timefunction,u) \in [\timefunction_0,\timefunctionboot) \times [- \rightu,\leftu]$:
\begin{align} \label{E:MAINWAVEESTIMATESINTERMEDIATE1}
\apriorimain(\timefunction,u) 
& 
\leq 
	C 
	\left(1 + \varsigma^{-1} \right) 
	\initialsmall^2, 
& \aprioripartial(\timefunction,u) 
&
\leq C
	\left(1 + \varsigma^{-1} \right) 
	\initialsmall^2, 
&
\apriorilower(\timefunction,u) 
& 
\leq C 
\left(1 + \varsigma^{-1} \right) 
\initialsmall^2.
\end{align}
We clarify that even though the constants
$C$ in 
\eqref{E:MAINWAVEESTIMATESINTERMEDIATE1}
are independent of $\varsigma$
and $\mathfrak{c}$,
the constants
on RHS~\eqref{E:MAINWAVEENERGYESTIMATESBLOWUP} can
depend on $\varsigma$ and, in view of definitions~\eqref{E:RENORMALIZEDTOPWAVEENERGY}--\eqref{E:RENORMALIZEDBELOWTOPWAVEENERGY}
and \eqref{E:INTEGRALOFINTEGRATINGFACTORS},
on $\mathfrak{c}$ as well.
We further clarify that our final choice of $\varsigma$ and $\mathfrak{c}$ will not be made until the
very end of the proof. The reason is that later on, we will use a downward induction scheme
to obtain the lower order estimates, and that scheme could in principle
require choosing $\varsigma$ to be smaller and $\mathfrak{c}$ to be larger at each step.
Moreover, for convenience, in the proof, we will set $\mathfrak{c} \eqdef \varsigma^{-2}$,
so that our final choice of $\mathfrak{c}$ will in fact be determined by choosing
and fixing $\varsigma \in (0,1]$ to be sufficiently small, where the final choice of
$\varsigma$ will not be made until the very end of the proof.

To prove \eqref{E:MAINWAVEESTIMATESINTERMEDIATE1}, 
we will show that there is a uniform $C > 0$
such that for every
$\varsigma \in (0,1]$,
$\mathfrak{c} \geq 1$,
sufficiently small $\fundbootsmall \geq 0$,
and
$(\timefunction,u) \in [\timefunction_0,\timefunctionboot) \times [- \rightu,\leftu]$,
the following estimates hold:
\begin{align}
	\begin{split} \label{E:MAINWAVEESTIMATESINTERMEDIATEF} 
	\apriorimain(\timefunction,u) 
	& 
	\leq C \left(1 + \varsigma^{-1}\right) \initialsmall^2 
		\\
	& \ \ 
		+ 
		\left\lbrace 
			\frac{\frac{4 \times 1.01}{1.99} + 4.13}{15.6} 
			+ 
			\frac{\frac{8(1.01)^2}{1.99}}{15.6 \times 7.8} 
			+ 
			\frac{4.13}{7.3} 
			+ 
			C \fundbootsmall 
			+ 
			C \varsigma 
			+ 
			\frac{C}{\mathfrak{c}} 
			\left(1 + \varsigma^{-1}\right) 
		\right\rbrace 
		\apriorimain(\timefunction,u) 
			\\
	& \ \ 
		+ 
		C 
		\left\lbrace 
			\fundbootsmall 
			+ 
			\varsigma 
			+ 
			\frac{1}{\mathfrak{c}} 
			\left(1 + \varsigma^{-1}\right) 
		\right\rbrace 
		\apriorilower(\timefunction,u)
		+
		C \apriorimain^{1/2}(\timefunction,u) \aprioripartial^{1/2}(\timefunction,u), 
	\end{split}	
			\\
	\begin{split} \label{E:MAINWAVEESTIMATESINTERMEDIATEG} 
	\aprioripartial(\timefunction,u) 
	& \leq 
		C \left(1 + \varsigma^{-1}\right) 
		\initialsmall^2 
		+ 
		C \left\lbrace 
			\fundbootsmall 
			+ 
			\varsigma 
			+ 
			\frac{1}{\mathfrak{c}} \left(1 + \varsigma^{-1}\right) 
		\right\rbrace 
		\apriorimain(\timefunction,u)  
		\\
	& \ \ 
	+ 
	C 
	\left\lbrace 
		\fundbootsmall 
		+ 
		\varsigma 
		+ 
		\frac{1}{\mathfrak{c}} \left(1 + \varsigma^{-1}\right) 
	\right\rbrace 
	\apriorilower(\timefunction,u),
	\end{split} 
		\\
\apriorilower(\timefunction,u) 
& \leq 
C 
\initialsmall^2 
+ 
C 
\apriorimain(\timefunction,u) 
+ 
\left\lbrace 
	\frac{1}{2} 
	+ 
	C \varsigma 
	+ 
	\frac{C}{\mathfrak{c}} 
	\left(1 + \varsigma^{-1}\right)
\right\rbrace
\apriorilower(\timefunction,u).
\label{E:MAINWAVEESTIMATESINTERMEDIATEH}
\end{align} 
Before proving \eqref{E:MAINWAVEESTIMATESINTERMEDIATEF}--\eqref{E:MAINWAVEESTIMATESINTERMEDIATEH},
we first show that these estimates imply \eqref{E:MAINWAVEESTIMATESINTERMEDIATE1}. 
To see this, we set (for convenience) $\mathfrak{c} \eqdef \varsigma^{-2}$. 
Also noting that
$
\frac{\frac{4 \times 1.01}{1.99} + 4.13}{15.6} + \frac{\frac{8(1.01)^2}{1.99}}{15.6 \times 7.8} + \frac{4.13}{7.3} 
< .995 < 1
$,
we see that for all sufficiently small 
$\varsigma > 0$ and $\fundbootsmall \geq 0$,
we can soak the second product on RHS~\eqref{E:MAINWAVEESTIMATESINTERMEDIATEF} 
back into LHS~\eqref{E:MAINWAVEESTIMATESINTERMEDIATEF} 
and soak the last factor on RHS~\eqref{E:MAINWAVEESTIMATESINTERMEDIATEH} back into 
LHS~\eqref{E:MAINWAVEESTIMATESINTERMEDIATEH},
thereby deducing: 
\begin{align} \label{E:HOWTOUSEMAINWAVEESTIMATESINTERMEDIATEF} 
	\apriorimain(\timefunction,u) 
	& \leq C \left(1 + \varsigma^{-1} \right) \initialsmall^2 
		+
		C 
		\left\lbrace 
			\fundbootsmall 
			+ 
			\varsigma 
			+ 
			\underbrace{\frac{1}{\mathfrak{c}} 
			\left(1 + \varsigma^{-1}\right)}_{\varsigma^2 + \varsigma}
		\right\rbrace 
		\apriorilower(\timefunction,u)
		+
		C \apriorimain^{1/2}(\timefunction,u) \aprioripartial^{1/2}(\timefunction,u)	
\end{align}
and:
\begin{align} \label{E:HOWTOUSEMAINWAVEESTIMATESINTERMEDIATEH}
\apriorilower(\timefunction,u) 
\leq 
C \initialsmall^2 
+ 
C \apriorimain(\timefunction,u).
\end{align}
Inserting \eqref{E:HOWTOUSEMAINWAVEESTIMATESINTERMEDIATEH} estimate into RHS~\eqref{E:MAINWAVEESTIMATESINTERMEDIATEG},
we find that for all sufficiently small $\varsigma > 0$ and $\varepsilon \geq 0$, we have:
\begin{align} \label{E:HOWTOUSEMAINWAVEESTIMATESINTERMEDIATEG} 
\aprioripartial(\timefunction,u) 
& \leq 
		C \left(1 + \varsigma^{-1}\right) 
		\initialsmall^2 
		+ 
		C \left\lbrace 
			\fundbootsmall 
			+ 
			\varsigma 
		\right\rbrace 
		\apriorimain(\timefunction,u).  
\end{align}
Next, using Young's inequality, we bound 
the last product on RHS~\eqref{E:HOWTOUSEMAINWAVEESTIMATESINTERMEDIATEF} 
as follows: 
$C \apriorimain^{1/2}(\timefunction,u) \aprioripartial^{1/2}(\timefunction,u) 
\leq 
\frac{1}{2} \apriorimain(\timefunction,u)
+
C
\aprioripartial(\timefunction,u) 
 $.
The term $\frac{1}{2} \apriorimain(\timefunction,u)$ can be absorbed back into LHS~\eqref{E:HOWTOUSEMAINWAVEESTIMATESINTERMEDIATEF},
which yields:
\begin{align} \label{E:SECONDHOWTOUSEMAINWAVEESTIMATESINTERMEDIATEF} 
	\apriorimain(\timefunction,u) 
	& \leq C \left(1 + \varsigma^{-1} \right) \initialsmall^2 
		+
		C 
		\left\lbrace 
			\fundbootsmall 
			+ 
			\varsigma 
			+ 
			\underbrace{\frac{1}{\mathfrak{c}} 
			\left(1 + \varsigma^{-1}\right)}_{\varsigma^2 + \varsigma}
		\right\rbrace 
		\apriorilower(\timefunction,u)
		+
		C \aprioripartial(\timefunction,u).
\end{align}
Inserting 
\eqref{E:HOWTOUSEMAINWAVEESTIMATESINTERMEDIATEH}--\eqref{E:HOWTOUSEMAINWAVEESTIMATESINTERMEDIATEG} into
\eqref{E:SECONDHOWTOUSEMAINWAVEESTIMATESINTERMEDIATEF},
we find that:
\begin{align} \label{E:THIRDHOWTOUSEMAINWAVEESTIMATESINTERMEDIATEF} 
	\apriorimain(\timefunction,u) 
	& \leq C \left(1 + \varsigma^{-1} \right) \initialsmall^2 
		+
		C 
		\left\lbrace 
			\fundbootsmall 
			+ 
			\varsigma 
		\right\rbrace 
		\apriorimain(\timefunction,u).
\end{align}
Hence, if $\varsigma > 0$ is sufficiently small, then for all
sufficiently small $\fundbootsmall \geq 0$,
we can absorb the last product on RHS~\eqref{E:THIRDHOWTOUSEMAINWAVEESTIMATESINTERMEDIATEF} back into the LHS.
This implies the desired bound
$
\apriorimain(\timefunction,u) 
\leq C \left(1 + \varsigma^{-1}\right) \initialsmall^2$.
Inserting this bound into RHSs~\eqref{E:HOWTOUSEMAINWAVEESTIMATESINTERMEDIATEH}--\eqref{E:HOWTOUSEMAINWAVEESTIMATESINTERMEDIATEG},
we find that
$
\aprioripartial(\timefunction,u) 
\leq 
C 
\left(1 + \varsigma^{-1}\right)
\initialsmall^2
$
and
$
\apriorilower(\timefunction,u) 
\leq 
C \left(1 + \varsigma^{-1}\right) \initialsmall^2
$.
We have therefore proved \eqref{E:MAINWAVEESTIMATESINTERMEDIATE1}.

\medskip

\noindent \emph{Proof of \eqref{E:MAINWAVEESTIMATESINTERMEDIATEF}--\eqref{E:MAINWAVEESTIMATESINTERMEDIATEH}}: 
It remains for us to prove
\eqref{E:MAINWAVEESTIMATESINTERMEDIATEF}--\eqref{E:MAINWAVEESTIMATESINTERMEDIATEH}.  
We set $N \eqdef \Ntop$. We fix any $(\timefunction,u) \in [\timefunction_0,\timefunctionboot]\times[-\rightu,\leftu]$, 
and we let $(\hat{\timefunction},\hat{u}) \in [\timefunction_0,\timefunction] \times [-\rightu,u]$. 
We evaluate the top-order integral inequality \eqref{E:TOPORDERWAVEL2CONTROLLINGINTEGRALINEQUALITY} 
at $(\hat{\timefunction},\hat{u})$ and multiply it by $\iota_{\apriorimain}^{-1}(\hat{\timefunction},\hat{u})$. 
Similarly, we evaluate \eqref{E:MAINWAVEPARTIALINTEGRALINEQUALITIES}--\eqref{E:MAINWAVEBELOWTOPINTEGRALINEQUALITIES} at $(\hat{\timefunction},\hat{u})$ and respectively multiply by $\iota_{\aprioripartial}^{-1}(\hat{\timefunction},\hat{u})$ and $\iota_{\apriorilower}^{-1}(\hat{\timefunction},\hat{u})$. We then obtain suitable bounds from the resulting products and then take $ \sup_{(\hat{\timefunction},\hat{u}) \in [\timefunction_0,\timefunction]\times[-\rightu,u]}$. The left-hand sides of the resulting inequalities are, by definition,
equal to the left-hand sides of \eqref{E:MAINWAVEESTIMATESINTERMEDIATEF}--\eqref{E:MAINWAVEESTIMATESINTERMEDIATEH},
while our ``suitable bounds,'' which we derive below, will yield the right-hand sides of
\eqref{E:MAINWAVEESTIMATESINTERMEDIATEF}--\eqref{E:MAINWAVEESTIMATESINTERMEDIATEH}.

We now prove \eqref{E:MAINWAVEESTIMATESINTERMEDIATEF}.
We will explain how to handle several representative terms on RHS~\eqref{E:MAINWAVEESTIMATESINTERMEDIATEF},
including the most difficult terms. The remaining terms can be handled using similar
or simpler arguments, and we omit the details.
As our first example, we consider the term 
$
C
\int_{\timefunction' = \timefunction_0}^{\hat{\timefunction}} 
			\frac{1}{|\timefunction'|^{4/3}} 
			\left\lbrace
				\int_{\timefunction'' = \timefunction_0}^{\timefunction'}
					\left(
						\hypersurfacecontrolVortVort_N^{1/2} + \hypersurfacecontrolDivGradEnt_N^{1/2}
					\right)
					(\timefunction'',u) 
				\, \mathrm{d} \timefunction'' 
			\right\rbrace^2 
		\, \mathrm{d}\timefunction' 
$ 
generated by the 
$3^{\text{rd}}$ line of RHS~\eqref{E:ERRORTOPORDERWAVEESTIMATES}. 
Using the already proved a priori estimates \eqref{E:MAINTOPORDERENERGYESTIMATESMODIFIEDFLUIDVARIABLESBLOWUP} for 
$\hypersurfacecontrolVortVort_{\Ntop}$ and $\hypersurfacecontrolDivGradEnt_{\Ntop}$,
we bound this term as follows:
\begin{align} 
\begin{split} \label{E:MAINWAVEESTIMATESINTERMEDIATE3}
	&  
	C 
	\int_{\timefunction' = \hat\timefunction_0}^{\hat\timefunction} 
		\frac{1}{|\timefunction'|^{4/3}} 
		\left\lbrace 
			\int_{\timefunction'' = \hat\timefunction_0}^{\timefunction'} 
				\left(\hypersurfacecontrolVortVort_N^{1/2} + \hypersurfacecontrolDivGradEnt_N^{1/2}\right)(\timefunction'',\hat{u}) 
			\, \mathrm{d} \timefunction'' 
		\right\rbrace^2
	\, \mathrm{d}\timefunction'  
			\\
	& \lesssim 
	\initialsmall^2 
	\int_{\timefunction' =  \timefunction_0}^{\hat\timefunction} 
		\frac{1}{|\timefunction'|^{4/3}}
		\left\lbrace \int_{\timefunction'' = \hat\timefunction_0}^{\timefunction'} 
			|\timefunction''|^{-8.55} 
			\, \mathrm{d} \timefunction'' 	
		\right\rbrace^2 \, \mathrm{d}\timefunction' 
			\\
	&  \lesssim 
		\initialsmall^2 
		\int_{\timefunction' =  \timefunction_0}^{\hat\timefunction} 
			\frac{1}{|\timefunction'|^{4/3}}
			|\timefunction'|^{-15.1} 
		\, \mathrm{d} \timefunction' 
			\\
	& \lesssim \initialsmall^2|\hat{\timefunction}|^{-(15 + 13/30)} 
	\lesssim  
	\initialsmall^2|\hat{\timefunction}|^{-15.6}.
\end{split}
\end{align}
Multiplying \eqref{E:MAINWAVEESTIMATESINTERMEDIATE3} 
by $\iota_{\apriorimain}^{-1}(\hat{\timefunction},\hat{u})$ and taking 
$\sup_{(\hat{\timefunction},\hat{u}) \in [\timefunction_0,\timefunction]\times[-\rightu,u]}$,
we obtain:
\begin{align} \label{E:MAINWAVEESTIMATESINTERMEDIATE4}
C \sup_{(\hat{\timefunction},\hat{u}) \in [\timefunction_0,\timefunction]\times[-\rightu,u]} 
	\iota_{\apriorimain}^{-1}(\hat{\timefunction},\hat{u})  
	\int_{\timefunction' = \hat\timefunction_0}^{\hat\timefunction} 
		\frac{1}{|\timefunction'|^{4/3}} 
		\left\lbrace 
			\int_{\timefunction'' = \hat\timefunction_0}^{\timefunction'} 
				\left(\hypersurfacecontrolVortVort_N^{1/2} + \hypersurfacecontrolDivGradEnt_N^{1/2}\right)(\timefunction'',\hat{u}) 
			\, \mathrm{d} \timefunction'' 
		\right\rbrace^2 
	\, \mathrm{d}\timefunction'   
	& 
	\leq 
	C \initialsmall^2,
\end{align}
which is $\leq \mbox{RHS~\eqref{E:MAINWAVEESTIMATESINTERMEDIATEF}}$ as desired.
The term 
$
\int_{\timefunction' = \timefunction_0}^{\hat{\timefunction}} 
		\frac{1}{|\timefunction'|^{4/3}} 
		\left\lbrace \int_{\timefunction'' = \timefunction_0}^{\timefunction'} 
			\frac{1}{|\timefunction''|^{1/2}} \left(\hypersurfacecontrolVortVort_{\leq N-1}^{1/2} 
			+ 
			\hypersurfacecontrolDivGradEnt_{\leq N-1}^{1/2}\right)(\timefunction'',u) 
		\, \mathrm{d} \timefunction'' \right\rbrace^2 
	\, \mathrm{d}\timefunction' 
$ 
generated by the $4^{\text{th}}$ line of RHS~\eqref{E:ERRORTOPORDERWAVEESTIMATES} can be handled using similar arguments,
this time with the help of the already proven estimates 
\eqref{E:MAINL2BELOWTOPORDERESTIMATESMODIFIEDFLUIDBLOWUP} and \eqref{E:MAINL2BELOWTOPORDERESTIMATESVORTANDENTROPYGRADIENTREGULAR}.
Similarly,
the following estimate holds for the term generated by the first term on RHS~\eqref{E:ERRORTOPORDERWAVEESTIMATES}:
\begin{align} \label{E:MAINWAVEESTIMATESINTERMEDIATE5}
C 
\sup_{(\hat{\timefunction},\hat{u}) \in [\timefunction_0,\timefunction]\times[-\rightu,u]} 
\iota_{\apriorimain}^{-1} 
\left(1+\varsigma^{-1}\right) 
\initialsmall^2 
\frac{1}{|\hat{\timefunction}|^{3/2}} 
& 
\leq
C
\left(1 + \varsigma^{-1} \right)
\initialsmall^2,
\end{align} 
which is $\leq \mbox{RHS~\eqref{E:MAINWAVEESTIMATESINTERMEDIATEF}}$ as desired.

We now handle the first term on RHS~\eqref{E:TOPORDERWAVEL2CONTROLLINGINTEGRALINEQUALITY},
i.e., the one multiplied by the boxed constant 
$\boxed{\left\lbrace \frac{4 \times 1.01}{1.99} + 4.13 \right\rbrace}$. 
Evaluating the term at 
$(\hat{\timefunction},\hat{u})$,
multiplying it by 
$\iota_{\apriorimain}^{-1}(\hat{\timefunction},\hat{u})$, 
multiplying and dividing the integrand by $|\timefunction'|^{15.6}$, 
taking 
$
\sup_{\timefunction' \in [\timefunction_0,\hat{\timefunction}]} 
\hypersurfacecontrolwave_N(\timefunction',\hat{u}) |\timefunction'|^{15.6}$, 
and pulling the sup-ed out quantity outside of the integral, we deduce: 
\begin{align}
\begin{split} \label{E:MAINWAVEESTIMATESINTERMEDIATE6}
	&  
		\left\lbrace 
			\frac{4 \times 1.01}{1.99} + 4.13 
		\right\rbrace 
		\iota_{\apriorimain}^{-1}(\hat{\timefunction},\hat{u}) 
		\int_{\timefunction' =  \timefunction_0}^{\hat{\timefunction}} 
			\frac{1}{|\timefunction'|}  \hypersurfacecontrolwave_N (\timefunction',\hat{u}) 
		\, \mathrm{d} \timefunction' 
		\\
	& \leq   
		\left\lbrace 
			\frac{4\times 1.01}{1.99} + 4.13 
		\right\rbrace  
		\iota_{\apriorimain}^{-1}(\hat{\timefunction},\hat{u})  
		\times \sup_{(\timefunction',u') \in [\timefunction_0,\hat{\timefunction}]\times[-\rightu,\hat{u}]} 
		\left\lbrace 
			\hypersurfacecontrolwave_N(\timefunction',\hat{u}) |\timefunction'|^{15.6}
		\right\rbrace 
		\times 
		\int_{\timefunction' = \timefunction_0}^{\hat{\timefunction}}  
			\frac{1}{|\timefunction'|^{16.6}}  
		\, \mathrm{d} \timefunction' 
			\\
	& \leq 
		\frac{1}{15.6}  
		\left\lbrace 
			\frac{4\times 1.01}{1.99} + 4.13
		\right\rbrace \iota_{\apriorimain}^{-1}(\hat{\timefunction},\hat{u}) 
		\times \sup_{(\timefunction',u') \in [\timefunction_0,\hat{\timefunction}]\times[-\rightu,\hat{u}]}  
		\left\lbrace 
			\hypersurfacecontrolwave_N(\timefunction',\hat{u}) |\timefunction'|^{15.6}
		\right\rbrace  
		\times |\hat{\timefunction}|^{-15.6} 
			\\
	& 
	\leq  
	\frac{1}{15.6} 
	\left\lbrace 
		\frac{4\times 1.01}{1.99} + 4.13  
	\right\rbrace \apriorimain(\hat{\timefunction},\hat{u}) 
		\\
	& 
	\leq
	\frac{1}{15.6}  
	\left\lbrace 
		\frac{4 \times 1.01}{1.99} + 4.13 
	\right\rbrace 
	\apriorimain(\timefunction,u).
\end{split}
\end{align}
Taking $\sup_{(\hat{\timefunction},\hat{u}) \in [\timefunction_0,\timefunction]\times[-\rightu,u]}$
of \eqref{E:MAINWAVEESTIMATESINTERMEDIATE6},
we find see that
$\sup_{(\hat{\timefunction},\hat{u}) \in [\timefunction_0,\timefunction]\times[-\rightu,u]}
\mbox{LHS~\eqref{E:MAINWAVEESTIMATESINTERMEDIATE6}}
\leq
\mbox{RHS~\eqref{E:MAINWAVEESTIMATESINTERMEDIATE6}}
$,
and the terms on RHS~\eqref{E:MAINWAVEESTIMATESINTERMEDIATE6} are
part of the ``main terms'' located on the second line of RHS~\eqref{E:MAINWAVEESTIMATESINTERMEDIATEF}.
We clarify that to obtain the third inequality in \eqref{E:MAINWAVEESTIMATESINTERMEDIATE6}, 
we multiplied by 
$
1
=
\mathfrak{l}^{\mathfrak{c}}(\timefunction')
\mathfrak{l}^{\mathfrak{c}}(\timefunction') 
e^{\mathfrak{c}u'}
e^{-\mathfrak{c}u'}$ 
under the sup,
used the monotonicity of 
$\mathfrak{l}^{\mathfrak{c}}(\timefunction')$ and $e^{\mathfrak{c}u'}$ 
to bound these two factors by $\mathfrak{l}^{\mathfrak{c}}(\hat{\timefunction})$ and $e^{\mathfrak{c}\hat{u}}$ (and so the remaining sup-ed quantity is $\sup_{(\timefunction',u')\in [\timefunction_0,\hat{\timefunction}],[-\rightu,\hat{u}]} \{|\timefunction'|^{15.6} \mathfrak{l}^{-\mathfrak{c}} (\timefunction')e^{-\mathfrak{c}u'}\hypersurfacecontrolwave_{\Ntop}(\timefunction',u')\} = \apriorimain(\hat{\timefunction},\hat{u})$), pulled the two factors 
$\mathfrak{l}^{\mathfrak{c}}(\hat{\timefunction})$ and $e^{\mathfrak{c}\hat{u}}$
out of the sup, and then note that these two factors and
$|\hat{\timefunction}|^{-15.6}$ multiply together to exactly cancel
$\iota_{\apriorimain}^{-1}(\hat{\timefunction},\hat{u})$. 

We can handle the second term on RHS~\eqref{E:TOPORDERWAVEL2CONTROLLINGINTEGRALINEQUALITY} 
(which is multiplied by $\boxed{\frac{8\times (1.01)^2}{1.99}}$)
using similar arguments, but we have to integrate twice in time.
We find that the corresponding term is: 
\begin{align} 
\begin{split} \label{E:MAINWAVEESTIMATESINTERMEDIATE7}
	& 
	\leq
	\boxed{\frac{8\times (1.01)^2}{1.99}}
	\iota_{\apriorimain}^{-1}(\hat{\timefunction},\hat{u})  
	\int_{\timefunction' =  \timefunction_0}^{\hat{\timefunction}} 
		\frac{1}{|\timefunction'|}  
		\hypersurfacecontrolwave_{ N}^{1/2}(\timefunction',\hat{u}) 
		\int_{\timefunction'' = \timefunction_0}^{\timefunction'} 
			\frac{1}{|\timefunction''|} \hypersurfacecontrolwave_N^{1/2}(\timefunction'',\hat{u}) 
		\, \mathrm{d}\timefunction'' 
	\, \mathrm{d} \timefunction' 
		\\
	& 
	\leq 
	\frac{1}{15.6 \times 7.8} 
	\left\lbrace 
		\frac{8\times (1.01)^2}{1.99} 
	\right\rbrace 
	\apriorimain(\hat{\timefunction},\hat{u}) 
	\leq 
	\frac{1}{15.6 \times 7.8} 
	\left\lbrace 
		\frac{8 \times (1.01)^2}{1.99}
	\right\rbrace 
	\apriorimain(\timefunction,u).  
\end{split}
\end{align}
The terms on RHS~\eqref{E:MAINWAVEESTIMATESINTERMEDIATE7} are also
part of the ``main terms'' located on the second line of RHS~\eqref{E:MAINWAVEESTIMATESINTERMEDIATEF}.

We can handle the third term on RHS~\eqref{E:TOPORDERWAVEL2CONTROLLINGINTEGRALINEQUALITY} 
(which is multiplied by $\boxed{4.13}$)
using similar arguments, involving only one integration in time,
to deduce:
\begin{align} 
\begin{split} \label{E:MAINWAVEESTIMATESINTERMEDIATE8}
	&
	\boxed{4.13}  
	\iota_{\apriorimain}^{-1}(\hat{\timefunction},\hat{u})  
	\frac{1}{|\hat{\timefunction}|^{1/2}} 
	\hypersurfacecontrolwave_N^{1/2}(\hat{\timefunction},\hat{u})   
	\int_{\timefunction' = \hat{\timefunction}_0}^{\hat{\timefunction}} \frac{1}{|\timefunction'|^{1/2}} 
	\hypersurfacecontrolwave_N^{1/2}(\timefunction',\hat{u}) 
	\, \mathrm{d} \timefunction' 
		\\
	& \leq   
	\left\lbrace 
		\frac{4.13}{7.3}
	\right\rbrace 
	\apriorimain(\hat{\timefunction},\hat{u}) 
	\leq 
	\left\lbrace 
		\frac{4.13}{7.3}
	\right\rbrace 
	\apriorimain(\timefunction,u).
\end{split}
\end{align}
The terms on RHS~\eqref{E:MAINWAVEESTIMATESINTERMEDIATE8} provide the last contribution to
the ``main terms'' located on the second line of RHS~\eqref{E:MAINWAVEESTIMATESINTERMEDIATEF}.

Using the same arguments we used to prove
\eqref{E:MAINWAVEESTIMATESINTERMEDIATE3}--\eqref{E:MAINWAVEESTIMATESINTERMEDIATE8},
we can bound the contribution of the three $C_*$-multiplied terms on 
RHS~\eqref{E:TOPORDERWAVEL2CONTROLLINGINTEGRALINEQUALITY} 
by $\leq C \apriorimain^{1/2}(\timefunction,u) \aprioripartial^{1/2}(\timefunction,u)$,
which in turn is bounded by the last term on RHS~\eqref{E:MAINWAVEESTIMATESINTERMEDIATEF}.

The remaining terms on RHS~\eqref{E:MAINWAVEESTIMATESINTERMEDIATEF},
which are all found in $\Errortop(\hat{\timefunction},\hat{u})$, 
are less dangerous than the three main terms we treated in \eqref{E:MAINWAVEESTIMATESINTERMEDIATE6}--\eqref{E:MAINWAVEESTIMATESINTERMEDIATE8}
because they are either 
\textbf{I)} critical with respect to the energy blow-up rates but feature a small factor of $C \fundbootsmall$,
\textbf{II)} sub-critical\footnote{By a critical term, we mean that by inserting the desired estimates 
of Prop.\,\ref{P:APRIORIL2ESTIMATESWAVEVARIABLES} into it, 
one derives an estimate that is exactly compatible with the blow up rates of Prop.\,\ref{P:APRIORIL2ESTIMATESWAVEVARIABLES} in terms of powers of $|\timefunction|^{-1}$. 
By sub-critical term, we mean that the estimate one derives after inserting 
the desired estimate of Prop.\,\ref{P:APRIORIL2ESTIMATESWAVEVARIABLES} is less singular with respect to powers of $|\timefunction|^{-1}$.} with respect to the blow-up rates,
or \textbf{III)} controlled by the already proven a priori estimates for the transport variables 
from Prop.\,\ref{P:MAINHYPERSURFACEENERGYESTIMATESFORTRANSPORTVARIABLES} 
(it turns out that all these terms are sub-critical too).
The type \textbf{I} terms can be handled as in
\eqref{E:MAINWAVEESTIMATESINTERMEDIATE6}--\eqref{E:MAINWAVEESTIMATESINTERMEDIATE8},
but they are much less delicate because we do not have to be careful about the size of the constants;
such terms contribute to the $C \fundbootsmall$-multiplied terms on the second line of RHS~\eqref{E:MAINWAVEESTIMATESINTERMEDIATEF}.
The type \textbf{III} terms have already been adequately controlled 
in Prop.\,\ref{P:MAINHYPERSURFACEENERGYESTIMATESFORTRANSPORTVARIABLES}
and contribute only to the
term $C \left(1 + \varsigma^{-1} \right) \initialsmall^2$ on RHS~\eqref{E:HOWTOUSEMAINWAVEESTIMATESINTERMEDIATEF},
as in
\eqref{E:MAINWAVEESTIMATESINTERMEDIATE4}--\eqref{E:MAINWAVEESTIMATESINTERMEDIATE5}.
The type \textbf{II} terms can be handled using arguments that rely only on the multiplicative factors 
$\mathfrak{l}^{\mathfrak{c}}(\timefunction)$ and $e^{\mathfrak{c} u}$
in
\eqref{E:TIMEONLYGRONWALLMULTIPLICATIONFACTORTOPFULLANDPARTIAL}--\eqref{E:GRONWALLMULTIPLICATIONFACTORTOPFULLANDPARTIAL}.
We will handle three representative type \textbf{II} terms:
the term on the next-to-last line of RHS~\eqref{E:ERRORTOPORDERWAVEESTIMATES} 
featuring a triple integral, 
the $u'$-integral term on the fifth-from-last line of RHS~\eqref{E:ERRORTOPORDERWAVEESTIMATES}, 
and the $\timefunction'$ integral featuring the lower-order term $\hypersurfacecontrolwave_{[1,N-1]}$ on the last line
of RHS~\eqref{E:ERRORTOPORDERWAVEESTIMATES}.
The remaining terms on RHS~\eqref{E:ERRORTOPORDERWAVEESTIMATES} can be handled using similar
arguments, where we use the estimates of Prop.\,\ref{P:MAINHYPERSURFACEENERGYESTIMATESFORTRANSPORTVARIABLES}
and argue as in
\eqref{E:MAINWAVEESTIMATESINTERMEDIATE3}--\eqref{E:MAINWAVEESTIMATESINTERMEDIATE4}
to handle error terms that involve the (already bounded) quantities
$
\hypersurfacecontrolVortVort_M
$,
$
\hypersurfacecontrolDivGradEnt_M
$,
$
\hypersurfacecontrolVortVort_M
$,
and
$
\hypersurfacecontrolDivGradEnt_M
$; 
we omit the details.

We will show that the contribution of the first representative term can be bounded as follows:
\begin{align}
\begin{split} \label{E:TRIPPLEINTEGRALERRORTERMGRONWALARGUMENT}
	& C 
	\iota_{\apriorimain}^{-1}(\hat{\timefunction},\hat{u}) 
	\int_{\timefunction' = \timefunction_0}^{\hat{\timefunction}} 
		\frac{1}{|\timefunction'|}  \hypersurfacecontrolwave_N^{1/2}(\timefunction',\hat{u}) 
		\int_{\timefunction'' = \timefunction_0}^{\timefunction'}  
			\frac{1}{|\timefunction''|} 
			\int_{\timefunction''' = \timefunction_0}^{\timefunction''} 
				\frac{1}{|\timefunction'''|^{1/2}} \hypersurfacecontrolwave_N^{1/2}(\timefunction''',\hat{u}) 
			\, \mathrm{d} \timefunction''' 
		\, \mathrm{d} \timefunction'' 
	\, \mathrm{d}\timefunction'
		\\
	& \leq 
		\frac{C}{\mathfrak{c}} 
		\iota_{\apriorimain}^{-1}(\hat{\timefunction},\hat{u}) 
		\mathfrak{l}^{\frac{\mathfrak{c}}{2}}(\hat{\timefunction}) 
		\int_{\timefunction' = \timefunction_0}^{\hat{\timefunction}}  
			\frac{1}{|\timefunction'|}  
			\hypersurfacecontrolwave_N^{1/2}(\timefunction',\hat{u})  
			\int_{\timefunction'' = \timefunction_0}^{\timefunction'}  
				\frac{1}{|\timefunction''|}  
				\sup_{(\timefunction''',u')\in[\timefunction_0,\timefunction'']\times[-\rightu,\hat{u}]} 
				\left\lbrace 
					\mathfrak{l}^{-\frac{\mathfrak{c}}{2}}(\timefunction''') 
					\hypersurfacecontrolwave_N^{1/2}(\timefunction''',u')
				\right\rbrace \,  \mathrm{d} \timefunction'' 
		\, \mathrm{d}\timefunction' 
		\\
	& \leq 
	\frac{C}{\mathfrak{c}} 
	\iota_{\apriorimain}^{-1}(\hat{\timefunction},\hat{u}) 
	\mathfrak{l}^{\frac{\mathfrak{c}}{2}}(\hat{\timefunction}) 
		\\
& \ \ \ \
	\times
	\sup_{(\timefunction''',u')\in[\timefunction_0,\hat{\timefunction}]\times[-\rightu,\hat{u}]} 
	\left\lbrace 
		|\timefunction'''|^{7.8}
		\mathfrak{l}^{-\frac{\mathfrak{c}}{2}}(\timefunction''') 
		\hypersurfacecontrolwave_N^{1/2}(\timefunction''',u')
	\right\rbrace
	\times
	\sup_{(\timefunction'''',u'')\in[\timefunction_0,\hat{\timefunction}]\times[-\rightu,\hat{u}]} 
	\left\lbrace 
		|\timefunction''''|^{7.8}
		\hypersurfacecontrolwave_N^{1/2}(\timefunction'''',u'')
	\right\rbrace
		\\
& \ \ \ \
	\times
	\int_{\timefunction' = \timefunction_0}^{\hat{\timefunction}}   
		\frac{1}{|\timefunction'|^{8.8}}  
		\int_{\timefunction'' = \timefunction_0}^{\timefunction'} 
			\frac{1}{|\timefunction''|^{8.8}} 
		\, \mathrm{d} \timefunction'' 
	\, \mathrm{d} \timefunction' 
		\\
	& \leq  
	\frac{C}{\mathfrak{c}} 
	\iota_{\apriorimain}^{-1}(\hat{\timefunction},\hat{u}) 
	\mathfrak{l}^{\mathfrak{c}}(\hat{\timefunction}) 
	\sup_{(\timefunction'',u')\in[\timefunction_0,\hat{\timefunction}]\times[-\rightu,\hat{u}]} 
	\left\lbrace 
		|\timefunction''|^{15.6}
		\mathfrak{l}^{-\mathfrak{c}}(\timefunction'')
		\hypersurfacecontrolwave_N(\timefunction'',u')
	\right\rbrace
	\int_{\timefunction' = \timefunction_0}^{\hat{\timefunction}}   
		\frac{1}{|\timefunction'|^{16.6}}  
	\, \mathrm{d} \timefunction' 
		\\
	& \leq   
		\frac{C}{\mathfrak{c}} 
		\iota_{\apriorimain}^{-1}(\hat{\timefunction},\hat{u}) 
		|\hat{\timefunction}|^{-15.6}
		\mathfrak{l}^{\mathfrak{c}}(\hat{\timefunction})
		e^{\mathfrak{c} \hat{u}}
		\sup_{(\timefunction'',u') \in [\timefunction_0,\hat{\timefunction}] \times [-\rightu,\hat{u}]} 
		\left\lbrace
			|\timefunction''|^{15.6}
			\mathfrak{l}^{-\mathfrak{c}}(\timefunction'')
			e^{-\mathfrak{c}u'}
			\hypersurfacecontrolwave_N(\timefunction'',u')
		\right\rbrace  
			\\
	& \leq   
	\frac{C}{\mathfrak{c}} 
	\apriorimain(\hat{\timefunction},\hat{u}) 
	\leq  
	\frac{C}{\mathfrak{c}} \apriorimain(\timefunction,u).
\end{split}
\end{align}
Taking $\sup_{(\hat{\timefunction},\hat{u}) \in [\timefunction_0,\timefunction]\times[-\rightu,u]}$
of \eqref{E:TRIPPLEINTEGRALERRORTERMGRONWALARGUMENT},
we find see that
$\sup_{(\hat{\timefunction},\hat{u}) \in [\timefunction_0,\timefunction]\times[-\rightu,u]}
\mbox{LHS~\eqref{E:TRIPPLEINTEGRALERRORTERMGRONWALARGUMENT}}
\leq
\mbox{RHS~\eqref{E:TRIPPLEINTEGRALERRORTERMGRONWALARGUMENT}}
$,
which is $\leq \mbox{RHS~\eqref{E:MAINWAVEESTIMATESINTERMEDIATEF}}$ as desired.
We now justify the sequence of inequalities in \eqref{E:TRIPPLEINTEGRALERRORTERMGRONWALARGUMENT}.
The first inequality follows from multiplying and dividing the $\mathrm{d} \timefunction'''$ integrand by 
$\mathfrak{l}^{\frac{\mathfrak{c}}{2}}(\timefunction''')$, 
pulling out 
$\sup_{(\timefunction''',u')\in[\timefunction_0,\timefunction'']\times[-\rightu,\hat{u}]} 
\left\lbrace 
	\mathfrak{l}^{-\frac{\mathfrak{c}}{2}}(\timefunction''') 
	\hypersurfacecontrolwave_N^{1/2}(\timefunction''',u')
\right\rbrace
$ 
from the $\mathrm{d} \timefunction'''$ integral, 
and using \eqref{E:INTEGRALOFINTEGRATINGFACTORS} to 
bound the remaining inner-most time integral
$
\int_{\timefunction''' = \timefunction_0}^{\timefunction''} 
		\frac{\mathfrak{l}^{\frac{\mathfrak{c}}{2}}(\timefunction''')}{|\timefunction'''|^{1/2}}
\, \mathrm{d} \timefunction''' 
$
by $\leq \frac{C}{\mathfrak{c}} \mathfrak{l}^{\frac{\mathfrak{c}}{2}}(\hat{\timefunction})$.
The second inequality follows from multiplying and dividing the 
$\mathrm{d} \timefunction''$
integrand by $|\timefunction''|^{7.8}$ 
and pulling out a sup-ed quantity from the integral,
and from multiplying and dividing the 
$\mathrm{d} \timefunction'$
integrand by $|\timefunction'|^{7.8}$ 
and pulling out a sup-ed quantity from the integral.
The remaining inequalities follow from straightforward integration,
the definitions of the quantities involved,
and the monotonicity of various factors.

Next, using \eqref{E:INTEGRALOFINTEGRATINGFACTORS} and arguments similar to but simpler
than the ones we used to prove \eqref{E:TRIPPLEINTEGRALERRORTERMGRONWALARGUMENT},
we bound the term generated by the $u'$-integral on the fifth-from-last line of \eqref{E:ERRORTOPORDERWAVEESTIMATES}
as follows:
\begin{align}
\begin{split} \label{E:UINTEGRALERRORTERMGRONWALLARGUMENT}
	&
	C 
	\left(1 + \varsigma^{-1} \right) 
	\iota_{\apriorimain}^{-1}(\hat{\timefunction},\hat{u}) 
	\int_{u' = -\rightu}^{\hat{u}} 
		\hypersurfacecontrolwave_{N}(\hat{\timefunction},u') 
	\, \mathrm{d} u' 
		\\
	&
	\leq 	
	\left(1 + \varsigma^{-1} \right) 
	\iota_{\apriorimain}^{-1}(\hat{\timefunction},\hat{u}) 
	\sup_{(\timefunction',u') \in [\timefunction_0,\hat{\timefunction}]\times[-\rightu,\hat{u}]}
	\{ e^{-\mathfrak{c}u'} \hypersurfacecontrolwave_{N}(\timefunction',u')\} 
	\int_{u' = -\rightu}^{\hat{u}} 
		e^{\mathfrak{c} u'} 
	\, \mathrm{d} u' 
		\\
	& 
	\leq 
	\frac{C}{\mathfrak{c}} 
	\left(1 + \varsigma^{-1} \right) 
	\iota_{\apriorimain}^{-1}(\hat{\timefunction},\hat{u}) 
	e^{\mathfrak{c} \hat{u}}  
	\sup_{(\timefunction',u') \in [\timefunction_0,\hat{\timefunction}]\times[-\rightu,\hat{u}]}
	\left\lbrace 
		e^{-\mathfrak{c}u'} \hypersurfacecontrolwave_{N}(\timefunction',u')
	\right\rbrace
		\\
	& 
	\leq 
		\frac{C}{\mathfrak{c}} 
		\iota_{\apriorimain}^{-1}(\hat{\timefunction},\hat{u}) 
		|\hat{\timefunction}|^{-15.6}
		\mathfrak{l}^{\mathfrak{c}}(\hat{\timefunction}) 
		e^{\mathfrak{c} \hat{u}}
		\sup_{(\timefunction',u')\in[\timefunction_0,\hat{\timefunction}]\times[-\rightu,\hat{u}]} 
		\left\lbrace 
			|\timefunction'|^{15.6}
			\mathfrak{l}^{-\mathfrak{c}}(\timefunction')
			e^{-\mathfrak{c}u'}
			\hypersurfacecontrolwave_N(\timefunction',u')
		\right\rbrace  
		\\
	& 
	\leq 
	\frac{C}{\mathfrak{c}} 
	\left(1 + \varsigma^{-1} \right)
	\apriorimain(\hat{\timefunction},\hat{u}) 
	\leq \frac{C}{\mathfrak{c}} 
	\left(1 + \varsigma^{-1} \right) 
	\apriorimain(\timefunction,u), 
\end{split}
\end{align}
which is $\leq \mbox{RHS~\eqref{E:MAINWAVEESTIMATESINTERMEDIATEF}}$ 
as desired.

Similarly, we bound the term generated by the 
$\hypersurfacecontrolwave_{[1,N-1]}$ integral 
on the last line of RHS~\eqref{E:ERRORTOPORDERWAVEESTIMATES} 
as follows:
\begin{align}
\begin{split} \label{E:CONTROBUTIONOFBELOWTOPINGRONWALARGUMENTFORTOP}
	&
	C 
	\iota_{\apriorimain}^{-1}(\hat{\timefunction},\hat{u}) 
	\int_{\timefunction' = \timefunction_0}^{\hat{\timefunction}} 
		\frac{1}{|\timefunction'|^{5/2}} 
		\hypersurfacecontrolwave_{[1,N-1]}(\timefunction',\hat{u}) 
	\, \mathrm{d} \timefunction'
		\\
& = C 
		|\hat{\timefunction}|^{15.6} 
		\mathfrak{l}^{-\mathfrak{c}}(\hat{\timefunction})
		e^{-\mathfrak{c} \hat{u}} 
		\int_{\timefunction' = \timefunction_0}^{\hat{\timefunction}} 
			\frac{1}{|\timefunction'|^{5/2}} 
			\hypersurfacecontrolwave_{[1, N-1]}(\timefunction',\hat{u}) 
			\mathfrak{l}(\timefunction')^{\mathfrak{c}}
			\mathfrak{l}^{-\mathfrak{c}}(\timefunction')
\, \mathrm{d} \timefunction' 
		\\
& \leq 
C 
|\hat{\timefunction}|^{13.6} 
\mathfrak{l}^{-\mathfrak{c}}(\hat{\timefunction})
e^{-\mathfrak{c} \hat{u}} 
\sup_{(\timefunction',u') \in [\timefunction_0,\hat{\timefunction}] \times[-\rightu,\hat{u}]} 
\left\lbrace 
	\mathfrak{l}^{-\mathfrak{c}} (\timefunction')
	\hypersurfacecontrolwave_{[1,N-1]}(\timefunction',u') 
\right\rbrace
\int_{\timefunction' = \timefunction_0}^{\hat{\timefunction}} 
	\frac{1}{|\timefunction'|^{1/2}} 
	\mathfrak{l}^{\mathfrak{c}}(\timefunction') 
\, \mathrm{d} \timefunction' 
	\\
& \leq 
	\frac{C}{\mathfrak{c}}
	|\hat{\timefunction}|^{13.6}
	e^{-\mathfrak{c}\hat{u}} 
	\sup_{(\timefunction',u') \in [\timefunction_0,\hat{\timefunction}]\times[-\rightu,\hat{u}]} 
	\left\lbrace 
		\mathfrak{l}^{-\mathfrak{c}}(\timefunction')
		\hypersurfacecontrolwave_{[1,N-1]}(\timefunction',u') 
	\right\rbrace
		\\
& \leq 
	\frac{C}{\mathfrak{c}}
	\sup_{(\timefunction',u') \in [\timefunction_0,\hat{\timefunction}]\times[-\rightu,\hat{u}]} 
	\left\lbrace 
		|\timefunction'|^{13.6}
		\mathfrak{l}^{-\mathfrak{c}}(\timefunction')
		e^{-\mathfrak{c} u'} 
		\hypersurfacecontrolwave_{[1,N-1]}(\timefunction',u') 
	\right\rbrace
	\\
	& 
	\leq 
	\frac{C}{\mathfrak{c}}\apriorilower(\hat{\timefunction},\hat{u}) 
	\leq 
	\frac{C}{\mathfrak{c}}\apriorilower(\timefunction,u),
\end{split}
\end{align}
which is $\leq \mbox{RHS~\eqref{E:MAINWAVEESTIMATESINTERMEDIATEF}}$ as desired.
This completes our proof of \eqref{E:MAINWAVEESTIMATESINTERMEDIATEF}.

The estimate \eqref{E:MAINWAVEESTIMATESINTERMEDIATEG} can be proved via
similar arguments that start with evaluating both sides of inequality \eqref{E:MAINWAVEPARTIALINTEGRALINEQUALITIES} 
at $(\hat{\timefunction},\hat{u})$, multiplying the inequality by 
$\iota_{\aprioripartial}^{-1}(\hat{\timefunction},\hat{u})$, 
and then taking $\sup_{(\hat{\timefunction},\hat{u}) \in [\timefunction_0,\timefunction]\times[-\rightu,u]}$. 
The key difference between the 
estimates \eqref{E:MAINWAVEESTIMATESINTERMEDIATEF} and \eqref{E:MAINWAVEESTIMATESINTERMEDIATEG} 
is that the critical boxed-constant-multiplied terms and $C_*$-multiplied terms appearing
on RHS~\eqref{E:TOPORDERWAVEL2CONTROLLINGINTEGRALINEQUALITY}
are \emph{absent} from RHS~\eqref{E:MAINWAVEPARTIALINTEGRALINEQUALITIES}.
Consequently, using arguments similar to the ones given above,
we find that the terms
$
	\left\lbrace 
			\frac{\frac{4 \times 1.01}{1.99} + 4.13}{15.6} 
			+ 
			\frac{\frac{8(1.01)^2}{1.99}}{15.6 \times 7.8} 
			+ 
			\frac{4.13}{7.3} 
	\right\rbrace 
	\apriorimain(\timefunction,u)
$
and 
$C \apriorimain^{1/2}(\timefunction,u) \aprioripartial^{1/2}(\timefunction,u)$ 
from RHS~\eqref{E:MAINWAVEESTIMATESINTERMEDIATEF} 
are \emph{absent} from RHS~\eqref{E:MAINWAVEESTIMATESINTERMEDIATEG}.

The estimate \eqref{E:MAINWAVEESTIMATESINTERMEDIATEH} can be proved via
similar arguments that start with evaluating both sides of inequality \eqref{E:MAINWAVEBELOWTOPINTEGRALINEQUALITIES}
at $(\hat{\timefunction},\hat{u})$, multiplying the inequality by 
$\iota_{\apriorilower}^{-1}(\hat{\timefunction},\hat{u})$, 
and then take the supremum $\sup_{(\hat{\timefunction},\hat{u}) \in [\timefunction_0,\timefunction]\times[-\rightu,u]}$. 
Using arguments similar to the ones given above,
we find that the terms $\Errorsubcriticalarg{N-1}(\timefunction,u)$ on RHS~\eqref{E:MAINWAVEBELOWTOPINTEGRALINEQUALITIES}
generate terms that can be bounded by
$\leq
C 
\initialsmall^2 
+ 
\left\lbrace 
	C \varsigma 
	+ 
	\frac{C}{\mathfrak{c}} 
	\left(1 + \varsigma^{-1}\right)
\right\rbrace
\apriorilower(\timefunction,u)
$,
which is $\leq \mbox{RHS~\eqref{E:MAINWAVEESTIMATESINTERMEDIATEH}}$ as desired.
To handle the term generated by the remaining term on RHS~\eqref{E:MAINWAVEBELOWTOPINTEGRALINEQUALITIES}
(i.e., the double time-integral involving $\hypersurfacecontrolwave_N^{1/2}$),
we can bound it using straightforward arguments based on
multiplying and dividing by $|\timefunction'|^{6.8}$ and
$|\timefunction''|^{7.8}$
in the two time integrals:
\begin{align}
\begin{split} \label{E:HANDLINGDERIVATIVELOSINGTERMSINMAINWAVEESTIMATESINTERMEDIATEH}
	&
	C 
	\iota_{\apriorilower}^{-1}(\hat{\timefunction},\hat{u}) 
	\int_{\timefunction' = \timefunction_0}^{\hat{\timefunction}} 
			\frac{1}{|\timefunction'|^{1/2}} \hypersurfacecontrolwave_{[1,N-1]}^{1/2}(\timefunction',\hat{u}) 
			\int_{\timefunction'' = \timefunction_0}^{\timefunction'} 
				\frac{1}{|\timefunction''|^{1/2}} 
				\hypersurfacecontrolwave_N^{1/2}(\timefunction'',\hat{u}) 
			\, \mathrm{d} \timefunction'' 
	\, \mathrm{d} \timefunction' 
		\\
& = C 
		|\hat{\timefunction}|^{13.6} 
		\mathfrak{l}^{-\mathfrak{c}}(\hat{\timefunction})
		e^{-\mathfrak{c} \hat{u}}
		\int_{\timefunction' = \timefunction_0}^{\hat{\timefunction}} 
			\frac{1}{|\timefunction'|^{1/2}} \hypersurfacecontrolwave_{[1,N-1]}^{1/2}(\timefunction',\hat{u}) 
			\int_{\timefunction'' = \timefunction_0}^{\timefunction'} 
				\frac{1}{|\timefunction''|^{1/2}} 
				\hypersurfacecontrolwave_N^{1/2}(\timefunction'',\hat{u}) 
			\, \mathrm{d} \timefunction'' 
	\, \mathrm{d} \timefunction' 
			\\
& 
\leq 
C 
|\hat{\timefunction}|^{13.6} 
\mathfrak{l}^{-\mathfrak{c}}(\hat{\timefunction})
e^{-\mathfrak{c} \hat{u}} 
\sup_{(\timefunction',u') \in [\timefunction_0,\hat{\timefunction}]\times[-\rightu,\hat{u}]} 
\left\lbrace 
	|\timefunction'|^{6.8} \hypersurfacecontrolwave_{[1,N-1]}^{1/2}(\timefunction',u') 
\right\rbrace 
	\\
& \ \ \ \ 
\times
\sup_{(\timefunction'',u'') \in [\timefunction_0,\hat{\timefunction}]\times[-\rightu,\hat{u}]} 
\left\lbrace 
	|\timefunction''|^{7.8} \hypersurfacecontrolwave_N^{1/2}(\timefunction'',u'') 
\right\rbrace 
	\\
& \ \ \ \  
\times
\int_{\timefunction' = \timefunction_0}^{\hat{\timefunction}} 
		\frac{1}{|\timefunction'|^{7.3}}
		\int_{\timefunction'' = \timefunction_0}^{\timefunction'} 
			\frac{1}{|\timefunction''|^{8.3}} 
		\, \mathrm{d} \timefunction'' 
\, \mathrm{d} \timefunction'  
	\\
& \leq
	C 
	\mathfrak{l}^{-\mathfrak{c}}(\hat{\timefunction})
	e^{-\mathfrak{c} \hat{u}}
	\sup_{(\timefunction',u') \in [\timefunction_0,\hat{\timefunction}]\times[-\rightu,\hat{u}]} 
	\left\lbrace 
		|\timefunction'|^{6.8} \hypersurfacecontrolwave_{[1,N-1]}^{1/2}(\timefunction',u') 
	\right\rbrace 
		\\
	& \ \ \ \ 
	\times
	\sup_{(\timefunction'',u'') \in [\timefunction_0,\hat{\timefunction}]\times[-\rightu,\hat{u}]} 
	\left\lbrace 
		|\timefunction''|^{7.8} \hypersurfacecontrolwave_N^{1/2}(\timefunction'',u'') 
	\right\rbrace 
		\\
& \leq
	C 
	\sup_{(\timefunction',u') \in [\timefunction_0,\hat{\timefunction}]\times[-\rightu,\hat{u}]} 
	\left\lbrace 
		|\timefunction'|^{6.8} 
		\mathfrak{l}^{-\frac{\mathfrak{c}}{2}}(\timefunction')
		e^{-\frac{\mathfrak{c}}{2} u'}
		\hypersurfacecontrolwave_{[1,N-1]}^{1/2}(\timefunction',u') 
	\right\rbrace 
		\\
	& \ \ \ \
	\times
	\sup_{(\timefunction'',u'') \in [\timefunction_0,\hat{\timefunction}]\times[-\rightu,\hat{u}]} 
	\left\lbrace 
		|\timefunction''|^{7.8} 
		\mathfrak{l}^{-\frac{\mathfrak{c}}{2}}(\timefunction'')
		e^{-\frac{\mathfrak{c}}{2} u''} 
		\hypersurfacecontrolwave_N^{1/2}(\timefunction'',u'') 
	\right\rbrace 
	\\
& 
\leq 
C
\apriorilower^{1/2}(\timefunction,u)
\apriorimain^{1/2}(\timefunction,u)
\leq 
\frac{1}{2}
\apriorilower(\timefunction,u)
+
C 
\apriorimain(\timefunction,u),
\end{split}
\end{align}
where to obtain the last inequality, we used Young's inequality.
We finally observe that 
$\mbox{RHS~\eqref{E:HANDLINGDERIVATIVELOSINGTERMSINMAINWAVEESTIMATESINTERMEDIATEH}}
\leq \mbox{RHS~\eqref{E:MAINWAVEESTIMATESINTERMEDIATEH}}$,
and we note that the inequality \eqref{E:HANDLINGDERIVATIVELOSINGTERMSINMAINWAVEESTIMATESINTERMEDIATEH}
is the only one that contributes the term
$\frac{1}{2}
\apriorilower(\timefunction,u)
$
to RHS~\eqref{E:MAINWAVEESTIMATESINTERMEDIATEH}.
We have therefore proved \eqref{E:HANDLINGDERIVATIVELOSINGTERMSINMAINWAVEESTIMATESINTERMEDIATEH},
which completes our proof of \eqref{E:MAINWAVEESTIMATESINTERMEDIATEF}--\eqref{E:MAINWAVEESTIMATESINTERMEDIATEH}.
This also yields the desired bounds for
$\totalcontrolwave_{[1,\Ntop]}$
and
$\totalcontrolwave_{[1,\Ntop-1]}$,
aside from the issue that, as we noted below
\eqref{E:MAINWAVEESTIMATESINTERMEDIATE1},
the final choice of $\mathfrak{c}$ is not made
until the end of the proof.

\medskip

\noindent \textbf{Estimates for}
$\totalcontrolwave_{[1,\Ntop-2]}$, $\totalcontrolwave_{[1,\Ntop-3]}$, $\cdots$, $\totalcontrolwave_1$.
We now explain how to derive the a priori estimates 
\eqref{E:MAINWAVEENERGYESTIMATESBLOWUP}--\eqref{E:MAINWAVEENERGYESTIMATESREGULAR}
for 
$\totalcontrolwave_{[1,\Ntop-2]}$, $\totalcontrolwave_{[1,\Ntop-3]}$, $\cdots$, 
$\totalcontrolwave_1$ via downward induction, starting with $\totalcontrolwave_{[1,\Ntop-2]}$. 

Unlike our analysis of the strongly coupled triple 
$\totalcontrolwave_{[1,\Ntop]}$,
$\totalcontrolwavepartial_{[1,\Ntop]}$, 
and $\totalcontrolwave_{[1,\Ntop-1]}$, we can derive the estimate for $\totalcontrolwave_{[1,\Ntop-2]}$ using only 
the integral inequality \eqref{E:MAINWAVEBELOWTOPINTEGRALINEQUALITIES}, 
the already proven vorticity and entropy estimates 
provided by Prop.\,\ref{P:MAINHYPERSURFACEENERGYESTIMATESFORTRANSPORTVARIABLES}, 
and our already proven bounds for $\totalcontrolwave_{[1,\Ntop-1]}$
(more precisely, the already proven \eqref{E:MAINWAVEESTIMATESINTERMEDIATE1} for $\apriorilower$).
To begin, we define an analog of \eqref{E:GRONWALLMULTIPLICATIONFACTORBELOWTOP}: 
$\iota_{\widetilde{\apriorilower}}(\timefunction,u) 
\eqdef 
|\timefunction|^{-11.6} 
\mathfrak{l}^{\mathfrak{c}} (\timefunction)
e^{\mathfrak{c} u}$, 
as well as an analog of \eqref{E:RENORMALIZEDBELOWTOPWAVEENERGY}: 
$\widetilde{\apriorilower}(\timefunction,u) \eqdef  \sup_{(\timefunction',u') \in [\timefunction_0,\timefunction]\times[-\rightu,u]} 
\left\lbrace \iota_{\widetilde{H}}^{-1}(\timefunction',u') \hypersurfacecontrolwave_{[1,\Ntop - 2]} \right\rbrace$.  
Note that compared to our definition \eqref{E:GRONWALLMULTIPLICATIONFACTORBELOWTOP} of $\iota_{\apriorilower}$,
we have reduced the power of $|\timefunction|^{-1}$ by two in our definition of $\iota_{\widetilde{\apriorilower}}$. 
As before, we will prove \eqref{E:MAINWAVEENERGYESTIMATESBLOWUP} for $K = 2$ by showing that 
there is a uniform $C > 0$ such that for all sufficiently small $\varsigma \in (0,1]$ and $\fundbootsmall \geq 0$, 
we have:
\begin{align}	\label{E:TWOBELOWTOPORDERMAINWAVEESTIMATESWEIGHTEDQUANTITYBOUND}
	\widetilde{\apriorilower}(\timefunction,u) 
	& \leq 
	C \left(1 + \varsigma^{-1} \right) \initialsmall^2.
\end{align}
To prove \eqref{E:TWOBELOWTOPORDERMAINWAVEESTIMATESWEIGHTEDQUANTITYBOUND}, we will show that:
\begin{align} \label{E:TWOBELOWTOPMAINWAVEENERGYESTIMATESBLOWUP}
\widetilde{\apriorilower}(\timefunction,u) 
	& 
	\leq
	C \left(1 + \varsigma^{-1} \right) \initialsmall^2 
	+
	\left\lbrace 
	\frac{1}{2} 
	+ 
	C \varsigma 
	+ 
	\frac{C}{\mathfrak{c}} 
	\left(1 + \varsigma^{-1}\right)
	\right\rbrace
	\widetilde{\apriorilower}(\timefunction,u),
\end{align}
where as before, we set $\mathfrak{c} \eqdef \varsigma^{-2}$.
Once we have proved \eqref{E:TWOBELOWTOPMAINWAVEENERGYESTIMATESBLOWUP},
then if $\varsigma \in (0,1]$ and $\fundbootsmall \geq 0$ are sufficiently small,
we can absorb all terms on RHS~\eqref{E:TWOBELOWTOPMAINWAVEENERGYESTIMATESBLOWUP}
except for $C \left(1 + \varsigma^{-1} \right) \initialsmall^2$
back into the left, thereby arriving 
at the desired bound \eqref{E:TWOBELOWTOPORDERMAINWAVEESTIMATESWEIGHTEDQUANTITYBOUND}.


It remains for us to prove \eqref{E:TWOBELOWTOPMAINWAVEENERGYESTIMATESBLOWUP}.
To proceed, we set $N = \Ntop -1$, multiply both sides of \eqref{E:MAINWAVEBELOWTOPINTEGRALINEQUALITIES} by 
$\iota_{\widetilde{\apriorilower}}^{-1}$, and evaluate the resulting expression at $(\hat{\timefunction},\hat{u}) \in [\timefunction_0,\timefunction]\times[-\rightu,u]$. 
Using same arguments we used to prove \eqref{E:MAINWAVEESTIMATESINTERMEDIATEH}
(in particular using the already proven estimates of Prop.\,\ref{P:MAINHYPERSURFACEENERGYESTIMATESFORTRANSPORTVARIABLES}),
we find that the terms $\Errorsubcriticalarg{N-1}(\timefunction,u)$ on RHS~\eqref{E:MAINWAVEBELOWTOPINTEGRALINEQUALITIES}
generate terms that we can bound by
$
\leq
	C \initialsmall^2 
	+
	\left\lbrace 
	C \varsigma 
	+ 
	\frac{C}{\mathfrak{c}} 
	\left(1 + \varsigma^{-1}\right)
	\right\rbrace
	\widetilde{\apriorilower}(\timefunction,u)
$. 
We now handle the remaining term, i.e., 
the term generated by the first term on RHS~\eqref{E:MAINWAVEBELOWTOPINTEGRALINEQUALITIES},
which is a double time-integral involving 
the above-present-order factor $\hypersurfacecontrolwave_{\Ntop-1}^{1/2}$.
Multiplying and dividing by $|\timefunction'|^{5.8}$ and
$|\timefunction''|^{6.8}$
in the two time integrals,
arguing as in the proof of \eqref{E:HANDLINGDERIVATIVELOSINGTERMSINMAINWAVEESTIMATESINTERMEDIATEH},
and using the already proven bound \eqref{E:MAINWAVEESTIMATESINTERMEDIATE1}
for $\apriorilower$,
we bound this term as follows:
\begin{align}
\begin{split} \label{E:TWOBELOWTOPHANDLINGDERIVATIVELOSINGTERMSINMAINWAVEESTIMATESINTERMEDIATEH}
	&
	C 
	\iota_{\widetilde{\apriorilower}}^{-1}(\hat{\timefunction},\hat{u}) 
	\int_{\timefunction' = \timefunction_0}^{\hat{\timefunction}} 
			\frac{1}{|\timefunction'|^{1/2}} \hypersurfacecontrolwave_{[1,\Ntop-2]}^{1/2}(\timefunction',\hat{u}) 
			\int_{\timefunction'' = \timefunction_0}^{\timefunction'} 
				\frac{1}{|\timefunction''|^{1/2}} 
				\hypersurfacecontrolwave_{\Ntop - 1}^{1/2}(\timefunction'',\hat{u}) 
			\, \mathrm{d} \timefunction'' 
	\, \mathrm{d} \timefunction' 
		\\
& = C 
		|\hat{\timefunction}|^{11.6} 
		\mathfrak{l}^{-\mathfrak{c}}(\hat{\timefunction})
		e^{-\mathfrak{c} \hat{u}}
		\int_{\timefunction' = \timefunction_0}^{\hat{\timefunction}} 
			\frac{1}{|\timefunction'|^{1/2}} \hypersurfacecontrolwave_{[1,\Ntop-2]}^{1/2}(\timefunction',\hat{u}) 
			\int_{\timefunction'' = \timefunction_0}^{\timefunction'} 
				\frac{1}{|\timefunction''|^{1/2}} 
				\hypersurfacecontrolwave_{\Ntop - 1}^{1/2}(\timefunction'',\hat{u}) 
			\, \mathrm{d} \timefunction'' 
	\, \mathrm{d} \timefunction'
			\\
& 
\leq 
C 
|\hat{\timefunction}|^{11.6} 
\mathfrak{l}^{-\mathfrak{c}}(\hat{\timefunction})
e^{-\mathfrak{c} \hat{u}} 
	\\
& \ \ \ \
\times
\sup_{(\timefunction',u') \in [\timefunction_0,\hat{\timefunction}]\times[-\rightu,\hat{u}]} 
\left\lbrace 
	|\timefunction'|^{5.8} \hypersurfacecontrolwave_{[1,\Ntop-2]}^{1/2}(\timefunction',u') 
\right\rbrace 
	\\
& \ \ \ \
\times
\sup_{(\timefunction'',u'') \in [\timefunction_0,\hat{\timefunction}]\times[-\rightu,\hat{u}]} 
\left\lbrace 
	|\timefunction''|^{6.8} \hypersurfacecontrolwave_{\Ntop-1}^{1/2}(\timefunction'',u'') 
\right\rbrace 
	\\
& \ \ \ \ 
\times
\int_{\timefunction' = \timefunction_0}^{\hat{\timefunction}} 
		\frac{1}{|\timefunction'|^{6.3}}
		\int_{\timefunction'' = \timefunction_0}^{\timefunction'} 
			\frac{1}{|\timefunction''|^{7.3}} 
		\, \mathrm{d} \timefunction'' 
\, \mathrm{d} \timefunction'  
	\\
& \leq
	C 
	\mathfrak{l}^{-\mathfrak{c}}(\hat{\timefunction})
	e^{-\mathfrak{c} \hat{u}}
	\sup_{(\timefunction',u') \in [\timefunction_0,\hat{\timefunction}]\times[-\rightu,\hat{u}]} 
	\left\lbrace 
		|\timefunction'|^{5.8} \hypersurfacecontrolwave_{[1,\Ntop-2]}^{1/2}(\timefunction',u') 
	\right\rbrace 
		\\
	& \ \ \ \
	\times
	\sup_{(\timefunction'',u'') \in [\timefunction_0,\hat{\timefunction}]\times[-\rightu,\hat{u}]} 
	\left\lbrace 
		|\timefunction''|^{6.8} \hypersurfacecontrolwave_{\Ntop-1}^{1/2}(\timefunction'',u'') 
	\right\rbrace 
	\\
& \leq
	C 
	\sup_{(\timefunction',u') \in [\timefunction_0,\hat{\timefunction}]\times[-\rightu,\hat{u}]} 
	\left\lbrace 
		|\timefunction'|^{5.8} 
		\mathfrak{l}^{-\frac{\mathfrak{c}}{2}}(\timefunction')
		e^{- \frac{\mathfrak{c}}{2} u'}
		\hypersurfacecontrolwave_{[1,\Ntop-2]}^{1/2}(\timefunction',u') 
	\right\rbrace 
		\\
& \ \ \ \
	\times
	\sup_{(\timefunction'',u'') \in [\timefunction_0,\hat{\timefunction}]\times[-\rightu,\hat{u}]} 
	\left\lbrace 
		|\timefunction''|^{6.8} 
		\mathfrak{l}^{-\frac{\mathfrak{c}}{2}}(\timefunction'')
		e^{-\frac{\mathfrak{c}}{2} u''} 
		\hypersurfacecontrolwave_{\Ntop-1}^{1/2}(\timefunction'',u'') 
	\right\rbrace 
	\\
& 
\leq 
C
\widetilde{\apriorilower}^{1/2}(\hat{\timefunction},\hat{u})
\apriorilower^{1/2}(\hat{\timefunction},\hat{u}) 
\leq 
\frac{1}{2}
\widetilde{\apriorilower}(\timefunction,u)
+
C 
\apriorilower(\timefunction,u)
	\\
&
\leq
\frac{1}{2}
\widetilde{\apriorilower}(\timefunction,u)
+
C \left(1 + \varsigma^{-1} \right) \initialsmall^2,
\end{split}
\end{align}
which is $\leq \mbox{RHS~\eqref{E:TWOBELOWTOPMAINWAVEENERGYESTIMATESBLOWUP}}$ as desired.
We have therefore proved \eqref{E:TWOBELOWTOPORDERMAINWAVEESTIMATESWEIGHTEDQUANTITYBOUND},
which in particular implies
\eqref{E:MAINWAVEENERGYESTIMATESBLOWUP} for $K = 2$.

The desired bounds \eqref{E:MAINWAVEENERGYESTIMATESBLOWUP}--\eqref{E:MAINWAVEENERGYESTIMATESREGULAR} 
for $\totalcontrolwave_{[1,\Ntop-3]}$, $\cdots$, $\totalcontrolwave_1$ 
can be derived by downward induction
based on an argument that is very similar to the one we used to prove the bound
\eqref{E:TWOBELOWTOPORDERMAINWAVEESTIMATESWEIGHTEDQUANTITYBOUND} for $\totalcontrolwave_{[1,\Ntop-2]}$. 
The only difference is that we define an analogous multiplicative factor 
$\iota_{\widetilde{\apriorilower},P}(\timefunction,u) \eqdef |\timefunction|^{-P} \mathfrak{l}^{\mathfrak{c}} (\timefunction)e^{\mathfrak{c} u}$, where $P= 9.6$ for $\totalcontrolwave_{[1,\Ntop-3]}$, 
$P = 7.6$ for $\totalcontrolwave_{[1,\Ntop-4]}$,  
$P = 5.6$ for $\totalcontrolwave_{[1,\Ntop-5]}$,  
$P = 3.6$ for $\totalcontrolwave_{[1,\Ntop-6]}$, 
$P = 1.6$ for $\totalcontrolwave_{[1,\Ntop-7]}$, 
and $P = 0$ for 
$\totalcontrolwave_{[1,\Ntop-8]}$. 
We stress that these latter estimates (i.e., \eqref{E:MAINWAVEENERGYESTIMATESREGULAR}) 
\emph{do not involve any singular factor of} $|\timefunction|^{-1}$. 

\hfill $\qed$

\subsection{Proof of Prop.\,\ref{P:APRIORIL2ESTIMATESACOUSTICGEOMETRY}}
\label{SS:PROOFOFENERGYESTIMATESFORACOUSTICGEOMETRY}
In this short section, we prove Prop.\,\ref{P:APRIORIL2ESTIMATESACOUSTICGEOMETRY},
which yields our energy estimates for the acoustic geometry on the rough foliations,
thereby completing our proof of the energy estimates.

To proceed, we note that \eqref{E:MINVALUEOFMUONFOLIATION} implies 
$\| f \|_{L^2\left(\hypthreearg{\timefunction}{[-\rightu,\leftu]}{\muxmulevelsetvalue} \right)} \leq 
|\timefunction|^{-1} 
\| \upmu f \|_{L^2\left(\hypthreearg{\timefunction}{[-\rightu,\leftu]}{\muxmulevelsetvalue} \right)}$. 
Hence, the desired bound \eqref{E:MAINL2CHITOPORDERBLOWUP} follows from inserting the 
already proven estimates of
Props.\,\ref{P:APRIORIL2ESTIMATESWAVEVARIABLES} and
\ref{P:MAINHYPERSURFACEENERGYESTIMATESFORTRANSPORTVARIABLES} 
for 
$\totalcontrolwave_M$, 
$\hypersurfacecontrolVortVort_M$, 
$\hypersurfacecontrolDivGradEnt_M$, 
$\hypersurfacecontrolVort_M$ 
and
$\hypersurfacecontrolGradEnt_M$
(for the relevant values of $M$)
into RHS~\eqref{E:TOPORDERL2ESTIMATEMUCHI} 
and integrating in $\timefunction$. 
Similarly, \eqref{E:MAINACOUSTGEOMETRYBELOWTOPORDERBLOWUP}--\eqref{E:MAINACOUSTGEOMETRYBELOWTOPORDERREGULAR} follow from inserting the estimates \eqref{E:MAINWAVEENERGYESTIMATESBLOWUP}--\eqref{E:MAINWAVEENERGYESTIMATESREGULAR} into RHS~\eqref{E:PRELIMINARYEIKONALWITHOUTL}. 

\hfill $\qed$

\section{Improvements of the fundamental quantitative $L^{\infty}$ bootstrap assumptions} 
\label{S:IMPROVEMENTSOFFUNDAMENTALQUANTITATIVEBOOTSTRAPASSUMPTIONS}
In this short section, 
we derive $L^{\infty}$ estimates that yield an improvement of the fundamental quantitative bootstrap assumptions
stated in Sect.\,\ref{SSS:FUNDAMENTALQUANTITATIVE}.  
The results follow easily from the 
Sobolev embedding estimate \eqref{E:H2LINFINITYFUNDAMENTALTHEOREMOFCALCULUSPLUSSOBOLEVEMBEDDINGONROUGHTORUS}
and the non-degenerate energy estimates we have already derived in Prop.\,\ref{P:APRIORIL2ESTIMATESWAVEVARIABLES}.

\begin{proposition}[Improvement of the fundamental quantitative bootstrap assumptions]
\label{P:IMPROVEMENTOFFUNDAMENTALQUANTITATIVEBOOTSTRAPASSUMPTIONS}
		Under the parameter-size assumptions of Sect.\,\ref{SS:PARAMETERSIZEASSUMPTIONS}, 
		the initial data assumptions of 
		Sects.\,\ref{SSS:QUANTITATIVEASSUMPTIONSONDATAAWAYFROMSYMMETRY}--\ref{SSS:LOCALIZEDDATAASSUMPTIONSFORMUANDDERIVATIVES},
		and the bootstrap assumptions of 
		Sects.\,\ref{S:BOOTSTRAPEVERYTHINGEXCEPTENERGIES} and \ref{SS:BOOTSTRAPASSUMPTIONSFORTHEWAVEENERGIES},
		there exists a constant $C > 0$ such that the following estimates hold for 
		$(\timefunction,u) \in [\timefunction_0,\timefunctionboot) \times [- \rightu,\leftu]$:
\begin{align} \label{E:RECOVEREDLINFINITYESTIMATE}
	\left\| 
		\tander^{[1,\Ntop - 10]} \wavearray 
	\right\|_{L^{\infty}\left(\twoargroughtori{\timefunction,u}{\muxmulevelsetvalue}\right)}, 
		\,
	\left\| \tander^{\leq \Ntop - 11} 
		(\vortrenormalized,\GradEnt) 
	\right\|_{L^{\infty}\left(\twoargroughtori{\timefunction,u}{\muxmulevelsetvalue}\right)},
		\,
	\left\| \tander^{\leq \Ntop - 12} 
		(\VortVort,\DivGradEnt) 
	\right\|_{L^{\infty}\left(\twoargroughtori{\timefunction,u}{\muxmulevelsetvalue}\right)}
	& \leq C \initialsmall.
\end{align}

In particular, if $\initialsmall$ is small enough such that $C \initialsmall < \fundbootsmall$,
then \eqref{E:RECOVEREDLINFINITYESTIMATE} yields a strict improvement of the
fundamental quantitative bootstrap assumptions
\eqref{E:FUNDAMENTALQUANTITATIVEBOOTWAVEANDTRANSPORT}--\eqref{E:FUNDAMENTALQUANTITATIVEBOOTMODIFIEDFLUIDVARS}.

\end{proposition}

\begin{proof}
The estimate \eqref{E:RECOVEREDLINFINITYESTIMATE} for
$
\left\| \tander^{\leq \Ntop - 11} 
		(\vortrenormalized,\GradEnt) 
	\right\|_{L^{\infty}\left(\twoargroughtori{\timefunction,u}{\muxmulevelsetvalue}\right)}
$
follows from 
the Sobolev embedding estimate \eqref{E:H2LINFINITYSOBOLEVEMBEDDINGROUGHTORUS},
the definitions
\eqref{E:TORIVORTICITYL2CONTROLLINGQUANTITY}--\eqref{E:TORIENTROPYGRADIENTL2CONTROLLINGQUANTITY}
of $\toricontrolVort_N(\timefunction,u)$ and $\toricontrolGradEnt_N(\timefunction,u)$,
and the rough tori energy estimates \eqref{E:MAINTORIL2VORTGRADENTBELOWTOPORDERREGULAR}.
Similarly, the estimate \eqref{E:RECOVEREDLINFINITYESTIMATE}
for
$
\left\| \tander^{\leq \Ntop - 12} 
		(\VortVort,\DivGradEnt) 
\right\|_{L^{\infty}\left(\twoargroughtori{\timefunction,u}{\muxmulevelsetvalue}\right)}
$
follows from 
the Sobolev embedding estimate \eqref{E:H2LINFINITYSOBOLEVEMBEDDINGROUGHTORUS},
the definitions
\eqref{E:TORIMODIFIEDVORTICITYL2CONTROLLINGQUANTITY}--\eqref{E:TORIDIVGRADENTL2CONTROLLINGQUANTITY}
of $\toricontrolVortVort_N(\timefunction,u)$ and $\toricontrolDivGradEnt_N(\timefunction,u)$,
and the rough tori energy estimates \eqref{E:MAINTORIL2MODIFIEDFLUIDVARIABLESBELOWTOPORDERREGULAR}.

To prove the estimate \eqref{E:RECOVEREDLINFINITYESTIMATE} for
$	
	\left\| 
		\tander^{[1,\Ntop - 10]} \wavearray 
	\right\|_{L^{\infty}\left(\twoargroughtori{\timefunction,u}{\muxmulevelsetvalue}\right)}
$,
we use the estimate \eqref{E:H2LINFINITYFUNDAMENTALTHEOREMOFCALCULUSPLUSSOBOLEVEMBEDDINGONROUGHTORUS}, 
the data-assumption \eqref{E:SMALLDATAOFPSIONINITIALROUGHTORI}, 
Lemma~\ref{L:COERCIVENESSOFL2CONTROLLINGQUANITIES},
and the estimate \eqref{E:MAINWAVEENERGYESTIMATESREGULAR}
to conclude that
$\left \|\tander^{[1,\Ntop-10]} \wavearray \right\|_{L^{\infty}\left(\twoargroughtori{\timefunction,u}{\muxmulevelsetvalue}\right)}^2 
\lesssim 
\initialsmall^2 
+ 
\totalcontrolwave_{[1,\Ntop-8]}(\timefunction,u) 
\lesssim 
\initialsmall^2$
as desired.

\end{proof}

\section{Existence up to the singular boundary at fixed $\muxmulevelsetvalue$ via continuation criteria}
\label{S:EXISTENCEUPTOSINGULARBOUNDARYATFIXEDKAPPA}
In this section, we prove our first main theorem,
Theorem~\ref{T:EXISTENCEUPTOTHESINGULARBOUNDARYATFIXEDKAPPA},
which shows that at fixed $\muxmulevelsetvalue \in [0,\muxmulevelsetvalue_0]$, 
the solution exists on the domain $\twoargMrough{[\timefunction_0,0],[- \rightu,\leftu]}{\muxmulevelsetvalue}$,
which in particular contains the torus $\twoargmumuxtorus{0}{-\muxmulevelsetvalue}$,
a subset of the singular boundary (see Sect.\,\ref{SS:SINGULARBOUNDARYANDCREASE}).
To prove the theorem, we rely mainly on results that we have already established,
though we also rely on standard continuation criteria,
which we prove independently in Prop.\,\ref{P:CONTINUATIONCRITERIA}.
Roughly, given the results we have already proved, 
the continuation criteria allow us to continue the solution classically
as long as $\upmu$ has not vanished.
We highlight that our setup guarantees 
(see in particular \eqref{E:MINVALUEOFMUONFOLIATION})
that within $\twoargMrough{[\timefunction_0,0],[- \rightu,\leftu]}{\muxmulevelsetvalue}$, 
the vanishing of $\upmu$ can happen only
along the top boundary $\hypthreearg{0}{[- \rightu,\leftu]}{\muxmulevelsetvalue}$, 
i.e., when $\timefunctionarg{\muxmulevelsetvalue} = 0$.

\subsection{Existence up to the singular boundary torus $\twoargmumuxtorus{0}{-\muxmulevelsetvalue}$}
\label{S:EXISTENCEUPTOTHESINGULARBOUNDARYATFIXEDKAPPA}
In this section, we state and prove our first main theorem, which concerns fixed $\muxmulevelsetvalue \in [0,\muxmulevelsetvalue_0]$.

\begin{theorem}[Existence in a region containing the singular boundary torus $\twoargmumuxtorus{0}{-\muxmulevelsetvalue}$]
	\label{T:EXISTENCEUPTOTHESINGULARBOUNDARYATFIXEDKAPPA}
	Fix any of the compactly supported admissible simple isentropic plane symmetric ``background'' solutions $\RRiemannPS$
	from Def.\,\ref{AD:ADMISSIBLEBACKGROUND} (recall that 
	$\LRiemann$, $v^2$, $v^3$, $\Ent$, 
	$\vortrenormalized$,
	$\GradEnt$
	$\VortVort$
	and
	$\DivGradEnt$ vanish for these background solutions).
	
Let
$(\RRiemann,\LRiemann,v^2,v^3,\Ent) \big|_{\Sigma_0}
\eqdef 
\left(\RRiemannpertinitial, \LRiemannpertinitial, \vtwopertinitial, \vthreepertinitial, \spertinitial \right)
$
be perturbed fluid data on the flat Cartesian hypersurface $\Sigma_0$,
as in \eqref{E:PERTURBEDBONAFIDEDATA},
and let $u|_{\Sigma_0} = - x^1$ be the initial condition of the eikonal function, 
as in \eqref{E:EIKONALEQUATION} and \eqref{AE:DATAFOREIKONALEQUATIONINPLANESYMMETRY}.
Let
$
\left(\vortrenormalizedpertinitial^i,
\GradEntpertinitial^i,
\VortVortpertinitial^i,
\DivGradEntpertinitial 
\right)_{i=1,2,3}
$
respectively denote the initial data on $\Sigma_0$
of 
$
\left(
\vortrenormalized^i,
\GradEnt^i,
\VortVort^i,
\DivGradEnt
\right)_{i=1,2,3}
$. Note that these data are determined by
$\left(\RRiemannpertinitial, \LRiemannpertinitial, \vtwopertinitial, \vthreepertinitial, \spertinitial \right)$,
the compressible Euler equations
\eqref{E:BVIEVOLUTION}--\eqref{E:BENTROPYEVOLUTION},
definition~\eqref{E:ALMOSTRIEMANNINVARIANTS},
and Def.\,\ref{D:HIGHERORDERFLUIDVARIABLES}.
Also recall that these data and the data of the eikonal function
determine the data of all the acoustic geometry on $\Sigma_0$; 
see Remark~\ref{R:DATAOFEIKONAFUNCTIONQUANTITIES}.

Assume the following:
	\begin{itemize}
		\item $\Ntop \geq 24$.
		\item The quantity
			$
			\mathring{\Delta}_{\Sigma_0^{[-\farrightu,\leftu]}}^{\Ntop+1}
			$
			defined in \eqref{E:PERTURBATIONSMALLNESSINCARTESIANDIFFERENTIALSTRUCTURE}
			is sufficiently small, where
			$
			\mathring{\Delta}_{\Sigma_0^{[-\farrightu,\leftu]}}^{\Ntop+1}
			$
			is a Sobolev norm of the perturbation of the fluid data away from the background solution,
			where $\farrightu > 0$ and $\leftu > 0$ are parameters from Sect.\,\ref{SS:PARAMETERSIZEASSUMPTIONS}.
		\item Recall that in Appendix~\ref{A:OPENSETOFDATAEXISTS} (see Prop.\,\ref{P:CAUCHYSTABILITYANDEXISTENCEOFOPENSETS}), 
			we showed that the smallness of 
			$
			\mathring{\Delta}_{\Sigma_0^{[-\farrightu,\leftu]}}^{\Ntop+1}
			$ 
			implies that the parameter-size assumptions of Sect.\,\ref{SS:PARAMETERSIZEASSUMPTIONS} hold
			and that the fluid variable and acoustic geometry data induced on
			$\twoargroughtori{\timefunction_0,u}{\muxmulevelsetvalue}$,
			$\hypthreearg{\timefunction_0}{[- \rightu,\leftu]}{\muxmulevelsetvalue}$,
			and $\nullhyparg{- \rightu}^{[0,\frac{4}{\mathring{\updelta}_*}]}$
			satisfy the assumptions stated in 
			Sects.\,\ref{SSS:QUANTITATIVEASSUMPTIONSONDATAAWAYFROMSYMMETRY}--\ref{SSS:LOCALIZEDDATAASSUMPTIONSFORMUANDDERIVATIVES},
			where $\rightu > 0$ and $\mathring{\updelta}_* > 0$ are parameters from Sect.\,\ref{SS:PARAMETERSIZEASSUMPTIONS}.
		\item In particular, Prop.\,\ref{P:CAUCHYSTABILITYANDEXISTENCEOFOPENSETS} implies that the
				parameter $\initialsmall$ from Sect.\,\ref{SS:PARAMETERSIZEASSUMPTIONS} can be chosen to satisfy
				$\initialsmall = \mathcal{O}(\mathring{\Delta}_{\Sigma_0^{[-\farrightu,\leftu]}}^{\Ntop+1})$,
				where the implicit constants in ``$\mathcal{O}(\cdot)$'' depend on the background solution.
		\item $\muxmulevelsetvalue \in [0,\muxmulevelsetvalue_0]$, where $\muxmulevelsetvalue_0$ is the parameter from Sect.\,\ref{SS:PARAMETERSIZEASSUMPTIONS}.
	\end{itemize}
	Then the following conclusions hold.
	
	\medskip
	
\noindent \underline{\textbf{The rough time function and classical existence relative to the geometric coordinates}}.

\begin{itemize}
	\item There exists a rough time function 
	$\timefunctionarg{\muxmulevelsetvalue} = \timefunctionarg{\muxmulevelsetvalue}(t,u,x^2,x^3)$ 
	(constructed in Sect.\,\ref{S:ROUGHTIMEFUNCTIONANDROUGHSUBSETS})
		with range $[\timefunction_0,0] = [-\mupositive,0]$ (recall that $\timefunction_0 = - \mupositive$),
		and we denote its level set portions,
		viewed as subsets of geometric coordinate space $\mathbb{R} \times \mathbb{R} \times \mathbb{T}^2$,
		as follows:
		$\hypthreearg{\timefunction}{[u_1,u_2]}{\muxmulevelsetvalue}
		= 
		\lbrace 
			(t,u,x^2,x^3) 
			\ | \
			 \timefunctionarg{\muxmulevelsetvalue}(t,u,x^2,x^3) = \timefunction,
				\,
			 u_1 \leq u \leq u_2,
				\,
			(x^2,x^3) \in \mathbb{T}^2
		\rbrace
	$.
		More precisely,
		$\timefunctionarg{\muxmulevelsetvalue}$ is
		defined on the portion 
		$\twoargMrough{[\timefunction_0,0],[- \rightu,\leftu]}{\muxmulevelsetvalue}
			= \bigcup_{\timefunction \in [\timefunction_0,0]} \hypthreearg{\timefunction}{[- \rightu,\leftu]}{\muxmulevelsetvalue}
		$
		of the maximal classical development 
		of the data with respect to the differential structure of the geometric coordinates $(t,u,x^2,x^3)$.
	\item The change of variables map 
		$\CHOVgeotorough{\muxmulevelsetvalue}(t,u,x^2,x^3) = (\timefunctionarg{\muxmulevelsetvalue},u,x^2,x^3)$
		is a diffeomorphism from $\twoargMrough{[\timefunction_0,0],[- \rightu,\leftu]}{\muxmulevelsetvalue}$
		onto its image $[-\mupositive,0] \times [- \rightu,\leftu] \times \mathbb{T}^2$
		satisfying $\| \CHOVgeotorough{\muxmulevelsetvalue} \|_{C_{\textnormal{geo}}^{2,1}(\twoargMrough{[\timefunction_0,0],[- \rightu,\leftu]}{\muxmulevelsetvalue})} \leq C$.
		Moreover, $\geop{t} \timefunctionarg{\muxmulevelsetvalue} \approx 1$ on $\twoargMrough{[\timefunction_0,0],[- \rightu,\leftu]}{\muxmulevelsetvalue}$.
	\item The fluid variables
	$\wavearray$,
	$\vortrenormalized^i$,
	$\GradEnt^i$,
	$\VortVort^i$,
	and
	$\DivGradEnt$, 
	the eikonal function $u$,
	$\upmu$, $\Lunit^i$,
	the Cartesian coordinate functions
	$(t,x^1,x^2,x^3)$ 
	and all of the auxiliary quantities constructed out of these quantities
	exist classically
	with respect to the geometric coordinates $(t,u,x^2,x^3)$
	on $\twoargMrough{[\timefunction_0,0],[- \rightu,\leftu]}{\muxmulevelsetvalue}$.
	In particular, with respect to the geometric coordinates,
	the compressible Euler equations \eqref{E:BVIEVOLUTION}--\eqref{E:BENTROPYEVOLUTION} are satisfied,
	and the equations of Theorem~\ref{T:GEOMETRICWAVETRANSPORTSYSTEM} are also satisfied.
	\item The H\"{o}lder estimates of Lemma~\ref{L:CONTINUOUSEXTNESION},
				the $L^{\infty}$ estimates of Prop.\,\ref{P:IMPROVEMENTOFAUXILIARYBOOTSTRAP}
				with $\fundbootsmall$ replaced by $C \initialsmall$,
				and the energy estimates of 
				Props.\,\ref{P:APRIORIL2ESTIMATESWAVEVARIABLES},
				\ref{P:MAINHYPERSURFACEENERGYESTIMATESFORTRANSPORTVARIABLES},
				\ref{P:ROUGHTORIENERGYESTIMATES},
				and \ref{P:APRIORIL2ESTIMATESACOUSTICGEOMETRY}
				hold with $\timefunctionboot = 0$, i.e., they
				hold on $\twoargMrough{[\timefunction_0,0],[- \rightu,\leftu]}{\muxmulevelsetvalue}$.
	\end{itemize}

\medskip

\noindent \underline{\textbf{The behavior of $\upmu$ and properties of $\Upsilon$}}.

\begin{itemize}
\item For $\timefunction \in [-\mupositive,0]$,
		 we have: 
		\begin{align} \label{E:MAINRESULTSMINVALUEOFMUONFOLIATION}
			\min_{\hypthreearg{\timefunction}{[- \rightu,\leftu]}{\muxmulevelsetvalue}} \upmu
			& = - \timefunction.
		\end{align}
		Moreover, within $\hypthreearg{\timefunction}{[- \rightu,\leftu]}{\muxmulevelsetvalue}$,
		the minimum value $- \timefunction$
		in \eqref{E:MAINRESULTSMINVALUEOFMUONFOLIATION} is achieved by $\upmu$ precisely on the
		set $\twoargmumuxtorus{-\timefunction}{-\muxmulevelsetvalue}$ from definition~\ref{E:MUXMUTORI},
		which is a $C^{1,1}$-embedded torus.
		In particular, in $\twoargMrough{[\timefunction_0,0],[- \rightu,\leftu]}{\muxmulevelsetvalue}$,
		$\upmu$ vanishes precisely along the torus $\twoargmumuxtorus{0}{-\muxmulevelsetvalue}$,
		which is a subset of $\hypthreearg{0}{[- \rightu,\leftu]}{\muxmulevelsetvalue}$
		that is contained in the singular boundary (see Sect.\,\ref{SS:SINGULARBOUNDARYANDCREASE}).
\item  On $\twoargMrough{[\timefunction_0,0],[- \rightu,\leftu]}{\muxmulevelsetvalue}$,
		the change of variables map 
		$\Upsilon(t,u,x^2,x^3) = (t,x^1,x^2,x^3)$ is an injection onto its image in Cartesian coordinate space
		satisfying $\| \Upsilon \|_{C_{\textnormal{geo}}^{3,1}(\twoargMrough{[\timefunction_0,0],[- \rightu,\leftu]}{\muxmulevelsetvalue})} \leq C$.
		In particular, $\Upsilon$ is a homeomorphism from the compact set
		$\twoargMrough{[\timefunction_0,0],[- \rightu,\leftu]}{\muxmulevelsetvalue}$ onto its image.
	\item With $d_{\textnormal{geo}} \Upsilon$ denoting the Jacobian matrix of $\Upsilon$,
		we have:
		\begin{align} \label{E:FIXEDKAPPAMAINTHEOREMGEOTOCARTESIANJACOBIANDETERMINANTESTIMATE}
			\mbox{\upshape det} d_{\textnormal{geo}} \Upsilon
			& \approx
				- 
				\upmu.
		\end{align}
		Hence, on $\twoargMrough{[\timefunction_0,0],[- \rightu,\leftu]}{\muxmulevelsetvalue} \backslash \twoargmumuxtorus{0}{-\muxmulevelsetvalue}$,
		$\Upsilon$ is a diffeomorphism.
	\item For $\mulevelsetvalue \in [0,\mupositive]$, 
			$\Upsilon\left(\twoargmumuxtorus{\mulevelsetvalue}{-\muxmulevelsetvalue} \right)$
			is an embedded two-dimensional $C^{1,1}$ torus in Cartesian coordinate space.
			In particular, the restriction of $\Upsilon$ to $\twoargmumuxtorus{\mulevelsetvalue}{-\muxmulevelsetvalue}$ 
			is a diffeomorphism from $\twoargmumuxtorus{\mulevelsetvalue}{-\muxmulevelsetvalue}$ onto its image
			$\Upsilon\left(\twoargmumuxtorus{\mulevelsetvalue}{-\muxmulevelsetvalue} \right)$.
\end{itemize}

\medskip

\noindent \underline{\textbf{A description of the singular and regular behavior with respect to Cartesian coordinates}}.
	\begin{itemize}
		\item (Region without singularities).
		On the subset
		$\Upsilon\left(\twoargMrough{[\timefunction_0,0],[- \rightu,\leftu]}{\muxmulevelsetvalue} \backslash \twoargmumuxtorus{0}{-\muxmulevelsetvalue} \right)$
		of Cartesian coordinate space,
		the solution exists classically with respect to the Cartesian coordinates.
	\item (The fluid singularity).
		The following lower bound holds in 
		$\Upsilon\left(\twoargMrough{[\timefunction_0,0],[- \interestingu,\interestingu]}{\muxmulevelsetvalue} \right)$: 
		\begin{align} \label{E:FIXEDKAPPAMAINTHEOREMBLOWUPLOWERBOUND}
		|X \almostRiemann_{(+)}| 
		&
		\geq
		\frac{\mathring{\updelta}_*}{\upmu |\bar{\Speed}_{;\LogDensity} + 1|},
		\end{align}
		where $\mathring{\updelta}_* > 0$ is the data-parameter from \eqref{E:DELTASTARDEF},
		$\bar{\Speed}_{;\LogDensity} \eqdef \Speed_{;\LogDensity}(\LogDensity = 0,\Ent=0)$
		is $\Speed_{;\LogDensity}$ evaluated at the trivial solution, 
		$\bar{\Speed}_{;\LogDensity} + 1$ is a \underline{non-zero} constant
		by assumption, and the $\Sigma_t$-tangent vectorfield $X$ has Euclidean length satisfying
		$\sqrt{\sum_{a=1}^3 (X^a)^2} = 1 + \mathcal{O}(\mathring{\upalpha})$,
		where $\mathring{\upalpha}$ is the small parameter from Sect.\,\ref{SS:PARAMETERSIZEASSUMPTIONS}.
		In particular, if $q \in \Upsilon\left(\twoargmumuxtorus{0}{-\muxmulevelsetvalue} \right)$,
		then since 
		$\twoargmumuxtorus{0}{-\muxmulevelsetvalue} \in \twoargMrough{[\timefunction_0,0],[- \frac{\interestingu}{2},\frac{\interestingu}{2}]}{\muxmulevelsetvalue}$ 
		by \eqref{E:IMPROVEDLEVELSETSTRUCTUREANDLOCATIONOFMIN}, 
		and since $\upmu = 0$ along 
		$\Upsilon\left(\twoargmumuxtorus{0}{-\muxmulevelsetvalue} \right)$,
		it follows that
		$|X \almostRiemann_{(+)}|(q') \to \infty$ as $q' \rightarrow q$ in
		$\Upsilon\left(\twoargMrough{[\timefunction_0,0],[- \rightu,\leftu]}{\muxmulevelsetvalue} 
		\backslash \twoargmumuxtorus{0}{-\muxmulevelsetvalue} \right)$.
		Similarly, the following lower bounds hold in
		$\Upsilon\left(\twoargMrough{[\timefunction_0,0],[- \interestingu,\interestingu]}{\muxmulevelsetvalue} \right)$, 
		where $\LogDensity$ is the logarithmic density (see \eqref{E:LOGDENS}):
		\begin{align} \label{E:FIXEDKAPPAMAINTHEOREMDENSITYANDVELOCITYBLOWUPLOWERBOUNDS}
		|X \LogDensity| 
		&
		\geq
		\frac{\mathring{\updelta}_*}{4 \upmu |\bar{\Speed}_{;\LogDensity} + 1|},
		&
		|X v^1| 
		&
		\geq
		\frac{\mathring{\updelta}_*}{4 \upmu |\bar{\Speed}_{;\LogDensity} + 1|}.
		\end{align}
	\item (Regular behavior along the characteristics).
	The derivatives of
	$\wavearray$,
	$\vortrenormalized^i$,
	$\GradEnt^i$
	up to order $\Ntop - 11$
	with respect to the vectorfields in the 
	$\nullhyparg{u}$-tangent commutation set $\Tanset$ defined in \eqref{E:COMMUTATIONVECTORFIELDS}
	and the derivatives of
	$\VortVort^i$
	and
	$\DivGradEnt$
	up to order $\Ntop - 12$
	with respect to the elements of $\Tanset$ are $L^{\infty}$-bounded
	on $\Upsilon\left(\twoargMrough{[\timefunction_0,0],[- \rightu,\leftu]}{\muxmulevelsetvalue} \right)$.
	Finally, for $\alpha = 0,1,2,3$ and $A=2,3$, 
	the derivatives of 
	$\gfour_{ab} \Yvf{A}^a \partial_{\alpha} v^b$
	up to order $\Ntop - 12$
	with respect to the elements of $\Tanset$
	are $L^{\infty}$-bounded on
	$\Upsilon\left(\twoargMrough{[\timefunction_0,0],[- \rightu,\leftu]}{\muxmulevelsetvalue} \right)$.
	\end{itemize}
\end{theorem}

\begin{proof}
	The standard local well-posedness results and Cauchy stability provided by Prop.\,\ref{P:CAUCHYSTABILITYANDEXISTENCEOFOPENSETS}
	imply that there exists a $\timefunction_{Local} \in (3\timefunction_0/4,0)$
	such that the solution variables 
$\wavearray$,
$\vortrenormalized^i$,
$\GradEnt^i$,
$\VortVort^i$,
$\DivGradEnt$, 
$u$, 
and
$\timefunctionarg{\muxmulevelsetvalue}$
are classical solutions on
$\twoargMrough{[\timefunction_0,\timefunction_{Local}),[- \rightu,\leftu]}{\muxmulevelsetvalue}$
and such that all of the bootstrap assumptions from Sect.\,\ref{S:BOOTSTRAPEVERYTHINGEXCEPTENERGIES}
hold on $\twoargMrough{[\timefunction_0,\timefunction_{Local}),[- \rightu,\leftu]}{\muxmulevelsetvalue}$.
Let $\timefunction_{Max}$ be the supremum over all such $\timefunction_{Local}$.
Then the solution exists classically 
and satisfies all the bootstrap assumptions from Sect.\,\ref{S:BOOTSTRAPEVERYTHINGEXCEPTENERGIES}
on $\twoargMrough{[\timefunction_0,\timefunction_{Max}),[- \rightu,\leftu]}{\muxmulevelsetvalue}$.
If it were true that $\timefunction_{Max} < 0$, then
Prop.\,\ref{P:CONTINUATIONCRITERIA} (which we prove independently in the next section)
would imply that there is a $\Delta > 0$ 
with $\timefunction_{Max} + \Delta < 0$
such that
the solution exists classically 
and satisfies the bootstrap assumptions from Sect.\,\ref{S:BOOTSTRAPEVERYTHINGEXCEPTENERGIES}
on $\twoargMrough{[\timefunction_0,\timefunction_{Max} + \Delta),[- \rightu,\leftu]}{\muxmulevelsetvalue}$,
which is impossible in view of the definition of $\timefunction_{Max}$. Hence, 
$\timefunction_{Max} = 0$.

Aside from \eqref{E:FIXEDKAPPAMAINTHEOREMDENSITYANDVELOCITYBLOWUPLOWERBOUNDS} and 
the results concerning the boundedness of the quantities $\gfour_{ab} \Yvf{A}^a \partial_{\alpha} v^b$, 
the remaining conclusions of the theorem now follow from
\eqref{E:MINVALUEOFMUONFOLIATION} and the statements just below it,
\eqref{E:LOWERBOUNDONMAGNITUDEOFXRPLUS},
Lemma~\ref{L:CONTINUOUSEXTNESION},
Props.\,\ref{P:IMPROVEMENTOFAUXILIARYBOOTSTRAP},
\ref{P:HOMEOMORPHICANDDIFFEOMORPHICEXTENSIONOFCARTESIANCOORDINATES},
\ref{P:APRIORIL2ESTIMATESWAVEVARIABLES},
\ref{P:MAINHYPERSURFACEENERGYESTIMATESFORTRANSPORTVARIABLES},
\ref{P:ROUGHTORIENERGYESTIMATES},
\ref{P:APRIORIL2ESTIMATESACOUSTICGEOMETRY},
and \ref{P:IMPROVEMENTOFFUNDAMENTALQUANTITATIVEBOOTSTRAPASSUMPTIONS}
with $0$ in the role of $\timefunctionboot$.

The lower bounds stated in \eqref{E:FIXEDKAPPAMAINTHEOREMDENSITYANDVELOCITYBLOWUPLOWERBOUNDS}
follow from \eqref{E:FIXEDKAPPAMAINTHEOREMBLOWUPLOWERBOUND}
and the estimates 
$|\upmu X \LogDensity|
= 
\frac{1}{2}
\left\lbrace
	1 + \mathcal{O}(\mathring{\upalpha})
\right\rbrace
|\upmu X \RRiemann|
+
\mathcal{O}(\initialsmall)
$
and
$|\upmu X v^1| 
= 
\frac{1}{2}
\left\lbrace
	1 + \mathcal{O}(\mathring{\upalpha})
\right\rbrace
|\upmu X \RRiemann|
+
\mathcal{O}(\initialsmall)
$,
which follow from
\eqref{E:BACKGROUNDSOUNDSPEEDISUNITY},
\eqref{E:ALMOSTRIEMANNINVARIANTS}
(see also Remark~\ref{R:HOWWEESTIMATEDENSITYANDV1}),
\eqref{E:SCHEMATICSTRUCTUREOFXSMALL} for $\Speed - 1$,
and the estimates of 
Prop.\,\ref{P:IMPROVEMENTOFAUXILIARYBOOTSTRAP}
and Cor.\,\ref{C:IMPROVEAUX}.

Finally, we show that the derivatives of 
$\gfour_{ab} \Yvf{A}^a \partial_{\alpha} v^b$
up to order $\Ntop - 12$
with respect to the elements of $\Tanset$
are $L^{\infty}$-bounded on
$\Upsilon\left(\twoargMrough{[\timefunction_0,0],[- \rightu,\leftu]}{\muxmulevelsetvalue} \right)$.
We first consider the case $\alpha = 0$. 
We start by using \eqref{E:CARTESIANPARTIALTTOCOMMUTATORS} to deduce that
$\gfour_{ab} \Yvf{A}^a \partial_t v^b
= - 
\frac{\Lunit^1X^1 + \Lunit^2X^2 + \Lunit^3X^3}{\Speed^2} 
\gfour_{ab} \Yvf{A}^a X v^b
+
\mathrm{Error}
$,
where $\mathrm{Error}$ involves only $\nullhyparg{u}$-tangential derivatives of $v$.
From Lemma~\ref{L:SCHEMATICSTRUCTUREOFVARIOUSTENSORSINTERMSOFCONTROLVARS}
and Prop.\,\ref{P:IMPROVEMENTOFAUXILIARYBOOTSTRAP}, it follows that
the derivatives of 
$
\mathrm{Error}
$
up to order $\Ntop - 12$
with respect to the elements of $\Tanset$
are $L^{\infty}$-bounded on
$\Upsilon\left(\twoargMrough{[\timefunction_0,0],[- \rightu,\leftu]}{\muxmulevelsetvalue} \right)$.
Next, we use \eqref{E:SPECIFICVORTICITYDEF} and 
\eqref{E:ACOUSTICALMETRIC}--\eqref{E:INVERSEACOUSTICALMETRIC} 
to compute that relative to the Cartesian coordinates, we have
(where $\upepsilon_{ijk}$ is the fully antisymmetric symbol normalized by $\upepsilon_{123}=1$):
$X v^b = X^d \partial_d v^b = X^d \partial_b v^d 
+ 
\upepsilon_{dbe} X^d (\Flatcurl v)^e
= 
\Speed^2
X_d \partial_b v^d 
+ 
\exp(\LogDensity)
\upepsilon_{dbe} X^d \vortrenormalized^e
$.
Again using \eqref{E:ACOUSTICALMETRIC}, we find that:
\begin{align} 
\begin{split} \label{E:VORTICITYTERMINGOODCONTRACTIONSOFBADDERIVATIVE}
- 
\frac{\Lunit^1 X^1 + \Lunit^2 X^2 + \Lunit^3 X^3}{\Speed^2} 
\gfour_{ab} \Yvf{A}^a X v^b
&
=
- 
\frac{\Lunit^1 X^1 + \Lunit^2 X^2 + \Lunit^3 X^3}{\Speed^2} 
X_d \Yvf{A} v^d 
	\\
& \ \
- 
\exp(\LogDensity)
\frac{\Lunit^1X^1 + \Lunit^2X^2 + \Lunit^3X^3}{\Speed^2} 
\gfour_{ab}
\Yvf{A}^a
\upepsilon_{dbe} X^d \vortrenormalized^e.
\end{split}
\end{align}
Finally, using Lemma~\ref{L:SCHEMATICSTRUCTUREOFVARIOUSTENSORSINTERMSOFCONTROLVARS}
and Prop.\,\ref{P:IMPROVEMENTOFAUXILIARYBOOTSTRAP}, we see that
the derivatives of the products on RHS~\eqref{E:VORTICITYTERMINGOODCONTRACTIONSOFBADDERIVATIVE}
up to order $\Ntop - 12$
with respect to the elements of $\Tanset$
are $L^{\infty}$-bounded on
$\Upsilon\left(\twoargMrough{[\timefunction_0,0],[- \rightu,\leftu]}{\muxmulevelsetvalue} \right)$.
We have therefore proved the desired result in the case $\alpha = 0$.
To handle the cases $\alpha = 1,2,3$, we use a similar argument
that relies on the identities \eqref{E:CARTESIANPARTIAL1TOCOMMUTATORS}--\eqref{E:CARTESIANPARTIAL3TOCOMMUTATORS}
in place of \eqref{E:CARTESIANPARTIALTTOCOMMUTATORS}.

\end{proof}

\subsection{Continuation criteria}
\label{SS:CONTINUATIONCRITERIA}
In the next proposition, we provide the continuation criteria that are needed for the proof of
Theorem~\ref{T:EXISTENCEUPTOTHESINGULARBOUNDARYATFIXEDKAPPA}.
Roughly, the proposition shows that in the solution regime under study,
if $\timefunctionboot < 0$,
then the solution can be continued beyond
$\twoargMrough{[\timefunction_0,\timefunctionboot],[- \rightu,\leftu]}{\muxmulevelsetvalue}$
as a classical solution with respect to the Cartesian and geometric coordinates.
Central to the proof is the estimate \eqref{E:MINVALUEOFMUONFOLIATION},
which in particular shows that $\upmu$ is strictly positive on
$\twoargMrough{[\timefunction_0,\timefunctionboot],[- \rightu,\leftu]}{\muxmulevelsetvalue}$
whenever $\timefunctionboot < 0$, i.e., no shocks are present in
$\twoargMrough{[\timefunction_0,\timefunctionboot],[- \rightu,\leftu]}{\muxmulevelsetvalue}$.
The proposition is rather standard, except for some aspects concerning the regularity of
the solution with respect to the geometric coordinates.

\begin{proposition}[\textbf{Continuation criteria}]
\label{P:CONTINUATIONCRITERIA}
Assume the following:
\begin{itemize}
	\item The assumptions of Theorem~\ref{T:EXISTENCEUPTOTHESINGULARBOUNDARYATFIXEDKAPPA} hold.
	\item $\timefunctionboot < 0$.
	\item The rough time function $\timefunctionarg{\muxmulevelsetvalue}$,
the
solution variables
$\wavearray$,
$\vortrenormalized^i$,
$\GradEnt^i$,
$\VortVort^i$,
$\DivGradEnt$, 
$u$,
etc.\ are classical solutions on $\twoargMrough{[\timefunction_0,\timefunctionboot),[- \rightu,\leftu]}{\muxmulevelsetvalue}$.
	\item The bootstrap assumptions of Sect.\,\ref{S:BOOTSTRAPEVERYTHINGEXCEPTENERGIES}
		and \ref{SS:BOOTSTRAPASSUMPTIONSFORTHEWAVEENERGIES} hold on 
		$\twoargMrough{[\timefunction_0,\timefunctionboot),[- \rightu,\leftu]}{\muxmulevelsetvalue}$.
\end{itemize}

Then there exists a $\Delta \in (0,|\timefunctionboot|)$ 
such that the 
rough time function $\timefunctionarg{\muxmulevelsetvalue}$,
the
solution variables
$\wavearray$,
$\vortrenormalized^i$,
$\GradEnt^i$,
$\VortVort^i$,
$\DivGradEnt$, 
$u$, 
and all of the
other geometric quantities defined
throughout the article can be uniquely extended 
(where the solution variables are classical solutions)
to a strictly larger region of the form
$\twoargMrough{[\timefunction_0,\timefunctionboot + \Delta),[- \rightu,\leftu]}{\muxmulevelsetvalue}$
on which all of the bootstrap assumptions of Sects.\,\ref{S:BOOTSTRAPEVERYTHINGEXCEPTENERGIES} 
and \ref{SS:BOOTSTRAPASSUMPTIONSFORTHEWAVEENERGIES} hold.

\end{proposition}

\begin{proof}
	Throughout this proof, we allow the small positive numbers 
	$m > 0$ and $\Delta > 0$ to vary from line to line,
	sometimes silently shrinking them as necessary.

	\medskip
	
	\noindent \textbf{Step 1: Extension to $\twoargMrough{[\timefunction_0,\timefunctionboot],[- \rightu,\leftu]}{\muxmulevelsetvalue}$}.
	Since $\timefunctionboot = - \upmuboot < 0$ by assumption,
	Lemmas~\ref{L:DIFFEOMORPHICEXTENSIONOFROUGHCOORDINATES} and \ref{L:CONTINUOUSEXTNESION},
	Prop.\,\ref{P:HOMEOMORPHICANDDIFFEOMORPHICEXTENSIONOFCARTESIANCOORDINATES},
	and the bootstrap assumptions
	imply that
	the quantities
$\wavearray$,
$\vortrenormalized^i$,
$\GradEnt^i$,
$\VortVort^i$,
$\DivGradEnt$, 
$\Upsilon$,
$t$, $x^1$, $x^2$, $x^3$,
$u$, 
and all the
other geometric quantities defined
throughout the article 
extend from $\twoargMrough{[\timefunction_0,\timefunctionboot),[- \rightu,\leftu]}{\muxmulevelsetvalue}$
to the compact set
$\twoargMrough{[\timefunction_0,\timefunctionboot],[- \rightu,\leftu]}{\muxmulevelsetvalue}$
as classical solutions relative to the geometric coordinates $(t,u,x^2,x^3)$
such that $\Upsilon$ is a diffeomorphism on $\twoargMrough{[\timefunction_0,\timefunctionboot],[- \rightu,\leftu]}{\muxmulevelsetvalue}$.
The same results yields that 
these quantities extend to the compact set
$[\timefunction_0,\timefunctionboot] \times [- \rightu,\leftu] \times \mathbb{T}^2$
as classical solutions relative to the rough adapted coordinates $(\timefunctionarg{\muxmulevelsetvalue},u,x^2,x^3)$,
and to the compact set
$\Upsilon\left(\twoargMrough{[\timefunction_0,\timefunctionboot],[- \rightu,\leftu]}{\muxmulevelsetvalue} \right)$
as classical solutions relative to the Cartesian coordinates $(t,x^1,x^2,x^3)$.
Moreover, in view of the energy estimates of Sect.\,\ref{S:STATEMENTOFALLL2ESTIMATESANDBOOTSTRAPASSUMPTIONSFORWAVEENERGIES},
it is a standard result 
that for all $(\timefunction,u) \in [\timefunction_0,\timefunctionboot] \times [- \rightu,\leftu]$,
the extended quantities enjoy the same Sobolev and Lebesgue regularity
(i.e., the corresponding norms are all finite)
with respect to the geometric coordinates
on 
$\hypthreearg{\timefunction}{[- \rightu,\leftu]}{\muxmulevelsetvalue}$,
$\twoargroughtori{\timefunction,u}{\muxmulevelsetvalue}$,
and
$\nullhypthreearg{\muxmulevelsetvalue}{u}{[\timefunction_0,\timefunction]}$
as the data  
on $\hypthreearg{\timefunction_0}{[- \rightu,\leftu]}{\muxmulevelsetvalue}$,
$\twoargroughtori{\timefunction_0,u}{\muxmulevelsetvalue}$,
and
$\nullhypthreearg{\muxmulevelsetvalue}{- \rightu}{[\timefunction_0,\timefunctionboot]}$,
and that relative to all of the corresponding Sobolev and Lebesgue function space topologies on these surfaces
that we have used throughout the paper,
the solution is continuous with respect to $(\timefunction,u)$
on $[\timefunction_0,\timefunctionboot] \times [- \rightu,\leftu]$;
we refer readers to \cite[Section~2.7]{jS2008c} 
for the main ideas behind the proof of these 
``propagation of regularity'' and ``continuity-in-norm'' results
in the context of the relativistic 
Euler equations coupled to N\"{o}rdstrom's theory of gravity.

Furthermore, by \eqref{E:C21BOUNDFORCONSTANTTIMEFUNCTIONGRAPH}--\eqref{E:LEVELSETSOFTIMEFUNCTIONAREAGRAPH},
there is a $\timefunctionboot$-dependent 
function $\Cartesiantisafunctiononlevelsetsofroughtimefunctionarg{\timefunctionboot}{\muxmulevelsetvalue}$ on $[- \rightu,\leftu] \times \mathbb{T}^2$ 
such that relative to the geometric coordinates, 
we have:
\begin{align*}
\hypthreearg{\timefunctionboot}{[- \rightu,\leftu]}{\muxmulevelsetvalue} 
& = 
\left\lbrace 
	\left(\Cartesiantisafunctiononlevelsetsofroughtimefunctionarg{\timefunctionboot}{\muxmulevelsetvalue}(u,x^2,x^3),u,x^2,x^3 \right) 
	\ | \   
	(u,x^2,x^3)
	\in
	[- \rightu,\leftu] \times \mathbb{T}^2
\right\rbrace,
\end{align*}
and such that
$\| \Cartesiantisafunctiononlevelsetsofroughtimefunctionarg{\timefunctionboot}{\muxmulevelsetvalue} \|_{C^{2,1}([- \rightu,\leftu] \times \mathbb{T}^2)} \leq C$.
In addition, using 
\eqref{E:SIZEOFRTRANS},
\eqref{E:RTRANSNORMSMALLFACTOR},
\eqref{E:HYPNORMALSIZE},
\eqref{E:MINVALUEOFMUONFOLIATION},
\eqref{E:BOUNDSONLMUINTERESTINGREGION},
and \eqref{E:POINTWISEBOUNDTANGENTIALANDTRANSVERSALDERIVATIVESOFRTRANSNORMSMALLFACTOR},
and our crucial assumption that $\timefunctionboot < 0$,
we see that
$\gfour(\hypnormalarg{\muxmulevelsetvalue},\hypnormalarg{\muxmulevelsetvalue}) < 0$ on $\hypthreearg{\timefunctionboot}{[- \rightu,\leftu]}{\muxmulevelsetvalue}$, 
i.e., that $\hypthreearg{\timefunctionboot}{[- \rightu,\leftu]}{\muxmulevelsetvalue}$ is $\gfour$-spacelike.
Also using Prop.\,\ref{P:HOMEOMORPHICANDDIFFEOMORPHICEXTENSIONOFCARTESIANCOORDINATES},
we further deduce that the hypersurface $\Upsilon\left(\hypthreearg{\timefunctionboot}{[- \rightu,\leftu]}{\muxmulevelsetvalue} \right)$
in Cartesian coordinate space is $C^{2,1}$.

\medskip

\noindent \textbf{Step 2: Extending all quantities -- except for the rough time function --
beyond $\twoargMrough{[\timefunction_0,\timefunctionboot],[- \rightu,\leftu]}{\muxmulevelsetvalue}$}.
Let $Q$ denote the quantities
$\wavearray$,
$\vortrenormalized^i$,
$\GradEnt^i$,
$\VortVort^i$,
$\DivGradEnt$, 
$\Upsilon$,
$t$, $x^1$, $x^2$, $x^3$,
$u$, 
etc.\ from the beginning of Step~1. 
Note that the rough time function $\timefunctionarg{\muxmulevelsetvalue}$ is not among
these quantities. 
To extend $Q$ beyond $\twoargMrough{[\timefunction_0,\timefunctionboot],[- \rightu,\leftu]}{\muxmulevelsetvalue}$,
we will use the Cauchy stability arguments given in Appendix~\ref{A:OPENSETOFDATAEXISTS};
while it is not essential for us to use the results of Appendix~\ref{A:OPENSETOFDATAEXISTS} here,
it allows us to avoid treating the characteristic initial value problem for the compressible Euler equations,
which would have involved inessential technical complications.
Specifically, in Step~3 of the proof of Prop.\,\ref{P:CAUCHYSTABILITYANDEXISTENCEOFOPENSETS},
we used Cauchy stability to show that $Q$ exists classically in the 
(geometric coordinate) region $\CSregion_{Small}^{[0,5 \blowuptimePS]}$ depicted in Fig.\,\ref{F:CAUCHYSTABILITYREGION}.
Moreover, in the proof of Lemma~\ref{L:CONTROLOFDATAFORROUGHTORIENERGYESTIMATES},
we showed that we can extend the rough time function $\timefunctionarg{\muxmulevelsetvalue}$
into a subset of $\CSregion_{Small}^{[0,5 \blowuptimePS]}$
so that it is defined on a region 
$
\twoargMrough{[\timefunction_0,\timefunctionboot],[-U_*,\leftu]}{\muxmulevelsetvalue}
$
containing $\hypthreearg{\timefunctionboot}{[- (\rightu + \epsilon),\leftu]}{\muxmulevelsetvalue}$
for all sufficiently small $\epsilon > 0$
(see Footnote~\ref{FN:USTARISBIGGERTHANRIGHTU})
and such that the estimate \eqref{E:C21GEOESTIMATEFORROUGHTIMEFUNCTIONINSMALLCACUHYSTABILITYREGION} holds
in 
$
\twoargMrough{[\timefunction_0,\timefunctionboot],[-U_*,\leftu]}{\muxmulevelsetvalue}
$.
In particular, combining these results with
\eqref{E:DATAASSUMPTIONSIZEOFCARTESIANT},
\eqref{E:PERTURBEDBLOWUPDELTAISCLOSETOBACKGROUNDONE},
and the results we derived in Step~1 of the present proof, 
we see that there is an $\epsilon > 0$ such that 
$\hypthreearg{\timefunctionboot}{[- (\rightu + \epsilon),\leftu]}{\muxmulevelsetvalue}
\subset
\twoargMrough{[\timefunction_0,\timefunctionboot],[- (\rightu + \epsilon),\leftu]}{\muxmulevelsetvalue}
\subset
\twoargMrough{[\timefunction_0,\timefunctionboot],[- \rightu,\leftu]}{\muxmulevelsetvalue}
\cup
\CSregion_{Small}^{[0,5 \blowuptimePS]}
$.
Shrinking $\epsilon$ if necessary, and using a standard compactness argument, we can assume that
$\Upsilon$ is a diffeomorphism on $\twoargMrough{[\timefunction_0,\timefunctionboot],[- (\rightu + \epsilon),\leftu]}{\muxmulevelsetvalue}$.
Hence, the quantities $Q$ and the rough time function $\timefunctionarg{\muxmulevelsetvalue}$ extend to
$\twoargMrough{[\timefunction_0,\timefunctionboot],[- (\rightu + \epsilon),\leftu]}{\muxmulevelsetvalue}$
as classical solutions with respect to the geometric coordinates, and similarly with
respect to the Cartesian coordinates 
(on $\Upsilon \left( \twoargMrough{[\timefunction_0,\timefunctionboot],[- (\rightu + \epsilon),\leftu]}{\muxmulevelsetvalue} \right)$).
Moreover, by continuity and the results of Step~1 of the present proof, 
we see that
$
\hypthreearg{\timefunctionboot}{[- (\rightu + \epsilon),\leftu]}{\muxmulevelsetvalue}
$
is a $\gfour$-spacelike hypersurface portion if
$\epsilon > 0$ is small enough.

We now use standard local well-posedness 
(see \cite[Section~2.7]{jS2008c} for the main ideas behind the analysis in the context of the relativistic Euler equations
coupled to N\"{o}rdstrom's theory of gravity)
for the compressible Euler equations and the eikonal equation \eqref{E:EIKONALEQUATION}
relative to the Cartesian coordinates. That is, starting from the data on the $\gfour$-spacelike hypersurface portion
$
\Upsilon\left(\hypthreearg{\timefunctionboot}{[- (\rightu + \epsilon),\leftu]}{\muxmulevelsetvalue} \right)
$
in Cartesian coordinate space,
we consider the
corresponding local solution
$\wavearray$,
$\vortrenormalized^i$,
$\GradEnt^i$,
$\VortVort^i$,
and
$\DivGradEnt$ to the equations of Theorem~\ref{T:GEOMETRICWAVETRANSPORTSYSTEM} 
and the solution $u$ to the eikonal equation \eqref{E:EIKONALEQUATION}.
This ``extended'' solution is classical and enjoys the same regularity as the solution in
$\twoargMrough{[\timefunction_0,\timefunctionboot],[- (\rightu + \epsilon),\leftu]}{\muxmulevelsetvalue}$;
we will discuss this in more detail below.
For the extended solution, we let 
$\mathscr{D}_{\epsilon}^+ 
= \mathscr{D}_{\epsilon}^+\left(\Upsilon\left(\hypthreearg{\timefunctionboot}{[- (\rightu + \epsilon),\leftu]}{\muxmulevelsetvalue} \right) \right)$ denote the future domain of dependence in Cartesian coordinate space 
$\R_t \times \R_{x^1} \times \mathbb{T}^2$
of the set 
$
\Upsilon\left(\hypthreearg{\timefunctionboot}{[- (\rightu + \epsilon),\leftu]}{\muxmulevelsetvalue} \right)
$
with respect to the acoustical metric $\gfour$.
For each small $\delta > 0$, let $\mathscr{D}_{\epsilon;\delta}^+$ denote the subset of $\mathscr{D}_{\epsilon}^+$
consisting of the points in $\mathscr{D}_{\epsilon}^+$ that can be joined to
$
\Upsilon\left(\hypthreearg{\timefunctionboot}{[- (\rightu + \epsilon),\leftu]}{\muxmulevelsetvalue} \right)
$
by a $C^1$ curve in $\mathscr{D}_{\epsilon}^+$ that has length $\leq \delta$ with respect to the 
standard Euclidean metric on $\mathbb{R}_t \times \mathbb{R}_{x^1} \times \mathbb{T}^2$;
$\mathscr{D}_{\epsilon;\delta}^+$ is, in particular, a 
neighborhood of 
$
\Upsilon\left(\hypthreearg{\timefunctionboot}{[- (\rightu + \epsilon),\leftu]}{\muxmulevelsetvalue} \right)
$
in $\mathscr{D}_{\epsilon}^+$ such that the lateral boundaries of $\mathscr{D}_{\epsilon;\delta}^+$
contain $\gfour$-null hypersurface portions, one of which is a portion of 
$\nullhyparg{\leftu}$; see Fig.\,\ref{F:EXTENSIONOFBOOTSTRAPSURFACE} for a picture of the setup in geometric coordinate space.
Since $\Upsilon$ is a diffeomorphism on $\twoargMrough{[\timefunction_0,\timefunctionboot],[- (\rightu + \epsilon),\leftu]}{\muxmulevelsetvalue}$,
a standard compactness argument yields that if 
$\delta$ and $\epsilon > 0$ are small enough, 
then $\Upsilon^{-1}$ is a diffeomorphism from $\mathscr{D}_{\epsilon;\delta}^+$ 
onto its image in geometric coordinate space $\R_t \times \R_u \times \T^2$. 
In particular, by using $\Upsilon^{-1}$ to change variables to geometric coordinates after extending in Cartesian coordinates, 
we see that the quantities $Q$
-- but not yet the rough time function --
can be extended as classical solutions relative to the geometric coordinates
to a larger region $\extendedtwoargMrough{m}{\muxmulevelsetvalue}$ 
of the following form for some sufficiently small $m > 0$:
\begin{align} \label{E:CONTINUATIONVERTICALLYTRANSLATEDFOLIATION}
\extendedtwoargMrough{m}{\muxmulevelsetvalue}
& \eqdef
	\twoargMrough{[\timefunction_0,\timefunctionboot],[- \rightu,\leftu]}{\muxmulevelsetvalue}
	\cup
	\newregionextendedtwoargMrough{m}{\muxmulevelsetvalue},
	\\
\begin{split} \label{E:NEWPARTCONTINUATIONVERTICALLYTRANSLATEDFOLIATION}
\newregionextendedtwoargMrough{m}{\muxmulevelsetvalue}
& 
\eqdef
\bigcup
\left\lbrace 
	\left(\Cartesiantisafunctiononlevelsetsofroughtimefunctionarg{\timefunctionboot}{\muxmulevelsetvalue}(u,x^2,x^3) + m',u,x^2,x^3 \right) 
	\ | \   
	(m',u,x^2,x^3)
	\in
	[0,m]
	\times
	[- \rightu,\leftu] 
	\times 
	\mathbb{T}^2
\right\rbrace
	\\
& \subset \Upsilon^{-1} \left(\mathscr{D}_{\epsilon;\delta}^+ \right),
\end{split}
\end{align}
such that each of the following hypersurfaces:
\begin{align} \label{E:CONTINUATIONVERTICALLYTRANSLATEDHYPERSURFACES}
\left\lbrace 
	\left(\Cartesiantisafunctiononlevelsetsofroughtimefunctionarg{\timefunctionboot}{\muxmulevelsetvalue}(u,x^2,x^3) + m',u,x^2,x^3 \right) 
	\ | \   
	(u,x^2,x^3)
	\times
	[- \rightu,\leftu] 
	\times 
	\mathbb{T}^2
\right\rbrace
\end{align}
is $\gfour$-spacelike and such that
$\Upsilon$ is a diffeomorphism from 
$
\extendedtwoargMrough{m}{\muxmulevelsetvalue}
$
onto its image in Cartesian coordinate space.
We clarify that in carrying out this argument, 
we have taken into account that the $\gfour$-null boundary of
$\extendedtwoargMrough{m}{\muxmulevelsetvalue}$
on which $u \equiv \leftu$ 
lies in the future domain of dependence (with respect to the acoustical metric $\gfour$)
of $\hypthreearg{\timefunctionboot}{[- \rightu,\leftu]}{\muxmulevelsetvalue}$ (without needing to extend to $u \in [- (\rightu + \epsilon),\leftu]$);
this is tied to the fact that by construction,
the level sets of $u$ are right-moving in Cartesian coordinate space,
as is shown in Fig.\,\ref{F:INFININTEDENSITYOFCHARACTERISTICSONSINGULARBOUNDARY}.
In contrast, the $\gfour$-null boundary of $\extendedtwoargMrough{m}{\muxmulevelsetvalue}$
on which $u \equiv -\rightu$
does not lie in the future domain of dependence of 
$\hypthreearg{\timefunctionboot}{[- \rightu,\leftu]}{\muxmulevelsetvalue}$.
It does, however, lie in the future domain of dependence of 
$\hypthreearg{\timefunctionboot}{[- (\rightu + \epsilon),\leftu]}{\muxmulevelsetvalue}$
whenever $\epsilon > 0$ and $m$ is small compared to $\epsilon$
(this is the reason we are working with the extended hypersurface 
$\hypthreearg{\timefunctionboot}{[- (\rightu + \epsilon),\leftu]}{\muxmulevelsetvalue}$). Moreover, thanks to the bound 
$\| \Cartesiantisafunctiononlevelsetsofroughtimefunctionarg{\timefunctionboot}{\muxmulevelsetvalue} \|_{C^{2,1}([- \rightu,\leftu] \times \mathbb{T}^2)} \leq C$ from Step~1,
arguments similar to the ones we used in 
the proof of Lemma~\ref{L:DIFFEOMORPHICEXTENSIONOFROUGHCOORDINATES}
can be used to show that $\extendedtwoargMrough{m}{\muxmulevelsetvalue}$
is quasi-convex in the sense of item 6 in the statement of that lemma.

\begin{center}
	\begin{figure}  
		\begin{overpic}[scale=.8, grid = false, tics=5, trim=-.5cm -1cm -1cm -.5cm, clip]{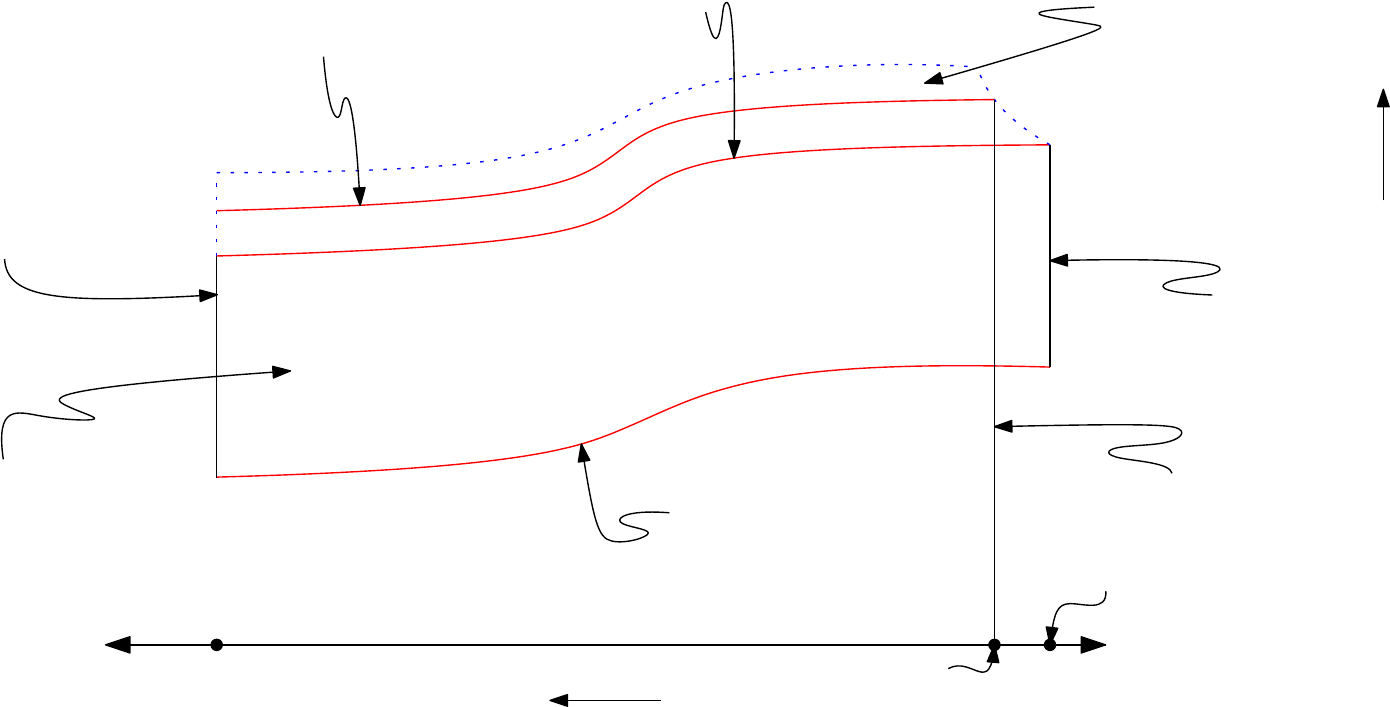}
			\put (75,51) {$\Upsilon^{-1}(\mathscr{D}_{\epsilon;\delta}^+)$}
			\put (0,36) {$\nullhypthreearg{\muxmulevelsetvalue}{\leftu}{[\timefunction_0,\timefunctionboot]}$}
			\put (0,50) {$\{(\Cartesiantisafunctiononlevelsetsofroughtimefunctionarg{\timefunctionboot}{\muxmulevelsetvalue}(u,x^2,x^3) + m,u,x^2,x^3)\}$}
			\put (83,32) {$\nullhypthreearg{\muxmulevelsetvalue}{- (\rightu + \epsilon)}{[\timefunction_0,\timefunctionboot]}$}
			\put (47,18) {$\hypthreearg{\timefunction_0}{[- (\rightu + \epsilon),\leftu]}{\muxmulevelsetvalue}$}
			\put (43,54) {$\hypthreearg{\timefunctionboot}{[-  (\rightu + \epsilon),\leftu]}{\muxmulevelsetvalue}$}
			\put (95,40) {$t$}
			\put (17,7) {$\leftu$}
			\put (60,7) {$-\rightu$}
			\put (73,15) {$- (\rightu + \epsilon)$}
			\put (0,20) {$\twoargMrough{[\timefunction_0,\timefunctionboot],[- (\rightu + \epsilon),\leftu]}{\muxmulevelsetvalue}$}
			\put (80,19) {$\nullhyparg{-\rightu}^{[0,\Cartesiantisafunctiononlevelsetsofroughtimefunctionarg{\timefunctionboot}{\muxmulevelsetvalue} + m]}$}
			\put (42,5) {$u$}
		\end{overpic}
		  \caption{The extension to $\newregionextendedtwoargMrough{m}{\muxmulevelsetvalue}$ in geometric coordinates}
  \label{F:EXTENSIONOFBOOTSTRAPSURFACE}
	\end{figure}
\end{center}

Given the energy estimate framework we established in the bulk of the paper
(which yields $L^2$estimates up to top-order under commutations with the elements of the 
$\nullhyparg{u}$-tangent commutation set $\Tanset$),
it is also standard (see \cite[Section~2.7]{jS2008c}) that,
as in Step 1, the solution enjoys the same regularity with respect to the geometric coordinates
in the extended region
$
\extendedtwoargMrough{m}{\muxmulevelsetvalue}
$
as it does in $\twoargMrough{[\timefunction_0,\timefunctionboot],[- \rightu,\leftu]}{\muxmulevelsetvalue}$,
where in $\newregionextendedtwoargMrough{m}{\muxmulevelsetvalue}$, 
regularity is measured on the spacelike hypersurfaces \eqref{E:CONTINUATIONVERTICALLYTRANSLATEDHYPERSURFACES} in geometric coordinate space
as well as constant-$u$ $\gfour$-null hypersurface portions contained in
$
\extendedtwoargMrough{m}{\muxmulevelsetvalue}
$.
Similarly, the solution enjoys the same regularity with respect to the Cartesian coordinates
in 
$
\Upsilon\left(\extendedtwoargMrough{m}{\muxmulevelsetvalue} \right)
$
as it does in $\Upsilon\left( \twoargMrough{[\timefunction_0,\timefunctionboot],[- \rightu,\leftu]}{\muxmulevelsetvalue} \right)$.

\medskip

\noindent \textbf{Step 3: Extending the rough time function
beyond $\twoargMrough{[\timefunction_0,\timefunctionboot],[- \rightu,\leftu]}{\muxmulevelsetvalue}$}.
We first note that 
Lemmas~\ref{L:DIFFEOMORPHICEXTENSIONOFROUGHCOORDINATES},
\ref{L:CONTINUOUSEXTNESION},
and
\ref{L:CHOVFROMROUGHCOORDINATESTOMUWEGIGHTEDXMUCOORDINATES} 
and the relation \eqref{E:SIMPLERELATIONSHIPBETWEENCHOVMAPS}
imply that the map 
$\CHOVgeotomumuxmu(t,u,x^2,x^3) = (\upmu,\muX \upmu,x^2,x^3)$
defined in \eqref{E:CHOVFROMGEOMETRICCOORDINATESTOMUWEGIGHTEDXMUCOORDINATES}
is a $C_{\textnormal{geo}}^{1,1}$ diffeomorphism from
$\twoargMrough{[\timefunction_0,\timefunctionboot],[-\interestingu,\interestingu]}{\muxmulevelsetvalue}$
onto a set containing $[\upmuboot,\mupositive] \times \lbrace - \muxmulevelsetvalue \rbrace \times \mathbb{T}^2$
such that 
$
\lbrace \upmuboot \rbrace \times \lbrace - \muxmulevelsetvalue \rbrace \times \mathbb{T}^2
\subset
\CHOVgeotomumuxmu\left(\hypthreearg{\timefunctionboot}{(-\interestingu,\interestingu)}{\muxmulevelsetvalue} \right)$.
By exploiting the compactness of 
$
\extendedtwoargMrough{m}{\muxmulevelsetvalue}
$,
we can, recalling that $\geop{t} \upmu < 0$ 
when $|u| \leq \interestingu$ by
\eqref{E:BOUNDSONGEOMETRICTDERIVATIVEMUINTERESTINGREGION}
and shrinking $m$ if necessary, 
assume that
$\CHOVgeotomumuxmu$ is a $C_{\textnormal{geo}}^{1,1}$ diffeomorphism from $\extendedtwoargMrough{m}{\muxmulevelsetvalue} \cap \lbrace |u| \leq \interestingu \rbrace$
onto a set containing 
$
\lbrace \upmuboot \rbrace \times \lbrace - \muxmulevelsetvalue \rbrace \times \mathbb{T}^2
$
in its interior, and that there is a $\Delta \in (0,\upmuboot)$ such that
$
[\upmuboot-\Delta, \upmuboot+\Delta] \rbrace \times [- \muxmulevelsetvalue - \Delta, - \muxmulevelsetvalue + \Delta] \times \mathbb{T}^2
$
is contained in the interior of
$
\CHOVgeotomumuxmu\left( \extendedtwoargMrough{m}{\muxmulevelsetvalue} \cap \lbrace |u| \leq \interestingu \rbrace \right)
$.
In particular, for $\mulevelsetvalue \in [\upmuboot - \Delta,\upmuboot + \Delta]$,
the $\upmu$-adapted tori $\twoargmumuxtorus{\mulevelsetvalue}{-\muxmulevelsetvalue}$
defined in \eqref{E:MUXMUTORI} are two-dimensional $C_{\textnormal{geo}}^{1,1}$ submanifolds
contained in the interior of $\extendedtwoargMrough{m}{\muxmulevelsetvalue} \cap \lbrace |u| \leq \interestingu \rbrace$,
and the hypersurface portion 
$
\datahypfortimefunctiontwoarg{-\muxmulevelsetvalue}{[-\upmuboot - \Delta,-\upmuboot + \Delta]}
\eqdef
\datahypfortimefunctionarg{-\muxmulevelsetvalue}
\cap
\left\lbrace (t,u,x^2,x^3) \in \R \times \R \times \T^2 
		\ | \ 
		\upmuboot - \Delta \leq \upmu(t,u,x^2,x^3) \leq \upmuboot + \Delta
\right\rbrace
$
is a three-dimensional $C_{\textnormal{geo}}^{1,1}$ submanifold contained in
the interior of
$\extendedtwoargMrough{m}{\muxmulevelsetvalue} \cap \lbrace |u| \leq \interestingu \rbrace$
(see \eqref{E:MUXMUTORI} for the definition of $\datahypfortimefunctionarg{-\muxmulevelsetvalue}$,
and compare with the alternate definition
\eqref{E:TRUNCATEDLEVELSETSOFMUXMU} of $\datahypfortimefunctiontwoarg{-\muxmulevelsetvalue}{[-\upmuboot - \Delta,-\upmuboot + \Delta]}$,
which will eventually agree with the definition given above).
Moreover, the estimate \eqref{E:BAMUTRANSVERSALCONVEXITY}, 
which by continuity holds in $\extendedtwoargMrough{m}{\muxmulevelsetvalue} \cap \lbrace |u| \leq \interestingu \rbrace$
with different constants,
implies that the vectorfield $\Wtransarg{\muxmulevelsetvalue}$ is transversal to 
$
\datahypfortimefunctiontwoarg{-\muxmulevelsetvalue}{[-\upmuboot - \Delta,-\upmuboot + \Delta]}
$.
Considering also that $\Wtransarg{\muxmulevelsetvalue}$ is tangent to the lower boundary
$\hypthreearg{\timefunctionboot}{[- \rightu,\leftu]}{\muxmulevelsetvalue}$ of $\newregionextendedtwoargMrough{m}{\muxmulevelsetvalue}$,
we see that these results are sufficient to allow us to extend the construction of the rough time function 
$\timefunctionarg{\muxmulevelsetvalue}$ 
(see Sect.\,\ref{SS:BASICCONSTRUCTIONSFORROUGHTIMEFUNCTION}, 
Lemma~\ref{L:FLOWMAPFORGENERATOROFROUGHTIMEFUNCTION},
and
Lemma~\ref{L:ODESOLUTIONSTHATARESMOOTHERTHANTHEDATAHYPERSURFACE})
from $\twoargMrough{[\timefunction_0,\timefunctionboot],[- \rightu,\leftu]}{\muxmulevelsetvalue}$ to all of
$\extendedtwoargMrough{m}{\muxmulevelsetvalue}$, i.e.,
so that $\timefunctionarg{\muxmulevelsetvalue}$
is defined on $\extendedtwoargMrough{m}{\muxmulevelsetvalue}$ (where we have perhaps shrunk $m$ if necessary)
and satisfies $\| \timefunctionarg{\muxmulevelsetvalue} \|_{W_{\textnormal{geo}}^{3,\infty}(\mbox{\upshape int}(\extendedtwoargMrough{m}{\muxmulevelsetvalue}))} \leq C$.
From this bound, the quasi-convexity of $\extendedtwoargMrough{m}{\muxmulevelsetvalue}$
mentioned in Step 2, and the Sobolev embedding result \eqref{E:SOBOELVEMBEDDINGRELYINGONQUASICONVEXITY}
-- which also holds for the domain $\extendedtwoargMrough{m}{\muxmulevelsetvalue}$ thanks to its quasi-convexity --
we also find that 
$\| \timefunctionarg{\muxmulevelsetvalue} \|_{C_{\textnormal{geo}}^{2,1}(\extendedtwoargMrough{m}{\muxmulevelsetvalue})} \leq C$.
Considering also the estimate $\geop{t} \timefunctionarg{\muxmulevelsetvalue} \approx 1$,
which holds in $\extendedtwoargMrough{m}{\muxmulevelsetvalue}$ by \eqref{E:PARTIALTIMEDERIVATIVEOFROUGHTIMEFUNCTIONISAPPROXIMATELYUNITY} 
and continuity,
we further deduce that the map $\CHOVgeotorough{\muxmulevelsetvalue}(t,u,x^2,x^3) = (\timefunctionarg{\muxmulevelsetvalue},u,x^2,x^3)$
is a diffeomorphism from
$\extendedtwoargMrough{m}{\muxmulevelsetvalue}$ onto its image,
and that there is a $\Delta > 0$ (perhaps smaller than before)
such that the image set contains
$
[\timefunction_0,\timefunctionboot + \Delta] \times [- \rightu,\leftu] \times \mathbb{T}^2
$.
That is, 
$\| \timefunctionarg{\muxmulevelsetvalue} \|_{C_{\textnormal{geo}}^{2,1}(\twoargMrough{[\timefunction_0,\timefunctionboot + \Delta],[- \rightu,\leftu]}{\muxmulevelsetvalue})} 
\lesssim 1$, 
and
$\CHOVgeotorough{\muxmulevelsetvalue}$
is a diffeomorphism from
$\twoargMrough{[\timefunction_0,\timefunctionboot + \Delta],[- \rightu,\leftu]}{\muxmulevelsetvalue}$ onto its image.

\medskip

\noindent \textbf{Step 4: Propagation of regularity on the foliation induced by the extended $\timefunctionarg{\muxmulevelsetvalue}$}.
We have constructed the solution on the extended region 
$\twoargMrough{[\timefunction_0,\timefunctionboot + \Delta],[- \rightu,\leftu]}{\muxmulevelsetvalue}$
as well as the foliation
$\left\lbrace \hypthreearg{\timefunction}{[- \rightu,\leftu]}{\muxmulevelsetvalue} \right\rbrace_{\timefunction 
\in [\timefunction_0,\timefunctionboot + \Delta]}$ of it.
We can now argue as in Step 2 to deduce that 
the solution enjoys the same regularity with respect to the geometric coordinates
in the extended region
$\twoargMrough{[\timefunction_0,\timefunctionboot + \Delta],[- \rightu,\leftu]}{\muxmulevelsetvalue}$
as it does in $\twoargMrough{[\timefunction_0,\timefunctionboot],[- \rightu,\leftu]}{\muxmulevelsetvalue}$,
where, as in the bulk of the paper, regularity is measured on the hypersurfaces
$
\hypthreearg{\timefunction}{[- \rightu,\leftu]}{\muxmulevelsetvalue}
$
(which, by continuity, are $\gfour$-spacelike for $\Delta$ sufficiently small),
on constant $u$ $\gfour$-null hypersurface portions
$\nullhypthreearg{\muxmulevelsetvalue}{u}{[\timefunction_0,\timefunction_0 + \Delta]}$,
and on the rough tori $\twoargroughtori{\timefunction,u}{\muxmulevelsetvalue}$.
 
\medskip


\medskip
\noindent \textbf{Step 5: The bootstrap assumptions hold for the extended solution}.
The results we proved throughout the paper have yielded, 
on $\twoargMrough{[\timefunction_0,\timefunctionboot),[- \rightu,\leftu]}{\muxmulevelsetvalue}$,
strict improvements of all the bootstrap assumptions of
Sects.\,\ref{S:BOOTSTRAPEVERYTHINGEXCEPTENERGIES}
and \ref{SS:BOOTSTRAPASSUMPTIONSFORTHEWAVEENERGIES};
see Sect.\,\ref{SS:SUMMARYOFIMPROVEMENTOFBOOTSTRAPASSUMPTIONS}
for a description of the results that yield the improvements.
Hence, by exploiting the continuity guaranteed by local well-posedness,
we conclude that all the bootstrap assumptions 
of Sects.\,\ref{S:BOOTSTRAPEVERYTHINGEXCEPTENERGIES} and \ref{SS:BOOTSTRAPASSUMPTIONSFORTHEWAVEENERGIES}
also hold on
$\twoargMrough{[\timefunction_0,\timefunctionboot + \Delta),[- \rightu,\leftu]}{\muxmulevelsetvalue}$.

\end{proof}

\section{Developments of the data, the singular boundary and the crease, and a new time function}
\label{S:DEVELOPMENTSOFTHEDATASINGULARBOUNDARYNEWTIMEFUNCTION}
In Sect.\,\ref{S:MAINRESULTS}, we will state and prove Theorem~\ref{T:DEVELOPMENTANDSTRUCTUREOFSINGULARBOUNDARY},
which is our main theorem on the behavior of the solution
up to the singular boundary. 
To prove Theorem~\ref{T:DEVELOPMENTANDSTRUCTUREOFSINGULARBOUNDARY},
we will amalgamate some of our prior results that we derived at fixed 
$\muxmulevelsetvalue \in [0,\muxmulevelsetvalue_0]$.
In this section, we carry out some of these tasks by
constructing a new region $\MInteresting$
that contains the portion of the singular boundary featured in Theorem~\ref{T:DEVELOPMENTANDSTRUCTUREOFSINGULARBOUNDARY}.
We also construct a corresponding time function $\newtimefunction$ that foliates $\MInteresting$,
and our construction is such that the singular boundary portion of interest 
(including the crease) is contained in the level set $\lbrace \newtimefunction = 0 \rbrace$,
which forms the top boundary of $\MInteresting$.
Finally, we derive fundamental properties of $\MInteresting$ and $\newtimefunction$.

Actually, as part of our construction in this section (see Sect.\,\ref{SS:SINGULARBOUNDARYANDCREASE}), 
we rigorously \emph{define} the singular boundary and crease. 
In Theorem~\ref{T:DEVELOPMENTANDSTRUCTUREOFSINGULARBOUNDARY}, 
we reveal the behavior of the solution up to these sets, and the results of the theorem
will justify our definitions.

\subsection{Definitions of developments}
\label{SS:DEFINITIONSFORNEWDEVELOPMENTS}
For each $\muxmulevelsetvalue \in [0,\muxmulevelsetvalue_0]$,
Theorem~\ref{T:EXISTENCEUPTOTHESINGULARBOUNDARYATFIXEDKAPPA}
yields the development
$
\twoargMrough{[\timefunction_0,0],[- \rightu,\leftu]}{\muxmulevelsetvalue}
$
of the data,
which contains the torus
$
\twoargmumuxtorus{0}{-\muxmulevelsetvalue}
$,
which is a subset of the singular boundary.
Using the
$
\twoargMrough{[\timefunction_0,0],[- \rightu,\leftu]}{\muxmulevelsetvalue}
$
as building blocks,
we now define other developments,
including
$\MInteresting$,
on which we will derive
refined estimates
revealing the detailed structure of the singular boundary.
In Prop.\,\ref{SS:PROPERTIESOFSINGULARDEVELOPMENTANDCREASE}, we will 
exhibit various key properties of $\MInteresting$.

	\begin{center}
	\begin{figure}  
		\begin{overpic}[scale=.36, grid = false, tics=5, trim=-.5cm -1cm -1cm -.5cm, clip]{Minterestingmainresults.pdf}
			\put (40,10) {$\Sigma_0$}
			\put (61,18.5) {$(x^2,x^3) \in \mathbb{T}^2$}
			\put (78,28) {$t$}
			\put (62.5,25) {$u \in\mathbb{R}$}
			\put (15,50) {$\MLeft$}
			\put (45,55) {$\MSingular$}
			\put (65,69) {$\MRight$}
			\put (30,25) {$\inthyp{\timefunction}{[- \rightu,\leftu]}$}
			\put (30,100) {$\inthyp{0}{[- \rightu,\leftu]}$}
			\put (34,34) {$\twoargmumuxtorus{-\timefunction}{0}$}
			\put (57,43) {$\twoargmumuxtorus{-\timefunction}{-\muxmulevelsetvalue_0}$}
			\put (26,81) {$\twoargmumuxtorus{0}{0}$}
			\put (49,89.5) {$\twoargmumuxtorus{0}{-\muxmulevelsetvalue}$}
			\put (34.5,80) {$\mathcal{B}^{[0,\muxmulevelsetvalue_0]}$}
			\put (43,35) {$\datahypfortimefunctiontwoarg{0}{[\timefunction_0,0]}$}
			\put (64,48) {$\datahypfortimefunctiontwoarg{-\muxmulevelsetvalue_0}{[\timefunction_0,0]}$}
	\end{overpic}
		\caption{The region $\MInteresting$ featured in Theorems~\ref{T:ABBREVIATEDSTATEMENTOFMAINRESULTS} and 
		\ref{T:DEVELOPMENTANDSTRUCTUREOFSINGULARBOUNDARY}}
	\label{F:MINTERESTINGDEVELOPMENT}
	\end{figure}
\end{center}

\begin{definition}[The development $\MInteresting$ and constituent subsets] 
	\label{D:DEVELOPMENTOFDATA}
	We define the following subsets of geometric coordinate space (see Fig.\,\ref{F:MINTERESTINGDEVELOPMENT}),
	where the subset $\twoargMrough{I,J}{\muxmulevelsetvalue}$ is defined in \eqref{E:TRUNCATEDMROUGH},
	the hypersurface portion 
	$\datahypfortimefunctiontwoarg{-\muxmulevelsetvalue}{[\timefunction_0,0]}$ is defined in \eqref{E:TRUNCATEDLEVELSETSOFMUXMU},
	and to obtain the second equality in \eqref{E:SINGULARDEVELOPMENT}, we used 
	\eqref{E:CLOSEDIMPROVEMENTLEVELSETSTRUCTUREOFMUXEQUALSMINUSKAPPA}:
	\begin{subequations}
	\begin{align} 
	\MLeft
	& \eqdef
	\twoargMrough{[\timefunction_0,0],[\interestingu,\leftu]}{0}
	\cup
	\left(
		\twoargMrough{[\timefunction_0,0],[-\interestingu,\interestingu]}{0}
		\cap
		\lbrace \muX \upmu > 0 \rbrace
	\right),
		 \label{E:LEFTDEVELOPMENT} 
			\\
	\MSingular
	&
	\eqdef
	\bigcup_{\muxmulevelsetvalue \in [0,\muxmulevelsetvalue_0]}
		\datahypfortimefunctiontwoarg{-\muxmulevelsetvalue}{[\timefunction_0,0]}
		= \bigcup_{(\mulevelsetvalue,\muxmulevelsetvalue) \in [0,\mupositive] \times [0,\muxmulevelsetvalue_0]} \twoargmumuxtorus{\mulevelsetvalue}{-\muxmulevelsetvalue},
		 \label{E:SINGULARDEVELOPMENT} 
			\\
\MRight
& \eqdef 
	\twoargMrough{[\timefunction_0,0],[- \rightu,-\interestingu]}{\muxmulevelsetvalue_0}
	\cup
	\left(
		\twoargMrough{[\timefunction_0,0],[-\interestingu,\interestingu]}{\muxmulevelsetvalue_0}
		\cap
		\lbrace \muX \upmu < - \muxmulevelsetvalue_0 \rbrace
	\right),
	\label{E:RIGHTDEVELOPMENT}
		\\
	\MInteresting
	&
	\eqdef
	\MLeft
	\cup
	\MSingular
	\cup
	 \MRight.
	\label{E:INTERESTINGDEVELOPMENTOFDATA}
	\end{align}
	\end{subequations}

	
\end{definition}

	\begin{remark}[$|u| \leq \frac{1}{2}\interestingu$ in $\MSingular$]
		\label{R:SINGULARREGIONHASSMALLUVALUES}
		Note that by \eqref{E:MUXMUKAPPALEVELSETLOCATION}, 
		we have $|u| \leq \frac{1}{2}\interestingu$ 
		in $\MSingular$.
		In our subsequent analysis, we often silently use this fact.
	\end{remark}
	
	\begin{remark}[The smoothness of $\MLeft$ and $\MRight$]
	\label{R:SMOOTHNESSOFMLEFTANDMRIGHT}
		The two subsets $\twoargMrough{[\timefunction_0,0],[\interestingu,\leftu]}{0}$ and
		$\twoargMrough{[\timefunction_0,0],[-\interestingu,\interestingu]}{0}
		\cap
		\lbrace \muX \upmu > 0 \rbrace$,
		whose union defines $\MLeft$,
		join smoothly across their common boundary;
		see Prop.\,\ref{P:PROPERTIESOFMSINGULARANDCREASE}.
		Analogous remarks apply for
		$\MRight$.
	\end{remark}

\subsection{The singular boundary and the crease}
\label{SS:SINGULARBOUNDARYANDCREASE}
We are now ready to define the singular boundary and the crease.
In Prop.\,\ref{SS:PROPERTIESOFSINGULARDEVELOPMENTANDCREASE}, we derive
some crucial properties that these sets enjoy.

\begin{definition}[$\mulevelsettwoarg{\mulevelsetvalue}{[0,\muxmulevelsetvalue_0]}$, the singular boundary portion $\mathcal{B}^{[0,\muxmulevelsetvalue_0]}$, and the crease $\partial_- \mathcal{B}^{[0,\muxmulevelsetvalue_0]}$]
\label{D:SINGULARBOUNDARYPORTION}
For each fixed $\mulevelsetvalue \in [0,\mupositive]$,
we define the set $\mulevelsettwoarg{\mulevelsetvalue}{[0,\muxmulevelsetvalue_0]}$ as follows,
where $\twoargmumuxtorus{\mulevelsetvalue}{-\muxmulevelsetvalue}$ is the $\upmu$-adapted torus defined in \eqref{E:MUXMUTORI}:
\begin{align} \label{E:LEVELSETSOFMUINSINGULARGETION}
	\mulevelsettwoarg{\mulevelsetvalue}{[0,\muxmulevelsetvalue_0]}
	& \eqdef 
	\bigcup_{\muxmulevelsetvalue \in [0,\muxmulevelsetvalue_0]} \twoargmumuxtorus{\mulevelsetvalue}{-\muxmulevelsetvalue}.
\end{align}

	We define the \textbf{singular boundary} portion $\mathcal{B}^{[0,\muxmulevelsetvalue_0]}$ as follows:
	\begin{align} \label{E:SINGULARBOUNDARYPORTION}
		\mathcal{B}^{[0,\muxmulevelsetvalue_0]}
		\eqdef
		\mulevelsettwoarg{0}{[0,\muxmulevelsetvalue_0]}
		=
		\bigcup_{\muxmulevelsetvalue \in [0,\muxmulevelsetvalue_0]} \twoargmumuxtorus{0}{-\muxmulevelsetvalue}.
	\end{align}
	
	We define the \textbf{crease} $\partial_- \mathcal{B}^{[0,\muxmulevelsetvalue_0]}$ as follows:
	\begin{align} \label{E:CREASE}
		\partial_- \mathcal{B}^{[0,\muxmulevelsetvalue_0]} 
		\eqdef
		\twoargmumuxtorus{0}{0}.
	\end{align}
\end{definition}

We note that since $\twoargmumuxtorus{0}{-\muxmulevelsetvalue} \subset \datahypfortimefunctiontwoarg{-\muxmulevelsetvalue}{[\timefunction_0,0]}$,
it follows that $\mathcal{B}^{[0,\muxmulevelsetvalue_0]} \subset \MSingular$ and
$\partial_- \mathcal{B}^{[0,\muxmulevelsetvalue_0]} \subset \MSingular$.
We also note that from definitions~\eqref{E:SINGULARDEVELOPMENT} and \eqref{E:LEVELSETSOFMUINSINGULARGETION}
and \eqref{E:CLOSEDIMPROVEMENTLEVELSETSTRUCTUREOFMUXEQUALSMINUSKAPPA}, it follows
that $\mulevelsettwoarg{\mulevelsetvalue}{[0,\muxmulevelsetvalue_0]}$
is the portion of the level set 
$\lbrace (t,u,x^2,x^3) \ | \ \upmu(t,u,x^2,x^3) = \mulevelsetvalue  \rbrace$ in $\MSingular$.

\subsection{The structure of $\MLeft$, 
$\MRight$,
$\MSingular$,
$\mathcal{B}^{[0,\muxmulevelsetvalue_0]}$, and $\partial_- \mathcal{B}^{[0,\muxmulevelsetvalue_0]}$}
\label{SS:PROPERTIESOFSINGULARDEVELOPMENTANDCREASE}
In the next proposition, we derive key properties of the sets that we defined
in Sect.\,\ref{S:DEVELOPMENTSOFTHEDATASINGULARBOUNDARYNEWTIMEFUNCTION}.

\begin{proposition}[The structure of $\MLeft$, 
$\MRight$,
$\MSingular$,
$\mathcal{B}^{[0,\muxmulevelsetvalue_0]}$, and $\partial_- \mathcal{B}^{[0,\muxmulevelsetvalue_0]}$]
\label{P:PROPERTIESOFMSINGULARANDCREASE}
Assume the hypotheses and conclusions of Theorem~\ref{T:EXISTENCEUPTOTHESINGULARBOUNDARYATFIXEDKAPPA}
for $\muxmulevelsetvalue \in [0,\muxmulevelsetvalue_0]$.
Let
$\twoargroughtori{\timefunction,u}{\muxmulevelsetvalue}$,
$\nullhypthreearg{\muxmulevelsetvalue}{u}{I}$,
$\twoargmumuxtorus{\mulevelsetvalue}{-\muxmulevelsetvalue}$,
and
$\datahypfortimefunctiontwoarg{-\muxmulevelsetvalue}{I}$
be the sets defined in 
\eqref{E:ROUGHTORI},
\eqref{E:NULLHYPERSURFACEROUGHTRUNCATED},
\eqref{E:TRUNCATEDLEVELSETSOFMUXMU},
and \eqref{E:MUXMUTORI}
respectively.
Then the following conclusions hold (see Fig.\,\ref{F:MINTERESTINGDEVELOPMENT}).

\medskip

\noindent \underline{\textbf{Differential-topological structure of} $\MLeft$}.
\begin{itemize}
\item The left lateral boundary of the set
	$\MLeft$ defined in \eqref{E:LEFTDEVELOPMENT}
	is $\nullhypthreearg{0}{\leftu}{[\timefunction_0,0]}$,
	which is a smooth hypersurface with boundary components
	$\twoargroughtori{0,\leftu}{0}$
	and $\twoargroughtori{\timefunction_0,\leftu}{0}$,
	each of which are $C^{2,1}$ graphs over $\mathbb{T}^2$.
\item The right lateral boundary of $\MLeft$
	is
	$\datahypfortimefunctiontwoarg{0}{[\timefunction_0,0]}$,
	which is a $C^{1,1}$ hypersurface with boundary components
	$\twoargmumuxtorus{0}{0}$
	and
	$\twoargmumuxtorus{-\timefunction_0}{0}$,
	each of which are
	$C^{1,1}$ graphs over $\mathbb{T}^2$.
	Moreover, 
	$\datahypfortimefunctiontwoarg{0}{[\timefunction_0,0]}
	=
	\lbrace
		(t,u,x^2,x^3) 
			\,
		\ | \
		(t,x^2,x^3) \in \mathscr{H}_{0}^{[0,\mupositive]},
			\,
		u = h^{(0)}(t,x^2,x^3)
	\rbrace
	$,
	where $h^{(0)}$ is the $C^{1,1}$ function on 
	$\mathscr{H}_{0}^{[0,\mupositive]}
	=
	\lbrace 
		(t,x^2,x^3) \in \mathbb{R} \times \mathbb{T}^2 
		\ | \
		\Cartesiantisafunctiononmumxtoriarg{\mupositive,0}(x^2,x^3) \leq t \leq \Cartesiantisafunctiononmumxtoriarg{0,0}(x^2,x^3)
	\rbrace
	$
	from Cor.\,\ref{C:QUANTITATIVECONTROLOFEMBEDDINGSONCLOSURESOFTHEIRDOMAINS}.
\item	The top boundary of $\MLeft$
	is the $C^{2,1}$ hypersurface
	$
	\lbrace
				(t,u,x^2,x^3)
					\ | \
					(x^2,x^3) \in \mathbb{T}^2,
						\,
					\Eikonalisafunctiononmumuxtoriarg{0}{0}(x^2,x^3) \leq u \leq \leftu,
						\,
					\mbox{and }
					t = \Cartesiantisafunctiononlevelsetsofroughtimefunctionarg{0}{0}(u,x^2,x^3)
	\rbrace
	$,
	where $\Eikonalisafunctiononmumuxtoriarg{0}{0}$ is the $C^{1,1}$ 
	function from Lemma~\ref{L:CHOVFROMROUGHCOORDINATESTOMUWEGIGHTEDXMUCOORDINATES}
	and
	$\Cartesiantisafunctiononlevelsetsofroughtimefunctionarg{0}{0}$ is the $C^{2,1}$ function from Lemma~\ref{L:DIFFEOMORPHICEXTENSIONOFROUGHCOORDINATES}.
	The two boundary components of this hypersurface are
	$\twoargroughtori{0,\leftu}{0}$,
	which is a $C^{2,1}$ graph over $\mathbb{T}^2$,
	and
	$\twoargmumuxtorus{0}{0}$,
	which is a $C^{1,1}$ graph over $\mathbb{T}^2$.
\item	The bottom boundary of $\MLeft$ has an analogous structure:
	it is the $C^{2,1}$ hypersurface
	$
	\lbrace
				(t,u,x^2,x^3)
					\ | \
					(x^2,x^3) \in \mathbb{T}^2,
						\,
					\Eikonalisafunctiononmumuxtoriarg{\timefunction_0}{0}(x^2,x^3) \leq u \leq \leftu,
						\,
					\mbox{and }
					t = \Cartesiantisafunctiononlevelsetsofroughtimefunctionarg{\timefunction_0,0}(u,x^2,x^3)
	\rbrace
	$,
	which has the $C^{1,1}$ boundary component
	$\twoargmumuxtorus{\timefunction_0}{0}$
	and
	the $C^{2,1}$ boundary component
	$\twoargroughtori{\timefunction_0,\leftu}{0}$,
	each of which are graphs over $\mathbb{T}^2$.
\item $\MLeft$ is not closed, but it contains all of its limit points except for the points in its right lateral boundary
		$\datahypfortimefunctiontwoarg{0}{[\timefunction_0,0]}$.
\end{itemize}

\medskip

\noindent \underline{\textbf{Differential-topological structure of} $\MRight$}.
\begin{itemize}
\item The left lateral boundary of the set $\MRight$ defined in \eqref{E:RIGHTDEVELOPMENT}
	is equal to
	$\datahypfortimefunctiontwoarg{-\muxmulevelsetvalue_0}{[\timefunction_0,0]}$,
	which is a $C^{1,1}$ hypersurface with boundary components
	$\twoargmumuxtorus{0}{-\muxmulevelsetvalue_0}$
	and
	$\twoargmumuxtorus{-\timefunction_0}{\muxmulevelsetvalue_0}$,
	each of which are
	$C^{1,1}$ graphs over $\mathbb{T}^2$.
	Moreover, 
	$\datahypfortimefunctiontwoarg{-\muxmulevelsetvalue_0}{[\timefunction_0,0]}
	=
	\lbrace
		(t,u,x^2,x^3) 
			\,
		\ | \
		(t,x^2,x^3) \in \mathscr{H}_{\muxmulevelsetvalue_0}^{[0,\mupositive]},
			\,
		u = h^{(\muxmulevelsetvalue_0)}(t,x^2,x^3)
	\rbrace
	$,
	where $h^{(\muxmulevelsetvalue_0)}$ is the $C^{1,1}$ function on 
	$\mathscr{H}_{\muxmulevelsetvalue_0}^{[0,\mupositive]}
	=
	\lbrace 
		(t,x^2,x^3) \in \mathbb{R} \times \mathbb{T}^2 
		\ | \
		F_{\mupositive,\muxmulevelsetvalue_0}(x^2,x^3) \leq t \leq F_{0,\muxmulevelsetvalue_0}(x^2,x^3)
	\rbrace
	$
	from Cor.\,\ref{C:QUANTITATIVECONTROLOFEMBEDDINGSONCLOSURESOFTHEIRDOMAINS}.
\item The right lateral boundary of
	$\MRight$
	is $\nullhypthreearg{\muxmulevelsetvalue_0}{- \rightu}{[\timefunction_0,0]}$,
	which is a smooth hypersurface with boundary components
	$\twoargroughtori{0,- \rightu}{\muxmulevelsetvalue_0}$
	and $\twoargroughtori{\timefunction_0,- \rightu}{\muxmulevelsetvalue_0}$,
	each of which are $C^{2,1}$ graphs over $\mathbb{T}^2$.
\item	The top boundary of $\MRight$
	is the $C^{2,1}$ hypersurface
	$
	\lbrace
				(t,u,x^2,x^3)
					\ | \
					(x^2,x^3) \in \mathbb{T}^2,
						\,
					 - \rightu \leq u \leq \Eikonalisafunctiononmumuxtoriarg{0}{\muxmulevelsetvalue_0}(x^2,x^3),
						\,
					\mbox{and }
					t = \Cartesiantisafunctiononlevelsetsofroughtimefunctionarg{0}{\muxmulevelsetvalue_0}(u,x^2,x^3)
	\rbrace
	$,
	where $\Eikonalisafunctiononmumuxtoriarg{0}{\muxmulevelsetvalue_0}$ 
	is the $C^{1,1}$ function from Lemma~\ref{L:CHOVFROMROUGHCOORDINATESTOMUWEGIGHTEDXMUCOORDINATES}
	and
	$\Cartesiantisafunctiononlevelsetsofroughtimefunctionarg{0}{\muxmulevelsetvalue_0}$ 
	is the $C^{2,1}$ function from Lemma~\ref{L:DIFFEOMORPHICEXTENSIONOFROUGHCOORDINATES}.
	The two boundary components are
	$\twoargroughtori{0,- \rightu}{\muxmulevelsetvalue_0}$,
	which is a $C^{2,1}$ graph over $\mathbb{T}^2$,
	and
	$\twoargmumuxtorus{0}{-\muxmulevelsetvalue_0}$,
	which is a $C^{1,1}$ graph over $\mathbb{T}^2$.
\item	The bottom boundary of $\MRight$ has an analogous structure:
	it is the $C^{2,1}$ hypersurface
	$
	\lbrace
				(t,u,x^2,x^3)
					\ | \
					(x^2,x^3) \in \mathbb{T}^2,
						\,
					\Eikonalisafunctiononmumuxtoriarg{\timefunction_0}{\muxmulevelsetvalue_0}(x^2,x^3) \leq u \leq \leftu,
						\,
					\mbox{and }
					t = \Cartesiantisafunctiononlevelsetsofroughtimefunctionarg{\timefunction_0,\muxmulevelsetvalue_0}(u,x^2,x^3)
	\rbrace
	$,
	which has the $C^{1,1}$ boundary component
	$\twoargmumuxtorus{\timefunction_0}{\muxmulevelsetvalue_0}$
	and
	the $C^{2,1}$ boundary component
	$\twoargroughtori{\timefunction_0,- \rightu}{\muxmulevelsetvalue_0}$,
	each of which are graphs over $\mathbb{T}^2$.
 \item $\MRight$ is not closed, but it contains all of its limit points except for the points in its left lateral boundary
		$\datahypfortimefunctiontwoarg{-\muxmulevelsetvalue_0}{[\timefunction_0,0]}$.
\end{itemize}

\medskip

\noindent \underline{\textbf{A diffeomorphism onto} $\MSingular$}.
Let $\extendedembeddatahypersurface$
be the map from $[0,\mupositive] \times [0,\muxmulevelsetvalue_0] \times \mathbb{T}^2$ 
into geometric coordinate space $\lbrace (t,u,x^2,x^3) \ | \ t,u \in \mathbb{R}, \, (x^2,x^3) \in \mathbb{T}^2 \rbrace$
defined by:
\begin{align} \label{E:EMBEDDINGOFSINGULARREGION}
	\extendedembeddatahypersurface(\mulevelsetvalue,\muxmulevelsetvalue,x^2,x^3)
	& 
	\eqdef 
	\embeddatahypersurfacearg{\muxmulevelsetvalue}(\mulevelsetvalue,x^2,x^3)
	=
	\left(\Cartesiantisafunctiononmumxtoriarg{\mulevelsetvalue}{-\muxmulevelsetvalue}(x^2,x^3),
	\Eikonalisafunctiononmumuxtoriarg{\mulevelsetvalue}{-\muxmulevelsetvalue}(x^2,x^3),x^2,x^3 \right),
\end{align}
where $\embeddatahypersurfacearg{\muxmulevelsetvalue}$, 
$\Cartesiantisafunctiononmumxtoriarg{\mulevelsetvalue}{-\muxmulevelsetvalue}$,
and 
$\Eikonalisafunctiononmumuxtoriarg{\mulevelsetvalue}{-\muxmulevelsetvalue}$ are the 
functions from \eqref{E:EMBEDDATAHYPERSURFACE}
and Cor.\,\ref{C:QUANTITATIVECONTROLOFEMBEDDINGSONCLOSURESOFTHEIRDOMAINS}.
Then the following results hold.

\begin{itemize}
\item For $(\mulevelsetvalue,\muxmulevelsetvalue) \in [0,\mupositive] \times [0,\muxmulevelsetvalue_0]$, 
we have
$\extendedembeddatahypersurface(\lbrace \mulevelsetvalue \rbrace \times \lbrace \muxmulevelsetvalue \rbrace \times \mathbb{T}^2)
= 
\twoargmumuxtorus{\mulevelsetvalue}{-\muxmulevelsetvalue}
$.
\item
$\extendedembeddatahypersurface$ is a $C^{1,1}$ diffeomorphism from
$[0,\mupositive] \times [0,\muxmulevelsetvalue_0] \times \mathbb{T}^2$
onto the set $\MSingular$ defined in \eqref{E:SINGULARDEVELOPMENT}
such that the following estimate holds:
 \begin{align} \label{E:C11BOUNDFOREMBEDDINGOFSINGULARREGION}
	\| \extendedembeddatahypersurface \|_{C^{1,1}([0,\mupositive] \times [0,\muxmulevelsetvalue_0] \times \mathbb{T}^2)}
	& \leq C.
 \end{align}
\item In particular, the map 
	$\embeddingofsingularboundaryintogeometriccoordinatespace$ defined by:
	\begin{align} \label{E:EMBEDDINGOFSINGULARBOUNDARYINTOGEOMETRICCOORDINATESPACE}
		\embeddingofsingularboundaryintogeometriccoordinatespace(\muxmulevelsetvalue,x^2,x^3)
		& \eqdef \extendedembeddatahypersurface(0,\muxmulevelsetvalue,x^2,x^3)
	\end{align}
	is a $C^{1,1}$ diffeomorphism from $[0,\muxmulevelsetvalue_0] \times \mathbb{T}^2$ onto 
	the singular boundary portion $\mathcal{B}^{[0,\muxmulevelsetvalue_0]}$ defined in \eqref{E:SINGULARBOUNDARYPORTION}
	such that: 
	\begin{align} \label{E:EMBEDDINGOFSINGULARBOUNDARYINTOGEOMETRICCOORDINATESPACEIMAGEOFKAPPAEQUALSCONSTANTISZEROMUMUXTORUS}
		\embeddingofsingularboundaryintogeometriccoordinatespace(\lbrace \muxmulevelsetvalue \rbrace \times \mathbb{T}^2)
		& = \twoargmumuxtorus{0}{-\muxmulevelsetvalue}.
	\end{align}
\item For $(\mulevelsetvalue,\muxmulevelsetvalue) \in [0,\mupositive] \times [0,\muxmulevelsetvalue_0]$ the torus
	$\twoargmumuxtorus{\mulevelsetvalue}{-\muxmulevelsetvalue}$ is a $C^{1,1}$ graph over $\mathbb{T}^2$
	that is $\gfour$-spacelike.
\item There exists a constant $C > 1$ such that for
$(\mulevelsetvalue,\muxmulevelsetvalue) \in [0,\mupositive] \times [0,\muxmulevelsetvalue_0]$, we have:
\begin{subequations}
\begin{align}
	- C 
		& 
		\leq
		\min_{(x^2,x^3) \in \mathbb{T}^2}
		\frac{\partial}{\partial \mulevelsetvalue} \Cartesiantisafunctiononmumxtoriarg{\mulevelsetvalue}{-\muxmulevelsetvalue}(x^2,x^3)
		\leq
		\max_{(x^2,x^3) \in \mathbb{T}^2}
		\frac{\partial}{\partial \mulevelsetvalue} \Cartesiantisafunctiononmumxtoriarg{\mulevelsetvalue}{-\muxmulevelsetvalue}(x^2,x^3)
		\leq - \frac{1}{C},
			 \label{E:TMSINGULAREMBEDDINGFUNCTIONDECREASINGINLAMBDA} \\
		- C
		& \leq 
			\min_{(x^2,x^3) \in \mathbb{T}^2}
			\frac{\partial}{\partial \muxmulevelsetvalue} \Eikonalisafunctiononmumuxtoriarg{\mulevelsetvalue}{-\muxmulevelsetvalue}(x^2,x^3)
			\leq
			\max_{(x^2,x^3) \in \mathbb{T}^2}
			\frac{\partial}{\partial \muxmulevelsetvalue} \Eikonalisafunctiononmumuxtoriarg{\mulevelsetvalue}{-\muxmulevelsetvalue}(x^2,x^3)
			\leq 
			- \frac{1}{C}.
		\label{E:TMSINGULAREMBEDDINGFUNCTIONINCREASINGINKAPPA}
\end{align}
\end{subequations}
\item For each fixed $\mulevelsetvalue \in [0,\mupositive]$,
the map $(\muxmulevelsetvalue,x^2,x^3) \rightarrow \left(\Eikonalisafunctiononmumuxtoriarg{\mulevelsetvalue}{-\muxmulevelsetvalue}(x^2,x^3),x^2,x^3 \right)$
is a $C^{1,1}$ diffeomorphism from
$[0,\muxmulevelsetvalue_0] \times \mathbb{T}^2$ onto its image,
which is: 
\begin{align} \label{E:DOMAINOFGRAPHFORMULEVELSETINSINGULARREGION}
\domainofgraphofmulevelsetinsingularregion{[0,\muxmulevelsetvalue_0]}{\mulevelsetvalue}
\eqdef
\left\lbrace
(u,x^2,x^3) 
 \ | \
(x^2,x^3) \in \mathbb{T}^2,
						\,
\Eikonalisafunctiononmumuxtoriarg{\mulevelsetvalue}{-\muxmulevelsetvalue_0}(x^2,x^3) \leq u \leq \Eikonalisafunctiononmumuxtoriarg{\mulevelsetvalue}{0}(x^2,x^3)
\right\rbrace.
\end{align}

\end{itemize}

\medskip

\noindent \underline{\textbf{$t$ is a function of $(u,x^2,x^3)$ along the level sets of
$\upmu$ in $\MSingular$}}.
Recall that $\mulevelsettwoarg{\mulevelsetvalue}{[0,\muxmulevelsetvalue_0]}$ is the set defined in 
\eqref{E:LEVELSETSOFMUINSINGULARGETION}, 
and that $\mulevelsettwoarg{\mulevelsetvalue}{[0,\muxmulevelsetvalue_0]}$
is the portion of the level set  
$\lbrace (t,u,x^2,x^3) \ | \ \upmu(t,u,x^2,x^3) = \mulevelsetvalue  \rbrace$ in $\MSingular$.
Then the following results hold.

\begin{itemize}
\item
Recall that the set $\domainofgraphofmulevelsetinsingularregion{[0,\muxmulevelsetvalue_0]}{\mulevelsetvalue}$ is defined in \eqref{E:DOMAINOFGRAPHFORMULEVELSETINSINGULARREGION}.
Then for each $\mulevelsetvalue \in [0,\mupositive]$,
there exists a function $\tisafunctiononlevelsetsofmu{\mulevelsetvalue} : \domainofgraphofmulevelsetinsingularregion{[0,\muxmulevelsetvalue_0]}{\mulevelsetvalue} \rightarrow \mathbb{R}$
such that relative to the geometric coordinates, we have:
\begin{align} \label{E:GRAPHSTRUCTUREOFLEVELSETSOFMUINSINGULARGETION}
	\mulevelsettwoarg{\mulevelsetvalue}{[0,\muxmulevelsetvalue_0]}
	& =
	\left\lbrace
				\left(\tisafunctiononlevelsetsofmu{\mulevelsetvalue}(u,x^2,x^3),u,x^2,x^3 \right) 
					\ | \
					(u,x^2,x^3) \in \domainofgraphofmulevelsetinsingularregion{[0,\muxmulevelsetvalue_0]}{\mulevelsetvalue}
	\right\rbrace.
\end{align}
\item There exists a $C > 0$ such that following estimate holds for $\mulevelsetvalue \in [0,\mupositive]$:
	\begin{align} \label{E:C21FORGRAPHFUNCTIONOFLEVELSETSOFMUINSINGULARGETION}
		\| \tisafunctiononlevelsetsofmu{\mulevelsetvalue} \|_{C^{2,1}\left(\domainofgraphofmulevelsetinsingularregion{[0,\muxmulevelsetvalue_0]}{\mulevelsetvalue}\right)} \leq C.
	\end{align}
	In particular, $\mulevelsettwoarg{\mulevelsetvalue}{[0,\muxmulevelsetvalue_0]}$ is a $3$-dimensional $C_{\textnormal{geo}}^{2,1}$ 
	submanifold-with-boundary in geometric coordinate space.
\item The boundary of $\mulevelsettwoarg{\mulevelsetvalue}{[0,\muxmulevelsetvalue_0]}$ in geometric coordinate space 
$\mathbb{R}_t \times \mathbb{R}_u \times \mathbb{T}^2$ satisfies
$\partial \mulevelsettwoarg{\mulevelsetvalue}{[0,\muxmulevelsetvalue_0]} = \twoargmumuxtorus{\mulevelsetvalue}{0} \cup \twoargmumuxtorus{\mulevelsetvalue}{-\muxmulevelsetvalue_0}$,
the $\upmu$-adapted tori
$\twoargmumuxtorus{\mulevelsetvalue}{-\muxmulevelsetvalue}$ are $C^{1,1}$ graphs over $\mathbb{T}^2$,
as indicated in \eqref{E:GRAPHDESCRIPTIONOFMUXMUTORUS}.
\end{itemize}

\medskip

\noindent \underline{\textbf{Differential-topological structure of} $\MSingular$}.
\begin{itemize}
	\item (\textbf{Quasi-convexity}) 
					$\MSingular$ is quasi-convex.
					That is, there is a constant $C > 0$ such that
					every pair of points
					$p_1,p_2 \in \MSingular$
					are connected by a $C_{\textnormal{geo}}^1$ curve in $\MSingular$
					whose length with respect to the standard flat Euclidean metric on 
					geometric coordinate space
					$\mathbb{R}_t \times \mathbb{R}_u \times \mathbb{T}^2$
					is $\leq C \mbox{\upshape dist}_{\mbox{\upshape flat}}(p_1,p_2)$.
	\item (\textbf{Sobolev embedding}).
					There is a constant $C > 0$ 
					such that the following Sobolev embedding result holds for scalar functions $f$ on 
					$\mbox{\upshape int}(\MSingular)$:
					\begin{align} \label{E:MSINGULARSOBOELVEMBEDDINGRELYINGONQUASICONVEXITY}
						\| f \|_{C_{\textnormal{geo}}^{0,1}(\MSingular)}
						& 
						\leq
						C
						\| f \|_{W_{\textnormal{geo}}^{1,\infty}(\mbox{\upshape int}(\MSingular))}.
					\end{align}
	\item Let $\extendedembeddatahypersurface^{-1}$ 
	denote the inverse function of the map $\extendedembeddatahypersurface$ from
	\eqref{E:EMBEDDINGOFSINGULARREGION}, i.e., 
	$\extendedembeddatahypersurface^{-1}(t,u,x^2,x^3) = \left(\upmu, - \muX \upmu,x^2,x^3 \right)$.
		Then $\extendedembeddatahypersurface^{-1}$
		is a $C_{\textnormal{geo}}^{1,1}$ diffeomorphism from 
		$\MSingular$
		onto $[0,\mupositive] \times [0,\muxmulevelsetvalue_0] \times \mathbb{T}^2$.
\item 
		$
		\extendedembeddatahypersurface^{-1}(\twoargmumuxtorus{\mulevelsetvalue}{-\muxmulevelsetvalue}) 
		= \lbrace \mulevelsetvalue \rbrace \times \lbrace - \muxmulevelsetvalue \rbrace \times \mathbb{T}^2
		$.
	\item The following estimates hold:
			\begin{subequations}
			\begin{align}	\label{E:C11BOUNDCHOVGEOTOMUXMUCOORDINSONSINGULARREGION}
				\| \extendedembeddatahypersurface^{-1} \|_{C_{\textnormal{geo}}^{1,1}(\MSingular)} 
				& \leq C,
					\\
				\| \upmu \|_{C_{\textnormal{geo}}^{2,1}(\MSingular)} 
				& \leq C.
				\label{E:C21BOUNDMUSINGULARREGION}
			\end{align}
			\end{subequations}
	\item There exists a constant $C > 0$ such that the following estimates hold on $\MSingular$:
		\begin{align} 
			-\frac{1}{C} 
			& \leq \mbox{\upshape det} d_{\textnormal{geo}} \extendedembeddatahypersurface^{-1}
			\leq - C,
				\label{E:JACOBIANDETERMINANTESTIMATEFORINVERSEEMBEDDINGOFSINGULARREGION} 
					\\
	- 
	\frac{9}{8}
	\mathring{\updelta}_*
	&
	\leq
	\min_{\MSingular} \geop{t} \upmu
	\leq 
	\max_{\MSingular} \geop{t} \upmu
	\leq 
	- 
	\frac{7}{8}
	\mathring{\updelta}_*
				 \label{E:GEOPTMUISNEGATIVEONMSINGULAR} 
					\\
		\frac{\secondtransversalderivativemulowerbound}{2}
		& 
		\leq \min_{\MSingular} \geop{u} \muX \upmu
		\leq \max_{\MSingular} \geop{u} \muX \upmu
		\leq \frac{2}{\secondtransversalderivativemulowerbound}.
				\label{E:GEOPUMUXMUISNEGATIVEONMSINGULAR}
		\end{align}
	\item The two lateral boundaries of $\MSingular$
		are the $C_{\textnormal{geo}}^{1,1}$ embedded hypersurfaces
		$\datahypfortimefunctiontwoarg{0}{[\timefunction_0,0]}$
		and
		$\datahypfortimefunctiontwoarg{-\muxmulevelsetvalue_0}{[\timefunction_0,0]}$
		mentioned above, which have $C_{\textnormal{geo}}^{1,1}$ boundaries.
	\item The top boundary of $\MSingular$
	is
	$\mulevelsettwoarg{0}{[0,\muxmulevelsetvalue_0]}$,
	which is equal to the singular boundary portion $\mathcal{B}^{[0,\muxmulevelsetvalue_0]}$ defined in \eqref{E:SINGULARBOUNDARYPORTION}.
	It is a $C_{\textnormal{geo}}^{2,1}$ embedded hypersurface with the
	boundary components 
	$\twoargmumuxtorus{0}{0}$
	and
	$\twoargmumuxtorus{0}{-\muxmulevelsetvalue_0}$,
	which are $C^{1,1}$ graphs over $\mathbb{T}^2$.
	\item The bottom boundary of $\MSingular$
	is
	$\mulevelsettwoarg{\mupositive}{[0,\muxmulevelsetvalue_0]}$,
	and it is a $C_{\textnormal{geo}}^{2,1}$ embedded hypersurface with the
	boundary components 
	$\twoargmumuxtorus{\mupositive}{0}$
	and
	$\twoargmumuxtorus{\mupositive}{-\muxmulevelsetvalue_0}$,
	which are $C^{1,1}$ graphs over $\mathbb{T}^2$.
\end{itemize}

\end{proposition}

\begin{proof}

\noindent \textbf{Proof of the properties $\MLeft$ and $\MRight$}:
We give the proof only for $\MLeft$ since the properties of $\MRight$
can be derived using similar arguments.
The fact that the left lateral boundary of $\MLeft$
is $\nullhypthreearg{0}{\leftu}{[\timefunction_0,0]}$
follows trivially from definition \eqref{E:LEFTDEVELOPMENT}.
The regularity and structure of its boundary components,
namely
$\twoargroughtori{0,\leftu}{0}$
and $\twoargroughtori{\timefunction_0,\leftu}{0}$,
follow from
the estimate \eqref{E:C21BOUNDFORINVERSECHOVGEOTOROUGH}
and the fact that 
$\twoargroughtori{0,\leftu}{0}
=
\InverseCHOVgeotorough{\muxmulevelsetvalue}\left( \lbrace 0 \rbrace \times \lbrace \leftu \rbrace \times \mathbb{T}^2 \right)$
and
$\twoargroughtori{\timefunction_0,\leftu}{0}
=
\InverseCHOVgeotorough{\muxmulevelsetvalue}\left( \lbrace \timefunction_0 \rbrace \times \lbrace \leftu \rbrace \times \mathbb{T}^2 \right)$.

The fact that the right lateral boundary of $\MLeft$
is $\datahypfortimefunctiontwoarg{0}{[\timefunction_0,0]}$
follows from definition \eqref{E:LEFTDEVELOPMENT},
\eqref{E:MUXMUKAPPALEVELSETLOCATION},
and 
Lemmas~\ref{L:DIFFEOMORPHICEXTENSIONOFROUGHCOORDINATES}
and
Lemma~\ref{L:CHOVFROMROUGHCOORDINATESTOMUWEGIGHTEDXMUCOORDINATES} with $\muxmulevelsetvalue = 0$,
which in particular imply that
$
\datahypfortimefunctiontwoarg{0}{[\timefunction_0,0]}
=
(^{(0)}\mathscr{T})^{-1}
\circ
(\CHOVroughtomumuxmu{0})^{-1}
\left([0,\mupositive] \times \lbrace 0 \rbrace \times \mathbb{T}^2 \right)
\subset 
\mathcal{M}^{(0)}_{[\timefunction_0,0],[-\frac{\interestingu}{2},\frac{\interestingu}{2}]}
$.
We established the regularity and structure of 
$\datahypfortimefunctiontwoarg{0}{[\timefunction_0,0]}$
and its 
boundary components
$\twoargmumuxtorus{0}{0}$
and
$\twoargmumuxtorus{-\timefunction_0}{0}$
in Cor.\,\ref{C:QUANTITATIVECONTROLOFEMBEDDINGSONCLOSURESOFTHEIRDOMAINS} (with $\muxmulevelsetvalue = 0$).

Next, we note that the arguments given in the previous paragraph,
together with \eqref{E:PHIINVERSEIMAGEOFTORUSISTORUSCONTAINEDINROUGHTIMEFUNCTIONLEVELSET},
\eqref{E:GRAPHDESCRIPTIONOFMUXMUTORUS},
and
\eqref{E:LEVELSETSOFTIMEFUNCTIONAREAGRAPH},
imply that the top boundary of $\MLeft$
is the subset of 
$\hypthreearg{0}{[- \rightu,\leftu]}{0}
=					\lbrace
							(t,u,x^2,x^3)
							\ | \
							t = \Cartesiantisafunctiononlevelsetsofroughtimefunctionarg{0}{0}(u,x^2,x^3),
								\,
							(u,x^2,x^3) \in [- \rightu,\leftu] \times \mathbb{T}^2
						\rbrace
$
in which $\Eikonalisafunctiononmumuxtoriarg{0}{0}(x^2,x^3) \leq u \leq \leftu$.
In particular, the two boundary components of the subset under consideration
are
$\hypthreearg{0}{[- \rightu,\leftu]}{0} \cap \lbrace u = \leftu \rbrace = \twoargroughtori{0,\leftu}{0}$
and 
$\hypthreearg{0}{[- \rightu,\leftu]}{0} \cap \lbrace (t,u,x^2,x^3) \ | \ 
u 
= \Eikonalisafunctiononmumuxtoriarg{0}{0}(x^2,x^3) \rbrace = \twoargmumuxtorus{0}{0}$,
where to obtain the last identity, we used
\eqref{E:LASTSLICETORIAREGRAPHSABOVEFLATTORIINGEOMETRICCOORDINATES}.
We derived the $C^{2,1}$ regularity of $\Cartesiantisafunctiononlevelsetsofroughtimefunctionarg{0}{0}$
in \eqref{E:C21FORGRAPHFUNCTIONOFLEVELSETSOFMUINSINGULARGETION}.
We derived the $C^{2,1}$ regularity of
$\twoargroughtori{0,\leftu}{0}$
and the $C^{1,1}$ regularity of $\twoargmumuxtorus{0}{0}$ earlier in the proof.
We have therefore established the claimed properties of 
$\MLeft$.
As we mentioned earlier,
the regularity and structure of the 
bottom boundary of $\MRight$
can be derived using similar arguments.

\medskip

\noindent \textbf{Proof of the properties of $\extendedembeddatahypersurface$ and $\embeddingofsingularboundaryintogeometriccoordinatespace$}:
To obtain these results, we will study the following map:
\begin{align} \label{E:MAPFROMGEOMETRICTOMUMUXMUCOORDINATESONMSINGULAR}
\widetilde{\breve{\mathscr{M}}}(t,u,x^2,x^3) 
& \eqdef \left(\upmu,-\muX \upmu,x^2,x^3 \right)
\end{align}
on the domain $\MSingular$ in geometric coordinate space. 
Note the sign difference of $-\muX \upmu$ on RHS~\eqref{E:MAPFROMGEOMETRICTOMUMUXMUCOORDINATESONMSINGULAR} compared to
the definition~\eqref{E:CHOVFROMGEOMETRICCOORDINATESTOMUWEGIGHTEDXMUCOORDINATES}
of $\CHOVgeotomumuxmu$; we inserted the minus sign on 
RHS~\eqref{E:MAPFROMGEOMETRICTOMUMUXMUCOORDINATESONMSINGULAR} to facilitate the discussion in parts of this proof. 
We will show that $\widetilde{\breve{\mathscr{M}}}$ is a diffeomorphism from
$\MSingular$ onto $[0,\mupositive] \times [0,\muxmulevelsetvalue_0] \times \mathbb{T}^2$,
and our proof will show that the desired map $\extendedembeddatahypersurface$ is equal to $\widetilde{\breve{\mathscr{M}}}^{-1}$. 
To proceed, we first compose the maps $\CHOVgeotorough{\muxmulevelsetvalue}$ and $\CHOVroughtomumuxmu{\muxmulevelsetvalue}$
(recall \eqref{E:SIMPLERELATIONSHIPBETWEENCHOVMAPS}, which states that 
$\CHOVroughtomumuxmu{\muxmulevelsetvalue} \circ \CHOVgeotorough{\muxmulevelsetvalue}(t,u,x^2,x^3) 
= \CHOVgeotomumuxmu(t,u,x^2,x^3) 
= 
(\upmu,\muX \upmu,x^2,x^3)$)
and use
Lemmas~\ref{L:DIFFEOMORPHICEXTENSIONOFROUGHCOORDINATES} and \ref{L:CHOVFROMROUGHCOORDINATESTOMUWEGIGHTEDXMUCOORDINATES}
and \eqref{E:CLOSEDIMPROVEMENTLEVELSETSTRUCTUREOFMUXEQUALSMINUSKAPPA}
to deduce that for each fixed $\muxmulevelsetvalue \in [0,\muxmulevelsetvalue_0]$, 
$\CHOVgeotomumuxmu$ is a $C_{\textnormal{geo}}^{1,1}$ diffeomorphism on a 
subset of $\twoargMrough{[\timefunction_0,0],[- \rightu,\leftu]}{\muxmulevelsetvalue}$ 
containing
$\datahypfortimefunctiontwoarg{-\muxmulevelsetvalue}{[\timefunction_0,0]} 
= 
\bigcup_{\mulevelsetvalue \in [0,\mupositive]} \twoargmumuxtorus{\mulevelsetvalue}{-\muxmulevelsetvalue}$.
Hence, using \eqref{E:PHIINVERSEIMAGEOFTORUSCROSSMUINTERVALISTORUSCROSSINTERVALCONTAINEDININTERESTINGREGION}
with $\upmuboot = 0$ (which is allowable by Theorem~\ref{T:EXISTENCEUPTOTHESINGULARBOUNDARYATFIXEDKAPPA}),
definition~\eqref{E:SINGULARDEVELOPMENT},
and Lemma~\ref{L:GEOTOMUMUXCOORDINATESJACOBIANDETERMINANTRATIOINEVOLUTIONBOUND},
we see (accounting for the minus sign difference between $\CHOVgeotomumuxmu$ and $\widetilde{\breve{\mathscr{M}}}$)
that $\widetilde{\breve{\mathscr{M}}}$ is a diffeomorphism from
$\MSingular$ onto $[0,\mupositive] \times [0,\muxmulevelsetvalue_0] \times \mathbb{T}^2$
whose Jacobian determinant satisfies the bound 
$
\mbox{\upshape det} \widetilde{\breve{\mathscr{M}}}
\approx 1
$
on $\MSingular$ and
such that for $(\mulevelsetvalue,\muxmulevelsetvalue) \in [0,\mupositive] \times [0,\muxmulevelsetvalue_0]$, 
we have
$\widetilde{\breve{\mathscr{M}}}(\twoargmumuxtorus{\mulevelsetvalue}{-\muxmulevelsetvalue}) 
= 
\lbrace \upmu \rbrace \times \lbrace \muxmulevelsetvalue \rbrace \times \mathbb{T}^2$.
Also using the H\"{o}lder estimates provided by Lemma~\ref{L:CONTINUOUSEXTNESION} and Rademacher's theorem, 
we see that $\| \widetilde{\breve{\mathscr{M}}} \|_{W_{\textnormal{geo}}^{2,\infty}(\mbox{\upshape int}(\MSingular))} \lesssim 1$.
From \eqref{E:BAMUTORI}, 
Cor.\,\ref{C:QUANTITATIVECONTROLOFEMBEDDINGSONCLOSURESOFTHEIRDOMAINS},
and \eqref{E:PHIINVERSEIMAGEOFTORUSISTORUSCONTAINEDINROUGHTIMEFUNCTIONLEVELSET},
it follows that the inverse map $\widetilde{\breve{\mathscr{M}}}^{-1}$ is precisely the map $\extendedembeddatahypersurface$ 
defined in \eqref{E:EMBEDDINGOFSINGULARREGION}.
Thus, with $\extendedembeddatahypersurface$ denoting the inverse function of $\widetilde{\breve{\mathscr{M}}}$, 
we can use these estimates and
differentiate the identity $\extendedembeddatahypersurface \circ \widetilde{\breve{\mathscr{M}}}(t,u,x^2,x^3) = (t,u,x^2,x^3)$
up to two times to deduce that
$\| \extendedembeddatahypersurface \|_{W^{2,\infty}((0,\mupositive) \times (0,\muxmulevelsetvalue_0) \times \mathbb{T}^2)}
\lesssim 1
$. Since $(0,\mupositive) \times (0,\muxmulevelsetvalue_0) \times \mathbb{T}^2$ is convex,
we further deduce from Sobolev embedding that
$\| \extendedembeddatahypersurface \|_{C^{1,1}([0,\mupositive] \times [0,\muxmulevelsetvalue_0] \times \mathbb{T}^2)}
\lesssim 1
$,
which yields \eqref{E:C11BOUNDFOREMBEDDINGOFSINGULARREGION}.
The properties of the $\embeddingofsingularboundaryintogeometriccoordinatespace$
defined in \eqref{E:EMBEDDINGOFSINGULARBOUNDARYINTOGEOMETRICCOORDINATESPACE},
including \eqref{E:EMBEDDINGOFSINGULARBOUNDARYINTOGEOMETRICCOORDINATESPACEIMAGEOFKAPPAEQUALSCONSTANTISZEROMUMUXTORUS},
follow from above arguments.

\medskip
\noindent \textbf{Proof of the quasi-convexity of $\MSingular$ and \eqref{E:MSINGULARSOBOELVEMBEDDINGRELYINGONQUASICONVEXITY}}:
First, we note that Lemma~\ref{L:GEOTOMUMUXCOORDINATESJACOBIANDETERMINANTRATIOINEVOLUTIONBOUND} 
(in particular the Jacobian estimate \eqref{E:GEOTOMUMUXCOORDINATESJACOBIANDETERMINANTRATIOINEVOLUTIONBOUND}),
the estimates of Prop.\,\ref{P:IMPROVEMENTOFAUXILIARYBOOTSTRAP},
and the convexity of $[0,\mupositive] \times [-\muxmulevelsetvalue_0,0] \times \mathbb{T}^2$
imply that for every pair of points 
$q_1, q_2 \in [0,\mupositive] \times [-\muxmulevelsetvalue_0,0] \times \mathbb{T}^2$,
we have the following estimates:
$
 \mbox{\upshape dist}_{\mbox{\upshape flat}}(q_1,q_2)
\approx
\mbox{\upshape dist}_{\mbox{\upshape flat}}\left(\InverseCHOVgeotomumuxmu(q_1),\InverseCHOVgeotomumuxmu(q_2) \right)
\in
\MSingular
$,
where $\mbox{\upshape dist}_{\mbox{\upshape flat}}(q_1,q_2)$ is the
standard Euclidean distance between $q_1$ and $q_2$ in the flat space 
$\mathbb{R} \times \mathbb{R} \times \mathbb{T}^2$,
and
$
\mbox{\upshape dist}_{\mbox{\upshape flat}}\left(\InverseCHOVgeotomumuxmu(q_1),\InverseCHOVgeotomumuxmu(q_2) \right)
$
is the standard Euclidean distance between $\InverseCHOVgeotomumuxmu(q_1)$ and $\InverseCHOVgeotomumuxmu(q_1)$ in the flat space 
$\mathbb{R}_t \times \mathbb{R}_u \times \mathbb{T}^2$.
From this bound, the convexity of $[0,\mupositive] \times [-\muxmulevelsetvalue_0,0] \times \mathbb{T}^2$,
and the estimates of Prop.\,\ref{P:IMPROVEMENTOFAUXILIARYBOOTSTRAP},
we conclude that $\MSingular = \InverseCHOVgeotomumuxmu\left( [0,\mupositive] \times [-\muxmulevelsetvalue_0,0] \times \mathbb{T}^2 \right)$
is quasi-convex in the sense stated in the proposition.
From this quasi-convexity, it is a standard result 
(see, for example, \cite[Theorem~7]{pHpKhT2008}),
that the Sobolev embedding result
\eqref{E:MSINGULARSOBOELVEMBEDDINGRELYINGONQUASICONVEXITY}
holds on $\MSingular$,
where the constant $C$ 
on RHS~\eqref{E:MSINGULARSOBOELVEMBEDDINGRELYINGONQUASICONVEXITY}
depends on the constant (a different one, also called $C$) in the 
definition of quasi-convexity.

\medskip
\noindent \textbf{Proof of 
\eqref{E:JACOBIANDETERMINANTESTIMATEFORINVERSEEMBEDDINGOFSINGULARREGION}--\eqref{E:GEOPUMUXMUISNEGATIVEONMSINGULAR}}:
Since $\extendedembeddatahypersurface^{-1} = \widetilde{\mathcal{M}}$,
the estimate \eqref{E:JACOBIANDETERMINANTESTIMATEFORINVERSEEMBEDDINGOFSINGULARREGION} follows from 
the fact that on $\MSingular$, we have 
$
\mbox{\upshape det} \widetilde{\breve{\mathscr{M}}}
\approx 1
$, 
as we showed above.

\eqref{E:GEOPTMUISNEGATIVEONMSINGULAR} follows from the estimate \eqref{E:BOUNDSONLMUINTERESTINGREGION}
(which holds on 
$\mathcal{M}^{(\muxmulevelsetvalue)}_{[\timefunction_0,0],[-\frac{1}{2}\interestingu,\frac{1}{2}\interestingu]}$
for each $\muxmulevelsetvalue \in [0,\muxmulevelsetvalue_0]$),
the definition $\MSingular$ of \eqref{E:SINGULARDEVELOPMENT},
and the fact that 
$\datahypfortimefunctiontwoarg{-\muxmulevelsetvalue}{[\timefunction_0,0]}
\subset \mathcal{M}^{(\muxmulevelsetvalue)}_{[\timefunction_0,0],[-\frac{1}{2}\interestingu,\frac{1}{2}\interestingu]}
$
by \eqref{E:MUXMUKAPPALEVELSETLOCATION}.
From similar reasoning, based on the estimate \eqref{E:MUTRANSVERSALCONVEXITY}, 
we conclude \eqref{E:GEOPUMUXMUISNEGATIVEONMSINGULAR}.

\medskip
\noindent \textbf{Proof of \eqref{E:TMSINGULAREMBEDDINGFUNCTIONDECREASINGINLAMBDA}--\eqref{E:TMSINGULAREMBEDDINGFUNCTIONINCREASINGINKAPPA}}:
We define the vectorfields $\muderivativevectorfield$ and $\muxmuderivativevectorfield$ as follows:
\begin{align}
\muderivativevectorfield
& \eqdef 
\frac{1}
{\geop{t} \upmu 
	- 
	\frac{(\geop{t} \muX \upmu) \geop{u} \upmu}
		{\geop{u} \muX \upmu}
} 
\left\lbrace 
	\geop{t} 
	- 
	\frac{\geop{t} \muX \upmu}
	{\geop{u} \muX \upmu} 
	\geop{u} 
\right\rbrace,
	 \label{E:PARTIALDERIVATIVEWITHRESPECTTOMUATFIXEDMUXMU} 
		\\
\muxmuderivativevectorfield 
& \eqdef 
\frac{-1}
{\geop{u} \muX \upmu 
	- 
	\frac{(\geop{u} \upmu) \geop{t} \muX \upmu}
	{\geop{t} \upmu}
} 
\left\lbrace 
	\geop{u} 
	- 
	\frac{\geop{u} \upmu}
	{\geop{t} \upmu} 
	\geop{t}
\right\rbrace.
	\label{E:PARTIALDERIVATIVEWITHRESPECTTOMUXMUATFIXEDMU}
\end{align}	
From definitions~\eqref{E:PARTIALDERIVATIVEWITHRESPECTTOMUATFIXEDMUXMU}--\eqref{E:PARTIALDERIVATIVEWITHRESPECTTOMUXMUATFIXEDMU} 
and straightforward computations, we find that:
\begin{align} \label{E:RELATIONSPARTIALDERIVATIVEWITHRESPECTTOMUATFIXEDMUXMU}
	\muderivativevectorfield \upmu 
	& = 1, 
	& &
	-
	\muderivativevectorfield \muX \upmu 
	= 
	\muderivativevectorfield x^2 
	= 
	\muderivativevectorfield x^3 
	= 0. 
\end{align}
Hence,
$
\muderivativevectorfield$ is the partial derivative with respect to $\upmu$ in the coordinate system $(\upmu,- \muX \upmu,x^2,x^3)$
(corresponding to RHS~\eqref{E:MAPFROMGEOMETRICTOMUMUXMUCOORDINATESONMSINGULAR})
on the region $[0,\mupositive] \times [0,\muxmulevelsetvalue_0] \times \mathbb{T}^2$.
Similarly, we compute that:
\begin{align} \label{E:RELATIONSPARTIALDERIVATIVEWITHRESPECTTOMUXMUATFIXEDMU}
	-
	\muxmuderivativevectorfield \muX \upmu 
	& = 1, 
	& &
	\muxmuderivativevectorfield \upmu 
	= 
	\muxmuderivativevectorfield x^2 
	= 
	\muxmuderivativevectorfield  x^3 
	= 0,
\end{align}
and thus
$\muxmuderivativevectorfield$ is the partial derivative with respect to 
$-\muX \upmu$ in the coordinate system $(\upmu,- \muX \upmu,x^2,x^3)$.
From \eqref{E:MAPFROMGEOMETRICTOMUMUXMUCOORDINATESONMSINGULAR}
and the inverse function theorem,
we compute that the $2 \times 2$ upper left-hand block of the matrix 
$[d_{\textnormal{geo}} \widetilde{\breve{\mathscr{M}}}]^{-1}$
is equal to
$
\begin{pmatrix} 
				\muderivativevectorfield t  & \muxmuderivativevectorfield  t \\
			 \muderivativevectorfield u & \muxmuderivativevectorfield  u
\end{pmatrix}
$.
Using \eqref{E:PARTIALDERIVATIVEWITHRESPECTTOMUATFIXEDMUXMU}--\eqref{E:PARTIALDERIVATIVEWITHRESPECTTOMUXMUATFIXEDMU},
\eqref{E:MUTRANSVERSALCONVEXITY},
and
\eqref{E:BOUNDSONGEOMETRICTDERIVATIVEMUINTERESTINGREGION},
and the fact that
$
\widetilde{\breve{\mathscr{M}}}$
is a diffeomorphism from $\MSingular$
(see also Remark~\ref{R:SINGULARREGIONHASSMALLUVALUES})
onto $[0,\mupositive] \times [0,\muxmulevelsetvalue_0] \times \mathbb{T}^2$,
we deduce that the diagonal entries of the $2 \times 2$ upper left-hand block of the matrix 
$[d_{\textnormal{geo}} \widetilde{\breve{\mathscr{M}}}]^{-1}$
satisfy the following estimates on $[0,\mupositive] \times [0,\muxmulevelsetvalue_0] \times \mathbb{T}^2$,
where $C > 1$:
\begin{align} \label{E:MONOTONICITYESTIMATES}
-C & \leq \muderivativevectorfield t \leq - \frac{1}{C},
& 
&
-C \leq \muxmuderivativevectorfield  u \leq - \frac{1}{C}.
\end{align}
Since $\extendedembeddatahypersurface = \widetilde{\breve{\mathscr{M}}}^{-1}$,
we see that estimates in \eqref{E:MONOTONICITYESTIMATES} are precisely 
\eqref{E:TMSINGULAREMBEDDINGFUNCTIONDECREASINGINLAMBDA}--\eqref{E:TMSINGULAREMBEDDINGFUNCTIONINCREASINGINKAPPA}.

\medskip
\noindent \textbf{Proof of the diffeomorphism property of the map $(\muxmulevelsetvalue,x^2,x^3) \rightarrow \left(\Eikonalisafunctiononmumuxtoriarg{\mulevelsetvalue}{-\muxmulevelsetvalue}(x^2,x^3),x^2,x^3 \right)$}:
The $\muxmulevelsetvalue$-monotonicity of $\Eikonalisafunctiononmumuxtoriarg{\mulevelsetvalue}{-\muxmulevelsetvalue}(x^2,x^3)$ provided by 
\eqref{E:TMSINGULAREMBEDDINGFUNCTIONINCREASINGINKAPPA}
yields the diffeomorphism property of the map
$(\muxmulevelsetvalue,x^2,x^3) \rightarrow \left(\Eikonalisafunctiononmumuxtoriarg{\mulevelsetvalue}{-\muxmulevelsetvalue}(x^2,x^3),x^2,x^3 \right)$
onto the set $\domainofgraphofmulevelsetinsingularregion{[0,\muxmulevelsetvalue_0]}{\mulevelsetvalue}$ defined in \eqref{E:DOMAINOFGRAPHFORMULEVELSETINSINGULARREGION}.

\medskip
\noindent \textbf{Proof of \eqref{E:C11BOUNDCHOVGEOTOMUXMUCOORDINSONSINGULARREGION} and \eqref{E:C21BOUNDMUSINGULARREGION}}:
Since the estimates of Prop.\,\ref{P:IMPROVEMENTOFAUXILIARYBOOTSTRAP} imply that
$\| \upmu \|_{W_{\textnormal{geo}}^{3,\infty}(\mbox{\upshape int}(\MSingular))} \lesssim 1$,
we can use \eqref{E:MSINGULARSOBOELVEMBEDDINGRELYINGONQUASICONVEXITY} to conclude 
\eqref{E:C21BOUNDMUSINGULARREGION}.
From \eqref{E:C21BOUNDMUSINGULARREGION} and \eqref{E:MAPFROMGEOMETRICTOMUMUXMUCOORDINATESONMSINGULAR},
it follows that
$\| \widetilde{\breve{\mathscr{M}}} \|_{C_{\textnormal{geo}}^{1,1}(\MSingular)} \leq C$. 
Since $\extendedembeddatahypersurface^{-1} = \widetilde{\breve{\mathscr{M}}}$,
we conclude \eqref{E:C11BOUNDCHOVGEOTOMUXMUCOORDINSONSINGULARREGION}.

\medskip

\noindent \textbf{Proof of the existence of $\tisafunctiononlevelsetsofmu{\mulevelsetvalue}$ and 
the properties of $\tisafunctiononlevelsetsofmu{\mulevelsetvalue}$ and
$\mulevelsettwoarg{\mulevelsetvalue}{[0,\muxmulevelsetvalue_0]}$}.
In view of \eqref{E:LEVELSETSOFMUINSINGULARGETION}
and the fact that 
$\extendedembeddatahypersurface(\lbrace \mulevelsetvalue \rbrace \times \lbrace \muxmulevelsetvalue \rbrace \times \mathbb{T}^2)
= 
\twoargmumuxtorus{\mulevelsetvalue}{-\muxmulevelsetvalue}
$,
we observe that the function $\tisafunctiononlevelsetsofmu{\mulevelsetvalue}$ is the first component of the composition of the
map $(\muxmulevelsetvalue,x^2,x^3) \rightarrow \extendedembeddatahypersurface(\mulevelsetvalue,\muxmulevelsetvalue,x^2,x^3)$
(the map has the domain \eqref{E:DOMAINOFGRAPHFORMULEVELSETINSINGULARREGION} and the
first component is denoted by $\Cartesiantisafunctiononmumxtoriarg{\mulevelsetvalue}{-\muxmulevelsetvalue}$ on RHS~\eqref{E:EMBEDDINGOFSINGULARREGION})
with the inverse of the map 
$(\muxmulevelsetvalue,x^2,x^3) \rightarrow \left(\Eikonalisafunctiononmumuxtoriarg{\mulevelsetvalue}{-\muxmulevelsetvalue}(x^2,x^3),x^2,x^3 \right)$.
This fact yields \eqref{E:GRAPHSTRUCTUREOFLEVELSETSOFMUINSINGULARGETION}.
To deduce the bound \eqref{E:C21FORGRAPHFUNCTIONOFLEVELSETSOFMUINSINGULARGETION},
we differentiate the identity $\upmu\left(\tisafunctiononlevelsetsofmu{\mulevelsetvalue}(u,x^2,x^3),u,x^2,x^3 \right) = \mulevelsetvalue$,
which holds on $\domainofgraphofmulevelsetinsingularregion{[0,\muxmulevelsetvalue_0]}{\mulevelsetvalue}$,
and use the chain rule,
\eqref{E:BOUNDSONGEOMETRICTDERIVATIVEMUINTERESTINGREGION},
and
\eqref{E:C21BOUNDMUSINGULARREGION}.
Also using the bounds 
$
\sup_{(\mulevelsetvalue,\muxmulevelsetvalue) \in [0,\mupositive] \times [0,\muxmulevelsetvalue_0]} 
\| \Cartesiantisafunctiononmumxtoriarg{\mulevelsetvalue}{-\muxmulevelsetvalue}\|_{C^{1,1}(\mathbb{T}^2)}
\leq C
$
and
$
\sup_{(\mulevelsetvalue,\muxmulevelsetvalue) \in [0,\mupositive] \times [0,\muxmulevelsetvalue_0]} 
\| \Eikonalisafunctiononmumuxtoriarg{\mulevelsetvalue}{-\muxmulevelsetvalue} \|_{C^{1,1}(\mathbb{T}^2)}
\leq C
$
implied by \eqref{E:C11BOUNDFOREMBEDDINGOFSINGULARREGION},
we conclude that for each fixed $\mulevelsetvalue \in [0,\mupositive]$,
the set $\mulevelsettwoarg{\mulevelsetvalue}{[0,\muxmulevelsetvalue_0]}$ is a $C_{\textnormal{geo}}^{2,1}$
embedded hypersurface with the $C_{\textnormal{geo}}^{1,1}$ boundary
components equal to
$\twoargmumuxtorus{\mulevelsetvalue}{0}$ 
and
$\twoargmumuxtorus{\mulevelsetvalue}{-\muxmulevelsetvalue_0}$.

\medskip
\noindent \textbf{Proof of the properties of the boundary $\MSingular$}:
Since $\extendedembeddatahypersurface$ is a diffeomorphism from
$[0,\mupositive] \times [0,\muxmulevelsetvalue_0] \times \mathbb{T}^2$ onto $\MSingular$,
and since diffeomorphisms map boundaries to boundaries,
it follows that the boundary of $\MSingular$
is the union of four sets:
$\extendedembeddatahypersurface([0,\mupositive] \times \lbrace 0 \rbrace \times \mathbb{T}^2)$,
$\extendedembeddatahypersurface(\lbrace 0 \rbrace \times [0,\muxmulevelsetvalue_0] \times \mathbb{T}^2)$,
$\extendedembeddatahypersurface([0,\mupositive] \times \lbrace \muxmulevelsetvalue_0 \rbrace \times \mathbb{T}^2)$,
and
$\extendedembeddatahypersurface(\lbrace 0 \rbrace \times \upmu_0 \times \mathbb{T}^2)$.
In view of the form \eqref{E:MAPFROMGEOMETRICTOMUMUXMUCOORDINATESONMSINGULAR}
of the map $\widetilde{\breve{\mathscr{M}}}$ (which is equal to $\extendedembeddatahypersurface^{-1}$),
we conclude (recalling that $\timefunction_0 = - \mupositive$)
that the boundary of $\MSingular$ is the union of the four sets
$\datahypfortimefunctiontwoarg{0}{[\timefunction_0,0]}$,
$\mulevelsettwoarg{0}{[0,\muxmulevelsetvalue_0]} = \mathcal{B}^{[0,\muxmulevelsetvalue_0]}$,
$\datahypfortimefunctiontwoarg{-\muxmulevelsetvalue_0}{[\timefunction_0,0]}$,
and $\mulevelsettwoarg{\mupositive}{[0,\muxmulevelsetvalue_0]}$,
as is stated in the proposition. 
The regularity properties of these boundary portions was derived earlier in the proof.
Finally, using 
\eqref{E:GEOPTMUISNEGATIVEONMSINGULAR}--\eqref{E:GEOPUMUXMUISNEGATIVEONMSINGULAR},
we find that the four sets mentioned above are respectively
(see Fig.\,\ref{F:MINTERESTINGDEVELOPMENT})
the left lateral boundary of $\MSingular$,
the top boundary of $\MSingular$,
the right lateral boundary of $\MSingular$,
and
the bottom boundary of $\MSingular$.

This completes our proof of the proposition.

\end{proof}

\subsection{The character of 
$\mulevelsettwoarg{\mulevelsetvalue}{[0,\muxmulevelsetvalue_0]}$,
$\twoargmumuxtorus{\mulevelsetvalue}{-\muxmulevelsetvalue}$,
 and $\partial_- \mathcal{B}^{[0,\muxmulevelsetvalue_0]}$}
\label{SS:CAUSALSTRUCTUREOFSINGULARBOUNDARYANDCREASE}
In this section, we study the character of various submanifolds of geometric coordinate space, i.e.,
whether they are $\gfour$-timelike, null, or spacelike.
The singular boundary is degenerate for reasons discussed in Remark~\ref{R:ACOUSTICALMETRICDEGENERACIESALONGSINGULARBOUNDARY}. 
Hence, we postpone our investigation of the character of the singular boundary until
Prop.\,\ref{P:DESCRIPTIONOFSINGULARBOUNDARYINCARTESIANSPACE},
where we describe how it is embedded into the physical Cartesian coordinate space
equipped with the acoustical metric $\gfour$.

We start with the following simple lemma, 
which provides various identities involving the gradient 
vectorfield of $\upmu$.

\begin{lemma}[Identities involving $\Dfour^{\#}\upmu$]
	\label{L:IDENTITIESFORGRADIENTVECTORFIELDOFINVERSEFOLIATIONDENSITY}
	Let $\Dfour \upmu$ denote the gradient one-form of $\upmu$, and let  
	$\Dfour^{\#} \upmu$ denote the $\gfour$-dual vectorfield of the gradient one-form. Then
	the following identities hold:
	\begin{subequations}
	\begin{align}  \label{E:GRADIENTVECTORFIELDOFMU}
		\upmu \Dfour^{\#} \upmu
		& = 
		- 
		(\upmu \Lunit \upmu + \muX \upmu) \Lunit
		-
		(\Lunit \upmu) \muX
		+
		\upmu \angDuparg{\#} \upmu,
		\\
	\upmu
	\gfour(\Dfour^{\#} \upmu, \Dfour^{\#} \upmu)
	& = 
	\upmu
		(\gfour^{-1})^{\alpha \beta}
		(\partial_{\alpha} \upmu)
		\partial_{\beta} \upmu
	= 
			-
			2 (\Lunit \upmu) \muX \upmu
			- 
			\upmu 
			\left\lbrace
				(\Lunit \upmu)^2
				-
				|\angD \upmu|_{\gtorus}^2
			\right\rbrace.
			\label{E:MUWEGHTEDSPACETIMEGRADMUNORMSQUARED}
	\end{align}
	\end{subequations}
	
	In particular, if $0 < \upmu' \leq \mupositive$,
	$q \in \lbrace  (t,u,x^2,x^3) \ | \ \upmu(t,u,x^2,x^3) = \upmu' \rbrace\cap \lbrace |u| \leq \interestingu \rbrace$,
	and if $\muX \upmu|_q \leq 0$,
	then since \eqref{E:BOUNDSONLMUINTERESTINGREGION} and \eqref{E:POINTWISEBOUNDROUGHGRADIENTOFMU}
	imply that $\mbox{RHS~\eqref{E:MUWEGHTEDSPACETIMEGRADMUNORMSQUARED}} < 0$ at $q$, it follows that 
	$\lbrace  (t,u,x^2,x^3) \ | \ \upmu(t,u,x^2,x^3) = \upmu' \rbrace$ is $\gfour$-spacelike 
	(i.e., $\gfour(\Dfour^{\#} \upmu, \Dfour^{\#} \upmu) < 0$)
	at $q$.
	

\end{lemma}

\begin{proof}
\eqref{E:GRADIENTVECTORFIELDOFMU}--\eqref{E:MUWEGHTEDSPACETIMEGRADMUNORMSQUARED} 
follow from a straightforward computation based on the fact that
$(\Dfour^{\#})^{\alpha} \eqdef (\gfour^{-1})^{\alpha \beta} \partial_{\beta}$ 
and the identity
$(\gfour^{-1})^{\alpha \beta} 
= - \Lunit^{\alpha} \Lunit^{\beta} - X^{\alpha}\Lunit^{\beta} - \Lunit^{\alpha} X^\beta + (\gtorus^{-1})^{\alpha \beta}$, 
which follows from \eqref{E:SMOOTHTORUSINVERSEMETRICINTERMSOFINVERSESIGMATMETRICANDX}. 
\end{proof}

In the next lemma, we exhibit the character of various submanifolds of geometric coordinate space.

\begin{lemma}[The character of 
$\mulevelsettwoarg{\mulevelsetvalue}{[0,\muxmulevelsetvalue_0]}$,
$\twoargmumuxtorus{\mulevelsetvalue}{-\muxmulevelsetvalue}$,
$\partial_- \mathcal{B}^{[0,\muxmulevelsetvalue_0]}$]
\label{L:CAUSALSTRUCTUREOFTORIANDLEVELSETSOFMUINGEOMETRICCOORDINATES}
Assume the hypotheses and conclusions of Theorem~\ref{T:EXISTENCEUPTOTHESINGULARBOUNDARYATFIXEDKAPPA}
for $\muxmulevelsetvalue \in [0,\muxmulevelsetvalue_0]$. 
Recall that for $(\mulevelsetvalue,\muxmulevelsetvalue) \in [0,\mupositive] \times [0,\muxmulevelsetvalue_0]$,
the $\upmu$-adapted torus $\twoargmumuxtorus{\mulevelsetvalue}{-\muxmulevelsetvalue}$
defined in \eqref{E:MUXMUTORI} is contained in $\MSingular$ (see \eqref{E:SINGULARDEVELOPMENT}),
and that for $\mulevelsetvalue \in [0,\mupositive]$, the set $\mulevelsettwoarg{\mulevelsetvalue}{[0,\muxmulevelsetvalue_0]}$ 
defined in \eqref{E:LEVELSETSOFMUINSINGULARGETION} is the $\mulevelsetvalue$-level set of $\upmu$ in $\MSingular$.
Then the following results hold.

\noindent \underline{\textbf{The character of} 
$\twoargmumuxtorus{\mulevelsetvalue}{-\muxmulevelsetvalue}$
\textbf{and}
$\partial_- \mathcal{B}^{[0,\muxmulevelsetvalue_0]}$
}.

\begin{itemize}
	\item For $(\mulevelsetvalue,\muxmulevelsetvalue) \in [0,\mupositive] \times [0,\muxmulevelsetvalue_0]$,
			the torus $\twoargmumuxtorus{\mulevelsetvalue}{-\muxmulevelsetvalue}$ is a $2$-dimensional, $\gfour$-spacelike
			submanifold.
	\item  In particular, the crease 
			$\partial_- \mathcal{B}^{[0,\muxmulevelsetvalue_0]} = \twoargmumuxtorus{0}{0}$ 
			is a $2$-dimensional, $\gfour$-spacelike
			submanifold.
\end{itemize}

\medskip

\noindent \underline{\textbf{The character of} 
$\mulevelsettwoarg{\mulevelsetvalue}{[0,\muxmulevelsetvalue_0]}$
\textbf{and}
$\mathcal{B}^{[0,\muxmulevelsetvalue_0]}$
}.
\begin{itemize}
	\item For $0 < \mulevelsetvalue \leq \mupositive$,
			$\mulevelsettwoarg{\mulevelsetvalue}{[0,\muxmulevelsetvalue_0]}$
			is a $3$-dimensional, $\gfour$-spacelike
			submanifold-with-boundary.
\end{itemize}

\end{lemma}

\begin{remark}[Acoustical metric degeneracies along the singular boundary]
	\label{R:ACOUSTICALMETRICDEGENERACIESALONGSINGULARBOUNDARY}
	Note that Lemma~\ref{L:CAUSALSTRUCTUREOFTORIANDLEVELSETSOFMUINGEOMETRICCOORDINATES} does not address
	the causal structure of the singular boundary portion
	$\mathcal{B}^{[0,\muxmulevelsetvalue_0]}
	=
	\mulevelsettwoarg{0}{[0,\muxmulevelsetvalue_0]}$
	(see definition~\ref{E:SINGULARBOUNDARYPORTION}, and recall that $\upmu$ vanishes along $\mulevelsettwoarg{0}{[0,\muxmulevelsetvalue_0]}$),
	viewed as a subset of geometric coordinate space.
	We have avoided discussing the causal structure of $\mathcal{B}^{[0,\muxmulevelsetvalue_0]}$
	because relative to the geometric coordinates $(t,u,x^2,x^3)$,
	some components of the acoustical metric $\gfour$ degenerate along $\mulevelsettwoarg{0}{[0,\muxmulevelsetvalue_0]}$.
	In particular, using Lemma~\ref{L:BASICPROPERTIESOFVECTORFIELDS} and Lemma~\ref{L:COMMUTATORSTOCOORDINATES},
	one can check that along $\mathcal{B}^{[0,\muxmulevelsetvalue_0]}$, all vectors 
	$V \in \mbox{\upshape span} \left\lbrace \Lunit, \geop{u}  \right\rbrace$ are $\gfour$-null, i.e,
	they satisfy $\gfour(V,V) = 0$.
	In contrast, even along $\mathcal{B}^{[0,\muxmulevelsetvalue_0]}$,
	the Cartesian component matrix $\lbrace \gfour_{\alpha \beta} \rbrace_{\alpha,\beta=0,1,2,3}$ 
	of the acoustical metric is a non-degenerate $4 \times 4$ Lorentzian matrix
	(in fact, by \eqref{E:SPLITMETRICINTOMINKOWSKIANDREMAINDERPART}--\eqref{E:METRICPERTURBATIONVANISHESATTRIVIALPSISOLUTION}
	and the estimates of Prop.\,\ref{P:IMPROVEMENTOFAUXILIARYBOOTSTRAP}, the matrix is close to the standard Minkowski matrix).
	This discrepancy is caused by the fact that the change of variables map $\Upsilon(t,u,x^2,x^3) = (t,x^1,x^2,x^3)$
	has a non-injective Jacobian matrix along $\mathcal{B}^{[0,\muxmulevelsetvalue_0]}$;
	see Prop.\,\ref{E:JACOBIANDETERMINANTBOUNDCHOVGEOMETRICTOCARTESIANONWECAREABOUT}.
	We refer to Prop.\,\ref{P:DESCRIPTIONOFSINGULARBOUNDARYINCARTESIANSPACE} 
	for a description of the structure of $\Upsilon(\mathcal{B}^{[0,\muxmulevelsetvalue_0]})$,
	that is, the structure of the image of $\mathcal{B}^{[0,\muxmulevelsetvalue_0]}$
	into Cartesian coordinate space under the map $\Upsilon$,
	which is shown in the proposition to be injective on all of $\MInteresting$.
\end{remark}

\begin{proof}[Proof of Lemma~\ref{L:CAUSALSTRUCTUREOFTORIANDLEVELSETSOFMUINGEOMETRICCOORDINATES}]
We already showed in Prop.\,\ref{P:PROPERTIESOFMSINGULARANDCREASE}
that $\twoargmumuxtorus{\mulevelsetvalue}{-\muxmulevelsetvalue}$ is a $2$-dimensional manifold
and that $\mulevelsettwoarg{\mulevelsetvalue}{[0,\muxmulevelsetvalue_0]}$is a $3$-dimensional, $\gfour$-spacelike
submanifold-with-boundary.

\medskip
\noindent \textbf{Proof that $\twoargmumuxtorus{\mulevelsetvalue}{-\muxmulevelsetvalue}$ is $\gfour$-spacelike}:
Consider the map $\widetilde{\breve{\mathscr{M}}}(t,u,x^2,x^3) = \left(\upmu,-\muX \upmu,x^2,x^3 \right)$
from \eqref{E:MAPFROMGEOMETRICTOMUMUXMUCOORDINATESONMSINGULAR}.
In the proof of Prop.\,\ref{P:PROPERTIESOFMSINGULARANDCREASE}, we showed that
$\widetilde{\breve{\mathscr{M}}}$ is a diffeomorphism from $\MSingular$ onto 
$[0,\mupositive] \times [0,\muxmulevelsetvalue_0] \times \mathbb{T}^2$.
Using
Lemma~\ref{L:COMMUTATORSTOCOORDINATES},
Lemma~\ref{L:SCHEMATICSTRUCTUREOFVARIOUSTENSORSINTERMSOFCONTROLVARS}
and the $L^{\infty}$ estimates of Prop.\,\ref{P:IMPROVEMENTOFAUXILIARYBOOTSTRAP}
(with $\fundbootsmall$ replaced by $C \initialsmall$, which Theorem~\ref{T:EXISTENCEUPTOTHESINGULARBOUNDARYATFIXEDKAPPA} allows for)
we compute that: 
\begin{align} \label{E:JACOBIANMATRIXFORCHOVGEOTOMUMINUSMUXMUCOORDINATES}
	d_{\textnormal{geo}} \widetilde{\breve{\mathscr{M}}}
	& =
	\begin{pmatrix} 
			\geop{t} \upmu & \geop{u} \upmu & * & *  
				\\
			- \geop{t} \muX \upmu & - \geop{u} \muX \upmu & * & * 
				\\
			0 & 0 & 1 & 0 
				\\
			0 & 0 & 0 & 1
		\end{pmatrix},
\end{align}
where $d_{\textnormal{geo}} \widetilde{\breve{\mathscr{M}}}$ is the Jacobian matrix of $\widetilde{\breve{\mathscr{M}}}$ and ``$*$'' denotes 
quantities that are bounded in magnitude by $\mathcal{O}(\initialsmall)$.
From \eqref{E:JACOBIANMATRIXFORCHOVGEOTOMUMINUSMUXMUCOORDINATES}, 
\eqref{E:MUTRANSVERSALCONVEXITY},
and \eqref{E:BOUNDSONGEOMETRICTDERIVATIVEMUINTERESTINGREGION},
we deduce that
$\mbox{\upshape det} d_{\textnormal{geo}} \widetilde{\breve{\mathscr{M}}} \approx 1$.
It follows that:
\begin{align} \label{E:INVERSEJACOBIANMATRIXFORCHOVGEOTOMUMINUSMUXMUCOORDINATES}
	[d_{\textnormal{geo}} \widetilde{\breve{\mathscr{M}}}]^{-1}
	& =
	\frac{1}{\mbox{\upshape det} d_{\textnormal{geo}} \widetilde{\breve{\mathscr{M}}}}
	\begin{pmatrix} 
			- \geop{u} \muX \upmu  & - \geop{u} \upmu & * & *  
				\\
			 \geop{t} \muX \upmu & \geop{t} \upmu  & * & * 
				\\
			0 & 0 & 1 & 0 
				\\
			0 & 0 & 0 & 1
		\end{pmatrix}.
\end{align}
Fix any $(\mulevelsetvalue,\muxmulevelsetvalue) \in [0,\mupositive] \times [0,\muxmulevelsetvalue_0]$.
Since $\upmu \equiv \mulevelsetvalue$ and $\muX \upmu \equiv - \muxmulevelsetvalue$ along $\twoargmumuxtorus{\mulevelsetvalue}{-\muxmulevelsetvalue}$,
the last two columns of \eqref{E:INVERSEJACOBIANMATRIXFORCHOVGEOTOMUMINUSMUXMUCOORDINATES} 
are the components of vectorfields  
$V_{(2)} = * \geop{t} + * \geop{u} + \geop{x^2}$ and $V_{(3)} = * \geop{t} + * \geop{u} + \geop{x^3}$
that, when restricted to $\twoargmumuxtorus{\mulevelsetvalue}{-\muxmulevelsetvalue}$,
span the tangent space of $\twoargmumuxtorus{\mulevelsetvalue}{-\muxmulevelsetvalue}$.
Using
Lemma~\ref{L:BASICPROPERTIESOFVECTORFIELDS},
Lemma~\ref{L:COMMUTATORSTOCOORDINATES},
\eqref{E:SMOOTHTORIGABEXPRESSION},
Lemma~\ref{L:SCHEMATICSTRUCTUREOFVARIOUSTENSORSINTERMSOFCONTROLVARS}
and the $L^{\infty}$ estimates of Prop.\,\ref{P:IMPROVEMENTOFAUXILIARYBOOTSTRAP},
we compute that
$\gfour(V_{(A)},V_{(B)}) = \updelta_{AB} + \mathcal{O}(\mathring{\upalpha})$,
where $\updelta_{AB}$ is the Kronecker delta. From this estimate, it easily follows that
 $\twoargmumuxtorus{\mulevelsetvalue}{-\muxmulevelsetvalue}$ is $\gfour$-spacelike, 
even in the case of the crease (in which $\mulevelsetvalue = \muxmulevelsetvalue = 0$).

\medskip
\noindent \textbf{Proof that $\mulevelsettwoarg{\mulevelsetvalue}{[0,\muxmulevelsetvalue_0]}$ is $\gfour$-spacelike when $\mulevelsetvalue > 0$}:	
Fix any $\mulevelsetvalue \in (0, \mupositive]$.
Note that $\upmu \equiv \mulevelsetvalue$ along $\mulevelsettwoarg{\mulevelsetvalue}{[0,\muxmulevelsetvalue_0]}$,
that $\muX \upmu \leq 0$ along $\mulevelsettwoarg{\mulevelsetvalue}{[0,\muxmulevelsetvalue_0]}$ 
(by \eqref{E:LEVELSETSOFMUINSINGULARGETION} and the fact that $\muX \upmu \equiv - \muxmulevelsetvalue$ 
along $\twoargmumuxtorus{\mulevelsetvalue}{-\muxmulevelsetvalue}$),
and that $\Lunit \upmu < 0$ along $\mulevelsettwoarg{\mulevelsetvalue}{[0,\muxmulevelsetvalue_0]}$
(by Remark~\ref{R:SINGULARREGIONHASSMALLUVALUES} and \eqref{E:BOUNDSONLMUINTERESTINGREGION}).
Using these results, the estimate $|\angD \upmu|_{\gtorus} \lesssim \fundbootsmall$ obtained in the proof of
\eqref{E:POINTWISEBOUNDROUGHGRADIENTOFMU}, and the identity
\eqref{E:MUWEGHTEDSPACETIMEGRADMUNORMSQUARED}, we conclude that along $\mulevelsettwoarg{\mulevelsetvalue}{[0,\muxmulevelsetvalue_0]}$,
we have $\gfour(\Dfour^{\#} \upmu, \Dfour^{\#} \upmu) < 0$. This implies that the $\gfour$-normal to 
$\mulevelsettwoarg{\mulevelsetvalue}{[0,\muxmulevelsetvalue_0]}$ is $\gfour$-timelike, which is the desired result.
	
\end{proof}

\subsection{A new time function and related geometric objects}
\label{SS:NEWTIMEFUNCTION}
In this section, on $\MInteresting$, 
we define a time function $\newtimefunction$ 
as well as several related geometric objects.

\subsubsection{Definitions}
\label{SSS:DEFINITIONNEWTIMEFUNCTION}

\begin{definition}[$\newtimefunction$, 
${^{(Interesting)}\mathscr{T}}(t,u,x^2,x^3) = (\newtimefunction,u,x^2,x^3)$,
$\tisafunctionalonglevelsetsofnewtimefunctionarg{\timefunction}$,
$\levelsetgeneratornewtimefunction$, and $\partialderivativewithrespecttonewtimefunction$]
\label{D:LEVELSETGENERATORFORNEWTIMEFUNCTIONANDNEWTIMEFUNCTION}
Let $\MInteresting 
	= \MLeft
	\cup
	\MSingular
	\cup
	 \MRight$
	be the set defined in \eqref{E:INTERESTINGDEVELOPMENTOFDATA}
	and depicted in Fig.\,\ref{F:MINTERESTINGDEVELOPMENT}.
	
	\medskip
	
	\noindent \underline{\textbf{Definition of} $\newtimefunction$}.
	On $\MInteresting$, we define the scalar function $\newtimefunction$ as follows:
	\begin{align} \label{E:NEWTIMEFUNCTION}
		\newtimefunction(t,u,x^2,x^3)
		& \eqdef
		\begin{cases}
			\timefunctionarg{0}(t,u,x^2,x^3), & \mbox{in } 
			\MLeft,
			\\
	- \upmu(t,u,x^2,x^3), & \mbox{in } 
			\MSingular,
				\\
	\timefunctionarg{\muxmulevelsetvalue_0}(t,u,x^2,x^3), & \mbox{in } 
		 \MRight.
	\end{cases}
	\end{align}
	
	\medskip
	\noindent \underline{\textbf{Definition of level set portions of} $\newtimefunction$}.
	For $- \rightu \leq u_1 \leq u_2 \leq \leftu$ and $\timefunction \in [\timefunction_0,0]$,
	we define:
	\begin{align} \label{E:NEWTIMEFUNCTIONLEVELSETPORTION}
		\inthyp{\timefunction}{[- \rightu,\leftu]} 
		&
		\eqdef 
		\left\lbrace
			(t,u,x^2,x^3) \ | \ 
			(u,x^2,x^3) \in [u_1,u_2] \times \mathbb{T}^2,
				\,
			\newtimefunction(t,u,x^2,x^3) = \timefunction
			\right\rbrace.
	\end{align}
	
	\medskip
	
	\noindent \underline{\textbf{Definition of} ${^{(Interesting)}\mathscr{T}}$ and 
	$\tisafunctionalonglevelsetsofnewtimefunctionarg{\timefunction}$}.
	We define the map    
	$\InterestingCHOV : \MInteresting \rightarrow [\timefunction_0,0] \times [- \rightu,\leftu] \times \mathbb{T}^2$
	as follows:
	\begin{align} \label{E:CHOVFROMGEOTOINTERESTINGCOORDS}
		\InterestingCHOV(t,u,x^2,x^3) 
		& = (\newtimefunction,u,x^2,x^3).
	\end{align}
	
	Next, for each $\timefunction \in [\timefunction_0,0]$, 
	we define the function 
	$\tisafunctionalonglevelsetsofnewtimefunctionarg{\timefunction} : [- \rightu,\leftu] \times \mathbb{T}^2 \rightarrow \mathbb{R}$
	as follows:
	\begin{align} \label{E:DEFININGFUNCTIONCARTESIANTISAGRAPHALONGLEELSETSOFNEWTIMEFUNCTION}
		\tisafunctionalonglevelsetsofnewtimefunctionarg{\timefunction}(u,x^2,x^3)
		& \eqdef
		\begin{cases}
			\Cartesiantisafunctiononlevelsetsofroughtimefunctionarg{\timefunction}{0}(u,x^2,x^3)
			& \mbox{if } 
			\Eikonalisafunctiononmumuxtoriarg{-\timefunction}{0}(x^2,x^3) \leq u \leq \leftu
				\\
		\tisafunctiononlevelsetsofmu{-\timefunction}(u,x^2,x^3)
		& \mbox{if } 
		\Eikonalisafunctiononmumuxtoriarg{-\timefunction}{-\muxmulevelsetvalue_0}(x^2,x^3) \leq u \leq \Eikonalisafunctiononmumuxtoriarg{-\timefunction}{0}(x^2,x^3)
				\\
		\Cartesiantisafunctiononlevelsetsofroughtimefunctionarg{\timefunction}{\muxmulevelsetvalue_0}(u,x^2,x^3)
			& \mbox{if } 
			- \rightu \leq u \leq \Eikonalisafunctiononmumuxtoriarg{-\timefunction}{-\muxmulevelsetvalue_0}(x^2,x^3),
	\end{cases}
	\end{align}
	where $\Cartesiantisafunctiononlevelsetsofroughtimefunctionarg{\timefunction}{\muxmulevelsetvalue}$
	is the function from \eqref{E:LEVELSETSOFTIMEFUNCTIONAREAGRAPH},
	$\Eikonalisafunctiononmumuxtoriarg{\mulevelsetvalue}{-\muxmulevelsetvalue}$ is the function on $\mathbb{T}^2$ from
	\eqref{E:LASTSLICETORIAREGRAPHSABOVEFLATTORIINGEOMETRICCOORDINATES},
	and $\tisafunctiononlevelsetsofmu{\mulevelsetvalue}$ is the function from \eqref{E:GRAPHSTRUCTUREOFLEVELSETSOFMUINSINGULARGETION}.
	
	\medskip
	\noindent \underline{\textbf{Definition of} $\levelsetgeneratornewtimefunction$ and $\partialderivativewithrespecttonewtimefunction$}.
	In $\MInteresting$, we define
	the vectorfields $\levelsetgeneratornewtimefunction$ 
	and $\partialderivativewithrespecttonewtimefunction$
	as follows, where $\phi$ is the cut-off function from Def.\,\ref{D:WTRANSANDCUTOFF}:
	\begin{subequations}
	\begin{align} \label{E:LEVELSETGENERATORFORNEWTIMEFUNCTION}
		\levelsetgeneratornewtimefunction
		& \eqdef
		\begin{cases}
			\muX
			-
			\muX^A
			\left(
			\geop{x^A} 
			- 
			\frac{\geop{x^A} \timefunctionarg{0}}{\geop{t}\timefunctionarg{0}} \geop{t}
			\right),
			& \mbox{in } 
			\MLeft,
					\\
			\left(
			\muX
			-
			\frac{\muX \upmu}{\Lunit \upmu} \Lunit
			\right)
			-
			\left(
			\muX^A
			-
			\frac{\muX \upmu}{\Lunit \upmu} \Lunit^A
			\right)
			\left(
			\geop{x^A} 
			- 
			\frac{\geop{x^A} \upmu}{\geop{t} \upmu} \geop{t}
			\right),
			& \mbox{in } 
				\MSingular,
			\\
			\left(
			\muX 
			+ 
			\phi 
			\frac{\muxmulevelsetvalue_0}{\Lunit \upmu} \Lunit
			\right)
			-
			\left(
			\muX^A 
			+ 
			\phi 
			\frac{\muxmulevelsetvalue_0}{\Lunit \upmu} \Lunit^A
			\right)
			\left(
			\geop{x^A} 
			- 
			\frac{\geop{x^A} \timefunctionarg{\muxmulevelsetvalue_0}}{\geop{t}\timefunctionarg{\muxmulevelsetvalue_0}} \geop{t}
			\right),
	& \mbox{in } 
		 		\MRight,
	\end{cases}
		\\
	\partialderivativewithrespecttonewtimefunction
		& \eqdef
		\begin{cases}
			\frac{1}{\geop{t} \timefunctionarg{0}} \geop{t},
			& \mbox{in } 
			\MLeft,
			\\
			-
			\frac{1}{\geop{t} \upmu} \geop{t},
			& \mbox{in } 
			\MSingular,
			\\
			\frac{1}{\geop{t} \timefunctionarg{\muxmulevelsetvalue_0}} \geop{t},
	& \mbox{in } 
		 	\MRight.
		\end{cases}
			 \label{E:PARTIALDERIVATIVEWITHRESPECTTOTIMEFUNCTION}
	\end{align}
	\end{subequations}
		
\end{definition}

\begin{remark}[The regularity of $\newtimefunction$ and the connection to the causal structure of $\lbrace \upmu = 0 \rbrace$]
	\label{R:NEWTIMEFUNCTIONISC11ANDNOTBETTERANDCONNECTIONTOCAUSALSTRUCTUREOFMUZEROLEVELSET}
	In Prop.\,\ref{P:INTERESTINGREGIONFOLIATEDBYINTERESTINGTIMFEFUNCTION}, we show that
	$\newtimefunction \in C^{1,1}(\MInteresting)$.
	That regularity is optimal in the sense that generally, $\newtimefunction \notin C^2(\MInteresting)$.
	The reason is that $\timefunction$ and $- \upmu$
	generally agree \underline{only} to first-order 
	along $\datahypfortimefunctiontwoarg{0}{[\timefunction_0,0]}$,
	which is the common boundary of the regions $\MLeft$ and $\MSingular$
	in the piecewise-defined definition \eqref{E:NEWTIMEFUNCTION} of $\newtimefunction$;
	the first-order agreement along $\datahypfortimefunctiontwoarg{0}{[\timefunction_0,0]}$
	follows from \eqref{E:INITIALCONDITIONFORROUGHTIMEFUNCTION} and
	\eqref{E:GRADIENTOFTIMEFUNCTIONAGREESWITHGRADIENTOFMUALONGMUXMUEQUALSMINUSKAPPHYPERSURFACE},
	while \eqref{E:IVPFORROUGHTTIMEFUNCTION} implies that 
	$\Dfour \Wtransarg{\muxmulevelsetvalue} \timefunctionarg{\muxmulevelsetvalue} \equiv 0$
	along $\datahypfortimefunctiontwoarg{0}{[\timefunction_0,0]}$, which does not generally hold
	for $\upmu$. It follows that the second derivatives of $\newtimefunction$ can
	jump across $\datahypfortimefunctiontwoarg{0}{[\timefunction_0,0]}$.
	
	Moreover, it is generally \underline{impossible} to modify the construction of $\newtimefunction$ and 
	$\MInteresting$
	to enforce $C^2$ agreement of $\newtimefunction$ with $-\upmu$
	across $\datahypfortimefunctiontwoarg{0}{[\timefunction_0,0]}$
	in a manner such that the zero level set of the new $\newtimefunction$ still contains
	the singular boundary portion $\mathcal{B}^{[0,\muxmulevelsetvalue_0]}$.
	The reason is that if $\newtimefunction$ were $C^2$
	and agreed with $-\upmu$ along $\datahypfortimefunctiontwoarg{0}{[\timefunction_0,0]}$,
	then by Taylor expanding $\newtimefunction$ starting from
	the crease $\twoargmumuxtorus{0}{0}$ 
	(i.e., the subset of $\datahypfortimefunctiontwoarg{0}{[\timefunction_0,0]} \cap \lbrace \newtimefunction = 0 \rbrace$  
	along which $\upmu = \muX \upmu = 0$),
	and using the estimate the estimates
	\eqref{E:MUTRANSVERSALCONVEXITY}, 
	\eqref{E:BOUNDSONLMUINTERESTINGREGION},
	and \eqref{E:POINTWISEBOUNDROUGHGRADIENTOFMU} for $\upmu$
	(which by Taylor expansion imply corresponding estimates for $\newtimefunction$),
	one could prove the following result: 
	on near-zero level sets of the new $\newtimefunction$,
	in a region near the crease $\twoargmumuxtorus{0}{0}$ with $\upmu > 0$,
	we would have: 
	$\upmu \gfour(\Dfour^{\#} \newtimefunction, \Dfour^{\#} \newtimefunction) > 0$.
	This would imply that the level sets of $\newtimefunction$ are spacelike near the crease, i.e.,
	it could not be used as a time function.
	The main idea of the proof is that formally,
	along the level set $\lbrace \upmu = 0 \rbrace$,
	RHS~\eqref{E:MUWEGHTEDSPACETIMEGRADMUNORMSQUARED} becomes positive 
	as one passes through the crease in the direction of increasing $u$
	because $\muX \upmu$ becomes positive (thanks to \eqref{E:MUTRANSVERSALCONVEXITY}) 
	while $\Lunit \upmu$ is strictly negative 
	everywhere near the crease (thanks to \eqref{E:BOUNDSONLMUINTERESTINGREGION}).
	
\end{remark}

\subsubsection{Properties of $\newtimefunction$ and related quantities}
\label{SSS:PROPERTIESNEWTIMEFUNCTION}
In the next proposition, we derive some fundamental properties of the quantities from
Def.\,\ref{D:LEVELSETGENERATORFORNEWTIMEFUNCTIONANDNEWTIMEFUNCTION}
as well as implications of these properties for the structure of $\MInteresting$
and the behavior of $\upmu$ on $\MInteresting$.

\begin{proposition}[Properties of $\newtimefunction$,
$\InterestingCHOV$,
$\tisafunctionalonglevelsetsofnewtimefunctionarg{\timefunction}$,
$\partialderivativewithrespecttonewtimefunction$,
$\levelsetgeneratornewtimefunction$,
and $\MInteresting$]
	\label{P:INTERESTINGREGIONFOLIATEDBYINTERESTINGTIMFEFUNCTION}
	Let $\newtimefunction$, $\tisafunctionalonglevelsetsofnewtimefunctionarg{\timefunction}$,
	$\InterestingCHOV$,
	$\partialderivativewithrespecttonewtimefunction$,
	and
	$\levelsetgeneratornewtimefunction$
	be as in Def.\,\ref{D:LEVELSETGENERATORFORNEWTIMEFUNCTIONANDNEWTIMEFUNCTION}.
	Then these quantities enjoy the following properties on $\MInteresting$.
	
	\medskip
	
	\noindent \underline{\textbf{Properties of} $\newtimefunction$}.
	\begin{itemize}
\item The following estimate holds:
	\begin{align} \label{E:NEWTIMEFUNCTIONISC11}
		\| \newtimefunction \|_{C_{\textnormal{geo}}^{1,1}(\MInteresting)}
		& \leq C.
	\end{align}
	\item The following estimate holds:
	\begin{align} \label{E:TDERIVATIVEOFNEWTIMEFUNCTIONISAPPROXIMATELYUNITY}
		\geop{t} \newtimefunction
		& \approx 1,
		&&
		\mbox{on } \MInteresting.
	\end{align}
	\item For $\timefunction \in [\timefunction_0,0]$, the
	level set portions $\inthyp{\timefunction}{[- \rightu,\leftu]}$
	defined in \eqref{E:NEWTIMEFUNCTIONLEVELSETPORTION}
	are $\gfour$-spacelike, except along the singular boundary portion
	$\mathcal{B}^{[0,\muxmulevelsetvalue_0]} \subset \inthyp{0}{[- \rightu,\leftu]} \cap \MSingular$.
	\end{itemize}
	
	\medskip
	\noindent \underline{\textbf{The behavior of $\upmu$ on $\MInteresting$}}.
	For $\timefunction \in [\timefunction_0,0]$,
		we have: 
		\begin{align} \label{E:MINVALUEOFMUONNEWTIMEFUNCTIONLEVELSETS}
			\min_{\inthyp{\timefunction}{[- \rightu,\leftu]}} \upmu
			& = - \timefunction.
		\end{align}
		Moreover, within $\inthyp{\timefunction}{[- \rightu,\leftu]}$,
		the minimum value of $- \timefunction$
		in \eqref{E:MINVALUEOFMUONNEWTIMEFUNCTIONLEVELSETS} is achieved by $\upmu$ precisely on the
		set $\mulevelsettwoarg{- \timefunction}{[0,\muxmulevelsetvalue_0]}
		\eqdef 
		\bigcup_{\muxmulevelsetvalue \in [0,\muxmulevelsetvalue_0]} \twoargmumuxtorus{- \timefunction}{-\muxmulevelsetvalue}$
		from definition~\eqref{E:LEVELSETSOFMUINSINGULARGETION}.
	
	\medskip
	
	\noindent \underline{\textbf{Properties of} $\tisafunctionalonglevelsetsofnewtimefunctionarg{\timefunction}$}.
	\begin{itemize}
	\item
	For each $\timefunction \in [\timefunction_0,0]$, we have:
	\begin{align}
		\| \tisafunctionalonglevelsetsofnewtimefunctionarg{\timefunction} \|_{C^{1,1}([- \rightu,\leftu] \times \mathbb{T}^2)}
		& \leq C.
		\label{E:GRAPHDEFININGFUNCTIONFORLEVELSETSOFNEWTIMEFUNCTIONISC11}
	\end{align}
	\item The level set portions $\inthyp{\timefunction}{[- \rightu,\leftu]}$
	defined in \eqref{E:NEWTIMEFUNCTIONLEVELSETPORTION}
	have the following graph structure:
	\begin{align} \label{E:CARTESIANTISAGRAPHALONGLEELSETSOFNEWTIMEFUNCTION}
		\inthyp{\timefunction}{[- \rightu,\leftu]}
		& = 
		\left\lbrace
			(t,u,x^2,x^3) 
				\ | \
				(u,x^2,x^3) \in [- \rightu,\leftu] \times \mathbb{T}^2,
					\,
				t = \tisafunctionalonglevelsetsofnewtimefunctionarg{\timefunction}(u,x^2,x^3)
		\right\rbrace.
	\end{align}
	\end{itemize}
	
	\medskip
	
	\noindent \underline{\textbf{Properties of} $\InterestingCHOV$}.
	\begin{itemize}
	\item
	The change of variables map $\InterestingCHOV(t,u,x^2,x^3) = (\newtimefunction,u,x^2,x^3)$
	is a diffeomorphism from $\MInteresting$ onto $[\timefunction_0,0] \times [- \rightu,\leftu] \times \mathbb{T}^2$.
	\item The following estimate holds:
	\begin{align} \label{E:INTERESTINGCHOVMAPISC11}
		\| \InterestingCHOV \|_{C_{\textnormal{geo}}^{1,1}(\MInteresting)}
		& \leq C.
	\end{align}
	\end{itemize}
	
	\medskip
	
	\noindent \underline{\textbf{Properties of} $\partialderivativewithrespecttonewtimefunction$}
	\begin{itemize}
	\item In $\MInteresting$, we have:
	\begin{align} \label{E:PARTIALTDERIVATIVEFORNEWTIMEFUNCTIONISINFACTPARTIALTDERIVATIVE}
		\partialderivativewithrespecttonewtimefunction \newtimefunction 
		& = 1,
		\qquad
		\partialderivativewithrespecttonewtimefunction u
		=
		\partialderivativewithrespecttonewtimefunction x^2 
		=
		\partialderivativewithrespecttonewtimefunction x^3
		=
		0.
	\end{align}
	In particular, 
	$\partialderivativewithrespecttonewtimefunction$
	is the partial derivative with respect to $\newtimefunction$ in the coordinate system $(\newtimefunction,u,x^2,x^3)$.
	\item The following estimate holds:
		\begin{align} \label{E:PARTIALDERIVATIVEWTIHRESPECTTONEWTIMEFUNCTIONISLIPSCHITZ}
		\| \partialderivativewithrespecttonewtimefunction \|_{C_{\textnormal{geo}}^{0,1}(\MInteresting)}
		& \leq C.
	\end{align}
	\end{itemize}
	
	\medskip
	
	\noindent \underline{\textbf{Properties of} $\levelsetgeneratornewtimefunction$}.
	\begin{itemize}
		\item In $\MInteresting$, we have:
	\begin{align} \label{E:LEVELSETGENOFNEWTIMEFUNCTIONISPARTIALUDERIVATIVE}
		\levelsetgeneratornewtimefunction u 
		& = 1,
		\qquad
		\levelsetgeneratornewtimefunction \newtimefunction 
		=
		\levelsetgeneratornewtimefunction x^2 
		=
		\levelsetgeneratornewtimefunction x^3
		=
		0.
	\end{align}
	In particular, 
	$\levelsetgeneratornewtimefunction$
	is the partial derivative with respect to $u$ in the coordinate system $(\newtimefunction,u,x^2,x^3)$.
	\item The following estimate holds:
		\begin{align} \label{E:LEVELSETGENERATOROFNEWTIMEFUNCTIONISLIPSCHITZ}
		\| \levelsetgeneratornewtimefunction \|_{C_{\textnormal{geo}}^{0,1}(\MInteresting)}
		& \leq C.
	\end{align}
	\item Every $u$-parametrized integral curve of $\levelsetgeneratornewtimefunction$ 
	is defined on the interval $[- \rightu,\leftu]$,
	and the images of these integral curves are contained in $\MInteresting$.
	\item For each fixed $\timefunction \in [\timefunction_0,0]$,
	every $u$-parametrized integral curve of $\levelsetgeneratornewtimefunction$ 
	in the level set $\inthyp{\timefunction}{[- \rightu,\leftu]}$
	intersects the torus $\twoargmumuxtorus{-\timefunction}{0} \subset \datahypfortimefunctiontwoarg{0}{[\timefunction_0,0]}$
	in precisely one point. 
	Moreover, $\twoargmumuxtorus{-\timefunction}{0} \subset \inthyp{\timefunction}{[-\frac{\interestingu}{2},\frac{\interestingu}{2}]}$.
	
	\end{itemize}
	
	\medskip
	
	\noindent \underline{\textbf{Properties of} $\MInteresting$}
	The following results hold.
	\begin{itemize}
		\item \begin{align} \label{E:INTERESTINGREGIONFOLIATEDBYINTERESTINGTIMFEFUNCTION}
				\MInteresting
				& =	
				\bigcup_{\timefunction \in [\timefunction_0,0]} \inthyp{\timefunction}{[- \rightu,\leftu]}.
			\end{align}
		\item 
			For every pair of points $q_1, q_2 \in [\timefunction_0,0] \times [- \rightu,\leftu] \times \mathbb{T}^2$,
			we have:
			\begin{align} \label{E:INTERESTINGCOMPARABLEDISTANCES}
						\mbox{\upshape dist}_{\mbox{\upshape flat}}
						\left(\InverseInterestingCHOV(q_1),\InverseInterestingCHOV(q_2) \right)
						& \approx
						\mbox{\upshape dist}_{\mbox{\upshape flat}}(q_1,q_2),
					\end{align}
					where on both sides of \eqref{E:INTERESTINGCOMPARABLEDISTANCES}, 
					$\mbox{\upshape dist}_{\mbox{\upshape flat}}(A,B)$ is the
					standard Euclidean distance between $A$ and $B$ in the flat space 
					$\mathbb{R} \times \mathbb{R} \times \mathbb{T}^2$.
		\item (\textbf{Quasi-convexity}) 
					$\MInteresting$ is quasi-convex.
					That is, there is a constant $C > 0$ such that
					every pair of points
					$p_1,p_2 \in \MInteresting$
					are connected by a $C_{\textnormal{geo}}^1$ curve in $\MInteresting$
					whose length with respect to the standard flat Euclidean metric on geometric coordinate space
					$\mathbb{R} \times \mathbb{R} \times \mathbb{T}^2$
					is $\leq C \mbox{\upshape dist}_{\mbox{\upshape flat}}(p_1,p_2)$.
		\item (\textbf{Sobolev embedding})
			There is a constant $C > 0$ 
			such that the following Sobolev embedding result holds for scalar functions $f$ on 
			$\mbox{\upshape int}(\MInteresting)$:
			\begin{align} \label{E:SOBOELVEMBEDDINGINTERESTINGREGIONRELYINGONQUASICONVEXITY}
						\| f \|_{C_{\textnormal{geo}}^{0,1}(\MInteresting)}
						& \leq
						C
						\| f \|_{W_{\textnormal{geo}}^{1,\infty}(\mbox{\upshape int}(\MInteresting))}.
			\end{align}
	\end{itemize}
	
\end{proposition} 

\begin{proof}
	\noindent \textbf{Proof of \eqref{E:INTERESTINGREGIONFOLIATEDBYINTERESTINGTIMFEFUNCTION}}:
	From definition~\eqref{E:LEFTDEVELOPMENT},
	we see that 
	$
	\MLeft 
	\subset
	\twoargMrough{[\timefunction_0,0],[- \rightu,\leftu]}{0}
	$.
	Hence, since
	$
	\twoargMrough{[\timefunction_0,0],[- \rightu,\leftu]}{0}
	$
	is foliated by the level sets of $\timefunctionarg{0}$,
	it follows from definition \eqref{E:NEWTIMEFUNCTION}
	that 
	$
	\MLeft 
	$
	is foliated by the level sets of $\newtimefunction$,
	which have the range $[\timefunction_0,0]$.
	Using similar reasoning based on definition \eqref{E:RIGHTDEVELOPMENT},
	we see that
	$
	\MRight 
	$
	is foliated by the level sets of $\newtimefunction$.
	Moreover, from Proposition~\ref{P:PROPERTIESOFMSINGULARANDCREASE},
	we see that $\MSingular$ is foliated by the level sets of
	$\upmu$, which have the range $[0,\mupositive] = [0,-\timefunction_0]$ in
	$\MSingular$.
	Hence, from definition~\eqref{E:NEWTIMEFUNCTION}, we see that
	$\MSingular$
	is foliated by the level sets of $\newtimefunction$,
	which have the range $[\timefunction_0,0]$.
	From these facts and definition~\eqref{E:INTERESTINGDEVELOPMENTOFDATA},
	we conclude \eqref{E:INTERESTINGREGIONFOLIATEDBYINTERESTINGTIMFEFUNCTION}.
	
	\medskip
	\noindent \textbf{Proof of \eqref{E:CARTESIANTISAGRAPHALONGLEELSETSOFNEWTIMEFUNCTION}}:
	In Prop.\,\ref{P:PROPERTIESOFMSINGULARANDCREASE}, we showed that the
	top boundary of $\MLeft$, which is contained in the level set
	$\lbrace \timefunctionarg{0} = 0 \rbrace$,
	is the hypersurface
	$
	\lbrace
				(t,u,x^2,x^3)
					\ | \
					(x^2,x^3) \in \mathbb{T}^2,
						\,
					\Eikonalisafunctiononmumuxtoriarg{0}{0}(x^2,x^3) \leq u \leq \leftu,
						\,
					\mbox{and }
					t = \Cartesiantisafunctiononlevelsetsofroughtimefunctionarg{0}{0}(u,x^2,x^3)
	\rbrace
	$.
	The same arguments used in the proof and definition~\eqref{E:NEWTIMEFUNCTION}
	also imply that for each fixed $\timefunction \in [\timefunction_0,0]$,
	we have
	$
	\lbrace \newtimefunction = \timefunction \rbrace \cap \MLeft
	=
	\lbrace
				(t,u,x^2,x^3)
					\ | \
					(x^2,x^3) \in \mathbb{T}^2,
						\,
					\Eikonalisafunctiononmumuxtoriarg{-\timefunction}{0}(x^2,x^3) \leq u \leq \leftu,
						\,
					\mbox{and }
					t = \Cartesiantisafunctiononlevelsetsofroughtimefunctionarg{-\timefunction}{0}(u,x^2,x^3)
	\rbrace
	$.
	Similarly, for each fixed $\timefunction \in [\timefunction_0,0]$,
	we have
	$
	\lbrace \newtimefunction = \timefunction \rbrace \cap \MRight
	=
	\lbrace
				(t,u,x^2,x^3)
					\ | \
					(x^2,x^3) \in \mathbb{T}^2,
						\,
					- \rightu \leq u \leq \Eikonalisafunctiononmumuxtoriarg{-\timefunction}{\muxmulevelsetvalue_0}(x^2,x^3),
						\,
					\mbox{and }
					t = \Cartesiantisafunctiononlevelsetsofroughtimefunctionarg{-\timefunction}{-\muxmulevelsetvalue_0}(u,x^2,x^3)
	\rbrace
	$.
	Moreover, from Prop.\,\ref{P:PROPERTIESOFMSINGULARANDCREASE}
	(in particular \eqref{E:DOMAINOFGRAPHFORMULEVELSETINSINGULARREGION}--\eqref{E:GRAPHSTRUCTUREOFLEVELSETSOFMUINSINGULARGETION})
	and definition~\eqref{E:NEWTIMEFUNCTION}, we see that for each fixed $\timefunction \in [\timefunction_0,0]$, we have
	$
	\lbrace \newtimefunction = \timefunction \rbrace \cap \MSingular
	=
	\left\lbrace
				\left(\tisafunctiononlevelsetsofmu{-\timefunction}(u,x^2,x^3),u,x^2,x^3 \right) 
					\ | \
					(x^2,x^3) \in \mathbb{T}^2,
						\,
					\Eikonalisafunctiononmumuxtoriarg{-\timefunction}{-\muxmulevelsetvalue_0}(x^2,x^3) 
					\leq 
					u 
					\leq \Eikonalisafunctiononmumuxtoriarg{-\timefunction}{0}(x^2,x^3),
	\right\rbrace
	$.
	From these facts and definition~\eqref{E:DEFININGFUNCTIONCARTESIANTISAGRAPHALONGLEELSETSOFNEWTIMEFUNCTION},
	we conclude \eqref{E:CARTESIANTISAGRAPHALONGLEELSETSOFNEWTIMEFUNCTION}.
	
	In $\MLeft$, the level sets of $\newtimefunction$ and $\timefunctionarg{0}$ agree.
	Hence, \eqref{E:LEVELSETSOFTIMEFUNCTIONAREAGRAPH} implies that in
	$\lbrace \newtimefunction = \timefunction \rbrace \cap \MLeft$,
	$t = \Cartesiantisafunctiononlevelsetsofroughtimefunctionarg{\timefunction}{0}(u,x^2,x^3)$.

	\medskip
	\noindent \textbf{Proof of \eqref{E:TDERIVATIVEOFNEWTIMEFUNCTIONISAPPROXIMATELYUNITY}}:
	From definitions \eqref{E:LEFTDEVELOPMENT}, \eqref{E:RIGHTDEVELOPMENT}, and \eqref{E:NEWTIMEFUNCTION},
	and Lemma~\ref{L:CONTINUOUSEXTNESION}
	in the cases $\muxmulevelsetvalue = 0$ and $\muxmulevelsetvalue = \muxmulevelsetvalue_0$,
	we see that 
	$\newtimefunction \in C_{\textnormal{geo}}^{2,1}(\mbox{\upshape cl}(\MLeft))$
	and
	$\newtimefunction \in C_{\textnormal{geo}}^{2,1}(\mbox{\upshape cl}(\MRight))$,
	where $\mbox{\upshape cl}$ denotes set closure in geometric coordinate space.
	Similarly, from definitions \eqref{E:SINGULARDEVELOPMENT} and \eqref{E:NEWTIMEFUNCTION} and the estimate
	\eqref{E:C21BOUNDMUSINGULARREGION}, we see that
	$\newtimefunction \in C_{\textnormal{geo}}^{2,1}(\MSingular)$.
	Moreover,
	Prop.\,\ref{P:PROPERTIESOFMSINGULARANDCREASE}
	yields that
	$
	\MLeft
	$
	and
	$\MSingular$
	have the common $C^{1,1}$ boundary
	$\datahypfortimefunctiontwoarg{0}{[\timefunction_0,0]}$,
	and from the above observations, definition \eqref{E:NEWTIMEFUNCTION}, and Def.\,\ref{D:ROUGHTIMEFUNCTION},
	we see that 
	$\newtimefunction$ and its first partial derivatives with respect to the geometric coordinates $(t,u,x^2,x^3)$
	are continuous across this common boundary. 
	Similar arguments yield that
	$
	\MRight
	$
	and
	$\MSingular$
	have the common $C^{1,1}$ boundary
	$\datahypfortimefunctiontwoarg{-\muxmulevelsetvalue_0}{[\timefunction_0,0]}$,
	and that $\newtimefunction$ and its first partial derivatives with respect to the geometric coordinates $(t,u,x^2,x^3)$
	are continuous across this common boundary.
	From these facts, definition~\eqref{E:CHOVFROMGEOTOINTERESTINGCOORDS}, 
	and Rademacher's theorem,
	it follows that 
	$\newtimefunction, \, \InterestingCHOV \in W^{2,\infty}(\mbox{\upshape int}(\MInteresting)) 
	\cap C^1(\MInteresting)$,
	and that the following estimates hold:
	\begin{align} 	\label{E:FIRSTW2INFTYESTIMATESFORNEWTIMEFUNCTIONANDCHOV}
		\| \newtimefunction \|_{W_{\textnormal{geo}}^{2,\infty}(\mbox{\upshape int}(\MInteresting))},
			\,
		\| \InterestingCHOV \|_{W_{\textnormal{geo}}^{2,\infty}(\mbox{\upshape int}(\MInteresting))}
		& \leq C,
			\\
	\| \newtimefunction \|_{C^1(\MInteresting)},	
			\,
		\|  \InterestingCHOV \|_{C^1(\MInteresting)}
		& \leq C.
		\label{E:FIRSTC1ESTIMATESFORNEWTIMEFUNCTIONANDCHOV}
	\end{align}
	Next, we use definition \eqref{E:NEWTIMEFUNCTION} 
	and the estimate \eqref{E:CLOSEDVERSIONKEYJACOBIANDETERMINANTESTIMATECHOVGEOTOROUGH}
	to deduce that
	$\geop{t} \newtimefunction|_{\MLeft} \approx 1$ 
	and
	$\geop{t} \newtimefunction|_{\MRight} \approx 1$. 
	Similarly, since definitions \eqref{E:SINGULARDEVELOPMENT}, \eqref{E:NEWTIMEFUNCTION}, and Def.\,\ref{D:ROUGHTIMEFUNCTION} 
	imply that
	$\newtimefunction|_{\datahypfortimefunctiontwoarg{-\muxmulevelsetvalue}{[\timefunction_0,0]}}
	=
	\timefunctionarg{\muxmulevelsetvalue}|_{\datahypfortimefunctiontwoarg{-\muxmulevelsetvalue}{[\timefunction_0,0]}}
	$,
	we deduce that
	$\geop{t} \newtimefunction|_{\MSingular} \approx 1$.
	From these bounds and \eqref{E:FIRSTC1ESTIMATESFORNEWTIMEFUNCTIONANDCHOV}, 
	we conclude \eqref{E:TDERIVATIVEOFNEWTIMEFUNCTIONISAPPROXIMATELYUNITY}.
	
	\medskip
	\noindent \textbf{Proof that $\InterestingCHOV$ is a diffeomorphism and proof of 
	\eqref{E:GRAPHDEFININGFUNCTIONFORLEVELSETSOFNEWTIMEFUNCTIONISC11}}:
	From the definition~\eqref{E:CHOVFROMGEOTOINTERESTINGCOORDS} of $\InterestingCHOV$,
	the estimates 
	\eqref{E:TDERIVATIVEOFNEWTIMEFUNCTIONISAPPROXIMATELYUNITY}
	and
	\eqref{E:FIRSTC1ESTIMATESFORNEWTIMEFUNCTIONANDCHOV},
	and the inverse function theorem,
	we see that
	$\InterestingCHOV$ is a local diffeomorphism on
	$\MInteresting$.
	Also using the graph structure of $\inthyp{\timefunction}{[- \rightu,\leftu]}$ from 
	\eqref{E:CARTESIANTISAGRAPHALONGLEELSETSOFNEWTIMEFUNCTION},
	we see that 
	$\InterestingCHOV$ is injective on $\MInteresting$
	and that $\InterestingCHOV(\MInteresting) = [\timefunction_0,0] \times [- \rightu,\leftu] \times \mathbb{T}^2$.
	That is, 
	$\InterestingCHOV$ is a global diffeomorphism from 
	$\MInteresting$ onto $[\timefunction_0,0] \times [- \rightu,\leftu] \times \mathbb{T}^2$.
	
	Next, using \eqref{E:TDERIVATIVEOFNEWTIMEFUNCTIONISAPPROXIMATELYUNITY},
	\eqref{E:FIRSTW2INFTYESTIMATESFORNEWTIMEFUNCTIONANDCHOV},
	and \eqref{E:FIRSTC1ESTIMATESFORNEWTIMEFUNCTIONANDCHOV},
	we deduce that the inverse map 
	satisfies 
	$
	\| \InverseInterestingCHOV \|_{C^1([\timefunction_0,0] \times [- \rightu,\leftu] \times \mathbb{T}^2)}
	\leq 
	C
	$
	and
	$
	\| \InverseInterestingCHOV \|_{W^{2,\infty}((\timefunction_0,0) \times (- \rightu,\leftu) \times \mathbb{T}^2)}
	\leq 
	C
	$.
	Thanks to the convexity of $(\timefunction_0,0) \times (- \rightu,\leftu) \times \mathbb{T}^2)$,
	standard Sobolev embedding also yields
	$
	\| ({^{Interesting)}\mathscr{T}})^{-1} \|_{C^{1,1}([\timefunction_0,0] \times [- \rightu,\leftu] \times \mathbb{T}^2)}
	\leq 
	C
	$.
	Since $\tisafunctionalonglevelsetsofnewtimefunctionarg{\timefunction}$ is the first component function of
	$({^{Interesting)}\mathscr{T}})^{-1}$,
	we conclude \eqref{E:GRAPHDEFININGFUNCTIONFORLEVELSETSOFNEWTIMEFUNCTIONISC11}.
	
	\medskip
	\noindent \textbf{Proof of \eqref{E:MINVALUEOFMUONNEWTIMEFUNCTIONLEVELSETS} and related properties of $\upmu$}:
	From Defs.\,\ref{D:SINGULARBOUNDARYPORTION} and
	\ref{D:LEVELSETGENERATORFORNEWTIMEFUNCTIONANDNEWTIMEFUNCTION},
	it follows that
	$\inthyp{\timefunction}{[- \rightu,\leftu]} \cap \MLeft 
	= 
	\hypthreearg{\timefunction}{[- \rightu,\leftu]}{0}
	\cap \MLeft
	$,
	$\inthyp{\timefunction}{[- \rightu,\leftu]} \cap \MRight 
	= 
	\hypthreearg{\timefunction}{[- \rightu,\leftu]}{\muxmulevelsetvalue_0}
	\cap \MRight
	$,
	and 
	$\inthyp{\timefunction}{[- \rightu,\leftu]} \cap \MSingular
	=
	\mulevelsettwoarg{- \timefunction}{[0,\muxmulevelsetvalue_0]}
	$. 
	Since 
	$\inthyp{\timefunction}{[- \rightu,\leftu]}  
	= 
	(\inthyp{\timefunction}{[- \rightu,\leftu]} \cap \MLeft)
	\cup
	(\inthyp{\timefunction}{[- \rightu,\leftu]} \cap \MRight)
	\cup
	(\inthyp{\timefunction}{[- \rightu,\leftu]} \cap \MSingular)
	$,
	\eqref{E:MINVALUEOFMUONNEWTIMEFUNCTIONLEVELSETS}
	and the results stated just below \eqref{E:MINVALUEOFMUONNEWTIMEFUNCTIONLEVELSETS}
	follow from \eqref{E:MINVALUEOFMUONFOLIATION},
	the results stated just below \eqref{E:MINVALUEOFMUONFOLIATION},
	and \eqref{E:LEVELSETSOFMUINSINGULARGETION}
	(which in particular shows that $\upmu \equiv - \timefunction$ along $\mulevelsettwoarg{- \timefunction}{[0,\muxmulevelsetvalue_0]}$).
	
	\medskip
	\noindent \textbf{Proof of \eqref{E:INTERESTINGCOMPARABLEDISTANCES},
	of the quantitative quasi-convexity of $\MInteresting$, and of
	\eqref{E:SOBOELVEMBEDDINGINTERESTINGREGIONRELYINGONQUASICONVEXITY}}:
	Thanks to the estimate \eqref{E:TDERIVATIVEOFNEWTIMEFUNCTIONISAPPROXIMATELYUNITY},
	the estimate
	$
	\| \InverseInterestingCHOV \|_{C^1([\timefunction_0,0] \times [- \rightu,\leftu] \times \mathbb{T}^2)}
	\leq 
	C
	$ proved above,
	and the convexity of $[\timefunction_0,0] \times [- \rightu,\leftu] \times \mathbb{T}^2$,
	the same arguments given in
	the proof of Lemma~\ref{L:DIFFEOMORPHICEXTENSIONOFROUGHCOORDINATES}
	yield 
	\eqref{E:INTERESTINGCOMPARABLEDISTANCES},
	the quasi-convexity of $\MInteresting$,
	and the Sobolev embedding result
	\eqref{E:SOBOELVEMBEDDINGINTERESTINGREGIONRELYINGONQUASICONVEXITY}.
	
	\medskip
	\noindent \textbf{Proof of \eqref{E:NEWTIMEFUNCTIONISC11}
	and
	\eqref{E:INTERESTINGCHOVMAPISC11}}:
	These estimates follow from
	\eqref{E:SOBOELVEMBEDDINGINTERESTINGREGIONRELYINGONQUASICONVEXITY}
	and \eqref{E:FIRSTW2INFTYESTIMATESFORNEWTIMEFUNCTIONANDCHOV}.

	\medskip
	\noindent \textbf{Proof that the $\inthyp{\timefunction}{[- \rightu,\leftu]}$ 
	are $\gfour$-spacelike for $\timefunction \in [\timefunction_0,0]$, except along $\mathcal{B}^{[0,\muxmulevelsetvalue_0]}$}:
	This follows from definition~\eqref{E:NEWTIMEFUNCTION},
	Lemma~\ref{L:CAUSALSTRUCTUREOFTORIANDLEVELSETSOFMUINGEOMETRICCOORDINATES}
	(note that $\MSingular \backslash \mathcal{B}^{[0,\muxmulevelsetvalue_0]} 
	= \bigcup_{\mulevelsetvalue \in (0,\mupositive]} \mulevelsettwoarg{\mulevelsetvalue}{[0,\muxmulevelsetvalue_0]}$),
	and \eqref{E:HYPNORMALSIZE}, which shows that the vectorfield
	$\hypnormalarg{\muxmulevelsetvalue}$ (which is $\gfour$-orthogonal to $\hypthreearg{\timefunction}{[- \rightu,\leftu]}{\muxmulevelsetvalue}$),
	is $\gfour$-timelike in regions where $\upmu > 0$ 
	(as is the case in $\MLeft \cup \MRight$).

\medskip
	\noindent \textbf{Proof of \eqref{E:PARTIALTDERIVATIVEFORNEWTIMEFUNCTIONISINFACTPARTIALTDERIVATIVE}
	and \eqref{E:LEVELSETGENOFNEWTIMEFUNCTIONISPARTIALUDERIVATIVE}}:
	These identities
	are straightforward to verify from Def.\,\ref{D:LEVELSETGENERATORFORNEWTIMEFUNCTIONANDNEWTIMEFUNCTION}.

\medskip
	\noindent \textbf{Proof of \eqref{E:PARTIALDERIVATIVEWTIHRESPECTTONEWTIMEFUNCTIONISLIPSCHITZ}
	and
	\eqref{E:LEVELSETGENERATOROFNEWTIMEFUNCTIONISLIPSCHITZ}}:	
	Since we have already shown that
	$\levelsetgeneratornewtimefunction$
	and
	$\partialderivativewithrespecttonewtimefunction$
	are coordinate partial derivatives in the coordinate system
	$(\newtimefunction,u,x^2,x^3)$,
	the
	estimates 
	\eqref{E:PARTIALDERIVATIVEWTIHRESPECTTONEWTIMEFUNCTIONISLIPSCHITZ}
	and
	\eqref{E:LEVELSETGENERATOROFNEWTIMEFUNCTIONISLIPSCHITZ}
	follow from
	\eqref{E:INTERESTINGCHOVMAPISC11}
	and the estimate
	$$
	\| \InverseInterestingCHOV \|_{C^{1,1}([\timefunction_0,0] \times [- \rightu,\leftu] \times \mathbb{T}^2)}
	\leq 
	C
	$$noted above.
	
	\medskip
	\noindent \textbf{Proof of the remaining properties of $\levelsetgeneratornewtimefunction$}:
	Since \eqref{E:LEVELSETGENOFNEWTIMEFUNCTIONISPARTIALUDERIVATIVE} shows that
	$\levelsetgeneratornewtimefunction$
	is the partial derivative with respect to $u$ in the coordinate system $(\newtimefunction,u,x^2,x^3)$,
	it trivially follows that every
	$u$-parametrized integral curve of $\levelsetgeneratornewtimefunction$ 
	that starts in the level set $\inthyp{\timefunction}{[- \rightu,\leftu]}$ (for some $\timefunction \in [\timefunction_0,0]$)
	is defined on the interval $[- \rightu,\leftu]$
	and remains in $\inthyp{\timefunction}{[- \rightu,\leftu]}$.
	
	Finally, we show that every integral curve in $\inthyp{\timefunction}{[- \rightu,\leftu]}$
	from the previous paragraph 
	must intersect
	$\twoargmumuxtorus{-\timefunction}{0}$ in a unique point.
	To this end, we note that
	Prop.\,\ref{P:PROPERTIESOFMSINGULARANDCREASE} implies that
	$\extendedembeddatahypersurface^{-1}$
	is a $C_{\textnormal{geo}}^{1,1}$ diffeomorphism from 
	$\MSingular$
	onto $[0,\mupositive] \times [0,\muxmulevelsetvalue_0] \times \mathbb{T}^2$.
	In particular, considering the form \eqref{E:EMBEDDINGOFSINGULARREGION} of $\extendedembeddatahypersurface$,
	and
	using \eqref{E:IMPROVEDLEVELSETSTRUCTUREANDLOCATIONOFMIN} and the fact
	that 
	$\newtimefunction|_{\datahypfortimefunctiontwoarg{0}{[\timefunction_0,0]}}
	=
	\timefunctionarg{0}|_{\datahypfortimefunctiontwoarg{0}{[\timefunction_0,0]}}
	$,
	we see that
	$\twoargmumuxtorus{-\timefunction}{0} 
	\subset 
	\datahypfortimefunctiontwoarg{0}{[\timefunction_0,0]}
	\cap
	\inthyp{\timefunction}{[\frac{-\interestingu}{2},\frac{\interestingu}{2}]}
	\subset
	\MSingular
	$,
	and that along $\twoargmumuxtorus{-\timefunction}{0}$,
	$u$ is a $C^{1,1}$ function of $(x^2,x^3) \in \mathbb{T}^2$
	satisfying $|u| \leq \frac{\interestingu}{2}$.
	Combining these results and using that
	$\levelsetgeneratornewtimefunction u = 1$ and $\levelsetgeneratornewtimefunction x^2 = \levelsetgeneratornewtimefunction x^3 = 0$,
	we see that every integral curve from the previous paragraph must intersect
	$\twoargmumuxtorus{-\timefunction}{0}$. The uniqueness of the point follows from
	the fact that 
	$\muX \upmu|_{\twoargmumuxtorus{-\timefunction}{0}} = 0$
	and the fact that when $|u| \leq \interestingu$, we have the estimate
	$\levelsetgeneratornewtimefunction \muX \upmu \approx 1$;
	this estimate follows from definition~\eqref{E:LEVELSETGENERATORFORNEWTIMEFUNCTION},
	Lemma~\ref{L:COMMUTATORSTOCOORDINATES},
	Lemma~\ref{L:SCHEMATICSTRUCTUREOFVARIOUSTENSORSINTERMSOFCONTROLVARS},
	\eqref{E:MUTRANSVERSALCONVEXITY},
	and the estimates of Lemma~\ref{L:LINFTYESTIMATESFORROUGHTIMEFUNCTIONANDDERIVATIVES} and Prop.\,\ref{P:IMPROVEMENTOFAUXILIARYBOOTSTRAP}.
	
	\end{proof}

\section{Homeomorphism and diffeomorphism properties of $\Upsilon$ on $\MInteresting$ and
a description of the singular boundary in Cartesian coordinate space}
\label{S:HOMEOANDDIFFEOANDSINGULARBOUNDARYINCARTESIAN}
In this section, we reveal the homeomorphism and diffeomorphism properties of 
the change of variables map $\Upsilon(t,u,x^2,x^3) = (t,x^1,x^2,x^3)$ on 
the region $\MInteresting$.
We also reveal how the singular boundary $\mathcal{B}^{[0,\muxmulevelsetvalue_0]}$ is embedded in Cartesian coordinate space under $\Upsilon$,
i.e., we exhibit various properties of the set $\Upsilon(\mathcal{B}^{[0,\muxmulevelsetvalue_0]})$,
including its structure as a $\gfour$-null hypersurface in the Cartesian coordinate differential structure.
We refer to Remark~\ref{R:NONUNIQUENESSOFINTEGRALCURVESOFLUNIT} for a discussion of interesting
degeneracies in the vectorfield $\Lunit$ that occur along $\Upsilon(\mathcal{B}^{[0,\muxmulevelsetvalue_0]})$.

\subsection{Homeomorphism and diffeomorphism properties of $\Upsilon$ on $\MInteresting$}
\label{SS:HOMEOANDDIFFEOANDSINGULARBOUNDARY}
In the next proposition, we reveal the homeomorphism and diffeomorphism properties of 
the change of variables map $\Upsilon$. The proposition is crucial for translating results that we have
derived with respect to the geometric coordinates on $\MInteresting$ into results
with respect to the Cartesian coordinates on $\Upsilon(\MInteresting)$.
This will become apparent in the proof of Theorem~\ref{T:DEVELOPMENTANDSTRUCTUREOFSINGULARBOUNDARY}.

\begin{proposition}[Homeomorphism and diffeomorphism properties of $\Upsilon$ on $\MInteresting$]
	\label{P:CHOVGEOMETRICTOCARTESIANISINJECTIVEONREGIONWECAREABOUT}
	Assume the hypotheses and conclusions of Theorem~\ref{T:EXISTENCEUPTOTHESINGULARBOUNDARYATFIXEDKAPPA}
	for $\muxmulevelsetvalue \in [0,\muxmulevelsetvalue_0]$. Recall that $\MInteresting$ is the set defined in \eqref{E:INTERESTINGDEVELOPMENTOFDATA}
	and depicted in Fig.\,\ref{F:MINTERESTINGDEVELOPMENT}.
	Then the change of variables map $\Upsilon(t,u,x^2,x^3) = (t,x^1,x^2,x^3)$ enjoys the following properties.
	
	\begin{itemize}
	 \item $\Upsilon$ is continuous, injective map on the compact set $\MInteresting$
		defined in \eqref{E:INTERESTINGDEVELOPMENTOFDATA}.
		In particular, $\Upsilon$ is a homeomorphism from $\MInteresting$ onto its image.
		\item The following estimates hold on $\MInteresting$:
		\begin{align}
			\| \Upsilon \|_{C_{\textnormal{geo}}^{3,1}(\MInteresting)} 
			& \leq C,
				\label{E:C31BOUNDCHOVGEOMETRICTOCARTESIANONWECAREABOUT} 
					\\
		\mbox{\upshape det} \frac{\partial \Upsilon(t,u,x^2,x^3)}{\partial (t,u,x^2,x^3)} 
		& = \upmu \frac{\Speed^2}{X^1}
			= 
			- \left\lbrace 1 + \mathcal{O}(\mathring{\upalpha}) \right\rbrace \upmu.
			\label{E:JACOBIANDETERMINANTBOUNDCHOVGEOMETRICTOCARTESIANONWECAREABOUT}
		\end{align}
	\item $\Upsilon$ is a global diffeomorphism on the subset
		$\MInteresting \backslash \mathcal{B}^{[0,\muxmulevelsetvalue_0]}$,
		i.e., it is a diffeomorphism away from the singular boundary.
	\end{itemize}
\end{proposition}

\begin{proof}
	Using Lemma~\ref{L:COMMUTATORSTOCOORDINATES} and
	Prop.\,\ref{P:IMPROVEMENTOFAUXILIARYBOOTSTRAP}, we compute that
	$\| \Upsilon \|_{W_{\textnormal{geo}}^{4,\infty}(\MInteresting)} \lesssim 1$.
	From this estimate and
	\eqref{E:SOBOELVEMBEDDINGINTERESTINGREGIONRELYINGONQUASICONVEXITY},
	we further deduce that
	$\| \Upsilon \|_{C_{\textnormal{geo}}^{3,1}(\MInteresting)} \lesssim 1$ as desired.
	
	\eqref{E:JACOBIANDETERMINANTBOUNDCHOVGEOMETRICTOCARTESIANONWECAREABOUT}
	follows from same the arguments we used to prove \eqref{E:GEOTOCARTESIANJACOBIANDETERMINANTESTIMATE}.
	From these facts, 
	the fact that $\upmu$ is positive on $\MInteresting \backslash \mathcal{B}^{[0,\muxmulevelsetvalue_0]}$,
	and the inverse function theorem, 
	we deduce that
	$\Upsilon$ is a local diffeomorphism on 
	$\MInteresting \backslash \mathcal{B}^{[0,\muxmulevelsetvalue_0]}$.
	
	The rest of proof is similar to the proof of
	Prop.\,\ref{P:HOMEOMORPHICANDDIFFEOMORPHICEXTENSIONOFCARTESIANCOORDINATES}.
	We will silently use 
	the following results from Prop.\,\ref{P:INTERESTINGREGIONFOLIATEDBYINTERESTINGTIMFEFUNCTION}:
	the vectorfield $\levelsetgeneratornewtimefunction$
	is the partial derivative with respect to $u$ in the coordinate system $(\newtimefunction,u,x^2,x^3)$
	and the vectorfield
	$\partialderivativewithrespecttonewtimefunction$
	is the partial derivative with respect to $\newtimefunction$ in the coordinate system $(\newtimefunction,u,x^2,x^3)$.
	Moreover, we will use the notation ``$*$'' to denote
	any quantity that is pointwise bounded in magnitude by $\mathcal{O}(\mathring{\upalpha})$.
	
	To complete the proof of the proposition, we must show that
	$\Upsilon$ is injective on $\MInteresting$.
	In view of the diffeomorphism properties of $\InterestingCHOV$ shown in
	Prop.\,\ref{P:INTERESTINGREGIONFOLIATEDBYINTERESTINGTIMFEFUNCTION},
	we see that it suffices to show that $\Upsilon \circ \InverseInterestingCHOV$,
	which maps $(\newtimefunction,u,x^2,x^3) \rightarrow (t,x^1,x^2,x^3)$,
	is injective on the domain
	$[\timefunction_0,0] \times [- \rightu,\leftu] \times \mathbb{T}^2$.
	
	As an intermediate step, we will show that the map
	$(\newtimefunction,u,x^2,x^3) \rightarrow (\newtimefunction,x^1,x^2,x^3)$
	is injective on the domain
	$[\timefunction_0,0] \times [- \rightu,\leftu] \times \mathbb{T}^2$
	and is a diffeomorphism away from the crease.
	Below we will show that the map
	$(\newtimefunction,x^1,x^2,x^3) \rightarrow (t,x^1,x^2,x^3)$
	is also injective, which will complete the proof.
	To achieve the intermediate step, we will show that
	for
	every fixed
	$(\timefunction,x^2,x^3) \in [\timefunction_0,0] \times \mathbb{T}^2$,
	the map $u \rightarrow x^1(\timefunction,u,x^2,x^3)$
	is strictly decreasing (here we stress that $\timefunction$ denotes a fixed value of $\newtimefunction$)
	on the domain $[- \rightu,\leftu]$
	with $\levelsetgeneratornewtimefunction x^1 < 0$ away from
	$\InterestingCHOV\left( \partial_- \mathcal{B}^{[0,\muxmulevelsetvalue_0]}  \right)$,
	i.e., away from the image of the crease under $\InterestingCHOV$.
	To this end, we first use 
	Lemma~\ref{L:COMMUTATORSTOCOORDINATES},
	Lemma~\ref{L:SCHEMATICSTRUCTUREOFVARIOUSTENSORSINTERMSOFCONTROLVARS},
	\eqref{E:LEVELSETGENERATORFORNEWTIMEFUNCTION}, 
	and the estimates of 
	Lemma~\ref{L:DIFFEOMORPHICEXTENSIONOFROUGHCOORDINATES}
	and
	Props.\,\ref{P:IMPROVEMENTOFAUXILIARYBOOTSTRAP} and
	\ref{P:INTERESTINGREGIONFOLIATEDBYINTERESTINGTIMFEFUNCTION}
	to compute the following, where 
	$\phi$ is the cut-off from Definition~\ref{D:WTRANSANDCUTOFF}:
	\begin{align} \label{E:LEVELSETOFNEWTIMEFUNCTIONGENERATORAPPLIEDTOCOORDINATEX1}
		\levelsetgeneratornewtimefunction x^1
		& =
		\begin{cases}
			\upmu(-1 + *)
			& \mbox{in } 
			\InterestingCHOV\left( \MLeft \right)
			\\
			-(1 + *) \upmu
			-
			(1 + *)
			\frac{\muX \upmu}{\Lunit \upmu}
			& \mbox{in } 
			\InterestingCHOV\left( \MSingular \right)
			\\
			-(1 + *) \upmu
			+ 
			\phi 
			(1 + *)
			\frac{\muxmulevelsetvalue_0}{\Lunit \upmu} 
	& \mbox{in } 
		 	\InterestingCHOV\left( \MRight \right).
	\end{cases}
	\end{align}
	Next, we recall that by Prop.\,\ref{P:SHARPCONTROLOFMUANDDERIVATIVES}
	and definition \eqref{E:SINGULARDEVELOPMENT},
	in $\InterestingCHOV\left( \MSingular \right)$,
	we have $\Lunit \upmu \approx - 1$
	and $\muX \upmu < 0$ except along $\datahypfortimefunctiontwoarg{0}{[\timefunction_0,0]}$,
	where $\muX \upmu$ vanishes.
	We also recall that in $\InterestingCHOV\left( \MInteresting \right)$,
	$\upmu$ is positive, except on $\InterestingCHOV\left( \partial_- \mathcal{B}^{[0,\muxmulevelsetvalue_0]} \right)$,
	where $\partial_- \mathcal{B}^{[0,\muxmulevelsetvalue_0]} \subset \MSingular$.
	From these facts and \eqref{E:LEVELSETOFNEWTIMEFUNCTIONGENERATORAPPLIEDTOCOORDINATEX1},
	we see that
	$\levelsetgeneratornewtimefunction x^1 < 0$, 
	except along 
	$\InterestingCHOV\left( \partial_- \mathcal{B}^{[0,\muxmulevelsetvalue_0]} \right)$,
	where $\levelsetgeneratornewtimefunction x^1|_{\InterestingCHOV\left( \partial_- \mathcal{B}^{[0,\muxmulevelsetvalue_0]} \right)} = 0$
	(recall that $\partial_- \mathcal{B}^{[0,\muxmulevelsetvalue_0]} = \twoargmumuxtorus{0}{0} \subset \lbrace \newtimefunction = 0 \rbrace$,
	and that $\twoargmumuxtorus{0}{0}$ 
	is precisely the subset of $\MInteresting$ along which $\upmu = \muX \upmu = 0$).
	Moreover, Prop.\,\ref{P:INTERESTINGREGIONFOLIATEDBYINTERESTINGTIMFEFUNCTION} implies that
	every $u$-parameterized integral curve of $\levelsetgeneratornewtimefunction$
	in $\InterestingCHOV\left(\inthyp{0}{[- \rightu,\leftu]} \right)$
	(i.e., in the level set $\lbrace \newtimefunction = 0 \rbrace$)
	intersects the crease
	at a unique value of $u \in [-\frac{\interestingu}{2},\frac{\interestingu}{2}]$.
	In total, we have shown that
	along every $u$-parameterized integral curve of 
	$\levelsetgeneratornewtimefunction$
	in
	$[\timefunction_0,0] \times [- \rightu,\leftu] \times \mathbb{T}^2$,
	we have $\levelsetgeneratornewtimefunction x^1 < 0$, except 
	at possibly one point on the integral curve.
	From this fact and the mean value theorem, we conclude that
	at every fixed
	$(\timefunction,x^2,x^3) \in [\timefunction_0,0] \times \mathbb{T}^2$,
	the map $u \rightarrow x^1(\timefunction,u,x^2,x^3)$
	is strictly decreasing.
	We have therefore shown that the map
	$(\newtimefunction,u,x^2,x^3) \rightarrow (\newtimefunction,x^1,x^2,x^3)$
	is injective on the domain
	$[\timefunction_0,0] \times [- \rightu,\leftu] \times \mathbb{T}^2$
	and is a $C^1$ diffeomorphism away from 
	$\InterestingCHOV\left( \partial_- \mathcal{B}^{[0,\muxmulevelsetvalue_0]} \right)$.
	
	Next, for use below,
	we use \eqref{E:PARTIALDERIVATIVEWITHRESPECTTOTIMEFUNCTION},
	Lemma~\ref{L:COMMUTATORSTOCOORDINATES},
	Lemma~\ref{L:SCHEMATICSTRUCTUREOFVARIOUSTENSORSINTERMSOFCONTROLVARS},
	and the estimates of
	Props.\,\ref{P:IMPROVEMENTOFAUXILIARYBOOTSTRAP} 
	and
	Prop.\,\ref{P:INTERESTINGREGIONFOLIATEDBYINTERESTINGTIMFEFUNCTION},
	to compute that in $\InterestingCHOV\left(\MInteresting \right)$,
	we have
	$
	\partialderivativewithrespecttonewtimefunction x^1
	\approx 
	\geop{t} x^1
	\approx
	\Lunit x^1
	- 
	\Lunit^A \geop{x^A} x^1
	\approx
	\Lunit^1
	\approx 1
	$.
	
	The rest of the proof now mirrors
	the proof of the injectivity of $\Upsilon$ on $\twoargMrough{[\timefunction_0,\timefunctionboot],[- \rightu,\leftu]}{\muxmulevelsetvalue}$
	provided by Prop.\,\ref{P:HOMEOMORPHICANDDIFFEOMORPHICEXTENSIONOFCARTESIANCOORDINATES},
	where the estimate 
	$
	\partialderivativewithrespecttonewtimefunction x^1
	\approx 1
	$
	plays the role of the estimate \eqref{E:ROUGHTIMEDERIVATIVEOFCARTESIANX1ISAPPROXIMATELYUNITY} used in that proof.
	We will sketch the details.
	Specifically, the argument requires that we show that
	$
	\frac{1}{\partial_t \newtimefunction} > 0
	$
	on $\MInteresting$,
	except possibly when $\newtimefunction = 0$,
	where $\partial_t$ is the Cartesian partial time derivative vectorfield.
	To prove this result,
	we first use 
	\eqref{E:LSMALLDEF},
	\eqref{E:XSMALL},
	\eqref{E:XSMALLINTERMSOFLSMALLANDVELOCITY},
	\eqref{E:WTRANSDEF},
	\eqref{E:IVPFORROUGHTTIMEFUNCTION},
	\eqref{E:CARTESIANPARTIALTTOCOMMUTATORS}, 
	Lemma~\ref{L:COMMUTATORSTOCOORDINATES},
	Lemma~\ref{L:SCHEMATICSTRUCTUREOFVARIOUSTENSORSINTERMSOFCONTROLVARS},
	\eqref{E:NEWTIMEFUNCTION},
	and the estimates of 
	Lemma~\ref{L:DIFFEOMORPHICEXTENSIONOFROUGHCOORDINATES}
	and
	Props.\,\ref{P:IMPROVEMENTOFAUXILIARYBOOTSTRAP} 
	and
	\ref{P:INTERESTINGREGIONFOLIATEDBYINTERESTINGTIMFEFUNCTION}
	to compute that
	$\partial_t = \Lunit + (1 + *) X + * \Yvf{2} + * \Yvf{3}$
	and that the following estimates hold:
	\begin{align} \label{E:PARTIALDERIVATIVEOFCARTESIANTWRTNEWTIMEFUNCTIONX1X2X3FIXED}
		\frac{1}{\partial_t \newtimefunction}
		\approx
		& 
		\begin{cases}
			 1, & \mbox{in } 
			\MLeft,
			\\
			\frac{1}{1  - \frac{\muX \upmu}{\upmu}}
			=
			\frac{\upmu}{\upmu - \muX \upmu},
			& \mbox{in } 
			\MSingular,
				\\
			\frac{1}{1 + \frac{\muxmulevelsetvalue_0 \phi}{\upmu}}
			=
			\frac{\upmu}{\upmu + \muxmulevelsetvalue_0 \phi},
			& \mbox{in } 
		 	\MRight.
	\end{cases}
	\end{align}
	We now recall that $\muX \upmu|_{\datahypfortimefunctiontwoarg{-\muxmulevelsetvalue}{[\timefunction_0,0]}} = - \muxmulevelsetvalue$
	and that within $\MInteresting$,
	$\upmu$ can vanish only on the singular boundary $\mathcal{B}^{[0,\muxmulevelsetvalue_0]} \subset \inthyp{0}{[- \rightu,\leftu]}$,
	i.e., that $\upmu$ can vanish only when $\newtimefunction = 0$.
	From these facts, 
	definitions \eqref{E:LEFTDEVELOPMENT}--\eqref{E:RIGHTDEVELOPMENT},
	and \eqref{E:PARTIALDERIVATIVEOFCARTESIANTWRTNEWTIMEFUNCTIONX1X2X3FIXED},
	follows that
	$
	\frac{1}{\partial_t \newtimefunction} > 0
	$
	except possibly when $\newtimefunction = 0$,
	which is the desired result.
	This concludes the proof of the proposition.
	
\end{proof}

\subsection{Description of the singular boundary in Cartesian coordinate space}
\label{SS:DESCRIPTIONOFSINGULARBOUNDARYINCARTESIANSPACE}
In the next proposition, we reveal how the singular boundary $\mathcal{B}^{[0,\muxmulevelsetvalue_0]}$
is embedded in Cartesian coordinate space, i.e., we exhibit various properties of $\Upsilon(\mathcal{B}^{[0,\muxmulevelsetvalue_0]})$.
Among the main conclusions is that $\Upsilon(\mathcal{B}^{[0,\muxmulevelsetvalue_0]})$ is ruled, 
in a degenerate sense made clear in the proposition and Remark~\ref{R:NONUNIQUENESSOFINTEGRALCURVESOFLUNIT}, 
by integral curves of the $\gfour$-null vectorfield $\Lunit$.

\begin{proposition}[Description of the singular boundary in Cartesian coordinate space]
\label{P:DESCRIPTIONOFSINGULARBOUNDARYINCARTESIANSPACE}
Assume the hypotheses and conclusions of Theorem~\ref{T:EXISTENCEUPTOTHESINGULARBOUNDARYATFIXEDKAPPA}
for $\muxmulevelsetvalue \in [0,\muxmulevelsetvalue_0]$.
Let 
$\Upsilon$ be the change of variables map from geometric to Cartesian coordinates defined in \eqref{E:CHOVGEOTOCARTESIAN},
and let 
$\embeddingofsingularboundaryintogeometriccoordinatespace$ be the $C^{1,1}$ diffeomorphism defined in 
\eqref{E:EMBEDDINGOFSINGULARBOUNDARYINTOGEOMETRICCOORDINATESPACE}.
Recall that $\mathcal{B}^{[0,\muxmulevelsetvalue_0]}$ is a portion of the singular boundary in geometric coordinate space
and that $\partial_- \mathcal{B}^{[0,\muxmulevelsetvalue_0]}$ denotes the crease, 
viewed as subsets of geometric coordinate space.
Recall also that in Prop.\,\ref{P:PROPERTIESOFMSINGULARANDCREASE}, 
we showed that $\embeddingofsingularboundaryintogeometriccoordinatespace$
is a diffeomorphism from $[0,\muxmulevelsetvalue_0] \times \mathbb{T}^2$ onto $\mathcal{B}^{[0,\muxmulevelsetvalue_0]}$.
Then the following conclusions hold.

\medskip
\noindent \underline{\textbf{A homeomorphism onto $\Upsilon(\mathcal{B}^{[0,\muxmulevelsetvalue_0]})$ and a diffeomorphism onto
$\Upsilon(\mathcal{B}^{[0,\muxmulevelsetvalue_0]} \backslash \partial_- \mathcal{B}^{[0,\muxmulevelsetvalue_0]})$}}.
$\Upsilon \circ \embeddingofsingularboundaryintogeometriccoordinatespace$ is a $C^{1,1}$ injective map 
from $[0,\muxmulevelsetvalue_0] \times \mathbb{T}^2$ onto $\Upsilon(\mathcal{B}^{[0,\muxmulevelsetvalue_0]})$
- the image of $\mathcal{B}^{[0,\muxmulevelsetvalue_0]}$ in Cartesian coordinate space - 
that satisfies:
\begin{align} \label{E:JACOBIANMATRIXFOREMBEDDINGSINGUULARBOUNDARYINCARTESIANSPACE}
\frac{\partial [\Upsilon \circ \embeddingofsingularboundaryintogeometriccoordinatespace](\muxmulevelsetvalue,x^2,x^3)}{\partial (\muxmulevelsetvalue,x^2,x^3)}
& = 
\frac{\partial (t,x^1,x^2,x^3)}{\partial (\muxmulevelsetvalue,x^2,x^3)}
=
\begin{pmatrix} 
\frac{\muxmulevelsetvalue}
{(\geop{u} \upmu) \geop{t} \muX \upmu
-
(\geop{t} \upmu) \geop{u} \muX \upmu 
	}				& * & *
					\\
\frac{\left\lbrace
	-1 + *
\right\rbrace
	\muxmulevelsetvalue}
{(\geop{u} \upmu) \geop{t} \muX \upmu
-
(\geop{t} \upmu) \geop{u} \muX \upmu 
} 
	& * & *
						\\
				0 & 1 & *  
					\\
				0 & * & 1  
\end{pmatrix},
\end{align}
where $*$ denotes terms of size $\mathcal{O}(\mathring{\upalpha})$
and the denominator terms satisfy 
$
(\geop{u} \upmu) \geop{t} \muX \upmu
-
(\geop{t} \upmu) \geop{u} \muX \upmu 
\approx 1
$.
Moreover, on $(0,\muxmulevelsetvalue_0] \times \mathbb{T}^2$, 
$\Upsilon \circ \embeddingofsingularboundaryintogeometriccoordinatespace$ is an embedding, i.e., 
the differential of $\Upsilon \circ \embeddingofsingularboundaryintogeometriccoordinatespace$ is injective on 
$(0,\muxmulevelsetvalue_0] \times \mathbb{T}^2$,
and $\Upsilon \circ \embeddingofsingularboundaryintogeometriccoordinatespace$ 
is a diffeomorphism from $(0,\muxmulevelsetvalue_0] \times \mathbb{T}^2$ onto its image
$\Upsilon(\mathcal{B}^{[0,\muxmulevelsetvalue_0]} \backslash \partial_- \mathcal{B}^{[0,\muxmulevelsetvalue_0]})$
in Cartesian coordinate space.

Furthermore, the map 
$(z,x^2,x^3) \rightarrow [\Upsilon \circ \embeddingofsingularboundaryintogeometriccoordinatespace](\sqrt{z},x^2,x^3)$
is a $C^{1,1/2}$
diffeomorphism from $[0,\muxmulevelsetvalue_0^2] \times \mathbb{T}^2$ onto its image
$\Upsilon(\mathcal{B}^{[0,\muxmulevelsetvalue_0]})$
in Cartesian coordinate space, i.e., 
$\Upsilon(\mathcal{B}^{[0,\muxmulevelsetvalue_0]})$ is 
a $C^{1,1/2}$ embedded submanifold-with-boundary of Cartesian coordinate space.

\medskip
\noindent \underline{\textbf{A description of the $\gfour$-spacelike embedded tori $\Upsilon(\twoargmumuxtorus{0}{-\muxmulevelsetvalue})$}}.
For each fixed $\muxmulevelsetvalue \in [0,\muxmulevelsetvalue_0]$,
the map $(x^2,x^3) \rightarrow [\Upsilon \circ \embeddingofsingularboundaryintogeometriccoordinatespace](\muxmulevelsetvalue,x^2,x^3)$
is a diffeomorphism from $\mathbb{T}^2$ onto $\Upsilon(\twoargmumuxtorus{0}{-\muxmulevelsetvalue})$,
where $\twoargmumuxtorus{0}{-\muxmulevelsetvalue}$ is the $\upmu$-adapted torus defined in \eqref{E:MUXMUTORI}.
In particular, the differential of this map is injective,
and $\Upsilon(\twoargmumuxtorus{0}{-\muxmulevelsetvalue})$ is an embedded $C^{1,1}$ graph over 
$\mathbb{T}^2$ in Cartesian coordinate space such that at each $q \in \twoargmumuxtorus{0}{-\muxmulevelsetvalue}$,
$\Upsilon(\twoargmumuxtorus{0}{-\muxmulevelsetvalue})$
is spacelike with respect to $\gfour|_{\Upsilon(q)}$ at the point $\Upsilon(q)$.

\medskip
\noindent \underline{\textbf{The causal structure of 
$\Upsilon(\mathcal{B}^{[0,\muxmulevelsetvalue_0]} \backslash \partial_- \mathcal{B}^{[0,\muxmulevelsetvalue_0]})$}}.
The vectorfield $\nullgeneratorofsingularboundary$ defined by:
\begin{align} \label{E:NULLGENERATOROFSINGULARBOUNDARYINGEOMETRICCOORDINATES}
	\nullgeneratorofsingularboundary
	& \eqdef
	\Lunit
	-
	\frac{\Lunit \upmu}{\muX \upmu} \muX
\end{align}
enjoys the following properties: 
\begin{itemize}
\item $\nullgeneratorofsingularboundary$ 
	is well-defined on and tangent to $\mathcal{B}^{[0,\muxmulevelsetvalue_0]} \backslash \partial_- \mathcal{B}^{[0,\muxmulevelsetvalue_0]}$,
	viewed as a subset of geometric coordinate space (which is the target of $\embeddingofsingularboundaryintogeometriccoordinatespace$).
\item The integral curves of $\nullgeneratorofsingularboundary$ 
	foliate $\mathcal{B}^{[0,\muxmulevelsetvalue_0]} \backslash \partial_- \mathcal{B}^{[0,\muxmulevelsetvalue_0]}$
	and satisfy $\nullgeneratorofsingularboundary u > 0$ and $\nullgeneratorofsingularboundary t > 0$
	along $\mathcal{B}^{[0,\muxmulevelsetvalue_0]} \backslash \partial_- \mathcal{B}^{[0,\muxmulevelsetvalue_0]}$.
	In particular, the integral curves are transversal to the characteristics $\nullhyparg{u}$.
\item Along $\mathcal{B}^{[0,\muxmulevelsetvalue_0]} \backslash \partial_- \mathcal{B}^{[0,\muxmulevelsetvalue_0]}$, 
	$\gfour(\nullgeneratorofsingularboundary,V) = 0$ holds for
	\underline{every} vectorfield $V = V^t \geop{t} + V^u \geop{u} + V^2 \geop{2} + V^3 \geop{3}$
	on $\mathcal{B}^{[0,\muxmulevelsetvalue_0]} \backslash \partial_- \mathcal{B}^{[0,\muxmulevelsetvalue_0]}$, regardless of whether
	$V$ is tangent to $\mathcal{B}^{[0,\muxmulevelsetvalue_0]} \backslash \partial_- \mathcal{B}^{[0,\muxmulevelsetvalue_0]}$.
	In particular, on $\mathcal{B}^{[0,\muxmulevelsetvalue_0]} \backslash \partial_- \mathcal{B}^{[0,\muxmulevelsetvalue_0]}$, 
	we have that $\gfour(\nullgeneratorofsingularboundary,\nullgeneratorofsingularboundary) = 0$.
\item For every $q \in \mathcal{B}^{[0,\muxmulevelsetvalue_0]} \backslash \partial_- \mathcal{B}^{[0,\muxmulevelsetvalue_0]}$, 
	the pushforward vectorfield
	$
	[d_{\textnormal{geo}} \Upsilon(q)] \cdot \nullgeneratorofsingularboundary(q)
	$,
	which is tangent to $\Upsilon(\mathcal{B}^{[0,\muxmulevelsetvalue_0]} \backslash \partial_- \mathcal{B}^{[0,\muxmulevelsetvalue_0]})$, 
	is equal to $\Lunit|_{\Upsilon(q)} = [\Lunit^{\alpha} \partial_{\alpha}]|_{\Upsilon(q)}$.
\item For every $q \in \mathcal{B}^{[0,\muxmulevelsetvalue_0]} \backslash \partial_- \mathcal{B}^{[0,\muxmulevelsetvalue_0]}$,
	$\Lunit|_{\Upsilon(q)}$ is $\gfour|_{\Upsilon(q)}$-orthogonal to the tangent space of
	$\Upsilon(\mathcal{B}^{[0,\muxmulevelsetvalue_0]} \backslash \partial_- \mathcal{B}^{[0,\muxmulevelsetvalue_0]})$
	at $\Upsilon(q)$. 
\item $\nullgeneratorofsingularboundary$ is the unique such vectorfield with the above properties.
\end{itemize}
In particular, in the differential structure on spacetime induced by the Cartesian coordinates, 
there exist integral curves of the $\gfour$-null vectorfield $\Lunit = \Lunit^{\alpha}\partial_{\alpha}$ 
that foliate $\Upsilon(\mathcal{B}^{[0,\muxmulevelsetvalue_0]})$
and that are everywhere $\gfour$-orthogonal to $\Upsilon(\mathcal{B}^{[0,\muxmulevelsetvalue_0]} \backslash \partial_- \mathcal{B}^{[0,\muxmulevelsetvalue_0]})$.
Considering also that the $\Upsilon(\twoargmumuxtorus{0}{-\muxmulevelsetvalue})$ are $\gfour$-spacelike,
we see that $\Upsilon(\mathcal{B}^{[0,\muxmulevelsetvalue_0]} \backslash \partial_- \mathcal{B}^{[0,\muxmulevelsetvalue_0]})$
is a $\gfour$-null hypersurface in the Cartesian coordinate differential structure.

\end{proposition}

\begin{remark}[Non-uniqueness of the integral curves of $\Lunit^{\alpha} \partial_{\alpha}$ along $\Upsilon(\mathcal{B}^{[0,\muxmulevelsetvalue_0]})$]
\label{R:NONUNIQUENESSOFINTEGRALCURVESOFLUNIT}
Note that along $\Upsilon(\mathcal{B}^{[0,\muxmulevelsetvalue_0]})$, 
even though the scalar functions $\muX \Lunit^{\beta} = \upmu X \Lunit^{\beta}$ remain bounded,
the non-$\upmu$-weighted quantities 
$X \Lunit^{\beta} = X^{\alpha} \partial_{\alpha} \Lunit^{\beta}$
can blow up, due to the vanishing of $\upmu$ there.
Since $X^{\alpha} \partial_{\alpha}$ is a non-degenerate (i.e., everywhere non-zero and bounded) 
vectorfield in the Cartesian differential structure,
it follows that the Cartesian partial derivatives
$\partial_{\alpha} \Lunit^{\beta}$ can blow up along $\Upsilon(\mathcal{B}^{[0,\muxmulevelsetvalue_0]})$.
Hence, in the Cartesian differential structure, the vectorfield $\Lunit$ does not have sufficient regularity
to ensure uniqueness of its integral curves up to $\Upsilon(\mathcal{B}^{[0,\muxmulevelsetvalue_0]})$,
i.e., standard uniqueness theorems would require Lipschitz regularity for $\Lunit^{\alpha}$.
This lack of uniqueness is the mechanism that allows for 
the existence of integral curves of $\Lunit$ that foliate $\Upsilon(\mathcal{B}^{[0,\muxmulevelsetvalue_0]})$.

More precisely, consider any fixed point $p_0 \in \mathcal{B}^{[0,\muxmulevelsetvalue_0]} \backslash \partial_- \mathcal{B}^{[0,\muxmulevelsetvalue_0]}$,
and let $u_0$ denote the eikonal function evaluated at $p_0$
(i.e., the $u$-coordinate of $p_0$ in geometric coordinates).
By Prop.\,\ref{P:DESCRIPTIONOFSINGULARBOUNDARYINCARTESIANSPACE}, 
there exists an interval $I$ of $u$-values containing $u_0$ and a 
unique integral curve 
$\upgamma_{p_0}: I \to \mathcal{B}^{[0,\muxmulevelsetvalue_0]} \backslash \partial_- \mathcal{B}^{[0,\muxmulevelsetvalue_0]}$ 
of $\nullgeneratorofsingularboundary$ in geometric coordinate space satisfying:
\begin{align} \label{E:INTEGRALCURVETANGENTTOB}
	\frac{\mathrm{d}}{\mathrm{d} u} \upgamma_{p_0}(u) 
	& 
	= 
	\nullgeneratorofsingularboundary \circ \upgamma_p(u) & & \upgamma_p(u_0) = p_0.
\end{align} 
It also follows from Prop.\,\ref{P:DESCRIPTIONOFSINGULARBOUNDARYINCARTESIANSPACE} 
that the pushforward of the tangent vector $\frac{\mathrm{d}}{\mathrm{d} u} \upgamma_{p_0}(u_0)$ under 
$\Upsilon$ is precisely $\Lunit|_{\Upsilon(p_0)}$. 
Similarly, since $p_0 \in \nullhyparg{u_0}$, 
$\Lunit$ is tangent to $\nullhyparg{u_0}$, and $\Lunit t = 1$,
there is an interval $J$ of $t$-values and a $t$-parameterized integral curve 
$\Lambda: J \to \nullhyparg{u_0}$ of $\Lunit$ in geometric coordinate space,
that is \emph{tangent to $\nullhyparg{u_0}$} such that $p_0 = \Lambda(t_0)$ for some $t_0 \in J$
(here, $t_0$ is the Cartesian time function evaluated at $p_0$).
Hence, in Cartesian coordinate space, there are two distinct integral curves through $\Upsilon(p_0)$, 
namely $\Upsilon \circ \upgamma$ and $\Upsilon \circ \Lambda$,
each with the same tangent vector $\Lunit|_{\Upsilon(p_0)}$; 
see Fig.\,\ref{F:NONUNIQUENESSOFINTEGRALCURVESOFL}.
\end{remark}

\begin{center}
	\begin{figure}[ht]  
		\begin{overpic}[scale=.7, grid = false, tics=5, trim=-.5cm -1cm -1cm -.5cm, clip]{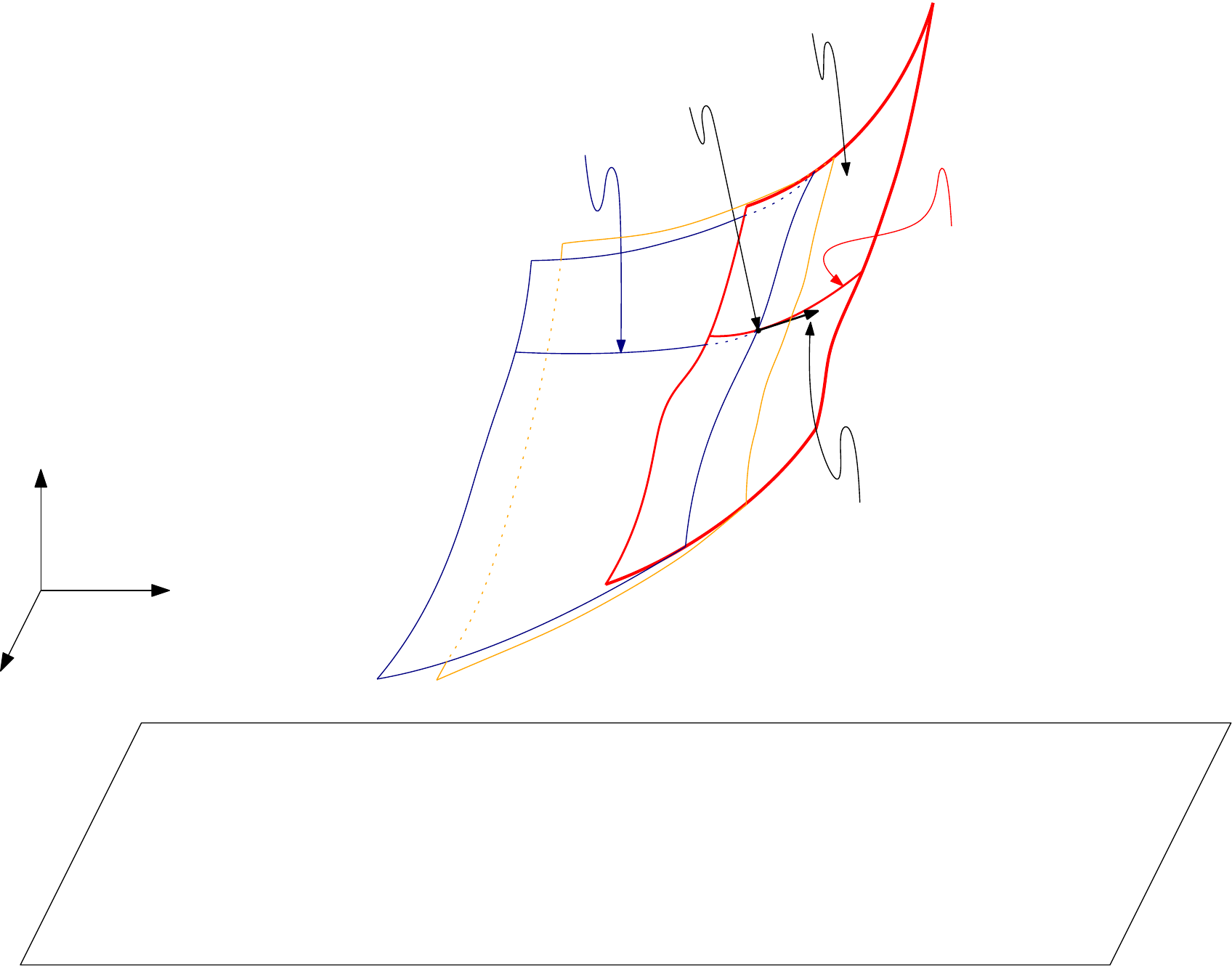}
			\put (40,10) {$\Sigma_0$}
			\put (4,25) {$(x^2,x^3) \in \mathbb{T}^2$}
			\put (4,40) {$t$}
			\put(51,71) {$\Upsilon(p_0)$}
			\put (10,34) {$x^1 \in\mathbb{R}$}
			\put (64,36) {$[\Lunit^{\alpha} \partial_{\alpha}]|_{\Upsilon(p_0)} = [d_{\textnormal{geo}} \Upsilon]|_{p_0} \cdot \upgamma_{p_0}(u_0)$}
			\put (58,76) {$\Upsilon(\mathcal{B}^{[0,\muxmulevelsetvalue_0]})$}
			\put (73,58) {$\Upsilon \circ \upgamma|_{p_0}$}
			\put (29,24.5) {$\nullhyparg{u_0}$}
			\put (38,25) {$\nullhyparg{u_1}$, \mbox{where} $u_1 < u_0$}
			\put (43,67) {$\Upsilon \circ \Lambda$}
	\end{overpic}
		\caption{Non-uniqueness of integral curves of $\Lunit$ along $\Upsilon(\mathcal{B}^{[0,\muxmulevelsetvalue_0]})$ in Cartesian coordinates.}
	\label{F:NONUNIQUENESSOFINTEGRALCURVESOFL}
	\end{figure}
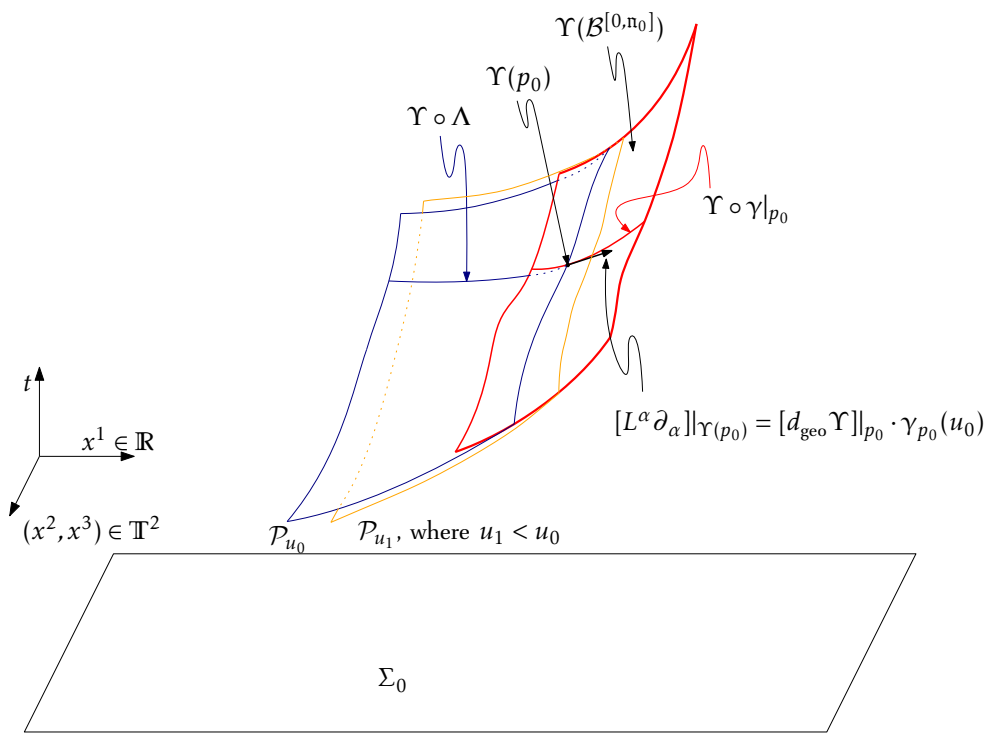
\end{center}

\begin{proof}[Proof of Prop.\,\ref{P:DESCRIPTIONOFSINGULARBOUNDARYINCARTESIANSPACE}]
	The regularity and injectivity properties of
	$\Upsilon \circ \embeddingofsingularboundaryintogeometriccoordinatespace$
	follow from Props.\,\ref{P:PROPERTIESOFMSINGULARANDCREASE} and \ref{P:CHOVGEOMETRICTOCARTESIANISINJECTIVEONREGIONWECAREABOUT}.
	At the end of the proof, we will prove \eqref{E:JACOBIANMATRIXFOREMBEDDINGSINGUULARBOUNDARYINCARTESIANSPACE}
	and the denominator term estimate
	$
(\geop{u} \upmu) \geop{t} \muX \upmu
-
(\geop{t} \upmu) \geop{u} \muX \upmu 
\approx 1
$.
	Taking these for granted for the time being, we
	deduce from the explicit form of \eqref{E:JACOBIANMATRIXFOREMBEDDINGSINGUULARBOUNDARYINCARTESIANSPACE}
	that the differential of $\Upsilon \circ \embeddingofsingularboundaryintogeometriccoordinatespace$ 
	(with respect to $(\muxmulevelsetvalue,x^2,x^3)$)
	is injective on $(0,\muxmulevelsetvalue_0] \times \mathbb{T}^2$
	and thus $\Upsilon \circ \embeddingofsingularboundaryintogeometriccoordinatespace$ 
	is a diffeomorphism from $(0,\muxmulevelsetvalue_0] \times \mathbb{T}^2$ onto its image,
	as is desired.

\medskip
\noindent \textbf{Proof of the properties of the map 
$(z,x^2,x^3) \rightarrow [\Upsilon \circ \embeddingofsingularboundaryintogeometriccoordinatespace](\sqrt{z},x^2,x^3)$}:
We set $z \eqdef \muxmulevelsetvalue^2$.
From the previous paragraph,
\eqref{E:JACOBIANMATRIXFOREMBEDDINGSINGUULARBOUNDARYINCARTESIANSPACE}
(note that the first column of the matrix on RHS~\eqref{E:JACOBIANMATRIXFOREMBEDDINGSINGUULARBOUNDARYINCARTESIANSPACE}
features a linear factor of $\muxmulevelsetvalue$), 
and the chain rule,
it follows that the map
$
(z,x^2,x^3) \rightarrow  [\Upsilon \circ \embeddingofsingularboundaryintogeometriccoordinatespace](\sqrt{z},x^2,x^3)$
is injective on $[0,\muxmulevelsetvalue_0^2] \times \mathbb{T}^2$
and that
$
\frac{\partial [\Upsilon \circ \embeddingofsingularboundaryintogeometriccoordinatespace](\sqrt{z},x^2,x^3)}{\partial (z,x^2,x^3)}
=
\begin{pmatrix} 
\frac{1}{2}
{(\geop{u} \upmu) \geop{t} \muX \upmu
-
\frac{1}{2}
(\geop{t} \upmu) \geop{u} \muX \upmu 
	}				& * & *
					\\
\frac{\left\lbrace
	- \frac{1}{2} + *
\right\rbrace
}
{(\geop{u} \upmu) \geop{t} \muX \upmu
-
(\geop{t} \upmu) \geop{u} \muX \upmu 
} 
	& * & *
						\\
				0 & 1 & *  
					\\
				0 & * & 1  
\end{pmatrix}
$
and that
$
\frac{\partial [\Upsilon \circ \embeddingofsingularboundaryintogeometriccoordinatespace](\sqrt{z},x^2,x^3)}{\partial (z,x^2,x^3)}
$
can be expressed as $C^{0,1}$ function of $(\muxmulevelsetvalue,x^2,x^3)$.
Since the map $z \rightarrow \sqrt{z}$ is $C^{0,1/2}$, it follows that
$
\frac{\partial [\Upsilon \circ \embeddingofsingularboundaryintogeometriccoordinatespace](\sqrt{z},x^2,x^3)}{\partial (z,x^2,x^3)}
$
can be expressed as a $C^{0,1/2}$ function of $(z,x^2,x^3)$.
Since the differential of the map
$
(z,x^2,x^3)
\rightarrow
[\Upsilon \circ \embeddingofsingularboundaryintogeometriccoordinatespace](\sqrt{z},x^2,x^3)
$
is clearly injective, we conclude that the map
is a $C^{1,1/2}$ diffeomorphism from $[0,\muxmulevelsetvalue_0^2] \times \mathbb{T}^2$
onto $\Upsilon(\mathcal{B}^{[0,\muxmulevelsetvalue_0]})$, as is desired.

\medskip

\noindent \textbf{Proof of the properties of $\Upsilon(\twoargmumuxtorus{0}{-\muxmulevelsetvalue})$}:
	From the explicit form of \eqref{E:JACOBIANMATRIXFOREMBEDDINGSINGUULARBOUNDARYINCARTESIANSPACE},
	we see that the last two columns of the matrix on RHS~\eqref{E:JACOBIANMATRIXFOREMBEDDINGSINGUULARBOUNDARYINCARTESIANSPACE}
	are linearly independent when $\mathring{\upalpha}$ is sufficiently small (even if $\muxmulevelsetvalue = 0$).
	It follows that for any $\muxmulevelsetvalue \in [0,\muxmulevelsetvalue_0]$, the map 
	$(x^2,x^3) \rightarrow [\Upsilon \circ \embeddingofsingularboundaryintogeometriccoordinatespace](\muxmulevelsetvalue,x^2,x^3)$
	is a diffeomorphism from $\mathbb{T}^2$ onto $\Upsilon(\twoargmumuxtorus{0}{-\muxmulevelsetvalue})$
	Since Lemma~\ref{L:CAUSALSTRUCTUREOFTORIANDLEVELSETSOFMUINGEOMETRICCOORDINATES} shows that each
	$\twoargmumuxtorus{0}{-\muxmulevelsetvalue}$ is $\gfour$-spacelike (in the differential structure of the geometric coordinates), 
	and since $\Upsilon|_{\twoargmumuxtorus{0}{-\muxmulevelsetvalue}}$ is a diffeomorphism,
	we immediately conclude that $\Upsilon(\twoargmumuxtorus{0}{-\muxmulevelsetvalue})$ is $\gfour$-spacelike.

\noindent \textbf{Proof of the properties of $\nullgeneratorofsingularboundary$}:
 The uniqueness statement made about $\nullgeneratorofsingularboundary$ follows from the fact that
	the differential of $\Upsilon \circ \embeddingofsingularboundaryintogeometriccoordinatespace$ is injective 
	and our assumption that 
	$
	[d_{\textnormal{geo}} \Upsilon(q)] \cdot \nullgeneratorofsingularboundary(q)
	$
	is equal to
	$
	[\Lunit^{\alpha} \partial_{\alpha}]|_{\Upsilon(q)}
	$.
	
Next, in view of \eqref{E:NULLGENERATOROFSINGULARBOUNDARYINGEOMETRICCOORDINATES},
we see that $\nullgeneratorofsingularboundary$ is well-defined at points where $\muX \upmu < 0$,
a condition that is satisfied on $\mathcal{B}^{[0,\muxmulevelsetvalue_0]} \backslash \partial_- \mathcal{B}^{[0,\muxmulevelsetvalue_0]}$
by \eqref{E:SINGULARBOUNDARYPORTION}--\eqref{E:CREASE}
and the fact that $\muX \upmu|_{\twoargmumuxtorus{0}{-\muxmulevelsetvalue}} \equiv - \muxmulevelsetvalue$.
We also note that $\nullgeneratorofsingularboundary \upmu = 0$
(i.e., $\nullgeneratorofsingularboundary$ is tangent to the singular boundary).
Furthermore, using Lemma~\ref{L:BASICPROPERTIESOFVECTORFIELDS}, 
$\muX \upmu|_{\twoargmumuxtorus{0}{-\muxmulevelsetvalue}} \equiv - \muxmulevelsetvalue$,
\eqref{E:IMPROVEDLEVELSETSTRUCTUREANDLOCATIONOFMIN},
and \eqref{E:BOUNDSONLMUINTERESTINGREGION},
we find that along $\mathcal{B}^{[0,\muxmulevelsetvalue_0]} \backslash \partial_- \mathcal{B}^{[0,\muxmulevelsetvalue_0]}$,
we have
$\nullgeneratorofsingularboundary t = 1 > 0$
and $\nullgeneratorofsingularboundary u = 
-
\frac{\Lunit \upmu}{\muX \upmu} > 0$,
as desired.

Next, using Lemma~\ref{L:BASICPROPERTIESOFVECTORFIELDS}, 
we compute that 
$\gfour(\nullgeneratorofsingularboundary,\Lunit) 
= 
\upmu
\frac{\Lunit \upmu}{\muX \upmu}
$,
$\gfour(\nullgeneratorofsingularboundary,\muX) 
= 
- \upmu
+
\upmu^2
\frac{\Lunit \upmu}{\muX \upmu}
$,
and
$\gfour(\nullgeneratorofsingularboundary,\geop{x^A}) 
= 0
$
for $A=2,3$. Note that on $\mathcal{B}^{[0,\muxmulevelsetvalue_0]} \backslash \partial_- \mathcal{B}^{[0,\muxmulevelsetvalue_0]}$,
where $\upmu = 0$, all of these inner products vanish.
Since $\lbrace \Lunit, \muX, \geop{2}, \geop{3} \rbrace$ spans the tangent space of 
$\mathcal{B}^{[0,\muxmulevelsetvalue_0]} \backslash \partial_- \mathcal{B}^{[0,\muxmulevelsetvalue_0]}$
(in the differential structure of the geometric coordinates),
it follows that along $\mathcal{B}^{[0,\muxmulevelsetvalue_0]} \backslash \partial_- \mathcal{B}^{[0,\muxmulevelsetvalue_0]}$,
$\nullgeneratorofsingularboundary$ is $\gfour$-orthogonal to \underline{every} vectorfield,
as we claimed in the proposition.

\medskip
\noindent \textbf{Proof of properties tied to foliations of 
$\mathcal{B}^{[0,\muxmulevelsetvalue_0]} \backslash \partial_- \mathcal{B}^{[0,\muxmulevelsetvalue_0]}$ and
$\Upsilon\left(\mathcal{B}^{[0,\muxmulevelsetvalue_0]} \backslash \partial_- \mathcal{B}^{[0,\muxmulevelsetvalue_0]} \right)$}:
Using \eqref{E:MUTRANSVERSALCONVEXITY}, \eqref{E:BOUNDSONLMUINTERESTINGREGION}, and 
$\muX \upmu|_{\twoargmumuxtorus{0}{-\muxmulevelsetvalue}} \equiv - \muxmulevelsetvalue$,
we find that $\nullgeneratorofsingularboundary \muX \upmu > 0$.
Hence, since \eqref{E:SINGULARBOUNDARYPORTION} implies that the level sets of $\muX \upmu$ foliate 
$\mathcal{B}^{[0,\muxmulevelsetvalue_0]} = \lbrace (t,u,x^2,x^3) \ | \ \upmu(t,u,x^2,x^3) = 0, 
	\,
- \muxmulevelsetvalue_0 \leq \muX \upmu(t,u,x^2,x^3) \leq 0, \, (x^2,x^3) \in \mathbb{T}^2 \rbrace$,
and since $\nullgeneratorofsingularboundary \muX \upmu > 0$ 
implies that $\nullgeneratorofsingularboundary$ is transversal to the level sets of $\muX \upmu$, 
we deduce that the integral curves of $\nullgeneratorofsingularboundary$
foliate $\mathcal{B}^{[0,\muxmulevelsetvalue_0]} \backslash \partial_- \mathcal{B}^{[0,\muxmulevelsetvalue_0]}$. 
Hence, the pushforward of
$
\nullgeneratorofsingularboundary|_{\mathcal{B}^{[0,\muxmulevelsetvalue_0]} \backslash \partial_- \mathcal{B}^{[0,\muxmulevelsetvalue_0]}}
$
by $\Upsilon$ is a vectorfield tangent to
$
\Upsilon(\mathcal{B}^{[0,\muxmulevelsetvalue_0]} \backslash \partial_- \mathcal{B}^{[0,\muxmulevelsetvalue_0]})
$, which is foliated by the integral curves.
We now come to the key point: for any
$q \in \mathcal{B}^{[0,\muxmulevelsetvalue_0]} \backslash \partial_- \mathcal{B}^{[0,\muxmulevelsetvalue_0]}$,
the pushforward of 
$
\nullgeneratorofsingularboundary|_q
$
by $\Upsilon(q)$
is $[\Lunit^{\alpha} \partial_{\alpha}]|_{\Upsilon(t,u,x^2,x^3)}$.
The reason is that the pushforward of
$
\muX|_q
$
by $\Upsilon$
is 
$[\upmu X^{\alpha} \partial_{\alpha}]|_{\Upsilon(q)}$,
which, in the Cartesian differential structure, 
vanishes along $\Upsilon(\mathcal{B}^{[0,\muxmulevelsetvalue_0]})$, 
where $\upmu \equiv 0$.

It remains for us to prove \eqref{E:JACOBIANMATRIXFOREMBEDDINGSINGUULARBOUNDARYINCARTESIANSPACE}.
We first complement the vectorfields $\muderivativevectorfield$ and 
	$\muxmuderivativevectorfield$ defined in
	\eqref{E:PARTIALDERIVATIVEWITHRESPECTTOMUATFIXEDMUXMU}--\eqref{E:PARTIALDERIVATIVEWITHRESPECTTOMUXMUATFIXEDMU}
	with the following pair of vectorfields ($A=2,3$):
	\begin{align} \label{E:PARTIALDERIVATIVEWITHRESPECTTOXAATFIXEDMUANDMUXMU}
		\derivativevectorfieldatfixemuandmuxmu{A}
		& \eqdef 
			\geop{x^A}
			-
			\left(\geop{x^A} \upmu \right) 
			\muderivativevectorfield
			+
			\left(\geop{x^A} \muX \upmu \right) 
			\muxmuderivativevectorfield.
\end{align}
As in \eqref{E:RELATIONSPARTIALDERIVATIVEWITHRESPECTTOMUATFIXEDMUXMU}, 
we compute (recalling that $\muxmuderivativevectorfield \muX \upmu =  -1$) that:
\begin{align} \label{E:RELATIONSPARTIALDERIVATIVEWITHRESPECTTOMUATFIXEDMUANDMUXMU}
	\derivativevectorfieldatfixemuandmuxmu{2} x^2
	& = 1, 
	& &
	\derivativevectorfieldatfixemuandmuxmu{2} \upmu 
	= 
	-
	\derivativevectorfieldatfixemuandmuxmu{2} \muX \upmu 
	= 
	\muderivativevectorfield x^3 
	= 0.
\end{align}
It follows from \eqref{E:RELATIONSPARTIALDERIVATIVEWITHRESPECTTOMUATFIXEDMUANDMUXMU} that
$
\derivativevectorfieldatfixemuandmuxmu{2}$ 
is the partial derivative with respect to $x^2$ in the coordinates 
$(\muxmulevelsetvalue,x^2,x^3)$ on 
$[0,\muxmulevelsetvalue_0] \times \mathbb{T}^2$,
i.e., on the domain of $\Upsilon \circ \embeddingofsingularboundaryintogeometriccoordinatespace$.
Similarly, 
$
\derivativevectorfieldatfixemuandmuxmu{3}$ 
is the partial derivative with respect to $x^3$ in the coordinates 
$(\muxmulevelsetvalue,x^2,x^3)$ on 
$[0,\muxmulevelsetvalue_0] \times \mathbb{T}^2$.
Hence, $\lbrace \muxmuderivativevectorfield, \derivativevectorfieldatfixemuandmuxmu{2}, \derivativevectorfieldatfixemuandmuxmu{3}\rbrace$ are
the coordinate partial derivative vectorfields on $[0,\muxmulevelsetvalue_0] \times \mathbb{T}^2$ 
and using
\eqref{E:PARTIALDERIVATIVEWITHRESPECTTOMUXMUATFIXEDMU},
\eqref{E:RELATIONSPARTIALDERIVATIVEWITHRESPECTTOMUXMUATFIXEDMU},
and
\eqref{E:PARTIALDERIVATIVEWITHRESPECTTOXAATFIXEDMUANDMUXMU},
we calculate that the Jacobian matrix on LHS~\eqref{E:JACOBIANMATRIXFOREMBEDDINGSINGUULARBOUNDARYINCARTESIANSPACE} 
can be expressed as:
\begin{align} 
\begin{split} \label{E:PROOFSTEPJACOBIANMATRIXFOREMBEDDINGSINGUULARBOUNDARYINCARTESIANSPACE}
\frac{\partial [\Upsilon \circ \embeddingofsingularboundaryintogeometriccoordinatespace](\muxmulevelsetvalue,x^2,x^3)}{\partial (\muxmulevelsetvalue,x^2,x^3)}
& = 
\frac{\partial (t,x^1,x^2,x^3)}{\partial (\muxmulevelsetvalue,x^2,x^3)}
	\\
&
=
\begin{pmatrix} 
\frac{- 
	\frac{\geop{u} \upmu}
	{\geop{t} \upmu}}
{\geop{u} \muX \upmu 
	- 
	\frac{(\geop{u} \upmu) \geop{t} \muX \upmu}
	{\geop{t} \upmu}
} 
& \derivativevectorfieldatfixemuandmuxmu{2} t  & \derivativevectorfieldatfixemuandmuxmu{3} t
					\\
\frac{1}
{\geop{u} \muX \upmu 
	- 
	\frac{(\geop{u} \upmu) \geop{t} \muX \upmu}
	{\geop{t} \upmu}
} 
\left\lbrace 
	\frac{\upmu}{\Speed^2}{X^1} 
	- 
	\frac{\geop{u} \upmu}
	{\geop{t} \upmu} 
	\left(
		\frac{\Lunit^1 X^1 + \Lunit^2 X^2 + \Lunit^3 X^3}{X^1}
	\right)
\right\rbrace	& \derivativevectorfieldatfixemuandmuxmu{2} x^1  & \derivativevectorfieldatfixemuandmuxmu{3} x^1
						\\
				0 & 1 & 0  
					\\
				0 & 0 & 1  
\end{pmatrix}.
\end{split}
\end{align}
\eqref{E:JACOBIANMATRIXFOREMBEDDINGSINGUULARBOUNDARYINCARTESIANSPACE} now follows from
\eqref{E:PARTIALDERIVATIVEWITHRESPECTTOMUATFIXEDMUXMU}--\eqref{E:PARTIALDERIVATIVEWITHRESPECTTOMUXMUATFIXEDMU},
\eqref{E:PARTIALDERIVATIVEWITHRESPECTTOXAATFIXEDMUANDMUXMU},
\eqref{E:PROOFSTEPJACOBIANMATRIXFOREMBEDDINGSINGUULARBOUNDARYINCARTESIANSPACE},	
Lemma~\ref{L:COMMUTATORSTOCOORDINATES},
Lemma~\ref{L:SCHEMATICSTRUCTUREOFVARIOUSTENSORSINTERMSOFCONTROLVARS},
the $L^{\infty}$ estimates of Prop.\,\ref{P:IMPROVEMENTOFAUXILIARYBOOTSTRAP},
\eqref{E:MUTRANSVERSALCONVEXITY}, \eqref{E:BOUNDSONLMUINTERESTINGREGION}, 
the fact that
$\muX \upmu|_{\twoargmumuxtorus{0}{-\muxmulevelsetvalue}} \equiv - \muxmulevelsetvalue$,
and the fact that by \eqref{E:MUXINTERMSOFGEOMETRICCOORDINATEVECTORFIELDS}, 
$\muX = \geop{u}$ along the singular boundary (where $\upmu \equiv 0$).

We have therefore proved the proposition.

\end{proof}

\section{The main results}
\label{S:MAINRESULTS}
In this section, we state and prove the main theorem of the paper.
The theorem provides an assimilated version of results we have already proved.

\begin{theorem}[The development and structure of the singular boundary]
	\label{T:DEVELOPMENTANDSTRUCTUREOFSINGULARBOUNDARY}
	Fix any of the compactly supported admissible simple isentropic plane symmetric ``background'' solutions $\RRiemannPS$
	from Def.\,\ref{AD:ADMISSIBLEBACKGROUND} (recall that 
	$\LRiemann$, $v^2$, $v^3$, $\Ent$, 
	$\vortrenormalized$,
	$\GradEnt$
	$\VortVort$
	and
	$\DivGradEnt$ vanish for these background solutions).
	Let
$(\RRiemann,\LRiemann,v^2,v^3,\Ent) \big|_{\Sigma_0}
\eqdef 
\left(\RRiemannpertinitial, \LRiemannpertinitial, \vtwopertinitial, \vthreepertinitial, \spertinitial \right)
$
be perturbed fluid data on the flat Cartesian hypersurface $\Sigma_0$,
as in \eqref{E:PERTURBEDBONAFIDEDATA},
and let $u|_{\Sigma_0} = - x^1$ be the initial condition of the eikonal function, 
as in \eqref{E:EIKONALEQUATION} and \eqref{AE:DATAFOREIKONALEQUATIONINPLANESYMMETRY}.
Assume the hypotheses and conclusions of Theorem~\ref{T:EXISTENCEUPTOTHESINGULARBOUNDARYATFIXEDKAPPA}.
	In particular, assume that $\Ntop \geq 24$,
	and that the quantity
			$
			\mathring{\Delta}_{\Sigma_0^{[-\farrightu,\leftu]}}^{\Ntop+1}
			$
			defined in \eqref{E:PERTURBATIONSMALLNESSINCARTESIANDIFFERENTIALSTRUCTURE}
			is sufficiently small, where
			$
			\mathring{\Delta}_{\Sigma_0^{[-\farrightu,\leftu]}}^{\Ntop+1}
			$
			is a Sobolev norm of the perturbation of the fluid data away from the background solution.
	Then the corresponding solution 
	exhibits the following properties.

\medskip

\noindent \underline{\textbf{Classical existence with respect to the geometric coordinates on $\MInteresting$}}.
\begin{itemize}
	\item There exists a compact region $\MInteresting$ in geometric coordinate space 
		$\mathbb{R}_t \times \mathbb{R}_u \times \mathbb{T}^2$,
		which is defined in Def.\,\ref{D:DEVELOPMENTOFDATA},
		depicted in Fig.\,\ref{F:MINTERESTINGDEVELOPMENT},
		and which has the properties revealed by Props.\,\ref{P:PROPERTIESOFMSINGULARANDCREASE} 
		and \ref{P:INTERESTINGREGIONFOLIATEDBYINTERESTINGTIMFEFUNCTION}.
		$\MInteresting$ is contained in 
		$\bigcup_{\muxmulevelsetvalue \in [0,\muxmulevelsetvalue_0]} \twoargMrough{[\timefunction_0,0],[- \rightu,\leftu]}{\muxmulevelsetvalue}$,
		where the $\twoargMrough{[\timefunction_0,0],[- \rightu,\leftu]}{\muxmulevelsetvalue}$ 
		are the developments from Theorem~\ref{T:EXISTENCEUPTOTHESINGULARBOUNDARYATFIXEDKAPPA}.
		Moreover, the singular boundary portion $\mathcal{B}^{[0,\muxmulevelsetvalue_0]}$,
		described below,
		is contained in the top boundary of $\MInteresting$.
	\item The fluid solution wave variables $\wavearray$ (see \eqref{E:ARRAYOFWAVEVARIABLES}), 
			the eikonal function $u$, 
			$\upmu$, $\Lunit^i$,
			and all of the auxiliary quantities constructed out of these quantities
			exist classically
			with respect to the geometric coordinates $(t,u,x^2,x^3)$
			on all of $\MInteresting$,
			\textbf{including the singular boundary portion $\mathcal{B}^{[0,\muxmulevelsetvalue_0]}$ described below}.
			In particular, with respect to the geometric coordinates,
			the fluid variables are
			solutions to equations \eqref{E:INTROTRANSPORTVI}--\eqref{E:INTROBS} \underline{and} the equations of
			Theorem~\ref{T:GEOMETRICWAVETRANSPORTSYSTEM}
			on $\MInteresting$.
	\item The following quantities extend as solutions to the compact set 
	$\MInteresting$
	as elements of the following spacetime H\"{o}lder spaces\footnote{Actually, thanks to the estimates offered by 
	Prop.\,\ref{P:IMPROVEMENTOFAUXILIARYBOOTSTRAP},
	the solution enjoys additional regularity in directions tangent to the characteristics
	compared to what we have stated here; we have stated simpler, sub-optimal regularity
	conclusions only to avoid cluttering the presentation. \label{FN:SUBOPTIMALREGULARITYSTATEMENTS}} 
	with respect to the geometric coordinates,
	and their corresponding spacetime H\"{o}lder norms on $\MInteresting$
	are bounded by $\leq C$:
	\begin{itemize}
		\item $\wavearray, \, \vortrenormalized^i, \, \GradEnt^i, \, \VortVort^i, \, \DivGradEnt 
			\in C_{\textnormal{geo}}^{3,1}(\MInteresting)$
		\item $\Upsilon \in C_{\textnormal{geo}}^{3,1}(\MInteresting)$
		\item $\Lunit^i, \, \upmu \in C_{\textnormal{geo}}^{2,1}(\MInteresting)$
	\end{itemize}
	\item $\MInteresting$ is foliated by the
		level sets of a time function $\newtimefunction$,
		which satisfies the bound 
		$\| \newtimefunction \|_{C_{\textnormal{geo}}^{1,1}(\MInteresting)} \leq C$
		and has the range $[\timefunction_0,0] \eqdef [-\mupositive,0]$
		on $\MInteresting$.
		That is, 
		$\MInteresting = \cup \timefunction_{\in [\timefunction_0,0]}
			\inthyp{\timefunction}{[- \rightu,\leftu]}
		$,
		where for $\timefunction \in [\timefunction_0,0]$,
		$
		\inthyp{\timefunction}{[- \rightu,\leftu]} 
		\eqdef 
		\left\lbrace
			(t,u,x^2,x^3) \ | \ 
			(u,x^2,x^3) \in [- \rightu,\leftu] \times \mathbb{T}^2,
				\,
			\newtimefunction(t,u,x^2,x^3) = \timefunction
			\right\rbrace
		$.
	\item The $L^{\infty}$ estimates of Prop.\,\ref{P:IMPROVEMENTOFAUXILIARYBOOTSTRAP}
				hold on $\MInteresting$
				with $\fundbootsmall$ replaced by $C \initialsmall$.
				Moreover, on each development $\twoargMrough{[\timefunction_0,0],[- \rightu,\leftu]}{\muxmulevelsetvalue}$
				with $\muxmulevelsetvalue \in [0,\muxmulevelsetvalue_0]$,
				the solution enjoys the energy estimates guaranteed by
				Theorem~\ref{T:EXISTENCEUPTOTHESINGULARBOUNDARYATFIXEDKAPPA}.
	\item For $\timefunction \in [\timefunction_0,0]$,
		 we have: 
		\begin{align} \label{E:MAINRESULTSMINVALUEOFMUONFOLIATION}
			\min_{\inthyp{\timefunction}{[- \rightu,\leftu]}} \upmu
			& = - \timefunction.
		\end{align}
		Moreover, within $\inthyp{\timefunction}{[- \rightu,\leftu]}$,
		the minimum value of $- \timefunction$
		in \eqref{E:MAINRESULTSMINVALUEOFMUONFOLIATION} is achieved by $\upmu$ precisely on the
		set $\mulevelsettwoarg{- \timefunction}{[0,\muxmulevelsetvalue_0]}
		\eqdef 
		\bigcup_{\muxmulevelsetvalue \in [0,\muxmulevelsetvalue_0]} \twoargmumuxtorus{- \timefunction}{-\muxmulevelsetvalue}$
		from definition~\eqref{E:LEVELSETSOFMUINSINGULARGETION},
		which is a three-dimensional $C_{\textnormal{geo}}^{2,1}$ embedded manifold contained in 
		$\inthyp{\timefunction}{[- \rightu,\leftu]}$
		with $C_{\textnormal{geo}}^{1,1}$ boundary components
		$\twoargmumuxtorus{- \timefunction}{0}$
		and
		$ \twoargmumuxtorus{- \timefunction}{-\muxmulevelsetvalue_0}$.
		In particular, in $\MInteresting$,
		$\upmu$ vanishes precisely along $\mulevelsettwoarg{0}{[0,\muxmulevelsetvalue_0]}$,
		which by Def.\,\ref{D:SINGULARBOUNDARYPORTION} 
		is equal to the singular boundary portion $\mathcal{B}^{[0,\muxmulevelsetvalue_0]}$
		and which is contained in $\inthyp{\timefunction}{[- \rightu,\leftu]}$.
	\item The change of variables map 
		$\InterestingCHOV(t,u,x^2,x^3) 
		= (\newtimefunction,u,x^2,x^3)$ defined in \eqref{E:CHOVFROMGEOTOINTERESTINGCOORDS}
		is a diffeomorphism from $\MInteresting$
		onto its image $[\timefunction_0,0] \times [- \rightu,\leftu] \times \mathbb{T}^2$
		satisfying $\| \InterestingCHOV \|_{C_{\textnormal{geo}}^{1,1}(\MInteresting)} \leq C$.
	\item On $\MInteresting$,
		the change of variables map 
		$\Upsilon(t,u,x^2,x^3) = (t,x^1,x^2,x^3)$ is an injection onto its image in Cartesian coordinate space
		satisfying $\| \Upsilon \|_{C_{\textnormal{geo}}^{3,1}(\MInteresting)} \leq C$.
		In particular, $\Upsilon$ is a homeomorphism from the compact set
		$\MInteresting$ onto its image.
		Moreover, $\Upsilon$ is a diffeomorphism on
		$\MInteresting \backslash \mathcal{B}^{[0,\muxmulevelsetvalue_0]}$.
\end{itemize}		

\medskip

\noindent \underline{\textbf{The geometric coordinate description of the singular boundary}}.
\begin{itemize}
\item The \textbf{singular boundary} portion
$\mathcal{B}^{[0,\muxmulevelsetvalue_0]} = \bigcup_{\muxmulevelsetvalue \in [0,\muxmulevelsetvalue_0]} \twoargmumuxtorus{0}{-\muxmulevelsetvalue}$ 
from Def.\,\ref{D:SINGULARBOUNDARYPORTION}.
is contained in
$ 
\partial \MInteresting
$,
and $\mathcal{B}^{[0,\muxmulevelsetvalue_0]}$ is a $3$-dimensional 
$C^{2,1}$-embedded submanifold-with-boundary of geometric coordinate space.
Its two boundary components are its future boundary $\twoargmumuxtorus{0}{-\muxmulevelsetvalue_0}$
and its past boundary $\partial_- \mathcal{B}^{[0,\muxmulevelsetvalue_0]} = \twoargmumuxtorus{0}{0}$,
which refer to as the \textbf{crease} (see definition \eqref{E:CREASE}).
\item The boundary components $\twoargmumuxtorus{0}{-\muxmulevelsetvalue_0}$ and $\partial_- \mathcal{B}^{[0,\muxmulevelsetvalue_0]}$
	are $C^{1,1}$ embedded $2$-dimensional tori in geometric coordinate space 
	that are spacelike with respect to $\gfour$.
\end{itemize}

\medskip

\noindent \underline{\textbf{The Cartesian coordinate description of the singularity formation in
 $\Upsilon\left(\MInteresting \right)$}}.
\begin{itemize}
	\item On $\MInteresting$,
		the change of variables map 
		$\Upsilon(t,u,x^2,x^3) = (t,x^1,x^2,x^3)$ is an injection onto its image 
		$\Upsilon\left(\MInteresting \right)$ in Cartesian coordinate space 
		$\mathbb{R}_t \times \mathbb{R}_{x^1} \times \mathbb{T}^2$
		verifying $\| \Upsilon \|_{C_{\textnormal{geo}}^{2,1}(\MInteresting)} \leq C$.
		In particular, $\Upsilon$ is a homeomorphism from the compact set
		$\MInteresting$ onto its image.
		Moreover, on $\MInteresting \backslash \mathcal{B}^{[0,\muxmulevelsetvalue_0]}$
		(where $\mathcal{B}^{[0,\muxmulevelsetvalue_0]}$ 
		is the singular boundary portion from Def.\,\ref{D:SINGULARBOUNDARYPORTION}),
		$\Upsilon$ is a diffeomorphism.
 	\item On $\Upsilon\left(\MInteresting \backslash \mathcal{B}^{[0,\muxmulevelsetvalue_0]} \right)$,
		the solution exists classically with respect to the Cartesian coordinates.
\item The following lower bound holds in 
		$\Upsilon\left(\MInteresting \cap \lbrace (t,u,x^2,x^3) \ | \ |u| \leq \interestingu \rbrace\right)$:
		\begin{align} \label{E:MAINTHEOREMBLOWUPLOWERBOUND}
		|X \almostRiemann_{(+)}| 
		&
		\geq
		\frac{\mathring{\updelta}_*}{\upmu |\bar{\Speed}_{;\LogDensity} + 1|},
		\end{align}
		where $\mathring{\updelta}_* > 0$ is the data-parameter from \eqref{E:DELTASTARDEF},
		$\bar{\Speed}_{;\LogDensity} \eqdef \Speed_{;\LogDensity}(\LogDensity = 0,\Ent=0)$
		is $\Speed_{;\LogDensity}$ evaluated at the trivial solution, 
		$\bar{\Speed}_{;\LogDensity} + 1$ is a \underline{non-zero} constant
		by assumption, and the $\Sigma_t$-tangent vectorfield $X$ has Euclidean length satisfying
		$\sqrt{\sum_{a=1}^3 (X^a)^2} = 1 + \mathcal{O}(\mathring{\upalpha})$,
		where $\mathring{\upalpha}$ is the small parameter from Sect.\,\ref{SS:PARAMETERSIZEASSUMPTIONS}.
		In particular, if $q \in \Upsilon\left(\mathcal{B}^{[0,\muxmulevelsetvalue_0]} \right)$,
		then since 
		$\mathcal{B}^{[0,\muxmulevelsetvalue_0]} 
		\subset
		\MInteresting \cap \left\lbrace (t,u,x^2,x^3) \ | \ |u| \leq \frac{\interestingu}{2} \right\rbrace$ 
		by \eqref{E:SINGULARBOUNDARYPORTION} and \eqref{E:IMPROVEDLEVELSETSTRUCTUREANDLOCATIONOFMIN}, 
		and since $\upmu = 0$ along $\Upsilon\left(\mathcal{B}^{[0,\muxmulevelsetvalue_0]} \right)$,
		it follows that
		$|X \almostRiemann_{(+)}|(q') \to \infty$ as $q' \rightarrow q$ in
		$\Upsilon\left(\MInteresting \backslash \mathcal{B}^{[0,\muxmulevelsetvalue_0]} \right)$.
		Similarly, the following lower bounds hold in
		$\Upsilon\left(\MInteresting \cap \lbrace (t,u,x^2,x^3) \ | \ |u| \leq \interestingu \rbrace\right)$, 
		where $\LogDensity$ is the logarithmic density (see \eqref{E:LOGDENS}):
		\begin{align} \label{E:MAINTHEOREMDENSITYANDVELOCITYBLOWUPLOWERBOUNDS}
		|X \LogDensity| 
		&
		\geq
		\frac{\mathring{\updelta}_*}{4 \upmu |\bar{\Speed}_{;\LogDensity} + 1|},
		&
		|X v^1| 
		&
		\geq
		\frac{\mathring{\updelta}_*}{4 \upmu |\bar{\Speed}_{;\LogDensity} + 1|}.
		\end{align}
	\item (Regular behavior\footnote{Here we have only highlighted some of the quantities that remain
	$L^{\infty}$-bounded on
	$\Upsilon\left(\MInteresting \right)$. We refer to 
	Prop.\,\ref{P:IMPROVEMENTOFAUXILIARYBOOTSTRAP} for more comprehensive results. 
	\label{FN:MAINTHEOREMONLYHIGHLIGHTEDSOMELINFINITYBOUNDEDQUANTITIES}} along the characteristics).
	The derivatives of
	$\wavearray$,
	$\vortrenormalized^i$,
	$\GradEnt^i$
	up to order $\Ntop - 11$
	with respect to the vectorfields in the 
	$\nullhyparg{u}$-tangent commutation set $\Tanset$ defined in \eqref{E:COMMUTATIONVECTORFIELDS}
	and the derivatives of
	$\VortVort^i$
	and
	$\DivGradEnt$
	up to order $\Ntop - 12$
	with respect to the elements of $\Tanset$ are $L^{\infty}$-bounded
	on $\Upsilon\left(\MInteresting \right)$.
	Finally, for $\alpha = 0,1,2,3$ and $A=2,3$, 
	the derivatives of 
	$\gfour_{ab} \Yvf{A}^a \partial_{\alpha} v^b$
	up to order $\Ntop - 12$
	with respect to the elements of $\Tanset$
	are $L^{\infty}$-bounded on
	$\Upsilon\left(\MInteresting \right)$.
	\end{itemize}

\medskip

\noindent \underline{\textbf{The Cartesian coordinate description of the singular boundary and the crease}}.
\begin{itemize}	
	\item $\Upsilon(\mathcal{B}^{[0,\muxmulevelsetvalue_0]})$
		is a $C^{1,1/2}$ embedded $\gfour$-null hypersurface in Cartesian coordinate space
		that is foliated by the integral curves of
		$\Lunit$; see Prop.\,\ref{P:DESCRIPTIONOFSINGULARBOUNDARYINCARTESIANSPACE}
		for a more detailed description.
	\item For $\muxmulevelsetvalue \in [0,\muxmulevelsetvalue_0]$,
		 $\Upsilon(\twoargmumuxtorus{0}{-\muxmulevelsetvalue}) \subset \Upsilon(\mathcal{B}^{[0,\muxmulevelsetvalue_0]})$ 
		is an embedded $C^{1,1}$ graph over 
		$\mathbb{T}^2$ in Cartesian coordinate space such that at each $q \in \twoargmumuxtorus{0}{-\muxmulevelsetvalue}$,
		$\Upsilon(\twoargmumuxtorus{0}{-\muxmulevelsetvalue})$
		is spacelike with respect to $\gfour|_{\Upsilon(q)}$ at the point $\Upsilon(q)$.
	\item In particular, considering the case $\muxmulevelsetvalue = 0$ in the previous point,
		we see that the image of the crease under $\Upsilon$, namely $\Upsilon(\twoargmumuxtorus{0}{0})$,
		is an embedded $C^{1,1}$ graph over $\mathbb{T}^2$ in Cartesian coordinate space
		that is $\gfour$-spacelike at each of its points.
\end{itemize}

\end{theorem}

\begin{proof}
		\hfill
		
		\medskip
		\noindent \textbf{Proof that $\MInteresting \subset 
		\bigcup_{\muxmulevelsetvalue \in [0,\muxmulevelsetvalue_0]} \twoargMrough{[\timefunction_0,0],[- \rightu,\leftu]}{\muxmulevelsetvalue}$}:
		This result
		follows from definition \eqref{E:INTERESTINGDEVELOPMENTOFDATA}
		and the fact that 
		$\datahypfortimefunctiontwoarg{-\muxmulevelsetvalue}{[\timefunction_0,0]}
		\subset 
		\twoargMrough{[\timefunction_0,0],[- \rightu,\leftu]}{\muxmulevelsetvalue}
		$,
		which in turn follows from the decomposition
		$
		\twoargMrough{[\timefunction_0,0],[- \rightu,\leftu]}{\muxmulevelsetvalue}
		= 
		\bigcup_{\timefunction \in [\timefunction_0,0]} \hypthreearg{\timefunction}{[- \rightu,\leftu]}{\muxmulevelsetvalue}
		$,
		the fact that $\twoargmumuxtorus{-\timefunction}{-\muxmulevelsetvalue} \subset \hypthreearg{\timefunction}{[- \rightu,\leftu]}{\muxmulevelsetvalue}$,
		and \eqref{E:CLOSEDIMPROVEMENTLEVELSETSTRUCTUREOFMUXEQUALSMINUSKAPPA} with $\timefunctionboot, \upmuboot = 0$.
		In the remainder of the proof, we will silently use
		the fact that $\MInteresting \subset 
		\bigcup_{\muxmulevelsetvalue \in [0,\muxmulevelsetvalue_0]} \twoargMrough{[\timefunction_0,0],[- \rightu,\leftu]}{\muxmulevelsetvalue}$
		and the following consequence of
		Theorem~\ref{T:EXISTENCEUPTOTHESINGULARBOUNDARYATFIXEDKAPPA}:
		at fixed $\muxmulevelsetvalue \in [0,\muxmulevelsetvalue_0]$, all results proved in the paper 
		prior to Theorem~\ref{T:EXISTENCEUPTOTHESINGULARBOUNDARYATFIXEDKAPPA}
		hold with $\timefunctionboot = 0$.
		
		\medskip
		\noindent \textbf{Proof of classical existence with respect to the geometric coordinates on $\MInteresting$}:
		Theorem~\ref{T:EXISTENCEUPTOTHESINGULARBOUNDARYATFIXEDKAPPA} yields 
		that for fixed $\muxmulevelsetvalue \in [0,\muxmulevelsetvalue_0]$,
		the solution exists classically with the respect to the geometric coordinates on
		$\twoargMrough{[\timefunction_0,0],[- \rightu,\leftu]}{\muxmulevelsetvalue}$.
		Since $\MInteresting \subset 
		\bigcup_{\muxmulevelsetvalue \in [0,\muxmulevelsetvalue_0]} \twoargMrough{[\timefunction_0,0],[- \rightu,\leftu]}{\muxmulevelsetvalue}$,
		we immediately conclude classical existence with respect to the geometric coordinates on $\MInteresting$.
		
		\medskip
		\noindent \textbf{Proof of the H\"{o}lder bounds}:
		To derive the H\"{o}lder bounds 
		$\| \wavearray \|_{C_{\textnormal{geo}}^{3,1}(\MInteresting)} \leq C$,
		$\| \vortrenormalized^i \|_{C_{\textnormal{geo}}^{3,1}(\MInteresting)} \leq C$,
		etc., 
		we first use Lemma~\ref{L:CONTINUOUSEXTNESION} with $\timefunctionboot = 0$
		(i.e., we use the lemma on $\twoargMrough{[\timefunction_0,0],[- \rightu,\leftu]}{\muxmulevelsetvalue}$),
		the fact that
		$\MInteresting 
		\subset 
		\bigcup_{\muxmulevelsetvalue \in [0,\muxmulevelsetvalue_0]} \twoargMrough{[\timefunction_0,0],[- \rightu,\leftu]}{\muxmulevelsetvalue}$,
		and Rademacher's theorem to deduce that
		$\| \wavearray \|_{W_{\textnormal{geo}}^{4,\infty}(\MInteresting)} \leq C$,
		$\| \vortrenormalized^i \|_{W_{\textnormal{geo}}^{4,\infty}(\MInteresting)} \leq C$,
		etc. From these bounds and the Sobolev embedding result \eqref{E:SOBOELVEMBEDDINGINTERESTINGREGIONRELYINGONQUASICONVEXITY},
		we conclude the desired H\"{o}lder bounds.
		
		\medskip
		\noindent \textbf{Proof of the properties of the time function $\newtimefunction$ and the map
		$\InterestingCHOV(t,u,x^2,x^3)$}:
		We derived these results in Prop.\,\ref{P:INTERESTINGREGIONFOLIATEDBYINTERESTINGTIMFEFUNCTION}.
		
		\medskip
		\noindent \textbf{Proof of the properties of $\Upsilon$ and 
		classical existence with respect to the Cartesian coordinates on 
		the domain $\Upsilon\left(\MInteresting \backslash \mathcal{B}^{[0,\muxmulevelsetvalue_0]} \right)$}:
		We derived the properties of $\Upsilon$ 
		in Prop.\,\ref{P:CHOVGEOMETRICTOCARTESIANISINJECTIVEONREGIONWECAREABOUT}.
		Since the proposition in particular shows that $\Upsilon$ is a global diffeomorphism 
		$\MInteresting \backslash \mathcal{B}^{[0,\muxmulevelsetvalue_0]}$,
		and since we have already shown classical existence with respect to the geometric coordinates on
		$\MInteresting$, we conclude that the solution exists classically 
		with respect to the Cartesian coordinates on 
		the domain $\Upsilon\left(\MInteresting \backslash \mathcal{B}^{[0,\muxmulevelsetvalue_0]} \right)$
		in Cartesian coordinate space.
		
		\medskip
		\noindent \textbf{Proof of the regular behavior along the characteristics}:
		These results follow from the last item stated in Theorem~\ref{T:EXISTENCEUPTOTHESINGULARBOUNDARYATFIXEDKAPPA}.
		
		\medskip
		\noindent \textbf{Proof of the properties of $\mathcal{B}^{[0,\muxmulevelsetvalue_0]}$, $\twoargmumuxtorus{0}{-\muxmulevelsetvalue_0}$,
		and $\partial_- \mathcal{B}^{[0,\muxmulevelsetvalue_0]} = \twoargmumuxtorus{0}{0}$}:
		We derived these results in
		in Prop.\,\ref{P:PROPERTIESOFMSINGULARANDCREASE}.
		
		\medskip
		\noindent \textbf{Proof of \eqref{E:MAINRESULTSMINVALUEOFMUONFOLIATION} and related properties of $\upmu$}:
		In Prop.\,\ref{P:INTERESTINGREGIONFOLIATEDBYINTERESTINGTIMFEFUNCTION}, we proved
		\eqref{E:MAINRESULTSMINVALUEOFMUONFOLIATION} and showed that 
		 within $\inthyp{\timefunction}{[- \rightu,\leftu]}$,
		the minimum value of $- \timefunction$
		in \eqref{E:MAINRESULTSMINVALUEOFMUONFOLIATION} is achieved by $\upmu$ precisely on the
		set $\mulevelsettwoarg{- \timefunction}{[0,\muxmulevelsetvalue_0]}
		\eqdef 
		\bigcup_{\muxmulevelsetvalue \in [0,\muxmulevelsetvalue_0]} \twoargmumuxtorus{- \timefunction}{-\muxmulevelsetvalue}$.
		
		\medskip
		\noindent \textbf{Proof of the lower bounds \eqref{E:MAINTHEOREMBLOWUPLOWERBOUND} and 
		\eqref{E:MAINTHEOREMDENSITYANDVELOCITYBLOWUPLOWERBOUNDS}}:
		These estimates follow from 
		the estimates
		\eqref{E:FIXEDKAPPAMAINTHEOREMBLOWUPLOWERBOUND}
		and
		\eqref{E:FIXEDKAPPAMAINTHEOREMDENSITYANDVELOCITYBLOWUPLOWERBOUNDS},
		which hold in 
		$
		\bigcup_{\muxmulevelsetvalue \in [0,\muxmulevelsetvalue_0]}
		\Upsilon\left(\twoargMrough{[\timefunction_0,0],[- \interestingu,\interestingu]}{\muxmulevelsetvalue} \right)
		$,
		and the fact that
		$\Upsilon\left(\MInteresting \cap \lbrace (t,u,x^2,x^3) \ | \ |u| \leq \interestingu \rbrace\right)
		\subset
		\bigcup_{\muxmulevelsetvalue \in [0,\muxmulevelsetvalue_0]}
		\Upsilon\left(\twoargMrough{[\timefunction_0,0],[- \interestingu,\interestingu]}{\muxmulevelsetvalue} \right)
		$
		(since 
		 $\MInteresting \subset 
		\bigcup_{\muxmulevelsetvalue \in [0,\muxmulevelsetvalue_0]} \twoargMrough{[\timefunction_0,0],[- \rightu,\leftu]}{\muxmulevelsetvalue}$).
		
		\medskip
		\noindent \textbf{Proof of the Cartesian coordinate description of the singular boundary and the crease}:
		We derived these results in Prop.\,\ref{P:DESCRIPTIONOFSINGULARBOUNDARYINCARTESIANSPACE}.
		
	\end{proof}

\appendix

 \section{Simple isentropic plane-symmetric solutions} \label{A:PS}
In this appendix, we show that there exists a large family of
shock-forming simple isentropic plane-symmetric solutions whose 
data satisfy the assumptions stated in Sect.\,\ref{S:ASSUMPTIONSONTHEDATA}. 
As we will see, simple isentropic plane-wave solutions are characterized by
$\initialsmall =0$, where $\initialsmall$ is the data-size parameter
featured in the assumptions of Sect.\,\ref{S:ASSUMPTIONSONTHEDATA}.
In Appendix~\ref{A:OPENSETOFDATAEXISTS}, we combine the 
results of this appendix with Cauchy stability arguments
to show that there exists an open set of data satisfying the assumptions of Sect.\,\ref{S:ASSUMPTIONSONTHEDATA}. 
By ``open set,'' we mean open relative to the topologies corresponding to the norm
$\mathring{\Delta}_{\Sigma_0^{[-\farrightu,\leftu]}}^{\Ntop+1}$ on $\Sigma_0$ defined in
\eqref{E:PERTURBATIONSMALLNESSINCARTESIANDIFFERENTIALSTRUCTURE}, where $\Ntop \geq 24$.

\subsection{Plane symmetric and simple plane-symmetric solutions} 
\label{SS:PSCONSTRUCTIONS} 
We begin with a quick presentation of the isentropic compressible Euler equations in plane-symmetry. 
By \emph{isentropic plane-symmetric solutions}, 
we mean those solutions with the following properties: $\LogDensity$ and $v^1$ are functions 
of only $t$ and $x^1$, ``the symmetry breaking fluid variables'' satisfy $v^2, v^3 \equiv 0$, and $\Ent$ is constant.
For convenience, we will assume that $\Ent \equiv 0$, though that is not essential for our main results,
i.e., we could have easily handled solutions with $\Ent \equiv \Ent_0$, where $\Ent_0$ is a constant.
It is straightforward to see that for such solutions, the fluid vorticity $\Flatcurl v$ also vanishes.
In Appendix~\ref{A:OPENSETOFDATAEXISTS}, we will view our isentropic plane-symmetric solutions as ``background'' solutions
in three spatial dimensions with trivial dependence on the $(x^2,x^3)$ coordinates. However, in this appendix,
to shorten the presentation, we will completely suppress the variables $(x^2,x^3)$ and instead view the 
plane-symmetric solutions as solutions in $1+1$ dimensions.

Throughout this appendix, we adorn symbols related to the background plane-symmetric solutions by a ``$PS$.'' 
For example, 
we denote the logarithmic density by $\LogDensityPS$, and we set $\velocityPS \eqdef v^1$. 
Although the eikonal function also depends on the background solution through the dependence
of the coefficients of the eikonal equation \eqref{E:EIKONALEQUATION} on the fluid,
we will continue to denote it by ``$u$'' instead of ``$u^{\text{PS}}$.''
This is consistent with the point of view we take in Appendix~\ref{A:OPENSETOFDATAEXISTS}
during our discussion of Cauchy stability, 
where for convenience, 
we will ``fix'' geometric coordinate space $\mathbb{R}_t \times \mathbb{R}_u \times \mathbb{T}_{x^2,x^3}^2$
and consider families of solutions that exist with respect to the geometric coordinates on
a common domain. Note that, although for general solutions, the maps $(t,u,x^2,x^3) \rightarrow x^i(t,u)$ from geometric coordinates to 
a Cartesian spatial coordinate depend on the fluid solution, to avoid clutter, in this appendix,
we will not adorn $x^1$ with a ``$PS$'' subscript.
However, we do adorn the corresponding change of variables map: $\PSUpsilon(t,u) \eqdef (t,x^1)$.

\subsubsection{The Riemann invariants in isentropic plane-symmetry}
\label{SSS:RIEMANNINVARIANTSISENTROPICPLANESYMMETRY} 
In plane-symmetry with $\Ent \equiv 0$, 
we define the Riemann invariants to be the following
functions of $\velocityPS$ and $\LogDensityPS$:
\begin{align}\label{AE:RIEMANNINVARIANTS} 
\RRiemannPS 
& \eqdef 
	\velocityPS 
	+ 
	\FPS(\LogDensityPS), 
& \LRiemannPS 
& \eqdef 
	\velocityPS
	- 
	\FPS(\LogDensityPS),
\end{align}
where $\FPS = \FPS(z)$ is defined to be the solution to the following ODE initial value problem:
\begin{align} \label{AE:FFORRIEMANNINVARIANT}
\FPS(0) & = 0, 
& \frac{\mathrm{d}}{\mathrm{d} z} \FPS(z)
& = \SpeedPS(z),
\end{align}
where we recall that by our assumption \eqref{E:BACKGROUNDSOUNDSPEEDISUNITY}, we have $\SpeedPS(0) = 1$.
On RHS~\eqref{AE:FFORRIEMANNINVARIANT} and throughout,
$\SpeedPS(z) := \Speed(z,0)$, where 
$\Speed(z,0)$ is the speed of sound (see \eqref{E:SOUNDSPEED}), 
viewed as a function of the logarithmic density $z$
and the entropy $\Ent$, evaluated at $\Ent = 0$.
We note that $\velocityPS$ and $\LogDensityPS$
can respectively be expressed in terms of the Riemann invariants as follows:
\begin{align}\label{AE:RIEMANNTOWAVE}
\velocityPS 
& = 
\frac{1}{2} \left(\RRiemannPS + \LRiemannPS \right), 
&
\LogDensityPS & = (\FPS)^{-1} \circ \left\lbrace \frac{1}{2} \left(\RRiemannPS - \LRiemannPS \right)\right\rbrace,
\end{align}
where $(\FPS)^{-1}$ is the inverse function of $\FPS$.
Note that by \eqref{AE:FFORRIEMANNINVARIANT} and \eqref{E:BACKGROUNDSOUNDSPEEDISUNITY}, 
$(\FPS)^{-1}$
is well-defined and smooth in a neighborhood of $0$. 
Note that the Riemann invariants \eqref{AE:RIEMANNINVARIANTS} 
agree with the almost Riemann invariants defined in \eqref{E:ALMOSTRIEMANNINVARIANTS}.
Much like in the bulk of the paper, when we are deriving estimates for the fluid,
it is understood that all fluid variables are to be viewed as functions of the Riemann invariants
via \eqref{AE:RIEMANNTOWAVE}.

\subsubsection{The compressible Euler equations in terms of the Riemann invariants in isentropic plane-symmetry}
\label{SSS:COPMRESSIBLEEULERINTERMSOFRIEMANNINVARIANTS}
Due to the isentropic plane-symmetry, the Riemann invariants are in fact invariant along the characteristics.
More precisely, it is straightforward to verify
that for smooth isentropic plane-symmetric solutions,  
the compressible Euler equations
\eqref{E:BVIEVOLUTION}--\eqref{E:BENTROPYEVOLUTION}
are equivalent to the following $2 \times 2$ system
of quasilinear transport equations:
\begin{align} \label{AE:COMPRESSIBLEEULEREQUATIONSFORRIEMANNINVARIANTS}
\LunitPS \RRiemannPS & = 0,
&
\uLunitPS \LRiemannPS & = 0, 
\end{align}
where: 
\begin{align} \label{AE:NULLVECTORFIELDSPLANESYMMETRY}
\LunitPS 
& \eqdef \partial_t + (\velocityPS + \SpeedPS) \partial_{1},
& \uLunitPS & \eqdef \partial_t 
+ 
(\velocityPS - \SpeedPS) \partial_{1}.
\end{align} 
We denote the initial data for
\eqref{AE:COMPRESSIBLEEULEREQUATIONSFORRIEMANNINVARIANTS}
by:
\begin{align} \label{AE:DATACOMPRESSIBLEEULEREQUATIONSFORRIEMANNINVARIANTS}
	\RRiemannPS|_{\Sigma_0} & = \dataRRiemannPS,
	&
	\LRiemannPS|_{\Sigma_0} & = \dataLRiemannPS.
\end{align}
In the rest of the appendix, we will consider only
\emph{simple} isentropic plane-symmetric solutions, defined to be solutions with $\LRiemannPS \equiv 0$.
From \eqref{AE:COMPRESSIBLEEULEREQUATIONSFORRIEMANNINVARIANTS}, it follows that such solutions
arise from initial data with $\dataLRiemannPS \equiv 0$.

A consequence of standard methods going back
Riemann's famous work \cite{bR1860} is that any compactly supported initial datum $\dataRRiemannPS$ launches 
a shock-forming simple isentropic plane-symmetric solution to \eqref{AE:COMPRESSIBLEEULEREQUATIONSFORRIEMANNINVARIANTS}. In the forthcoming subsections, we will construct ``admissible" (see Def.\,\ref{AD:ADMISSIBLEBACKGROUND}) 
compactly supported data along $\Sigma_0$ that launch solutions
whose perturbations we study in Appendix~\ref{A:OPENSETOFDATAEXISTS}. 
To this end, we find it useful to revisit our construction of the acoustical geometry
so that we can capitalize on the many simplifications that occur in simple isentropic plane-symmetry.

\subsection{The acoustic geometry and explicit solution formulas in simple isentropic plane-symmetry}
\label{SS:ACOUSTICGEOMETRYINISENTROPICPLANESYMMETRY} 

\subsubsection{Definitions and identities}
We begin by defining the eikonal function in plane-symmetry to be the solution to the following transport equation initial value problem:
\begin{subequations}
\begin{align}
\LunitPS u & = 0, \label{AE:EIKONALEQUATIONINPLANESYMMETRY} 
	\\
u |_{\Sigma_0} & = - x^1. \label{AE:DATAFOREIKONALEQUATIONINPLANESYMMETRY} 
\end{align}
\end{subequations}
Note that the initial condition stated in \eqref{AE:DATAFOREIKONALEQUATIONINPLANESYMMETRY} is the same as the one
\eqref{E:EIKONALEQUATION} we assumed in the bulk of the paper.
It is straightforward to check (cf.\ \eqref{E:ACOUSTICALMETRICINDOUBLENULLFRAME})
that in plane-symmetry, 
$
(\gfour^{-1})^{\alpha \beta} \partial_{\alpha} u \partial_{\beta} u
= 
- (\LunitPS u) \uLunitPS u
$,
and that
$u$ solves
\eqref{AE:EIKONALEQUATIONINPLANESYMMETRY}--\eqref{AE:DATAFOREIKONALEQUATIONINPLANESYMMETRY}
if and only if it solves \eqref{E:EIKONALEQUATION}. In particular, in plane-symmetry,
the (fully nonlinear) eikonal equation \eqref{E:EIKONALEQUATION} 
is equivalent to \eqref{AE:EIKONALEQUATIONINPLANESYMMETRY}, 
which is linear in $u$ because the operator $\LunitPS$ can be defined by
equation \eqref{AE:NULLVECTORFIELDSPLANESYMMETRY}, which does not depend on $u$.

We now define the inverse foliation density as follows:
\begin{align} \label{AE:INVERSEFOLIATIONDENSITYINPLANESYMMETRY}
	\muPS 
	& \eqdef 
	- 
	\frac{1}{\SpeedPS \partial_{1} u}. 
\end{align}
One can check that in isentropic plane-symmetry, the quantity ``$\muPS$'' defined by \eqref{AE:INVERSEFOLIATIONDENSITYINPLANESYMMETRY}
is equal to the quantity ``$\upmu$'' defined in \eqref{E:MUDEF}.

Straightforward calculations yield that in isentropic plane-symmetry, 
with $\LunitPS$ as in \eqref{AE:NULLVECTORFIELDSPLANESYMMETRY},
we have:
\begin{align} \label{AE:IDENTITYNULLCVECTORFIELDSINPLANESYMMETRY}
	\uLunitPS 
	& = \LunitPS + 2 \XPS,
	&
	\XPS 
	& = - \SpeedPS \p_1,
	&
	\muXPS 
	& \eqdef \muPS \XPS,
\end{align}
where, given our construction of $u$ in 
\eqref{AE:EIKONALEQUATIONINPLANESYMMETRY}--\eqref{AE:DATAFOREIKONALEQUATIONINPLANESYMMETRY},
$\LunitPS$, $\XPS$, $\muXPS$ coincide with the vectorfields defined 
in Def.\,\ref{D:COMVECTORFIELDS},
while $\uLunitPS$ coincides with the vectorfield defined in Def.\,\ref{D:ULUNIT}. 

Next, using \eqref{AE:DATAFOREIKONALEQUATIONINPLANESYMMETRY} and definition~\ref{AE:INVERSEFOLIATIONDENSITYINPLANESYMMETRY},
we compute that the following identities hold on $\Sigma_0$:
\begin{align}\label{E:PSMUINITIAL}
\left. \muPS \right|_{\Sigma_0} 
& = 
\frac{1}{\SpeedPS|_{\Sigma_0}}\qquad \Longrightarrow \qquad \left. (\SpeedPS\muPS) \right|_{\Sigma_0} = 1.
\end{align}

As in the bulk of the paper, we define $(t,u)$ to be the \emph{geometric coordinates}, 
and we denote the corresponding geometric coordinate vectorfields
as $\left\lbrace\geop{t},\geop{u}\right\rbrace$.
In plane-symmetry, $\LunitPS x^2 = \LunitPS x^3 = \XPS x^2 = \XPS x^2 = 0$ and thus, by Lemma~\ref{L:COMMUTATORSTOCOORDINATES}, we have:
\begin{align} \label{AE:PLANESYMMETRYCOMMUTATORVECTORFIELDSAREGEOCOORDINATEVECTORFIELDS}
	\LunitPS 
	& = \geop{t}, 
	& 
	\muXPS & = \geop{u}.
\end{align}
It follows that:
\begin{align} \label{AE:LANDMUXCOMMUTEINPLANESYMMETRY}
[\LunitPS, \muXPS] 
& = 0.
\end{align} 
We will often silently use 
\eqref{AE:PLANESYMMETRYCOMMUTATORVECTORFIELDSAREGEOCOORDINATEVECTORFIELDS}--\eqref{AE:LANDMUXCOMMUTEINPLANESYMMETRY}
in the rest of this appendix.

\subsubsection{Explicit solution formula in geometric coordinates}
\label{SSS:PSEXPLICITFLUIDSOLUTION}
From \eqref{AE:PLANESYMMETRYCOMMUTATORVECTORFIELDSAREGEOCOORDINATEVECTORFIELDS}, we see 
that the transport equation \eqref{AE:COMPRESSIBLEEULEREQUATIONSFORRIEMANNINVARIANTS} for
$\RRiemannPS$ takes the following form in geometric coordinates:
\begin{align} \label{AE:RRIEMANNTRIVIALTRANSPORTEQUATIONINGEOMETRICCOORDINATES}
	\geop{t} \RRiemannPS(t,u) 
	& = 0.
\end{align}
From \eqref{AE:RRIEMANNTRIVIALTRANSPORTEQUATIONINGEOMETRICCOORDINATES} and \eqref{AE:DATACOMPRESSIBLEEULEREQUATIONSFORRIEMANNINVARIANTS},
we see (recalling that $\dataLRiemannPS \equiv 0$ by assumption)
that in geometric coordinates, the solution to \eqref{AE:RRIEMANNTRIVIALTRANSPORTEQUATIONINGEOMETRICCOORDINATES} is:
\begin{align} \label{AE:SOLUIONRRIEMANNTRIVIALTRANSPORTEQUATIONINGEOMETRICCOORDINATES}
	\RRiemannPS(t,u) 
	& = \dataRRiemannPS(u).
\end{align}

\subsubsection{The evolution equation for $\muPS$}
The following lemma provides the evolution equation for $\muPS$. 

\begin{lemma}[Transport equation satisfied by $\muPS$]
\label{AL:PSTRANSPORTFORMU}
For simple isentropic plane-symmetric solutions, 
the inverse foliation density satisfies the following transport equation,
where 
$
\SpeedPSdot = \SpeedPSdot(\LogDensity) \eqdef \frac{d}{d \LogDensity} \SpeedPS(\LogDensity)
$
is the derivative of the speed of sound with respect to the logarithmic density:
\begin{align} \label{AE:LMUPS}
 \LunitPS (\SpeedPS \muPS)
	& = \PSLmusourcetermfunction,
		\\
\PSLmusourcetermfunction 
& \eqdef 
- 
\frac{1}{2} 
\left\lbrace\frac{\SpeedPSdot}{\SpeedPS} 
	+
	1 
\right\rbrace 
\muXPS \RRiemannPS.
\label{AE:PSLMUSOURCETERMFUNCTION}
\end{align}
Moreover, 
\begin{align}\label{AE:LLUPS}
\LunitPS \LunitPS \muPS 
& = 
\LunitPS \LunitPS (\SpeedPS \muPS)
= 0.
\end{align}

Finally, we have the following identities: 
\begin{subequations}
\begin{align} \label{AE:USEFULIDENTITYFORPSLMUSOURCETERMFUNCTION}
	\PSLmusourcetermfunction
	& = 
			\muXPS
			\antiderivativePSLmusourcetermfunction,
				\\
\antiderivativePSLmusourcetermfunction
=
\antiderivativePSLmusourcetermfunction[\RRiemannPS]
& \eqdef 
			- \FPSdot
			\circ 
			(\FPS)^{-1}
			\circ
			\left(\frac{1}{2} \RRiemannPS \right)
			-
			\frac{1}{2} \RRiemannPS
			+ 
			\backgroundSpeedPS,
\label{AE:EXPLICITFORMUSEFULIDENTITYFORPSLMUSOURCETERMFUNCTION}
	\\
\frac{d}{d \RRiemannPS}
\antiderivativePSLmusourcetermfunction[\RRiemannPS]
& = 
- 
\frac{1}{2}
\frac{\FPStwodots}{\FPSdot} 
\circ 
(\FPS)^{-1}
\circ
\left(\frac{1}{2} \RRiemannPS \right)
-
\frac{1}{2},
\label{AE:DERIVATIVEOFUSEFULIDENTITYFORPSLMUSOURCETERMFUNCTION}
\end{align}
\end{subequations}
where $\FPSdot = \FPSdot(\LogDensity) \eqdef \frac{d}{d \LogDensity} \FPS(\LogDensity) = \Speed(\LogDensity)$,
$\backgroundSpeedPS \eqdef \FPSdot(\LogDensity)|_{\LogDensity = 0}$,
$(\FPS)^{-1}$ is the inverse of the map $\LogDensity \rightarrow \FPS(\LogDensity)$,
$\FPStwodots = \FPStwodots(\LogDensity) = \SpeedPSdot(\LogDensity) \eqdef \frac{d}{d \LogDensity} \Speed(\LogDensity)$,
$
\antiderivativePSLmusourcetermfunction'[\RRiemannPS]
\eqdef
\frac{d}{d \RRiemannPS}
\antiderivativePSLmusourcetermfunction[\RRiemannPS]
$,
and $\circ$ denotes the composition of functions.
\end{lemma}

\begin{remark}
	\label{AR:NORMALIZATIONFORANTIDERIVATIVEPSLMUSOURCETERM}
	Note that $\antiderivativePSLmusourcetermfunction[0] = 0$.
	We will use this basic fact in Sect.\,\ref{SS:EXAMPLESOFBONAFIDEDATALEADINGTOADMISSIBLE}.
\end{remark}

\begin{proof}[Proof of Lemma~\ref{AL:PSTRANSPORTFORMU}]
The transport equation 
\eqref{AE:LMUPS} follows from \eqref{E:MUTRANSPORT}, \eqref{E:IDENTITYFORMAINTERMDRIVINGTHESHOCK},
and our assumption that the solution is isentropic, simple, and plane-symmetric,
which in particular implies that $ \LunitPS \SpeedPS = 0$. 
Equation \eqref{AE:LLUPS} follows from differentiating 
\eqref{AE:LMUPS} with $\LunitPS$ and using \eqref{AE:LANDMUXCOMMUTEINPLANESYMMETRY} along with $\LunitPS \RRiemannPS = 0$.
\eqref{AE:USEFULIDENTITYFORPSLMUSOURCETERMFUNCTION}--\eqref{AE:DERIVATIVEOFUSEFULIDENTITYFORPSLMUSOURCETERMFUNCTION}
follow from the chain rule,
\eqref{AE:FFORRIEMANNINVARIANT},
and \eqref{AE:RIEMANNTOWAVE}.
\end{proof}

\subsubsection{Recalling the non-degeneracy condition}
As in \eqref{E:NONDEGENCONDITION}, we will assume the following 
non-degeneracy condition:
$
\frac{1}{\backgroundSpeedPS}
\left\lbrace
	\frac{\backgroundSpeedprimePS}{\backgroundSpeedPS}
	+ 
	1 
\right\rbrace
\neq 0
$,
where the factor in braces on the LHS of this relation
denotes the factor in braces on RHS~\eqref{AE:PSLMUSOURCETERMFUNCTION}
evaluated at the trivial solution $\LogDensityPS \equiv 0$
(which, in the present context, is equivalent to $\RRiemannPS \equiv 0$).
In view of our normalization assumption \eqref{E:BACKGROUNDSOUNDSPEEDISUNITY}
and \eqref{AE:DERIVATIVEOFUSEFULIDENTITYFORPSLMUSOURCETERMFUNCTION}, it follows that
the non-degeneracy condition is equivalent to:
\begin{align}	\label{AE:PLANESYMMETRYNONDEGENCONDITION}
	-
	2
	\frac{d}{d \RRiemannPS}
	\antiderivativePSLmusourcetermfunction[\RRiemannPS]|_{\RRiemannPS = 0}
	&
	=
	\frac{\backgroundSpeedprimePS}{\backgroundSpeedPS}
	+ 
	1 
	\neq 0.
\end{align}
Note that \eqref{AE:PLANESYMMETRYNONDEGENCONDITION}
is equivalent to RHS~\eqref{AE:DERIVATIVEOFUSEFULIDENTITYFORPSLMUSOURCETERMFUNCTION} being non-zero
when it is evaluated at $\RRiemannPS = 0$, and that it implies the
invertibility of the map
$\RRiemannPS \rightarrow \antiderivativePSLmusourcetermfunction[\RRiemannPS]$ 
in a neighborhood of the origin (i.e., near $\RRiemannPS = 0$).
As we discussed in Sect.\,\ref{SS:THEFACTORDRIVINGTHESHOCKFORMATION},
for any equation of state except for that of the Chaplygin gas,
there are always background densities $\overline{\varrho} > 0$
such that \eqref{AE:PLANESYMMETRYNONDEGENCONDITION} holds.

Next, we note that in simple isentropic plane-symmetry, the 
quantity $\PSLmusourcetermfunction$ defined in \eqref{AE:PSLMUSOURCETERMFUNCTION}
can be viewed as a function of $\RRiemannPS$ and $\muXPS \RRiemannPS$.
To emphasize this point of view, we use the notation 
$\PSLmusourcetermfunction[\RRiemannPS,\muXPS \RRiemannPS]$.
Moreover, by \eqref{AE:SOLUIONRRIEMANNTRIVIALTRANSPORTEQUATIONINGEOMETRICCOORDINATES},
in geometric coordinates, $\PSLmusourcetermfunction[\RRiemannPS,\muXPS \RRiemannPS]$ is a function of 
only $u$. To simplify the presentation, we will emphasize this point of view with the shorthand notation
$\PSLmusourcetermfunction(u)$, i.e.,
\begin{align} \label{AE:SHORTHANDNOTATIONFORMAINLMUSOURCETERM}
	\PSLmusourcetermfunction[\RRiemannPS,\muXPS \RRiemannPS](t,u)
	& 
	\eqdef
	\PSLmusourcetermfunction[\RRiemannPS(t,u),\muXPS \RRiemannPS(t,u)]
	=
	\PSLmusourcetermfunction[\dataRRiemannPS(u),\geop{u} \dataRRiemannPS(u)]
	\eqdef 
	\PSLmusourcetermfunction(u).
\end{align}

In a similar vein, $\SpeedPS = \FPSdot$ can be viewed as a function of $\RRiemannPS$,
where $\RRiemannPS$ is a function of $u$ alone, and we will use 
the following shorthand notation:
\begin{align} \label{AE:SHORTHANDNOTATIONFORSPEEDOFSOUND}
	\SpeedPS[\RRiemannPS](t,u)
	& 
	\eqdef
	\FPSdot[\RRiemannPS(t,u)]
	=
	\FPSdot[\dataRRiemannPS(u)]
	\eqdef 
	\SpeedPS(u).
\end{align}

\subsubsection{Explicit expressions for various solution variables}
\label{SSS:EXPLICITEXPRESSIONSFORVARIOUSSOLUTIONVARIABLES}
The following corollary is a straightforward consequence of the prior discussion in this appendix,
and we therefore omit the simple proof.

\begin{corollary}[Explicit expressions for $\muPS$, $\LunitPS\muPS$, and $\partial_1 \RRiemannPS$] 
\label{AC:PSEXPLICITEXPRESSIONSFORSOLUTION}
For simple isentropic plane-symmetric solutions, the following identities hold
relative to the geometric coordinates:
\begin{subequations}
\begin{align}
\LunitPS 
\left\lbrace 
	\SpeedPS(u) \muPS(t,u) 
\right\rbrace
& = \PSLmusourcetermfunction(u)
	=
	\frac{d}{du}
	\antiderivativePSLmusourcetermfunction[\dataRRiemannPS(u)]
	=
- 
\frac{1}{2}
\left\lbrace
	\frac{\FPStwodots}{\FPSdot} 
	\circ 
	(\FPS)^{-1}
	\circ
	\left(\frac{1}{2} \RRiemannPS \right)
	+
	1
\right\rbrace
\frac{d}{du} \dataRRiemannPS(u),
	\label{AE:LSPEEDTIMESMUPLANESYMMETRYWITHEXPLICITSOURCE} 
	\\
\SpeedPS(u) \muPS(t,u) 
& = 
	1
	+ 
	t \PSLmusourcetermfunction(u)
= 1 
	+
	t
	\frac{d}{du}
	\antiderivativePSLmusourcetermfunction[\dataRRiemannPS(u)]
= 1 
	+
	t
	\antiderivativePSLmusourcetermfunction'[\dataRRiemannPS(u)]
	\frac{d}{du}
	\dataRRiemannPS(u),
	\label{AE:SPEEDTIMESMUPLANESYMMETRYWITHEXPLICITSOURCE} 
\end{align}
\end{subequations}

\begin{align} \label{AE:MUXMUPLANESYMMETRYWITHEXPLICITSOURCE}
	\muXPS \muPS(t,u) 
	& = 
	\frac{t
	\frac{d^2}{du^2}
	\antiderivativePSLmusourcetermfunction[\dataRRiemannPS(u)]}{\FPSdot[\dataRRiemannPS(u)]}
	- 
	\left\lbrace
		\frac{\frac{d}{du} \FPSdot[\dataRRiemannPS(u)]}{\FPSdot[\dataRRiemannPS(u)]}
	\right\rbrace
	\left\lbrace
	1 
	+
	t
	\frac{d}{du}
	\antiderivativePSLmusourcetermfunction[\dataRRiemannPS(u)]
	\right\rbrace,
\end{align}

\begin{align} \label{AE:CLOSEDPARTIAL1DERIVATIVEOFRPLUSPS}
	[\partial_1 \RRiemannPS](t,u)
	& = - \frac{1}{\SpeedPS(u) \muPS(t,u)} \geop{u} \RRiemannPS(u)
	= 
	- \frac{\frac{d}{du} \dataRRiemannPS(u)}{1
	+ 
	t \PSLmusourcetermfunction(u)}
	= 
	- \frac{\frac{d}{du} \dataRRiemannPS(u)}{
	1 
	+
	t
	\antiderivativePSLmusourcetermfunction'[\dataRRiemannPS(u)]
	\frac{d}{du}
	\dataRRiemannPS(u)}.
\end{align}	

\end{corollary}

\subsection{The rough time function in simple isentropic plane-symmetry}
\label{SS:ROUGHTIMEFUNCTIONSIMPLEISENTROPICPLANESYMMETRY}
In simple isentropic plane-symmetry, 
using \eqref{AE:PLANESYMMETRYCOMMUTATORVECTORFIELDSAREGEOCOORDINATEVECTORFIELDS},
we can write the transport equation \eqref{E:TRANSPORTEQUATIONFORROUGHTIMEFUNCTION} 
for $\PStimefunctionarg{\muxmulevelsetvalue}(t,u)$ as follows:
\begin{align} \label{AE:PSTRANSPORTEQUATIONFORROUGHTIMEFUNCTION}
\geop{u} \PStimefunctionarg{\muxmulevelsetvalue}(t,u)
+ 
\phi(u)
\frac{\muxmulevelsetvalue}{\geop{t} \muPS(t,u)} \geop{t} \PStimefunctionarg{\muxmulevelsetvalue}(t,u)
& = 0,
\end{align}
where $\phi(u) = \psi\left(\frac{u}{\interestingu} \right)$ is the cut-off from Definition~\ref{D:WTRANSANDCUTOFF}, 
and the differential operator on LHS~\eqref{AE:PSTRANSPORTEQUATIONFORROUGHTIMEFUNCTION} is 
the rough adapted coordinate vectorfield $\roughgeop{u}$,
i.e., we have the following identities (see \eqref{E:ROUGHUDERIVATIVEINTERMSOFGEOMETRICVECTORFIELDSANDROUGHTIMEFUNCTION}):
\begin{align} \label{AE:PSROUGHADAPTEDUDERIVTATIVEINTERMSOFGEOMETRICPARTIALDERIVATIVES}
	\roughgeop{u}
	& 
	=
	\PSWtransarg{\muxmulevelsetvalue}
	= \geop{u}
			+
			\phi
			\frac{\muxmulevelsetvalue}{\geop{t} \muPS} \geop{t}.
\end{align}

We also note that the initial condition \eqref{E:INITIALCONDITIONFORROUGHTIMEFUNCTION} is equivalent to:
\begin{align} \label{AE:PSDATAFORFORROUGHTIMEFUNCTION}
\PStimefunctionarg{\muxmulevelsetvalue}|_{\lbrace \muXPS \muPS = - \muxmulevelsetvalue \rbrace}
& = 
- \muPS |_{\lbrace \muXPS \muPS = - \muxmulevelsetvalue \rbrace}.
\end{align}

\subsection{Admissible background simple isentropic plane-symmetric solutions} 
\label{SS:ADMISSIBLEDATA}
In this section, we construct a large class of ``admissible" simple isentropic plane-symmetric solutions to 
\eqref{AE:COMPRESSIBLEEULEREQUATIONSFORRIEMANNINVARIANTS} that,
on a rough hypersurface,
induce data satisfying the assumptions stated in Sect.\,\ref{S:ASSUMPTIONSONTHEDATA}. 
The main results are Theorem~\ref{T:PROPERTIESOFADMISSIBLESOLUTIONSLAUNCEDBYRESCALED} and
Cor.\,\ref{AC:SIMPLEISENTROPICPLANESYMMETRICSOLUTIONSSATISFYALLTHEASSUMPTIONS}.
We formalize the notion of ``admissible'' in Def.\,\ref{AD:ADMISSIBLEBACKGROUND}.

In Theorem~\ref{T:PROPERTIESOFADMISSIBLESOLUTIONSLAUNCEDBYRESCALED}, we will state
some Sobolev estimates for the solution that involve high order $\geop{u}$
derivatives. The availability of these estimates will simplify our discussion of 
Cauchy stability in Appendix~\ref{A:OPENSETOFDATAEXISTS}, although with additional effort,
it would have been possible for us to prove the needed Cauchy stability result without
such high order $\geop{u}$ estimates. We now define the corresponding Sobolev norms.
Given a scalar function $f = f(t,u)$ of the geometric coordinates,
$N \in \mathbb{N}$,
and real numbers $u_1 \leq u_2$, we define:
\begin{align}
 \| f \|_{H_u^N(\Sigma_t^{[u_1,u_2]})}
& \eqdef	
	\sqrt{
	\sum_{M = 0}^N
	\int_{u_1}^{u_2}
		\left|\left(\geop{u} \right)^M f(t,u) \right|^2
	\, \mathrm{d}u}.
\end{align}

\subsection{Construction of the initial data that lead to admissible shock-forming solutions}
\label{SS:EXAMPLESOFBONAFIDEDATALEADINGTOADMISSIBLE}
In this section, we will exhibit a large family of profiles that, upon rescaling,
yield initial data that launch admissible solutions.

\subsubsection{Assumptions on the ``seed'' profile}
\label{SSS:SEEDPROFILE}
To start, we fix any scalar ``seed profile'' $\mathring{\varphi}$ with the following properties
(it is straightforward to show that such functions exist):

\begin{itemize}
	\item $\mathring{\varphi} = \mathring{\varphi}(u)$ is compactly supported in
		an interval $[-\rightu,\leftu]$ of $u$-values, where $\rightu, \leftu > 1$.
	\item $\mathring{\varphi} \in H_u^{\Ntop+1}(\Sigma_0^{[-\rightu,\leftu]})$ for some integer $\Ntop \geq 24$.
	\item $\frac{d}{du} \mathring{\varphi}(u)$ 
		has a unique, non-degenerate minimum at $u=0$.
	\item Modifying $\mathring{\varphi}$ by multiplying it by a constant 
		and composing it with a linear map of the form $u \rightarrow z u$
		for some constant $z \in \mathbb{R}$ if necessary, we assume that 
		(the modified $\mathring{\varphi}$)
		satisfies
		$\frac{d}{du} \mathring{\varphi}(u)|_{u=0} = - 1$,
		$\frac{d^2}{du^2} \mathring{\varphi}(u)|_{u=0} = 0$,
		that there is a constant $\PSdatamuHessianTaylorcoefficient$ satisfying:
		\begin{align}  \label{AE:PSKEYMUHESSIANTAYLORCOEFFICIENT}
			\PSdatamuHessianTaylorcoefficient 
			& > 0 
		\end{align} 
		such that
		$\frac{d^3}{du^3} \mathring{\varphi}(0) = \PSdatamuHessianTaylorcoefficient$
		and such that
		$\frac{1}{2} \PSdatamuHessianTaylorcoefficient \leq \frac{d^3}{du^3} \mathring{\varphi}(u) \leq 
		2 \PSdatamuHessianTaylorcoefficient$ when $|u| \leq 1$,
		and that there is a constant $\PSLmunottoonegativeparameter$ satisfying:
		\begin{align} \label{AE:PSLOWERBOUNDONNULLDERIVATIVEOFMUINBORINGREGION}
			0 
			& 
			< 
			\PSLmunottoonegativeparameter 
			< 1, 
		\end{align}
		such that
		$\frac{d}{du} \mathring{\varphi}(u) > - \PSLmunottoonegativeparameter$ for $|u| \geq 1$.
\end{itemize}
See Fig.\,\ref{F:SEEDFUNCTION} for the graph (with $u$ increasing from right to left) 
of a representative seed profile.

\begin{center}
	\begin{figure}  
		\begin{overpic}[scale=.6, grid = false, tics=5, trim=-.5cm -1cm -1cm -.5cm, clip]{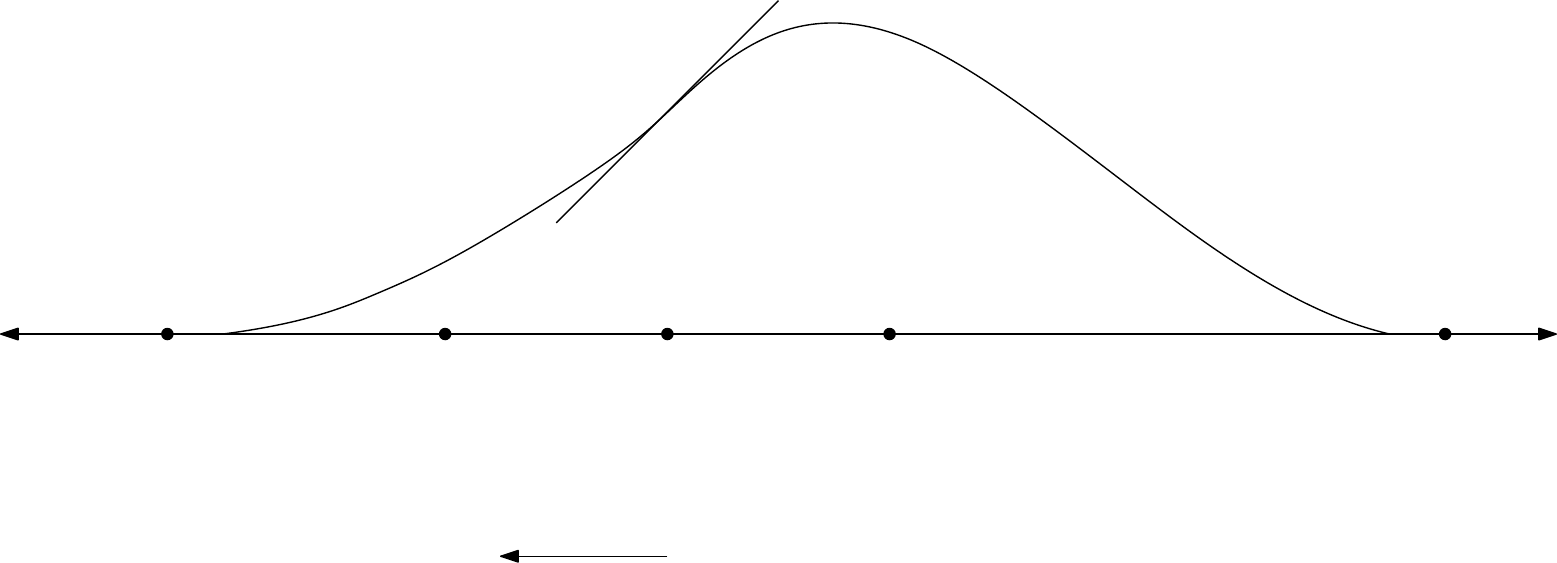}
			\put (9,15) {$u=\leftu$}
			\put (26,15) {$u=1$}
			\put (38,15) {$u=0$}
			\put (51,15) {$u=-1$}
			\put (82,15) {$u = - \rightu$}
			\put (36,7) {$u$}
			\put (65,32) {$\mathring{\varphi}(u)$}
		\end{overpic}
		\caption{The graph of a representative ``seed profile.''}
		\label{F:SEEDFUNCTION}
	\end{figure}
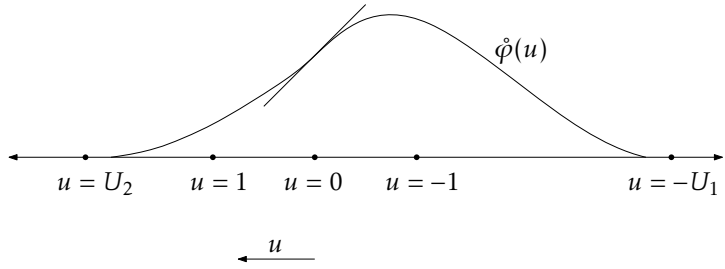
\end{center}

\subsubsection{Construction of one-parameter families of initial data for $\RRiemannPS$}
\label{SSS:ONEPARAMETERFAMILIESOFDATA}
Let $\mathring{\varphi}$ be as in Sect.\,\ref{SSS:SEEDPROFILE}.
Given a real parameter $\PSdataamplitude > 0$, we define:
\begin{align} \label{AE:RESCALEDSEEDDATA}
	\mathring{\varphi}_{\PSdataamplitude}(u) 
	& 
	\eqdef
	\PSdataamplitude \mathring{\varphi}(u),
		\\
	(\dataRRiemannPS)_{\PSdataamplitude}(u)
	& \eqdef \antiderivativePSLmusourcetermfunction^{-1}[\mathring{\varphi}_{\PSdataamplitude}(u)],
	\label{AE:RIEMANNINVARIANTRESCALEDSEEDDATA}
\end{align}
where $\antiderivativePSLmusourcetermfunction^{-1}$ is the inverse function of
of the map
$\RRiemannPS \rightarrow \antiderivativePSLmusourcetermfunction[\RRiemannPS]$ 
define by \eqref{AE:EXPLICITFORMUSEFULIDENTITYFORPSLMUSOURCETERMFUNCTION}.

Taking into account Remark~\ref{AR:NORMALIZATIONFORANTIDERIVATIVEPSLMUSOURCETERM} and \eqref{AE:PLANESYMMETRYNONDEGENCONDITION},
we deduce from Taylor expansions and the standard Sobolev calculus 
that if $\PSdataamplitude$ is sufficiently small, then the following conclusions
hold:
\begin{itemize}
	\item
		$(\dataRRiemannPS)_{\PSdataamplitude}$ is compactly supported in $[-\rightu,\leftu]$
		and satisfies the following bounds:
		\begin{subequations}
		\begin{align} \label{AE:PSRESCALEDDATALINFINITYBOUND} 
			\left\| \frac{d^M}{du^M} (\dataRRiemannPS)_{\PSdataamplitude} \right\|_{L^{\infty}(\Sigma_0^{[-\rightu,\leftu]})} 
			& \lesssim \PSdataamplitude,
			& 
			M = 0,1,2,3,4,
				\\
		\left\| 
			(\dataRRiemannPS)_{\PSdataamplitude} 
		\right\|_{H_u^{\Ntop+1}(\Sigma_0^{[-\rightu,\leftu]})} 
		& \lesssim \PSdataamplitude.
		&&
		\label{AE:PSRESCALEDDATASOBOLVBOUND} 
		\end{align}
		\end{subequations}
	\item $\frac{d}{du} \antiderivativePSLmusourcetermfunction[(\dataRRiemannPS)_{\PSdataamplitude}(u)]$ 
		has a unique, negative, non-degenerate minimum at $u=0$.
	\item There is a differentiable function $\PSthirdorderTaylorremaindercoefficientfunction: [-\rightu,\leftu] \rightarrow \mathbb{R}$ 	
		such that
		$
		\left\| 
			\PSthirdorderTaylorremaindercoefficientfunction
		\right\|_{C_{\textnormal{geo}}^1(\Sigma_0^{[-\rightu,\leftu]})} 
		\lesssim 1
		$ 
		and such that:
		\begin{align} \label{E:TAYLOREXPANSIONOFKEYMUTERMABOUTCREASE}
			\frac{d}{du} \antiderivativePSLmusourcetermfunction[(\dataRRiemannPS)_{\PSdataamplitude}(u)]
			& = - 
					\PSdataamplitude
					+
					\frac{1}{2} \PSdataamplitude \PSdatamuHessianTaylorcoefficient u^2
					+
					\PSthirdorderTaylorremaindercoefficientfunction(u) \PSdataamplitude u^3.
		\end{align}
		\item For $j=1,2,3$, there exists a continuous function 
		$\PSinitialkeymufunctioncoefficients_j : [-\rightu,\leftu] \rightarrow (0,\infty)$
		such that:
		\begin{align} \label{E:INITIALMUFUNCTIONSKEYEXPANSION}
			\PSinitialkeymufunctioncoefficients_j(u)
			& \approx 1,
			& & u \in [-\rightu,\leftu],
		\end{align}
		and such that for $M=0,1,2$, we have:
		\begin{subequations}
		\begin{align}
			\frac{d^M}{du^M}
			\left(
			\frac{d}{du} \antiderivativePSLmusourcetermfunction[(\dataRRiemannPS)_{\PSdataamplitude}(u)]
			+
			\PSdataamplitude
			\right)
			& = \PSinitialkeymufunctioncoefficients_{M+1}(u) \PSdataamplitude u^{2-M},
				\label{E:RESCALEDMUFUNCTIONANDDERIVATIVESEXPANSIONSNEAR0} 
			& \mbox{for } |u| \leq 1,
				\\
			\frac{d}{du} \antiderivativePSLmusourcetermfunction[(\dataRRiemannPS)_{\PSdataamplitude}(u)] 
			& \geq - \PSLmunottoonegativeparameter \PSdataamplitude,
			& \mbox{for } |u| \geq 1,
				\label{E:RESCALEDFIRSTDERIVATIVEOFMUFUNCTIONUNIFORMLYPOSITIVEAWAYFROMINTERESTINGREGION} 
		\end{align}
		\end{subequations}
\end{itemize}
where $\PSLmunottoonegativeparameter < 1$ is the positive constant on 
RHS~\eqref{E:RESCALEDFIRSTDERIVATIVEOFMUFUNCTIONUNIFORMLYPOSITIVEAWAYFROMINTERESTINGREGION}.
In the above relations, all implicit constants are independent of $\PSdataamplitude$
and $u$, though they depend on $\mathring{\varphi}$ and its derivatives.

\subsection{Admissible solutions}
\label{SS:ADMISSIBLESOLUTIONS} 
We now prove the main results of Appendix~\ref{A:PS}, i.e.,
we show that when $\PSdataamplitude$ is sufficiently small, 
the initial datum $(\dataRRiemannPS)_{\PSdataamplitude}$
launches a simple isentropic plane-symmetric solution
that satisfies all the assumptions we stated in stated in Sect.\,\ref{S:ASSUMPTIONSONTHEDATA}.

\subsubsection{Some preliminary definitions}
\label{SSS:PRELIMDEFSFORADMISSIBLESOLUTIONS}
Given any initial data function $\dataRRiemannPS = \dataRRiemannPS(u)$, as in 
\eqref{AE:DATACOMPRESSIBLEEULEREQUATIONSFORRIEMANNINVARIANTS},
we define:
\begin{subequations}
\begin{align} \label{AE:PSDELTASTARDEF}
	\blowupdeltaPS 
	&  
	= \blowupdeltaPS[\dataRRiemannPS] 
	\eqdef 
	\max_{u \in [- \rightu,\leftu]}
	[\PSLmusourcetermfunction(u)]_-
	= \max_{u \in [- \rightu,\leftu]}
		\left[
		\antiderivativePSLmusourcetermfunction'[\dataRRiemannPS(u)]
		\frac{d}{du} \dataRRiemannPS(u)
		\right]_-,
		\\
	\blowuptimePS
	& \eqdef \frac{1}{\blowupdeltaPS},
	\label{AE:PSBLOWUPTIME}
\end{align}
\end{subequations}
where $[z]_- \eqdef \max\lbrace -z, 0 \rbrace$.
In view of \eqref{AE:PSLMUSOURCETERMFUNCTION} and \eqref{AE:SHORTHANDNOTATIONFORMAINLMUSOURCETERM},
we see that for simple isentropic plane-wave solutions,
$\blowupdeltaPS$ coincides with the quantity defined in \eqref{E:DELTASTARDEF}.
Moreover, from Cor.\,\ref{AC:PSEXPLICITEXPRESSIONSFORSOLUTION}, we see that
$\blowuptimePS$ is the Cartesian time of first blowup of $\partial_1 \RRiemannPS$.

\subsubsection{The existence of admissible simple isentropic plane-symmetric solutions}
\label{SSS:EXISTENCEOFADMISSIBLESOLUTIONS}
We now state and prove the main results of this appendix. 

\begin{theorem}[The existence of admissible simple isentropic plane-symmetric solutions]
\label{T:PROPERTIESOFADMISSIBLESOLUTIONSLAUNCEDBYRESCALED}
Let $\mathring{\varphi}$ 
be a ``seed profile'' function with the properties stated in 
Sect.\,\ref{SS:EXAMPLESOFBONAFIDEDATALEADINGTOADMISSIBLE}.
In particular, $\mathring{\varphi}$ is supported in the interval $[-\leftu,\rightu]$ of $u$-values 
(where by \eqref{AE:DATAFOREIKONALEQUATIONINPLANESYMMETRY}, $u=-x^1$ along $\Sigma_0$),
and $\| \mathring{\varphi} \|_{H_u^{\Ntop+1}(\Sigma_0^{[-\rightu,\leftu]} )} < \infty$.
There exist small constants $\PSdataamplitude_0 > 0$
and $\PSinterestingusmallmultipleofamplitude > 0$,
depending on $\mathring{\varphi}$ 
(including the constant $\PSdatamuHessianTaylorcoefficient > 0$ from \eqref{AE:PSKEYMUHESSIANTAYLORCOEFFICIENT})
and satisfying\footnote{Our assumption $\PSinterestingusmallmultipleofamplitude < \PSLmunottoonegativeparameter$
ensures that \eqref{E:MU1BIGGERTHANMU0} is satisfied.} $\PSinterestingusmallmultipleofamplitude < \PSLmunottoonegativeparameter$
(where $\PSLmunottoonegativeparameter < 1$ is the positive constant on 
RHS~\eqref{AE:PSLOWERBOUNDONNULLDERIVATIVEOFMUINBORINGREGION})
such that if $0 < \PSdataamplitude \leq \PSdataamplitude_0$,
then the following conclusions hold,
where in all estimates,
the constants (including the implicit ones corresponding to ``$\lesssim$'') 
are independent of $\PSdataamplitude$, 
$\PSinterestingusmallmultipleofamplitude$, $t$, and $u$
(on the domains where the estimates are asserted to hold).

Let $(\dataRRiemannPS)_{\PSdataamplitude}$
be the initial datum,
defined by \eqref{AE:RESCALEDSEEDDATA}--\eqref{AE:RIEMANNINVARIANTRESCALEDSEEDDATA},
for the simple isentropic plane-symmetric compressible Euler equations,
i.e., for \eqref{AE:COMPRESSIBLEEULEREQUATIONSFORRIEMANNINVARIANTS} with $\LRiemannPS \equiv 0$.
Let $\timefunction_0$ be any negative real number such that: 
\begin{align} \label{AE:PSTIMEFUNCTION0CONSTRAINTS}
- \frac{\PSinterestingusmallmultipleofamplitude}{4}
& 
<  
\timefunction_0 
< 0,
\end{align}
and define $\interestingu > 0$ and $\muxmulevelsetvalue_0 > 0$ by:
\begin{align} \label{AE:PSINTERESTINGUISSMALLMULTIPLEOFAMPLITUDE}
	\interestingu
	& \eqdef \frac{\PSdataamplitude}{\PSdatamuHessianTaylorcoefficient^2},
		\\
	\muxmulevelsetvalue_0
	& \eqdef \frac{|\timefunction_0| \PSdataamplitude}{16 \PSdatamuHessianTaylorcoefficient},
	\label{AE:PSKAPPSOISMULTIPLEOFTIMEFUNCTION0TIMESPSIAMPLITUDE}
\end{align}
where $\PSdatamuHessianTaylorcoefficient > 0$ is the $\mathring{\varphi}$-dependent constant from Sect.\,\ref{SSS:SEEDPROFILE}.


\medskip
\noindent \underline{\textbf{Classical existence with respect to the geometric coordinates}}. 
With respect to the geometric coordinates $(t,u)$, the solution $\RRiemannPS$ is a function of $u$ alone
(i.e., $\RRiemannPS = \RRiemannPS(u)$),
vanishes on the complement of the region 
$\lbrace (t,u) \in \mathbb{R} \times \mathbb{R} \ | \ u \in [- \rightu,\leftu] \rbrace$,
and exists classically for $(t,u) \in [0,\blowuptimePS] \times \mathbb{R}$,
where the parameter $\blowupdeltaPS = \blowupdeltaPS[(\dataRRiemannPS)_{\PSdataamplitude}]$ defined by \eqref{AE:PSDELTASTARDEF}
satisfies:
\begin{align} \label{AE:PSBLOWUPTIMEINTERMSOFAMPITUDE}
	\blowupdeltaPS
	& 
	= 
	\PSdataamplitude,
\end{align}
and:
\begin{align} \label{AE:PSBLOWUPTIME}
\blowuptimePS
& \eqdef
	\frac{1}{\blowupdeltaPS}
	= \frac{1}{\PSdataamplitude}.
\end{align}
Similarly, the null vectorfield $\LunitPS$ and inverse foliation density $\muPS$ exist classically 
for $(t,u) \in [0,\blowuptimePS] \times \mathbb{R}$,
and on the complement of 
$\lbrace (t,u) \in \mathbb{R} \times \mathbb{R} \ | \ u \in [- \rightu,\leftu] \rbrace$,
we have $(\LunitPS)^1 = \LunitPS x^1 = 1$ and $\muPS = 1$.
Finally, $\muPS > 0$ on $\left([0,\blowuptimePS] \times \mathbb{R} \right) \backslash \lbrace (\blowuptimePS,0) \rbrace$,
and $\muPS(\blowuptimePS,0) = 0$.

\medskip
\noindent \underline{\textbf{Description of the crease relative to the geometric coordinates}}. 
Relative to the geometric coordinates $(t,u)$, 
the crease, which by definition is 
$\creasePS 
\eqdef 
\{ (t,u) \ | \ \muPS(t,u) = 0\} 
\cap \{(t,u) \ | \ \muXPS \muPS(t,u) = 0\} 
\cap 
\left([0,\blowuptimePS] 
\times 
[-\interestingu,\interestingu] \right)$,
is equal to the single point $(\blowuptimePS,0)$.

\medskip
\noindent \underline{\textbf{The Cartesian coordinate description of the singularity formation up to the crease}}. 
\begin{itemize}
	\item Let $\PSUpsilon(t,u) \eqdef (t,x^1)$ denote the change of variables map from geometric to Cartesian coordinates. 
	Then $\PSUpsilon$ is a homeo- (resp.\ diffeo)-morphism from 
	$[0,\blowuptimePS] \times \mathbb{R}$ 
	(resp.\ from $\left([0,\blowuptimePS] \times \mathbb{R} \right) \setminus \creasePS$) onto its image. 
	\item The solution $\RRiemannPS$ exists classically with respect to the Cartesian coordinates on the subset
	$\PSUpsilon\left(([0,\blowuptimePS] \times \mathbb{R}) \setminus \creasePS \right)$ of Cartesian coordinate space 
	$\mathbb{R}_t \times \mathbb{R}_{x^1}$. 
	\item There exists a past neighborhood $\mathcal{N}$ of $\PSUpsilon(\creasePS)$ in Cartesian coordinate space
		such that, with $\RRiemannPS = \RRiemannPS(t,u) = \RRiemannPS(u)$, we have:
	\begin{align} \label{AE:PSBLOWUPLOWERBOUND}
		\left|
			[\partial_1 \RRiemannPS] \circ \PSInverseUpsilon(t,x^1)
		\right| 
		=
		\frac{1}{\muPS}
		\left|\frac{1}{\SpeedPS} \muXPS \RRiemannPS \right| 
		\circ 
		\PSInverseUpsilon(t,x^1)
		& \gtrsim 
		\frac{1}{\muPS},
		&&
		(t,x^1) \in \mathcal{N}.
	\end{align}
		In particular, for any sequence of points 
		$\lbrace q_n \rbrace_{n \in \mathbb{N}} 
		\subset 
		\PSUpsilon\left(([0,\blowuptimePS] \times \mathbb{R}) \setminus \creasePS \right)$
		converging to the single point in $\PSUpsilon(\creasePS)$ 
		(which implies that $\PSInverseUpsilon(q_n) \rightarrow (\blowuptimePS,0)$),
		we have that $|\partial_1 \RRiemannPS| \circ \PSInverseUpsilon(q_n) \to \infty$ as $n \to \infty$.  
\end{itemize}

\medskip
\noindent \underline{\textbf{Estimates for the solution and acoustic geometry with respect to the geometric coordinates}}.
For $(t,u) \in [0,\blowuptimePS] \times \mathbb{R}$, the
speed of sound satisfies the following estimate:
\begin{align} \label{AE:PSSPEEDOFSOUNDESTIMATE}
	\SpeedPS(t,u)
	=
	\SpeedPS(u)
	 = 1 + \mathcal{O}(\PSdataamplitude).
\end{align}

For $t \in [0,\blowuptimePS]$, the following Sobolev estimates hold:
\begin{subequations}
\begin{align} \label{AE:PSSOBOLEVBOUNDSFORRRIEMANN}
	\left\| 
		\RRiemannPS 
	\right\|_{H_u^{\Ntop+1}(\Sigma_t)}
	& 
	\lesssim \PSdataamplitude,
		\\
	\left\| 
		(\LsmallPS)^1 
	\right\|_{H_u^{\Ntop+1}(\Sigma_t)}
	& 
	\lesssim \PSdataamplitude,
		\label{AE:PSSOBOLEVBOUNDSFORLSMALL} 
			\\
	\left\| 
		\geop{t} \muPS 
	\right\|_{H_u^{\Ntop}(\Sigma_t)}
	& 
	\lesssim \PSdataamplitude,
	\label{AE:PSSOBOLEVBOUNDSFORGEOPTDERIVATIVEOFINVERSEFOLIATIONDENSITY} 
		\\
	\left\| 
		\muPS 
	\right\|_{H_u^{\Ntop}(\Sigma_t)}
	& 
	\lesssim 1.
	\label{AE:PSSOBOLEVBOUNDSFORINVERSEFOLIATIONDENSITY} 
\end{align}
\end{subequations}

For $M = 0,1,2,3,4$ and $t \in [0,\blowuptimePS]$, the following estimates hold:
\begin{subequations}
\begin{align}  
	\left\| 
		(\muXPS)^M \RRiemannPS 
	\right\|_{L^{\infty}(\Sigma_t)} 
	& 
	\lesssim \PSdataamplitude,
		\label{AE:PSSMALLAMPLITUDESOLUTIONLINFINITYRRIEMANNTRANSVERSALDERIVATIVES}
			\\
	\left\| 
		(\muXPS)^M (\LsmallPS)^1
	\right\|_{L^{\infty}(\Sigma_t)} 
	& 
	\lesssim \PSdataamplitude.
	\label{AE:PSSMALLAMPLITUDELISMALLLINFINITYRRIEMANNTRANSVERSALDERIVATIVES}
\end{align}
\end{subequations}

Similarly, for $M=0,1,2,3$ and $t \in [0,\blowuptimePS]$, we have:
\begin{subequations} 
	\begin{align} \label{AE:LINFINITYINITIALROUGHHYPERSURFACELDERIVATIVEOFMUANDTRANSVERSALDERIVATIVES} 
	\left\| \LunitPS (\muXPS)^M \muPS \right\|_{L^{\infty}(\Sigma_t)}
	& 	= 
			\frac{1}{2} 
			\left\|
				(\muXPS)^M
				\left\lbrace
				\frac{1}{\SpeedPS}
				\left(
				\frac{(\SpeedPS)'}{\SpeedPS}
				+ 
				1
				\right)  
				\muXPS \RRiemannPS 
				\right\rbrace
			\right\|_{L^{\infty}(\Sigma_0)}
			\lesssim \PSdataamplitude,
			\\
	\begin{split}  \label{AE:LINFINITYMUANDTRANSVERSALDERIVATIVES}  
	\left\| (\muXPS)^M \muPS \right \|_{L^{\infty}(\Sigma_t)}
	& \leq
		\left\|
			 (\muXPS)^M
			\left\lbrace
				\frac{1}{\SpeedPS}
			\right\rbrace
		\right\|_{L^{\infty}(\Sigma_0)}
		+
		\frac{1}{2}
		\frac{1}{\blowupdeltaPS}
		\left\|
				 (\muXPS)^M
				\left\lbrace
				\frac{1}{\SpeedPS}
				\left(
				\frac{(\SpeedPS)'}{\SpeedPS}
				+ 
				1\right)  
				\muXPS \RRiemannPS 
				\right\rbrace
			\right\|_{L^{\infty}(\Sigma_0)}
			\\
		& \lesssim 1.
\end{split}
\end{align}

\end{subequations} 

\medskip

\medskip
\noindent \underline{\textbf{Estimates tied to the change of variables map $\PSUpsilon$}}. 
For $t \in [0,\blowuptimePS]$, the Cartesian spatial coordinate $x^1 = x^1(t,u)$ satisfies 
the following estimates:
\begin{align} \label{AE:PSSIZEOFCARTESIANX1}
		- 
		\leftu
		+
		t
		& 
		\leq
		\min_{\Sigma_t^{[- \rightu,\leftu]}} x^1
		\leq
		\max_{\Sigma_t^{[- \rightu,\leftu]}} x^1
		\leq 
		\rightu
		+
		t.
\end{align}

Moreover, 
\begin{align} \label{AE:PSIDENTITYFORDIFFERENTIALOFGEOTOCARTESIANCHOVMAP}
	d_{\textnormal{geo}} \PSUpsilon(t,u)
	&
	\eqdef
	\frac{\partial \PSUpsilon(t,u)}{\partial (t,u)}
	= \begin{pmatrix}
			1 & 0 
				\\
			(\LunitPS)^1 & - \SpeedPS \muPS
	\end{pmatrix},
\end{align}
and the following estimates hold for $t \in [0,\blowuptimePS]$:
\begin{subequations}
\begin{align} \label{AE:PSHIGHORDERESTIMATESFORUPSILON}
	\left\| 
		d_{\textnormal{geo}} \PSUpsilon
	\right\|_{C_{\textnormal{geo}}^{\Ntop-1}(\Sigma_t)}
	& \lesssim 
		1,
			\\
	\left\| 
		\geop{t} d_{\textnormal{geo}} \PSUpsilon
	\right\|_{H_u^{\Ntop}(\Sigma_t)},
		\,
	\left\| 
		\geop{t} d_{\textnormal{geo}} \PSUpsilon
	\right\|_{C_{\textnormal{geo}}^{\Ntop-1}(\Sigma_t)}
	& \lesssim 
		\PSdataamplitude.
		\label{AE:PSSOBOLEVESTIMATESFORGEOPTUPSILON}
\end{align}
\end{subequations}

\medskip
\noindent \underline{\textbf{Properties of the rough time functions $\PStimefunctionarg{\muxmulevelsetvalue}$ and the location of their level sets}}. 
We define: 
\begin{align} \label{AE:BIGDELTAPS} 
\PSBigDelta 
& 
\eqdef 
\frac{|\timefunction_0|}{16}
\blowuptimePS
=
\frac{|\timefunction_0|}{16 \PSdataamplitude}.
\end{align}
For $\timefunction \in [2 \timefunction_0,\frac{1}{2} \timefunction_0]$ and $\muxmulevelsetvalue \in [0,\muxmulevelsetvalue_0]$,
there exist functions 
$\PSCartesiantisafunctiononlevelsetsofroughtimefunctionarg{\timefunction}{\muxmulevelsetvalue}: 
\mathbb{R} \rightarrow [\blowuptimePS - 2\PSinterestingusmallmultipleofamplitude \blowuptimePS,\blowuptimePS - 2 \PSBigDelta]$,
depending on $\timefunction$ and $\muxmulevelsetvalue$ and strictly increasing with respect to
$\timefunction$, 
such that 
$\| \PSCartesiantisafunctiononlevelsetsofroughtimefunctionarg{\timefunction}{\muxmulevelsetvalue} \|_{C^3(\mathbb{R})} \lesssim 1$
and such that the rough time function $\PStimefunctionarg{\muxmulevelsetvalue}$ exists as a $C^3$ function of 
$(t,u)$ on the domain: 
\begin{align} \label{AE:PSROUGHTIMEFUNCTIONDOMAININGEOMETRICCOORDINATES}
\PStwoargMrough{[2 \timefunction_0,\frac{1}{2} \timefunction_0],\mathbb{R}}{\muxmulevelsetvalue}
=
\left\lbrace (t,u) \ | \ 
\PSCartesiantisafunctiononlevelsetsofroughtimefunctionarg{2 \timefunction_0}{\muxmulevelsetvalue}(u) \leq t \leq 
\PSCartesiantisafunctiononlevelsetsofroughtimefunctionarg{\frac{1}{2} \timefunction_0}{\muxmulevelsetvalue}(u),
	\,
u \in \mathbb{R}
\right\rbrace,
\end{align}
and such that relative to the geometric coordinates,
the level sets 
$\twoarghypPS{\timefunction}{\muxmulevelsetvalue} \eqdef \lbrace (t,u) \in \mathbb{R} \times \mathbb{R} 
\ | \ \PStimefunctionarg{\muxmulevelsetvalue}(t,u) = \timefunction \rbrace$ 
are the following graphical surfaces:
					\begin{align} \label{AE:LEVELSETSOFTIMEFUNCTIONAREAGRAPH}
						\twoarghypPS{\timefunction}{\muxmulevelsetvalue}
						& =
						\left\lbrace
							(t,u)
							\ | \
							t = \PSCartesiantisafunctiononlevelsetsofroughtimefunctionarg{\timefunction}{\muxmulevelsetvalue}(u),
								\,
							u \in \mathbb{R}
						\right\rbrace.
					\end{align}
In particular, 
\begin{align} \label{AE:ROUGHHYPSFORADMISSIBLESOLUTIONSCONTAINEDINCAUCHYSTABILITYREGION}
\bigcup_{\timefunction \in [2 \timefunction_0,\frac{1}{2} \timefunction_0]}
\twoarghypPS{\timefunction}{\muxmulevelsetvalue}
& \subset 
\bigcup_{t \in [\blowuptimePS - 2\PSinterestingusmallmultipleofamplitude \blowuptimePS,\blowuptimePS - 2 \PSBigDelta]} \Sigma_t.
\end{align} 
Moreover, the following estimates hold for $(t,u) \in \PStwoargMrough{[2 \timefunction_0,\frac{1}{2} \timefunction_0],(-\infty,\infty)}{\muxmulevelsetvalue}$:
\begin{align}
	\frac{15}{16} \blowupdeltaPS
	\leq
	\geop{t} \PStimefunctionarg{\muxmulevelsetvalue}(t,u)
	&
	\leq \frac{17}{16} \blowupdeltaPS,
		 \label{AE:PSGEOPTDERIVATIVEOFROUGHTIMEFUNCTIONKEYBOUND} 
			\\
	\left|
		\geop{u} \PStimefunctionarg{\muxmulevelsetvalue} 
	\right|
	&
	\leq 2 \muxmulevelsetvalue_0
			=
			\frac{|\timefunction_0| \PSdataamplitude}{8 \PSdatamuHessianTaylorcoefficient}.
				\label{AE:PSGEOPUDERIVATIVEOFROUGHTIMEFUNCTIONSIMPLEBOUND} 
\end{align}

\medskip
\noindent \underline{\textbf{Properties of $\PSCHOVgeotorough{\muxmulevelsetvalue}$ and $\PSCHOVJacobianroughtomumuxmu{\muxmulevelsetvalue}$}}. 
The change of variables map $\PSCHOVgeotorough{\muxmulevelsetvalue}$ defined by:
\begin{align} \label{AE:PSCHOVGEOTOROUGH}
	\PSCHOVgeotorough{\muxmulevelsetvalue}(t,u) 
	& 
	\eqdef (\PStimefunctionarg{\muxmulevelsetvalue},u)
\end{align}		
is a diffeomorphism from
$
\twoargMrough{[2 \timefunction_0,\frac{1}{2} \timefunction_0],\R}{\muxmulevelsetvalue}
$
onto its image, which is
$
[2 \timefunction_0,\frac{1}{2} \timefunction_0] \times \mathbb{R}
$,
and it satisfies: 
\begin{align}\label{AE:PSC3ESTIMATESFORCHOVMAPFROMGEOTOROUGH}
	\| 	\PSCHOVgeotorough{\muxmulevelsetvalue} \|_{C_{\textnormal{geo}}^3(\PStwoargMrough{[2 \timefunction_0,\frac{1}{2} \timefunction_0],(-\infty,\infty)}{\muxmulevelsetvalue})}
	& \lesssim 1, 
	&&
	\\
	\frac{15}{16} \blowupdeltaPS
	&
	\leq
	\mbox{\upshape det} \left( d_{\textnormal{geo}} \PSCHOVgeotorough{\muxmulevelsetvalue} \right)  = \geop{t} \PStimefunctionarg{\muxmulevelsetvalue} 
	\leq \frac{17}{16} \blowupdeltaPS,
	&&
	\text{on } \PStwoargMrough{[2\timefunction_0,\frac{1}{2}\timefunction_0],(-\infty,\infty)}{\muxmulevelsetvalue}. 
	\label{AE:DETERMINATOFPSCHOVMAPFROMGEOTOROUGH}
\end{align}

In addition, for every $\muxmulevelsetvalue \in [0,\muxmulevelsetvalue_0]$,
the Jacobian matrix 
$
\PSCHOVJacobianroughtomumuxmu{\muxmulevelsetvalue} 
\eqdef
\frac{\partial (\muPS,\muXPS \muPS)}{(\PStimefunctionarg{\muxmulevelsetvalue},u)}
$
(see also \eqref{E:JACOBIANMATRIXFORCHOVFROMROUGHCOORDINATESTOMUWEGIGHTEDXMUCOORDINATES})
is invertible for every 
$q \eqdef (\timefunction,u) \in [2 \timefunction_0,\frac{1}{2} \timefunction_0] \times [-\interestingu,\interestingu]$
and satisfies:
\begin{align} \label{AE:PSDATAJACOBIANDETERMINANTRATIOBOUND}
	\sup_{q_1,q_2 \in [2 \timefunction_0,\frac{1}{2} \timefunction_0] \times [-\interestingu,\interestingu]}
	\left|
		\PSInverseCHOVJacobianroughtomumuxmu{\muxmulevelsetvalue}(q_1) \PSCHOVJacobianroughtomumuxmu{\muxmulevelsetvalue}(q_2)
		-
		\mbox{\upshape ID} 
	\right|_{\mbox{\upshape}Euc}
	& 
	\leq 
	\frac{1}{4},
\end{align}
where 
$|\cdot|_{\mbox{\upshape}Euc}$ is the standard Frobenius norm on matrices
(equal to the square root of the sum of the squares of the matrix entries)
and $\mbox{\upshape ID}$ denotes the $2 \times 2$ identity matrix.

\medskip
\noindent \underline{\textbf{Properties of $\PSCHOVgeotomumuxmu$}}. 
We define the map $\PSCHOVgeotomumuxmu$ 
	from geometric coordinates to ``$(\upmu,\muX \upmu)$-space''
	and its Jacobian 
	$\PSCHOVJacobiangeotomumuxmu$ as follows
	(see also 
	\eqref{E:CHOVFROMGEOMETRICCOORDINATESTOMUWEGIGHTEDXMUCOORDINATES}--\eqref{E:JACOBIANMATRIXFORCHOVFROMGEOMETRICCOORDINATESTOMUWEGIGHTEDXMUCOORDINATES}):
	\begin{subequations}
	\begin{align} \label{E:PSCHOVFROMGEOMETRICCOORDINATESTOMUWEGIGHTEDXMUCOORDINATES}
		\PSCHOVgeotomumuxmu(t,u,x^2,x^3)
		& \eqdef (\upmu,\muX \upmu),
			\\
		\PSCHOVJacobiangeotomumuxmu(t,u)
		& \eqdef 
		\frac{\partial \CHOVgeotomumuxmu(t,u)}{\partial(t,u)}
		=
		\frac{\partial (\upmu,\muX \upmu)}{\partial(t,u)}.
			\label{E:PSJACOBIANMATRIXFORCHOVFROMGEOMETRICCOORDINATESTOMUWEGIGHTEDXMUCOORDINATES}
	\end{align}
	\end{subequations}

	Then there is a $C > 1$ such that for $\muxmulevelsetvalue \in [0,\muxmulevelsetvalue_0]$,
	$\PSCHOVJacobiangeotomumuxmu$
	is invertible for 
	$(t,u) 
\in
[\blowuptimePS - 2 \PSinterestingusmallmultipleofamplitude \blowuptimePS,\blowuptimePS] 
			\times [- \interestingu, \interestingu]$
	and satisfies the following bounds:
	\begin{align} 
	- C
	&
	\leq
	\mbox{\upshape det} \CHOVJacobiangeotomumuxmu
	\leq 
		- \frac{1}{C},
	&&
	\text{on }
	[\blowuptimePS - 2 \PSinterestingusmallmultipleofamplitude \blowuptimePS,\blowuptimePS] 
			\times [- \interestingu, \interestingu],
		\label{E:PSJACOBIANDETBOUNDCHOVGEOTOMUXMUCOORDS} 
\end{align}
\begin{align}
	\sup_{p_1,p_2 \in [\blowuptimePS - 2 \PSinterestingusmallmultipleofamplitude \blowuptimePS,\blowuptimePS] 
			\times [- \interestingu, \interestingu]}
	\left|
		\PSCHOVJacobiangeotomumuxmu(p_1) \PSInverseCHOVJacobiangeotomumuxmu(p_2) 
		-
		\mbox{\upshape ID} 
	\right|_{\mbox{\upshape}Euc}
	& \leq \frac{1}{4},
		\label{E:PSGEOTOMUMUXCOORDINATESJACOBIANDETERMINANTRATIOINEVOLUTIONBOUND}
\end{align}
where $|\cdot|_{\mbox{\upshape}Euc}$ is the standard Frobenius norm on matrices
(equal to the square root of the sum of the squares of the matrix entries)
and $\mbox{\upshape ID}$ denotes the $2 \times 2$ identity matrix.

\medskip
\noindent \underline{\textbf{Behavior of $\muPS$ in the interesting region}}.
The following estimates hold for $t \in [0,\blowuptimePS]$:
\begin{align} \label{AE:BOUNDSONLMUINTERESTINGREGION} 
	- 
	\frac{17}{16}
	\blowupdeltaPS
	\leq
	\min_{\Sigma_t^{[-\interestingu,\interestingu]}} \LunitPS \muPS(t,u)
	\leq 
	\max_{\Sigma_t^{[-\interestingu,\interestingu]}} \LunitPS \muPS(t,u)
	\leq 
	- \frac{15}{16} \blowupdeltaPS.
\end{align}

Moreover, for
$
(t,u) 
\in
[0,\blowuptimePS] 
			\times [-\interestingu,\interestingu]
$,
the following estimates hold:
\begin{align}
			 \muPS(t,u) 
				& 
				= 
					\left\lbrace
						1 + \mathcal{O}(\PSdataamplitude)
					\right\rbrace
					\frac{\PSdatamuHessianTaylorcoefficient}{2} u^2
					+
					\left\lbrace
						1 + \mathcal{O}(\PSdataamplitude)
					\right\rbrace
					\PSdataamplitude
					(\blowuptimePS - t), 
				\label{AE:PSMUTAYLOREXPANSIONININTERESTINGREGION} 
					\\
				\LunitPS \muPS(t,u) 
				& 
				= - 
					\left\lbrace
						1 + \mathcal{O}(\PSdataamplitude)
					\right\rbrace
					\PSdataamplitude, 
				\label{AE:PSLUNTMUTAYLOREXPANSIONININTERESTINGREGION} 
					\\
		\muXPS \muPS(t,u) 
		& 
			=
			\left\lbrace
				1 + \mathcal{O}(\PSdataamplitude)
			\right\rbrace
			\PSdatamuHessianTaylorcoefficient u
			+
			\mathcal{O}(\PSdataamplitude^2)
			(\blowuptimePS - t),
			\label{AE:PSMUFIRSTDERIVATIVETAYLOREXPANSIONININTERESTINGREGION}
				\\
			\LunitPS \muXPS \muPS(t,u)
				& 
				= 
					\mathcal{O}(\PSdataamplitude^2),
			\label{AE:PSLUNITMUXMUTAYLOREXPANSIONININTERESTINGREGION} 
					\\
		\muXPS \muXPS \muPS(t,u) 
		& = \left\lbrace
				1 + \mathcal{O}(\PSdataamplitude)
			\right\rbrace
			\PSdatamuHessianTaylorcoefficient
			+
			\mathcal{O}(\PSdataamplitude)
			(\blowuptimePS - t).
			\label{AE:PSMUSECONDDERIVATIVETAYLOREXPANSIONININTERESTINGREGION}
\end{align}

In addition, for $t \in [\blowuptimePS - \blowuptimePS \PSinterestingusmallmultipleofamplitude,\blowuptimePS]$
and $\muxmulevelsetvalue \in [0,\muxmulevelsetvalue_0]$, we have:\footnote{It might appear that there are fewer quantities in braces in \eqref{AE:PSMUTRANSVERSALCONVEXITY} compared
to \eqref{E:DATATASSUMPTIONMUTRANSVERSALCONVEXITY}, but this is not true. 
The reason is that in simple isentropic plane-symmetry,
some of the quantities in braces in \eqref{E:DATATASSUMPTIONMUTRANSVERSALCONVEXITY} are equal to each other.
\label{AFN:SOMEREDUNDANTTRANSVERSALCONVEXITYQUANTITIES}}
	\begin{align} \label{AE:PSLOCATIONOFXMUEQUALSMINUSKAPPA}
			\lbrace \muXPS \muPS = - \muxmulevelsetvalue \rbrace
			\cap
			\Sigma_t^{[-\interestingu,\interestingu]}
			\subset
			\Sigma_t^{[-\frac{1}{4}\interestingu,\frac{1}{4}\interestingu]},
					\\
	\min_{\Sigma_t^{[-\interestingu,\interestingu]}
			\backslash
			\Sigma_t^{[-\frac{1}{2}\interestingu,\frac{1}{2}\interestingu]}}
			|\muXPS \muPS + \muxmulevelsetvalue| 
			& \geq \frac{\PSdatamuHessianTaylorcoefficient \interestingu}{8},
			\label{AE:PSREGIONWHEREMUXMUKAPPALEVELSETISNOTLOCATED}
	\end{align}

\begin{equation}  \label{AE:PSMUTRANSVERSALCONVEXITY}
\begin{split}
			\frac{\PSdatamuHessianTaylorcoefficient}{2}
			& \leq 
			\min_{\Sigma_t^{[-\interestingu,\interestingu]}}
				\left\lbrace
				\PSWtransarg{\muxmulevelsetvalue} \PSWtransarg{\muxmulevelsetvalue} \muPS,
					\,
				\PSWtransarg{\muxmulevelsetvalue} \muXPS \muPS,
					\,
				\muXPS \muXPS \muPS,
					\,
				\geop{u} \muXPS \muPS - \frac{(\geop{u} \muPS) \geop{t} \muXPS \muPS}{\geop{t} \muPS}
				\right\rbrace
					\\
			&
			\leq 
			\max_{\Sigma_t^{[-\interestingu,\interestingu]}} 
			\left\lbrace
				\PSWtransarg{\muxmulevelsetvalue} \PSWtransarg{\muxmulevelsetvalue} \muPS,
					\,
				\PSWtransarg{\muxmulevelsetvalue} \muXPS \muPS,
					\,
				\muXPS \muXPS \muPS,
					\,
				\geop{u} \muXPS \muPS - \frac{(\geop{u} \muPS) \geop{t} \muXPS \muPS}{\geop{t} \muPS}
				\right\rbrace
			\leq 
			2 \PSdatamuHessianTaylorcoefficient.
\end{split}
\end{equation}

\medskip
\noindent \underline{\textbf{$\muPS$ is uniformly positive away from the interesting region}}.
With $\PSLmunottoonegativeparameter < 1$ denoting the positive constant on 
RHS~\eqref{E:RESCALEDFIRSTDERIVATIVEOFMUFUNCTIONUNIFORMLYPOSITIVEAWAYFROMINTERESTINGREGION},
we have the following estimate:
\begin{align} \label{AE:PSMUISLARGEAWAYFROMINTERESTINGREGION}
	\min_{\lbrace (t,u) \ | \ t \in [0,\blowuptimePS], \,  |u| \geq \interestingu \rbrace} \muPS(t,u)
	& 
	\geq \frac{1}{2} \PSLmunottoonegativeparameter.
\end{align}
\end{theorem}

\begin{remark}[Generalizations of Theorem~\ref{T:PROPERTIESOFADMISSIBLESOLUTIONSLAUNCEDBYRESCALED}] 
\label{AR:GENERALIZATIONSOFTHEOREMPROPERTIESOFADMISSIBLESOLUTIONSLAUNCEDBYRESCALED}
Before proving the theorem, we first make a series of remarks on how it could be extended.

\begin{itemize}
	\item In Theorem~\ref{T:PROPERTIESOFADMISSIBLESOLUTIONSLAUNCEDBYRESCALED},
		we chose to follow the solution up to the Cartesian time of first blowup
		because we believe that the results could help prepare the reader for the more
		difficult analysis in the bulk of the paper.
		However, we could have ``stopped the analysis'' before then; our main
		goal in the theorem was to construct the rough time functions $\PStimefunctionarg{\muxmulevelsetvalue}$
		and to describe the state of the solution near the hypersurfaces
		$\twoarghypPS{\timefunction_0}{\muxmulevelsetvalue}$ so that we can use these
		results in Appendix~\ref{A:OPENSETOFDATAEXISTS}, in our study of Cauchy stability.
	\item Our definition \eqref{AE:PSKAPPSOISMULTIPLEOFTIMEFUNCTION0TIMESPSIAMPLITUDE} of	
		$\muxmulevelsetvalue_0$ is such that $\muxmulevelsetvalue_0$ decreases as $|\timefunction_0| \downarrow 0$.
		This is highly non-optimal and is an artifact of our insistence
		(out of convenience) that the level sets 
		$\twoarghypPS{\timefunction_0}{\muxmulevelsetvalue}$
		should be contained in the rectangular shaped domain
		$[0,\blowuptimePS - 2 \PSBigDelta] \times \mathbb{R}$ in geometric coordinate space;
		by studying the solution on a larger (curved) subset of geometric coordinate space, 
		we could have shown that $\muxmulevelsetvalue_0$ can be chose to be independent of all sufficiently small $|\timefunction_0|$.
		In a similar vein, with modest additional effort, 
		we could have shown that $\interestingu$  
		can be chosen to be independent of all sufficiently small $\PSdataamplitude$.
	\item One could generalize Theorem~\ref{T:PROPERTIESOFADMISSIBLESOLUTIONSLAUNCEDBYRESCALED} 
		in a straightforward fashion to allow for much more general initial data 
		of simple isentropic plane-symmetric type.
		For example, one could consider seed profile functions $\mathring{\varphi}$ of multi-bump type,
		leading to the formation distinct shocks that are separated in space.
		One could also consider two-parameter families of rescaled seed profile functions of the form
		$\mathring{\varphi}_{\PSdataamplitude_1;\PSdataamplitude_2}(u) 
		\eqdef \PSdataamplitude_1 \mathring{\varphi}(\PSdataamplitude_2 u)$,
		which would allow one to produce small-amplitude shock-forming solutions 
		for initial data with small derivatives
		(i.e., when $\PSdataamplitude_1$ is small and $\PSdataamplitude_1 \PSdataamplitude_2$ is even smaller)  
		or large derivatives 
		(i.e., when $\PSdataamplitude_1$ is small and $\PSdataamplitude_1 \PSdataamplitude_2$ is large).
		Moreover, the assumption that $\mathring{\varphi}$ is compactly supported can easily be eliminated.
	\item One could also prove an analog of Theorem~\ref{T:PROPERTIESOFADMISSIBLESOLUTIONSLAUNCEDBYRESCALED}
		for isentropic plane-symmetric solutions that are not simple,
		i.e., when both Riemann invariants in the system
		\eqref{AE:COMPRESSIBLEEULEREQUATIONSFORRIEMANNINVARIANTS} are non-vanishing.
		In this way, one could produce shocks along the characteristics of $\LunitPS$
		(i.e., the blowup of $\partial_1 \RRiemann$, as in Theorem~\ref{T:PROPERTIESOFADMISSIBLESOLUTIONSLAUNCEDBYRESCALED})
		as well as along the characteristics of $\uLunitPS$
		(i.e., the blowup of $\partial_1 \LRiemann$).
	\item In Theorem~\ref{T:PROPERTIESOFADMISSIBLESOLUTIONSLAUNCEDBYRESCALED}, 
		we followed the solution up to the crease $\creasePS$. 
		With minor additional effort, 
		we could have followed the solution up to a portion of the singular boundary that contains a neighborhood of the  
		crease. We have omitted such results because we have already derived them away from plane-symmetry, in 
		Theorem~\ref{T:DEVELOPMENTANDSTRUCTUREOFSINGULARBOUNDARY}.
\end{itemize}
\end{remark}

\begin{center}
	\begin{figure}  
		\begin{overpic}[scale=.6, grid = false, tics=5, trim=-.5cm -1cm -1cm -.5cm, clip]{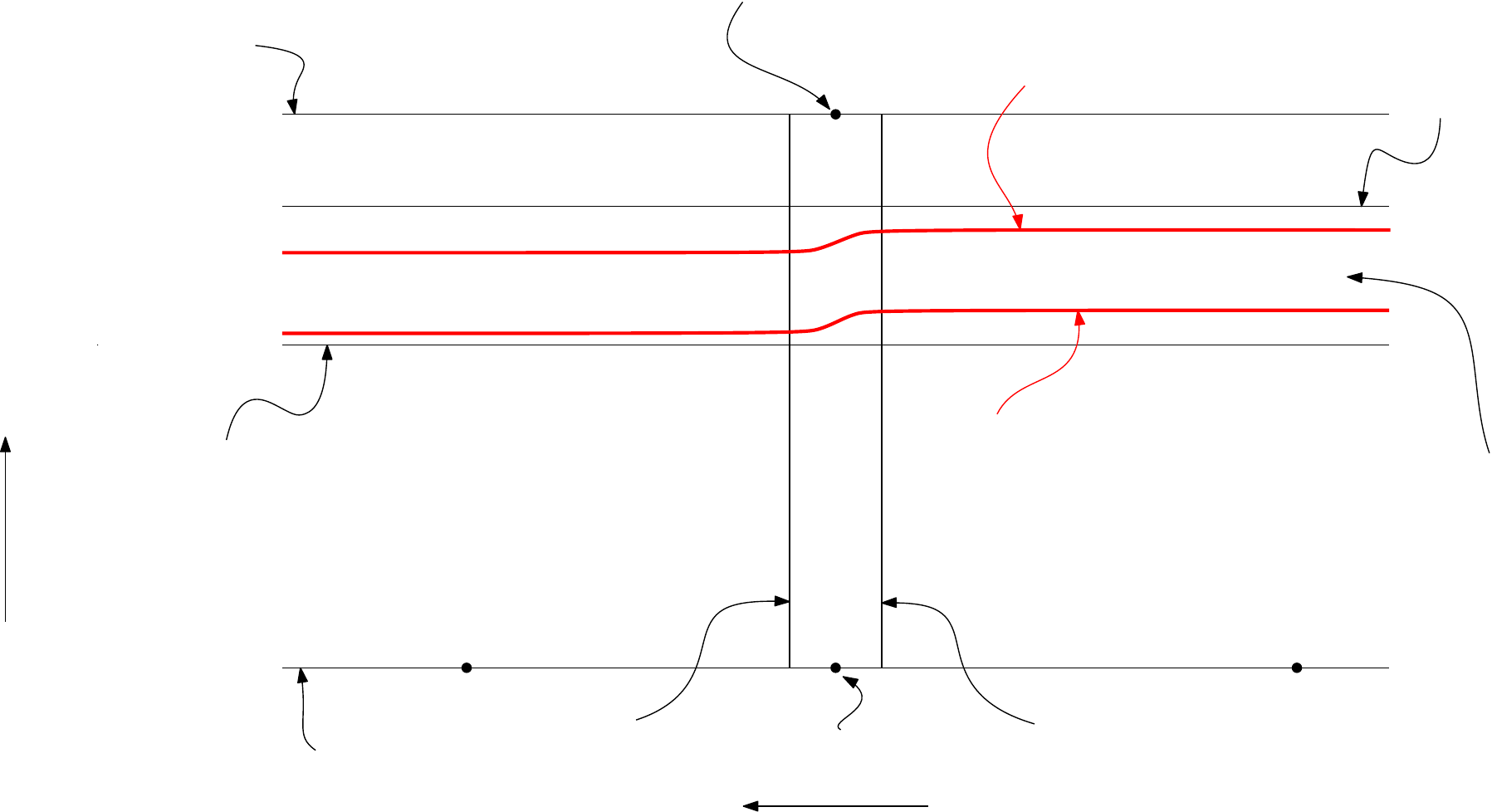}
			\put (22,6) {$\Sigma_0$}
			\put (85,25) {$\twoargMrough{[2 \timefunction_0,\frac{1}{2} \timefunction_0],\R}{\muxmulevelsetvalue}$}
			\put (58,26) {$\twoarghypPS{2\timefunction_0}{\muxmulevelsetvalue}$}
			\put (1,20) {$t$}
			\put (61,51) {$\twoarghypPS{\frac{1}{2}\timefunction_0}{\muxmulevelsetvalue}$}
			\put (39,56) {$\creasePS = (\blowuptimePS, 0)$}
			\put (91,49) {$\Sigma_{\blowuptimePS - 2\PSBigDelta}$}
			\put (13,55) {$\Sigma_{\blowuptimePS}$}
			\put (14,25) {$\Sigma_{\blowuptimePS - 2 \PSinterestingusmallmultipleofamplitude \blowuptimePS}$}
			\put (28,15) {$u = \leftu$}
			\put (79,15) {$u = - \rightu$}
			\put (36,6) {$\nullhyparg{\interestingu}^{[0,\blowuptimePS]}$}
			\put (66,7) {$\nullhyparg{-\interestingu}^{[0,\blowuptimePS]}$}
			\put (54,3) {$u$}
			\put (52,8) {$u=0$}
		\end{overpic}
	\end{figure}
\end{center}

\begin{proof}[Proof of Theorem~\ref{T:PROPERTIESOFADMISSIBLESOLUTIONSLAUNCEDBYRESCALED}]
In each step of the proof, we will silently adjust the smallness of the positive constants
$\PSdataamplitude_0$ and $\PSinterestingusmallmultipleofamplitude$
without explicitly mentioning it each time.

\medskip

\noindent \textbf{Proof of classical existence with respect to the geometric coordinates}: The facts that $\RRiemannPS = \RRiemannPS(u)$ and that $\RRiemannPS$
vanishes on the complement of the region 
$\lbrace (t,u) \in \mathbb{R} \times \mathbb{R} \ | \ u \in [- \rightu,\leftu] \rbrace$
follow from the evolution equation $\geop{t} \RRiemannPS(t,u) = 0$
(see \eqref{AE:COMPRESSIBLEEULEREQUATIONSFORRIEMANNINVARIANTS} and \eqref{AE:PLANESYMMETRYCOMMUTATORVECTORFIELDSAREGEOCOORDINATEVECTORFIELDS})
and our assumption that the initial data are supported
in the $u$-interval $[- \rightu,\leftu]$. The identity \eqref{AE:PSBLOWUPTIMEINTERMSOFAMPITUDE} follows from definition \eqref{AE:PSDELTASTARDEF},
\eqref{AE:USEFULIDENTITYFORPSLMUSOURCETERMFUNCTION}--\eqref{AE:DERIVATIVEOFUSEFULIDENTITYFORPSLMUSOURCETERMFUNCTION},
and the properties of $\antiderivativePSLmusourcetermfunction[(\dataRRiemannPS)_{\PSdataamplitude}(u)]$
described in Sect.\,\ref{SS:EXAMPLESOFBONAFIDEDATALEADINGTOADMISSIBLE}.
The properties of $\RRiemannPS$, $\LunitPS$, and $\muPS$ follow easily from
our assumptions on the support of $\mathring{\varphi}$,
the fact that relative to the geometric coordinates, $\RRiemannPS$ depends only on $u$,
equations \eqref{AE:RIEMANNTOWAVE} and \eqref{AE:NULLVECTORFIELDSPLANESYMMETRY},
the normalization assumption \eqref{E:BACKGROUNDSOUNDSPEEDISUNITY}
(which implies that $\SpeedPS = 1 + \smoothfunction(\RRiemannPS)$, where $\smoothfunction$ is smooth),
and the explicit solution formulas provided by Cor.\,\ref{AC:PSEXPLICITEXPRESSIONSFORSOLUTION}.

\medskip

\noindent \textbf{Proof of 
\eqref{AE:PSSPEEDOFSOUNDESTIMATE}--\eqref{AE:LINFINITYMUANDTRANSVERSALDERIVATIVES}:}
These estimates are straightforward consequences of the fact that in geometric coordinates, 
$\RRiemannPS$ depends only on $u$,
\eqref{AE:NULLVECTORFIELDSPLANESYMMETRY},
\eqref{AE:PLANESYMMETRYCOMMUTATORVECTORFIELDSAREGEOCOORDINATEVECTORFIELDS},
\eqref{AE:SPEEDTIMESMUPLANESYMMETRYWITHEXPLICITSOURCE},
the normalization assumption \eqref{E:BACKGROUNDSOUNDSPEEDISUNITY}
(which implies that $\SpeedPS = 1 + \smoothfunction(\RRiemannPS)$, where $\smoothfunction$ is smooth),
the data estimates 
\eqref{AE:PSRESCALEDDATALINFINITYBOUND}--\eqref{AE:PSRESCALEDDATASOBOLVBOUND},
and standard Sobolev calculus.

\medskip

\noindent \textbf{Proof \eqref{AE:PSSIZEOFCARTESIANX1}:}
Note that by \eqref{AE:NULLVECTORFIELDSPLANESYMMETRY}, 
$\geop{t} x^1 = \LunitPS x^1 = \velocityPS + \SpeedPS$.
Considering also that $\RRiemannPS$ vanishes on the complement of
$\lbrace (t,u) \in \mathbb{R} \times \mathbb{R} \ | \ u \in [- \rightu,\leftu] \rbrace$,
we see that for
$(t,u) \in  [0,\blowuptimePS] \times \left((-\infty,-\rightu] \cup [\leftu,\infty) \right)$, 
we have $\velocityPS(t,u) + \SpeedPS(t,u) = 1$.
Recalling also that $x^1(0,u) = - u$ (see \eqref{AE:DATAFOREIKONALEQUATIONINPLANESYMMETRY}), 
we see that for $t \in [0,\blowuptimePS]$,
we have $x^1(t,-\rightu) = \rightu + t$ 
and $x^1(t,\leftu) = - \leftu + t$.
Moreover, since \eqref{AE:IDENTITYNULLCVECTORFIELDSINPLANESYMMETRY} and 
\eqref{AE:PLANESYMMETRYCOMMUTATORVECTORFIELDSAREGEOCOORDINATEVECTORFIELDS} 
imply that
$\geop{u} x^1 = \muXPS x^1 = - \SpeedPS \muPS \leq 0$,
we conclude that for $(t,u) \in [0,\blowuptimePS] \times [-\rightu,\leftu]$,
we have $x^1(t,\leftu) \leq x^1(t,u) \leq x^1(t,-\rightu)$.
Combining these results, we conclude \eqref{AE:PSSIZEOFCARTESIANX1}.

\medskip

\noindent \textbf{Proof of \eqref{AE:PSIDENTITYFORDIFFERENTIALOFGEOTOCARTESIANCHOVMAP}
and \eqref{AE:PSHIGHORDERESTIMATESFORUPSILON}--\eqref{AE:PSSOBOLEVESTIMATESFORGEOPTUPSILON}:}
Since $\PSUpsilon(t,u) = (t,x^1)$,
\eqref{AE:PSIDENTITYFORDIFFERENTIALOFGEOTOCARTESIANCHOVMAP}
follows easily from 
\eqref{AE:NULLVECTORFIELDSPLANESYMMETRY},
\eqref{AE:IDENTITYNULLCVECTORFIELDSINPLANESYMMETRY},
and
\eqref{AE:PLANESYMMETRYCOMMUTATORVECTORFIELDSAREGEOCOORDINATEVECTORFIELDS}.

The estimates \eqref{AE:PSHIGHORDERESTIMATESFORUPSILON}--\eqref{AE:PSSOBOLEVESTIMATESFORGEOPTUPSILON}
follow from \eqref{AE:PSIDENTITYFORDIFFERENTIALOFGEOTOCARTESIANCHOVMAP},
the fact that relative to the geometric coordinates, $\RRiemannPS$ depends only on $u$,
the estimates \eqref{AE:PSSOBOLEVBOUNDSFORRRIEMANN}--\eqref{AE:PSSOBOLEVBOUNDSFORINVERSEFOLIATIONDENSITY},
and standard Sobolev calculus.

\medskip

\noindent \textbf{Proof of the remaining properties of $\PSUpsilon$:}
We now prove that the map $\PSUpsilon(t,u) = (t,x^1)$ is a diffeomorphism on $[0,\blowuptimePS) \times \mathbb{R}$.
First, using \eqref{AE:PSIDENTITYFORDIFFERENTIALOFGEOTOCARTESIANCHOVMAP},
we compute that $\mbox{\upshape det} \mathrm{d} \PSUpsilon = - \SpeedPS \muPS$. 
Also using \eqref{AE:PSSPEEDOFSOUNDESTIMATE} and the fact that
$\muPS > 0$ on $[0,\blowuptimePS) \times \mathbb{R}$,
we see that $\mbox{\upshape det} \mathrm{d} \PSUpsilon < 0$ on $[0,\blowuptimePS) \times \mathbb{R}$,
and since $\geop{u} x^1 = - \SpeedPS \muPS$ (by \eqref{AE:PSIDENTITYFORDIFFERENTIALOFGEOTOCARTESIANCHOVMAP}),
we see that for $t \in [0,\blowuptimePS)$,
the map $u \rightarrow x^1(t,u)$ is strictly decreasing for $u \in \mathbb{R}$.
We therefore find that $\PSUpsilon$ is injective on $[0,\blowuptimePS) \times \mathbb{R}$
and that it is a diffeomorphism on the same domain.
Moreover, since $\muPS(\blowuptimePS,u)$ vanishes only at the origin $u=0$
(where the crease is located),
the map $u \rightarrow x^1(\blowuptimePS,u)$ is strictly decreasing for $u \in \mathbb{R}$.
It follows that $\PSUpsilon$ is a homeomorphism on $[0,\blowuptimePS] \times \mathbb{R}$,
as is desired.

\medskip

\noindent \textbf{Proof of \eqref{AE:BOUNDSONLMUINTERESTINGREGION},
\eqref{AE:PSMUTAYLOREXPANSIONININTERESTINGREGION}--\eqref{AE:PSMUSECONDDERIVATIVETAYLOREXPANSIONININTERESTINGREGION},
\eqref{AE:PSLOCATIONOFXMUEQUALSMINUSKAPPA}--\eqref{AE:PSREGIONWHEREMUXMUKAPPALEVELSETISNOTLOCATED},
and \eqref{AE:PSMUTRANSVERSALCONVEXITY}:}
These estimates are straightforward to obtain via standard Sobolev calculus,
Taylor expansions,
the identities \eqref{AE:LSPEEDTIMESMUPLANESYMMETRYWITHEXPLICITSOURCE}--\eqref{AE:SPEEDTIMESMUPLANESYMMETRYWITHEXPLICITSOURCE},
the properties of the function $\antiderivativePSLmusourcetermfunction[(\dataRRiemannPS)_{\PSdataamplitude}(u)] = \mathring{\varphi}_{\PSdataamplitude}(u)$
stated in Sect.\,\ref{SS:EXAMPLESOFBONAFIDEDATALEADINGTOADMISSIBLE},
the estimates \eqref{AE:PSSPEEDOFSOUNDESTIMATE}--\eqref{AE:LINFINITYMUANDTRANSVERSALDERIVATIVES},
and the definition \eqref{AE:PSKAPPSOISMULTIPLEOFTIMEFUNCTION0TIMESPSIAMPLITUDE} of $\muxmulevelsetvalue_0$.

\medskip
\noindent \textbf{Proof of \eqref{AE:PSBLOWUPLOWERBOUND}:}
We first use \eqref{AE:CLOSEDPARTIAL1DERIVATIVEOFRPLUSPS},
\eqref{E:TAYLOREXPANSIONOFKEYMUTERMABOUTCREASE},
\eqref{AE:PLANESYMMETRYNONDEGENCONDITION},
and \eqref{AE:PSSMALLAMPLITUDESOLUTIONLINFINITYRRIEMANNTRANSVERSALDERIVATIVES}
to deduce that for $(t,u) \in [\blowuptimePS - 2 \PSBigDelta,\blowuptimePS]  \times [-\interestingu,\interestingu]$,
we have:
\begin{align} \label{AE:PROOFSTEPPSBLOWUPLOWERBOUND}
	[\partial_1 \RRiemannPS](t,u)
	& = 
	- 
	\frac{\frac{d}{du} (\dataRRiemannPS)_{\PSdataamplitude}(u)}{
	1 
	+
	t
	\antiderivativePSLmusourcetermfunction'[(\dataRRiemannPS)_{\PSdataamplitude}(u)]
	\frac{d}{du}
	(\dataRRiemannPS)_{\PSdataamplitude}(u)}
	=
	\frac{2 + \mathcal{O}(\PSdataamplitude)}{\frac{\backgroundSpeedprimePS}{\backgroundSpeedPS} + 1}
	\cdot
	\frac{\PSdataamplitude}{1 - t \PSdataamplitude \left\lbrace 1 - \frac{1}{2} \left[\PSdatamuHessianTaylorcoefficient + \mathcal{O}(\PSdataamplitude) \right] u^2 \right\rbrace}.
\end{align}
From \eqref{AE:PROOFSTEPPSBLOWUPLOWERBOUND}
and the fact that the crease has geometric coordinates $(t,u) = (\PSdataamplitude^{-1},0)$,
we conclude \eqref{AE:PSBLOWUPLOWERBOUND}.

\medskip
\noindent \textbf{Proof of \eqref{AE:PSMUISLARGEAWAYFROMINTERESTINGREGION}:}
This estimate follows from \eqref{E:RESCALEDFIRSTDERIVATIVEOFMUFUNCTIONUNIFORMLYPOSITIVEAWAYFROMINTERESTINGREGION},
\eqref{AE:SPEEDTIMESMUPLANESYMMETRYWITHEXPLICITSOURCE},
\eqref{AE:PSSPEEDOFSOUNDESTIMATE},
and \eqref{AE:PSBLOWUPTIME}.

\medskip

\noindent \textbf{An intermediate step - the map $(t,u) \rightarrow (\muPS,\muXPS \muPS)$ is a local diffeomorphism:}
Using 
\eqref{AE:PLANESYMMETRYCOMMUTATORVECTORFIELDSAREGEOCOORDINATEVECTORFIELDS},
\eqref{AE:PSTIMEFUNCTION0CONSTRAINTS},
\eqref{AE:PSINTERESTINGUISSMALLMULTIPLEOFAMPLITUDE},
\eqref{AE:PSBLOWUPTIME},
and
\eqref{AE:PSLUNTMUTAYLOREXPANSIONININTERESTINGREGION}--\eqref{AE:PSMUSECONDDERIVATIVETAYLOREXPANSIONININTERESTINGREGION},
we compute that for
$
(t,u) 
\in
[\blowuptimePS - 2 \PSinterestingusmallmultipleofamplitude \blowuptimePS,\blowuptimePS] 
			\times [- \interestingu, \interestingu]
$,
we have: 
\begin{align}	 \label{AE:PSEXPANSIONOFMAPFROMGEOMETRICTOMUXUMCOORDINATES}
	\frac{\partial(\muPS,\muXPS \muPS)}{\partial (t,u)}
	= \begin{pmatrix}
			- \PSdataamplitude & 0 
				\\
		0 & \PSdatamuHessianTaylorcoefficient
		\end{pmatrix}
		+
		\begin{pmatrix}
			\mathcal{O}(\PSdataamplitude^2) & \mathcal{O}(\PSdataamplitude) 
				\\
		 \mathcal{O}(\PSdataamplitude^2) & \mathcal{O}(\PSdataamplitude) + \mathcal{O}(\PSinterestingusmallmultipleofamplitude)
		\end{pmatrix}.
\end{align}
Using \eqref{AE:PSEXPANSIONOFMAPFROMGEOMETRICTOMUXUMCOORDINATES}, we deduce that if
	$
(t,u) 
\in
[\blowuptimePS - 2 \PSinterestingusmallmultipleofamplitude \blowuptimePS,\blowuptimePS] 
			\times [-\interestingu , \interestingu]
$,
and 
if $\PSdataamplitude$ 
and
$\PSinterestingusmallmultipleofamplitude$ are sufficiently small, 
then
$
\frac{\partial(\muPS,\muXPS \muPS)}{\partial (t,u)}
$
is invertible,
and moreover, that:
\begin{align} \label{AE:PSKEYESTIMATEFORINTERTIBILITYOFTUTOMUMUXMUMAP}
\max_{(t_1,u_1),
	\,
(t_2,u_2)
\in
[\blowuptimePS - 2 \PSinterestingusmallmultipleofamplitude \blowuptimePS,\blowuptimePS] 
			\times [-\interestingu , \interestingu]}
\left|
\left(\frac{\partial(\muPS,\muXPS \muPS)}{\partial (t,u)}|_{(t_1,u_1)} \right)^{-1}
\frac{\partial(\muPS,\muXPS \muPS)}{\partial(t,u)}|_{(t_2,u_2)}
-
\mbox{\upshape ID} 
\right|_{\mbox{\upshape}Euc}
& \leq \frac{1}{4},
\end{align}
where  
$|\cdot|_{\mbox{\upshape}Euc}$ is the standard Frobenius norm on matrices
(equal to the square root of the sum of the squares of the matrix entries)
and $\mbox{\upshape ID}$ denotes the $2 \times 2$ identity matrix.
These estimates, 
together with \eqref{AE:LLUPS}
and
\eqref{AE:LINFINITYINITIALROUGHHYPERSURFACELDERIVATIVEOFMUANDTRANSVERSALDERIVATIVES}--\eqref{AE:LINFINITYMUANDTRANSVERSALDERIVATIVES},
imply that the map $(t,u) \rightarrow (\muPS,\muXPS \muPS)$
is a $C^2$ diffeomorphism from the compact, convex set
$
[\blowuptimePS - 2 \PSinterestingusmallmultipleofamplitude \blowuptimePS,\blowuptimePS] 
			\times [-\interestingu,\interestingu ]
$
onto its image, which we denote by $\mathcal{I}$,
i.e.,
\begin{align} \label{AE:PSDESCRPTIONOFPORTIONSOFLEVELSETSOFMUXEQUALSMINUSKAPPAMUFOLIATIONS}
	\mathcal{I}
	& \eqdef
		\left\lbrace
			\left(\muPS(t,u), \muXPS \muPS(t,u) \right)
			\ | \ (t,u) \in 
			[\blowuptimePS - 2 \PSinterestingusmallmultipleofamplitude \blowuptimePS,\blowuptimePS] 
			\times 
			[-\interestingu,\interestingu]
		\right\rbrace.
\end{align}
Let $\widetilde{\mathcal{I}}$ be the following subset of $\mathcal{I}$, where
$\PSBigDelta$ is defined in \eqref{AE:BIGDELTAPS}:
\begin{align} \label{AE:PSSUBSETDESCRPTIONOFPORTIONSOFLEVELSETSOFMUXEQUALSMINUSKAPPAMUFOLIATIONS}
\widetilde{\mathcal{I}}
	& \eqdef
		\left\lbrace
			\left(\muPS(t,u), \muXPS \muPS(t,u) \right)
			\ | \ (t,u) \in 
			[\blowuptimePS - \PSinterestingusmallmultipleofamplitude \blowuptimePS,\blowuptimePS - 3 \PSBigDelta] 
			\times 
			[-\interestingu,\interestingu]
		\right\rbrace.
\end{align}

\medskip

\noindent \textbf{Proof of the properties of the $\PStimefunctionarg{\muxmulevelsetvalue}$ and the location of their level sets:}
Let $\widetilde{\mathcal{I}}$ 
be as in \eqref{AE:PSSUBSETDESCRPTIONOFPORTIONSOFLEVELSETSOFMUXEQUALSMINUSKAPPAMUFOLIATIONS}.
Considering 
\eqref{AE:PSMUTAYLOREXPANSIONININTERESTINGREGION} 
and
\eqref{AE:PSMUFIRSTDERIVATIVETAYLOREXPANSIONININTERESTINGREGION},
recalling that $\timefunction_0$ is allowed to be any negative real number satisfying
\eqref{AE:PSTIMEFUNCTION0CONSTRAINTS},
that $\muxmulevelsetvalue_0$ is defined by \eqref{AE:PSKAPPSOISMULTIPLEOFTIMEFUNCTION0TIMESPSIAMPLITUDE},
and that $\PSBigDelta$ is defined in \eqref{AE:BIGDELTAPS},
we see that 
if $\PSdataamplitude$ 
and
$\PSinterestingusmallmultipleofamplitude$ 
are sufficiently small,
then:
\begin{align} \label{AE:ESTIMATESFORPSDESCRPTIONOFPORTIONSOFLEVELSETSOFMUXEQUALSMINUSKAPPAMUFOLIATIONS}
\begin{split}
\left[\frac{1}{2}|\timefunction_0|,2 |\timefunction_0| \right]
\times
[- \muxmulevelsetvalue_0,0]
&
\subset
\left[\frac{|\timefunction_0|^2 \PSdataamplitude^2}{64 \PSdatamuHessianTaylorcoefficient^3} + \frac{1}{4}|\timefunction_0|,
\frac{1}{2} \PSinterestingusmallmultipleofamplitude \right]
\times
[- \muxmulevelsetvalue_0,0]
=
\left[\frac{4 \muxmulevelsetvalue_0^2}{\PSdatamuHessianTaylorcoefficient} + \frac{1}{4}|\timefunction_0|,
\frac{1}{2} \PSinterestingusmallmultipleofamplitude \right]
\times
[- \muxmulevelsetvalue_0,0]
	\\
& \subset
\bigcup_{\muxmulevelsetvalue \in [-\frac{\PSinterestingusmallmultipleofamplitude \PSdataamplitude}{2 \PSdatamuHessianTaylorcoefficient},\frac{\PSinterestingusmallmultipleofamplitude \PSdataamplitude}{2 \PSdatamuHessianTaylorcoefficient}]}
\left[\frac{4 \muxmulevelsetvalue^2}{\PSdatamuHessianTaylorcoefficient} + \frac{1}{4}|\timefunction_0|,
\frac{1}{2} \PSinterestingusmallmultipleofamplitude 
\right]
\times
\lbrace - \muxmulevelsetvalue \rbrace	
\subset
\widetilde{\mathcal{I}}.
\end{split}
\end{align}
In particular, \eqref{AE:ESTIMATESFORPSDESCRPTIONOFPORTIONSOFLEVELSETSOFMUXEQUALSMINUSKAPPAMUFOLIATIONS} shows
that for $\muxmulevelsetvalue \in [0,\muxmulevelsetvalue_0]$,
along the portion of the level set $\lbrace (t,u) \ | \ \muXPS \muPS(t,u) = - \muxmulevelsetvalue \rbrace$
that is contained in 
$
[\blowuptimePS - \PSinterestingusmallmultipleofamplitude \blowuptimePS,\blowuptimePS - 3 \PSBigDelta] 
			\times [-\interestingu,\interestingu]
$,
$\muPS$ ranges over an interval that contains
$
\left[\frac{1}{2}|\timefunction_0|,2 |\timefunction_0| \right]
$.

We now study the initial value problem 
for the rough time function $\PStimefunctionarg{\muxmulevelsetvalue}$ (see Definition~\ref{D:ROUGHTIMEFUNCTION}).
Fix $\muxmulevelsetvalue \in [0,\muxmulevelsetvalue_0]$, where $\muxmulevelsetvalue_0$ is defined by \eqref{AE:PSKAPPSOISMULTIPLEOFTIMEFUNCTION0TIMESPSIAMPLITUDE}.
From the diffeomorphism properties of the map 
$(t,u) \rightarrow (\muPS,\muXPS \muPS)$ established above,
\eqref{AE:PSDESCRPTIONOFPORTIONSOFLEVELSETSOFMUXEQUALSMINUSKAPPAMUFOLIATIONS}, and
\eqref{AE:ESTIMATESFORPSDESCRPTIONOFPORTIONSOFLEVELSETSOFMUXEQUALSMINUSKAPPAMUFOLIATIONS},
we see that for every $\timefunction \in [2 \timefunction_0,\frac{1}{2} \timefunction_0]$, 
there exists a unique point
$
q_{\timefunction} 
\in
\lbrace (t,u) \ | \ \muXPS \muPS(t,u) = - \muxmulevelsetvalue \rbrace
\cap
[\blowuptimePS - \PSinterestingusmallmultipleofamplitude \blowuptimePS,\blowuptimePS - 3 \PSBigDelta] 
\times [-\interestingu,\interestingu]
$
such that $\muPS(q_{\timefunction)} = - \timefunction$
and such that the map $\timefunction \rightarrow (t_{q_{\timefunction}},u_{q_{\timefunction}})$
is $C^2$, where
$(t_{q_{\timefunction}},u_{q_{\timefunction}})$ are the geometric coordinates of $q_{\timefunction}$.
Note that $q_{\timefunction}$ is a point on the ``initial'' data hypersurface for $\PStimefunctionarg{\muxmulevelsetvalue}(t,u)$
(see Definition~\ref{D:ROUGHTIMEFUNCTION}).
In view of Definition~\ref{D:ROUGHTIMEFUNCTION}, we see that
$\PStimefunctionarg{\muxmulevelsetvalue}(t_{q_{\timefunction}},u_{q_{\timefunction}}) = - \muPS(t_{q_{\timefunction}},u_{q_{\timefunction}}) = \timefunction$,
i.e., the initial value of $\PStimefunctionarg{\muxmulevelsetvalue}$ at $q_{\timefunction}$ is $\timefunction$.
Much like in the bulk of the paper, we will use the notation
$\PSdatahypfortimefunctiontwoarg{- \muxmulevelsetvalue}{[2 \timefunction_0,\frac{1}{2} \timefunction_0]}$
to denote the union of the points $q_{\timefunction}$ as $\timefunction$ 
varies over the interval $[2 \timefunction_0,\frac{1}{2} \timefunction_0]$,
i.e., 
\begin{align} \label{AE:PSRELEVANTDATAHYPERSURFACEPORTIONFORROUGHTIMEFUNCTION}
\begin{split}
\PSdatahypfortimefunctiontwoarg{- \muxmulevelsetvalue}{[2 \timefunction_0,\frac{1}{2} \timefunction_0]}
\eqdef
&
\lbrace (t,u) \ | - \frac{1}{2} \timefunction_0  \leq \muPS(t,u) \leq - 2 \timefunction_0 \rbrace 
\cap 
\lbrace (t,u) \ | \ \muXPS \muPS(t,u) = - \muxmulevelsetvalue \rbrace
	\\
& 
\cap
\left([\blowuptimePS - \PSinterestingusmallmultipleofamplitude \blowuptimePS,\blowuptimePS] \times [-\interestingu,\interestingu] \right).
\end{split}
\end{align}

Let $\upgamma_{q_{\timefunction}} : \mathbb{R} \rightarrow \mathbb{R} \times \mathbb{R}$ 
be the $u$-parameterized integral curve of $\PSWtransarg{\muxmulevelsetvalue} = \roughgeop{u}$ 
(recall that $\roughgeop{u} u = 1$)
that emanates from $q_{\timefunction}$,
where the target is geometric coordinate space, 
i.e., 
$\upgamma_{q_{\timefunction}}(u_{q_{\timefunction}}) = (t_{q_{\timefunction}},u_{q_{\timefunction}})$, where
$\upgamma_{q_{\timefunction}}(u)$ belongs to geometric coordinate space for $u$ belonging to the interval of existence of 
$\upgamma_{q_{\timefunction}}$.
Below we will show that for $\timefunction \in [2 \timefunction_0,\frac{1}{2} \timefunction_0]$,
we have:
\begin{align} \label{AE:PSINTEGRALCURVESOFROUGHGEOPCONTAINEDINFLATDEVELOPMENT}
\upgamma_{q_{\timefunction}}(\mathbb{R}) 
& \subset [\blowuptimePS - 2\PSinterestingusmallmultipleofamplitude \blowuptimePS,\blowuptimePS - 2 \PSBigDelta] \times \mathbb{R},
\end{align}
which shows in particular that the entire integral curve
is contained in the region of classical existence with respect to the geometric coordinates
and is temporally separated from the Cartesian time of first blowup by at least $2 \PSBigDelta$.
Moreover, \eqref{AE:PSMUTRANSVERSALCONVEXITY} implies that these integral curves
are transversal to 
$
\PSdatahypfortimefunctiontwoarg{- \muxmulevelsetvalue}{[2 \timefunction_0,\frac{1}{2} \timefunction_0]}
$.
It follows that the map $(\timefunction,u) \rightarrow (t,u)$,
is an injection from $[2 \timefunction_0,\frac{1}{2} \timefunction_0] \times \mathbb{R}$
onto a subset of $[0,\blowuptimePS] \times \mathbb{R}$,
where the image component function $t = \upgamma_{q_{\timefunction}}^0(u)$ 
is defined to be the Cartesian time coordinate of the point $\upgamma_{q_{\timefunction}}(u)$.
Considering also that \eqref{AE:PSTRANSPORTEQUATIONFORROUGHTIMEFUNCTION} implies
$\PStimefunctionarg{\muxmulevelsetvalue}$ is constant along the integral curves $u \rightarrow \upgamma_{q_{\timefunction}}(u)$,
we see that the map $(\timefunction,u) \rightarrow (t,u)$ is precisely the map from rough adapted coordinates
to geometric coordinates, and that its inverse is the map
$
	\PSCHOVgeotorough{\muxmulevelsetvalue}(t,u) 
	= 
	(\PStimefunctionarg{\muxmulevelsetvalue},u)
$
from \eqref{AE:PSCHOVGEOTOROUGH}.
In addition,
from this reasoning and \eqref{AE:PSINTEGRALCURVESOFROUGHGEOPCONTAINEDINFLATDEVELOPMENT},
we also conclude \eqref{AE:ROUGHHYPSFORADMISSIBLESOLUTIONSCONTAINEDINCAUCHYSTABILITYREGION}.

We now prove 
\eqref{AE:PSINTEGRALCURVESOFROUGHGEOPCONTAINEDINFLATDEVELOPMENT}.
First, using
\eqref{AE:PSSUBSETDESCRPTIONOFPORTIONSOFLEVELSETSOFMUXEQUALSMINUSKAPPAMUFOLIATIONS}--\eqref{AE:ESTIMATESFORPSDESCRPTIONOFPORTIONSOFLEVELSETSOFMUXEQUALSMINUSKAPPAMUFOLIATIONS},
we deduce that
for $\timefunction \in [2 \timefunction_0,\frac{1}{2} \timefunction_0]$, 
we have:
\begin{align} \label{E:PSINITIALCARTESIANTIMEVALUEALONGXMUEQUALSMINUSKAPPAISAWAYFROMSHOCK}
	t_{q_{\timefunction}} \in 
	[\blowuptimePS - \PSinterestingusmallmultipleofamplitude \blowuptimePS,\blowuptimePS - 3 \PSBigDelta].
\end{align}

Next, in view of \eqref{AE:PSROUGHADAPTEDUDERIVTATIVEINTERMSOFGEOMETRICPARTIALDERIVATIVES}, 
we see that if $0 \leq \muxmulevelsetvalue \leq \muxmulevelsetvalue_0$,
then as $u$ varies over $\mathbb{R}$, 
we can bound the total change in the Cartesian time coordinate $t$ along $\upgamma_{q_{\timefunction}}$, denoted by
$\Delta_{\upgamma_{q_{\timefunction}}}^0$, as follows: 
\begin{align}  \label{AE:ESTIMATEFORCARTESIANTDISPLACEMENTROUGHHYPPS}
\begin{split}	
	|\Delta_{\upgamma_{q_{\timefunction}}}^0|
	& \leq 
	\int_{u' = - \infty}^{\infty} 
		\left| 
			\roughgeop{u} t  
		\right|
	\, \mathrm{d} u' 
	= \int_{u' = - \infty}^{\infty} 
			\left|
				\frac{\muxmulevelsetvalue \phi}{\geop{t} \muPS} 
			\right|
		\, \mathrm{d} u'
	= \int_{u' = - \interestingu}^{\interestingu}  
			\left|
				\frac{\muxmulevelsetvalue \phi}{\geop{t} \muPS} 
			\right|
			\, \mathrm{d} u'
	\leq
	2
	\int_{|u'| \leq \frac{\PSdataamplitude}{b^2}}
		\frac{\muxmulevelsetvalue}{\PSdataamplitude}
	\, \mathrm{d} u'
		\\
	& = \frac{4 \muxmulevelsetvalue}{b^2} 
		\leq
		 \frac{4 \muxmulevelsetvalue_0}{b^2} 
		= 
	\frac{\PSdataamplitude |\timefunction_0|}{4 \PSdatamuHessianTaylorcoefficient^3}, 
\end{split}
\end{align}
where to obtain the second ``$=$'' and the next-to-last ``$\leq$'' on RHS~\eqref{AE:ESTIMATEFORCARTESIANTDISPLACEMENTROUGHHYPPS},
we used the properties of $\phi$ from Definition~\ref{D:WTRANSANDCUTOFF},
the estimate \eqref{AE:BOUNDSONLMUINTERESTINGREGION} for $\geop{t} \muPS = \LunitPS \muPS$,
and the definition \eqref{AE:PSINTERESTINGUISSMALLMULTIPLEOFAMPLITUDE} of $\interestingu$,
and to obtain the last ``$=$'' on RHS~\eqref{AE:ESTIMATEFORCARTESIANTDISPLACEMENTROUGHHYPPS},
we used definition \eqref{AE:PSKAPPSOISMULTIPLEOFTIMEFUNCTION0TIMESPSIAMPLITUDE}.
From \eqref{E:PSINITIALCARTESIANTIMEVALUEALONGXMUEQUALSMINUSKAPPAISAWAYFROMSHOCK}
and
\eqref{AE:ESTIMATEFORCARTESIANTDISPLACEMENTROUGHHYPPS},
we see that for any point 
$
\upgamma_{q_{\timefunction}}(u) \in \twoarghypPS{\timefunction}{\muxmulevelsetvalue}
$,
we can bound its Cartesian time coordinate $t$ as follows: 
\begin{align} \label{AE:ESTIMATEFORCARTESIANTIMEONROUGHHYPERSURFACES}
 \begin{split} 
	t
	& = 
	\upgamma_{q_{\timefunction}}^0(u)
	\in 
	[\blowuptimePS - \PSinterestingusmallmultipleofamplitude \blowuptimePS - |\Delta_{\upgamma_{q_{\timefunction}}}^0|,
	\blowuptimePS - 3 \PSBigDelta + |\Delta_{\upgamma_{q_{\timefunction}}}^0|]
		\\
	& 
	\subset
	[\blowuptimePS - \PSinterestingusmallmultipleofamplitude \blowuptimePS - \frac{\PSdataamplitude |\timefunction_0|}{4 \PSdatamuHessianTaylorcoefficient^3}|,
	\blowuptimePS - 3 \PSBigDelta + \frac{\PSdataamplitude |\timefunction_0|}{4 \PSdatamuHessianTaylorcoefficient^3}].
\end{split}
\end{align}
From
\eqref{AE:BIGDELTAPS}
and
\eqref{AE:ESTIMATEFORCARTESIANTIMEONROUGHHYPERSURFACES},
we conclude that
if $\PSdataamplitude$ 
and
$\PSinterestingusmallmultipleofamplitude$ 
are sufficiently small,
then the desired result \eqref{AE:PSINTEGRALCURVESOFROUGHGEOPCONTAINEDINFLATDEVELOPMENT} holds.

Next we prove \eqref{AE:PSGEOPTDERIVATIVEOFROUGHTIMEFUNCTIONKEYBOUND}--\eqref{AE:PSGEOPUDERIVATIVEOFROUGHTIMEFUNCTIONSIMPLEBOUND}. We first note that in simple isentropic plane-symmetry,
$\geop{t}$ commutes with $\geop{u} 
+ 
\phi(u)
\frac{\muxmulevelsetvalue}{\geop{t} \muPS(t,u)} \geop{t}$ 
and thus, by \eqref{AE:LSPEEDTIMESMUPLANESYMMETRYWITHEXPLICITSOURCE} and \eqref{AE:PSTRANSPORTEQUATIONFORROUGHTIMEFUNCTION},
$\geop{t} \PStimefunctionarg{\muxmulevelsetvalue}$ satisfies the following transport equation:
\begin{align} \label{AE:PSGEOPTCOMMUTEDTRANSPORTEQUATIONFORROUGHTIMEFUNCTION}
\left\lbrace
\geop{u} 
+ 
\phi(u)
\frac{\muxmulevelsetvalue}{\geop{t} \muPS(t,u)} \geop{t} 
\right\rbrace
\geop{t} \PStimefunctionarg{\muxmulevelsetvalue}(t,u)
& = 0.
\end{align} 
Moreover, since the same arguments used to prove \eqref{E:GRADIENTOFTIMEFUNCTIONAGREESWITHGRADIENTOFMUALONGMUXMUEQUALSMINUSKAPPHYPERSURFACE}
imply that
$
\geop{t} \PStimefunctionarg{\muxmulevelsetvalue}|_{\PSdatahypfortimefunctiontwoarg{- \muxmulevelsetvalue}{[2 \timefunction_0,\frac{1}{2} \timefunction_0]}}
= - 
\geop{t} \muPS|_{\PSdatahypfortimefunctiontwoarg{- \muxmulevelsetvalue}{[2 \timefunction_0,\frac{1}{2} \timefunction_0]}}
$,
we deduce from \eqref{AE:BOUNDSONLMUINTERESTINGREGION} and \eqref{AE:PSLOCATIONOFXMUEQUALSMINUSKAPPA}
the following ``data estimates:''
\begin{align} \label{AE:PSDATAESTIMATESFORGEOPTROUGHTIMEFUNCTION}
\frac{15}{16} \PSdataamplitude
&
=
\frac{15}{16} \blowupdeltaPS
\leq
\geop{t} \PStimefunctionarg{\muxmulevelsetvalue}|_{\PSdatahypfortimefunctiontwoarg{- \muxmulevelsetvalue}{[2 \timefunction_0,\frac{1}{2} \timefunction_0]}}
\leq
\frac{17}{16} \blowupdeltaPS = \frac{17}{16} \PSdataamplitude.
\end{align}
We emphasize that \eqref{AE:PSINTEGRALCURVESOFROUGHGEOPCONTAINEDINFLATDEVELOPMENT} implies that the
data hypersurfaces $\PSdatahypfortimefunctiontwoarg{- \muxmulevelsetvalue}{[2 \timefunction_0,\frac{1}{2} \timefunction_0]}$ and the integral curves of $\geop{u} + \phi(u) \frac{\muxmulevelsetvalue}{\geop{t}\muPS} \geop{t}$ emanating from them are contained in the region of classical existence. 
The desired bounds in \eqref{AE:PSGEOPTDERIVATIVEOFROUGHTIMEFUNCTIONKEYBOUND}
now follow from the transport equation \eqref{AE:PSGEOPTCOMMUTEDTRANSPORTEQUATIONFORROUGHTIMEFUNCTION}
for $\geop{t} \PStimefunctionarg{\muxmulevelsetvalue}$ and the data estimates \eqref{AE:PSDATAESTIMATESFORGEOPTROUGHTIMEFUNCTION}.
To derive \eqref{AE:PSGEOPUDERIVATIVEOFROUGHTIMEFUNCTIONSIMPLEBOUND},
we first use \eqref{AE:PSTRANSPORTEQUATIONFORROUGHTIMEFUNCTION}
to deduce the pointwise bound
$|\geop{u} \PStimefunctionarg{\muxmulevelsetvalue}(t,u)|
\leq
\phi(u)
\frac{\muxmulevelsetvalue}{|\geop{t} \muPS(t,u)|} |\geop{t} \PStimefunctionarg{\muxmulevelsetvalue}(t,u)|$.
Also using \eqref{AE:PSDATAESTIMATESFORGEOPTROUGHTIMEFUNCTION} and \eqref{AE:BOUNDSONLMUINTERESTINGREGION}
and the fact that $\phi \geq 0$ is supported on $\lbrace |u| \leq \interestingu \rbrace$ and bounded by $1$,
we further deduce that
$|\geop{u} \PStimefunctionarg{\muxmulevelsetvalue}(t,u)|
\leq 2 \muxmulevelsetvalue \leq 2 \muxmulevelsetvalue_0
$.
From this bound and \eqref{AE:PSKAPPSOISMULTIPLEOFTIMEFUNCTION0TIMESPSIAMPLITUDE}, 
we conclude \eqref{AE:PSGEOPUDERIVATIVEOFROUGHTIMEFUNCTIONSIMPLEBOUND}.

We now exhibit the diffeomorphism properties of 
the map $(\timefunction,u) \rightarrow (t,u)$ and its inverse $\PSCHOVgeotorough{\muxmulevelsetvalue}$.  
Straightforward calculations show that 
$\mbox{\upshape det} \frac{\widetilde{\partial} (t,u)}{\widetilde{\partial} (\timefunction,u)} = \roughgeop{\timefunction} t$. By the chain rule, we see that 
$\roughgeop{\timefunction} t= \frac{1}{\geop{t} \timefunctionarg{\muxmulevelsetvalue}}$ and thus the estimate \eqref{AE:PSGEOPTDERIVATIVEOFROUGHTIMEFUNCTIONKEYBOUND} implies that the matrix is $\frac{\widetilde{\partial} (t,u)}{\widetilde{\partial} (\timefunction,u)}$ is  invertible. In view of the injectivity established shortly after \eqref{AE:PSINTEGRALCURVESOFROUGHGEOPCONTAINEDINFLATDEVELOPMENT}, 
we conclude that the map $(\timefunction,u) \rightarrow (t,u)$
is a diffeomorphism from $[2 \timefunction_0,\frac{1}{2} \timefunction_0] \times \mathbb{R}$
onto its image $\PStwoargMrough{[2 \timefunction_0,\frac{1}{2} \timefunction_0],(-\infty,\infty)}{\muxmulevelsetvalue}$, 
and that its inverse map
$
	\PSCHOVgeotorough{\muxmulevelsetvalue}(t,u) 
	= 
	(\timefunctionarg{\muxmulevelsetvalue},u)
$
is a diffeomorphism from
$\PStwoargMrough{[2 \timefunction_0,\frac{1}{2} \timefunction_0],(-\infty,\infty)}{\muxmulevelsetvalue}$ onto $[2 \timefunction_0,\frac{1}{2} \timefunction_0] \times \mathbb{R}$. 
Moreover, from the identity 
$\mbox{\upshape det} \frac{\widetilde{\partial}(\PStimefunctionarg{\muxmulevelsetvalue},u)}{\widetilde{\partial}(t,u)} 
= 
\geop{t} \PStimefunctionarg{\muxmulevelsetvalue}$
and the estimate \eqref{AE:PSGEOPTDERIVATIVEOFROUGHTIMEFUNCTIONKEYBOUND},
we also conclude \eqref{AE:DETERMINATOFPSCHOVMAPFROMGEOTOROUGH}.

Now that we have shown that $\PSCHOVgeotorough{\muxmulevelsetvalue}$
is a diffeomorphism on $\PStwoargMrough{[2 \timefunction_0,\frac{1}{2} \timefunction_0],(-\infty,\infty)}{\muxmulevelsetvalue}$,
we can derive $C^3$ estimates for the map.
First, we recall that
$
\timefunctionarg{\muxmulevelsetvalue}|_{\PSdatahypfortimefunctiontwoarg{- \muxmulevelsetvalue}{[2 \timefunction_0,\frac{1}{2} \timefunction_0]}}
= - \muPS|_{\PSdatahypfortimefunctiontwoarg{- \muxmulevelsetvalue}{[2 \timefunction_0,\frac{1}{2} \timefunction_0]}}
$.
We also note that the same arguments used to prove \eqref{E:GRADIENTOFTIMEFUNCTIONAGREESWITHGRADIENTOFMUALONGMUXMUEQUALSMINUSKAPPHYPERSURFACE}
imply that
$
\geop{t} \PStimefunctionarg{\muxmulevelsetvalue}|_{\PSdatahypfortimefunctiontwoarg{- \muxmulevelsetvalue}{[2 \timefunction_0,\frac{1}{2} \timefunction_0]}}
= - 
\geop{t} \muPS|_{\PSdatahypfortimefunctiontwoarg{- \muxmulevelsetvalue}{[2 \timefunction_0,\frac{1}{2} \timefunction_0]}}
$,
and
$
\geop{u} \PStimefunctionarg{\muxmulevelsetvalue}|_{\PSdatahypfortimefunctiontwoarg{- \muxmulevelsetvalue}{[2 \timefunction_0,\frac{1}{2} \timefunction_0]}}
= - 
\geop{u} \muPS|_{\PSdatahypfortimefunctiontwoarg{- \muxmulevelsetvalue}{[2 \timefunction_0,\frac{1}{2} \timefunction_0]}}
$.
Hence,
we can commute equation \eqref{AE:PSTRANSPORTEQUATIONFORROUGHTIMEFUNCTION} up to $3$ times
with the elements of $\lbrace \geop{t}, \geop{u} \rbrace$, use the initial conditions on 
$\PSdatahypfortimefunctiontwoarg{- \muxmulevelsetvalue}{[2 \timefunction_0,\frac{1}{2} \timefunction_0]}$ noted above,
use
\eqref{AE:LLUPS}
and
\eqref{AE:LINFINITYINITIALROUGHHYPERSURFACELDERIVATIVEOFMUANDTRANSVERSALDERIVATIVES}--\eqref{AE:LINFINITYMUANDTRANSVERSALDERIVATIVES},
and use a standard argument based on Gr\"{o}nwall's inequality to deduce that:
\begin{align}\label{AE:PSC3ESTIMATESFORROUGHTIMEFUNCTION}
	\| \PStimefunctionarg{\muxmulevelsetvalue} \|_{C_{\textnormal{geo}}^3(\PStwoargMrough{[2 \timefunction_0,\frac{1}{2} \timefunction_0],(-\infty,\infty)}{\muxmulevelsetvalue})}
	& \lesssim 1.
\end{align}
From \eqref{AE:PSC3ESTIMATESFORROUGHTIMEFUNCTION}, we conclude \eqref{AE:PSC3ESTIMATESFORCHOVMAPFROMGEOTOROUGH}.
We clarify that the fact that $\PStimefunctionarg{\muxmulevelsetvalue}$ agrees with the $C^3$ function $\muPS$ up to first-order along
$\PSdatahypfortimefunctiontwoarg{- \muxmulevelsetvalue}{[2 \timefunction_0,\frac{1}{2} \timefunction_0]}$
would allow us to obtain the estimate \eqref{AE:PSC3ESTIMATESFORROUGHTIMEFUNCTION}
using only the bound
$
\| \muPS \|_{C_{\textnormal{geo}}^3(\PStwoargMrough{[2 \timefunction_0,\frac{1}{2} \timefunction_0],(-\infty,\infty)}{\muxmulevelsetvalue})}
\lesssim 1
$
and a $C^2$ bound on the embedding of the hypersurface portion
$\PSdatahypfortimefunctiontwoarg{- \muxmulevelsetvalue}{[2 \timefunction_0,\frac{1}{2} \timefunction_0]}$ into
geometric coordinate space;
this is the same phenomenon we encountered in the bulk of the paper in Lemma~\ref{L:ODESOLUTIONSTHATARESMOOTHERTHANTHEDATAHYPERSURFACE}.

Next, we note that the function
$\PSCartesiantisafunctiononlevelsetsofroughtimefunctionarg{\timefunction}{\muxmulevelsetvalue}$
from the statement of the theorem is the first component of the map
$u \rightarrow \PSInverseCHOVgeotorough{\muxmulevelsetvalue}(\timefunction,u)$. Moreover, 
\eqref{AE:PSINTEGRALCURVESOFROUGHGEOPCONTAINEDINFLATDEVELOPMENT} implies that
$\PSCartesiantisafunctiononlevelsetsofroughtimefunctionarg{\timefunction}{\muxmulevelsetvalue}(\mathbb{R})
\in [\blowuptimePS - 2\PSinterestingusmallmultipleofamplitude \blowuptimePS,\blowuptimePS - 2\PSBigDelta]$. 
Since \eqref{AE:PSC3ESTIMATESFORCHOVMAPFROMGEOTOROUGH} and the diffeomorphism properties of
$\PSCHOVgeotorough{\muxmulevelsetvalue}$ imply that 
$
\| \PSInverseCHOVgeotorough{\muxmulevelsetvalue} \|_{C^3([2 \timefunction_0,\frac{1}{2} \timefunction_0] \times \mathbb{R})}
\lesssim 1
$,
we conclude that for $\timefunction \in [2 \timefunction_0,\frac{1}{2} \timefunction_0]$,
$\| \PSCartesiantisafunctiononlevelsetsofroughtimefunctionarg{\timefunction}{\muxmulevelsetvalue} \|_{C^3(\mathbb{R})} \lesssim 1$
and that the level sets $\twoarghypPS{\timefunction}{\muxmulevelsetvalue}$ 
are the graphical surfaces in \eqref{AE:LEVELSETSOFTIMEFUNCTIONAREAGRAPH}.

\medskip
\noindent \textbf{Proof of \eqref{AE:PSDATAJACOBIANDETERMINANTRATIOBOUND}:}
Using the chain rule relation
$
\frac{\partial(\muPS,\muXPS \muPS)}{\partial (\timefunction,u)}
=
\frac{\partial(\muPS,\muXPS \muPS)}{\partial (t,u)}
\cdot
\left(\frac{\partial (\timefunction,u)}{\partial (t,u)} \right)^{-1}
$
and the estimates
\eqref{AE:PSEXPANSIONOFMAPFROMGEOMETRICTOMUXUMCOORDINATES}
and
\eqref{AE:PSGEOPTDERIVATIVEOFROUGHTIMEFUNCTIONKEYBOUND}--\eqref{AE:PSGEOPUDERIVATIVEOFROUGHTIMEFUNCTIONSIMPLEBOUND},
we compute that for $(\timefunction,u) \in [2 \timefunction_0,\frac{1}{2} \timefunction_0] \times [-\interestingu,\interestingu]$, 
we have:
\begin{align}	 \label{AE:PSEXPANSIONOFMAPFROMROUGHADAPTEDTOMUXUMCOORDINATES}
	\frac{\partial(\muPS,\muXPS \muPS)}{\partial (\timefunction,u)}
	= \begin{pmatrix}
			- 1 & 0 
				\\
		0 & \PSdatamuHessianTaylorcoefficient
		\end{pmatrix}
		+
		\mathcal{E},
\end{align}
where the entries of the ``error matrix''
$\mathcal{E} \eqdef 
	\begin{pmatrix}
			\mathcal{E}_{11} & \mathcal{E}_{12}
				\\
		\mathcal{E}_{21} & \mathcal{E}_{22}
		\end{pmatrix}
$
satisfy the following estimates: 
\begin{align} \label{AE:PSERRORMATRIXENTRIESBOUNDEXPANSIONOFMAPFROMROUGHADAPTEDTOMUXUMCOORDINATES}
	|\mathcal{E}_{11}|
	& 
	\leq \frac{1}{8},
	&
	|\mathcal{E}_{12}|,
		\,
	|\mathcal{E}_{21}|,
		\,
	|\mathcal{E}_{22}|
	&
	= 
	\mathcal{O}(\PSdataamplitude)
	+
	\mathcal{O}(\PSinterestingusmallmultipleofamplitude).
\end{align}
From \eqref{AE:PSEXPANSIONOFMAPFROMROUGHADAPTEDTOMUXUMCOORDINATES} and 
\eqref{AE:PSERRORMATRIXENTRIESBOUNDEXPANSIONOFMAPFROMROUGHADAPTEDTOMUXUMCOORDINATES},
we deduce that if
$
(\timefunction,u)
\in
[2 \timefunction_0,\frac{1}{2} \timefunction_0] \times [-\interestingu,\interestingu]
$,
and 
if $\PSdataamplitude$ 
and
$\PSinterestingusmallmultipleofamplitude$ are sufficiently small, 
then the Jacobian matrix
$
\PSCHOVJacobianroughtomumuxmu{\muxmulevelsetvalue}(\timefunction,u) 
\eqdef
\frac{\partial(\muPS,\muXPS \muPS)}{\partial (\timefunction,u)}
$
is invertible,
and that for every pair of points 
$(\timefunction_1,u_1), (\timefunction_2,u_2) 
\in
[2 \timefunction_0,\frac{1}{2} \timefunction_0] \times [-\interestingu,\interestingu]
$,
we have
$
\left|
\PSInverseCHOVJacobianroughtomumuxmu{\muxmulevelsetvalue}(\timefunction_1,u_1) \PSCHOVJacobianroughtomumuxmu{\muxmulevelsetvalue}(\timefunction_2,u_2)
-
\mbox{\upshape ID} 
\right|_{\mbox{\upshape}Euc}
\leq \frac{1}{4}
$,
which is the desired bound \eqref{AE:PSDATAJACOBIANDETERMINANTRATIOBOUND}.

\noindent \textbf{Proof of the properties of $\PSCHOVgeotomumuxmu$ and \eqref{E:PSCHOVFROMGEOMETRICCOORDINATESTOMUWEGIGHTEDXMUCOORDINATES}--\eqref{E:PSJACOBIANMATRIXFORCHOVFROMGEOMETRICCOORDINATESTOMUWEGIGHTEDXMUCOORDINATES}:}
These results follow from the ``intermediate step'' mentioned above, including
the estimates 
\eqref{AE:PSEXPANSIONOFMAPFROMGEOMETRICTOMUXUMCOORDINATES}--\eqref{AE:PSKEYESTIMATEFORINTERTIBILITYOFTUTOMUMUXMUMAP}.

\end{proof}

\begin{corollary}[The data assumptions of Sect.\,\ref{SS:ASSUMPTIONSONDATA} are satisfied]
	\label{AC:SIMPLEISENTROPICPLANESYMMETRICSOLUTIONSSATISFYALLTHEASSUMPTIONS}
	The solutions provided by Theorem~\ref{T:PROPERTIESOFADMISSIBLESOLUTIONSLAUNCEDBYRESCALED}
	induce data on $\hypthreearg{\timefunction_0}{[- \rightu,\leftu]}{\muxmulevelsetvalue}$
	that satisfy all the assumptions stated in Sect.\,\ref{SS:ASSUMPTIONSONDATA},
	where 
	$\overline{\varrho} > 0$ is a fixed constant density (see Def.\,\ref{D:LOGDENS}), 
	$\rightu$, 
	$\leftu$, 
	$\interestingu$,
	$\timefunction_0$, and $\muxmulevelsetvalue_0$ are as in the statement of the theorem,
	and the parameters 
	$\PSdataalpha$,
	$\blowupdeltaPS$,
	$\transversalsizedeltaPS$,
	$\initialsmall$,
	$\PSmutransversalHessiansize$, 
	and
	$\PSboringregionmupositive$
	can be chosen to satisfy the following relations, where the implicit constants can depend on the 
	function $\mathring{\varphi}$ from the statement of the theorem:
		\begin{itemize}
		\item $\PSdataalpha \lesssim \PSdataamplitude$
		\item $\blowupdeltaPS = \PSdataamplitude$
		\item  $\transversalsizedeltaPS \lesssim \PSdataamplitude$
		\item $\initialsmall = 0$
		\item $\PSmutransversalHessiansize 
		= 
		\min\lbrace 
			\frac{\PSdatamuHessianTaylorcoefficient}{2}, 
			\frac{1}{2 \PSdatamuHessianTaylorcoefficient}
		\rbrace$ 
		(note that by \eqref{AE:PSTIMEFUNCTION0CONSTRAINTS}--\eqref{AE:PSKAPPSOISMULTIPLEOFTIMEFUNCTION0TIMESPSIAMPLITUDE},
		the assumption \eqref{E:KAPPA0ISSMALLERTHANM2RIGHTUOVER16} is satisfied)
	\item $\PSboringregionmupositive = \frac{1}{2} \PSLmunottoonegativeparameter$
	\end{itemize}
	
\end{corollary}

\begin{proof}
	All aspects of the corollary, aside from 
	the estimate \eqref{E:LINFINITYINITIALROUGHHYPERSURFACEMUANDTRANSVERSALDERIVATIVES},
	can be deduced by comparing the results provided by
	Theorem~\ref{T:PROPERTIESOFADMISSIBLESOLUTIONSLAUNCEDBYRESCALED} 
	with the assumptions stated in Sect.\,\ref{SS:ASSUMPTIONSONDATA},
	and using the following facts: 
	\begin{itemize}
		\item The solution vanishes on the complement of the region
		$\lbrace (t,u) \in \mathbb{R} \times \mathbb{R} \ | \ u \in [- \rightu,\leftu] \rbrace$,
			and in particular is trivial along the null hypersurface
			$\nullhyparg{- \rightu}^{[0,\infty)}$.
		\item By \eqref{AE:ROUGHHYPSFORADMISSIBLESOLUTIONSCONTAINEDINCAUCHYSTABILITYREGION},
			the hypersurface $\twoarghypPS{\timefunction_0}{\muxmulevelsetvalue}$ 
			is contained in the region of classical existence relative to the geometric coordinates 
			(where the estimates of the theorem hold).
		\item Relative to the geometric coordinates, 
			$\RRiemannPS$ depends only on $u$.
		\item Along $\twoarghypPS{\timefunction_0}{\muxmulevelsetvalue}$,
			$t$ is a function of $u$ (see \eqref{AE:LEVELSETSOFTIMEFUNCTIONAREAGRAPH}).
		\end{itemize}
	The estimate \eqref{E:LINFINITYINITIALROUGHHYPERSURFACEMUANDTRANSVERSALDERIVATIVES} 
	(with $\initialsmall = 0$, 
	$\blowupdeltaPS$ in the role of $\mathring{\updelta}_*$, $\muXPS$ in the role of $\muX$, and $\muPS$ in the role of $\upmu$)
	follows from the facts noted above,
	\eqref{AE:EXPLICITFORMUSEFULIDENTITYFORPSLMUSOURCETERMFUNCTION},
	\eqref{AE:LSPEEDTIMESMUPLANESYMMETRYWITHEXPLICITSOURCE}--\eqref{AE:SPEEDTIMESMUPLANESYMMETRYWITHEXPLICITSOURCE},
	and the fact that in the context of Theorem~\ref{T:PROPERTIESOFADMISSIBLESOLUTIONSLAUNCEDBYRESCALED},
	we have $\LunitPS = \geop{t}$, 
	$\muXPS = \geop{u}$, 
	and
	$0 \leq t \leq \frac{1}{\blowupdeltaPS}$.
	
\end{proof}

\begin{definition}[Admissible background solutions] \label{AD:ADMISSIBLEBACKGROUND}
	We refer to the compactly supported
	solutions furnished by Theorem~\ref{T:PROPERTIESOFADMISSIBLESOLUTIONSLAUNCEDBYRESCALED}
	as ``admissible background solutions.''
	We consider their parameters $\PSdataalpha$,
	$\blowupdeltaPS$,
	etc.\
	to be the ones guaranteed by Cor.\,\ref{AC:SIMPLEISENTROPICPLANESYMMETRICSOLUTIONSSATISFYALLTHEASSUMPTIONS}.
\end{definition}

\subsection{The Cauchy stability region} \label{SS:REGIONOFCAUCHYSTABILITY}
Fix any of the admissible ``background'' (shock-forming) simple isentropic plane-symmetric 
solutions from Def.\,\ref{AD:ADMISSIBLEBACKGROUND}. Recall that these solutions are supported in the strip
$\lbrace (t,u) \in \mathbb{R} \times \mathbb{R} \ | \ u \in [-\rightu,\leftu] \rbrace$.
In this section, we discuss the behavior of the solution in a ``Cauchy stability region,''
which is a large region close to the singularity where the background solution exists classically.
In particular, we describe Fig.\,\ref{F:CAUCHYSTABILITYREGION}, which
will guide our discussion in Appendix~\ref{A:OPENSETOFDATAEXISTS}.
In Appendix~\ref{A:OPENSETOFDATAEXISTS}, we will use Cauchy stability arguments to 
show that perturbations of the background solution
stay close to it, in all relevant norms, 
in the Cauchy stability region.
From the Cauchy stability estimates, it will follow that there exist open sets data -- without symmetry --
satisfying the assumptions stated in
Sect.\,\ref{SS:ASSUMPTIONSONDATA}.
In the rest of this section (and also in Appendix~\ref{A:OPENSETOFDATAEXISTS}), 
we view the background solutions 
as solutions in three spatial dimensions that are independent of the torus coordinates $(x^2,x^3)$.

Let $\blowupdeltaPS$ denote the quantity \eqref{E:DELTASTARDEF} 
evaluated at the (shock-forming) background solution.
In Theorem~\ref{T:PROPERTIESOFADMISSIBLESOLUTIONSLAUNCEDBYRESCALED},
we showed that relative to the geometric coordinates, 
the first singular point for the background solution
occurs at $(t,u) = (\blowuptimePS,0)$,
where $\blowuptimePS  = \frac{1}{\blowupdeltaPS}$.
We define:
\begin{align} \label{AE:FARRIGHTU}
	\farrightu 
	& \eqdef \rightu + \frac{18}{\blowupdeltaPS}.
\end{align}

For $0 \leq t_1 \leq t_2 \leq \blowuptimePS - 2\PSBigDelta$,
where $\PSBigDelta > 0$ is the small constant defined in \eqref{AE:BIGDELTAPS},
we define (see Def.\,\ref{D:ACOUSTICSUBSETSOFSPACETIME}):
\begin{subequations}
\begin{align} \label{E:MAINCAUCHYSTABILITYREGION}
	\CSregion_{Main}^{\blowuptimePS - 2\PSBigDelta} 
	& \eqdef 
		\bigcup_{t \in [0,\blowuptimePS - 2\PSBigDelta]}
		\Sigma_t^{[- \rightu,\leftu]},
			\\
\CSregion_{Small}^{[t_1,t_2]}
& \eqdef
	\bigcup_{t \in [t_1,t_2]}
	\Sigma_t^{[3t -\farrightu,- \rightu]},
		\label{E:SMALLCAUCHYSTABILITYREGION} 
			\\
\begin{split}
\csspacelikehypersurface^{[t_1,t_2]} 
&
\eqdef 
\lbrace (t,u,x^2,x^3) \ | \ 3 t - u = \farrightu, t \in [t_1,t_2], \mbox{ and } (x^2,x^3) \in \mathbb{T}^2 \rbrace
	\label{E:FARRIGHTSPACELIKEBOUNDARYCAUCHYSTABILITYREGION} 
		\\
& =
	\bigcup_{t \in [t_1,t_2]} \ell_{t,3t - \farrightu}.
\end{split}
\end{align}
\end{subequations}
Note that $\csspacelikehypersurface^{[0,5 \blowuptimePS]}$ is the right boundary
of $\CSregion_{Small}^{[0,5 \blowuptimePS]}$; see Fig.\,\ref{F:CAUCHYSTABILITYREGION}.

Consider the ``Cauchy stability region''
$\CSregion^{\blowuptimePS;\PSBigDelta}$
depicted relative to the geometric coordinates $(t,u,x^2,x^3)$ in Fig.\,\ref{F:CAUCHYSTABILITYREGION}.
We decompose this region into various sub-regions:
\begin{subequations}
\begin{align}
\CSregion^{\blowuptimePS;2\PSBigDelta} 
& 
\eqdef 
\CSregion_{Main}^{[0,\blowuptimePS - 2\PSBigDelta]} 
\cup
\CSregion_{Small}^{[0,5 \blowuptimePS]},
	\label{E:CAUCHYSTABILITYREGIONDECOMPOSTEDINTOMAINREGIONANDSMALLREGION} \\
\CSregion_{Small}^{[0,5 \blowuptimePS]} 
& = 
\CSregion_{Small}^{[0,\blowuptimePS - 2\PSBigDelta]} 
\cup
\CSregion_{Small}^{[\blowuptimePS - 2\PSBigDelta,5 \blowuptimePS]},
	 \label{E:CAUCHYSTABILITYSOLUTIONISSMALLREGION} \\
\widehat{\CSregion}^{[0,\blowuptimePS - 2\PSBigDelta]}
& 
= 
\CSregion_{Main}^{[0,\blowuptimePS - 2\PSBigDelta]} 
\cup
\CSregion_{Small}^{[0,\blowuptimePS - 2\PSBigDelta]}.
	\label{E:BEFORESHOCKTIMECAUCHYSTABILITYREGION}
\end{align}
\end{subequations}

\begin{center}
\begin{figure}  
\begin{overpic}[scale=.6, grid = false, tics=5, trim=-.5cm -1cm -1cm -.5cm, clip]{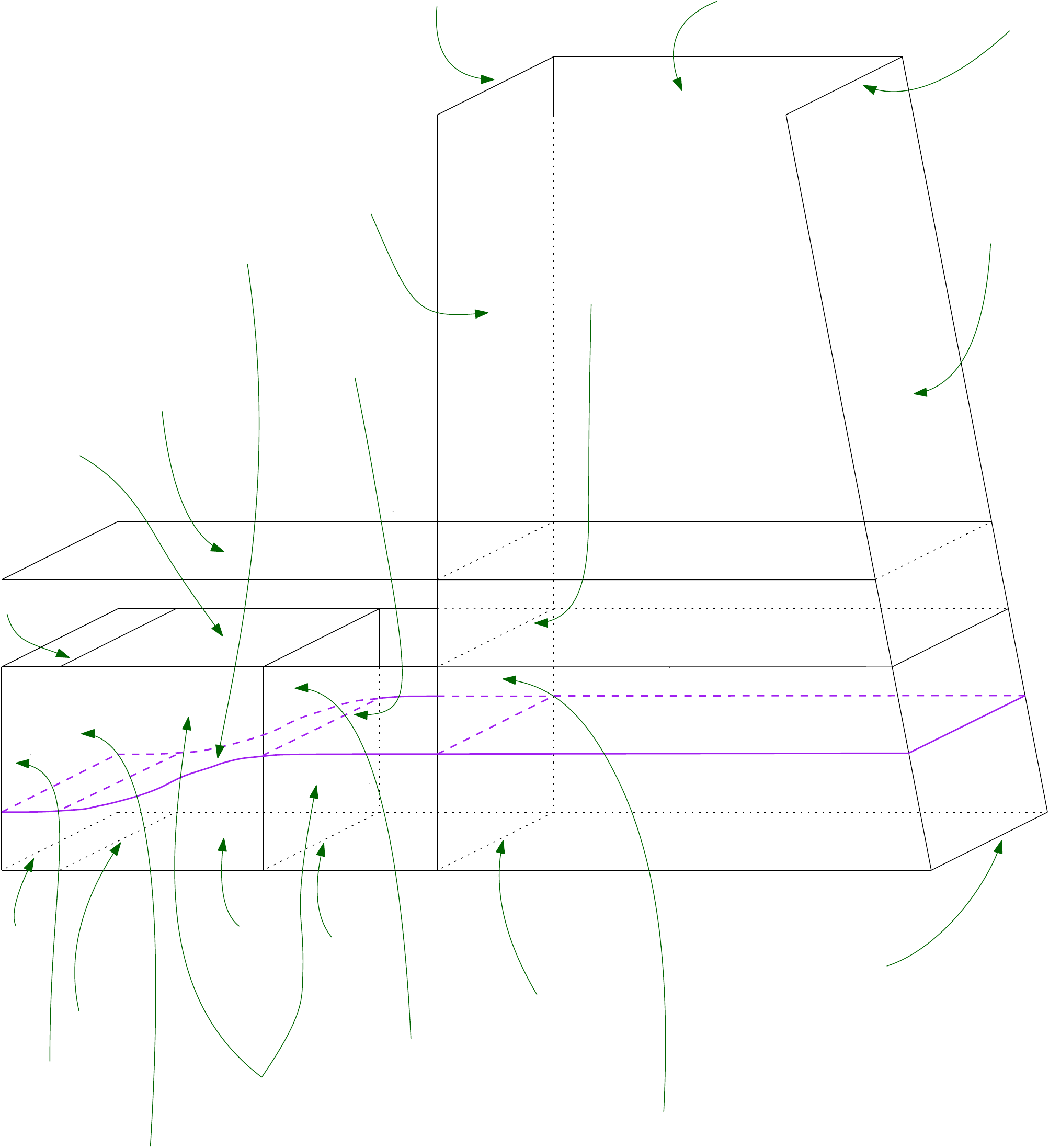}
	\put (1,21) {$\ell_{0,\leftu}$}
	\put (7,14) {$\ell_{0,\interestingu}$}
	\put (28,19) {$\ell_{0,-\interestingu}$}
	\put (44,15) {$\ell_{0,- \rightu}$}
	\put (68,18) {$\ell_{0,-\farrightu}$}
	\put (-3,49) {$\ell_{\blowuptimePS - 2\PSBigDelta,\interestingu}$}
	\put (46,74) {$\ell_{\blowuptimePS - 2\PSBigDelta,- \rightu}$}
	\put (33,100) {$\ell_{5\blowuptimePS,- \rightu}$}
	\put (80,98) {$\ell_{5\blowuptimePS,15 \blowuptimePS -\farrightu}$}
	\put (54,101) {$\Sigma_{5 \blowuptimePS}^{[15 \blowuptimePS - \farrightu,- \rightu]}$}
	\put (17.5,19) {$\Sigma_0^{[- \rightu,\leftu]}$}
	\put (11,68) {$\Sigma_{\blowuptimePS}^{[- \rightu,\leftu]}$}
	\put (1,64) {$\Sigma_{\blowuptimePS - 2\PSBigDelta}^{[- \rightu,\leftu]}$}
	\put (14,78) {$\threearghypPS{\timefunction_0/2}{[- \rightu,\leftu]}{\muxmulevelsetvalue}$}
	\put (26,69) {$\PStwoargroughtori{\timefunction_0/2,-\interestingu}{\muxmulevelsetvalue}$}
	\put (29,82) {$\nullhyparg{- \rightu}^{5 \blowuptimePS}$}
	\put (51,3) {$\nullhyparg{- \rightu}^{[0,\blowuptimePS - 2\PSBigDelta]}$}
	\put (33,8) {$\nullhyparg{-\interestingu}^{[0,\blowuptimePS - 2\PSBigDelta]}$}
	\put (0,7) {$\nullhyparg{\leftu}^{[0,\blowuptimePS - 2\PSBigDelta]}$}
	\put (10,0) {$\nullhyparg{\interestingu}^{[0,\blowuptimePS - 2\PSBigDelta]}$}
	\put (16,6) {$\CSregion_{Main}^{[0,\blowuptimePS - 2\PSBigDelta]}$}
	\put (53,34) {$\CSregion_{Small}^{[0,\blowuptimePS - 2\PSBigDelta]}$}
	\put (55,63) {$\CSregion_{Small}^{[0,5 \blowuptimePS]}$}
	\put (79,79) {$\csspacelikehypersurface^{[0,5 \blowuptimePS]}$}
	\put(54,22.5){\vector(-1,0){7}}
	\put(56,22){$u$}
	\put(-2,34){\vector(0,1){7}}
	\put (-2.5,31) {$t$}
\end{overpic}
\caption{The Cauchy stability region $\CSregion^{\blowuptimePS;2\PSBigDelta}$, not drawn to scale}
\label{F:CAUCHYSTABILITYREGION}
\end{figure}
\end{center}
Note that \eqref{AE:PSBLOWUPTIME} and \eqref{AE:FARRIGHTU}
imply that the top boundary of $\CSregion_{Small}^{[0,5 \blowuptimePS]}$,
namely $\Sigma_{5 \blowuptimePS}^{[15 \blowuptimePS - \farrightu,- \rightu]}$
has $u$-width equal to $3 \blowuptimePS = \frac{3}{\blowupdeltaPS}$.

The background solutions provided by Theorem~\ref{T:PROPERTIESOFADMISSIBLESOLUTIONSLAUNCEDBYRESCALED} are smooth
in $\CSregion^{\blowuptimePS;2\PSBigDelta}$ and are trivial in the sub-region $\CSregion_{Small}^{[0,\blowuptimePS - 2\PSBigDelta]}$, 
that is, both Riemann invariants identically vanish in $\CSregion_{Small}^{[0,\blowuptimePS - 2\PSBigDelta]}$.
This ensures, in particular, that the surface portions
$\csspacelikehypersurface^{[0,5 \blowuptimePS]}$
are $\gfour$-spacelike with respect to the acoustical metric of the background.
To see this, we compute that the future-directed unit normal to 
$\csspacelikehypersurface^{[0,5 \blowuptimePS]}$
is 
$3 \Transport - \frac{1}{\upmu} \Lunit = 2 \Lunit + 3 X$,
where we have used \eqref{E:BISLPLUSX} and the fact 
that $\upmu \equiv 1$ in the trivial region
$\CSregion_{Small}^{[0,5 \blowuptimePS]}$
(which contains $\csspacelikehypersurface^{[0,5 \blowuptimePS]}$).
We can therefore use Lemma~\ref{L:BASICPROPERTIESOFVECTORFIELDS} to 
compute that $\gfour(2 \Lunit + 3 X,2 \Lunit + 3 X) = -3$,
which indeed implies that $3 \Transport - \frac{1}{\upmu} \Lunit$ is 
$\gfour$-timelike along $\csspacelikehypersurface^{[0,5 \blowuptimePS]}$.

To help prepare the reader for
Appendix~\ref{A:OPENSETOFDATAEXISTS}, we now further discuss some aspects of Fig.\,\ref{F:CAUCHYSTABILITYREGION} that
follow from the conclusions of Theorem~\ref{T:PROPERTIESOFADMISSIBLESOLUTIONSLAUNCEDBYRESCALED}.
Let $\PStimefunctionarg{\muxmulevelsetvalue}$ denote the rough time function of the background,
and let $\threearghypPS{\timefunction}{[u_1,u_2]}{\muxmulevelsetvalue}$ denote
the level set portions
$\lbrace \PStimefunctionarg{\muxmulevelsetvalue} = \timefunction \rbrace \cap \lbrace u \in [u_1,u_2] \rbrace$,
and let 
$
\PStwoargroughtori{\timefunction,u'}{\muxmulevelsetvalue} 
= 
\lbrace \PStimefunctionarg{\muxmulevelsetvalue} = \timefunction \rbrace \cap \lbrace u = u' \rbrace
$
denote the background rough tori.
The background solution, though smooth,
is ``about'' to form a shock in $\CSregion_{Main}^{\blowuptimePS - 2\PSBigDelta}$ near $u=0$. 
Moreover, if $\timefunction_0 < 0$ and $\muxmulevelsetvalue_0$ are sufficiently small
as in the statement of the theorem,
then for $\muxmulevelsetvalue \in [0,\muxmulevelsetvalue_0]$,
the background rough time functions
$\PStimefunctionarg{\muxmulevelsetvalue}$
are defined on a subset of $\CSregion^{\blowuptimePS;2\PSBigDelta}$
such that for $\timefunction \in [2 \timefunction_0,\frac{\timefunction_0}{2}]$,
the level set portions
$\threearghypPS{\timefunction}{[-\farrightu,\leftu]}{\muxmulevelsetvalue}$
are contained in $\CSregion^{\blowuptimePS;2\PSBigDelta}$,
i.e., they are temporally separated from
$\Sigma_{\blowuptimePS}$ by a distance at least equal to $2\PSBigDelta$.
Note that for the background solution, if we followed it all the way to the first singular point
(which is contained in $\threearghypPS{0}{[-\farrightu,\leftu]}{0}$),
we would have
$\threearghypPS{0}{[-\farrightu,\leftu]}{0} = \Sigma_{\blowuptimePS}^{[-\farrightu,\leftu]}$. 
This is because in plane-symmetry, by \eqref{AE:PSTRANSPORTEQUATIONFORROUGHTIMEFUNCTION} with $\muxmulevelsetvalue = 0$,
$\PStimefunctionarg{0}$ can be expressed as a function of $t$ alone.

\section{The existence of an open set of data satisfying the assumptions}
\label{A:OPENSETOFDATAEXISTS}
By Cor.\,\ref{AC:SIMPLEISENTROPICPLANESYMMETRICSOLUTIONSSATISFYALLTHEASSUMPTIONS},
there exists a large family of isentropic plane-symmetric initial data on $\Sigma_0$
such that the corresponding background solutions
(and corresponding parameters) induce
data on the background-solution-dependent rough hypersurface
$\threearghypPS{\timefunction_0}{[-\rightu,\leftu]}{\muxmulevelsetvalue}$
and null hypersurface portion
$\nullhyparg{- \rightu}^{[0,\frac{5}{\blowupdeltaPS}]}$
that satisfy the assumptions in Sects.\,\ref{SS:ASSUMPTIONSONDATA}--\ref{SSS:LOCALIZEDDATAASSUMPTIONSFORMUANDDERIVATIVES}
and the parameter-size assumptions of Sect.\,\ref{SS:PARAMETERSIZEASSUMPTIONS}
with $\initialsmall = 0$.
In the following proposition, 
we show that if one perturbs 
-- without symmetry, irrotationality, or isentropicity assumptions --
the background initial data on $\Sigma_0$,
then the corresponding perturbed solutions also induce data that
satisfy the assumptions in 
Sects.\,\ref{SSS:QUANTITATIVEASSUMPTIONSONDATAAWAYFROMSYMMETRY}--\ref{SSS:LOCALIZEDDATAASSUMPTIONSFORMUANDDERIVATIVES}
and the parameter-size assumptions of Sect.\,\ref{SS:PARAMETERSIZEASSUMPTIONS}
with $\initialsmall$ non-negative but small.
In conjunction with Theorem~\ref{T:DEVELOPMENTANDSTRUCTUREOFSINGULARBOUNDARY},
this shows that our main results hold for open sets of solutions.

\begin{remark}[We can choose the smallness of $|\timefunction_0|$]
		\label{R:INITIALROUGHHYPERSURFACEISCLOSETOSINGULARITY}
		Theorem~\ref{T:PROPERTIESOFADMISSIBLESOLUTIONSLAUNCEDBYRESCALED} implies that
		we can choose and fix the parameter
		$\timefunction_0 = - \mupositive < 0$
		to be as close to $0$ as we want. 
		For perturbed solutions, the smallness of $|\timefunction_0|$
		corresponds to assuming that their initial rough hypersurfaces 
		$\hypthreearg{\timefunction_0}{[- \rightu,\leftu]}{\muxmulevelsetvalue}$
		are close to the singularity of the background solution.
		While such smallness is not essential for our analysis of perturbed solutions,
		it is helpful because it allows us to simplify the proofs of various estimates 
		in the bulk of the paper, i.e., 
		it allows us to exploit that we only have to control perturbed solutions for
		$|\timefunction_0|$ amounts of rough time.
\end{remark}

\begin{proposition}[Cauchy stability and the existence of open sets of data satisfying our assumptions]
	\label{P:CAUCHYSTABILITYANDEXISTENCEOFOPENSETS}
	Fix any of the admissible ``background'' (shock-forming) simple isentropic plane-symmetric solutions
	from Def.\,\ref{AD:ADMISSIBLEBACKGROUND}.
	Recall (see Cor.\,\ref{AC:SIMPLEISENTROPICPLANESYMMETRICSOLUTIONSSATISFYALLTHEASSUMPTIONS}) that 
	$\overline{\varrho}$, 
	$\rightu$,
	$\leftu$,	
	$\interestingu$, 
	$\blowupdeltaPS$,
	$\transversalsizedeltaPS$,
	$\muxmulevelsetvalue_0$, 
	$|\timefunction_0| = - \timefunction_0$,
	$\PSdataalpha$,
	$\blowupdeltaPS$,
	$\transversalsizedeltaPS$,
	$\PSmutransversalHessiansize$, 
	and
	$\PSboringregionmupositive$
	are \textbf{positive} parameters associated to the background solution
	and that
	$\farrightu \eqdef \rightu + \frac{18}{\blowupdeltaPS}$ (see \eqref{AE:FARRIGHTU}).
	Recall that 
	we can choose and fix
	$\mupositive$
	and
	$|\timefunction_0|$
	to be as small as we want
	(see Remark~\ref{R:INITIALROUGHHYPERSURFACEISCLOSETOSINGULARITY}).
	Let $\mathring{\Delta}_{\Sigma_0^{[-\farrightu,\leftu]}}^{\Ntop+1}$ 
	be the norm of the bona fide data perturbation
	defined in \eqref{E:PERTURBATIONSMALLNESSINCARTESIANDIFFERENTIALSTRUCTURE}.
	If $\PSBigDelta > 0$ is the small constant defined in \eqref{AE:BIGDELTAPS},
	then for all sufficiently small
	$\mathring{\Delta}_{\Sigma_0^{[-\farrightu,\leftu]}}^{\Ntop+1} > 0$
	(where the required smallness depends on the background solution),
	the perturbed fluid solution,
	the eikonal function $u$,
	and all of the auxiliary geometric quantities constructed out of $u$ 
	(such as $\upmu$, $\Lunit^i$, and $\upchi$)
	exist classically in the Cauchy stability region
	$\CSregion^{\blowuptimePS;\PSBigDelta} = \CSregion_{Main}^{\PSBigDelta} \cup \CSregion_{Small}$ 
	(which we view to be a fixed subset of geometric coordinate space)
	defined in Sect.\,\ref{SS:REGIONOFCAUCHYSTABILITY}.
	Moreover, for $\muxmulevelsetvalue \in [0,\muxmulevelsetvalue_0]$
	each perturbed time function $\timefunctionarg{\muxmulevelsetvalue}$
	exists classically on a subset of
	$\CSregion^{\PSBigDelta}$ and has a range containing $[2 \timefunction_0,\timefunction_0/2]$
	such that for $\timefunction \in [2 \timefunction_0,\timefunction_0/2]$,
	the level sets
	$\hypthreearg{\timefunction}{[-\rightu,\leftu]}{\muxmulevelsetvalue}$
	are contained in $\CSregion_{Main}^{\PSBigDelta}$.
	In particular, the perturbed solution exists classically on
	$\twoargMrough{[2 \timefunction_0,\timefunction_0/2],[-\rightu,\leftu]}{\muxmulevelsetvalue}$.
	
	Moreover, for the perturbed solution, 
	we define $\mathring{\updelta}_*$
	by \eqref{E:DELTASTARDEF},
	we define
	$\mupositive$, 
	$\boringregionmupositive$,
	$\mathring{\updelta}$,
	$\secondtransversalderivativemulowerbound$
	by:\footnote{The parameters that serve as lower bounds are defined to be half the value of the corresponding
	background parameters, while parameters that serve as upper bounds are defined to be twice the value of the corresponding
	background parameter. \label{FN:PARAMETERSSOMEAREHALFTHEBACKGROUNDVALUESANDSOMEARETWICE}} 
	\begin{align} \label{E:PERTURBEDPARAMETERSARECLOSETOBACKGROUNDONES}
	\frac{\mupositive^{\text{PS}}}{\mupositive},
		\,
	\frac{\boringregionmupositive^{\text{PS}}}{\boringregionmupositive},
		\,
	\frac{\mathring{\updelta}}{\blowupdeltaPS},
		\,
	\frac{\mathring{\upalpha}}{\PSdataalpha},
		\,
	\frac{\secondtransversalderivativemulowerboundPS}{\secondtransversalderivativemulowerbound}
	= 2,
	\end{align}
	and we define the remaining parameters to be the same as for the background solution.
	Then the following estimate holds:
	\begin{align} \label{E:PERTURBEDBLOWUPDELTAISCLOSETOBACKGROUNDONE}
	\mathring{\updelta}_* 
	= 
	\mathring{\updelta}_*^{\text{PS}} 
	+ 
	\mathcal{O}(\mathring{\Delta}_{\Sigma_0^{[-\farrightu,\leftu]}}^{\Ntop+1}),
	\end{align}
	and for these parameters, the
	perturbed solutions induce data on the perturbed rough hypersurface
	$\hypthreearg{\timefunction_0}{[- \rightu,\leftu]}{\muxmulevelsetvalue}$,
	the null hypersurface portion
	$\nullhyparg{- \rightu}^{[0,4 \mathring{\updelta}_*]}$,
	and the perturbed rough tori $\twoargroughtori{\timefunction_0,u}{\muxmulevelsetvalue}$
	that satisfy all of the assumptions of
	Sects.\,\ref{SSS:QUANTITATIVEASSUMPTIONSONDATAAWAYFROMSYMMETRY}--\ref{SSS:LOCALIZEDDATAASSUMPTIONSFORMUANDDERIVATIVES}
	with $\initialsmall = \mathcal{O}(\mathring{\Delta}_{\Sigma_0^{[-\farrightu,\leftu]}}^{\Ntop+1})$,
	where the implicit constants in ``$\mathcal{O}(\cdot)$''
	depend on the background solution.
\end{proposition}

\begin{proof}[Proof sketch] \hfill

\noindent \textbf{Overview of the main ideas of the proof}.
Because the background solutions
satisfy the assumptions of
Sects.\,\ref{SS:ASSUMPTIONSONDATA}--\ref{SSS:LOCALIZEDDATAASSUMPTIONSFORMUANDDERIVATIVES}
with $\initialsmall = 0$,
most aspects of the proposition follow from standard arguments based on Cauchy stability.
The only non-standard aspects are the following, which we flesh out in Steps 1-5 below:
\begin{itemize}
	\item (\textbf{Estimates - without derivative loss - on flat spacelike hypersurfaces and null hypersurfaces})
		We need to show that in $\CSregion^{\blowuptimePS;\PSBigDelta}$, 
		we can control the solution with respect to the geometric coordinates $(t,u,x^2,x^3)$
		up to top-order, i.e., without losing a derivative relative to the data norm
		$\mathring{\Delta}_{\Sigma_0^{[-\farrightu,\leftu]}}^{\Ntop+1}$
		defined in \eqref{E:PERTURBATIONSMALLNESSINCARTESIANDIFFERENTIALSTRUCTURE}.
		This is essentially a much easier version of the proofs of the energy estimates of
		Props.\,\ref{P:APRIORIL2ESTIMATESWAVEVARIABLES},
	\ref{P:MAINHYPERSURFACEENERGYESTIMATESFORTRANSPORTVARIABLES},
	\ref{P:ROUGHTORIENERGYESTIMATES},
	and
	\ref{P:APRIORIL2ESTIMATESACOUSTICGEOMETRY},
	except for there is one new conceptually important aspect of the
	tori energy estimates 
	(i.e., the analog of the estimates of Prop.\,\ref{P:ROUGHTORIENERGYESTIMATES}), 
	described below.
	The main simplification is that $\upmu$ is uniformly positive in 
	$\CSregion^{\blowuptimePS;\PSBigDelta}$
	(see \eqref{E:MUISUNIFORMLYPOSITIVEINCAUCHYSTABILITYREGION}),
	and this positivity allows one to use standard arguments based on
	Gr\"{o}nwall's inequality to show that the energies and null fluxes grow at most exponentially.
	\item (\textbf{Estimates - without derivative loss - on tori})
		We need to adequately control the solution up to top order on various tori.
		This control does not directly come from the energy estimates described 
		in the previous step.
		We obtain the desired control by combining averaging arguments based
		on Chebychev's inequality
		with elliptic-hyperbolic identities in the spirit of \eqref{E:INTEGRALIDENTITYFORELLIPTICHYPERBOLICCURRENT}.
\end{itemize}

Throughout this proof, we view 
$\CSregion^{\blowuptimePS;\PSBigDelta}$ to be a fixed subset of geometric coordinate space.
We will tacitly assume that $\mathring{\Delta}_{\Sigma_0^{[-\farrightu,\leftu]}}^{\Ntop+1}$ is sufficiently small.
We will also freely use notation defined in Sect.\,\ref{SS:REGIONOFCAUCHYSTABILITY}.

\medskip

\noindent \textbf{Step 1: Standard Cauchy stability with respect to the Cartesian coordinates
 in $\CSregion^{\blowuptimePS;\PSBigDelta}$}
By Cauchy stability, 
if $\mathring{\Delta}_{\Sigma_0^{[-\farrightu,\leftu]}}^{\Ntop+1}$ is small enough,
then the perturbed fluid solutions 
vary continuously with respect to the fluid data in the Cauchy stability 
region $\CSregion^{\blowuptimePS;\PSBigDelta}$ depicted in Fig.\,\ref{F:CAUCHYSTABILITYREGION}
(see also \eqref{E:CAUCHYSTABILITYREGIONDECOMPOSTEDINTOMAINREGIONANDSMALLREGION}).
Similarly, the eikonal function $u$, 
which solves the fully nonlinear transport equation \eqref{E:EIKONALEQUATION},
varies continuously in $\CSregion^{\blowuptimePS;\PSBigDelta}$ 
with respect to the fluid data (here we are viewing $u$ as a function of the Cartesian coordinates).
Cauchy stability results of this type are standard
if one measures continuity using topologies corresponding to Sobolev spaces 
respect to the Cartesian coordinates on portions of the hypersurfaces $\Sigma_t$ of constant Cartesian time
that are contained in $\CSregion^{\blowuptimePS;\PSBigDelta}$; 
see \cite{hR2013}*{Proposition~9.17} for a detailed proof in the context of Einstein's equations.
We emphasize that the present \textbf{Step 1 yields the continuous dependence of the solution 
only with respect to standard Sobolev norms defined through
Cartesian partial derivatives}.
In particular, this ``Cartesian approach'' does not yield the top-order regularity of the solution
with respect to the geometric coordinates; in the remaining steps, we will obtain the desired
top-order regularity in $\CSregion^{\blowuptimePS;\PSBigDelta}$.
We clarify that even though the eikonal function $u = u(t,x^1,x^2,x^3)$ 
varies slightly with the solution, domain of dependence considerations imply that
$\CSregion^{\blowuptimePS;\PSBigDelta}$ is always a development of the portion of the initial data
on the fixed subset $\Sigma_0^{[-\farrightu,\leftu]}$ of $\Sigma_0$.

Next, we highlight that all of the auxiliary geometric quantities constructed out of $u$, 
such as $\upmu$, $\Lunit^i$, and $\upchi$,
also vary continuously\footnote{All of these quantities can be expressed in terms of the fluid
variables, $u$, and their Cartesian coordinate partial derivatives.}
with respect to the fluid data.
In particular, since the background solution's inverse foliation density $\muPS$ is strictly positive in
$\CSregion^{\blowuptimePS;\PSBigDelta}$, the perturbed solution satisfies:
\begin{align} \label{E:MUISUNIFORMLYPOSITIVEINCAUCHYSTABILITYREGION}
	\upmu|_{\CSregion^{\blowuptimePS;\PSBigDelta}} 
	& \gtrsim 1,
\end{align}
where the implicit constants in \eqref{E:MUISUNIFORMLYPOSITIVEINCAUCHYSTABILITYREGION}
depend on the parameters of the background solution.
Similarly, since $\muPS|_{\CSregion_{Small}^{[0,5 \blowuptimePS]}} = 1$,
it is a standard consequence of Cauchy stability that the perturbed solution satisfies:
\begin{align} \label{E:MUISALMOSTUNITYINTRIVIALCAUCHYSTABILITYREGION}
	\upmu|_{\CSregion_{Small}^{[0,5 \blowuptimePS]}} 
	& = 1 + \mathcal{O}(\mathring{\Delta}_{\Sigma_0^{[-\farrightu,\leftu]}}^{\Ntop+1}).
\end{align}

Moreover, for $\muxmulevelsetvalue \in [0,\muxmulevelsetvalue_0]$ (where $\muxmulevelsetvalue_0$ is as in Theorem~\ref{T:PROPERTIESOFADMISSIBLESOLUTIONSLAUNCEDBYRESCALED}),
the perturbed rough time functions $\timefunctionarg{\muxmulevelsetvalue}$,
which solve the initial value problem
\eqref{E:TRANSPORTEQUATIONFORROUGHTIMEFUNCTION}--\eqref{E:INITIALCONDITIONFORROUGHTIMEFUNCTION},
also vary continuously (on the subset of $\CSregion^{\blowuptimePS;\PSBigDelta}$ where they are defined) 
with respect to the fluid data.
Furthermore, if $\mathring{\Delta}_{\Sigma_0^{[-\farrightu,\leftu]}}^{\Ntop+1}$ 
is small, then, like the map $\PSUpsilon(t,u) = (t,x^1)$ from Theorem~\ref{T:PROPERTIESOFADMISSIBLESOLUTIONSLAUNCEDBYRESCALED},
the perturbed change of variables map
$\Upsilon(t,u,x^2,x^3) = (t,x^1,x^2,x^3)$ 
defined in \eqref{E:CHOVGEOTOCARTESIAN} is a diffeomorphism on
$\CSregion^{\blowuptimePS;\PSBigDelta}$.
Hence, from the chain rule and the standard Sobolev calculus, it follows that for the perturbed solution,
we can view the following quantities as functions of the geometric coordinates $(t,u,x^2,x^3)$,
and for a neighborhood of the data of the background solution, 
they all vary continuously with respect to the fluid data  
in various Sobolev and Lebesgue norms corresponding to the geometric coordinates:
all of the fluid variables, 
the eikonal function,
all of the auxiliary quantities, such as $\upmu$, $\Lunit^i$, $\upchi$, etc.\
constructed out of the eikonal function (the auxiliary quantities solve various geometric PDEs),
the rough time functions,
and the Cartesian coordinates 
(which we view to be functions of the geometric coordinates, where those functions depend on the perturbed solution).
In particular, all the perturbed solutions
exist and are smooth on the common domain $\CSregion^{\blowuptimePS;\PSBigDelta}$
in geometric coordinate space.
Similarly, if $\mathring{\Delta}_{\Sigma_0^{[-\farrightu,\leftu]}}^{\Ntop+1}$
is sufficiently small, then
for $\muxmulevelsetvalue \in [0,\muxmulevelsetvalue_0]$,
the perturbed change of variables map
$
		\CHOVgeotorough{\muxmulevelsetvalue}(t,u,x^2,x^3) 
		= 
		(\timefunctionarg{\muxmulevelsetvalue},u,x^2,x^3)
$
defined in \eqref{E:CHOVGEOTOROUGH}
is a diffeomorphism from a subset of $\CSregion^{\blowuptimePS;\PSBigDelta}$
onto a region that contains
$[2 \timefunction_0,\timefunction_0/2] \times [-\farrightu,\leftu] \times \mathbb{T}^2$.

\begin{quote}
	In total, these Cauchy stability arguments yield all of the conclusions of
	Prop.\,\ref{P:CAUCHYSTABILITYANDEXISTENCEOFOPENSETS}, except for the following
	top-order energy estimates, which we will derive in the remaining steps:
	\begin{itemize}
		\item The $L^2$ estimates at the highest derivative level in
		\eqref{E:TANGENTIALL2NORMSOFWAVEVARIABLESSMALLALONGINITIALROUGHHYPERSURFACE}--\eqref{E:TANGENTIALL2NORMSOFMODIFIEDFLUIDVARIABLESSMALLALONGINITIALROUGHHYPERSURFACE},
	\eqref{E:WAVESARESMALLONINITIALNULLHYERSURFACE}--\eqref{E:MODIFIEDFLUIDVARIABLESARESMALLONINITIALNULLHYERSURFACE},
		\eqref{E:SMALLDATAOFPSIONINITIALROUGHTORI}--\eqref{E:SMALLDATAOFMODIFIEDFLUIDVARIABLESONINITIALROUGHTORI},
		and
		\eqref{E:L2DATASMALLNESSFOREIKONALFUNCTIONQUANTITIES}--\eqref{E:L2DATASMALLNESSFORCHIUPTOTOPORDER}.
	\end{itemize}
\end{quote}

\medskip

\noindent \textbf{Step 2: Estimates -- without derivative loss -- relative to the geometric coordinates in 
$\widehat{\CSregion}^{[0,\blowuptimePS - 2 \PSBigDelta]}$.}	
Note that \eqref{E:MAINCAUCHYSTABILITYREGION}--\eqref{E:SMALLCAUCHYSTABILITYREGION} 
and \eqref{E:BEFORESHOCKTIMECAUCHYSTABILITYREGION} imply that
$\widehat{\CSregion}^{[0,\blowuptimePS - 2 \PSBigDelta]}
= 
\bigcup_{t \in [0,\blowuptimePS - 2 \PSBigDelta]} \Sigma_t^{[-\farrightu + 3t,\leftu]}$.
In particular, 
$\widehat{\CSregion}^{[0,\blowuptimePS - 2 \PSBigDelta]}$
is foliated by portions of the flat hypersurfaces $\Sigma_t$
that are bounded by level sets of $u$.
Hence, we can control the solution in $\widehat{\CSregion}^{[0,\blowuptimePS - 2 \PSBigDelta]}$
by constructing
energies and null fluxes as in Sect.\,\ref{S:BASICINGREDIENTSFORL2ANALYSIS},
but instead of using the rough time functions and rough hypersurfaces,
we use the Cartesian time function $t$ and the spacelike hypersurface portions
$\Sigma_t^{[-\farrightu + 3t,u]}$ and $\csspacelikehypersurface^{[0,t]}$
(note that $\Sigma_t^{[-\farrightu + 3t,u]} \cup \csspacelikehypersurface^{[0,t]}$ is Lipschitz
and piecewise smooth, which is sufficient regularity for applying the divergence theorem on 
$\widehat{\CSregion}^{[0,\blowuptimePS - 2 \PSBigDelta]}$ to obtain $L^2$ estimates for the solution).
We refer to these as ``flat geometric'' energies and null fluxes.
The same arguments used in the proofs of
Props.\,\ref{P:APRIORIL2ESTIMATESWAVEVARIABLES},
\ref{P:MAINHYPERSURFACEENERGYESTIMATESFORTRANSPORTVARIABLES},
\ref{P:ROUGHTORIENERGYESTIMATES},
and
\ref{P:APRIORIL2ESTIMATESACOUSTICGEOMETRY}
show that analogous estimates 
also hold\footnote{To derive estimates for the acoustic geometry quantities $\upmu$, $\Lunit^i$, $\upchi$, etc.,
we in particular need to control their data on $\Sigma_0$.
To this end, we recall that $u|_{\Sigma_0} = - x^1$. Hence, $u|_{\Sigma_0}$ is $C^{\infty}$, 
and it is straightforward to use the eikonal equation \eqref{E:EIKONALEQUATION} 
and the definitions of $\upmu$, $\Lunit^i$, $\upchi$, etc.\
to ``solve'' for their initial data
on $\Sigma_0$ in terms of the fluid variable data on $\Sigma_0$;
the size of the corresponding initial data functions can then be controlled
using the standard Sobolev calculus.
\label{FN:DATAOFEIKONALFUNCTIONQUANTITIES}}  
for the flat geometric energies and null fluxes,
where the role of the smallness parameter $\initialsmall$
is now played by
$\mathring{\Delta}_{\Sigma_0^{[-\farrightu,\leftu]}}^{\Ntop+1}$.
The proof is dramatically easier in the present context because 
the Cartesian time function $t$ is much easier to control and because
$\upmu$ is uniformly positive by
\eqref{E:MUISUNIFORMLYPOSITIVEINCAUCHYSTABILITYREGION}. In particular, the
energy estimates can be derived using the standard version of
Gr\"{o}nwall's inequality, as opposed to the very technical
arguments used in e.g.\ Sect.\,\ref{SSS:PROOFOFAPRIORIL2ESTIMATESWAVEVARIABLES}.

There is, however, one new detail of significance that we now describe.
In carrying out the above arguments, which involve using Gr\"{o}nwall's inequality on
the sub-regions $\widehat{\CSregion}_{[0,t],[-\farrightu,u]}$ defined by:
\begin{align} \label{E:CAUCHYSTABILITYSUBREGIONS}
	\widehat{\CSregion}_{[0,t],[-\farrightu,u]} 
	& 
	\eqdef \bigcup_{t' \in [0,t]} \Sigma_{t'}^{[-\farrightu + 3t,u]}, 
\end{align}
we must
control the top-order
derivatives of $\vortrenormalized$ and $\GradEnt$ 
by using a Cartesian time function-analog of 
the integral identity \eqref{E:INTEGRALIDENTITYFORELLIPTICHYPERBOLICCURRENT}
on $\widehat{\CSregion}_{[0,t],[-\farrightu,u]}$.
Since the right boundary of $\widehat{\CSregion}_{[0,t],[-\farrightu,u]}$
is the $\gfour$-\emph{spacelike} hypersurface 
$\csspacelikehypersurface^{[0,t]}$
 -- in contrast to the $\gfour$-\emph{null} right boundary of the domain 
$\twoargMrough{[\timefunction_1,\timefunction_2),[u_1,u_2]}{\muxmulevelsetvalue}$
featured in \eqref{E:INTEGRALIDENTITYFORELLIPTICHYPERBOLICCURRENT} --
the new integral identity features additional $\csspacelikehypersurface^{[0,t]}$-integrals
and tori-integrals.
More precisely, for $(t,u) \in [0,\blowuptimePS - 2 \PSBigDelta] \times [-\farrightu,\leftu]$,
the same arguments that we used to prove \eqref{E:INTEGRALIDENTITYFORELLIPTICHYPERBOLICCURRENT}
can be used to prove the following similar identity,
where we have suppressed the volume and area forms to simplify the presentation,
and in practice, 
the role of the vectorfield $\SigmatTan$ is played by $\tander^{\Ntop} \vortrenormalized$ and $\tander^{\Ntop} \GradEnt$:
\begin{align}	
\begin{split} \label{E:FLATCAUCHYSTABILITYINTEGRALIDENTITYFORELLIPTICHYPERBOLICCURRENT}
		&
		\int_{\widehat{\CSregion}_{[0,t],[-\farrightu,u]}}
			\ellipticCoerciveQuadratic[\pmb{\partial} \SigmatTan,\pmb{\partial} \SigmatTan]
		+
		\int_{\ell_{t,u}}
			\frac{1}{4 \upmu}
			|\SigmatTan|_g^2
		-
		\int_{\ell_{t,-\farrightu + 3t}}
			\frac{1}{4 \upmu}
			|\SigmatTan|_g^2
		+ 
		\int_{\ell_{t,-\farrightu + 3t}}
			\frac{1}{4 (\upmu - \frac{1}{3})}
			|\SigmatTan|_g^2
			\\
		& =	
			\int_{\ell_{0,u}}
				\frac{1}{4 \upmu}
				|\SigmatTan|_g^2
			-
			\int_{\ell_{0,-\farrightu}}
				\frac{1}{4 \upmu}
				|\SigmatTan|_g^2
			+
			\int_{\ell_{0,-\farrightu}}
				\frac{1}{4 (\upmu - \frac{1}{3})}
				|\SigmatTan|_g^2
				    \\
		& \ \
			+
			\int_{\Sigma_t^{[-\farrightu + 3t,u]}}
				\cdots
			+
			\int_{\csspacelikehypersurface^{[0,t]}}
			-
			\int_{\Sigma_0^{[-\farrightu,u]}}
				\cdots
			+
			\int_{\widehat{\CSregion}_{[0,t],[-\farrightu,u]}}	
			\cdots.
	\end{split}
	\end{align}
	In \eqref{E:FLATCAUCHYSTABILITYINTEGRALIDENTITYFORELLIPTICHYPERBOLICCURRENT},
	``$\cdots$'' denotes error integrands that are
	similar to the ones on RHS~\eqref{E:INTEGRALIDENTITYFORELLIPTICHYPERBOLICCURRENT},
	but are much simpler to control because they depend on the Cartesian time function 
	(as opposed to the rough time function) and its derivatives;
	the error integrals of the ``$\cdots$'' terms can be controlled by the flat energies
	mentioned above, much like in the proof of Prop.\,\ref{P:ELLIPTICHYPERBOLICINTEGRALINEQUALITIES}.
	We point out that the Cauchy stability arguments from Step 1 yield that
	$\csspacelikehypersurface^{[0,t]}$ is $\gfour$-spacelike for the perturbed solution
	(because it is also spacelike for the background solution).
	
	The key feature of the identity \eqref{E:FLATCAUCHYSTABILITYINTEGRALIDENTITYFORELLIPTICHYPERBOLICCURRENT}
	is that 
	the tori integrals
	$	-
		\int_{\ell_{t,-\farrightu + 3t}}
			\frac{1}{4 \upmu}
			|\SigmatTan|_g^2
		+ 
		\int_{\ell_{t,-\farrightu + 3t}}
			\frac{1}{4 (\upmu - \frac{1}{3})}
			|\SigmatTan|_g^2
	$
	on LHS~\eqref{E:FLATCAUCHYSTABILITYINTEGRALIDENTITYFORELLIPTICHYPERBOLICCURRENT}
	sum to yield \emph{positive definite} control of
	$
	\int_{\ell_{t,-\farrightu + 3t}}
		|\SigmatTan|_g^2
	$,
	thus allowing us to propagate the ``torus regularity''
	corresponding to the tori integrals on
	the RHS of the definition
	\eqref{E:PERTURBATIONSMALLNESSINCARTESIANDIFFERENTIALSTRUCTURE}
	for $\mathring{\Delta}_{\Sigma_0^{[-\farrightu,\leftu]}}^{\Ntop+1}$
	(note that the integrals
	$\int_{\ell_{0,u}} \cdots$ and $\int_{\ell_{0,-\farrightu}} \cdots$
	on RHS~\eqref{E:FLATCAUCHYSTABILITYINTEGRALIDENTITYFORELLIPTICHYPERBOLICCURRENT}
	are data-terms that are controlled by $\mathring{\Delta}_{\Sigma_0^{[-\farrightu,\leftu]}}^{\Ntop+1}$).
	This coercive control allows us to prove an analog of
	Prop.\,\ref{P:ROUGHTORIENERGYESTIMATES} on the tori $\ell_{t,u}$
	using \underline{only} our smallness assumption on the bona fide data norm
	$\mathring{\Delta}_{\Sigma_0^{[-\farrightu,\leftu]}}^{\Ntop+1}$.
	We highlight that this coercive control stands
	in contrast to the \emph{rough tori} control
	provided by the identity \eqref{E:INTEGRALIDENTITYFORELLIPTICHYPERBOLICCURRENT}, where the 
	integral 
	$\int_{\twoargroughtori{\timefunction_2,u_1}{\muxmulevelsetvalue}}
				\CurrentboundaryerrorperfectRderivative[\SigmatTan,\SigmatTan]
			\, \volroughtorus$
	appears on the \emph{right-hand} side with a positive (unfavorable) sign,
	i.e., one needs a \emph{new estimate} to control 
	$\int_{\twoargroughtori{\timefunction_2,u_1}{\muxmulevelsetvalue}}
				\CurrentboundaryerrorperfectRderivative[\SigmatTan,\SigmatTan]
			\, \volroughtorus$.
	We derive the needed new estimate in Lemma~\ref{L:CONTROLOFDATAFORROUGHTORIENERGYESTIMATES}.
	However, we stress that, logically speaking, the proof of
	Lemma~\ref{L:CONTROLOFDATAFORROUGHTORIENERGYESTIMATES} can be completed only \emph{after}
	one has derived the Cauchy stability estimates of Prop.\,\ref{P:CAUCHYSTABILITYANDEXISTENCEOFOPENSETS},
	which are independent of Lemma~\ref{L:CONTROLOFDATAFORROUGHTORIENERGYESTIMATES}.
	We also clarify that the tori integrands
	$\frac{1}{4 \upmu} |\SigmatTan|_g^2$
	in \eqref{E:FLATCAUCHYSTABILITYINTEGRALIDENTITYFORELLIPTICHYPERBOLICCURRENT}
	reflect the fact that
	the integrand
	$
	\CurrentboundaryerrorperfectRderivative[\SigmatTan,\SigmatTan]
	$
	defined in \eqref{E:PERFECTRDERIVAIVEERRORPUTANGENTCURRENTCONTRACTEDAGAINSTVECTORFIELD}
	takes a simplified form when we use foliations by level sets of $t$
	instead of the rough time functions.
	Similarly, the tori integrands 
	$\frac{1}{4 (\upmu - \frac{1}{3})} |\SigmatTan|_g^2$
	in \eqref{E:FLATCAUCHYSTABILITYINTEGRALIDENTITYFORELLIPTICHYPERBOLICCURRENT}
	are analogs of the integrand
	$
	\CurrentboundaryerrorperfectRderivative[\SigmatTan,\SigmatTan]
	$
	defined in \eqref{E:PERFECTRDERIVAIVEERRORPUTANGENTCURRENTCONTRACTEDAGAINSTVECTORFIELD},
	but now corresponding to the time function $3 t - u$, whose $\farrightu$ level set defines 
	$\csspacelikehypersurface^{[0,5 \blowuptimePS]}$
	(see \eqref{E:FARRIGHTSPACELIKEBOUNDARYCAUCHYSTABILITYREGION}).

\medskip

\noindent \textbf{Step 3: Estimates -- without derivative loss -- relative to the geometric coordinates in 
$\CSregion_{Small}^{[0,5 \blowuptimePS]}$.}
Thanks to the positivity of $\upmu$ in $\CSregion_{Small}^{[0,5 \blowuptimePS]}$
guaranteed by \eqref{E:MUISALMOSTUNITYINTRIVIALCAUCHYSTABILITYREGION},
we can treat this region using essentially the same arguments we used in Step 2,
using foliations of $\CSregion_{Small}^{[0,5 \blowuptimePS]}$ 
by portions of flat hypersurfaces $\Sigma_t$ and 
portions of null hypersurfaces $\nullhyparg{u}$.

\medskip

\noindent \textbf{Recap}: We have now shown that
on the entire Cauchy stability region $\CSregion^{\blowuptimePS;\PSBigDelta}$,
the fluid variables' flat geometric energies and null fluxes
are bounded up to top order 
by $\lesssim \left(\mathring{\Delta}_{\Sigma_0^{[-\farrightu,\leftu]}}^{\Ntop+1} \right)^2$,
and that the same result holds for the acoustic geometry quantities (such as $\upmu$, $\Lunit^i$, and $\upchi$).
Since $\nullhyparg{- \rightu}^{[0,\frac{4}{\mathring{\updelta}_*}]} \subset \CSregion^{\blowuptimePS;\PSBigDelta}$,
this shows in particular that the data estimates
\eqref{E:WAVESARESMALLONINITIALNULLHYERSURFACE}--\eqref{E:MODIFIEDFLUIDVARIABLESARESMALLONINITIALNULLHYERSURFACE}
hold with
$\initialsmall \lesssim \mathring{\Delta}_{\Sigma_0^{[-\farrightu,\leftu]}}^{\Ntop+1}$.

\medskip
\noindent \textbf{What remains to be accomplished:}
The previous steps have provided estimates up to top order for the solution with respect to the geometric coordinates
on portions of the hypersurfaces $\Sigma_t$ of constant Cartesian time, null hypersurface portions of the form
$\nullhyparg{u}^{[0,t]}$,
and smooth tori $\ell_{t,u}$.
It remains for us to control the solution along suitable rough hypersurfaces
$\hypthreearg{\timefunction'}{[- \rightu,\leftu]}{\muxmulevelsetvalue}$,
null hypersurface portions of the form $\nullhypthreearg{\muxmulevelsetvalue}{u}{[\timefunction_1,\timefunction_2]}$,
and 
$
\twoargroughtori{\timefunction,u}{\muxmulevelsetvalue}
$.
Our arguments will rely on 
the previously derived estimates,
averaging methods, 
and some additional energy estimates
that are similar to, but much simpler than, the ones derived in the bulk of the paper.

\medskip

\noindent \textbf{Step 4: Estimates -- without derivative loss -- on a special rough hypersurface,
a special null hypersurface, and special rough tori via Chebychev's inequality.}
	Recall that our main goal is show that the perturbed solution induces
	data on the perturbed rough hypersurface
	$\hypthreearg{\timefunction_0}{[- \rightu,\leftu]}{\muxmulevelsetvalue}$,
	the null hypersurface portions $\nullhyparg{- \rightu}^{[0,\frac{4}{\mathring{\updelta}_*}]}$,
	and the rough tori $\twoargroughtori{\timefunction_0,u}{\muxmulevelsetvalue}$
	that satisfy all of the assumptions of
	Sects.\,\ref{SSS:QUANTITATIVEASSUMPTIONSONDATAAWAYFROMSYMMETRY}--\ref{SSS:LOCALIZEDDATAASSUMPTIONSFORMUANDDERIVATIVES}.
	In the present Step 4, we will combine the results of the previous steps 
	with averaging arguments based on Chebychev's inequality to show that
	the desired estimates hold on nearby rough hypersurfaces and rough tori. 
	Then, in Step 5,
	we will derive energy estimates, starting from the data on these nearby rough hypersurfaces and tori,
	showing that the estimates hold
	on
	$\hypthreearg{\timefunction_0}{[- \rightu,\leftu]}{\muxmulevelsetvalue}$,
	$\nullhyparg{- \rightu}^{[0,\frac{4}{\mathring{\updelta}_*}]}$,
	and $\twoargroughtori{\timefunction_0,u}{\muxmulevelsetvalue}$.

To proceed, we note that the standard Cauchy stability results yielded by Step 1 imply that
for $\muxmulevelsetvalue \in [0,\muxmulevelsetvalue_0]$, 
the region $\twoargMrough{[2 \timefunction_0,\timefunction_0/2],[- \farrightu,\leftu]}{\muxmulevelsetvalue}$
is contained in the region 
$\CSregion^{\blowuptimePS;\PSBigDelta}$ defined in \eqref{E:CAUCHYSTABILITYREGIONDECOMPOSTEDINTOMAINREGIONANDSMALLREGION}
(recall that $\timefunction_0 < 0$).
Hence, by Fubini's theorem, for any non-negative function $F$, we have
the following estimate, where in what follows, we suppress the area and volume forms to simplify the presentation:
\begin{align} \label{E:FUBINIPLUSBOUND}
\int_{\timefunction' = 2 \timefunction_0}^{\timefunction_0/2}
\int_{\hypthreearg{\timefunction'}{[-\farrightu,\leftu]}{\muxmulevelsetvalue}}
	F
\, d \timefunction'
&
=
\int_{u' = - \farrightu}^{\leftu}
\int_{\nullhypthreearg{\muxmulevelsetvalue}{u'}{[2 \timefunction_0,\timefunction_0/2]}}
	F
\, d u'
=
\int_{\twoargMrough{[2 \timefunction_0,\timefunction_0/2],[- \farrightu,\leftu]}{\muxmulevelsetvalue}}
	F
\leq 
\int_{\CSregion^{\blowuptimePS;\PSBigDelta}}
	F.
\end{align}
We now let (see Def.\,\ref{D:STRINGSOFCOMMUTATIONVECTORFIELDS} regarding the notation):
\begin{align} 
\begin{split} \label{E:SQUARESOFINTEGRANDSTHATARECONTROLLEDBYENERGIESANDELLIPTICESTIMATES}
	F
	& \eqdef
		|\comdersmall^{[1,\Ntop+1];1} \wavearray|^2
		+
		|\tander^{\leq \Ntop} (\vortrenormalized,\GradEnt)|^2
		+
		|\tander^{\leq \Ntop} (\VortVort,\DivGradEnt)|^2
		+
		|\pmb{\partial}  \tander^{\Ntop} (\vortrenormalized,\GradEnt)|^2
			\\
	&
		+
		\sum_{a=1}^3 |\comdersmall^{[1,\Ntop];1} \Lunit^a|^2
		+
		|\tandersmall^{[1,\Ntop]} \upmu|^2
		+
		|\tandersmall^{\leq \Ntop} \mytr_{\gtorus} \upchi|^2
		+
		|\tandersmall^{\leq \Ntop} \mytr_{\gtorus} \upchi|^2
		+
		|\angLie_{\tander}^{\Ntop} \upchi|_{\gtorus}^2.
\end{split}
\end{align}
Note that $F$ is precisely the quantity that we have controlled
in $\CSregion^{\blowuptimePS;\PSBigDelta}$ by
our flat geometric energies and null fluxes
and the integral identity \eqref{E:FLATCAUCHYSTABILITYINTEGRALIDENTITYFORELLIPTICHYPERBOLICCURRENT}.
In particular, the arguments provided by Steps 2 and 3 imply\footnote{Actually, aside from
the term $|\pmb{\partial}  \tander^{\Ntop} (\vortrenormalized,\GradEnt)|^2$, the estimates
from Steps 2 and 3 show that various energies of the terms in $F$ on 
portions of the hypersurfaces $\Sigma_t$ and $\nullhyparg{u}$
are bounded by $\lesssim \left(\mathring{\Delta}_{\Sigma_0^{[-\farrightu,\leftu]}}^{\Ntop+1} \right)^2$.
These hypersurface estimates imply the spacetime estimate \eqref{E:FLASTSPACETIMEINTEGRALOFEVERYTHINGSMALLINCSREGION},
where the implicit constants in \eqref{E:FLASTSPACETIMEINTEGRALOFEVERYTHINGSMALLINCSREGION}
depend on the size of the region $\CSregion^{\blowuptimePS;\PSBigDelta}$
(which is compact with dimensions controlled by the background solution).
\label{FN:HYPERSURFACEESTIMATESIMPLYSPACETIMEESTIMATES}}
the following \emph{spacetime} integral estimate:
\begin{align} \label{E:FLASTSPACETIMEINTEGRALOFEVERYTHINGSMALLINCSREGION}
\int_{\CSregion^{\blowuptimePS;\PSBigDelta}}
	F
& \lesssim 
\left(\mathring{\Delta}_{\Sigma_0^{[-\farrightu,\leftu]}}^{\Ntop+1} \right)^2.
\end{align}

We now use \eqref{E:FUBINIPLUSBOUND}
and \eqref{E:FLASTSPACETIMEINTEGRALOFEVERYTHINGSMALLINCSREGION} to deduce that
$
\int_{\timefunction' = 2 \timefunction_0}^{\timefunction_0/2}
\int_{\hypthreearg{\timefunction'}{[-\farrightu,\leftu]}{\muxmulevelsetvalue}}
	F
\, d \timefunction'
\lesssim 
\left(\mathring{\Delta}_{\Sigma_0^{[-\farrightu,\leftu]}}^{\Ntop+1} \right)^2$.
From this bound and Chebychev's inequality, 
we see that there must exist a number $\timefunction_*$ satisfying:
\begin{align} \label{E:LOCATIONOFTIMEFUNCTIONSTAR}
	\timefunction_* \in [2 \timefunction_0, (3/2) \timefunction_0]
\end{align}
such that the following estimate holds:
\begin{align} \label{E:ENERGIESAREBOUNDEDONUNKNOWNHYPERSURFACE}
\int_{\hypthreearg{\timefunction_*}{[-\farrightu,\leftu]}{\muxmulevelsetvalue}}
	F
& \lesssim
\frac{1}{|\timefunction_0|}
\left(\mathring{\Delta}_{\Sigma_0^{[-\farrightu,\leftu]}}^{\Ntop+1} \right)^2
\lesssim 
\left(\mathring{\Delta}_{\Sigma_0^{[-\farrightu,\leftu]}}^{\Ntop+1} \right)^2.
\end{align}
Using Fubini's theorem again,
\eqref{E:ENERGIESAREBOUNDEDONUNKNOWNHYPERSURFACE}, 
and the fact that $[- \rightu - \frac{1}{\blowupdeltaPS},- \rightu] \subset [-\farrightu,\leftu]$,
we further deduce that:
\begin{align} \label{E:FUBINIFORENERGIESAREBOUNDEDONUNKNOWNHYPERSURFACE}
\int_{u' = - \rightu - \frac{1}{\blowupdeltaPS}}^{- \rightu}
\int_{\roughtori{\timefunction_*,u'}}
	F
\, \mathrm{d} u'
& =
\int_{\hypthreearg{\timefunction_*}{[- \rightu - \frac{1}{\blowupdeltaPS},- \rightu]}{\muxmulevelsetvalue}}
	F
\leq
\int_{\hypthreearg{\timefunction_*}{[-\farrightu,\leftu]}{\muxmulevelsetvalue}}
	F
\lesssim 
\left(\mathring{\Delta}_{\Sigma_0^{[-\farrightu,\leftu]}}^{\Ntop+1} \right)^2.
\end{align}
From \eqref{E:FUBINIFORENERGIESAREBOUNDEDONUNKNOWNHYPERSURFACE}
and Chebychev's inequality,
we find that there is a $u_* \in [- \rightu - \frac{1}{\blowupdeltaPS},- \rightu]$
such that:
\begin{align} \label{E:ENERGIESAREBOUNDEDONUNKNOWNTORUS}
		\int_{\roughtori{\timefunction_*,u_*}}
			|\tander^{\leq \Ntop} (\Omega,\GradEnt)|^2
			&
			\lesssim
			\frac{1}{\blowupdeltaPS}
			\left(\mathring{\Delta}_{\Sigma_0^{[-\farrightu,\leftu]}}^{\Ntop+1} \right)^2
			\lesssim
			\left(\mathring{\Delta}_{\Sigma_0^{[-\farrightu,\leftu]}}^{\Ntop+1} \right)^2.
\end{align}
From \eqref{E:ENERGIESAREBOUNDEDONUNKNOWNTORUS},
the bound
$
\int_{\hypthreearg{\timefunction_*}{[-\farrightu,\leftu]}{\muxmulevelsetvalue}}
			|\pmb{\partial} \tander^{\leq \Ntop} \Omega|^2
			\lesssim 
			\left(\mathring{\Delta}_{\Sigma_0^{[-\farrightu,\leftu]}}^{\Ntop+1} \right)^2
$
implied by \eqref{E:ENERGIESAREBOUNDEDONUNKNOWNHYPERSURFACE},
the fact that $[- \rightu - \frac{1}{\blowupdeltaPS},\leftu] \subset [-\farrightu,\leftu]$,
and fundamental theorem of calculus-type 
arguments similar to the ones we used to prove
\eqref{E:ROUGHTORUSINNULLHYPERSURFACEL2FUNDAMENTALTHEOREMOFCALCULUSESTIMATE}
-- but based on the identity \eqref{E:IDENTITYROUGHEIKONALDERIVATIVEOFROUGHTORUSINTEGRAL} 
and Gr\"{o}nwall's inequality with respect to $u$ --
we further deduce that:
\begin{align} \label{E:ENERGIESAREBOUNDEDONALLUNKNOWNTORIALONGROUGHHYPERSURFACE}
		\sup_{u \in [- \rightu - \frac{1}{\blowupdeltaPS},\leftu]}
		\int_{\roughtori{\timefunction_*,u}}
			|\tander^{\leq \Ntop} (\Omega,\GradEnt)|^2
			&
			\lesssim
			\left(\mathring{\Delta}_{\Sigma_0^{[-\farrightu,\leftu]}}^{\Ntop+1} \right)^2.
\end{align}
Moreover, the same arguments yield:
\begin{align} \label{E:WAVEVARIABLEENERGIESAREBOUNDEDONALLUNKNOWNTORIALONGROUGHHYPERSURFACE}
		\sup_{u \in [- \rightu - \frac{1}{\blowupdeltaPS},\leftu]}
		\int_{\roughtori{\timefunction_*,u}}
			|\tander^{[1, \Ntop]} \wavearray|^2
			&
			\lesssim
			\left(\mathring{\Delta}_{\Sigma_0^{[-\farrightu,\leftu]}}^{\Ntop+1} \right)^2.
\end{align}

Similarly, since 
\eqref{E:FUBINIPLUSBOUND}
and \eqref{E:FLASTSPACETIMEINTEGRALOFEVERYTHINGSMALLINCSREGION}
yield
$
\int_{u' = - \rightu - \frac{1}{\blowupdeltaPS}}^{- \rightu}
\int_{\nullhypthreearg{\muxmulevelsetvalue}{u'}{[2 \timefunction_0,\timefunction_0/2]}}
	F
\, \mathrm{d} u'
\lesssim
\left(\mathring{\Delta}_{\Sigma_0^{[-\farrightu,\leftu]}}^{\Ntop+1} \right)^2
$,
we can again use Chebychev's inequality
and fundamental theorem of calculus-type
arguments similar to the ones we used to prove
\eqref{E:ROUGHTORUSINNULLHYPERSURFACEL2FUNDAMENTALTHEOREMOFCALCULUSESTIMATE}
to deduce that there exists
a $U_* > 0$ such that:\footnote{Note that \eqref{E:MINUSUSTARTLOCATION} implies that
$U_* > \rightu$. We will use this basic fact in Step~2 of our proof of Prop.\,\ref{P:CONTINUATIONCRITERIA}.
\label{FN:USTARISBIGGERTHANRIGHTU}}
\begin{align} \label{E:MINUSUSTARTLOCATION}
	- U_* 
	& \in [- \rightu - \frac{1}{\blowupdeltaPS}, - \rightu - \frac{1}{2 \blowupdeltaPS}]
\end{align}
such that:
\begin{align} \label{E:ENERGIESAREBOUNDEDONUNKNOWNNULLHYPERSURFACE}
	\int_{\nullhypthreearg{\muxmulevelsetvalue}{-U_*}{[2 \timefunction_0,\timefunction_0/2]}}
		F
	& \lesssim
	\left(\mathring{\Delta}_{\Sigma_0^{[-\farrightu,\leftu]}}^{\Ntop+1} \right)^2
\end{align}
and:
\begin{align} \label{E:ENERGIESAREBOUNDEDONALLUNKNOWNTORIALONGNULLHYPERSURFACE}
		\sup_{\timefunction \in [2 \timefunction_0,\timefunction_0/2]}
		\int_{\roughtori{\timefunction,-U_*}}
			|\tander^{\leq \Ntop} (\Omega,\GradEnt)|^2
			&
			\lesssim
			\left(\mathring{\Delta}_{\Sigma_0^{[-\farrightu,\leftu]}}^{\Ntop+1} \right)^2.
\end{align}

\noindent \textbf{Step 5: The desired geometric energy estimates -- without derivative loss -- 
on the rough hypersurfaces, null hypersurfaces, and rough tori
relative to the geometric coordinates in 
$\twoargMrough{[\timefunction_*,\timefunction_0],[-U_*,\leftu]}{\muxmulevelsetvalue}$.}
In total, the arguments given in Steps 1-4 have shown that for $\muxmulevelsetvalue \in [0,\muxmulevelsetvalue_0]$,
the bona fide initial data on 
$\Sigma_0^{[-\farrightu,\leftu]}$
induce data on the rough hypersurface portion
$
\hypthreearg{\timefunction_*}{[-U_*,\leftu]}{\muxmulevelsetvalue}
$,
the null hypersurface portion
$\nullhypthreearg{\muxmulevelsetvalue}{-U_*}{[2 \timefunction_0,\timefunction_0/2]}$,
the rough tori $\twoargroughtori{\timefunction_*,u}{\muxmulevelsetvalue}$ for $u \in [-U_*,\leftu]$,
and the rough tori $\twoargroughtori{\timefunction,-U_*}{\muxmulevelsetvalue}$ for $\timefunction \in [2 \timefunction_0,\timefunction_0/2]$,
such that on these surfaces, 
all of the energies and null fluxes (up to top-order) 
defined in Sect.\,\ref{SS:FUNDAMENTALL2CONTROLLINGQUANTITIES}
are bounded by
$\lesssim \left(\mathring{\Delta}_{\Sigma_0^{[-\farrightu,\leftu]}}^{\Ntop+1} \right)^2$.
Starting from these ``data estimates''
(including the ones on the rough tori provided by 
\eqref{E:ENERGIESAREBOUNDEDONALLUNKNOWNTORIALONGROUGHHYPERSURFACE},
\eqref{E:WAVEVARIABLEENERGIESAREBOUNDEDONALLUNKNOWNTORIALONGROUGHHYPERSURFACE},
and \eqref{E:ENERGIESAREBOUNDEDONALLUNKNOWNTORIALONGNULLHYPERSURFACE}),
and using the same arguments we used in the proofs of
\eqref{E:LOSSOFONEDERIVATIVEL2ESTIMATESFORWAVEVARIABLESONROUGHTORIINTERMSOFDATAANDCONTROLLING},
Props.\,\ref{P:APRIORIL2ESTIMATESWAVEVARIABLES},
\ref{P:MAINHYPERSURFACEENERGYESTIMATESFORTRANSPORTVARIABLES},
\ref{P:ROUGHTORIENERGYESTIMATES},
and
\ref{P:APRIORIL2ESTIMATESACOUSTICGEOMETRY},
we can derive the same geometric energy estimates on the region
$
\twoargMrough{[\timefunction_*,\timefunction_0/2],[-U_*,\leftu]}{\muxmulevelsetvalue}
$,
i.e., we can bound the geometric energies up to top order on
$
\hypthreearg{\timefunction}{[-U_*,\leftu]}{\muxmulevelsetvalue}
$
for $\timefunction \in [\timefunction_*,\timefunction_0/2]$,
the geometric null fluxes up to top order on 
$
\nullhypthreearg{\muxmulevelsetvalue}{u}{[\timefunction_*,\timefunction_0/2]}
$
for $u \in [-U_*,\leftu]$, 
and the geometric rough tori energies (as in Prop.\,\ref{P:ROUGHTORIENERGYESTIMATES})
on $\twoargroughtori{\timefunction,u}{\muxmulevelsetvalue}$ up to top order
for $(\timefunction,u) \in [\timefunction_*,\timefunction_0/2] \times [-U_*,\leftu]$;
all of these quantities are bounded by $\lesssim \left(\mathring{\Delta}_{\Sigma_0^{[-\farrightu,\leftu]}}^{\Ntop+1} \right)^2$,
e.g.,
\begin{align} \label{E:INCSSMALLREGIONENERGIESAREBOUNDEDONALLTORIINTEGRALSALONDATANULLHYPERSURFACEUPTOBOOSTRAPTIME}
		\sup_{(\timefunction,u) \in [\timefunction_*,\timefunction_0/2] \times [-U_*,\leftu]}
		\int_{\roughtori{\timefunction,u}}
			|\tander^{\leq \Ntop} (\Omega,\GradEnt)|^2
		\, \volroughtorus
		&
		\lesssim
		\left(\mathring{\Delta}_{\Sigma_0^{[-\farrightu,\leftu]}}^{\Ntop+1} \right)^2.
\end{align}

\begin{remark}[The proof of \eqref{E:ENERGIESAREBOUNDEDONALLUNKNOWNTORIALONGNULLHYPERSURFACE} 
	does not rely on Lemma~\ref{L:CONTROLOFDATAFORROUGHTORIENERGYESTIMATES}]
	\label{R:ROUGHTORIESTIMATESONPUSTARAREINDEPENDENTOFLEMMACONTROLOFDATAFORROUGHTORIENERGYESTIMATES}
	The proof of \eqref{E:INCSSMALLREGIONENERGIESAREBOUNDEDONALLTORIINTEGRALSALONDATANULLHYPERSURFACEUPTOBOOSTRAPTIME}
	relies in particular on the integral identity \eqref{E:INTEGRALIDENTITYFORELLIPTICHYPERBOLICCURRENT}
	with 
	$\timefunction_1 = \timefunction_*$,
	$\timefunction_2 \in [\timefunction_*,\timefunction_0/2]$,
	$u_1 = - U_*$,
	and
	$u_2 \in [-U_*,\leftu]$.
	The rough tori $L^2$ estimates 
	\eqref{E:WAVEVARIABLEENERGIESAREBOUNDEDONALLUNKNOWNTORIALONGROUGHHYPERSURFACE}
	and
	\eqref{E:ENERGIESAREBOUNDEDONALLUNKNOWNTORIALONGNULLHYPERSURFACE}
	are needed to control the corresponding rough tori integrals
	$
	\int_{\twoargroughtori{\timefunction_2,- \rightu}{\muxmulevelsetvalue}}
				\cdots
	$,
	$
			\int_{\twoargroughtori{\timefunction_*,u_2}{\muxmulevelsetvalue}}
				\cdots
	$,
	and
	$
	\int_{\twoargroughtori{\timefunction_*,- \rightu}{\muxmulevelsetvalue}}
				\cdots
$
	on RHS~\eqref{E:INTEGRALIDENTITYFORELLIPTICHYPERBOLICCURRENT}.
	In particular, the rough tori estimates \eqref{E:ENERGIESAREBOUNDEDONALLUNKNOWNTORIALONGNULLHYPERSURFACE}
	provide an analog of the estimates of Lemma~\ref{L:CONTROLOFDATAFORROUGHTORIENERGYESTIMATES}
	that are relevant for the region under study here.
	We stress that our proof of \eqref{E:ENERGIESAREBOUNDEDONALLUNKNOWNTORIALONGNULLHYPERSURFACE}
	given above is independent of Lemma~\ref{L:CONTROLOFDATAFORROUGHTORIENERGYESTIMATES};
	this is important for the logic of the paper.
\end{remark}

Since $\timefunction_0 \in (\timefunction_*,\timefunction_0/2)$,
particular cases of these bounds are the ones along
$
\hypthreearg{\timefunction_0}{[-U_*,\leftu]}{\muxmulevelsetvalue}
$,
which, in view of the fact that $[- \rightu,\leftu] \subset [-U_*,\leftu]$,
imply the data bounds
\eqref{E:TANGENTIALL2NORMSOFWAVEVARIABLESSMALLALONGINITIALROUGHHYPERSURFACE}--\eqref{E:TANGENTIALL2NORMSOFMODIFIEDFLUIDVARIABLESSMALLALONGINITIALROUGHHYPERSURFACE} 
and
\eqref{E:L2DATASMALLNESSFOREIKONALFUNCTIONQUANTITIES}--\eqref{E:L2DATASMALLNESSFORCHIUPTOTOPORDER}
with
$\initialsmall \lesssim \mathring{\Delta}_{\Sigma_0^{[-\farrightu,\leftu]}}^{\Ntop+1}$.
In particular, the analog of Prop.\,\ref{P:ROUGHTORIENERGYESTIMATES} yields:
\begin{align} \label{E:ENERGIESAREBOUNDEDONALLTORIINTEGRALSALONDATAROUGHHYPERSURFACE}
		\max_{u \in [-U_*,\leftu]}
		\int_{\twoargroughtori{\timefunction_0,u}{\muxmulevelsetvalue}}
		|\tander^{\leq \Ntop} (\Omega,\GradEnt)|^2
		&
		\lesssim
		\left(\mathring{\Delta}_{\Sigma_0^{[-\farrightu,\leftu]}}^{\Ntop+1} \right)^2,
\end{align}
which, in view of the fact that $[- \rightu,\leftu] \subset [-U_*,\leftu]$,
implies \eqref{E:SMALLDATAOFVORTICITYANDENTORPYGRADIENTONINITIALROUGHTORI}
with $\initialsmall \lesssim \mathring{\Delta}_{\Sigma_0^{[-\farrightu,\leftu]}}^{\Ntop+1}$.
Similarly, 
from the data estimate \eqref{E:WAVEVARIABLEENERGIESAREBOUNDEDONALLUNKNOWNTORIALONGROUGHHYPERSURFACE},
the same arguments we used to prove
\eqref{E:LOSSOFONEDERIVATIVEL2ESTIMATESFORWAVEVARIABLESONROUGHTORIINTERMSOFDATAANDCONTROLLING},
and the geometric energy and null-flux estimates,
we conclude the rough tori $L^2$ estimates \eqref{E:SMALLDATAOFPSIONINITIALROUGHTORI} for $\tander^{[1, \Ntop]} \wavearray$.
Finally, we will derive the rough tori $L^2$ estimates \eqref{E:SMALLDATAOFMODIFIEDFLUIDVARIABLESONINITIALROUGHTORI} 
for $\tander^{\leq \Ntop-1}(\VortVort,\DivGradEnt)$.
To this end, we note that the geometric energy estimates mentioned above
imply the following bounds:
$
\max_{u \in [-\rightu,\leftu]}
\int_{\nullhypthreearg{\muxmulevelsetvalue}{u}{[\timefunction_0,\timefunction_0/2]}}
	|\tander^{\leq \Ntop-1}(\VortVort,\DivGradEnt)|^2
\lesssim
\left(\mathring{\Delta}_{\Sigma_0^{[-\farrightu,\leftu]}}^{\Ntop+1} \right)^2
$.
From this bound
and arguments similar to the ones we used to prove 
\eqref{E:ENERGIESAREBOUNDEDONALLUNKNOWNTORIALONGNULLHYPERSURFACE},
based on Chebychev's inequality
and fundamental theorem of calculus-type estimates
(cf.\ \eqref{E:ROUGHTORUSINNULLHYPERSURFACEL2FUNDAMENTALTHEOREMOFCALCULUSESTIMATE}),
we find that for any $u \in [-\rightu,\leftu]$, we have
$
\max_{\timefunction \in [\timefunction_0,\timefunction_0/2]}
		\int_{\twoargroughtori{\timefunction,u}{\muxmulevelsetvalue}}
			|\tander^{\leq \Ntop - 1} (\VortVort,\DivGradEnt)|^2
			\lesssim
			\left(\mathring{\Delta}_{\Sigma_0^{[-\farrightu,\leftu]}}^{\Ntop+1} \right)^2
$.
In particular, this bound implies
\eqref{E:SMALLDATAOFMODIFIEDFLUIDVARIABLESONINITIALROUGHTORI} 
with $\initialsmall \lesssim \mathring{\Delta}_{\Sigma_0^{[-\farrightu,\leftu]}}^{\Ntop+1}$.

We have therefore derived 
\eqref{E:TANGENTIALL2NORMSOFWAVEVARIABLESSMALLALONGINITIALROUGHHYPERSURFACE}--\eqref{E:TANGENTIALL2NORMSOFMODIFIEDFLUIDVARIABLESSMALLALONGINITIALROUGHHYPERSURFACE},
	\eqref{E:WAVESARESMALLONINITIALNULLHYERSURFACE}--\eqref{E:MODIFIEDFLUIDVARIABLESARESMALLONINITIALNULLHYERSURFACE},
		\eqref{E:SMALLDATAOFPSIONINITIALROUGHTORI}--\eqref{E:SMALLDATAOFMODIFIEDFLUIDVARIABLESONINITIALROUGHTORI},
		and
		\eqref{E:L2DATASMALLNESSFOREIKONALFUNCTIONQUANTITIES}--\eqref{E:L2DATASMALLNESSFORCHIUPTOTOPORDER}
		with $\initialsmall \lesssim \mathring{\Delta}_{\Sigma_0^{[-\farrightu,\leftu]}}^{\Ntop+1}$,
		thereby completing our proof sketch of Prop.\,\ref{P:CAUCHYSTABILITYANDEXISTENCEOFOPENSETS}.

\end{proof}

\section{Notation} \label{A:NOTATION}
For the reader's convenience, in this appendix, we have gathered some of the notation and conventions used throughout the paper.

\begin{center}
	\begin{tabular}{ | m{3cm} | m{10cm}| m{3cm} | } 
		\hline
		\multicolumn{3}{|c|}{\textbf{Cartesian coordinate space and derivatives}} \\
		\hline 
		Symbols & Descriptions & Reference \\ 
		\hline
		$\R \times \R \times \T^2$ & The ambient $(1+3)$-dimensional spacetime. The spatial domain is $\R\times\T^2$ & 
			Sect.\,\ref{SS:FIRSTVERSIONOFEQUATIONS} \\
		\hline
		$\{x^\alpha\}_{\alpha = 0,1,2,3}$ & The Cartesian coordinates on 
		$\R\times\R\times\T^2$, relative to which the Minkowski metric $m$ has components 
		$m_{\alpha\beta} = \diag(-1,1,1,1)$. We often use the notation $t \eqdef x^0$ for Cartesian time. & Sect.\,\ref{SS:FIRSTVERSIONOFEQUATIONS} \\ 
		\hline
		$\{\p_\alpha\}_{\alpha = 0,1,2,3}$ & Cartesian coordinate partial derivative vectorfields. Note that $\{\p_2,\p_3\}$ can be extended to a globally defined positively oriented frame on $\T^2$. & Sect.\,\ref{SS:FIRSTVERSIONOFEQUATIONS} \\
		\hline
		$\pmb{\partial} f$ & The array of Cartesian coordinate spacetime partial derivatives of a function $f$, i.e.,
		$\pmb{\partial} f \eqdef (\partial_t f, \partial_1, \partial_2 f, \partial_3 f)$ & Sect.\,\ref{SSS:INTRONEWFORMULATION} \\
		\hline
		$\partial f$ & The array of Cartesian coordinate spatial partial derivatives of $f$, i.e., 
		$\p f \eqdef (\partial_1, \partial_2 f, \partial_3 f)$. & 
			Sect.\,\ref{SSS:INTRONEWFORMULATION} \\
		\hline
		\multicolumn{3}{|C{16cm}|}{Lowercase Greek ``spacetime'' indices, such as $\alpha$, vary over $0,1,2,3$,
		while Lowercase Latin ``spatial'' indices, such as $a$, vary over $1,2,3$. 
		We use Einstein's summation convention in that repeated indices are summed over their respective ranges.} \\
		\hline
		$\Sigma_t$ & The hypersurface constant Cartesian time $t$: $\Sigma_t \eqdef \{(x^0,x^1,x^2,x^3) \in \R\times\R\times\T^2 \ | \ x^0 \equiv t\}$ &  Sect.\,\ref{SS:FIRSTVERSIONOFEQUATIONS} \\ 
		\hline 
		$\Flatcurl V$ & The Euclidean curl of a $\Sigma_t$-tangent vectorfield $V$, with components $(\Flatcurl V)^k \eqdef \upepsilon_{ijk}\updelta^{jl} \p_{l} V^k$, where $\upepsilon_{ijk}$ denotes the fully antisymmetric symbol normalized by $\upepsilon_{123} = 1$ and $\updelta^{jl}$ is the Kronecker delta. & Def.\,\ref{D:EUCLIDEANDIVERGENCEANDCURL} \\
		\hline 
		$\Flatdiv V$ & The Euclidean divergence  of a $\Sigma_t$-tangent vectorfield $V$, $\Flatdiv V \eqdef \p_a V^a$. & 
			Def.\,\ref{D:EUCLIDEANDIVERGENCEANDCURL} \\
		\hline

	\end{tabular}
\end{center}

\begin{center}
	\begin{tabular}{ | m{3cm} | m{10cm}| m{3cm} | } 
		\hline
		\multicolumn{3}{|c|}{\textbf{Fluid variables and some geometric objects}} \\
		\hline 
		Symbols & Descriptions & Reference \\ 
		\hline 
		$v^i,\Density,\Ent, p$ & $\{v^i\}_{i=1,2,3}$ is the $\Sigma_t$-tangent fluid velocity, $\Density$ is the fluid density, and $\Ent$ is the entropy. $p = p(\Density,\Ent)$ denotes the equation of state. & Sect.\,\ref{SS:FIRSTVERSIONOFEQUATIONS} \\
		\hline
		$\Transport$ & The material vectorfield. & \eqref{E:MATERIALDERIVATIVEVECOTRFIELD} \\
		\hline
		$\overline{\varrho}$ & A fixed positive constant ``background density.'' & \eqref{E:BACKGROUNDDENSITY} \\
		\hline 
		$\LogDensity, \Speed$ & $\LogDensity = \ln\left(\frac{\Density}{\overline{\Density}}\right)$ denotes the logarithmic density. The speed of sound is $\Speed(\LogDensity,\Ent)
		= \sqrt{(\overline{\varrho})^{-1} \exp(-\LogDensity) p;_{\LogDensity}}$,
		where $p;_{\LogDensity} \eqdef \tfrac{\p p}{\p \LogDensity}$ denotes the derivative of the equation of state 
		with respect to the logarithmic density at fixed $\Ent$. & Sect.\,\ref{SSS:LOGDENSITYASSUMPTIONSONEOSANDNORMALIZATIONS} \\
		\hline
		$\RRiemann$, $\LRiemann$ & The almost Riemann invariants. & \eqref{E:ALMOSTRIEMANNINVARIANTS} \\
		\hline
		$\Omega^i$, $\GradEnt^i$, $\VortVort^i$, $\DivGradEnt$ & $\{\Omega^i\}_{i=1,2,3}$ denotes the specific vorticity, $\{\GradEnt^i\}_{i=1,2,3}$ denotes the entropy gradient, and $\{\VortVort^i\}_{i=1,2,3},\DivGradEnt$ denote the modified fluid variables. & Def.\,\ref{D:HIGHERORDERFLUIDVARIABLES} \\
		\hline 
		$\wavearray$, $\wavearraypartial$ & The array and partial array of wave variables. & \eqref{E:ARRAYOFWAVEVARIABLES}--\eqref{E:PARTIALWAVEARRAY} \\
		\hline
		$\gfour$, $\gfour^{-1}$ & The acoustical metric and inverse metric of spacetime. & Def.\,\ref{D:ACOUSTICALMETRICDEF}. \\
		\hline
		$ G_{\alpha\beta}^{\iota}$, $\vec{G}_{\alpha\beta}$ & $\wavearray_\iota$-derivatives of $\gfour_{\alpha\beta}$ and the array of $\wavearray_\iota$-derivatives of $\gfour_{\alpha\beta}$. & Def.\,\ref{E:DERIVATIVESOFMETRICWRTFLUID} \\
		\hline
		$ D\wavearray$, $\vec{G}_{V_1V_2}\diamond D\wavearray$ & Differential operators involving $\wavearray$. & Def.\,\ref{D:DERIVATIVESOFARRAYS} \\
		\hline
		$\Box_{\gfour}$ & The covariant wave operator of $\gfour$. & Def.\,\ref{D:COVWAVEOP} \\
		\hline 
		$ \mathfrak{Q}^{(\gfour)}, \mathfrak{Q}_{\alpha\beta}$ & Standard $\gfour$-null forms. & Def.\,\ref{D:STANDARDNULLFORMS} \\
		\hline
	\end{tabular}
\end{center}

\begin{center}
	\begin{tabular}{ | m{3cm} | m{10cm}| m{3cm} | } 
		\hline
		\multicolumn{3}{|c|}{\textbf{The acoustic geometry}} \\
		\hline 
		Symbols & Descriptions & Reference \\ 
		\hline 
		$u$, $\upmu$ & The eikonal function and inverse foliation density. & Def.\,\ref{D:EIKONALFUNCTIONANDMU} \\
		\hline
		$\Sigma_t$, $\nullhyparg{u}$, $\ell_{t,u}$, $\Sigma_{t}^{[u_1,u_2]}$, $  \nullhyparg{u}^{[t_1,t_2]}$ & Acoustic regions and truncated acoustic regions of spacetime. & Def.\,\ref{D:ACOUSTICSUBSETSOFSPACETIME} \\
		\hline
		$g$, $\gtorus$ & First fundamental forms of $\Sigma_t$ and $\ell_{t,u}$ induced by $\gfour$, respectively. & Def.\,\ref{D:FIRSTFUNDAMENTALFORMS} \\
		\hline 
		$\Dfour$, $\angrmD$, $\angdiv$, $\angLap$ & Levi-Civita connections of $g$, $\gtorus$ and associated differential operators. & Def.\,\ref{D:CONNECTIONSANDDIFFERENTIALOPERATORS} \\
		\hline
		$\Riemfour$, $\Ricfour$ & The Riemann and Ricci curvature tensors of $\gfour$. & \eqref{E:ACOUSTICALCURVATURETENSOR}--\eqref{E:SPACETIMERICCITENSOR} \\
		\hline 
		$(t,u,x^2,x^3)$ & The geometric coordinate system. & Def.\,\ref{D:GEOMETRICCOORDIANTESANDPARTIALDERIVATIVEVECTORFIELDS}. \\
		\hline
		$ \{ \geop{t},\geop{u},\geop{x^2},\geop{x^A}\}$ &  The coordinate partial derivative vectorfields in geometric coordinates.  & Def.\,\ref{D:GEOMETRICCOORDIANTESANDPARTIALDERIVATIVEVECTORFIELDS}. \\
		\hline
		$\frac{\p^{\vec{\alpha}}}{\p(t,u,x^2,x^3)}$ & Multi-index notation in geometric coordinates. & \eqref{E:GEOMETRICCOORDINATEMULTIINDEXDIFFERENTIALOPERATOR} \\
		\hline 
		\multicolumn{3}{|C{16cm}|}{Uppercase Latin indices such as $A$ vary over 2,3. We use Einstein's summation convention in that repeated indices are summed over their respective ranges.} \\
		\hline
		\multicolumn{3}{|C{16cm}|}{The contraction of a one-form $\upxi$ and $\geop{x^A}$ is denoted by $\upxi_A \eqdef \upxi_\alpha (\geop{x^A})^{\alpha}$ for $A = 2,3$. A similar convention is used for higher order contractions, i.e. 
		$\gtorus_{AB} = \gtorus(\geop{x^A},\geop{x^B})$ for $A,B = 2,3$.} \\
		\hline
		$\Lgeo$, $\Lunit$, $X$, $\muX$, $\Yvf{A}$ & The important acoustic vectorfields. & Def.\,\ref{D:COMVECTORFIELDS} \\
		\hline 
		$\Lsmall^i$, $\Xsmall^i$, $\Yvfsmall{A}^i$ & $L^{\infty}$-small ``Error parts'' of the acoustic vectorfields. & \eqref{E:LSMALLDEF}, \eqref{E:XSMALL}, \eqref{E:YSMALL} \\
		\hline 
		$\comder$, $\tander$, $\tanderY$ & Sets of the commutation vectorfields, $\nullhyparg{u}$-tangent vectorfields, and $\ell_{t,u}$-tangent vectorfields, respectively. & \eqref{E:COMMUTATIONVECTORFIELDS}	\\
		\hline
		$\comder^{N;M}$, $\tander^N$, $\tanderY^N$, $\comdersmall^{N;M}$, $\tandersmall^N$, $\comderdoublesmall^{N;M}$, $\mathfrak{P}^{(N)}, \mathfrak{Y}^{(N)}$ & Strings of commutation vectorfields. & Def.\,\ref{D:STRINGSOFCOMMUTATIONVECTORFIELDS} \\
		\hline
		$\uLunit$ & A $\gfour$-null vectorfield normalized by $\gfour(\uLunit,\Lunit) = -2$ & \eqref{E:ULUNIT} \\
		\hline
		$\Sigmatproject, \smoothtorusproject, \Sigmatproject \xi, \smoothtorusproject \xi, \slashed{\upxi}$ & $\Sigma_t$ and $\ell_{t,u}$-projection tensorfields. Projections of spacetime tensors $\xi$ onto $\Sigma_t$ and $\ell_{t,u}$. & Def.\,\ref{D:PROJECTIONTENSORFIELDSANDTANGENCYTOHYPERSURFACES} \\
		\hline
		$\SigmatLie_Z \upxi, \angLie_Z \upxi$ & $\Sigma_t$-projected and $\ell_{t,u}$-projected Lie derivatives. & Def.\,\ref{D:PROJECTEDLIEDERIVATIVES} \\
		\hline
		$ \angLie_{\comder}^{N;M} \xi, \angLie_{\tander}^{N}\xi, \angLie_{\tanderY}^{N} \xi$. & Strings of $\ell_{t,u}$-projected Lie derivatives. & Def.\,\ref{D:STRINGSOFCOMMUTATIONVECTORFIELDS} \\
		\hline 
		$ \mytr_{\gfour} \xi, \mytr_{\gtorus}\xi$ & Trace of spacetime and $\ell_{t,u}$-tangent tensors. & Def.\,\ref{D:TRACEOFTENSORS} \\
		\hline
		$|\xi|_{\gfour}, |\xi|_g, |\xi|_{\gtorus}$ & Pointwise norms of tensorfields. & Def.\,\ref{D:POINTWISESEMINORMS} \\
		\hline
		$k, \upchi, \zeta$ & The second fundamental form of $\Sigma_t$, the null second fundamental form of $\ell_{t,u}$, and the $\ell_{t,u}$-tangent one form. & Def.\,\ref{D:SECONDFUNDAMENTALFORMSANDZETAONEFORM} \\
		\hline
		$\controlvars,\badcontrolvars$ & The controlling quantities. & Def.\,\ref{D:CONTROLVARS} \\
		\hline
		$ \angrmD \vec{x}$ & The array of $\ell_{t,u}$-projected spatial coordinate one-forms $\angrmD \vec{x}= (\angrmD x^1, \angrmD x^2, \angrmD x^3)$ & \eqref{E:ANDGULARDIFFERENTIALOFCARTESIANSPATIALCOORDIANTES} \\
		\hline
		$\deform{Z}$ & The deformation tensor of a spacetime vectorfield $Z$ with respect to $\gfour$. & \eqref{E:DEFORMATIONTENSORDEF} \\
		\hline 
		$\mathring{\Delta}_{\Sigma_0^{[u_1,u_2]}}^{\Ntop+1}$ & Norm of the data perturbation on $\Sigma_0$. & \eqref{E:PERTURBATIONSMALLNESSINCARTESIANDIFFERENTIALSTRUCTURE} \\
		\hline
		$\| f \|_{H_{\textnormal{Cartesian}}^N(\Sigma_0^{[u_1,u_2]} )}$, $\| f \|_{L_{\textnormal{Cartesian}}^2(\ell_{0,u})}$ & Sobolev and $L^2$ norms on $\Sigma_0^{[u_1,u_2]}$ and $\ell_{0,u}$. & \eqref{E:CARTESIANSOBOLEVNORMONSIGMA0}--\eqref{E:CARTESIANSOBOLEVNORMONELLTU}\\
		\hline
	\end{tabular}
\end{center}

\begin{center}
	\begin{tabular}{ | m{3cm} | m{10cm}| m{3cm} | } 
		\hline
		\multicolumn{3}{|c|}{\textbf{The rough time function and associated regions of spacetime}} \\
		\hline 
		Symbols & Descriptions & Reference \\ 
		\hline 
		$\phi = \phi(u)$ & The cutoff function which is identically equal to $1$ 
		when $|u| \leq  \frac{1}{2} \interestingu$. & \eqref{E:CUTOFFFUNCTION} \\
		\hline
		$\Wtransarg{\muxmulevelsetvalue}$ & The rough transversal vectorfield $\Wtransarg{\muxmulevelsetvalue} = \muX + \phi \frac{\muxmulevelsetvalue}{\Lunit \upmu} \Lunit$. & \eqref{E:WTRANSDEF} \\
		\hline 
		$\mulevelsetarg{\mulevelsetvalue}$, $\datahypfortimefunctionarg{-\muxmulevelsetvalue}$, 
		$\twoargmumuxtorus{\mulevelsetvalue}{-\muxmulevelsetvalue}$ & Level sets of $\upmu, \muX \upmu$, and the $\upmu$-adapted tori. & Def.\,\ref{D:LEVELSETSOFMUANDXMUANDMUMUXTORI} \\
		\hline
		$\timefunctionarg{\muxmulevelsetvalue}$ & The rough time function. & Def.\,\ref{D:ROUGHTIMEFUNCTION} \\
		\hline
		$\timefunction_0, - \mupositive$ & The largest-in-magnitude value of $\timefunctionarg{\muxmulevelsetvalue}$. & \eqref{E:RELATIONBETWEENINITIALROUGHTIMESLICEANDSMALLMUVALUE} \\
		\hline
		$(\timefunctionarg{\muxmulevelsetvalue},u,x^2,x^3)$ & The rough adapted coordinates. & Def.\,\ref{D:ROUGHCOORDINATESANDPARTIALDERIVATIVES} \\
		\hline
		$\{\roughgeop{\timefunctionarg{\muxmulevelsetvalue}},\roughgeop{u},\roughgeop{x^2},\roughgeop{x^3}\}$ & The rough adapted coordinate partial derivative vectorfields. & Def.\,\ref{D:ROUGHCOORDINATESANDPARTIALDERIVATIVES} \\
		\hline
		$\frac{\widetilde{\p}^{\vec{\alpha}}}{\widetilde{\p}(\timefunctionarg{\muxmulevelsetvalue},u,x^2,x^3)}$ & Multi-index notation in the rough adapted coordinates. & \eqref{E:ROUGHADAPATEDCOORDINATEMULTIINDEXDIFFERENTIALOPERATOR} \\
		\hline
		$\hypthreearg{\timefunction}{I}{\muxmulevelsetvalue}$, $\twoargroughtori{\timefunction,u}{\muxmulevelsetvalue}$ & The rough hypersurfaces and rough tori. & Def.\,\ref{D:TRUNCATEDROUGHSUBSETS} \\
		\hline
		$\nullhypthreearg{\muxmulevelsetvalue}{u}{J}$, 
		$\twoargMrough{I,J}{\muxmulevelsetvalue}$ & Truncated rough foliations of spacetime. & Def.\,\ref{D:TRUNCATEDROUGHSUBSETS} \\
		\hline
		$\mulevelsettwoarg{\mulevelsetvalue}{I}$, 
		$\datahypfortimefunctiontwoarg{-\muxmulevelsetvalue}{I}$ & Truncated level sets of $\upmu$ and $\muX \upmu$. & Def.\,\ref{D:TRUNCATEDLEVELSETSOFMUANDXMUANDMUMUXTORI}. \\
		\hline			$\smallneighborhoodofcreasetwoarg{[\timefunction_0,\timefunctionboot]}{\muxmulevelsetvalue}$ & Region of spacetime for which there is especially sharp control of $\upmu$. Specifically, the 
		estimate \eqref{E:WIDETILDELMUISALMOSTMINUSONEINSMALLNEIGHBORHOOD} holds on it. & \eqref{E:SMALLNEIGHBORHOOD} \\
		\hline
	\end{tabular}
\end{center}

\begin{center}
	\begin{tabular}{ | m{3cm} | m{10cm}| m{3cm} | } 
		\hline
		\multicolumn{3}{|c|}{\textbf{Relationships between the different coordinate systems}} \\
		\hline 
		Symbols & Descriptions & Reference \\ 
		\hline 
		$\Upsilon$ & The map $(t,u,x^2,x^3) \mapsto (t,x^1,x^2,x^3)$. & \eqref{E:CHOVGEOTOCARTESIAN} \\
		\hline
		$\CHOVgeotorough{\muxmulevelsetvalue}$ & The map $(t,u,x^2,x^3) \mapsto (\timefunctionarg{\muxmulevelsetvalue},u,x^2,x^3)$. & \eqref{E:CHOVGEOTOROUGH} \\
		\hline
		$\CHOVroughtomumuxmu{\muxmulevelsetvalue}, \CHOVJacobianroughtomumuxmu{\muxmulevelsetvalue}$ & The map $(\timefunctionarg{\muxmulevelsetvalue},u,x^2,x^3) \mapsto (\upmu,\muX \upmu, x^2,x^3)$ and its Jacobian. & \eqref{E:CHOVFROMROUGHCOORDINATESTOMUWEGIGHTEDXMUCOORDINATES}--\eqref{E:JACOBIANMATRIXFORCHOVFROMROUGHCOORDINATESTOMUWEGIGHTEDXMUCOORDINATES} \\
		\hline
	\end{tabular}
\end{center}
 
 \begin{center}
 	\begin{tabular}{ | m{3cm} | m{10cm}| m{3cm} | } 
 		\hline
 		\multicolumn{3}{|c|}{\textbf{The rough acoustic geometry}} \\
 		\hline 
 		Symbols & Descriptions & Reference \\ 
 		\hline 
 		$\gtorusroughfirstfund, \hypg$ & The first fundamental forms of $\twoargroughtori{\timefunction,u}{\muxmulevelsetvalue}$ and $\hypthreearg{\timefunction}{[-\farrightu,\leftu]}{\muxmulevelsetvalue}$. & Def.\,\ref{D:ROUGHFIRSTFUNDS} \\
 		\hline
 		$\argLrough{\muxmulevelsetvalue}$ & The rough null vectorfield. & \eqref{E:LROUGH} \\
 		\hline
 		$\FlowmapLrougharg{\muxmulevelsetvalue}$ & The $\timefunction_0$-normalized flow map of $\argLrough{\muxmulevelsetvalue}$. & Lemma~\ref{L:PROPERTIESOFFLOWMAPOFWIDETILDEL} \\
 		\hline 
 		$\Roughtoritangentvectorfieldarg{\muxmulevelsetvalue}$, $\Rtransarg{\muxmulevelsetvalue}$, 
		$\Rtransunitarg{\muxmulevelsetvalue}$, $\hypnormalarg{\muxmulevelsetvalue}$, $\hypunitnormalarg{\muxmulevelsetvalue}$ & Several geometric vectorfields adapted to the rough foliations. & Def.\,\ref{D:GEOMETRICVECTORFIELDSADAPTEDTOROUGHFOLIATIONS}\\
 		\hline
 		$\gtorusCOV_A^B$ & Matrix governing the relationship between $\geop{x^A}, \roughgeop{x^A}$ and $\gtorus,\gtorusroughfirstfund$. & \eqref{E:CHOVCOEFFICIENTSSMOOTHANGULARDERIVATIVESINTERMSOFROUGHONESANDL} \\
 		\hline
 		$\Rtransnormsmallfactorarg{\muxmulevelsetvalue} $ & Small, non-negative factor related to the size of 
		$\Rtransarg{\muxmulevelsetvalue}$. & \eqref{E:RTRANSNORMSMALLFACTOR}. \\
 		\hline
 		$\Riemtorus$, $\Rictorus$, $\Scalartorus$ & The Riemann and Ricci curvature tensors of $\gtorusroughfirstfund$. The scalar and Gauss curvatures of $\gtorusroughfirstfund$. & \eqref{E:RIEMANNCURVATURETENSOROFROUGHTORUS}--\eqref{E:2DGAUSSCURVATUREISTWICESCALARCURVATURE}. \\
\hline
 	\end{tabular}
\end{center}

 \begin{center}
	\begin{tabular}{ | m{3cm} | m{10cm}| m{3cm} | } 
		\hline
		\multicolumn{3}{|c|}{\textbf{Norms, volume forms, and $L^2$ ingredients}} \\
		\hline 
		Symbols & Descriptions & Reference \\ 
		\hline 
		$\| f\|_{W_{\textnormal{geo}}^{m,\infty}\left(\twoargMrough{I,J}{\muxmulevelsetvalue}\right)}$, $\| f \|_{C^{m,1}_{\textnormal{geo}}\left(\twoargMrough{I,J}{\muxmulevelsetvalue}\right)} $ & $L^\infty$-type Sobolev and H\"{o}lder norms of $\twoargMrough{I,J}{\muxmulevelsetvalue}$ in the $(t,u,x^2,x^3)$ coordinate system. & Def.\,\ref{D:ESSENTIALSUPNORMTYPENORMSANDHOLDERNORMS} \\
		\hline
		 $\| f \|_{C^{m,1}_{\textnormal{rough}}\left(I\times J \times \T^2\right)} $ & $L^\infty$-type H\"{o}lder norms in the $(\timefunctionarg{\muxmulevelsetvalue},u,x^2,x^3)$ coordinate system. & \eqref{E:CK1ROUGHNORMS} \\
		\hline
		$	\volMRoughCoordinates,\volPuRoughCoordinates,\volRoughHypersurface, \volroughtorus$ & Volume forms induced on the rough foliations by $\gfour$ in the rough adapted coordinates used to define $L^2$ norms. & Def.\,\ref{D:ROUGHVOLFORMS} \\
		\hline
		$\|\xi\|_{L^2(\twoargroughtori{\uptau,u}{\muxmulevelsetvalue})}$, $\| \xi\|_{L^2(\nullhypthreearg{u}{J}{\muxmulevelsetvalue})} $, $\|\xi\|_{L^2(\hypthreearg{\timefunction}{I}{\muxmulevelsetvalue})}$, $\|\xi\|_{L^2(\twoargMrough{I,J}{\muxmulevelsetvalue})}$ & $L^2$ norms relative to the volume forms 	$	\volMRoughCoordinates,\volPuRoughCoordinates,\volRoughHypersurface, \volroughtorus$ & Def.\,\ref{D:GEOMETRICL2NORMS} \\
		\hline
		$\volcanonical_{\hypg}, \volcanonical_{\gfour}$ & The canonical area and volume forms induced by $\hypg$ and $\gfour$, respectively. & Def.\,\ref{D:CANONICALVOLFORMSINROUGHADAPATEDCOORDINATES} \\
		\hline
		$\enmomem[f]$ & The energy-momentum tensor $\enmomem_{\alpha\beta}[f] = \Dfour_\alpha \Dfour_\beta f - \frac{1}{2} \gfour_{\alpha\beta} \gfour^{-1}(\Dfour f, \Dfour f)$. & \eqref{E:DEFOFENERGYMOMENTUM} \\
		\hline
		$\Jen{Z}[f]$ & The energy current vectorfield $\Jenarg{Z}{\alpha}[f] = \enmomem^{\alpha\beta}[f]Z_\beta$. & \eqref{E:WAVEEQUATIONENERGYCURRENT} \\
		\hline
		$\multipliervectorfield$ & The multiplier vectorfield $(1 + 2 \upmu) \Lunit + 2 \muX$. & \eqref{E:MULTIPLIERVECTORFIELD} \\
		\hline
		$\mathbb{E}_{(\textnormal{Wave})}, \mathbb{F}_{(\textnormal{Wave})}$, $\mathbb{E}_{(\textnormal{Transport})}$, $\mathbb{F}_{(\textnormal{Transport})}$ & The energies and null-fluxes for the wave and transport variables, respectively. & Def.\,\ref{D:WAVEANDTRANSPORTENERGIESANDNULLFLUXES} \\
		\hline
		$\spacetimeintegralcontrolwave[f]$ & A coercive spacetime integral used in the wave equation energy estimates. & \eqref{E:COERCIVESPACETIMEWAVENERGYINTEGRAL} \\
		\hline
		${^{(\multipliervectorfield)}\mathfrak{B}}[f], {^{(\multipliervectorfield)}\mathfrak{B}}_{(1)}[f]$, $\cdots, {^{(\multipliervectorfield)}\mathfrak{B}}_{(6)}[f]$ & Bulk terms arising from $\upmu \enmomem^{\alpha\beta}[f] \deformarg{\multipliervectorfield}{\alpha}{\beta}[f]$ in the fundamental energy identity for wave equations. & \eqref{E:WAVEENERGYIDENTITYBULKERRORTERM}-- \eqref{E:WAVEENERGYIDENTITYBULKERRORTERM6} \\
		\hline
		$\hypersurfacecontrolwave_N, \spacetimeintegralcontrolwave_N$, $\totalcontrolwave_N$ & The $N$-th order wave-controlling quantities for all the wave variables $\{\RRiemann,\LRiemann,v^2,v^3,\Ent\}$. & \eqref{E:WAVESPACELIKEANDNULLHYPERSURFACEL2CONTROLLINGQUANTITY}--\eqref{E:WAVETOTALL2CONTROLLINGQUANTITY} \\
		\hline
		$\hypersurfacecontrolwavepartial_N, \spacetimeintegralcontrolwavepartial_N$, $\totalcontrolwavepartial_N$ & The $N$-th order wave-controlling quantities for all the partial wave variables $\{\LRiemann,v^2,v^3,\Ent\}$. & \eqref{E:WAVEPARTIALSPACELIKEANDNULLHYPERSURFACEL2CONTROLLINGQUANTITY}--\eqref{E:WAVEPARTIALL2CONTROLLINGQUANTITY} \\
		\hline
		$\hypersurfacecontrolVort_N,\hypersurfacecontrolGradEnt_N$ & The $N$-th order specific vorticity and entropy gradient controlling quantities. & \eqref{E:VORTICITYL2CONTROLLINGQUANTITY}--\eqref{E:ENTROPYGRADIENTL2CONTROLLINGQUANTITY} \\
		\hline	
		 $\toricontrolVort_N, \toricontrolGradEnt_N$ & The $N$-th order specific vorticity and entropy gradient controlling quantities \emph{on the rough tori}. & \eqref{E:TORIVORTICITYL2CONTROLLINGQUANTITY}--\eqref{E:TORIENTROPYGRADIENTL2CONTROLLINGQUANTITY} \\
		 \hline
		 $\hypersurfacecontrolVortVort_N, \hypersurfacecontrolDivGradEnt_N$ & The $N$-th order controlling quantities for the modified fluid variables. & \eqref{E:MODIFIEDVORTICITYVORTICITYL2CONTROLLINGQUANTITY}--\eqref{E:MODIFIEDDIVGRADENTL2CONTROLLINGQUANTITY} \\
		 \hline
$\toricontrolVortVort_N, \toricontrolDivGradEnt_N$ & The $N$-th order controlling quantities 
for the modified fluid variables 
\emph{on the rough tori}. & \eqref{E:TORIMODIFIEDVORTICITYL2CONTROLLINGQUANTITY}--\eqref{E:TORIDIVGRADENTL2CONTROLLINGQUANTITY} 
\\
		 \hline
		 \multicolumn{3}{|C{16cm}|}{We use the following summation conventions:
		 $\hypersurfacecontrolwave_{[N_1,N_2]}= \sum_{M = N_1}^{N_2} \hypersurfacecontrolwave_M$, $\hypersurfacecontrolVort_{\leq N} = \sum_{M = 0}^N \hypersurfacecontrolVort_M$, and similarly for the other controlling quantities.}\\
		 \hline
\end{tabular}
\end{center} 		

 \begin{center}
	\begin{tabular}{ | m{3cm} | m{10cm}| m{3cm} | } 
		\hline
		\multicolumn{3}{|c|}{\textbf{Bootstrap assumptions}} \\
		\hline 
		Symbols & Descriptions & Reference \\ 
		\hline 
		$\timefunctionboot$ & The bootstrap rough-time. & \eqref{E:BOOTSTRAPTIME} \\
		\hline
		$\Cartesiantisafunctiononmumxtoriarg{\mulevelsetvalue}{-\muxmulevelsetvalue}, \, \Eikonalisafunctiononmumuxtoriarg{\mulevelsetvalue}{-\muxmulevelsetvalue}$ & The $\upmu$-adapted torus $\twoargmumuxtorus{\mulevelsetvalue}{-\muxmulevelsetvalue}$ is a graph over $(x^2,x^3) \in \T^2$: $\twoargmumuxtorus{\mulevelsetvalue}{-\muxmulevelsetvalue} = \{( \Cartesiantisafunctiononmumxtoriarg{\mulevelsetvalue}{-\muxmulevelsetvalue}(x^2,x^3), \Eikonalisafunctiononmumuxtoriarg{\mulevelsetvalue}{-\muxmulevelsetvalue}(x^2,x^3), x^2,x^3) \ | \ (x^2,x^3) \in \T^2\}$. & \eqref{E:BAMUTORI} \\
		\hline
		$\embeddatahypersurfacearg{\muxmulevelsetvalue}$ & The mapping defined by $(\mulevelsetvalue,x^2,x^3) \in (\upmu_{Boot},\upmu_0]\times \T^2 \mapsto$ $( \Cartesiantisafunctiononmumxtoriarg{\mulevelsetvalue}{-\muxmulevelsetvalue}(x^2,x^3), \Eikonalisafunctiononmumuxtoriarg{\mulevelsetvalue}{-\muxmulevelsetvalue}(x^2,x^3), x^2,x^3)$. & \eqref{E:EMBEDDATAHYPERSURFACE} \\
		\hline 
		\multicolumn{3}{|C{16cm}|}{
		(\textbf{Quantitative improvement of bootstrap assumptions}) By this, we mean 
		that some quantity $Q$ 
		was assumed to satisfy $A_1 \leq Q \leq A_2$ in the bootstrap assumptions
		(where $A_1,A_2$ are real numbers),
		and we derive the improved bound $B_1 \leq Q \leq B_2$, 
		where $A_1 < B_1 \leq B_2 < A_2$.} \\
		\hline
		\multicolumn{3}{|C{16cm}|}{(\textbf{From soft to quantitative improvement of bootstrap assumptions}) By this, we mean that in the bootstrap assumptions,
			we assumed that some function $Q$ belonged to some function space and had a finite norm in that space,
			and our improvement is a quantitative estimate for the norm of $Q$.} \\
		\hline 
		\multicolumn{3}{|C{16cm}|}{(\textbf{Extension to the closure improvement of bootstrap assumptions})
			By this, we mean that our bootstrap assumptions involved an assumption on the ``open-at-the-top'' domain
			$\twoargMrough{[\timefunction_0,\timefunctionboot),[-\farrightu,\leftu]}{\muxmulevelsetvalue}$,
			and we derive an improved result showing that the assumption holds
			on the closed domain $\twoargMrough{[\timefunction_0,\timefunctionboot],[-\farrightu,\leftu]}{\muxmulevelsetvalue}$.} \\
		\hline
 		\end{tabular}
 \end{center}

 \begin{center}
	\begin{tabular}{ | m{3cm} | m{10cm}| m{3cm} | } 
		\hline
		\multicolumn{3}{|c|}{\textbf{Embeddings and flow maps}} \\
		\hline 
		Symbols & Descriptions & Reference \\ 
		\hline
		$\domainforembeddingdatahypfortimefunctiontwoarg{\muxmulevelsetvalue}{(\upmuboot,\mupositive]}$ & A subset of $\R \times \T^2$ which is diffeomorphic to $\datahypfortimefunctiontwoarg{-\muxmulevelsetvalue}{[\timefunction_0,\timefunctionboot)}$. & \eqref{E:DOMAINOFEMBEDDINGFORXMUEQUALSMINUSKAPPAYSURFACE} \\
		\hline
		$\embeddingdatahypfortimefunctionarg{\muxmulevelsetvalue}, \scalarembeddingdatahypfortimefunctionarg{\muxmulevelsetvalue}$ & The function $\embeddingdatahypfortimefunctionarg{\muxmulevelsetvalue}(t,x^2,x^3) = (t,\scalarembeddingdatahypfortimefunctionarg{\muxmulevelsetvalue}(t,x^2,x^3),x^2,x^3)$ is an embedding from $\domainforembeddingdatahypfortimefunctiontwoarg{\muxmulevelsetvalue}{(\upmuboot,\mupositive]} \to \twoargMrough{[\timefunction,\timefunctionboot),[-\frac{3}{4} \rightu,\frac{3}{4}\rightu]}{\muxmulevelsetvalue}$ whose image is $\datahypfortimefunctiontwoarg{-\muxmulevelsetvalue}{[\timefunction_0,\timefunctionboot)}$ & Lemma~\ref{L:XMUISMINUSCAPPISAGRAPH} \\
		\hline
		$\flowmapWtransargtwoarg{\muxmulevelsetvalue}{\Delta u}$ & The flow map of $\Wtrans$. & \eqref{E:FLOWMAPFORGENERATOROFROUGHTIMEFUNCTION} \\
		\hline
		$\composedflowmapdiffeoarg{\muxmulevelsetvalue}, \domaincomposedflowmapdiffeoarg{\muxmulevelsetvalue}$ & $\composedflowmapdiffeoarg{\muxmulevelsetvalue}(\Delta u,t,x^2,x^3) = \flowmapWtransargtwoarg{\muxmulevelsetvalue}{\Delta u} \circ \embeddingdatahypfortimefunctionarg{\muxmulevelsetvalue}(t,x^2,x^3)$, and $\domaincomposedflowmapdiffeoarg{\muxmulevelsetvalue}$ is its domain. & Lemma~\ref{L:FLOWMAPFORGENERATOROFROUGHTIMEFUNCTION}. \\
		\hline
		$\flowmapWtransMatrixarg{\muxmulevelsetvalue}$ & The $4\times 4$ matrix-valued  function on 
		$\domainforembeddingdatahypfortimefunctiontwoarg{\muxmulevelsetvalue}{(\upmuboot,\mupositive]}$
		whose first column is
		$(0,1,0,0)^{\top}$
		and whose last three columns are the Jacobian $d_{(t,x^2,x^3)} \embeddingdatahypfortimefunctionarg{\muxmulevelsetvalue}$. & \eqref{E:MATRIXARISINGINSTUDYOFFLOWMAPOFWTRANS} \\
		\hline
		$\ambient$ & An ``ambient" spacetime function of class $W_{\textnormal{geo}}^{3,\infty}(\twoargMrough{(\timefunction_0,\timefunctionboot),(-\farrightu,\leftu)}{\muxmulevelsetvalue})$, satisfying the constraint $(\Wtransarg{\muxmulevelsetvalue} \ambient)|_{\datahypfortimefunctiontwoarg{-\muxmulevelsetvalue}{[\timefunction_0,\timefunctionboot)}}
		= 0$ & Lemma~\ref{L:ODESOLUTIONSTHATARESMOOTHERTHANTHEDATAHYPERSURFACE} \\
		\hline
		$\Cartesiantisafunctiononlevelsetsofroughtimefunctionarg{\timefunction}{\muxmulevelsetvalue}$ & Describes the Cartesian time function $t$ as a function of $(u,x^2,x^3)$ on the rough hypersurfaces. That is, $	\hypthreearg{\timefunction}{[-\farrightu,\leftu]}{\muxmulevelsetvalue}
		 =
		\left\lbrace
		(t,u,x^2,x^3)
		\ | \
		t = \Cartesiantisafunctiononlevelsetsofroughtimefunctionarg{\timefunction}{\muxmulevelsetvalue}(u,x^2,x^3),
		\,
		(u,x^2,x^3) \in [-\farrightu,\leftu] \times \mathbb{T}^2
		\right\rbrace$. & \eqref{E:LEVELSETSOFTIMEFUNCTIONAREAGRAPH} \\
		\hline
		$\embeddingofmuadapatedtorinCartesianspace{\muxmulevelsetvalue}$ & The map defined by $\embeddingofmuadapatedtorinCartesianspace{\muxmulevelsetvalue}(x^2,x^3) =
		\Upsilon
		\circ
		\left(\Cartesiantisafunctiononmumxtoriarg{\mulevelsetvalue}{-\muxmulevelsetvalue}(x^2,x^3),
		\Eikonalisafunctiononmumuxtoriarg{\mulevelsetvalue}{-\muxmulevelsetvalue}(x^2,x^3),
		x^2,
		x^3 
		\right)$, which is a diffeomorphism from $\T^2$ onto 	$\Upsilon\left(\twoargmumuxtorus{\mulevelsetvalue}{-\muxmulevelsetvalue} \right)$. & \eqref{E:MUADAPATEDTORIINCARTESIANSPACE} \\
		\hline
	\end{tabular}
\end{center}

 \begin{center}
	\begin{tabular}{ | m{3cm} | m{10cm}| m{3cm} | } 
		\hline
		\multicolumn{3}{|c|}{\textbf{The elliptic--hyperbolic integral identities}} \\
		\hline 
		Symbols & Descriptions & Reference \\ 
		\hline
		$\Nullhypersurfaceproject$ & The projection tensorfield 
		$\Nullhypersurfaceproject_{\beta}^{\ \alpha} 
		= 
		\updelta_{\beta}^{\ \alpha} + \frac{1}{2} \uLunit^{\alpha} \Lunit_{\beta}$ onto $\nullhyparg{u}$. & 
			\eqref{E:NULLHYPERSURFACEPROJECTION} \\
		\hline
		 $\hfour, \hfour^{-1}$ & Spacetime Riemannian metric and inverse metrics. & \eqref{E:RIEMANNIANACOUSTICALMETRIC}--\eqref{E:INVERSERIEMANNIANACOUSTICALMETRIC} \\
		 \hline
		$\Nullhypersurfacemetric, \Nullhypersurfaceinversemetric$ & Positive semi-definite $\binom{0}{2}, \binom{2}{0}$ tensorfields which restrict to Riemannian metric and inverse metrics on $\nullhyparg{u}$. & \eqref{E:NULLHYPERSURFACESRIEMANNIANMETRIC}--\eqref{E:NULLHYPERSURFACESINVERSERIEMANNIANMETRIC} \\
		\hline
		$|Z V|_{g}$ & The $g$-norm of the $\Sigma_t$-tangent vectorfield $ZV^a$, where $Z$ is an arbitrary spacetime vectorfield and $V^a$ is a $\Sigma_t$-tangent vectorfield. & Def.\,\ref{D:DERIVATIVEOFATENSORFIELD} \\
		\hline
		$|\xi|_{\hfour}$ & The $\hfour$-norm of a $\binom{m}{n}$ spacetime tensorfield. & 
		Def.\,\ref{D:RIEMANNIANACOUSTICALMETRICSPACETIMENORMOFTENSORFIELDS} \\
		\hline
		$\ellipticCoerciveQuadratic[\pmb{\partial}V,\pmb{\partial}V]$ & The coercive elliptic-hyperbolic quadratic form. &
		\eqref{E:NULLHYPERSURFACEADAPTEDCOERCIVEQUADRATICFORM} \\
		\hline
		$\ehcurrent[V,\pmb{\partial}V]$ & The characteristic current, also referred to as the elliptic-hyperbolic current. & \eqref{E:PUTANGENTELLIPTICHYPERBOLICCURRENT} \\
		\hline
		$\weight$ & An arbitrary weight function. & Lemma~\ref{L:COVARIANTDIVERGENCEIDENTITYFORELLIPTICHYPERBOLICCURRENT} \\
		\hline
		$\mathfrak{J}_{(\textnormal{Antisymmetric})}[\pmb{\partial} \SigmatTan,$ $\pmb{\partial} \SigmatTan]$, $\mathfrak{J}_{(\textnormal{Div})}[\pmb{\partial} \SigmatTan,\pmb{\partial} \SigmatTan]$ & Error terms arising in the covariant divergence identity for $\ehcurrent$ that are \emph{quadratic} in $\pmb{\partial} V^a$. & \eqref{E:ANTISYMMETRICNULLCURRENTSPACETIMERRORTERM}--\eqref{E:DIVERGENCENULLCURRENTSPACETIMERRORTERM} \\
		\hline
		$\mathfrak{J}_{(\pmb{\partial} \weight)}[V,\pmb{\partial} \SigmatTan]$, $\mathfrak{J}_{(\textnormal{Absorb-1})}[V,\pmb{\partial} \SigmatTan]$, $\mathfrak{J}_{(\textnormal{Absorb-2})}[V,\pmb{\partial} \SigmatTan]$, $\mathfrak{J}_{(\textnormal{Null Geometry})}[V,$ $\pmb{\partial} \SigmatTan]$ &  Error terms arising in the covariant divergence identity for $\ehcurrent$ that are \emph{linear} in $\pmb{\partial} V^a$. & \eqref{E:DERIVATIVEOFWEIGHTNULLCURRENTSPACETIMERRORTERM}--\eqref{E:DERIVATIVESOFNULLGEOMETRYNULLCURRENTSPACETIMERRORTERM} \\
		\hline
		$\CurrentboundaryerrorperfectRderivative[\SigmatTan,\SigmatTan]$ & Term arising in the boundary terms for the characteristic current identity which is a \emph{perfect} $\Rtransarg{\muxmulevelsetvalue}$ derivative. It is positive definite in $\Sigma_t$-tangent vectorfields in the sense that $\CurrentboundaryerrorperfectRderivative[\SigmatTan,\SigmatTan] \approx \frac{1}{\upmu - \phi \frac{\muxmulevelsetvalue}{\Lunit \upmu}} |V|_g^2$. & \eqref{E:PERFECTRDERIVAIVEERRORPUTANGENTCURRENTCONTRACTEDAGAINSTVECTORFIELD} \\
		\hline
		$\Currentboundaryerrorhavetocontrolprincipal[\SigmatTan,\pmb{\partial} \SigmatTan]$ & Principal order error terms arising in the boundary terms for the characteristic current identity. & \eqref{E:PRINCIPALERRORTERMHAVETOCONTROLKEYIDPUTANGENTCURRENTCONTRACTEDAGAINSTVECTORFIELD} \\
		\hline
		$\Currentboundaryerrorhavetocontrollowerorder[\SigmatTan,\SigmatTan]$ & Lower order error terms arising in the boundary terms for the characteristic current identity. & \eqref{E:LOWERORDERERRORTERMHAVETOCONTROLKEYIDPUTANGENTCURRENTCONTRACTEDAGAINSTVECTORFIELD} \\
		\hline
		$\EllipticHyperbolicCurrentIntegralIdentityTotalSpacetimeErrorTerm[\SigmatTan,\pmb{\partial} \SigmatTan]$ & Spacetime bulk error term arising in the elliptic-hyperbolic integral identity. & \eqref{E:ELLIPTICHYPERBOLICINTEGRALIDENTITYBULKERRORTERM} \\
		\hline	
	\end{tabular}
\end{center}

\begin{center}
	\begin{tabular}{ | m{3cm} | m{10cm}| m{3cm} | } 
		\hline
		\multicolumn{3}{|c|}{\textbf{Error terms in the geometric wave equations}} \\
		\hline 
		Symbols & Descriptions & Reference \\ 
		\hline
		$\HarmlessWave{N}$ & $N$-th order harmless wave error terms. & \eqref{E:HARMLESSWAVE} \\
		\hline
		$\vec{\mathfrak{G}},\mathfrak{G}_{\iota}$ & For $\wavearray = (\Psi_0,\Psi_1,\Psi_2,\Psi_3,\Psi_4) 
		= (\RRiemann,\LRiemann, v^2,v^3,\Ent)$, $\vec{\mathfrak{G}} = ( \mathfrak{G}_0,\cdots,\mathfrak{G}_4)$ is the vector array of the inhomogeneous  terms in the covariant wave equations $\upmu \Box_{\gfour}\Psi_{\iota} = \mathfrak{G}_{\iota}$ 
		& Prop.\,\ref{P:MOSTDIFFICULTWAVETERMS} \\
		\hline
	\end{tabular}
\end{center}

\begin{center}
	\begin{tabular}{ | m{3cm} | m{10cm}| m{3cm} | } 
		\hline
		\multicolumn{3}{|c|}{\textbf{Avoiding derivative loss when controlling the acoustic geometry}} \\
		\hline 
		Symbols & Descriptions & Reference \\ 
		\hline
		$\mathfrak{A}$ & Error term arising in the key identity verified by $\upmu \Ricfour_{LL}$. The identity requires the use of the geometric wave equations \eqref{E:VELOCITYWAVEEQUATION}--\eqref{E:ENTROPYWAVEEQUATION}. & \eqref{E:AINHOMRIC} \\
		\hline
		$\mathfrak{B}$ & Error term arising from key identity verified by $\Ricfour_{LL}$. The identity does not require that the  geometric wave equations \eqref{E:VELOCITYWAVEEQUATION}--\eqref{E:WAVEEQUATIONFORALMOSTRIEMANNINVARIANTS} are satisfied. & \eqref{E:BINHOMRIC} \\
		\hline
		$\fullymodquant{\tander^N}$ & Fully modified version of $\upmu \tander^N \mytr_{\gtorus} \upchi$ that satisfies the favorable transport equation \eqref{E:TRANSPORTEQUATIONFORFULLYMODIFIEDQUANTITY}. & \eqref{E:FULLYMODIFIEDQUANTITY} \\
		\hline
		$\fullymodquantinhom$ & Inhomogeneous term in the modified quantity $\fullymodquant{\tander^N} = \upmu \tander^N \mytr_{\gtorus}\upchi + \tander^N \fullymodquantinhom$. & \eqref{E:MODIFIEDQUANTITYINHOM} \\
		\hline
		$\partialmodquant{\tander^N}$ & Partially modified version of $\tander^N \mytr_{\gtorus}\upchi$ that satisfies the favorable transport equation \eqref{E:TRANSPORTEQUATIONFORPARTIALMODIFIEDQUANTITY}. & \eqref{E:PARTIALMODIFIEDQUANTITY} \\
		\hline
		$\partialmodquantinhom{\tander^N}, \widetilde{\mathfrak{X}}$ & $N$-th and $0$-th order inhomogeneous terms in the partially modified quantity $\partialmodquant{\tander^N}$. & \eqref{E:PARTIALMODIFIEDQUANTITYINHOM}--\eqref{E:PARTIALMODIFIEDQUANTITYINHOMZEROORDER} \\
		\hline
		$^{(\tander^{N-1})}\mathfrak{B}$ & Error terms arising in the transport equation for $\partialmodquant{\tander^N}$. & \eqref{E:COMMUTEDPARTIALMODIFIEDQUANTITYINHOM} \\
		\hline
		$|\widetilde{\upxi}|_{\gtorusroughfirstfund}$ & Norms of $\twoargroughtori{\timefunction,u}{\muxmulevelsetvalue}$-tensors with respect to $\gtorusroughfirstfund$. & Def.\,\ref{D:POINTWISESEMINORMWITHRESPECTTOFIRSTFUNDOFROUGHTORI} \\
		\hline
		$\mytr_{\gtorusroughfirstfund} \widetilde{\upxi}$ & The $\gtorusroughfirstfund$-trace of a type-$\binom{0}{2}$ tensorfield on $\twoargroughtori{\timefunction,u}{\muxmulevelsetvalue}$: 
		$\mytr_{\gtorusroughfirstfund} \widetilde{\upxi} = (\gtorusroughinversefirstfund)^{\alpha \beta}
\upxi_{\alpha \beta}.$ & \eqref{E:TRACEOFROUGHTORITANGENT02TENSORS}\\
		\hline
		$\widetilde{\angrmD} \varphi$ & The rough differential of a scalar function. & Def.\,\ref{D:ROUGHTORUSDIFFERENTIAL} \\
		\hline
		$\roughangdiv \xi$ & The rough divergence of a $\twoargroughtori{\timefunction,u}{\muxmulevelsetvalue}$-tangent tensorfield $\xi$. & Def.\,\ref{D:CONNECTIONSANDDIFFERENTIALOPERATORSONROUGHTORI} \\
		\hline
		$\{e_A\}_{A=2,3}$, $\{f_A\}_{A=2,3}$ & The frames obtained from applying Gram--Schmidt to $\{\geop{x^A}\}_{A=2,3}$ and 
		$\{\roughgeop{x^A}\}_{A = 2,3}$ with $\gtorus$, $\gtorusroughfirstfund$ respectively. & 
		Def.\,\ref{D:ORTHONORMALFRAMESONSMOOTHANDROUGHTORI} \\
		\hline
		$\COVframe_{AB}$, $\COVL_A$ & Change of frame coefficients between $\{e_A\}_{A=2,3}$ and $\{f_A\}_{A=2,3}$. & Lemma~\ref{L:CHANGEOFORTHONORMALFRAMES} \\
		\hline
	\end{tabular}
\end{center}

\begin{center}
	\begin{tabular}{ | m{3cm} | m{10cm}| m{3cm} | } 
		\hline
		\multicolumn{3}{|c|}{\textbf{Developments of the data, the singular boundary and crease, and a new time function}} \\
		\hline 
		Symbols & Descriptions & Reference \\ 
		\hline
		$\MLeft$ & The portion of the development sandwiched between $\nullhypthreearg{U_2}{[\timefunction_0,0]}{\muxmulevelsetvalue}$ and  $\datahypfortimefunctiontwoarg{0}{[\timefunction_0,0]}$. & \eqref{E:LEFTDEVELOPMENT} \\
		\hline
		$\MRight$ & The portion of the development sandwiched between $\datahypfortimefunctiontwoarg{-\muxmulevelsetvalue_0}{[\timefunction_0,0]}$ and $\nullhypthreearg{-\farrightu}{[\timefunction_0,0]}{\muxmulevelsetvalue}$. & \eqref{E:RIGHTDEVELOPMENT} \\
		\hline
		$\MSingular$ & The portion of the development sandwiched between $\datahypfortimefunctiontwoarg{0}{[\timefunction_0,0]}$ and $\datahypfortimefunctiontwoarg{-\muxmulevelsetvalue_0}{[\timefunction_0,0]}$. & \eqref{E:SINGULARDEVELOPMENT} \\
		\hline
		$\MInteresting$ & The union of $\MLeft,\MSingular,\MRight$. & \eqref{E:INTERESTINGDEVELOPMENTOFDATA} \\
		\hline
		$\mathcal{B}^{[0,\muxmulevelsetvalue_0]}$ & The singular boundary portion of the development, given by $\mathcal{B}^{[0,\muxmulevelsetvalue_0]} = \bigcup_{\muxmulevelsetvalue \in [0,\muxmulevelsetvalue_0]} \twoargmumuxtorus{0}{-\muxmulevelsetvalue}$. & \eqref{E:SINGULARBOUNDARYPORTION} \\
		\hline
		 $\partial_- \mathcal{B}^{[0,\muxmulevelsetvalue_0]}$ & The crease, given by $\partial_- \mathcal{B}^{[0,\muxmulevelsetvalue_0]} 
		 \eqdef
		 \twoargmumuxtorus{0}{0}$. & \eqref{E:CREASE} \\
		 \hline
		 $\extendedembeddatahypersurface$ & The map given by $\extendedembeddatahypersurface(\mulevelsetvalue,\muxmulevelsetvalue,x^2,x^3) = (\Cartesiantisafunctiononmumxtoriarg{\mulevelsetvalue}{-\muxmulevelsetvalue}(x^2,x^3),\Eikonalisafunctiononmumuxtoriarg{\mulevelsetvalue}{-\muxmulevelsetvalue}(x^2,x^3),x^2,x^3)$, which is a $C^{1,1}$ diffeomorphism from
		 $[0,\mupositive] \times [0,\muxmulevelsetvalue_0] \times \mathbb{T}^2$
		 onto $\MSingular$. & \eqref{E:EMBEDDINGOFSINGULARREGION} \\
		 \hline
		 $\domainofgraphofmulevelsetinsingularregion{[0,\muxmulevelsetvalue_0]}{\mulevelsetvalue}$ & The $\mulevelsetvalue$-level sets of $\upmu$ in the singular development as a graph. & \eqref{E:DOMAINOFGRAPHFORMULEVELSETINSINGULARREGION} \\
		 \hline
		 $\tisafunctiononlevelsetsofmu\mulevelsetvalue$ & The graph of the Cartesian $t$ as a function of $(u,x^2,x^3)$ along the $\mulevelsetvalue$-level sets of $\upmu$ in the singular development. & \eqref{E:GRAPHSTRUCTUREOFLEVELSETSOFMUINSINGULARGETION} \\
		 \hline
		 $\newtimefunction$ & The time function whose level sets foliate $\MInteresting$. & \eqref{E:NEWTIMEFUNCTION} \\
		 \hline
		 $(\newtimefunction,u$, $x^2,x^3)$ & The interesting coordinate system. & Def.\,\ref{D:LEVELSETGENERATORFORNEWTIMEFUNCTIONANDNEWTIMEFUNCTION} \\
		 \hline
		 $\InterestingCHOV$ & The map $(t,u,x^2,x^3) \mapsto (\newtimefunction,u,x^2,x^3)$. & \eqref{E:CHOVFROMGEOTOINTERESTINGCOORDS} \\
		 \hline
		 $\tisafunctionalonglevelsetsofnewtimefunctionarg{\timefunction}$ & The function on $[-\farrightu,\leftu]\times\T^2$ whose graph is the Cartesian $t$ on $\MInteresting$. & \eqref{E:DEFININGFUNCTIONCARTESIANTISAGRAPHALONGLEELSETSOFNEWTIMEFUNCTION}\\
		 \hline
		 $\levelsetgeneratornewtimefunction$ & The vectorfield that is the coordinate partial derivative with respect to $u$ in the interesting coordinate system. & \eqref{E:LEVELSETGENERATORFORNEWTIMEFUNCTION} \\
		 \hline
		 $\partialderivativewithrespecttonewtimefunction$ & The vectorfield that is the coordinate partial derivative with respect to $\newtimefunction$ in the interesting coordinate system. & \eqref{E:PARTIALDERIVATIVEWITHRESPECTTOTIMEFUNCTION} \\
		 \hline
	\end{tabular}
\end{center}

\bibliographystyle{amsalpha}

\bibliography{JBib}

\end{document}